%% file: main_algebra.tex
\documentclass[12pt,oneside,reqno]{amsbook}

\input{preambleThAV.tex}

\begin{document}
\pagestyle{plain}

\frontmatter

\include{parts/front/TitlePage}          
\include{parts/front/Copyright}          
\include{parts/front/Dedication}         
\include{parts/front/Ackno}              
\include{parts/front/Abstract}           
\include{parts/front/Contents}           
\mainmatter
\bibliographystyle{plain}                     

\pagestyle{headings}
\include{parts/front/Introduction}    

\part{The Pre-Categorical Phase: Set Theory and Fundamental Algebra Concepts}    
\input{parts/AlgebraNC/AlgebraNC1}    
\input{parts/AlgebraNC/AlgebraNC2}    

\input{parts/AlgebraNC/AlgebraNC3}    
\input{parts/AlgebraNC/AlgebraNC4}    
\input{parts/AlgebraNC/AlgebraNC5}    
\input{parts/AlgebraNC/AlgebraNC6}    

\part{The 1st Categorical Phase: Systematic Structures and Imbedding of Abelian Categories}\label{ChCatPers}  
\input{parts/AlgebraCat/AlgebraCatS1}    
\input{parts/AlgebraCat/AlgebraCatS2}    
\input{parts/AlgebraCat/AlgebraCatS3}    
\input{parts/AlgebraCat/AlgebraCatS4}    
\input{parts/AlgebraCat/AlgebraCatS5}    
\input{parts/AlgebraCat/AlgebraCatS6}    
\input{parts/AlgebraCat/AlgebraCatS7}    

\part{The 2nd Categorical Phase: Homological Structures and Triangulation of Derived Categories}\label{ChCatPers2}  
\input{parts/AlgebraCat/AlgebraCatS8}     
\input{parts/AlgebraCat/AlgebraCatS9}     
\input{parts/AlgebraCat/AlgebraCatS10}    
\input{parts/AlgebraCat/HomDim}           

\input{parts/AlgebraCat/AlgebraCatS11}    
\input{parts/AlgebraCat/AlgebraCatS12}    
\input{parts/AlgebraCat/AlgebraCatS13}    
\input{parts/AlgebraCat/AlgebraCatS14}    

\part{The Miscellaneous Phase: Hints on Application}       
\input{parts/AlgebraM/RepsTheoryI}    
\input{parts/AlgebraM/RepsTheoryII}   
\input{parts/AlgebraM/RepsTheoryIII}  
\input{parts/AlgebraM/GeomAnaI}       
\input{parts/AlgebraM/GeomAnaII}      
\input{parts/AlgebraM/GeomAnaIII}     


\backmatter
\bibliographystyle{amsplain}
\input{parts/references/RefsAlg}                    
{\footnotesize
\printindex 
}

\include{parts/front/Vita}                   
\end{document}

%% file: preambleThAV.tex
\usepackage{amsmath,amsthm,amsfonts,amssymb,amsrefs,amscd,graphicx,nicefrac}
\allowdisplaybreaks  

\usepackage[titletoc]{appendix}
\usepackage[active]{srcltx}
\usepackage{dsfont,amscd}
\usepackage{mathrsfs}
\usepackage{amstext}
\usepackage[capposition=bottom]{floatrow}
\usepackage{setspace}
\usepackage{color,xcolor}
\usepackage{enumerate}
\usepackage{fancyhdr,url}
\usepackage[top=0.5in, bottom=0.5in, left=0.5in, right=0.5in]{geometry}

\usepackage{comment} 
\usepackage[english]{babel}
\usepackage{todonotes}
\usepackage{tikz,tikz-cd,adjustbox}
\usepackage{etoolbox}
\usepackage{xfrac,soul}

\usepackage[inline, shortlabels]{enumitem} 

\usepackage{tabu,tabularx}

\usepackage{imakeidx}
\makeindex[title=Index of Concepts,columns=3]

\usepackage[hidelinks]{hyperref}
\usepackage{setspace}\setdisplayskipstretch{}  
\definecolor{OliveGreen}{rgb}{0,0.6,0}

\newcommand{\Natural}{\mathbb{N}}
\newcommand{\Integer}{\mathbb{Z}}
\newcommand{\Rational}{\mathbb{Q}}
\newcommand{\Real}{\mathbb{R}}
\newcommand{\Complex}{\mathbb{C}}
\newcommand{\Cardinal}{\mathbb{K}}
\newcommand{\Ordinal}{\mathbb{O}}

\DeclareMathOperator{\Ob}{Ob} 
\DeclareMathOperator{\Mor}{Mor} \DeclareMathOperator{\mor}{mor}
\DeclareMathOperator{\Tot}{Tot} 
\DeclareMathOperator{\Tor}{Tor} 
\DeclareMathOperator{\Ext}{Ext}

\DeclareMathOperator{\Top}{Top} 
\DeclareMathOperator{\Sets}{Sets} 
 \DeclareMathOperator{\Spec}{Spec}

\DeclareMathOperator{\im}{im}

\DeclareMathOperator{\rank}{rank}
\DeclareMathOperator{\coker}{coker}
\DeclareMathOperator{\coim}{coim}
\DeclareMathOperator{\dom}{dom}
\DeclareMathOperator{\cod}{cod}
\DeclareMathOperator{\diam}{diam}

\DeclareMathOperator{\dist}{dist}
\DeclareMathOperator{\Conv}{Conv}

\DeclareMathOperator{\sign}{sign}

\DeclareMathOperator{\Int}{Int}
\DeclareMathOperator{\Cl}{Cl}
\DeclareMathOperator{\FreeSpan}{FreeSpan}
\DeclareMathOperator{\Span}{Span}

\newcommand {\sr}{\stackrel}

\newcommand {\eqv}{\equiv}

\newcommand {\ol}{\overline}

\newcommand {\td}{\tilde}

\newcommand {\ra}{\rightarrow}
\newcommand {\Ra}{\Rightarrow}
\newcommand {\la}{\leftarrow}
\newcommand {\La}{\Leftarrow}

\newcommand {\Lra}{\Leftrightarrow}
\newcommand {\ral}{\longrightarrow}
\newcommand {\lal}{\longleftarrow}

\newcommand {\ub}{\underbrace}

\newcommand {\ep}{\epsilon}
\newcommand {\vep}{\varepsilon}
\newcommand {\vphi}{\varphi}
\newcommand {\del}{\partial}
\newcommand {\Ld}{\Lambda}
\newcommand {\ld}{\lambda}

\newcommand {\al}{\alpha}

\newcommand{\wt}{\widetilde}
\newcommand{\wh}{\widehat}

\newcommand{\modf}{\!\!\!\mod}

\newcommand {\bfr}{\begin{flushright}}
\newcommand {\efr}{\end{flushright}}
\newcommand {\bfl}{\begin{flushleft}}
\newcommand {\efl}{\end{flushleft}}

\newcommand {\nn} {\nonumber}

\newcommand {\txt}{\textrm}
\newcommand {\bd}{\begin{document}}
\newcommand {\ed}{\end{document}}

\newcommand {\bt}{\begin{tikzcd}}
\newcommand {\et}{\end{tikzcd}}
\newcommand {\be}{\begin{equation}}
\newcommand {\ee}{\end{equation}}
\newcommand {\bea}{\begin{eqnarray}}
\newcommand {\eea}{\end{eqnarray}}
\newcommand {\ba}{\begin{aligned}}
\newcommand {\ea}{\end{aligned}}

\newcommand{\bbib}{\begin{thebibliography}{99}}
\newcommand{\ebib}{\end{thebibliography}}
\newcommand {\bab}{\begin{abstract}}
\newcommand {\eab}{\end{abstract}}

\newcommand {\bc}{\begin{center}}
\newcommand {\ec}{\end{center}}
\newcommand {\bit}{\begin{itemize}}
\newcommand {\eit}{\end{itemize}}
\newcommand{\blue}{\textcolor{blue}}
\newcommand{\red}{\textcolor{red}}
\newcommand{\magenta}{\textcolor{magenta}}

\def\A{\mathcal{A}}\def\B{\mathcal{B}}\def\C{\mathcal{C}}\def\D{\mathcal{D}}\def\E{\mathcal{E}}\def\F{\mathcal{ F}}\def\G{\mathcal{G}}\def\H{\mathcal{H}}\def\I{\mathcal{I}}
\def\J{\mathcal{J}}\def\K{\mathcal{K}}\def\L{\mathcal{L}}\def\M{\mathcal{M}}\def\N{\mathcal{N}}\def\O{\mathcal{ O}}\def\P{\mathcal{P}}\def\Q{\mathcal{Q}}\def\R{\mathcal{R}}
\def\S{\mathcal{S}}\def\T{\mathcal{T}}\def\U{\mathcal{U}}\def\V{\mathcal{V}}\def\W{\mathcal{W}}\def\X{\mathcal{ X}}\def\Y{\mathcal{Y}}\def\Z{\mathcal{Z}}
\def\fp{\mathfrak{p}} \def\fq{\mathfrak{q}} \def\fm{\mathfrak{m}} \def\fn{\mathfrak{n}}
\def\mfp{\mathfrak{p}} \def\mfq{\mathfrak{q}} \def\mfm{\mathfrak{m}} \def\mfn{\mathfrak{n}}
\def\abar{\bar{a}}\def\bbar{\bar{b}}\def\cbar{\bar{c}}\def\dbar{\bar{d}}\def\ebar{\bar{e}}\def\fbar{\bar{f}}
\def\gbar{\bar{g}}\def\ibar{\bar{i}}\def\jbar{\bar{j}}\def\kbar{\bar{k}}\def\lbar{\bar{l}}
\def\mbar{\bar{m}}\def\nbar{\bar{n}}\def\obar{\bar{o}}\def\pbar{\bar{p}}\def\qbar{\bar{q}}\def\rbar{\bar{r}}
\def\sbar{\bar{s}}\def\tbar{\bar{t}}\def\ubar{\bar{u}}\def\vbar{\bar{v}}\def\wbar{\bar{w}}
\def\xbar{\bar{x}}\def\ybar{\bar{y}}\def\zbar{\bar{z}}
\def\Abar{\bar{A}}\def\Bbar{\bar{B}}\def\Cbar{\bar{C}}\def\Dbar{\bar{D}}\def\Ebar{\bar{E}}\def\Fbar{\bar{F}}
\def\Gbar{\bar{G}}\def\Hbar{\bar{H}}\def\Ibar{\bar{I}}\def\Jbar{\bar{J}}\def\Kbar{\bar{K}}\def\Lbar{\bar{L}}
\def\Mbar{\bar{M}}\def\Nbar{\bar{N}}\def\Obar{\bar{O}}\def\Pbar{\bar{P}}\def\Qbar{\bar{Q}}\def\Rbar{\bar{R}}
\def\Sbar{\bar{S}}\def\Tbar{\bar{T}}\def\Ubar{\bar{U}}\def\Vbar{\bar{V}}\def\Wbar{\bar{W}}
\def\Xbar{\bar{X}}\def\Ybar{\bar{Y}}\def\Zbar{\bar{Z}}
\def\0bar{\bar{0}}\def\1bar{\bar{1}}\def\2bar{\bar{2}}\def\3bar{\bar{3}}\def\4bar{\bar{4}}\def\5bar{\bar{5}}
\def\6bar{\bar{6}}\def\7bar{\bar{7}}\def\8bar{\bar{8}}\def\9bar{\bar{9}}
\def\cprod{\mathop{\prod\hspace{-0.457cm}\times}}
\def\tprod{\mathop{\prod\hspace{-0.457cm}\otimes}}
\def\dsum{\mathop{\sum\hspace{-0.458cm}\oplus}}
\def\Tr{\txt{Tr}}\def\tr{\txt{tr}}\def\str{\txt{Str}}
\def\bde{\begin{description}}\def\ede{\end{description}}

\newtheorem{thm}{Theorem}[section]
\newtheorem{lmm}[thm]{Lemma}
\newtheorem{crl}[thm]{Corollary}
\newtheorem{prp}[thm]{Proposition}
\newtheorem{algorithm}[thm]{Algorithm}
\newtheorem{caution}[thm]{Caution}
\newtheorem{claim}[thm]{Claim}
\newtheorem{fact}[thm]{Fact}
\newtheorem{solution}[thm]{Solution}
\newtheorem{conjecture}[thm]{Conjecture}
\newtheorem{assumption}[thm]{Assumption}

\theoremstyle{definition}
\newtheorem{dfn}{Definition}[section]
\newtheorem{dfns}{Definitions}[section]
\newtheorem{example}{Example}[section]
\newtheorem{examples}{Examples}[section]
\newtheorem{rmk}{Remark}[section]
\newtheorem{rmks}{Remarks}[section]
\newtheorem{problem}{Problem}[section]
\newtheorem{problems}{Problems}[section]
\newtheorem{algorithms}{Algorithms}[section]
\newtheorem{question}{Question}[section]
\newtheorem{rquestion}{Research question}[section]
\newtheorem{questions}{Questions}[section]
\newtheorem{rquestions}{Research questions}[section]
\newtheorem{recall}{Recall}[section]
\newtheorem{recalls}{Recalls}[section]
\newtheorem{exercise}{Exercise}[section]
\newtheorem{exercises}{Exercises}[section]
\newtheorem{note}{Note}[section]
\newtheorem{notes}{Notes}[section]
\newtheorem{notation}{Notation}[section]
\newtheorem{convention}{Convention}[section]
\newtheorem{claims}{Claims}[section]
\newtheorem{facts}{Facts}[section]
\newtheorem{assumptions}{Assumptions}
\newtheorem{summary}{Summary}
\newtheorem{consequence}{Consequence}
\newtheorem{consequences}{Consequences}

\theoremstyle{definition}\newtheorem*{thm*}{Theorem}
\newtheorem*{lmm*}{Lemma}                     \newtheorem*{lmms*}{Lemmas}
\newtheorem*{dfn*}{Definition}                \newtheorem*{dfns*}{Definitions}
\newtheorem*{crl*}{Corollary}
\newtheorem*{prp*}{Proposition}
\newtheorem*{example*}{Example}               \newtheorem*{examples*}{Examples}
\newtheorem*{rmk*}{Remark}                    \newtheorem*{rmks*}{Remarks}
\newtheorem*{problem*}{Problem}               \newtheorem*{problems*}{Problems}
\newtheorem*{algorithm*}{Algorithm}           \newtheorem*{algorithms*}{Algorithms}
\newtheorem*{question*}{Question}             \newtheorem*{questions*}{Questions}
\newtheorem*{rquestion*}{Research question}
\newtheorem*{rquestions*}{Research questions}
\newtheorem*{recall*}{Recall}                 \newtheorem*{recalls*}{Recalls}
\newtheorem*{exercise*}{Exercise}             \newtheorem*{exercises*}{Exercises}
\newtheorem*{note*}{Note}                     \newtheorem*{notes*}{Notes}
\newtheorem*{notation*}{Notation}
\newtheorem*{convention*}{Convention}         \newtheorem*{conventions*}{Conventions}
\newtheorem*{caution*}{Caution}
\newtheorem*{claim*}{Claim}                   \newtheorem*{claims*}{Claims}
\newtheorem*{fact*}{Fact}                     \newtheorem*{facts*}{Facts}
\newtheorem*{assumption*}{Assumption}         \newtheorem*{assumptions*}{Assumptions}
\newtheorem*{summary*}{Summary}
\newtheorem*{consequence*}{Consequence}       \newtheorem*{consequences*}{Consequences}

\numberwithin{section}{chapter}
\numberwithin{equation}{chapter}

\makeatletter
\let\sv@thm\@thm
\def\@thm{\let\indent\relax\sv@thm}
\makeatother

\makeatletter
\let\@afterindenttrue\@afterindentfalse
\makeatother

\setul{1.5pt}{.4pt}

\AtBeginDocument{\renewcommand{\contentsname}{Table of contents}}        

\makeatletter
\patchcmd{\@tocline}
  {\hfil}
  {\leaders\hbox{\,.\,}\hfil}
  {}{}
\makeatother

\AtBeginDocument{\renewcommand{\bibname}{References}}

\makeatletter
\renewcommand{\l@chapter}{\@tocline{0}{8pt plus1pt}{0pt}{}{\bfseries}}
\makeatother

\DeclareMathOperator{\bigutimes}{\textstyle\bigcup\bigutimeskern\times}
\newcommand{\bigutimeskern}{%
  \mkern-17.5mu
  \mathchoice{}{}{\mkern0.2mu}{\mkern0.5mu}%
  }
\DeclareMathOperator{\utimes}{\cup\utimeskern\times}
\newcommand{\utimeskern}{%
  \mkern-13mu
  \mathchoice{}{}{\mkern0.2mu}{\mkern0.5mu}%
}

\newenvironment{myquote}[1]%
  {\list{}{\leftmargin=#1\rightmargin=#1}\item[]}%
  {\endlist}

\makeatletter
\@namedef{subjclassname@2020}{%
  \textup{2020} Mathematics Subject Classification}
\makeatother
\newcommand\blfootnote[1]{%
  \begingroup
  \renewcommand\thefootnote{}\footnote{#1}%
  \addtocounter{footnote}{-1}%
  \endgroup
}

%% file: parts/front/TitlePage.tex
\makeatletter
\patchcmd{\@maketitle}{\newpage}{}{}{}  
\makeatother
\title{Basic Set Theory and Algebra: Hints on Representation, Topology, Geometry, Analysis}
\author{\large ~\\[4em]
Earnest Akofor \\[7em]


(Low prerequisite Algebra notebook\blfootnote{
\begin{center}
{\scriptsize
\begin{align}
\txt{\shortstack[l]{
03E20:  Other classical set theory (including functions,\\
~\hspace{1cm} relations, and set algebra)\\
03E30:  Axiomatics of classical set theory and its fragments\\
03E25:  Axiom of choice and related propositions\\
03E75:  Applications of set theory\\
06A05:  Total orders \\
06A06:  Partial orders, general\\
18A30:  Limits and colimits\\
18A05:  Definitions and generalizations in theory of categories\\
18E05:  Preadditive, additive categories\\
18E10:  Abelian categories, Grothendieck categories\\
18E20:  Categorical embedding theorems\\
18E35:  Localization of categories, calculus of fractions\\
18G05:  Projectives and injectives (category-theoretic aspects)\\
18G10:  Resolutions; derived functors (category-theoretic aspects)\\
18G15:  Ext and Tor, generalizations, Künneth formula \\
~\hspace{1cm} (category-theoretic aspects)}}
~~~~
\txt{\shortstack[l]{
18G20:  Homological dimension (category-theoretic aspects)\\
18G80:  Derived categories, triangulated categories\\
20F29:  Representations of groups as automorphism groups \\
~\hspace{1cm} of algebraic systems\\
16W22:  Actions of groups and semigroups; invariant theory \\
~\hspace{1cm} (associative rings and algebras)\\
15A21:  Canonical forms, reductions, classification\\
16G10:  Representations of associative Artinian rings\\
54B30:  Categorical methods in general topology\\
54B40:  Presheaves and sheaves in general topology\\
18F60:  Categories of topological spaces and continuous \\
~\hspace{1cm} mappings\\
46E25:  Rings and algebras of continuous, differentiable or\\
~\hspace{1cm} analytic functions\\~\vspace{1cm}~}}\nn
\end{align}
}
\end{center}
})

%
%
%

}



\maketitle

%% file: parts/front/Copyright.tex
\thispagestyle{empty}
{
\begin{center}
  \vspace*{3in}
  \copyright \hspace{5mm} Copyright 2020 Earnest Akofor\\[1.5em]

\end{center}
\begin{myquote}{0.2in}
{\flushleft These} notes may be used in any form (including edited/modified forms) by an individual to jump-start or improve their studies. The notes may not be used for commercial purposes or propaganda.
\end{myquote}
}
\newpage

%% file: parts/front/Dedication.tex
\vskip 2cm
\begin{center}
\textbf{}\\[1.5em]
\end{center}
\bc In memory of the lost ones. \ec

\vspace{1cm}

\bc For the school I promised, longed, and slaved in vain to give my people. \ec

\vspace{1cm}

\bc To my primary school teacher Samson Mbamgah. \ec

\vspace{1cm}

\bc To my parents Esther Nkimi and Barnabas Ngocho. \ec

\vspace{1cm}

\newpage 

%% file: parts/front/Ackno.tex
\vskip 2cm
{\doublespacing
\begin{center}
\textbf{Acknowledgments}\\[1.5em]
\end{center}
\begin{myquote}{0.4in}
I am grateful to all real instructors whose classes I attended at Syracuse University.

{\flushleft I} am equally grateful for the efforts of all selfless persons who try to share (especially on the internet in various ways) every little they have, without which I could not have learned much, let alone found something to write down.
\end{myquote}
}
\newpage 

%% file: parts/front/Abstract.tex
\newpage
{\doublespacing
\thispagestyle{empty}
  \begin{center}
    {\large Abstract}\\[1em]
  \end{center}
{
\begin{myquote}{0.9in}
In these self-contained low prerequisite introductory notes we first present (in part 1) basic concepts of set theory and algebra without explicit category theory.

{\flushleft We} then present (in part 2) basic category theory involving a somewhat detailed discussion of system limits and the exact imbedding of abelian categories.

{\flushleft This} is followed (in part 3) by a discussion of localization, homological algebra, and generalizations of additive and abelian categories such as triangulated and derived categories.

{\flushleft Based} on the concepts of category theory from parts 2 and 3, (in part 4) we provide hints  for constructive discussions on familiar mathematics such as representation theory and topological geometry/analysis (i.e., topology-based geometry/analysis).

{\flushleft If} events permit, the notes will be revised/updated regularly.
\end{myquote}
}
}
{\let\thefootnote\relax}
\newpage

%% file: parts/front/Contents.tex
\tableofcontents
\listoffigures
\listoftables

%% file: parts/front/Introduction.tex
\setcounter{chapter}{-1}
\chapter{Introduction}\label{Introduction}

\section{On the frequent use of ``word$_1$/word$_2$'' (in singular first person)}
Throughout, ``word$_1$/word$_2$'' is a shorthand for ``word$_1$ and/or word$_2$'' (i.e., ``word$_1$ and word$_2$'' if the words convey different relevant meanings, or, ``word$_1$ or word$_2$'' if both words convey the same relevant meaning). Due to ambiguity in the English language caused by the multiplicity of words whose meanings overlap differently in various contexts, I will very often deliberately use ``word$_1$/word$_2$'' simply for convenience.

If the reader completely ignores context (which can be reasonably safe in many cases) and thus assumes that every word has a different meaning then ``word$_1$/word$_2$'' becomes simply ``word$_1$ and word$_2$''. However, in a phrase of the form ``... a word$_1$/word$_2$ ...'' it will often (but not always) be the case that ``word$_1$/word$_2$'' means ``word$_1$ or word$_2$''.

Therefore, it is best for the reader to simply get comfortable with the dual meaning of ``word$_1$/word$_2$''. Fortunately enough, the appropriate meaning of ``word$_1$/word$_2$'' will often be immediately clear either from preceding words, or from following words, or from both preceding and following words.

\section{Motivation and objectives (in singular first person)}\label{MotObjSec}

These are merely my personal notes on some basic aspects of (algebra-oriented) mathematics, and not a specialized textbook on any particular area of the vast subject of mathematics. While a ``graduate student'' at Syracuse University from 2015-2020, I took notes as carefully as possible during several mathematics course lectures given by various instructors at different times. As usual, these lecture notes contained interesting ideas/outlines (\blue{footnote}\footnote{The missing details were often about results from set theory and category theory, which are central areas of mathematics that nevertheless are (for yet mysterious reasons) often either only briefly mentioned in depending courses or ignored altogether.}) for the somewhat tedious proofs of certain well-known standard results in algebra, examples of which are the K\"unneth formulas in algebraic topology. Since I was more interested in the proof details than in the results themselves or their uses, I decided afterwards to undertake the challenge of efficiently reconstructing the omitted connections that help explain such intriguing results.

The notes should be self-contained as far as concepts/proofs are concerned and therefore best for personal study but, of course, (along with needed supplements/exercises) can also be easily used to teach a course. The notes assume nothing beyond high school mathematics/logic along with a certain degree of maturity and, of course, the willingness to think/reason in a flexible manner towards intelligent abstraction (\blue{footnote}\footnote{Here, ``\emph{abstraction}'', which is essentially ``\emph{efficient generalization via symbolic (rather than concrete) representation}'' is done in a carefully guided way, so that it is not merely ``\emph{nonsense}''.}). As in an old saying, ``\emph{it is best to make a fisherman out of a serious fish consumer}'', my main goal is to provide stage setting guidance/tools, the best way I can, on how to learn/do algebra more efficiently and independently, but not necessarily to attempt the obviously unattainable task of ``teaching'' algebra itself. If someone successfully learns/understands these notes then, I believe, such a person can easily follow/understand an endless lot of the so-called ``advanced stuff'' in mathematics. In any case, my greatest hope is that the notes can help a gap hopper to healthily progress into a bridge builder.

The notes present algebra concepts, hopefully, in a ``purified'' manner in the sense that they avoid (as much as possible) equally important but often distracting/uncertain details such as historical accounts, award of credit, credibility of sources, applicability, and more. The notes are intended to be an efficient attempt to reveal what I believe makes algebra work (or lies at the heart of algebra), notably, the roles of set theory and category theory, but neither who is credited for what piece of algebra, nor what algebra can do. Therefore, the notes are best for persons who by nature like/enjoy doing or sharing mathematics even on their own (i.e., even with little or no appreciation/encouragement/support from others).

Recognize however that it is only through application (e.g., working on exercises from mathematics texts/literature) that one can really begin to appreciate the strength/significance of these fundamentals. The emphasis here, quite contrary to common practice, is on first getting a fairly extensive framework of the fundamentals before attempting to apply/evaluate them. That is, the ``motto'' is to ``get the idea first'' (which is by no means a trivial task, and so might require patience/endurance depending on one's natural capability or background), and then skillful/enjoyable exercise is likely to follow afterwards.

\ul{Important note}: Names of persons and all references are used strictly as reasonably convenient guiding labels (for easily/efficiently navigating through the literature), and no crediting whatsoever is implied, understood, or even recognized by me. To me, relevant/important/great persons are simply those who are \ul{each} consciously/willfully doing their very best in whatever situation they happen to find themselves (irrespective of whether or not any of their efforts are ever going to be appreciated). This is what I call justice.

\section{About errors and revision (in singular first person)}
These notes have no deliberate/intentional mistakes, but of course they might have typos (including unforseen ambiguities/incompletes) that can easily be misunderstood/misinterpreted by some readers and so become real errors. In that case, feel free to create an improved version of the notes for yourself. If circumstances permit, the notes will be revised/updated regularly.

\ul{Important note}: In seemingly rare occasions, I find myself assuming incorrectly that some well known fact is trivial (probably because of false impressions from handwaving arguments by instructors and writers), and thus give an inadequate or even incorrect proof for that fact. I am highly confident that this version of the notes is free from such mistakes, but one must never be too certain. Therefore, if you find such an unfortunate mistake in the notes, then please send me a brief email, using \textcolor{blue}{eakofor@gmail.com}, simply indicating the particular result whose use/proof is problematic; e.g., an email  with subject ``\emph{Problem in notes}'' and body ``\emph{Check the proof of $\cdots$ on page(s) $\cdots$.}''.

Also, for brevity, I sometimes cannot avoid intentionally making a formal statement, i.e., an appealing claim that is false in general but which is true under certain conditions that should be reasonably clear or obvious to the reader. Therefore, a reader who has recognized a formal statement might consider it an easy exercise to attempt to make sense of the statement.

\section{Scientific view of the notebook}
The \emph{nature} of a thing (physical \emph{object} or \emph{system}) refers to those aspects/properties of it that do not depend on any particular/individual viewpoint on (observation of) it. \emph{Science} is an attempt to find the nature of a thing by putting together various viewpoints/perceptions of it, including its different compositions in terms of smaller things (e.g., its structural \emph{appearances}), its roles as part of larger things (e.g., its \emph{origins}), its effects through interaction with other things (i.e., its \emph{environment}), etc.

\emph{Mathematics} is a simplified$\slash$mechanical pseudo-programming language that serves as a means of description and communication of science. After ignoring overlap, mathematics consists roughly of three major parts: \emph{analysis} (object \emph{size} measurements), \emph{geometry} (object \emph{shape} classifications), \emph{algebra} (object background$\slash$\emph{reference}$\slash$framework categorizations).

Among many others, an important application of mathematics is found in statistics: \emph{Statistics} is a method of scientific investigation that involves the learning of \emph{characteristics} of a \emph{population} (a system of objects) from \emph{data} collected through a \emph{study} on a \emph{sample} (an accessible portion) of the population. \emph{Probability}, which is a branch of analysis that concerns the quantification of randomness, is a mathematical tool for avoiding \emph{biased samples} (samples that do not sufficiently represent the population) in statistics.

These notes partly resemble a typical scientific study, with basic algebra concepts/results/proofs as the population of interest, rigor/abstraction/consistency as characteristics of interest, and references (i.e., selected relevant sources from the mathematics literature) as a pseudo-random sample from the population of interest. Therefore, in these notes, our main concern is not with any specific area or application of \emph{mathematical science}, but rather with some basic concepts that underly advanced studies in (algebra-oriented) mathematics. The basic idea is to unhide/expose the ingredient that makes mathematics outstanding as an instructional tool. This ingredient is in the form of a genuine \ul{attempt} to reveal as much of the truth about nature as possible through \emph{mathematical rigor}: At every step, we \ul{attempt} to build an \emph{unambiguous} conceptual foundation based on which certain results can then be \emph{unambiguously} established.

\section{Organization and structural summary of contents}
The notes are divided into four parts as shown in Figure \ref{NotesStruct}, in which components of the same part have the same color.

\begin{figure}[H]
{\scriptsize
\bea\label{NotesStrDiag}\adjustbox{scale=1.1}{
\bt[row sep=small, column sep=tiny]
\shortstack{
\textbf{KEY} \\~\\
\fbox{\begin{minipage}{1cm}\centering
\textcolor{OliveGreen}{Part 1} \\ \textcolor{blue}{Part 2} \\ \textcolor{magenta}{Part 3} \\ Part 4
\end{minipage}}
}
&
&&&&
&&&&
\fbox{\begin{minipage}{2cm}\centering
Representable \\ systems
\end{minipage}}\ar[from=d]\ar[dllll,bend right]\ar[drrrr,bend left]
&&&&
 \\ 
&
&&&&
\fbox{\begin{minipage}{1cm}\centering
\textcolor{OliveGreen}{Classes}
\end{minipage}}\ar[ddllll,bend right,"\txt{orders}"description]\ar[dd,"\txt{axioms}"']\ar[rrrr]
&&&&
\fbox{\begin{minipage}{1.5cm}\centering
\textcolor{blue}{Categories}
\end{minipage}}\ar[ddllll,hook,dashed,Rightarrow,"\txt{Imbed}"description]\ar[dd]\ar[rrrr]\ar[ddrrrr]
&&&&
\fbox{\begin{minipage}{2cm}\centering
Topological \\ systems
\end{minipage}}\ar[dd,"\txt{analogy}"]
 \\ 
&
&&&&
&&&&
&&&&
 \\ 
&
\fbox{\begin{minipage}{1.25cm}\centering
\textcolor{OliveGreen}{Number} \\ \textcolor{OliveGreen}{systems}
\end{minipage}}\ar[d]
&&&&
\fbox{\begin{minipage}{0.7cm}\centering
\textcolor{OliveGreen}{Sets}
\end{minipage}}\ar[from=llll]\ar[d]
&&&&
\fbox{\begin{minipage}{1.5cm}\centering
\textcolor{blue}{Additive} \\ \textcolor{blue}{categories}
\end{minipage}}\ar[d]\ar[drrrr,"\txt{axioms}"description]
&&&&
\fbox{\begin{minipage}{2cm}\centering
Measurable \\ systems
\end{minipage}}
  \\ 
&
\fbox{\begin{minipage}{1.25cm}
\textcolor{OliveGreen}{Natural} \\ \textcolor{OliveGreen}{Cardinal} \\ \textcolor{OliveGreen}{Ordinal}
\end{minipage}}\ar[rrrr]
&&&&
\fbox{\begin{minipage}{1.3cm}\centering
\textcolor{OliveGreen}{Groups} \\ \textcolor{OliveGreen}{Rings} \\ \textcolor{OliveGreen}{Modules} \\ \textcolor{OliveGreen}{Algebras}
\end{minipage}}\ar[d]
&&&&
\fbox{\begin{minipage}{1.5cm}\centering
\textcolor{blue}{Abelian} \\ \textcolor{blue}{categories}
\end{minipage}}\ar[dllll,hook,dashed,Rightarrow,"\txt{Imbed}"description]\ar[d,"\txt{Localize}"]\ar[rrrr]
&&&&
\fbox{\begin{minipage}{2cm}\centering
\textcolor{magenta}{Triangulated} \\ \textcolor{magenta}{categories}
\end{minipage}}\ar[d,"\txt{Localize}"]
  \\ 
&
&&&&
\fbox{\begin{minipage}{1.2cm}\centering
\textcolor{OliveGreen}{Modules}
\end{minipage}}
&&&&
\fbox{\begin{minipage}{1.5cm}\centering
\textcolor{magenta}{Derived} \\ \textcolor{magenta}{categories}
\end{minipage}}\ar[rrrr]
&&&&
\fbox{\begin{minipage}{2cm}\centering
\textcolor{magenta}{Triangulated} \\ \textcolor{magenta}{Derived} \\ \textcolor{magenta}{categories}
\end{minipage}}
\et}\nn
\eea}
\caption{~~Structural summary of the main contents}
\label{NotesStruct}
\end{figure}

Part 1 introduces set theory and algebra concepts without explicit category theory. Both the construction and arithmetic of natural, cardinal, and ordinal numbers are discussed to some extent. Also introduced right from the beginning of part 1 is the concept of classes, a generalization of sets, which serves as an indication that there are equally useful objects/structures beyond relatively familiar set-based objects/structures such as groups, rings, modules, algebras, topological spaces, etc. These other useful objects/structures include (class-based) graphs and categories.

Part 2 introduces basic initial objects/structures (e.g., graphs, categories, functors) in category theory along with commonly used descriptive tools/operations such as diagrams, systems, and limits of systems. Because of their central role, abelian categories, along with their exact imbedding into module categories are discussed somewhat extensively.

Part 3 similarly introduces basic homological objects/structures (e.g., the long exact homology sequence, derived functors, homological functors) in category theory along with commonly used descriptive tools/operations such as localization, resolution, and triangulation in categories. In particular, generalizations of additive and abelian categories such as triangulated, derived, and triangulated derived categories are also discussed.

Part 4 is application oriented and provides hints for constructive categorical discussions on familiar part 1 based mathematics with the help of the ideas/methods introduced in part 2 and part 3. Part 4 is further divided into two main parts, namely, (A) representation theory and (B) geometry/analysis. Part 4(A) discusses symmetry operations, finitely generated modules, and semisimple rings/modules. Part 4(B) discusses some applications of topology, concepts of measure theory, and topological analysis.

\newpage
\section{Prerequisite and recommendable attitudes}
As already mentioned in section \ref{MotObjSec}, the notes are (i) self-contained as far as concepts/proofs are concerned, and (ii) assume nothing beyond high school mathematics/logic along with a certain degree of maturity and the willingness to be flexible with abstraction. Some college or graduate level mathematics might be helpful, but the notes are designed with the assumption that the reader has no streamlined previous knowledge. Therefore, it is recommendable for a reader with prior knowledge to attempt to view these notes as somehow a trimmed/cleaned up review of their knowledge.

It is possible that our notes, just like anything else, are not perfect for every single reader. Therefore, if necessary, it is also recommendable for a reader to freely use their own intuition to fill-in/account for whatever they consider to be relevant but missing details. Such details can include simple exercises, alternative viewpoints, and additional sources of information/knowledge.

Finally, any result that is stated but not explicitly proved in the notes is either (1) reasonably obvious, as it is often the case with corollaries, or (2) not strictly relevant for the rest of the discussion. Therefore, the reader can freely ignore those very few/rare unproved results of the type (2), which are included only for general knowledge and as side information that might be useful to readers with wider areas of interest than the notes provide.

\section{Terminology, notation, conventions (important)}
The word ``object'' will be our preferred synonym for the word ``thing''. Given a \emph{collection} of objects (or simply a \emph{collection}), each object in the collection will be called a \emph{member} of the collection. A collection is itself an object, but in general, not of the same type as its own members. The specification/defining rules for an object are called its \emph{structure}. Depending on its structure type, an object may be called a \emph{set}, \emph{class}, \emph{category}, \emph{graph}, etc (to be discussed later), where classes include sets, categories include classes, graphs include categories, and so on. In particular, a \emph{set} (an object in the category of sets) is a special type of collection that is also an object of the same type as its own members. A member of a set is called an \emph{element} of the set. A rule between two collections (which are not necessarily distinct) that uses the members of one collection to label/mark some of the members of the other collection is called a \emph{map} (or \emph{mapping}, \emph{indexing}, \emph{function}, \emph{operation}) between the said collections. A rule (e.g., a map) between two objects that preserves structure is called a \emph{morphism}. A morphism between categories is called a \emph{functor}.

Capital letters such as $A,B,\cdots$ will denote sets (and later, objects in a category and functors between categories), meanwhile "curly" capital letters such as $\A,\B,\cdots$ will denote \emph{families} or collections of collections (and later, categories).  Lower case letters will denote members of a collection (and later, morphisms in a category), i.e., an arbitrary member of a collection $A$ will be denoted by $a$.

Membership and non-membership in a collection will be denoted by $\in$ and $\not\in$ respectively. For example, if $C$ is a collection based on an entry requirement $\phi$, written as $C:=\{c~|~\txt{$c$ satisfies $\phi$}\}$ or $C:=\{c:\txt{$c$ satisfies $\phi$}\}$, then ``$c\in C$'' means ``$c$ satisfies $\phi$'' or ``$c$ belongs in $C$'' or ``$c$ is a member of $C$'', etc, and similarly ``$c\not\in C$'' means ``$c$ does not satisfy $\phi$'' or ``$c$ does not belong in $C$'' or ``$c$ is not a member of $C$'', etc.

The symbols $\{~\}$, $(~)$, $[~]$ and more will be used to denote/enclose specifications for collections; $\{~\}$ for \emph{unindexed collections} (collections whose members are not distinguished by \emph{any labels} found on them), $(~)$ for \emph{indexed collections} (collections whose members are distinguished by \emph{specific labels} found on them), and $[~]$ for \emph{partitioning/equivalence collections} (unindexed collections that arise automatically through some type of equivalence relation between their members). (\blue{footnote}\footnote{Labels are also called ``indices'', but for us ``indices'' on the members of a collection do not automatically imply that the collection is an indexed collection.}). Thus, a \emph{labeled} collection of the form $(C_\al)_{\al\in A}$ is unambiguously an indexed collection. However, a \emph{labeled} collection of the form $\{C_\al\}_{\al\in A}$ can be considered to be indexed or unindexed and, in situations where a distinction between the two is essential, the correct option will be clear from the context. (\blue{footnote}\footnote{This is deliberately done in order to prevent more serious difficulties, including poor readability, that can arise due to overuse of $(~)$. For example, depending on the context, in an indexed collection $(F_x)_{x\in X}$, if we replace $F_x$ with $F(x)$, then the result $(F(x))_{x\in X}$ might look better in the alternative form $\{F(x)\}_{x\in X}$.})

Inclusion (or containment) between collections will be denoted by $\subset$, i.e., if $A,B$ are collections, then $A\subset B$ (or $B\supset A$) means ``$A$ is included (or contained) in $B$'' or ``$B$ includes (or contains) $A$'' or ``every member of $A$ is a member of $B$'', etc. Concerning the strictness of inclusion, our preferred inclusion symbols are $\subset$ for non-strict inclusion, and $\subsetneq$ for strict inclusion. However, when dealing with \emph{partially ordered} collections (collections with a weak or partial means of comparison between them), we will for convenience also use $\subseteq$ (in place of $\subset$) for non-strict inclusion. That is $\subset,\subseteq$ will each denote ``non-strict inclusion'' while $\subsetneq$ will denote ``strict inclusion''.

Given collections $X$ and $Y$, if $f:X\ra Y$, $x\mapsto f(x)$ is a \emph{map} (\emph{mapping}, \emph{indexing}, \emph{function}, \emph{operation}) between them (i.e., a \emph{labeling} rule that assigns to each object $x\in X$ a \emph{single} object $f(x)\in Y$ called the \emph{value} of $f$ at $x$), then the value of $f$ at $x$ is denoted by $f(x)$ or $f_x$ or $f^x$ depending on convenience, where the labeled subcollection $\{f(x)\}_{x\in X}:=\{f(x)~|~x\in X\}$ of $Y$ is called the \emph{image} of $f$ and the map $f$ itself is sometimes viewed as the indexed collection $\big(f(x)\big)_{x\in X}$ when the identity of the target collection $Y$ is understood/unimportant. The same notation applies to all labeled collections (whether indexed or not). For example, we can write $\C=\{C(\al)\}_{\al\in A}$ or $\C=\{C_\al\}_{\al\in A}$ or $\C=\{C^\al\}_{\al\in A}$. For \emph{finite/countable} collections (collections with at most as many members as there are the natural numbers $\Natural$ or integers $\Integer$), we will often write $\A=\{A_i\}_{i\in C}$, for a finite/countable set $C$. That is, we will mostly (but not always) use $\al,\beta,\gamma,\cdots$ for arbitrary labels and mostly (but not always) use $i,j,k,l,\cdots$ for finite/countable labels. In cases where the finite/countable set $C\subset\Integer:=\{\cdots,-1,0,1,\cdots\}$, we will often (but not always) use $m,n,r,s,t,\cdots$, e.g., $\C=\{A_n\}_{n=1}^\infty$.

Any specialized notation that deviates from the above basic rules/conventions will be explicitly stated (along with any associated terminology) in definitions. Such specialized notation is common in all sections/chapters/parts. Note that all sections/chapters/parts adopt and use notation introduced in earlier sections/chapters/parts, often without any explicit reminder.

We will sometimes abbreviate certain frequent phrases and words (to be defined later if necessary) as in the following table (where the list of symbols is only a typical sample that is not meant to be exhaustive).

\begin{table}[H]
  \centering
\adjustbox{scale=0.8}{
\begin{tabularx}{17cm}{|||c|X|||c|X|}
    \hline
    \textbf{Symbol} & \textbf{Phrase}
    & \textbf{Symbol} & \textbf{Phrase}\\
     \hline\hline
\hline &&&\\
$\exists$ &  there exists
  &    
 s.t. &    such that
   \\  
\hline &&&\\
$\nexists$  &   there does not exist
  &    
 a.e. &    almost everywhere
   \\  
\hline &&&\\
    $\exists!$  &    there exists a unique
  &    
 a.e.f &    all except finitely many
   \\  
\hline &&&\\
$:=$  &  defined to be
  &    
 wlog &   without loss of generality
   \\  
\hline &&&\\
$\Ra$  &    implies
  &    
 nbd &    neighborhood
   \\  
\hline &&&\\
$\substack{\txt{iff} \\~\\ (\txt{or}~\iff)}$  &    if and only if
  &    
 resp. &    respectively
   \\  
   \hline &&&\\
$|$  & restricted to, ~divides
  &    
 equiv. &    equivalently
   \\  
\hline &&&\\
$:$   &  restricted to, ~such that
  &    
 const. &  constant
   \\  
  \hline &&&\\
$\ll$  &    much less than
  &    
 hty &    homotopy
   \\  
\hline &&&\\
$\gg$  &    much greater than
  &    
 endo &    endomorphism
   \\  
\hline &&&\\
$\approx$  & approximately equal to, cardinally equivalent to
  &    
 mono &    monomorphism
   \\  
\hline &&&\\
$\simeq$  &  homotopic to, homotopy equivalent to
  &    
 epi &    epimorphism
   \\  
\hline &&&\\
$\cong$  &  isomorphic to
  &    
 iso &    isomorphism
   \\  
\hline\hline
\end{tabularx}}
  \caption{\textbf{Common math abbreviations for certain phrases and words}}\label{PhrAbbTab}
\end{table}

The symbol $\cong$ will denote \ul{``isomorphism''} (of all types of objects), the symbol $\sim$ (among others) will denote an \ul{``equivalence relation''}, the symbol $\simeq$ (among others) will denote \ul{``(homotopy) equivalence of objects''} (e.g., chain complexes/morphisms) in a category, and the symbol $\approx$ will denote either \ul{``approximate equality''} of numbers and numerical estimates or \ul{``cardinal equivalence''} of classes of objects.

In categorical diagrams, we will use the following conventions for morphisms (arrows):
\bit[leftmargin=0.5cm]
\item Morphism~~ \bt X\ar[r,"f"]& Y\et
\item Induced (implied/automatic/existential) morphism~~ \bt X\ar[r,dashed,"f"]& Y\et ~~or~~ \bt X\ar[r,dotted,"f"]& Y\et
\item Monic morphism~~ \bt X\ar[r,hook,"f"]& Y\et
\item Epic morphism~~ \bt X\ar[r,two heads,"f"]& Y\et
\item Monic epic morphism~~ $\bt X\ar[r,hook,two heads,"f"]& Y\et$  ~~(will also be called a ``unitary morphism'')
\eit
In an \ul{abelian category}, to be defined later, consider a sequence (of objects and arrows) of the form
\[\bt \cdots\ar[r]& A\ar[r,"f"]& B\ar[r,"g"]& C\ar[r]&\cdots\et.\]
If the sequence is \ul{exact at $B$}, ``exactness'' to be defined later, we will indicate this by putting a tail on the arrow leaving $B$ as follows:
\[\bt \cdots\ar[r]& A\ar[r,"f"]& B\ar[r,tail,"g"]& C\ar[r]&\cdots\et\]
A sequence is \ul{exact} if it is exact at every object. In particular, we will frequently encounter exact versions of sequences of the form
\[\bt 0\ar[r]& A\ar[r,"f"]& B\ar[r,"g"]& C\ar[r]&0.\et\]
which by our conventions above will precisely take the form
\[\bt 0\ar[r]& A\ar[r,hook,"f"]& B\ar[r,tail,two heads,"g"]& C\ar[r]&0.\et\]
However, the latter exact sequence will very often simply be written without the tail as follows:
\[\bt 0\ar[r]& A\ar[r,hook,"f"]& B\ar[r,two heads,"g"]& C\ar[r]&0.\et\]

%% file: parts/AlgebraNC/AlgebraNC1.tex
\chapter{Basic Set Theory I: Tools and Axioms}

For this chapter, if needed, additional sources of reading include for example \cite{enderton1977,goldrei1996,gaillard2012,lewin1991}. As usual, by necessity, we will initially take for granted a few intuitively obvious properties of the natural (or counting) numbers $\Natural:=\{0,1,2,\cdots\}$ before actually establishing them later. These properties include \ul{linear/total/consecutive ordering} $0<1<2<\cdots$, (\blue{footnote}\footnote{Here, $a<b$ (or $b>a$) denotes ``$a$ is less than $b$'' (or ``$b$ is greater than $a$'') and $a\leq b$ (or $b\geq a$) denotes ``$a$ is at most $b$'' (or ``$b$ is at least $a$'').}), \ul{well-ordering} (i.e., the existence of a smallest member for each nonempty subcollection of the natural numbers $\Natural$), and \ul{induction} (i.e., the observation that given any statements $S_0,S_1,\cdots$, we can conclude that they are all true provided (i) $S_0$ is true and (ii) $S_{k+1}$ is true whenever either $S_k$ is true or $S_0,S_1,...,S_k$ are all true). Moreover, given any collection, we will call it a \ul{finite} collection if it has as many members as the collection of numbers $\{1,2,\cdots,n\}$ for some natural number $n\in\Natural$, otherwise, we will call it an \ul{infinite} collection.

\section{Classes, Maps, Relations, Cardinality}
\begin{rmk}[\textcolor{blue}{\index{Russell's paradox}{Russell's paradox}}]\label{RusParRmk}
Naively, a \emph{set} is simply \ul{any} collection of objects (or things) called \emph{elements} of the set (where a set $\emptyset$ with no elements exists and is called the \emph{empty set}), and a subcollection of objects from a given set is called a \emph{subset} of the set. This allows for the idea/notion of ``a set of all sets, or \emph{universal set}, $u$'' (suggesting that a set, such as $u$, can be a member of itself), which leads to a contradiction as follows. Let $y:=\{x~|~x\not\in x\}\subset u$ be the subset of $u$ consisting of all sets $x$ that are not members of themselves, i.e., such that $x\not\in x$. Then $y\in y$ if and only if $y\not\in y$ (a contradiction). Consequently, the axioms (i.e., set defining rules) of set theory must be constructed such that a set of all sets does not exist, and so a set is not just \ul{any} collection of objects. It turns out (with details in \cite{enderton1977,goldrei1996} not of immediate relevance to us) that, in a consistent theory of sets, a set cannot be a member of itself (a condition that automatically rules out the existence of a set of all sets).
\end{rmk}

The following is an interpretation of Russell's paradox and other set-theoretic paradoxes (\blue{footnote}\footnote{Here, a paradox is a claim/statement that is appealing in simplicity/tidiness but also vague/ambiguous enough to be self-contradicting. Concerning paradoxes in general, there need not be a contest or conflict between reason/logic (i.e., our wishes/expectations) and reality (i.e., what we actually get). Sound reason/logic is by construction at least a step towards (if not already a foundation for) understanding reality, but to someone who is either confused or has too little knowledge, any form of reason/logic might as well look like a step away from reality instead. It might be a fact that we cannot fully understand reality, but it should not be any less of a fact that we cannot live/exist without some understanding of reality.}) as consequences of poorly designed entry or membership requirements for sets.
\begin{rmk}[\blue{Interpretation of Russell's paradox, \index{Berry's paradox}{Berry's paradox}}]
If a set of all sets $u$ existed, then as a collection we would write it as follows:
\bea
u:=\{x~|~\txt{$x$ is a set}\}.\nn
\eea
However, the entry requirement into $u$, namely, ``$x$ is a set'', is not well defined particularly because the word ``set'' is yet undefined. This is the main reason for Russell's paradox, i.e., it is simply the consequence of a poorly designed entry requirement (which is vague/ambiguous in the sense that it has more than one inequivalent interpretations that combine to yield a contradiction). The entry rule for $u$ is vague/ambiguous because its interpretation allows for two incompatible assumptions/outcomes: (i) a set can be a member of itself, and (ii) any collection of sets is a set (as a subset of $u$). But (i) and (ii) together imply that the collection $y:=\{x\in u:x\not\in x\}$ of sets that do not contain themselves is a set (which leads to a contradiction as we have seen in Remark \ref{RusParRmk}). Therefore, the construction of a consistent theory of sets must include rules for precise/unambiguous entry requirements for sets.

From (\cite[p.5]{enderton1977}), another situation revealing a set paradox (called \ul{Berry's paradox}) as a result of a poorly designed entry requirement for a set  is as follows. Consider the collection
\bea
v:=\{x~|~\txt{$x$ is a positive integer whose definition can be typed in one line of text}\}.\nn
\eea
Again, the trouble here is that the word ``integer'' (a set related object) is yet undefined, because we have not yet defined the phrase ``set of integers''. The entry rule for $v$ is vague/ambiguous because its interpretation allows for two incompatible assumptions/outcomes: (i) $v$ is a nonempty finite collection (because ``1'' or ``the least positive integer'' is in $v$ and only finitely many lines of text are possible, as the typing machine's alphabet is finite), and (ii) if $x\in v$ then $x+1\in v$ as well (because we can always make the definition of $x+1$ even shorter than that of $x=(x+1)-1$). By (i) $v$ has a greatest element $l$ (thus $l+1\not\in v$). But by (ii), $l+1\in v$ (a contradiction). Therefore, a proper set-theoretic definition of the set of integers (along with precise/unambiguous subset specification rules) is necessary.
\end{rmk}

As a consequence of the above discussion alone (even while we have not yet fully defined a set), there exist collections (such as the collection of all sets and the collection of all nonempty sets) that do not form sets as we know them so far. Such collections are called \emph{proper classes}, where a \emph{class} is any (logically) meaningful collection. Just as a (logical) set of all sets does not exist, a (logical) class of all classes does not exist. Similarly, and ultimately, a (logical) collection of all collections does not exist. Hence, in the absolute sense, there is no such \emph{thing} as \emph{everything}, i.e., \emph{everything} (whatever it is) is not a \emph{thing} or \emph{object}: Indeed, \emph{everything} includes \emph{the collection of all things}, and therefore cannot be a \emph{thing}.

\begin{dfn}[\textcolor{blue}{
\index{Class}{{Class}},
\index{Class structure}{Class structure},
\index{Empty class}{Empty class},
\index{Subclass}{Subclass},
\index{Superclass}{Superclass},
\index{Powerclass}{Powerclass},
\index{Equal! classes}{Equal classes},
\index{Distinct classes}{Distinct classes},
\index{Proper subclass}{Proper subclass},
\index{Proper superclass}{Proper superclass},
\index{Union! of classes}{Union of classes},
\index{Intersection of! classes}{Intersection of classes},
\index{Relative complement of a class}{Relative complement of a class},
\index{Complement of a class}{Complement of a class},
\index{Disjoint! classes}{Disjoint classes},
\index{Disjoint! collection}{Disjoint collection},
\index{Disjoint! union of classes}{Disjoint union of classes}}]
A \ul{class} is a collection of objects (of arbitrary type) that might have one or more properties/characteristics in common (which we will call the \ul{class structure}). In order to avoid Russell's paradox, (i) we will assume a class cannot be an object of itself, i.e., if $C$ is a class then $C\not\in C$ (\blue{footnote}\footnote{In particular, a class of all classes does not exist.}), and (ii) any collection of classes will be assumed to be ``sufficiently'' smaller than the collection of all classes (unless stated otherwise).

For the convenience of expression, we introduce empty classes, where an \ul{empty class}, denoted by $\emptyset$, is a class containing no objects whatsoever.

Let $X,Y$ be classes. The class $X$ is a \ul{subclass} of $Y$, written $X\subset Y$ (or equivalently, $Y$ is a \ul{superclass} of $X$, written $Y\supset X$), if every object of $X$ is an object of $Y$ (i.e., if $x\in X$ then $x\in Y$). The \ul{powerclass} of $X$ is the collection $\P(X):=\{X':X'\subset X\}$ of all subclasses of $X$. (\blue{footnote}\footnote{We are assuming the class $X$ is such that its powerclass can be considered a class as well.}). The classes $X$ and $Y$ are \ul{equal} (written $X=Y$) if $X\subset Y$ and $Y\subset X$. If the classes $X$ and $Y$ are not equal, then they are called \ul{distinct} or \ul{unequal} (written $X\neq Y$). The class $X$ is a \ul{proper subclass} of $Y$, written $X\subsetneq Y$ (or equivalently, $Y$ is a \ul{proper superclass} of $X$, written $Y\supsetneq X$), if $X\subset Y$ and $X\neq Y$ (i.e., $X\subset Y$ and there exists an object $y\in Y$ such that $y\not\in X$).

Let $\C$ be a collection of classes (i.e., a class of classes). The \ul{union} of the collection $\C$ is the class
\bea
\textstyle \bigcup\C:=\{c~|~c\in C~\txt{for some}~C\in\C\}.~~~~\txt{(\blue{footnote}\footnotemark).}\nn
\eea
\footnote{We are assuming the collection of classes $\C$ is such that its union can be considered a class as well.}
The \ul{intersection} of the collection $\C$ is the class
\bea
\textstyle\bigcap\C:=\{c~|~c\in C~\txt{for every}~C\in \C\}.\nn
\eea
In particular, if $X,Y$ are classes, their union is $X\cup Y:=\{z~|~z\in X~\txt{or}~z\in Y\}$, their intersection is $X\cap Y:=\{z~|~z\in X~\txt{and}~z\in Y\}$, and the \ul{relative complement} of $Y$ in $X$ is the class $Y^{c_X}:=\{x\in X~|~x\not\in Y\}$ (which is also denoted by $X-Y$ or $X\backslash Y$). If $X\supset Y$ is a superclass, the \ul{complement} of $Y$ in $X$ is the class $Y^c:=Y^{c_X}=X-Y=X\backslash Y$. In particular, if $X,Y\subset Z$ are subclasses, then we see that $X-Y=X\cap Y^{c_Z}=X\cap(Z-Y)$.

Two classes $C,C'$ are \ul{disjoint} if $C\cap C'=\emptyset$. A collection of classes $\C$ is a \ul{disjoint collection} if any two distinct/unequal classes $C,C'\in\C$, $C\neq C'$, are disjoint. If $\C$ is a disjoint collection of classes, the disjointness of $\C$ is sometimes emphasized (if necessary) in the union $\bigcup\C$ by alternatively writing the union as $\bigsqcup\C$, called a \ul{disjoint union} (where of course $\bigsqcup\C=\bigcup\C$, except when an entirely different meaning is attached to $\bigsqcup$, as done in Definition \ref{DisjUnDef}).
\end{dfn}

It is immediate that unions and intersections of classes as specified above satisfy \index{DeMorgan's laws}{\ul{DeMorgan's laws}} for sets (Lemma \ref{DemogLaws}), which states that if $X$ is a class of objects, $\C\subset\P(X)$ a collection of subclasses of $X$, and $\C^c:=\{X-C:C\in\C\}\subset\P(X)$ the collection of complements of these subclasses in $X$, then
\[
\textstyle\left(\bigcup\C\right)^c=\bigcap\C^c~~~~\txt{and}~~~~\left(\bigcap\C\right)^c=\bigcup\C^c.
\]
We will therefore use these basic laws in proofs (e.g., the proof of Theorem \ref{CSBthm}) often without explicit mention.

\begin{dfn}[\textcolor{blue}{
\index{Map (mapping, indexing, function, operation)}{{Map (mapping, indexing, function, operation)}},
\index{Value of a map}{Value of a map},
\index{Restriction! of a map}{Restriction of a map},
\index{Equal! maps}{Equal maps},
\index{Distinct maps}{Distinct maps},
\index{Extension of a map}{Extension of a map},
\index{Domain}{Domain},
\index{Codomain}{Codomain},
\index{Image of! a set}{Image of a set},
\index{Image of! (range of) a map}{Image (range) of a map},
\index{Preimage of a! set}{Preimage of a set},
\index{Fiber of a map}{Fiber of a map},
\index{Section of a map}{Section of a map},
\index{Local! section}{Local section},
\index{Constant map}{Constant map},
\index{Invariant! class}{Invariant class},
\index{Fixed point}{Fixed point},
\index{Point-invariant class}{Point-invariant class},
\index{Fix of a map}{Fix of a map},
\index{Finite! map}{Finite map},
\index{Identity! map}{Identity map},
\index{Injective! map}{Injective map},
\index{Surjective map}{Surjective map},
\index{Bijective map}{Bijective map},
\index{Inverse! map}{Inverse map},
\index{Permutation (symmetry)}{Permutation (symmetry)}}]~\\~
Let $X,Y$ be classes of objects. A \ul{map} (mapping, indexing, function, operation) $f$ from $X$ to $Y$ (denoted by $f:X\ra Y$, or $X\sr{f}{\ral}Y$, or just $X\ra Y$ when a label is not necessary) is a rule that assigns to each object $x$ of $X$ a \emph{single} object $f(x)$ of $Y$ (with the assignment rule denoted by $x\mapsto f(x)$) called the \ul{value of $f$} at $x$ (or the \ul{image of $x$} under $f$). The map $f$ from $X$ to $Y$ that assigns to each object $x$ of $X$ a \emph{single} object $f(x)$ of $Y$ is often briefly specified in the form
\bea
\bt[column sep=small] f:X\ar[r] & Y,~~x\mapsto f(x).\et\nn
\eea
If $Z\subset X$ is a subclass, the map $Z\ra Y$, $z\mapsto f(z)$ is called the \ul{restriction} of $f$ to $Z$, and denoted by
\bea
\bt[column sep=small] f|_Z:Z\subset X\ar[r] & Y,~~z\mapsto f(z).\et\nn
\eea
A map $g:X'\ra Y'$ is \ul{equal} to $f:X\ra Y$, written $g=f$ (otherwise, $g$ and $f$ are \ul{distinct} or \ul{unequal}, written $g\neq f$), if (i) $X'=X$ and (ii) $g(x)=f(x)$ for all $x\in X=X'$. A map $g:X'\ra Y'$ is an \ul{extension} of $f:X\ra Y$ if $f=g|_X$, i.e., (i) $X\subset X'$ and (ii) $g(x)=f(x)$ for all $x\in X\subset X'$.

In the map $f:X\ra Y$, the class $X$ is called the \ul{domain} of $f$, denoted by $\dom(f)$. The target class $Y$ is called the \ul{codomain} of $f$, denoted by $\cod(f)$. The \ul{image (or range) of $f$} is the subclass
\bea
\im f:=\{f(x):x\in A\}\subset Y~~~~(\txt{also ambiguously written as $f(X)$ })\nn
\eea
of all objects of $Y$ that are each assigned to an object of $X$. If $A\subset X$, the \ul{image of $A$} under $f$ is
\bea
\im(f|_A):=\{f(x):x\in A\}\subset Y~~~~(\txt{also ambiguously written as $f(A)$ }).~~~~(\blue{\txt{footnote}}\footnotemark).\nn
\eea
\footnotetext{\magenta{Caution}: If it happens that both $A\subset X$ and $A\in X$ are true statements, then it is clear that $f(A)$ does not have a unique meaning. However, with a bit of extra care, the actual meaning of $f(A)$ will often be clear from the context.}If $B\subset Y$, the \ul{preimage of $B$} under $f$ is the subclass
\bea
\textstyle f^{-1}(B):=\{x\in X:f(x)\in B\}=\bigcup\{A\subset X:f(A)\subset B\}\subset X~~~~(\blue{\txt{footnote}}\footnotemark)\nn
\eea
\footnotetext{Given any maps $f:X\ra Y$, $g:Y\ra Z$ and any subclasses $A\subset X$, $B\subset Y$, $\{B_i\subset Y\}_{i\in I}$, $C\subset Z$, we always have
\bit[leftmargin=0cm]
\item[](1) $(fg)^{-1}(C)=g^{-1}(f^{-1}(C))$, ~~(2) $f^{-1}\big(f(A)\big)\cap A=A$, ~~(3) $f^{-1}(B)=f^{-1}\big(B\cap f(X)\big)$, ~~(4) $f(f^{-1}(B))=B\cap f(X)$,
\item[](5) $f^{-1}(\bigcup_i B_i)=\bigcup_i f^{-1}(B_i)$, ~~(6) $f^{-1}(\bigcap_i B_i)=\bigcap_i f^{-1}(B_i)$, ~~(7) $f(\bigcup_i B_i)=\bigcup_i f(B_i)$, ~~(8) $f(\bigcap_i B_i)\subset\bigcap_i f(B_i)$,
\item[] where (4), (5), and (7) imply that ~~$B\cap(\bigcup_i B_i)=\bigcup_i(B\cap B_i)$.
\eit
}of all objects of $X$ whose images lie in $B$. If $y\in Y$, the \ul{fiber} of $f$ at $y$ is {\small $f^{-1}(y):=f^{-1}(\{y\})$}. A map of the form $s:Y\ra X$, $y\mapsto s(y)\in f^{-1}(y)$, if it exists, is called a \ul{(global) section} of $f$. (\blue{footnote}\footnote{That is, a (global) section of $f:X\ra Y$ is a map $s:Y\ra X$ with values in the respective (or corresponding) fibers of $f$.}). A map of the form $s:B\subset Y\ra X$, $y\mapsto s(y)\in f^{-1}(y)$, i.e., a section of $f|_{f^{-1}(B)}:f^{-1}(B)\subset X\ra B\subset Y$, is called a \ul{local section} of $f$.

The map $f:X\ra Y$ is \ul{constant} if there exists an object $y_0\in Y$ such that $f(x)=y_0$ for all $x\in X$, i.e., the image $f(X)=\{y_0\}$ consists of a single object. The map $f$ is \ul{finite} if the image $f(X)$ is finite (i.e., a finite collection of objects). The map $f$ is \ul{identity}, written $f=id_X$, if (i) $X\subset Y$ and (ii) $f(x)=x$ for all $x\in X$. A subclass $A\subset X\cap Y$ is an \ul{invariant class} of $f:X\ra Y$ if $f(A)\subset A$ (i.e., $f|_A:A\mapsto A$). In particular, an object $x\in X\cap Y$ is called a \ul{fixed point} of $f:X\ra Y$ if $\{x\}\subset X\cap Y$ is an invariant class of $f$ (i.e., $f(x)=x$). A subclass $A\subset X\cap Y$ is a \ul{point-invariant class} of $f:X\ra Y$ (written $A\subset Fix(f):=\{x\in X\cap Y:f(x)=x\}$) if every object $x\in A$ is a fixed point of $f$.

The map $f:X\ra Y$ is \ul{injective} (or one-to-one) if it \emph{uniquely} assigns to each object $x\in X$ a \emph{single} object $f(x)\in Y$, (\blue{footnote}\footnote{\magenta{Caution}: The phrase ``A rule that \emph{uniquely} assigns to each $x\in X$ a \emph{single} $f(x)\in Y$ ... '' is possibly also written as ``A rule that assigns to each $x\in X$ a \emph{unique} $f(x)\in Y$ ... ''. The word ``unique'' in the second phrase can also be interpreted as ``single'', which is a point of confusion since strictly speaking it only implies the map is well-defined but not that it is injective. Therefore, one must learn to extract the actual meaning of the word ``unique'' from the context. For us, the phrase ``\magenta{For each $x\in X$ there exists a \emph{unique} $y=y_x\in Y$ ...} '' is the same as ``\magenta{For each $x\in X$ there exists a \emph{single} $y=y_x\in Y$ ...} '' and will precisely mean that ``\magenta{We have a \emph{well-defined} map $X\ra Y,~x\mapsto y_x$ ...} ''.}) in the sense that each object $y\in Y$ is assigned to at most a \emph{single} object $x\in X$, i.e., if $f(x)=f(x')$ then $x=x'$, or, if $x\neq x'$ then $f(x)\neq f(x')$. Briefly, an injective map is specified as
\bea
\bt[column sep=small] f:X\ar[r,hook] & Y,~~x\mapsto f(x).\et\nn
\eea
The map $f:X\ra Y$ is \ul{surjective} (or onto) if $f(X)=Y$. Briefly, a surjective map is specified  as
\bea
\bt[column sep=small] f:X\ar[r,two heads] & Y,~~x\mapsto f(x).\et\nn
\eea
The map $f:X\ra Y$ is \ul{bijective} if it is both injective and surjective. Briefly, a bijective map is specified as
\bea
\bt[column sep=small] f:X\ar[r,hook,two heads] & Y,~~x\mapsto f(x).\et\nn
\eea
If the map $f:X\ra Y$ is bijective, then the map $f^{-1}:Y\ra X$, $f(x)\mapsto x$ is called the \ul{inverse map} of $f$.

A bijective map of the form $f:X\ra X$ (i.e., a bijective map from a class to itself) is called a \ul{permutation} (or \ul{symmetry}) of $X$.
\end{dfn}

It is often convenient to view a map $f$ as a diagram
{\footnotesize $\bt[column sep=small]\dom(f)\ar[r,"f"]&\cod(f).\et$} An injective (resp. surjective) map is also called an \ul{injection} (resp. \ul{surjection}), while a bijective map is also called a \ul{bijection}.

\begin{rmk}[\textcolor{blue}{\index{Well-defined map}{Well-defined map}}]
A map $f:X\ra Y$, $x\mapsto f(x)$ must be \ul{well-defined} in the sense that for each $x\in X$, $f(x)$ is a single object (as the definition requires). Equivalently, no two objects of $Y$ can be assigned to the same object of $X$, i.e., for all $x,x'\in X$, if $x=x'$ then $f(x)=f(x')$, or, if $f(x)\neq f(x')$ then $x\neq x'$.
\end{rmk}

\begin{dfn}[\textcolor{blue}{\index{Composition! of maps}{Composition of maps}}]
Given maps {\footnotesize $X\sr{f}{\ral}Y$} and {\footnotesize $Y\sr{g}{\ral}Z$}, their composition is the map
{\small\bea
\left(X\sr{g\circ f}{\ral}Z\right)=\left(X\sr{f}{\ral}Y\sr{g}{\ral}Z\right),~x\mapsto (g\circ f)(x):=g\big(f(x)\big).\nn
\eea}
\end{dfn}

Whenever no confusion is possible, the composition $g\circ f$ will often be written simply as $gf$.

\begin{rmk}[\textcolor{blue}{Nature and existence of sections}]
Consider a map $f:X\ra Y$. Observe that a map $s:Y\ra X$ is a section of $f$ iff $f\circ s=id_Y$ (which shows that if a section of $f$ exists, then it is injective and also makes $f$ surjective). Moreover, a section $s:Y\ra X$ of $f$ exists iff (i) $f$ is surjective (i.e., all fibers of $f$ are nonempty) and (ii) given any $y\in Y$, we can make a unique selection or choice $s(y)\in f^{-1}(y)$.
\end{rmk}

\begin{thm}[\textcolor{blue}{Characterization of a bijection}]
Let $f:X\ra Y$ and $g:Y\ra X$ be maps of classes. Then $f$ is a bijection with inverse map $g$ (equivalently, $g$ is a bijection with inverse map $f$) $\iff$ $g\circ f=id_X:X\sr{f}{\ral}Y\sr{g}{\ral}X$ and $f\circ g=id_Y:Y\sr{g}{\ral}X\sr{f}{\ral}Y$.
\end{thm}
\begin{proof}
($\Ra$): Assume $f$ (resp. $g$) is a bijection with inverse map $g$ (resp. $f$). Then it is clear by direct substitution and evaluation that $g\circ f=id_X$ and $f\circ g=id_Y$. ($\La$): Conversely, assume $g\circ f=id_X$ and $f\circ g=id_Y$. Then $g\circ f=id_X$ implies $g$ is surjective and $f$ is injective, while $f\circ g=id_Y$ implies $f$ is surjective and $g$ is injective, i.e., $f$ and $g$ are bijections. It then follows that $g$ is the inverse map of $f$ and $f$ is the inverse map of $g$.
\end{proof}

\begin{dfn}[\textcolor{blue}{
\index{Indexing (Indexed or Cartesian collection)}{Indexing (Indexed collection or Cartesian collection)} of classes,
\index{Cartesian product}{Cartesian product},
\index{Cartesian power}{Cartesian power},
\index{Ordered collection}{Ordered collection},
\index{Cartesian pair (ordered pair)}{Cartesian pair (Ordered pair)}}]
Let $I$ be a class, and $\C$ a collection of classes. An \ul{$I$-indexing} (\ul{$I$-indexed collection} or \ul{$I$-cartesian collection}) of classes in $\C$ is a map {\footnotesize$C_I=(C_i)_{i\in I}:I\ra\C$, $i\mapsto C_i$}. (\blue{footnote}\footnote{The indexed collection $C_I=(C_i)_{i\in I}$, being a map $I\ra\C$, is not the same thing as its image $C_I(I)=\{C_i\}_{i\in I}$ which, unlike $C_I$, is an unindexed collection. However, the indexed collection $(C_i)_{i\in I}$ is sometimes (without any warning) also imprecisely denoted by the unindexed (hence in general smaller) collection $\{C_i\}_{i\in I}$, in which case the brackets $\{,\}$ must accordingly be understood as representing the braces $(,)$.}). Let {\footnotesize$\bigcup C_I:=\bigcup_{i\in I} C_i:=\{x~|~x\in C_i~\txt{for some}~i\in I\}$}. The \ul{cartesian product} of $C_I$ is the class of $I$-indexings (i.e., $I$-indexed or $I$-cartesian collections) of objects of $C_I$ given by
\bea
\textstyle\prod C_I:=\prod_{i\in I}C_i:=\big\{\txt{maps}~f:I\ra\bigcup C,~i\mapsto f(i)\in C_i\big\}=\big\{(c_i)_{i\in I}~|~c_i\in C_i~\txt{for each}~i\in I\big\}.\nn
\eea
If there is a single class $B\in\C$ such that $C_i=B$ for all $i\in I$ (hence $\bigcup_{i\in I}C_i=B$), then we write
\bea
\textstyle B^I:=\prod_{i\in I}B:=\prod_{i\in I}C_i=\big\{\txt{maps}~f:I\ra B\big\}.~~~~~~~~~~~~\txt{(The \ul{$I$-cartesian power} of $B$)}\nn
\eea
If $I=\{0,1,\cdots,n-1\}$, then the $I$-cartesian collection  $C_I=(C_0,C_1,...,C_{n-1})$ is called an \ul{$n$-tuple}, and its cartesian product is the class of $n$-tuples given by
\bea
\prod_{i=0}^{i=n-1}C_i=C_0\times C_1\times\cdots\times C_{n-1}:=\{(c_0,c_1,...,c_{n-1})~|~c_i\in C_i\},\nn
\eea
and we write the cartesian power $B^I=B^{\{0,1,\cdots,n-1\}}$ as $B^n=\{(c_0,c_1,...,c_{n-1})~|~c_i\in B\}$.

If $I$ is an ordered class (Definition \ref{OrderDef}), then an $I$-cartesian collection is also called an \ul{ordered collection}. In particular (with $2:=\{0,1\}$ viewed as an ordered class) a $2$-tuple (i.e., a $2$-cartesian collection, also called a \ul{cartesian pair}) is called an \ul{ordered pair}.
\end{dfn}

\begin{dfn}[\blue{\index{Disjoint! union}{Disjoint union} of an indexing of classes}]\label{DisjUnDef}
Let $\C$ be a collection of classes, $I$ a class, and $C_I=(C_i)_{i\in I}:I\ra\C,~i\mapsto C_i$ an $I$-indexing of classes. The disjoint union of $C_I$ is the union class
\begin{align}
\textstyle \bigsqcup C_I:=\bigsqcup_{i\in I}C_i:=\bigcup_{i\in I}C_i\times\{i\}=\bigcup_{i\in I}\big\{(c_i,i):c_i\in C_i\big\}.~~~(\txt{\blue{footnote}}\footnotemark)\nn
\end{align}
\footnotetext{With $\C_i:=\{C\times\{i\}:C\in\C\}$ for each $i\in I$, we have $\bigsqcup C_I=\bigcup D_I$, where ~$D_I=(D_i)_{i\in I}:I\ra\bigcup_{i\in I}\C_i,~i\mapsto D_i:=C_i\times\{i\}=\{(c_i,i):c_i\in C_i\}\in\C_i$.}
\end{dfn}

\begin{dfn}[\textcolor{blue}{
\index{Partition of a class}{Partition of a class},
\index{Cover of a class}{Cover of a class}}]
Let $C$ be a class. A collection $\C\subset\P(C)$ of subclasses of $\C$ is a partition of $C$ if the subclasses in $\C$ are
\bit
\item[(i)] nonempty (i.e., if $A\in \C$ then $A\neq\emptyset$),
\item[(ii)] disjoint (i.e., for any $A,B\in\C$, if $A\neq B$ then $A\cap B=\emptyset$), and
\item[(iii)] cover $C$ (i.e., $C=\bigcup\C:=\{\txt{all objects from members of $\C$}\}$).
\eit
Any collection $\C\subset\P(C)$ of subclasses of $C$ that satisfies (iii) above is a \ul{cover} (or \ul{covering}) of $C$.
 \end{dfn}

\begin{dfn}[\blue{\index{Graph of a map}{Graph of a map}, \index{Submap}{Submap}}]
Let $f:X\ra Y$ and $f':X'\ra Y'$ be maps of classes. The \ul{graph} of $f$ is the subclass $G(f):=\{(x,f(x)):x\in X\}\subset X\times Y$. $f$ is a \ul{submap} of $f'$, written $f\subset f'$, if $G(f)\subset G(f')$, i.e., $\dom f\subset\dom f'$ and $f'|_{\dom f}=f$.
\end{dfn}

\begin{dfn}[\textcolor{blue}{
\index{Relation! between classes}{{Relation} between classes},
\index{Domain! of a relation}{Domain of a relation},
\index{Image of! a relation}{Image of a relation},
\index{Opposite! relation}{Opposite relation},
\index{Relation! on a class}{Relation on a class},
\index{Subrelation}{Subrelation},
\index{Reflexive relation}{Reflexive relation},
\index{Symmetric! relation}{Symmetric relation},
\index{Antisymmetric relation}{Antisymmetric relation},
\index{Linear! (total) relation}{Linear (or total) relation},
\index{Transitive! relation}{Transitive relation},
\index{Strict! relation}{Strict relation},
\index{Strictly! antisymmetric relation}{Strictly antisymmetric relation},
\index{Equivalence! relation}{Equivalence relation},
\index{Equivalence! class}{Equivalence class},
\index{Quotient! class}{Quotient class},
\index{Partial! order}{Partial order},
\index{Strict! partial order}{Strict partial order},
\index{Partially ordered class}{Partially ordered class},
\index{Minimal! object}{Minimal object},
\index{Maximal! object}{Maximal object},
\index{Lower bound}{Lower bound},
\index{Upper bound}{Upper bound},
\index{Strict! lower bound}{Strict lower bound},
\index{Strict! upper bound}{Strict upper bound},
\index{Least! object}{Least object},
\index{Greatest! object}{Greatest object},
\index{Infimum (glb)}{Infimum (glb)},
\index{Supremum (lub)}{Supremum (lub)},
\index{Infimum (glb) property}{Infimum (glb) property},
\index{Supremum (lub) property}{Supremum (lub) property},
\index{Linearly ordered class (chain, list, line, linear class)}{Linearly-ordered class (chain, list, line, linear class)},
\index{Well-ordered! class}{Well-ordered class}}]\label{OrderDef}
Let $A,B$ be classes. A \ul{(binary) relation between} $A$ and $B$ is a subclass $R\subset A\times B$. The \ul{domain}, \ul{image}, and \ul{opposite relation} of a relation $R\subset A\times B$ are respectively defined as
\bea
&&\dom R:=\{a\in A~|~(a,b)\in R~\txt{for some}~b\in B\},~~~~\im R:=\{b\in B~|~(a,b)\in R~\txt{for some}~a\in A\},\nn\\
&&R^{op}:=\{(b,a)\in B\times A:(a,b)\in R\}\subset B\times A,~~~~\txt{i.e., ~$(a,b)\in R \iff (b,a)\in R^{op}$}.\nn
\eea
If $B=A$, a relation between $A$ and itself (or a relation within $A$) $R\subset A\times A$ is called a \ul{relation on} $A$. (\blue{footnote}\footnote{Relations of the form $R\subset A\times A$ suffice because $A\times B\subset(A\cup B)\times(A\cup B)$, and so a relation between $A$ and $B$ is a relation on $A\cup B$.}). Given a relation $R\subset A\times B$ between $A$ and $B$ (resp. $R\subset A\times A$ on $A$), every pair of subsets $A'\subset A$, $B'\subset B$ (resp. every subset $A'\subset A$) inherits the \ul{subrelation} $R':=R\cap(A'\times B')\subset A'\times B'$ between $A'$ and $B'$ (resp. $R':=R\cap(A'\times A')\subset A'\times A'$ on $A'$).

For notational simplicity, we have the following abbreviations: If $R\subset A\times A$,
\bit
\item $(a,b)\in R$ is written as ~$aRb$ (or as ~$bR^{op}a$),
\item $(a,b),(b,c)\in R$ (i.e., $aRb$ and $bRc$) is written as ~$aRbRc$.
\eit
We will need the following properties (also inherited by any subrelation). A relation $R\subset A\times A$ is
\bit
\item \ul{reflexive} ~if ~$aRa$ ~~(for all ~$a\in A$).
\item \ul{symmetric} ~if ~$aRb$ ~$\Ra$~ $bRa$ ~~(for all ~$a,b\in A$).
\item \ul{antisymmetric} ~if ~$aRbRa$ ~$\Ra$~ $a=b$ ~~(for all ~$a,b\in A$).
\item \ul{transitive} ~if ~$aRbRc$ ~$\Ra$~ $aRc$ ~~(for all ~$a,b,c\in A$).
\item \ul{strict} ~if ~$aRb$ ~$\Ra$~ $a\neq b$ ~~(for all ~$a,b\in A$).
\item \ul{strictly antisymmetric} ~if $aRbRa$ is impossible ~~(for all ~$a,b\in A$).
\item \ul{linear (or total)} ~if ~$aRb$ ~or ~$bRa$ ~~(for all ~$a,b\in A$).
\eit
A relation $\sim~\subset A\times A$ is an \ul{equivalence relation} on $A$ if it is reflexive, symmetric, and transitive, i.e.,
\bea
(i)~a\sim a,~~~~(ii)~a\sim b~~\Ra~~b\sim a,~~~~(iii)~a\sim b\sim c~~\Ra~~a\sim c.\nn
\eea
Given an equivalence relation $\sim$ on $A$ (i.e., an equivalence relation $\sim~\subset A\times A$), the \ul{equivalence class} of $a\in A$ is the subclass ~$[a]~:=~\{b\in A:~(a,b)\in R\}=\{b\in A:~b\sim a\}\subset A$. The \ul{(equivalence-) quotient class} associated with $\sim$ is the collection of $\sim$-equivalence classes $A_\sim:={A\times A\over\sim}:=\{[a]:a\in A\}$.

A relation $\leq~\subset A\times A$ (with opposite $\geq~\subset A\times A$) is a \ul{partial order}, making $(A,\leq)$ a \ul{partially ordered class}, if it is reflexive, antisymmetric, and transitive, i.e., for all $a,b\in A$,
\bea
(i)~a\leq a,~~~~(ii)~a\leq b\leq a~~\Ra~~a=b,~~~~(iii)~a\leq b\leq c~~\Ra~~a\leq c.\nn
\eea
A relation $<~\subset A\times A$ (with opposite $>~\subset A\times A$) is a \ul{strict partial order}, making $(A,<)$ a \ul{strictly partially ordered class}, if it is strict and transitive, i.e., (i) $a<b$ $\Ra$ $a\neq b$, and (ii) $a<b<c$ $\Ra$ $a<c$. Specifying a partially ordered class $(A,\leq)$ is equivalent to specifying an associated strictly partially ordered class $(A,<)$ and vice versa (see the \blue{footnote}\footnote{\ul{A class is partially ordered iff strictly partially ordered}: Associated with any partial order $\leq~\subset A\times A$ is the unique strict partial order given by ``$a<b$ if $a\leq b$ and $a\neq b$''. Similarly, associated with any strict partial order $<~\subset A\times A$ is the unique partial order $\leq~\subset A\times B$ given by ``$a\leq b$ if $a<b$ or $a=b$''.}).

Let $(A,\leq)$ be a partially ordered class, $B\subset A$, and $l,u\in A$. Then $l$ is a \ul{lower bound} of $B$, written $l\leq B$ or $B\geq l$, if $l\leq b$ for all $b\in B$ (resp. $l$ is a \ul{strict lower bound} of $B$, written $l<B$ or $B>l$, if $l<b$ for all $b<B$). Similarly, $u$ is an \ul{upper bound} of $B$, written $B\leq u$ (or $u\geq B$), if $b\leq u$ for all $b\in B$ (resp. $u$ is a \ul{strict upper bound} of $B$, written $B<u$ or $u>B$, if $b<u$ for all $b\in B$).

Let $(A,\leq)$ be a partially ordered class and $a,a'\in A$. Then $a$ and $a'$ are \ul{comparable objects} if $a\leq a'$ or $a'\leq a$. An object $m\in A$ is \ul{minimal} if for all $a\in A$, $a\leq m$ $\Ra$ $a=m$ (i.e., $m$ has no strict lower bound). An object $m\in A$ is \ul{maximal} if for all $a\in A$, $m\leq a$ $\Ra$ $m=a$ (i.e., $m$ has no strict upper bound).

Let $(A,\leq)$ be a partially ordered class and $B\subset A$. A \ul{least object} of $B$ is an object $l\in B$ that is a lower bound $l\leq B$ (written $l\in\min_A B$, or $l=\min_AB$ if unique). A \ul{greatest object} of $B$ is an object $g\in B$ that is an upper bound $g\geq B$ (written $g\in\max_AB$, or $g=\max_AB$ if unique). (If it exists) an \ul{infimum} of $B$ is a greatest lower bound (glb) of $B$, i.e., $\inf B\in\max_A\{l\in A:l\leq B\}$ or $\inf B:=\max_A\{l\in A:l\leq B\}$ if unique. (If it exists) a \ul{supremum} of $B$ is a least upper bound (lub) of $B$, i.e., $\sup B\in\min_A\{u\in A:B\leq u\}$ or $\sup B:=\min_A\{u\in A:B\leq u\}$ if unique.

A partially ordered class $(A,\leq)$ has the \ul{infimum (glb) property} (resp. \ul{supremum (lub) property}) if every nonempty subclass of $A$ with a lower bound (resp. upper bound) has an infimum (resp. supremum).

A partially ordered class $(A,\leq)$ is \ul{linearly ordered} (or \ul{totally ordered}, or a \ul{chain}, or a \ul{list}, or a \ul{line}, or a \ul{linear class}) if the relation $\leq$ is a linear (or total) relation. A linearly-ordered class $(A,\leq)$ is \ul{well-ordered} if every nonempty subclass $\emptyset\neq B\subset A$ contains a least object.
\end{dfn}

\begin{notation}
Let $(A,\leq)$ be a partially ordered class. If $B,C\subset A$, we will write $B\leq C$ if $b\leq c$ for all $b\in B,c\in C$ (resp. $B<C$ if $b<c$ for all $b\in B,c\in C$).
\end{notation}

\begin{dfn}[\blue{
\index{Initial segment}{Initial segment},
\index{Strict! initial segment}{Strict initial segment},
\index{Lower section}{Lower section},
\index{Strict! lower section}{Strict lower section},
\index{Upper section}{Upper section},
\index{Strict! upper section}{Strict upper section},
\index{Lower segment}{Lower segment},
\index{Strict! lower segment}{Strict lower segment},
\index{Upper segment}{Upper segment},
\index{Strict! upper segment}{Strict upper segment}
}]
Let $(A,\leq_A)$ and $(B,\leq_B)$ be partially ordered classes. Then $A$ is an \ul{initial segment} of $B$, written $A\preceq B$, if (i) $A\subset B$, (ii) $\leq_A~=~\leq_B|_A~:=~\leq_B\cap(A\times A)$, and (iii) $A<B\backslash A$. An initial segment $A\preceq B$ is a \ul{strict initial segment}, written $A\prec B$, if $A\subsetneq B$.

Given $A'\subset A$, let $L_A(A'):=\{l\in A:l\leq A'\}$ be the class of lower bounds (resp. $L_A^\ast(A'):=\{l\in A:l<A'\}$ the class of strict lower bounds) of $A'$ and $U_A(A'):=\{u\in A:A'\leq u\}$ the class of upper bounds (resp. $U_A^\ast(A'):=\{u\in A:A'<u\}$ the class of strict upper bounds) of $A'$. The \ul{lower section} at $a\in A$ is $L_A(a):=L_A(\{a\})$ (resp. the \ul{strict lower section} at $a\in A$ is $L^\ast_A(a):=L^\ast_A(\{a\})$ ). The \ul{upper section} at $a\in A$ is $U_A(a):=U_A(\{a\})$ (resp. the \ul{strict upper section} at $a\in A$ is $U^\ast_A(a):=U^\ast_A(\{a\})$ ).

If $(A,\leq)$ is a linearly-ordered class, then a lower section (resp. strict lower section) of $A$ will be called a \ul{lower segment} (resp. \ul{strict lower segment}), and similarly, an upper section (resp. strict upper section) of $A$ will be called an \ul{upper segment} (resp. \ul{strict upper segment}).
\end{dfn}

\begin{dfn}[\blue{
\index{Directed! class}{Directed class},
\index{Lower-directed class}{Lower-directed class},
\index{Filter in a class}{Filter in a class}}]
Let $(A,\leq)$ be a partially ordered class. The class $A$ is \ul{(upper-) directed} if any two objects $a,a'\in A$ have a common upper bound $u\geq\{a,a'\}$ (equivalently, every finite subclass $\{a_i\}_{i=1}^n\subset A$ has an upper bound $u\geq\{a_i\}_{i=1}^n$). The class $A$ is \ul{lower-directed} if any two objects $a,a'\in A$ have a common lower bound $l\leq\{a,a'\}$ (equivalently, every finite subclass $\{a_i\}_{i=1}^n\subset A$ has a lower bound $l\leq\{a_i\}_{i=1}^n$).

A subclass $F\subset A$ is a \ul{filter} if (i) $(F,\leq)$ is lower-directed (i.e., any two objects $f,f'\in F$ have a common lower bound $l\leq\{f,f'\}$) and (ii) for any $f\in F$, $U_A(f)\subset F$.
\end{dfn}

\begin{dfn}[\blue{Product of partially ordered classes}]
Let $I$ be a class, $\C$ a class of partially ordered classes, and $A_I:=(A_i)_{i\in I}=\big((A_i,\leq_{A_i})\big)_{i\in I}:I\ra\C$ an $I$-indexing in $\C$. The product of $A_I$ is the partially ordered class
\[
\textstyle\prod A_I=\prod_{i\in I}A_i=\prod_{i\in I}(A_i,\leq_{A_i})=\left(\prod A_I,\leq\right)=\left(\prod_{i\in I} A_i,\leq\right),\nn
\]where for all $(a_i)_{i\in I},(a'_i)_{i\in I}\in\prod A_I$, we have ``$(a_i)_{i\in I}\leq(a'_i)_{i\in I}$ if $a_i\leq_{A_i}a'_i$ for all $i\in I$''. (\blue{footnote}\footnote{Given a directed class $D$ and $n\geq 1$, its finite cartesian power $D^n$ is also a directed class (by induction on $n$), since for any {\scriptsize $(a,a'),(b,b')\in D\times D$}, there exist $c,c'\in D$ such that $a,b\leq c$ and $a',b'\leq c'$ in $D$, and so $(a,a'),(b,b')\leq(c,c')$ in $D\times D$.})
\end{dfn}

\begin{prp}[\textcolor{blue}{Outgoing maps, Partitions, and Equivalence relations are the same thing}]
Let $A$ be a class. (i) A map of classes $f:A\ra B$ is given $\iff$ (ii) a partition of $A$ is given, $\iff$ (iii) an equivalence relation $\sim$ on $A$ is given.
\end{prp}
\begin{proof}
{\flushleft(i)$\Ra$(ii):} Given a map $f:A\ra B$, the collection $\left\{f^{-1}(b):b\in B\right\}$ contains a partition of $A$.

{\flushleft(ii)$\Ra$(iii):} Given a partition $\C\subset\P(A)$ of $A$, define a relation $\sim$ on $A$ by $a\sim b$ if $a$ and $b$ belong to the same member of the partition $\C$. Then $\sim$ is an equivalence relation on $A$.

{\flushleft(iii)$\Ra$(i):} Given an equivalence relation $\sim$ on $A$, let $[a]:=\{a'\in A:a'\sim a\}$ for $a\in A$. Then we get a map
\[
f:A\ra B,~a\mapsto [a],~~~~B:=\big\{[a]:a\in A\big\}. \qedhere
\]
\end{proof}

\begin{thm}[\textcolor{blue}{\index{Cantor-Schr\"oder-Bernstein theorem}{Cantor-Schr\"oder-Bernstein theorem}: \cite[page 147]{enderton1977}}]\label{CSBthm}
Let $X,Y$ be classes, and $f:X\ra Y$, $g:Y\ra X$ injective maps. Then there exists a bijective map $h:X\ra Y$.
\end{thm}
\begin{proof}
If $f$ or $g$ is surjective the result holds trivially. So, assume $f,g$ are not surjective. Form the following sequences of classes:
\bea
X_0:=X-g(Y),~~Y_0:=f(X_0),~~X_n:=g(Y_{n-1}):=g\big(f(X_{n-1})\big),~~Y_n:=f(X_n)=fg(Y_{n-1}),~~\txt{for}~~n\geq1.\nn
\eea
Observe that for all $n\geq 1$, $X_0\cap X_n=(X-g(Y))\cap g(Y_{n-1})=\emptyset$, and so the injectivity of $f$ implies $Y_0\cap Y_n=f(X_0)\cap f(X_n)=f(X_0\cap X_n)=\emptyset$. Also, for all $i,j\geq 0$, the injectivity of $f,g$ implies
\bea
&&X_i\cap X_j=g(Y_{i-1})\cap g(Y_{j-1})=g(Y_{i-1}\cap Y_{j-1})=g(f(X_{i-1})\cap f(X_{j-1}))=gf(X_{i-1}\cap X_{j-1}),\nn\\
&&Y_i\cap Y_j=f(X_i)\cap f(X_j)=f(X_i\cap X_j)=fg(Y_{i-1}\cap Y_{j-1}),\nn
\eea
and so by induction and symmetry, $X_i\cap X_j=\emptyset$, $Y_i\cap Y_j=\emptyset$ whenever $i\neq j$.

Now define a map $h:X\ra Y$ as follows:
\bea
\textstyle h(x):=\left\{
        \begin{array}{ll}
          f(x), & \txt{if}~~x\in \bigcup X_i:=X_0\cup\bigcup X_{i+1}=X_0\cup\bigcup g(Y_i)\\
          g^{-1}(x), & \txt{if}~~x\in X-\bigcup X_i=g(Y)-\bigcup g(Y_i)=g(Y-\bigcup Y_i)
        \end{array}
      \right.\nn
\eea
To check \ul{injectivity} of $h$, let $x,x'\in X$ be distinct (i.e., $x\neq x'$). If $x,x'\in\bigcup X_i$ or $x,x'\in X-\bigcup X_i$, then $h(x)\neq h(x')$ by injectivity of $f$ and $g$. So, assume wlog that $x\in\bigcup X_i$ and $x'\in X-\bigcup X_i$. Then
\bea
&&\textstyle h(x)=f(x)\in f\left(\bigcup X_i\right)=\bigcup f(X_i)=\bigcup Y_i,~~~~h(x')=g^{-1}(x')\in g^{-1}\left(X-\bigcup X_i\right)=Y-\bigcup Y_i,\nn\\
&&\textstyle~~\Ra~~h(x)\neq h(x').\nn
\eea
To check \ul{surjectivity} of $h$, let $y\in Y$. If $y\in \bigcup Y_i=\bigcup f(X_i)=f\left(\bigcup X_i\right)=h\left(\bigcup X_i\right)$, it is clear that $y\in h(X)$. So, assume $y\in Y-\bigcup Y_i$. Then
\[
\textstyle g(y)\in g(Y-\bigcup Y_i)=X-\bigcup X_i~~\Ra~~h(g(y))=g^{-1}\big(g(y)\big)=y,~~\Ra~~y\in h(X). \qedhere
\]

\end{proof}
\begin{crl}[\textcolor{blue}{Choice based equivalent of the Cantor-Schr\"oder-Bernstein theorem}]\label{CSBcrl}
Let $X,Y$ be classes, and $f:X\ra Y$, $g:Y\ra X$ surjective maps. Then there exists a bijective map $h:X\ra Y$.
\end{crl}
\begin{proof}
If a surjective map $u:A\ra B$ exists, we get an injective map $u':B\ra A$, $b\mapsto u'(b)\in u^{-1}\big(b\big)$, assuming we can always make a unique selection or \ul{choice} $u'(b)\in u^{-1}\big(b\big)$. Also, if an injective map $v:A\ra B$ exists, we get a surjective map $v':B\ra A$ given by $v'|_{f(A)}:=v^{-1}$ and $v'|_{B\backslash f(A)}:=\txt{const.}$

Hence, a surjective (resp. injective) map $X\ra Y$ exists $\iff$ an injective (resp. surjective) map $X\la Y$ exists.
\end{proof}

\begin{dfn}[\textcolor{blue}{
\index{Cardinal! order}{Cardinal order},
\index{Cardinal! equivalence}{Cardinal equivalence},
\index{Cardinality}{{Cardinality}},
\index{Finite! class}{Finite class},
\index{Infinite class}{Infinite class},
\index{Countable! class}{Countable class},
\index{At most countable class}{At most countable class},
\index{Uncountable class}{Uncountable class}}]
\ul{Cardinal order} is the relation $\lesssim$ on classes defined as follows: Let $C, D$ be classes. We say $C$ is \ul{cardinally less} than $ D$, written $|C|\leq| D|$ or $C\lesssim D$, if there exists an injective map $C\hookrightarrow D$ (\blue{footnote}\footnote{Or equivalently, if a surjection $C\twoheadleftarrow D$ exists (when choice is allowed).}).

\ul{Cardinal equivalence} is the equivalence relation $\approx$ on classes given by ``$C\approx D$ if there exists a bijection $g:C\ra D$''. The \ul{cardinality} $|C|$ of $C$ is a measure of the number of objects of $C$, and represented by the equivalence class ~$|C|:=[C]_\approx :=\big\{ D~|~ D\approx C\big\}$.

The class $C$ is \ul{finite} if $C\approx\{1,2,\cdots,n\}$ for some natural number $n\in\Natural$ (in which case we by convention simply write $|C|=n$), otherwise $C$ is \ul{infinite}. Also, we will call the class $C$ \ul{countable} if $|C|=|\Natural|$, or in some cases, \ul{(at most) countable} if $|C|\leq|\Natural|$, otherwise $C$ is \ul{uncountable}.
\end{dfn}

\begin{lmm}
Cardinal order is a partial order. (\blue{footnote}\footnote{Later, an axiom (called axiom of choice) will make cardinal order a linear order on classes of sets.})
\end{lmm}
\begin{proof}
Let $C,D$ be classes. (i) Cardinal order $\lesssim$ is trivially reflexive. (ii) By Theorem \ref{CSBthm}, we have $|C|=|D|$ $\iff$ $|C|\leq|D|$ and $|D|\leq|C|$ (i.e., $C\approx D$ $\iff$ $C\lesssim D$ and $D\lesssim C$), and so $\lesssim$ is antisymmetric. (iii) Also, $\lesssim$ is transitive since the composition of injective maps is injective.
\end{proof}

In the discussion under the class of sets to follow and beyond, many of those concepts that do not rely on any special property of sets can be applied to classes in general (i.e., not only to classes of sets).

\section{The Class of Sets}\label{ClassOfSets}
The following definition presents a list of axioms that satisfactorily defines the class of sets. However, as usual, even in the proofs of those major results requiring all of these axioms, some intermediate results might not depend on some of the axioms. Also, while we will briefly refer to the full list of axioms as ZF-axioms (instead of ZFC-axioms), it is actually known as ZFC-axioms, where $C$ stands for ``\emph{including the axiom of choice}'', while ZF-axioms usually stands for the axioms without the axiom of choice.

There exist several reformulations (motivated by paradoxes, cardinality questions, and more) that shorten the list of fundamental axioms for sets. However, the ZF-axioms are  sufficient for our purpose.

\begin{dfn}[\textcolor{blue}{(\cite{enderton1977,goldrei1996})
\index{The class ``Sets''}{{Sets}},
\index{Zermelo-Fraenkel (ZF) axioms}{Zermelo-Fraenkel (ZF) axioms},
\index{Elements of a set}{Elements of a set},
\index{Empty set}{Empty set},
\index{Equal! sets}{Equal sets},
\index{Subset}{Subset},
\index{Intersection of! sets}{Intersection},
\index{Set! difference (Relative complement)}{Set difference (Relative complement)}, \index{Complement}{Complement},
\index{Union! of sets}{Union},
\index{Powerset}{Powerset},
\index{Map! of sets}{Map of sets},
\index{Singleton set}{Singleton set},
\index{Well-ordered! set}{Well-ordered set},
\index{Least! element}{Least element},
\index{Inductive! set}{Inductive set},
\index{Natural! numbers $\Natural$}{Natural numbers $\Natural$},
\index{Countable! set}{Countable set}}]
Sets are a collection of objects that satisfy a list of defining rules called \ul{Zermelo-Fraenkel (ZF) axioms}, where objects that form a set are called \ul{elements} of the set. According to the ZF-axioms, the class of sets (denoted by \ul{Sets}) has the following defining properties: Let $x,y$ be any sets (not necessarily distinct);
\begin{enumerate}[leftmargin=0.9cm]
\item[(0)] \textbf{Content}: Every element of a set is itself a set. (\blue{footnote}\footnote{This redundant, but important, observation follows directly from ZF6 and ZF4.}).

\item (ZF1) \textbf{Extensionality (equality of sets)}: $x=y$ ~$\iff$~ $x\subset y$ and $y\subset x$. (\blue{footnote}\footnote{Here, we assume the meaning of inclusion $\subset$ (which was discussed earlier under classes) is already clear to the reader.})

\item (ZF2) \textbf{Empty set}: There exists a set $\emptyset$ (called \ul{empty set}) that has no elements. (\blue{footnote}\footnote{This axiom is redundant because it is deducible from ZF4 as a subset of the infinite set postulated in ZF7.}).

(Note: Axiom ZF1 implies the empty set $\emptyset$ is unique, since by definition $\emptyset\subset x$ for every set $x$.)

\item (ZF3) \textbf{Pairing}: There exists a set $z$ whose only elements are $x$ and $y$, i.e., such that $z=\{x,y\}$.

    (Note: For any set $x$, the set $\{x\}:=\{x,x\}$ is called a \ul{singleton} set.)

\item (ZF4) \textbf{Separation (subsets)}: Any subcollection, say via some specification rule $\phi$, of elements  $z:=\{a\in x~|~\phi(a)\}$ of a set $x$ is itself a set called a \ul{subset} of $x$ (written $z\subset x$).

[[ For example, we have \ul{intersection} $x\cap y:=\{z\in x~|~z\in y\}=\{z\in y~|~z\in x\}$ and \ul{difference} (or \ul{relative-complement}) $x-y=x\backslash y:=\{z\in x~|~z\not\in y\}$ of sets. If $y\subset x$, the difference $x-y=x\backslash y:=\{z\in x~|~z\not\in y\}$ is also denoted by $y^c$, and called the \ul{complement} of $y$ in $x$. ]]

\item (ZF5) \textbf{Powerset}: The collection $\P(x)=2^x:=\{y~|~y\subset x\}$ of all subsets of $x$ is itself a set (called the \ul{powerset} of $x$).

\item (ZF6) \textbf{Union}: The \ul{union} $\cup x=\bigcup_{z\in x}z:=\{t\in z~|~z\in x\}=\{t~|~t\in z~\txt{for some}~z\in x\}$ of all elements of $x$ (i.e., the collection whose members are the elements of elements of $x$) is itself a set.

(Note: The union of all sets contains all sets, and so is not a set.)

\item (ZF7) \textbf{Infinity}: There exists a set $s$ with infinitely many elements (called an \ul{inductive set});
\bea
&&s:=\{\emptyset\}\cup\big\{a\cup\{a\}~|~a\in s\big\}=\Big\{\emptyset,\big\{\emptyset\big\},\big\{\emptyset,\{\emptyset\}\big\},\cdots\Big\}\supset\{0,1,2,\cdots\}=\Natural,\nn\\
&&\big(\txt{where, in general, a set $s$ is called \ul{inductive} if (i) $\emptyset\in s$, and (ii) ~$a\in s$ ~$\Ra$~ $a\cup\{a\}\in s$}~\big)\nn
\eea
i.e., such that its cardinality $|s|\geq|\Natural|=\infty$, where we define the set of \ul{natural numbers} by
\bea
\textstyle\Natural:=\bigcap\{s~|~s~~\txt{an inductive set}\}~~~~~~~~\txt{(also itself clearly an inductive set)}.\nn
\eea

(Note that a set $x$ is called a \ul{countable set} if there exists a bijection $f:\Natural\ra x$.)

\item (ZF8) \textbf{Replacement (image of a map)}: If $c\subset\txt{Sets}$ is any subclass of sets, and $f:x\ra c$ any map, the image $f(x):=\{f(a)~|~a\in x\}\subset c$ is a set. (Equivalently, for any map of classes $f:c\ra c'$ and any subclass $c_1\subset c$ that forms a set, its image $f(c_1)\subset c'$ is a set, and thus gives a \ul{map of sets} $f|_{c_1}:c_1\ra f(c_1)$.)

\item (ZF9) \textbf{Regularity or foundation}: If $x\neq\emptyset$, there exists an element $z\in x$ such that $x\cap z=\emptyset$ (i.e., a nonempty set is disjoint from at least one of its elements (\blue{footnote}\footnote{In particular, $z\not\in z$.}) ). That is (from Section \ref{ParOrdZorn}), every nonempty set $x$ contains an element $z$ that is $\in$-minimal in the sense that $z$ is a \ul{minimal element} of the \ul{poset} $(x,<):=(x,\in)$, i.e., for every $a\in x$, $a\not\in z$.

    (In particular, a set $x$ is not a member of itself (\blue{footnote}\footnote{It follows that if $x\in y$, then $y\not\subset x$ (and so $x\in y$ $\Ra$ $x\neq y$).}): Indeed, $x$ is the only member of the pair $\{x\}=\{x,x\}$, and so $\{x\}\cap x=\emptyset$ by regularity, which implies $x\not\in x$.)

\item (ZFC) \textbf{Choice or well-ordering}: If the elements of $x$ are all nonempty (disjoint) sets, then their union $\cup x$ has at least one subset $z\subset\cup x$ containing exactly one member of each element $a\in x$ (i.e., $z\cap a$ is a single element set for each $a\in x$), or equivalently, there exists a map
\bea
\textstyle f:x\ra\cup x:=\bigcup_{a\in x}a,~~a\mapsto f(a)\in a~~\left(\txt{i.e.,}~~\prod_{a\in x}a\neq\emptyset\right).\nn
\eea
Equivalently (via \ul{Zorn's lemma}, i.e., Theorem \ref{ZornLemma}) for every set $s$, there exists a linear order $\leq~\subset s\times s:=\{(a,b):a,b\in s\}$ that \ul{well-orders} $s$ in the sense that every nonempty subset $\emptyset\neq s'\subset s$ has a \ul{least element} (i.e., an element $m\in s'$ such that $m\leq a$ for all $a\in s'$).
\end{enumerate}
\end{dfn}

The statements of most of the following results are either (i) immediately clear from definitions or (ii) implicitly/explicitly contain their proofs, which are therefore not given separately.

\begin{lmm}[\textcolor{blue}{The cartesian pair of sets is representable as a set}]
Let $x,y$ be sets. The \ul{ordered pair} $(x,y):\{0,1\}\ra\txt{Sets}$, $0\mapsto x$, $1\mapsto y$ is uniquely determined by the set $\big\{\{x\},\{x,y\}\big\}$, because $\big\{\{x\},\{x,y\}\big\}=\big\{\{x'\},\{x',y'\}\big\}$ $\iff$ $x=x'$ and $y=y'$. Thus, we can define the ordered pair of $x,y$ as (\blue{footnote}\footnote{In particular, for any set $x$, ~$(x,x)=\{\{x\}\}\neq\{x\}\neq x$,~ even though for a special element $0$ (called ``zero'') in a ring or module later we will sometimes briefly/imprecisely write $\{0\}$ simply as $0$ (which, in the case of a ring or module that contains $\Natural$, strictly speaking equals $1:=0\cup\{0\}=\emptyset\cup\{\emptyset\}=\{\emptyset\}=\{0\}$).})
\bea
(x,y):=\big\{\{x\},\{x,y\}\big\}.\nn
\eea
\end{lmm}

\begin{thm}[\textcolor{blue}{A set-valued map on a set is representable as a set}]
If $x$ is a set and $c\subset\txt{Sets}$ a subclass of sets, then any map $f:x\ra c$ can be uniquely represented by its \index{Graph of a map of sets}{\ul{graph}}
\bea
G(f):=\{(a,f(a))~|~a\in x\}~\subset~x\times f(x),\nn
\eea
which is a set, as the image of the map ~$F:x\ra Sets,~a\mapsto(a,f(a))$.
\end{thm}
\begin{proof}
For any map $g:x\ra c$, $G(f)=G(g)$ $\iff$ $f=g$. This is because for all $a,b\in x$, we have $(a,f(a))=(b,g(b))$ $\iff$ $a=b$ and $f(a)=g(b)=g(a)$.
\end{proof}

In the following immediate corollaries of the above result, we have also considered alternative set-representations for the indexed collections involved.

\begin{crl}[\textcolor{blue}{A set-indexing of sets is representable as a set}]
Let $x,y,x',y'$ be sets.
\begin{enumerate}[leftmargin=0.7cm]
\item[a.] An \ul{$n$-tuple} $(x_1,...,x_n):\{1,2,\cdots,n\}\ra Sets,~i\mapsto x_i$ of sets $x_i$ can be represented either (i) directly by its graph $\{(i,x_i)\}_{i=1}^n$ or (ii) by the set defined inductively as
{\footnotesize
\bea
(x_1,...,x_n):=\big((x_1,...,x_{n-1}),x_n\big)=\big(((x_1,...,x_{n-2}),x_{n-1}),x_n\big)=\cdots=\big(((...(x_1,x_2),...,x_{n-2}),x_{n-1}),x_n\big).\nn
\eea}

\item[b.] For any \ul{$\Natural$-indexing} of sets $(x_i)_{i\in\Natural}:\Natural\ra Sets,~i\mapsto x_i$, the collection $\{(x_1,...,x_i):i\in\Natural\}$ is a set, as the image of the map $f:\Natural\ra\txt{Sets},~i\mapsto(x_1,...,x_i)$, and so we can represent $(x_i)_{i\in\Natural}$ either (i) directly by its graph $\{(i,x_i)\}_{i\in\Natural}$ or (ii) by the image ~$(x_i)_{i\in\Natural}:=\im f=f(\Integer)=\{(x_1,...,x_i):i\in\Natural\}$~ (\blue{footnote}\footnote{Indeed, as maps, $(x_1,...,x_i)=(y_1,...,y_j)$ $\iff$ $i=j$ and $x_k=y_k$ for all $k=1,...,i$. Equivalently, observe that a map is completely determined by its restrictions on the members of a covering of its domain.}).

\item[c.] If $I$ is a set, then for any \ul{$I$-indexing} of sets $(x_i)_{i\in I}:I\ra Sets,~i\mapsto x_i$, the collection $\big\{(x_i)_{i\in C}~|~C\subset I~\txt{countable}\big\}$ is a set, as the image of the map $f:\{\txt{countable}~C\subset I\}\ra\txt{Sets},~C\mapsto(x_i)_{i\in C}$, and so we can represent $(x_i)_{i\in I}$ either (i) directly by its graph $\{(i,x_i)\}_{i\in I}$ or (ii) by the image ~$(x_i)_{i\in I}:=\im f=\big\{(x_i)_{i\in C}~|~C\subset I~\txt{countable}\big\}$ (\blue{footnote}\footnote{Indeed, as maps, $(x_i)_{i\in C}=(x'_i)_{i\in C'}$ $\iff$ $C=C'$ and $x_i=x_i'$ for all $i\in C=C'$. Equivalently, observe that a map is completely determined by its restrictions on the members of a covering of its domain.}).
\end{enumerate}
\end{crl}

\begin{crl}[\textcolor{blue}{The cartesian product of a set-indexing of sets is representable as a set}]\label{CartProdSets}
If $x,y$ are sets, then their cartesian product $x\times y$ is a set, since
{\small\[
x\times y=\{(a,b)~|~a\in x,b\in y\}=\Big\{\big\{\{a\},\{a,b\}\big\}~\big|~a\in x,b\in y\Big\}~\subset~\P\big(\P(x)\cup\P(y)\cup\P(x\cup y)\big).
\]}If $I$ is a set, then for any \ul{set-indexing of sets}  $(x_i)_{i\in I}$ (which we know to be itself a set), the cartesian product $\prod_{i\in I}x_i$ is itself a set, because it is cardinally smaller than the powerset $\P\left(I\times\bigcup_i x_i\right)$, via the fact that each map of sets can be uniquely viewed as a set (its graph) as follows:
{\footnotesize\[
\textstyle \prod_{i\in I}x_i=\left\{\txt{maps}~f:I\ra\bigcup_i x_i,~i\mapsto f(i)\in x_i\right\}\subset\{\txt{maps}~f:I\ra\bigcup_i x_i\}\thickapprox\left\{\txt{graphs}~G\left(f:I\ra\bigcup_i x_i\right)\right\}\lesssim \P\left(I\times\bigcup_i x_i\right).
\]}
\end{crl}

\begin{thm}[\textcolor{blue}{A \index{Set-union of sets}{set-union of sets} is a set}]~
Let $I$ be a set and $(x_i)_{i\in I}:I\ra\txt{Sets},~i\mapsto x_i$ an $I$-indexing of sets. Then the un-indexed collection of sets $\{x_i~|~i\in I\}$ is a set, as the image of a map $I\ra Sets$. Therefore, the \ul{set-union of sets} ~$\bigcup \{x_i~|~i\in I\}=\bigcup_{i\in I}x_i$ ~is a set.
\end{thm}

\begin{note}
The intersection $\bigcap c$ of any subclass of sets ~$c\subset\txt{Sets}$ ~is itself a set (since it is a subset of a member of $c$). In particular, the intersection of all sets is equal to $\emptyset$. On the other hand, the union of a subclass of sets is not necessarily a set. In particular, the union of all sets contains all sets and so is (not a set but) a proper class.
\end{note}

The proofs of the following useful identities are brief and follow directly from definitions.
\begin{lmm}[\textcolor{blue}{\index{DeMorgan's laws}{DeMorgan's laws}, \index{Complement}{Complement}}]\label{DemogLaws}
Given any collection of subsets $\{S_\al\}_{\al\in A}$ of a set $S$,
\bit
\item[] {\small $\left(\bigcup_{\al\in A} S_\al\right)^c=\bigcap_{\al\in A} S_\al^c$~~ and ~~ $\left(\bigcap_{\al\in A} S_\al\right)^c=\bigcup_{\al\in A} S_\al^c$},
\eit
where as usual, the \ul{complement} of a subset $U\subset S$ is the subset ~$U^c:=\{x\in S:x\not\in U\}$. (\blue{footnote}\footnote{Recall that for (relative) convenience, the complement $U^c\subset S$ is also written as $S\backslash U$ or as $S-U$.})
\end{lmm}

\section{Partial Ordering of Sets: Zorn's Lemma and Choice}\label{ParOrdZorn}
Recall that our preferred set inclusion operations are $\subset$ for non-strict inclusion and $\subsetneq$ for strict inclusion. In a poset we (for convenience) sometimes \ul{also} use $\subseteq$ (in place of $\subset$) for non-strict inclusion. That is $\subset,\subseteq$ each stand for ``non-strict inclusion'' while $\subsetneq$ stands for ``strict inclusion''.

We have already encountered most of this section's definitions in our discussion on classes. However, the definitions are repeated here using a more explicit/expanded format for the purposes of reminding, reviewing, and connecting with the literature on the centrally important special case of sets.

\begin{dfn}[\textcolor{blue}{
\index{Relation! between sets}{Relation between sets},
\index{Domain! of a relation}{Domain of a relation},
\index{Image of! a relation}{Image of a relation},
\index{Opposite! relation}{Opposite relation},
\index{Relation! on a set}{Relation on a set},
\index{Related elements}{Related elements},
\index{Intermediate element}{Intermediate element},
\index{Subrelation}{Subrelation}}]
Let $X,Y,A$ be sets. A \ul{(binary) relation between} $X$ and $Y$ is any subset $R\subset X\times Y$. The \ul{domain} and \ul{image} of a relation $R\subset X\times Y$ are
\bea
\dom R:=\{x\in X~|~(x,y)\in R~\txt{for some}~y\in Y\},~~\im R:=\{y\in Y~|~(x,y)\in R~\txt{for some}~x\in X\}.\nn
\eea
The \ul{opposite relation} of a relation $R\subset X\times Y$ is the relation
\[
R^{op}:=\{(y,x)\in Y\times X:(x,y)\in R\}\subset Y\times X,~~~~\txt{i.e., ~$(x,y)\in R \iff (y,x)\in R^{op}$}.
\]
A \ul{relation on} $A$ is any subset $R\subset A\times A$. (\blue{footnote}\footnote{Relations of the form $R\subset A\times A$ suffice because $X\times Y\subset(X\cup Y)\times(X\cup Y)$, and so a relation between $X$ and $Y$ is a relation on $X\cup Y$.}). If $a,b\in A$, we say ``$a$ is \ul{related} to $b$ (through $R$)'', written $aR b$ or $bR^{op}a$, if $(a,b)\in R$. For any $a,b,c\in A$, we say ``$c$ lies between $a$ and $b$ (with respect to $R$)'', or ``$c$ is \ul{intermediate} to $a$ and $b$ (with respect to $R$)'', written $aR cR b$, if $aR c$ and $cRb$.

Given a relation $R\subset X\times Y$ between $X$ and $Y$ (resp. $R\subset A\times A$ on $A$), every pair of subsets $X'\subset X$, $Y'\subset Y$ (resp. every subset $A'\subset A$) inherits the \ul{subrelation} $R':=R\cap(X'\times Y')\subset X'\times Y'$ between $X'$ and $Y'$ (resp. $R':=R\cap(A'\times A')\subset A'\times A'$ on $A'$).

\end{dfn}
\begin{dfn}[\textcolor{blue}{
\index{Reflexive relation}{Reflexive relation},
\index{Symmetric! relation}{Symmetric relation},
\index{Antisymmetric relation}{Antisymmetric relation},
\index{Transitive! relation}{Transitive relation},
\index{Strict! relation}{Strict relation},
\index{Linear! (total) relation}{Linear (total) relation}}]
Let $A$ be a set and $R\subset A\times A$. The relation $R$ is (1) \ul{reflexive} if $aRa$ for all $a\in A$, (2) \ul{symmetric} if $aRb$ implies $b R a$ for all $a,b\in A$, (3) \ul{antisymmetric} if $aRbRa$ implies $a=b$ for all $a,b\in A$, (4) \ul{transitive} if $aRbRc$ implies $aRc$ for all $a,b,c\in A$, (5) \ul{strict} if $aRb$ implies $a\neq b$ for all $a,b\in A$, (6) \ul{linear} (\ul{total}) if $aRb$ or $bRa$ for all $a,b\in A$.
\end{dfn}

Note that (unless the relation $R\subset A\times A$ is symmetric) if $a\in A$ is related to $b\in A$, it does not follow that $b$ is related to $a$. Similarly (unless the relation, if transitive, is symmetric), if $c$ lies between $a$ and $b$, it does not follow that $c$ lies between $b$ and $a$.

\begin{dfn}[\textcolor{blue}{\index{Equivalence! relation}{Equivalence relation}}]
An \ul{equivalence relation} is a relation that is reflexive, symmetric, and transitive. That is, a relation $\sim~\subset A\times A$ such that for all $a,b,c\in A$,
\bea
(i)~a\sim a,~~~~(ii)~a\sim b~~\Ra~~b\sim a,~~~~(iii)~a\sim b\sim c~~\Ra~~a\sim c.\nn
\eea
\end{dfn}

\begin{dfn}[\textcolor{blue}{
\index{Strict! partial order}{Strict partial order},
\index{Strict! poset}{Strict poset},
\index{Partial! order}{Partial order},
\index{Poset}{{Poset}},
\index{Product of! posets}{Product of posets}}]
Let $P$ be a set. A relation $<~\subset P\times P$ (with opposite relation $>~\subset P\times P$) is a \ul{strict partial order} on $P$ (making $P=(P,<)$ a \ul{strict partially ordered set} or \ul{strict poset}) if it is strict and transitive, i.e., for all $a,b,c\in P$,
\bea
(i)~a<b~~\Ra~~a\neq b~~~~\txt{and}~~~~(ii)~a<b<c~~\Ra~~a<c.\nn
\eea
A relation $\leq~\subset P\times P$ (with opposite relation $\geq~\subset P\times P$) is a \ul{partial order} on $P$ (making $P=(P,\leq)$ a \ul{partially ordered set} or \ul{poset}) if it is reflexive, antisymmetric, and transitive, i.e., for all $a,b,c\in P$,
\bea
(i)~a\leq a,~~~~(ii)~a\leq b\leq a~~\Ra~~a=b,~~~~(iii)~a\leq b\leq c~~\Ra~~a\leq c.\nn
\eea
Given a strict poset $(P,<)$, we get a unique poset via ``$a\leq b$ $\iff$ $a<b$ or $a=b$''. Likewise, given a poset $(P,\leq)$, we get a unique strict poset via ``$a<b$ $\iff$ $a\leq b$ and $a\neq b$''. Consequently, specifying the poset $(P,\leq)$ is equivalent to specifying the associated strict poset $(P,<)$, and vice versa.

If $a\leq b$ (or $b\geq a$) we say ``\ul{$a$ is less than $b$}'' (or ``\ul{$b$ is greater than $a$}''). Also, if $a\leq b$ and $a\neq b$, we write $a<b$ (or $b>a$) and say ``\ul{$a$ is strictly less than $b$}'' (or ``\ul{$b$ is strictly greater than $a$}'').

Given a family of posets $\{P_i\}_{i\in I}$, their \ul{product} is the poset
\bea
\textstyle\prod_{i\in I}P_i=\{(p_i)_{i\in I}:p_i\in P_i\},~~~~\txt{where}~~~~(p_i)_{i\in I}\leq (p'_i)_{i\in I}~~\iff~~p_i\leq p'_i~~\txt{for all}~~i\in I.\nn
\eea
\end{dfn}
Note that in a poset $P=(P,\leq)$, every subset $S\subset P$ is itself a poset $(S,\leq)$.

\begin{dfn}[\textcolor{blue}{
\index{Comparable elements}{Comparable elements},
\index{Linearly or totally ordered set (chain, list, line)}{Linearly or totally ordered set (chain, list, line) in a poset}}]
Let $(P,\leq)$ be a poset. Two elements $a,b\in P$ are \ul{comparable} if $a\leq b$ or $b\leq a$. A set $C\subset P$ is \ul{linearly ordered} (or \ul{totally ordered}, or a \ul{chain}, or a \ul{list}, or a \ul{line}) if every two elements in $C$ are \ul{comparable}, i.e., for any $a,b\in C$, we have $a\leq b$ or $b\leq a$ (equivalently, the order $\leq|_C~:=~\leq\cap(C\times C)$ on $C$ is linear).
\end{dfn}

\begin{dfn}[\textcolor{blue}{
\index{Upper bound}{Upper bound},
\index{Lower bound}{Lower bound},
\index{Strict! upper bound}{Strict upper bound},
\index{Strict! lower bound}{Strict lower bound},
\index{Directed! set}{Directed set}}]
Let $(P,\leq)$ be a poset and $A\subset P$. An element $u\in P$ is an \ul{upper bound} of $A$ (written $A\leq u$) if $a\leq u$ for all $a\in A$. An element $l\in P$ is a \ul{lower bound} of $A$ (written $l\leq A$) if $l\leq a$ for all $a\in A$.

An upper bound $u\geq A$ (i.e., $A\leq u$) is a \ul{strict upper bound}, written $A<u$, if $u\not\in A$. Similarly, a lower bound $l\leq A$ is a \ul{strict lower bound}, written $l<A$, if $l\not\in A$.

The set $A\subset P$ (or the associated poset $(A,\leq)\subset(P,\leq)$ ) is a \ul{directed set} if any two elements of $A$ have a common upper bound, i.e., for any $a,b\in A$, there exists $c\in A$ such that $a,b\leq c$ (i.e., $a\leq c$ and $b\leq c$). (\blue{footnote}\footnote{Given a directed set $D$ and $n\geq 1$, its finite cartesian power $D^n$ is also a directed set (by induction on $n$), since for any {\scriptsize $(a,a'),(b,b')\in D\times D$}, there exist $c,c'\in D$ such that $a,b\leq c$ and $a',b'\leq c'$ in $D$, and so $(a,a'),(b,b')\leq(c,c')$ in $D\times D$.})
\end{dfn}

\begin{dfn}[\textcolor{blue}{\index{Minimal! element}{Minimal element}, \index{Maximal! element}{Maximal element}, \index{Least! element}{Least element}, \index{Greatest! element}{Greatest element}}]
Let $(P,\leq)$ be a poset. An element $m\in P$ is \ul{minimal} if for all $x\in P$,~ $x\leq m$ implies $m=x$ (i.e., $m$ cannot be strictly greater than an element that is comparable with $m$). An element $m\in P$ is \ul{maximal} if for all $x\in P$,~ $m\leq x$ implies $m=x$ (i.e., $m$ cannot be strictly less than an element that is comparable with $m$).

An element $y\in P$ is a \ul{least} (resp. \ul{greatest}) element if $y\leq x$ (resp. $x\leq y$) for all $x\in P$. (\blue{footnote}\footnote{A least (resp. largest) element is a minimal (resp. maximal) element, but in general a minimal (resp. maximal) element is not a least (resp. largest) element. In a linearly ordered set, least = minimal (resp. largest = maximal).})
\end{dfn}

\begin{dfn}[\textcolor{blue}{\index{Well-ordered! set}{Well-ordered set}}]
Let $(P,\leq)$ be a poset. A subset $W\subset P$ (or the associated poset $(W,\leq)\subset(P,\leq)$ ) is a \ul{well-ordered} set if (i) $(W,\leq)$ is a chain and (ii) every nonempty subset $\emptyset\neq A\subset W$ contains a least element.
\end{dfn}

\begin{dfn}[\textcolor{blue}{
\index{Greatest! lower bound (glb or infimum)}{Greatest lower bound (glb or infimum)},
\index{glb property}{glb property}}]
Let $(P,\leq)$ be a poset, and $\emptyset\neq A\subset P$ a nonempty subset with a lower bound. An \ul{infimum} of $A$, written $\inf(A)$ or $\txt{glb}_P(A)$, is a greatest lower bound (\ul{glb}) of $A$, i.e., a greatest element in the set of all lower bounds  of $A$. Equivalently,
\bea
L_P(A):=\{l\in P~|~l\leq A\}\leq\inf(A)\leq A.\nn
\eea

A poset $(P,\leq)$ has the \ul{glb property} if every nonempty subset of $P$ with a lower bound has an infimum.
\end{dfn}

\begin{dfn}[\textcolor{blue}{\index{Least! upper bound (lub or supremum)}{Least upper bound (lub or supremum)}, \index{lub property}{lub property}}]
Let $(P,\leq)$ be a poset, and $\emptyset\neq A\subset P$ a nonempty subset with an upper bound. A \ul{supremum} of $A$, written $\sup(A)$ or $\txt{lub}_P(A)$, is a least upper bound (\ul{lub}) of $A$, i.e., a least element in the set of all upper bounds of $A$. Equivalently,
\bea
A\leq\sup(A)\leq U_P(A):=\{u\in P~|~A\leq u\}.\nn
\eea

A poset $(P,\leq)$ has the \ul{lub property} if every nonempty subset of $P$ with an upper bound has a supremum.
\end{dfn}

\begin{notation}
Let $(A,\leq)$ be a poset. If $B,C\subset A$, we will write $B\leq C$ if $b\leq c$ for all $b\in B,c\in C$ (resp. $B<C$ if $b<c$ for all $b\in B,c\in C$).
\end{notation}

\begin{dfn}[\textcolor{blue}{
\index{Initial segment}{Initial segment},
\index{Strict! initial segment}{Strict initial segment}}]
Let $(P,\leq)$ be a poset and $C,D\subset P$ chains. Then $C$ is an \ul{initial segment} of $D$ (written $C\preceq D$ or $C\leq_{\txt{seg}}D$) if $C\subseteq D$ and $C<D\backslash C$ (i.e., $C<d$ for all $d\in D\backslash C$). If $C\preceq D$ and $C\subsetneq D$ (i.e., $C\subsetneq D$ and $C<D\backslash C$), we say $C$ is a \ul{strict initial segment} of $D$, written $C\prec D$.
\end{dfn}

Note that the empty chain $\emptyset$ is a (strict) initial segment of every (nonempty) chain. Also, the collection $\C$ of all chains in $P$ gives posets $(\C,\subseteq)$ and $(\C,\preceq)$.

\begin{dfn}[\textcolor{blue}{
\index{Lower section}{Lower section},
\index{Strict! lower section}{Strict lower section},
\index{Lower segment}{Lower segment},
\index{Strict! lower segment}{Strict lower segment}}]
Let $(P,\leq)$ be a poset and $x\in P$. The \ul{lower section} at $x$ is the set of lower bounds $L_x:=L_P(x):=\{y\in P:y\leq x\}$ of $x$. Similarly, the \ul{strict lower section} at $x$ is the set of strict lower bounds $L^\ast_x:=L^\ast_P(x):=\{y\in P:y<x\}$ of $x$. (\blue{footnote}\footnote{In applications, we will sometimes define the lower section to be strict in the sense ~$L_x:=\{y\in P:y<x\}$.~ When both types of lower sections (i.e., strict and non-strict) are necessary, then we will make a distinction.})

If $(A,\leq)$ is a linearly-ordered set, then a lower section (resp. strict lower section) of $A$ will be called a \ul{lower segment} (resp. \ul{strict lower segment}).
\end{dfn}

\begin{dfn}[\textcolor{blue}{\index{Regular! poset}{Regular poset}}]
A \ul{regular poset} is a poset in which every (well-ordered) chain has an upper bound.
\end{dfn}

\begin{lmm}\label{ZrnPreLmm}
Let $(P,\leq)$ be a poset. The collections $\C=\C(P):=\{\txt{chains}~C\subset P\}$, $\L=\L(P):=\{L_x~|~x\in P\}$, and $\L^\ast=\L^\ast(P):=\{L^\ast_x~|~x\in P\}$ give posets $(\C,\subseteq)$, $(\L,\subseteq)$, and $(\L^\ast,\subseteq)$, as subsets of the poset $(\P(P),\subseteq)$. Moreover, the following hold. For all $x,y\in P$,
\bit[leftmargin=0.7cm]
\item[(a)] $x\leq y$ $\iff$ $L_x\subseteq L_y$ (and ~$x<y$ $\iff$ $L_x\subsetneq L_y$).
\item[(b)] $x$ is maximal in $P$ $\iff$ $L_x$ is maximal in $\L$.
\item[(c)] Assume $P$ is regular. If there exists a maximal chain $C\in\C=(\C,\subseteq)$, then $P$ has a maximal element in $C$ (an upper bound of $C$ in $P$).
\item[(d)] $(\C,\subseteq)$ is a regular poset. Also, $(P,\leq)$ is regular $\iff$ $(\L,\subseteq)$ is regular.
\eit
\end{lmm}
\begin{proof}
(a),(b) are immediate consequences of the definitions. For (c) let $C\leq u\in P$. Then $C\subset C\cup\{u\}$ $\Ra$ $C=C\cup\{u\}$, $\Ra$ $u\in C$. Moreover for any $x\in P$, $u\leq x$ $\Ra$ $C\leq u\leq x$, $\Ra$ (as before) $x\in C$, $\Ra$ $x\leq u$, i.e., $x=u$, and so $u$ is maximal in $P$. Finally, for (d), observe that unions give upper bounds.
\end{proof}

\begin{thm}[\textcolor{blue}{\index{Zorn's lemma}{Zorn's lemma}: \cite{gaillard2012,lewin1991}}]\label{ZornLemma}
A nonempty regular poset has a maximal element.
\end{thm}
\begin{proof}
Let $(P,\leq)$ be a nonempty regular poset (i.e., every (well-ordered) chain in $P$ has an upper bound). For any sets $A\subset B$ in $P$, let ${}_{A<}B=U^\ast_A(B):=\{b\in B:A<b\}$ be the set of strict upper bounds of $A$ in $B$, and $B_{<A}=L^\ast_A(B):=\{b\in B:b<A\}$  be the set of strict lower bounds of $A$ in $B$. That is, ~$B_{<A}~<~A~<~{}_{A<}B$.

Suppose $P$ has no maximal element. Then by Lemma \ref{ZrnPreLmm}(c), every chain $C\subset P$ is not maximal in $(\C(P),\subseteq)$, and so its set of strict upper bounds ${}_{C<}P$ is nonempty. Thus (by the axiom of choice) we have a strict upper bound selection function
\bea
u:\{\txt{nonempty chains in}~P\}\ra P,~~C\mapsto u(C)\in{}_{C<}P=U^\ast_C(P).~~~~(\txt{\blue{footnote}}\footnotemark).\nn
\eea
\footnotetext{(Note I): Only the restriction $u|_{WO(P)}$ of $u$ to well-ordered chains $WO(P)$ in $P$ is essential. (Note II): If $W\in\C(P)$ is a well-ordered chain and $w\in W$, then because (i) $w\in W\cap\big({}_{(W_{<w})}P\big)\neq\emptyset$ and (ii) $W$ is well-ordered, the subset $W\cap\big({}_{(W_{<w})}P\big)\subset{}_{(W_{<w})}P$ of strict upper bounds of $W_{<w}$ contains a least element $l\geq w$ (hence $l=w$), otherwise, $l<w$ $\Ra$ $l\in W_{<w}<l$ (a contradiction). Therefore, in problems where we are only concerned with the restriction $u|_{WO(P)}$, it is possible to choose $u|_{WO(P)}$ such that for any well-ordered chain $W\in\C(P)$, if desired, we can set $u(W_{<w}):=w$ for all $w\in W$.}
Consider the associated \emph{nonempty} poset $(\W,\subseteq)$ of those well-ordered chains in $P$ given by
\bea
\W:=\left\{W\in \C(P):~W~\txt{is well-ordered, and}~u(W_{<w})=w~\txt{for all}~w\in W\right\}.~~~~\txt{(\blue{footnote}\footnotemark)}.\nn
\eea
\footnotetext{Note that every finite chain $F\in\C(P)$ is well-ordered.}
Then for each $W\in\W$, we have $W\cup\{u(W)\}\in \W$ as well, because ~{\small$u\left((W\cup\{u(W)\})_{<u(W)}\right)=u\left(W_{<u(W)}\right)=u(W)$} implies
~{\small $u\big((W\cup\{u(W)\})_{<w}\big)=
\left\{
  \begin{array}{ll}
   u(W) , & \txt{if}~~w=u(W) \\
    u\left(W_{<w}\right)=w, & \txt{if}~~w\in W
  \end{array}
\right\}=w$}.

For any $W_1,W_2\in \W$, let $L:=\bigcup\{C\in \W:C\preceq W_1,W_2\}\in\W$ be the union of all common \emph{initial segments} of $W_1$ and $W_2$. Then we see that $L\preceq W_1,W_2$, i.e., $L\in\W$, and that $L$ is maximal in $\W$ with respect to this property (i.e., wrt $\preceq W_1,W_2$). If $L\prec W_1,W_2$, then there exist $w_1\in W_1$, $w_2\in W_2$ such that $L=(W_1){}_{<w_1}=(W_2){}_{<w_2}$, $L\cup\{w_1\}\preceq W_1$, and $L\cup\{w_2\}\preceq W_2$, which imply
\begin{align}
&w_1=u\left((W_1){}_{<w_1}\right)=u(L)=u\left((W_2){}_{<w_2}\right)=w_2,~~\Ra~~L\cup\{u(L)\}\preceq W_1,W_2,\nn\\
&~~\Ra~~L\subsetneq L\cup\{u(L)\}\in\W,\nn
\end{align}
contradicting maximality (wrt $\preceq W_1,W_2$) of $L$ in $\W$. Thus, $L=W_1$ or $L=W_2$, i.e.,
\bea
W_1\preceq W_2~~\txt{or}~~W_2\preceq W_1,~~~~\Ra~~~~W_1\subseteq W_2~~\txt{or}~~W_2\subseteq W_1,\nn
\eea
which shows $\W$ is linearly ordered, hence a chain, in $(\C(P),\subseteq)$.

Let {\small $W_0:=\bigcup\{W\in\W\}$}. Then $W_0\in \W$ (\blue{footnote}\footnote{
\ul{$W_0$ is a chain}: If $a,b\in W_0$, let $a\in W_a,b\in W_b$. Since $W_a\subseteq W_b$ or $W_b\subseteq W_a$, $a,b\in W_a$ or $a,b\in W_b$. Hence $a\leq b$ or $b\leq a$.\\
\ul{$W_0$ is well-ordered}: If $\emptyset\neq B\subset W_0$, we have $B=\bigcup_{W\in\W}B\cap W$. Since $(\W,\subseteq)$ is a chain, we have $\min(B\cap W)=\min(B\cap W')$ for all $W,W'\in\W$ such that $B\cap W\neq\emptyset$ and $B\cap W'\neq\emptyset$. Thus $B$ has a least element $\min B=\min(B\cap W)$ for all $W\in\W$ such that $B\cap W\neq\emptyset$. \ul{$W_0\in\W$}: Given $w\in W_0$, let $w\in W\in\W$. Then $u\big((W_0){}_{<w}\big)=u(W_{<w})=w$.}), and so {\small $W_0\subsetneq W_0\cup\{u(W_0)\}\in \W$} (a contradiction).
\end{proof}

\begin{rmk}\label{ZornLmmRmk}
The above proof only requires that ``every \ul{well-ordered} chain has an upper bound''. This observation will be useful in the proof of Zermelo's well-ordering theorem (Corollary \ref{ZermWoThm}).
\end{rmk}

\begin{crl}[\textcolor{blue}{\index{Hausdorff's maximality principle}{Hausdorff's maximality principle I}}]
Every poset contains a maximal chain.
\end{crl}
\begin{proof}
Let $(P,\leq)$ be a poset. The poset of chains in $P$, namely, $(\C(P),\subseteq)$ is regular by Lemma \ref{ZrnPreLmm}(d), and so has a maximal element by Zorn's lemma.
\end{proof}

\begin{crl}[\textcolor{blue}{Hausdorff's maximality principle II}]
In a poset every chain is contained in a maximal chain. ({\small Hence, by Lemma \ref{ZrnPreLmm}(c), a regular poset $(P,\leq)$ has a maximal element $\iff$ it has a maximal chain $C\in(\C(P),\subseteq)$.})
\end{crl}
\begin{proof}
Let $(P,\leq)$ be a poset and $C\in\C(P)$. Then the poset of all chains in $P$ containing $C$, namely $\big(\C_C(P),\subseteq\big)$, where $\C_C(P):=\{C'\in\C(P):C\subset C'\}$, is regular by Lemma \ref{ZrnPreLmm}(d), and so contains a maximal element $C_m$. But $C_m$ is also maximal in $\C(P)$, since every chain $C''\in\C(P)$ containing $C_m$ also contains $C$ and so lies in $\C_C(P)$.
\end{proof}

\begin{crl}[\textcolor{blue}{In a poset, every maximal element is the greatest element of a maximal chain}]
Let $(P,\leq)$ be a poset. If $x\in P$ is a maximal element, then there exists a maximal chain $C_x\in\C(P)$ such that $C_x\leq x\in C_x$ (i.e., $x$ is a maximal element of $C_x$).
\end{crl}
\begin{proof}
By the preceding result, the chain $\{x\}$ is contained in a maximal chain $C_x\in\C(P)$. Since $x$ is maximal in $P$, it follows that $C_x\leq x\in C_x$.
\end{proof}

\begin{crl}[\textcolor{blue}{\index{Zermelo's well-ordering theorem}{Zermelo's well-ordering theorem}}]\label{ZermWoThm}
Every set $S$ can be well-ordered.
\end{crl}
\begin{proof}
Let $(P,\leq)$ be the poset with $P:=\{\txt{well-orderings}~(A,\leq_A)~|~A\subset S\}$ consisting of all possible well-orderings of subsets of $S$, and (with ~$\leq_B|_A~:=~\leq_B\cap(A\times A)$~ we have)
{\small\bea
(A,\leq_A)\leq (B,\leq_B)~~~~\txt{iff}~~~~\leq_B\big|_A~=~\leq_A~~\txt{and}~~A\preceq B~~(\txt{i.e., $A$ is an initial segment of $B$}).\nn
\eea}By Remark \ref{ZornLmmRmk}, to conclude that $P$ has a maximal element by Zorn's lemma, it is enough to show that every well-ordered chain in $P$ has an upper bound. Consider a well-ordered chain {\small $W:=\{(A_\ld,\leq_{A_\ld})\}_{\ld\in\Ld}\subset P$}.

Let $A:=\bigcup_{\ld\in\Ld} A_\ld\subset S$ and define a relation $\leq_A~\subset A\times A$ by $\leq_A|_{A_\ld}~:=~\leq_{A_\ld}$, for each $\ld\in\Ld$. \ul{$(A,\leqslant_A)$ is a chain}: If $a,b\in A$, let $a\in A_{\ld_a},b\in A_{\ld_b}$. Since $A_{\ld_a}\subseteq A_{\ld_b}$ or $A_{\ld_b}\subseteq A_{\ld_a}$, $a,b\in A_{\ld_a}$ or $a,b\in A_{\ld_b}$. Thus $a\leq b$ or $b\leq a$. \ul{$(A,\leqslant_A)$ is well-ordered}: If $\emptyset\neq B\subset A$, we have $B=\bigcup_{\ld\in\Ld}B\cap A_\ld$. Since $\{A_\ld\}_{\ld\in\Ld}$ is a chain, we have $\min(B\cap A_\ld)=\min(B\cap A_{\ld'})$ for all $\ld,\ld'\in\Ld$ such that $B\cap A_\ld\neq\emptyset$ and $B\cap A_{\ld'}\neq\emptyset$. Thus $B$ has a least element $\min B:=\min(B\cap A_\ld)$ for all $\ld$ such that $B\cap A_\ld\neq\emptyset$. (\blue{footnote}\footnote{\ul{Alternative proof that $(A,\leqslant_A)$ is well-ordered}: Assume wlog that $W:=\{A_\ld\}_{\ld\in\Ld}$ via an injective indexing $\Ld\ra P,~\ld\mapsto A_\ld$. Order $\Ld$ by declaring ``$\ld\geq\ld'$ if $A_\ld\leq A_{\ld'}$''. Then by construction, every nonempty subset of $\Ld$ contains a greatest element (i.e., $\Ld$ is reverse-well-ordered since $W$ is well-ordered). \\
Define the poset $(A,\leq_A)$ as before, which is a chain as we have already shown. Let $\emptyset\neq C\subset A$. For each $c\in C$, let $c\in A_{\ld_c}$. Define $\Ld_C:=\{\ld_c:c\in C\}$. Then $\Ld_C$ contains a greatest element, say $g$. Since $\ld_c\leq g$ for all $c\in C$, we have $A_g\leq A_{\ld_c}$ for all $c\in C$. Thus, if $c_l\in A_g\cap C$ is the least element of $A_g\cap C$, then $c_l\in C$ is the least element in $C=\bigcup_{c\in C}A_{\ld_c}\cap C$.}). It follows that $(A,\leq_A)\in P$ is an upper bound of $W$ in $P$.

By Zorn's lemma (or Remark \ref{ZornLmmRmk}) $P$ has a maximal element $(M,\leq_M)\in P$. Suppose $M\neq S$, i.e., there exists $e\in S\backslash M$. Then $M_e:=M\cup\{e\}$ can be well-ordered as $(M_e,\leq_{M_e})$ by declaring that (i) $\leq_{M_e}\big|_M~:=~\leq_M$ and (ii) $m\leq_{M_e}e$ (i.e., $m<_{M_e}e$) for all $m\in M$. Hence $(M,\leq_M)<(M_e,\leq_e)\in P$, contradicting the maximality of $(M,\leq_M)$ in $P$.
\end{proof}

\begin{crl}[\textcolor{blue}{Zorn's lemma is equivalent to the axiom of choice for sets}]
\end{crl}
\begin{proof}
Zorn's lemma implies Corollary \ref{ZermWoThm}. Corollary \ref{ZermWoThm} in turn implies the \ul{axiom of choice for sets}, i.e., if $I$ is a set and $(s_i)_{i\in I}$ an indexing of nonempty (disjoint) sets, then there exists a map
\bea
\textstyle f:I\ra\bigcup_{i\in I}s_i,~~i\mapsto f(i):=\min_R s_i~\in~s_i,~~~~\txt{i.e.,}~~~~\prod_{i\in I}s_i\neq\emptyset,\nn
\eea
where $\min s_i$ is based on some well-ordering $R\subset S\times S$ of the set $S:=\bigcup_{i\in I}s_i$.

Conversely, the axiom of choice for sets implies Zorn's lemma (by the proof of Zorn's Lemma).
\end{proof} 

%% file: parts/AlgebraNC/AlgebraNC2.tex
\chapter{Basic Set Theory II: Classical Number Systems}
For this chapter, if needed, additional sources of reading include for example \cite{enderton1977,goldrei1996,rudin1976,gallian2013,cohn1982}. Also, in this chapter, it will often not be necessary to assume knowledge of all of the axioms of set theory (i.e., the ZFC-axioms) given in the previous chapter.
\section{Natural Numbers}\label{NatNumSec}
\subsection{Properties and arithmetic of natural numbers}\label{NatNumArth}
\begin{dfn}[\textcolor{blue}{
\index{Natural! numbers $\Natural$}{Natural numbers $\Natural$},
\index{Order! on $\Natural$}{Order on $\Natural$},
\index{Successor function}{Successor function},
\index{Predecessor function}{Predecessor function}}]
Recall from the set theory axioms that in the set of natural numbers $\Natural=\{0,1,2,\cdots\}$, which is inductive, we have
\bea
0:=\emptyset~~\subset~~ 1:=0\cup\{0\}~~\subset~~ 2:=1\cup\{1\}~~\subset~~ 3:=2\cup\{2\}~~\subset~~\cdots.\nn
\eea
The map $f_+:\Natural\ra\Natural$, $n\mapsto n^+:=n\cup\{n\}$, called the \ul{successor function} of $\Natural$, is injective by construction. (\blue{footnote}\footnote{(i) By construction, it is clear that for any natural number $n\in\Natural$, we have $n\cap\{n\}=\emptyset$, i.e., $n\not\in n$. (ii) Also, observe that $n\in n\cup\{n\}$ and $n\subset n\cup\{n\}$. So, if $n\cup\{n\}=n'\cup\{n'\}$, then $n\in n'\cup\{n'\}$ (and $n\subset n'\cup\{n'\}$) and $n'\in n\cup\{n\}$ (and $n'\subset n\cup\{n\}$). In particular, ``$n\in n'$ or $n=n'$'' and ``$n'\in n$ or $n=n'$'', $\Ra$ $n=n'$ (which shows $f_+$ is injective).}). Consequently, we also get a well defined map $f_-:\Natural\backslash\{0\}\ra\Natural,~n^+\mapsto(n^+)^-:=n$ (call it the \ul{predecessor function} of $\Natural$). By definition, we have $f_-\circ f_+=id_\Natural$, and so for any $n\in\Natural$, by applying $f_+$ we get $f_+f_-f_+(n)=f_+(n)$. This means $f_+\circ f_-|_{f_+(\Natural)}=id_{f_+(\Natural)}$, and so by the surjectivity of $f_+$ onto $\Natural\backslash\{0\}$, we get $f_+\circ f_-=id_{\Natural\backslash\{0\}}$. That is, $f_-\circ f_+=id_\Natural$ and $f_+\circ f_-=id_{\Natural\backslash\{0\}}$. Also, for each $n\in\Natural$, we have both $n\in n^+$ and $n\subset n^+$.

We define a \ul{linear order} relation $\leq$ on $\Natural$ by ``$m<n$ if $m\in n$'' (i.e., ``$m\leq n$ if $m\in n$ or $m=n$''). With respect to this ordering,
\bea
n^-<n<n^+,~~~~\txt{for all}~~~~n\in\Natural,~~(\txt{with ~$0^-\not\in\Natural$~ as an extra symbol}).\nn
\eea
\end{dfn}

\begin{lmm}[\textcolor{blue}{\index{Well-ordering principle for $\Natural$}{Well-ordering principle for $\Natural$}}]
Every nonempty subset $A\subset\Natural$ has a least element.
\end{lmm}
\begin{proof}
Suppose $A$ has no least element. Let $B:=\{b\in\Natural~|~\txt{for all}~n\in\Natural,~~n\leq b~\Ra~n\not\in A\}=\{b\in\Natural~|~L_b\cap A=\emptyset\}\subset\Natural\backslash A$.
\bit[leftmargin=1cm]
\item Since $A$ has no least element, it follows that $0\not\in A$. Thus, $0\in B$, since ~$n\leq 0~~\Ra~~n=0\not\in A$.
\item Let $b\in B$. If $n\leq b^+$, then either $n\leq b$ or $n=b^+$.  If $n\leq b$, then $n\not\in A$ (and so $b^+\in B$). If $n=b^+$, then again $n\not\in A$ (otherwise, if $n=b^+\in A$, then $b^+$ is a least element of $A$, since every $m\leq b$ is not in $A$, a contradiction), and so $b^+\in B$.
\eit
This shows $B$ is an inductive subset of $\Natural$, and so $B=\Natural$, which implies $A=\emptyset$ (a contradiction).
\end{proof}

\begin{crl}[\textcolor{blue}{\index{First principle of induction}{First principle of induction}}]\label{PrincInd1}
Let $A\subset\Natural$ and $a\in A$. If $A$ has the property
\bea
\label{FrstIndHypEq}\txt{``for all ~$n\geq a,~~n\in A~~\Ra~~n^+\in A$''},
\eea
then $A$ contains every natural number $n\geq a$.
\end{crl}
\begin{proof}
Suppose the set $B:=\{b\in\Natural~|~a\leq b\not\in A\}\neq\emptyset$. Then by well-ordering, $B$ has a least element $m>a$ (since $B\subset\Natural\backslash A$), i.e., $m^-\geq a$. Thus, by hypotheses (\ref{FrstIndHypEq}) (which imply $m\not\in A$ $\Ra$ $m^-\not\in A$),
\bea
a<m\not\in A~~\Ra~~a\leq m^-\not\in A,~~\Ra~~m^-\in B,\nn
\eea
which contradicts the minimality of $m$ in $B$.
\end{proof}

\begin{crl}[\textcolor{blue}{\index{Second (strong) principle of induction}{Second (strong) principle of induction}}]\label{PrincInd2}
Let $A\subset\Natural$ and $a\in A$. If $A$ has the property
\bea
\label{SecIndHypEq} \txt{``for all ~$n\geq a,~~~[a,n):=\{k\in\Natural~|~a\leq k<n\}\subset A~~\Ra~~n\in A$''},
\eea
then $A$ contains every natural number $n\geq a$.
\end{crl}
\begin{proof}
Suppose the set $B:=\{b\in\Natural~|~a\leq b\not\in A\}\neq\emptyset$. Then by well-ordering, $B$ has a least element $m>a$ (since $B\subset\Natural\backslash A$), i.e., $m^-\geq a$. Thus, by hypotheses (\ref{SecIndHypEq}),
\bea
m\not\in A~~\Ra~~k\not\in A~~\txt{for some}~~a\leq k<m,~~\Ra~~a\leq k\not\in A,~~\Ra~~k\in B,\nn
\eea
which contradicts the minimality of $m$ in $B$.
\end{proof}

\begin{lmm}[\textcolor{blue}{\index{Recursion function}{Recursion function} of $\Natural$: \cite[Recursion theorem, p.73]{enderton1977}}]\label{RecurThm}
Recall the successor function $f_+:\Natural\ra\Natural,~n\mapsto n^+:=n\cup\{n\}$ of $\Natural$.

Let $X$ be a set, $x_0\in X$, and $F:X\ra X$ any map. Then there exists a unique map $h=h(x_0,F):\Natural\ra X$ (call it the \ul{$(x_0,F)$-recursion function} of $\Natural$) such that
\bea
\label{InductEq}h(0)=x_0,~~~~(h\circ f_+)(n):=h(n^+)=F\big(h(n)\big),~~~~\txt{for all}~~n\in\Natural.
\eea
\bea\bt
(\Natural,0)\ar[d,dashed,"h"]\ar[rr,"f_+"] && (\Natural,0^+)\ar[d,dashed,"h"]\\
(X,x_0)\ar[rr,"F"] && (X,F(x_0))
\et~~~~~~~~h\circ f_+=F\circ h.\nn
\eea
\end{lmm}
\begin{proof}
Given a map $h_A:A\subset\Natural\ra X$, we will say $h_A$ is an $h$-approximation (i.e., an approximation to $h$) if $h_A$ satisfies the following:
\bit
\item[(i)] If $0\in A$, then $h_A(0)=x_0$.
\item[(ii)] For any $n\in\Natural$, if $n^+\in A$, then ~$n\in A$~ and ~$(h_A\circ f_+)(n):=h_A(n^+)=F\big(h_A(n)\big)$.
\eit
Let ~$P:=\{\txt{$h$-approximations}~h_A:A\subset\Natural\ra X\}$,~ and define
\bea
\textstyle H:=\bigcup P=\bigcup_{h_A\in P}h_A:=\bigcup_{h_A\in P}\big\{\big(a,h_A(a)\big)~|~a\in A\big\}.\nn
\eea
Let $M:=\dom H:=\bigcup_{h_A\in P} A$. It remains to verify the following properties:
\begin{enumerate}[leftmargin=0.9cm]
\item \ul{$H$ is a map $h_M:M\subset\Natural\ra X$.}~ We need to show that all $h$-approximations agree at intersections, i.e., for any $h_A,h_{A'}\in P$, we have ~$h_A|_{A\cap A'}=h_{A'}|_{A\cap A'}$. Equivalently, we need to show the following set is inductive (and so equals $\Natural$):
\bea
&&S:=\{n\in\Natural~|~(n,x)\in H~\txt{for at most one}~x\in X\}\nn\\
&&~~~~=\big\{n\in\Natural~|~\txt{the set}~X_n:=\{h_A(n):h_A\in P\}~\txt{contains at most one element}\big\}.\nn
\eea
If $(0,x),(0,x')\in H$, then $x=h_A(0)$ and $x'=h_{A'}(0)$ for some $h_A,h_A'\in P$, and so
\bea
x=h_A(0)=x_0=h_{A'}(0)=x',~~\Ra~~0\in S.\nn
\eea
Now, let $n\in S$, i.e., if $(n,x),(n,x')\in H$ then $x=x'$ (definition of $S$). We need to show $n^+\in S$. Suppose $(n^+,x),(n^+,x')\in H$. Then $x=h_{A_+}(n^+)$, $x'=h_{A_+'}(n^+)$ for some $h_{A_+},h_{A_+'}\in P$. Thus,
\bea
x=h_{A_+}(n^+)=F(h_{A_+}(n)),~~~~x'=h_{A_+'}(n^+)=F(h_{A_+'}(n)),\nn
\eea
Since $n\in S$ and $(n,h_{A_+}(n)),(n,h_{A_+'}(n))\in H$, it follows that $h_{A_+}(n)=h_{A_+'}(n)$, and so
\bea
x=h_{A_+}(n^+)=F(h_{A_+}(n))=F(h_{A_+'}(n))=h_{A_+'}(n^+)=x',~~\Ra~~n^+\in S.\nn
\eea

\item \ul{$H=h_M$ is an $h$-approximation, i.e., $H=h_M\in P$.}~ (i) If $0\in M$, then $0\in A$ for some $h_A\in P$, and so $H(0)=h_A(0)=x_0$. (ii) For any $n\in\Natural$, if $n^+\in M$, then $n^+\in A$ for some $h_A\in P$, and so
\bea
n\in A\subset M~~\txt{and}~~H(n^+)=h_A(n^+)=F\big(h_A(n)\big)=F(H(n)).\nn
\eea

\item \ul{$M=\Natural$.}~ It suffices to show $M$ is inductive. The constant map $h_{\{0\}}=\{(0,x_0)\}:\{0\}\ra X,~0\mapsto x_0$ is in $P$, and so $0\in M$. Suppose $n\in M$. We need to show $n^+\in M$. Define
\bea
H_n:=H\cup\big\{\big(n^+,F(H(n))\big)\big\},~~\Ra~~M_n:=\dom H_n=M\cup\{n^+\}.\nn
\eea
(i) It is clear that $H_n(0)=H(0)=x_0$. (ii) For any $k\in\Natural$, if $k^+\in M_n$, then $k^+\in M$ or $k^+=n^+$. If $k^+\in M$, then $k\in M\subset M_n$ and {\small $H_n(k^+)=H(k^+)=F(H(k))=F(H_n(k))$}. If $k^+=n^+$ then (by injectivity of $f_+$) $k=n\in M\subset M_n$ and {\small $H_n(k^+)=H(k^+)=H(n^+)=F\big(H(n)\big)=F\big(H(k)\big)=F\big(H_n(k)\big)$}.

It follows that $H_n\in P$, and so $M_n=M$ (by the maximality, wrt $\subset$, of $H$ in $P$), i.e., $n^+\in M$.

\item \ul{$H=h_M$ is unique.}~ Let $h,h':\Natural\ra X$ be maps satisfying (\ref{InductEq}). It suffices to show the set
\bea
S:=\{n\in\Natural~|~h(n)=h'(n)\}\nn
\eea
is inductive. We have $h(0)=x_0=h'(0)$, and so $0\in S$. Let $n\in S$, i.e., $h(n)=h'(n)$. Then
\[
h(n^+)=F(h(n))=F(h'(n))=h(n^+). \qedhere
\]
\end{enumerate}
\end{proof}

\begin{dfn}[\textcolor{blue}{Arithmetic:
\index{Addition! in $\Natural$}{Addition in $\Natural$},
\index{Multiplication! in $\Natural$}{Multiplication in $\Natural$},
\index{Exponentiation in $\Natural$}{Exponentiation in $\Natural$}}]~\\~
Fix $m\in\Natural$. Our intent is to transform $m$ using (i.e., \ul{add to $m$}, \ul{multiply $m$ by}, or \ul{exponentiate $m$ by}) a variable $n\in\Natural$ in a recursive manner, leading to a desired binary/arithmetic operation (i.e., a map of the form) $\Natural\times\Natural\ra\Natural$.

\bit[leftmargin=0.5cm]
\item In Lemma \ref{RecurThm}, if we set $X:=\Natural$, $x_0:=m$, and $F:=f_+$, then there exists a unique map $~(~h~=~)~A_m:\Natural\ra\Natural$ (``\ul{adding to $m$}'') such that
\bea
\label{AddDefEq1}A_m(0)=m,~~~~A_m(n^+)=\big(A_m(n)\big)^+,~~~~\txt{i.e.,}~~~~A_m\circ f_+=f_+\circ A_m.
\eea

\bea\bt
(\Natural,0)\ar[d,dashed,"A_m"]\ar[rr,"f_+"] && (\Natural,0^+)\ar[d,dashed,"A_m"]\\
(\Natural,m)\ar[rr,"f_+"] && (\Natural,m^+)
\et~~~~~~~~A_m\circ f_+=f_+\circ A_m.\nn
\eea
We define \ul{addition} $+:\Natural\times\Natural\ra\Natural$, $(m,n)\mapsto m+n$ on $\Natural$ by
\bea
\label{AddDefEq2}m+n:=A_m(n),~~~~\txt{for all}~~m,n\in\Natural.
\eea
It follows that for any $m,n\in\Natural$, we have ~$m+0=m$ ~and ~$m+n^+=(m+n)^+$,~ and so
\bea
\label{AddDefEq3} m+0=m~~~~\txt{and}~~~~m+n=(m+n^-)^+=(m+n^+)^-.
\eea
Therefore, ~$n+1=(n+1^-)^+=(n+0)^+=n^+=f_+(n)$,~ and so
\bea
\label{AddDefEq4} n+0=n,~~~n+1=f_+(n),~~~n+2=f_+^2(n),~~~n+3=f_+^3(n),~~~\cdots,~~~n+m=f_+^m(n),~~~\cdots
\eea
where $f_+^m=f_+^{m^-}\circ f_+$ is the $m$-fold composition of the successor function $f_+$.

\item Similarly, in Lemma \ref{RecurThm}, if we set $X:=\Natural$, $x_0:=0$, and $F:=A_m:\Natural\ra\Natural$, $n\mapsto n+m$, then there exists a unique map $~(~h~=~)~M_m:\Natural\ra\Natural$ (``\ul{multiplying $m$ by}'') such that
\bea
\label{MultDefEq1}M_m(0)=0,~~~~M_m(n^+)=M_m(n)+m,~~~~\txt{i.e.,}~~~~M_m\circ f_+=A_m\circ M_m.
\eea
\bea\bt
(\Natural,0)\ar[d,dashed,"M_m"]\ar[rr,"f_+"] && (\Natural,0^+)\ar[d,dashed,"M_m"]\\
(\Natural,0)\ar[rr,"A_m"] && (\Natural,A_m(0))
\et~~~~~~~~M_m\circ f_+=A_m\circ M_m.\nn
\eea
We define \ul{multiplication} $\cdot:\Natural\times\Natural\ra\Natural$, $(m,n)\mapsto m\cdot n$ on $\Natural$ by
\bea
\label{MultDefEq2}m\cdot n:=M_m(n),~~~~\txt{for all}~~m,n\in\Natural.
\eea
It follows that for any $m,n\in\Natural$, we have ~$m\cdot 0=0$ ~and ~$m\cdot n^+=m\cdot n+m$,~ and so
\bea
\label{MultDefEq3} m\cdot 0=m~~~~\txt{and}~~~~m\cdot n=m\cdot n^-+m.
\eea
Therefore, $m\cdot 1=m\cdot 1^-+m=m\cdot 0+m=A_m(0)=m$, $m\cdot 2=m\cdot 1+m=m+m=(A_m)^2(0)$, and so
\bea
\label{MultDefEq4} m\cdot 0=0,~~m\cdot 1=m,~~~m\cdot 2=m+m,~~~m\cdot 3=m+m+m,~~~\cdots,~~~m\cdot n=(A_m)^n(0),~~~\cdots
\eea
where $(A_m)^n=(A_m)^{n^-}\circ A_m$ is the $n$-fold composition of ``adding to $m$''.

\item Finally, in Lemma \ref{RecurThm}, if we set $X:=\Natural$, $x_0:=1$, and $F:=M_m:\Natural\ra\Natural,~n\mapsto m\cdot n$, then there exists a unique map $~(~h~=~)~E_m:\Natural\ra\Natural$ (``\ul{exponentiating $m$ by}'') such that
\bea
\label{ExpDefEq1} E_m(0)=1,~~~~E_m(n^+)=m\cdot E_m(n),~~~~\txt{i.e.,}~~~~E_m\circ f_+=M_m\circ E_m.
\eea
\bea\bt
(\Natural,0)\ar[d,dashed,"E_m"]\ar[rr,"f_+"] && (\Natural,0^+)\ar[d,dashed,"E_m"]\\
(\Natural,1)\ar[rr,"M_m"] && (\Natural,M_m(1))
\et~~~~~~~~E_m\circ f_+=M_m\circ E_m.\nn
\eea
We define \ul{exponentiation} $\Natural\times\Natural\ra\Natural$, $(m,n)\mapsto m^n$ on $\Natural$ by
\bea
\label{ExpDefEq2}m^n:=E_m(n),~~~~\txt{for all}~~m,n\in\Natural.
\eea
It follows that for any $m,n\in\Natural$, we have ~$m^0=1$ ~and ~$m^{n^+}=m^n\cdot m$,~ and so
\bea
\label{ExpDefEq3} m^0=1~~~~\txt{and}~~~~m^n=m^{n^-}\cdot m.
\eea
Therefore, $m^1=m^{1^-}\cdot m=m^0\cdot m=m$, $m^2=m^1\cdot m=m\cdot m=(M_m)^2(1)$, and so
\bea
\label{ExpDefEq4} m^0=1,~~m^1=m,~~~m^2=m\cdot m,~~~m^3=m\cdot m\cdot m,~~~\cdots,~~~m^n=(M_m)^n(1),~~~\cdots
\eea
where $(M_m)^n=(M_m)^{n^-}\circ M_m$ is the $n$-fold composition of ``multiplying $m$ by''.
\eit
\end{dfn}

\subsection{Integer and Rational numbers}\label{IntRatNumSec}
\begin{dfn}[\textcolor{blue}{
\index{Integers $\Integer$}{Integers $\Integer$},
\index{Rational! numbers $\Rational$}{Rational numbers $\Rational$},
\index{Intrinsic density of rationals}{Intrinsic density of rationals},
\index{Absolute value}{Absolute value}: \cite{enderton1977}}]
The set of \ul{integers} is the quotient set
\bea
\label{IntDefEq}\textstyle\Integer:={\Natural\times\Natural\over\sim}:=\big\{[(a,b)]~|~a,b\in\Natural\big\}=\big\{[(n,0)]~|~n\in\Natural\big\}\cup\big\{-[(n,0)]:=[(0,n)]~|~n\in\Natural\big\},
\eea
where the equivalence relation $\sim$ on $\Natural\times\Natural$ is given by
\bea
(a,b)\sim(a',b')~~~~\txt{if}~~~~a+b'=a'+b,\nn
\eea
with \ul{addition} $+:\Integer\times\Integer\ra\Integer,~(x,x')\mapsto x+x'$, \ul{multiplication} $\cdot:\Integer\times\Integer\ra\Integer,~(x,x')\mapsto x\cdot x'$, ``\ul{zero}'' (i.e., $0_\Integer$), and ``\ul{one}'' (i.e., $1_\Integer$) respectively given by
\[
[(a,b)]+[(c,d)]:=[(a+c,b+d)],~~[(a,b)]\cdot[(c,d)]:=[(ac+bd,ad+bc)],~~0_\Integer:=[(0,0)],~~1_\Integer:=[(1,0)].
\]

The above operations on $\Integer$ are well-defined (as maps) because if $(a,b)\sim(a',b')$ and $(c,d)\sim(c',d')$, i.e., if $a+b'=a'+b$ and $c+d'=c'+d$, then
{\small\begin{align}
\label{IntWDeq1}&a+b'+c+d'=a'+b+c'+d,~~\Ra~~(a+c,b+d)\sim(a'+c',b'+d'),\\
& (a+b')c=(a'+b)c,~~(a'+b)d=(a+b')d,~~~~a'(c+d')=a'(c'+d),~~b'(c'+d)=b'(c+d'),\nn\\
\label{IntWDeq2}&~~\Ra~~a'c'+b'd'+ad+bc=ac+bd+a'd'+b'c',~~\Ra~~[(ac+bd,ad+bc)]\sim[(a'c'+b'd',a'd'+b'c')].
\end{align}}We define a \ul{linear order} relation $\leq$ on $\Integer$ by
\bea
[(a,0)]\leq[(b,0)]~~\txt{and}~~-[(b,0)]\leq-[(a,0)]~~~~\txt{if}~~~~a\leq b.\nn
\eea
\ul{Subtraction} on $\Integer$ is the operation ~$\Integer\times\Integer\ra\Integer,~(m,n)\mapsto m-n:=m+(-n)$, where $-n\in\Integer$ (the \ul{negative} or \ul{additive inverse} of $n$) is the unique integer such that $n+(-n)=0$. \ul{Distance} (or \ul{separation}) on $\Integer$ is given by the map $\Integer\times\Integer\ra\Integer,~(m,n)\mapsto|m-n|$, where $|m-n|$ is called the \emph{distance between} $m$ and $n$. The \ul{natural numbers} are contained in $\Integer$ via the injection ~$\Natural\hookrightarrow\Integer$, $n\mapsto[(n,0)]$.

Now consider the equivalence relation $\sim$ on $\Integer\times\Integer$ (\blue{footnote}\footnote{Note that, although we are using the same symbol $\sim$, this new equivalence relation is entirely different from the equivalence relation that was just used above to construct $\Integer$ from $\Natural$.}) given by
\bea
\textstyle(a,b)\sim (a',b')~~~~\txt{if}~~~~b\neq 0,~~b'\neq 0,~~\txt{and}~~ab'=a'b.\nn
\eea
Then, in terms of the equivalence classes ${a\over b}:=[(a,b)]$, the set of \ul{rational numbers} (\ul{fractions of integers}) is
\bea
\label{RatDefEq}\textstyle \Rational:={\Integer\times\Integer\over\sim}:=\left\{{a\over b}~|~a,b\in \Integer,~b\neq 0\right\}
\eea
with \ul{addition} $+:\Rational\times\Rational\ra\Rational,~(x,x')\mapsto x+x'$, \ul{multiplication} $\cdot:\Rational\times\Rational\ra\Rational,~(x,x')\mapsto x\cdot x'$, ``\ul{zero}'' (i.e., $0_\Rational$), and ``\ul{one}'' (i.e., $1_\Rational$) respectively given by
\bea
\textstyle{a\over b}+{a'\over b'}:={ab'+a'b\over bb'},~~~~{a\over b}\cdot{a'\over b'}:={aa'\over bb'},~~~~0_\Rational:={0\over 1_\Integer},~~~~1_\Rational:={1_\Integer\over 1_\Integer},\nn
\eea
which are well-defined (as maps), in the sense that
\bea
\label{RatWDeq}\textstyle {a\over b}\sim{c\over d}~~\txt{and}~~{a'\over b'}\sim{c'\over d'}~~~~\Ra~~~~{a\over b}+{a'\over b'}\sim{c\over d}+{c'\over d'}~~\txt{and}~~{a\over b}\cdot{a'\over b'}\sim{c\over d}\cdot{c'\over d'},
\eea
because ${a\over b}\sim{c\over d}$ and ${a'\over b'}\sim{c'\over d'}$ together imply $ad=cb$ and $a'd'=c'b'$, and so
{\small\begin{align}
(ab'+a'b)dd'=(adb'd'+a'd'bd)=(cbb'd'+c'b'bd)=(cd'+c'd)bb',~~(aa')(dd')=ada'd'=cbc'b'=(cc')(bb').\nn
\end{align}}We define a \ul{linear order} relation on $\Rational$ by
\bea
\textstyle{a\over b}\leq{a'\over b'}~~~~\txt{if}~~~~b>0,~~b'>0,~~\txt{and}~~ab'\leq a'b.\nn
\eea
This ordering of rational numbers is \ul{intrinsically dense} in the sense that between any two distinct rational numbers ${a\over b}<{a'\over b'}$ lies the distinct rational number ${1\over 2}\left({a\over b}+{a'\over b'}\right)={ab'+a'b\over 2bb'}$, i.e.,
\bea
\textstyle{a\over b}<{a'\over b'}~~~~\Ra~~~~{a\over b}<{1\over 2}\left({a\over b}+{a'\over b'}\right)<{a'\over b'}.\nn
\eea
The \ul{integers} are contained in $\Rational$ via the injection ~$\Integer\hookrightarrow\Rational$, $a\mapsto{a\over 1_\Integer}$. We define the \ul{absolute value} $|~|:\Rational\ra\Rational$ of rational numbers by
\bea
|q|:=
\left\{
  \begin{array}{ll}
    q, & q\geq 0 \\
    -q, & q\leq 0
  \end{array}
\right\},~~~~\txt{for}~~~~q\in\Rational.\nn
\eea
\ul{Subtraction} on $\Rational$ is the operation {\small $\Rational\times\Rational\ra\Rational,~(q,p)\mapsto q-p:=q+(-p)$}, where $-p\in\Rational$ (the \ul{negative} or \ul{additive inverse} of $p$) is the unique rational such that $p+(-p)=0$. \ul{Distance} (or \ul{separation}) on $\Rational$ is given by the map {\small $\Rational\times\Rational\ra\Rational,~(q,p)\mapsto|q-p|$}, where $|q-p|$ is called the \emph{distance between} $q$ and $p$. \ul{Division} on $\Rational$ is the operation {\small $\Rational\times(\Rational\backslash\{0\})\ra\Rational,~(q,p)\mapsto {q\over p}:=q\cdot p^{-1}$}, where $p^{-1}\in \Rational$ (the \ul{reciprocal} or \ul{multiplicative inverse} of $p$) is the unique rational such that $p\cdot p^{-1}=1$. In a rational number of the form ${q\over p}$ (called a \ul{fractional rational number} or a \ul{rational fraction}), we call $q$ the \ul{numerator} and $p$ the \ul{denominator} of ${q\over p}$.
\end{dfn}

\subsection{Real and Complex numbers}\label{ReComNumSec}
One construction of the real numbers uses rational Cauchy sequences as follows.
\begin{dfn}[\textcolor{blue}{\index{Rational! Cauchy sequences}{Rational Cauchy sequences}, \index{Real numbers $\Real$}{Real numbers $\Real$}: \cite{enderton1977}}]
A \ul{rational Cauchy sequence} (i.e., a Cauchy sequence in $\Rational$) is a map $s:\Natural\ra\Rational$, $i\mapsto s_i$ such that ``$|s_i-s_j|\ra 0$ as $i,j\ra\infty$'' (or $\lim_{i,j\ra\infty}|s_i-s_j|=0$), in the sense that for any $\vep\in \Rational$, $\vep>0$, there exists $n=n_\vep\in\Natural$ such that (\blue{footnote}\footnote{Note that a Cauchy sequence is \ul{bounded} because if $|s_i-s_j|<\vep$ for all $i,j\geq n=n_\vep$, then
\bea
|s_i-s_n|<\vep~~\txt{or}~~s_n-\vep<s_i<s_n+\vep~~~~\txt{for all}~~i\geq n.\nn
\eea})
\bea
|s_i-s_j|<\vep~~~~\txt{for all}~~i,j\geq n.\nn
\eea
Let $C(\Rational^\Natural):=\big\{\txt{Cauchy sequence}~s:\Natural\ra\Rational\big\}$. Define an equivalence relation $\sim$ on $C(\Rational^\Natural)$ by
\bea
\textstyle s\sim s'~~~\txt{if}~~~\lim_{i\ra\infty}|s_i-s_i'|=0,\nn
\eea
where $\lim_{i\ra\infty}|s_i-s_i'|=0$ means ``for any $\vep\in\Rational$, $\vep>0$, there exists $n=n_\vep\in\Natural$ such that $|s_i-s_i'|<\vep$ for all $i\geq n$''. The set of \ul{real numbers} is the quotient set
\bea
\label{RealDefEq1}\textstyle\Real:={C(\Rational^\Natural)\over\sim}:=\big\{[s]~|~s\in C(\Rational^\Natural)\big\},
\eea
with \ul{addition} $+:\Real\times\Real\ra\Real,~(x,x')\mapsto x+x'$, \ul{multiplication} $\cdot:\Real\times\Real\ra\Real,~(x,x')\mapsto x\cdot x'$, ``\ul{zero}'', and ``\ul{one}'' respectively given by (well-defined operations)
\bea
[s]+[s']:=[s+s'],~~~~[s]\cdot[s']:=[s\cdot s'],~~~~0_\Real:=\big[c_{0_\Rational}\big],~~~~1_\Real:=\big[c_{1_\Rational}\big],\nn
\eea
where
\bit
\item sequences are added and multiplied pointwise as follows:
\bea
s+s':\Natural\ra\Rational,~i\mapsto s_i+s'_i,~~~~\txt{and}~~~~s\cdot s':\Natural\ra\Rational,~i\mapsto s_i\cdot s'_i\nn
\eea
\item for any $q\in\Rational$, $c_q$ denotes the constant sequence at $q$ given by
\bea
c_q:\Natural\ra\Rational,~~i\mapsto q.\nn
\eea
\eit
As maps, addition and multiplication in $\Real$ are well-defined (i.e., $s\sim t$ and $s'\sim t'$ together imply $s+s'\sim t+t'$ and $s\cdot s'\sim t\cdot t'$) because if $|s_i-t_i|\ra 0$ and $|s'_i-t'_i|\ra 0$, then
\[
|s_i+s_i'-(t_i+t_i')|\leq|s_i-t_i|+|s_i'-t_i'|\ra 0~~~~\txt{and}~~~~|s_is_i'-t_it_i'|\leq|s_i-t_i||s_i'|+|t_i||s_i'-t_i'|\ra 0.
\]
We define a \ul{linear order} relation on $\Real$ by
{\small\[
[s]\leq[s']~~~~\txt{if}~~~~[s]=[s']~~\txt{or}~~\txt{``there exists $\vep\in\Rational$, $\vep>0$, and $n=n_\vep\in\Natural$ such that}~~\vep+s_i<s'_i~~\txt{for all}~~i\geq n\txt{''}.
\]}The \ul{rational numbers} are contained in $\Real$ via the injection ~$\Rational\hookrightarrow\Real$, $q\mapsto[c_q]$, where $c_q:\Natural\ra\Rational$, $i\mapsto q$ is the constant sequence at $q\in\Rational$. As for rational numbers, we define the \ul{absolute value} $|~|:\Real\ra\Real$ of real numbers by
\bea
|r|:=
\left\{
  \begin{array}{ll}
    r, & r\geq 0 \\
    -r, & r\leq 0
  \end{array}
\right\},~~~~\txt{for}~~~~r\in\Real.\nn
\eea
\ul{Subtraction} on $\Real$ is the operation ~$\Real\times\Real\ra\Real,~(r,s)\mapsto r-s:=r+(-s)$, where $-s\in\Real$ (the \ul{negative} or \ul{additive inverse} of $s$) is the unique real number such that $s+(-s)=0$. \ul{Distance} (or \ul{separation}) on $\Real$ is given by the map $\Real\times\Real\ra\Real,~(r,s)\mapsto|r-s|$, where $|r-s|$ is called the \emph{distance between} $r$ and $s$. \ul{Division} on $\Real$ is the operation ~$\Real\times(\Real\backslash\{0\})\ra\Real,~(r,s)\mapsto {r\over s}:=r\cdot s^{-1}$, where $s^{-1}\in\Real$ (the \ul{reciprocal} or \ul{multiplicative inverse} of $s$) is the unique real number such that $s\cdot s^{-1}=1$. In a real number of the form ${r\over s}$ (called a \ul{fractional real number} or a \ul{real fraction}), we call $r$ the \ul{numerator} and $s$ the \ul{denominator} of ${r\over s}$.
\end{dfn}

Another construction of the real numbers uses special subsets of $\Rational$ (called Dedekind cuts) that generalize strict lower segments in $\Rational$ as follows.
\begin{dfn}[\textcolor{blue}{\index{Dedekind cut}{Dedekind cut}, \index{Real numbers $\Real$}{Real numbers $\Real$}: \cite{enderton1977}}]
Given $q\in \Rational$, consider the set of strict lower bounds $L_q:=\{p\in\Rational~|~p<q\}$. A set $c\subset\Rational$ is a \ul{Dedekind cut} if it satisfies the following:
\begin{enumerate}
\item $c$ is nonempty and proper, i.e., $\emptyset\neq c\subsetneq\Rational$.
\item $L_q\subset c$ for all $q\in c$, i.e., ~$\bigcup_{q\in c}L_q~\subset~c$.
\item $c$ has no greatest element, i.e., for every $q\in c$, there exists $q'\in c$ such that $q<q'$.
\end{enumerate}
(Note that for any rational $q\in\Rational$, the set $L_q$ is clearly a Dedekind cut). The set of \ul{real numbers} is
\bea
\label{RealDefEq2} \Real:=\{\txt{Dedekind cuts}~~c\subset\Rational\},~~~~~~~\txt{where}~~~\Rational\subset\Real~~~\txt{via}~~~\Rational\hookrightarrow\Real,~q\mapsto L_q,
\eea
as a linearly ordered set via a \ul{linear order} relation given by
\bea
c\leq c'~~~~\txt{if}~~~~c\subset c'~~~~\txt{(or $c\leq q'$ for some $q'\in c'$)},\nn
\eea
with \ul{addition} $+:\Real\times\Real\ra\Real,~(x,x')\mapsto x+x'$, ``\ul{zero}'' ($0_\Real$), \ul{multiplication} $\cdot:\Real\times\Real\ra\Real,~(x,x')\mapsto x\cdot x'$, and ``\ul{one}'' ($1_\Real$) respectively given by
{\small\begin{align}
& c+c':=\{q+q':~q\in c,~q'\in c'\},~~~~0_\Real:=L_{0_\Rational},\nn\\
&c\cdot c':=
\left\{
  \begin{array}{ll}
    0_\Real\cup\{qq':~0\leq q\in c,~0\leq q'\in c'\}, & \txt{if}~~c,c'\geq0\\
    |c|\cdot|c'|, & \txt{if}~~c,c'\leq0 \\
    -\big(|c|\cdot|c'|\big), & \txt{if one of $c,c'$ is $\geq0$ and the other $\leq0$}
  \end{array}
\right\},~~~~1_\Real:=L_{1_\Rational},\nn
\end{align}}
where for any real number $c\in\Real$, its \ul{absolute value} (\blue{footnote}\footnote{As before, the absolute value for real numbers is a map $|~|:\Real\ra\Real$.}) (in terms of its \ul{additive inverse}, i.e., the unique real number $-c$ satisfying $c+(-c)=0_\Real$, which is $-c=\big\{-q~|~q\in\Rational\backslash c\big\}$) is given by
\bea
|c|:=
\left\{
  \begin{array}{ll}
    c, & c\geq 0 \\
    -c, & c\leq 0
  \end{array}
\right\}=c\cup(-c).\nn
\eea
\end{dfn}

\begin{thm}[\textcolor{blue}{\index{lub (supremum) property of $\Real$}{lub (supremum) property of the ordering of $\Real$}: \cite[Theorem 2.2, p.13]{goldrei1996}}]
(With respect to our standard order on $\Real$) every nonempty subset $A\subset\Real$ with an upper bound has a \ul{least upper bound} (\ul{lub} or \ul{supremum}), denoted by $\txt{lub}(A)$ or $\sup(A)$.
\end{thm}
\begin{proof}
With the Dedekind cut representation of real numbers, let $A\subset\Real$ be any nonempty set with an upper bound $u\in\Real$, i.e., $A\leq u$, equivalently,
\bea
\textstyle\al:=\bigcup A:=\bigcup_{c\in A}c=\{q\in\Rational:q\in c~\txt{for some}~c\in A\}~\subset~u.\nn
\eea
It is clear that (i) $\al$ contains every element of $A$ and (ii) $\al$ is contained in every upper bound of $A$. Thus, it suffices to show $\al\in\Real$ (so that $\al$ is then a lub of $A$).
\begin{enumerate}
\item $\al\neq\emptyset$ by hypotheses, and $\al\subsetneq\Rational$ because $A$ is bounded.
\item If $q\in\al$, then $q\in c$ for some $c\in A$, and so $L_q\subset c\subset\al$.
\item If $q\in\al$, then $q\in c$ for some $c\in A$, and so $q<q'$ for some $q'\in c\subset\al$.
\end{enumerate}
It follows that $\al\in\Real$.
\end{proof}

\begin{crl}[\textcolor{blue}{\index{Archimedian property of $\Integer$ in $\Real$}{Archimedian property of $\Integer$ in $\Real$}, \index{Density of $\Rational$ in $\Real$}{Density of $\Rational$ in $\Real$}: \cite{rudin1976}}]~
Let $x,y\in\Real$.
\bit
\item[(a)]\ul{$\Integer$ is Archimedean in $\Real$}: If $x>0$, then there exists a positive integer $n>0$ such that
$nx>y$.
\item[(b)]\ul{$\Rational$ is dense in $\Real$}: If $x<y$, then there exists a rational number $q\in\Rational$ such that $x<q<y$.
\eit
\end{crl}
\begin{proof}
{\flushleft(a)} Let $A:=\{nx~|~n\in\Natural\}$. Suppose on the contrary that $nx\leq y$ for all $n\in\Natural$. Then $A$ is bounded above by $y$ (i.e., $A\leq y$), and thus by the lub property, $\al:=\sup A\leq y$ exists.

Since~ $x>0~~\Ra~~-x<0,~~\Ra~~\al-x<\al$, we see that $\al-x$ is not an upper bound of $A$, i.e., $A\not\leq \al-x$. This means there is some~ $m x\in A$ (where $m\in\Natural$) such that
\bea
\al-x<mx,~~~~\Ra~~~~\al<(m+1)x\in A,\nn
\eea
which is a contradiction, since $\al$ is an upper bound of $A$.
{\flushleft(b)} Let $x<y$. Since $y-x>0$, by part (a), there is a natural number $n\geq 1$ such that $n(y-x)>1$, i.e., such that $nx+1<ny$. Also, $nx\in\Real$ and so between $nx$ and $nx+1$ lies an integer $m$ such that $nx<m\leq nx+1$ (since consecutive integers are separated by a distance of $1$). Hence,
\[
nx<m\leq nx+1<ny,~~~~\Ra~~~~x<m/n<y. \qedhere
\]
\end{proof}

\begin{crl}[\textcolor{blue}{\index{glb (infimum) property of $\Real$}{glb (infimum) property of the ordering of $\Real$: \cite{rudin1976}}}]
(With respect to our order on $\Real$) every nonempty subset $A\subset\Real$ with a lower bound has a \ul{greatest lower bound} (\ul{glb} or \ul{infimum}), denoted by $\txt{glb}(A)$ or $\inf(A)$.
\end{crl}
\begin{proof}
$A\subset\Real$ has a lower bound $\iff$ $-A:=\{-a~|~a\in A\}$ has an upper bound. Hence, $\inf(A)$ exists, and is given by ~$\inf(A)=-\sup(-A)$.
\end{proof}

\begin{lmm}[\textcolor{blue}{
\index{Supremum criterion}{Supremum (lub) criterion for $\Real$},
\index{Infimum criterion}{Infimum (glb) criterion for $\Real$}}]\label{SupInfCrit}
Let $A\subset\Real$ be a set with an upper bound (resp. a lower bound), and $x\in \Real$. Then $x=\sup A$ iff for any $\vep>0$, there exists $a_\vep\in A$ such that $A\leq x<a_\vep+\vep$ (resp. $x=\inf A$ iff for any $\vep>0$, there exists $a_\vep\in A$ such that $a_\vep-\vep<x\leq A$).
\end{lmm}
\begin{proof}
It suffices to prove the supremum case only (because we can obtain the infimum case by replacing the pair $A,x$ with the pair $-A,-x$). ($\Ra$): Assume $x=\sup A$. Pick any $\vep>0$. Then $A\leq x$ and $A\not\leq x-\vep$ (otherwise $x$ is not the least upper bound). Thus, there exists $a=a_\vep\in A$ such that $x-\vep<a_\vep$. That is, $A\leq x<a_\vep+\vep$. ($\La$): Assume for all $\vep>0$, there exists $a_\vep\in A$ such that $A\leq x<a_\vep+\vep$. Suppose $x\neq\sup A$, i.e., $x>\sup A$. Let $\vep_0:=x-\sup A$. Then for any $a\in A$, $a+\vep_0=x+(a-\sup A)\leq x$. Hence
\[
A\leq x<a_{\vep_0}+\vep_0\leq x~~~~ \txt{(a contradiction)}. \qedhere
\]
\end{proof}

\begin{rmk}[\blue{\index{Selection map notation}{Selection map notation}}]
Let $\C$ be a collection of classes, $I$ a class, and $i\in I$ an object. Define a map ~$some_i:\C\ra\bigcup\C,~C\mapsto some_i(C)$~ that simply picks out one member from a class based on a selection rule that depends on $i$.

Let $A\subset\Real$ be a set with a lower bound (resp. an upper bound) and $x\in\Real$. Then $x=\inf A$ $\iff$ for any $\vep> 0$, $some_\vep(A-\vep)<x\leq A$ ~( resp. $x=\sup A$ $\iff$ for any $\vep>0$, $A\leq x<some_\vep(A+\vep)$ ),~ where $A-\vep:=\{a-\vep:a\in A\}$ and $A+\vep:=\{a+\vep:a\in A\}$.
\end{rmk}

\begin{dfn}[\blue{\index{Unconditional sum in $\Real$}{Unconditional sum in $\Real$}}]
Let $A$ be a set (or any class of sets), and consider maps $p:A\ra[0,\infty)$, $s:A\ra\Real$, $s_+:=s|_{s^{-1}[0,\infty)}$, and $s_-:=-s|_{s^{-1}(-\infty,0)}$. We define the \ul{unconditional sum} of $s$ (based on that of $p$) to be the unique possibly infinite number given by the following:
\begin{align}
&\textstyle \sum p:=\sum_{a\in A}p(a):=\sup\left\{\sum_{a\in F}p(a):\txt{finite}~F\subset A\right\}~\in~[0,\infty]:=[0,\infty)\cup\{|\Natural|\},\nn\\
&\textstyle\sum s:=\sum s_+-\sum s_-~\in~\ol{\Real}:=\Real\cup\{|\Natural|\},~~~~~~\txt{provided}~~~\sum s_+<\infty~~~\txt{or}~~~\sum s_-<\infty,\nn
\end{align}
where the inequality $\sum p\leq|\Natural|$ follows from the archimedean property of $\Natural$ (or of $\Integer$) in $\Real$.
\end{dfn}

\begin{dfn}[\textcolor{blue}{\index{Ring! of polynomials over $\Real$}{Ring of polynomials over $\Real$}, \index{Complex numbers $\Complex$}{Complex numbers $\Complex$}, \index{Modulus of a complex number}{Modulus}}]
The \ul{ring of polynomials} in one variable (\blue{footnote}\footnote{In accordance with more general polynomial rings we will encounter later whose elements are not necessarily functions in the current sense (i.e., in order to suppress the $\Real\ra\Real$ function character of its elements), this polynomial ring is often simply written as $\Real[x]:=\{f:=f_0+f_1x+\cdots+f_nx^n~|~f_i\in\Real,n\in\Natural\}$, with $x$ called an ``indeterminate/variable''. Note that the pointwise/function version of multiplication in (\ref{CompNoDefEq0}) is not valid for the general version $R[x]$ of the polynomial ring, while the scale-wise/coefficient-wise multiplication in (\ref{CompNoDefEq01}) remains valid. Also, we have chosen to work with $R[\ast]$, instead of $R[x]$, here only because we do not yet have a precise definition/meaning for $R[x]$, which will only come later after we have introduced basic abstract algebra operations.}) over $\Real$ is the set ~$\Real[\ast]:=\{f:\Real\ra\Real,~x\mapsto f_0+f_1x+\cdots+f_nx^n~|~f_i\in\Real,n\in\Natural\}\subset\Real^\Real$,~ along with \emph{addition} $+:\Real[\ast]\times \Real[\ast]\ra\Real[\ast],~(f,g)\mapsto f+g$ and \emph{multiplication} $\cdot:\Real[\ast]\times \Real[\ast]\ra\Real[\ast],~(f,g)\mapsto f\cdot g$ respectively given by
\bea
\label{CompNoDefEq0}(f+g)(x)=f(x)+g(x),~~~~(f\cdot g)(x):=f(x)g(x),~~~~\txt{for all}~~~~x\in\Real,
\eea
or more explicitly, with ~$(f\cdot g)_i:=\sum_{j+k=i}f_jg_k$,~ we define
\bea
\label{CompNoDefEq01}\textstyle\left(\sum f_ix^i\right)+\left(\sum g_ix^i\right):=\sum(f_i+g_i)x^i,~~~~\left(\sum f_ix^i\right)\cdot\left(\sum g_ix^i\right):=\sum (f\cdot g)_ix^i,
\eea
as well as ``\emph{zero}'' $0_{\Real[\ast]}$ and ``\emph{one}'' $1_{\Real[\ast]}$ respectively given by
\bea
0_{\Real[\ast]}:\Real\ra\Real,~x\mapsto 0_{\Real},~~~~1_{\Real[\ast]}:\Real\ra\Real,~x\mapsto 1_{\Real}.\nn
\eea
Define an equivalence relation $\sim$ on $\Real[\ast]$ as follows: For all $f,g\in\Real[\ast]$,
\bea
f\sim g~~\txt{if}~~\txt{there exists}~~h\in \Real[\ast]~~\txt{such that}~~ f(x)-g(x)=(1+x^2)h(x)~~\txt{for all}~~x\in\Real.\nn
\eea

The set of \ul{complex numbers} $\Complex$ is the set quotient set
\bea
\label{CompNoDefEq}\textstyle\Complex:={\Real[\ast]\over\sim}:=\big\{[f]~|~f\in \Real[\ast]\big\},~~~~[f]:=\{g\in\Real[\ast]: g\sim f\},
\eea
with \ul{addition} $\Complex\times\Complex\ra\Complex,~(z,z')\mapsto z+z'$, \ul{multiplication} $\Complex\times\Complex\ra\Complex,~(z,z')\mapsto z\cdot z'$, ``\ul{zero}'', and ``\ul{one}'' respectively given by
\bea
[f]+[g]:=[f+g],~~~~[f]\cdot[g]:=[f\cdot g],~~~~0_\Complex:=[x\mapsto 1+x^2],~~~~1_\Complex:=[1_{\Real[\ast]}].\nn
\eea
The complex number ~$i:=[id_\Real]=[x\mapsto x]$~ is called the \ul{imaginary unit}, and it satisfies
\bea
i^2=-1,~~~~\txt{where}~~~~1:=1_{\Complex},\nn
\eea
from which it follows that every complex number $[f]$ can be written as
\bea
[f]=f_0+if_0',~~~~\txt{for some}~~~~f_0,f_0'\in\Real.\nn
\eea
That is, $\Complex=\{a+ib~|~a,b\in\Real\}\approx\Real^2$ via $a+ib\mapsto(a,b)$, with addition and multiplication satisfying
\bea
(a+ib)+(a'+ib')=(a+a')+i(b+b'),~~~~(a+ib)\cdot(a'+ib')=(aa'-bb')+i(ab'+ba'),\nn
\eea
while ``zero'' and ``one'' are respectively given by ~$0_\Complex=0+i0=0$~ and ~$1_\Complex=1+i0$.

Consider any complex number $z=a+ib\in\Complex$. The \ul{conjugate} of $z$ is the complex number $\ol{z}:=a-ib$. The \ul{real part}, and the \ul{imaginary part}, of $z$ are respectively the real numbers
\bea
\textstyle\txt{Re}(z):=a={z+\ol{z}\over 2},~~~~\txt{and}~~~~\txt{Im}(z):=b={z-\ol{z}\over 2i}.\nn
\eea
The \ul{modulus} (as a map $|~|:\Complex\ra\Complex$) of $z$ is the nonnegative real number $|z|\geq 0$ whose square is given by
\bea
|z|^2:=z\ol{z}=a^2+b^2.\nn
\eea
The real numbers $\Real$ are contained in the complex numbers $\Complex$ via the injection ~$\Real\hookrightarrow\Complex$, $x\mapsto x+i0$.

\ul{Subtraction} on $\Complex$ is the operation ~$\Complex\times\Complex\ra\Complex,~(z,w)\mapsto z-w:=z+(-w)$, where $-w\in\Complex$ (the \ul{negative} or \ul{additive inverse} of $w$) is the unique complex number such that $w+(-w)=0$. \ul{Distance} (or \ul{separation}) on $\Complex$ is given by the map $\Complex\times\Complex\ra\Complex,~(z,w)\mapsto|z-w|$, where $|z-w|$ is called the \emph{distance between} $z$ and $w$. \ul{Division} on $\Complex$ is the operation ~$\Complex\times(\Complex\backslash\{0\})\ra\Complex,~(z,w)\mapsto {z\over w}:=z\cdot w^{-1}$, where $w^{-1}\in\Complex$ (the \ul{reciprocal} or \ul{multiplicative inverse} of $w$) is the unique complex number such that $w\cdot w^{-1}=1$. In a complex number of the form ${z\over w}$ (called a \ul{fractional complex number} or a \ul{complex fraction}), we call $z$ the \ul{numerator} and $w$ the \ul{denominator} of ${z\over w}$.
\end{dfn}

\section{Cardinal Numbers}\label{CardNumSec}
Recall that if $A,B$ are classes, then the cardinal ordering $|A|\leq|B|$ means there exists an injection $A\hookrightarrow B$. The classes we defined earlier were not necessarily classes of sets. In this section and the next, however, we are only interested in classes of sets. So, we will make the following definition.

\begin{dfn}[\blue{\index{Sets-cardinality}{Sets-cardinality}}]
Let $A$ be a set (or a class of sets if desired). The \ul{Sets-cardinality} of $A$ is the class of sets $|A|_{\txt{Sets}}:=\{B\in\txt{Sets}~|~B\approx A\}$. As usual, if $A$ is finite, i.e., $A\approx n$ for some $n\in\Natural$, then we simply write $|A|=n$ (i.e., $|A|=|n|$, where $|n|:=n$).
\end{dfn}
Henceforth, ``cardinality'' $|A|$ will mean ``Sets-cardinality'' $|A|_{\txt{Sets}}$, unless specified otherwise or the distinction is irrelevant.

\subsection{Cantor's theorem and the arithmetic of cardinals}
\begin{thm}[\textcolor{blue}{\index{Cantor's theorem}{Cantor's theorem}: \cite[Theorem 6.5, p.148]{goldrei1996}}]\label{CantorThm}
For any set $X$, we have $|X|<|\P(X)|$, i.e, the cardinality of $\P(X)$ is strictly greater than that of $X$. In particular, every element of the infinite set ~$\{\Natural,\P(\Natural),\P\P(\Natural),\P\P\P(\Natural),\cdots\}\approx\Natural$ ~has a distinct cardinality.
\end{thm}
\begin{proof}
We have $|X|\leq|\P(X)|$ due to the injective map $X\ra\P(X)$, $x\mapsto\{x\}$. Suppose $|X|=|\P(X)|$, i.e., there exists a bijection $f:X\ra\P(X)$. Define $A:=\{x\in X~|~x\not\in f(x)\}$. Since $f$ is surjective, there exists $x_0\in X$ such that $f(x_0)=A$. We have either (i) $x_0\in A$ or (ii) $x_0\not\in A$. If $x_0\in A=f(x_0)$, then $x_0\not\in f(x_0)=A$, a contradiction. If $x_0\not\in A=f(x_0)$, then $x_0\in A=f(x_0)$, a contraction.
\end{proof}

\begin{dfn}[\textcolor{blue}{
\index{Cardinal! numbers $\Cardinal$}{Cardinal numbers $\Cardinal$},
\index{Cardinal! arithmetic}{Cardinal arithmetic}}]
A subclass of sets $\kappa\subset Sets$ is a \ul{cardinal number}, written $\kappa\in\Cardinal$, if there exists a set $A$ with cardinality equal to $\kappa$, i.e., with $|A|=\kappa$.

Note that every $0\neq\kappa\in\Cardinal$ is a \ul{proper class of sets} (i.e., a class of sets that is not a set). We define \ul{cardinal arithmetic} operations (i.e., addition, multiplication, and exponentiation in $\Cardinal$ respectively) as follows (\textcolor{blue}{\cite[p.139]{enderton1977}}):
\bea
&&\kappa+\kappa':=|A\cup A'|,~~~~\txt{for any \ul{disjoint} sets $A,A'$ such that}~~|A|=\kappa,~|A'|=\kappa',\nn\\
&&\kappa\cdot\kappa':=|A\times A'|,~~~~\txt{for any sets $A,A'$ such that}~~|A|=\kappa,~|A'|=\kappa',\nn\\
&&\kappa^{\kappa'}:=|A^{A'}|,~~~~\txt{for any sets $A,A'$ such that}~~|A|=\kappa,~|A'|=\kappa'.\nn
\eea
\end{dfn}
Some basic properties (not of immediate relevance to us) of the arithmetic of cardinals are given in \textcolor{blue}{\cite[Theorem 6I, p.142]{enderton1977}}.

\begin{lmm}[\textcolor{blue}{The ordering in $\Cardinal$ is linear: \cite[Theorem 6M, p.151]{enderton1977}}]\label{CarOrdLin}
For any sets $A,B$ we have $|A|\leq |B|$ or $|B|\leq|A|$. (This requires the axiom of choice.)
\end{lmm}
\begin{proof}
Let the set ~$P:=\{\txt{bijections}~f_{UV}:U\ra V~|~U\subset A,~V\subset B\}$~ be ordered by inclusion, i.e.,
\bea
f_{UV}\leq f_{U'V'}~~~~\txt{if}~~~~f_{UV}:=Graph(f_{UV})\subseteq f_{U'V'}:=Graph(f_{U'V'}).\nn
\eea
If $\F=(f_{U_\ld V_\ld})_{\ld\in \Ld}$ is a nonempty chain in $P$, then with $U:=\bigcup_\ld U_\ld\subset A$ and $V:=\bigcup_\ld V_\ld\subset B$, the union
\[
\textstyle f_{UV}:=\bigcup\F=\bigcup_\ld f_{U_\ld V_\ld}:U\ra V~~\txt{is in $P$, because}\nn
\]
{\small\begin{align}
&\textstyle (u,v)\in f_{UV}~\iff~(u,v)\in f_{U_\ld V_\ld}~\txt{for some}~\ld\in\Ld,~~\txt{and}~~u\neq u'~\txt{in}~U~\iff~u\neq u'~\txt{in}~ U_\ld~\txt{for some}~\ld\nn\\
&\textstyle~~\iff~~f_{U_\ld V_\ld}(u)\neq f_{U_\ld V_\ld}(u')~\txt{in}~ V_\ld\subset V ~~\iff~~f_{UV}(u)\neq f_{UV}(u')~\txt{in}~U,\nn\\
&~~\Ra~~\txt{$f$ is bijective}.\nn
\end{align}}It follows from Zorn's lemma that $P$ has a maximal element $F_{CD}$. Suppose $C\neq A$ and $D\neq B$. Let $a\in A-C$, $b\in B-D$. Then we get a contradiction in the form
\bea
F_{CD}<G_{(C\cup\{a\})(D\cup\{b\})}:=F_{CD}\cup\{(a,b)\}\in P.\nn
\eea
Hence, $C=A$ (i.e., we have an injective map $A\ra B$) or $D=B$ (i.e., we have an injective map $B\ra A$).
\end{proof}

\begin{dfn}[\textcolor{blue}{\index{Characteristic! function}{Characteristic function of a subset}, \index{Subfinite partition of a set}{Subfinite partition of a set}}]
Let $A$ be a set, $B\subset A$ a subset, and $0\neq n\in\Natural$. The \ul{characteristic function} of $B$ with respect to $A$ is the map
\bea
\chi_B:A\ra 2=\{0,1\},~~a\mapsto\left\{
                                  \begin{array}{ll}
                                    1, &\txt{if}~~a\in B \\
                                    0, &\txt{if}~~a\in A-B
                                  \end{array}
                                \right\},~~~~~~~~\chi_\emptyset:=0:A\ra 2,~a\mapsto 0.\nn
\eea
An \ul{$n$-subfinite partition} of $A$ is a finite disjoint cover of $A$ of the form $\{A_0,A_1,...,A_{n-1}\}$ (hence containing a partition of $A$ of cardinality $k\leq n$) in the sense we have a disjoint union $A=A_0\sqcup A_1\sqcup\cdots\sqcup A_{n-1}$.
\end{dfn}

\begin{lmm}[\textcolor{blue}{Some basic cardinal equivalences}]\label{CantorLmm}
Let $A$ be a set, $0\neq n\in\Natural$, and
\[
\textstyle\P_n(A):=\{\txt{ordered $n$-subfinite partitions $(A_0,...,A_{n-1})$ of $A$}\}.~~\txt{In particular $\P_1(A)=\P(A)$}.
\]
Then the following are true:
\begin{enumerate}[leftmargin=0.9cm]
\item[(1)] $\P(A)=\P_1(A)\approx \P_2(A)\approx 2^A:=\big\{\txt{maps}~f:A\ra 2=\{0,1\}\big\}$.
\item[(2)] $\P_n(A)\approx n^A:=\big\{\txt{maps}~f:A\ra n=\{0,1,\cdots,n-1\}\big\}$.
\item[(3)] $\Real\approx n^{\Natural}\approx\P_n(\Natural)\approx\P(\Natural)\approx 2^\Natural$ ~~(for all ~$0\neq n\in\Natural$).
\item[(4)] $|\P(A)|=2^{|A|}$. (From \textcolor{blue}{\cite[p.141]{enderton1977}})
\end{enumerate}
\end{lmm}
\begin{proof}
\begin{enumerate}[leftmargin=0.7cm]
\item[(1)] We have the injective map $\chi:\P(A)\ra 2^A,~B\mapsto\chi_B$, with $\chi_B$ the characteristic function.

The map $\chi$ is surjective because every map $f:A\ra 2$ can be written as $f=\chi_{f^{-1}(1)}$.
\item[(2)] This is a straightforward extension of the proof of (1) above, i.e., the injective map
\[
\textstyle\chi:\P_n(A)\ra n^A,~(A_0,...,A_{n-1})\mapsto\sum_{i=0}^{n-1}i\chi_{A_i} 
\]
is also surjective because (with $\chi_\emptyset=0$) any map $f:A\ra n$ can be written as
\[
\textstyle f=\sum_{i=0}^{n-1}i\chi_{f^{-1}(i)}=\sum_{i=0}^{n-1}\chi_{\bigcup_{j=i+1}^nf^{-1}(j)},~~~~f^{-1}(n):=\emptyset.
\]

\item[(3)] Let $\Real_+:=(0,\infty):=\{x\in\Real:x>0\}$. Then $\Real\approx\Real_+$ via the map $h:\Real_+\ra\Real,~x\mapsto {1\over 2}(x-{1\over x})$, with inverse map $h^{-1}:\Real\ra\Real_+:x\mapsto x+(x^2+1)^{1/2}$. It thus suffices to prove $\Real_+\approx n^{\Natural}$.

($\Ra$): We can obtain an expression (called \index{Decimal expansion}{\ul{base-$n$ decimal expansion}}) for any positive real number as follows: Let $x > 0$ be any positive real number. Let $N_x\in\Natural$ be the largest number such that $N_x\leq x$, i.e., $N_x+1>x$ ($N_x$ exists due to the archimedean property of $\Natural$ in $\Real$, and the well-ordering of $\Natural$). Let $c_{-1}\in n$ be the largest number such that $N_x+{c_{-1}\over n}\leq x$, $c_{-2}\in n$ the largest number such that $N_x+{c_{-1}\over n}+{c_{-2}\over n^2}\leq x$, and so on (which each exist by induction). Having obtained the numbers $N_x,c_{-1},...,c_{-(k-1)}\in n$ this way, the next number $c_{-k}\in n$ is the largest number such that $N_x+c_{-1}/n+...+c_{-k}/n^k\leq x$.

Now let $k_x\in\Natural$ be the largest number such that $n^{k_x+1}>N_x$ (i.e., $n^{k_x}\leq N_x$). Next let $c_{k_x}\in n$ be the largest number such that $c_{k_x}n^{k_x}\leq N_x$, $c_{k_x-1}\in n$ the largest number such that $c_{k_x}n^{k_x}+c_{k_x-1}n^{k_x-1}\leq N_x$, and so on, until we obtain a relation $c_{k_x}n^{k_x}+c_{k_x-1}n^{k_x-1}+\cdots+c_1n^1+c_0n^0\leq N_x$. Observe that
{\small\bea
N_x<n^{k_x+1}=nn^{k_x}=(n-1)n^{k_x}+n^{k_x}=(n-1)n^{k_x}+(n-1)n^{k_x-1}+\cdots+(n-1)n+(n-1)+1.\nn
\eea}Moreover, if a nonzero natural number {\small $u\in n[n]:=\{a_{k_x}n^{k_x}+a_{k_x-1}n^{k_x-1}+\cdots+a_1n^1+a_0n^0~|~a_0,...,a_{k_x}\in n\}$}, then so does $u-1\in n[n]$. It thus follows that {\small $N_x=c_{k_x}n^{k_x}+c_{k_x-1}n^{k_x-1}+\cdots+c_1n+c_0\in n[n]$}.

Let $E:=\{x_k:k\in\Natural\}$, where the ``$x$-approximating'' numbers $x_k\leq x$ are given by
{\small\begin{align}
&\label{decimal1}x_k:=N_x+c_{-1}n^{-1}+\cdots+c_{-k}n^{-k}~~~~~(k=0,1,2,...)\nn\\
&~~~~=c_{k_x}n^{k_x}+c_{k_x-1}n^{k_x-1}+\cdots+c_1n+c_0+c_{-1}n^{-1}+\cdots+c_{-k}n^{-k}~~~~(k=0,1,2,...).
\end{align}}Since $E\leq x$, the supremum $\al:=\sup E=c_{k_x}n^{k_x}+c_{k_x-1}n^{k_x-1}+\cdots+c_1n+c_0+c_{-1}n^{-1}+c_{-2}n^{-2}+\cdots$ exists. If $\al<x$, then by the construction of $E$, we can choose $x_k$ such that $\al<x_k\leq x$ (a contraction). Hence, $x=\sup E=c_{k_x}n^{k_x}+c_{k_x-1}n^{k_x-1}+\cdots+c_1n+c_0+c_{-1}n^{-1}+c_{-2}n^{-2}+\cdots$.

The following associated unique sequence of coefficients is the \ul{base-$n$ decimal expansion} of $x$.
\bea
\label{decimal2}\wt{x}~:=~...~000c_{k_x}c_{k_x-1}...c_2c_1c_0\cdot c_{-1}c_{-2}~...~\in~n^\Natural.
\eea
This defines an injective map ~$\Real_+\hookrightarrow n^\Natural$, $x\mapsto\wt{x}$.

($\La$): Conversely, any map $f:\Natural\ra n$ is a sequence $\wt{x}_f$ of the form (\ref{decimal2}), for which the associated set $E_f:=\{x^f_k:k\in\Natural\}$ of real numbers given by (\ref{decimal1}) is bounded above by a sum of the form ~$C+n\sum_{k=0}^\infty n^{-k}$ (which converges as a \index{Geometric series}{\ul{geometric series}} (\blue{footnote}\footnote{Given a complex number $z\in\Complex$, a sum of the form $f(z):=\sum_{k\in\Natural}z^k:=\lim_{n\ra\infty}\sum_{k=0}^nz^k$ is called a \ul{geometric series}, where the limit is said to \ul{exist} or \ul{make sense} (meaning it \ul{is finite}) iff the following is true: With $f_n(z):=\sum_{k=0}^nz^k$, for any $\vep>0$, there exists $N_\vep$ such that $|f(z)-f_n(z)|<\vep$ for all $n\geq N_\vep$.})) for a constant $C\geq 0$, and so (\ref{decimal2}) is the base-$n$ decimal expansion of the unique real number $x_f:=\sup E_f$. This defines an injective map $n^\Natural\hookrightarrow\Real_+,~f\mapsto x_f$.

Hence, by the Cantor-Schr\"{o}der-Bernstein theorem, we have a bijection $\Real_+\ra n^\Natural$.

\item[(4)] Since $\P(A)\approx 2^A$, we have $|\P(A)|=|2^A|\sr{(s)}{=}2^{|A|}$, where step (s) is the exponentiation of cardinals. \qedhere
\end{enumerate}
\end{proof}

\subsection{The arithmetic of infinite cardinals}
\begin{lmm}[\blue{\cite[Lemma 6R, p.162]{enderton1977}}]
If $\kappa\geq|\Natural|$ is an infinite cardinal number, then ~$\kappa\cdot\kappa=\kappa$.
\end{lmm}
\begin{proof}
First, observe that we have the disjoint union $\Natural\times\Natural=\bigcup_{k\in\Natural}F_k$, where {\footnotesize $F_k:=\{(i,j):i+j=k\}=\{(i,k-i):i=0,1,...,k\}$} has cardinality $|F_k|=k+1$. So, we get an injective map {\footnotesize $f:\Natural\times\Natural\ra\Natural$} given by
{\small\[
\textstyle f|_{F_k}:F_k\ra \wt{F}_k:=\left\{i+{k(k+1)\over 2}:i=0,1,...,k\right\},~(i,i-k)\mapsto i+{k(k+1)\over 2},~~~~\txt{for each}~~k\in\Natural.
\]}That is, $\Natural\times\Natural\approx\Natural$, or $|\Natural|\cdot|\Natural|=|\Natural|$. (\blue{footnote}\footnote{If we consider an infinite set such as $\ol{\Real}_+:=[0,\infty)=\{r\in \Real:r\geq 0\}$, we also have a disjoint union $\ol{\Real}_+\times\ol{\Real}_+=\bigcup_{r\in\ol{\Real}_+}D_r$, where $D_r:=\{(a,b):a+b=r\}=\{(tr,(1-t)r):0\leq t\leq 1\}\approx[0,1]:=\{t\in\Real:0\leq t\leq 1\}$, but $|D_r|$ is not finite (although "compact" when we put a ``topology'' on $\Real$ later).})

Let $B$ be any set such that $|B|=\kappa$. As a poset under inclusion $\subseteq$, consider the set
\bea
P:=\{\txt{bijections}~h:C^2\ra C,~\txt{for infinite subsets}~C\subset B\}~\subset~ M(B),\nn
\eea
where ~$M(B):=\{\txt{maps}~h\subset B^2\times B\}=\{\txt{maps}~h:U\ra V~|~U\in \P(B^2),~V\in\P(B)\}$. $P$ is nonempty because every infinite subset of $B$ contains a set of cardinality $|\Natural|$. Let $\C=\{h_\ld:C_\ld^2\ra C_\ld\}_{\ld\in\Ld}$ be a nonempty chain in $P$. Then as shown in the proof if Lemma \ref{CarOrdLin}, $g:=\bigcup\C=\bigcup_\ld h_\ld:\bigcup_\ld C_\ld^2\ra\bigcup_\ld C_\ld$ is \ul{bijective}. Let ~$D:=\im g=\bigcup_\ld \im h_\ld=\bigcup_\ld C_\ld$, which is infinite since $\C$ is nonempty. It is clear that $\bigcup_\ld C_\ld^2\subset D^2$. Also, if $(a,a')\in D^2$, let $a\in\im h_{\ld_a}=C_{\ld_a}$ and $a'\in\im h_{\ld_{a'}}=C_{\ld_{a'}}$. Since $\C$ is a chain, wlog, assume $h_{\ld_{a'}}\subset h_{\ld_a}$, i.e., $a,a'\in \im h_{\ld_a}=C_{\ld_a}$. Then $(a,a')\in C_{\ld_a}^2\subset\bigcup_\ld C_\ld^2$, and so $D^2\subset \bigcup_\ld C_\ld^2$. Hence $D^2=\bigcup_\ld C_\ld^2$, and so $g$ is a bijection $D^2\ra D$, i.e., $g\in P$.

By Zorn's lemma, $P$ has a maximal element $f:A\times A\ra A$. Let $\ld:=|A|\geq|\Natural|$. Then by construction, $\ld\cdot\ld=\ld$. Moreover, $|B-A|<\ld$, otherwise, if $\ld=|A|\leq|B-A|$, then $B-A$ contains a subset $A'\subset B-A$ of cardinality $|A'|=|A|=\ld$, and so
\bea
&&|(A\sqcup A')\times(A\sqcup A')|=|(A\times A)\sqcup(A\times A')\sqcup(A'\times A)\sqcup(A'\times A')|\nn\\
&&~~~~=\ld\cdot\ld+\ld\cdot\ld+\ld\cdot\ld+\ld\cdot\ld=\ld+\ld+\ld=3\cdot\ld=2\cdot\ld=\ld+\ld=|A\sqcup A'|,\nn\\
&&~~\Ra~~(A\sqcup A')^2\sr{F}{\approx} A\sqcup A'~~~~\txt{(a contradiction of the maximality of $f\subsetneq F\in P$)}.\nn
\eea
Hence, $|B|=|A\sqcup(B-A)|=|A|+|B-A|=\ld+|B-A|\leq\ld+\ld=2\cdot\ld\leq\ld\cdot\ld=\ld$, i.e., $\kappa=\ld$.
\end{proof}

\begin{crl}[\blue{Absorption law of cardinal arithmetic: \cite[p.164]{enderton1977}}]
Let $\kappa,\kappa'\in\Cardinal$. If $\min\{\kappa,\kappa'\}\neq 0$ and $\max\{\kappa,\kappa'\}=\infty$, then ~$\kappa+\kappa'=\kappa\cdot\kappa'=\max(\kappa,\kappa')$.
\end{crl}

\begin{rmk}[\blue{The arithmetic of infinities (i.e., infinite cardinals)}]
Cantor's theorem shows that we have many different$\slash$inequivalent sizes of infinity i.e., cardinalities $\infty_A:=|A|\geq|\Natural|$ of infinite sets $A$. Let $A$ and $B$ be two infinite sets, with cardinalities/infinities $\infty_A:=|A|$ and $\infty_B:=|B|$. If $\infty_A\leq\infty_B$, then because of cardinal equivalences (due to the preceding result)
\bea
A\approx A\sqcup A\approx A\times A,~~~~B\approx A\sqcup B\approx B\sqcup B\approx A\times B\approx B\times B,\nn
\eea
we have the following ``simplified'' basic rules for adding and multiplying infinities:
\bea
\infty_A\cdot\infty_A=\infty_A+\infty_A=\infty_A\leq\infty_B=\infty_A+\infty_B=\infty_B+\infty_B=\infty_A\cdot\infty_B=\infty_B\cdot\infty_B\nn
\eea
\end{rmk}

\begin{crl}\label{FinSubCaCrl}
Let $A$ be a set. If $|A|=\infty$ (i.e., $|A|\geq|\Natural|$), then the set of finite subsets $FS(A):=\{B\subset A:|B|<\infty\}$ of $A$ has cardinality ~$|FS(A)|=|A|$.
\end{crl}
\begin{proof}
Observe that $FS(A)=\bigcup_{n=0}^\infty A(n)$, where $A(n):=\{B\subset A:|B|\leq n\}$ has cardinality $|A(n)|\leq|A^n|=|A|^n$. Therefore, $|A|\leq|FS(A)|=\left|\bigcup_{n=0}^\infty A(n)\right|=\sup\{ |A(n)|:n\in\Natural\} \leq \sup\{|A|^n:n\in\Natural\}=\sup\{|A|:n\in\Natural\}=|A|$.
\end{proof}

\subsection{Countability of sets (a concretizing digression)}~\\~
This section is mostly repeating (in a more concrete way) things we have already encountered.
\begin{dfn}[\textcolor{blue}{\index{Finite! set}{Finite set}, \index{Infinite set}{Infinite set}, \index{Countable! set}{Countable set}, \index{At most countable set}{At most countable set}, \index{Uncountable set}{Uncountable set}}]
Let $A$ be a set. We say $A$ is
\begin{enumerate}
\item \ul{finite} if $A\approx n$ for some $n\in\Natural$ (i.e., $|A|<|\Natural|$).
\item \ul{infinite} if $A$ is not finite (i.e., $|A|\not<|\Natural|$).
\item \ul{countable} (or \ul{countably-infinite}) if $A\approx\Natural$ (i.e., $|A|=|\Natural|$).
\item \ul{at most countable} if $A$ is finite or countable (i.e., $|A|\leq|\Natural|$).
\item \ul{uncountable} (or \ul{uncountably-infinite}) if $A$ is neither finite nor countable (i.e., $|A|\not\leq|\Natural|$).
\end{enumerate}
Countable sets are sometimes called enumerable, or denumerable, sets.
\end{dfn}

\begin{dfn}[\textcolor{blue}{\index{Sequence in a set}{Sequence in a set}, \index{Terms of a sequence}{Terms of a sequence}}]
A \ul{sequence} in a set $A$ (or a sequence of elements of $A$) is a map $f:\Natural\ra A$, i.e., an element of the set $A^\Natural$. With $a_n:=f_n:=f(n)$, for $n\in \Natural$, we also denote the sequence $f$ by $(a_n)\subset A$, or $\{a_n\}\subset A$, or $a_0,a_1,\cdots\in A$, or just $a_n\in A$. The values of $f$, i.e., the elements $a_n$, are called the \ul{terms} of the sequence.
\end{dfn}

Since every countable set $C$ is the image of an injective map $\Natural\ra C$, we
may regard $C$ as the image of a sequence with \ul{distinct terms}. That is, the elements of any countable set $C$ can be \ul{enumerated} as ~$C=\{c_0,c_1,\cdots\}$, ~where ~$i\neq j$ ~$\Ra$~ $c_i\neq c_j$.

\begin{thm}[\textcolor{blue}{Countable subsets: \cite{rudin1976}}]\label{infinite-subcount}
An infinite subset of a countable set is countable.
\end{thm}
\begin{proof}
Let $C$ be countable and $E\subset C$ infinite. Consider an enumeration $C=\{c_n\}=\{c_0,c_1,...\}$. Then $E$ can also be enumerated as ~$E=\{e_k\}:=\{c_{n_k}\}=\{c_{n_0},c_{n_1},\cdots\}$,~ where
\bea
n_0:=\min\{n\geq0:~c_n\in E\},~~~~n_k:=\min\{n>n_{k-1}:~c_n\in E\},~~\txt{for $k\geq 1$}.\nn
\eea
The map $f:\Natural\ra E,~i\mapsto c_{n_i}$ is bijective by construction, and so $E$ is countable.
\end{proof}

\begin{thm}[\textcolor{blue}{Countable unions}]
If $\{C_i\}_{i\in\Natural}$ are countable sets, their union $C=\bigcup_{n=0}^\infty C_n$ is countable.
\end{thm}
\begin{proof}
For each $n\in \Natural$, consider an enumeration {\small $C_n=\{c_{nk}\}_{k\in \Natural}$}. Then the union $C$ can be written as
\bea
C=\{c_{ij}\}_{i,j\in \Natural}.\nn
\eea
If the sets $C_n$ are disjoint, then a bijection between $C$ and $\Natural$ is given by a map $f:\Natural\ra\{c_{ij}\}$ of the form
{\footnotesize\bea\bt[row sep=tiny,column sep=tiny]
f(1)\ar[r,draw=none] &~& f(2) \ar[r,draw=none] & f(3)\ar[r,draw=none] &~& f(4) \ar[r,draw=none] & f(5) \ar[r,draw=none] & f(6)\ar[r,draw=none] &~& f(7) \ar[r,draw=none] & f(8) \ar[r,draw=none] & f(9) \ar[r,draw=none] & f(10)\ar[r,draw=none] &~& f(11) \ar[r,draw=none] & ...\ar[r,draw=none] &~& \cdots\\
c_{11}\ar[r,draw=none] &~& c_{21} \ar[r,draw=none] & c_{12}\ar[r,draw=none] &~& c_{31} \ar[r,draw=none] & c_{22} \ar[r,draw=none] & c_{13}\ar[r,draw=none] &~& c_{41} \ar[r,draw=none] & c_{32} \ar[r,draw=none] & c_{23} \ar[r,draw=none] & c_{14}\ar[r,draw=none] &~& c_{51}\ar[r,draw=none] & ...\ar[r,draw=none] &~& \cdots
\et\nn
\eea}where the grouping of terms is based on the sets of pairs $I_2:=\{(i,j)\}_{i+j=2}$, $I_3:=\{(i,j)\}_{i+j=3}$, $I_4:=\{(i,j)\}_{i+j=4}$, $I_5:=\{(i,j)\}_{i+j=5}$, $I_6:=\{(i,j)\}_{i+j=6}$, $\cdots$

If some of the sets $C_n$ overlap, then (by Theorem \ref{infinite-subcount}) $C$ will be countable as an infinite subset of the countable set $\{c'_{ij}\}$ in which the repeating elements in $C=\{c_{ij}\}$ are replaced by distinct ones.
\end{proof}

\begin{crl}[\textcolor{blue}{At most countable unions}]
Let $I$ be a set and $\{A_i~|~i\in I\}$ a collection of sets. If both $I$ and the $A_i$'s are at most countable, then the union ~$A:=\bigcup_{i\in I}A_i$~ is at most countable.
\end{crl}

\begin{thm}[\textcolor{blue}{Countable products}]
If $A$ is a countable set, then $A^n$ (for $0\neq n\in\Natural$) is countable.
\end{thm}
\begin{proof}
The claim holds for $n=1$. If the claim holds for $n-1$, then it also holds for $n$ since
\bea
\textstyle A^n=A\times A^{n-1}=\left\{(a,a'):a\in A,~a'\in A^{n-1}\right\}=\bigcup_{a\in A}\left\{(a,a'):a'\in A^{n-1}\right\}
=\bigcup_{a\in A}~\{a\}\times A^{n-1}.\nn
\eea
Hence, by induction, the claim holds for all $n$.
\end{proof}
\begin{crl}[\textcolor{blue}{\index{Countability of $\Rational$}{Countability of $\Rational$}}]
The set of rationals $\Rational$ is countable.
\end{crl}
\begin{proof}
$\Integer\approx\Natural$ via the map $f:\Integer\ra\Natural$ given by $f|_{\Natural}:i\mapsto 2i$ and $f|_{\Integer\backslash\Natural}:i\mapsto -2i-1=2|i|-1$. This implies $\Integer\approx\Natural\times\Natural\approx\Integer\times\Integer$. Hence $\Rational:={\Integer\times\Integer\over\sim}=\left\{{a\over b}:a,b\in\Integer,~b\neq 0\right\}$ is countable as an infinite set whose cardinality satisfies $|\Integer|\leq|\Rational|\leq|\Integer\times\Integer|$.
\end{proof}

We have already seen, from Lemma \ref{CantorLmm} and Theorem \ref{CantorThm}, that $\Real\approx 2^\Natural\approx\P(\Natural)>\Natural$ (where the strict inequality is due to Cantor's theorem, i.e., Theorem \ref{CantorThm}), and so we know $\Real$ is uncountable. The following result gives an alternative argument that does not rely on Cantor's theorem.

\begin{thm}[\textcolor{blue}{\cite{rudin1976}}]
The set of $\{0,1\}$-sequences ~$2^\Natural=\big\{f:\Natural\ra 2=\{0,1\}\big\}$~ is uncountable.
\end{thm}
\begin{proof}
Let $E\subset 2^\Natural$ be countable. Then $E=\{s(1),s(2),...\}$, where for each $i$, $s(i)=\big\{s(i)_k\big\}_{k\in \Natural}$ is a sequence of $0$'s and $1$'s. Consider the sequence $s=\{s_k\}\in 2^\Natural$ given by
~$s_k:=\left\{
      \begin{array}{ll}
        0, & \hbox{if}~s(k)_k=1, \\
        1, & \hbox{if}~s(k)_k=0.
      \end{array}
    \right\}$.~ Then $s\neq s(i)$ for all $i$, i.e., $s$ differs from any given sequence in $E$ (by at least one term), and so $s\not\in E$. Therefore $E$ is a proper subset of $2^\Natural$. That is, every countable subset of $2^\Natural$ is a proper subset of $2^\Natural$. It follows that $2^\Natural$ is uncountable (otherwise, $2^\Natural$ would be a proper subset of itself, a contradiction).
\end{proof}

\begin{crl}[\textcolor{blue}{\index{Uncountability of $\Real$}{Uncountability of $\Real$}}]
The set of real numbers $\Real\approx 2^\Natural$ is uncountable.
\end{crl}

\begin{caution}[\textcolor{blue}{A practical convention for countability}]
We have named a set $A$ countable iff $A$ has the same cardinality as $\Natural$. However, for brevity in later chapters, ``\emph{countable}'' will mean ``\emph{at most countable}'' (i.e., ``\emph{finite or countable}''). Even when this is not explicitly mentioned, the appropriate convention in use will be obvious from the context.
\end{caution}

\section{Ordinal Numbers}\label{OrdNumSec}
\subsection{Linearly-ordered sets: Transfinite induction, basic properties, and arithmetic}
\begin{dfn}[\blue{
\index{Lower segment}{Lower segment},
\index{Minimum element}{Minimum element},
\index{Maximum element}{Maximum element},
\index{Successor}{Successor},
\index{Predecessor}{Predecessor},
\index{Limit! element}{Limit element},
\index{Dense! order}{Dense order}}]
Let $(A,\leq)$ be a linearly-ordered set, $a\in A$, and $\emptyset\neq B\subset A$. The \ul{(strict) lower segment} of $A$ at $a$ is the set ~$[A,a):=L_a:=L_A(a):=Seg_A(a):=\{a'\in A:a'<a\}$.

(If it exists) the \ul{minimum element} $\min B$ of $B$ is the least element of $B$, i.e., $\min B\in B$ and $\min B\leq b$ for all $b\in B$. (If it exists) the \ul{maximum element} $\max B$ of $B$ is the greatest element of $B$, i.e., $\max B\in B$ and $\max B\geq b$ for all $b\in B$.

(If it exists) the \ul{predecessor} of the element $a\in A$ is $Pred_A(a):=\max\{a'\in A:a'<a\}$. (If it exists) the \ul{successor} of the element $a\in A$ is $Succ_A(a):=\min\{a'\in A:a<a'\}$. (\blue{footnote}\footnote{If $(A,\leq)$ is well-ordered, then it is clear that for every $a\in A$, its successor $Succ_A(a)$ exists.}).

The element $a\in A$ is a \ul{(lower-) limit element} of $A$ if it is ``lower-dense'' in the sense that (i) $a\neq\min A$ and (ii) $Pred_A(a)$ does not exist. Equivalently, $a\in A$ is a limit element of $A$ if (i) $[A,a)\neq\emptyset$, i.e., some $a'<a$, and (ii) for any $a'<a$ we have $a'<c<a$ for some $c\in A$. (\blue{footnote}\footnote{In particular, if $(A,\leq)$ is well-ordered, then for any $\emptyset\neq B\subset A$, $\min B$ (even if a limit element of $A$) is not a limit element of $B=(B,\leq)\subset(A,\leq)$.}).

The order $\leq$ is a \ul{dense order} if for any $a\in A$, both $Pred_A(a)$ and $Succ_A(a)$ do not exist (i.e., for all $a,a'\in A$, if $a<a'$ then $a<c<a'$ for some $c\in A$).
\end{dfn}

\begin{dfn}[\textcolor{blue}{\index{Inductive! subset of a linearly-ordered set}{Inductive subset of a linearly-ordered set}}]
Let $A$ be a set and $(A,\leq)$ a linear-ordering of $A$. If $a\in A$, let ~$[A,a):=L_a:=L_A(a):=Seg_A(a):=\{a'\in A~|~a'<a\}$. ~A nonempty subset $\emptyset\neq B\subset A$ is \ul{$(\leqslant)$-inductive} (i.e., \ul{inductive with respect to $\leqslant$}) if for all $a\in A$, ~$[A,a)\subset B$ ~$\Ra$~ $a\in B$. (\blue{footnote}\footnote{If $(A,\leq)$ is well-ordered, then every inductive subset $B\subset A$ contains $m:=\min A$, since $[A,m)=\emptyset\subset B$.})
\end{dfn}

The direction ($\Ra$) in the following result is a generalization of the \emph{strong principle of induction} for $\Natural$.
\begin{thm}[\textcolor{blue}{\index{Transfinite induction principle for! well-ordered sets}{Transfinite induction principle for well-ordered sets}: \cite[p.174]{enderton1977}}]\label{TFIPforLOS}
Let $(A,\leq)$ be a linearly-ordered set. Then $(A,\leq)$ is well-ordered $\iff$ every $(\leq)$-inductive subset $B\subset A$ equals $A$.
\end{thm}
\begin{proof}
{\flushleft($\Ra$):} Assume $(A,\leq)$ is well-ordered. Suppose we have a $(\leq)$-inductive subset $B\subsetneq A$. Let $m:=\min(A-B)$. Then $[A,m)\subset B$. But this implies $m\in B$ (by inductivity of $B$), which is a contradiction.
{\flushleft($\La$):} Assume every $(\leq)$-inductive subset of $A$ equals $A$. Let $\emptyset\neq C\subset A$. Let $B:=[A,C):=\{a\in A~|~a<C\}$ be the set of strict lower bounds of $C$ in $A$ (where it is clear that $B\cap C=\emptyset$). Then, either $B$ is $(\leq)$-inductive or $B$ is not $(\leq)$-inductive. If $B$ is $(\leq)$-inductive, then $B=A$, and so $C=\emptyset$ (a contradiction). Therefore $B$ is not $(\leq)$-inductive. That is, there exists $a\in A$ such that $[A,a)\subset B=[A,C)$ but $a\not\in B$, i.e., $a\geq c$ for some $c\in C$. Since $[A,a)\subset B$ and $B\cap C=\emptyset$, it follows that $c\not<a$ (i.e., $c\not<a\geq c$), and so $c=a$ (by the linear-ordering of $A$). Since $[A,c)\cap C=[A,a)\cap C\subset B\cap C=\emptyset$, we see that $c=\min C$.
\end{proof}

\begin{dfn}[\blue{
\index{Order-invariant map}{Order-invariant map},
\index{Order-symmetric map}{Order-symmetric map},
\index{Order-imbedding map}{Order-imbedding map},
\index{Order! isomorphism}{Order isomorphism},
\index{Order-isomorphic sets}{Order-isomorphic sets}}]
Let $(A,\leq_A)$ and $(B,\leq_B)$ be partially ordered sets. A map $f:A\ra B$ is \ul{order-invariant} (or \ul{order-preserving} -- \blue{footnote}\footnote{\magenta{Caution}: In the literature, ``order-preserving'' might stand for something stronger (e.g., in \cite{goldrei1996}). In our definition, a more direct interpretation of ``order-preserving'' involving the map alone (without any thought to its invertibility) is preferred.}) if for all $a,a'\in A$, $a<a'$ $\Ra$ $f(a)<f(a')$. (\blue{footnote}\footnote{If the set $A$ is linearly ordered, then an order-invariant map $f:A\ra B$ is automatically injective.}). A map $f:A\ra B$ is \ul{order-symmetric} if for all $a,a'\in A$, $a<a'$ $\iff$ $f(a)<f(a')$. A map $f:A\ra B$ is \ul{order-imbedding} if it is (i) injective and (ii) order-symmetric. An \ul{order isomorphism} is a surjective order-imbedding map (i.e., an order-invariant bijection whose inverse map is also order-invariant). If an order isomorphism $f:A\ra B$ exists, we say $A$ and $B$ are \ul{order-isomorphic} sets, written $A{\cong_o}B$ (or just $A\cong B$ for brevity).
\end{dfn}

\begin{prp}[\blue{Sectional representation of a partially ordered set}]
Let $(X,\leq)$ be a partially ordered set. Let $Y:=\{\ol{L}_x:x\in X\}=(Y,\subseteq)$ be the set of lower sections $[X,x]:=\ol{L}_x:=\ol{L}_X(x):=\{x'\in X:x'\leq x\}$ of $X$ ordered by inclusion. Then the map $f:X\ra Y,~x\mapsto \ol{L}_x$ is an order-isomorphism, i.e., $X\cong_oY$.
\end{prp}
\begin{proof}
Observe that the surjective map $f$ is injective (if $\ol{L}_{x}=\ol{L}_{x'}$, then $\ol{L}_{x}\subset \ol{L}_{x'}$ and $\ol{L}_{x'}\subset \ol{L}_{x}$, which implies $x\leq x'$ and $x'\leq x$, and so $x=x'$). Also, for all $x,x'\in X$, we have ~$x<x'$ $\iff$ $\ol{L}_x\subsetneq \ol{L}_{x'}$ (because $x<x'$ implies both $\ol{L}_x\subset \ol{L}_{x'}$ and $x'\in \ol{L}_{x'}\backslash \ol{L}_x$, while $\ol{L}_x\subsetneq \ol{L}_{x'}$ implies both $x\leq x'$ and $x\neq x'$).
\end{proof}

\begin{prp}[\blue{Strict segmental representation of a linearly ordered set}]
Let $(X,\leq)$ be a linearly ordered set. Let $Y:=\{L_x:x\in X\}=(Y,\subseteq)$ be the set of strict lower segments $[X,x):=L_x:=L_X(x):=Seg_X(x):=\{x'\in X:x'<x\}$ of $X$ ordered by inclusion. Then the map $f:X\ra Y,~x\mapsto L_x$ is an order-isomorphism, i.e., $X\cong_oY$.
\end{prp}
\begin{proof}
Observe that the surjective map $f$ is injective (if $x\neq x'$ then $x<x'$ or $x'<x$, which implies $L_x\subsetneq L_{x'}$ or $L_{x'}\subsetneq L_x$, and so $L_x\neq L_{x'}$). Also, for all $x,x'\in X$, we have ~$x<x'$ $\iff$ $L_x\subsetneq L_{x'}$ (because $x<x'$ implies both $L_x\subset L_{x'}$ and $x'\in L_{x'}\backslash L_x$, while $L_x\subsetneq L_{x'}$ implies both $x'\not <x$ and $x\neq x'$).
\end{proof}

\begin{thm}[\blue{\cite[Theorem 7.7, p.195]{goldrei1996}}]\label{OrderInvThm}
If $(X,\leq)$ is a well-ordered set, then every order-invariant map $f:X\ra X$ satisfies $x\leq f(x)$ for all $x\in X$.
\end{thm}
\begin{proof}
Let $A:=\{x\in X:x\leq f(x)\}$. We need to show $A=X$. It suffices (by Theorem \ref{TFIPforLOS}) to show $A\subset X$ is inductive. Let $x_0:=\min X$. Then $x_0\leq f(x_0)$, and so $x_0\in A$. Fix $x\in X$ such that $[X,x):=L_x:=\{x'\in X:x'<x\}\subset A$, i.e., $x'<x$ $\Ra$ $x'\leq f(x')$. Since $f$ is order-invariant, $x'<x$ $\Ra$ $x'\leq f(x')<f(x)$, i.e., $x'<x$ $\Ra$ $x'<f(x)$. Thus, if $f(x)<x$, i.e., $f(x)\in[X,x)$, then $f(x)<f(x)$, a contradiction. Therefore $x\leq f(x)$, i.e., $x\in A$. This shows $A$ is inductive, and since $X$ is well-ordered, it follows that $A=X$.
\end{proof}
\begin{crl}[\blue{\cite[Theorem 7.8, p.196]{goldrei1996}}]
If $(X,\leq)$ is a well-ordered set, then every order isomorphism $f:X\ra X$ is the identity map $id_X$.
\end{crl}
\begin{proof}
The preceding result for both $f$ and $f^{-1}$ gives {\small $x\leq f(x)\leq f^{-1}f(x)=x$} for all $x\in X$.
\end{proof}

\begin{crl}[\blue{\cite[Exercise 7.58, p.197]{goldrei1996}}]
If $(X,\leq_X)$ and $(Y,\leq_Y)$ are well-ordered sets, then any order isomorphisms $f,g:X\ra Y$ are equal.
\end{crl}
\begin{proof}
By the preceding result, $g^{-1}\circ f:X\sr{f}{\ral}Y\sr{g^{-1}}{\ral}X$ is the identity, and so $f=g$.
\end{proof}

\begin{crl}[\blue{\cite[Exercise 7.59, p.197]{goldrei1996}}]
Let $(X,\leq)$ be a well-ordered set, $x\in X$, and $L_x:=[X,x):=\{x'\in X:x'<x\}$ the strict lower segment at $x$. Then there is no order-invariant map $f:X\ra L_x$. In particular, $X$ is not order-isomorphic to $L_x$.
\end{crl}
\begin{proof}
Suppose $f:X\ra L_x=[X,x)\subset X$ is an order-invariant map. Then, by Theorem \ref{OrderInvThm}, $x'\leq f(x')$ for all $x'\in X$. In particular, $x\leq f(x)<x$ (a contradiction).
\end{proof}

\begin{crl}[\blue{\cite[Exercise 7.60, p.197]{goldrei1996}}]
Let $(Z,\leq)$ be a well-ordered set, $z,z'\in Z$, and $L_z$,$L_{z'}$ the respective lower segments at $z,z'$. If $L_z\cong_o L_{z'}$, then $z=z'$.
\end{crl}
\begin{proof}
Set $(X,\leq):=(L_{\max\{z,z'\}},\leq)$ and $x:=\min\{z,z'\}$ in the preceding result. Then we see that if $z\neq z'$ (i.e., $\min\{z,z'\}<\max\{z,z'\}$), then $L_{\max\{z,z'\}}=X\not\cong_oL_x=L_{\min\{z,z'\}}$, and so $L_z\not\cong_oL_{z'}$.
\end{proof}

\begin{thm}[\blue{\cite[Theorem 7.5, p.197]{goldrei1996}}]\label{WellOrdSetThm}
Let $(X,\leq_X)$ and $(Y,\leq_Y)$ be well-ordered sets. Then exactly one of the following holds:
\bit
\item[1.] $X$ is order-isomorphic to an lower segment of $Y$.
\item[2.] $X$ is order-isomorphic to $Y$.
\item[3.] $Y$ is order-isomorphic to an lower segment of $X$.
\eit
\end{thm}
\begin{proof}
Consider the relation $R:=\{(x,y)\in X\times Y:L_x\cong_oL_y\}$, where $L_x\subset X$ and $L_y\subset Y$ are the strict lower segments at $x\in X$ and $y\in Y$ respectively.  The domain and image of $R$ are $\dom R:=\{x\in X:(x,y)\in R~\txt{for some}~y\in Y\}$ and $\im R:=\{y\in Y:(x,y)\in R~\txt{for some}~x\in X\}$. We have $(\min X,\min Y)\in R$, and so $R\neq\emptyset$.

If $x\in\dom X$, let $y,y'\in Y$ be such that $(x,y),(x,y')\in R$, i.e., $L_x\cong L_y$ and $L_x\cong L_{y'}$. Then $L_y\cong L_{y'}$ implies $y=y'$. This means $R$ is a \ul{map} $R:\dom R\ra Y,~x\mapsto y_x$ via $L_x\cong L_{y_x}$. If $x\neq x'$, then
\begin{align}
L_{y_x}\cong L_x\not\cong L_{x'}\cong L_{y_{x'}}~~\Ra~~L_{y_x}\not\cong L_{y_{x'}},~~\Ra~~y_x\neq y_{x'}.\nn
\end{align}
That is, the map $R:\dom R\ra Y,~x\mapsto y_x$ is \ul{injective}.

To show the map $R:\dom R\ra Y,~x\mapsto y_x$ is order-imbedding, let $x,x'\in \dom R$. Then we have order-isomorphisms $f:L_x\ra L_{y_x}$ and $f':L_{x'}\ra L_{y_{x'}}$. If $x<x'$, then $L_x\subsetneq L_{x'}$, and so
\bea
f'|_{L_x}:L_x\subsetneq L_{x'}\ra L_{f'(x)}\subsetneq L_{y_{x'}},~~\Ra~~L_{f'(x)}\cong L_x\cong L_{y_x},~~\Ra~~y_x=f'(x)<y_{x'}.\nn
\eea
That is, $R:\dom R\ra Y,~x\mapsto y_x$ is \ul{order-imbedding}. Hence $R:\dom R\ra\im R,~x\mapsto y_x$ is an \ul{order-isomorphism}.

Also, observe from above that if $x'\in\dom R$, i.e., some order-isomorphism $f':L_{x'}\ra L_{y_{x'}}$ exists, then for any $x<x'$, we get the order-isomorphism $f'|_{L_x}:L_x\subsetneq L_{x'}\ra L_{f'(x)}\subsetneq L_{y_{x'}}$, and so $L_{x'}\subset\dom R$. That is, if $x\in\dom R$ then $L_x\subset\dom R$. It follows that ~$\dom R=\bigcup_{a\in A\subset X}L_a$,~ for some subset $A\subset X$. By symmetry, we also have ~$\im R=\bigcup_{b\in B\subset Y}L_b$,~ for some subset $B\subset Y$.

Hence either $\dom R=X$ or $\dom R=L_m$, where $m:=\min(X\backslash\dom R)$. Similarly, either $\im R=Y$ or $\im R=L_n$, where $n:=\min(Y\backslash\im R)$. It is now clear that at least one of 1., 2., 3. holds.

Finally, by the fact that each of $X$ and $Y$ cannot be order-isomorphic to one of its lower segments, it follows that at most one of 1., 2., 3. can hold. Hence exactly one of 1., 2., 3. holds.
\end{proof}

\begin{dfn}[\blue{
\index{Addition! of ordered sets}{Addition of ordered sets},
\index{Order-sum}{Order-sum},
\index{Multiplication! of ordered sets}{Multiplication of ordered sets},
\index{Order-product}{Order-product},
\index{Exponentiation of ordered sets}{Exponentiation of ordered sets},
\index{Order-exponent}{Order-exponent}}]
Let $(A,\leq_A)$ and $(B,\leq_B)$ be linearly-ordered sets. Then their \ul{order-sum} (\blue{footnote}\footnote{As done with natural numbers (where $m+n:=A_m(n)$ = ``adding to m'' a variable $n$), $A+_oB$ is to be interpreted as ``adding to $A$'' a variable $B$.}) $A\substack{+_o}B=(A\substack{+_o}B,\leq)$ is (up to order isomorphism) the set $A\substack{+_o}B:=A\sqcup B:=(A\times\{0\})\cup(B\times\{1\})$ as a linearly-ordered set with order given by the following: For all $a,a'\in A$, $b,b'\in B$,
\bit
\item[] (i) $(a,0)<(a',0)$ if $a<_Aa'$, (ii) $(a,0)<(b,1)$, and (iii) $(b,1)<(b',1)$ if $b<_Bb'$.
\eit
With the disjoint union $A\times B=\bigcup_{b\in B}A\times\{b\}=\bigsqcup_{b\in B}A\times\{b\}$, viewing this as the sum $+_o$ (over $B$) of $B$-indexed copies of $A$, we introduce a linear-order on $A\times B$ as follows: For all $a,a'\in A$, $b,b'\in B$,
\bit
\item[] (i) $(a,b)<(a',b)$ if $a<_Aa'$, and (ii) $(a,b)<(a',b')$ if $b<_Bb'$.
\eit
(That is, for all $(a,b),(a',b')\in A\times B$, ~$(a,b)<(a',b')$~ if (i) $b=b'$ and $a<_Aa'$ or (ii) $b<_Bb'$). The \ul{order-product} (\blue{footnote}\footnote{As done with natural numbers (where $mn:=M_m(n)$ = ``multiplying $m$ by'' a variable $n$), $A\times_oB$ is to be interpreted as ``multiplying $A$ by'' a variable $B$. Also, if $B$ is finite then, we can consider an inductive definition of the multiplication as follows: $A\times_oB:=\big(A\times_o(B-\{\txt{one element}\})\big)+_oA\cong_oA+_o\big(A\times_o(B-\{\txt{one element}\})\big)$.}) $A\substack{\times_o}B=(A\substack{\times_o}B,\leq)$ of $A$ and $B$ is the set $A\times B$ as a linearly-ordered set with the above-defined order. That is, $A\substack{\times_o}B~\substack{\cong_o}~(+_o)_{b\in B}A~\substack{\cong_o}~(+_o)_{b\in B}A\times\{b\}$, where in $A\times\{b\}$, ``$(a,b)<(a',b)$ if $a<_Aa'$''.
\end{dfn}

By construction, $\substack{+_o}$ and $\substack{\times_o}$ are \ul{order-associative} (or associative up to order-isomorphism), i.e.,
\[
(A\substack{+_o}B)\substack{+_o}C~\substack{\cong_o}~A\substack{+_o}(B\substack{+_o}C)~~~~\txt{and}~~~~(A\substack{\times_o}B)\substack{\times_o}C~\substack{\cong_o}~A\substack{\times_o}(B\substack{\times_o}C),
\]
but are in general not \ul{order-commutative} (i.e., not commutative up to order-isomorphism), i.e.,
\[
A\substack{+_o}B\not\!\substack{\cong_o}B\substack{+_o}A~~~\txt{and}~~~~A\substack{\times_o}B\not\!\substack{\cong_o}B\substack{\times_o}A.
\]Also, by construction $A\times(B+C)\cong (A\times B)+(A\times C)$, while $(A+B)\times C\not\cong(A\times C)+(B\times C)$ in general. As done in these equations, for convenience, we will often write $A+B$, $A\times B$, $A\cong B$ instead of $A\substack{+_o}B$, $A\substack{\times_o}B$, $A\substack{\cong_o}B$ respectively.

Based on the multiplication introduced above, if $B$ is finite, we can define \ul{exponentiation} of $A$ by $B$ inductively as follows: $A_o^B:=A_o^{B-\{\txt{one element}\}}\times_oA\cong_oA\times_oA_o^{B-\{\txt{one element}\}}$. In particular, $A_o^n:=A_o^{n-1}\times_oA\cong_oA\times_oA_o^{n-1}$. For a general definition of order exponentiation, we need to put a suitable linear order (i.e., one that is consistent, up to order-isomorphism, with that considered for multiplication above) on the Cartesian product $A^B:=\prod_{b\in B}A$. On $A^n$ we will therefore consider the following linear order: For each {\small $1\leq i\leq n$, ~$(a_0,a_1,...,a_i,a_{i+1},...,a_{n-1})< (a'_0,a'_1,...,a'_i,a_{i+1},...,a_{n-1})$ if $a_i<a'_i$}. Similarly, on {\small $A^B=\prod_{b\in B}A=\{(a_b)_{b\in B}:a_b\in A\}$} we will consider the following linear order: For each $\ol{b}\in B$,
\bea
(a_b)_{b\in B}<(a'_b)_{b\in B}~~\txt{if}~~\txt{$a_{\ol{b}}<a'_{\ol{b}}$ and $a_b=a'_b$ for all $b>\ol{b}$}.\nn
\eea
Up to order-isomorphism, the resulting linearly-ordered set {\small $A_o^B=(A_o^B,\leq)$} is the \ul{order-exponent} of $A$ by $B$.

\begin{dfn}[\blue{\index{Ordered naturals $\omega$}{Ordered naturals $\omega$}}]
We define the well-ordered set $\omega=(\omega,<)$ such that (i) $\omega:=\Natural$ and (ii) $<$ is given (for all $m,n\in\omega$) by ``$m<n$ if $m\in n$'' (i.e., ``$m\leq n$ if $m\in n$ or $m=n$'').

That is, ~$\omega=(\omega,<):=(\Natural,\in)$ as an ordered set.
\end{dfn}

\begin{rmk}[\blue{Order-commutativity in $\omega$}]
If $m,n\in\omega$, then it is clear by induction that as a set,
\bea
&&m+n=(m+n-1)\cup\{m+n-1\}=m\cup\{m+0,m+1,...,m+n-1\}\nn\\
&&~~~~\substack{\cong_o}(m\times\{0\})\cup(\{0,1,...,n-1\}\times\{1\})=(m\times\{0\})\cup(n\times\{1\})\nn\\
&&~~~~=m\substack{+_o}n.\nn
\eea
That is, up to order-isomorphism, we have ~$m+n=m\substack{+_o}n$. It follows that $\substack{+_o}|_{\omega^2}:\omega^2\ra\omega,~(m,n)\ra m\substack{+_o}n$ is order-commutative, since ~$m\substack{+_o}n~\substack{\cong_o}~m+n=n+m~\substack{\cong_o}~n\substack{+_o}m$. Similarly,
{\small\bea
m\substack{\times_o}n=\ub{m\substack{+_o}m\substack{+_o}\cdots\substack{+_o}m}_n~\substack{\cong_o}~ \ub{m+m+\cdots+m}_n=m\cdot n=n\cdot m~\substack{\cong_o}~n\substack{\times_o}m\nn
\eea}shows ~$\substack{\times_o}|_{\omega^2}:\omega^2\ra\omega,~(m,n)\ra m\substack{\times_o}n$~ is order-commutative.
(\blue{footnote}\footnote{By \cite[Exercise 7.27, p.179]{goldrei1996}, $\omega~\substack{\cong_o}~2\substack{+_o}\omega~\not\!\substack{\cong_o}~\omega\substack{+_o}2$. By \cite[Exercise 7.36, p.183]{goldrei1996}, $\omega~\substack{\cong_o}~2\substack{\times_o}\omega~\not\!\substack{\cong_o}~\omega\substack{\times_o}2~\substack{\cong_o}~\omega\substack{+_o}\omega$.}). It follows by induction that
\bea
m_o^n:=m_o^{n-1}\times_om\cong_om^{n-1}\cdot m=m^n~~\txt{and}~~~m_o^n:=m_o^{n-1}\times_om\cong_om\times_om_o^{n-1}.\nn
\eea
\end{rmk}

\subsection{The construction of ordinal numbers}~\\~
Let $A$ be a set and $(A,\leq)$ a linear-ordering of $A$. As before, if $a\in A$, we write ~$[A,a):=L_a:=L_A(a):=Seg_A(a):=\{a'\in A~|~a'<a\}$. Recall that a nonempty subset $\emptyset\neq B\subset A$ is \ul{$(\leqslant)$-inductive} (i.e., \ul{inductive with respect to $\leqslant$}) if for all $a\in A$, ~$[A,a)\subset B$ ~$\Ra$~ $a\in B$. (\blue{footnote}\footnote{If $(A,\leq)$ is well-ordered, then every inductive subset $B\subset A$ contains $m:=\min A$, since $[A,m)=\emptyset\subset B$.}). Recall the \ul{transfinite induction principle for well-ordered sets} from Theorem \ref{TFIPforLOS} (on page \pageref{TFIPforLOS}): A linearly-ordered set $(A,\leq)$ is well-ordered $\iff$ every $(\leq)$-inductive subset $B\subset A$ equals $A$.

\begin{dfn}[\textcolor{blue}{\index{Ordinal! numbers $\Ordinal$}{Ordinal numbers $\Ordinal$}, \index{Transitive! set}{Transitive set}: \cite[p.189]{enderton1977}, \cite[p.204]{goldrei1996}}]\label{Ordinals}
A set $\al$ is an \ul{ordinal number}, written $\al\in\Ordinal$, if $\al$ has the following properties:
\bit
\item[(i)] $\al$ is well-ordered by $\in$, in the sense that $(\al,<):=(\al,\in)$ is a well-ordering of $\al$ (where for $a,b\in\al$, we say ``$a<b$ if $a\in b$'').
\item[(ii)] $\al$ is a \ul{transitive set}, in the sense that ~$a\in\al$ ~$\Ra$~ $a\subset\al$ ~(or equivalently, ~$a\in b\in\al$ ~$\Ra$~ $a\in\al$).
\eit
The collection of all \ul{ordinal numbers} (\blue{footnote}\footnote{The \ul{ordinal numbers} $\Ordinal$ are a \ul{proper class of sets} (Theorem \ref{BonFortThm}) that is meant to be a generalization of the natural numbers $\Natural$.}), denoted by $\Ordinal$, is viewed as the partially ordered class $(\Ordinal,<):=(\Ordinal,\in)$, i.e., if $\al,\beta\in\Ordinal$, then ``$\al\leq\beta$ if $\al\in\beta$ or $\al=\beta$''.
\end{dfn}

\begin{lmm}[\blue{Elements of ordinals are lower segments}]
Let $\al\in\Ordinal$. If $a\in\al$ (hence $a\subset\al$), then
\[
Seg_\al(a):=\{a'\in\al:a'<a\}=\{a'\in\al:a'\in a\subset\al\}=a.
\]
\end{lmm}

From here on, as done above, if $(A,\leq)$ is linearly-ordered, we will preferably use the more explicit notation $Seg_A(a)$ to denote the lower segment at $a\in A$. Therefore it is worth recalling that {\small $Seg_A(a):=L_A(a):=L_a:=[A,a):=\{a'\in A:a'<a\}$}.

\begin{thm}[\textcolor{blue}{\index{Transfinite induction principle for! ordinals}{Transfinite induction principle for ordinals}: \cite[Theorem 8.1, p.224]{goldrei1996}}]
Let $\al\in\Ordinal$ and $A\subset\al$. If for all $a\in\al$, ~$a\subset A$ ~$\Ra$~ $a\in A$, then $A=\al$.
\end{thm}
\begin{proof}
This follows immediately from Theorem \ref{TFIPforLOS} (transfinite induction for linearly-ordered sets).
\end{proof}

\begin{lmm}[\blue{\cite[Theorem 8.2, p.206]{goldrei1996}}]\label{OrdinalEqlt}
Let $\al,\beta\in\Ordinal$. If $\al\cong_o\beta$, then $\al=\beta$.
\end{lmm}
\begin{proof}
Let $f:\al\ra\beta$ be an order-isomorphism. Suppose $f\neq id_\al$. Then $A:=\{a\in\al:f(a)\neq a\}\neq\emptyset$. Let $m:=\min A$. Then $f|_{Seg_\al(m)}:Seg_\al(m)\subset\al\ra Seg_\beta(f(m))\subset\beta$, i.e., we have an order-isomorphism $f|_m:m\subset\al\ra f(m)\subset\beta$. For every $a'\in m$ (i.e., $a'<m$), $a'\not\in A$, and so $f(a')=a'$. That is,
\[
f(m)=Seg_\beta(f(m))=Seg_\al(m)=m~~~~\txt{(a contradiction since $m\in A$)}.\qedhere
\]
\end{proof}

\begin{dfn}[\textcolor{blue}{\index{Successor function}{Successor function}}]
The \ul{successor function} of sets is the map
\[
f_+:Sets\ra Sets,~s\mapsto s^+:=s\cup\{s\}.
\]
\end{dfn}
If $T\subset Sets$ is any class of sets such that $(T,<):=(T,\in)$ is a transitive poset (i.e., for all $a,b,c\in T$, $a\in b\in c$ $\Ra$ $a\in c$) then the restriction $f_+|_T:T\ra Sets$ is injective, since for any $s,s'\in T$,
\bea
&&s\cup\{s\}=s'\cup\{s'\}~~\iff~~s\cup\{s\}\subset s'\cup\{s'\}~~\txt{and}~~s'\cup\{s'\}\subset s\cup\{s\},\nn\\
&&~~\Ra~~s\in s'\cup\{s'\}~~\txt{and}~~s'\in s\cup\{s\}~~\Ra~~(s\in s'~~\txt{or}~~s=s')~~\txt{and}~~(s'\in s~~\txt{or}~~s'=s)\nn\\
&&~~\Ra~~s=s'.\nn
\eea
In particular, the following result shows that (i) $f_+|_\Ordinal:\Ordinal\ra Sets$~ is injective and (ii) $\im f_+|_\Ordinal\subset\Ordinal$ (but unlike for natural numbers, we will see that ~$\im f_+|_\Ordinal\subsetneq\Ordinal\backslash\{0\}$~ due to the existence of ``limit ordinals'').

\begin{lmm}[\textcolor{blue}{Transitivity, Succession, Heredity in $\Ordinal$}]\label{FirstOrdPts}
Let $\al,\beta,\gamma\in\Ordinal$. Then the following hold.
\bit
\item[] (a) $\al\not\in\al$.~~ (b) $\al\in\beta\in\gamma$ implies $\al\in\gamma$.~~ (c) $\al^+\in\Ordinal$.~~ (d) $\al\subset\Ordinal$.~~ (e) $\al\leq\beta$ $\iff$ $\al\subset\beta$.
\eit
(Part (e) equivalently means ~$\al\in\beta$ $\iff$ $\al\subsetneq\beta$.)
\end{lmm}
\begin{proof}
{\flushleft (a)}: Since $\al$ is well-ordered (hence linearly-ordered) by $\in$, if $\al\in\al$, then $\al<\al$, $\Ra$ $\al\neq\al$ (a contradiction).
{\flushleft (b)}: If $\al\in\beta\in\gamma$, then $\al\in\beta\subset\gamma$, and so $\al\in\gamma$.
{\flushleft(c)} (i): If $a\in\al^+:=\al\cup\{\al\}$, then either $a\in\al$ or $a=\al$ (but not both, by part (a)). Thus, for any $a,b\in \al^+=\al\cup\{\al\}$, we get either $a<b$ (i.e., $a\in b$) or $a=b$ or $a>b$ (i.e., $b\in a$), and so $\al^+$ is \ul{linearly-ordered} by $\in$. If $\emptyset\neq C\subset \al^+$, then it is also clear that $C=(\al\cap C)\cup(C\cap\{\al\})$ has a least element which is either $\min(C\cap\al)$, if $C\cap\al\neq\emptyset$, or $\al$ otherwise. Hence $\al^+$ is \ul{well-ordered} by $\in$. (\blue{footnote}\footnote{A subset $B\subset\al^+$ is inductive if for any $a\in \al^+$, $[\al^+,a):=\{a'\in\al^+:a'<a\}=\{a'\in\al^+:a'\in a\}=a\subset B~~\Ra~~a\in B$.}).

(ii): If $a\in \al^+$ (i.e., $a\in\al$ or $a=\al$), then $a\subset\al$ (since $\al$ is an ordinal), i.e., $a\subset\al\subset\al^+$, and so $a\subset\al^+$. Hence $\al^+$ is an ordinal.
{\flushleft(d)} To see that $\al\subset\Ordinal$ (i.e., if $a\in\al$, then $a\in\Ordinal$), we need to show (i) $a\in\al$ is well-ordered by $\in$ and (ii) $a'\in a$ $\Ra$ $a'\subset a$. For (i), $a\in\al$ $\Ra$ $a\subset\al$, $\Ra$ $a$ is well-ordered by $\in$, since $\al$ is. For (ii), let $a''\in a'\in a\in\al$. Then $a\subset\al$ $\Ra$ $a'\in\al$ ($\Ra$ $a'\subset\al$, $\Ra$ $a''\in\al$). This implies $a'',a',a\in\al$ satisfy $a''\in a'\in a$, and so $a''\in a$ (by the transitivity of $\in$ in $\al$).

{\flushleft (e)}: If $\al\leq\beta$ (i.e., $\al\in\beta$ or $\al=\beta$) then it is clear that $\al\subset\beta$. Conversely, assume $\al\subset\beta$. Then by Theorem \ref{WellOrdSetThm}, either $\al\cong_oSeg_\beta(a)=a$ for some $a\in\beta$, or $\al\cong_o\beta$, and so by Lemma \ref{OrdinalEqlt}, either $\al=a\in\beta$ or $\al=\beta$.
\end{proof}

\begin{crl}\label{WellOrdUnCrl1}
Let $A\subsetneq\Ordinal$ be a set. If $(A,<):=(A,\in)$ is well-ordered (\blue{footnote}\footnote{We will also show that ``every set of ordinals is well-ordered by $\in$'', but we are not assuming the result at this point.}), then $lub_\Ordinal(A)=\bigcup A\in\Ordinal$. In particular, for any ascending $\omega$-chain of ordinals $\al_0,\al_1,\cdots$, we have $\bigcup_{i\in\omega}\al_i\in\Ordinal$.
\end{crl}
\begin{proof}
(i) Since $A$ is \ul{linearly-ordered} by $\in$, so is {\small $\bigcup A:=\bigcup_{\al\in A}\al$}. If {\small $\emptyset\neq C\subset \bigcup A$}, then for any $c\in C$, let $c\in\al_c\in A$. Then $\emptyset\neq A':=\{\al_c\}_{c\in C}\subset A$ has a least element $\al_{c_m}:=\min A'$, i.e., $\al_{c_m}\leq\al_c$ for all $c\in C$. Let $m:=\min\al_{c_m}$. Then (because for all $c\in C$, $m\in\al_{c_m}\subset\al_c$ $\Ra$ $m,c\in\al_c$) $m\leq c$ for all $c\in C$, and so $m=\min C$. This shows $\bigcup A$ is well-ordered by $\in$. (ii) If $a\in \bigcup A$, then with $a\in\al_a\in A$, we get $a\subset\al_a\in A$, and so $a\subset\bigcup A$. This proves $\bigcup A\in\Ordinal$.

It is clear that $A\leq\bigcup A$ in $\Ordinal$ (i.e., $\al\leq\bigcup A$ for all $\al\in A$). If $\beta\in\Ordinal$ is an upper bound of $A$, i.e., $A\leq\beta$, suppose $\beta<\bigcup A$, i.e., $\beta\in\bigcup A$. Then $\beta\in \al$ (i.e., $\beta<\al$) for some $\al\in A$, which is a contradiction. This proves $\bigcup A=lub_\Ordinal(A)$, the least upper bound of $A$ in $\Ordinal$.

If $\al_0\subset\al_1\subset\cdots$ is an ascending $\omega$-chain of ordinals, let $A:=\{\al_i\}_{i\in\omega}$. Then by Lemma \ref{FirstOrdPts}(e), $\al_i\leq\al_j$ in $\Ordinal$ iff $\al_i\subset\al_j$, iff $i\leq j$. Thus, by the well-ordering of $\omega$, $(A,\in)$ is well-ordered. Hence $\bigcup_{i\in\omega}\al_i\in\Ordinal$.
\end{proof}

\begin{thm}[\textcolor{blue}{\index{Burali-Forti theorem (or paradox)}{Burali-Forti theorem (or paradox)}: \cite[p.194]{enderton1977}}]\label{BonFortThm}
There is no set to which every ordinal number belongs.
\end{thm}
\begin{proof}
Suppose $A$ is a set such that $\Ordinal\subset A$. Then $\Ordinal$ is a set (as a subset), and so by its properties from Lemma \ref{FirstOrdPts}, $\Ordinal$ an ordinal number, i.e., $\Ordinal\in\Ordinal$ (a contradiction).
\end{proof}

With the properties of ordinals we have seen so far, we can now list a few countable ordinals.
\begin{dfn}[\textcolor{blue}{\index{Integrating map}{Integrating map}, \index{Integral! of a set}{Integral of a set}}]
The \ul{integrating map} is the map ~$\omega\times Sets\ra Sets$, $(n,s)\mapsto s[n]$ ~given (for any set $s$ and any ~$n\in\omega=\Natural=\{\emptyset,\emptyset^+,\emptyset^{++},\cdots\}=\{0,1,2,\cdots\}$) by
\bea
s[0]:=s,~~~~s[n^+]:=s[n]\cup\big\{s[n],s[n]^+,s[n]^{++},\cdots\big\}=s[n][1].\nn
\eea
We will call the set $s[n]$ the \ul{$n$-integral} of the set $s$. (\blue{footnote}\footnote{Note that ~$s[n+m]=s[n][m]=s[m][n]$ ~for all $n,m\in\omega$.}). Note that
\bea
s^+=s\cup\{s\},~~s^{++}=s\cup\{s,s^+\},~~s^{+++}=s\cup\{s,s^+,s^{++}\},~~...,~~s^{+++\cdots}=s\cup\{s,s^+,s^{++},\cdots\}=s[1].\nn
\eea
\end{dfn}

\begin{examples}[\blue{Examples of countable ordinal numbers}]
The following are some countable ordinal numbers (each has a successor but need not have a predecessor, i.e., limit ordinals such as $\omega$ exist):
{\small\begin{align}
&\Ordinal\supset\big\{\emptyset,\emptyset^+,\emptyset^{++},\cdots,\omega,\omega^+,\omega^{++},\cdots,\omega[1],\omega[1]^+,\omega[1]^{++},\cdots,\omega[2],\omega[2]^+,\omega[2]^{++},\cdots\big\}\nn\\
&~~~~=\overbrace{\ub{\omega\cup\{\omega,\omega^+,\omega^{++},\cdots\}}_{\omega[1]}\cup\big\{\omega[1],\omega[1]^+,\omega[1]^{++},\cdots\big\}}^{\omega[2]}\cup\big\{\omega[2],\omega[2]^+,\omega[2]^{++},\cdots\big\}\cup\cdots,\nn\\
&~~~~=\omega[0]\cup\omega[1]\cup\omega[2]\cup\omega[3]\cup\cdots=0[1]\cup0[2]\cup0[3]\cup0[4]\cup\cdots=0[0]\cup0[1]\cup0[2]\cup0[3]\cup0[4]\cup\cdots,\nn
\end{align}}where ~$0[0]=0=\emptyset$, ~$0[1]=0[0]\cup\{0[0],0[0]^+,0[0]^{++},\cdots\}=\omega=\omega[0]$, ~and ~$0[n^+]=0[1][n]=\omega[n]$.
\end{examples}

\subsection{Further properties of ordinal numbers}
\begin{thm}[\textcolor{blue}{The class of ordinals $\Ordinal$ is ``well-ordered by $\in$'': \cite[Theorem 8.3, p.207]{goldrei1996}}]\label{WellOrdSubThm}
Every \ul{set} of ordinals $A\subsetneq\Ordinal$ is well-ordered. That is, for any set $A\subset\Ordinal$, the poset $(A,<):=(A,\in)$ is well-ordered.
\end{thm}
\begin{proof}
{\flushleft \ul{$(A,<):=(A,\in)$ is linearly-ordered}}: Let $a,a'\in A$ be distinct, i.e., $a\neq a'$. We need to show that $a\in a'$ or $a'\in a$. By Theorem \ref{WellOrdSetThm}, either $a\cong_oSeg_{a'}(c)=c$ for some $c\in a'$, or $a'\cong_oSeg_a(c')=c'$ for some $c'\in a$. Therefore, by Lemma \ref{OrdinalEqlt}, $a=c\in a'$ or $a'=c'\in a$.

{\flushleft \ul{$(A,<):=(A,\in)$ is well-ordered}}: Let $\emptyset\neq C\subset A$. We need to show $\min_{(A,\in)}C$ exists. Let $c\in C$. Define $C':=\{c'\in C:c'\in c\}=c\cap C$. Either (i) $C'=\emptyset$ or (ii) $C'\neq\emptyset$.
(i) If $C'=\emptyset$, then $c=\min_{(A,\in)}C$ since $(A,\in)$ is linearly-ordered (and so $C'=\emptyset$ means $c\in c'$ for every $c'\in C\backslash\{c\}$). (ii) If $C'\neq\emptyset$ then as a nonempty subset $\emptyset\neq C'\subset c$ of the ordinal $c$, $C'$ has a least element $m:=\min C'$. Since $(A,\in)$ is linearly-ordered, $m=\min C'=\min C$ (because for any $c''\in C\backslash C'$, we have $c\leq c''$, and so $m<c\leq c''$).
\end{proof}

\begin{crl}[\blue{\cite[Theorem 8.4, p.209]{goldrei1996}}]\label{WellOrdUnCrl2}
For any \ul{set} of ordinals $A\subsetneq\Ordinal$, we have $lub_\Ordinal(A)=\bigcup A\in\Ordinal$. That is, the union of any \ul{set} of ordinals is the least upper bound (or supremum) of the set in $\Ordinal$.
\end{crl}
\begin{proof}
By Theorem \ref{WellOrdSubThm}, $(A,\in)$ is well-ordered, and so $\bigcup A\in\Ordinal$ by Corollary \ref{WellOrdUnCrl1}.
\end{proof}

\begin{notation}[\textcolor{blue}{Initial segments}]
Let $(A,\leq_A)$ and $(B,\leq_B)$ be linearly-ordered sets. Write $A\leq_{\txt{seg}}B$ if $A$ is an \ul{initial segment} of $B$ (i.e., $A\subset B$ and $A<b$ for all $b\in B-A$). Similarly, write $A<_{\txt{seg}}B$ if $A$ is a \ul{strict initial segment} of $B$ (i.e., $A\subsetneq B$ and $A<b$ for all $b\in B-A$).
\end{notation}
\begin{crl}
Let $\al\in\Ordinal$. Every initial segment $s\leq_{\txt{seg}}\al$ (resp. strict initial segment $s<_{\txt{seg}}\al$) of $\al$ is an ordinal, i.e., $s\in\Ordinal$.
\end{crl}
\begin{proof}
Since $s\subset\al$ is an initial segment (resp. $s\subsetneq\al$ is a strict initial segment), for each $\al'\in s$, we have $\al'=Seg_\al(\al')\subset s$, and so $s=\bigcup_{\al'\in s}Seg_\al(\al')=\bigcup_{\al'\in s}\al'\in\Ordinal$.
\end{proof}

\begin{lmm}[\blue{\cite[Theorem 8.5, p.212]{goldrei1996}}]\label{OrdSuffLmm}
Let $(X,\leq)$ be a well-ordered set. Then for each $x\in X$, there exists a unique ordinal $\al_x\in\Ordinal$
such that ~$Seg_X(x)\cong_o\al_x$.
\end{lmm}
\begin{proof}
Let $A:=\{x\in X:Seg_X(x)\cong_o\al~\txt{for some}~\al\in\Ordinal\}$. Then $A$ is nonempty since $\min X\in A$. We need to show $A\subset X$ is inductive, hence that $A=X$. Pick any $x\in X$ such that $Seg_X(x)\subset A$. We need to show $x\in A$. For each $x'<x$, $Seg_X(x')\cong \al_{x'}$ for an ordinal $\al_{x'}\in\Ordinal$, which is unique by Lemma \ref{OrdinalEqlt}. Thus, we get a \ul{map} $f:Seg_X(x)\ra\Ordinal,~x'\mapsto\al_{x'}$, whose image $\im f:=\{{\al'}:=\al_{x'}\in\Ordinal~|~Seg_X(x')\cong {\al'}~\txt{for some}~x'<x\}$ is a set by the axiom of replacement. Moreover, as before, $f$ is \ul{injective} because $x'\neq x''$ $\Ra$ $Seg_X(x')\not\cong  Seg_X(x'')$ $\Ra$ $\al_{x'}\not\cong \al_{x''}$ $\Ra$ $\al_{x'}\neq\al_{x''}$.

We will prove that $\im f\in\Ordinal$, i.e., $(i)$ $\im f\subset\Ordinal$ is well-ordered (which is already the case by Theorem \ref{WellOrdSubThm}) and (ii) ${\al'}\in \im f$ $\Ra$ ${\al'}\subset\im f$. For (ii), let $h:Seg_X(x')\ra{\al'}\in\im f$ be an order isomorphism. If $\al''\in{\al'}\in\im f$, then we get an order isomorphism
\bea
&&h^{-1}|_{\al''}:{\al''}=Seg_{\al'}({\al''})\subsetneq{\al'}\ra Seg_X(h^{-1}({\al''}))\subsetneq Seg_X(x'),~~~~\Ra~~x'':=h^{-1}({\al''})<x'<x,\nn\\
&&~~\txt{or}~~h|_{Seg_X(x'')}:Seg_X(x'')\subsetneq Seg_X(x')\ra{\al''}\subsetneq{\al'},~~\Ra~~{\al''}\in\im f,\nn
\eea
which proves $\al'\subset\im f$, and so $\im f\in\Ordinal$.

To prove $f:Seg_X(x)\ra\im f$ is an \ul{order isomorphism}, if $x_1'<x_2'<x$, consider order isomorphisms
$h_1:Seg_X(x_1')\ra\al_1'$ and $h_2:Seg_X(x_2')\ra\al_2'$. Then because of the induced order isomorphism
\[
h_2|_{Seg_X(x_1')}:Seg_X(x_1')\subsetneq Seg_X(x_2')\ra Seg_{\al_2'}(h_2(x_1'))=h_2(x_1')~\subsetneq\al_2',~\Ra~h_2(x_1)<\al_2',
\]
the uniqueness of $\al_1'$ implies ~$\al_1'=h_1(x_1')<\al_2'$.

Therefore, with $\al_x:=\im f$, we get an order isomorphism $f:Seg_X(x)\ra\al_x$, and so $x\in A$. That is, $A\subset X$ is inductive. Hence $A=X$.
\end{proof}

\begin{thm}[\textcolor{blue}{\cite[Theorem 8.6, p.214]{goldrei1996}}]\label{OrdSuffThm}
Let $(X,\leq)$ be a well-ordered set. Then there exists a unique ordinal $\al_X\in\Ordinal$ such that $X\cong_o\al_X$.
\end{thm}
\begin{proof}
Let $Y:=X+_o1:=(X\times\{0\})\cup(1\times\{1\})=\{(x,0):x\in X\}\cup\{(0,1)\}$ be the order-sum of $X$ and $1=\{0\}\subset\omega$. Then with $y:=(0,1)\in Y$, we have $X\cong_oSeg_Y(y)$. Hence the result follows from Lemma \ref{OrdSuffLmm}.
\end{proof}

\begin{thm}[\textcolor{blue}{\index{Hartogs' theorem}{Hartogs' theorem}: \cite[p.195]{enderton1977},  \cite[Theorem 8.18, p.254]{goldrei1996}}]\label{HartSetThm}
For any set $X$, there exists an ordinal $\al\in\Ordinal$ such that $|\al|\not\leq|X|$ ~(i.e., such that there is \ul{no} injective map ~$\al\hookrightarrow X$). (\blue{footnote}\footnote{The axiom of choice (hence Lemma \ref{CarOrdLin} on the linearity of cardinal order) is not assumed here, otherwise, the result follows immediately from Cantor's theorem (Theorem \ref{CantorThm}) and Theorem \ref{OrdSuffThm}.}).
\end{thm}
\begin{proof}
Let $V:=\{\beta\in\Ordinal:|\beta|\leq|X|\}$. If $\beta\in{V}$, i.e., $|\beta|\leq|X|$, then we have an injective map $f=f_\beta:\beta\hookrightarrow X$. Let $Y=Y_\beta:=\im f=(Y,\leq_Y)$, ordered such that for all $y,y'\in Y$, ``$y\leq_Yy$ if $f^{-1}(y)\leq f^{-1}(y')$''. Then it is clear that $f:\beta\ra Y\subset X$ is an order-isomorphism (hence $(Y,\leq_Y)$ is well-ordered). Thus, we have the nonempty set of well-orderings on subsets of $X$ given by
\bea
W:=\{R\subset X\times X~|~(Y,\leq_Y):=(Y,R|_{Y\times Y})~\txt{is well-ordered for some}~Y\subset X\}.\nn
\eea
It is clear that for each $R\in W$, we have a subset $Y_R\subset X$ and a well-ordered set $(Y_R,\leq_R):=(Y_R,R|_{Y_R\times Y_R})\cong_o\beta_R$ (by Theorem \ref{OrdSuffThm}) for a unique ordinal $\beta_R$ (where it is clear that $\beta_R\in {V}$, i.e., $|\beta_R|=|Y_R|\leq|X|$). This gives a \ul{map} ~$F:W\ra{V},~R\mapsto\beta_R$.~ It is also clear that for each $\beta\in{V}$, there is $R_\beta\in W$ and $Y_\beta\subset X$ such that ~$(Y_\beta,\leq_{Y_\beta}):=(Y_\beta,R_\beta|_{Y_\beta\times Y_\beta})\cong_o\beta$.~ This implies the map $F$ is \ul{surjective}. Hence ${V}$ is a set (by the axiom of replacement).

(i) As a set of ordinals, ${V}\subset\Ordinal$ is well-ordered by Theorem \ref{WellOrdSubThm}. (ii) If $\beta\in V$, then for any $\beta'\in\beta\subset\Ordinal$, we know $\beta'\subset\beta$, and so $|\beta'|\leq|\beta|\leq|X|$, i.e., $\beta'\in V$.  That is, if $\beta\in V$, then $\beta\subset V$, which shows $V\in\Ordinal$.

With $\al:=V\in\Ordinal$, we have $|\al|\not\leq|X|$, otherwise $\al\in\al$ (a contradiction). (\blue{footnote}\footnote{The obtained ordinal $\al$ is the least ordinal such that $|\al|\not\leq|X|$: Indeed, if $\gamma<\al$, then $\gamma\in\al$, and so $|\gamma|\leq|X|$.}).
\end{proof}

\begin{crl}[\textcolor{blue}{\cite[Theorem 8.19, p.256]{goldrei1996}}]
For every ordinal $\al\in\Ordinal$, there exists an ordinal $\al'\in\Ordinal$ such that $|\al|<|\al'|$, i.e., there is no cardinal of greatest cardinality. (\blue{footnote}\footnote{The axiom of choice (hence the linearity of cardinal order from Lemma \ref{CarOrdLin}) is not assumed, otherwise, the result follows immediately from Theorems \ref{CantorThm} and \ref{OrdSuffThm}.}).
\end{crl}
\begin{proof}
Let $\al\in\Ordinal$. Then by Hartogs' theorem, there is an ordinal $\al'\in\Ordinal$ such that $|\al'|\not\leq|\al|$. Since $\{\al,\al'\}\subset\Ordinal$ is linearly-ordered, either $\al<\al'$, or $\al=\al'$, or $\al>\al'$. But $|\al'|\not\leq|\al|$ implies $\al'\not\leq\al$, and so $\al<\al'$, i.e., $\al\in\al'$. This implies $|\al|\leq|\al'|$, and since $|\al'|\not\leq|\al|$ implies $|\al'|\neq|\al|$, we get $|\al|<|\al'|$.
\end{proof}

\begin{thm}[\textcolor{blue}{\cite[Theorem 6M, p.151]{enderton1977}}]
The following are equivalent:
\bit[leftmargin=1cm]
\item[(a)] \ul{Cardinal comparability}: For any sets $A,B$ we have $|A|\leq |B|$ or $|B|\leq|A|$.
\item[(b)] \ul{Axiom of choice}: For any set $I$ and nonempty sets $(x_i)_{i\in I}$, there is a map {\footnotesize $f:I\ra\bigcup_{i\in I} x_i,~i\mapsto f(i)\in x_i$}.
\item[(c)] \ul{Zorn's lemma}: A poset in which every (well-ordered) chain has an upper bound has a maximal element.
\eit
\end{thm}
\begin{proof}
{\flushleft\ul{(a)$\Ra$(b)}}: Assume (a) holds. Then for any set $A$, Hartogs' theorem implies there exists an ordinal $\al\in\Ordinal$ such that $|A|\leq|\al|$. This means there is an injective map $f:A\hookrightarrow\al$. Consider the ordering $(A,\leq_f)$ of $A$ given by ``$a\leq_f a'$ if $f(a)\leq f(a')$, for all $a,a'\in A$''. Then $f:A\ra \im f$ is an order-isomorphism (wrt $\leq_f$), and so $(A,\leq_f)$ is a well-ordering of $A$.

Given a set $I$ and an $I$-indexing of sets $(x_i)_{i\in I}$, consider a well-ordering $R\subset X\times X$ on the union $X:=\bigcup_{i\in I}x_i$. Then we get a map $f:I\ra \bigcup_{i\in I}x_i,~i\mapsto\min_R x_i\in x_i$, i.e., $\prod_{i\in I}x_i\neq\emptyset$.

{\flushleft\ul{(b)$\Ra$(c)}}: This is the usual proof of Zorn's lemma.

{\flushleft\ul{(c)$\Ra$(a)}}: This is the proof of Lemma \ref{CarOrdLin}.
\end{proof}

\begin{dfn}[\blue{
\index{Countable! ordinal}{Countable ordinal},
\index{Uncountable ordinal}{Uncountable ordinal},
\index{Smallest uncountable ordinal}{Smallest uncountable ordinal},
\index{Limit! ordinal}{Limit ordinal},
\index{Initial ordinals}{Initial ordinals}}]
An ordinal $\al\in\Ordinal$ is a \ul{countable ordinal} if $\al\approx\Natural$, i.e., $|\al|=|\Natural|$. Recall that the countable ordinal $\omega$ is defined to be $\Natural$ viewed as a subset/element of $\Ordinal$. An ordinal $\al\in\Ordinal$ is an \ul{uncountable ordinal} if $|\omega|<|\al|$. The \ul{smallest uncountable ordinal} is denoted by $\omega_1$. (\blue{footnote}\footnote{By \cite[Theorem 8.7, p.215]{goldrei1996}, $\omega_1=\{\txt{countable ordinals}~\al\in\Ordinal\}$. Therefore, $\bigcup\omega_1=lub_\Ordinal(\omega_1)=\omega_1$.}).

Let $\al\in\Ordinal$. Then $\al$ is a \ul{limit ordinal} if $\al$ has no predecessor in $\Ordinal$ (i.e., $\al$ is a limit element in any set of ordinals $A\subsetneq\Ordinal$ such that $\al\in A$ and $\al\subset A$), that is, for any ordinal $\beta<\al$ there exists an ordinal $\beta'$ such that $\beta<\beta'<\al$ (equivalently, for any ordinal $\beta<\al$, we have $\beta^+<\al$). The ordinal $\al$ is an \ul{initial ordinal} if (i) it is infinite, i.e., $|\omega|\leq|\al|$ and (ii) for every $\al'<\al$, we have $|\al'|<|\al|$. (\blue{footnote}\footnote{Every initial ordinal $\al$ is a limit ordinal since (i) $0<\al$ and (i) for any ordinal $\beta<\al$, we have $\beta<\beta^+<\al$.}).

For any $\al\in\Ordinal$, the \ul{$\al$th initial ordinal} is denoted by $\omega_\al$ (with cardinality $\aleph_\al:=|\omega_\al|$), where $\omega_0:=\omega$, and $\omega_1$ is the smallest uncountable ordinal. (\blue{footnote}\footnote{Since the collection $\{\aleph_\al\}_{\al\in\Ordinal}=\{|\omega_\al|\}_{\al\in\Ordinal}$ contains all infinite cardinals (on the class of sets), for any infinite cardinal $\kappa\in\Cardinal\backslash\Natural$, there is an ordinal $\al_{\kappa}\in\Ordinal$ such that $\kappa=\aleph_{\al_{\kappa}}$. This is also \cite[Theorem 8A, p.213]{enderton1977}.})
\end{dfn}

\begin{rmk}[\blue{\index{Continuum hypotheses}{Continuum hypotheses}}]
It is clear that $|\omega|<|\omega_1|\leq 2^{|\omega|}$ (and more generally, for any ordinal $\al$, $|\omega_\al|<|\omega_{\al^+}|\leq 2^{|\omega_\al|}$). The anticipated equality $|\omega_1|=2^{|\omega|}$ is known as the \ul{continuum hypothesis}. A related concern (called \ul{generalized continuum hypothesis}) is whether there exist infinite cardinals $\kappa,\kappa'\in\Cardinal\backslash\Natural$ such that $\kappa<\kappa'<2^\kappa$ (equivalently, whether we have $|\omega_{\al^+}|=2^{|\omega_\al|}$ for all $\al\in\Ordinal$). As a passing note, according to the literature, the equality $|\omega_1|=2^{|\omega|}$ has been shown to be independent of (i.e., to be neither provable nor disprovable using) the ZF/ZFC axioms, and therefore encodes separate axioms for improved set theories.
\end{rmk}

\begin{lmm}[\blue{Initial ordinals represent the infinite cardinals}]
Let $\al,\beta\in\Ordinal_\infty:=\{\eta\in\Ordinal:|\omega|\leq|\eta|\}$ be infinite ordinals. Let $\al'\in\Ordinal$ with $|\al|<|\al'|$. Define {\small $\al_0:=\min\{\gamma\in\Ordinal:\gamma\leq\al,~|\gamma|=|\al|\}$} and {\small $\al_1:=\min\{\gamma\in\Ordinal:\gamma\leq\al',~|\al|<|\gamma|\}$}. The following hold:
\begin{enumerate}[leftmargin=1cm]
\item[(1)] If $\al,\beta$ are initial ordinals, then $\al\leq\beta$ $\iff$ $|\al|\leq|\beta|$. (\blue{footnote}\footnote{That is, on initial ordinals, cardinal ordering coincides with the ordering in $\Ordinal$.})
\item[(2)] $\al_0$ is the unique initial ordinal such that ~$|\al|=|\al_0|$ (i.e., such that $\al\approx\al_0$).
\item[(3)] $\al_1$ is the least initial ordinal such that ~$|\al|<|\al_1|$.
\item[(4)] There is a chain of initial ordinals {\small $\I_\al=\{\al_\ld\}_{\ld\in\Ordinal}:\al_0<\al_1<\al_2<\cdots$} satisfying the following:
\bit
\item[(a)] For each $\ld\in\Ordinal$, (1) above implies $\al_\ld=(\al_\ld)_0$, and we define $\al_{\ld^+}:=(\al_\ld)_1$.
\item[(b)] For each $\ld\in\Ordinal$, there is no initial ordinal $q$ such that ~$\al_\ld<q<\al_{\ld^+}$.
\item[(c)] $\I_{\al}=\I_{\beta}$ $\iff$ $|\al|=|\beta|$ ~(i.e., the chain is unique up to cardinal equivalence).
\eit
\end{enumerate}
\end{lmm}
\begin{proof}
{\flushleft (1)}: Assume $\al,\beta$ are initial ordinals. If $\al<\beta$, then $|\al|<|\beta|$ (since $\beta$ is initial). If $|\al|<|\beta|$ then $\al<\beta$ (since $\beta\leq\al$ implies $|\beta|\leq|\al|$). Since $\al,\beta$ are both initial, we have $\al=\beta$ $\iff$ $|\al|=|\beta|$.

{\flushleft (2)}: By definition, $\al_0\leq\al$ and $|\al_0|=|\al|$, and (the minimality of $\al_0$ means) for any $x\in\Ordinal$, if $x<\al_0$ then $\al<x$ or $|x|\neq |\al|$ (which is the case $\iff$ $|x|<|\al|=|\al_0|$). Hence $\al_0$ is an initial ordinal (which is unique by its definition as a least element).

{\flushleft (3)}: By definition, $\al_1\leq\al'$ and $|\al|<|\al_1|$, and (the minimality of $\al_1$ means) for any $x\in\Ordinal$, if $x<\al_1$ then $\al'<x$ or $|x|\leq|\al|$ (which is the case $\iff$ $|x|\leq|\al|<|\al_1|$). Hence $\al_1$ is an initial ordinal.

{\flushleft (4)}: This follows by transfinite induction based on (1), (2), (3) above.
\end{proof}

\subsection{Transfinite recursion and the arithmetic of ordinals}~\\~
The arithmetic (addition, multiplication, exponentiation, etc) of ordinals is based on their arithmetic as linearly-ordered sets. That is, we can directly add and multiply ordinals as linearly-ordered sets, while noting that addition and multiplication each lack commutativity in general. Since we still have associativity (up to order isomorphism), we can proceed to establish appropriate arithmetic operations using recursion procedures exactly as we did for natural numbers. (\blue{footnote}\footnote{Even in a non-associative number system based on objects (of course other than linearly-ordered sets) whose addition and multiplication each lack both commutativity and associativity, although (unlike in the case of arithmetic in $\Natural$) there will be a general lack of uniqueness in the expression of a given product using sums and products of smaller ordinals, a convenient way to guarantee/ensure uniqueness of arithmetic expressions for such numbers is to use recursion functions (just as we did with the natural numbers) to pick out self-consistent expression patterns/formats/rules for addition, multiplication, and exponentiation of those numbers. Since there is no associativity, we need to fix a conventional rule of inductively ``adding, multiplying, and exponentiating incrementally (say) on the right'' (as opposed to doing so from the left, or mixing left and right).}).

\begin{dfn}[\blue{
\index{Ordinal! element}{Ordinal element},
\index{Ordinal! set}{Ordinal set},
\index{Ordinal! map}{Ordinal map}}]
Let $X$ be a set. An element $x\in X$ is an \ul{ordinal element} of $X$ if $x\subsetneq X$ (not to be confused with an element of $X\cap\Ordinal$). Let $Ord(X)$ denote the set of ordinal elements of $X$. The set $X$ is an \ul{ordinal set} if $Ord(X)=X$, i.e., every element of $X$ is an ordinal element.

A map of sets $f:X\ra Y$ is an \ul{ordinal map} if $f(Ord(X))\subset Ord(f(X))$, i.e., $f(x)\subsetneq f(X)$ for all $x\in Ord(X)$, such that ~$\im f|_x:=\{f(x'):x'\in x\}=f(x)$~ for all $x\in Ord(X)$.
\end{dfn}

In particular, if $\al,\beta\in\Ordinal$ and $f:\al\ra\beta$ is an ordinal map, then ~$\im(f|_{Seg_\al(\al')})=\im(f|_{\al'})=f(\al')$~ for all $\al'\in\al$. That is, the pointwise value $f(\al')\in f(\al)$ associated with $\al'\in\al$ (which by convention is always assumed unless indicated otherwise), can also be viewed as the image value $f(\al')\subsetneq f(\al)$ associated with $\al'\subsetneq\al$ (which happens to be a basic property of the ordinals).

Therefore, it is reasonable to try defining arithmetic operations (addition, multiplication, exponentiation) for ordinals in terms of ordinal maps (including realizing those operations themselves as bivariate maps $\Ordinal\times\Ordinal\ra\Ordinal$ that are ordinal in at least one argument).

\begin{dfn}[\blue{\index{Descending chain}{Descending chain}}]
Let $(A,\leq)$ be a linearly-ordered set. A descending chain in $A$ is a map $f:\omega\ra A$ such that $f(n^+)<f(n)$ for all $n\in\omega$.
\end{dfn}

\begin{thm}[\blue{\cite[Theorem 7B, p.173]{enderton1977}}]\label{OrdTransRec1}
Let $(A,\leq)$ be a linearly-ordered set. $A$ is well-ordered $\iff$ $A$ has no descending chain.
\end{thm}
\begin{proof}
($\Ra$): Assume $A$ is well-ordered. Suppose $f:\omega\ra A$ is a descending chain. Then $f(\omega)\neq\emptyset$ but $\min f(\omega)$ does not exist (a contradiction).

($\La$): Assume $A$ has no descending chains. Suppose $A$ is not well-ordered. Then there exists a nonempty set $\emptyset\neq B\subset A$ such that $\min B$ does not exist, i.e., for each $b\in B$, there exists $b'\in B$ such that $b'<b$. Fix any $b_0\in B$. Let $b_1<b_0$ in $B$, $b_2<b_1$ in $B$, and so on. Then we obtain a descending chain $f:\omega\ra A,~n\mapsto b_n$ (a contradiction).
\end{proof}

\begin{notation}[\textcolor{blue}{Recall: Initial segments}]
Let $(A,\leq_A)$ and $(B,\leq_B)$ be linearly-ordered sets. Write $A\leq_{\txt{seg}}B$ if $A$ is an \ul{initial segment} of $B$ (i.e., $A\subset B$ and $A<b$ for all $b\in B-A$). Similarly, write $A<_{\txt{seg}}B$ if $A$ is a \ul{strict initial segment} of $B$ (i.e., $A\subsetneq B$ and $A<b$ for all $b\in B-A$).
\end{notation}

\begin{dfn}[\blue{
\index{Segment maps}{Segment maps},
\index{Segmentation map}{Segmentation map}}]
Let $(A,\leq)$ be a linearly-ordered set, $\S\subset Sets$ a class of sets, and $Seg_A:=\{Seg_A(a)~|~a\in A\}\subset\P(A)$ the set of segments of $A$. A \ul{segment map} from $A$ to $\S$ is a map of the form $f:Seg_A(a)\ra\S$ for some $a\in A$. The \ul{set of segment maps} from $A$ to $\S$ is
\bea
Seg(A,\S):=\{\txt{maps}~f:Seg_A(a)\ra \S:~a\in A\}\subset\P(A\times\S).\nn
\eea
The segmentation map of $A$ is the order isomorphism
\[
\pi=\pi_A:(A,\leq)\ra(Seg_A,\subseteq),~a\mapsto Seg_A(a).
\]
\end{dfn}

\begin{thm}[\blue{\index{Transfinite recursion formula for! well-ordered sets}{Transfinite recursion} formula for well-ordered sets: \cite[Theorem 7B, p.173]{enderton1977}}]\label{WellOrdTransRec}
Let $(A,\leq)$ be a well-ordered set, $\C,\S\subset Sets$ two classes with $\P(A,\S)\subset\C$ (or just $Seg(A,\S)\subset\C$), and $h:\C\ra \S$ a map. Then there exists a unique map ~$f:A\ra \S,~a\mapsto h(f|_{Seg_A(a)})$. (\blue{footnote}\footnote{In the class $\C$ (the domain of $h:\C\ra\S$) we are treating maps $f\in\C$ as sets (say as their graphs $G(f)$ ) and so $h(f|_{Seg_A(a)})$ is really $h\big(G(f|_{Seg_A(a)})\big)$.})

{\footnotesize\bea\bt
(Seg_A,Seg_A(a))\ar[d,dashed,"f|"]\ar[rrrr,"\pi^{-1}","\cong_o"'] &&&& (A,a)\ar[d,dashed,"f"] \\
 (Seg(A,\S),f|_{Seg_A(a)})\ar[rr,hook] && (\C,f|_{Seg_A(a)})\ar[rr,"h"] && (\S,f(a))
\et~~~~~~~~f\circ\pi^{-1}=h\circ f|\nn
\eea}

\end{thm}
\begin{proof}
If {\footnotesize $f_B:B\leq_{seg} A\ra \S$} is a map (where $B\leq_{seg} A$ means $B$ is an initial segment of $A$), we will say $f_B$ is an $f$-approximation (i.e., an approximation to $f$) if $f_B$ satisfies {\footnotesize $f_B(b)=h(f_B|_{Seg_A(b)})$} for all $b\in B$.

$f_B:B\leq_{seg} A\ra \S$ is unique (in particular, if $f:=f_A$ exists, then it is unique). To see this consider any map $g_B:B\leq_{seg}A\ra \S$ that also satisfies $g_B(b)=h(g_B|_{Seg_A(b)})$ for all $b\in B$. Let $C:=\{c\in B:f_B(c)=g_B(c)\}$. Then $a_0:=\min A\in C$ (so $C\neq\emptyset$), since $f_B(a_0)=h(f_B|_{Seg_A(a_0)})=h(f_B|_\emptyset)=h(\emptyset)=h(g_B|_\emptyset)=h(g_B|_{Seg_A(a_0)})=g_B(a_0)$. Let $c\in C$ be such that $Seg_A(c)\subset C$, i.e., $f_B(c')=g_B(c')$ for all $c'\in Seg_A(c)$ (or equivalently, $f_B|_{Seg_A(c)}=g_B|_{Seg_A(c)}$). Then $f_B(c)=h(f_B|_{Seg_A(c)})=h(g_B|_{Seg_A(c)})=g_B(c)$, and so $c\in C$. This shows $C\subset B$ is inductive and so $C=B$.

The uniqueness of the $f_B$'s implies for any $B,B'\subset A$, we have
\[
f_B|_{B\cap B'}=f_{B'}|_{B\cap B'}.
\]

Let ~$P:=\{\txt{$f$-approximations}~f_B:B\leq_{seg} A\ra \S\}$.~ Also let
\bea
\textstyle F:=\bigcup P=\bigcup_{f_B\in P}f_B:=\bigcup_{f_B\in P}\big\{\big(b,f_B(b)\big)~|~b\in B\big\}~~~~\txt{and}~~~~M:=\dom F:=\bigcup_{f_B\in P}B.\nn
\eea
Then it remains to verify the following properties:
\begin{enumerate}[leftmargin=0.9cm]
\item \ul{$F$ is a map $f_M:M\leqslant_{seg} A\ra \S$.}~ We need to show that all $f$-approximations agree at intersections, i.e., for any $f_B,f_{B'}\in P$, we have ~$f_B|_{B\cap B'}=f_{B'}|_{B\cap B'}$ (which already was proved at the beginning). Equivalently, this means the following set is inductive (and so equals $A$):
\bea
&&S:=\{a\in A~|~(a,x)\in F~\txt{for at most one}~x\in \S\}\nn\\
&&~~~~=\big\{a\in A~|~\txt{the set}~X_a:=\{f_B(a):f_B\in P\}~\txt{contains at most one element}\big\}.\nn
\eea

\item \ul{$F=f_M$ is an $f$-approximation, i.e., $F=f_M\in P$.}~ If $a\in M$, let $a\in B_a\leq_{seg}A$ with $f_{B_a}\in P$. Then
\bea
f_M(a)=f_{B_a}(a)=h(f_{B_a}|_{Seg_A(a)})=h(f_M|_{Seg_A(a)}).\nn
\eea

\item \ul{$M=A$}: (\blue{footnote}\footnote{Alternatively, suppose $M\subsetneq A$. Let $a:=\min(A\backslash M)$. Then $Seg_A(a)\subset M$. So, we get the map $F_a:=F\cup\{(a,h(f_M|_{Seg_A(a)}))\}\in P$ with domain larger than that of $F$ (a contradiction).}). It suffices to show $M\subset A$ is inductive. Set $f_{\{a_0\}}(a_0):=h(\emptyset)=h(f|_{\{a_0\}}|_\emptyset)=h(f|_{\{a_0\}}|_{Seg_A(a_0)})$. Then $a_0\in M$. Let $a\in A$ be such that $Seg_A(a)\subset M$, i.e., for all $a'<a$ we have $a'\in B_{a'}$ for some $B_{a'}\leq_{seg}A$ such that $f_{B_{a'}}\in P$. Define $F_a:=F\cup\{(a,h(f_M|_{Seg_A(a)}))\}$. Then $F_a\in P$. Since $\dom F$ is maximal in $P$, we have $M\cup\{a\}=\dom F_a\subset\dom F=M$, and so $a\in M$.
\item \ul{$F=f_M$ is unique.}~ This part follows from the discussion at the start of the proof. \qedhere
\end{enumerate}
\end{proof}

The following result is a special case of Theorem \ref{WellOrdTransRec}, and so its proof follows the same steps. The only ``simplification'' comes from the fact that if $\al\in\Ordinal$ and $\al'\in\al$, then $Seg_{\al}(\al')=\al'$.
\begin{thm}[\blue{\index{Transfinite recursion formula for! ordinals}{Transfinite recursion} formula for ordinals: \cite[Theorem 8.9, p.228]{goldrei1996}}]\label{OrdTransRec2}
Let $\al\in\Ordinal$, $\S:=Sets$ the class of sets, and $h:\Ordinal\times \S\ra \S$ a map. Then there exists a unique map ~$f:\al\ra \S,~{\al'}\mapsto h({\al'},f|_{\al'})$. (\blue{footnote}\footnote{If $h(\al',f|_{\al'}):=\im(f|_{\al'})$, then $f$ becomes an ordinal map, i.e., $f(\al')=h(\al',f|_{\al'})=\im(f|_{\al'})$ for $\al'\in\al$.})
{\footnotesize\bea\bt
(Seg_\al,\al')\ar[d,dashed,"f|"]\ar[rrrr,"\pi^{-1}","\cong_o"'] &&&& (\al,{\al'})\ar[d,dashed,"f"] \\
 (Seg(\al,\S),f|_{\al'})\ar[rr,hook,"D"] && (\Ordinal\times \S,({\al'},f|_{\al'}))\ar[rr,"h"] && (\S,f({\al'}))
\et~~~~f=f\circ\pi^{-1}=h\circ D\circ f|\nn
\eea}
\end{thm}
\begin{proof}
If $f_\beta:\beta\leq_{seg}\al\ra\S$ is a map (where of course $\beta\leq_{seg}\al$ implies $\beta\in\Ordinal$), we will say $f_\beta$ is an $f$-approximation (i.e., an approximation to $f$) if $f_\beta$ satisfies $f_\beta(\beta')=h(\beta',f_\beta|_{\beta'})$ for all $\beta'\in\beta$. Since $\beta\leq_{seg}\al$ $\iff$ $\beta\in\al^+$, we will use the two expressions interchangeably.

$f_\beta:\beta\leq_{seg}\al\ra \S$ is unique (in particular, if $f:=f_\al$ exists, then it is unique). To see this consider any map $g_\beta:\beta\leq_{seg}\al\ra \S$ that also satisfies $g_\beta(\beta')=h(\beta',g_\beta|_{\beta'})$ for all $\beta'\in\beta$. Let $B:=\{b\in\beta:f_\beta(b)=g_\beta(b)\}$. Then $0\in B$ (so $B\neq\emptyset$), since $f_\beta(0)=h(0,f_\beta|_0)=h(0,\emptyset)=h(0,g_\beta|_0)=g_\beta(0)$. Let $b\in\beta$ be such that $Seg_\beta(b)=b\subset\beta$, i.e., $f_\beta(b')=g_\beta(b')$ for all $b'\in b$ (or equivalently, $f_\beta|_b=g_\beta|_b$). Then $f_\beta(b)=h(b,f_\beta|_b)=h(b,g_\beta|_b)=g_\beta(b)$, and so $b\in B$. This shows $B\subset\beta$ is inductive and so $B=\beta$.

The uniqueness of the $f_\beta$'s implies for any $\gamma,\gamma'\in\Ordinal$, we have
\[
f_\gamma|_{\gamma\cap\gamma'}=f_{\gamma'}|_{\gamma\cap\gamma'}.
\]

Let ~$P:=\{\txt{$f$-approximations}~f_\beta:\beta\leq_{seg}\al\ra \S\}=\{\txt{$f$-approximations}~f_\beta:\beta\in\al^+\ra \S\}$.~ Also let
\bea
\textstyle F:=\bigcup P=\bigcup_{f_\beta\in P}f_\beta:=\bigcup_{f_\beta\in P}\big\{\big(b,f_\beta(b)\big)~|~b\in \beta\big\}~~~~\txt{and}~~~~M:=\dom F:=\bigcup_{f_\beta\in P}\beta.\nn
\eea
Then it remains to verify the following properties:
\begin{enumerate}[leftmargin=0.9cm]
\item \ul{$F$ is a map $f_M:M\leqslant_{seg}\al\ra \S$.}~ We need to show that all $f$-approximations agree at intersections, i.e., for any $f_{\beta},f_{\beta'}\in P$, we have ~$f_{\beta}|_{\beta\cap \beta'}=f_{\beta'}|_{\beta\cap \beta'}$ (which already was proved at the beginning). Equivalently, this means the following set is inductive (and so equals $\al$):
\bea
&&S:=\{a\in\al^+~|~(a,x)\in F~\txt{for at most one}~x\in \S\}\nn\\
&&~~~~=\big\{a\in\al^+~|~\txt{the set}~X_a:=\{f_{\beta}(a):f_{\beta}\in P\}~\txt{contains at most one element}\big\}.\nn
\eea

\item \ul{$F=f_M$ is an $f$-approximation, i.e., $F=f_M\in P$.}~ If $a\in M$, let $a\in\beta_a\leq_{seg}\al$ with $f_{\beta_a}\in P$. Then
\bea
f_M(a)=f_{\beta_a}(a)=h(a,f_{\beta_a}|_a)=h(a,f_M|_a).\nn
\eea

\item \ul{$M=\al$}: (\blue{footnote}\footnote{Alternatively, suppose $M\subsetneq \al$. Let $a:=\min(\al\backslash M)$. Then $a=Seg_\al(a)\subset M$. We get the map $F_a:=F\cup\{(a,h(a,f_M|_a))\}\in P$ with domain larger than that of $F$ (a contradiction).}). It suffices to show $M\subset\al$ is inductive. Set $f_{\{0\}}(0):=h(0,\emptyset)=h(0,f|_{\{0\}}|_0)$. Then $0\in M$. Let $a\in\al$ be such that $a=Seg_\al(a)\subset M$, i.e., for all $a'\in a$ we have $a'\in \beta_{a'}$ for some $\beta_{a'}\leq_{seg}\al$ such that $f_{\beta_{a'}}\in P$. Define $F_a:=F\cup\{(a,h(a,f_M|_a))\}$. Then $F_a\in P$. Since $\dom F$ is maximal in $P$, we have $M\cup\{a\}=\dom F_a\subset\dom F=M$, and so $a\in M$.
\item \ul{$F=f_M$ is unique.}~ This part follows from the discussion at the start of the proof. \qedhere
\end{enumerate}
\end{proof}

We will now define arithmetic operations for the ordinals. Fix $\al\in\Ordinal$. As done before for natural numbers (\blue{footnote}\footnote{If the number system were non-associative, we would need a new ``right-increments'' (as opposed to ``left-increments'', or a mixing of ``left-increments'' and ``right-increments'') convention involved in multiplication via incremental right-addition and in exponentiation via incremental right-multiplication.}), our intent is to transform $\al$ using (i.e., \ul{add to $\al$}, \ul{multiply $\al$ by}, or \ul{exponentiate $\al$ by}) a variable $\beta\in\Ordinal$ in a recursive manner, leading to a desired binary/arithmetic operation (i.e., a map of the form $\Ordinal\times\Ordinal\ra\Ordinal$).

{\flushleft $\bullet$} \ul{Construction of addition}: To define addition, in Theorem \ref{OrdTransRec2}, define $h(\beta,g^\beta)$ for any map ~$g^\beta:\dom g^{\beta}\supset\beta\ra Sets$~ as follows: (\blue{footnote}\footnote{The explicit requirement $h(\beta,g^\beta):=\im(g^\beta|_\beta)$ in the 3rd option is precisely what will make ``adding to $\al$'' $A_\al:\Ordinal\ra\Ordinal$ an ordinal map, i.e., $A_\al(\beta)=h(\beta,A_\al|_\beta):=\im(A_\al|_\beta)$. The same argument holds for multiplication and exponentiation.})
\bea
h(\beta,g^\beta):=\left\{
                     \begin{array}{ll}
                      \al, & \txt{if $\beta=0$} \\
                       h(\beta^-,g^{\beta^-})^+, & \txt{if $\beta$ has a predecessor $\beta^-$} \\
                       \bigcup_{\beta'\in\beta}h(\beta',g^{\beta'}):=\bigcup_{\beta'\in\beta}\im(g^{\beta'}|_{\beta'}), &\txt{if $\beta$ has no predecessor}
                     \end{array}
                   \right.
\eea
Then there exists a unique map ( $f=$ ) $A_\al:\Ordinal\ra\Ordinal$ (``adding to $\al$'') such that $A_\al(\beta)=h(\beta,A_\al|_\beta)$.

{\footnotesize\bea\bt
(Seg_\beta,\beta')\ar[d,dashed,"A_\al|"]\ar[rrrr,"\pi^{-1}","\cong_o"'] &&&& (\beta,{\beta'})\ar[d,dashed,"A_\al"] \\
 (Seg(\beta,Sets),A_\al|_{\beta'})\ar[rr,hook,"D"] && (\Ordinal\times Sets,({\beta'},A_\al|_{\beta'}))\ar[rr,"h"] && (Sets,A_\al({\beta'}))
\et~~~~A_\al=A_\al\circ\pi^{-1}=h\circ D\circ A_\al|\nn
\eea}

We define \ul{addition} $+:\Ordinal\times\Ordinal\ra\Ordinal,~(\al,\beta)\ra\al+\beta$ of ordinals by
\bea
\al+\beta:=A_\al(\beta),~~~~\txt{for all}~~~~\al,\beta\in\Ordinal.
\eea
Therefore, with $g^\beta:=A_\al|_\beta$,
{\footnotesize\begin{align}
&\textstyle\al+\beta:=A_\al(\beta)=
h(\beta,A_\al|_\beta):=\left\{
                     \begin{array}{ll}
                      \al , & \txt{if $\beta=0$} \\
                       h(\beta^-,A_\al|_{\beta^-})^+, & \txt{if $\beta^-$ exists} \\
                       \bigcup_{\beta'\in\beta}h(\beta',A_\al|_{\beta'})
:=\bigcup_{\beta'\in\beta}\im((A_\al|_{\beta'})|_{\beta'})=\bigcup_{\beta'\in\beta}\im(A_\al|_{\beta'}), &\txt{if no $\beta^-$ exists}
                     \end{array}
                   \right\}\nn\\
&~~~~=\left\{
                     \begin{array}{ll}
                      \al , & \txt{if $\beta=0$} \\
                       (\al+\beta^-)^+, & \txt{if $\beta^-$ exists} \\
                       \bigcup_{\beta'\in\beta}~\al+\beta', &\txt{if no $\beta^-$ exists}
                     \end{array}
                   \right\}.
\end{align}}

{\flushleft $\bullet$} \ul{Construction of multiplication}: To define multiplication, we do likewise while making use of addition as defined above. In Theorem \ref{OrdTransRec2}, define $h(\beta,g^\beta)$ for any map ~$g^\beta:\dom g^{\beta}\supset\beta\ra Sets$~ as follows: (\blue{footnote}\footnote{If the number system were non-associative, we would need to write this particular instance/application/use of $A_\al$ as $A^{\txt{right}}_\al$ to emphasize that placing the additive increments on the right (as opposed to the left, or mixing left and right) is essential.})
\bea
h(\beta,g^\beta):=\left\{
                     \begin{array}{ll}
                      0, & \txt{if $\beta=0$} \\
                       A_\al\big(h(\beta^-,g^{\beta^-})\big), & \txt{if $\beta$ has a predecessor $\beta^-$} \\
                       \bigcup_{\beta'\in\beta}h(\beta',g^{\beta'}):=\bigcup_{\beta'\in\beta}\im(g^{\beta'}|_{\beta'}), &\txt{if $\beta$ has no predecessor}
                     \end{array}
                   \right.
\eea
Then there exists a unique map ( $f=$ ) $M_\al:\Ordinal\ra\Ordinal$ (``multiplying $\al$ by'') such that $M_\al(\beta)=h(\beta,M_\al|_\beta)$.

{\footnotesize\bea\bt
(Seg_\beta,\beta')\ar[d,dashed,"M_\al|"]\ar[rrrr,"\pi^{-1}","\cong_o"'] &&&& (\beta,{\beta'})\ar[d,dashed,"M_\al"] \\
 (Seg(\beta,Sets),M_\al|_{\beta'})\ar[rr,hook,"D"] && (\Ordinal\times Sets,({\beta'},M_\al|_{\beta'}))\ar[rr,"h"] && (Sets,M_\al({\beta'}))
\et~~~~M_\al=M_\al\circ\pi^{-1}=h\circ D\circ M_\al|\nn
\eea}

We define \ul{multiplication} $\cdot=\times:\Ordinal\times\Ordinal\ra\Ordinal,~(\al,\beta)\ra\al\times\beta$ of ordinals by
\bea
\al\cdot\beta=\al\times\beta:=M_\al(\beta),~~~~\txt{for all}~~~~\al,\beta\in\Ordinal.
\eea
Therefore, with $g^\beta:=M_\al|_\beta$,
{\footnotesize\begin{align}
&\textstyle\al\times\beta:=M_\al(\beta)=
h(\beta,M_\al|_\beta):=\left\{
                     \begin{array}{ll}
                      0, & \txt{if $\beta=0$} \\
                       A_\al\big(h(\beta^-,M_\al|_{\beta^-})\big), & \txt{if $\beta^-$ exists} \\
                       \bigcup_{\beta'\in\beta}h(\beta',M_\al|_{\beta'})
:=\bigcup_{\beta'\in\beta}\im((M_\al|_{\beta'})|_{\beta'})=\bigcup_{\beta'\in\beta}\im(M_\al|_{\beta'}), &\txt{if no $\beta^-$ exists}
                     \end{array}
                   \right\}\nn\\
&~~~~=\left\{
                     \begin{array}{ll}
                      0, & \txt{if $\beta=0$} \\
                       (\al\times\beta^-)+\al, & \txt{if $\beta^-$ exists} \\
                       \bigcup_{\beta'\in\beta}~\al\times\beta', &\txt{if no $\beta^-$ exists}
                     \end{array}
                   \right\}.
\end{align}}

{\flushleft $\bullet$} \ul{Construction of exponentiation}: We similarly use multiplication to define exponentiation. In Theorem \ref{OrdTransRec2}, define $h(\beta,g^\beta)$ for any map ~$g^\beta:\dom g^{\beta}\supset\beta\ra Sets$~ as follows: (\blue{footnote}\footnote{If the number system were non-associative, we would need to write this instance/application/use of $M_\al$ as $A^{\txt{right}}_\al$ to emphasize that placing the multiplicative increments on the right (as opposed to the left, or mixing left and right) is essential.})
\bea
h(\beta,g^\beta):=\left\{
                     \begin{array}{ll}
                      1, & \txt{if $\beta=0$} \\
                       M_\al\big(h(\beta^-,g^{\beta^-})\big), & \txt{if $\beta$ has a predecessor $\beta^-$} \\
                       \bigcup_{\beta'\in\beta}h(\beta',g^{\beta'}):=\bigcup_{\beta'\in\beta}\im(g^{\beta'}|_{\beta'}), &\txt{if $\beta$ has no predecessor}
                     \end{array}
                   \right.
\eea
Then there exists a unique map ( $f=$ ) $E_\al:\Ordinal\ra\Ordinal$ (``exponentiating $\al$ by'') such that $E_\al(\beta)=h(\beta,E_\al|_\beta)$.

{\footnotesize\bea\bt
(Seg_\beta,\beta')\ar[d,dashed,"E_\al|"]\ar[rrrr,"\pi^{-1}","\cong_o"'] &&&& (\beta,{\beta'})\ar[d,dashed,"E_\al"] \\
 (Seg(\beta,Sets),E_\al|_{\beta'})\ar[rr,hook,"D"] && (\Ordinal\times Sets,({\beta'},E_\al|_{\beta'}))\ar[rr,"h"] && (Sets,E_\al({\beta'}))
\et~~~~E_\al=E_\al\circ\pi^{-1}=h\circ D\circ E_\al|\nn
\eea}

We define \ul{exponentiation} $\Ordinal\times\Ordinal\ra\Ordinal,~(\al,\beta)\ra\al^\beta$ of ordinals by
\bea
\al^\beta:=E_\al(\beta),~~~~\txt{for all}~~~~\al,\beta\in\Ordinal.
\eea
Therefore, with $g^\beta:=E_\al|_\beta$,
{\footnotesize\begin{align}
&\textstyle\al^\beta:=E_\al(\beta)=
h(\beta,E_\al|_\beta):=\left\{
                     \begin{array}{ll}
                      1 , & \txt{if $\beta=0$} \\
                       M_\al\big(h(\beta^-,E_\al|_{\beta^-})\big), & \txt{if $\beta^-$ exists} \\
                       \bigcup_{\beta'\in\beta}h(\beta',E_\al|_{\beta'})
:=\bigcup_{\beta'\in\beta}\im((E_\al|_{\beta'})|_{\beta'})=\bigcup_{\beta'\in\beta}\im(E_\al|_{\beta'}), &\txt{if no $\beta^-$ exists}
                     \end{array}
                   \right\}\nn\\
&~~~~=\left\{
                     \begin{array}{ll}
                      1, & \txt{if $\beta=0$} \\
                       \al^{\beta^-}\times\al, & \txt{if $\beta^-$ exists} \\
                       \bigcup_{\beta'\in\beta}~\al^{\beta'}, &\txt{if no $\beta^-$ exists}
                     \end{array}
                   \right\}.
\end{align}}

A lot of details on the basic properties of ordinal arithmetic (which are not of immediate interest to us) can be found in \cite[Sections 8.3/8.4, pages 218/231]{goldrei1996} and in \cite[pp.227-240]{enderton1977}. In particular,  for any $\al,\al',\beta,\beta'\in\Ordinal$, we have the following:
{\footnotesize\[\bt[column sep=small, row sep=tiny]
\al(\beta+\beta')=\al\beta+\al\beta', &
(\al+\al')\beta\neq\al\beta+\al'\beta, &
\al^{\beta+\beta'}=\al^\beta\al^{\beta'}, &
(\al^\beta)^{\beta'}=\al^{\al\beta}, &
(\al\al')^\beta\neq\al^\beta\al'^\beta \\
(1) & (2) & (3) & (4) & (5)
\et\]
}where (3) and (4) are \cite[Theorem 8R on p.238 and Theorem 8S on p.239]{enderton1977} respectively.

\begin{convention}
Due to the lack of commutativity of multiplication in particular, we use the following convention for expressing repeated sums as products (\blue{footnote}\footnote{This convention is directly based on consistency with our definition for the product of linearly-ordered sets.}): For any $\al,\beta\in\Ordinal$ and $n\in\Natural$,
{\footnotesize\[\bt[column sep=small]
\al+\al:=\al\times2~\neq~2\times\al, &
\cdots, &
\overbrace{\al+\al+\cdots+\al}^{n~\txt{times}}:=\al\times n~\neq~n\times\al, &
\cdots, &
\sum_{\beta'\in\beta}\al:=\al\times\beta~\neq~\beta\times\al.
\et\]}With this convention, the first few countable ordinal numbers are often listed as follows:
{\scriptsize\[
0,~ 1, \cdots,~ \omega,~ \omega+1,~ \omega+2,~ \cdots,~ \omega+\omega=\omega\times 2,~ \cdots,~ \omega\times 3,~ \cdots,~ \omega\times\omega=\omega^2,~ \cdots,~ \omega^3,~ \cdots,~ \omega^\omega,~ \cdots,~ \omega^{\omega^\omega},~ \cdots,~ \omega^{\omega^{\omega^\omega}},~\cdots
\]}where $\omega:=\sup\{0,1,2,\cdots\}$, $\omega^\omega:=\sup\{1,\omega,\omega^2,\omega^3,\cdots\}$, and so on.
\end{convention}

\begin{dfn}[\blue{\index{Epsilon numbers}{Epsilon numbers}, \index{Fixed point}{Fixed point}}]
The \ul{epsilon numbers} are the elements of the subclass $\E:=\{\vep\in\Ordinal:\omega^\vep=\vep\}\subset\Ordinal$ consisting of \ul{fixed points} (\blue{footnote}\footnote{Recall: Given a map $f:X\ra Y$, a point $x\in X\cap Y$ is a \ul{fixed point} of $f$ if $f(x)=x$.}) of the map $E_\omega:\Ordinal\ra\Ordinal,\al\ra\omega^\al$. For $\al\in\Ordinal$, the $\al$th epsilon number is denoted by $\vep_\al$. In particular, $\vep_0$ denotes the smallest epsilon number.
\end{dfn}
\begin{lmm}
The smallest epsilon number satisfies ~$\vep_0=\sup\{\omega,\omega^\omega,\omega^{\omega^\omega},\cdots\}$.
\end{lmm}
\begin{proof}
Let $S:=\{\omega,\omega^\omega,\omega^{\omega^\omega},\cdots\}$ and $\omega^S:=\{\omega^s:s\in S\}=S\backslash\{\omega\}$. Then $\omega^{\sup S}\sr{(a)}{=}\sup(\omega^S)=\sup S$, where step (a) holds by the order-imbedding property of the map $E_\omega:\Ordinal\ra\Ordinal,\al\mapsto\omega^\al$ (along with uniqueness of the supremum). This shows $\sup S\geq\vep_0$. On the other hand, since $\vep_0\geq 1$, by iterating as $\vep_0=\omega^{\vep_0}=\omega^{\omega^{\vep_0}}=\omega^{\omega^{\omega^{\vep_0}}}=\cdots$, we see that $\vep_0\geq s$ for all $s\in S$, and so $\vep_0\geq\sup S$.
\end{proof}

\begin{dfn}[\blue{
\index{Cumulative hierarchy}{Cumulative hierarchy},
\index{Ordinal! rank}{Ordinal rank} of a set}]
The \ul{cumulative hierarchy} is the collection of classes of sets $\{\C_\al\}_{\al\in\Ordinal}$ defined recursively as follows: For $\al\in\Ordinal$,
\[
\textstyle\C_\al:=\left\{
          \begin{array}{ll}
            \emptyset, & \txt{if $\al=0$} \\
            \P(\C_{\al^-}), & \txt{if $\al$ has a predecessor $\al^-$} \\
            \bigcup\{\C_{\al'}:\al'\in\al\}, & \txt{if $\al$ has no predecessor}
          \end{array}
        \right.
\]
Let $A$ be a set. The \ul{ordinal rank} of $A$ is the least ordinal $\al$ such that $A\in\C_\al$, i.e.,
\[
\textstyle\rank_\Ordinal(A):=\min_\Ordinal\{\al\in\Ordinal:A\in\C_\al\}\sr{(s)}{=}\min_\Ordinal\{\al\in\Ordinal:A\in\C_{\al^+}=\P(\C_\al)\}=\min_\Ordinal\{\al\in\Ordinal:A\subset\C_\al\},
\]where step (s) holds because if $\al$ is a limit ordinal (i.e., if $\al$ has no predecessor), then $A\in\C_\al=\bigcup_{\al'\in\al}\C_{\al'}$ $\Ra$ $A\in\C_{\al'}$ for some $\al'\in\al$, and so $\al$ cannot be minimal wrt $A\in\C_\al$.
\end{dfn}

Let $\al\in\Ordinal$. Then a set $A_\al\in\C_\al$ $\iff$ $A_\al\subset\C_{\al'}$ for some $\al'\in\al$, and so the ordinal rank
\[
\rank_\Ordinal(A_\al)=\min_\Ordinal\{\beta\in\Ordinal:A_\al\subset\C_\beta\cap\C_{\al'}~\txt{for some}~\al'\in\al\}=\min\{\al'\in\al:A_\al\subset\C_{\al'}\}\nn
\]exists since $\al\subset\Ordinal$ is well-ordered.

In general, the rank of a set (if it exists) is unique since $\Ordinal$ is well-ordered, and so with $\C:=\bigcup\{\C_\al:\al\in\Ordinal\}$, we have a map ~$\rank_\Ordinal:\C\subset Sets\ra\Ordinal,~A\mapsto\rank_\Ordinal(A)$.~ It is shown in  \cite[Theorem 8.11, p.226]{goldrei1996} that $\C=Sets$, and so every set has a rank (i.e., every set $A$ lies in some member $\C_{\rank_\Ordinal(A)^+}$ of the cumulative hierarchy). The cumulative hierarchy thus gives an explicit way of listing all sets using the ordinals:
\bea
&&\C_0=\emptyset=0,~~~~\C_1=\P(0)=\{0\}=1,~~~~\C_2=\P(1)=\{0,1\}=2,~~~~
\C_3=\P(2)=3\cup\{\{1\}\},\nn\\
&&~~\C_4=\P(3\cup\{\{1\}\})=4\cup\{\txt{$12$ more elements}\},~~\cdots\nn
\eea

%% file: parts/AlgebraNC/AlgebraNC3.tex
\chapter{Groups, Rings, Modules, and Algebras}

For this chapter, if needed, additional sources of reading include for example \cite{gallian2013,artin1991,rotman2010,lang2002,cohn1982,cohn1989}. Many of the concepts for sets that follow (such as binary operations, identity, inverse, etc) can also be applied to classes in general.

\section{Underlying Constructions}
\begin{dfn}[\textcolor{blue}{
\index{Binary operation (Multiplication)}{Binary operation (Multiplication)},
{\index{Multiplicative! set}{Multiplicative set}},
\index{Product of! elements}{Product of elements},
\index{Action!}{Action},
\index{Left! action}{Left action},
\index{Right! action}{Right action},
\index{Center of a multiplicative set}{Center of a multiplicative set},
\index{Commutative! operation}{Commutative operation},
{\index{Commutative! set}{Commutative set}},
\index{Associative! operation}{Associative operation},
{\index{Associative! set (Semigroup)}{Associative set (Semigroup)}}}]~\\~
Let $S$ be a set. A \ul{binary operation} or \ul{multiplication} on $S$, making $S=[S,\cdot]$ a \ul{multiplicative set}, is a map of the form ~$\cdot:S\times S\ra S,~(a,b)\mapsto a\cdot b$ or $ab$. The image $ab$ of $(a,b)$ is called a \ul{product} of $a$ and $b$.

More generally, given any sets $A,B,C$, a map $A\times B\ra C,~(a,b)\mapsto ab$ is also called a \ul{\emph{multiplication}}, a possible justification being that it can be considered a restriction of a map of the form $(A\cup B\cup C)\times(A\cup B\cup C)\ra A\cup B\cup C$. In this sense, a multiplication of the form $A\times B\ra B$ (i.e., with $C=B$) is called an \ul{\emph{action}} (or a \ul{\emph{left action}}) of $A$ on $B$. Similarly, a multiplication of the form $A\times B\ra A$ is called a  \ul{\emph{right action}} of $B$ on $A$.


Let $S=[S,\cdot]$ be a multiplicative set. The \ul{center} of $S$ is $Z(S):=\{c\in S:cs=sc~\txt{for all}~s\in S\}$. The binary operation of $S$ is \ul{commutative}, making $S=\{S,\cdot\}$ a \ul{commutative set}, if $ab=ba$ for any $a,b\in S$ (i.e., if $Z(S)=S$). The binary operation of $S$ is \ul{associative}, making $S=(S,\cdot)$ an \ul{associative set}, if $(ab)c=a(bc)\sr{\txt{written}}{\eqv}abc$ for any $a,b,c\in S$. An associative set is also called a \ul{semigroup}.
\end{dfn}

\begin{prp}[\textcolor{blue}{\index{General associativity}{General associativity}}]
Let $S=(S,\cdot)$ be an associative set. Given any finite list of elements $a_1,...,a_n\in S$, their product $a_1a_2...a_n$ is independent of the way brackets are used.
\end{prp}
\begin{proof}
We proceed by induction on $n$. The cases of $n=1,2$ are trivially true. The case of $n=3$ is true by associativity (Assoc). \emph{Induction hypothesis} (IH): Assume the statement is true for each of the partial products $a_1a_2...a_k$ for all $1\leq k<n$. Then by repeatedly applying the induction hypothesis and associativity, we get
\begin{align}
(a_1a_2\cdots a_{n-1})a_n&\sr{\txt{IH}}{=}((a_1a_2\cdots a_{n-2})a_{n-1})a_n\sr{\txt{Assoc}}{=}(a_1a_2\cdots a_{n-2})(a_{n-1}a_n)\nn\\
&\sr{\txt{IH}}{=}((a_1a_2\cdots a_{n-3}) a_{n-2})(a_{n-1}a_n)\sr{\txt{Assoc}}{=}(a_1a_2\cdots a_{n-3}) (a_{n-2}(a_{n-1}a_n))\nn\\
&\sr{\txt{IH}}{=}(a_1a_2\cdots a_{n-3}) (a_{n-2}a_{n-1}a_n)\nn\\
&~~\vdots\hspace{2cm} \vdots\hspace{2cm}\vdots\nn\\
&\sr{\txt{IH}}{=}a_1(a_2a_3\cdots a_n).\nn \qedhere
\end{align}
\end{proof}
To see how the above proof works more explicitly, observe that for a fixed positive integer $N$, all brackets can be removed by progressively partitioning each $n\leq N$ and then applying the induction hypothesis followed by associativity (IH+Assoc) as in Table \ref{EqCoeffTab1}:
\begin{center}
\adjustbox{scale=0.8}{
\begin{minipage}{23cm}
\begin{table}[H]
  \centering
\begin{tabular}{l|l|l}
  $n$ & partial product & bracket removal \\
\hline
3 &  $a_1a_2a_3$                    &$(a_1a_2)a_3=a_1(a_2a_3)$\\
 &                                    &$~~~~(2,1)=(1,2),~~~~\txt{(Assoc)}$\\&&\\
4 &  $a_1a_2a_3a_4$                 &$(a_1a_2a_3)a_4=(a_1a_2)(a_3a_4)=a_1(a_2a_3a_4)$\\
 &                                    &$~~~~(3,1)=(2,2)=(1,3),~~~~\txt{(IH+Assoc)}$\\&&\\
5 &  $a_1a_2a_3a_4a_5$              &$(a_1a_2a_3a_4)a_5=(a_1a_2a_3)(a_4a_5)=(a_1a_2)(a_3a_4a_5)=a_1(a_2a_3a_4a_5)$\\
 &                                    &$~~~~(4,1)=(3,2)=(2,3)=(1,4),~~~~\txt{(IH+Assoc)}$\\&&\\
6 &  $a_1a_2a_3a_4a_5a_6$            &$~~~~(5,1)=(4,2)=(3,3)=(2,4)=(1,5),~~~~\txt{(IH+Assoc)}$\\&&\\
 &  $~~\vdots~~\hspace{1cm}$         &\hspace{1cm}~~\vdots~~\hspace{2cm}~~\vdots\\&&\\
N &  $a_1a_2a_3a_4\cdots a_N$         &$~~~~(N-1,1)=(N-2,2)=(N-3,3)=\cdots=(3,N-3)=(2,N-2)=(1,N-1)~~~~\txt{(IH+Assoc)}$\\&&\\
N+1 &  $a_1a_2a_3a_4\cdots a_{N+1}$   &$~~~~(N,1)=(N-1,2)=(N-2,3)=\cdots=(3,N-2)=(2,N-1)=(1,N)~~~~\txt{(IH+Assoc)}$\\
  \hline
\end{tabular}
  \caption{}\label{EqCoeffTab1}
\end{table}
\end{minipage}}
\end{center}

\begin{dfn}[\textcolor{blue}{
\index{Left! identity}{Left identity},
\index{Right! identity}{Right identity},
\index{Identity!}{Identity},
{\index{Identity! set}{Identity set}},
{\index{Associative! identity set (Monoid)}{Associative identity set (Monoid)}}}]
Let $S=[S,\cdot]$ be a multiplicative set. An element $l\in S$ is a \ul{left identity} if $la=a$ for all $a\in S$. An element $r\in S$ is a \ul{right identity} if $ar=a$ for all $a\in S$. An element $e\in S$ is an \ul{identity}, making $S=[S,\cdot,e]$ an \ul{identity set}, if it is both a left identity and a right identity (i.e., if $ea=ae=a$ for all $a\in S$).

An \ul{associative identity set} $S=(S,\cdot,e)$ is also called a \ul{monoid}.
\end{dfn}

\begin{note}
(i) The identity $e$ in an identity set $S=[S,\cdot,e]$ is unique. Indeed, if $e'\in S$ is any identity, then ~$e=ee'=e'e=e'$. (ii) Moreover, if both a left identity $l$ and a right identity $r$ exist, then $l=lr=r$ (and so give an identity $e:=l=r$).
\end{note}

\begin{notation}[\blue{Symbols for identity elements}]
For all identity sets $S=[S,\cdot,e]$ that we will encounter, the identity element $e=e_S$ will often also be denoted by $1=1_S$, or by $0=0_S$, or by $id=id_S$, depending on the context. Concerning the present random choice of the symbol $e$ to denote \ul{identity}, the reader may consider $e$ to stand for ``\ul{exceptional element}''.
\end{notation}

\begin{dfn}[\textcolor{blue}{
\index{Left! infinity}{Left infinity},
\index{Right! infinity}{Right infinity},
\index{Infinity!}{Infinity},
{\index{Infinity! set}{Infinity set}}}]
Let $S=[S,\cdot]$ be a multiplicative set. An element $l\in S$ is a \ul{left infinity} if $la=l$ for all $a\in S$. An element $r\in S$ is a \ul{right infinity} if $ar=r$ for all $a\in S$. An element $\infty\in S$ is an \ul{infinity}, making $S=[S,\cdot,\infty]$ an \ul{infinity set}, if it is both a left infinity and a right infinity (i.e., if $\infty a=a\infty=\infty$ for all $a\in S$).
\end{dfn}

\begin{note}
(i) The infinity $\infty$ in an infinity set $S=[S,\cdot,\infty]$ is unique. Indeed, if $\infty'\in S$ is any infinity, then ~$\infty=\infty\infty'=\infty'\infty=\infty'$. (ii) Moreover, if both a left infinity $l$ and a right identity $r$ exist, then $l=lr=r$ (and so give an infinity $\infty:=l=r$).
\end{note}

\begin{dfn}[\textcolor{blue}{
\index{Left! inverse}{Left inverse},
\index{Right! inverse}{Right inverse},
\index{Inverse!}{Inverse},
\index{Invertible element}{Invertible element},
{\index{Group}{Group}},
\index{Inverse! operation}{Inverse operation},
{\index{Abelian group}{Abelian group}}}]
Let $S=[S,\cdot,e]$ be an identity set and $a\in S$. An element $a_l\in S$ is a \ul{left inverse} of $a$ if $a_la=e$. An element $a_r\in S$ is a \ul{right inverse} of $a$ if $aa_r=e$. An element $a^{-1}\in S$ is an \ul{inverse} of $a$ (making $a$ \ul{invertible}) if it is both a left inverse of $a$ and a right inverse of $a$ (i.e., if $a^{-1}a=aa^{-1}=e$).

A \ul{group} $G=\big[G,\cdot,e,(~)^{-1}\big]$ is an identity set $G=[G,\cdot,e]$ in which every element $a\in G$ has an inverse $a^{-1}$ (equivalent to the existence of a map $G\ra G$, $a\mapsto a^{-1}$ called \ul{inverse operation}). An \ul{abelian group} $A=\langle A,\cdot,e,(~)^{-1}\rangle$ is a commutative associative group (i.e., a commutative associative identity set in which every element has an inverse). Here, the limiters $\langle,\rangle$ jointly express commutativity and associativity.
\end{dfn}

\begin{note}
(i) Inverses in an associative identity set (monoid) $S=(S,\cdot,e)$ are unique. Indeed, if $s\in S$ and $s',s''\in S$ are inverses of $s$, then ~$s'=s'e=s'(ss'')=(s's)s''=es''=s''$. (ii) Moreover, if $s\in S$ has both a left inverse $u$ and a right inverse $v$, then $u=ue=u(sv)=(us)v=ev=v$ (and so give an inverse $s^{-1}:=u=v$).
\end{note}

\begin{lmm}[\textcolor{blue}{\index{Inverse! of a product}{Inverse of a product} in an associative identity set (monoid)}]
Let $S=(S,\cdot,e)$ be an associative identity set. If $a,b\in S$ are invertible elements, then $ab$ is invertible, and $(ab)^{-1}=b^{-1}a^{-1}$.
\end{lmm}
\begin{proof}
Observe that $(ab)b^{-1}a^{-1}=a(bb^{-1})a^{-1}=aea^{-1}=aa^{-1}=e$ and $b^{-1}a^{-1}(ab)=b^{-1}(a^{-1}a)b=b^{-1}ea=a^{-1}a=e$. Hence, $(ab)b^{-1}a^{-1}=b^{-1}a^{-1}(ab)=e$.
\end{proof}

Because the most interesting groups are associative groups, unless specified otherwise, we will henceforth assume all groups are associative groups. Note that a binary operation does not have to be surjective. However, for any multiplicative set $S=[S,\cdot]$ with a left identity $l\in S$, or right identity $r\in S$, the binary operation $\cdot:S\times S\ra S$ is clearly surjective.

\begin{dfn}[\textcolor{blue}{
\index{Addition! in a multiplicative set}{Addition in a multiplicative set},
\index{One (or Unity)}{One (or Unity)},
\index{Zero!}{Zero},
{\index{Ring! on a multiplicative set}{Ring on a multiplicative set}},
\index{Subtraction in a ring}{Subtraction in a ring},
\index{Distributivity in a ring}{Distributivity in a ring},
\index{Unital! ring}{Unital ring},
\index{Associative! ring}{Associative ring},
\index{Commutative! ring}{Commutative ring},
\index{Center of a ring}{Center of a ring}}]
Let $R=[R,\cdot]$ be a multiplicative set. A binary operation $+:R\times R\ra R$, $(a,b)\mapsto a+b$ is called an \ul{addition}, making $R=[R,\cdot,+]$ a \ul{ring}, if the following hold:
\begin{enumerate}[leftmargin=1cm]
\item If the multiplicative set $[R,\cdot]$ has an identity element, it is denoted by $1$ and called ``\ul{one}'' or ``\ul{unity}'' or ``\ul{multiplicative identity}''.
\item $\langle R,+\rangle$ is an abelian group $\langle R,+,0,-\rangle$ with (i) identity element denoted by $0$ and called ``\ul{zero}'' or ``\ul{additive identity}'', and (ii) the inverse operation $R\ra R$, $a\ra -a:=a^{-1}$ equivalently expressed as a binary operation $R\times R\ra R,~(a,b)\mapsto a-b:=a+(-b)$ called \ul{subtraction}.
\item \ul{Distributivity}:~ $a(b+c)=ab+ac$ ~~and~~ $(a+b)c=ac+bc$,~~ for all $a,b,c\in R$.
\end{enumerate}
That is, a ring $R=[R,\cdot,+]$ is a set $R$ with two binary operations, namely, \ul{multiplication} $\cdot:R\times R\ra R$, $(a,b)\mapsto ab$, and \ul{addition} $+:R\times R\ra R$, $(a,b)\mapsto a+b$ such that the following hold.
\begin{enumerate}[leftmargin=1cm]
\item $[R,\cdot]$ is a multiplicative set with identity element (if it exists) denoted by $1$ and called ``\ul{one}'' or ``\ul{unity}'' or ``\ul{multiplicative identity}''.
\item $\langle R,+\rangle$ is an abelian group $\langle R,+,0,-\rangle$ with (i) identity element denoted by $0$ and called ``\ul{zero}'' or ``\ul{additive identity}'', and (ii) the inverse operation $R\ra R$, $a\ra -a:=a^{-1}$ equivalently expressed as a binary operation $R\times R\ra R,~(a,b)\mapsto a-b:=a+(-b)$ called \ul{subtraction}.
\item $a(b+c)=ab+ac$ ~~and~~ $(a+b)c=ac+bc$,~~ for all $a,b,c\in R$. Here, we say multiplication is \ul{distributive over addition}.
\end{enumerate}

A ring is a \ul{unital ring} if it has a multiplicative identity (i.e., a unity). A ring is an \ul{associative ring} if multiplication is associative. A ring is a \ul{commutative ring} if multiplication is commutative. The \ul{center} of a ring $R$ is its center wrt multiplication, i.e., $Z(R):=\{c\in R:cr=rc~\txt{for all}~r\in R\}$.
\end{dfn}

\begin{rmks}[\textcolor{blue}{Basic deductions}]
Let $R=[R,\cdot,+]$ be a ring. Then for all $r,s\in R$, we have (i) $r0=0r=0$ (because $0+0=0$), and (ii) $(-r)s=-(rs)$, because $rs+(-r)s=(r+(-r))s=0s=0$.
\end{rmks}

\begin{dfn}[\textcolor{blue}{
\index{Scalar multiplication (Scaling)}{Scalar multiplication (Scaling)},
{\index{Module! (on an abelian group) over a ring}{Module (on an abelian group) over a ring}},
\index{Module}{Module},
\index{Left! module}{Left module},
\index{Right! module}{Right module},
\index{Additive! abelian group}{Additive abelian group},
\index{Subtraction in a module}{Subtraction in a module},
\index{Distributivity in a module}{Distributivity in a module},
\index{Associative! module}{Associative module},
\index{Unital! module}{Unital module}}]
Let $R=[R,\cdot,+,0,-]$ be a ring and $M=\langle M,+,0,-\rangle=\langle M,+_M,0_M,-_M\rangle$ an abelian group. A map
\bea
{}_s\cdot:R\times M\ra M,~(r,m)\mapsto r\cdot_sm=rm~~~~(~\txt{resp.}~~\cdot_s:M\times R\ra M,~(r,m)\mapsto r\cdot_sm=rm~)\nn
\eea
is called a \ul{scalar multiplication} or \ul{scaling}, making $M={}_RM$ a \ul{(left) $R$-module} (resp. making $M=M_R$ a \ul{right $R$-module}) if the following hold. (\blue{footnote}\footnote{Because of left-right symmetry (or commutativity) of addition both in the module and in the ring, (i) the distinction between a left-module and a right-module is necessary only when associativity is considered, after which (ii) the distinction is again necessary only if the ring is \ul{not} commutative.})
\bit
\item[(1)] The abelian group operation $+=+_M$ is now called \ul{addition}, which makes $M=\langle M,+,0,-\rangle$ an ``\ul{additive abelian group}'' with (i) identity element $0=0_M$ called ``\ul{zero}'' or ``\ul{additive identity}'', and (ii) the inverse operation $M\ra M$, $m\mapsto -m:=m^{-1}$ equivalently expressed as a binary operation $M\times M\ra M$, $(m,m')\mapsto m-m':=m+(-m')$ called \ul{subtraction}.
\item[(2)] Distributivity: $(r+r')m=rm+r'm$~~ and ~~ $r(m+m')=rm+rm'$ for all $r,r'\in R$, $m,m'\in M$.~\\~
~~~~( resp. $m(r+r')=mr+mr'$~~ and ~~ $(m+m')r=mr+m'r$ for all $r,r'\in R$, $m,m'\in M$.  )
\eit
That is, given a ring $R=[R,\cdot,+,0,-]$, a (left) $R$-module $M={}_RM$ ( resp. a right $R$-module $M=M_R$ ) is a set $M$ with two operations, namely, \ul{scalar multiplication}
\bea
\cdot_s:R\times M\ra M,~(r,m)\mapsto rm~~~~(\txt{resp. }~{}_s\cdot:M\times R\ra M,~(m,r)\mapsto mr~)\nn
\eea
and \ul{addition} $+=+_M:M\times M\ra M$, $(m,m')\mapsto m+m'$, such that the following hold.
\begin{enumerate}
\item[(1)] $\langle M,+\rangle=\langle M,+_M\rangle$ is an abelian group $\langle M,+,0,-\rangle=\langle M,+_M,0_M,-_M\rangle$ with (i) identity element denoted by $0=0_M$ and called ``\ul{zero}'' or ``\ul{additive identity}'', and (ii) the inverse operation $M\ra M$, $m\mapsto -m=-_Mm:=m^{-1}$ equivalently expressed as a binary operation $M\times M\ra M$, $(m,m')\mapsto m-m':=m+(-m')$ called \ul{subtraction}.
\item[(2)] \ul{Distributivity}: $(r+r')m=rm+r'm$~~ and ~~ $r(m+m')=rm+rm'$ for all $r,r'\in R$, $m,m'\in M$.~\\~
~~~~( resp. $m(r+r')=mr+mr'$~~ and ~~ $(m+m')r=mr+m'r$ for all $r,r'\in R$, $m,m'\in M$.  )
\end{enumerate}
The $R$-module $M={}_RM$ is \ul{associative} if, in addition to (1),(2) above, the condition (3) below also holds.
\bit
\item[(3)] Associativity: $r(r'm)=(rr')m$, ~~ for all $r,r'\in R$, $m\in M$.~\\~
~~~~( resp. $(mr)r'=m(rr')$, ~~ for all $r,r'\in R$, $m\in M$. )
\eit
The $R$-module $M={}_RM$ is \ul{unital} if, in addition to (1),(2) above, the condition (4) below also holds.
\bit
\item[(4)] Unitality: The ring $R$ is unital, and $1m=m$ for all $m\in M$ (resp. $m1=m$ for all $m\in M$).
\eit
\end{dfn}

\begin{rmks}[\textcolor{blue}{Basic deductions}]
Let $M={}_RM$. Then for all $r\in R$, $m\in M$, we have (i) $0_Rm=r0_M=0_M$ (because $0_R+0_R=0_R$, $0_M+0_M=0_M$) and (ii) $(-r)m=-(rm)$, because $rm+(-r)m=(r+(-r))m=0_Rm=0_M$.
\end{rmks}
Consequently, it is no longer necessary to distinguish between $0_R$ and $0_M$, or between subtraction in $R$ and subtraction in $M$.

\begin{lmm}[\textcolor{blue}{\index{Abelian means additive}{Abelian means additive}}]
An abelian group is precisely a (unital associative) $\Integer$-module.
\end{lmm}
\begin{proof}
($\Ra$) If $G=\langle G,\cdot,e,(~)^{-1}\rangle=\langle G,+_G,0_G,-_G\rangle$ is an abelian group, where for every $g,h\in G$ we define
\bea
g+_Gh:=gh,~~~~g-_Gh:=g+_G(-_Gh):=g+_Gh^{-1}=gh^{-1},\nn
\eea
then the map $\Integer\times G\ra G$, $(n,g)\mapsto n\ast g:=g^n$ is a scalar multiplication since
\bea
&&(n+m)\ast g=g^{n+m}=g^ng^m=(n\ast g)(m\ast g)=n\ast g+_Gm\ast g,\nn\\
&&n\ast(g+_Gh)=n\ast(gh)=(gh)^n=g^nh^n=(n\ast g)(n\ast h)=n\ast g+_Gn\ast h,\nn\\
&&(nm)\ast g=g^{nm}=(g^n)^{m}=(n\ast g)^m=m\ast(n\ast g),\nn\\
&&1\ast g=g^{1}=g.\nn
\eea
It is now clearly neither necessary to distinguish (notation-wise) between addition in $\Integer$ and the operation ~$\cdot~=~+_G$~ of ~$G$, nor to distinguish (notation-wise) between multiplication in $\Integer$ and scalar multiplication $\ast$.

($\La$) Conversely, it is clear that a (unital associative) $\Integer$-module is an abelian group.
\end{proof}

\begin{dfn}[\textcolor{blue}{
{\index{Bimodule}{Bimodule}},
\index{Associative! bimodule}{Associative bimodule},
\index{Left-associative bimodule}{Left-associative bimodule},
\index{Right-associative bimodule}{Right-associative bimodule},
\index{Unital! bimodule}{Unital bimodule},
\index{Left-unital bimodule}{Left-unital bimodule},
\index{Right-unital bimodule}{Right-unital bimodule}}]
Let $R=[R,\cdot_R]$, $S=[S,\cdot_S]$ be rings and $M=\langle M,+\rangle$ an abelian group. Then $M$ is an \ul{RS-bimodule} (written $M={}_RM_S$) if the following holds.
\bit
\item[(a)] $M$ is both a left $R$-module and a right $S$-module (i.e., $M={}_RM=M_S$).
\eit

The RS-bimodule is \ul{associative} if, in addition to (a) above, the condition (b) below also holds.
\bit
\item[(b)] Associativity: $_RM$ and $M_S$ are both associative, and ~$r(ms)=(rm)s$ for all $r\in R$, $m\in M$, $s\in S$.~\\~
~~~~(More terminology: ${}_RM_S$ is \ul{left-associative} if ${}_RM$ is associative. ${}_RM_S$ is \ul{right-associative} if $M_S$ is associative.)
\eit

The RS-bimodule is \ul{unital} if, in addition to (a) above, the condition (c) below also holds.
\bit
\item[(c)] Unitality: $_RM$ and $M_S$ are both unital (i.e., $R,S$ are unital and $1_Rm=m=m1_S$ for all $m\in M$).~\\~
~~~~(More terminology: ${}_RM_S$ is \ul{left-unital} if ${}_RM$ is unital. ${}_RM_S$ is \ul{right-unital} if $M_S$ is unital.)
\eit
\end{dfn}

\begin{note}[\blue{Exclusion of $0$}]
If $K$ is a ring or module, we will often write $K\backslash\{0\}$ simply as $K\backslash 0$ (which is technically incorrect/unacceptable in classical set theory).
\end{note}
\begin{note}[\blue{A module over a commutative ring is a bimodule}]
If $R$ is a commutative ring and $M={}_RM$ an $R$-module (whether associative or not) with scalar multiplication {\small $R\times M\ra M,~(r,m)\mapsto rm$}, then $M$ is also a right $R$-module (hence an $R$-bimodule {\small $M={}_RM=M_R={}_RM_R$}) with scalar multiplication {\small $\ast:M\times R\ra M,~(m,r)\mapsto m\ast r:=rm$} (i.e., the same operation), because for all $r,r'\in R$, $m,m'\in M$,
\bit
\item[] $m\ast(r+r')=(r+r')m=rm+r'm=m\ast r+m\ast r'$,
\item[] $(m+m')\ast r=r(m+m')=rm+rm'=m\ast r+m'\ast r$, ~and if associative,
\item[] $m\ast(rr')=(rr')m=(r'r)m=r'(rm)=(m\ast r)\ast r'$.
\eit
\end{note}
Because of the above note, we can refer to a module over a commutative ring as a commutative module.
\begin{dfn}[\blue{\index{Commutative! module}{Commutative module}}]
An $R$-module $M$ is a \ul{commutative $R$-module} if $R$ is a commutative ring. In particular, an $R$-algebra (to be defined next) is an example of a commutative $R$-module.
\end{dfn}

\begin{dfn}[\blue{\index{Opposite! ring}{Opposite ring}}]
Let $R$ be a ring. The \ul{opposite ring} $R^{op}:=(R,\ast)$ of $R$ is the same set $R$ as a ring wrt a new multiplication (in terms of multiplication $R\times R\ra R,~(r,s)\mapsto rs$ in $R$) given by
\bea
\ast:R\times R\ra R,~(r,s)\mapsto r\ast s:=sr.\nn
\eea
Note that $R^{op}$-modules ${}_{R^{op}}M$ are precisely right $R$-modules $M_R$. Similarly, $R$-modules ${}_RM$ are precisely right $R^{op}$-modules $M_{R^{op}}$.
\end{dfn}

\begin{dfn}[\textcolor{blue}{
{\index{Algebra! (on a ring) over a commutative ring}{Algebra (on a ring) over a commutative ring}},
\index{Unital! algebra}{Unital algebra},
\index{Associative! algebra}{Associative algebra},
\index{Graded!-commutative algebra}{Graded-commutative algebra}}]
Let $R$ be a commutative ring and $A$ a ring. Then $A$ is an \ul{$R$-algebra} if the following hold.
\bit
\item[(i)] $A$ is an $R$-module (i.e., $A={}_RA$).
\item[(ii)] $R$ is \ul{central} in $A$, in the sense that ~$(ra)a'=(ar)a'$, for all $r\in R$, $a,a'\in A$.
\eit
An $R$-algebra $A={}_RA$ is \ul{unital} if (i) $A$ is a unital ring and (ii) $_RA$ is a unital module.  An $R$-algebra $A={}_RA$ is \ul{associative} if (i) $A$ is an associative ring and (ii) $_RA$ is an associative module.

An $R$-algebra $_RA$ is called \ul{graded-commutative} if there exists a map ~$\al:A\times A\ra R$~ such that
\bea
ab=\al(a,b)ba,~~~~\txt{for all}~~a,b\in A,\nn
\eea
(in which case, we say $A$ is $\al$-commutative). Note that ~$\al(a,b)\al(b,a)ab=ab$,~ for all ~$a,b\in A$.
\end{dfn}
It is worth noting that a ring (as defined earlier) is precisely a $\Integer$-algebra.

\begin{notation*}[\textcolor{blue}{
\index{Multiplication! of sets}{Multiplication of sets},
\index{Addition! of sets}{Addition of sets}}]
Let $G$ be a group and $M$ an $R$-module. If $A,A'\subset G$, $C\subset R$, and $N,N'\subset M$, we define
\bea
&&\textstyle AA':=\{aa':a\in A,a'\in A'\}\subset G,~~~~N+N':=\{n+n':~n\in N,~n'\in N'\}\subset M,\nn\\
&&\textstyle CN:=\left\{\txt{finite sums}~\sum\limits_{i=1}^kc_in_i:c_i\in C,~n_i\in N,~k\geq 1\right\}\subset M,\nn
\eea
where if $0\in C$, then we can more conveniently write (with a.e.f. meaning ``all except finitely many'')
\[
\textstyle CN=\{\sum_{n\in N}c_nn:c_n\in C,~c_n=0~~\txt{for a.e.f}~~n\in N\}.
\]
\end{notation*}
If $M$ is an $R$-module, then in general equality does not hold in $RM\subseteq M$ (i.e., scalar multiplication is not surjective in general). With the following simplification, we will then assume $RM=M$ by default.

\begin{assumption*}[\textcolor{blue}{Major simplification}]
Observe that a ring $R$ is a bimodule over itself, i.e., $R={}_RR=R_R={}_RR_R$. Also, the most interesting modules (rings and algebras included) are associative unital modules. Therefore, unless specified otherwise, we will henceforth assume all modules (rings and algebras included) are associative unital modules.
\end{assumption*}

Note however that nonassociative structures are also relevant. For example, we have Lie algebras defined as follows:
\begin{dfn}[\blue{
\index{Lie! bracket}{Lie bracket},
\index{Lie! algebra}{Lie algebra},
\index{Graded! Lie bracket}{Graded Lie bracket},
\index{Graded! Lie algebra}{Graded Lie algebra}}]
Let $R$ be a commutative ring and $L$ an $R$-module. A \ul{Lie bracket} on $L$ (making $L$ a \ul{Lie $R$-algebra}) is a binary operation $\cdot:L\times L\ra L,~(a,b)\ra a\cdot b$ satisfying the following: For all $a,b,c\in L$ and $r\in R$, we have
\begin{enumerate}
\item\ul{bilinearity}: ~$(ra)\cdot b=r(a\cdot b)$ ~and~ $a\cdot (rb)=r(a\cdot b)$
\item\ul{biadditivity}: ~$(a+b)\cdot c=a\cdot c+b\cdot c$ ~and~ $a\cdot(b+c)=a\cdot b+a\cdot c$
\item\ul{graded commutativity (antisymmetry)}: ~$a\cdot b=-b\cdot a$
\item\ul{Jacobi-associativity (or Jacobi identity)}: ~$a\cdot(b\cdot c)=(a\cdot b)\cdot c+b\cdot (a\cdot c)$
\end{enumerate}
Similarly, given a map ~$\deg:R\cup L\ra\Integer$~ with $\deg|_R=\txt{constant}=0$,~ a \ul{graded Lie bracket} on $L$ (making $L$ a \ul{graded Lie $R$-algebra}) is a binary operation $\cdot:L\times L\ra L,~(a,b)\ra a\cdot b$ satisfying the following: For all $a,b,c\in L$ and $r\in R$,
\begin{enumerate}
\item\ul{bilinearity}: ~$(ra)\cdot b=r(a\cdot b)$ ~and~ $a\cdot (rb)=r(a\cdot b)$
\item\ul{biadditivity}: ~$(a+b)\cdot c=a\cdot c+b\cdot c$ ~and~ $a\cdot(b+c)=a\cdot b+a\cdot c$
\item\ul{graded commutativity}: ~$a\cdot b=-(-1)^{\deg(a)\deg(b)}b\cdot a$
\item\ul{graded Jacobi-associativity}: ~$a\cdot(b\cdot c)=(a\cdot b)\cdot c+(-1)^{\deg(a)\deg(b)}b\cdot (a\cdot c)$
\end{enumerate}
If $\deg=\txt{constant}=0$, the associated graded Lie algebra becomes the Lie algebra.
\end{dfn}

\begin{dfn}[\textcolor{blue}{
\index{Multiplicative! inverse in a ring}{Multiplicative inverse in a ring},
\index{Left! unit}{Left unit},
\index{Right! unit}{Right unit},
\index{Unit (Unitary element)}{Unit (Unitary element)},
{\index{Group! of units}{Group of units}},
{\index{Division! (or Divisible) ring}{Division (or Divisible) ring}},
{\index{Field}{Field}}}]
Let $R=(R,\cdot,1,+,0,-)$ be a ring. An element $u\in R$ is a \ul{left unit} if it has a \ul{multiplicative left inverse}, i.e., $Ru=R$ or there exists $l\in R$ such that $lu=1$ (resp. a \ul{right unit} if it has a \ul{multiplicative right inverse}, i.e., $uR=R$ or there exists $r\in R$ such that $ur=1$). The set of left units in $R$ is denoted by $U_l(R)$ and the set of right units by $U_r(R)$.

An element $u\in R$ is a \ul{unit} if it is both a left unit and a right unit (i.e., it has a \ul{multiplicative inverse}, i.e. there exists $v\in R$, denoted by $u^{-1}$, such that $uv=vu=1$). The set of all units $U(R):=\big\{u\in R:~\txt{$u^{-1}$ exists}\big\}=U_l(R)\cap U_r(R)$ is a group (called \ul{group of units} of $R$) with multiplication as operation and identity element $1$, i.e., $U(R)=\big(U(R),\cdot,1\big)$.

A ring $D=(D,\cdot,1,+,0,-)$ is a \ul{division ring} (or \ul{divisible ring}) if $U(D)=D\backslash\{0\}$ (i.e., if every nonzero element of $D$ is a unit, or equivalently, $rR=Rr=R$ for all every $r\in R\backslash\{0\}$). A commutative division ring $F=\langle F,\cdot,1,+,0,-\rangle$ is called a \ul{field}.
\end{dfn}

\begin{notation}[\blue{Symbols for fields}]
We will often denote an arbitrary field by $k$. Also, in extended discussions (e.g., whole subsections/sections) involving mostly classical fields like $\Real$ or $\Complex$, such fields will be denoted by $\mathbb{F}$.
\end{notation}

\begin{dfn}[\textcolor{blue}{\index{Division! (or Divisible) module}{Division (or Divisible) module}}]
An $R$-module $M$ is a \ul{division module} (or \ul{divisible module}) if for any $r\in R\backslash 0$, we have $M\subset rM$ (i.e., $rM=M$).
\end{dfn}

Once we have discussed \emph{localization} (to come later), it will not be difficult to create/invent various examples of division rings and division modules.

\begin{dfn}[\textcolor{blue}{
\index{Vector space over a field}{Vector space over a field},
\index{Linear! transformation}{Linear transformation}}]
Let $F$ be a field. An $F$-module $V=\langle V,+,0,-\rangle={}_FV$ is called a \ul{vector space over $F$} (or \ul{$F$-space}, or \ul{$F$-vector space}). A map between $F$-vector spaces $f:V\ra V'$ is \ul{$F$-linear} (and hence called an \ul{$F$-linear transformation}) if for all $a\in F$, $v,u\in V$, we have (i) $f(av)=af(v)$ and (ii) $f(v+u)=f(v)+f(u)$.
\end{dfn}
Linear transformations are examples of structure preserving maps (generally called ``homomorphisms'') between groups, rings, modules, algebras, etc that we will study in subsequent sections.

\begin{dfn}[\textcolor{blue}{
\index{Polynomial ring (Polynomial algebra}{Polynomial ring (Polynomial algebra)},
\index{Commutative! polynomial ring}{Commutative polynomial ring}}]
Let $R$ be a commutative ring. A ``\ul{polynomial $R$-ring} $R\langle X\rangle$ over a \ul{variable set} (or set of variables) $X$'' is understood to be the \ul{smallest ring containing} (also called the \ul{ring generated by}) a subset $X\subset A$ of an $R$-algebra $A={}_RA=A_R={}_RA_R$ that is variable in the sense that $A$ need only belong to a desired/given class of $R$-algebras (e.g., commutative $R$-algebras, graded $R$-algebras, graded-commutative $R$-algebras, etc).

The polynomial $R$-ring $R\langle X\rangle$ is \ul{commutative}, denoted by $R[X]$, if $X\subset A$ for $A$ in the class of commutative $R$-algebras. In particular, the polynomial $R$-ring in one variable is the set
\bea
\textstyle R\langle x\rangle=R[x]:=\big\{f:=f_0+f_1x+\cdots+f_nx^n~|~f_i\in R,n\in\Natural\big\}=\bigcup_{n\in\Natural}(R+Rx+\cdots+Rx^n)\nn
\eea
along with addition and multiplication given, with ~$(f\cdot g)_i:=\sum_{j+k=i}f_jg_k$, respectively by
\bea
\textstyle\left(\sum f_ix^i\right)+\left(\sum g_ix^i\right):=\sum(f_i+g_i)x^i,~~~~\left(\sum f_ix^i\right)\cdot\left(\sum g_ix^i\right):=\sum (f\cdot g)_ix^i,\nn
\eea
as well as ``zero'' and ``one'' given respectively by ~$0_{R[x]}:= 0_{R}$ ~and ~$1_{R[x]}:=1_{R}$.

For any set of variables $X$, we can define the \ul{commutative polynomial $R$-ring} over $X$ recursively as follows: Consider a well-ordering {\small $X=\{x_\al\}_{\al\in\Ordinal}$}, and let {\small $X_\al:=\{x_{\al'}\}_{\al'\in\al}$}. Then {\small $R[X]=\bigcup_{\al\in\Ordinal}R[X_\al]$}, where
\bea
&&\textstyle R[X_\al]=\left\{
           \begin{array}{ll}
             (R[X_{\al^-}])[x_\al], & \txt{if $\al$ has a predecessor $\al^-$} \\
             \bigcup_{\al'\in\al} R[X_{\al'}], & \txt{if $\al$ has no predecessor}
           \end{array}
         \right\},\nn\\
&&\textstyle (R[X_{\al^-}])[x_\al]=\sum_{i=0}^\infty R[X_{\al^-}]x_\al^i:=\bigcup_{i\in\Natural}\Big(R[X_{\al^-}]+R[X_{\al^-}]x_\al+\cdots+R[X_{\al^-}]x_\al^i\Big).\nn
\eea
Also, for any set of variables $X$, with
\bea
X^{\ast 0}:=\{1_R\},~~~~X^{\ast 1}:=X,~~~~X^{\ast i}:=\{x_1x_2\cdots x_i~|~x_1,...,x_i\in X\}~~\txt{for}~~i\geq 1,\nn
\eea
the general \ul{polynomial $R$-ring} over $X$ can be written as
\bea
\textstyle R\langle X\rangle=\sum_{i=0}^\infty RX^{\ast i}:=\bigcup_{i=0}^\infty\Big(RX^{\ast 0}+RX^{\ast 1}+\cdots+RX^{\ast i}\Big).\nn
\eea
\end{dfn}

\begin{dfn}[\textcolor{blue}{
\index{Matrix! groups/rings}{Matrix groups/rings},
\index{Matrix! addition}{Matrix addition},
\index{Matrix! multiplication}{Matrix multiplication},
\index{Matrix! scalar multiplication}{Matrix scalar multiplication,
\index{Matrix! module}{Matrix module},
\index{Matrix! algebra}{Matrix algebra},
\index{Matrix! transpose}{Matrix transpose}}}]
Let $R$ be a ring and $0\neq m,n\in\Natural$. The \ul{group/ring of $m\times n$ matrices} (\blue{footnote}\footnote{As usual, ~$m\times n=\{0,1,\cdots,n-1\}\times\{0,1,\cdots,m-1\}\approx nm=\{0,1,\cdots,mn-1\}$ ~is the cartesian product of sets.}) over $R$ is
\bea
M_{(m,n)}(R):=R^{m\times n}=\{\txt{maps}~f:m\times n\ra R\}=\left\{[f_{ij}]=[f_{ij}]_{i\in m,j\in n}~|~f_{ij}\in R\right\},~~~~M_n(R):=M_{(n,n)}(R),\nn
\eea
with \ul{addition} $+:R^{m\times n}\times R^{m\times n}\ra R^{m\times n}$, \ul{multiplication} $\cdot:R^{m\times t}\times R^{t\times n}\ra R^{m\times n}$, ``\ul{zero}'' in the abelian group $R^{m\times n}$,  and ``\ul{one}'' in the ring $R^{n\times n}$ given respectively by
\bea
\textstyle[f_{ij}]+[g_{ij}]:=[f_{ij}+g_{ij}],~~~~[f_{ij}]\cdot[g_{ij}]:=\big[\sum_kf_{ik}g_{kj}\big],~~~~0_{M(R)}:= [0_R],~~~~1_{M(R)}:=[1_{R}\delta_{ij}].\nn
\eea
If $S$ is another ring such that $R$ is an $S$-module, then $M_{(m,n)}(R):=R^{m\times n}$ is an $S$-module with \ul{scalar multiplication} $S\times M_{(m,n)}(R)\ra M_{(m,n)}(R)$ given by
\bea
\al\cdot[f_{ij}]:=[\al f_{ij}],~~~~\txt{for}~~~~\al\in S,~~[f_{ij}]\in M_{(m,n)}(R),\nn
\eea
and so if $S$ is commutative (which makes $R$ an $S$-algebra), then the matrix ring $M_n(R):=R^{n\times n}$ clearly becomes an $S$-algebra.
\end{dfn}
The matrix relations above are often written in expanded form as follows: \ul{Matrix addition} is given by
\[
\left[
  \begin{array}{cccc}
   f_{11}  & f_{12} & \cdots & f_{1n} \\
   f_{21}  & f_{22} & \cdots & f_{2n} \\
    \vdots & \vdots &  & \vdots \\
   f_{m1}  & f_{m2} & \cdots & f_{mn} \\
  \end{array}
\right]+\left[
  \begin{array}{cccc}
   g_{11}  & g_{12} & \cdots & g_{1n} \\
   g_{21}  & f_{22} & \cdots & g_{2n} \\
    \vdots & \vdots &  & \vdots \\
   g_{m1}  & g_{m2} & \cdots & g_{mn} \\
  \end{array}
\right]:=
\left[
  \begin{array}{cccc}
   f_{11}+g_{11}  & f_{21}+g_{21} & \cdots & f_{1n}+g_{1n} \\
   f_{21}+g_{21}  & f_{22}+g_{22} & \cdots & f_{2n}+g_{2n} \\
    \vdots & \vdots &  & \vdots \\
   f_{m1}+g_{m1}  & f_{2n}+g_{2n} & \cdots & f_{mn}+g_{mn} \\
  \end{array}
\right],
\]
\ul{matrix multiplication} is given by
\[
\left[
  \begin{array}{cccc}
   f_{11}  & f_{12} & \cdots & f_{1t} \\
   f_{21}  & f_{22} & \cdots & f_{2t} \\
    \vdots & \vdots &  & \vdots \\
   f_{m1}  & f_{m2} & \cdots & f_{mt} \\
  \end{array}
\right]\cdot\left[
  \begin{array}{cccc}
   g_{11}  & g_{12} & \cdots & g_{1n} \\
   g_{21}  & f_{22} & \cdots & g_{2n} \\
    \vdots & \vdots &  & \vdots \\
   g_{t1}  & g_{t2} & \cdots & g_{tn} \\
  \end{array}
\right]:=
\left[
  \begin{array}{cccc}
   \sum f_{1k}g_{k1}  & \sum f_{1k}g_{k2} & \cdots & \sum f_{1k}g_{kn} \\
   \sum f_{2k}g_{k1}  & \sum f_{2k}g_{k2} & \cdots & \sum f_{2k}g_{kn} \\
    \vdots & \vdots &  & \vdots \\
   \sum f_{mk}g_{k1}  & \sum f_{mk}g_{k2} & \cdots & \sum f_{mk}g_{kn} \\
  \end{array}
\right],
\]
and \ul{matrix scalar multiplication} is given by
\[
\al\cdot\left[
  \begin{array}{cccc}
   f_{11}  & f_{12} & \cdots & f_{1n} \\
   f_{21}  & f_{22} & \cdots & f_{2n} \\
    \vdots & \vdots &  & \vdots \\
   f_{m1}  & f_{m2} & \cdots & f_{mn} \\
  \end{array}
\right]:=\left[
  \begin{array}{cccc}
   \al f_{11}  & \al f_{12} & \cdots & \al f_{1n} \\
   \al f_{21}  & \al f_{22} & \cdots & \al f_{2n} \\
    \vdots & \vdots &  & \vdots \\
   \al f_{m1}  & \al f_{m2} & \cdots & \al f_{mn} \\
  \end{array}
\right].
\]
The \ul{transpose} of a matrix $f=[f_{ij}]\in M_{(m,n)}(R)$ is the matrix
$f^T=[f^T_{ij}]:=[f_{ji}]\in M_{(n,m)}(R)$. That is,
\[
\left[
  \begin{array}{cccc}
   f_{11}  & f_{12} & \cdots & f_{1n} \\
   f_{21}  & f_{22} & \cdots & f_{2n} \\
    \vdots & \vdots &  & \vdots \\
   f_{m1}  & f_{m2} & \cdots & f_{mn} \\
  \end{array}
\right]^T:=
\left[
  \begin{array}{cccc}
   f_{11}  & f_{21} & \cdots & f_{m1} \\
   f_{12}  & f_{22} & \cdots & f_{m2} \\
    \vdots & \vdots &  & \vdots \\
   f_{1n}  & f_{2n} & \cdots & f_{mn} \\
  \end{array}
\right].
\]

\begin{examples}
Common examples of rings $R$ (which also provide examples of groups, modules, and algebras) include subsets of integers $R_1\subset\Integer$, of rationals $R_2\subset\Rational$, of real numbers $R_3\subset\Real$, of complex numbers $R_4\subset\Complex$, of polynomials $R_5\subset R_i[x_1,...,x_n]$ for $i\in\{1,\cdots,4\}$, and of matrices $R_6\subset M_n(R_k)$ for $k=1,...,5$.
\end{examples}

\section{Commutative Ring Concepts (a concretizing digression)}
\begin{dfn}[\textcolor{blue}{
\index{Division! in a commutative ring}{Division in a commutative ring},
\index{Divisor (factor) of an element}{Divisor (factor) of an element},
\index{Multiple of an element}{Multiple of an element},
\index{Zerodivisor}{{Zerodivisor}},
\index{Integral! domain (ID)}{{Integral domain (ID)}}}]
Let $R=\langle R,\cdot,1,+,0,-\rangle$ be a commutative ring and $m\in R$. Then an element $d\in R$ divides $m$ (or $d$ is a \ul{divisor or factor} of $m$, or $m$ is a \ul{multiple} of $d$), written $d|m$, if $m=dm'$ for some $m'\in R$ (equivalently, $m\in Rd$, or $Rm\subset Rd$).

An element $z\in R$ is a \ul{zerodivisor} if (i) $z\neq 0$ and (ii) $0=zz'$ for some nonzero element $z'\neq 0$. An \ul{integral domain} (ID) is a commutative ring $R$ with no zerodivisors (equivalently, for all $a,b\in R$, if $ab=0$ then $a=0$ or $b=0$).
\end{dfn}

\begin{dfn}[\textcolor{blue}{
\index{Division! in an ID}{Division in an ID},
\index{Irreducible! element}{{Irreducible element}},
\index{Prime! element}{{Prime element}}}]
Let $R$ be an ID and $x\in R$. Then $x$ is \ul{irreducible} or ``\ul{weakly prime}'' if $x$ is a nonunit that cannot be written as a product of nonunits, i.e., (i) $x\not\in U(R)$ and (ii) if $x=ab$ then $a\in U(R)$ or $b\in U(R)$. (Since $0=00$, it is clear that an irreducible element is nonzero.)

Let $Irred(R)$ denote the set of irreducible elements of $R$.

An element $p\in R$ is \ul{prime} or ``\ul{strongly irreducible}'' if (i) $p\not\in \{0\}\cup U(R)$ and (ii) if $p|ab$ then $p|a$ or $p|b$. Let $Prime(R)\subset Irred(R)$ denote the set of prime elements of $R$. (Note that a prime $p$ is irreducible, because $p=ab$ $\Ra$ $p|ab$, $\Ra$ $p|a$ or $p|b$, $\Ra$ $a=\al p$ or $b=\beta p$, $\Ra$ $p=(\al p)b$ or $p=a(\beta p)$, $\Ra$ $1=\al b$ or $1=a\beta$, $\Ra$ $a\in U(R)$ or $b\in U(R)$.)
\end{dfn}

\begin{dfn}[\textcolor{blue}{
\index{Fully! factorable element in an ID}{Fully factorable element in an ID},
\index{Fully! factorable domain (FFD)}{{Fully factorable domain (FFD)}}}]
Let $R$ be an ID and $x\in R\backslash\big(\{0\}\cup U(R)\big)$. The nonzero nonunit $x$ is \ul{fully factorable}, or \ul{fully factorizable}, if we can write
\bea
x=x_1x_2\cdots x_k,~~~~\txt{for a \ul{finite} set of irreducibles}~~\{x_1,...,x_k\}\subset Irred(R).\nn
\eea
An ID, $R$, is a \ul{fully factorable domain (FFD)} if every $x\in R\backslash\big(\{0\}\cup U(R)\big)$ is fully factorable.
\end{dfn}

\begin{dfn}[\textcolor{blue}{
\index{Associates in an ID}{{Associates} in an ID},
{\index{Unique factorization domain (UFD)}{Unique factorization domain (UFD)}}}]
Let $R$ be an ID. Elements $a,b\in R$ are \ul{associates} (denoted $a\sim b$) if $a=bu$, for some $u\in U(R)$. The ID, $R$, is a \ul{unique factorization domain} (UFD) if the following hold (where (1) means $R$ is a FFD).
\begin{enumerate}
\item \emph{Nonzero non-units can be factored into irreducibles:} If $x\in R\backslash(0\cup U(R))$, then we can write
\bea
x=x_1x_2\cdots x_k,~~~~\txt{for some (not necessarily distinct)}~~~~x_1,...,x_k\in Irred(R).\nn
\eea

\item \emph{The factoring into irreducibles is unique up to associates:} If $x\in R\backslash(0\cup U(R))$ is given by
\bea
x=x_1x_2\cdots x_k=x'_1x'_2\cdots x'_{k'},~~~~\txt{for}~~~~x_i,x_i'\in Irred(R),\nn
\eea
then $k=k'$, and after reordering the irreducible factors, $x_i\sim x'_i$ (i.e., associates) for all $i=1,...,k$.
\end{enumerate}
\end{dfn}
Note that the relation $\sim$ on the ID, $R$, above is an equivalence relation, with the equivalence class of $a\in R$ given by $[a]=U(R)a=\{ua:u\in U(R)\}$.

\begin{dfn}[\textcolor{blue}{
\index{Degree map}{{Degree} map},
\index{Division! algorithm}{{Division algorithm}},
\index{Euclidean domain (ED)}{{Euclidean domain (ED)}},
\index{Quotient! (after division)}{Quotient (after division)},
\index{Remainder (after division)}{Remainder (after division)}}]
Let $R$ be an $ID$. A \ul{degree map} on $R$ is a map of the form $\sigma:R\backslash\{0\}\ra\Natural=\{0,1,2,\cdots\}$. We say a degree map $\sigma$ induces the \ul{division algorithm} on $R$ (or the division algorithm holds in $R$ with respect to $\sigma$), making $R=\langle R,\sigma\rangle$ a \ul{Euclidean domain (ED)}, if the following holds: For any $a,b\in R$ with $a\neq 0$, there exist $q,r\in R$ such that
\[
b=qa+r=\txt{quotient}\times a+\txt{remainder},~~~~\txt{with}~~r=0~~\txt{or}~~\sigma(r)<\sigma(a),
\]
where $q$ is called the \ul{quotient} and $r$ the \ul{remainder} (of $b$ wrt $a$). In other words, given any nonzero element $a\in R\backslash\{0\}$, if we define $[R,a):=\{0\}\cup\{r\in R:\sigma(r)<\sigma(a)\}$, then
\[
R=Ra+[R,a).\nn
\]
\end{dfn}

\begin{dfn}[\textcolor{blue}{
\index{Greatest! common divisor (gcd) in an ID}{Greatest common divisor ({gcd}) in an ID},
\index{Least! common multiple (lcm) in an ID}{Least common multiple ({lcm}) in an ID},
\index{Relatively! prime elements in an ID}{{Relatively prime elements} in an ID}}]\label{GcdLcmDef}
Let $R$ be an ID and $a_1,...,a_n\in R$. An element $d\in R$ is a \ul{greatest common divisor} (gcd) of $a_1,....,a_n$, denoted by $d\in gcd(a_1,...,a_n)$, if the following hold.
\bit
\item[(a)] $d$ is a common divisor of $a_1,...,a_n$:~ That is, $d|a_i$ for all $i$. Equivalently, $Ra_1+\cdots+Ra_n\subset Rd$.
\item[(b)] Every common divisor of $a_1,...,a_n$ is a divisor of $d$:~ That is, for any $d'\in R$, if $d'|a_i$ for all $i$, then $d'|d$. Equivalently, if $Ra_1+\cdots+Ra_n\subset Rd'$, then $Rd\subset Rd'$. In particular,
\bea
\label{GCDFormulaEq}\textstyle Ra_1+\cdots+Ra_n\subset Rd\subset\bigcap\limits_{d'\in R}\left\{Rd':~Ra_1+\cdots+Ra_n\subset Rd'\right\}.
\eea
\eit
Similarly, an element $m\in R$ is a \ul{least common multiple} (lcm) of $a_1,....,a_n$, denoted by $m\in lcm(a_1,...,a_n)$, if the following hold.
\bit
\item[(a)] $m$ is a common multiple of $a_1,...,a_n$:~ That is, $a_i|m$ for all $i$. Equivalently, $Rm\subset Ra_1\cap\cdots\cap Ra_n$.
\item[(b)] Every common multiple of $a_1,...,a_n$ is a multiple of $m$:~ That is, for any $m'\in R$, if $a_i|m'$ for all $i$, then $m|m'$. Equivalently, if $Rm'\subset Ra_1\cap\cdots\cap Ra_n$, then $Rm'\subset Rm$. In particular,
\bea
\label{LCMFormulaEq}\textstyle Ra_1\cap\cdots\cap Ra_n\supset Rm\supset\bigcup\limits_{m'\in R}\left\{Rm':~Rm'\subset Ra_1\cap\cdots\cap Ra_n\right\}.
\eea
\eit
We say $a_1,\cdots,a_n$ are \ul{relatively prime} if $1\in gcd(a_1,...,a_n)$, i.e., $1$ is a gcd of $a_1,...,a_n$.
\end{dfn}
To see why condition (b) is reasonable in the definition of the gcd above, observe that if $R$ is a FFD (fully factorable domain) and {\small $d=x_1x_2\cdots x_k=\prod\{x_1,...,x_k\}_{rp}$}, {\small $d'=x_1'x_2'\cdots x_{k'}'=\prod\{x_1',...,x_{k'}'\}_{rp}$} (for irreducibles $x_i,x_i'$) are any two common divisors of $a_1,...,a_n$, then another common divisor of $a_1,...,a_n$ is
\bea
\label{GcdFormula}\textstyle d''=d\Cup d':=\prod\left(\{x_1,x_2,...,x_k\}_{rp}\Cup\{x_1',x_2',...,x_{k'}'\}_{rp}\right),
\eea
where (i) the subscript $rp$ in $\{...\}_{rp}$ means if a factor $x_i$ repeats as $x_i=x_{i_0}\sim x_{i_1}\sim\cdots\sim x_{i_{k_i}}$ then it is treated as the pair $(x_i,k_i)$ in order to note its \ul{multiplicity or power} $k_i$, and (ii) $\Cup$ is the union taken in such a way that for any repeated factor only the \ul{higher} of its two possible multiplicities is retained, up to associates (i.e., associates represent the same element). It is clear that $d|d''$ and $d'|d''$. That is, the set
\[
CD(a_1,...,a_n):=\{\txt{common divisors $d$ of $a_1,...,a_n$}\}
\]
is a directed set with respect to the order relation given by ``$d\leq d'$ if $d|d'$''.

Similarly, to see why condition (b) is reasonable in the definition of the lcm above, observe that if $R$ is a FFD (fully factorable domain) and {\small $m=x_1x_2\cdots x_k=\prod\{x_1,...,x_k\}_{rp}$}, {\small $m'=x_1'x_2'\cdots x_{k'}'=\prod\{x_1',...,x_{k'}'\}_{rp}$} (for irreducibles $x_i,x_i'$) are any two common multiples of $a_1,...,a_n$, then another common multiple is
\bea
\label{LcmFormula}\textstyle m''=m\Cap m':=\prod\left(\{x_1,x_2,...,x_k\}_{rp}\Cap\{x_1',x_2',...,x_{k'}'\}_{rp}\right),
\eea
where (i) the subscript $rp$ in $\{...\}_{rp}$ has the same meaning as before, and (ii) $\Cap$ is the intersection taken in such a way that for any repeated factor only the \ul{lower} of its two possible multiplicities is retained, up to associates. It is clear that $m''|m$ and $m''|m'$. That is, the set
\[
CM(a_1,...,a_n):=\{\txt{common multiples $m$ of $a_1,...,a_n$}\}
\]
is a directed set with respect to the order relation given by ``$m\leq m'$ if $m'|m$''.

The following is an immediate consequence of unique factorization:

\begin{lmm}[\textcolor{blue}{In a UFD, unique gcd's and lcm's exist}]
Let $R$ be a UFD and $a_1,...,a_n\in R$. Then $gcd(a_1,...,a_n)$ and $lcm(a_1,...,a_n)$ exist and are unique up to associates.
\end{lmm}
\begin{proof}
It is enough to do what follows for $n=2$ and then conclude by induction, since by associativity
\[
gcd(a_1,\cdots,a_n)=gcd(gcd(a_1,\cdots,a_{n-1}),a_n)~~\txt{and}~~lcm(a_1,\cdots,a_n)=lcm(lcm(a_1,\cdots,a_{n-1}),a_n).
\]
Factorize each $a_i$ into irreducibles as $a_i=x_{i1}x_{i2}\cdots x_{in_i}=\prod\{x_{i1},x_{i2},\cdots, x_{in_i}\}_{rp}$. Then
\[
\textstyle gcd(a_1,...,a_n)=\prod\big(\Cap_{i=1}^n\{x_{i1},x_{i2},\cdots, x_{in_i}\}_{rp}\big),~~~~lcm(a_1,...,a_n)=\prod\big(\Cup_{i=1}^n\{x_{i1},x_{i2},\cdots, x_{in_i}\}_{rp}\big),
\]
where the expressions/operations $\{...\}_{rp}$, $\Cup$, $\Cap$ have the same meaning as in (\ref{GcdFormula}) and (\ref{LcmFormula}).
\end{proof}

\section{Homomorphisms: Structure preserving maps}
A structured object, or ``structure'' for short, (such as a group, ring, module, algebra, or category to come later) by definition consists of a list of auxiliary objects and (binary) operations involving the auxiliary objects.

\begin{dfn}[\textcolor{blue}{
\index{Morphism in a class of objects}{Morphism in a class of objects},
\index{Endomorphism!}{Endomorphism},
\index{Homomorphism! in a class of classes}{Homomorphism in a class of classes}, 
\index{Monomorphism}{Monomorphism},
\index{Epimorphism}{Epimorphism}}]\label{MorDef0}
Let $\C$ be a class. A \ul{morphism in $\C$} (or \ul{$\C$-morphism}) is a \ul{rule} $f:A\ra B$, between objects $A,B\in\C$, that respects/preserves the class structure of $\C$ (in a specific/precise sense). A $\C$-morphism of the form $f:A\ra A$ (i.e., from a class to itself) is called a \ul{$\C$-endomorphism}.

Let $\C$ be a class of classes and $A,B\in\C$. A map $f:A\ra B$ is a \ul{$\C$-homomorphism} if it is a $\C$-morphism. An injective (resp. surjective) $\C$-homomorphism is called a \ul{$\C$-monomorphism} (resp. \ul{$\C$-epimorphism}). (\blue{footnote}\footnote{These concepts will be defined later for a general category based on a class whose objects are not necessarily maps.})
\end{dfn}
Classes of sets we have seen include the class \index{Class of! groups (Groups)}{\ul{Groups} (Gp)} of groups, the class \index{Class of! rings (Rings)}{\ul{Rings}} (Rng) of rings, the class \index{Class of! modules (Modules or Mod)}{\ul{Modules} (Mod)} of modules, and the class \index{Class of! algebras (Algebras or Alg)}{\ul{Algebras} (Alg)} of algebras. Given rings $R$ and $S$, the class of \ul{left $R$-modules} is denoted by \index{Class of! left $R$-modules ($R$-mod or $_R$Mod)}{$R$-Mod or $_R$Mod}, the class of \ul{right $R$-modules} by \index{Class of! right $R$-modules (Mod-$R$ or Mod$_R$)}{Mod-$R$ or Mod$_R$}, and the class of \ul{$RS$-bimodules} by \index{Class of! RS-bimodules ($R$-mod-S or $_R$Mod$_S$)}{$R$-Mod-$S$ or $_R$Mod$_S$}. The \index{Class of! abelian groups ($\Integer$-mod or Ab)}{class of \ul{abelian groups}}, $\Integer$-Mod, is also denoted by Ab. The class of \ul{topological spaces}, denoted by  \index{Class of! topological spaces (Top)}{Top}, will be considered in chapter \ref{GeomAnaI}.

\begin{dfn}[\textcolor{blue}{
\index{Group! homomorphism}{Group homomorphism},
\index{Kernel of! a group homomorphism}{Kernel of a group homomorphism}}]\label{MorDef1}
A map between groups $f:G\ra H$ is a \ul{group homomorphism} if $f(gg')=f(g)f(g')$, for all $g,g'\in G$. The \ul{kernel} of $f$ is ~$\ker f:=f^{-1}(e_H)=\{g\in G:f(g)=e_H\}$.
\end{dfn}

\begin{dfn}[\textcolor{blue}{
\index{Ring! homomorphism}{Ring homomorphism},
\index{Kernel of! a ring homomorphism}{Kernel of a ring homomorphism}}]\label{MorDef2}
A map between rings $f:R\ra S$ is a \ul{ring homomorphism} if (i) $f(r+r')=f(r)+f(r')$, (ii) $f(rr')=f(r)f(r')$, and (iii) $f(1_R)=1_S$, for all $r,r'\in R$. The \ul{kernel} of $f$ is ~$\ker f:=f^{-1}(0_S)=\{r\in R:f(r)=0_S\}$.
\end{dfn}

\begin{dfn}[\textcolor{blue}{
\index{Module! homomorphism}{Module homomorphism},
\index{Kernel of! a module homomorphism}{Kernel of a module homomorphism},
\index{Linear! map}{Linear map},
\index{Linear! transformation}{Linear transformation}}]\label{MorDef3}
A map between $R$-modules $f:M\ra N$ is an \ul{$R$-module homomorphism} (or \ul{$R$-homomorphism}) if (i) $f(m+m')=f(m)+f(m')$, and (ii) $f(rm)=rf(m)$, for all $m,m'\in M$, $r\in R$. The \ul{kernel} of $f$ is ~$\ker f:=f^{-1}(0_N)=\{m\in M:f(m)=0_N\}$.

Let $M,N$ be $R$-modules and $A\subset M$ a subset. A map $f:A\subset M\ra N$ is \ul{$R$-linear} if it satisfies (i) and (ii) above for all $m,m'\in A$, $r\in R$ such that $m+m',rm\in A$. Under certain conditions, ensuring well-definedness, a linear map $f:A\subset M\ra N$ immediately extends to an $R$-homomorphism in the form $h_f:RA\subset M\ra N,~\sum r_aa\mapsto \sum r_af(a)$. In particular, this is always true whenever $R$ is a field, in which case every $R$-linear map $f:A\subset M\ra N$ is equivalent to an $R$-homomorphism $h_f:RA\subset M\ra N$ (a \ul{linear transformation} between the $R$-vector spaces $RA$ and $N$).
\end{dfn}

\begin{dfn}[\textcolor{blue}{
\index{Algebra! homomorphism}{Algebra homomorphism},
\index{Kernel of! an algebra homomorphism}{Kernel of an algebra homomorphism}}]\label{MorDef4}
A map between $R$-algebras $f:A\ra B$ is an $R$-algebra \ul{homomorphism} if (i) $f(a+a')=f(a)+f(a')$, (ii) $f(aa')=f(a)f(a')$, (iii) $f(1_A)=1_B$, and (iv) $f(ra)=rf(a)$, for all $a,a'\in A$, $r\in R$. The \ul{kernel} of $f$ is ~$\ker f:=f^{-1}(0_B)=\{a\in A:f(a)=0_B\}$.
\end{dfn}

\section{Isomorphisms: Reversible structure preserving maps}
\begin{dfn}[\textcolor{blue}{
\index{Isomorphism}{Isomorphism},
\index{Isomorphic! classes}{Isomorphic classes},
\index{Inverse! of an isomorphism}{Inverse of an isomorphism},
\index{Automorphism}{Automorphism}}]
Let $\C$ be a class of classes and $A,B\in\C$. A $\C$-homomorphism $f:A\ra B$ is a \ul{$\C$-isomorphism} (i.e., an \ul{isomorphism in $\C$}) if (i) $f$ is a bijection and (ii) the bijection-inverse $\ol{f}:B\ra A$ is a $\C$-homomorphism. In this case, the bijection inverse $\ol{f}:B\ra A$ is called the \ul{$\C$-inverse} of $f$ and written $f^{-1}:B\ra A$. If a $\C$-isomorphism $f:A\ra B$ exists, we say $A$ and $B$ are \ul{$\C$-isomorphic classes} (written $A\cong B$). A $\C$-isomorphism of the form $f:A\ra A$ (i.e., from a class to itself) is called a \ul{$\C$-automorphism}.
\end{dfn}

\begin{dfn}[\textcolor{blue}{
\index{Group! isomorphism}{Group isomorphism},
\index{Isomorphic! groups}{Isomorphic groups}}]
A group isomorphism is an isomorphism in the class Groups. If a group isomorphism $f:G\ra H$ exists, we say $G,H$ are isomorphic, written $G\cong H$.
\end{dfn}

\begin{dfn}[\textcolor{blue}{
\index{Ring! isomorphism}{Ring isomorphism},
\index{Isomorphic! rings}{Isomorphic rings}}]
A ring isomorphism is an isomorphism in the class Rings. If a ring isomorphism $f:R\ra S$ exists, we say $R,S$ are isomorphic, written $R\cong S$.
\end{dfn}

\begin{dfn}[\textcolor{blue}{
\index{Module! isomorphism}{Module isomorphism},
\index{Isomorphic! modules}{Isomorphic modules}}]
An $R$-module isomorphism (or $R$-isomorphism) is an isomorphism in the class $R$-Mod. If an $R$-module isomorphism $f:M\ra N$ exists, we say $M,N$ are $R$-isomorphic, written $M\cong_RN$.
\end{dfn}

\begin{dfn}[\textcolor{blue}{
\index{Algebra! isomorphism}{Algebra isomorphism},
\index{Isomorphic! algebras}{Isomorphic algebras}}]
An $R$-algebra isomorphism is an isomorphism in the class $R$-Alg. If an $R$-algebra isomorphism $f:A\ra B$ exists, we say $A,B$ are $R$-isomorphic, written $A\cong_RB$.
\end{dfn}

\begin{thm}
(i) A group homomorphism $f:G\ra H$ is a group isomorphism $\iff$ it is bijective $\iff$ $\ker f=\{0\}$ and $f(G)=H$. (ii) The same is true if ``group'' is replaced with ``ring'', or ``$R$-module'', or ``$R$-algebra''. (\magenta{footnote}\footnote{\magenta{Caution}: Note however that this is not true for isomorphisms $f:A\ra B$ in classes/categories of sets (classes) where the bijection-inverse $\ol{f}:B\ra A$ is not automatically a homomorphism.})
\end{thm}
\begin{proof}
This follows directly from the definitions.
\end{proof}

%% file: parts/AlgebraNC/AlgebraNC4.tex
\chapter{Subobjects, Quotient Objects, and Related Structures}
\section{Subobjects and Substructures}
\begin{dfn}[\textcolor{blue}{
\index{Subgroup}{Subgroup},
\index{Generated! subgroup}{Generated subgroup},
\index{Generator of! a subgroup}{Generator of a subgroup},
\index{Finitely generated! group}{Finitely generated group},
\index{Cyclic! group}{Cyclic group},
\index{Order! of a group}{Order of a group},
\index{Order! of a group element}{Order of a group element},
\index{Left! coset}{Left coset},
\index{Right! coset}{Right coset},
\index{Normal subgroup}{Normal subgroup},
\index{Simple! group}{Simple group}}]
Let $G=[G,\cdot,e,()^{-1}]$ be a group. A subset $H\subset G$ is a \ul{subgroup}, denoted $H\leq G$, if $H$ is itself a group under the operations of $G$, i.e., (i) $e\in H$, (ii) $a,b\in H$ $\Ra$ $ab\in H$, (iii) $a\in H$ $\Ra$ $a^{-1}\in H$ (or equivalently, $a,b\in H$ $\Ra$ $ab^{-1}\in H$).

If $A\subset G$, the \ul{subgroup generated} by $A$, $\langle A\rangle\leq G$, (with $A$ called a \ul{generator} of $\langle A\rangle$) is the smallest subgroup containing $A$, i.e.,
\bea
\textstyle\langle A\rangle:=\bigcap\{H\leq G:A\subset H\}=\big\{b_1b_2\cdots b_n~|~b_i\in\{e\}\cup A\cup A^{-1},~n\geq 0\big\},~~~~A^{-1}:=\{a^{-1}:a\in A\}.\nn
\eea
The group $G$ is \ul{finitely generated} if $G=\langle A\rangle$ for a finite set $A\subset G$.

The group $G$ is called \ul{cyclic} if it is generated by one element (i.e., $G=\langle g\rangle=\{g^k:k\in\Integer\}$ for some element $g\in G$). The \ul{order} of $G$ is its cardinality $|G|$. The \ul{order of an element} $g\in G$ is $|g|:=|\langle g\rangle|$.

If $H\leq G$, the sets of \ul{left cosets} and \ul{right cosets} of $H$ in $G$ are, respectively, given by the partitions
{\small\[
G/H:=(G/H)_l:=\{gH~|~g\in G\}~~(\txt{left cosets}),~~~~(G/H)_r:=\{Hg~|~g\in G\}~~(\txt{right cosets}).~~~~(\textcolor{blue}{footnote}\footnotemark).
\]}
\footnotetext{\blue{Properties of cosets}: Let $H\leq G$.
{\flushleft (i)} For any $g,g'\in G$, ~~$g'\in gH$ $\iff$ $g^{-1}g'\in H$, $\iff$ $g^{-1}g'H=H$, $\iff$ $gH=g'H$ ($\iff$ $Hg^{-1}=Hg'^{-1}$).
{\flushleft (ii)} For any $g,g'\in G$, ~~$g\in Hg'$ $\iff$ $gg'^{-1}\in H$, $\iff$ $Hgg'^{-1}=H$, $\iff$ $Hg=Hg'$ ($\iff$ $g^{-1}H=g'^{-1}H$).
{\flushleft (iii) By definition,} $\{gH:g\in G\}=\{Hg:g\in H\}$ $\iff$ for any $g\in G$ we have $gH=Hg'$ for some $g'=g'_g\in G$.
{\flushleft (iv) Left-right symmetry:} For any $g,g'\in G$, ~$gH=Hg'$ $\iff$ $g'\in gH=Hg'\ni g$, $\iff$ $gH=g'H=Hg'=Hg$.

{\flushleft (v) Disjointness:} For any $g,g'\in G$, ~$gH\cap g'H\neq\emptyset$ (i.e., some $g''\in gH\cap g'H$) $\iff$ $gH=g''H=g'H$.}A subgroup $N\leq G$ is a \ul{normal subgroup}, denoted by $N\lhd G$, if $(G/N)_l=(G/N)_r$ (or equivalently, $gN=Ng$ for all $g\in G$, or equivalently, $gNg^{-1}\subset N$ for all $g\in G$).  The group $G$ is called \ul{simple} if its only normal subgroups are $\{e\}$ and $G$.

\end{dfn}

\begin{lmm}[\blue{\index{Index of a group}{Index of a subgroup}}]
The set of left cosets $G/H:=\{aH:a\in G\}$ and the set of right cosets $(G/H)_r:=\{Ha:a\in G\}$ have the same cardinality (called the \ul{index}, $[G:H]$, of $H$ in $G$).
\end{lmm}
\begin{proof}
The map $f:G/H\ra (G/H)_r,~aH\mapsto Ha^{-1}$ is a well defined bijection, since $aH=bH$ $\iff$ $a^{-1}bH=H$, $\iff$ $a^{-1}b\in H$, $\iff$ $Ha^{-1}b=H$, $\iff$ $Ha^{-1}=Hb^{-1}$.
\end{proof}

\begin{lmm}
Every subgroup of index $2$ is normal.
\end{lmm}
\begin{proof}
Let $G$ be a group and $H\leq G$ a subgroup with index $[G:H]=2$. The left cosets are $\{H,aH\}$ for some $a\in G$. The right cosets are $\{H,Hb\}$ for some $b\in G$. Since the cosets partition $G$, we have~ $aH=G\backslash H=Hb$.~ Hence every left coset is a right coset, and so $H$ is normal.
\end{proof}

\begin{lmm}
Let $H,K\leq G$ be subgroups. (i)  If $H$ or $K$ is normal, then $HK=KH\subset G$ is a subgroup. (ii) If both $H$ and $K$ are normal, then the subgroup $HK\leq G$ is normal.
\end{lmm}
\begin{proof}
{\flushleft (i)} If $K$ is normal, then for any $hk,h'k'\in HK$, we have $hk(h'k')^{-1}=hkk'^{-1}h'^{-1}\in K\subset HK$. Similarly, if $H$ is normal, then for any $hk,h'k'\in HK$, we have $(hk)^{-1}h'k'=k^{-1}h^{-1}h'k'\in H\subset HK$.
{\flushleft (ii)} By (i), $HK\subset G$ is a subgroup, and for any $g\in G$,~ $gHKg^{-1}=gHg^{-1}~gKg^{-1}~\subset~HK$.
\end{proof}

\begin{lmm}
Let $G$ be a group and $H,K\leq G$ subgroups. Then $H\subset K$ $\iff$ $aH\subset bK$ for some $a,b\in G$ (i.e., some coset of $H$ is contained in some coset of $K$).
\end{lmm}
\begin{proof}
If $H\subset K$, it is clear that some coset of $H$ is contained in some coset of $K$. Conversely, let $a,b\in G$ such that $aH\subset bK$. Then $a\in aH\subset bK$ $\Ra$ $aK=bK$, $\Ra$ $H\subset a^{-1}bK=K$.
\end{proof}

\begin{dfn}[\textcolor{blue}{
\index{Subring}{Subring},
\index{Generated! subring}{Generated subring},
\index{Generator of! a subring}{Generator of a subring},
\index{Finitely generated! ring}{Finitely generated ring},
\index{Left! ideal}{Left ideal},
\index{Right! ideal}{Right ideal},
\index{Ideal (two-sided)}{(two-sided) Ideal},
\index{Generated! ideal}{Generated ideal},
\index{Generator of! an ideal}{Generator of an ideal},
\index{Principal ideal!}{Principal ideal},
\index{Principal ideal! ring}{Principal ideal ring},
\index{Principal ideal! domain (PID)}{Principal ideal domain (PID)},
\index{Prime! ideal}{Prime ideal},
\index{Maximal! ideal}{Maximal ideal},
\index{Minimal! ideal}{Minimal ideal},
\index{Local! ring}{Local ring}}]
Let $R$ be a ring. A subset $S\subset R$ is a \ul{subring}, denoted $S\leq R$, if $S$ is itself a ring under the operations of $R$, i.e., (i) $(S,+)\leq(R,+)$ is a subgroup, (ii) $a,b\in S$ $\Ra$ $ab\in S$, (iii) $1_R\in S$ (or more generally if desired, there exists $e\in S$ such that $se=es=s$ for all $s\in S$). If $C\subset R$, the \ul{subring generated} by $C$, $S(C)\leq R$, (with $C$ called a \ul{generator} of $S(C)$) is the smallest subring containing $C$, i.e.,
\begin{align}
&\textstyle S(C):=\bigcap\big\{S\leq R~|~C\subset S\big\}=\Integer\langle 1_R,C\rangle=\Integer 1_R\langle C\rangle=\big\{n_0b_0+\cdots+n_kb_k~|~n_i\in\Integer,~b_i\in C^{\ast i},~k\in\Natural\big\},\nn\\
&~~~~=\Integer 1_R+\Integer C+\Integer C^{\ast 2}+\Integer C^{\ast 3}+\cdots,~~C^{\ast 0}:=\{1_R\},~~C^{\ast i}:=\big\{c_1c_2\cdots c_i~|~c_j\in C\big\}~~\txt{for}~~i\geq 1.\nn
\end{align}
The ring $R$ is \ul{finitely generated (as a ring)} if $R=S(C)$ for a finite set $C\subset R$.

A subset $I\subset R$ is a \ul{left ideal} of $R$, denoted $I\lhd_lR$ or $_RI\lhd R$ (resp. \ul{right ideal} of $R$, denoted $I\lhd_rR$ or $I_R\lhd R$) if (i) $(I,+)\leq(R,+)$ is a subgroup and (ii) $RI\subset I$ (resp. $IR\subset I$). A subset $I\subset R$ is a \ul{(two-sided) ideal} of $R$, denoted $I\lhd R$, if $I$ is both a left and right ideal of $R$ (i.e., $(I,+)\leq (R,+)$ is a subgroup and $RI,IR\subset I$).

If $C\subset R$, the \ul{left ideal generated} by $C$, $I_l(C)\lhd_l R$, is the smallest left ideal (resp. \ul{right ideal generated} by $C$, $I_r(C)\lhd_r R$, is the smallest right ideal) containing $C$, i.e.,
\bea
\textstyle I_l(C):=\bigcap\big\{I\lhd_l R~|~C\subset I\big\}=RC,~~~~I_r(C):=\bigcap\big\{I\lhd_r R~|~C\subset I\big\}=CR.\nn
\eea
If $C\subset R$, the \ul{ideal generated} by $C$, $I(C)\lhd R$, (with $C$ called a \ul{generator} of $I(C)$) is the smallest ideal containing $C$, i.e.,
\bea
\textstyle I(C):=\bigcap\big\{I\lhd R~|~C\subset I\big\}=RCR.\nn
\eea

A commutative ring $R$ is a \ul{principal ideal ring} if every ideal $I\lhd R$ is of the form (called \ul{principal ideal}) $I=Ra$ for some $a\in R$ (i.e., every ideal is generated by a single element). A \ul{principal ideal domain} (PID) is a principal ideal integral domain.

A proper ideal $P\lhd R$, $P\neq R$, is a \ul{prime ideal} if for all ideals $I,J\lhd R$, if $IJ\subset P$ then $I\subset P$ or $J\subset P$ (i.e., for all $a,b\in R$,~ $aRb\subset P$ implies $a\in P$ or $b\in P$). A proper ideal $M\lhd R$, $M\neq R$, is a \ul{maximal ideal} if for any ideal $I\lhd R$, if $I\supset M$ then $I=M$ or $I=R$ (i.e., $RaR+M=R$ for any $a\not\in M$). A nonzero proper ideal $0\neq M\lhd R$, $M\neq R$, is a \ul{minimal ideal} or \ul{simple ideal} if for any ideal $I\lhd R$, if $I\subset M$ then $I=0$ or $I=M$ (i.e., $RaR=I$ for all $0\neq a\in I$).

If a ring $R$ has a unique maximal ideal $\mfm$, then $R=(R,\mfm)$ is called a \ul{local ring}.
\end{dfn}

\begin{lmm}
A proper ideal $P\vartriangleleft R\neq P$ is prime (\blue{footnote}\footnote{Prime in the sense that for all ideals $I,J\vartriangleleft R$, ~$IJ\subset P$ ~$\Ra$~ $I\subset P$ or $J\subset P$.}) if and only if for all $a,b\in R$,
\bea
aRb\subset P~~\Ra~~a\in P~~\txt{or}~~b\in P.\nn
\eea
\end{lmm}
 \begin{proof}
Assume $P$ is prime, i.e., for all $I,J\lhd R$, $IJ\subset P~\Ra~I\subset P$ or $J\subset P$. Then for any $a,b\in R$,
\bea
aRb\subset P~~\Ra~~RaR~RbR\subset RPR=P~~\Ra~~RaR\subset P~~\txt{or}~~RbR\subset P,~~\Ra~~a\in P~~\txt{or}~~b\in P.\nn
\eea
Conversely, assume $P$ is not prime, i.e., there exist $I,J\lhd R$ such that $IJ\subset P$ but $I\not\subset P$ and $J\not\subset P$. Let $a\in I\backslash P$,~ $b\in J\backslash P$. Then $aRb=(aR)(Rb)\subset IJ\subset P$ but $a\not\in P$ and $b\not\in P$.
\end{proof}

\begin{lmm}[\blue{A maximal ideal is a prime ideal}]
\end{lmm}
 \begin{proof}
  Let $M\lhd R$ be maximal. Let $I,J\lhd R$ be such that $IJ\subset M$. Suppose on the contrary that $I\not\subset M$ and $J\not\subset M$. Then ~$R=RR=(I+M)(J+M)~\subset~IJ+M=M$~ (a contradiction).
 \end{proof}

Directly from the definitions (i.e., Definition \ref{GcdLcmDef} for gcd and lcm, and the definition of a PID), we have the following immediate result.

\begin{prp}[\textcolor{blue}{\index{gcd and lcm in a PID}{gcd and lcm in a PID}}]
If $R$ is a PID, the gcd and lcm (for any $a_1,...,a_n\in R$) satisfy
\bea
Ra_1+\cdots+Ra_n~=~R~\txt{gcd}(a_1,...,a_n)~~~~\txt{and}~~~~Ra_1\cap\cdots\cap Ra_n~=~R~\txt{lcm}(a_1,...,a_n).\nn
\eea
\end{prp}
\begin{proof}
These results are immediate consequences of the relations (\ref{GCDFormulaEq}) and (\ref{LCMFormulaEq}).
\end{proof}

\begin{prp}[\textcolor{blue}{Every ED is a PID}]
If $R$ is a ED (Euclidean domain), then $R$ is a PID.
\end{prp}
\begin{proof}
Let $R$ be a ED with degree map $\sigma:R\backslash\{0\}\ra\Natural$, and let $I\vartriangleleft R$. If $I=0$, then $I=R0$, and thus $I$ is principal. So assume $I\neq 0$. Choose $a\in I\backslash 0$ such that $\sigma(a)$ is minimal in $I$, i.e., $\sigma(a)=\min\{\sigma(c):c\in I\backslash 0\}$. If $b\in I$, we can write
\bea
b=qa+r,~~\txt{for some}~~q,r\in R,~~r=0~~\txt{or}~~\sigma(r)<\sigma(a).\nn
\eea
Now if $r\neq 0$, then $r=b-qa=b+(-q)a\in I$ implies $\sigma(a)\leq \sigma(r)$, which is a contradiction. Thus $b=qa\in Ra$, i.e., $I\subset Ra$. Since $a\in I$, we also have $Ra\subset I$. Hence $I=Ra$.
\end{proof}

\begin{thm}[\textcolor{blue}{Division algorithm for $\Integer$}]
Let $b\in\Integer$, $b>0$. For each $a\in\Integer$, there exist $q,r\in\Integer$ such that
~$a=qb+r$, with $q,r\in\Integer$ unique and ~$0\leq r<b$. That is, $\Integer=\Integer b+\Integer\cap[0,b)$.
\end{thm}
\begin{proof}
{\flushleft\ul{Existence}}: Given $a\in R$, let $S:=\{a-kb:~k\in\Integer,~a-kb\geq 0\}=(a-\Integer b)\cap\Natural=(a+\Integer b)\cap\Natural$. Then $S$ is nonempty, since (i) $a\geq 0$ $\Ra$ $a=a-0\cdot b\in S$, and (ii) $a<0$ $\Ra$ $a-ab=-a(b-1)\in S$. Thus by the well-ordering of $\Natural$, $S$ has a least element $r$, and so there is $q\in\Integer$ such that $r=a-qb$, i.e., $a=qb+r$. If $r\geq b$, then $r-b=a-(q+1)b\geq 0$ implies $r>r-b\in S$ (a contradiction). Therefore ~$0\leq r<b$.

{\flushleft\ul{Uniqueness}}: Let $a=qb+r,~0\leq r<b,$ and $a=q'b+r',~0\leq r'<b$. Assume wlog that $r\leq r'$. Then
\bea
qb+r=q'b+r'~~\Ra~~0\leq r'-r=(q-q')b\leq r'<b,\nn
\eea
that is, (i) $b$ divides $r'-r$ and (ii) $0\leq r'-r\leq r'<b$. Therefore, $r'-r=0$ and $q-q'=0$.
\end{proof}

\begin{crl}[\textcolor{blue}{The ring $\Integer$ is a ED, and hence a PID}]
\end{crl}

\begin{lmm}[\textcolor{blue}{\textcolor{blue}{Maximal ideals in a PID}}]
Let $R$ be a PID. Then every \ul{nonzero} maximal ideal has the form $Rp$ for a prime element (equivalently, an irreducible element) $p\in R$.
\end{lmm}
\begin{proof}
($\Ra$): Assume $0\neq Rx\lhd R$ is maximal. We show $x$ is prime. It is clear that $x\neq 0$. Since $Rx\lhd R$ is maximal, for any $r\in R$, $r\not\in Rx$ $\iff$ $Rr+Rx=R$. Thus, for any $a,b\in R$, if $a\not\in Rx$ and $b\not\in Rx$ then {\small $R=RR=(Ra+Rx)(Rb+Rx)=Rab+Rx$}, i.e., if $a\not\in Rx$ and $b\not\in Rx$ then $ab\not\in Rx$. Equivalently, for any $a,b\in R$, if $ab\in Rx$ then $a\in Rx$ or $b\in Rx$. Therefore $x$ is prime.

($\La$): Conversely, let $p\in R$ be prime (hence irreducible). We show $Rp\lhd R$ is maximal. If $Rp\subset I\subset R$ for some ideal $I\vartriangleleft R$, then $I=Ra$ for some $a\in R$ (as $R$ is a PID). Thus ~$Rp\subset Ra$, which implies
\bit
\item[] $p=ra$ ~for some ~$r\in R$, ~which in turn implies ~$r\in U(R)$ ~or ~$a\in U(R)$.
\eit
If $r\in U(R)$, then $I=Ra=R(r^{-1}p)=Rp$. If $a\in U(R)$, then $I=Ra=R$. Hence $Rp$ is maximal.
\end{proof}

\begin{crl}[\textcolor{blue}{In a PID, irreducible = prime}]
Let $R$ be a PID. Then $x\in R$ is irreducible iff prime.
\end{crl}

\begin{lmm}[\textcolor{blue}{In an ID, prime factorization is unique}]~
If $R$ is an ID and ~$p_1p_2\cdots p_r=q_1q_2\cdots q_s$ (where $p_i,q_i$ are prime), then (i) $r=s$, and (ii) after reordering, $p_i$ is an associate of $q_i$, for $1\leq i\leq r$.
\end{lmm}
\begin{proof}
We use induction on $r+s$ (may also use induction on $\max(r,s)$). If $r+s=2$, the result is immediate since we must have $r=s=1$. So assume $r+s>2$. Then
\bea
p_1|q_1q_2...q_s~~\Ra~~p_1|q_i~~\txt{for some}~~i.\nn
\eea
Without loss of generality, let $i=1$. Then $q_1=up_1$ for some $u\in U(R)$. Thus
\bea
p_1p_2...p_r=p_1uq_2...q_s,~~\Ra~~p_2p_3...p_r=uq_2q_3...q_s.\nn
\eea
Since $(r-1)+(s-1)=r+s-2<r+s$, the result follows by the induction hypothesis.
\end{proof}

\begin{lmm}[\textcolor{blue}{In a UFD, irreducible = prime}]~
Let $R$ be a UFD. Then $x\in R$ is irreducible iff prime.
\end{lmm}
\begin{proof}
{\flushleft ($\Ra$)}: Assume $x$ is irreducible. Let $a,b\in R$ such that $x|ab$, i.e., $ab=xc$ for some $c\in R$. Write $a=x_1x_2...x_l$, $b=y_1y_2...y_m$, $c=z_1z_2...z_n$, where $x_i,y_i,z_i$ are all irreducible. Then
\[
x(z_1z_2...z_n)=x_1x_2...x_ly_1y_2...y_m.
\]
By the uniqueness of factorization, $1+n=l+m$, and $x$ is an associate of $x_i$ (for some $i$) or of $y_j$ (for some $j$). This implies $x|a$ or $x|b$, and so $x$ is prime.
{\flushleft ($\La$)}: We already know that a prime is irreducible.
\end{proof}

\begin{lmm}[\textcolor{blue}{Character of FF elements of an ID}]\label{FFCharLmm}
Let $R$ be an ID. Then every $a\in R\backslash(0\cup U(R))$ is FF (i.e., fully factorable) $\iff$ there is no infinite sequence of principal ideals~ $Ra_1\subsetneq Ra_2\subsetneq Ra_3\subsetneq \cdots$.
\end{lmm}
\begin{proof}
{\flushleft($\Ra$)} Assume every $a\in R\backslash(0\cup U(R))$ is FF. Suppose on the contrary that an infinite sequence~ $Ra_1\subsetneq Ra_2\subsetneq Ra_3\subsetneq \cdots$~ exists. Then for each $i$, $a_i\in Ra_{i+1}$, and so $a_i=b_{i+1}a_{i+1}$ for some $b_{i+1}\not\in U(R)$. Therefore, ~$a_1=b_2a_2=b_2b_3a_3=b_2b_3b_4a_4\cdots$, ~and we get an infinite number of factors (as the process does not stop), and so $a_1$ is not FF (a contradiction).
{\flushleft($\La$)} Assume no such sequence of principal ideals exists. Suppose on the contrary that $a\in R\backslash(0\cup U(R))$ is not FF. Then $a=b_1a_1$, where $a_1,b_1\not\in U(R)$, and (without loss of generality,) $a_1$ is not FF. Similarly, $a_1=b_2a_2$, where $a_2,b_2\not\in U(R)$, and $a_2$ is not FF. Continuing this way, we obtain an infinite sequence~ $Ra\subsetneq Ra_1\subsetneq Ra_2\subsetneq \cdots$ (a contradiction).
\end{proof}

\begin{thm}[\textcolor{blue}{Every PID is a UFD}]
If $R$ is a PID, then $R$ is a UFD.
\end{thm}
\begin{proof}
Since $R$ is a PID, an element is irreducible if and only if it is prime. Moreover, uniqueness of factorization into irreducibles follows from uniqueness of factorization into primes. Thus, it remains only to show that every $a\in R\backslash(0\cup U(R))$ has a factorization into primes, and that every such $a$ is FF.

By Lemma \ref{FFCharLmm}, suppose on the contrary that some $a\not\in(\{0\}\cup U(R))$ is not FF (into primes), i.e., we have an infinite sequence $Ra_1\subsetneq Ra_2\subsetneq Ra_3\subsetneq \cdots$. Let $I:=\bigcup_{i=1}^\infty Ra_i$. Then it is clear that $I\vartriangleleft R$. Since $R$ is a PID, $I=Rb$ for some $b\in R$. That is, $b\in Rb=\bigcup_{i=1}^\infty Ra_i$, which implies $b\in Ra_n$ for some $n$, and so $Rb\subset Ra_n\subsetneq Ra_{n+1}\subsetneq Ra_{n+2}\subsetneq...\subset Rb$,
which is a contradiction. Hence every $a\in R\backslash(0\cup U(R))$ is FF.
\end{proof}

\begin{summary*}[\textcolor{blue}{ED's $\subset$ PID's $\subset$ UFD's $\subset$ FFD's $\subset$ ID's}]
\end{summary*}

\begin{dfn}[\textcolor{blue}{
\index{Submodule}{Submodule},
\index{Minimal! submodule (Simple submodule)}{Minimal submodule (Simple submodule)},
\index{Maximal! submodule}{Maximal submodule},
\index{Generated! submodule}{Generated submodule},
\index{Generator of! a submodule}{Generator of a submodule},
\index{Finitely generated! module}{Finitely generated module},
\index{Cyclic! module}{Cyclic module},
\index{Simple! module}{Simple module},
\index{Semisimple module}{Semisimple module},
\index{Socle of a module}{Socle of a module},
\index{Decomposable module}{Decomposable module},
\index{Indecomposable module}{Indecomposable module},
\index{Completely decomposable (Completely reducible) module}{Completely decomposable (Completely reducible) module},
\index{Noetherian module}{Noetherian module},
\index{Artinian module}{Artinian module},
\index{Left! noetherian ring}{Left noetherian ring},
\index{Right! noetherian ring}{Right noetherian ring},
\index{Left! artinian ring}{Left artinian ring},
\index{Right! artinian ring}{Right artinian ring},
\index{Simple! ring}{Simple ring},
\index{Left-semisimple ring}{Left-semisimple ring},
\index{Right-semisimple ring}{Right-semisimple ring},
\index{Semisimple ring}{Semisimple ring}}]
Let $M$ be an $R$-module. A subset $N\subset M$ is a \ul{submodule}, denoted by ${}_RN\subset{}_RM$ or $N\subset_RM$, if $N$ is itself an $R$-module under the operations of $M$, i.e., (i) $(N,+)\leq(M,+)$ is a subgroup nd (i) $RN\subset N$.

A nonzero proper submodule $\{0\}\neq{}_RN\subsetneq {}_RM$ is a \ul{minimal submodule} or \ul{simple submodule} if it cannot properly contain a \emph{nonzero submodule}, i.e., if $\{0\}\subset{}_RX\subset{}_RN\subsetneq{}_RM$, then $X=\{0\}$ or $X=N$ (or equivalently, $Rn=N$ for any $0\neq n\in N$). A proper submodule ${}_RN\subsetneq {}_RM$ is a \ul{maximal submodule} if it cannot be properly contained in a \emph{proper submodule}, i.e., if ${}_RN\subset{}_RX\subset{}_RM$, then $X=N$ or $X=M$ (or equivalently, $Rm+N=M$ for any $m\not\in N$).

If $C\subset M$, the \ul{submodule generated} by $C$, ${}_R\langle C\rangle\subset{}_RM$, (with $C$ called a \ul{generator} of $\langle C\rangle$) is the smallest submodule containing $C$, i.e., ~$\langle C\rangle:=\bigcap\big\{{}_RN\subset{}_RM~|~C\subset N\big\}=RC$. The module $M$ is called \ul{cyclic} if it is generated by one element (i.e., $M=Rm$ for some element $m\in M$). The module $M$ is \ul{finitely generated} if $M=\langle C\rangle$ for a finite set $C\subset M$, in which case, if $C=\{c_1,...,c_t\}$, then
\bea
M=Rc_1+Rc_2+\cdots+Rc_t.\nn
\eea

The module $M$ is called \ul{simple} if (i) $M\neq\{0\}$ and (ii) its only submodules are $\{0\}$ and $M$ (i.e., $Rm=M$ for any $0\neq m\in M$).
(\blue{footnote}\footnote{Note that $M$ is simple $\iff$ ${}_R\{0\}\subset{}_RM$ is a maximal submodule. Also, it is clear that a simple module $M$ is necessarily cyclic, i.e., $M=Rm$ for some $0\neq m\in M$. Thus, as for any cyclic module, submodules of $M=Rm$ are $\{Im:\txt{for left ideals}~{}_RI\lhd R\}$.
}).
The module $M$ is \ul{semisimple} if $M=\sum_{i\in I}S_i:=\left\langle\bigcup_{i\in I}S_i\right\rangle=\left\{m_1+\cdots+m_n~|~m_j\in\bigcup S_i,~n\in\Natural\right\}$ for simple submodules ${}_RS_i\subset{}_RM$. (\blue{footnote}\footnote{The sum  of all simple submodules of $M$ is called the \ul{socle} of $M$, written $Soc(M):=\sum\{\txt{simple submodules}~S\subset M\}$. Thus, $M$ is semisimple iff $Soc(M)=M$.}).

The module $M$ is \ul{decomposable} if $M$ a direct sum of two nonzero submodules in the sense $M=N_1\oplus N_2$ (meaning $M=N+N'$, $N\cap N'=\{0\}$) for nonzero submodules $\{0\}\neq N_1,N_2\subset M$, otherwise $M$ is \ul{indecomposable} (i.e., for any nonzero submodules $\{0\}\neq N_1,N_2\subset M$, we have $M\neq N_1\oplus N_2$). (\blue{footnote}\footnote{The ``\ul{disjoint sum}'' $\oplus$ here is equivalent to the \ul{direct sum} $\oplus$, or coproduct $\coprod$, defined in section \ref{MdPdSec}.}). The module $M$ is \ul{completely decomposable} (or \ul{completely reducible}) if every submodule ${}_RN\subset{}_RM$ is a direct summand of $M$ in the sense $M=N\oplus N'$ for some submodule ${}_RN'\subset{}_RM$.

The module $M$ is \ul{noetherian} (resp. \ul{artinian}) if every increasing sequence of submodules $N_1\subset N_2\subset\cdots$ (resp. decreasing sequence of submodules $N_1\supset N_2\supset\cdots$) of $M$ terminates or stabilizes, i.e., there exists $t$ such that $N_i=N_t$ for all $i\geq t$. A ring $R$ is \ul{left-noetherian} (resp. \ul{left-artinian}) if the module ${}_RR$ is noetherian (resp. artinian). Similarly, a ring $R$ is \ul{right-noetherian} (resp. \ul{right-artinian}) if the module ${}_RR$ is noetherian (resp. artinian). A ring $R$ is a \ul{noetherian ring} (resp. an \ul{artinian ring}) if it is both left and right noetherian (resp. both left and right artinian). (\blue{footnote}\footnote{A simple module is clearly both noetherian and artinian.}).

A ring $R$ is \ul{simple} if its only ideals are $\{0\}$ and $R$ (i.e., $RaR=R$ for all $0\neq a\in R$). A ring $R$ is \ul{left-semisimple} (resp. \ul{right-semisimple}) if the module $_RR$ ($R_R$) is semisimple. A ring is a \ul{semisimple ring} if it is both left-semisimple and right-semisimple.
\end{dfn}

\begin{note*}
Let $R$ be a ring. The submodules of ${}_RR$ (resp. $R_R$) are precisely the left (resp. right) ideals of $R$. Similarly, the submodules of ${}_RR_R$ are precisely the (two-sided) ideals of $R$.
\end{note*}

\begin{dfn}[\textcolor{blue}{\index{Subalgebra}{Subalgebra}, \index{Generated! subalgebra (``Polynomial ring/algebra'')}{Generated subalgebra (``Polynomial ring/algebra'')}, \index{Generator of! a subalgebra}{Generator of a subalgebra}, \index{Finitely generated! algebra}{Finitely generated algebra}}]
Let $R$ be a commutative ring and $A$ an $R$-algebra. A subset $B\subset A$ is a \ul{subalgebra}, denoted by ${}_RB\leq {}_RA$ or $B\leq_RA$, if $B$ is itself an $R$-algebra under the operations of $A$, i.e., (i) $(B,\cdot,+)\leq(A,\cdot,+)$ is a subring and (ii) ${}_RB\subset{}_RA$ is an $R$-submodule.

If $C\subset A$, the \ul{subalgebra generated} (or \ul{``polynomial ring/algebra'' generated}) by $C$, ${}_RS(C)\leq{}_RA$, (with $C$ called a \ul{generator} of $S(C)$) is the smallest subalgebra containing $C$, i.e.,
\begin{align}
&\textstyle S(C):=\bigcap\big\{{}_RB\leq{}_RA~|~C\subset B\big\}=R\langle1_A,C\rangle=R1_A\langle C\rangle=\big\{r_0b_0+\cdots+r_kb_k~|~r_i\in R,~b_i\in C^{\ast i},~k\in\Natural\big\},\nn\\
&~~~~=R1_A+RC+RC^{\ast 2}+RC^{\ast 3}+\cdots,~~C^{\ast 0}:=\{1_A\},~~C^{\ast i}:=\big\{c_1c_2\cdots c_i~|~c_j\in C\big\},~~i\geq 1.\nn
\end{align}
The algebra $A$ is \ul{finitely generated} if it is finitely generated as a ring, i.e., $A=S(C)$ for a finite set $C\subset A$.
\end{dfn}

\begin{note*}
By construction, all concepts for rings and modules also apply to algebras.
\end{note*}

\section{Linear algebra of modules (a concretizing digression)}\label{LinAlgModSec1}
\begin{dfn}[\textcolor{blue}{
\index{Span of a subset}{Span of a subset},
\index{Spanning! set}{Spanning set},
\index{Linearly independent set}{Linearly independent set},
\index{Linear! combination}{Linear combination},
\index{Basis}{Basis},
\index{Free! module (Module with a basis)}{Free module (Module with a basis)}}]
Let $M$ be an $R$-module and $X\subset M$ a subset. The \ul{$R$-span} (or just the \ul{span}) of $X$ is the submodule generated by $X$, i.e., $\Span_RX:=RX$. The set $X$ is a \ul{spanning set} for $M$ if $\Span_RX:=M$, i.e., $X$ generates $M$. The set $X$ is \ul{linearly independent} (or just \ul{independent}) if for every finite subset $\{x_1,...,x_n\}\subset X$,
\bea
(\txt{for every}~~r_1,...,r_n\in R)~~~~r_1x_1+\cdots r_nx_n=0~~\Ra~~r_1=\cdots=r_n=0,\nn
\eea
where the sum $\sum_{i=1}^nr_ix_i$ is called a \ul{linear combination} of $x_1,...,x_n$.

 A \ul{basis} for $M$ (making $M$ a \ul{free module} or a \ul{module with a basis}) is a linearly independent spanning set for $M$.
\end{dfn}

\begin{thm}[\textcolor{blue}{Basis cardinality theorem I}]\label{SThCardThm1}
Let $F$ be a free R-module with two bases $X$ and $Y$.
\bit[leftmargin=0.7cm]
\item[(1)] $X$ is finite $\iff$ $Y$ is finite (i.e., $X$ is infinite $\iff$ $Y$ is infinite).
\item[(2)] If $X$ is infinite, then $|X|=|Y|$.
\eit
\end{thm}
\begin{proof}
{\flushleft(1)} By symmetry, assume wlog that $|X|$ is finite (so we need to show $|Y|$ is also finite). Since $X$ is a finite set, we have $X\subset Ry_1+\cdots+Ry_n$ for some $y_1,...,y_n\in Y$. This implies $F=\sum_{x\in X}Rx=Ry_1+\cdots+Ry_n$. Since $Y\subset F$ and the elements of $Y$ are linearly independent, it follows that $Y=\{y_1,...,y_n\}$, because any element of $Y$ not in $\{y_1,...,y_n\}$ must be a linear combination of the $y_i$'s (contradicting linear independence of $Y$).
{\flushleft(2)} By part (1) it is clear that $X$ and $Y$ are both infinite. Let $FS(X)$ be the set of finite subsets of $X$, and  $FS(Y)$ likewise. Then we know (from Corollary \ref{FinSubCaCrl}) that $|FS(X)|=|X|$ and $|FS(Y)|=|Y|$. Since $Y\subset F=\Span_RX$, every $y\in Y$ can be written as $y=\sum_{i=1}^{n_y}r_i^yx_i^y$ for some $r_i^y\in R$ and $x_i^y\in X$. Since $Y$ is \ul{linearly independent}, it follows from part (1) that, the set $\{x_1^y,...,x_{n_y}^y\}$ determines at most a finite number of elements of $Y$ in its span $Rx_1^y+\cdots+Rx_{n_y}^y$. The cardinality of the set of coefficient-tuples {\small $C_X(Y):=\left\{(r_i^y)_{i=1}^{n_y}:y\in Y\right\}$} equals that of {\small $FS_Y(X):=\left\{(x_i^y)_{i=1}^{n_y}:y\in Y\right\}\subset FS(X)$} and is at least that of $Y$.
\[
|Y|=|FS(Y)|\leq|C_X(Y)|=FS_Y(X)\leq |FS(X)|=|X|.
\]
By symmetry we also have $|X|\leq|Y|$. Hence ~$|X|=|Y|$.
\end{proof}

\begin{dfn}[\textcolor{blue}{
\index{Sum (or span) of submodules}{Sum (or span) of submodules},
\index{Spanning! family}{Spanning family},
\index{Free! family of submodules}{Free family of submodules},
\index{Internal direct sum}{Internal direct sum},
\index{Complemented submodule}{Complemented submodule}}]
Let $M$ be an $R$-module and $\F:=\{M_i:i\in I\}$ a family of submodules of $M$. The \ul{sum} (or \ul{span}) of $\F$ is the span of the union $\bigcup_{i\in I} M_i$ (i.e., the submodule generated by $\bigcup_{i\in I} M_i$), that is,
\bea
\textstyle\Span_R\F:=\sum_{i\in I}M_i:=\Span_R(\bigcup_{i\in I}M_i)=\left\langle\bigcup_{i\in I}M_i\right\rangle=\left\{m_1+\cdots+m_n~|~m_j\in\bigcup M_i,~n\in\Natural\right\}.\nn
\eea
The family $\F$ is a \ul{spanning family} if $\sum_{i\in I} M_i=M$. The family $\F$ is a \ul{free family} if $M_i\cap\left(\sum_{j\neq i}M_j\right)=0$ for all $i\in I$ (in which case, the sum $\sum_{i\in I}M_i$ is called an \ul{internal direct sum}, and written as $\bigoplus_{i\in I}M_i$).

A submodule $N\subset M$ is \ul{complemented} in $M$ if there exists a submodule $N'\subset M$ such that $N\oplus N'=M$ (i.e., such that $N\cap N'=0$ and $N+N'=M$).

\end{dfn}
It is clear that a module is an internal direct sum iff it has a spanning free family of submodules.

\begin{thm}
Let $V$ be a $k$-vector space and $B\subset V$ a subset. The following are equivalent.
\begin{enumerate}[leftmargin=0.7cm]
\item[(a)] $B$ is a linearly independent spanning set for $V$.
\item[(b)] $B$ is a maximal linearly independent subset of $V$. (\blue{footnote}\footnote{That is, no linearly independent subset of $V$ contains $B$ as a proper subset. Equivalently, (1) $B$ is linearly independent, and (2) if $v\not\in B$, then $B\cup\{v\}$ is dependent.})
\item[(c)] $B$ is a minimal spanning set for $V$. (\blue{footnote}\footnote{That is, no proper subset of $B$ can span $V$. Equivalently, (1) $Span_kB=V$, and (2) if $C\subset B$ and $Span_kC=V$, then $C=B$.})
\end{enumerate}
\end{thm}
\begin{proof}~
\bit[leftmargin=0.5cm]
\item\ul{(a)$\Ra$(b)}: Let $B$ be an independent spanning set for $V$. Suppose $B\cup\{v\}$ is linearly independent for some $v\not\in B$. Since $B$ is a spanning set, we have
    \bea
    \textstyle v=\sum_{i=1}^nr_ib_i,~~\txt{for some}~~r_i\in k,~b_i\in B,~~\Ra~~0=(-1)b+\sum_{i=1}^nr_ib_i,\nn
    \eea
which is a contradiction.
\item\ul{(b)$\Ra$(c)}: Let $B$ be a maximal linearly independent set. To show $B$ spans $V$, let $v\in V$. If $v\in B$, then $v\in\txt{span}B$. If $v\not\in B$, then $B\cup\{v\}$ is dependent (since $B$ is a maximal linearly independent set), i.e.,
\bea
\textstyle r_0v+\sum_{i=1}^nr_ib_i=0~~~~\txt{for some $r_0,r_1...,r_n \in k$ that are not all $0$}.\nn
\eea
Since $B$ is linearly independent, $r_0\neq 0$ (otherwise $r_1=\cdots=r_n=0$ as well), and so
\bea
\textstyle v=-{1\over r_0}\sum_{i=1}^nr_ib_i~~~\in~\txt{span}B.\nn
\eea
To see that $B$ a minimal spanning set, suppose for some $b\in B$, the set $B\backslash\{b\}$ still spans $V$. Then
\bea
\textstyle b=\sum_{i=1}^nr_1b_1~~~~\txt{for some}~~~~r_i\in k,~~b_i\in B\backslash\{v\},~~~~\Ra~~~~0=(-1)b+\sum_{i=1}^nr_ib_i,\nn
\eea
which shows $B$ is not independent (a contradiction).
\item\ul{(c)$\Ra$(a)}: Let $B$ be a minimal spanning set for $V$. Suppose $B$ is dependent (i.e., not independent). Consider any relation $\sum_{i=1}^n r_ib_i=0$, where $r_i\in k$ and $b_i\in B$, such that (wlog) $r_1\neq 0$. Then
\bea
\textstyle b_1=-{1\over r_1}\sum_{i=2}^nr_ib_i\in \txt{span}(B\backslash\{b_1\}),\nn
\eea
and so $B\backslash\{b_1\}\subsetneq B$ spans $V$ (a contradiction since $B$ is a minimal spanning set). \qedhere
\eit
\end{proof}

\begin{thm}[\textcolor{blue}{Basis for a vector space}]
Every vector space has a basis.
\end{thm}
\begin{proof}
Let $V$ be a $k$-vector space. The case $V=0$ is trivial (since $V=0=\Span_k\emptyset$). So let $V\neq 0$. Then we get a \ul{nonempty} poset under inclusion $\subset$ given by $\L:=\{\txt{linearly independent subsets $X\subset V$}\}$. Let $\C\subset \L$ be a chain (linearly ordered subset). Let $U:=\bigcup_{X\in \C}X$. For any finite set $\{u_1,...,u_n\}\subset U$, let $u_i\in X_i\in \C$, where wlog $X_1\subset \cdots\subset X_n$ (since $\C$ is a chain). Then $\{u_1,...,u_n\}\subset X_n\in \C$. Thus, $\{u_1,...,u_n\}\subset U$ (and hence $U$) is linearly independent. Therefore, $U$ is an upper bound of $\C$ in $\L$.

By Zorn's Lemma, $\L$ has a maximal element $B$ (a basis for $V$ by the preceding theorem).
\end{proof}

\begin{thm}[\textcolor{blue}{Maximal submodule for a nonzero finitely generated module}]\label{FinGenMSub}
Every nonzero finitely generated module has a maximal submodule.
\end{thm}
\begin{proof}
Let $M$ be a nonzero finitely generated $R$-module. Let $\S:=\{\txt{proper submodules $_RN\subsetneq M$}\}$ as a poset under inclusion $\subset$. Then $\S$ is nonempty since $0\in \S$. Let $\C\subset \S$ be a chain (linearly ordered subset). Let $U:=\bigcup_{N\in \C}N$. If $x,y\in U$, let $x\in N_1$ and $y\in N_2$ for some $N_1,N_2\in\C$, where wlog $N_1\subset N_2$. Then for $r\in R$, $x,y\in N_2$ implies $x-y,rx\in N_2\subset U$, and so $U$ is a submodule.

To show $U\in \S$ (i.e., $U\subsetneq M$), observe that $M$ (being finitely generated) can be written as
\bea
M=Rx_1+Rx_2+\cdots +Rx_n,~~~~\txt{for some}~~x_1,...,x_n\in M.\nn
\eea
Suppose on the contrary that $U=M$. Then with $x_i\in N_i\in \C$ (for each $i=1,...,n$), where wlog $N_1\subset N_2\subset\cdots\subset N_n$, we get $x_1,...,x_n\in N_n\in \C\subset \S$. This implies
\bea
M=Rx_1+\cdots +Rx_n\subset N_n,~~\Ra~~N_n\not\in \S~~~~\txt{(a contradiction)}.\nn
\eea

It is clear by its definition that $U$ is an upper bound of $\C$ in $\S$. By Zorn's lemma, $\S$ contains a maximal element.
\end{proof}

\begin{crl}[\textcolor{blue}{Containing maximal submodule}]
If $M$ is a finitely generated module, then every proper submodule $W\subsetneq M$ is contained in a maximal submodule.
\end{crl}
\begin{proof}
Let $\S$ be the poset (under $\subset$) consisting of proper submodules of $M$ containing $W$. Then it is clear that we can proceed exactly as in the proof of Theorem \ref{FinGenMSub} (since the proof steps/conclusions are not affected by the containment of $W$ by the submodules).
\end{proof}

\begin{crl}[\textcolor{blue}{Containing maximal ideal}]
Let $R$ be a ring. Every proper ideal $I\vartriangleleft R$ is contained in a maximal ideal.
\end{crl}
\begin{proof}
Let $\I:=\{\txt{proper ideals $J\vartriangleleft R$ containing $I$}\}=\{J\vartriangleleft R:I\subset J,~1\not\in J\}$, as a poset under inclusion $\subset$. Let $\C\subset \I$ be a chain in $\I$, and let ~$U:=\bigcup_{J\in \C}J$.

It is clear that $I\subset U$ (so in particular, $U\neq\emptyset$). For $a,b\in U$ and $r\in R$, we have $a\in J_a$, $b\in J_b$ for some $J_a,J_b\in \C$ (where wlog, $J_a\subset J_b$), which implies $ra,ar,ab,a-b\in J_b\subset U$. Hence $U\vartriangleleft R$. Moreover, $1\not\in R$, i.e., $U$ is a proper ideal (otherwise $1\in J_i$ for some $i$, a contradiction). Thus $U$ is an upper bound of $\C$ in $\I$. By Zorn's Lemma, $\I$ contains a maximal element.
\end{proof}

\section{Quotient objects and Quotient structures}
\begin{dfn}[\textcolor{blue}{
\index{Quotient! group}{Quotient group},
\index{Natural! group homomorphism}{Natural group homomorphism}}]
Let $G=[G,\cdot,e,()^{-1}]$ be a group and $N\lhd G$ a normal subgroup. The \ul{quotient group} of $G$ by $N$ is the set of cosets
\bea
\textstyle{G\over N}:=\{gN~|~g\in G\}\nn
\eea
as a group with multiplication $\cdot:{G\over N}\times{G\over N}\ra{G\over N}$, inverse $()^{-1}:{G\over N}\ra{G\over N}$, and identity given respectively by
\bea
(gN)(g'N):=gg'N,~~~~(gN)^{-1}:=g^{-1}N,~~~~e_{G\over N}:=eN=N.\nn
\eea
The \ul{natural group homomorphism} associated with $N$ is the map ~$\pi:G\ra {G\over N}$, $g\mapsto gN$.
\end{dfn}

\begin{dfn}[\textcolor{blue}{
\index{Quotient! ring}{Quotient ring},
\index{Natural! ring homomorphism}{Natural ring homomorphism}}]
Let $R=[R,\cdot,+,0,1]$ be a ring and $I\lhd R$ an ideal. The \ul{quotient ring} of $R$ by $I$ is the set of cosets
\bea
\textstyle{R\over I}:=\{a+I~|~a\in R\}\nn
\eea
as a ring with multiplication {\small $\cdot:{R\over I}\times{R\over I}\ra{R\over I}$}, addition {\small $+:{R\over I}\times{R\over I}\ra{R\over I}$}, zero, and unity given respectively by
\bea
(a+I)(b+I):=ab+I,~~~~(a+I)+(b+I):=(a+b)+I,~~~~0_{R\over I}:=0+I=I,~~~~1_{R\over I}:=1+I.\nn
\eea
The \ul{natural ring homomorphism} associated with $I\lhd R$ is the map ~$\pi:R\ra {R\over I}$, $a\mapsto a+I$.
\end{dfn}

\begin{dfn}[\textcolor{blue}{
\index{Quotient! module}{Quotient module},
\index{Natural! $R$-homomorphism}{Natural $R$-homomorphism}}]
Let $M={}_RM$ be an $R$-module and ${}_RN\subset{}_RM$ a submodule. The \ul{quotient $R$-module} of $M$ by $N$ is the set of cosets
\bea
\textstyle{M\over N}:=\{m+N~|~m\in M\}\nn
\eea
as an $R$-module with scaling {\small $R\times{M\over N}\ra{M\over N}$}, addition {\small $+:{M\over N}\times{M\over N}\ra{M\over N}$}, and zero given respectively by
\bea
r(m+N):=rm+N,~~~~(m+N)+(m'+N):=(m+m')+N,~~~~0_{M\over N}:=0+N=N.\nn
\eea
The \ul{natural $R$-homomorphism} associated with $_RN\subset M$ is the map ~$\pi:M\ra {M\over N}$, $m\mapsto m+N$.
\end{dfn}

\begin{dfn}[\textcolor{blue}{
\index{Quotient! algebra}{Quotient algebra},
\index{Natural! algebra homomorphism}{Natural algebra homomorphism}}]
Let $A$ be an $R$-algebra and $I\lhd A$ an ideal. The \ul{quotient $R$-algebra} of $A$ by $I$ is the quotient ring and quotient $R$-module (\blue{footnote}\footnote{The center of $A/I$ is $Z(A/I)={C_I(A)+I\over I}$, where $C_I(A):=\{a\in A~|~aa'-a'a\in I~\txt{for all}~a'\in A\}\supset Z(A)$, which is in general larger than the \ul{induced center} ${Z(A)+I\over I}$.})
\bea
\textstyle{A\over I}:=\{a+I~|~a\in A\}~~~~\txt{viewed as an ~${R+I\over I}=R(1+I)$-algebra}.\nn
\eea
The \ul{natural $R$-algebra homomorphism} associated with $I\lhd A$ is the map ~$\pi:A\ra {A\over I}$, $a\mapsto a+I$.
\end{dfn}

\begin{examples}[\blue{
\index{Tensor! algebra}{Tensor algebra},
\index{Grassmann (exterior) algebra}{Grassmann (exterior) algebra},
\index{Clifford algebra}{Clifford algebra},
\index{Weyl algebras}{Weyl algebras},
\index{Moyal algebras}{Moyal algebras}}]
Let $R$ be a commutative ring and $M$ an $R$-module. The \ul{tensor $R$-algebra} $Te(M)$ of $M$ is the fully noncommutative polynomial $R$-algebra over $M$ given (with $M^i:=M\times M^{i-1}$) by the set
\[
\textstyle Te(M):=R\langle M\rangle:=\sum_{i\in\Natural}RM^i=\{\sum_{i\in\Natural}r_ib_i:r_i\in R,~b_i\in M^i,~r_i=0~\txt{a.e.f.}\},~~~~M^0:=R,~~M^1=M,\nn
\] as a ring with addition $+:Te(M)\times Te(M)\ra Te(M)$ and multiplication $Te(M)\times Te(M)\ra Te(M)$ given (along with zero $0_{T_e(R)}:=\sum_i0_{M^i}b_i$ and unity $1_{Te(M)}:=1_R=\sum_i1_R\delta_{0i}b_i$) respectively by
\[
\textstyle\sum_ir_ib_i+\sum_ir'_ib'_i:=\sum_i(r_ib_i+r'_ib'_i),~~~~\left(\sum_ir_ib_i\right)\left(\sum_ir'_ib'_i\right):=\sum_ic_i,~~~~c_i:=\sum_{j+k=i}r_jr_k'b_jb'_k\in RM^i,
\]
where the multiplication of $b_n:=m_1m_2\cdots m_n:=(m_1,...,m_n)\in M^n$ and $b_{n'}':=m_1'm_2'\cdots m_{n'}':=(m_1',...,m_{n'}')\in M^{n'}$ is given by concatenation of words as follows (\blue{foonote}\footnote{In the tensor product context, the finite sequence/string ~$b_n:=m_1m_2\cdots m_n:=(m_1,...,m_n)\in M^n$~ is often also written as $m_1\otimes m_2\otimes\cdots\otimes m_n$, but we will not do this, and instead reserve the symbol $\otimes$ for later use.}):
\[
b_nb_{n'}'=(m_1,...,m_n)(m_1',...,m_{n'}'):=(m_1,...,m_n,m_1',...,m_{n'}')\in M^{n+n'},\nn
\]
and each $M^i:=\{(m_1,...,m_i):m_j\in M\}$ is viewed as an $R$-module with addition and scalar multiplication
\[
(m_1,...,m_i)+(m'_1,...,m'_i):=(m_1+m'_1,...,m_i+m'_i),~~~~r(m_1,...,m_i):=(rm_1,...,rm_i).
\]

Consider the ideal $I:=Te(M) M(2)\lhd Te(M)$ generated by elements of the form $M(2):=\{a^2:=(a,a)\in M^2~|~a\in M\}\subset R\langle M\rangle$. The \ul{Grassmann $R$-algebra} (or \ul{exterior $R$-algebra}) of $M$ is
{\small\[
\textstyle\bigwedge(M):={Te(M)\over I}={R\langle M\rangle\over R\langle M\rangle\{a^2:a\in M\}}=R\langle\{a+I:a\in M\}\rangle,~~~~\txt{where}~~~~a^2+I=I~~\txt{for all}~~a\in M.
\]}Observe that $a^2+I=I$ for all $a\in M$ $\Ra$ $(a+b)^2+I=I$ for $a,b\in M$ $\Ra$ $ab+I=-ba+I$ for all $a,b\in M$, but the converse is not always true (\blue{footnote}\footnote{For an ideal $J\lhd R\langle M\rangle$ such that $ab+J=-ba+J$ for all $a,b\in M$, we also have $2a^2+J=J$ for all $a\in M$. However, in general we must have $2\in U(R)$ in order to conclude that $a^2+J=J$ for all $a\in M$.}).

Given a map $Q:M\ra R$, again with $Te(M):=R\langle M\rangle$, consider the ideal $I:=Te(M)M(Q)\lhd Te(M)$ generated by elements of the form $M(Q):=\{a^2-Q(a):a\in M\}\subset R\langle M\rangle$. The \ul{Clifford $R$-algebra} of $M$ (associated with $Q$) is
{\small\[
\textstyle Cl(M,Q):={Te(M)\over I}={R\langle M\rangle\over R\langle M\rangle\{a^2-Q(a):a\in M\}}=R\langle\{a+I:a\in M\}\rangle,~~~~\txt{where}~~~~a^2+I=Q(a)+I~~\txt{for all}~~a\in M.
\]}Observe that with ~$Q:=0:M\ra R,~m\mapsto 0$,~ we get ~$Cl(M,0)=\bigwedge(M)$.

Given a finite subset $X_n:=\{x_1,...,x_n,y_1,...,y_n\}\subset M$, let $I_n:=R\langle X_n\rangle D_n\lhd R\langle X_n\rangle\leq Te(M)$ be the ideal of the subring $R\langle X_n\rangle\leq Te(M)$ generated by the set $D_n:=\{x_iy_i-y_ix_i-1_R\}_{i=1,...,n}\cup\{x_iy_j-y_jx_i\}_{i\neq j}\subset R\langle X_n\rangle$. The $n$th \ul{Weyl $R$-algebra} is
{\small\[
\textstyle W_n:={R\langle X_n\rangle\over I_n}={R\langle X_n\rangle\over R\langle X_n\rangle D_n}=R\langle\{x_1+I_n,...,x_n+I_n,y_1+I_n,...,y_n+I_n\}\rangle,~~\txt{where}~~x_iy_j-y_jx_i+I_n=\delta_{ij}1_R+I_n.\nn
\]}

Given a finite subset $X_n:=\{x_1,...,x_n\}\subset M$ and elements $\{\theta_{ij}:1\leq i,j\leq n\}\subset R$, let $J_n:=R\langle X_n\rangle E_n\lhd R\langle X_n\rangle\leq Te(M)$ be the ideal of the subring $R\langle X_n\rangle\leq Te(M)$ generated by the set $E_n:=\{x_ix_j-x_jx_i-\theta_{ij}\}_{i=1,...,n}\subset R\langle X_n\rangle$. The $n$th \ul{Moyal $R$-algebra} is
\[
\textstyle M_n:={R\langle X_n\rangle\over J_n}={R\langle X_n\rangle\over R\langle X_n\rangle E_n}=R\langle\{x_1+J_n,...,x_n+J_n\}\rangle,~~\txt{where}~~x_ix_j-x_jy_i+J_n=\theta_{ij}+J_n.\nn
\]

\end{examples}

\section{Classical Polynomial Rings (a concretizing digression)}
\subsection{Division in polynomial rings}
\begin{dfn}[\textcolor{blue}{\index{Degree (of polynomial)}{Degree of a polynomial}, \index{Leading coefficient}{Leading coefficient} of a polynomial}]
Let $R$ be a commutative ring. Then a nonzero polynomial $f\in R[x]\backslash\{0\}$ can be written as
\bea
f=a_0+a_1x+a_2x^2+...+a_nx^n,~~~~\txt{with}~~a_n\neq 0,\nn
\eea
where $n$ is called the \ul{degree} of $f$, and $a_n$ is called the \ul{leading coefficient} of $f$ (denoted by $l(f):=a_n$).
\end{dfn}

\begin{lmm}[\textcolor{blue}{Division algorithm for polynomials}]
Let $R$ be a commutative ring and $f\in R[x]$. If $l(f)\in U(R)$ (in particular, $f\neq 0$), then every $g\in R[x]$ can be written as
\bea
g=qf+r,~~~~\txt{for some}~~q,r\in R[x],~~~~\txt{where either}~~r=0~~\txt{or}~~\deg r<\deg f.\nn
\eea
\end{lmm}
\begin{proof}
If $g=0$, or $\deg g<\deg f$, we can set $q=0$ and $r=g$. So let $\deg g\geq\deg f$ (in particular, $g\neq 0$ since $f\neq 0$). We will proceed by induction on $\deg g$. Let
\bea
g=b_0+b_1x+...+b_mx^m,~~~~f=a_0+a_1x+...+a_nx^n,~~~~m\geq n,~~a_n\in U(R).\nn
\eea
Then $h:=g-b_mx^{m-n}a_n^{-1}f\in R[x]$ has degree $\deg h<m$. Therefore, by induction hypotheses,
\bea
g-b_mx^{m-n}a_n^{-1}f=q_1f+r,~~~\txt{for some}~~~q_1,r\in R[x],~~\txt{where either}~~r=0~~\txt{or}~~\deg r<\deg f.\nn
\eea
Hence ~$g=[b_mx^{m-n}a_n^{-1}+q_1]f+r=qf+r$, ~where ~$q:=b_mx^{m-n}a_n^{-1}+q_1$.
\end{proof}

\begin{crl}[\textcolor{blue}{If $k$ is a field, then the ring $k[x]$ is a ED, and hence a PID}]
\end{crl}

\begin{lmm}
Let $R$ be a commutative ring, and $f,g\in R[x]\backslash\{0\}$ with leading coefficients $a$ and $b$ respectively. If $ab\neq 0$, then $\deg(fg)=\deg f+\deg g$.
\end{lmm}
\begin{proof}
$fg=\cdots+(ax^{\deg f})(bx^{\deg g})=\cdots+abx^{\deg f+\deg g}$.
\end{proof}

\begin{prp}
Let $R$ be an ID. Then $R[x]$ is an ID $\iff$ {\small $\deg(fg)=\deg f+\deg g$} for all $f,g\in R[x]$.
\end{prp}
\begin{proof}
($\Ra$) Assume $R$ be an ID. Let $f,g\in R[x]$, where $f=a_0+a_1x+\cdots+a_mx^m$ and $g=b_0+b_1x+\cdots+b_nx^n$. Then $fg=a_0b_0+\cdots+a_mb_nx^{m+n}$, where $a_mb_n\neq 0$ since $R$ is an ID. Hence $\deg (fg)=\deg f+\deg g$. ($\La$) Conversely, assume $\deg (fg)=\deg f+\deg g$ for all $f,g\in R[x]$. Then for any $f=a_0+a_1x+\cdots+a_mx^m$ and $g=b_0+b_1x+...+b_nx^n$ in $R[x]$, if $fg=0$ then $a_mb_n=0$, $\Ra$ $a_m=0$~ or $b_n=0$ (since $R$ is an ID),
$\Ra$ $f=0$ or $g=0$.
\end{proof}

\begin{thm}[\textcolor{blue}{Polynomial ID's}]
Let $R$ be a commutative ring.
\begin{enumerate}[leftmargin=0.9cm]
\item $R$ is an ID $\iff$ $R[x]$ is an ID ~~~~(~$\iff$~ $\deg(fg)=\deg f+\deg g$ for all $f,g\in R[x]$~).
\item $R$ is an ID $\iff$ $R[[x]]$ is an ID. (Here, $R[[x]]$ is the \index{Ring! of power series}{\ul{ring of power series}} over $R$ (see \blue{footnote}\footnote{Here $R[[x]]$, called the \ul{ring of power series} over $R$, is defined to be the set of sequences in $R$, $R^\Natural:=\{\txt{maps}~f:\Natural\ra R\}$,  along with addition, multiplication, zero, and unity given (pointwise for each $i\in\Natural$) by
\[
\textstyle (f+g)(i):=f(i)+g(i),~~~~ (fg)(i):=\sum_{j+k=i}f(j)g(k),~~~~0_{R[[x]]}(i)=0_R,~~~~1_{R[[x]]}(i)=1_R\delta_{0i}.\nn
\]
Given $f\in R^\Natural$, we formally write ~$f=\sum_{i=0}^\infty f(i)x^i:\Natural\ra R$ ~in analogy with elements of the polynomial ring $R[x]$.}).)
\end{enumerate}
\end{thm}
\begin{proof}
{\flushleft (1)} ($\Ra$): If $R$ is an ID, let $f,g\in R[x]\backslash\{0\}$, where $f=a_0+\cdots+a_nx^n$ and $g=b_0+\cdots+b_mx^m$ (with $a_n\neq 0$, $b_m\neq 0$). Then $fg=\cdots+a_nb_mx^{m+n}\neq 0$, since $a_nb_m\neq 0$. ($\La$): The converse is clear since $R\subset R[x]$ is a subring.

{\flushleft (2)} The proof is similar to that of part (1). ($\Ra$): If $R$ is an ID, let $f,g\in R[[x]]\backslash\{0\}$, where $f=a_ix^i+a_{i+1}x^{i+1}+\cdots$ and $g=b_jx^j+b_{j+1}x^{j+1}+\cdots$ (such that $a_i\neq 0$, $b_j\neq 0$ are minimal with this property). Then $fg=a_ib_jx^{i+j}+\cdots\neq 0$. ($\La$): The converse is clear since $R\subset R[[x]]$ is a subring.
\end{proof}

\begin{crl}[\textcolor{blue}{Polynomial ID's}]
Let $R$ be commutative. Then $R$ is an ID $\iff$ {\small $R[x_1,...,x_n]$} is an ID.
\end{crl}
\begin{proof}
($\Ra$) If $R$ is an ID, so is $R[x]$, and the rest follows by induction since {\footnotesize $R[x_1,...,x_n]=R[x_1,..,x_{n-1}][x_n]$}. ($\Ra$) This is clear because $R\subset R[x_1,...,x_n]$ is a subring.
\end{proof}

\subsection{Factorization of polynomials}
\begin{dfn}[\textcolor{blue}{\index{Primitive polynomial}{Primitive polynomial}}]
Let $R$ be an ID, and $f\in R[x]$. Then $f$ is \ul{primitive} (or \ul{$R$-irreducible}) if the only elements of $R$ that divide $f$ are units, i.e., for all $r\in R$, if $r|f$ (i.e., $f=rh$) then $r\in U(R)$. Equivalently, because $r|f$ $\iff$ $r$ divides every nonzero coefficient of $f$, $f\in R[x]$ is primitive $\iff$ $f$ has a coefficient that is irreducible in $R$. Equivalently, $f=a_0+a_1x+\cdots+a_nx^n\in R[x]$ is primitive $\iff$ any ~$\txt{gcd}(\txt{coefficients of $f$}):=gcd(a_0,...,a_n)\in U(R)$.
\end{dfn}

\begin{lmm}(\textcolor{blue}{Nonzero polynomials over a UFD are primitive up to scaling})
Let $R$ be a UFD and $R_0[x]$ the set of primitive polynomials in $R[x]$. Then ~$R[x]~=~R~R_0[x]$, i.e., for every $f\in R\backslash\{0\}$,~ we have $f=cf_0$ for some $c\in R$ and a primitive $f_0\in R[x]$.
\end{lmm}
\begin{proof}
Let $c:=\txt{gcd}(\txt{coefficients of $f$})$. Then $f=cf_0$, were it is clear that $f_0$ is primitive since all common divisors of coefficients of $f$ are already contained in $c$, and so $\txt{gcd}(\txt{coefficients of $f_0$})\in U(R)$.
\end{proof}

\begin{lmm}[\textcolor{blue}{A product of primitive polynomials is primitive}]
Let $R$ be a UFD. If $f,g\in R[x]$ are primitive, so is the product $fg$, i.e., $R_0[x]R_0[x]\subset R_0[x]$.
\end{lmm}
\begin{proof}
Let $f,g\in R[x]$ be primitive. If $fg$ is not primitive, then $fg=rh$, for some irreducible $r\in R\backslash(\{0\}\cup U(R))$. Since $r$ is prime in $R$ (hence prime in $R[x]$), $r|fg~\Ra~r|f~\txt{or}~r|g$ (a contradiction).
\end{proof}

\begin{dfn}[\textcolor{blue}{\index{Ring! of quotients for an ID}{Ring of quotients for an ID}}]
Let $R$ be an ID. (As done earlier in section \ref{IntRatNumSec} for $R=\Integer$ to obtain $\Rational=Q(\Integer)$) consider the equivalence relation $\sim$ on $R\times R$ given by
\bea
\textstyle(a,b)\sim (a',b')~~~~\txt{if}~~~~b\neq 0,~~b'\neq 0,~~\txt{and}~~ab'=a'b.\nn
\eea
Then, in terms of the equivalence classes ${a\over b}:=[(a,b)]$, the \ul{ring of quotients} (\ul{fractions of elements of $R$}) is
\bea
\label{ROQdefID}\textstyle Q(R):={R\times R\over\sim}=\left\{{a\over b}~|~a,b\in R,~b\neq 0\right\}
\eea
with addition, multiplication, ``zero'' (i.e., $0_{Q(R)}$), and ``one'' (i.e., $1_{Q(R)}$) given respectively by
\bea
\textstyle{a\over b}+{a'\over b'}:={ab'+a'b\over bb'},~~~~{a\over b}\cdot{a'\over b'}:={aa'\over bb'},~~~~0_{Q(R)}:={0\over 1_R},~~~~1_{Q(R)}:={1_R\over 1_R}.\nn
\eea
\end{dfn}
Note that for an arbitrary commutative ring $R$, its ring of quotients can also be constructed, but the presence of zerodivisors in $R$ requires a (slight) modification of the equivalence relation to be considered on $R\times R$, and this is done in chapter \ref{AlgebraCatS8}.

\begin{thm}[\textcolor{blue}{Gauss lemma I}] Let $R$ be a UFD. Let $f_0,g\in R[x]$, where $f_0$ is primitive. If $f_0|g$ in $Q(R)[x]$, then $f_0|g$ in $R[x]$. Equivalently, if $f_0\in R[x]$ is primitive, then
\bea
\label{gauss-lemma-eq1}\textstyle R[x]\cap\big(Q(R)[x]f_0\big)~\subset~R[x]f_0.
\eea
\end{thm}
\begin{proof}
Let $g=f_0h$ for some $h\in Q(R)[x]$. Since $h$ has finitely many coefficients, $ch\in R[x]$ for some $c\in R\backslash 0$ (where multiplication by $c$ clears denominators in $h$). Therefore,
\bea
cg=cf_0h=f(ch)=f_0l\in R[x],~~~~l:=ch\in R[x].\nn
\eea
Either (i) $c\in U(R)$, or (ii) $c\not\in U(R)$. Suppose $c\not\in U(R)$. Let $c=c_1c_2\cdots c_n$, with each $c_i$ irreducible in $R$ (hence irreducible in $R[x]$), such that $n$ is \ul{minimal} (i.e., just enough so that $l:=ch\in R[x]$). Then for each $1\leq i\leq n$, we have ${l\over c_i}={ch\over c_i}\not\in R[x]$.

Since $c_1|cg$ in $R[x]$ (and $cg=f_0l$), we also have $c_1|f_0l$ in $R[x]$. But $c_1$ is prime in $R[x]$, and $c_1\nmid f_0$ (since $f_0$ is primitive). Therefore $c_1|l$ in $R[x]$, and so
\bea
\textstyle c_2c_3\cdots c_nh={ch\over c_1}={l\over c_1}\in~R[x]~~~~\txt{(a contradiction since $n$ was minimal)}.\nn
\eea
Thus we conclude that  $c\in U(R)$. Hence $f_0|g$ in $R[x]$.
\end{proof}

\begin{thm}[\textcolor{blue}{Gauss lemma II}]
Let $R$ be a UFD. If (i) $g\in R[x]$ and (ii) $g=fh$ for some $f,h\in Q(R)[x]$, then there exists $\al\in Q(R)$ such that $\al f,\al^{-1}h\in R[x]$~ (so that $g=\al f ~\al^{-1}h$ in $R[x]$). (\blue{footnote}\footnote{That is (assuming $R$ is a UFD), if $g\in R[x]$ factorizes a certain way in $Q(R)[x]$, then $g$ also factorizes the same way in $R[x]$. Since the converse is also trivially true (because $R[x]\subset Q(R)[x]$), it follows that if $g\in R[x]$, then $g$ factorizes in $R[x]$ $\iff$ $g$ factorizes in $Q(R)[x]$.}).

Equivalently, this means that if $R$ is a UFD, then with $A\ast B:=\{ab:a\in A,~b\in B\}$,
\bea
\label{gauss-lemma-eq2}R[x]\cap\big(Q(R)[x]\ast Q(R)[x]\big)\subset R[x]\ast R[x].
\eea
\end{thm}
\begin{proof}
Choose $a\in R\backslash\{0\}$ such that $af\in R[x]$. Then
\bea
af=bf_0,~~~~\txt{where}~~b\in R~~\txt{and}~~f_0\in R[x]~~\txt{is primitive}.\nn
\eea
Now,
\bea
\textstyle g=fh={a\over b}f~{b\over a}h=f_0~{b\over a}h~~\left(f_0:={a\over b}f\right),~~\Ra~~f_0|g~~\txt{in}~~Q(R)[x].\nn
\eea
By the previous lemma, $f_0|g$ in $R[x]$, and so

\[\textstyle {b\over a}h~\in~R[x],~~~~\Ra~~~~\al f,\al^{-1}h~\in~R[x],~~~~\txt{where}~~~~\al:={a\over b}.\nn\qedhere\]
\end{proof}
\begin{proof}[Alternative proof]
Using the equivalent statement (\ref{gauss-lemma-eq2}), the result follows more easily from (\ref{gauss-lemma-eq1}) by the observation that for all $f\in R[x]$,
{\small\begin{align}
&R[x]\cap\big(Q(R)[x]f\big)=R[x]\cap\Big(Q(R)[x]cf_0\Big)=c\Big(R[x]\cap\big(Q(R)[x]f_0\big)\Big)
\sr{(\ref{gauss-lemma-eq1})}{\subset}~c(R[x]f_0)=R[x]cf_0=R[x]f,\nn\\
&\textstyle~~\Ra~~R[x]\cap\Big(Q(R)[x]\ast Q(R)[x]\Big)=R[x]\cap\Big(Q(R)[x]\ast R[x]\Big)=\bigcup_{f\in R[x]}R[x]\cap\Big(Q(R)[x]f\Big)\subset R[x]\ast R[x].\nn\qedhere
\end{align}}
\end{proof}

The above result implies that if $R$ is a UFD, then a primitive polynomial $f_0\in R[x]$ that is irreducible in $Q(R)[x]$ is also irreducible in $R[x]$.

\begin{crl}[\textcolor{blue}{Polynomial UFD's}]
If $R$ is a UFD, so is $R[x]$.
\end{crl}
\begin{proof}
Since in an ID (i) primes are irreducibles and (ii) factorization into primes is unique, it suffices to show that $f\in R[x]$ can be written as a finite product of primes. Let $f=cf_0$, where $c\in R$ and $f_0\in R[x]$ is primitive. Let
\[
 c=p_1p_2\cdots p_n,\nn
\]
where each $p_i\in R$ is irreducible (or prime since $R$ is a UFD), and
\[
f_0=g_1g_2\cdots g_l,~~~~\txt{(with $l$ finite as $\deg f_0$ is finite)}
\]
where each $g_i\in R[x]$ is primitive and irreducible in $R[x]$. It remains to show each $g_i$ is prime in $R[x]$.

Since $g_i\in R[x]$ is irreducible, we know $g_i\in Q(R)[x]$ is irreducible, hence $g_i$ is prime in $Q(R)[x]$ since $Q(R)[x]$ is a PID (as polynomials over fields). Since $R$ is a UFD and $g_i$ is primitive, for any $\gamma,\rho\in R[x]\subset Q(R)[x]$,
\bea
&& g_i|\gamma\rho~~~~\txt{in}~~~~Q(R)[x]\nn\\
&&~~\Ra~~g_i|\gamma~~\txt{or}~~g_i|\rho~~~~\txt{in}~~~~Q(R)[x]~~~~~~\txt{(since $g_i$ is prime in $Q(R)[x]$)}\nn\\
&&~~\Ra~~g_i|\gamma~~\txt{or}~~g_i|\rho~~~~\txt{in}~~~~R[x]~~~~~~\txt{(Gauss lemma)}.\nn
\eea
Hence $g_i$ is prime in $R[x]$.
\end{proof}

\begin{crl}[\textcolor{blue}{Polynomial UFD's}]
If $R$ is a UFD, so is $R[x_1,...,x_n]$.
\end{crl}
\begin{proof}
This follows by induction on $n$ since $R[x_1,...,x_n]=R[x_1,...,x_{n-1}][x_n]$.
\end{proof}

%% file: parts/AlgebraNC/AlgebraNC5.tex
\chapter{Basic Isomorphism and Correspondence Theorems}
Recall that in our definition of an $R$-algebra $A$, the subring $R\leq A$ is commutative and lies in the center of $A$, i.e., $R\subset Z(A)$. On the other hand, our rings $R$ in general (whether by themselves or on $R$-modules) are simply $\Integer$-algebras with no assumption of commutativity.

\section{Isomorphism and correspondence theorems for groups}
Here is a summary of the group isomorphism theorems we want to prove:
\begin{enumerate}
\item Given a group homomorphism $f:G\ra G'$, a normal subgroup $N\lhd G$, and the associated normal subgroup $M:=N\cap\ker f\lhd G$ with natural map $\pi_M:G\ra{G\over M},~g\mapsto gM$, we get (i) a group homomorphism $f_M:{G\over M}\ra G'$ with kernel $\ker f_M={\ker f\over M}$ such that $f=f_M\circ\pi_M:G\sr{\pi_M}{\ral}{G\over M}\sr{f_M}{\ral}G'$, and (ii) a group isomorphism ~${G\over\ker f}\cong\im f$.
\item Given subgroups $H,N\leq G$ with $N\lhd G$ normal, we get a group isomorphism ~${H\over H\cap N}\cong {HN\over N}$.
\item Given normal subgroups $N,N'\lhd G$, we get a group isomorphism ~${G/N\cap N'\over NN'/N\cap N'}\cong {G\over NN'}$.
\end{enumerate}
\begin{thm}[\textcolor{blue}{\index{First isomorphism theorem for! groups}{First isomorphism theorem for groups}}]
Let $f:G\ra G'$ be a group homomorphism, $N\lhd G$ a normal subgroup, $M:=N\cap\ker f$, and $\pi_M:G\ra{G\over M},~g\mapsto gM$ the natural homomorphism. Then there exists a group homomorphism $f_M:{G\over M}\ra G'$ such that (i) $f=f_M\circ\pi_M$ and (ii) $\ker f_M={\ker f\over M}={\ker f\over N\cap\ker f}$.
{\small\[
\textstyle\bt
G\ar[d,"f"]\ar[rr,"\pi_M"] && {G\over M}\ar[dll,dashed,"f_M"]\\
G'
\et~~~~~~~~\txt{In particular, with $N:=\ker f$, we get a group isomorphism ~~${G\over \ker f}\cong f(G)$}.\nn
\]}
\end{thm}
\begin{proof}
The map $f_M:{G\over M}\ra G',~gM\mapsto f(g)$ is well-defined, since
\bea
\textstyle gM=hM~~\Ra~~h^{-1}g\in M\subset\ker f,~~\Ra~~e_G=f(h^{-1}g)=f(h)^{-1}f(g),~~\Ra~~f(g)=f(h).\nn
\eea
Moreover, $f_M\circ\pi_M(g)=f_M(gM)=f(g)$, for all $g\in G$, and
\[
\textstyle \ker f_M=\big\{gM:e_{G'}=f_M(gM)=f(g)\big\}=\big\{gM:g\in\ker f\big\}={(\ker f)M\over M}={\ker f\over M}.\nn \qedhere
\]
\end{proof}

\begin{thm}[\textcolor{blue}{\index{Second isomorphism theorem for! groups}{Second isomorphism theorem for groups}}]
Let $H,N\leq G$ be subgroups with $N\lhd G$ normal. Then (i) $H\cap N\lhd H$, (ii) $N\lhd HN$, and (iii) ${H\over H\cap N}\cong{HN\over N}$.
\end{thm}
\begin{proof}
(i) For any $h\in H$, we have $h(H\cap N)h^{-1}=(hHh^{-1})\cap(hNh^{-1})=H\cap N$. (ii) For any $hn\in HN$, we have $hn(N)(hn)^{-1}=N$. (iii) The surjective group homomorphism $f:H\ra {HN\over N}$, $h\mapsto hN$ has kernel
\[
\textstyle \ker f=\{h\in H:hN=N\}=\{h\in H:h\in N\}=H\cap N. \qedhere
\]
\end{proof}

\begin{thm}[\textcolor{blue}{\index{Third isomorphism theorem for! groups}{Third isomorphism theorem for groups}}]
Let $N,N'\lhd G$ be normal subgroups. Then (i) $N\cap N',NN'\lhd G$, (ii) $N\cap N'\lhd NN'$, and (iii) ${G/N\cap N'\over NN'/N\cap N'}\cong{G\over NN'}$.

In particular, if $N\subset N'$, then ~${G/N\over N'/N}\cong{G\over N'}$.
\end{thm}
\begin{proof}
The surjective map $f:{G\over N\cap N'}\ra{G\over NN'}$, $gN\cap N'\mapsto gNN'$
is well defined because
\bea
gN\cap N'=hN\cap N'~~\Ra~~h^{-1}g\in N\cap N'\subset NN',~~\Ra~~gNN'=hNN',\nn
\eea
and it is a group homomorphism with kernel
\[
\textstyle \ker f=\{gN\cap N':gNN'=NN'\}=\{gN\cap N':g\in NN'\}={NN'\over N\cap N'}.\nn \qedhere
\]
\end{proof}

\begin{thm}[\textcolor{blue}{\index{Correspondence theorem for! groups}{Correspondence theorem for groups}}]
Let $N\lhd G$ be a normal subgroup, $Sub_N(G):=\{\txt{subgroups of $G$ containing $N$}\}$, and $Sub\left({G\over N}\right):=\{\txt{subgroups of ${G\over N}$}\}$. Then there exists a bijection
\bea
\bt Sub_N(G)\ar[rr,leftrightarrow] && Sub\left({G\over N}\right).\et\nn
\eea
\end{thm}
\begin{proof}
Let $\pi_N:G\ra {G\over N},~g\mapsto gN$ be the natural group homomorphism. Consider the maps
\bea
\textstyle\Phi:Sub_N(G)\ra Sub\left({G\over N}\right),~H\mapsto\pi_N(H)={H\over N},~~~~\Psi:Sub\left({G\over N}\right)\ra Sub_N(G),~H'\mapsto\pi_N^{-1}(H').\nn
\eea
Then for all ~$H\in Sub_N(G)$ and $H'\in Sub\left({G\over N}\right)$,
\begin{align}
&\textstyle\Phi\circ\Psi(H')=\pi_N\left(\pi_N^{-1}(H')\right)=H'\cap\pi_N(G)\sr{(a)}{=}H',~~\txt{where step (a) holds because $\pi_N$ is surjective},\nn\\
&\textstyle\Psi\circ\Phi(H)=\pi_N^{-1}\big(\pi_N(H)\big)=\pi_N^{-1}\left({H\over N}\right)\sr{(b)}{=}HN\sr{(c)}{=}H,~~\txt{where step (c) holds because $N\subset H$},\nn
\end{align}
and step (b) is a general result obtained as follows:
{\small\begin{align}
&\textstyle \pi_N^{-1}\left({H\over N}\right)=\left\{g\in G:gN\in {H\over N}\right\}=\left\{g\in G:gN=hN,~\txt{some}~h\in H\right\}=\left\{g\in G:g\in hN,~\txt{some}~h\in H\right\}\nn\\
&\textstyle~~~~=\left\{g\in G:g\in HN\right\}=HN.\nn \qedhere
\end{align}}
\end{proof}

\begin{crl}
A normal subgroup $N\lhd G$ is maximal $\iff$ $G/N$ is a simple group.
\end{crl}

\section{Isomorphism and correspondence theorems for rings}
Here is a summary of the ring isomorphism theorems we want to prove:
\begin{enumerate}
\item Given a ring homomorphism $f:R\ra R'$, an ideal $I\lhd R$, and the associated ideal $J:=I\cap\ker f\lhd R$ with natural map $\pi_J:R\ra{R\over J},~r\mapsto r+J$, we get (i) a ring homomorphism $f_J:{R\over J}\ra R'$ with kernel $\ker f_J={\ker f\over J}$ such that $f=f_J\circ\pi_J:R\sr{\pi_J}{\ral}{R\over J}\sr{f_J}{\ral}R'$, and (ii) a ring isomorphism ${R\over\ker f}\cong\im f$.
\item Given a subring $S\leq R$ and an ideal $I\vartriangleleft R$, we get a ring isomorphism ~${S\over S\cap I}\cong {S+I\over I}$.
\item Given ideals $I,J\lhd R$ of a ring $R$, we get a ring isomorphism ~${R/I\cap J\over (I+J)/I\cap J}\cong {R\over I+J}$.
\end{enumerate}

\begin{thm}[\textcolor{blue}{\index{First isomorphism theorem for! rings}{First isomorphism theorem for rings}}]
Let $f:R\ra R'$ be a ring homomorphism, $I\lhd R$ an ideal, $J:=I\cap\ker f$, and $\pi_J:R\ra{R\over J},~r\mapsto r+J$ the natural ring homomorphism. Then there exists a ring homomorphism $f_J:{R\over J}\ra R'$ such that (i) $f=f_J\circ\pi_J$ and (ii) $\ker f_J={\ker f\over J}={\ker f\over I\cap\ker f}$.
{\small\[
\textstyle\bt
R\ar[d,"f"]\ar[rr,"\pi_J"] && {R\over J}\ar[dll,dashed,"f_J"]\\
R'
\et~~~~~~~~\txt{In particular, with $I:=\ker f$, we get a ring isomorphism ~~${R\over \ker f}\cong f(R)$}.\nn
\]}
\end{thm}
\begin{proof}
The map $f_J:{R\over J}\ra R',~r+J\mapsto f(r)$ is well defined, since
\bea
\textstyle r+J=s+J~~\Ra~~r-s\in J\subset\ker f~~\Ra~~0=f(r-s)=f(r)-f(s),~~\Ra~~f(r)=f(s).\nn
\eea
Moreover, $f_J\circ\pi_J(r)=f_J(r+J)=f(r)$, for all $r\in R$, and
\[
\textstyle \ker f_J=\big\{r+J:0=f_J(r+J)=f(r)\big\}=\big\{r+J:r\in\ker f\big\}={\ker f+J\over J}={\ker f\over J}.\nn \qedhere
\]
\end{proof}

\begin{thm}[\textcolor{blue}{\index{Second isomorphism theorem for! rings}{Second isomorphism theorem for rings}}]
Let $S\leq R$ be a subring and $I\lhd R$ an ideal. Then (i) $S\cap I\lhd S$, (ii) $I\lhd S+I\leq R$, and (iii) ${S\over S\cap I}\cong{S+I\over I}$.
\end{thm}
\begin{proof}
The surjective ring homomorphism $f:S\ra {S+I\over I}$, $s\mapsto s+I$ has kernel
\[
\textstyle \ker f=\{s\in S:s+I=I\}=\{s\in S:s\in I\}=S\cap I.\nn \qedhere
\]
\end{proof}

\begin{thm}[\textcolor{blue}{\index{Third isomorphism theorem for! rings}{Third isomorphism theorem for rings}}]
Let $I,J\lhd R$ be ideals (in a ring $R$). Then (i) $I\cap J,I+J\lhd R$, (ii) $I\cap J\lhd I+J$, and (iii) ${R/I\cap J\over(I+J)/I\cap J}\cong{R\over I+J}$.

In particular, if $I\subset J$, then ~${R/I\over J/I}\cong{R\over J}$.
\end{thm}
\begin{proof}
The surjective map $f:{R\over I\cap J}\ra{R\over I+J}$, $r+I\cap J\mapsto r+(I+J)$ is well defined because
\bea
r+I\cap J=s+I\cap J~~\Ra~~r-s\in I\cap J\subset I+J,~~\Ra~~r+(I+J)=s+(I+J),\nn
\eea
and it is a ring homomorphism with kernel
\[
\textstyle \ker f=\{r+I\cap J:r+(I+J)=I+J\}=\{r+I\cap J:r\in I+J\}={I+J\over I\cap J}.\nn \qedhere
\]
\end{proof}

\begin{thm}[\textcolor{blue}{\index{Correspondence theorem for! rings}{Correspondence theorem for rings}}]
Let $I\lhd R$ be an ideal (in a ring $R$), $Sub_I(R):=\{\txt{subrings of $R$ containing $I$}\}$, and $Sub\left({R\over I}\right):=\{\txt{subrings of ${R\over I}$}\}$. Then there exists a bijection
\bea
\bt Sub_I(R)\ar[rr,leftrightarrow] && Sub\left({R\over I}\right).\et\nn
\eea
\end{thm}
\begin{proof}
Let $\pi_I:R\ra {R\over I},~r\mapsto r+I$ be the natural ring homomorphism. Consider the maps
\bea
\textstyle\Phi:Sub_I(R)\ra Sub\left({R\over I}\right),~S\mapsto\pi_I(S)={S\over I},~~~~\Psi:Sub\left({R\over I}\right)\ra Sub_I(R),~S'\mapsto\pi_I^{-1}(S').\nn
\eea
Then for all ~$S\in Sub_I(R)$ and $S'\in Sub\left({R\over I}\right)$,
\begin{align}
&\Phi\circ\Psi(S')=\pi_I\left(\pi_I^{-1}(S')\right)=S'\cap\pi_I(R)\sr{(a)}{=}S',~~\txt{where step (a) holds because $\pi_I$ is surjective},\nn\\
&\textstyle\Psi\circ\Phi(S)=\pi_I^{-1}\big(\pi_I(S)\big)=\pi_I^{-1}\left({S\over I}\right)\sr{(b)}{=}S+I\sr{(c)}{=}S,~~\txt{where step (c) holds because $I\subset S$},\nn
\end{align}
and step (b) is obtained as before.
\end{proof}

\begin{crl}
An ideal $I\lhd R$ is maximal $\iff$ $R/I$ is a simple ring. In particular, if $R$ is commutative, then an ideal $I\lhd R$ is maximal $\iff$ $R/I$ is a field.
\end{crl}

\section{Isomorphism and correspondence theorems for modules}
Later, generalizations of the isomorphism theorems in this section will be considered for abelian categories.

Here is a summary of the module isomorphism theorems we want to prove:
\begin{enumerate}
\item Given an $R$-homomorphism $f:M\ra M'$, a submodule $N\subset M$, and the associated submodule $L:=N\cap\ker f\subset M$ with natural map $\pi_L:M\ra{M\over L},~m\mapsto m+L$, we get (i) an $R$-module homomorphism $f_L:{M\over L}\ra M'$ with kernel $\ker f_L={\ker f\over L}$ such that $f=f_L\circ\pi_L:M\sr{\pi_L}{\ral}{M\over L}\sr{f_L}{\ral}M'$, and (ii) an $R$-module isomorphism ${M\over\ker f}\cong\im f$.
\item Given submodules ${}_RN,{}_RN'\subset{}_RM$, we get an $R$-module isomorphism ~${N\over N\cap N'}\cong {N+N'\over N'}$.
\item Given submodules ${}_RN,{}_RN'\subset{}_RM$, we get an $R$-module isomorphism ~${M/N\cap N'\over (N+N')/N\cap N'}\cong {M\over N+N'}$.
\end{enumerate}

\begin{thm}[\textcolor{blue}{\index{First isomorphism theorem for! modules}{First isomorphism theorem for modules}}]
Let $f:M\ra M'$ be an $R$-homomorphism, ${}_RN\subset{}_RM$ a submodule, $L:=N\cap\ker f$, and $\pi_L:M\ra{M\over L},~m\mapsto m+L$ the natural $R$-homomorphism. Then there exists an $R$-homomorphism $f_L:{M\over L}\ra M'$ such that (i) $f=f_L\circ\pi_L$ and (ii) $\ker f_L={\ker f\over L}={\ker f\over N\cap\ker f}$.
{\small\[
\textstyle\bt
M\ar[d,"f"]\ar[rr,"\pi_L"] && {M\over L}\ar[dll,dashed,"f_L"]\\
M'
\et~~~~~~~~\txt{In particular, with $N:=\ker f$, we get an $R$-module isomorphism ~~${M\over \ker f}\cong f(M)$}.
\]}
\end{thm}
\begin{proof}
The map $f_L:{M\over L}\ra M',~m+L\mapsto f(m)$ is well defined, since
\bea
\textstyle m+L=m'+L~~\Ra~~m-m'\in L\subset\ker f~~\Ra~~0=f(m-m')=f(m)-f(m'),~~\Ra~~f(m)=f(m').\nn
\eea
Moreover, $f_L\circ\pi_L(m)=f_L(m+L)=f(m)$, for all $m\in M$, and
\[
\textstyle \ker f_L=\big\{m+L:0=f_L(m+L)=f(m)\big\}=\big\{m+L:m\in\ker f\big\}={\ker f+L\over L}={\ker f\over L}.\nn \qedhere
\]
\end{proof}

\begin{thm}[\textcolor{blue}{\index{Second isomorphism theorem for! modules}{Second isomorphism theorem for modules}}]
Let $M$ be an $R$-module, and ${}_RN,{}_RN'\subset{}_RM$ submodules. Then (i) ${}_R(N\cap N')\subset{}_RN$, (ii) ${}_RN'\subset{}_R(N+N')$, and (iii) ${N\over N\cap N'}\cong{N+N'\over N'}$.
\end{thm}
\begin{proof}
The surjective $R$-homomorphism $f:N\ra {N+N'\over N'}$, $n\mapsto n+N'$ has kernel
\[
\textstyle \ker f=\{n\in N:n+N'=N'\}=\{n\in N:n\in N'\}=N\cap N'.\nn \qedhere
\]
\end{proof}

\begin{thm}[\textcolor{blue}{\index{Third isomorphism theorem for! modules}{Third isomorphism theorem for modules}}]
Let $M$ be an $R$-module and ${}_RN,{}_RN'\subset{}_RM$ submodules. Then (i) ${}_R(N\cap N'),{}_R(N+N')\subset{}_RM$, (ii) ${}_R(N\cap N')\subset {}_R(N+N')$, and (iii) ${M/N\cap N'\over(N+N')/N\cap N'}\cong{M\over N+N'}$. In particular, if $N\subset N'$, then ~${M/N\over N'/N}\cong{M\over N'}$.
\end{thm}
\begin{proof}
The map $f:{M\over N\cap N'}\ra{M\over N+N'}$, $m+N\cap N'\mapsto m+(N+N')$ is well defined because
{\small\bea
m_1+N\cap N'=m_2+N\cap N'~~\Ra~~m_1-m_2\in N\cap N'\subset N+N',~~\Ra~~m_1+(N+N')=m_2+(N+N'),\nn
\eea}and it is a surjective $R$-homomorphism with kernel
\[
\textstyle \ker f=\{m+N\cap N':m+(N+N')=N+N'\}=\{m+N\cap N':m\in N+N'\}={N+N'\over N\cap N'}.\nn \qedhere
\]
\end{proof}

\begin{thm}[\textcolor{blue}{\index{Correspondence theorem for! modules}{Correspondence theorem for modules}}]
Let $M$ be an $R$-module, $_RV\subset{}_RM$ a submodule, $Sub_V(M):=\{\txt{submodules of $M$ containing $V$}\}$, and $Sub\left({M\over V}\right):=\{\txt{submodules of ${M\over V}$}\}$. Then there exists a bijection ~
\bt Sub_V(M)\ar[rr,leftrightarrow] && Sub\left({M\over V}\right).\et
\end{thm}
\begin{proof}
Let $\pi_V:M\ra {M\over V},~m\ra m+V$ be the natural $R$-homomorphism. Consider the maps
\bea
\textstyle\Phi:Sub_V(M)\ra Sub\left({M\over V}\right),~N\mapsto\pi_V(N)={N\over V},~~~~\Psi:Sub\left({M\over V}\right)\ra Sub_V(M),~N'\mapsto\pi_V^{-1}(N').\nn
\eea
Then for all ~$S\in Sub_V(M)$ and $S'\in Sub\left({M\over V}\right)$,
\begin{align}
&\Phi\circ\Psi(N')=\pi_V\left(\pi_V^{-1}(N')\right)=N'\cap\pi_V(M)\sr{(a)}{=}N',~~\txt{where step (a) holds because $\pi_V$ is surjective},\nn\\
&\textstyle\Psi\circ\Phi(N)=\pi_V^{-1}\big(\pi_V(N)\big)=\pi_V^{-1}\left({N\over V}\right)\sr{(b)}{=}N+V\sr{(c)}{=}N,~~\txt{where step (c) holds because $V\subset N$},\nn
\end{align}
where step (b) is obtained as before.
\end{proof}

\begin{crl}
A submodule $N\subset M$ is maximal if and only if $M/N$ is a simple module.
\end{crl}

\section{Isomorphism and correspondence theorems for algebras}
Here is a summary of the $R$-algebra isomorphism theorems we want to prove:
\begin{enumerate}
\item Given an $R$-algebra homomorphism $f:A\ra A'$, an ideal $I\lhd A$, and the associated ideal $J:=I\cap\ker f\lhd A$ with natural map $\pi_J:A\ra{A\over J},~a\mapsto a+J$, we get (i) an $R$-algebra homomorphism $f_J:{A\over J}\ra A'$ with kernel $\ker f_J={\ker f\over J}$ such that $f=f_J\circ\pi_J:A\sr{\pi_J}{\ral}{A\over J}\sr{f_J}{\ral}A'$, and (ii) an $R$-algebra isomorphism ${A\over\ker f}\cong\im f$.
\item Given an $R$-subalgebra ${}_RS\leq{}_RA$ and an ideal $I\lhd A$, we get an $R$-algebra isomorphism ~${S\over S\cap I}\cong {S+I\over I}$.
\item Given ideals $I,J\lhd A$ of an $R$-algebra $A$, we get an $R$-algebra isomorphism ~${A/I\cap J\over (I+J)/I\cap J}\cong {A\over I+J}$.
\end{enumerate}

\begin{thm}[\textcolor{blue}{\index{First isomorphism theorem for! algebras}{First isomorphism theorem for algebras}}]
Let $f:A\ra A'$ be an $R$-algebra homomorphism, $I\lhd A$ an ideal, $J:=I\cap\ker f$, and $\pi_J:A\ra{A\over J},~a\mapsto a+J$ the natural $R$-algebra homomorphism. Then there exists an $R$-algebra homomorphism $f_J:{A\over J}\ra A'$ such that (i) $f=f_J\circ\pi_J$ and (ii) $\ker f_J={\ker f\over J}={\ker f\over I\cap\ker f}$.
{\small\[
\textstyle\bt
A\ar[d,"f"]\ar[rr,"\pi_J"] && {A\over J}\ar[dll,dashed,"f_J"]\\
A'
\et~~~~~~~~\txt{In particular, with $I:=\ker f$, we get an $R$-algebra isomorphism ~~${A\over \ker f}\cong f(A)$}.\nn
\]}
\end{thm}
\begin{proof}
The map $f_J:{A\over J}\ra A',~a+J\mapsto f(a)$ is well defined, since
\bea
\textstyle a+J=b+J~~\Ra~~a-b\in J\subset\ker f~~\Ra~~0=f(a-b)=f(a)-f(b),~~\Ra~~f(a)=f(b).\nn
\eea
Moreover, $f_J\circ\pi_J(a)=f_J(a+J)=f(a)$, for all $a\in A$, and
\[
\textstyle \ker f_J=\big\{a+J:0=f_J(a+J)=f(a)\big\}=\big\{a+J:a\in\ker f\big\}={\ker f+J\over J}={\ker f\over J}.\nn \qedhere
\]
\end{proof}

\begin{thm}[\textcolor{blue}{\index{Second isomorphism theorem for! algebras}{Second isomorphism theorem for algebras}}]
Let $A$ be an $R$-algebra, ${}_RS\leq{}_RA$ a subalgebra and $I\lhd A$ an ideal. Then (i) $S\cap I\lhd S$, (ii) $I\lhd S+I$, and (iii) ${S\over S\cap I}\cong{S+I\over I}$.
\end{thm}
\begin{proof}
The surjective $R$-algebra homomorphism $f:S\ra {S+I\over I}$, $s\mapsto s+I$ has kernel
\[
\textstyle \ker f=\{s\in S:s+I=I\}=\{s\in S:s\in I\}=S\cap I.\nn \qedhere
\]
\end{proof}

\begin{thm}[\textcolor{blue}{\index{Third isomorphism theorem for! algebras}{Third isomorphism theorem for algebras}}]
Let $A$ be an $R$-algebra and $I,J\lhd A$ ideals. Then (i) $I\cap J,I+J\lhd A$, (ii) $I\cap J\lhd I+J$, and (iii) ${A/I\cap J\over(I+J)/I\cap J}\cong{A\over I+J}$.

In particular, if $I\subset J$, then ~${A/I\over J/I}\cong{A\over J}$.
\end{thm}
\begin{proof}
The map $f:{A\over I\cap J}\ra{A\over I+J}$, $a+I\cap J\mapsto a+(I+J)$ is well defined because
\bea
a+I\cap J=b+I\cap J~~\Ra~~a-b\in I\cap J\subset I+J,~~\Ra~~a+(I+J)=b+(I+J),\nn
\eea
and it is a surjective $R$-algebra homomorphism with kernel
\[
\textstyle \ker f=\{a+I\cap J:a+(I+J)=I+J\}=\{a+I\cap J:a\in I+J\}={I+J\over I\cap J}.\nn \qedhere
\]
\end{proof}

\begin{thm}[\textcolor{blue}{\index{Correspondence theorem for! algebras}{Correspondence theorem for algebras}}]
Let $A$ be an $R$-algebra, $I\lhd A$ an ideal, $Sub_I(A):=\{\txt{subalgebras of $A$ containing $I$}\}$, and $Sub\left({A\over I}\right):=\{\txt{subalgebras of ${A\over I}$}\}$. Then there exists a bijection
\bea
\bt Sub_I(A)\ar[rr,leftrightarrow] && Sub\left({A\over I}\right).\et\nn
\eea
\end{thm}
\begin{proof}
Let $\pi_I:A\ra {A\over I},~a\mapsto a+I$ be the natural $R$-algebra homomorphism. Consider the maps
\bea
\textstyle\Phi:Sub_I(A)\ra Sub\left({A\over I}\right),~S\mapsto\pi_I(S)={S\over I},~~~~\Psi:Sub\left({A\over I}\right)\ra Sub_I(A),~S'\mapsto\pi_I^{-1}(S').\nn
\eea
Then for all ~$S\in Sub_I(A)$ and $S'\in Sub\left({A\over I}\right)$,
\begin{align}
&\Phi\circ\Psi(S')=\pi_I\left(\pi_I^{-1}(S')\right)=S'\cap\pi_I(A)\sr{(a)}{=}S',~~\txt{where step (a) holds because $\pi_I$ is surjective},\nn\\
&\textstyle\Psi\circ\Phi(S)=\pi_I^{-1}\big(\pi_I(S)\big)=\pi_I^{-1}\left({S\over I}\right)\sr{(b)}{=}S+I\sr{(c)}{=}S,~~\txt{where step (c) holds because $I\subset S$},\nn
\end{align}
and step (b) is obtained as before.
\end{proof}

\begin{crl}
An ideal $I\lhd A$ of an $R$-algebra $A$ is maximal $\iff$ ${A\over I}$ is a simple ${R+I\over I}$-algebra. In particular, if $A$ is commutative, then an ideal $I\lhd A$ is maximal $\iff$ ${A\over I}$ is a field (and ${R+I\over I}$-vector space).
\end{crl}

\section{Smallest Subfield and Characteristic of a Field}
\begin{dfn}[\textcolor{blue}{\index{Prime! subfield}{Prime subfield}}]
The \ul{prime subfield} of a field $k$ is the smallest subfield $k_0\subset k$, i.e.,
\bea
\textstyle k_0:=\bigcap\{k'\subset k:~k'~\txt{a subfield}\}.\nn
\eea
\end{dfn}

\begin{prp}
Let $k$ is be a field. Then, either (i) {\small $k_0\cong\Rational$} or (ii) {\small $k_0\cong\Integer_p:={\Integer\over\Integer p}$}, for a prime $p\in\Integer$.
\end{prp}
\begin{proof}
Consider the ring homomorphism $\phi:\Integer\ra k,~n\mapsto n\cdot1_k:=\txt{sign}(n)(\overbrace{1_k+\cdots+1_k}^{n~\txt{times}})$. Then $\phi(\Integer)=\Integer 1_k$ is an ID (integral domain) since for $n,m\in\Integer$, we have
{\small\bea
0_k=(n\cdot 1_k)(m\cdot 1_k)=(mn)\cdot 1_k~~\Ra~~mn=0,~~\Ra~~m=0~~\txt{or}~~n=0,~~\Ra~~n\cdot 1_k=0_k~~\txt{or}~~m\cdot 1_k=0_k.\nn
\eea}Either $\ker\phi=0$ (i.e., $\phi$ is injective) or $\ker\phi\neq 0$ (i.e., $\phi$ is not injective).
{\flushleft (i)} Suppose $\ker\phi=0$. In this case, the smallest subfield of $k$ (which is by construction the smallest subfield of $k$ containing $1_k$) is the same as the smallest subfield of $k$ containing $\phi(\Integer)\cong{\Integer\over\ker\phi}=\Integer$,
\bea
~~\Ra~~k_0=Q(\phi(\Integer))=Q(\Integer\cdot 1_k)\cong Q(\Integer)=\Rational.\nn
\eea
{\flushleft (ii)} Suppose $\ker\phi\neq 0$. Since ${\Integer\over\ker\phi}\cong\phi(\Integer)$ and $\phi(\Integer)=\Integer\cdot 1_k$ is an ID, it follows that $\ker\phi\vartriangleleft \Integer$ is prime, i.e., $\ker\phi=\Integer p$ for a prime $p\in\Integer$. In this case, the smallest subfield of $k$ (which is by construction the smallest subfield of $k$ containing $1_k$) is the same as the smallest subfield of $k$ containing $\phi(\Integer)\cong{\Integer\over\Integer p}\cong\Integer_p$,
\[
~~\Ra~~k_0\cong \Integer_p. \qedhere
\]
\end{proof}

\begin{dfn}[\textcolor{blue}{\index{Characteristic! of a field}{Characteristic of a field}}]
If $k$ is a field, the \ul{characteristic} of $k$ is
$$char(k):=\left\{
                 \begin{array}{ll}
                   0, & \hbox{if the prime subfield}~~k_0\cong\Rational, \\
                   |p|, & \hbox{if the prime subfield}~~k_0\cong\Integer_p,~~\txt{for a prime}~~p\in\Integer.
                 \end{array}
               \right.
$$
\end{dfn}
Observe that if $k\subset K$ are fields, then $char(k)=char(K)$. Also, it is clear that the fields $\Rational$, $\Real$, $\Complex$ each have characteristic $0$.

%% file: parts/AlgebraNC/AlgebraNC6.tex
\chapter{Products, Coproducts, and Grading}
Like with a few other chapters, the material in this chapter, although small compared to that in a typical chapter in the notes, seems important enough to stand out on its own. This is essentially a summary of further basic concepts (for groups, rings, modules, and algebras) whose categorical generalizations will be encountered later. Such basic concepts are therefore worth introducing at this stage, mainly in order to improve on the vocabulary required by the subsequent discussion on categories.
\section{Products and Coproducts}
The ``product objects'' discussed in this section have alternative characterizing properties (called ``\emph{universal properties}'') that will be described in part \ref{ChCatPers}. Throughout this section we will let $I$ be any set.

\subsection{Products and coproducts of multiplicative sets}\label{GpPdSec}~\\~
The ``coproducts'' of multiplicative sets are also called their ``free products''.
\begin{dfn}[\blue{Products and coproducts of associative sets}]
Let $S,S'$ be associative sets (or just multiplicative sets where sensible). A map $f:S\ra S'$ is a \index{Homomorphism! of associative sets}{\ul{homomorphism of associative sets}} if $f(ab)=f(a)f(b)$ for all $a,b\in S$. An \index{Isomorphism of! associative sets}{\ul{isomorphism of associative sets}} is a bijective homomorphism of associative sets. If an isomorphism of associative sets $f:S\ra S'$ exists, we say $S$ and $S'$ are \index{Isomorphic! associative sets}{\ul{isomorphic associative sets}}, written $S\cong S'$.

An equivalence relation $\sim$ on $S$ is a \index{Congruence relation}{\ul{congruence relation}} if for all $a,b,c,d\in S$, $a\sim b$ and $c\sim d$ imply $ac\sim bd$. The \ul{quotient associative set} of $S$ associated with $\sim$ is the set of equivalence classes ${S\over\sim}:=\{[a]:a\in S\}$ as an associative set with multiplication ${S\over\sim}\times {S\over\sim}\ra {S\over\sim},~([a],[b])\mapsto[a][b]:=[ab]$. The \index{Congruence homomorphism of associative sets}{\ul{congruence homomorphism of associative sets}} associated with $\sim$ is the surjective map
\bea
\textstyle\pi_\sim:S\ra {S\over\sim},~s\mapsto [s].\nn
\eea

Let $A,B\subset S$. We write $AB:=\{ab:a\in A,b\in B\}$. The set $A$ is an \index{Associative! subset}{\ul{associative subset}} of $S$, written $A\leq S$, if $aa'\in A$ for all $a,a'\in A$. An associative subset $N\leq S$ is \index{Normal associative subset}{\ul{normal}}, written $N\vartriangleleft S$, if (i) $sN=Ns$ for all $a\in S$ and (ii) $NN=N$. (\blue{footnote}\footnote{If $N\lhd S$ is a normal associative subset, then we get a congruence relation on $S$ given (for $a,b\in S$) by $a\sim b$ if $aN=bN$.}). The \index{Quotient! associative set}{\ul{quotient associative set}} of $S$ by $N$ is the set
\bea
S/N:=\{sN:s\in S\}\subset\P(S)\nn
\eea
as an associative set with multiplication $\cdot:S/N\times S/N\ra S/N$ given by $(sN)\cdot(s'N):=ss'N$, for all $s,s'\in S$. The \index{Natural! homomorphism of associative sets}{\ul{natural homomorphism of associative sets}} associated with $N\vartriangleleft S$ is the surjective map
\bea
\pi_N:S\ra S/N,~s\mapsto sN.\nn
\eea

Let $\S:=(S_i)_{i\in {I}}$ be an indexing of associative sets $S_i$. The \index{Product of! associative sets}{\ul{product}} $\prod\S$ of the associative sets $\S$ is the smallest associative set (up to isomorphism of associative sets) having each member $S_i\in\S$ as a quotient associative set (or alternatively, the smallest associative set for which there exists a surjective homomorphism of associative sets $p_i:\prod\S\twoheadrightarrow S_i$ for each $i\in I$). Explicitly, the product of $\S$ is the set
\bea
\textstyle\prod\S:=\prod_{i\in I}S_i:=\big\{\txt{maps}~~s:{I}\ra\bigcup_i S_i,~i\mapsto s_i\in S_i\big\}=\big\{(s_i)_{i\in{I}}~|~s_i\in S_i\big\}\nn
\eea
as an associative set with \emph{multiplication} given \emph{pointwise} by
\bea
(s_i)_{i\in{I}}\cdot(s'_i)_{i\in{I}}:=(s_is_i')_{i\in{I}}.\nn
\eea

The \index{Coproduct! (free product) of associative sets}{\ul{coproduct} (or \ul{free product})} $\coprod\S=\coprod_{i\in I}S_i$  of the associative sets $\S=(S_i)_{i\in I}$ is the associative set generated by $\S$, in the sense that it is the smallest associative set (up to isomorphism of associative sets) containing every member of $\S$ as an associative subset.
\end{dfn}
\begin{example}
In particular, if $\S=(S,S')$, then with (i)~$T(S,S'):=\bigcup_{n\geq 1}(S\sqcup S')^n$ and (ii) for all $a,b\in S\sqcup S'$, ``$a\sim b$ if either both $a,b\in S$ or both $a,b\in S'$'',
\[
\textstyle S\ast S':=S\amalg S':=\coprod(S,S'):=\left\{(a_0,a_1,...,a_n)\in T(S,S'):~a_i\not\sim a_{i+1}~\txt{for all}~i,~n\in\Natural\right\}\nn
\]as an associative set with multiplication given by the following: For any $(a_1,...,a_m)\in(S\ast S')\cap(S\sqcup S')^m$ and $(b_1,...,b_n)\in (S\ast S')\cap(S\sqcup S')^n$,
{\small\[
(a_1,...,a_m)\cdot (b_1,...,b_n):=
\left\{
  \begin{array}{ll}
    (a_1,...,a_{m-1},a_mb_1,b_2,...,b_n)\in (S\ast S')\cap(S\sqcup S')^{m+n-1}, & \txt{if}~~a_m\sim b_1\\
    (a_1,...,a_m,b_1,...,b_n)\in (S\ast S')\cap(S\sqcup S')^{m+n}, & \txt{if}~~a_m\not\sim b_1
  \end{array}
\right\}.
\]}
The component associative sets $S,S'$ are included via the associative set homomorphisms
\[
S\hookrightarrow S\ast S',~s\mapsto (s)~~~~\txt{and}~~~~S'\hookrightarrow S\ast S',~s'\mapsto (s').
\]
\end{example}

\begin{dfn}[\blue{Products and coproducts of associative identity sets (monoids)}]
Let $\M:=(M_i)_{i\in {I}}$ be an indexing of associative identity sets (monoids) $M_i$. The \index{Product of! monoids}{\ul{product}} of the monoids $\M$ is the smallest monoid (up to isomorphism of associative sets) having each member $M_i\in\M$ as a quotient associative set. Explicitly, the product of $\M$ is the set
\bea
\textstyle\prod\M:=\prod_{i\in I}M_i:=\big\{\txt{maps}~~s:{I}\ra\bigcup_i M_i,~i\mapsto m_i\in M_i\big\}=\big\{(m_i)_{i\in{I}}~|~m_i\in M_i\big\}\nn
\eea
as a monoid with \emph{multiplication} and \emph{identity element} respectively given \emph{pointwise} by
\bea
(m_i)_{i\in{I}}\cdot(s'_i)_{i\in{I}}:=(m_im_i')_{i\in{I}},~~~~e_{\prod\M}:=(e_{M_i})_{i\in{I}}.\nn
\eea

The \index{Coproduct! (free product) of monoids}{\ul{coproduct} (\ul{free product})} $\coprod\M=\coprod_{i\in I}M_i$ of the monoids $\M=(M_i)_{i\in I}$ is their coproduct as associative sets with the identity element structure maintained: That is, it is the monoid generated by $\M$, in the sense that it is the smallest monoid (up to isomorphism of associative sets) containing every member of $\M$ as an associative subset.
\end{dfn}

\begin{example}
In particular, if $\M=(M,M')$, then with (i) $S:=M\backslash\{e_M\}$, $S':=M'\backslash\{e_{M'}\}$, $T(M,M'):=\bigcup_{n\geq 1}(S\sqcup S')^n$, and (ii) for all $a,b\in S\sqcup S'$, ``$a\sim b$ if either both $a,b\in S$ or both $a,b\in S'$'',
\[
\textstyle M\ast M':=M\amalg M':=\coprod(M,M'):=\left\{(a_0,a_1,...,a_n)\in T(M,M'):~a_i\not\sim a_{i+1}~\txt{for all}~i,~n\in\Natural\right\}
\]as an associative set with multiplication given by the following: For any $(a_1,...,a_m)\in(M\ast M')\cap(S\sqcup S')^m$ and $(b_1,...,b_n)\in (M\ast M')\cap(S\sqcup S')^n$,
{\scriptsize\[
(a_1,...,a_m)\cdot (b_1,...,b_n):=
\left\{
  \begin{array}{ll}
    (a_1,...,a_m,b_1,...,b_n)\in (M\ast M')\cap(S\sqcup S')^{m+n}, & \txt{if}~~a_m\not\sim b_1\\
    (a_1,...,a_{m-1},a_mb_1,b_2,...,b_n)\in (M\ast M')\cap(S\sqcup S')^{m+n-1}, & \txt{if}~~a_m\sim b_1~\txt{and}~a_mb_1\not\in\{e_M,e_{M'}\}\\
    (a_1,...,a_{m-1})\cdot(b_2,...,b_n)\in (M\ast M')\cap(S\sqcup S')^{m+n-2}, & \txt{if}~~a_m\sim b_1~\txt{and}~a_mb_1\in\{e_M,e_{M'}\}
  \end{array}
\right\}.
\]}

By construction, $M\ast M'$ is a monoid with identity element $e_{M\ast M'}:=()$, i.e., the empty tuple. The component monoids $M,M'$ are included via the associative set homomorphisms
\[
M\hookrightarrow M\ast M',~m\mapsto (m),~e_M\mapsto e_{M\ast M'}~~~~\txt{and}~~~~M'\hookrightarrow M\ast M',~m'\mapsto (m'),~e_{M'}\mapsto e_{M\ast M'}.
\]
\end{example}

\begin{dfn}[\blue{Products and coproducts of groups}]
Let $\G:=(G_i)_{i\in {I}}$ be an indexing of groups $G_i$. The \index{Product of! groups}{\ul{product}} of the groups $\G$ is the smallest group (up to isomorphism of groups) having each member $G_i\in\G$ as a quotient group. Explicitly, the product of $\G$ is the set
\bea
\textstyle\prod\G:=\prod_{i\in I}G_i:=\big\{\txt{maps}~~g:{I}\ra\bigcup_i G_i,~i\mapsto g_i\in G_i\big\}=\big\{(g_i)_{i\in{I}}~|~g_i\in G_i\big\}\nn
\eea
as a group with \emph{multiplication}, \emph{inverse}, and \emph{identity element} respectively given \emph{pointwise} by
\bea
(g_i)_{i\in{I}}\cdot(g'_i)_{i\in{I}}:=(g_ig_i')_{i\in{I}},~~~~(g_i)_{i\in{I}}{}^{-1}:=(g_i^{-1})_{i\in{I}},~~~~e_{\prod\G}:=(e_{G_i})_{i\in{I}}.\nn
\eea

The \index{Coproduct! (free product) of groups}{\ul{coproduct} (\ul{free product})} $\coprod\G=\coprod_{i\in I}G_i$ of the groups $\G=(G_i)_{i\in I}$ is their coproduct as associative sets with the group structure maintained: That is, it is the group generated by $\G$, in the sense that it is the smallest group (up to isomorphism of groups) containing every member of $\G$ as a subgroup.
\end{dfn}

\begin{example}
In particular, if $\G=(G,G')$, then with (i) $S:=G\backslash\{e_G\}$, $S':=G'\backslash\{e_{G'}\}$, $T(G,G'):=\bigcup_{n\geq 1}(S\sqcup S')^n$, and (ii) for all $a,b\in S\sqcup S'$, ``$a\sim b$ if either both $a,b\in S$ or both $a,b\in S'$'',
\[
\textstyle G\ast G':=G\amalg G':=\coprod(G,G'):=\left\{(a_0,a_1,...,a_n)\in T(G,G'):~a_i\not\sim a_{i+1}~\txt{for all}~i,~n\in\Natural\right\}
\]as an associative set with multiplication given by the following: For any $(a_1,...,a_m)\in(G\ast G')\cap(S\sqcup S')^m$ and $(b_1,...,b_n)\in (G\ast G')\cap(S\sqcup S')^n$,
{\scriptsize\[
(a_1,...,a_m)\cdot (b_1,...,b_n):=
\left\{
  \begin{array}{ll}
    (a_1,...,a_m,b_1,...,b_n)\in (G\ast G')\cap(S\sqcup S')^{m+n}, & \txt{if}~~a_m\not\sim b_1\\
    (a_1,...,a_{m-1},a_mb_1,b_2,...,b_n)\in (G\ast G')\cap(S\sqcup S')^{m+n-1}, & \txt{if}~~a_m\sim b_1~\txt{and}~a_mb_1\not\in\{e_G,e_{G'}\}\\
    (a_1,...,a_{m-1})\cdot(b_2,...,b_n)\in (G\ast G')\cap(S\sqcup S')^{m+n-2}, & \txt{if}~~a_m\sim b_1~\txt{and}~a_mb_1\in\{e_G,e_{G'}\}
  \end{array}
\right\}.
\]}

By construction, $G\ast G'$ is a group with identity element $e_{G\ast G'}:=()$, i.e., the empty tuple. (\blue{footnote}\footnote{The element-wise multiplication used here applies to the coproduct of identity sets in general, and not just of groups.}). The component groups $G,G'$ are included via the group homomorphisms
\[
G\hookrightarrow G\ast G',~g\mapsto (g),~e_G\mapsto e_{G\ast G'}~~~~\txt{and}~~~~G'\hookrightarrow G\ast G',~g'\mapsto (g'),~e_{G'}\mapsto e_{G\ast G'}.
\]
\end{example}

\subsection{Products and coproducts of rings, modules, and algebras}
\label{RgPdSec}\label{MdPdSec}\label{AgPdSec}
\begin{dfn}[\blue{Products and coproducts of rings}]
Let $\R:=(R_i)_{i\in {I}}$ be an indexing of rings $R_i$. The \index{Product of! rings}{\ul{product}} of the rings $\R$ is the smallest ring (up to isomorphism of rings) having each member $R_i\in\R$ as a quotient ring. Explicitly, the product of $\R$ is the set
\bea
\textstyle\prod\R:=\prod_{i\in I}R_i:=\big\{\txt{maps}~~r:{I}\ra\bigcup_i R_i,~i\mapsto r_i\in R_i\big\}=\big\{(r_i)_{i\in{I}}~|~r_i\in R_i\big\},\nn
\eea
as a ring with \emph{addition}, \emph{multiplication}, \emph{zero}, and \emph{unity} respectively given \emph{pointwise} by
\bea
(r_i)_{i\in{I}}+(r'_i)_{i\in{I}}:=(r_i+r_i')_{i\in{I}},~~(r_i)_{i\in{I}}\cdot(r'_i)_{i\in{I}}:=(r_ir_i')_{i\in{I}},~~0_{\prod\R}:=(0_{R_i})_{i\in{I}},~~1_{\prod\R}:=(1_{R_i})_{i\in{I}}.\nn
\eea

The \index{Coproduct! of rings}{\ul{coproduct}} $\coprod\R=\coprod_{i\in I}R_i$ of the rings $\R=(R_i)_{i\in I}$ is the ring generated by $\R$, in the sense that it is the smallest ring (up to isomorphism of rings) containing every member of $\R$ as a subring. (\blue{footnote}\footnote{This is related to the coproduct of the rings as monoids (assoc. id. sets) with respect to the multiplications in the rings.})
\end{dfn}

\begin{example}
In particular, if $\R=(R,R')$, consider the product of abelian groups ${M_{\Integer,R,R'}}:=$ ``$R\times R'=R\times\{0\}+\{0\}\times R'$ as a $\Integer$-module'' and its tensor algebra $Te({M_{\Integer,R,R'}}):=\Integer\langle{M_{\Integer,R,R'}}\rangle:=\sum_{i\in\Natural}\Integer{M_{\Integer,R,R'}}^i$. With $S:=(R\times\{0\})\backslash\{(1_R,0)\}$ and $S':=(\{0\}\times R')\backslash\{(0,1_{R'})\}$, and noting $(0,0)\in S\cap S'$, we have
{\small\begin{align}
&M_{\Integer,R,R'}=(S\sqcup\{(1_R,0)\})+(\{(0,1_{R'})\}\sqcup S')=(S+S')\sqcup\big(S+\{(0,1_{R'})\}\big)\sqcup\big(\{(1_R,0)\}+S'\big)\sqcup\big\{(1_R,1_{R'})\big\}\nn\\
&~~~~=\Integer(S+S')+\Integer\{(1_R,0),(0,1_{R'})\}=\Integer S+\Integer S'+\Integer(1_R,0)+\Integer(0,1_{R'}).\nn
\end{align}}Let (i)~$T(R,R'):=\bigcup_{n\geq 1}(S\sqcup S')^n$ and (ii) for all $a,b\in S\sqcup S'$, let ``$a\sim b$ if either both $a,b\in S$ or both $a,b\in S'$''. Given $(a_1,...,a_m),(b_1,...,b_n)\in T(R,R')\subset Te({M_{\Integer,R,R'}})$, define
{\scriptsize\[
D(a_1,...,a_m,b_1,...,b_n):=\left\{
  \begin{array}{ll}
    (a_1,...,a_m)\cdot (b_1,...,b_n)-(a_1,...,a_m,b_1,...,b_n), & \txt{if}~~a_m\not\sim b_1\\
    (a_1,...,a_m)\cdot (b_1,...,b_n)-(a_1,...,a_{m-1},a_mb_1,b_2,...,b_n), & \txt{if}~~a_m\sim b_1~\txt{and}~a_mb_1\not\in\{(1_R,0),(0,1_{R'})\}\\
    (a_1,...,a_m)\cdot (b_1,...,b_n)-(a_1,...,a_{m-1})\cdot(b_2,...,b_n), & \txt{if}~~a_m\sim b_1~\txt{and}~a_mb_1\in\{(1_R,0),(0,1_{R'})\}
  \end{array}
\right\}.
\]}Consider the subgroups $\wt{R}:=R\times\{0\}\subset R\times R'$ and $\wt{R}':=\{0\}\times R'\subset R\times R'$, and the ideal
\[
Q:=\left\langle \left\{D(a_1,...,a_m,b_1,...,b_n)~\big|~(a_1,...,a_m),(b_1,...,b_n)\in T(R,R')\right\}\right\rangle\lhd Te({M_{\Integer,R,R'}}).\nn
\]
With the subsring $\Integer\langle \wt{R}\ast\wt{R}'\rangle\leq Te({M_{\Integer,R,R'}})$ generated by the coproduct of multiplicative-monoids $\wt{R}\ast\wt{R}'$,
\begin{align}
\textstyle R\amalg R'&\textstyle:=\coprod(R,R'):={Te({M_{\Integer,R,R'}})\over Q}={\Integer\langle \wt{R}\ast\wt{R}'\rangle+Q\over Q}={\sum_{i\in\Natural}\Integer\langle \wt{R}\ast\wt{R}'\rangle\cap{M_{\Integer,R,R'}}^i+Q\over Q}=\sum_{i\in\Natural}{\Integer\langle \wt{R}\ast\wt{R}'\rangle\cap{M_{\Integer,R,R'}}^i+Q\over Q}\nn\\
&=\textstyle\big\{\sum_in_ic_i+Q:n_i\in\Integer,~c_i\in \Integer\langle \wt{R}\ast\wt{R}'\rangle\cap{M_{\Integer,R,R'}}^i,~n_i=0~\txt{a.e.f.}~i\in I\big\}\nn\\
&\textstyle={\Integer\langle S+S'\rangle+Q\over Q},\nn
\end{align} where $\Integer\langle S+S'\rangle\leq Te({M_{\Integer,R,R'}})$ is the subring generated by $S+S'\subset {M_{\Integer,R,R'}}$.

By construction, $R\amalg R'$ is a ring with unity $1_{R\amalg R'}:=()+Q$, i.e., the class of the empty tuple. The component rings $R,R'$ are included via the ring homomorphisms
{\small\[
\textstyle R\hookrightarrow R\amalg R',~r\mapsto \big((r,0)\big)+Q,~1_R\mapsto 1_{R\amalg R'}~~~~\txt{and}~~~~R'\hookrightarrow R\amalg R',~r'\mapsto \big((0,r')\big)+Q,~1_{R'}\mapsto 1_{R\amalg R'}.
\]}
\end{example}

\begin{dfn}[\blue{Products and coproducts (direct sums) of modules}]
Let $\M:=(M_i)_{i\in {I}}$ be an indexing of $R$-modules $M_i$. The \index{Product of! modules}{\ul{product}} of the $R$-modules $\M$ is the smallest $R$-module (up to isomorphism of $R$-modules) having each member $M_i\in\M$ as a quotient $R$-module. Explicitly, the product of $\M$ is the set
\bea
\textstyle\prod\M:=\prod_{i\in I}M_i:=\big\{\txt{maps}~~m:{I}\ra\bigcup_i M_i,~i\mapsto m_i\in M_i\big\}=\big\{(m_i)_{i\in{I}}~|~m_i\in M_i\big\},\nn
\eea
as an $R$-module with \emph{addition}, \emph{scalar multiplication}, and \emph{zero} respectively given \emph{pointwise} by
\bea
(m_i)_{i\in{I}}+(m'_i)_{i\in{I}}:=(m_i+m_i')_{i\in{I}},~~~~r\cdot(m_i)_{i\in{I}}:=(rm_i)_{i\in{I}},~~~~0_{\prod\M}:=(0_{M_i})_{i\in{I}}.\nn
\eea
The \index{Coproduct! (Direct sum) of modules}{\ul{coproduct} (\ul{direct sum})} of the $R$-modules $\M=(M_i)_{i\in I}$ is the $R$-module generated by $\M$, in the sense that it is the smallest $R$-module (up to isomorphism of $R$-modules) containing every member of $\M$ as a submodule. Explicitly (up to isomorphism), the coproduct (\ul{direct sum}) of the $R$-modules $\M$ is the \emph{submodule} of $\prod\M$ given by
{\small\begin{align}
\textstyle\coprod\M~~\txt{or}~~\bigoplus\M:=\big\{\txt{finite maps}~~m:{I}\ra\bigcup_i M_i,~i\mapsto m_i\in M_i\big\}=\big\{(m_i)_{i\in{I}}\in\prod\M~|~m_i=0_{M_i}~~\txt{a.e.f.}~~i\in{I}\big\},\nn
\end{align}}where ``a.e.f.'' means ``\emph{for all except finitely many}''. In particular, $\bigoplus_{i\in I} R$ will mean the coproduct of copies of $R={}_RR$ or $R={}_\Integer R$ (i.e., $R$ is viewed as an $R$-module or abelian group).

If $(M_1,...,M_n)$ is a finite indexing of $R$-modules, then $\prod_{i=1}^nM_i\cong\bigoplus_{i=1}^nM_i$, which is also written as
\[
M_1\times M_2\times\cdots\times M_n\cong M_1\oplus M_2\oplus\cdots\oplus M_n.
\]
\end{dfn}

\begin{dfn}[\blue{Products and coproducts of algebras}]
Let $\A:=(A_i)_{i\in {I}}$ be an indexing of $R$-algebras $A_i$. The \index{Product of! algebras}{\ul{product}} of the $R$-algebras $\A$ is the smallest $R$-algebra (up to isomorphism of $R$-algebras) having each member $A_i\in\A$ as a quotient $R$-algebra. Explicitly, the product of $\A$ is the product of rings
\bea
\textstyle\prod\A:=\prod_{i\in I}A_i:=\big\{\txt{maps}~~a:{I}\ra\bigcup_i A_i,~i\mapsto a_i\in A_i\big\}=\big\{(a_i)_{i\in{I}}~|~a_i\in A_i\big\},\nn
\eea
as an $R$-algebra with \emph{scalar multiplication} given \emph{pointwise} by
\bea
r\cdot(a_i)_{i\in{I}}:=(ra_i)_{i\in{I}}.\nn
\eea

The \index{Coproduct! of algebras}{\ul{coproduct}} $\coprod\A=\coprod_{i\in I}A_i$ of the $R$-algebras $\A=(A_i)_{i\in I}$ is the $R$-algebra generated by $\A$, in the sense that it is the smallest $R$-algebra (up to isomorphism of $R$-algebras) containing every member of $\A$ as a subalgebra.
\end{dfn}

\begin{example}
In particular, if $\A=(A,A')$, consider the product of $R$-modules ${M_{R,A,A'}}:=$ ``$A\times A'=A\times\{0\}+\{0\}\times A'$ as an $R$-module'' and its tensor algebra $Te({M_{R,A,A'}}):=R\langle{M_{R,A,A'}}\rangle:=\sum_{i\in\Natural}R{M_{R,A,A'}}^i$. With $S:=(A\times\{0\})\backslash\{(1_A,0)\}$ and $S':=(\{0\}\times A')\backslash\{(0,1_{A'})\}$, and noting $(0,0)\in S\cap S'$, we have
{\small\begin{align}
&M_{R,A,A'}=(S\sqcup\{(1_A,0)\})+(\{(0,1_{A'})\}\sqcup S')=(S+S')\sqcup\big(S+\{(0,1_{A'})\}\big)\sqcup\big(\{(1_A,0)\}+S'\big)\sqcup\big\{(1_A,1_{A'})\big\}\nn\\
&~~~~=R(S+S')+R\{(1_A,0),(0,1_{A'})\}=R S+R S'+R(1_A,0)+R(0,1_{A'}).\nn
\end{align}}Let (i)~$T(A,A'):=\bigcup_{n\geq 1}(S\sqcup S')^n$ and (ii) for all $a,b\in S\sqcup S'$, let ``$a\sim b$ if either both $a,b\in S$ or both $a,b\in S'$''. Given $(a_1,...,a_m),(b_1,...,b_n)\in T(A,A')\subset Te({M_{R,A,A'}})$, define
{\scriptsize\[
D(a_1,...,a_m,b_1,...,b_n):=\left\{
  \begin{array}{ll}
    (a_1,...,a_m)\cdot (b_1,...,b_n)-(a_1,...,a_m,b_1,...,b_n), & \txt{if}~~a_m\not\sim b_1\\
    (a_1,...,a_m)\cdot (b_1,...,b_n)-(a_1,...,a_{m-1},a_mb_1,b_2,...,b_n), & \txt{if}~~a_m\sim b_1~\txt{and}~a_mb_1\not\in\{(1_A,0),(0,1_{A'})\}\\
    (a_1,...,a_m)\cdot (b_1,...,b_n)-(a_1,...,a_{m-1})\cdot(b_2,...,b_n), & \txt{if}~~a_m\sim b_1~\txt{and}~a_mb_1\in\{(1_A,0),(0,1_{A'})\}
  \end{array}
\right\}.
\]}Consider the $R$-submodules $\wt{A}:=A\times\{0\}\subset A\times A'$ and $\wt{A}':=\{0\}\times A'\subset A\times A'$, and the ideal
\[
Q:=\left\langle \left\{D(a_1,...,a_m,b_1,...,b_n)~\big|~(a_1,...,a_m),(b_1,...,b_n)\in T(A,A')\right\}\right\rangle\lhd Te({M_{R,A,A'}}).\nn
\]
With the subsring $R\langle \wt{A}\ast\wt{A}'\rangle\leq Te({M_{R,A,A'}})$ generated by the coproduct of multiplicative-monoids $\wt{A}\ast\wt{A}'$,
\begin{align}
\textstyle A\amalg A'&\textstyle:=\coprod(A,A'):={Te({M_{R,A,A'}})\over Q}={R\langle \wt{A}\ast\wt{A}'\rangle+Q\over Q}={\sum_{i\in\Natural}R\langle \wt{A}\ast\wt{A}'\rangle\cap{M_{R,A,A'}}^i+Q\over Q}=\sum_{i\in\Natural}{R\langle \wt{A}\ast\wt{A}'\rangle\cap{M_{R,A,A'}}^i+Q\over Q}\nn\\
&=\textstyle\big\{\sum_ir_ic_i+Q:r_i\in R,~c_i\in R\langle \wt{A}\ast\wt{A}'\rangle\cap{M_{R,A,A'}}^i,~r_i=0~\txt{a.e.f.}~i\in I\big\}\nn\\
&\textstyle={R\langle S+S'\rangle+Q\over Q},\nn
\end{align} where $R\langle S+S'\rangle\leq Te({M_{R,A,A'}})$ is the subring generated by the subset $S+S'\subset {M_{R,A,A'}}$.

By construction, $A\amalg A'$ is a ring with unity $1_{A\amalg A'}:=()+Q$, i.e., the class of the empty tuple. The component $R$-algebras $A,A'$ are included via the $R$-algebra homomorphisms
{\small\[
\textstyle A\hookrightarrow A\amalg A',~a\mapsto \big((a,0)\big)+Q,~1_A\mapsto 1_{A\amalg A'}~~~~\txt{and}~~~~A'\hookrightarrow A\amalg A',~a'\mapsto \big((0,a')\big)+Q,~1_{A'}\mapsto 1_{A\amalg A'}.
\]}
\end{example}

\section{Vector Space Dimension and Rank of a Free IBN-module}
\begin{thm}[\blue{\index{Steinitz replacement theorem}{Steinitz replacement theorem}}]
Let $V$ be a vector space, $S=\{v_1,...,v_n\}$ a finite spanning set for $V$, and $T=\{w_1,...,w_m\}$ a finite linearly independent set in $V$. Then $n\geq m$.
\end{thm}
\begin{proof}
The set $S\cup\{w_1\}=\{w_1,v_1,...,v_n\}$ is linearly dependent (since $\Span S=V$ and $w_1\neq 0$ by the linear independence of $T$) and so we can remove some $v_i$, say wlog $v_i=v_1$, and still have a spanning set $S_1:=\{w_1,v_2,...,v_n\}$. (\blue{footnote}\footnote{If $\al_1w_1+\al v_1+\al_2v_2+\cdots+\al_nv_n=0$ for scalars $\{\al,\al_1,...,\al_n\}$, where wlog $\al\neq 0$, then $v_1=-{1\over\al}\left(\al_1w_1+\al_2v_2+\cdots+\al_nv_n\right)$. }). Similarly, in the linearly dependent union $S_1\cup\{w_2\}$, we can wlog remove $v_2$ and still have a spanning set $S_2:=\{w_1,w_2,v_3,...,v_n\}$. Proceeding this way (wlog), at the $t$th step ($1\leq t\leq n$), we get a spanning set $S_t:=\{w_1,\cdots,w_t,v_{t+1},\cdots,v_n\}$. Suppose $n<m$. Then at step $t:=n$, $S_n:=\{w_1,\cdots,w_n\}\subsetneq T$ spans $V$, and so spans the remaining vectors $\{w_{n+1},...,w_m\}\subset T$ (a contradiction of the linear independence of $T$). Hence $n\geq m$.
\end{proof}

\begin{crl}[\textcolor{blue}{\index{Dimension of a vector space}{Dimension of a vector space}}]
Let $V$ be a $k$-vector space (for a field $k$). Then all bases of $V$ have the same cardinality (called the \ul{dimension} of $V$, and denoted by $\dim_kV$).
\end{crl}
\begin{proof}
Let $X,Y\subset V$ two bases for $V$. By Theorem \ref{SThCardThm1}, it is enough to assume $X,Y$ are finite sets. Let $X=\{x_1,...,x_n\}$ and $Y=\{y_1,...,y_m\}$. Since the bases are linearly independent spanning sets, the preceding theorem implies $n\geq m$ and $m\geq n$. Hence $n=m$.
\end{proof}

\begin{crl}
Let $V,W$ be $k$-vector spaces (for a field $k$). Then the following hold:
\bit
\item[(1)] If $V\cong W$, then ~$\dim_kV=\dim_kW$.
\item[(2)] If $V,W$ are finite-dimensional, then ~$V\cong W$ $\iff$ $\dim_kV=\dim_kW$.
\item[(3)] $k^n\cong k^m$ $\iff$ $m=n$ (where $k^n:=k\times k^{n-1}$).
\eit
\end{crl}
\begin{proof}
A $k$-isomorphism $f:V\ra W$ maps a spanning set to a spanning set and a linearly independent set to a linearly independent set.
\end{proof}

\begin{dfn}[\textcolor{blue}{
\index{Invariant! basis number (IBN) ring}{Invariant basis number (IBN) ring},
\index{IBN-module}{IBN-module}
}]
A ring $R$ is an \ul{invariant basis number (IBN) ring} if for any positive integers $n,m$, an isomorphism ~$R^n\cong R^m$ implies $n=m$, ~where $R^n:=R^{n-1}\times R$. An \ul{IBN-module} is a module over an IBN-ring.
\end{dfn}

\begin{dfn}[\textcolor{blue}{\index{Rank! of a free IBN-module}{Rank of a free IBN-module}}]
Let $R$ be an IBN ring, and $F$ a free R-module. The \ul{rank} of $F$ is the cardinality of any basis for $F$. (\blue{footnote}\footnote{The \ul{rank} of a free IBN-module is analogous to the \ul{dimension} of a vector space.})
\end{dfn}

Next, we will show in Theorem \ref{SThCardThm2} that every commutative ring is an IBN ring.

\begin{rmk}[\blue{${R\over I}$-module correspondence}]\label{RoverIrmk}
Let $R$ be a ring, $I\lhd R$ an ideal, and $M$ an $R$-module. If $IM=0$, then the map
${R\over I}\times M\ra M,~(r+I,m)\mapsto rm$ is well-defined (\blue{footnote}\footnote{If $(r+I,m)=(r'+I,m')$, then $r+I=r'+I$ or $r-r'\in I$ and $m=m'$, and so $rm=rm'=(r-r')m'+r'm'=0+r'm'$.}), thereby making $M$ an ${R\over I}$-module. Moreover, because the original $R$-scalar multiplication $R\times M\ra M,~(r,m)\mapsto rm$ has exactly the same effect on $M$ as the ${R\over I}$-scalar multiplication, we have an automatic bijective correspondence
\bea
\textstyle \left\{\txt{$R$-submodules ${}_RN\subset M$}\right\}~~\sr{id}{\longleftrightarrow}~~\left\{\txt{${R\over I}$-submodules ${}_{R\over I}N\subset M$}\right\}.\nn
\eea
\end{rmk}

\begin{thm}[\textcolor{blue}{Basis cardinality theorem II}]\label{SThCardThm2}
Let $F$ be a free R-module with two bases $X$ and $Y$. If $R$ is commutative then $|X|=|Y|$. (Hence every commutative ring is an IBN ring).
\end{thm}
\begin{proof}
By Theorem \ref{SThCardThm1}, $|X|$ is infinite $\iff$ $|Y|$ is infinite, in which case $|X|=|Y|$. So, we can assume $m:=|X|<\infty$ and $n:=|Y|<\infty$. By hypotheses, we have $F\cong R^m\cong R^n$. Let $M\vartriangleleft R$ be a maximal ideal. Recall $MR=M$ since $1\in R$. Let $\al:R^m\ra R^n$ be an isomorphism. Then
    \bea
    \al(M^m)=\al(MR^m)=M\al(R^m)=MR^n=M^n,\nn
    \eea
    and so the restriction {\small $\al|_{M^m}:M^m\ra M^n$} is an isomorphism. Thus, because the kernel of the map {\small $f:R^m\ra{R^n\over M^n}$} is $M^m$, it follows by the 1st isomorphism theorem that $\al$ induces an isomorphism of $R$-modules
    \bea
    \textstyle\ol{\al}:{R^m\over M^m}\ra {R^n\over M^n},~~(r_1,...,r_m)+M^m\mapsto\al(r_1,...,r_m)+M^n,~~~~\txt{or}~~~~\ol{\al}:\left({R\over M}\right)^m\ra \left({R\over M}\right)^n.\nn
    \eea
   Since $M\left({R\over M}\right)^m=0$ and $M\left({R\over M}\right)^n=0$, $\ol{\al}$ is also an isomorphism of $R/M$-modules. When $R$ is commutative, $k:=R/M$ is a field, and so $\ol{\al}$ is an isomorphism of $k$-vector spaces, i.e., $|X|=m=n=|Y|$. \qedhere
\end{proof}

\section{Central Subring Characterization of an Algebra}
\begin{thm}
Let $R$ be a commutative ring and $A$ a ring. Then $A$ is an $R$-algebra $\iff$ there exists a ring homomorphism $f:R\ra Z(A)\subset A$. (Equivalently, $A$ is an $R$-algebra $\iff$ there exists an ideal $I\vartriangleleft R$ and a subring $S\leq Z(A)\subset A$ such that $S\cong {R\over I}$.)
\end{thm}
\begin{proof}
{\flushleft ($\Ra$)}: If $A$ is an $R$-algebra, we get the ring homomorphism $f:R\ra Z(A)\subset A,~r\mapsto r1_A$. ($\La$): Conversely, if $f:R\ra Z(A)\subset A$ is a ring homomorphism, then $A$ is an $R$-algebra with respect to the operation ~$\cdot:R\times A\ra A,~(r,a)\mapsto r\cdot a:=f(r)a$.
\end{proof}

\begin{rmk}[\blue{In an $R$-algebra $A$, we can write $R\subset Z(A)$}]
In the above theorem, the ideal $I=\ker f\vartriangleleft R$ satisfies $I\cdot A=0$, since $I\cdot S\cong I\cdot{R\over I}=0$. So, by Remark \ref{RoverIrmk}, $R$-submodules of $A$ are precisely ${R\over I}$-submodules of $A$. Hence, we can assume without loss of generality that $R\subset Z(A)$, i.e., we can set $I=0$ in the actual situation ${R\over I}\subset Z(A)$.
\end{rmk}

\begin{crl}
Let $R$ be a commutative ring and $A$ a ring. Then $A$ is an $R$-algebra $\iff$ there exists an injective ring homomorphism $f:R\hookrightarrow Z(A)\subset A$.
\end{crl}

\begin{crl}
Given a ring $A$, the $R$-algebra structures on $A$ correspond to the subrings of $Z(A)$.
\end{crl}

\section{Graded Abelian Groups, Rings, Modules, and Algebras}
Given two classes of objects $\I$ and $\C$, an \index{Grading! (Sprectrum or Spectral decomposition)}{\ul{$\I$-grading} (or \ul{$\I$-spectrum} or \ul{$\I$-spectral decomposition})} $F$ of $\C$ (making $\C=(\C,F)$ an \index{Graded! class}{\ul{$\I$-graded class}}) is a map $F:\I\ra\P(\C)$ such that $\bigcup_{i\in\I}F(i)=\C$, and if necessary, together with regularity conditions on the indexed collection $\big(\C_i\big)_{i\in\I}:=\big(F(i)\big)_{i\in \I}$. If $\I=\Integer$, we simply say $\C$ is \ul{graded}. If $\I=\Integer^n$ (for some $n\in\Natural$) we say $\C=\bigcup_{i\in\Integer^n}\C_i=\bigcup_{i_j\in\Integer}\C_{i_1,...,i_n}$ is \index{Multiply graded class}{\ul{$n$-graded} (or \ul{$n$-multigraded}, or \ul{$n$-multiply graded})}. In particular, if $\I=\Integer^2$, we say $\C=\bigcup_{i,j\in\Integer}\C_{i,j}$ is \index{Bigraded! class}{\ul{bigraded}}.

The following are examples of commonly used gradings of abelian groups, rings, modules, and algebras. One will notice that they all involve the direct sum $\oplus$ of abelian groups (i.e., $\Integer$-modules). For more on grading and applications, see chapter \ref{AlgebraCatS12}.

Even though the objects (abelian groups, rings, modules, algebras) are assumed to be associative in the following definitions, it is not difficult to see that the concepts therein also apply (with appropriate modifications where necessary) to non-associative objects such as Lie algebras.

\begin{dfn}[\textcolor{blue}{\index{Grading! on an abelian group}{Grading on an abelian group}, \index{Graded! abelian group}{Graded abelian group}, \index{Homogeneous components}{Homogeneous components}, \index{Trivial grading}{Trivial grading}, \index{Graded! group homomorphism}{Graded group homomorphism}}]
Let $G=(G,e)$ be an abelian group and $I$ a set. An \ul{$I$-grading} $G_I$ on $G$ (making $G=(G,G_I)$ an \ul{$I$-graded abelian group}) is an indexed collection of abelian groups $G_I:=(G_i)_{i\in I}$ such that ~{\small $G=\bigoplus G_I=\bigoplus_{i\in I}G_i$}. The abelian groups $G_i$ are called the \ul{homogeneous components} of $G$ with respect to the grading $G_I$. The grading $G_I$ is \ul{trivial} if there exists $i\in I$ such that $G_j=\{e\}$ for all $j\neq i$ (i.e., $G_i=G$).

A group homomorphism {\small $f:G=\bigoplus_{i\in I}G_i\ra H=\bigoplus_{i\in I}H_i$} is \ul{$I$-graded} (\blue{footnote}\footnote{A graded group homomorphism is meant to be a homomorphism in the class of graded groups.}) if $f(G_i)\subset H_i$ for all $i\in I$.
\end{dfn}

\begin{dfn}[\textcolor{blue}{\index{Graded! ring}{Graded ring}, \index{Graded! ring homomorphism}{Graded ring homomorphism}}]
Let $R$ be a ring and $I=(I,e)$ a group (or any associative identity set, i.e., monoid) with identity element $e$. Then $R$ is $I$-graded if the following hold:
\bit
\item[(i)] The abelian group $(R,+)$ is $I$-graded.
\item[(ii)] The given $I$-grading $(R,+)=\bigoplus_{i\in I}(R_i,+)$ satisfies $R_iR_j\subset R_{i\cdot j}$ for all $i,j\in I$.
\item[(iii)] $1_R\in R_e$. (If $R$ has no unity or preferred element, we may replace $I$ with any associative set.)
\eit
A ring homomorphism $f:R=\bigoplus_{i\in I}R_i\ra S=\bigoplus_{i\in I}S_i$ is \ul{$I$-graded} (\blue{footnote}\footnote{A graded ring homomorphism is meant to be a homomorphism in the class of graded rings.}) if
\bea
f(R_i)\subset S_i~~~~\txt{for all}~~i\in I.\nn
\eea
\end{dfn}

\begin{dfn}[\textcolor{blue}{\index{Graded! module}{Graded module}, \index{Graded! module homomorphism}{Graded module homomorphism}}]
Let $M$ be an $R$-module and $I=(I,e)$ a group (or any associative identity set, i.e., monoid) with identity element $e$ . Then $M$ is $I$-graded if the following hold:
\bit
\item[(i)] The ring $R$ is $I'$-graded with respect to a subgroup $I'\leq I$.
\item[(ii)] The abelian group $(M,+)$ is $I$-graded.
\item[(iii)] The given gradings $(R,+)=\bigoplus_{i'\in I'}(R_{i'},+)$ and $(M,+)=\bigoplus_{i\in I}(M_i,+)$ satisfy
    \[
    R_{i'}M_i\subset M_{i'\cdot i}~~\txt{for}~~i'\in I',i\in I,~~\big(\txt{in addition to}~~\{R_{i'}R_{j'}\subset R_{i'\cdot j'}\}_{i',j'\in I'}~~\txt{and}~~1_R\in R_e\big).
    \]
\eit
Let $Z_I(I'):=\{i\in I:ij=ji~\txt{for all}~j\in I'\}\leq I$ be the subgroup by elements of $I$ that commute with all elements of $I'$. An $R$-homomorphism {\small $f:M=\bigoplus_{i\in I}M_i\ra N=\bigoplus_{i\in I}N_i$} is \ul{$I$-graded of degree $d\in Z_I(I')$} if
\bea
f(M_i)\subset N_{d\cdot i},~~~~\txt{for all}~~~~i\in I.\nn
\eea
If $f$ is $I$-graded of degree $e$, we simply say $f$ is \ul{$I$-graded}. (\blue{footnote}\footnote{A graded $R$-module homomorphism is meant to be a homomorphism in the class of graded $R$-modules.})
\end{dfn}

\begin{dfn}[\textcolor{blue}{\index{Graded! algebra}{Graded algebra}, \index{Graded! algebra homomorphism}{Graded algebra homomorphism}}]
Let $R$ be a commutative ring, $A$ an $R$-algebra, and $I=(I,e)$ a group (or any associative identity set, i.e., monoid) with identity element $e$. Then $A$ is $I$-graded if the following hold:
\bit
\item[(i)] The module $_RA$ is $I$-graded. ($R=\bigoplus_{i'\in I'}R_{i'}$ and $A=\bigoplus_{i\in I} A_i$ for a subgroup $I'\leq I$).
\item[(ii)] The ring $A$ is $I$-graded. ($A=\bigoplus_{i\in I} A_i$).
\eit
If $A,B$ are $I$-graded $R$-algebras and $i\in I$, then an $R$-algebra homomorphism $f:A\ra B$ is $I$-graded if it is $I$-graded as a ring homomorphism.
\end{dfn}

Note in the \ul{last} definition above (i.e., of an $I$-graded $R$-algebra) that the given gradings $R=\bigoplus_{i'\in I'}R_{i'}$ and $A=\bigoplus_{i\in I}A_i$ necessarily satisfy $1_R=1_A\in R_e\subset A_e$ along with
\bea
R_{i'}R_{j'}\subset R_{i'\cdot j'},~~R_{i'}A_i\subset A_{i'\cdot i},~~A_iA_j\subset A_{i\cdot j},~~~~\txt{for all}~~~~i',j'\in I',~~i,j\in I.\nn
\eea

In the \ul{last two} definitions above (i.e., of an $I$-graded $R$-module and of an $I$-graded $R$-algebra), if we set $I':=\{e\}\leq I$ in the $I'$-grading on $R$, then the homogeneous components of an $I$-graded $R$-module (or $I$-graded $R$-algebra) are themselves $R$-modules.

%% file: parts/AlgebraCat/AlgebraCatS1.tex
\chapter{Categories and Morphisms}\label{AlgebraCatS1}
From here on, if needed, additional sources of reading include \cite{cohn1991,mitchell1965,grothdk1957,rotman2009,weibel1994,gelfand-manin2010,gelfand-manin1999,maclane1963,cartan-elbg1956}.

\section{Graphs, Categories, Subcategories}
\begin{dfn}[\textcolor{blue}{
\index{Graph over a class of objects}{{Graph over a class of objects}},
\index{Vertex}{Vertex},
\index{Edge}{Edge},
\index{Path! in a graph}{Path},
\index{Connected! vertices}{Connected vertices},
\index{Connected! graph}{Connected graph},
\index{Cycle (Loop) in a graph}{Cycle (Loop)},
\index{Acyclic graph}{Acyclic graph},
\index{Tree}{Tree},
\index{Subgraph}{{Subgraph}},
\index{Subtree}{Subtree},
\index{Spanning! tree}{Spanning tree},
\index{Set! graph}{Set graph},
\index{Small! graph}{Small graph},
\index{Locally! small graph}{Locally small graph}}]\label{CatDef0}
Let $V$ be a class of objects. A \ul{graph} $G$ over $V$ consists of \ul{vertices} $v\in V$ (i.e., every object of $V$ is a vertex of $G$), \ul{edges} $E:=\bigcup\{R_i\subset V\times V: i\in\I,~\txt{for a class $\I$}\}$ (\blue{footnote}\footnote{That is, $E$ is a class whose objects are members $e\in R_i$ of relations $R_i\subset V\times V$ in a class of relations $\{R_i\subset V\times V:i\in\I\}$ between objects of $V$. The situation can also be described in terms of a map/indexing ~$R:\I\ra\P(V\times V),~i\mapsto R_i$.}) with an edge $e_{uv}$, written
{\footnotesize$\bt[column sep=small] u\ar[r,dash,"e_{uv}"]&v~,\et$} between $u,v\in V$ iff $(u,v)\in E$ (\blue{footnote}\footnote{That is, $(u,v)\in R_i\subset E$ for some $i\in\I$.}). Briefly, we write $G=(V,E)=\big(V_G,E_G\big)$. If $u,v\in V$, we will denote the collection of all edges between $u$ and $v$ by $E_G(u,v)\approx\{e_{i,uv}:=R_i|_{\{u\}\times\{v\}}~|~(u,v)\in R_i,~i\in\I\}$. (\blue{footnote}\footnote{Here, the edge $u\sr{e_{i,uv}}{^{\ul{~~~~~}}}v$ is the way in which the two objects $u,v\in V$ are related under the relation $R_i\subset V\times V$ from the indexing $R:\I\ra\P(V\times V),~i\mapsto R_i$ of all relevant such relations on which $G$ is based.}). Therefore, $E=\bigcup_{u,v\in V}E_G(u,v)=\bigsqcup_{u,v\in V}E_G(u,v)$.

A \ul{path} in $G$ with endpoints $u,v\in V$ is a collection of vertices $\gamma=\{u=v_0,v_1,...,v_k=v\}$ such that $(v_{i-1},v_i)\in E$. In which case, we say the vertices $u$ and $v$ are \ul{connected} by $\gamma$, or that $\gamma$ is a path \ul{between} $u$ and $v$. A graph $G$ is a \ul{connected graph} if every two vertices $u,v$ in $G$ are connected by a path. Let $P_G(u,v)$ denote all paths between $u$ and $v$. Note that $E_G(u,v)\subset P_G(u,v)$.

In a graph $G$, a path $\gamma=\{u=v_0,v_1,...,v_k=v\}$ is a \ul{cycle} (\ul{loop}) if its endpoints coincide, i.e., if $u=v$. A graph is \ul{acyclic} if it contains no cycles. A \ul{tree} is a connected acyclic graph (or equivalently, a graph in which every two vertices are connected by a unique path).

Given two graphs $G=(V,E)$ and $G'=(V',E')$, we say $G$ is a \ul{subgraph} of $G'$ (written $G\subset G'$) if $V\subset V'$ and $E\subset E'$. If a subgraph $T\subset G$ is a tree, we say $T$ is a tree in $G$. A \ul{subtree} is a subgraph of a tree that is itself a tree.

Given a graph $G$, a \ul{spanning tree} of $G$ is a subgraph $T\subset G$ such that (i) $T$ is a tree and (ii) $V_T=V_G$. A graph $G=(V,E)$ is a \ul{small graph} if $E$ (or $V$) is a set. A graph $G=(V,E)$ is a \ul{locally small graph} if $E_G(u,v)$ is a set for all $u,v\in V$.
\end{dfn}

\begin{thm}[\textcolor{blue}{\index{Connected! graph theorem}{Connected graph theorem}}]
A small graph has a spanning tree iff connected.
\end{thm}
\begin{proof}
Let $G$ be a small graph. ($\Ra$): If $G$ has a spanning tree, then it is clear that $G$ is connected since a tree is connected. ($\La$): Conversely, assume $G$ is a connected graph. If $G$ is finite, then by repeatedly eliminating loops (where a loop is eliminated by removing/excluding one edge from the loop, a process that clearly preserves connectedness of the graph) we obtain a spanning tree of $G$. If $G$ is infinite, we can proceed as follows.

Consider the set $\P:=\{\txt{Trees}~~T\subseteq G\}$ as a poset under inclusion $\subseteq$. If $\{T_\ld\}_{\ld\in\Ld}$ is a chain in $\P$ (i.e., $T_{\ld_1}\subseteq T_{\ld_2}$ if $\ld_1\leq\ld_2$), then $T':=\bigcup T_\ld$ is also in $\P$. Indeed, $T'\subseteq G$ and $T'$ cannot contain a loop, otherwise the loop will lie in some $T_\ld$, and so $T'$ is a tree. That is, every chain in $\P$ has an upper bound in $\P$, and so by Zorn's lemma, $\P$ has a maximal element $T$. Suppose $T=(V_T,E_T)$ is not a spanning tree of $G$. Then $V_G\backslash V_T\neq\emptyset$. Since $G$ is connected, there is an edge $e$ connecting $V_T$ to a vertex $v\in V_G\backslash V_T$. It follows that $T_1=\big(\{v\}\cup V_T,\{e\}\cup E_T\big)$ is a tree in $\P$ strictly containing $T$ (a contradiction).
\end{proof}

\begin{dfn}[\textcolor{blue}{\index{Simple! graph}{Simple graph}, \index{Multigraph}{Multigraph}}]
A graph $G=(V,E)$ is \ul{simple} if it has at most one edge between any two vertices (i.e., $|E_G(u,v)|\leq 1$ for all $u,v\in V$). If a graph is not simple, it is called a \ul{multigraph}.
\end{dfn}

\begin{dfn}[\textcolor{blue}{
\index{Directed! graph}{{Directed graph}},
\index{Orientation (Direction) in a graph}{Orientation (Direction)},
\index{Arrow in a graph}{Arrow},
\index{Directed! path}{Directed path},
\index{Direction-connected vertices}{Direction-connected vertices},
\index{Direction-connected graph}{Direction-connected graph},
\index{Directed! cycle (Directed loop)}{Directed cycle (Directed loop)},
\index{Direction-acyclic graph}{Direction-acyclic graph},
\index{Direction-tree}{Direction-tree}}]
A graph $G=(V,E)$ is a \ul{directed graph} if every edge $e\in E$ connecting vertices $u,v\in E$ is assigned a \ul{direction or orientation} $\sign(e)\in\{\pm1\}$, making $e$ an \ul{arrow} either from $u$ to $v$ (written $u\sr{e}{\ral}v$), or from $u$ to $v$ (written $u\sr{e}{\lal}v$ or $v\sr{e}{\ral}u$).

A \ul{directed path} in $G$ with endpoints $u,v\in V$ is a collection of vertices $\gamma=\{u=v_0,v_1,...,v_k=v\}$ such that $(v_{i-1},v_i)\in E$, $i=1,...,k$, are edges with the same orientation (in which case the vertices $u,v$ are \ul{direction-connected}). A directed graph $G$ is \ul{direction-connected} if every two vertices in $G$ are direction-connected (i.e., connected by a directed path).

A \ul{directed cycle} (or \ul{directed loop}) in $G$ is a directed path $\gamma=\{u=v_0,v_1,...,v_k=v\}$ whose endpoints coincide, i.e., with $u=v$. A graph is a \ul{direction-acyclic graph} if it contains no directed cycles (even though it might be non-acyclic in the sense it might contain cycles that are not directed). A \ul{direction-tree} is a direction-connected direction-acyclic graph, or equivalently, a directed graph in which every two vertices are connected by a unique directed path (even though it might be non-acyclic, hence not a tree).

We will write $dE_G(u,v)$ (resp. $dP_G(u,v)$ ) for the collection of arrows (resp. directed paths) from $u$ to $v$. Therefore, as before, $E=\bigcup_{u,v\in V}dE_G(u,v)$ and $dE_G(u,v)\subset dP_G(u,v)$.
\end{dfn}

\begin{dfn}[\textcolor{blue}{Recall:
\index{Morphism}{Morphism},
\index{Endomorphism! in a class}{Endomorphism}}]\label{MorphismDef}
Given a class $\C$, a \ul{morphism in $\C$} (or \ul{$\C$-morphism}) is a \ul{rule} $f:A\ra B$, between objects $A,B\in\C$, that respects/preserves the class structure of $\C$ (in a specific/precise sense). A \ul{$\C$-endomorphism} is a $\C$-morphism $f:A\ra A$ from an object $A\in\C$ to itself.
\end{dfn}

\begin{dfn}[\textcolor{blue}{
\index{Category over a class}{{Category over a class}},
\index{Morphism classes}{Morphism classes},
\index{Composition! of morphisms}{Composition of morphisms},
\index{Identity! morphism}{Identity morphism},
\index{Non-associative category}{Non-associative category},
\index{Small! category}{Small category},
\index{Locally! small category}{Locally small category},
\index{Domain! of a morphism}{Domain of a morphism},
\index{Codomain of a morphism}{Codomain of a morphism},
\index{Simple! category}{Simple category},
\index{Opposite! category}{{Opposite category}},
\index{Subcategory}{Subcategory},
\index{Generated! subcategory}{Generated subcategory}
}]\label{CatDef1}
A category {\small $\C=\big(\C,\Mor_\C,\circ\big)=\big(\Ob\C,\Mor\C,\circ\big)=\big(\Ob\C,\Mor_\C(\Ob\C,\Ob\C),\circ\big)$} is a directed graph consisting of the following.
\begin{enumerate}[leftmargin=0.9cm]
\item[(1)] A \ul{class of objects} ~$\Ob\C$ as vertices of $\C$. (\blue{footnote}\footnote{If the objects $\Ob\C$ form a set, then the resulting category $\C$ is called a \ul{small category}.}: \ul{Small category})
\item[(2)] A \ul{class of morphisms} ~{\small$\Mor_\C(A,B)=\left\{f:A\ra B~~\txt{or}~~A\sr{f}{\ral}B\right\}$} for every $A,B\in\Ob\C$, forming the edges ~{\small $\Mor\C =\Mor_\C(\Ob\C,\Ob\C):=\bigcup_{A,B\in\Ob\C}\Mor_\C(A,B)$}~ of $\C$,~ such that the following exist/hold. (\blue{footnote}\footnote{If the classes $\Mor_\C(A,B)$ are all sets, the category $\C$ is called a \ul{locally small} category.}: \ul{Locally small category})
\item[(3)] \ul{An associative composition operation}: Given by ~{\small$\circ :\Mor_\C(A,B)\times \Mor_\C(B,C)\ra \Mor_\C(A,C)$},
{\small\bea
\left(A\sr{f}{\ral}B,B\sr{g}{\ral}C\right)~\longmapsto~A\sr{gf}{\ral}B~~=~~A\sr{f}{\ral}B\sr{g}{\ral}C,~~~~\txt{for all}~~~~A,B,C\in\Ob\C,\nn
\eea}where associativity means {\footnotesize$\bt[column sep=small] (hg)f=h(gf):A\ar[r,"f"]&B\ar[r,"g"]&C\ar[r,"h"]&D,\et$} and the composition, which we have written simply as $gf$, is sometimes also explicitly written as $g\circ f$ when the emphasis is essential.

(If associativity of composition is dropped, then we get a \ul{non-associative category}.)

\item[(4)] \ul{Identity morphisms}:~ {\small$\big\{1_A=id_A\in \Mor_\C(A,A):A\in\Ob\C\big\}$} such that $1_Bf=f=f1_A$ for any morphism {\small$A\sr{f}{\ral}B$} (i.e. for any {\small$f\in \Mor(A,B)$}).
\item[(5)] \ul{Pairwise disjointness}:~ {\small$\Mor_\C(A,B)\cap \Mor_\C(A',B')=\emptyset$}, unless $A=A'$  and $B=B'$.
\end{enumerate}

{\flushleft If} $f\in \Mor_\C(A,B)$, we define the \ul{domain} and \ul{codomain} of $f$ as $\dom f:=A$ and $\cod f:=B$.

{\flushleft The} category $\C$ is a \ul{simple category} if it is simple as a graph, i.e., if $|\Mor_\C(A,B)|\leq 1$ for all $A,B\in\Ob\C$. The \ul{opposite category} of $\C$ is the category
{\small\bea
\C^{op}=\big(\C^{op},\Mor_{\C^{op}},\circ\big)=\big(\Ob\C^{op},\Mor\C^{op},\circ\big)=\big(\Ob\C^{op},\Mor_{\C^{op}}(\Ob\C^{op},\Ob\C^{op}),\circ\big)\nn
\eea}with objects $\Ob\C^{op}:=\Ob\C$, morphisms given (for all $A,B\in\Ob\C^{op}=\Ob\C$) by \ul{new-arrow assignments}
\bea
\txt{\small $\Mor_{\C^{op}}(A,B):=\Mor_\C(B,A)^{op}:=\left\{A\sr{f^{op}}{\ral}B~:~f\in \Mor_\C(B,A)\right\}$},\nn
\eea
and composition in $\C^{op}$ given in terms of composition in $\C$ by
\bea
f^{op}\circ g^{op}:=(g\circ f)^{op},~~~~\txt{i.e.,}~~~~f^{op}\circ_{\C^{op}} g^{op}:=(g\circ_\C f)^{op}.\nn
\eea

Given categories $\C=\big(\Ob\C,\Mor\C,\circ_\C\big)$  and $\D=\big(\Ob\D,\Mor\D,\circ_\D\big)$, we say $\C$ is a \ul{subcategory} of $\D$ (written $\C\subset\D$) if  the following hold:
\bit[leftmargin=0.9cm]
\item[(1)] The objects of $\C$ are objects of $\D$. That is, ~$\Ob\C\subset\Ob\D$.
\item[(2)] The morphisms of $\C$ are morphisms of $\D$. That is, ~$\Mor\C\subset\Mor\D$~ in the sense that for all $C,C'\in\Ob\C$, we have ~$\Mor_\C(C,C')\subset \Mor_\D(C,C')$.
\item[(3)] Composition in $\C$ is the same as composition in $\D$. That is, ~$\circ_\D|_{\Mor\C}=\circ_\C$.
\eit

Given a diagram $\D=(V,E)\subset\C$ in a category $\C$, the \ul{subcategory generated} by $\D$ is the smallest subcategory $\langle\D\rangle\subset\C$ of $\C$ containing $\D$, i.e.,
\[
\textstyle\langle\D\rangle:=\bigcap\{\txt{subcategories}~\C'\subset\C:~\D\subset\C'\}=\langle(V,E)\rangle=(V,\langle E\rangle),
\]
where ~$\langle E\rangle:=\{\txt{finite compositions of members of $E$}\}$.
\end{dfn}

\begin{notation}[\blue{Homomorphism classes}]
Let $\C$ be a category. Whenever (i) the objects $C,C'\in\Ob\C$ are classes and (ii) the morphisms $\Mor_\C(C,C')$ are maps, we will often (but not always) write $Hom_\C(C,C')$ instead of $\Mor_\C(C,C')$.
\end{notation}

\begin{rmk*}[\textcolor{blue}{Important simplifying assumption}]
Most of our major results will require/involve locally small categories (while general intermediate results will hold for general categories). Consequently, unless stated otherwise, we will henceforth assume that any given category is locally small. However, since many results about locally small categories extend to general categories, identifying and keeping track of such results would be a commendable exercise.
\end{rmk*}

\section{Selected Examples of Categories}
\begin{dfn}[\textcolor{blue}{\index{Filtered! category}{Filtered category}}]
A category $\C$ is called \ul{filtered} if it is
\bit[leftmargin=0.9cm]
\item[(1)] \index{Directed! category}{\ul{Directed}}: For any objects $A,B\in\Ob\C$, there exists $C\in\Ob\C$ ``bounding'' $A$ and $B$ in the sense
\bea
\Mor_\C(A,C)\neq\emptyset~~\txt{and}~~\Mor_\C(B,C)\neq\emptyset,\nn
\eea
or equivalently, there exist morphisms $f:A\ra C$, $g:B\ra C$, as in the diagram
\adjustbox{scale=0.7}{%
\bt[row sep=tiny] A\ar[dr,dashed,"f"] &  \\
     & C\\
    B\ar[ur,dashed,"g"'] & \et}
\item[(2)] \index{Equalizing category}{\ul{Equalizing}}: For any ``parallel'' morphisms $f,g\in \Mor_\C(A,B)$,
there exists a morphism $h\in \Mor_\C(B,C)$ ``equalizing'' $f,g$ in the sense that ~$h\circ f=h\circ g\in \Mor_\C(A,C)$.
\bea\bt
h\circ f=h\circ g~:~A\ar[r,shift left,"f"]\ar[r,shift right,"g"'] & B \ar[r,dashed,"h"] & C
\et\nn
\eea
\eit
\end{dfn}
\begin{note}[\textcolor{blue}{
\index{Poset category (Poset as a category)}{Poset category (Poset as a category)},
\index{Directed! set category (Directed set as a category)}{Directed set category (Directed set as a category)}
}]
We will always view a \ul{poset} (resp. \ul{directed set}) $\I=(\I,\leq)$ as a \ul{small category} (resp. \ul{filtered small category}) with \ul{objects} elements of the set $\I$, and \ul{morphisms} between objects $i,j\in\Ob\I:=\I$ given by
\bea
\Mor_\I(i,j):=
\left\{
  \begin{array}{ll}
    \left\{\txt{unique arrow}~~i\sr{\kappa_{ij}}{\ral}j\right\}, &~\txt{if}~~i\leq j \\
    \emptyset, & \txt{otherwise}
  \end{array}
\right.\nn
\eea
\end{note}
In the above note, we can also replace the poset (resp. directed set) with any partially ordered class (resp. directed class).

\begin{examples*}[\textcolor{blue}{\index{Examples of categories}{Some familiar categories}}]
Associated with each subclass of sets $\S\subset\txt{\ul{Sets}}$ is the category $\S=(\S,Hom_\S,\circ)$. Similarly, associated with any subclasses  $\G\subset\txt{\ul{Groups}}$ of groups, $\R\subset\txt{\ul{Rings}}$ of rings, $\M\subset\txt{\ul{Mod}}$ of modules, $\A\subset\txt{\ul{Alg}}$ of algebras, $\T\subset\txt{\ul{Top}}$ of topological spaces (\blue{footnote}\footnote{The category of topological spaces, Top, will be considered in chapter \ref{GeomAnaI}.})) are the categories
{\footnotesize\bea
\G=(\G,Hom_\G,\circ),~~\R=(\R,Hom_\R,\circ),~~\M=(\M,Hom_\M,\circ),~~\A=(\A,Hom_\A,\circ),~~ \T=(\T,Hom_\T,\circ),~~\cdots~),\nn
\eea}
where for sets $S,S'\in\S$, groups $G,G'\in\G$, rings $R,R'\in\R$, modules $M,M'\in\M$, algebras $A,A'\in\A$, spaces $X,X'\in\T$, $\cdots$, morphisms between them are given by the following.
\begin{enumerate}[leftmargin=0.9cm]
\item[$\ast$] $Hom_\S(S,S')$ is a set of maps ~$h:S\ra S'$.
\item[$\ast$] $Hom_\G(G,G')$ is a set of group homomorphisms ~$h:G\ra G'$.
\item[$\ast$] $Hom_\R(R,R')$ is a set of ring homomorphisms ~$h:R\ra R'$.
\item[$\ast$] $Hom_\M(M,M')$ is a set of module homomorphisms ~$h:M\ra M'$.
\item[$\ast$] $Hom_\A(A,A')$ is a set of algebra homomorphisms ~$h:A\ra A'$.
\item[$\ast$] $Hom_\T(X,X')$ is a set of homeomorphisms ~$h:X\ra X'$.
(\blue{footnote}\footnote{The category of topological spaces, Top, will be considered in chapter \ref{GeomAnaI}.})
\item[$\ast$] and so on.
\end{enumerate}
\end{examples*}

{\flushleft In }particular, for any ring $R$, the category of (left) $R$-modules is denoted by
\bea
\txt{$R$-mod $=(R\txt{-mod},Hom_R,\circ)$}.\nn
\eea
Similarly, the category of right $R$-modules is denoted by
\bea
\txt{mod-$R$ $=(\txt{mod-}R,Hom_R,\circ)$}.\nn
\eea
Of course, if $R$ is commutative, then the above left-right distinction is not necessary.

The category of abelian groups, $\Integer\txt{-mod}=(\Integer\txt{-mod},Hom_\Integer,\circ)$, is also denoted by
\bea
Ab=(Ab,Hom_{Ab},\circ).\nn
\eea

Notice that each of the categories given above is locally small, because (by Corollary \ref{CartProdSets}) given two sets $A$ and $B$, the collection $B^A$ of maps from $A$ to $B$ forms a set.

\section{Morphisms in a Category: Morphism Types and Properties}
\begin{dfn}[\textcolor{blue}{
\index{Monic! morphism (Monomorphism)}{Monic morphism (Monomorphism)},
\index{Epic! morphism (Epimorphism)}{Epic morphism (Epimorphism)},
\index{Restriction! of a morphism}{Restriction of a morphism},
\index{Unitary morphism (Unitary equivalence)}{{Unitary morphism (Unitary equivalence of objects)}},
\index{Unitarily equivalent objects}{Unitarily equivalent objects},
\index{Left! inverse}{Left inverse},
\index{Right! inverse}{Right inverse},
\index{Inverse! of a morphism}{Inverse of a morphism},
\index{Split! monomorphism}{Split monomorphism},
\index{Split! epimorphism}{Split epimorphism},
\index{Isomorphism}{Isomorphism},
\index{Isomorphic! objects}{Isomorphic objects},
\index{Endomorphism!}{Endomorphism},
\index{Automorphism}{Automorphism},
\index{Equivalence! in a category}{{Equivalence in a category}},
\index{Equivalence! of objects}{Equivalence of objects},
\index{Equivalent objects}{Equivalent objects},
\index{Pseudo-inverse}{Pseudo-inverse},
\index{Equality in a category}{Equality in a category},
\index{Equality of objects}{Equality of objects},
\index{Equal! objects}{Equal objects},
\index{Category of! morphisms}{Category of morphisms},
\index{Weak! isomorphism of objects}{Weak isomorphism of objects},
\index{Weakly isomorphic objects}{Weakly isomorphic objects}}]\label{CatDef2}
Fix a category {\small $\C=\big(\C,\Mor_\C,\circ\big)=\big(\Ob\C,\Mor\C,\circ\big)=\big(\Ob\C,\Mor_\C(\Ob\C,\Ob\C),\circ\big)$}, with morphisms $\Mor\C=\Mor_\C(\Ob\C,\Ob\C)$ $:=$ $\bigcup_{A,B\in\Ob\C}\Mor_\C(A,B)$.

A morphism $f\in \Mor_\C(A,B)\subset\Mor\C$ is \ul{monic} or a \ul{monomorphism}, written $f:A\hookrightarrow B$ (resp. \ul{epic} or an \ul{epimorphism}, written $f:A\twoheadrightarrow B$) if it is \ul{left-cancelling} in the sense $f\circ g=f\circ h~\Ra~g=h$ (resp. \ul{right-cancelling} in the sense $g\circ f=h\circ f~\Ra~g=h$), for all $g,h\in\Mor\C$. If $g$ is any morphism and $f$ is a monomorphism, then the composition $g\circ f$ is called the \ul{restriction} of $g$ to $f$ (or to $\cod f$), and denoted by $g|_f$ (or by $g|_{\cod f}$). A morphism $f\in \Mor_\C(A,B)\subset\Mor\C$ is \ul{unitary} (or a \ul{unitary equivalence} of objects, making $A$ and $B$ \ul{unitarily equivalent} or $A\sim_uB$), written $\bt[column sep=small] f:A\ar[r,hook,two heads]& B,\et$ if it is both monic and epic.

Given a morphism $f\in \Mor_\C(A,B)$, a \ul{left-inverse} (resp. \ul{right-inverse}) of $f$ is a morphism $g\in \Mor_\C(B,A)$ such that $g\circ f=id_A:A\sr{f}{\ral}B\sr{g}{\ral}A$ (resp. $f\circ g=id_B:B\sr{g}{\ral}A\sr{f}{\ral}B$). A (two-sided) \ul{inverse} of $f$ is a morphism $g\in \Mor_\C(B,A)$ that is both a left-inverse and a right-inverse of $f$, i.e., $g\circ f=id_A:A\sr{f}{\ral}B\sr{g}{\ral}A$ and $f\circ g=id_B:B\sr{g}{\ral}A\sr{f}{\ral}B$. (\blue{footnote}\footnote{By Remark \ref{CaUnInvRmk}, if both a left inverse and a right inverse exist, they are equal. Hence, if an inverse exists, it is unique.})

{\flushleft A} monomorphism (resp. epimorphism) is \ul{split} if it has a left-inverse (resp. right-inverse).
A morphism $f\in \Mor_\C(A,B)$ is an \ul{isomorphism}, written {\small $\bt[column sep=tiny] f:A\ar[rr]&\cong&B\et$} or {\small $\bt A\ar[r,"f","\cong"']&B\et$}, (making $A,B$ \ul{isomorphic objects}, written $A\cong B$) if it has an inverse $g\in \Mor_\C(B,A)$. Recall (from Definition \ref{MorphismDef}) that an \ul{endomorphism} is a morphism $f\in \Mor_\C(A,A)$ from an object $A$ to itself. An isomorphism $f\in \Mor_\C(A,A)$ is called an \ul{automorphism} of $A$.

{\flushleft An} \ul{equivalence in the category} $\C$ is an equivalence relation ~$\sim$ ~$\subset$ ~$\Mor\C\times\Mor\C$~ on $\Mor\C$. A morphism $f\in \Mor_\C(A,B)$ in $\Mor\C$ is a \ul{$\sim$-equivalence of objects} (or a $\sim$-induced equivalence of objects) if there exists a morphism $g\in \Mor_\C(B,A)$ (called \ul{$\sim$ pseudo-inverse} of $f$) such that $g\circ f\sim id_H:A\ra A$ and $f\circ g\sim id_B:B\ra B$. If a $\sim$-equivalence of objects $f\in \Mor_\C(A,B)$ exists, we say $A,B\in\Ob\C$ are \ul{$\sim$-equivalent objects}, written $A\sim B$.

{\flushleft In} particular, \ul{equality in the category} $\C$ is an equality relation $=$ (as a subclass of $\Mor\C\times\Mor\C$) on $\Mor\C$. A morphism $f\in \Mor_\C(A,B)$ in $\Mor\C$ that is an \ul{equality of objects} is precisely an isomorphism (in which case the pseudo-inverse of $f$ is an inverse of $f$). If an \ul{equality} (or \ul{isomorphism}) of objects $f\in \Mor_\C(A,B)$ exists, we say $A,B$ are \ul{equal objects} (or \ul{isomorphic objects}), written $A=B$ (or $A\cong B$).

{\flushleft Now} assume that the morphisms $\Mor\C$ of $\C$ form another category (the \ul{category of $\C$-morphisms})
{\small\bea
Cat(\Mor\C)=\big(Cat(\Mor\C),\Mor_{Cat(\Mor\C)},\circ\big)=\big(\Ob(Cat(\Mor\C)),\Mor(Cat(\Mor\C)),\circ\big)\nn
\eea}with objects $\Ob(Cat(\Mor\C)):=\Mor\C$ and morphisms arbitrarily written as
{\small\bea
\textstyle\Mor(Cat(\Mor\C))=\Mor_{Cat(\Mor\C)}(\Mor\C,\Mor\C)=\bigsqcup\limits_{f,g\in\Mor\C}\Mor_{Cat(\Mor\C)}(f,g),\nn
\eea}where a morphism in $Cat(\Mor\C)$, written ~$\al:f\ra g$~ or ~$\big(\dom f\sr{f}{\ral}\cod f\big)\sr{\al}{\ral}\big(\dom g\sr{g}{\ral}\cod g\big)$,~ is a pair of $\C$-morphisms ~$\al_1:\dom f\ra \dom g$~ and ~$\al_2:\cod f\ra \cod g$~ such that ~$\al_2f=g\al_1$. (\blue{footnote}\footnote{By construction, $\al$ is an isomorphism in $Cat(\Mor\C)$ $\iff$ $\al_1$ and $\al_2$ are isomorphisms in $\C$.})
\bea\bt
\dom f\ar[d,"f"]\ar[rr,"\al_1"] && \dom g\ar[d,"g"]\\
\cod f\ar[rr,"\al_2"] && \cod g
\et\nn
\eea
Then isomorphism $f\cong g$ (of objects in $Cat(\Mor\C)$ ) of morphisms of $\C$ is of course (i.e., by definition) an equivalence relation on $\Mor\C$, and hence an equivalence in $\C$. A morphism $f\in \Mor_\C(A,B)$ in $\Mor\C$ is a \ul{$\cong$-weak isomorphism of objects} in $\C$ if it is a \ul{$\cong$-equivalence of objects} (or a $\cong$-induced equivalence of objects), i.e., if there exists a morphism $g\in \Mor_\C(B,A)$ (a \ul{$\cong$-pseudo-inverse} of $f$) such that $g\circ f\cong id_H:A\ra A$ and $f\circ g\cong id_B:B\ra B$. If a weak isomorphism of objects $f\in \Mor_\C(A,B)$ exists, we say $A,B\in\Ob\C$ are \ul{weakly isomorphic objects}, written $A\cong_w B$.
\end{dfn}

\begin{rmk}[\blue{Uniqueness of the inverse}]\label{CaUnInvRmk}
A morphism $f:A\ra B$ has an inverse iff it has both a left-inverse and a right-inverse: Indeed if $f$ has an inverse, the inverse is clearly both a left-inverse and a right-inverse. Conversely, if $g:B\ra A$ is a left-inverse of $f$ and $h:B\ra A$ a right-inverse of $f$, then $fg=1_B$ and $hf=1_A$ imply $g=1_Ag=(hf)g=h(fg)=h1_B=h$. (This proof also shows that an inverse, when it exists, is unique.)
\end{rmk}

\begin{notes}
{\flushleft (i)} Even though it is so in many categories, in an arbitrary category, a unitary morphism $\bt[column sep=small] f:A\ar[r,hook,two heads]& B\et$ is not necessarily an isomorphism {\small $\bt[column sep=small] f:A\ar[r,draw=none,"\cong"description]&B.\et$} On the other hand, an isomorphism {\small $\bt[column sep=small] f:A\ar[r,draw=none,"\cong"description]&B\et$} is always a monic epic morphism $\bt[column sep=small] f:A\ar[r,hook,two heads]& B.\et$
{\flushleft (ii)} In a non-associative category, an inverse (if it exists) may not be unique.
{\flushleft (iii)} Except for educational reasons, it is not strictly necessary in practice to distinguish between equality, weak-isomorphism, and isomorphism.
\end{notes}

\begin{question}
It is clear that a morphism with a left-inverse (resp. right-inverse) inverse is monic (resp. epic). Is it true that a monic morphism (resp. an epic morphism) has a left inverse (resp. right-inverse)?
\end{question}

%% file: parts/AlgebraCat/AlgebraCatS2.tex
\chapter{Functors and Equivalences in Categories}\label{AlgebraCatS2}

In this chapter we will not be making a full introduction of diagrams, but to properly understand the definition of a natural transformation (and the composition of natural transformations in a category of functors) we need the basic idea of a commutative diagram in a category.
\begin{dfn}[\blue{
\index{Diagram in a category}{Diagram in a category},
\index{Diagram as a graph}{Diagram as a graph},
\index{Image of! a diagram}{Image of a diagram},
\index{Commutative! diagram}{Commutative diagram}}]
Let $\I,\C$ be categories (each of which is by definition a directed graph with objects as vertices and morphisms as edges). An \ul{$\I$-diagram} in $\C$ (written $D:\I\dashrightarrow\C$) is a map of the form $D:\Ob\I\sqcup\Mor\I\ra\Ob\C\sqcup\Mor\C$ such that its restrictions on objects and morphisms give separate maps
\[
D:\Ob\I\ra\Ob\C,~i\mapsto D(i)~~~~\txt{and}~~~~D:\Mor\I\ra\Mor\C,~\big(i\sr{\kappa}{\ral}j\big)\mapsto\big(D(i)\sr{D(\kappa)}{\ral}D(j)\big).
\]
The \ul{diagram} $D$ is equivalently the directed graph $D=(V,E)$ with vertices {\small $V:=\big(D(i)\big)_{i\in\Ob\I}$} and edges {\small $E:=\big(D(\kappa)\big)_{\kappa\in\Mor\I}$} (which must be viewed as indexings since they fully represent $D$ as a map). The following are worth noting:
\begin{enumerate}[leftmargin=0.9cm]
\item In typical applications we will encounter, $\I$ will often be a poset viewed as a category.
\item The directed graph $D=(V,E)$ representing the diagram $D:\I\dashrightarrow\C$ is not a subgraph of $\C$, even though its vertices are indexed objects of $\C$ and its edges are indexed morphisms of $\C$.
\item The \ul{image of the diagram} $D$ in $\C$ is the subgraph $D'=(V',E')\subset\C$ given by the unindexed collections $V':=\{D(i)\}_{i\in\Ob\I}\subset\Ob\C$ and $E':=\{D(\kappa)\}_{\kappa\in\Mor\I}\subset\Mor\C$.
\item A subgraph of $\C$ has the general form {\small $D'=(V',E')=\big(V'\subset\Ob\C,E'\subset \Mor_\C(V',V')\big)$}, where $\Mor_\C(V',V'):=\bigcup_{u,v\in V'}\Mor_\C(u,v)$, and it is a simple subgraph iff {\small $\big|E'\cap \Mor_\C(u,v)\big|\leq 1$} for all $u,v\in V'$.
\item We will see later that a subgraph of $\C$ is equivalently the image of a functor $\D:\I\ra\C$, for some category $\I$.
\end{enumerate}
The diagram $D=(V,E)$ is a \ul{commutative diagram} if for any pair of vertices $u,v\in V$, any two directed paths between them {\small $\gamma=\{f_i:v_i\ra v_{i+1}~|~1\leq i\leq n,~v_1=u,~v_{n+1}=v\}\subset E$} and {\small $\gamma'=\{f'_i:v'_i\ra v'_{i+1}~|~1\leq i\leq n',~v'_1=u,~v'_{n'+1}=v\}\subset E$}, or more explicitly,
{\footnotesize\[
 \gamma:=\big(u=v_1\sr{f_1}{\ral}v_2\sr{f_2}{\ral}\cdots\sr{f_n}{\ral}v_{n+1}=v\big)~~~~\txt{and}~~~~ \gamma':=\big(u=v'_1\sr{f'_1}{\ral}v'_2\sr{f'_2}{\ral}\cdots\sr{f'_{n'}}{\ral}v'_{n'+1}=v\big),
\]}give the same composition in the sense that (as morphisms in $\C$)
\[
f_nf_{n-1}\cdots f_1=f'_{n'}f'_{n'-1}\cdots f'_1:u\ra v.\nn
\]
\end{dfn}

\section{Functors between categories, Natural transformations}
\begin{dfn}[\textcolor{blue}{
\index{Functor}{{(Covariant) functor}},
\index{Endofunctor}{Endofunctor},
\index{Faithful! functor}{Faithful functor},
\index{Full functor}{Full functor},
\index{Fully! faithful functor}{Fully faithful functor},
\index{Dense! functor}{Dense functor},
\index{Morphism of! functors (Natural transformation)}{Morphism of functors (Natural transformation)},
\index{Isomorphism of! functors (Natural isomorphism)}{Isomorphism of functors (Natural isomorphism)},
\index{Isomorphic! functors}{Isomorphic functors},
\index{Category of! functors}{{Category of functors}},
\index{Contravariant functor (Cofunctor)}{\textcolor{blue}{Contravariant functor (Cofunctor)}},
\index{Morphism of! contravariant functors}{Morphism of contravariant functors},
\index{Isomorphism of! contravariant functors}{Isomorphism of contravariant functors},
\index{Isomorphic! contravariant functors}{Isomorphic contravariant functors},
\index{Opposite! functor of a category}{{Opposite functor of a category}},
\index{Category of! contravariant functors}{{Category of contravariant functors}}
}]\label{CatDef3}
Let $\C,\D$ be categories. A \ul{(covariant) functor} from $\C$ to $\D$ (written $F:\C\ra\D$ or $\C\sr{F}{\ral}\D$) is a map of the form
\bit[leftmargin=0.9cm]
\item $F:\Ob\C\sqcup\Mor\C\ra \Ob\D\sqcup\Mor\D$ ~whose restrictions to objects and morphisms give parallel maps
{\small\[
F:\Ob\C\ra\Ob\D,~C\mapsto F(C)~~~~\txt{and}~~~~F:\Mor\C\ra\Mor\D,~\Big(C\sr{f}{\ral}C'\Big)~\mapsto~\Big(F(C)\sr{F(f)}{\ral}F(C')\Big),\nn
\]}
\eit
and the map on morphisms has the following properties
(\blue{footnote}\footnote{
It is also enough to briefly specify the functor in the form $F:\C\ra\D,~\big(\dom f\sr{f}{\ral}\cod f\big)\mapsto\big(F(\dom f)\sr{f}{\ral}F(\cod f)\big)$, i.e., as a map on morphisms only. Indeed, the map on objects is automatically accounted for as follows: With $id_{\Ob\C}:=\{id_C:C\in\Ob\C\}$, the requirement ``$F(id_C)=id_{F(C)}$ for all $C\in\Ob\C$'' ensures the restriction $F:id_{\Ob\C}\ra id_{\Ob\D},~id_C\mapsto id_{F(C)}$ is a map that is clearly equivalent to the map on objects $F:\Ob\C\ra \Ob\D,~C\mapsto F(C)$.

Nevertheless, the description of various functors and their properties requires the full understanding of the functor $F:\C\ra\D$ as explicitly composed of two parallel maps $F:\Ob\C\ra\Ob\D$ and $F:\Mor\C\ra\Mor\D$. It is also worth noting that the map on morphisms $F:\Mor\C\ra\Mor\D$ is equivalently given by maps of the form $F:\Mor_\C(C,C')\ra \Mor_\D(F(C),F(C')),~f\mapsto F(f)$, for each pair of objects $C,C'\in\Ob\C$, in which context it is clearly no longer necessary to explicitly specify the domain and codomain of $f$ in the assignment rule $f\mapsto F(f)$.
}):
For all $C\in\Ob\C$ and $f,f'\in\Mor\C$,
\bit[leftmargin=0.9cm]
\item $F(id_C)=id_{F(C)}$~~ and ~~$F(f\circ f')=F(f)\circ F(f')$.
\eit

{\flushleft An} \ul{endofunctor} is a functor $F:\C\ra\C$ from a category $\C$ to itself.

{\flushleft The} functor $F:\C\ra\D$ is \ul{faithful} (resp. \ul{full}) if it is injective (resp. surjective) on morphisms, i.e., for each $C,C'\in\Ob\C$, the following map is injective (resp. surjective):
\begin{align}
{\small\txt{$F:\Mor_\C(C,C')\ra \Mor_\C(FC,FC'),~~\Big(C\sr{f}{\ral}C'\Big)\mapsto\Big(F(C)\sr{F(f)}{\ral}F(C')\Big)$}}.\nn
\end{align}
A functor that is both full and faithful is said to be \ul{fully faithful}.

The functor $F:\C\ra\D$ is \ul{dense} if it is surjective on isomorphism classes of objects, i.e., the induced map on isomorphism classes of objects $F:{\Ob\C\over\cong}\ra{\Ob\D\over\cong}$, $[C]\mapsto[F(C)]$ is surjective (equivalently, for any $D\in\Ob\D$, we have $D\cong F(C)$ for some $C\in\Ob\C$).

{\flushleft Let} $\I,\C$ be categories. (\blue{footnote}\footnote{(Instead of $\C,\D$) we are using $\I,\C$ simply so that the domain category $\I$ (while arbitrary) resembles an index set/class.}). Given functors $F,F':\I\ra\C$, a \ul{morphism of functors} (or \ul{natural transformation}) $\tau:F\ra F'$ is a collection of morphisms in $\C$
{\small\bea
\tau=\left\{F(i)\sr{\tau_i}{\ral}F'(i)\right\}_{i\in\Ob\I}~\subset~\Mor\C=\Mor_\C(\Ob\C,\Ob\C)\nn
\eea} such that for any $i,j\in\Ob\I$ and any $\kappa_{ij}\in \Mor_\I(i,j)$, we have a commutative diagram
{\small\bea
\left(F\sr{\tau}{\ral}F'\right)\left(i\sr{\kappa_{ij}}{\ral}j\right)~~:=~~
\bt
F_i\ar[d,"\tau_i"]\ar[rr,"f_{ij}"]&&F_j\ar[d,"\tau_j"]\\
F'_i\ar[rr,"f_{ij}'"]&&F'_j
\et~~:=~~
\bt
F(i)\ar[d,"\tau_i"]\ar[rr,"F(\kappa_{ij})"]&&F(j)\ar[d,"\tau_j"]\\
F'(i)\ar[rr,"F'(\kappa_{ij})"]&&F'(j)
\et~~~~~~~~\tau_j\circ f_{ij}=f'_{ij}\circ \tau_i.\nn
\eea}
A morphism of functors $\tau=\{\tau_i\}:F\ra F'$ is an \ul{isomorphism of functors} (or \ul{natural isomorphism}) if each $\tau_i:F_i\ra F'_i$ is an isomorphism in $\C$ (in which case, we say  $F,F'$ are \ul{isomorphic}, written $F\cong F'$).

{\flushleft The} \ul{category of functors from $\I$ to $\C$} is the category
\[
\C^\I=(\C^\I,\Mor_{\C^\I},\circ)=(\Ob\C^\I,\Mor\C^\I,\circ)=\big(\Ob\C^\I,\Mor_{\C^\I}(\Ob\C^\I,\Ob\C^\I),\circ\big)\nn
\]
in which (i) objects $\Ob\C^\I$ are functors $\I\sr{F}{\ral}\C$, (ii) morphisms are morphisms of functors
\[
\textstyle\tau~\in~\bigsqcup\limits_{F,F'\in\Ob\C^\I}\Mor_{\C^\I}(F,F'),
\]
and (iii) the composition of two morphisms of functors {\small $F\sr{\tau}{\ral} F'$} and {\small $F'\sr{\tau'}{\ral}F''$}, which is of course formally expressed as {\small $F\sr{\tau'\circ \tau}{\ral}F''~=~F\sr{\tau}{\ral}F'\sr{\tau'}{\ral}F''$} with componentwise form
{\small\begin{align}
\tau'\circ \tau&=\big\{F'_i\sr{\tau'_i}{\ral}F''_i~|~i\in\Ob\I\big\}\circ\big\{F_i\sr{\tau_i}{\ral}F'_i~|~i\in\Ob\I\big\}\nn\\
&:=\big\{F_i\sr{\tau'_i\circ \tau_i}{\ral}F''_i~|~i\in\Ob\I\big\}=\big\{F_i\sr{\tau_i}{\ral}F'_i\sr{\tau'_i}{\ral}F''_i~|~i\in\Ob\I\big\},\nn
\end{align}}
has the following property: For any $i,j\in\I$ and any $\kappa_{ij}\in \Mor_\I(i,j)$, we have a commutative diagram
{\small\[
\left(F\sr{\tau}{\ral}F'\sr{\tau'}{\ral}F''\right)\left(i\sr{\kappa_{ij}}{\ral}j\right)~~:=~~\bt
F_i\ar[d,"\tau_i"]\ar[rr,"f_{ij}"]&&F_j\ar[d,"\tau_j"]\\
F'_i\ar[d,"\tau'_i"]\ar[rr,"{f'_{ij}}"]&&F'_j\ar[d,"\tau'_j"]\\
F''_i\ar[rr,"f''_{ij}"]&&F''_j
\et~~~~~~~~\tau_j\circ f_{ij}=f'_{ij}\circ \tau_i,~~~~\tau'_j\circ f'_{ij}=f''_{ij}\circ \tau'_i.\nn
\]}

Diagrammatically (simply to aid our intuition) the morphism of functors $F\sr{\tau}{\ral}F'$ looks like

\bc\bt
i\ar[d,"id_i"]\ar[r,"\kappa_{ij}"]&j\ar[d,"id_j"]\ar[r,draw=none,"\in"description]&\Mor\I\ar[rr,"F"]&~\ar[d,"\tau=\{\tau_i\}"]&\Mor\C\ar[r,draw=none,"\ni"description]&F\big(i\sr{\kappa_{ij}}{\ral}j\big)\ar[r,draw=none,"="description]&F_i\ar[d,"\tau_i"]\ar[r,"f_{ij}"]&F_j\ar[d,"\tau_j"]\\
i\ar[r,"\kappa_{ij}"]&j\ar[r,draw=none,"\in"description]&\Mor\I\ar[rr,"F'"]&\mathop{~}\limits^{\substack{~\\~\\~}}&\Mor\C\ar[r,draw=none,"\ni"description]&F'\big(i\sr{\kappa_{ij}}{\ral}j\big)\ar[r,draw=none,"="description]&F'_i\ar[r,"f'_{ij}"]&F'_j\\
\et\ec
and similarly the composition of two morphisms of functors $F\sr{\tau}{\ral}F'$ and $F'\sr{\tau'}{\ral}F''$ looks like
\bc\bt
i\ar[d,"id_i"]\ar[r,"\kappa_{ij}"]&j\ar[d,"id_j"]\ar[r,draw=none,"\in"description]&\Mor\I\ar[rr,"F"]&~\ar[d,"\tau=\{\tau_i\}"]&\Mor\C\ar[r,draw=none,"\ni"description]&F\big(i\sr{\kappa_{ij}}{\ral}j\big)\ar[r,draw=none,"="description]&F_i\ar[d,"\tau_i"]\ar[r,"f_{ij}"]&F_j\ar[d,"\tau_j"]\\
i\ar[d,"id_i"]\ar[r,"\kappa_{ij}"]&j\ar[d,"id_j"]\ar[r,draw=none,"\in"description]&\Mor\I\ar[rr,"F'"]&\mathop{~}\limits^{\substack{~\\~\\~}}\ar[d,"\tau'=\{\tau'_i\}"]&\Mor\C\ar[r,draw=none,"\ni"description]&F'\big(i\sr{\kappa_{ij}}{\ral}j\big)\ar[r,draw=none,"="description]&F'_i\ar[d,"\tau'_i"]\ar[r,"f'_{ij}"]&F_j\ar[d,"\tau'_j"]\\
i\ar[r,"\kappa_{ij}"]&j\ar[r,draw=none,"\in"description]&\Mor\I\ar[rr,"F''"]&\mathop{~}\limits^{\substack{~\\~\\~}}&\Mor\C\ar[r,draw=none,"\ni"description]&F''\big(i\sr{\kappa_{ij}}{\ral}j\big)\ar[r,draw=none,"="description]&F'_i\ar[r,"f''_{ij}"]&F''_j\\
\et\ec
A \ul{contravariant functor} (\ul{cofunctor}) from $\C$ to $\D$ (written {\small $F:\C\ra\D$} or {\small $\C\sr{F}{\ral}\D$}) is a map of the form
\bit[leftmargin=0.9cm]
\item $F:\Ob\C\sqcup\Mor\C\ra \Ob\D\sqcup\Mor\D$ ~whose restrictions to objects and morphisms give parallel maps
{\small\[
F:\Ob\C\ra\Ob\D,~C\mapsto F(C)~~~~\txt{and}~~~~F:\Mor\C\ra\Mor\D,~\Big(C\sr{f}{\ral}C'\Big)~\mapsto~\Big(F(C)\sr{F(f)}{\lal}F(C')\Big),\nn
\]}
\eit
and the map on morphisms has the following properties
(\blue{footnote}\footnote{
As before, it is also enough to briefly specify the cofunctor as {\scriptsize $F:\C\ra\D,~\big(\dom f\sr{f}{\ral}\cod f\big)\mapsto\big(F(\dom f)\sr{f}{\lal}F(\cod f)\big)$}.
}):
For all $C\in\Ob\C$ and $f,f'\in\Mor\C$,
\bit[leftmargin=0.9cm]
\item $F(id_C)=id_{F(C)}$~~ and ~~$F(f\circ f')=F(f')\circ F(f)$.
\eit

{\flushleft A} \ul{contravariant endofunctor} is a contravariant functor $F:\C\ra\C$ from a category $\C$ to itself.

{\flushleft Let} $\I,\C$ be categories. Given cofunctors $F,F':\I\ra\C$, a \ul{morphism of cofunctors} (or \ul{conatural transformation}) $\tau:F\ra F'$ is a collection of morphisms in $\C$
{\small\bea
\tau=\Big\{F(i)\sr{\tau^i}{\ral}F'(i)\Big\}_{i\in\Ob\I}~\subset~\Mor\C=\Mor_\C(\Ob\C,\Ob\C)\nn
\eea} such that for any $i,j\in\I$ and any $\kappa_{ij}\in \Mor_\I(i,j)$, we have a commutative diagram
{\small\bea
\left(F\sr{\tau}{\ral}F'\right)\left(i\sr{\kappa_{ij}}{\ral}j\right)~~:=~~
\bt
F^i\ar[d,"\tau^i"]\ar[from=rr,"f^{ij}"']&&F^j\ar[d,"\tau^j"]\\
F'{}^i\ar[from=rr,"f'{}^{ij}"']&&F'{}^j
\et~~:=~~
\bt
F(i)\ar[d,"\tau^i"]\ar[from=rr,"F(\kappa_{ij})"']&&F(j)\ar[d,"\tau^j"]\\
F'(i)\ar[from=rr,"F'(\kappa_{ij})"']&&F'(j)
\et~~~~~~~~\tau^i\circ f^{ij}=f'{}^{ij}\circ \tau^j.\nn
\eea}
A morphism of cofunctors $\tau=\{\tau^i\}:F\ra F'$ is an \ul{isomorphism of cofunctors} (or \ul{conatural isomorphism}) if each $\tau^i:F^i\ra F'{}^i$ is an isomorphism in $\C$ (in which case, we say  $F,F'$ are \ul{isomorphic}, written $F\cong F'$).

{\flushleft The} \ul{opposite functor of  a category} $\C$ is the cofunctor $F_{op}$ from $\C$ to the opposite category $\C^{op}$ given by
\bea
F_{op}:\C\ra\C^{op},~~\big(C\sr{f}{\ral}C'\big)~\mapsto~\big(C\sr{f^{op}}{\lal}C'\big).\nn
\eea
Observe that any cofunctor $F:\I\ra\C$ can be uniquely expressed as a composition
{\small\bea
&&F=\wt{F}\circ F_{op}:\I\sr{F_{op}}{\ral}\I^{op}\sr{\wt{F}}{\ral}\C,~~\big(C\sr{f}{\ral}C'\big)~\sr{F_{op}}{\longmapsto}~
\big(C\sr{f^{op}}{\lal}C'\big)~\sr{\wt{F}}{\longmapsto}~
\big(F(C)\sr{F(f)}{\lal}F(C')\big)\nn\\
&&~~\txt{where}~~~~\wt{F}:\I^{op}\ra\C,~\big(C\sr{f}{\ral}C'\big)~\mapsto~\big(F(C)\sr{F(f^{op})}{\ral}F(C')\big)
\hspace{1.1cm}
\adjustbox{scale=0.7}{\bt
\I\ar[d,"F"]\ar[rr,"F_{op}"] && \I^{op}\ar[dll,dashed,"\wt{F}"]\\
\C\nn
\et}\nn
\eea}That is, a \ul{cofunctor} $F:\I\ra\C$ is  equivalently a functor $\wt{F}:\I^{op}\ra\C$ (or a functor $\wt{F}:\I\ra\C^{op}$).

The \ul{category of cofunctors} from $\I$ to $\C$ is thus the category of functors from $\I^{op}$ to $\C$, i.e., the category
\[
\textstyle \C^{\I^{op}}=\big(\C^{\I^{op}},\Mor_{\C^{\I^{op}}},\circ\big)=\big(\Ob \C^{\I^{op}},\Mor \C^{\I^{op}},\circ\big)=\big(\Ob\C^{\I^{op}},\Mor_{\C^{\I^{op}}}(\Ob\C^{\I^{op}},\Ob\C^{\I^{op}}),\circ\big).
\]
\end{dfn}

\begin{notation*}[\textcolor{blue}{Notation}]
It will often be convenient/sufficient to specify a functor $F$ from $\C$ to $\D$ in the form
{\small\[
F:\C\ra\D,~~~~\big(C\sr{f}{\ral}C'\big)~\mapsto~\big(C_\ast\sr{f_\ast}{\ral}C'_\ast\big)~:=~\Big(F(C)\sr{F(f)}{\ral}F(C')\Big).\nn
\]}Similarly, it will often be convenient/sufficient to specify a cofunctor $F$ from $\C$ to $\D$ in the form
{\small\[
F:\C\ra\D,~~~~\big(C\sr{f}{\ral}C'\big)~\mapsto~ \big(C'{}^\ast\sr{f^\ast}{\ral}C^\ast\big)~:=~\Big(F(C')\sr{F(f)}{\ral}F(C)\Big).\nn
\]}
\end{notation*}

\begin{rmk}[\textcolor{blue}{
\index{Faithfulness of a subcategory}{{Faithfulness of a subcategory}},
\index{Full subcategory}{Full subcategory},
\index{Dense! subcategory}{Dense subcategory}}]
Let $\C\subset\D$ be a subcategory. As the image of the inclusion map $\C\hookrightarrow\D$ (being a faithful functor) the subcategory $\C\subset\D$ is \ul{always faithful}. The subcategory $\C\subset\D$ is said to be \ul{full} (resp. \ul{dense}) if the inclusion map $\C\hookrightarrow\D$ is a full functor (resp. dense functor).
\end{rmk}

\begin{rmk}[\textcolor{blue}{``Functors'' include ``Cofunctors''}]
Because any given cofunctor $F:\C\ra\D$ can be understood in terms of a unique functor $\wt{F}:\C^{op}\ra\D$ through a composition $F=\wt{F}\circ F_{op}:\C\sr{F_{op}}{\ral}\C^{op}\sr{\wt{F}}{\ral}\D$, all concepts introduced for functors automatically carry over to cofunctors. Therefore, we will often (but not always) silently introduce concepts for functors (without explicit mention of the same concepts for cofunctors) but apply them to cofunctors (with notation and conventions \ul{unambiguously/canonically} adjusted accordingly for consistency).
\end{rmk}

\section{Categories of categories, Equivalences in a category}

\begin{rmk}[\textcolor{blue}{Vagueness of ``equivalence of categories''}]
As explained in Definition \ref{CatDef4} below, weak isomorphism (a very special type of equivalence by Definition \ref{CatDef2}) of categories is vaguely known as ``equivalence of categories''. However, it is apparent from Definition \ref{CatDef2} that there unavoidably exist many different notions of equivalences of categories, some of which may be known by different names (e.g., homotopy equivalence, to come later, of ``objects'' that are in a sense ``categories'').
\end{rmk}

\begin{dfn}[\textcolor{blue}{
\index{Category of! categories}{{Category of categories}},
\index{Weak! isomorphism (``Equivalence'') of categories}{Weak isomorphism (``Equivalence'') of categories},
\index{Weakly isomorphic (``Equivalent'') categories}{Weakly isomorphic (``Equivalent'') categories},
\index{(Strong) Isomorphism of categories}{(Strong) Isomorphism of categories},
\index{(Strongly) Isomorphic categories}{(Strongly) Isomorphic categories},
\index{Endofunctor}{Endofunctor},
\index{Autofunctor}{Autofunctor}
}]\label{CatDef4}
A category
\bea
\C\C=\big(\C\C,\Mor_{\C\C},\circ\big)=\big(\Ob\C\C,\Mor\C\C,\circ\big)=\big(\Ob\C\C,\Mor_{\C\C}(\Ob\C\C,\Ob\C\C),\circ\big)\nn
\eea
is a \ul{category of categories} if (i) its objects {\footnotesize $\C\in\Ob\C\C$} are categories, (ii) its morphisms {\footnotesize $F\in \bigsqcup\limits_{\C,\D\in\Ob\C\C}\Mor_{\C\C}(\C,\D)$} are functors, and (iii) the composition of two functors
{\footnotesize\[
F_1:\C_1\ra\C_2,~~\big(C_1\sr{f_1}{\ral}C_1'\big)~\mapsto~\Big(F_1(C_1)\sr{F_1(f_1)}{\ral}F_1(C_1')\Big)~~~~\txt{and}~~~~F_2:\C_2\ra\C_3,~~\big(C_2\sr{f_2}{\ral}C_2'\big)~\mapsto~\Big(F_2(C_2)\sr{F_2(f_2)}{\ral}F_2(C_2')\Big)\nn
\]}
is their composition as maps, which is explicitly given by
{\footnotesize\[
F_2\circ F_1:\C_1\sr{F_1}{\ral}\C_2\sr{F_2}{\ral}\C_3,~~\Big(C_1\sr{f_1}{\ral}C_1'\Big)\sr{F_1}{\longmapsto}
\Big(F_1(C_1)\sr{F_1(f_1)}{\ral}F_1(C_1')\Big)\sr{F_2}{\longmapsto}
\Big(F_2(F_1(C_1))\sr{F_2(F_1(f_1))}{\ral}F_2(F_1(C_1'))\Big).\nn
\]}
{\flushleft \ul{Weak isomorphism (or ``equivalence'') of categories}} in a category of categories $\C\C$ is the equivalence of objects in $\Ob\C\C$ induced by isomorphism of functors in $\Mor\C\C$. That is, a functor $F:\C\ra\D$ in $\Mor\C\C$ is a \ul{weak isomorphism (``equivalence'') of categories} (making $\C,\D$ \ul{weakly isomorphic (or ``equivalent'') categories}, written $\C\cong_w\D$) if there exists a functor $G:\D\ra\C$ (a pseudo-inverse of $F$ in $\C\C$) such that $G\circ F\cong id_\C:\C\ra\C$ and $F\circ G\cong id_\D:\D\ra\D$.

{\flushleft A} functor $F:\C\ra\D$ is a \ul{(strong) isomorphism of categories} (making $\C,\D$ \ul{(strongly) isomorphic categories}, written $\C\cong\D$) if there exists a functor $G:\D\ra\C$ (an inverse of $F$ in $\C\C$) such that $G\circ F=id_\C:\C\ra\C$ and $F\circ G=id_\D:\D\ra\D$.

Recall that an \ul{endofunctor} is a functor $F\in \Mor_{\C\C}(\C,\C)$ from a category $\C$ to itself. An isomorphism of categories $F\in \Mor_{\C\C}(\C,\C)$ is called an \ul{autofunctor} of $\C$.
\end{dfn}

\begin{thm}[\textcolor{blue}{Criterion for weak iso (``equivalence'') of categories: \cite[Theorem 7, p.71]{gelfand-manin2010}}]
A functor $F:\C\ra\D$ is a weak isomorphism (an ``equivalence'') of categories $\iff$ it is faithful, full, and dense.
\end{thm}
\begin{proof}
($\Ra$): Assume $F:\C\ra\D$ is an equivalence. Then there is $G:\D\ra\C$ such that $G\circ F\cong id_{\C}$, $F\circ G\cong id_{\D}$, i.e., we have isomorphisms of functors $\al:G\circ F\ra id_\C$,~ $\beta:F\circ G\ra id_{\D}$, defined by the following commutative diagram (in which $f\al_C=\al_{C'}GF(f)$ and $g\beta_{{D}}=\beta_{{D'}}FG(g)$):
{\small
\begin{figure}[H]
\centering
\adjustbox{scale=0.9}{
\begin{tikzcd}
 \ar[r] & GF(C)\ar[d,"\al_C"]\ar[rr,"GF(f)"] && GF(C')\ar[d,"\al_{C'}"]\ar[r] &~\\
 \ar[r] & C\ar[rr,"f"] && C'\ar[r] &~
\end{tikzcd}
\hspace{1cm}
\begin{tikzcd}
 \ar[r] & FG(D)\ar[d,"\beta_D"]\ar[rr,"FG(g)"] && FG(D')\ar[d,"\beta_{D'}"]\ar[r] &~\\
 \ar[r] & D\ar[rr,"g"] && D'\ar[r] &~
\end{tikzcd}}
\caption{}\label{dg1.4ss}
\end{figure}}

\bit[leftmargin=0.9cm]
\item \ul{$F$ is faithful}: For $f_1,f_2\in \Mor_\C(C,{C'})$, let $F(f_1)=F(f_2)$ in $\Mor_{\D}(F(C),F({C'}))$. Then
\bea
\al_{C'}^{-1}f_1\al_C=GF(f_1)=GF(f_2)=\al_{C'}^{-1}f_2\al_C,~~\Ra~~f_1=f_2.\nn
\eea
(Similarly, $G$ is faithful)
\item \ul{$F$ is full}: Let $g\in \Mor_{\D}(F(C),F({C'}))$. Then $G(g)\in \Mor_\C\big(GF(C),GF({C'})\big)$, which implies~ $f:=\al_{C'}G(g)\al_C^{-1}\in \Mor_\C(C,{C'})$. Therefore,
    \bea
    &&f\al_C=\al_{C'}GF(f)~~\Ra~~[\al_{C'}G(g)\al_C^{-1}]\al_C=\al_{C'}GF(f),~~\Ra~~G(g)=GF(f),\nn\\
    &&~~\Ra~~g=F(f)=F\big(\al_{C'}G(g)\al_C^{-1}\big),~~~~\txt{since $G$ is faithful}.\nn
    \eea
\item \ul{$F$ is dense}: Let ${D}\in Ob~\D$. Then it is clear that ${D}\cong FG({D})=F(C)$, where $C:=G({D})\in Ob~\C$.
\eit
{\flushleft($\La$):} Assume $F$ is faithful, full, and dense, i.e., $F:\Mor_\C(C,{C'})\ra \Mor_{\D}(F(C),F({C'}))$ is bijective, and for any ${D}\in Ob~\D$, we have an isomorphism $\phi_{{D}}:{D}\ra F(C_{{D}})$ for some $C_{{D}}\in Ob~\C$.

    We need to find a functor $G:\D\ra \C$ such that $GF\cong id_\C$, $FG\cong id_{\D}$, i.e., for all $f\in \Mor_\C(C,{C'})$, $g\in \Mor_{\D}({D},{D'})$ we need to find isomorphisms and relations given (as in Figure \ref{dg1.4ss}) by
    \bea
   \label{zzeq1} \al_C:GF(C)\ra C,~~\beta_{{D}}:FG({D})\ra {D},~~~~f\al_C=\al_{C'}GF(f),~~~~g\beta_{{D}}=\beta_{{D'}}FG(g).
    \eea
Consider the above given information (on the denseness and bijectivity $F$) in diagram form as follows:
\begin{figure}[H]
\centering
\adjustbox{scale=0.9}{
\begin{tikzcd}
 \ar[r] & F(C_D)\ar[rr,dashed,"FG(g)"] && F(C_{D'})\ar[r] &~\\
 \ar[r] & D\ar[u,"\phi_D"]\ar[rr,"g"] && D'\ar[u,"\phi_{D'}"]\ar[r] &~\\
  & & (a) & &
\end{tikzcd}
\hspace{1cm}
\begin{tikzcd}
 \ar[r] & C_{F(C)}\ar[rr,dashed,"GF(f)"] && C_{F(C')}\ar[r] &~\\
 \ar[r] & C\ar[u,"\psi_C",":=F^{-1}(\phi_{F(C)})"']\ar[rr,"f"] && C'\ar[u,"\psi_{C'}",":=F^{-1}(\phi_{F(C')})"']\ar[r] &~\\
 & & (b) & &
\end{tikzcd}}
\caption{}\label{dg2.4ss}
\end{figure}
\bit[leftmargin=0.9cm]
\item Based on Fig \ref{dg2.4ss} (a), define a functor $G:\D\ra\C$ by
\bea
&&G:\Big({D}\sr{g}{\ral}{D'}\Big)\mapsto \Big(G(D)\sr{G(g)}{\ral}G(D')\Big):= \Big(C_{{D}}\sr{G(g)}{\ral}C_{{D'}}\Big),~~~~G(g):=F^{-1}\big(\phi_{{D'}}g\phi_{{D}}{}^{-1}\big),\nn
\eea
which implies
\bea
\label{zzeq2}FG({D})=F(C_{{D}})~\cong~{D},~~~~\phi_{{D'}}g=FG(g)\phi_{{D}}~~~~\txt{for all}~~~~{D},{D'}\in~Ob~\D.
\eea
\item Next, for any $C\in Ob~\C$, consider an isomorphism
\bea
\phi_{F(C)}:F(C)\ra F(C_{F(C)})=F\big(GF(C)\big).\nn
\eea
Since ~$F:\Mor_\C\big(C,C_{F(C)}\big)\ra \Mor_{\D}\big(F(C),F(C_{F(C)})\big)$~ is bijective, there is a unique isomorphism ~$\psi_C:C\ra C_{F(C)}=GF(C)$~ such that ~$\phi_{F(C)}=F(\psi_C)$,~ and so we also have
\bea
\label{zzeq3}~~~~GF(C)~\cong~C~~~~\txt{for all}~~~~C\in Ob~\C.
\eea
Moreover, for any $f\in \Mor_\C(C,{C'})$, we have $g:=F(f)\in \Mor_{\D}(F(C),F({C'}))$, and so with ${D}:=F(C)$ and ${D'}:=F({C'})$ in (\ref{zzeq2}),
{\small
\bea
&&\phi_{F({C'})}g=FG(g)\phi_{F(C)}~~\Ra~~F(\psi_{C'})F(f)=FG(F(f))F(\psi_C),~~\Ra~~F(\psi_{C'}f)=F\big(GF(f)\psi_C\big)\nn\\
\label{zzeq4}&&~~\Ra~~\psi_{C'}f=GF(f)\psi_C,~~~~\txt{for all}~~~~C,{C'}\in Ob~\C.
\eea}
\item Hence, with $\al_C:=\psi_C^{-1}$, $\beta_{{D}}:=\phi_{{D}}^{-1}$, we get isomorphisms of functors $\al:GF\ra id_\C$, $\beta:FG\ra id_{\D}$. \qedhere
\eit
\end{proof}

\begin{dfn}[\textcolor{blue}{\index{Morita equivalence of rings}{Morita equivalence of rings}}]
Two rings $R,S$ are \ul{Morita-equivalent} if there exists an equivalence of categories $F:\txt{$R$-mod}\ra\txt{$S$-mod}$ (i.e., for all $R$-modules $M,M'$ the restriction $F:\Mor_R(M,M')\ra \Mor_S(F(M),F(M'))$ is bijective, and for every $S$-module $N$ there is an $R$-module $M$ such that $N\cong F(M)$ ).
\end{dfn}

\begin{dfn}[\textcolor{blue}{
\index{Imbedding of a category}{Imbedding of a category},
\index{Imbedded (sub)category}{Imbedded (sub)category}}]
A functor $F:\C\ra\D$ is called an \ul{imbedding} (making $F(\C)\subset\D$ an \ul{imbedded (sub)category}) if it is injective as a map. That is, $F$ is (i) injective on objects, and also (ii) injective on morphisms (i.e., faithful and injective on objects).
\end{dfn}

\begin{dfn}[\textcolor{blue}{
\index{Higher category}{Higher category},
\index{Higher morphism}{Higher morphism},
\index{Weak! higher category}{Weak higher category}}]
Given $n\geq 1$, an \ul{$n$-category}
\bea
\C^n=(\Ob\C^n,\Mor\C^n,\circ^n)=\Big(\Ob\C^n,(\Mor\C_1,...,\Mor\C_n),(\circ_1,...,\circ_n)\Big)\nn
\eea
is a category in the form of a hierarchical sequence of categories $\C_1,\cdots,\C_n$ defined inductively as follows:
\bit
\item[(i)] There is a starting category ~$\C_1:=(\Ob\C^n,\Mor\C_1,\circ_1)$, whose objects are the objects of $\C^n$. The morphisms of $\C_1$ are called \ul{$1$-morphisms} of $\C^n$.
\item[(ii)] For any $2\leq k\leq n$, we have a (diagram in a) category ~$\C_k:=(\Mor\C_{k-1},\Mor\C_k,\circ_k)$, whose objects are the morphisms of $\C_{k-1}$. The morphisms of $\C_k$ are called \ul{$k$-morphisms} of $\C^n$.
\eit

If the associativity and identity relations (as given in Definition \ref{CatDef1}) only hold up to isomorphism (i.e., if the equalities $=$ are replaced with isomorphisms $\cong$), the resulting $n$-category is called a \ul{weak $n$-category}.
\end{dfn}

For example, a diagram in some $3$-category $\C^3$ might have the following form:
\bea\bt
A\ar[rr,"f"] &\ar[dd,Rightarrow,"F"']& A'\ar[rr,"f'"]  && C\ar[rr,"h"] &\ar[dd,Rightarrow,"G"]\ar[dr,Rightarrow,"H"']& C'\ar[dd,"h'"] \\
     &\ar[rrrr,Rightarrow,dashed,shift right=2,"\eta"]&    &&    &~& ~\ar[dl,Rightarrow,"L"] \\
B\ar[rr,"g"'] &~& B'  && D\ar[rr,"l"'] &~& D' \\
\et\nn
\eea
where $A,A',B,B',C,C',D,D'$ are objects of $\C^3$, $f,f',g,h,h',l$ are $1$-morphisms, $F,G,H,L$ are $2$-morphisms, and $\eta$ is a $3$-morphism.

\begin{note}[\textcolor{blue}{Interpretation of a higher category}]
Recall that a category of categories is a category whose objects are themselves categories. A higher category on the other hand is a category in which certain morphisms are ``upgraded'' to play a double role in the sense they are also viewed as ``objects'' (in addition to being morphisms) so that they also have morphisms between them, and this upgrading can be repeated on subsequent morphisms between morphisms to reach any desired number of upgrading levels.
\end{note}

For a further reading on higher categories, see for example \cite{lurie2009}.

\section{Quotients, Products, and Coproducts of Categories}
\begin{dfn}[\textcolor{blue}{\index{Multiplicative! class}{Multiplicative class},
\index{Congruence relation on a class}{Congruence relation on a class},
\index{Quotient! multiplicative class}{Quotient multiplicative class},
\index{Congruence relation on a category}{Congruence relation on a category},
\index{Quotient! category}{{Quotient category}},
\index{Congruence functor}{Congruence functor}
}]
Let $M=[M,\cdot]$ be a \ul{multiplicative class} (i.e., a class $M$ with a product $\cdot:M\times M\ra M,~(a,b)\mapsto ab$). A \ul{congruence relation} on $M$ is an equivalence relation $\sim$ on $M$ such that for all $a,a',b,b'\in M$, if $a\sim a'$ and $b\sim b'$ then $ab\sim a'b'$. That is, if $[a]=[a']$ and $[b]=[b']$ then $[ab]=[a'b']$, where $[m]:=\{m'\in M:m'\sim m\}$ is the equivalence class of $m\in M$.

Given a congruence relation $\sim$ on $M$, the collection of equivalence classes ${M\over\sim}:=\{[a]:a\in M\}$ is a multiplicative class (called \ul{quotient multiplicative class}) with binary operation ${M\over\sim}\times{M\over\sim}\ra{M\over\sim}$,
\bea
([a],[b])\mapsto [a][b]:=[ab],\nn
\eea
which is well-defined because if $[a]=[a']$ and $[b]=[b']$, then ~$[a][b]=[ab]\sr{\txt{congruence}}{=}[a'b']=[a'][b']$.

Let $\C=\big(\C,\Mor_\C,\circ\big)=\big(\Ob\C,\Mor\C,\circ\big)=\big(\Ob\C,\Mor_\C(\Ob\C,\Ob\C),\circ\big)$ be a category. A \ul{congruence relation on $\C$} is an equivalence relation $\sim$ on $\Mor\C$ (through an equivalence relation $\sim_{AB}$ on $\Mor_\C(A,B)$ for each pair of objects $A,B\in\Ob\C$) such that for all $f,f',g,g'\in\Mor\C$, if $f\sim f'$ and $g\sim g'$ then $f\circ g\sim f'\circ g'$.

Given a congruence relation $\sim$ on $\C$, the \ul{quotient category} of $\C$ associated with $\sim$ is the category ${\C\over\sim}=\left({\C\over\sim},\Mor_{\C\over\sim},\circ\right)=\left(\Ob{\C\over\sim},\Mor{\C\over\sim},\circ\right)=\left(\Ob{\C\over\sim},\Mor_{\C\over\sim}(\Ob{\C\over\sim},\Ob{\C\over\sim}),\circ\right)$  given by the following:
\bit[leftmargin=0.9cm]
\item[(1)] $\Ob{\C\over\sim}:=\Ob\C$ (i.e., the same objects as $\C$).
\item[(2)] $\Mor_{\C\over\sim}(A,B):={\Mor_\C(A,B)\over\sim}:={\Mor_\C(A,B)\over\sim_{AB}}$.
\item[(3)] The composition of $[g]\in \Mor_{\C\over\sim}(A,B)$ and $[f]\in \Mor_{\C\over\sim}(B,C)$ is given by
\bea
[f]\circ[g]:=[f\circ g]\in \Mor_{\C\over\sim}(A,C),\nn
\eea
while the identity morphisms are the same since the objects have not changed.
\eit
Given a congruence relation $\sim$ on $\C$, the \ul{congruence functor} associated with $\sim$ is the obvious functor
\[
\textstyle \Pi:\C\ra{\C\over\sim},~\big(C\sr{f}{\ral}C'\big)\mapsto\big(C\sr{[f]}{\ral}C'\big).\nn
\]
\end{dfn}

\begin{dfn}[\textcolor{blue}{
\index{Product of! categories}{Product of categories},
\index{Coproduct! of categories}{Coproduct of categories}}]
As done with previously seen structures (multiplicative sets, groups, rings, modules, algebras), if $(\C_i)_{i\in I}$ is an indexing of categories, then their \ul{product} $\prod_{i\in I}\C_i$ is the smallest category (up to isomorphism of categories) having each member $\C_i$ as a quotient category (or alternatively, the smallest category for which there exists a mapwise-surjective functor $F_i:\prod_{i\in I}\C_i\ra\C_i$ for each $i\in I$). The following is an explicit construction of the said product.

Given categories $\C,\D$, their \ul{product category} is the category $\C\times\D=(\C\times\D,\Mor_{\C\times\D},\circ)=(\Ob(\C\times\D),\Mor(\C\times\D),\circ)=\Big(\Ob(\C\times\D),\Mor_{\C\times\D}\big(\Ob(\C\times\D),\Ob(\C\times\D)\big),\circ\Big)$ such that the following hold:
\bit[leftmargin=0.9cm]
\item[(1)] \ul{Objects} are cartesian (ordered) pairs of $\C$-objects with $\D$-objects:
\bea
\Ob(\C\times\D):=\Ob\C\times\Ob\D=\{(C,D):C\in\Ob\C,D\in\Ob\D\}.\nn
\eea
\item[(2)] \ul{Morphisms} are cartesian (ordered) pairs of $\C$-morphisms with $\D$-morphisms:
\bea
&&\textstyle \Mor_{\C\times\D}\big(\Ob(\C\times\D),\Ob(\C\times\D)\big):=\Mor_\C(\Ob\C,\Ob\C)\times \Mor_{\D}(\Ob\D,\Ob\D)\nn\\
&&\textstyle~~~~=\big\{(f,g):f\in \Mor_\C(\Ob\C,\Ob\C),g\in \Mor_\D(\Ob\D,\Ob\D)\big\}\nn\\
&&~~~~=\left\{(C,D)\sr{(f,g)}{\ral}(C',D')=\left(C\sr{f}{\ral}C',D\sr{g}{\ral}D'\right)~\Big|~\substack{f\in \Mor_\C(C,C'),~g\in \Mor_\D(D,D')\\ C,C'\in\Ob\C,~D,D'\in\Ob\D~~~~~~}\right\}\nn\\
&&\textstyle~~~~=\bigsqcup\limits_{(C,D),(C',D')\in\Ob(\C\times\D)}\Mor_\C(C,C')\times \Mor_\D(D,D').\nn
\eea
where the \ul{composition of morphisms} and the \ul{identity morphisms} are given respectively by
\bea
(f,g)\circ(f',g'):=(f\circ f',g\circ g'),~~~~id_{(C,D)}:=(id_C,id_D),~~\txt{for all}~~(C,D)\in\Ob(\C\times\D).\nn
\eea
\eit

{\flushleft It} is clear that given an indexing (indexed or cartesian collection) of categories $(\C_i)_{i\in I}$, we can generalize the above procedure to obtain a \ul{product category} denoted by $\prod(\C_i)_{i\in I}$ or $\prod_{i\in I}\C_i$, as the category
{\footnotesize\begin{align}
\textstyle\prod_i\C_i&=\textstyle\Big(\prod_i\C_i,\Mor_{\prod_i\C_i},\circ\Big)=\Big(\Ob(\prod_i\C_i),\Mor(\prod_i\C_i),\circ\Big)\nn\\
&=\textstyle\Big(\Ob(\prod_i\C_i),\Mor_{\prod_i\C_i}\big(\Ob(\prod_i\C_i),\Ob(\prod_i\C_i)\big),\circ\Big)\nn
\end{align}}such that the following hold:
\bit[leftmargin=0.9cm]
\item[(1)] \ul{Objects} are indexings of objects from the respective $\C_i$-objects:
\[
\textstyle\Ob\left(\prod_i\C_i\right):=\prod_i\Ob\C_i=\left\{\txt{maps}~C_I:I\ra\bigsqcup_i\Ob\C_i,~i\mapsto C_i\in\Ob\C_i\right\}=\left\{C_I=(C_i)_{i\in I}:C_i\in\Ob\C_i\right\}.
\]
\item[(2)] \ul{Morphisms} are indexings of morphisms from the respective $\C_i$-morphisms:
\bea
&&\textstyle \Mor_{\prod_i\C_i}\Big(\Ob(\prod_i\C_i),\Ob(\prod_i\C_i)\Big):=\prod_i \Mor_\C(\Ob\C_i,\Ob\C_i)\nn\\
&&\textstyle~~~~=\left\{\txt{maps}~f_I:I\ra\bigcup_i\Mor_\C(\Ob\C_i,\Ob\C_i),~i\mapsto f_i\in \Mor_\C(\Ob\C_i,\Ob\C_i)\right\}\nn\\
&&\textstyle~~~~=\big\{(f_i)_{i\in I}:f_i\in \Mor_{\C_i}(\Ob\C_i,\Ob\C_i)\big\}
=\Big\{(C_i)_{i\in I}\sr{(f)_{i\in I}}{\ral}(C'_i)_{i\in I}=\big(C_i\sr{f_i}{\ral}C'_i\big)_{i\in I}\Big\}_{\substack{f_i\in \Mor_{\C_i}(C_i,C'_i)\\ C_i,C'_i\in\Ob\C_i~~~~~}}\nn\\
&&\textstyle~~~~=\prod_i\bigsqcup\limits_{C_i,C'_i\in\Ob\C_i}\Mor_{\C_i}(C_i,C'_i)=\bigsqcup\limits_{C_I,C'_I\in\Ob(\prod_i\C_i)}\prod_i \Mor_{\C_i}(C_i,C'_i),\nn
\eea
where the \ul{composition of morphisms} and the \ul{identity morphisms} are given respectively by
\bea
\textstyle (f_i)_{i\in I}\circ(f'_i)_{i\in I}:=(f_i\circ f'_i)_{i\in I},~~~~id_{(C_i)_{i\in I}}:=\left(id_{C_i}\right)_{i\in I},~~\txt{for all}~~(C_i)_{i\in I}\in\Ob(\prod_i\C_i).\nn
\eea
\eit

The \ul{coproduct} $\coprod_{i\in I}\C_i$ of the categories $(\C_i)_{i\in I}$ is the category generated by $(\C_i)_{i\in I}$, in the sense it is the smallest category (up to isomorphism of categories) containing each $\C_i$ as a subcategory. Compare $\coprod_{i\in I}\C_i$ with the subcategory $\coprod'_{i\in I}\C_i$ of $\prod_{i\in I}\C_i$ with objects and morphisms defined as follows:
\bit[leftmargin=0.9cm]
\item[(1)] \ul{Objects} of $\coprod'_{i\in I}\C_i$ are those object-indexings $(C_i)_{i\in I}\in\Ob\prod_{i\in I}\C_i$ such that only a finite number of them are distinct (i.e. $C_i$ are distinct only for finitely many $i$).
\item[(2)] \ul{Morphisms} of $\coprod'_{i\in I}\C_i$ are those morphism-indexings $(f_i)_{i\in I}\in\Mor\prod_{i\in I}\C_i$ such that only a finite number of them are distinct (i.e. $f_i$ are distinct only for finitely many $i$).
\eit
\end{dfn}

\begin{dfn}[\textcolor{blue}{
\index{Multifunctor}{Multifunctor},
\index{Bifunctor}{Bifunctor},
\index{Symmetric! bifunctor}{Symmetric bifunctor},
\index{Associative! bifunctor}{Associative bifunctor},
\index{Associative! bifunctor system}{Associative bifunctor system}
}] Given a collection of categories $\C_1,...,\C_k,\D$, a \ul{$k$-multifunctor} is a functor of the form
\bea
F:\C_1\times\cdots\times\C_k\ra\D,~~\Big((C_1,...,C_k)\sr{(f_1,...,f_k)}{\ral}(C_1',...,C_k')\Big)\mapsto \Big(F(C_1,...,C_k)\sr{F(f_1,...,f_k)}{\ral}F(C_1',...,C_k')\Big),\nn
\eea
where the expression $(C_1,...,C_k)\sr{(f_1,...,f_k)}{\ral}(C_1',...,C_k')$ has the following expanded/componentwise form:
\bea
(C_1,...,C_k)\sr{(f_1,...,f_k)}{\ral}(C_1',...,C_k')~:=~\big(C_1\sr{f_1}{\ral}C_1',\cdots,C_k\sr{f_k}{\ral}C_k'\big).\nn
\eea
A \ul{bifunctor} is a $2$-multifunctor {\small $F:\C_1\times\C_2\ra\D$}. A bifunctor $F:\C\times\C\ra\D$ is \ul{symmetric} if $F(C,C')\cong F(C',C)$ in $\D$, for any $C,C'\in\Ob\C$. A bifunctor of the form {\small $F:\C\times\C\ra\C$} is \ul{associative} if
\bea
F\big(F(A,B),C\big)\cong F\big(A,F(B,C)\big)~~~~\txt{for all}~~~~A,B,C\in\Ob\C,\nn
\eea
i.e., the following diagram \ul{commutes up to isomorphism}:
\bea\bt
\C\times\C\times\C\ar[d,"id_\C\times F"]\ar[rr,"F\times id_\C"] && \C\times\C\ar[d,"F"]\\
\C\times\C\ar[rr,"F"] && \C
\et~~~~~~~~F\circ(F\times id_\C)\cong F\circ(id_\C\times F)\nn
\eea
which equivalently means the associated $3$-multifunctors $F_{(1,2)}(-,-,-):=F(-,F(-,-))=F\circ(id_\C\times F)$ and $F_{(2,1)}(-,-,-):=F(F(-,-),-)=F\circ(F\times id_\C)$ are isomorphic, i.e., $F_{(1,2)}\cong F_{(2,1)}$.

A system of bifunctors {\small $F_{ij}^k:\C_i\times\C_j\ra\C_k$} is an \ul{associative bifunctor system} if for all $i,j,k,r,s,t$,
{\small\[
F_{rk}^t\big(F_{ij}^r(C_i,C_j),C_k\big)\cong F_{is}^t\big(C_i,F_{jk}^s(C_j,C_k)\big)~~~\txt{for all}~~~C_i\in\Ob\C_i,~C_j\in\Ob\C_j,~C_k\in\Ob\C_k,\nn
\]}i.e., for any given $i,j,k,r,s,t$ the following diagram \ul{commutes up to isomorphism}:
\bea\bt
\C_i\times\C_j\times\C_k\ar[d,"id_{\C_i}\times F_{jk}^s"]\ar[rr,"F_{ij}^r\times id_{\C_k}"] && \C_r\times\C_k\ar[d,"F_{rk}^t"]\\
\C_i\times\C_s\ar[rr,"F_{is}^t"] && \C_t
\et~~~~~~~~F_{rk}^t\circ(F_{ij}^r\times id_{\C_k})\cong F_{is}^t\circ(id_{\C_i}\times F_{jk}^s).\nn
\eea
\end{dfn}

\begin{rmk}[\blue{
\index{Partial! multifunctor}{Partial multifunctor},
\index{Partial! bifunctor}{Partial bifunctor}
}]
In some applications, it is sufficient to consider a map of the form $F:\C_1\times\cdots\times\C_k\ra\D$ that is only ``partially functorial'' (making $F$ a \ul{partial $k$-multifunctor}) in the sense that each ``restriction'' $F_i:=F|_{C_i}$ (defined to be $F$ with all arguments, except the $i$th one, fixed) is a functor. A partial $2$-multifunctor would be a \ul{partial bifunctor}.

For brevity in our subsequent discussion, except in very specific/precise situations, we will not be making this distinction, and so use ``multifunctor'' (e.g., ``bifunctor'') where ``partial multifunctor'' (e.g., ``partial bifunctor'') is more useful/appropriate. Of course a multifunctor is a partial multifunctor but the reverse might not hold in general. Therefore, if ``multifunctor'' is used where ``partial multifunctor'' also works, then we might lose some generality but without a serious loss of correctness.
\end{rmk}

\begin{question}
Is a partial multifunctor equivalent (in an absolute/canonical way) to a multifunctor?
\end{question}

%% file: parts/AlgebraCat/AlgebraCatS3.tex
\chapter{Diagrams, Systems, Processes, Functoriality, and Limits}\label{AlgebraCatS3}

\section{Subobjects, Quotient objects, Restrictions}
\begin{dfn}[\textcolor{blue}{
\index{Subobject}{Subobject in a category},
\index{Containing object}{Containing object},
\index{Quotient! object}{Quotient object in a category},
\index{Numerator object}{Numerator object},
}]
Let $\C$ be a category and $A,B,C\in\Ob\C$. Then $A$ is a \ul{subobject} of $B$ (making $B$ a \ul{containing object} of $A$), written ~$A\subset B$~ or ~$A\subset_\C B$~ or ~$A\in Mono_\C(A,B)$,~ if there exists a monomorphism $m:A\hookrightarrow B$. We say $C$ is a \ul{quotient} object of $B$ (making $B$ a \ul{numerator object} of $C$), written ~$C\subset_{op}B$~ or ~$C\subset_{\C^{op}} B$~ or ~$C\in Epi_\C(B,C)=Mono_{\C^{op}}(C,B)$,~ if there exists an epimorphism $e:B\twoheadrightarrow C$. (\blue{footnote}\footnote{The notation here is based on (i) the fact that the categories $\C,\C^{op}$ have the same objects and (ii) the observation that a morphism $f\in\Mor\C$ is a mono (resp. an epi) in $\C$ $\iff$ $f^{op}\in\Mor\C^{op}$ is an epi (resp. a mono) in $\C^{op}$.})
\end{dfn}
The notions here of a \emph{subobject} and a \emph{quotient object} are categorical generalizations of the same notions for the objects seen earlier (groups, rings, modules, algebras). Since an object is only unique up to isomorphism, so are subobjects and quotient objects (viewed as objects in the earlier seen category $Cat(\Mor\C)$ whose objects are the morphisms of $\C$).
\begin{rmk}
If $A\subset B$ and $B\subset A$ it does \ul{not} necessarily follow that $A\cong B$. (See for example \cite{chih2014}.)
\end{rmk}

\begin{question}
Let $\C$ be a category. Suppose an imbedding $F:\C\ra Sets$ exists, and $F$ maps monomorphisms to monomorphisms. If two objects $C,C'\in\Ob\C$ are such that $C\subset C'$ and $C'\subset C$, does it follow that $C\cong C'$?
\end{question}

\begin{rmk}[\textcolor{blue}{Subobjects and quotient objects in $R$-mod}]~
If $A,B,C$ are $R$-modules, then it is clear that (i) $A$ is a subobject of $B$ $\iff$ $A\cong A'$ for a submodule $A'\subset B$, and (ii) $C$ is a quotient object of $B$ $\iff$ $C\cong C'$  for a quotient module $C':={B\over A'}$ of $B$ (i.e., for a submodule $A'\subset B$).
\end{rmk}

\begin{dfn}[\textcolor{blue}{
\index{Submorphism}{Submorphism},
\index{Quotient! morphism}{Quotient morphism},
\index{Restriction! of a morphism}{Restriction of a morphism}
}]
Let $\C$ be a category and $f,g$ morphisms in $\C$. Then $f$ is a \ul{submorphism} (resp. \ul{quotient morphism}) of $g$, written $f\subset g$ (resp. $f\subset_{op}g$), if $f$ is a subobject (resp. quotient object) of $g$ in the category $Cat(\Mor\C)$.

$f$ is a \ul{restriction} of $g$, written $f=g|_{\dom f,\cod f}$, if (i) there exist monomorphisms $u:\dom f\hookrightarrow \dom g$ and $v:\cod f\hookrightarrow \cod g$ (i.e., $\dom f\subset \dom g$ and $\cod f\subset \cod g$) satisfying (ii) $v\circ f=g\circ u$, i.e., the following diagram commutes.
\bea\adjustbox{scale=0.8}{\bt
\dom f\ar[r,hook,"u"]\ar[d,"f"] & \dom g\ar[d,"g"]\\
\cod f\ar[r,hook,"v"] & \cod g
\et}\nn\eea

\end{dfn}

\begin{question}
Let $\C$ be a category and $f,g\in\Mor\C$. Are the statements  ``$f$ is a submorphism of $g$'' and ``$f$ is a restriction of $g$'' equivalent?
\end{question}

\begin{dfn}[\textcolor{blue}{
\index{Subfunctor}{Subfunctor}, 
\index{Quotient! functor}{Quotient functor}}]
Let $\C,\D$ be categories and $F,G:\C\ra\D$ functors. Then $F$ is a \ul{subfunctor} (resp. \ul{quotient functor}) of $G$, written $F\subset G$ (resp. $F\subset_{op}G$), if $F\subset G$ as a subobject (resp. $F\subset_{op}G$ is a quotient object) in the category of functors $\D^\C$, in the sense there exists a monic (resp. epic) morphism of functors $\eta:F\hookrightarrow G$ (resp. $\delta:G\twoheadrightarrow F$), i.e., for any morphism $C\sr{f}{\ral}C'$ in $\C$, we have the following commutative diagram in $\D$.
\[\adjustbox{scale=0.8}{\bt
C\ar[d,"f"] & F(C)\ar[r,hook,"\eta_C"]\ar[d,"{F(f)}"] & G(C)\ar[d,"{G(f)}"]\\
C' & F(C')\ar[r,hook,"\eta_{C'}"] & G(C')
\et}~~~\left(~\txt{resp.}~~
\adjustbox{scale=0.8}{\bt
C\ar[d,"f"] & G(C)\ar[r,two heads,"\delta_C"]\ar[d,"{G(f)}"] & F(C)\ar[d,"{F(f)}"]\\
C'          & G(C')\ar[r,two heads,"\delta_{C'}"]            & F(C')
\et}\right)
\]

Thus, equivalently, ~$F\subset G$~ if ~$F(f)=G(f)|_{\dom F(f),\cod F(f)}$~ for all ~$f\in\Mor\C$. (\blue{footnote}\footnote{For any subfunctor ~$F\subset G:\C\ra\D$,~ it is clear that $F(C)\subset G(C)$ in $\D$ for all $C\in\Ob\C$.})
\end{dfn}

\section{Diagrams, Systems, Morphisms of systems, Categories of systems}
\begin{dfn}[\textcolor{blue}{
\index{Diagram in a category}{Diagram in a category},
\index{Functoral diagram (Functor)}{Functoral diagram (Functor)},
\index{Cofunctoral diagram (Cofunctor)}{Cofunctoral diagram (Cofunctor)},
\index{Diagram as a graph}{{Diagram as a graph}},
\index{Image of! a diagram}{Image of a diagram},
\index{Simple! diagram}{Simple diagram},
\index{Commutative! diagram}{Commutative diagram},
\index{Simple! category}{Simple category},
\index{Commutative! category}{Commutative category},
\index{System! in a category}{{System in a category}},
\index{System! objects}{System objects},
\index{System! transition morphisms}{System transition morphisms},
\index{Cosystem in a category}{{Cosystem in a category}},
\index{Trivial system}{Trivial system},
\index{Isosystem}{Isosystem},
\index{Morphism of! systems}{Morphism of systems},
\index{Morphism of! cosystems}{Morphism of cosystems},
\index{Category of! systems}{{Category of systems}},
\index{Category of! cosystems}{Category of cosystems},
\index{Multisystem}{Multisystem},
\index{Bisystem}{Bisystem}}]\label{CatDiag}
Let $\I,\C$ be categories (each of which is by definition a directed graph with objects as vertices and morphisms as edges). An \ul{$\I$-diagram} in $\C$ (written $D:\I\dashrightarrow\C$) is a map of the form $D:\Ob\I\sqcup\Mor\I\ra\Ob\C\sqcup\Mor\C$ such that its restrictions on objects and morphisms give separate maps
\[
D:\Ob\I\ra\Ob\C,~i\mapsto D(i)~~~~\txt{and}~~~~D:\Mor\I\ra\Mor\C,~\big(i\sr{\kappa}{\ral}j\big)\mapsto\big(D(i)\sr{D(\kappa)}{\ral}D(j)\big).
\]
The diagram $D:\I\dashrightarrow\C$ is \ul{functorial or natural}, making it a functor or natural transformation $D:\I\ra\C$, (resp. \ul{cofunctorial or conatural}, making it a cofunctor or conatural transformation,) if it satisfies the axioms of a functor (resp. cofunctor) as follows: For all $i\in\Ob\I$ and $\kappa,\kappa'\in\Mor\I$,
{\flushleft (i)} $D(id_i)=id_{D(i)}$~ and ~(ii) $D(\kappa\circ\kappa')=D(\kappa)\circ D(\kappa')$ ~$\Big($resp. $D(\kappa\circ\kappa')=D(\kappa')\circ D(\kappa)$ $\Big)$.

The \ul{diagram} $D$ is equivalently the directed graph $D=(V,E)$ with vertices {\small $V:=\big(D(i)\big)_{i\in\Ob\I}$} and edges {\small $E:=\big(D(\kappa)\big)_{\kappa\in\Mor\I}$} (which must be viewed as indexings since they fully represent $D$ as a map). The following are worth noting:
\begin{enumerate}[leftmargin=0.9cm]
\item In typical applications we will encounter, (i) $\I$ will often be a poset viewed as a category and (ii) the diagram $D:\I\dashrightarrow\C$ will often be a functor.
\item The directed graph $D=(V,E)$ representing the diagram $D:\I\dashrightarrow\C$ is not a subgraph of $\C$, even though its vertices are indexed objects of $\C$ and its edges are indexed morphisms of $\C$.
\item The \ul{image of the diagram} $D$ in $\C$ is the subgraph $D'=(V',E')\subset\C$ given by the unindexed collections $V':=\{D(i)\}_{i\in\Ob\I}\subset\Ob\C$ and $E':=\{D(\kappa)\}_{\kappa\in\Mor\I}\subset\Mor\C$.
\item A subgraph of $\C$ has the general form {\small $D'=(V',E')=\big(V'\subset\Ob\C,E'\subset \Mor_\C(V',V')\big)$}, where $\Mor_\C(V',V'):=\bigcup_{u,v\in V'}\Mor_\C(u,v)$, and it is a simple subgraph iff {\small $\big|E'\cap \Mor_\C(u,v)\big|\leq 1$} for all $u,v\in V'$.
\item As shown in Proposition \ref{SubgFuntRep}, a subgraph of $\C$ is equivalently the image of a functor $\D:\I\ra\C$, for some category $\I$.
\end{enumerate}

A diagram is a \ul{simple diagram} if it is a simple graph. The diagram $D=(V,E)$ is a \ul{commutative diagram} if for any two vertices $u,v\in V$, any two directed paths between them, {\small $\gamma:=\big(u=v_1\sr{f_1}{\ral}v_2\sr{f_2}{\ral}\cdots\sr{f_n}{\ral}v_{n+1}=v\big)$} and {\small $\gamma':=\big(u=v'_1\sr{f'_1}{\ral}v'_2\sr{f'_2}{\ral}\cdots\sr{f'_{n'}}{\ral}v'_{n'+1}=v\big)$} in $D$ (i.e., with $f_i,f_i'\in E$), have the same composition, i.e.,
\bea
\txt{Comp}(\gamma)~:=~f_n\circ f_{n-1}\circ\cdots\circ f_1=f'_{n'}\circ f'_{n'-1}\circ\cdots\circ f'_1~=:~\txt{Comp}(\gamma').~~~~\txt{(\blue{footnote}\footnotemark)}.\nn
\eea
\footnotetext{Such a relation can of course be accordingly modified (to get a weaker form of commutativity) if an application requires it.}A \ul{simple category} (resp. \ul{commutative category}) is a category that is simple (resp. commutative) as a graph. (\blue{footnote}\footnote{\ul{Notes}: (i) A commutative diagram is necessarily simple since edges that share the same vertices and direction must be equal. (ii) A simple diagram, however, is clearly not necessarily commutative. (iii) In a commutative diagram in a category $\C$ the edges of any directed loop are clearly isomorphisms in $\C$.}). We will now consider a practical class of functoral diagrams, namely, systems (or functors in general). (\magenta{footnote}\footnote{\magenta{Caution}: In the literature a ``diagram'' is simply a functor (which is already what we will call a ``system''). But for us, a diagram does not have to be a functor, and instead, \ul{``system''} really means a \ul{``functorial diagram or functor (to be assumed commutative/faithful whenever necessary)''}.}).

Let $\I,\C$ be categories (with $\I$ commutative if necessary). A functor (faithful if necessary)
\[
S:\I\ra\C,~~\big(i\sr{\kappa_{ij}}{\ral}j\big)~\mapsto~ \big(S_i\sr{s_{ij}}{\ral}S_j\big)~:=~\big(S(i)\sr{S(\kappa_{ij})}{\ral}S(j)\big)~~~~~~\txt{(\blue{footnote}\footnotemark)}
\]
\footnotetext{(i) When necessary, assume $\I$ is commutative and/or $S$ is faithful. (ii) For us, it is enough to assume $\I$ is commutative because most of our applications will require it. However, the commutativity condition can be relaxed/weakened if an exotic application requires it.}is called an \ul{$\I$-system} (or \ul{system over $\I$}) in $\C$ with \ul{system objects} $\{S_i:=S(i)\}_{i\in\Ob\I}$ and \ul{system transition morphisms} $\{s_{ij}:=S(\kappa_{ij})\}_{i,j\in\Ob\I}$. A (faithful) contravariant functor
\[
S:\I\ra\C,~~\big(i\sr{\kappa_{ij}}{\ral}j\big)~\mapsto~ \big(S^j\sr{s^{ji}}{\ral}S^i\big)~:=~\big(S(j)\sr{S(\kappa_{ij})}{\ral}S(i)\big)
\]
is called an \ul{$\I$-cosystem} (or \ul{cosystem over $I$}) in $\C$. A system or cosystem $S:\I\ra\C$ is a \ul{trivial system} (or a system with no transition morphisms) if $\I=\Ob\I$ is a \ul{trivial category} in the sense $\Mor_\I(i,i)=\{id_i\}$ and $\Mor_\I(i,j)=\emptyset$ for all $i\neq j$ in $\Ob\I=\I$ (in which case $S$ is just a collection of objects $\{S(i)\in\Ob\C\}_{i\in\I}$ in $\C$). A system or cosystem $S:\I\ra\C$ is an \ul{isosystem} if all transition morphisms $S(\kappa_{ij})$ are isomorphisms.

Given systems $S,S':\I\ra\C$, a morphism of functors $f:S\ra S'$ is called a \ul{morphism of systems}. That is, a morphism of systems $f:S\ra S'$ is a collection of morphisms in $\C$
{\small\bea
f=\Big\{S_i\sr{f_i}{\ral}S_i'\Big\}_{i\in\Ob\I}~\subset~\Mor\C=\Mor_\C(\Ob\C,\Ob\C)\nn
\eea} such that for any $i,j\in\I$ and and $\kappa_{ij}\in \Mor_\I(i,j)$, we have a commutative diagram

{\small\[
\left(S\sr{\tau}{\ral}S'\right)\left(i\sr{\kappa_{ij}}{\ral}j\right)~~=~~
\bt
S_i\ar[d,"f_i"]\ar[rr,"s_{ij}"]&&S_j\ar[d,"f_j"]\\
S'_i\ar[rr,"s_{ij}'"]&&S'_j
\et~~:=~~
\bt
S(i)\ar[d,"f_i"]\ar[rr,"S(\kappa_{ij})"]&&S(j)\ar[d,"f_j"]\\
S'(i)\ar[rr,"S'(\kappa_{ij})"]&&S'(j)
\et~~~~~~~~f_j\circ s_{ij}=s'_{ij}\circ f_i.
\]}

Given cosystems {\small$S,S':\I\ra\C$}, a morphism of cofunctors {\small$f:S\ra S'$} is called a \ul{morphism of cosystems}. That is, a morphism of cosystems {\small$f:S\ra S'$} is a collection of morphisms in $\C$
{\small\bea
f=\Big\{S^i\sr{f^i}{\ral}S'{}^i\Big\}_{i\in\Ob\I}~\subset~\Mor\C=\Mor_\C(\Ob\C,\Ob\C)\nn
\eea} such that for any $i,j\in\I$ and and $\kappa_{ij}\in \Mor_\I(i,j)$, we have a commutative diagram

{\small\[
\left(S\sr{\tau}{\ral}S'\right)\left(i\sr{\kappa_{ij}}{\ral}j\right)~~=~~
\bt
S^i\ar[d,"f^i"]\ar[from=rr,"s^{ji}"']&&S^j\ar[d,"f^j"]\\
S'{}^i\ar[from=rr,"s'{}^{ji}"']&&S'{}^j
\et~~:=~~
\bt
S(i)\ar[d,"f^i"]\ar[from=rr,"S(\kappa_{ij})"']&&S(j)\ar[d,"f^j"]\\
S'(i)\ar[from=rr,"S'(\kappa_{ij})"']&&S'(j)
\et~~~~~~~~f^i\circ s^{ji}=s'{}^{ji}\circ f^j.
\]}

The \ul{category of $\I$-systems in $\C$} is the category
\[
\textstyle\C^\I=(\C^\I,\Mor_{\C^\I},\circ)=(\Ob\C^\I,\Mor\C^\I,\circ)=\Big(\Ob\C^\I,\Mor_{\C^\I}(\Ob\C^\I,\Ob\C^\I),\circ\Big)
\]
whose objects $\Ob\C^\I$ are $\I$-systems in $\C$, its morphisms are morphisms of systems, and the composition {\small $\big(S\sr{g\circ f}{\ral}S''\big)~=~\big(S\sr{f}{\ral}S'\sr{g}{\ral}S''\big)$} of two morphisms of systems {\small $S\sr{f}{\ral} S'$}, {\small $S'\sr{g}{\ral}S''$} is their composition as morphisms of functors. The \ul{category of $\I$-cosystems in $\C$} (i.e., \ul{category of $\I^{op}$-systems in $\C$}) is the category
\[
\textstyle \C^{\I^{op}}=(\C^{\I^{op}},\Mor_{\C^{\I^{op}}},\circ)=(\Ob\C^{\I^{op}},\Mor\C^{\I^{op}},\circ)=\Big(\Ob\C^{\I^{op}},\Mor_{\C^{\I^{op}}}(\Ob\C^{\I^{op}},\Ob\C^{\I^{op}}),\circ\Big)
\]
whose objects $\Ob\C^{\I^{op}}$ are $\I$-cosystems (i.e., $\I^{op}$-systems) in $\C$ , its morphisms are morphisms of cosystems, and the composition {\small $\big(S\sr{g\circ f}{\ral}S''\big)~=~\big(S\sr{f}{\ral}S'\sr{g}{\ral}S''\big)$} of two morphisms of cosystems {\small $S\sr{f}{\ral} S'$}, {\small $S'\sr{g}{\ral}S''$} is their composition as morphisms of contravariant functors.

Given categories $\I_1,...,\I_k,\C$, and $\I:=\I_1\times\cdots\times\I_k$, an $\I$-system $S\in\C^\I$ is a \ul{$k$-multisystem} in $\C$. A \ul{bisystem} is a $2$-multisystem.
\end{dfn}

\begin{prp}[\blue{A subgraph is precisely the image of a functor}]\label{SubgFuntRep}
Let $\C$ be a category and $G=(V,E)\subset\C$ a subgraph of $\C$. Then there exists a category $\I$ and a functor $F:\I\ra\C$ such that ~$\im F=G$. (Conversely, it is clear that the image of a functor is a subgraph.)
\end{prp}
\begin{proof}
For each $v\in V$, let {\small $E_v:=\{e\in E:\dom e=v\}$, $E^v:=\{e\in E:\cod e=v\}$, $W_v:=(\{v\}\times E_v)\cup(E^v\times\{v\})$, and $W:=\bigsqcup_{v\in V}W_v$}. Consider the directed graph $\I:=(W,E)$ in which each $e\in E$ is the unique path from $(\dom e,e)$ to $(e,\cod e)$ and there are no other paths except the identities (since they are self inverses). Consider the obvious composition operation on $\I$ based on the fact that each $e\in E$ satisfies $\{e\}=E\cap\Mor_\C(\dom e,\cod e)$. Then $\I=(W,E,\circ)$ is trivially a category (since distinct edges are isolated, hence non-composable), and $G$ is the image of the functor $F:\I\ra\C$ given by $F|_E:=id_E$ and $F|_{W_v}:=\txt{constant}=v$ for each $v\in V$.
\end{proof}

\begin{dfn}[\blue{\index{Induced! system functor}{Induced system functor}}]
Let $\I,\C,\D$ be categories and consider any functor
\[
F:\C\ra\D,~(C\sr{f}{\ral}C')\ra\big(F(C)\sr{F(f)}{\ral}F(C')\big).
\]
Then (by the functor properties) $F$ gives the following functor (the \ul{induced system functor}):
\[
F:\C^\I\ra\D^\I,~(S\sr{\eta}{\ral}S')\ra\big(F\circ S\sr{F\circ\eta}{\ral}F\circ S'\big).
\]
\end{dfn}

\section{Graded objects/morphisms and Categories of graded objects}
\begin{dfn}[\textcolor{blue}{\index{Graded! object (or Grading) in a category}{Graded object (or Grading) in a category}}]
Let $\I,\C$ be categories (with $\I$ commutative if desired). An \ul{$\I$-graded object} (or an \ul{$\I$-grading}) in $\C$ is a functor $S:\I\ra\C$ (i.e., an $\I$-system in $\C$). (\blue{footnote}\footnote{Note that an $\I$-graded object is an $\I$-indexed object, but of a special type because of the functorial property of the indexing operation, i.e., the grading operation is a functorial indexing operation (hence a functor, and not just any map).})
\end{dfn}

\begin{dfn}[\textcolor{blue}{
\index{Graded! morphism}{Graded morphism},
\index{Degree of a graded morphism}{Degree of a graded morphism},
\index{Category of! graded objects}{Category of graded objects}}]
Let $\I,\C$ be categories, $g:\I\ra\I$ an isomorphism, and $S,T\in\C^\I$ any $\I$-graded objects. An \ul{$\I$-graded morphism of degree $g$} from $S$ to $T$ (written $f:S\ra T$) is a morphism of systems $\wt{f}:S\ra T\circ g$ (i.e., $\wt{f}_i:S_i\ra T_{g(i)}$ for all $i\in I$). We write \ul{$d_f:=\deg f:=g$}. If $\deg f=id_\I$ (i.e., $\wt{f}:S\ra T$), we say $f$ is an \ul{$\I$-graded morphism}. The \ul{category of $\I$-graded objects} in $\C$ is the category $Gr(\C^\I)$ whose (i) objects are $\I$-systems $S\in\Ob\C^\I$, (ii) morphisms $f\in \Mor_{Gr(\C^\I)}(S,T)$ are $\I$-graded morphisms $f:S\ra T$, i.e.,
\bea
\textstyle \Mor_{Gr(\C^\I)}(S,T):=\bigcup\limits_{\txt{isos}~g:\I\ra\I}\Mor_{\C^\I}(S,T\circ g),\nn
\eea
and (iii) the composition of $\I$-graded morphisms $f'\circ f:S\sr{f}{\ral}S'\sr{f'}{\ral}S''$ is given by
\[
\wt{f'\circ f}:=(\wt{f'}\circ d_f)\circ\wt{f}:S\sr{\wt{f}}{\ral}S'\circ d_f\sr{\wt{f'}\circ d_f}{\ral}S''\circ d_f\circ d_{f'},~~\Ra~~~~d_{f'\circ f}=d_f\circ d_{f'}.\nn
\]
\end{dfn}

A further discussion of grading for specific applications is given in chapter \ref{AlgebraCatS12}.

\section{Directed systems, Linear systems, Directed Multisystems}
\begin{dfn}[\textcolor{blue}{
\index{Directed! system}{Directed system},
\index{Directed! cosystem}{Directed cosystem},
\index{Linear! system}{Linear system},
\index{Linear! cosystem}{Linear cosystem},
\index{Sequence in a category}{Sequence in a category},
\index{Cosequence in a category}{Cosequence in a category},
\index{Directed! multisystem}{Directed multisystem},
\index{Directed! bisystem}{Directed bisystem},
\index{Mixed variant multisystem}{Mixed-variant multisystem}
}]
Let $\I,\C$ be categories. An $\I$-system (resp. $\I$-cosystem) $S:\I\ra\C$ is said to be \ul{directed} if $\I$ is a directed set (directed class) category, i.e., the following hold:
\bit
\item[(i)] $\Ob\I=(\Ob\I,\leq)$ is a directed set (class), i.e., for any $i,j\in\Ob\I$, there is $k\in\Ob\I$ such that $i,j\leq k$.
\item[(ii)] For each $i,j\in\Ob\I$,~~
$\Mor_\I(i,j)=
\left\{
  \begin{array}{ll}
    \left\{\txt{one morphism}~~i\sr{\kappa_{ij}}{\ral}j\right\}, & \txt{if}~~i\leq j, \\\vspace{-0.3cm}~\\
    \emptyset, & \txt{otherwise}.
  \end{array}
\right.$
\eit
A \ul{linear system} (resp. \ul{linear cosystem}) is a directed $\I$-system (resp. $\I$-cosystem) $S:\I\ra\C$ for which the directed set $\Ob\I=(\Ob\I,\leq)$ is linear, i.e., for all $i,j\in\Ob\I$, $i\leq j$ or $j\leq i$.

If $\I=\Integer$, each of the linear systems $\C^\Integer:=\{\txt{systems}~S:\Integer\ra\C\}$ (resp. linear cosystems $\C^{\Integer^{op}}:=\{\txt{cosystems}~S:\Integer\ra\C\}$) is called a \ul{sequence} (resp. \ul{cosequence}) of objets and morphisms in $\C$. A sequence (resp. cosequence) is denoted by
\bea
&&\bt
S=S_\ast=(S_\ast,s_\ast):~\cdots\ar[r,"s_{n+2}"]&S_{n+1}\ar[r,"s_{n+1}"]&S_n\ar[r,"s_n"]&S_{n-1}\ar[r,"s_{n-1}"]&\cdots
~~~~\txt{(sequence)}\et\nn\\
&&\bt
S=S^\ast=(S^\ast,s^\ast):~\cdots\ar[from=r,"s^{n+1}"']&S^{n+1}\ar[from=r,"s^n"']&S^n\ar[from=r,"s^{n-1}"']&S^{n-1}\ar[from=r,"s^{n-2}"']&\cdots
~~~~\txt{(cosequence)}\et\nn
\eea

Given linear set (i.e., linear directed set) categories $\I_1$, $\I_2$, $\cdots$, $\I_k$ and $\I:=\I_1\times\cdots\times\I_k$ (where $\Ob\I$ is of course viewed as a directed set with $(i_1,...,i_k)\leq(i_1',...,i_k')$ iff $i_j\leq i_j'$ for all $j=1,...,k$), a system $S\in\C^\I$ is a \ul{directed $k$-multisystem} in $\C$. A \ul{directed bisystem} is a directed $2$-multisystem $S\in\C^{\I_1\times\I_2}$. A directed bisystem $S_{\ast\ast}\in\C^{\Integer\times\Integer}$ can be displayed in the form

{\footnotesize
\bea\bt
 ~       &  \ar[from=dd,"s^v_{i+1,j-1}"']  ~                         &&  \ar[from=dd,"s^v_{i,j-1}"']  ~                        &&   \ar[from=dd,"{s^v_{i-1,j-1}}"']  ~                              && ~\\
 ~       &    ~                                                         &&    ~                                                      &&    ~                                                                 && ~\\
 ~\ar[r] & S_{i+1,j-1}\ar[rr,"s^h_{i+1,j-1}"]                        && S_{i,j-1}\ar[rr,"s^h_{i,j-1}"]                         && S_{i-1,j-1}\ar[rr,"s^h_{i-1,j-1}"] && ~\\
 ~       &     ~                                                        &&     ~                                                     &&     ~                                                           && ~\\
 ~\ar[r] & S_{i+1,j}\ar[uu,"s^v_{i+1,j}"']\ar[rr,"{s^h_{i+1,j}}"] && S_{ij}\ar[rr,"s^h_{ij}"]\ar[uu,"s^v_{ij}"'] && S_{i-1,j}\ar[uu,"{s^v_{i-1,j}}"']\ar[rr,"s^h_{i-1,j}"] && ~\\
 ~       &     ~                                                        &&     ~                                                     &&     ~                                                           && ~\\
 ~\ar[r] & S_{i+1,j+1}\ar[uu,"{s^v_{i+1,j+1}}"']\ar[rr,"s^h_{i+1,j+1}"] && S_{i,j+1}\ar[uu,"s^v_{i,j+1}"']\ar[rr,"s^h_{i,j+1}"]   && S_{i-1,j+1}\ar[uu,"s^v_{i-1,j+1}"']\ar[rr,"s^h_{i-1,j+1}"]    && ~\\
 ~       &   \ar[u]  ~                           &&    \ar[u] ~                         &&   \ar[u]  ~                              && ~
\et\nn
\eea}
\end{dfn}

In typical applications, we consider a multisystem $S\in\C^{\I_1\times\cdots\times\I_k}$ having \ul{mixed variance} in the sense that $\I_1,...,\I_k\in\{\I,\I^{op}\}$ for some directed set category $\I$. For example (say with $\I:=\Integer$),
\[
S^\ast{}_\ast\in\C^{\I^{op}\times\I},~~~~S_\ast{}^{\ast\ast}\in\C^{\I\times\I^{op}\times\I^{op}},~~~~S_{\ast\ast\ast}{}^\ast\in\C^{\I^3\times\I},~~~~ \cdots\nn
\]

\begin{rmk}[\textcolor{blue}{Images of diagrams and systems are not categories}]
The image of a diagram in a category $\C$ is not a subcategory of $\C$ because (from its definition) the diagram's image does not have to be closed under composition in $\C$. In particular, the image of an $\I$-system (a functor $\I\ra\C$) in $\C$ is not a subcategory because it could happen that it is not closed under composition in $\C$.

Indeed, if $F:\I\ra\C$ is a functor, morphisms $\kappa,\eta$ in $\I$ that are non-composable because $\dom\kappa\neq\cod\eta$ and $\cod\kappa\neq\dom\eta$ can nevertheless have composable images $F(\kappa),F(\eta)$ in $\C$ if for example $\cod F(\kappa):=F(\cod\kappa)=F(\dom\eta)=:\dom F(\eta)$, in which case, the composition $F(\eta)\circ F(\kappa):\dom F(\kappa)\ra\cod F(\eta)$ might not lie in $\im F$ (and so $\im F$ might not be closed under composition).
\end{rmk}

\begin{question}
We know that if a system (functor) is an imbedding then its image is a subcategory. If the image of a system (functor) happens to be a subcategory, does it follow that the system (functor) is an imbedding?
\end{question}

\section{Processes Between Categories of Systems: Limit Processes}
\begin{dfn}[\textcolor{blue}{
\index{Process between categories of systems}{{Process between categories of systems}}, 
\index{Inductive! process}{Inductive process},
\index{Induced! morphism of systems}{Induced morphism of systems},
\index{Canonical process}{Canonical process},
\index{Functorial (or natural) process}{Functorial (or natural) process},
\index{Cofunctorial (or conatural) process}{Cofunctorial (or conatural) process},
\index{Subprocess}{Subprocess}}]
Let $\C^\I$ (resp. $\D^\J$) be the category of systems in $\C$ over $\I$ (resp. the category of systems in $\D$ over $\J$). A \ul{process} between $\C^\I$ and $\D^\J$ (denoted $P:\C^\I\dashrightarrow\D^\J$) is a map $P:\Ob\C^\I\ra\Ob\D^\J$, $S\ra P(S)$ between their objects. The process $P$ is said to be \ul{inductive} if for every morphism of systems $S\sr{\tau}{\ral}S'$ in $\C^\I$, the process \ul{induces} (i.e., is by construction automatically accompanied by) \ul{at least one} morphism of systems $P(S)\sr{\tau_P}{\ral}P(S')$ in $\D^\J$. The process $P$ is said to be \ul{canonical} if (i) it is inductive and (ii) for every morphism of systems $S\sr{\tau}{\ral}S'$ in $\C^\I$, the process induces a unique (i.e., only one) morphism of systems $\tau_P:P(S)\ra P(S')$. That is, a \ul{canonical process} between $\C^\I$ and $\D^\J$ is precisely a map of the form
\bea
P:\Ob\C^\I\sqcup\Mor\C^\I\ra\Ob\D^\J\sqcup\Mor\D^\J,~~S\mapsto P(S),~~\big(S\sr{\tau}{\ral}S'\big)\mapsto\Big(P(S)\sr{P(\tau)}{\ral}P(S')\Big).\nn
\eea

A canonical process $P:\C^\I\dashrightarrow\D^\J$ is \ul{functorial or natural}, making it a functor or natural transformation $P:\C^\I\ra\D^\J$, (resp. \ul{cofunctorial or conatural}, making it a cofunctor or conatural transformation) if it satisfies the axioms of a functor (resp. cofunctor) as follows: For all $S\in\Ob\C^\I$ and $\tau,\tau'\in\Mor\C^\I$,
{\flushleft (i)} $P(id_S)=id_{P(S)}$~ and ~(ii) $P(\tau\circ\tau')=P(\tau)\circ P(\tau')$ ~$\Big($resp. $P(\tau\circ\tau')=P(\tau')\circ P(\tau)$ $\Big)$.

Given two processes $P,Q:\C^\I\dashrightarrow\D^\J$, $P$ is a \ul{subprocess} of $Q$ (written $P\subset Q$) if $P(S)\subset Q(S)$ for every $S\in\Ob\C^\I$.
\end{dfn}

\begin{rmk*}[\blue{The idea of system limits}]
Limits and colimits, as we shall see, are a useful class of functorial processes. Let $S:\I\ra\C$ be a system (resp. cosystem) in a category $\C$. A colimit (resp. limit) of $S$ in $\C$ is a \ul{single object} in $\C$ that acts as a morphism guiding center or gateway through which all simultaneous morphisms of a  \ul{specific type} must factor, or equivalently pass, if they are converging from objects of $S$ to an object of $\C$ (resp. diverging from an object of $\C$ to objects of $S$). Thus, a (co)limit of $S$ (approximately) summarizes a particular aspect of $S$. (\blue{footnote}\footnote{We of course also have \ul{mixed limit} types involving objects acting as morphism guiding centers or gateways for simultaneous morphisms with some converging from objects of $S$ and some diverging to objects of $S$.})

More specifically (based on the clear meanings of coproducts as encountered earlier), the colimit $\varinjlim S$ of a system $S:\I\ra\C$ is, up to isomorphism, the unique \ul{single object} in $\C$ that best (i) acts as a \ul{containing object} for all objects of $S$ and (ii) supports the interaction/transition morphisms between those objects. That is, the colimit ~$\varinjlim S$ = best ``\emph{containing object + interaction/transition support}'' for objects of $S$.

Similarly (based on the clear meanings of products as encountered earlier), the limit $\varprojlim S$ of a cosystem $S:\I\ra\C$ is, up to isomorphism, the unique \ul{single object} in $\C$ that best (i) acts as a \ul{numerator object} for all objects of $S$ and (ii) supports the interaction/transition morphisms between those objects. That is, the limit ~$\varprojlim S$ = best ``\emph{numerator object + interaction/transition support}'' for objects of $S$.
\end{rmk*}

\section{Colimits, Directed colimits, Coproducts}
In order to minimize confusion, we will place the ``\emph{system on the left}'' and its ``\emph{limit object on the right}'', and adopt the following conventions/interpretations (which seem to be most standard):
\bit[leftmargin=0.8cm]
\item[(i)] View a \index{Colimit}{\ul{colimit}} (\index{Injective! limit}{\ul{injective limit}} or \index{Inductive! limit}{\ul{inductive limit}}) as a \ul{right-limit} (or \index{Pushout}{\ul{``pushout''}} of the system), denoted by $\varinjlim$. Here, system objects ``\emph{inject}'', ``\emph{induct}'', or ``\emph{push}'' into the limit object. (\blue{footnote}\footnote{Special colimits in some categories will be called ``\emph{union}'', ``\emph{coproduct}'' (``\emph{direct sum}''), ``\emph{pushout}''.})
\bea
S(i)\sr{q_i}{\ral}\varinjlim S\nn
\eea
\item[(ii)] View a \index{Limit}{\ul{limit}} (\index{Projective! limit}{\ul{projective limit}}) as a \ul{left-limit} (or \index{Pullback}{\ul{``pullback''}} of the system), denoted by $\varprojlim$. Here, system objects ``\emph{project}'' or ``\emph{pull}'' out of the limit object. (\blue{footnote}\footnote{Special limits in some categories will be called ``\emph{intersection}'', ``\emph{product}'' (``\emph{direct product}''), ``\emph{pullback}''.})
\bea
S(i)\sr{p_i}{\lal}\varprojlim S\nn
\eea
\eit
\begin{note*}[\textcolor{blue}{Generality and the arrow under $\lim$}]
The left/right arrows under $\lim$ are only meant to give the (obviously false, but potentially convenient) impression that the ``\emph{system is always on the left}'' and its ``\emph{limit object is always on the right}'', so that (as explained above) it is helpful especially in the beginning to think of the \emph{colimit} as a \ul{\emph{right-limit}} or \ul{``\emph{pushout}''} of the system (hence $\varinjlim$) and the \emph{limit} as a \ul{\emph{left-limit}} or \ul{``\emph{pullback}''} of the system (hence $\varprojlim$).

Therefore, the arrow under $\lim$ does not imply the (co)limit is ``directed''. Instead, we will explicitly say so whenever a (co)limit is ``directed''. Also, ``directedness'' of a (co)limit will sometimes be obvious from the context if the system in question happens to be a directed system.
\end{note*}

Recall that we have two main types of systems, namely, \ul{contravariant-systems} (\ul{cosystems}) and \ul{covariant-systems} (\ul{systems}). In some traditional settings, one tends to take the \ul{colimit of a system} and to take the \ul{limit of a cosystem}. However, given any system (covariant or contravariant), one can take its limit, its colimit, and any other desired type of limit as well. Consequently, we will encounter many different types of limits in various applications.

\begin{dfn}[\textcolor{blue}{
\index{System-compatible collection}{System-compatible collection}, 
\index{Colimit of a system}{Colimit of a system}, 
\index{Weak! colimit}{Weak colimit}, 
\index{Directed! colimit of a directed system}{Directed colimit of a directed system}, 
\index{Coproduct! (Direct sum) of! a trivial system}{Coproduct (Direct sum) of a trivial system}, \index{Cokernel of a morphism}{Cokernel of a morphism} (\blue{footnote}\footnote{The cokernel makes sense only for additive/abelian categories, to be defined later, but it is mentioned here to emphasize the importance of colimits.})}]
Let $S:\I\ra\C$ be a system in a category $\C$. A collection of morphisms {\small $\Big\{S(i)\sr{f_i}{\ral} X_i,~X_i\sr{\psi_{ij}}{\ral}X_j~\big|~i,j\in\Ob\I\Big\}$} in $\C$ is $S$-\ul{compatible} (i.e., compatible with the system $S$) if for any $i,j\in\Ob\I$ and any $\kappa_{ij}\in \Mor_\I(i,j)$, we have a commutative diagram~~
\bt
S(i)\ar[d,"f_i"]\ar[r,"S(\kappa_{ij})"] & S(j)\ar[d,"f_j"']\\
X_i\ar[r,"\psi_{ij}"] & X_j
\et~~~~$\psi_{ij}\circ f_i=f_j\circ S(\kappa_{ij})$.

A \ul{colimit} of $S$ in $\C$ consists of an object $\varinjlim S\in\Ob\C$ and morphisms $\Big\{q_i:S(i)\ra\varinjlim S~\big|~i\in\Ob\I\Big\}$ such that the following hold.
\bit[leftmargin=0.9cm]
\item[(i)] \ul{Compatibility}: The collection {\small $\Big\{S(i)\sr{q_i}{\ral} \varinjlim S,~\varinjlim S\sr{id}{\ral}\varinjlim S~\big|~i\in\Ob\I\Big\}$} is $S$-compatible, i.e., for any $i,j\in\Ob\I$ and any $\kappa_{ij}\in \Mor_\I(i,j)$, we have a commutative diagram
\bea
\bt[row sep=tiny]
S(i)\ar[dr,"q_i"']\ar[rr,"S(\kappa_{ij})"] & & S(j)\ar[dl,"q_j"] \\
 & \varinjlim S &
\et~~~~~~~~q_i=q_j\circ S(\kappa_{ij})\nn
\eea

\item[(ii)] \ul{Universal property}: For any object $X\in\Ob\C$ and any $S$-compatible collection
{\footnotesize\[
\txt{{\small $\Big\{S(i)\sr{f_i}{\ral} X,~X\sr{id_X}{\ral}X~\big|~i\in\Ob\I\Big\}$}},~~\txt{i.e., for al $i,j,\kappa_{ij}$,}~~
\bt[row sep=tiny]
S(i)\ar[dr,"f_i"']\ar[rr,"S(\kappa_{ij})"] & & S(j)\ar[dl,"f_j"] \\
 & X &
\et~~~~f_i=f_j\circ S(\kappa_{ij}),\nn
\]}there exists a \ul{unique} morphism $f:\varinjlim S\ra X$ such that the following diagram commutes.
\bea\adjustbox{scale=0.9}{\bt
S(i)\ar[dr,"q_i"']\ar[rr,"S(\kappa_{ij})"]\ar[ddr,bend right,"f_i"'] &&  S(j)\ar[dl,"q_j"]\ar[ddl,bend left,"f_j"]\\
 & \varinjlim S\ar[d,dashed,"f"]& \\
 & X &
\et}~~~~~~~~f_i=f\circ q_i,~~\txt{for all}~~i\in \Ob\I.\nn
\eea
\eit

If the induced morphism $f$ is not required to be unique, then the colimit is called a \ul{weak colimit}. The colimit $\varinjlim:\C^\I\ra\C$ is called \ul{directed} if $\C^\I$ is a category of directed systems.

If $S:\I\ra\C$ is a trivial system (often written $\{S(i)\}_{i\in\I}$), then $\varinjlim S$ is called the \ul{direct sum} $\bigoplus_{i\in\I} S(i)$, or the \ul{coproduct} $\coprod_{i\in\I} S(i)$, of $S$.

\bea\adjustbox{scale=0.9}{\bt[column sep=tiny]
S(i)\ar[dr,"q_i"]\ar[ddr,bend right,"f_i"'] &&  S(j)\ar[dl,"q_j"']\ar[ddl,bend left,"f_j"]\\
 & \scalebox{0.7}{$\bigoplus\limits_{i\in\I} S(i)~\txt{or}~\coprod\limits_{i\in\I} S(i)$}\ar[d,dashed,"f"]& \\
 & X &
\et}~~~~~~~~f_i=f\circ q_i,~~\txt{for all}~~i\in \I.\nn
\eea

As we will see later in Definition \ref{KerAndCoker}, if $\C$ is an \ul{abelian category}, then given a morphism $f:A\ra B$ in $\C$, the \ul{cokernel} of $f$ (viewed as a system) satisfies
\bea
\lim_{\ral}\Big(A\sr{f}{\ral}B\Big)=\Big\{A\sr{0}{\ral}\coker f,~B\sr{c}{\ral}\coker f\Big\}.\nn
\eea

\end{dfn}

\begin{rmk}\label{ColimUnqRmk}
The colimit $\varinjlim S$ of a system $S:\I\ra\C$ is unique up to isomorphism in $\C$.
\end{rmk}
\begin{proof}
Let $\{q'_i:S(i)\ra X~|~i\in\I\}$ be another limit as in the diagram below:
\bea\adjustbox{scale=0.9}{\bt
S(i)\ar[dr,"q_i"']\ar[rr,"S(\kappa_{ij})"]\ar[ddr,bend right,"q'_i"'] &&  S(j)\ar[dl,"q_j"]\ar[ddl,bend left,"q'_j"]\\
 & \varinjlim S\ar[d,dashed,shift left,"q'"]& \\
 & X\ar[u,dashed,shift left,"q"] &
\et}~~~~~~~~q_i=q\circ q'_i,~~q'_i=q'\circ q_i,~~\txt{for all}~~i\in \Ob\I.\nn
\eea
Then $q_i=q\circ q'\circ q_i$, $q_i'=q'\circ q\circ q'_i$, and so by the uniqueness requirement in the universal property, we get $q'\circ q=id_X$ and $q\circ q'=id_{\varinjlim S}$.
\end{proof}

\section{Limits, Directed limits, Products}
\begin{dfn}[\textcolor{blue}{
\index{Cosystem-compatible collection}{Cosystem-compatible collection},
\index{Limit! of a cosystem}{Limit of a cosystem},
\index{Weak! limit}{Weak limit},
\index{Directed! limit of a directed cosystem}{Directed limit of a directed cosystem},
\index{Product of! a trivial system}{Product of a trivial system},
\index{Kernel of! a morphism}{Kernel of a morphism} (\blue{footnote}\footnote{The kernel makes sense only for additive/abelian categories, to be defined later, but it is mentioned here to emphasize the importance of limits.})}]
Let $S:\I\ra\C$ be a cosystem in a category $\C$. A collection of morphisms {\small $\Big\{S(i)\sr{f^i}{\lal} X^i,~X^i\sr{\psi^{ji}}{\lal}X^j~\big|~i,j\in\Ob\I\Big\}$} in $\C$ is $S$-\ul{compatible} (i.e., compatible with the cosystem $S$) if for any $i,j\in\Ob\I$ and any $\kappa_{ij}\in \Mor_\I(i,j)$, we have a commutative diagram~~
\bt
S(i) & S(j)\ar[l,"S(\kappa_{ij})"']\\
X^i\ar[u,"f^i"']  & X^j\ar[l,"\psi^{ji}"']\ar[u,"f^j"']
\et~~

A \ul{limit} of $S$ in $\C$ consists of an object $\varprojlim S\in\Ob\C$ and morphisms {\small $\Big\{S(i)\sr{p^i}{\lal}\varprojlim S~\big|~i\in\Ob\I\Big\}$} such that the following hold.
\bit[leftmargin=0.9cm]
\item[(i)] \ul{Compatibility}:~ The collection {\small $\Big\{S(i)\sr{p^i}{\lal} \varprojlim S,~\varprojlim S\sr{id}{\lal}\varprojlim S~\big|~i\in\Ob\I\Big\}$} is $S$-compatible, i.e., for any $i,j\in\Ob\I$ and any $\kappa_{ij}\in \Mor_\I(i,j)$, we have a commutative diagram

\bea
\bt[row sep=tiny]
S(i)&& S(j)\ar[ll,"S(\kappa_{ij})"']\\
    &\varprojlim S\ar[ul,"p^i"]\ar[ur,"p^j"'] &
\et~~~~~~p^i=S(\kappa_{ij})\circ p^j.\nn
\eea

\item[(ii)] \ul{Universal property}: For any object $X\in\Ob\C$ and any $S$-compatible collection
{\footnotesize
\[\Big\{S(i)\sr{f^i}{\lal} X,~X\sr{id_X}{\lal}X~\big|~i\in\Ob\I\Big\},~~\txt{i.e., for all $i,j,\kappa_{ij}$}~~\bt[row sep=tiny]
S(i)&& S(j)\ar[ll,"S(\kappa_{ij})"']\\
    &\varprojlim S\ar[ul,"f^i"]\ar[ur,"f^j"'] &
\et~~~~f^i=S(\kappa_{ij})\circ f^j,\nn
\]}there exists a \ul{unique} morphism $f:\varprojlim S\ra X$ such that the following diagram commutes.
\bea\adjustbox{scale=0.9}{\bt
S(i) &&  S(j)\ar[ll,"S(\kappa_{ij})"']\\
 & \varprojlim S\ar[ul,"p^i"]\ar[ur,"p^j"']& \\
 & X\ar[u,dashed,"f"]\ar[uul,bend left,"f^i"]\ar[uur,bend right,"f^j"'] &
\et}~~~~~~~~f^i=p^i\circ f,~~\txt{for all}~~i\in \Ob\I.\nn
\eea
\eit

If the induced morphism $f$ is not required to be unique, then the limit is called a \ul{weak limit}. The limit $\varprojlim:\C^{\I^{op}}\ra\C$ is called \ul{directed} if $\C^{\I^{op}}$ is a category of directed cosystems.

If $S:\I\ra\C$ is a trivial system $\{S(i)\in\Ob\C:i\in\I\}$, then $\varprojlim S$ is called the \ul{product} $\prod\limits_{i\in\I} S(i)$ of $S$.
\bea
\adjustbox{scale=0.9}{\bt[column sep=small]
S(i) &&  S(j)\\
 & \prod\limits_{i\in\I} S(i)\ar[ul,"p^i"']\ar[ur,"p^j"]& \\
 & X\ar[u,dashed,"f"]\ar[uul,bend left,"f^i"]\ar[uur,bend right,"f^j"'] &
\et}~~~~~~~~f^i=p^i\circ f,~~\txt{for all}~~i\in\I.\nn
\eea

As we will see later in Definition \ref{KerAndCoker}, if $\C$ is an abelian category, then given a morphism $f:A\ra B$ in $\C$, the \ul{kernel} of $f$ (viewed as a cosystem) satisfies
\bea
\lim_{\lal}\Big(B\sr{f}{\lal}A\Big)=\Big\{B\sr{0}{\lal}\ker f,~A\sr{k}{\lal}\ker f\Big\}.\nn
\eea
\end{dfn}

\begin{rmk}\label{LimUnqRmk}
The limit $\varprojlim S$ of a cosystem $S:\I\ra\C$ is unique up to isomorphism in $\C$. The proof is similar to that for the colimit of a system (Remark \ref{ColimUnqRmk}).
\end{rmk}

\begin{rmk}[\blue{Potentially useful observations}]
We have already seen that a monomorphism (resp. epimorphism) is precisely an epimorphism (resp. monomorphism) in the opposite category. Likewise, the colimit of a system (resp. limit of a cosystem) is precisely the limit of a cosystem (resp. colimit of a system) in the opposite category.
\end{rmk}

\begin{thm}
An equivalence of categories preserves colimits and limits.
\end{thm}
\begin{proof}
Let $E:\C\ra\D$ be an equivalence of categories (i.e., a faithful, full, and dense functor) and $S:\I\ra\C$ a system. Then the conclusion is clear  from the following diagrams of the colimit and limit setups between the two categories.
\[
\adjustbox{scale=0.7}{\bt
S(i)\ar[dr,"q_i"']\ar[rr,"S(\kappa_{ij})"]\ar[ddr,bend right,"f_i"'] &&  S(j)\ar[dl,"q_j"]\ar[ddl,bend left,"f_j"]\\
 & \varinjlim S\ar[d,dashed,"f"]& \\
 & X &
\et}~~\sr{E}{\longleftrightarrow}~~
\adjustbox{scale=0.7}{\bt
E\circ S(i)\ar[dr,"E(q_i)"']\ar[rr,"E\circ S(\kappa_{ij})"]\ar[ddr,bend right,"E(f_i)"'] &&  E\circ S(j)\ar[dl,"E(q_j)"]\ar[ddl,bend left,"E(f_j)"]\\
 & E(\varinjlim S)\ar[d,dashed,"E(f)"]& \\
 & E(X) &
\et}~~=~~
\adjustbox{scale=0.7}{\bt
E\circ S(i)\ar[dr,"E(q_i)"]\ar[rr,"E\circ S(\kappa_{ij})"]\ar[ddr,bend right,"E(f_i)"]\ar[dddr,bend right,"f'_i"'] &&  E\circ S(j)\ar[dl,"E(q_j)"']\ar[ddl,bend left,"E(f_j)"']\ar[dddl,bend left,"f'_j"]\\
 & E(\varinjlim S)\ar[d,dashed,"E(f)"]& \\
 & E(X) & \\
 &  X'\ar[u,"\cong","\vphi"']    &
\et}
\]

\[
\adjustbox{scale=0.7}{\bt
S(i) &&  S(j)\ar[ll,"S(\kappa_{ij})"']\\
 & \varprojlim S\ar[ul,"p^i"]\ar[ur,"p^j"']& \\
 & X\ar[u,dashed,"f"']\ar[uul,bend left,"f^i"]\ar[uur,bend right,"f^j"'] &
\et}~~\sr{E}{\longleftrightarrow}~~
\adjustbox{scale=0.7}{\bt
E\circ S(i) &&  E\circ S(j)\ar[ll,"E\circ S(\kappa_{ij})"']\\
 & E(\varprojlim S)\ar[ul,"E(p^i)"]\ar[ur,"E(p^j)"']& \\
 & E(X)\ar[u,dashed,"E(f)"']\ar[uul,bend left,"E(f^i)"]\ar[uur,bend right,"E(f^j)"'] &
\et}~~=~~
\adjustbox{scale=0.7}{\bt
E\circ S(i) &&  E\circ S(j)\ar[ll,"E\circ S(\kappa_{ij})"']\\
 & E(\varprojlim S)\ar[ul,"E(p^i)"']\ar[ur,"E(p^j)"]& \\
 & E(X)\ar[u,dashed,"E(f)"']\ar[uul,bend left,"E(f^i)"']\ar[uur,bend right,"E(f^j)"] &\\
 &  X'\ar[from=u,"\cong"',"\psi"]\ar[uuul,bend left,"f'{}^i"]\ar[uuur,bend right,"f'{}^j"']    &
\et}
\]
That is, ~$E\big(\varinjlim S\big)\cong \varinjlim E\circ S$~ and ~$E\big(\varprojlim S\big)\cong \varprojlim E\circ S$.
\end{proof}

\section{Directed Colimits/Limits over a Cofinal Subset}

\begin{dfn}[\textcolor{blue}{\index{Cofinal subset of a poset}{Cofinal subset of a poset}}]
Let $(P,\leq)$ be a poset. A subset $C\subset P$ is cofinal if every element of $P$ has an upper bound in $C$, i.e., for every $p\in P$ there exists $c\in C$ such that $p\leq c$.
\end{dfn}
If $(P,\leq)$ is directed and $C\subset P$ is cofinal, then $(C,\leq)$ is directed as well: Indeed for all $c,c'\in C$, there exist $p\in P$ and $c''\in C$ such that $c,c'\leq p\leq c''$, and so $c,c'\leq c''$.

\begin{lmm}[\textcolor{blue}{Limits over a cofinal subset}]
Let $S:\I\ra\C$ be a directed system, $T:\I\ra\C$ a directed cosystem, and $\J\subset\I$ a cofinal subset (hence $S|_\J,T|_\J:\J\subset\I\ra\C$ are directed as well). Then
\[
(i)~\varinjlim S|_{\J}\cong\varinjlim S~~~~\txt{and}~~~~(ii)~\varprojlim T|_{\J}\cong\varprojlim T.
\]
\end{lmm}
\begin{proof}
We prove (i) only, noting that (ii) is similar (except with arrows reversed in the diagram below). Let {\small $\Big\{q_j:S_j\ra \varinjlim_{j\in {\J}}S_j\Big\}_{j\in\J}$} be the limit over ${\J}$, i.e., the limit of {\small $S|_{\J}:\J\subset\I\ra\C$}. Let {\small $s_{ii'}:=S(\kappa_{ii'}):S(i)\ra S(i')$} be the system transition morphisms. Given any morphisms $\{f_i:S_i\ra X\}_{i\in\I}$ such that $f_i=f_{i'}s_{ii'}$ (for all $i'\geq i$), choose $j\geq i$ and $j'\geq i'$ (which is possible since ${\J}\subset {\I}$ is cofinal). Then
\bea
f_i=f_js_{ij},~~~~~f_{i'}=f_{j'}s_{i'j'}.\nn
\eea

\bea
\adjustbox{scale=0.9}{\bt
S_i\ar[ddddrr,bend right,"f_i"']\ar[dr,dashed,"s_{ij}"]\ar[rrrr,"s_{ii'}"] &&&& S_{i'}\ar[ddddll,bend left,"f_{i'}"]\ar[dl,dashed,"s_{i'j'}"']\\
  & S_j\ar[dddr,bend right,"f_j"]\ar[dr,"q_j"']\ar[rr,"s_{jj'}"] &  & S_{j'}\ar[dddl,bend left,"f_{j'}"']\ar[dl,"q_{j'}"]& \\
                      &                                   & \varinjlim_{j\in {\J}}S_j\ar[dd,dashed,"\theta"] &  & \\
                      &                                   &     &  & \\
                      &                                   & X    &  &
\et}\nn
\eea
By the universal property for $\varinjlim_{j\in {\J}}S_j$, there exists a unique morphism $\theta:\varinjlim_{j\in {\J}}S_j\ra X$ such that $f_j=\theta q_j$, $f_{j'}=\theta q_{j'}$. Thus, it follows that
\bea
f_i=f_js_{ij}=\theta q_js_{ij},~~~~~f_{i'}=f_{j'}s_{i'j'}=\theta q_{j'}s_{i'j'}.\nn
\eea
This shows the collection {\small$\Big\{q_js_{ij}:S_i\sr{s_{ij}}{\ral}S_j\sr{q_j}{\ral}\varinjlim_{j\in {\J}}S_j~\big|~i\in {\I},~j\in {\J},~j\geq i\Big\}$} is a limit of the system over ${\I}$. Since the limit over ${\I}$ is unique up to isomorphism, we see that $\varinjlim_{j\in J}S_j\cong \varinjlim_{i\in {\I}}S_i$.
\end{proof}

\section{Mor-stability and Functoriality of Colimits/Limits}
\begin{rmk}[\blue{
\index{Universal! morphism property}{Universal morphism property},
\index{UMP bijection}{UMP bijection}}]
Let $\C$ be a category, $\D\subset\C$ a subcategory, $C\in\Ob\C$, $D\in\Ob\D$, and $f\in\Mor_\C(C,D)$ a $\C$-morphism (i.e., a morphism in $\C$). For any $D'\in\Ob\D$, consider any subclass of morphisms $\Mor_{\C,f}(C,D')\subset\Mor_\C(C,D')$ that is (i) compatible with $f$ in a specific way and (ii) for any $\tau\in\Mor_\D(D,D')$, $\tau\circ f\in\Mor_{\C,f}(C,D')$. Let $\Mor(f,\D):=\bigcup_{D'\in\Ob\D}\Mor_{\C,f}(C,D')$.

We say the $\C$-morphism $f\in\Mor_\C(C,D)$ has a \ul{universal morphism property (UMP)} (wrt $\Mor(f,\D)$ and $\D$) if for any $\C$-morphism $g\in\Mor_{\C,f}(C,D')$, there exists a \textbf{\emph{unique}} $\D$-morphism $\wt{g}\in\Mor_\D(D,D')$ such that $\wt{g}\circ f=g$ (i.e., diagram \ul{\emph{commutes}}).
\[\bt
C\ar[d,"g"']\ar[rr,"f"] && D\ar[dll,dashed,"\exists!~{\wt{g}}"] \\
D' &&
\et\]

This gives a \textbf{\emph{well-defined}} map {\small $\phi_f:\Mor_{\C,f}(C,D')\ra\Mor_\D(D,D'),~g\mapsto\wt{g}$}. Moreover, $\phi_f$ is \ul{\emph{injective}} because if $\wt{g}_1=\wt{g}_2$, then $g_1=\wt{g}_1\circ f=\wt{g}_2\circ f=g_2$.
On the other hand, by definition/hypotheses, given a $\D$-morphism $\tau\in \Mor_\D(D,D')$, $\tau\circ f:C\ra D'$ is a $\C$-morphism (but not necessarily a $\D$-morphism) that is $f$-compatible. So, there is a unique $\D$-morphism $\wt{\tau\circ f}\in \Mor_\D(D,D')$ satisfying $\wt{\tau\circ f}\circ f=\tau\circ f$. By \ul{\emph{uniqueness}}, we have $\wt{\tau\circ f}=\tau$ (as in the following commutative diagram):
\[\bt
C\ar[d,"\tau\circ f"']\ar[rr,"f"] && D\ar[dll,shift right,"\tau"']\ar[dll,dashed,shift left,"\exists!~{\wt{\tau\circ f}}"] \\
D' &&
\et\]
That is, the map $\phi_f$ is \ul{\emph{surjective}}. Hence we have a bijection (which we will call the \ul{UMP bijection})
\begin{align}
\Mor_{\C,f}(C,D')\sr{\phi_f}{\longleftrightarrow}\Mor_\D(D,D').
\end{align}
\end{rmk}
The following result is a ``generalization'' of the above situation with surjectivity already built-in (but without explicit mention of the subcategory to which the limit object is possibly restricted).

\begin{prp}[\textcolor{blue}{\index{Mor-stability of colimits and limits}{Mor-stability of colimits and limits}}]\label{HomStabLim}
Let $\I,\C$ be categories, $S:\I\ra\C$ a system, and $X\in\Ob\C$. The limits (i.e., limit and colimit) of the system satisfy the following (\magenta{footnote}\footnote{\magenta{Usefulness of the results}: In practice, $\varinjlim S(i)$ or $\varprojlim S(i)$ (as well as the variable object $X$) is an object in a certain relevant/restricted subcategory $\D\subset\C$. In that case, the limits $\varprojlim \Mor_\C\big(S(i),X\big)$ and $\varprojlim \Mor_\C\big(X,S(i)\big)$ in Sets should really be the limits $\varprojlim \Mor_\D\big(S(i),X\big)$ and $\varprojlim \Mor_\D\big(X,S(i)\big)$ taken in a likewise relevant/restricted subcategory $\wt{D}\subset\txt{Sets}$.}):
\bea
(a)~\Mor_\C\left(\varinjlim S(i),X\right)\cong\varprojlim \Mor_\C\big(S(i),X\big),~~~~(b)~\Mor_\C\left(X,\varprojlim S(i)\right)\cong\varprojlim \Mor_\C\big(X,S(i)\big).\nn
\eea
\end{prp}
\begin{proof}
It is enough to prove (a) only, since (b) is similar (except with arrows reversed in the first diagram below). Recall that the colimit setup in $\C$ is given by the following commutative diagram:
\bea
\adjustbox{scale=0.9}{\bt
S(i)\ar[dr,"q_i"']\ar[rr,"S(\kappa_{ij})"]\ar[ddr,bend right,"f_i"'] &&  S(j)\ar[dl,"q_j"]\ar[ddl,bend left,"f_j"]\\
 & \varinjlim S\ar[d,dashed,"h"]& \\
 & X &
\et}~~~~~~~~f_i=h\circ q_i,~~\txt{for all}~~i\in \Ob\I.\nn
\eea
Let the class of all $S$-compatible collections of morphisms $(f_i:S(i)\ra X)_{i\in\Ob\I}$ in $\C$ be denoted by
\bea
&&\textstyle\txt{CMP}_\C(S,X):=\left\{f=\big(f_i:S(i)\ra X\big)_{i\in\Ob\I}\in\prod_{i\in\Ob\I}\Mor_\C(S(i),X):~\txt{$f$ is $S$-compatible}\right\}\nn\\
&&~~~~\subset \prod_{i\in\Ob\I}\Mor_\C(S(i),X).\nn
\eea
Then by the universal property of the colimit, we have a bijection ~$\Phi:~\Mor_\C\left(\varinjlim S,X\right)~\ra~\txt{CMP}_\C(S,X)$,
{\small\bea
\Phi:~\varinjlim S\sr{h}{\ral}X~~\longmapsto~~\Big(h\circ q_i:S(i)\sr{q_i}{\ral}\varinjlim S\sr{h}{\ral}X\Big)_{i\in\Ob\I}~~=~~\Big(S(i)\sr{h\circ q_i}{\ral}X\Big)_{i\in\Ob\I}.\nn
\eea}It remains to verify that, up to isomorphism (in any relevant/compatible subcategory of Sets), $\txt{CMP}_\C(S,X)$ is the same as ~$\varprojlim \Mor_\C\big(S(i),X\big)$. Consider the cosystem
\[
S^\ast:=\Mor_\C(S,X):\I\ra Sets,~(i\sr{\kappa_{ij}}{\ral}j)\mapsto\Big(\Mor_\C(S(i),X)\sr{S(\kappa_{ij})^\ast}{\lal}\Mor_\C(S(j),X)\Big).
\]

The relevant limit setup in Sets is contained in the following commutative diagram:
\bc
\adjustbox{scale=0.9}{\bt
\Mor_\C(S(i),X)\ar[from=dddrr,bend left,"g^i"]\ar[from=ddrr,bend left,"p^i"']\ar[from=drr,two heads,"p_0^i"']\ar[from=rrrr,bend right=10,"S(\kappa_{ij})^\ast:=(-)\circ S(\kappa_{ij})"'] &&&&
\Mor_\C(S(j),X)\ar[from=dddll,bend right,"g^j"']\ar[from=ddll,bend right,"p^j"]\ar[from=dll,two heads,"p_0^j"]  \\
 && \prod\limits_{i\in\Ob\I}\Mor_\C(S(i),X)\ar[from=d,hook,"q"'] &&\\
 &&\txt{CMP}_\C(S,X)\ar[from=d,dashed,"g"'] &&\\
 && U &&
\et}\ec
 Let $\big(g^i:U\ra \Mor_\C(S(i),X)\big)_{i\in\Ob\I}$ be $S^\ast:=\Mor_\C(S,X)$-compatible, i.e., $g^i=S(\kappa_{ij})^\ast\circ g^j$. Define $g:U\ra\prod_{i\in\Ob\I}\Mor_\C(S(i),X)$ by $g(u):=\big(g^i(u):S(i)\ra X\big)_{i\in\Ob\I}$ for $u\in U$. Then $g(u)$ is $S$-compatible (i.e., lies in $\txt{CMP}_\C(S,X)$) for each $u\in U$ (and so gives a well-defined map $g:U\ra\txt{CMP}_\C(S,X)$), since
\bea
&&g^i(u)=g^j(u)\circ S(\kappa_{ij})=S(\kappa_{ij})^\ast(g^j(u)),~~\txt{for all}~~u\in U,~~\iff~~g^i=S(\kappa_{ij})^\ast\circ g^j.\nn
\eea
With $p^i:=p_0^i\circ q$, where $p_0^i$ are the canonical projections and $q$ the inclusion, it is clear that $p^i\circ g=g^i$, for all $i\in\Ob\I$. Hence, $\txt{CMP}_\C(S,X)\cong\varprojlim \Mor_\C\big(S(i),X\big)$.
\end{proof}

\begin{prp}[\textcolor{blue}{\index{Functoriality of (co)limits}{Colimits and limits are functors}}]\label{LimFunct}
Let $\I,\C$ be categories and $\I_+\supsetneq\I$ a strict containing category of $\I$ (in the sense $\Ob\I_+\supsetneq\Ob\I$). Both (i) the colimit process $\varinjlim:\C^\I\ra\C^{\I_+}$ and (ii) the limit process $\varprojlim:\C^{\I^{op}}\ra\C^{\I^{op}_+}$ are covariant functors. (\blue{footnote}\footnote{Depending on practical convenience, we may sometimes ignore both the transition morphisms and the connecting morphisms associated with the colimit/limit object and so simply view either processes as follows; (i) colimit $\varinjlim:\C^\I\ra\C\cong\C^1$ and (ii) limit $\varprojlim:\C^{\I^{op}}\ra\C\cong\C^1$.})
\end{prp}
\begin{proof}
Consider systems $S,S',S'':\I\ra\C$ and morphisms between them. Then the colimit setup gives a commutative diagram:
\bea
\adjustbox{scale=0.9}{\bt
& S(i)\ar[dddd,near start,bend right=45,"S(\kappa_{ij})"']\ar[rr,"f_i"]\ar[ddr,"q_i"]&& S'(i)\ar[dddd,near start,bend right=45,"S'(\kappa_{ij})"']\ar[rr,"g_i"]\ar[ddr,"q'_i"] && S''(i)\ar[dddd,near start,bend right=45,"S''(\kappa_{ij})"']\ar[ddr,"q''_i"]  &  \\
 & &&  &&   &  \\
 &  &\varinjlim S\ar[rr,dashed,crossing over,"{\exists!~f~=~\varinjlim f_i}"]&& \varinjlim S'\ar[rr,dashed, crossing over,"{\exists!~g~=~\varinjlim g_i}"] && \varinjlim S''  &  \\
 & &&  && & \\
& S(j)\ar[rr,"f_j"]\ar[uur,"q_j"']&& S'(j)\ar[rr,"g_j"]\ar[uur,"q'_j"'] && S''(j)\ar[uur,"q''_j"'] &
\et}\nn
\eea
From the diagram, (i) is an immediate consequence of the universal property of the colimit, which shows
\[
\textstyle\varinjlim id_S=id_{\varinjlim S}~~~~\txt{and}~~~~\varinjlim(f_i\circ g_i)=\left(\varinjlim f_i\right)\circ\left(\varinjlim g_i\right).
\]
Similarly, consider cosystems $S,S',S'':\I\ra\C$ and morphisms between them. Then the limit setup gives a commutative diagram:
\bea
\adjustbox{scale=0.9}{\bt
& S(i)\ar[from=dddd,near end,bend left=45,"S(\kappa_{ij})"]\ar[from=rr,"f^i"']\ar[from=ddr,"p^i"']&& S'(i)\ar[from=dddd,near end,bend left=45,"S'(\kappa_{ij})"]\ar[from=rr,"g^i"']\ar[from=ddr,"p'{}^i"'] && S''(i)\ar[from=dddd,near end,bend left=45,"S''(\kappa_{ij})"]\ar[from=ddr,"p''_i"']  &  \\
 & &&  &&   &  \\
 &  &\varprojlim S\ar[from=rr,dashed,crossing over,"{\exists!~f~=~\varprojlim f^i}"']&& \varprojlim S'\ar[from=rr,dashed, crossing over,"{\exists!~g~=~\varprojlim g^i}"'] && \varprojlim S'' &  \\
 & &&  && & \\
& S(j)\ar[from=rr,"f^j"']\ar[from=uur,"p^j"]&& S'(j)\ar[from=rr,"g^j"']\ar[from=uur,"p'{}^j"] && S''(j)\ar[from=uur,"p''{}^j"]  &
\et}\nn
\eea
From the diagram, (ii) is an immediate consequence of the universal property of the limit, which shows
\[
\textstyle\varprojlim id_S=id_{\varprojlim S}~~~~\txt{and}~~~~\varprojlim(f^i\circ g^i)=\left(\varprojlim f^i\right)\circ\left(\varprojlim g^i\right). \qedhere
\]
\end{proof}

\section{Intersection, Union, Images, Preimages, Equalizers, Universality, Freeness}
Given a trivial system $S^0:\I\ra\C$ (one with no transition morphisms), a \ul{union} (resp. an \ul{intersection}) of $S^0$ is roughly speaking a minimal superobject or containing-object (resp. maximal subobject) of the objects $\{S^0(i)\}_{i\in\I}$ of $S^0$ with respect to the ordering $\subset$ (i.e., subobject or monic morphism relation) on the objects $\Ob\C$ of $\C$. In some applications, unions and intersections will only be considered on families of subobjects of a given object. In such cases, it will be possible to redefine union and intersection say through the coproduct or otherwise.

In the following definitions, if the underlying category $\C$ is ``\ul{exact}'' in the sense that every morphism $f$ factors as $f=m\circ e$ for a monomorphism $m$ and an epimorphism $e$, it might be enough in the universal property to require only a subclass of morphisms (e.g., monics only, epics only, etc, depending on the situation) to factor through the limit object of interest, the reason being that the rest of the morphisms would then automatically factor through the limit object as well. Under such conditions, it might then also be possible for the weak versions of the limits to be unique. (\blue{footnote}\footnote{Recall that a weak limit/colimit (i.e., a weak version of a limit/colimit) is one in which the universal property still requires the existence of factorizations as usual but not their uniqueness anymore.})

\begin{dfn}[\textcolor{blue}{
\index{Monic! limit of a system}{Monic limit of a system},
\index{Intersection of! objects}{{Intersection of objects}}
}]
Let $S:\I\ra\C$ be a system. A \ul{monic limit} of $S$ is a collection of monomorphisms $\{u_i:\varprojlim_{\txt{mo}}S\hookrightarrow S(i)~|~i\in I\}$ such that (i) $S(\kappa_{ij})\circ u_j=u_i$, and (ii) for any collection of morphism $\{x_i:X\ra S(i)~|~i\in I\}$ satisfying $S(\kappa_{ij})\circ x_j=x_i$, there exists a unique morphism $h:\varprojlim_{\txt{mo}}S\ra X$ such that $x_i=u_i\circ h$.

\bea
\adjustbox{scale=0.9}{\bt
S(i)\ar[rr,"S(\kappa_{ij})"] &&  S(j)\\
 & \varprojlim_{\txt{mo}}S\ar[ul,hook,"u_i"]\ar[ur,hook,"u_j"']& \\
 & X\ar[u,dashed,"h"]\ar[uul,bend left,"x_i"]\ar[uur,bend right,"x_j"'] &
\et}~~~~~~~~x_i=u_i\circ h,~~\txt{for all}~~i\in \Ob\I.\nn
\eea
Given a trivial system $S^0:\I\ra\C$, the \ul{intersection} of its objects (i.e., a ``largest common subobject of the objects of $S^0$'') with respect to a collection of monomorphisms {\small $M=\big\{S^0(i)\sr{m_i}{\hookrightarrow} M_{ij}\sr{m_j}{\hookleftarrow}S^0(j)\big\}_{i,j\in\I}\subset\Ob\C$} is a monic limit of the following form:
\[
\adjustbox{scale=0.8}{\bt
M_{ij} &&  S^0(j)\ar[ll,hook,"m_j"']\\
 & \bigcap_M S^0\ar[ul,hook,"u_i"]\ar[ur,hook,"u_j"']& \\
 & X\ar[u,dashed,"h"]\ar[uul,bend left,"x_i"]\ar[uur,bend right,"x_j"'] &
\et}~~~=~~~~
\adjustbox{scale=0.8}{\bt
 & M_{ij} & \\
S^0(i)\ar[ur,hook,"m_i"] &&  S^0(j)\ar[ul,hook,"m_j"']\\
 & \bigcap_M S^0\ar[ul,hook,"u_i"]\ar[ur,hook,"u_j"']& \\
 & X\ar[u,dashed,"h"]\ar[uul,bend left,"x_i"]\ar[uur,bend right,"x_j"'] &
\et}~~~~~~~~\shortstack[l]{$m_iu_i=u_jm_j$, ~for all ~$i,j\in \Ob\I$ \\~~\\ $x_i=u_i\circ h$, ~for all ~$i\in \Ob\I$}.
\]
Whenever $S^0$ is given as a trivial system of subobjects of a given object $C\in\Ob\C$, we will by default set $M_{ij}:=C$ for all $i,j\in\Ob\I$. On the other hand, if the collection of monomorphisms $M$ is not specified, we will by default assume it is constructed from the trivial system $S^0$ itself as follows: For each $i,j\in\I$,
{\small\[\textstyle M_{ij}:=M_{ij}(S^0):=\min_\subset\{C\in\Ob\C:S^0(i),S^0(j)\subset C\},~~\txt{(a smallest common upper bound wrt $\subset$)}\]}with the choice of inclusions $M(S^0):=\big\{S^0(i)\sr{m_i}{\hookrightarrow}M_{ij}(S^0)\sr{m_j}{\hookleftarrow}S^0(j)\big\}$ either arbitrary or determined by the intended application. For brevity, we will define
\[\textstyle S^0(i)\cap S^0(j):=S^0(i)\cap_{M_{ij}(S^0)}S^0(j)~~~~\txt{and}~~~~\bigcap S^0:=\bigcap_{M(S^0)}S^0.\]
\end{dfn}

\begin{note}
Suppose that it is enough to consider only monics in the morphism collections $\{x_i:X\ra S(i)~|~i\in I\}$. Then a monic weak limit (hence \ul{weak intersection}), whenever it exists, is unique. This is because monomorphisms cancel from the left, i.e., if $m$ is monic, then $mf=mf'$ ~$\Ra$~ $f=f'$.
\end{note}

\begin{dfn}[\textcolor{blue}{\index{Pairwise intersection system}{Pairwise intersection system} of a trivial system}]
Let $\C$ be a category, $S^0:\I\ra\C$ a trivial system, and {\small $M=\big\{S^0(i)\sr{m_i}{\hookrightarrow} M_{ij}\sr{m_j}{\hookleftarrow}S^0(j)\big\}_{i,j\in\I}\subset\Ob\C$}. Then, with a category $\J$ defined by $\Ob\J:=\I\sqcup\{(i,j):i\neq j\}_{i,j\in \Ob\I}$ and
\[
\Mor_\J\big((i,j),k\big):=
\left\{
  \begin{array}{ll}
    \{\ld_i\}, & ~\txt{if}~~i=k\neq j \\
    \{\ld_j\}, & ~\txt{if}~~i\neq k=j \\
    \emptyset, & ~\txt{if}~~i\neq k\neq j
  \end{array}
\right\}~~~~\txt{for all}~~i,j,k\in\I~~\txt{such that}~~i\neq j,
\]
the \ul{$M$-pairwise intersection system} of $S^0$ is the extension (extending system) {\small $S^0_M:\J\supset\I\ra\C$} for $S^0$ defined by $S^0_M|_\I:=S^0$ and
\bea
&&S^0_M\Big(\Mor_\J\big((i,j),i\big)\Big):=\Big\{\al_i:S^0(i)\cap_{M_{ij}} S^0(j)\hookrightarrow S^0(i)\Big\},\nn\\
&&S^0_M\Big(\Mor_\J\big((i,j),j\big)\Big):=\Big\{\al_j:S^0(i)\cap_{M_{ij}} S^0(j)\hookrightarrow S^0(j)\Big\},\nn
\eea
where for a fixed pair $(i,j)$, $S^0(i)\cap_{M_{ij}} S^0(j)$ is the monic limit given (in the usual way) by the following diagram:
\[\adjustbox{scale=0.7}{\bt
M_{ij} &&  S^0(j)\ar[ll,hook,"m_j"']\\
 & S^0(i)\cap_{M_{ij}} S^0(j)\ar[ul,hook,"u_i"]\ar[ur,hook,"u_j"']& \\
 & X\ar[u,dashed,"h"]\ar[uul,bend left,"x_i"]\ar[uur,bend right,"x_j"'] &
\et}~~=~~
\adjustbox{scale=0.7}{\bt
 & M_{ij} & \\
S^0(i)\ar[ur,hook,"m_i"] &&  S^0(j)\ar[ul,hook,"m_j"']\\
 & S^0(i)\cap_{M_{ij}} S^0(j)\ar[ul,hook,"u_i"]\ar[ur,hook,"u_j"']& \\
 & X\ar[u,dashed,"h"]\ar[uul,bend left,"x_i"]\ar[uur,bend right,"x_j"'] &
\et}~~~~~~~~\shortstack[l]{$m_iu_i=u_jm_j$, \\~~\\ $x_i=u_i\circ h$, $x_j=u_j\circ h$}.
\]
\end{dfn}

\begin{dfn}[\textcolor{blue}{
\index{Monic! colimit of a system}{Monic colimit of a system},
\index{Union! of objects}{{Union of objects}}
}]
Let $S:\I\ra\C$ be a system. A \ul{monic colimit} of $S$ is a collection of monomorphisms $\{u_i:S(i)\hookrightarrow\varinjlim_{\txt{mo}}S~|~i\in I\}$ such that (i) $u_j\circ S(\kappa_{ij})=u_i$, and (ii) for any collection of morphisms $\{x_i:S(i)\ra X~|~i\in I\}$ satisfying $x_j\circ S(\kappa_{ij})=x_i$, there exists a unique morphism $h:\varinjlim_{\txt{mo}}S\ra X$ such that $x_i=h\circ u_i$.

\bea
\adjustbox{scale=0.9}{\bt
S(i)\ar[dr,hook,"u_i"']\ar[rr,"S(\kappa_{ij})"]\ar[ddr,bend right,"x_i"'] &&  S(j)\ar[dl,hook,"u_j"]\ar[ddl,bend left,"x_j"]\\
 & \varinjlim_{\txt{mo}}S\ar[d,dashed,"h"]& \\
 & X &
\et}~~~~~~~~x_i=h\circ u_i,~~\txt{for all}~~i\in \Ob\I.\nn
\eea
Let $S^0:\I\ra\C$ be a trivial system and  $S^0_M:\J\supset\I\ra\C$ its $M$-pairwise intersection system. The \ul{union} of the objects of $S^0$ (i.e., a ``smallest object containing the objects of $S^0$ as subobjects) with respect to the collection of monomorphisms $M$ is the monic colimit of $S^0_M$: (\blue{footnote}\footnote{It is of course also possible (and perhaps advisable depending on the intended application) to define union with respect to an arbitrary collection of objects (independently of intersection). However, what we have so far is sufficient for our intended applications. Our approach assumes that intersections already exist, and so the formation of unions is then required to be mindful of the existing intersections. By symmetry, the process can be reversed; i.e., form unions first (or assume unions already exist), and then require the formation of intersections to be mindful (in some way) of the existing unions. Also, in some categories (e.g., abelian as we shall see) it is possible to require simultaneous/joint existence of union and intersection in a consistent way.})
\[
\adjustbox{scale=0.7}{\bt
S^0_M((i,j))\ar[dr,hook,"u_{ij}"']\ar[rr,hook,"\al_j"]\ar[ddr,bend right,"x_{ij}"'] &&  S^0_M(j)\ar[dl,hook,"u_j"]\ar[ddl,bend left,"x_j"]\\
 & \varinjlim_{mo}S^0_M\ar[d,dashed,"h"]& \\
 & X &
\et}~~=~~
\adjustbox{scale=0.7}{\bt
 & S^0(i)\cap_{M_{ij}} S^0(j)\ar[dl,hook,"\al_i"']\ar[dr,hook,"\al_j"]& \\
S^0(i)\ar[dr,hook,"u_i"']\ar[ddr,bend right,"x_i"'] &&  S^0(j)\ar[dl,hook,"u_j"]\ar[ddl,bend left,"x_j"]\\
 & \bigcup_M S^0\ar[d,dashed,"h"]& \\
 & X &
\et}~~~~\bt[row sep=small] \bigcup_M S^0~:=~\varinjlim_{mo}S^0_M \\ x_i=h\circ u_i,~~\txt{for all}~~i\in \Ob\I.\et\nn
\]
Just as with intersection, whenever $S^0$ is given as a trivial system of subobjects of a given object $C\in\Ob\C$, we will by default set $M_{ij}:=C$ for all $i,j\in\I$. Also, if the collection of monomorphisms $M$ is not specified, we will by default assume it is $M:=M(S^0)$, where $M(S^0)$ is defined as before. Again for brevity,
\[\textstyle \bigcup S^0:=\bigcup_{M(S^0)}S^0:=\varinjlim_{mo}S^0_{M(S^0)}.\]
\end{dfn}

\begin{note}
Suppose that it is enough to consider only monics in the morphism collections $\{x_i:X\ra S(i)~|~i\in I\}$. Then a monic weak colimit (hence \ul{weak union}) is \ul{not} necessarily unique. This is because a monomorphism might not cancel from the right, i.e., if $m$ is monic, $fm=f'm$ ~$\not\Ra$~ $f=f'$.
\end{note}

\begin{question}[\textcolor{blue}{Union of two sets}]
In the category Sets, consider two sets $A,B$: (i) Is it true that the union $A\cup B$ is the same (up to isomorphism) as the colimit $\varinjlim(A\hookleftarrow A\cap B\hookrightarrow B)$? (ii) Is it true that the intersection $A\cap B$ is the same (up to isomorphism) as the limit $\varprojlim(A\hookrightarrow A\cup B\hookleftarrow B)$?
\end{question}

\begin{dfn}[\textcolor{blue}{\index{Epic! limit of a cosystem}{Epic limit of a cosystem}}]
Let $S:\I\ra\C$ be a cosystem. An \ul{epic limit} of $S$ is a collection of epimorphisms $\{e^i:S(i)\twoheadleftarrow\varprojlim_{\txt{ep}}S~|~i\in I\}$ such that (i) $S(\kappa_{ij})\circ e^j=e^i$, and (ii) for any collection of morphisms $\{x^i:S(i)\la X~|~i\in I\}$ satisfying $S(\kappa_{ij})\circ x^j=x^i$, there exists a unique morphism $h:\varprojlim_{\txt{ep}}S\ra X$ such that $x^i=e^i\circ h$.

\bea
\adjustbox{scale=0.9}{\bt
S(i) &&  S(j)\ar[ll,"S(\kappa_{ij})"']\\
 & \varprojlim_{\txt{ep}}S\ar[ul,two heads,"e^i"]\ar[ur,two heads,"e^j"']& \\
 & X\ar[u,dashed,"h"]\ar[uul,bend left,"x^i"]\ar[uur,bend right,"x^j"'] &
\et}~~~~~~~~x^i=e^i\circ h,~~\txt{for all}~~i\in \Ob\I.\nn
\eea
\end{dfn}

\begin{note}
Suppose that it is enough to consider only epics in the morphism collections $\{x^i:S(i)\la X~|~i\in I\}$. Then an epic weak limit is \ul{not} necessarily unique (since an epimorphism might not cancel from the left, i.e., if $e$ is epic, $ef=ef'$ ~$\not\Ra$~ $f=f'$).
\end{note}

\begin{dfn}[\textcolor{blue}{\index{Epic! colimit of a cosystem}{Epic colimit of a cosystem}}]
Let $S:\I\ra\C$ be a cosystem. An \ul{epic colimit} of $S$ is a collection of epimorphisms $\{e^i:S(i)\twoheadrightarrow\varinjlim_{\txt{ep}}S~|~i\in I\}$ such that (i) $e^j=e^i\circ S(\kappa_{ij})$, and (ii) for any collection of morphisms $\{x^i:S(i)\ra X~|~i\in I\}$ satisfying $x^j=x^i\circ S(\kappa_{ij})$, there exists a unique morphism $h:\varinjlim_{\txt{ep}}S\ra X$ such that $x^i=h\circ e^i$.

\bea
\adjustbox{scale=0.9}{\bt
S(i)\ar[dr,two heads,"e^i"']\ar[from=rr,"S(\kappa_{ij})"']\ar[ddr,bend right,"x^i"'] &&  S(j)\ar[dl,two heads,"e^j"]\ar[ddl,bend left,"x^j"]\\
 & \varinjlim_{\txt{ep}}S\ar[d,dashed,"h"]& \\
 & X &
\et}~~~~~~~~x^i=h\circ e^i,~~\txt{for all}~~i\in \Ob\I.\nn
\eea
\end{dfn}

\begin{note}
Suppose that it is enough to consider only epics in the morphism collections $\{x^i:X\ra S(i)~|~i\in I\}$. Then an epic weak colimit (whenever it exists) is unique. This is because epimorphisms cancel from the right, i.e., if $e$ is epic, then $he=h'e$ ~$\Ra$~ $h=h'$.
\end{note}

The image and coimage of a morphism in the following definition are special cases of the \ul{image of a domain subobject} and \ul{coimage of a codomain quotient object} in the succeeding definitions.
\begin{dfn}[\textcolor{blue}{
\index{Image of! a morphism}{Image of a morphism},
\index{Epic! image}{Epic image},
\index{Coimage of a morphism}{Coimage of a morphism}}]
Let $\C$ be a category and $f:A\ra B$ a morphism in $\C$. An \ul{image} of $f$ is a monomorphism $m_f:\im f\hookrightarrow B$ (where $\im f$ is also denoted by $f(A)$) through which (i) $f$ factors as $f=m_f\al_f$ for a morphism $\al_f:A\ra\im f$, and (ii) for any morphism $m:X\hookrightarrow B$ through which $f$ factors as $f=mf_m$ for a morphism $f_m:A\ra X$, there exists a unique morphism $h:X\ra\im f$ such that diagram (a) below commutes. (\blue{footnote}\footnote{If the morphism $\al_f:A\ra\im f$ is an epimorphism, then we say $\im f$ is an \ul{epic image} of $f$.}).
\bea
\adjustbox{scale=0.9}{
(a)~~\bt
A\ar[ddr,bend right,dotted,"f_m"']\ar[dr,dotted,"\al_f"']\ar[rr,"f"] && B \\
 &\im f\ar[ur,hook,"m_f"'] & \\
 & X\ar[u,dashed,"h"]\ar[uur,bend right,"m"'] &
\et}\hspace{3cm}
(b)~~\bt
A\ar[ddr,bend right,"e"']\ar[dr,two heads,"e_f"']\ar[rr,"f"] && B \\
 &\coim f\ar[d,dashed,"h"]\ar[ur,dotted,"\beta_f"'] & \\
 & X\ar[uur,bend right,dotted,"f_e"'] &
\et\nn
\eea
A \ul{coimage} of $f$ is an epimorphism $e_f:A\twoheadrightarrow\coim f$ through which (i) $f$ factors as $f=\beta_fe_f$ for a morphism $\beta_f:\coim f\ra B$, and (ii) for any morphism $e:A\twoheadrightarrow X$ through which $f$ factors as $f=f_ee$ for a morphism $f_e:X\ra B$, there exists a unique morphism $h:\coim f\ra X$ such that diagram (b) above commutes.
\end{dfn}

The definition of a preimage below (as well as that of the intersection of two objects) belongs to a class of limits called pullback-limits that we will study in some detail later. Likewise, the definition of a postimage below (as well as that of the union of two objects) belongs to a class of colimits called pushout-limits. In other words, it is useful to remember that the preimage and two-object intersection are ``pullbacks'' meanwhile the postimage and two-object union are ``pushouts''. In the special case of an abelian category, all four and more will be intimately related to each other (including an isomorphism between the coimage and the image).

\begin{dfn}[\textcolor{blue}{
\index{Image of! a subobject}{{Image} of a subobject},
\index{Preimage of a subobject}{{Preimage} of a subobject}}]
Let $\C$ be a category, $f:A\ra B$ a morphism, $u:A'\hookrightarrow A$ a subobject, and $v:B'\hookrightarrow B$ a subobject.

An \ul{image} of $A'$ (wrt $f$) is a subobject $u':f(A')\hookrightarrow B$, together with a morphism $f':A'\ra f(A')$ such that (i) $u'f'=fu:A'\ra B$, and (ii) for any morphism $u'':X\ra B$ together with a morphism $f'':A'\ra X$ satisfying $u''f''=fu:A'\ra B$, there exists a unique morphism $h:X\ra f(A')$ such that diagram (a) below commutes: (\blue{footnote}\footnote{If we set $u:=id_A:A\ra A$, then we see that $f(A)=\im f$, i.e., this becomes just the definition of the image of $f$.})
\bea
(a)~~\adjustbox{scale=0.9}{\bt
   & A\ar[dr,"f"] & \\
A'\ar[ur,hook,"u"] & &  B\\
 & f(A')\ar[from=ul,dotted,"f'"]\ar[ur,hook,"u'"]& \\
 & X\ar[u,dashed,"h"]\ar[from=uul,dotted,bend right,"f''"']\ar[uur,bend right,"u''"'] &
\et}\hspace{3cm}(b)~~
\adjustbox{scale=0.9}{\bt
   & B & \\
A\ar[ur,"f"] & &  B'\ar[ul,hook,"v"']\\
 & f^{-1}(B')\ar[ul,hook,"v'"']\ar[ur,dotted,"f'"]& \\
 & X\ar[u,dashed,"h"]\ar[uul,bend left,"v''"]\ar[uur,bend right,dotted,"f''"'] &
\et}\nn
\eea
A \ul{preimage} of $B'$ (wrt $f$) is a subobject $v':f^{-1}(B')\hookrightarrow A$, together with a morphism $f':f^{-1}(B')\ra B'$ such that (i) $vf'=fv':f^{-1}(B')\ra B$, and (ii) for any morphism $v'':X\ra A$ together with a morphism $f'':X\ra B'$ satisfying $vf''=fv'':X\ra B$, there exists a unique morphism $h:X\ra f^{-1}(B')$ such that diagram (b) above commutes: (\blue{footnote}\footnote{This limit is an example of a ``pullback-limit'' (to be discussed later). We will see (after studying properties of pullbacks) that the preimage always exists.})
\end{dfn}

\begin{dfn}[\textcolor{blue}{
\index{Coimage of a quotient object}{{Coimage} of a quotient object},
\index{Postimage of a quotient object}{{Postimage} of a quotient object}}]
Let $\C$ be a category, $f:A\ra B$ a morphism, $u:A\twoheadrightarrow A'$ a quotient object, and $v:B\twoheadrightarrow B'$ a quotient object.

A \ul{coimage} of $B'$ (wrt $f$) is a quotient object $v':A\twoheadrightarrow f^c(B')$, together with a morphism $f':f^c(B')\ra B'$ such that (i) $f'v'=vf:A\ra B'$, and (ii) for any morphism $v'':A\ra X$ together with a morphism $f'':X\ra B'$ satisfying $f''v''=vf:A\ra B'$, there exists a unique morphism $h:f^c(B')\ra X$ such that diagram (a) below commutes: (\blue{footnote}\footnote{If we set $v:=id_B:B\ra B$, then we see that $f^c(B)=\coim f$, i.e., this becomes just the definition of the coimage of $f$.})
\bea
(a)~~
\adjustbox{scale=0.9}{\bt
   & B & \\
A\ar[ur,"f"] & &  B'\ar[from=ul,two heads,"v"]\\
 & f^c(B')\ar[from=ul,two heads,"v'"]\ar[ur,dotted,"f'"]& \\
 & X\ar[from=u,dashed,"h"']\ar[from=uul,bend right,"v''"']\ar[uur,bend right,dotted,"f''"'] &
\et}\hspace{3cm}
(b)~~\adjustbox{scale=0.9}{\bt
   & A\ar[dr,"f"] & \\
A'\ar[from=ur,two heads,"u"'] & &  B\\
 & f^{+1}(A')\ar[from=ul,dotted,"f'"]\ar[from=ur,two heads,"u'"']& \\
 & X\ar[from=u,dashed,"h"']\ar[from=uul,dotted,bend right,"f''"']\ar[from=uur,bend left,"u''"] &
\et}\nn
\eea
A \ul{postimage} of $A'$ (wrt $f$) is a quotient object $u':B\twoheadrightarrow f^{+1}(A')$, together with a morphism $f':A'\ra f^{+1}(A')$ such that (i) $f'u=u'f:A\ra f^{+1}(A')$, and (ii) for any morphism $u'':B\ra X$ together with a morphism $f'':A'\ra X$ satisfying $f''u=u''f:A\ra X$, there exists a unique morphism $h:f^{+1}(A')\ra X$ such that diagram (b) above commutes: (\blue{footnote}\footnote{This limit is an example of a ``pushout-limit'' (to be discussed later). We will see (after studying properties of pushouts) that the postimage always exists.})
\end{dfn}

\begin{questions}
Let $\C$ be a category, $f:A\ra B$, $g:B\ra C$ morphisms, $A'\hookrightarrow A$, $A\twoheadrightarrow A''$, $C'\hookrightarrow C$, and $C\twoheadrightarrow C''$.
(i) When is $(gf)(A')=g(f(A'))$? (ii) When is $(gf)^{-1}(C')=f^{-1}(g^{-1}(C'))$? (iii) When is $(gf)^c(C'')=f^c(g^c(C''))$? (iv) When is $(gf)^{+1}(A'')=g^{+1}(f^{+1}(A''))$? (v) Wherever isomorphism fails above, what partial results hold in its place?
\end{questions}
Some of the above questions will be answered later.

\begin{dfn}[\textcolor{blue}{\index{Equalizer}{Equalizer}, \index{Coequalizer}{Coequalizer}}]
Let $\C$ be a category and $f,g:A\ra B$ morphisms. Then a morphism $e:E\ra A$ is an equalizer of $f,g$ if (i) $fe=ge$, and (ii) for any morphism $e':E'\ra A$ satisfying $fe'=ge'$ there exists a morphism $h:E'\ra E$ such that diagram (b) below commutes.
\bea
\adjustbox{scale=0.9}{
(a)~~\bt
A\ar[rr,"{f,g}"] & &  B\\
 & E\ar[ul,"e"']& \\
 & E'\ar[u,dashed,"h"']\ar[uul,bend left,"e'"'] &
\et\hspace{3cm}
(b)~~\bt
A\ar[rr,"{f,g}"] & &  B\\
 & E\ar[from=ur,"c"']& \\
 & E'\ar[from=u,dashed,"h"]\ar[from=uur,bend left,"c'"] &
\et}\nn
\eea
Similarly, a morphism {\small $c:B\ra E$} is a \ul{coequalizer} of {\small $f,g:A\ra B$} if (i) $cf=cg$, and (ii) for any morphism $c':B\ra E'$ satisfying $c'f=c'g$ there exists a morphism $h:E\ra E'$ such that diagram (b) above commutes.
\end{dfn}

We might not be making explicit use of the following definition later, but it is a way to give meaning to universality.
\begin{dfn}[\textcolor{blue}{
\index{Base subcategory}{Base subcategory},
\index{Universal! (Limit) object}{Universal (Limit) object},
\index{Universal! (Limit) property}{Universal (Limit) property}}]
Let $\C$ be a category. A subcategory $\B\subset\C$ is a \ul{base subcategory} for $\C$ if for every object $C\in\Ob\C$, there exists a system $S:\I\ra\B\subset\C$ such that $\varprojlim S\cong C$ or $\varinjlim S\cong C$ (i.e., every object in $\C$ is the limit or colimit of a system in $\C$ with image in $\B$). (\blue{footnote}\footnote{The objects of the system $S$ lie in $\B$ but the limit or colimit we are referring to is taken in $\C$.}). An object $C\in\Ob\C$ is a \ul{$\B$-universal object} (or \ul{$\B$-limit object}) in $\C$ if $C\not\cong B$ for all $B\in\Ob\B$ (equivalently, $C$ has a \ul{$\B$-universal property} or \ul{$\B$-limit property} in the sense $C$ is the strict/nontrivial limit of a system in $\C$ with image in $\B$). (\blue{footnote}\footnote{If $R$ is a ring and $\C$ the full subcategory of semisimple modules in  $R$-mod, then $\B\subset\C$ might be the full subcategory of simple modules for example.})
\end{dfn}

\begin{dfn}[\textcolor{blue}{\index{Free! object in a category}{Free object in a category}}]
Let $\C$ be a category and $X$ a set. A free object over $X$ (or generated by $X$) in $\C$ is an object $\langle X\rangle\in\Ob\C$ such that given any imbedding (i.e., mapwise-injective functor) $J:\C\hookrightarrow Sets$, the following hold.
\bit
\item[(i)] There exists an injection $i:X\hookrightarrow J\big(\langle X\rangle\big)$.
\item[(ii)] Universal property: For every map $\theta:X\ra J(C)$, $C\in\Ob\C$, there exists a unique morphism $\theta':J(\langle X\rangle)\ra J(C)$ in Sets (or equivalently, $\theta'':=J^{-1}(\theta'):\langle X\rangle\ra C$ in $\C$) such that $\theta=\theta'\circ i$.
\bea\adjustbox{scale=0.9}{\bt
X\ar[d,"\theta"]\ar[rr,hook,"i"] && J\big(\langle X\rangle\big)\ar[dll,dashed,"\exists!~\theta'"]\\
J(C)
\et}~~~~~~~~\theta=\theta'|_X~:=~J(\theta')\circ i.\nn
\eea
\eit
\end{dfn}

(Note: If $\C\subset\txt{Sets}$ already, then we simply drop the imbedding $J$, i.e., $J$ becomes the identity functor.)

%% file: parts/AlgebraCat/AlgebraCatS4.tex
\chapter{Preadditive, Additive, and Abelian Categories}\label{AlgebraCatS4}

\section{Special Objects and the Additivity Property}
\begin{dfn}[\textcolor{blue}{
\index{Initial object}{Initial object},
\index{Final object}{Final object},
\index{Zero! object}{Zero object},
\index{Zero!}{Zero},
\index{Zero! map}{Zero map},
\index{Additive! category}{{Additive category}},
\index{Additive! functor}{{Additive functor}},
\index{Preadditive category}{Pre-additive category}
}]

Consider a category $\A=\big(\A,\Mor_\A,\circ\big)=\big(\Ob\A,\Mor\A,\circ\big)=\big(\Ob\A,\Mor_\A(\Ob\A,\Ob\A),\circ\big)$ with morphisims {\small $\Mor\C =\Mor_\A(\Ob\A,\Ob\A):=\bigsqcup_{A,B\in\Ob\A}\Mor_\C(A,B)$}.

An object $I\in \Ob\A$ is an \ul{initial object} if for every object $A\in\Ob\A$ there exists a unique morphism $q_{IA}\in \Mor_\A(I,A)$. An object $J\in\Ob\A$ is a \ul{final object} if for every object $A\in\Ob\A$ there exists a unique morphism $p_{AJ}\in \Mor_\A(A,J)$. An object $Z\in\Ob\A$ is a \ul{zero object} if it is both an initial object and a final object. If a zero object exists, it is necessarily unique up to isomorphism (by Remark \ref{ZeroObjUnq}), denoted by $0$, and called \ul{zero}. If a zero object $0\in\Ob\A$ exists, then for any $A,B\in\Ob\A$, the (associated) \ul{zero map} $0_{AB}\in \Mor_\A(A,B)$ is the composition $0_{AB}=q_{0B}\circ p_{A0}:A\sr{p_{A0}}{\ral}0\sr{q_{0B}}{\ral}B$.

The category $\A$ is called an \ul{additive category} if the following hold.
\begin{enumerate}[leftmargin=0.9cm]
\item[(1)] \ul{Zero object}: A zero object $0\in \Ob\A$ exists.
\item[(2)] \ul{Pre-additivity}: For every $A,B\in\Ob\A$, $\Mor_\A(A,B)$ is an abelian group $\Big(\Mor_\A(A,B),+,0_{AB}\Big)$, such that composition in $\A$ is distributive over addition as follows: For all $f,f',g,h\in\Mor\A$ as in the diagram
\adjustbox{scale=0.8}{
\bt
D\ar[r,"h"]& A\ar[r,shift left,"f"]\ar[r,shift right,"f'"']& B\ar[r,"g"]& C
\et}(i.e., {\small $f,f'\in \Mor_\A(A,B)$, $h\in \Mor_\A(B,C)$, $g\in \Mor_\A(D,A)$} ), we have
\bea
g\circ(f+f')=g\circ f+g\circ f',~~~~(f+f')\circ h=f\circ h+f'\circ h.\nn
\eea
\item[(3)] \ul{Finite products and finite sums}: For any objects $A,B\in\Ob\A$ (i.e., any trivial system $(A,B):2\ra\A$), there exists a product $A\times B:=\prod(A,B)$ and a coproduct $A\oplus B:=\coprod(A,B)$. (\blue{footnote}\footnote{Because of the pre-additivity property, we have $A\times B\cong A\oplus B$, as shown in Lemma \ref{FinPrCpr}.})
\end{enumerate}
If condition (3) above is dropped (in the sense products/coproducts are not required), then $\A$ is called a \ul{pre-additive category}.

Let $\A,\B$ be pre-additive categories. A functor (resp. cofunctor) $F:\A\ra\B$ is \ul{additive} (\blue{footnote}\footnote{Thinking of a ring $R$ as the pre-additive category $\C(R)$ with a single object $R$ and morphisms the elements $r\in R$ as left multiplication maps $r:R\ra R,~a\mapsto ra$ with composition the usual composition of maps, a ring homomorphism $f:R\ra R'$ is precisely an additive functor $F=F_f:\C(R)\ra\C(R')$. Hence, additive functors between pre-additive categories generalize ring homomorphisms (or equivalently, a pre-additive category is a generalization of a ring).}) if for any $A,A'\in\Ob\A$, the map ~$F:\Mor_\A(A,A')\ra \Mor_\B\big(F(A),F(A')\big)$~ satisfies
\bea
F(f+g)=F(f)+F(g),~~~~\txt{for all}~~~~f,g\in \Mor_\A(A,B).\nn
\eea
(That is a functor (resp. cofunctor) between pre-additive categories $F:\A\ra\B$ is additive if for any $A,A'\in\Ob\A$, the map $F:\Mor_\A(A,A')\ra \Mor_\B\big(F(A),F(A')\big)$ is a \ul{homomorphism of abelian groups} --
see Definition \ref{MorDef1}.)
\end{dfn}

\begin{rmk}[\textcolor{blue}{Uniqueness of the zero object}]\label{ZeroObjUnq}
If a zero object exists, it is unique up to isomorphism. To see this, pick any two zero objects $Z,Z'\in\Ob\A$. Since $Z$ is initial, we have a unique morphism $f=q_{ZZ'}:Z\ra Z'$, i.e., $\Mor_\A(Z,Z')=\{f:Z\ra Z'\}$. Since $Z$ is final we have a unique morphism $g=p_{Z'Z}:Z'\ra Z$, i.e., $\Mor_\A(Z',Z)=\{g:Z'\ra Z\}$. Thus, because $\Mor_\A(Z,Z)=\{id_Z:Z\ra Z\}$ and $\Mor_\A(Z',Z')=\{id_{Z'}:Z'\ra Z'\}$, it follows that $g\circ f=id_Z:Z\sr{f}{\ral}Z'\sr{g}{\ral}Z$ and $f\circ g=id_{Z'}:Z'\sr{g}{\ral}Z\sr{f}{\ral}Z'$.
\end{rmk}

\begin{lmm}[\textcolor{blue}{Equivalence of the finite product and finite coproduct}]\label{FinPrCpr}
Let $\A$ be an additive category and $A,B\in\Ob\A$ any objects in $\A$ (i.e., any trivial system $(A,B):2\ra\A$). Then $A\times B\cong A\oplus B$.
\end{lmm}
\begin{proof}
Let $\{p_1:A\times B\ra A,~p_2:A\times B\ra B\}$ be morphisms defining the product. Then (as the diagram below shows), there exist morphisms $q_1:\ra A\times B$, $q_2:B\ra A\times B$ such that
\bea
p_1q_1=1_A,~~~p_2q_1=0,~~~~~~~p_2q_2=1_B,~~~p_1q_2=0.\nn
\eea

\bea\adjustbox{scale=0.9}{\bt
A\ar[from=dd,"1_A"description]\ar[from=ddrr,bend right,near end,"0"description] &           & B\ar[from=ddll,bend left,near end,"0"description]\ar[from=dd,"1_B"description] \\
  & A\times B\ar[ul,"p_1"description]\ar[ur,"p_2"description] &   \\
A\ar[dddr,bend right,"f_1"description]\ar[ur,dashed,"q_1"description]\ar[dr,"q_1"description] &
& B\ar[dddl,bend left,"f_2"description]\ar[ul,dashed,"q_2"description]\ar[dl,"q_2"description]\\
  & A\times B\ar[from=uu,dashed,"1_{A\times B}=q_1p_1+q_2p_2"description] &  \\
  & &\\
  & X\ar[from=uu,dashed,"h:=f_1p_1+f_2p_2"description] &
\et}~~~~~~~~q_1p_1+q_2p_2=1_{A\times B}~~\txt{(by uniqueness)}.\nn
\eea
Thus, for any morphisms ~$f_1:A\ra X$, $f_2:B\ra X$, ~if we set ~$h:=f_1p_1+f_2p_2:A\times B\ra X$, ~then
\bit[leftmargin=0.9cm]
\item[] using pre-additivity we get (i)~$hq_1=f_1$ ~and ~$hq_2=f_2$,~ and (ii) for any $h'\in\Mor\A$,
\item[] if $h'q_1=f_1$ and $h'q_2=f_2$ then $h=f_1p_1+f_2p_2=h'q_1p_1+h'q_2p_2=h'1_{A\times B}=h'$ (i.e., $h$ is unique).
\eit
This shows ~$A\times B\cong A\oplus B$,~ where the coproduct $A\oplus B$ is given by the induced morphisms $\{q_1\ra A\times B,~q_2:B\ra A\times B\}$ which by construction satisfy
\[
q_1p_1+q_2p_2=1_{A\times B}\eqv 1_{A\oplus B}. \qedhere
\]
\end{proof}

\begin{rmk}[\textcolor{blue}{Non-additivity of Rings}]
The category $Rings$ is not pre-additive, because given rings $R$ and $S$, the hom-set $Hom_{Rings}(R,S)$ is not an abelian group (the sum of two ring homomorhisms is not a ring homomorphism). Moreover, the canonical projection $p:R\times S\ra R,~(r,s)\mapsto r$ is a ring homomorphism but the canonical inclusion $q:R\ra R\times S,~r\mapsto (r,0)$ is not a ring homomorphism, since $1_{R\times S}=(1_R,1_S)\neq(1_R,0)$. Therefore, in general, $R\times S\not\cong R\oplus S$ as rings (even though $R\times S\cong_\Integer R\oplus S$ as abelian groups).
\end{rmk}

\begin{dfn}[\textcolor{blue}{
\index{Infinity! object}{Infinity object},
\index{Infinity! category}{{Infinity category}},
\index{Left-infinity morphism}{Left-infinity morphism},
\index{Right-infinity morphism}{Right-infinity morphism},
\index{Infinity! morphism}{Infinity morphism},
\index{Composition-infinity category}{Composition-infinity category},
}]
Let $\C$ be a category. An object $\infty\in\Ob\C$ is an \ul{infinity object} (making $\C$ an \ul{infinity category}) if the coproduct $\coprod(\infty,C)\cong\infty$, whenever it exists, for all $C\in\Ob\C$. Consider the category of single-morphism systems $\C^{\I_2}$ (where $\I_2:=2=\{0,1\}$ as a poset), whose objects are of course morphisms of $\C$. A morphism $f\in\Mor\C$ is a \ul{left-infinity morphism} (resp. \ul{right-infinity morphism}) if for any $h\in\Mor\C$, we have $fh\cong f$ (resp. $hf\cong f$) in $\C^{\I_2}$ whenever the composition exists. A morphism $f\in\Mor\C$ is an \ul{infinity morphism} (making $\C$ a \ul{composition-infinity category}) if it is both a left-infinity morphism and a right-infinity morphism.
\end{dfn}

In a set, class, or category with more than one binary operation (e.g., a ring or module), an infinity is naturally defined with respect to a given binary operation. For example, in a ring, $0$ is a multiplicative infinity while also an additive identity, and there is no additive infinity in the ring.

\section{Kernels, Cokernels, and the Abelian Property}
\begin{dfn}[\textcolor{blue}{
\index{Right! annihilator of a morphism}{Right annihilator of a morphism},
\index{Kernel}{Kernel},
\index{Left! annihilator of a morphism}{Left annihilator of a morphism},
\index{Cokernel}{Cokernel},
\index{Abelian category}{{Abelian category}}
}]\label{KerAndCoker}
Consider a category $\A$ with a zero object $0\in\Ob\A$, where as usual
\bea
\A=\big(\A,\Mor_\A,\circ\big)=\big(\Ob\A,\Mor\A,\circ\big)=\big(\Ob\A,\Mor_\A(\Ob\A,\Ob\A),\circ\big)\nn
\eea
with morphisms ~{\small $\Mor\A =\Mor_\A(\Ob\A,\Ob\A):=\bigsqcup_{A,B\in\Ob\A}\Mor_\A(A,B)$}.

Let ~$A\sr{f}{\ral}B$~ be a morphism. A morphism $g:A'\ra A$ is a \ul{right annihilator} of $f$ if $f\circ g=0:A'\sr{g}{\ral}A\sr{f}{\ral}B$. A morphism $h:B\ra B'$ is a \ul{left annihilator} of $f$ if $h\circ f=0:A\sr{f}{\ral}B\sr{h}{\ral}B'$.

A \ul{kernel} of $f$ is an object $\ker f\in\Ob\A$ and a morphism $k_f:\ker f\ra A$ with the following properties.
\bit
\item[(i)] \ul{Right annihilator of $f$}:~ $f\circ k_f=0:\ker f\sr{k_f}{\ral}A\sr{f}{\ral}B$
\item[(ii)] \ul{Universal property}: For every morphism $K\sr{k}{\ral}A$ that is a right annihilator of $f$ (i.e., $f\circ k=0$), there exists a unique morphism $K\sr{h}{\ral}\ker f$ such that the following diagram commutes.
\bea\adjustbox{scale=0.9}{\bt
B &&  A\ar[ll,"f"']\\
 & \ker f\ar[ul,"0"]\ar[ur,"k_f"']& \\
 & K\ar[u,dashed,"h"]\ar[uul,bend left,"0"]\ar[uur,bend right,"k"'] &
\et}~~~~~~~~k=k_f\circ h.\nn
\eea
\eit

Note that $\ker f$ is a limit of the directed cosystem $S:\I\ra\A$ given by
{\small\begin{align}
& \Ob\I:=\{1,2\},~~\Mor_\I(1,1):=\{id_1\},~~\Mor_\I(2,2):=\{id_2\},~~\Mor_\I(1,2):=\{\tau\},~~\Mor_\I(2,1):=\emptyset,\nn\\
& S(1):=B,~~~~S(2):=A,~~~~S\Big(1\sr{\tau}{\ral}2\Big):=\Big(S(1)\sr{S(\tau)}{\lal}S(2)\Big):=\Big(B\sr{f}{\lal}A\Big).\nn
\end{align}}

A {cokernel} of $f$ is an object $\coker f\in\Ob\A$ and a morphism $c_f:B\ra\coker f$ with the following properties.
\bit
\item[(i)] \ul{Left annihilator of $f$}:~ $c_f\circ f=0:A\sr{f}{\ral}B\sr{c_f}{\ral}\coker f$.
\item[(ii)] \ul{Universal property}: For every morphism $B\sr{c}{\ral}C$ that is a left annihilator of $f$ (i.e., $c\circ f=0$), there exists a unique morphism ~$\coker f\sr{h}{\ral}C$ ~such that the following diagram commutes.
\bea\adjustbox{scale=0.9}{\bt
A\ar[dr,"0"']\ar[rr,"f"]\ar[ddr,bend right,"0"'] &&  B\ar[dl,"c_f"]\ar[ddl,bend left,"c"]\\
 & \coker f\ar[d,dashed,"h"]& \\
 & C &
\et}~~~~~~~~c=h\circ c_f.\nn
\eea
\eit

Note that $\coker f$ is a colimit of the directed system $S:\I\ra\A$ given by
{\small\begin{align}
& \Ob\I:=\{1,2\},~~\Mor_\I(1,1):=\{id_1\},~~\Mor_\I(2,2):=\{id_2\},~~\Mor_\I(1,2):=\{\tau\},~~\Mor_\I(2,1):=\emptyset,\nn\\
& S(1)=A,~~~~S(2)=B,~~~~S\Big(1\sr{\tau}{\ral}2\Big):=\Big(S(1)\sr{S(\tau)}{\ral}S(2)\Big):=\Big(A\sr{f}{\ral}B\Big).\nn
\end{align}}

The category $\A$ is an \ul{abelian category} if the following hold.
\begin{enumerate}
\item $\A$ is an additive category.

(That is, a zero object $0\in\Ob\A$ exists, morphism classes are abelian groups with addition over which composition of morphisms is distributive, and finite products of objects exist.)
\item Every morphism has a kernel and a cokernel.
\item Every monic morphism is a kernel.
\item Every epic morphism is a cokernel.
\end{enumerate}
\end{dfn}

\begin{notation}[\textcolor{blue}{Combined Kernel and Cokernel setup}]
Let $\A$ be an abelian category and $f:A\ra B$ a morphism in $\A$. (i) The kernel and cokernel of $f$ can be expressed as in the following diagram.
\[
\adjustbox{scale=0.9}{\begin{tikzcd}
 && K\ar[dll,dashed,"\al"']\ar[d,"k"]\ar[drr,"0"] && \\
  \ker f\ar[rr,"k_f"]  && A\ar[rr,"f"]\ar[drr,"0"'] && \ar[d,"c"]B\ar[rr,"c_f"] && \coker f\ar[dll,dashed,"\beta"] \\
    && && C &&
\end{tikzcd}}
\]
{\flushleft (ii)} For brevity in specific situations, we will sometimes write $\ker f$ (resp. $\coker f$) to mean the entire morphism $\ker f\sr{k_f}{\ral}\dom f$ (resp. $\cod f\sr{c_f}{\ral}\coker f$).
\end{notation}

\begin{rmk}
Because kernels and cokernels are limits of systems, if any of them exists it is unique up to isomorphism (as shown in the following proposition).
\end{rmk}

\begin{prp}
In an abelian category the following hold:
\begin{enumerate}[leftmargin=0.9cm]
\item Any kernel is monic, and any two kernels of the same morphism are isomorphic.
\item Any cokernel is epic, and any two cokernels of the same morphism are isomorphic.
\end{enumerate}
\end{prp}
\begin{proof}
Let $\A$ be an abelian category and $f:A\ra B$ a morphism in $\A$.
{\flushleft (1)} \ul{Any kernel of $f$ is monic}: (\blue{footnote}\footnote{In this proof, it is enough to assume $h=0$ due to the abelian property.}). Let $k:K\ra A$ be a kernel of $f$. Suppose $g,h:U\ra K$ satisfy $\bt kg=kh~:~U\ar[r,shift left,"g"]\ar[r,shift right,"h"'] & K\ar[r,"k"] & A\et$. Then because $\bt f(kg)=0~:~U\ar[r,shift left,"g"]\ar[r,shift right,"h"'] & K\ar[r,"k"] & A\ar[r,"f"] & B\et$, there exists a unique morphism $\theta:U\ra K$ such that $k\theta=kg$ ( = $kh$). By the uniqueness of $\theta$, we get $\theta=g=h$. Hence $k:K\ra A$ is monic.

\ul{Any two kernels of $f$ are isomorphic}: (\blue{footnote}\footnote{In this proof, it is enough to assume $h=0$ due to the abelian property.}). Let $K\sr{k}{\ral}A$, $K'\sr{k'}{\ral}A$ be kernels of $f$. Then there exists a unique $\al:K\ra K'$, and a unique $\al':K'\ra K$, such that $k=k'\al$ and $k'=k\al'$. Hence $\al\al'=id_{K'}$ and $\al'\al=id_K$ (since $k$ and $k'$ are monic).
{\flushleft (2)} \ul{Any cokernel of $f$ is epic}: Let $c:B\ra C$ be a cokernel of $f$. Suppose $g,h:B\ra V$ satisfy $\bt gc=hc~:~B\ar[r,"c"] & C\ar[r,shift left,"g"]\ar[r,shift right,"h"'] & V\et$. Then $\bt (gc)f=0~:~A\ar[r,"f"] & B\ar[r,"c"] & C\ar[r,shift left,"g"]\ar[r,shift right,"h"'] & V\et$ implies there exists a unique morphism $\theta:C\ra V$ such that $\theta c=gc$ ( = $hc$). By the uniqueness of $\theta$, we get $\theta=g=h$. Hence $c:B\ra C$ is epic.

\ul{Any two cokernels of $f$ are isomorphic}: Let $c:B\ra C$, $c:B\ra C'$ be cokernels of $f$. Then there exists a unique $\al:C\ra C'$, and a unique $\al':C'\ra C$, such that $c=\al' c'$ and $c'=\al c$. Hence $\al\al'=id_{C'}$ and $\al'\al=id_C$ (since $c$ and $c'$ are epic).
\end{proof}

\begin{crl}
In an abelian category $\A$ the following each hold: A morphism $f$ in $\A$ is
\bit[leftmargin=0.9cm]
\item[(1)] monic $\iff$ a kernel (also of its own cokernel, necessarily).
\item[(2)] epic $\iff$ a cokernel (also of its own kernel, necessarily).
\item[(3)] an isomorphism $\iff$ both monic and epic ($\iff$ $\ker f=0$ ~and ~$\coker f=0$).
\eit
\end{crl}
\begin{proof}
(1),(2) are clear, and so we will prove (3). If $f$ is an isomorphism it is clear that $f$ is both monic and epic. Conversely assume $f$ is both monic and epic (i.e., $\ker f=0$ and $\coker f=0$).

\bc\adjustbox{scale=0.9}{\bt
  &&  && B\ar[dll,dashed,"f''"']\ar[d,"id_B"] &&\\
  \ker f=0\ar[rr,"k_f=0"]  && A\ar[rr,hook,two heads,"f"] && B\ar[rr,"c_f=0"] && \coker f=0 \\
    && A\ar[from=u,"id_A"']\ar[from=urr,dashed,"f'"]&&  &&
\et}\ec
Since ``$f=\coker(k_f)$'' and $0=fk_f=id_Ak_f$, it follows that $id_A=f'f$ (for some $f':B\ra A$). Similarly, since ``$f=\ker(c_f)$'' and $0=c_ff=c_fid_B$, it follows that $id_B=ff''$ (for some $f'':B\ra A$).
\end{proof}

\begin{dfn}[\blue{
\index{Image}{Image},
\index{Coimage}{Coimage}}]
Let $\A$ be an abelian category and $f\in \Mor_\A(A,B)$ a morphism in $\A$. The \ul{image} (resp. \ul{coimage}) of $f$ is a subobject $\im f\hookrightarrow B$ (resp. quotient object $A\twoheadrightarrow\coim f$) defined as follows: With the kernel and cokernel $\bt\ker f\ar[r,hook,"k_f"]& A\ar[r,"f"]& B\ar[r,two heads,"c_f"]& \coker f\et,$
\bea
&&\textstyle\im f := \ker c_f=\ker\big(B\sr{c_f}{\ral}\coker f\big)~=~\txt{``$\ker\coker f$''}~=~kc_f,\nn\\
&&\textstyle\coim f := \coker k_f=\coker\big(\ker f\sr{k_f}{\ral} A\big)~=~\txt{``$\coker\ker f$''}~~\txt{or}~~{A\over\ker f}~=~ck_f.\nn
\eea
\bc\adjustbox{scale=0.9}{\bt
                        &                & \coim f && \im f\ar[dr,hook,"kc_f"] & \\
\ker f\ar[r,hook,"k_f"] & A\ar[ur,two heads,"ck_f"]\ar[rrrr,"f"] &         &&       &  B\ar[r,two heads,"c_f"] & \coker f
\et}\ec
\end{dfn}

\begin{rmk}
Let $\C$ be a category. Observe that any diagram in $\C$ of the form $A\sr{f}{\ral}B\sr{g}{\ral}C$ is equivalent to each of two triangular commutative diagrams as follows:
\bea
\adjustbox{scale=0.7}{\bt A\ar[d,"f"]\\ B\ar[d,"g"]\\ C\et}~~~~=~~~~
\adjustbox{scale=0.7}{\bt
A\ar[dd,"f"]\ar[dr,"g\circ f"] & \\
 & C \\
B\ar[ur,"g"'] &
\et}~~~~=~~~~
\adjustbox{scale=0.7}{\bt
B\ar[dd,"g"]\ar[from=dr,"f'"'] & \\
 & A \\
C\ar[from=ur,"g\circ f"] &
\et}
\nn\eea
Consequently, the explicit specification of a one-morphism (i.e., connected 2-object) colimit/limit as a functor is relatively convenient. The following definition/characterization is based on the above observation, where in each case we either replace $C$ (in the above diagrams) with the colimit ~$\varinjlim\big(A\sr{f}{\ral}B\big)$~ or replace $A$ (in the above diagrams) with the limit ~$\varprojlim\big(B\sr{g}{\ral}C\big)$.
\end{rmk}

\begin{dfn}[\blue{\index{Connected! 2-object limits}{Connected 2-object limits: Explicit form}}]
Let $\C$ be a category, $\C^{2\ast}:=Cat(\Mor\C)\subset\C^2$ the category of morphisms of $\C$, and $\C^{3\ast}$ the category of \ul{connected 3-object diagrams} (i.e., diagrams of the form $A\sr{f}{\ral}B\sr{g}{\ral}C$, $A\sr{f}{\lal}B\sr{g}{\ral}C$, $A\sr{f}{\ral}B\sr{g}{\lal}C$) in $\C$. A \ul{connected 2-object limit} over $\C$ is either a colimit $F:=\varinjlim$ as a functor (functorial process) of the form:
\bea
&&F=\varinjlim:\C^{2\ast}\ra\C^{3\ast},~~
\Big(f\sr{\eta}{\ral}g\Big)=\adjustbox{scale=0.9}{\bt
A\ar[d,"f"]\ar[r,"\eta_1"] & B\ar[d,"g"] \\
A'\ar[r,"\eta_2"] & B'
\et}\longmapsto \Big(F(f)\sr{F(\eta)}{\ral}F(g)\Big)=
\adjustbox{scale=0.7}{\bt
A\ar[d,"f"]\ar[r,"\eta_1"] & B\ar[d,"g"] \\
A'\ar[r,"\eta_2"]\ar[d,"F_f"] & B'\ar[d,"F_g"] \\
F(f)\ar[r,"F(\eta)"] & F(g)
\et}\nn
\eea
or a limit $G:=\varprojlim$ as a functor (functorial process) of the form:
\bea
&&G=\varprojlim:\C^{2\ast}\ra\C^{3\ast},~~
\Big(f\sr{\eta}{\ral}g\Big)=\adjustbox{scale=0.9}{\bt
A\ar[d,"f"]\ar[r,"\eta_1"] & B\ar[d,"g"] \\
A'\ar[r,"\eta_2"] & B'
\et}\longmapsto \Big(G(f)\sr{G(\eta)}{\ral}G(g)\Big)=
\adjustbox{scale=0.7}{\bt
G(f)\ar[r,"G(\eta)"]\ar[d,"G_f"] & G(g)\ar[d,"G_g"] \\
A\ar[d,"f"]\ar[r,"\eta_1"] & B\ar[d,"g"] \\
A'\ar[r,"\eta_2"] & B'
\et}\nn
\eea
In particular, if $\C$ is an abelian category, the kernel and cokernel are connected 2-object limits in the form
\bea
&&\ker:\C^{2\ast}\ra\C^{3\ast},~~
\Big(f\sr{\eta}{\ral}g\Big)~~=~~\adjustbox{scale=0.9}{\bt
A\ar[d,"f"]\ar[r,"\eta_1"] & B\ar[d,"g"] \\
A'\ar[r,"\eta_2"] & B'
\et}~~\mapsto~~
\adjustbox{scale=0.7}{\bt
\ker f\ar[d,hook,"k_f"]\ar[r,"k_{\eta}"] & \ker g\ar[d,hook,"k_g"] \\
A\ar[d,"f"]\ar[r,"\eta_1"] & B\ar[d,"g"] \\
A'\ar[r,"\eta_2"] & B'
\et}\nn\\
&&\coker:\C^{2\ast}\ra\C^{3\ast},~~
\Big(f\sr{\eta}{\ral}g\Big)~~=~~\adjustbox{scale=0.9}{\bt
A\ar[d,"f"]\ar[r,"\eta_1"] & B\ar[d,"g"] \\
A'\ar[r,"\eta_2"] & B'
\et}~~\mapsto~~
\adjustbox{scale=0.7}{\bt
A\ar[d,"f"]\ar[r,"\eta_1"] & B\ar[d,"g"] \\
A'\ar[d,two heads,"c_f"]\ar[r,"\eta_2"] & B'\ar[d,two heads,"c_g"] \\
\coker f\ar[r,"c_{\eta}"] & \coker g
\et}\nn
\eea
The image and coimage can therefore be viewed as (higher order) connected 2-object limits of the form
\bea
&&\im=\ker\circ\coker:\C^{2\ast}\sr{\coker}{\ral}\C^{3\ast}\sr{\ker}{\ral}\C^{4\ast},\nn\\
&&\coim=\coker\circ\ker:\C^{2\ast}\sr{\ker}{\ral}\C^{3\ast}\sr{\coker}{\ral}\C^{4\ast},\nn
\eea
where $\C^{4\ast}$ is the category of (commutative) connected $4$-object diagrams in $\C$.
\end{dfn}

\begin{rmk}[\blue{\index{Connected! 3-object limits}{Connected 3-object limits}}]
Let $\C$ be a category and $\C^{n\ast}$ the category of \ul{(commutative) connected $n$-object diagrams} in $\C$. A \ul{connected 3-object limit} over $\C$ is any colimit $F:=\varinjlim$ or limit $F:=\varprojlim$ as a functor of the form $F:\C^{3\ast}\ra\C^{4\ast}$. Preimages, equalizers/coequalizers (which we have already seen), and pullback/pushout limits (to be seen later) are examples of connected 3-object limits.
\end{rmk}

\section{The Isomorphism Theorem for Abelian Categories}
\begin{lmm}[\blue{
\index{Monic! left-invariance of kernel}{Monic left-invariance of kernel},
\index{Epic! right-invariance of cokernel}{Epic right-invariance of cokernel}
}]
Let $f:A\ra B$ be a morphism in an abelian category $\A$. If $m:B\hookrightarrow B'$ is a monic morphism, and $e:A'\twoheadrightarrow A$ an epic morphism, then from the definitions,
\bea
\ker(m\circ f)\cong\ker f,~~~~coker(f\circ e)\cong\coker f.\nn
\eea
\end{lmm}
\begin{proof}
Observe that $m\circ f\circ k_f=0$ $\iff$ $f\circ k_f=0$, and $c_f\circ f\circ e=0$ $\iff$ $c_f\circ f=0$.
\end{proof}

As we will see in the following theorem (Proposition \ref{UnitaryThm}), if we ignore the fact that a monic epic (i.e., unitary) morphism in an abelian category must be an isomorphism, such a result then becomes an immediate consequence of the resulting decomposition (\ref{UnitExpEq}), as discussed in the proofs of steps (i) and (ii) of the theorem.

\begin{prp}[\textcolor{blue}{\index{Unitarity theorem}{Unitarity theorem in an abelian category}}]\label{UnitaryThm}
Let $\A$ be an abelian category and $f\in \Mor_\A(A,B)$. Then a unitary equivalence $\bt[column sep=small]\coim f\ar[r,hook,two heads]&\im f\et$ exists, and $f$ has the expansion
\bea
\label{UnitExpEq}\bt f=kc_f\circ u_f\circ ck_f:A\ar[r,two heads,"ck_f"] & \coim f\ar[r,hook,two heads,"u_f"] & \im f\ar[r,hook,"kc_f"] & B.\et
\eea
Moreover, we have the following special cases:
\bit
\item[(i)] If $f$ is a monomorphism, then $ck_f=id_A$ (i.e., $\coim f=A$) and $u_f$ is an isomorphism.
\bea
\bt f=kc_f\circ u_f:A\ar[r,"u_f","\cong"'] & \im f\ar[r,hook,"kc_f"] & B.\et\nn
\eea
\item[(ii)] If $f$ is an epimorphism, then $kc_f=id_B$ (i.e., $\im f=A$) and $u_f$ is an isomorphism.
\bea
\bt f=u_f\circ ck_f:A\ar[r,two heads,"ck_f"] & \coim f\ar[r,"u_f","\cong"'] & B.\et\nn
\eea
\item[(iii)] If $f$ is both monic and epic, then it is \ul{clear} that  ~$f=u_f:A=\coim f\ra \im f=B$ ~is an isomorphism.
\item[(iv)] It is \ul{now also clear} that the unitary morphism ~$u_f:\coim f\ra\im f$ ~is an isomorphism.
\eit
\end{prp}
\begin{proof}
We first derive the decomposition (\ref{UnitExpEq}). Consider the following diagram:
\bea\adjustbox{scale=0.9}{\bt
\ker f\ar[rr,hook,"k_f"]&& A\ar[rrrrrr,bend left=50,"f"]\ar[rrrr,two heads,dashed,bend right=35,"\underline{f}"]\ar[rr,two heads,"ck_f"]&&\overbrace{\coker k_f}^{\coim f:=}\ar[rrrr,hook,dashed,bend left=35,"\ol{f}"]\ar[rr,hook,two heads,dashed,"{u}_f"]&&\overbrace{\ker c_f}^{\im f:=}\ar[rr,hook,"kc_f"]&&B\ar[rr,two heads,"c_f"]&&\coker f
\et}\nn
\eea
Since $f\circ k_f=0$, $ck_f\circ k_f=0$, there exists (def. of $\coker k_f$) a unique morphism $\ol{f}:\coim f\ra B$ such that
\bea
&&f=\ol{f}\circ ck_f,~~\Ra~~c_f\circ f=c_f\circ\ol{f}\circ ck_f=0,~~\Ra~~c_f\circ\bar{f}=0,~c_f\circ kc_f=0,\nn\\
&&~~\Ra~~\txt{a unique morphism ~$\coim f\sr{\ol{{u}}}{\ral}\im f$~ exists such that}~~\ol{f}=kc_f\circ\ol{{u}}.\nn
\eea
Similarly, because $c_f\circ f=0$, $c_f\circ kc_f=0$, there exists (def. of $\ker c_f$) a unique morphism $\underline{f}:A\ra\im f$ such that
\bea
&&f=kc_f\circ\underline{f},~~\Ra~~f\circ k_f=kc_f\circ\underline{f}\circ k_f=0,~~\Ra~~\underline{f}\circ k_f=0,~ck_f\circ k_f=0,\nn\\
&&~~\Ra~~\txt{a unique morphism ~$\coim f\sr{\underline{{u}}}{\ral}\im f$~ exists such that}~~\underline{f}=\underline{{u}}\circ ck_f.\nn
\eea
The two possible expansions of $f$ in terms of $\ol{{u}}$ and $\underline{{u}}$, namely
\bea
f=\ol{f}\circ ck_f=kc_f\circ\ol{{u}}\circ ck_f~~~~\txt{and}~~~~f=kc_f\circ\underline{f}=kc_f\circ\underline{{u}}\circ ck_f,\nn
\eea
imply ~$kc_f\circ\ol{{u}}\circ ck_f=kc_f\circ\underline{{u}}\circ ck_f$,~ and so ~$\ol{{u}}=\underline{{u}}$. It remains to show that $u_f:=\underline{u}=\ol{u}$ is unitary. Let
\[
K_f:=\ker f,~~~~C_f:=\coker f,~~~~CK_f:=\coim f,~~~~KC_f:=\im f.
\]
{\flushleft $\bullet$} To show $u_f$ is monic, consider any morphism $h:H\ra\coim f$ such that $u_f\circ h=0$ in the setup:

\bc\adjustbox{scale=0.9}{\bt
         &              &   K_{c_h\circ ck_f}\ar[ddl,dashed,"k_{f,h}'"']\ar[dd,hook,"\substack{k_{f,h}:=\\k_{c_h\circ ck_f}}"]    & &                 &                              &   \\
         &              &       & &                  &                              &   \\
  & K_f\ar[r,hook,"k_f"] & A\ar[dd,two heads,"ck_f"']\ar[ddrr,dashed,"\underline{f}"pos=0.3]\ar[rr,"f"] & & B\ar[r,two heads,"c_f"] & C_f &  \\
         &              &       & &                  &                              &   \\
         &     H\ar[r,"h"]         & \overbrace{CK_f}^{\coim f}\ar[dd,two heads,"c_h"']\ar[uurr,dashed,"\ol{f}"pos=0.7]\ar[rr,hook,two heads,"u_f",dashed] &  & \overbrace{KC_f}^{\im f}\ar[uu,hook,"kc_f"']                   &                              & \\
         &              &       & &                  &                              &   \\
         &              &   C_h\ar[uurr,dashed,bend right," u_f'"']\ar[uu,dashed, bend right=40,"c_h'"']    & &                  &                              &   \\
\et}\ec
The equalities $c_h\circ h=0$ and $u_f\circ h=0$ together imply the existence of $u_f'$. Thus, $f\circ k_f=0$ and
\bea
&& f\circ k_{f,h}=kc_f\circ u_f\circ ck_f\circ k_{f,h}=kc_f\circ u_f'\circ \overbrace{c_h\circ ck_f\circ k_{f,h}}^{=0}=0\nn
\eea
together imply the existence of $k_{f,h}'$. Next,
\bea
&& ck_f\circ k_{f,h}=ck_f\circ k_f\circ k_{f,h}'=0\circ k_{k,h}'=0,~~\txt{along with}~~(c_h\circ ck_f)\circ k_{f,h}=0\nn\\
&&~~\Ra~~ck_f=c_h'\circ(c_h\circ ck_f),~~\Ra~~c_h'\circ c_h=id_{\coim f},\nn\\
&&~~\Ra~~h=id_{\coim f}\circ h=(c_h'\circ c_h)\circ h=c_h'\circ 0=0.\nn
\eea

{\flushleft $\bullet$} Similarly, to show $u_f$ is epic, consider any morphism $h:\im f\ra H$ such that $h\circ u_f=0$ as in the setup:
\bc\adjustbox{scale=0.9}{\bt
         &              &       & &           C_{kc_f\circ k_h}       &                              &   \\
         &              &       & &                  &                              &   \\
  & K_f\ar[r,hook,"k_f"] & A\ar[dd,two heads,"ck_f"']\ar[ddrr,dashed,"\underline{f}"pos=0.3]\ar[rr,"f"] & & B\ar[r,two heads,"c_f"]\ar[uu,two heads,"\substack{c_{f,h}:=\\c_{kc_f\circ k_h}}"] & C_f\ar[uul,dashed,"c_{f,h}'"'] &  \\
         &              &       & &                  &                              &   \\
         &              & \overbrace{CK_f}^{\coim f}\ar[uurr,dashed,"\ol{f}"pos=0.7]\ar[rr,hook,two heads,"u_f",dashed]\ar[ddrr,dashed, bend right,"u_f'"'] &  & \overbrace{KC_f}^{\im f}\ar[dd,bend right=40,dashed,"k_h'"']\ar[uu,hook,"kc_f"']\ar[r,"h"]                   &H                              & \\
         &              &       & &                  &                              &   \\
         &              &       & &        K_h\ar[uu,hook,"k_h"']          &                              &   \\
\et}\ec
As before, the equalities $h\circ u_f=0$ and $h\circ k_h=0$ together imply the existence of $u_f'$. Thus, $c_f\circ f=0$ and
\bea
&& c_{f,h}\circ f=c_{f,h}\circ kc_f\circ u_f\circ ck_f=\overbrace{c_{f,h}\circ kc_f\circ k_h}^{=0}\circ u_f'\circ ck_f=0\nn
\eea
together imply the existence of $c_{f,h}'$. Next,
\bea
&& c_{f,h}\circ kc_f=c_{f,h}'\circ c_f\circ kc_f=c_{f,h}'\circ 0=0,~~\txt{along with}~~c_{f,h}\circ(kc_f\circ k_h)=0\nn\\
&&~~\Ra~~kc_f=(kc_f\circ k_h)\circ k_h'=kc_f\circ(k_h\circ k_h'),~~\Ra~~k_h\circ k_h'=id_{\im f},\nn\\
&&~~\Ra~~h=h\circ id_{\im f}=h\circ k_h\circ k_h'=0\circ k_h'=0.\nn
\eea
Hence, $u_f$ is both monic and epic (which means $\ol{f}=kc_f\circ u_f$ is monic and $\underline{f}=u_f\circ ck_f$ is epic).

\bit[leftmargin=0.9cm]
\item[(i)] If $f$ is monic (so $k_f=0$), then $ck_f=c_0=id_A$ (i.e., $\coim f=A$), and so $f=kc_f\circ u_f$.
\bea\adjustbox{scale=0.9}{\bt
\overbrace{\ker f}^{=0}\ar[rr,hook,"k_f=0"]&& A\ar[rrrrrr,bend left=50,"f"]\ar[rrrr,two heads,dashed,bend right=35,"\underline{f}=u_f","\cong"']\ar[rr,two heads,"ck_f=id_A"]&&\overbrace{\coker k_f}^{A=\coim f:=}\ar[rrrr,hook,dashed,bend left=35,"\ol{f}=f"]\ar[rr,hook,two heads,dashed,"{u}_f","\cong"']&&\overbrace{\ker c_f}^{\im f:=}\ar[rr,hook,"kc_f"]&&B\ar[rr,two heads,"c_f"]&&\coker f
\et}\nn
\eea
Thus, if $f=kc_f\circ u_f:A\sr{u_f}{\ral}\im f\sr{kc_f}{\ral} B$ is a monic morphism in an abelian category, then $u_f$ is an isomorphism (since every monic morphism is a kernel of its cokernel).

\item[(ii)] Similarly, if $f$ is epic (so $c_f=0$), then $kc_f=k_0=id_B$ (i.e., $\im f=B$), and so $f=u_f\circ ck_f$.
\bea\adjustbox{scale=0.9}{\bt
\ker f\ar[rr,hook,"k_f"]&& A\ar[rrrrrr,bend left=50,"f"]\ar[rrrr,two heads,dashed,bend right=35,"\underline{f}=f"]\ar[rr,two heads,"ck_f"]&&\overbrace{\coker k_f}^{\coim f:=}\ar[rrrr,hook,dashed,bend left=35,"\ol{f}=u_f","\cong"']\ar[rr,hook,two heads,dashed,"{u}_f","\cong"']&&\overbrace{\ker c_f}^{B=\im f:=}\ar[rr,hook,"kc_f=id_B"]&&B\ar[rr,two heads,"c_f=0"]&&\overbrace{\coker f}^{=0}
\et}\nn
\eea
Thus, if $f=u_f\circ ck_f:A\sr{ck_f}{\ral}\coim f\sr{u_f}{\ral} B$ is an epic morphism in an abelian category, then $u_f$ is an isomorphism (since every epic morphism is a cokernel of its kernel). \qedhere
\eit
\end{proof}

\begin{rmks}[\textcolor{blue}{
\index{KC-CK identities}{KC-CK identities},
\index{SES-decomposition of a morphism}{SES-decomposition of a morphism},
\index{Exactness! morphisms}{Exactness morphisms},
\index{Exactness! monomorphism}{Exactness monomorphism},
\index{Exactness! epimorphism}{Exactness epimorphism},
\index{Decomposition of a complex}{Decomposition of a complex}}]\label{SeqFactRmks}~
Let $f:A\ra B$ be a morphism in an abelian category $\A$. Consider the decomposition (\ref{UnitExpEq}) of $f$ derived in Proposition \ref{UnitaryThm}. (\blue{footnote}\footnote{In the discussion that follows, like in several others, it is not essential to know the unitary morphism $u_f$ is an isomorphism.})
\begin{enumerate}[leftmargin=0.9cm]
\item[(1)] \ul{KC-CK identities}: Recall that if $m$ is a monic morphism, and $e$ an epic morphism, then $\ker(m\circ f)\cong\ker f$ and $coker(f\circ e)\cong\coker f$. In particular, we get the following (i.e., ~$k_f=kck_f$~ and ~$c_f=ckc_f$~):
\bea
\ker f=\ker(kc_f\circ u_f\circ ck_f)\cong\ker(ck_f),~~~~\coker f=\coker(kc_f\circ u_f\circ ck_f)\cong\coker(kc_f).\nn
\eea
This shows the \emph{complex} (\blue{footnote}\footnote{To be defined later, a linear system in an abelian category is called a \emph{complex} if the composition of any two consecutive transition morphisms of the system is zero.}) ~$0\ra A\sr{f}{\ral}B\ra 0$ ~\ul{decomposes} into SES's (\blue{footnote}\footnote{To be defined later, ``SES'' means ``\emph{short exact sequence}'', which is a \emph{complex} of the form $0\ra A\sr{f}{\ral}B\sr{g}{\ral}C\ra 0$, which precisely means (i) $f$ is monic, (ii) $\im f\cong\ker g$, and (iii) $g$ is epic.}) as follows:
\bea\adjustbox{scale=0.9}{
\bt[row sep=small]
 & & & 0\ar[d] & & & \\
0\ar[r] & \ker f\ar[r,hook,"k_f"] & A\ar[r,two heads,"ck_f"] & \coim f\ar[dd,hook,two heads,"u_f"]\ar[r] & 0 & &\\
 & & &  & & & \\
 & & 0\ar[r] & \im f\ar[d]\ar[r,hook,"kc_f"] & B\ar[r,two heads,"c_f"] & \coker f\ar[r] & 0\\
  & & & 0 & & &
\et}\nn
\eea
\item[(2)]\ul{Exactness morphisms}: Given any two morphisms $\bt[column sep=small] A\ar[r,"f"] & B\ar[r,"g"] & C\et$, we have $g\circ f=0$ $\iff$ there exists a \ul{unique monomorphism} $\bt\im f\ar[r,hook,"m_{f|g}"] & \ker g\et$ that factors $\bt \im f\ar[r,hook,"kc_f"] & B\et$ as follows:
    \bea
    \bt kc_f=k_g\circ m_{f|g}:\im f\ar[r,hook,"m_{f|g}"] & \ker g\ar[r,hook,"k_g"] & B.\et\nn
    \eea
This follows because $0=g\circ f=g\circ(kc_f\circ u_f\circ ck_f)$ $\iff$ $g\circ kc_f=0$ (along with $g\circ k_g=0$).
\bea\adjustbox{scale=0.8}{\bt
                         &                           & \coim f\ar[r,hook,two heads,"u_f"] & \im f\ar[dd,hook,dashed,near start,"m_{f|g}"']\ar[dr,hook,"kc_f"] &                           & \coim g\ar[r,hook,two heads,"u_g"] & \im g\ar[dr,hook,"kc_g"] &   &  \\
                         & A\ar[ur,two heads,"ck_f"]\ar[rrr,crossing over,near start,"f"] &                                    &                          & B\ar[dr,two heads,"c_f"]\ar[ur,two heads,"ck_g"]\ar[rrr,"g"] &                                    &       & C\ar[dr,two heads,"c_g"] &  \\
\ker f\ar[ur,hook,"k_f"] &                           &                                    & \ker g\ar[ur,hook,"k_g"] &                           & \coker f \ar[uu,crossing over,two heads,dashed,near start,"e_{f|g}"']                          &       &   & \coker g
\et}\nn
\eea

Similarly, $g\circ f=0$ $\iff$ there exists a \ul{unique epimorphism} $\bt\coker f\ar[r,two heads,"e_{f|g}"] & \coim g\et$ that factors $\bt B\ar[r,two heads,"ck_g"] & \coim g\et$ as follows:
    \bea
    \bt ck_g=e_{f|g}\circ c_f:B\ar[r,two heads,"c_f"] & \coker f\ar[r,two heads,"e_{f|g}"] & \coim g.\et\nn
    \eea
This follows because $0=g\circ f=(kc_g\circ u_g\circ ck_g)\circ f$ $\iff$ $ck_g\circ f=0$ (along with $c_f\circ f=0$). We will call $m_{f|g}$ the \ul{exactness monomorphism} (resp. $e_{f|g}$ the \ul{exactness epimorphism}) of $A\sr{f}{\ral}B\sr{g}{\ral}C$ (when $g\circ f=0$).

\item[(3)] \ul{ Decomposition of a complex}: Given any three morphisms $\bt[column sep=small] A\ar[r,"f"] & B\ar[r,"g"] & C\ar[r,"h"]& D\et$ satisfying $g\circ f=0$ and $h\circ g=0$, the middle morphism $g$ decomposes as follows:
\bc\adjustbox{scale=0.8}{\bt
A\ar[r,"f"] & B\ar[rrrr,bend left=40,"g"]\ar[r,two heads,"c_f"]\ar[rr,bend right,two heads,"ck_g"'] & \coker f\ar[r,two heads,"e_{f|g}"] & \coim g\cong\im g\ar[rr,bend right,hook,"kc_g"']\ar[r,hook,"m_{g|h}"] & \ker h\ar[r,hook,"k_h"] & C\ar[r,"h"]& D
\et}\ec
Moreover, by part (1) above,
\[
\coker f\cong\coim g~~\iff~~\im f\cong\ker g.
\]
\end{enumerate}
\end{rmks}

\section{The Abelian Structure of $R$-mod: The Isomorphism Theorem}
\begin{rmk}[\textcolor{blue}{Abelian structure of $R$-mod}]~
\begin{enumerate}[leftmargin=0.9cm]
\item[(1)] In R-mod, the \ul{kernel} of $f:A\ra B$ is given by ~$k_f:\ker f\cong f^{-1}(0)\sr{i}{\hookrightarrow}A,~a\mapsto a$, where $f^{-1}(0):=\{a\in A:f(a)=0\}$.
\begin{proof}
Define $k_f:=i:=i_{f^{-1}(0)}:f^{-1}(0)\ra A,~a\mapsto a$ (i.e., the inclusion). If $k:K\ra A$ satisfies $f\circ k=0$, then $k(K):=\{k(x):x\in K\}\subset f^{-1}(0)$, and so we can set $h:=k$ to satisfy $k_f\circ h=i\circ k=k$.
\bea\adjustbox{scale=0.9}{\bt
B &&  A\ar[ll,"f"']\\
 & f^{-1}(0)\ar[ul,"0"]\ar[ur,hook,near start,"k_f:=i"]& \\
 & K\ar[u,dashed,"h:=k"']\ar[uul,bend left,"0"]\ar[uur,bend right,"k"'] &
\et}~~~~~~~~k_f\circ h=i\circ k=k.\nn
\eea
\end{proof}
\item[(2)] In R-mod, the \ul{cokernel} of $f:A\ra B$ is given by ~$c_f:B\sr{\pi}{\twoheadrightarrow}{B\over f(A)}\cong \coker f,~b\mapsto b+f(A)$, ~where ~$f(A):=\{f(a):a\in A\}$.
\begin{proof}
Define $c_f:=\pi:B\ra {B\over f(A)},~b\mapsto b+f(A)$ (i.e., the natural map). If $c:B\ra C$ satisfies $c\circ f=0$, then we get a well-defined map $h:{B\over f(A)},~b+f(A)\mapsto c(b)$, and so $h\circ c_f=h\circ\pi=c$.
\bea\adjustbox{scale=0.9}{\bt
A\ar[dr,"0"']\ar[rr,"f"]\ar[ddr,bend right,"0"'] &&  B\ar[dl,two heads,near end,"c_f:=\pi"']\ar[ddl,bend left,"c"]\\
 & {B\over f(A)}\ar[d,dashed,"h"]& \\
 & C &
\et}~~~~~~~~h\circ c_f=h\circ\pi=c.\nn
\eea
\end{proof}

\item[(3)] In R-mod, the \ul{image} of $f:A\ra B$ is given by ~$kc_f:\im f\cong f(A)\sr{i}{\hookrightarrow}B,~f(a)\mapsto f(a)$, where $f(A):=\{f(a):\in A\}\subset B$.
\begin{proof}
{\small $\im f:=\ker\left(A\twoheadrightarrow\coker f\right)\cong\ker\left(A\sr{\pi}{\twoheadrightarrow}{B\over f(A)},~a\mapsto a+f(A)\right)\sr{(1)}{\cong}\pi^{-1}(0)=f(A)$}.
\end{proof}

\item[(4)] In R-mod, the \ul{coimage} of $f:A\ra B$ is given by ~$ck_f:A\sr{\pi}{\twoheadrightarrow}{A\over f^{-1}(0)}\cong \coim f,~a\mapsto a+f^{-1}(0)$.
\begin{proof}
{\small $\coim f:=\coker\left(\ker f\hookrightarrow A\right)\cong\coker\left(f^{-1}(0)\sr{i}{\hookrightarrow}A,~a\mapsto a\right)\sr{(2)}{\cong}{\cod i\over i\big(f^{-1}(0)\big)}={A\over f^{-1}(0)}$}.
\end{proof}

\item[(5)] In R-mod, the \ul{exactness monomorphism} of $A\sr{f}{\ral}B\sr{g}{\ral}C$ (when $g\circ f=0$) is given by inclusion as
\bea
m_{f|g}:\im f\cong f(A)\sr{i}{\hookrightarrow}g^{-1}(0)\cong\ker g,~~f(a)\mapsto f(a),\nn
\eea
and so the associated \ul{exactness epimorphism} is given by a well defined map as follows:
\bea
\textstyle e_{f|g}:\coker f\cong{B\over f(A)}\sr{}{\twoheadrightarrow}{B\over g^{-1}(0)}\cong\coim g,~~b+f(A)\mapsto b+g^{-1}(0).\nn
\eea

\begin{proof}
This follows from the uniqueness (up to isomorphism) of the exactness monomorphism and the fact that both $\im f$ and $\ker g$ are given by inclusions into $B$.
\end{proof}
\end{enumerate}
\end{rmk}

\begin{crl}[\textcolor{blue}{Unitarity theorem in $R$-mod}]\label{IsomThmI}
In $R$-mod, for any morphism $f:A\ra B$, the unitary morphism \bt[column sep=small] u_f:\coim f\ar[r,hook,two heads] & \im f\et in the decomposition (\ref{UnitExpEq}) in Proposition \ref{UnitaryThm} immediately gives an isomorphism ~$u_f:{A\over f^{-1}(0)}\ra f(A),~a+f^{-1}(0)\mapsto f(a)$~ through which $f$ uniquely factors as
\bea
\textstyle\bt[column sep=small] f=\pi\circ u_f:A\ar[r,two heads,"\pi"] & {A\over f^{-1}(0)}\ar[r,"u_f"]& f(A)\ar[r,hook] & B.\et\nn
\eea
\end{crl}

\section{Additivity of the Colimit/Limit Functors}
\begin{prp}[\textcolor{blue}{\index{Additivity of! limits}{Limits as additive functors}}]\label{LimAddFunct}
Let $\I$ be a category and $\A$ an abelian category. Both the colimit $\varinjlim:\A^\I\ra\A$ and the limit $\varprojlim:\A^{\I^{op}}\ra\A$ are additive functors.
\end{prp}
\begin{proof}
{\flushleft \ul{Additivity of the colimit}}: Let $S,S':\I\ra\A$ be systems, and $f_i,h_i:S(i)\ra S'(i)$ morphisms of systems. (Ignore all other information such as $g_i:S'(i)\ra S''(i)$ in the diagram below):
\bea\adjustbox{scale=0.9}{\bt
 0\ar[r]& S(i)\ar[dddd,near start,bend right=45,"S(\kappa_{ij})"']\ar[rr,"{f_i,~h_i,~f_i+h_i}"]\ar[ddr,"q_i"]&& S'(i)\ar[dddd,near start,bend right=45,"S'(\kappa_{ij})"']\ar[rr,"g_i"]\ar[ddr,"q'_i"] && S''(i)\ar[dddd,near start,bend right=45,"S''(\kappa_{ij})"']\ar[ddr,"q''_i"] \ar[r] & 0 \\
 & &&  &&   &  \\
 & 0\ar[r,crossing over] &\varinjlim S\ar[rr,dashed,crossing over,"{\exists!~f~=~\varinjlim f_i}"]&& \varinjlim S'\ar[rr,dashed, crossing over,"{\exists!~g~=~\varinjlim g_i}"] && \varinjlim S'' \ar[r,crossing over] & 0 \\
 & &&  && & \\
0\ar[r]& S(j)\ar[rr,"{f_j,~h_j,~f_j+h_j}"]\ar[uur,"q_j"']&& S'(j)\ar[rr,"g_j"]\ar[uur,"q'_j"'] && S''(j)\ar[uur,"q''_j"'] \ar[r] & 0
\et}\nn
\eea
As the above diagram shows, there exist unique morphisms $\varinjlim(f_i)$, $\varinjlim(h_i)$, $\varinjlim(f_i+h_i)$ such that
\[
q_i'\circ f_i=\varinjlim(f_i)\circ q_i,~~~~q_i'\circ h_i=\varinjlim(h_i)\circ q_i,~~~~q_i'\circ(f_i+h_i)=\varinjlim(f_i+h_i)\circ q_i.
\]
Therefore, using the pre-additivity of $\A$, we see that
\begin{align}
&\varinjlim(f_i+h_i)\circ q_i=q_i'\circ(f_i+h_i)=q_i'\circ f_i+q_i'\circ h_i=\varinjlim(f_i)\circ q_i+\varinjlim(h_i)\circ q_i=\big(\varinjlim(f_i)+\varinjlim(h_i)\big)\circ q_i,\nn\\
&~~\Ra~~\varinjlim(f_i+h_i)=\varinjlim(f_i)+\varinjlim(h_i).\nn
\end{align}

{\flushleft \ul{Additivity of the limit}}: Let $S,S':\I\ra\A$ be cosystems, and $f^i,h^i:S(i)\la S'(i)$ morphisms of systems. (Ignore all other information such as $g^i:S'(i)\la S''(i)$ in the diagram below):
\bea\adjustbox{scale=0.9}{\bt
 0\ar[from=r]& S(i)\ar[from=dddd,near end,bend left=45,"S(\kappa_{ij})"]\ar[from=rr,"{f^i,~h^i,~f^i+h^i}"']\ar[from=ddr,"p^i"']&& S'(i)\ar[from=dddd,near end,bend left=45,"S'(\kappa_{ij})"]\ar[from=rr,"g^i"']\ar[from=ddr,"p'{}^i"'] && S''(i)\ar[from=dddd,near end,bend left=45,"S''(\kappa_{ij})"]\ar[from=ddr,"p''_i"'] \ar[from=r] & 0 \\
 & &&  &&   &  \\
 & 0\ar[from=r,crossing over] &\varprojlim S\ar[from=rr,dashed,crossing over,"{\exists!~f~=~\varprojlim f^i}"']&& \varprojlim S'\ar[from=rr,dashed, crossing over,"{\exists!~g~=~\varprojlim g^i}"'] && \varprojlim S'' \ar[from=r,crossing over] & 0 \\
 & &&  && & \\
0\ar[from=r]& S(j)\ar[from=rr,"{f^j,~h^j,~f^j+h^j}"']\ar[from=uur,"p^j"]&& S'(j)\ar[from=rr,"g^j"']\ar[from=uur,"p'{}^j"] && S''(j)\ar[from=uur,"p''{}^j"] \ar[from=r] & 0
\et}\nn
\eea
Then following the same steps as for the colimit, we see from the above diagram that
\[
\varprojlim(f^i+h^i)=\varprojlim(f^i)+\varprojlim(h^i). \qedhere
\]
\end{proof}

%% file: parts/AlgebraCat/AlgebraCatS5.tex
\chapter{Complexes, Homotopy, Homology, and Cohomology}\label{AlgebraCatS5}

\section{Types of complexes, Categories of complexes, Exact and Split sequences}
As usual, many of the terms in the following definition make sense in any additive category, and not just for the abelian category under which they are defined.

\begin{dfn}[\textcolor{blue}{
\index{Complex}{Complex in an abelian category},
\index{Cocomplex}{Cocomplex},
\index{Chain! complex}{Chain complex},
\index{Cochain complex}{Cochain complex},
\index{Boundary morphism (Differential)}{Boundary morphism (Differential)},
\index{Coboundary morphism (Codifferential)}{Coboundary morphism (Codifferential)},
\index{Morphism of! (co)complexes}{Morphism of (co)complexes},
\index{Morphism of! (co)chain complexes ((co)chain morphism)}{Morphism of (co)chain complexes ((co)chain morphism)},
\index{Category of! (co)complexes}{{Category of (co)complexes}},
\index{Category of! (co)chain complexes}{{Category of (co)chain complexes}},
\index{Exact! (co)chain complex (Exact sequence)}{Exact (co)chain complex (Exact sequence)},
\index{Short exact sequence (SES)}{Short exact sequence (SES)},
\index{Quotient! of two objects}{Quotient of two objects},
\index{Split! sequence}{Split sequence},
\index{Multicomplex}{Multicomplex},
\index{Morphism of! multicomplexes}{Morphism of multicomplexes},
\index{Category of! multicomplexes}{Category of multicomplexes},
\index{Bicomplex}{{Bicomplex}},
\index{Total! chain complex of a bicomplex}{Total chain complex of a bicomplex},
\index{Total! functor}{Total functor}}]\label{HomCohomDf0}

Let $\I$ be a linear set (i.e., linearly-ordered set) category and $\A$ an abelian category (or just an additive category where sensible). An $\I$-system (resp. $\I$-cosystem) in $\A$ explicitly written as
\bea
&&C:\I\ra\A,~~~~i\sr{\kappa_{ij}}{\ral}j~~\mapsto~~C_i\sr{\del^C_{ij}}{\ral}C_j~:=~C(i)\sr{C(\kappa_{ij})}{\ral}C(j),\nn\\
&&\Big(\txt{resp}.~~C:\I\ra\A,~~~~i\sr{\kappa_{ij}}{\ral}j~~\mapsto~~C_i\sr{\del^{ji}_C}{\lal}C_j~:=~C(i)\sr{C(\kappa_{ij})}{\lal}C(j)\Big),\nn
\eea
is called an \ul{$\I$-complex} (resp. \ul{$\I$-cocomplex}) in $\A$ if composition $\circ_\A$ of morphisms in $\A$ vanishes on $C(\Mor\I)\times C(\Mor\I)\subset\Mor\A\times\Mor\A$, i.e., $\circ_\A|_{C(\Mor\I)\times C(\Mor\I)}=0$, in the sense that the \ul{transition morphisms} $\del_{ij}=\del^C_{ij}$, now called \ul{boundary morphisms} or \ul{differentials}, (resp. $\del^{ij}=\del_C^{ij}$, now called \ul{coboundary morphisms} or \ul{codifferentials}) satisfy the following, which are often written as $\del_\ast^2=0$ (resp. $\del^\ast{}^2=0$):
\bea
&&\del^C_{ij}\circ \del^C_{i'j'}=0,~~\txt{whenever the composition is ``defined''},~~\txt{for all}~~i,j,i',j'\in\Ob\I.\nn\\
&&\Big(\txt{resp.}~~\del_C^{ij}\circ \del_C^{i'j'}=0,~~\txt{whenever the composition is ``defined''},~~\txt{for all}~~i,j,i',j'\in\Ob\I.\Big)\nn
\eea

\ul{A morphism of (co)complexes} is just the associated morphism of systems. We will denote the \ul{category of} \ul{$\I$-complexes} in $\A$ by $\A_0^\I\subset\A^\I$ (resp. \ul{category of $\I$-cocomplexes} in $\A$ by $\A_0^{\I^{op}}\subset\A^{\I^{op}}$).

A \ul{chain complex} (resp. \ul{cochain complex}) in $\A$ is a $\Integer$-complex (resp. $\Integer$-comcomplex) (\magenta{footnote}\footnote{\magenta{Caution}: A ``chain complex'' (i.e., $\Integer$-complex) here should not be confused with ``linear complex'' which would be a tautological reference to ``complex'' since a complex is already a linear system by definition. Confusion is possible because in our discussion of partially ordered sets/classes, ``chain'' and ``linear set/class'' had the same meaning, meanwhile ``chain complex'' here is simply a traditional phrase for ``$\Integer$-complex''.})
\begin{align}
&C_\ast:\Integer\ra\A,~~\big(n\sr{\kappa_n}{\ral}n-1\big)~\mapsto~\big(C_n\sr{\del^C_n}{\ral}C_{n-1}\big):=\big(C_\ast(n)\sr{C_\ast(\kappa_n)}{\ral}C_\ast(n-1)\big),~~~~(\txt{\magenta{footnote}}\footnotemark)\nn\\
&\Big(\txt{resp}.~~~C^\ast:\Integer\ra\A,~~\big(n\sr{\kappa_n}{\ral}n-1\big)~\mapsto~\big(C^n\sr{\del_C^{n-1}}{\lal}C^{n-1}\big):=\big(C^\ast(n)\sr{C^\ast(\kappa_n)}{\lal}C^\ast(n-1)\big)~\Big),\nn
\end{align}
\footnotetext{\magenta{Caution}: Here, our preferred ordering of $\Integer$ is the reverse/opposite of the standard ordering. This is done simply for convenience with the intended applications, but there is no inconsistency and no loss of generality.}where the \ul{boundary maps} or \ul{differentials} $\del^C_n$ (resp. \ul{coboundary maps} or \ul{codifferentials} $\del_C^n$) satisfy
\begin{align}
\del_\ast^2=0,~~\txt{i.e.,}~~\del^C_n\circ\del^C_{n+1}=0,~~\txt{for all}~~n\in\Integer~~\Big(\txt{resp.}~~\del^\ast{}^2=0,~~\txt{i.e.,}~~\del_C^n\circ\del_C^{n-1}=0,~~\txt{for all}~~n\in\Integer\Big).\nn
\end{align}
The chain (resp. cochain) complex is often more conveniently specified in the form of an explicit sequence, or ``chain'', as follows:
\bea
&&\adjustbox{scale=0.9}{\bt
C_\ast=(C_\ast,\del_\ast):~\cdots\ar[r,"\del_{n+2}^C"]&C_{n+1}\ar[r,"\del_{n+1}^C"]&C_n\ar[r,"\del_n^C"]&C_{n-1}\ar[r,"\del_{n-1}^C"]&\cdots
\et}\nn\\
&&\Big(\txt{resp.}~~\adjustbox{scale=0.9}{\bt
C^\ast=(C^\ast,\del^\ast):~\cdots\ar[from=r,"\del^{n+1}_C"']&C^{n+1}\ar[from=r,"\del^n_C"']&C^n\ar[from=r,"\del^{n-1}_C"']&C^{n-1}\ar[from=r,"\del^{n-2}_C"']&\cdots
\et}\Big)\nn
\eea
A \ul{morphism of chain complexes} (resp. \ul{morphism of cochain complexes}) is called a \ul{chain morphism} (\ul{cochain morphism}). We will denote the \ul{category of chain complexes} in $\A$ by $\A_0^\Integer\subset\A^\Integer$ (resp. \ul{category of cochain complexes} in $\A$ by $\A_0^{\Integer^{op}}\subset\A^{\Integer^{op}}$).

A chain complex $(C_\ast,\del_\ast)$ is \ul{exact} at degree $n$ if $\im\del_{n+1}\cong\ker\del_n$. Similarly, a cochain complex $(C^\ast,\del^\ast)$ is \ul{exact} at degree $n$ if $\ker\del^n\cong\im\del^{n-1}$. A chain complex (resp. cochain complex) is an \ul{exact sequence} or an \ul{exact complex} (resp. \ul{exact cocomplex}) if it is exact at all degrees $n\in\Integer$. (\blue{footnote}\footnote{In the literature, the word ``acyclic'' is also used in place of ``exact''.})

An exact complex in $\A$ of the form $0\ra A\sr{f}{\ral}B\sr{g}{\ral}C\ra 0$ is called a \ul{short exact sequence} (written \ul{SES}). Given objects $A,B,C\in\Ob\A$, we say $C$ is the \ul{quotient of $B$ by $A$}, written $C=B/A$ if $C=\coker(A\sr{m}{\ral}B)$ for a monomorphism $m:A\ra B$, or equivalently, we have a SES
\[
0\ra A\ra B\ra C\ra 0~~~~\big(\txt{hence also written}~~~~0\ra A\ra B\ra B/A\ra 0\big),\nn
\]
which also means that given any objects $A,B\in\A$, the expression $B/A$ alone denotes a SES
\[
0\ra A\ra B\ra B/A\ra 0.
\]

A sequence {\footnotesize $S=\left\{S_n\sr{s_n}{\ral}S_{n-1}\right\}$} $\in\A^\Integer$ is \ul{split at degree/index $n$} if in the SES ~$0\ra\ker s_n\sr{k_n}{\hookrightarrow}S_n\sr{\pi_n}{\twoheadrightarrow}\im s_n\ra 0$~ the monic morphism $k_n$ and epic morphism $\pi_n$ are each split (See the decomposition in the proof of Proposition \ref{UnitaryThm}). The sequence $S\in\A^\Integer$ is a \ul{split sequence} if it is split at each index $n\in\Integer$.

Given linear set categories $\I_1,...,\I_k$, a directed multisystem $C\in\A^{\I_1\times\cdots\times\I_k}$ (\blue{footnote}\footnote{Here, $\Ob(\I_1\times\cdots\times\I_k)$ is also viewed as a directed set with $(i_1,...,i_k)\leq(i_1',...,i_k')$ iff $i_j\leq i_j'$ for all $j=1,...,k$.}) is a \ul{$k$-multicomplex} in $\A$ if for each $j\in\{1,...,k\}$ and $\{a_{j'}\in\Ob\I_{j'}:j'\neq j\}$ the restriction
\bea
C_j:=C|_{\{a_1\}\times\cdots\times\{a_{j-1}\}\times\I_j\times\{a_{j+1}\}\times\cdots\times\{a_k\}}:\I_j\ra\A\nn
\eea
is a complex or cocomplex (i.e., $C_j\in \A_0^{\I_j}$ or $C_j\in \A_0^{\I_j^{op}}$). As for complexes and cocomplexes (i.e., $1$-multicomplexes), a \ul{morphism of $k$-multicomplexes} is a just the associated morphism of systems, and we will denote the \ul{category of $k$-multicomplexes} in $\A$ by $\A_0^{\I_1\times\cdots\times\I_k}\subset\A^{\I_1\times\cdots\times\I_k}$. (\blue{footnote}\footnote{By definition, a multicomplex can be mixed in the sense that some of its arguments can be covariant while others are contravariant. This is done for brevity, i.e., so as to avoid making three different definitions for covariant, contravariant and mixed multicomplexes. Moreover, there is no loss of generality due to the obvious covariant-contravariant symmetry.})

A \ul{bicomplex} is a $2$-multicomplex $C_{\ast\ast}\in\A_0^{\I_1\times\I_2}\subset\A^{\I_1\times\I_2}$.
\end{dfn}

\begin{dfn}[\textcolor{blue}{
\index{Total! chain complex of a bicomplex}{Total chain complex of a bicomplex},
\index{Total! functor}{Total functor}}]\label{HomCohomDf1}
Let $\A$ be an abelian category. The \ul{total chain complex} of a bicomplex $C_{\ast\ast}\in\A_0^{\Integer\times\Integer}\subset \A^{\Integer\times\Integer}$ (expressed in diagram form as)

\bea
\label{DoubleDiffEq}\adjustbox{scale=0.8}{\bt
 ~       &  \ar[from=dd,"\del^v_{i+1~j-1}"']  ~                         &&  \ar[from=dd,"\del^v_{i~j-1}"']  ~                        &&   \ar[from=dd,"{\del^v_{i-1~j-1}}"']  ~                              && ~\\
 ~       &    ~                                                         &&    ~                                                      &&    ~                                                                 && ~\\
 ~\ar[r] & C_{i+1~j-1}\ar[rr,"\del^h_{i+1~j-1}"]                        && C_{i~j-1}\ar[rr,"\del^h_{i~j-1}"]                         && C_{i-1~j-1}\ar[rr,"\del^h_{i-1~j-1}"] && ~\\
 ~       &     ~                                                        &&     ~                                                     &&     ~                                                           && ~\\
 ~\ar[r] & C_{i+1~j}\ar[uu,"\del^v_{i+1~j}"']\ar[rr,"{\del^h_{i+1~j}}"] && C_{i~j}\ar[rr,"\del^h_{i~j}"]\ar[uu,"\del^v_{i~j}"'] && C_{i-1~j}\ar[uu,"{\del^v_{i-1~j}}"']\ar[rr,"\del^h_{i-1~j}"] && ~\\
 ~       &     ~                                                        &&     ~                                                     &&     ~                                                           && ~\\
 ~\ar[r] & C_{i+1~j+1}\ar[uu,"{\del^v_{i+1~j+1}}"']\ar[rr,"\del^h_{i+1~j+1}"] && C_{i~j+1}\ar[uu,"\del^v_{i~j+1}"']\ar[rr,"\del^h_{i~j+1}"]   && C_{i-1~j+1}\ar[uu,"\del^v_{i-1~j+1}"']\ar[rr,"\del^h_{i-1~j+1}"]    && ~\\
 ~       &   \ar[u]  ~                           &&    \ar[u] ~                         &&   \ar[u]  ~                              && ~
\et}
\eea
is the chain complex ~$Tot(C_{\ast\ast})_\ast=\big(Tot(C_{\ast\ast})_\ast,\del_\ast\big)\in\A_0^\Integer\subset\A^\Integer$~ given \index{Formal expression}{\ul{formally}} (\blue{footnote}\footnote{We call a mathematical expression (involving mappings/operations) a \ul{formal expression} if (i) it is appealing in some way but (ii) it is not well defined as a mapping without some (relatively trivial/obvious) additional information.}) by
\bea
\textstyle Tot(C_{\ast\ast})_n:=\bigoplus\limits_{i+j=n}C_{ij},~~~~\del_n=\del^{Tot(C_{\ast\ast})}_n:=\sum\limits_{i+j=n}\left(\del^h_{ij}+(-1)^i\del^v_{ij}\right):Tot(C_{\ast\ast})_n\ra Tot(C_{\ast\ast})_{n-1}.\nn
\eea
Since the ordinary sum indeed makes sense between morphisms with the same destination, it follows from the diagram (\ref{DoubleDiffEq}) that a clearer (i.e., less formal) rewriting of the differential of the total complex is
\[
\textstyle\del^{Tot(C_{\ast\ast})}_n:=\sum\limits_{i+j=n}\left(\del^h_{i+1~j}+(-1)^i\del^v_{i~j+1}\right)
=\sum\limits_{i+j=n+1}\left(\del^h_{ij}+(-1)^{i-1}\del^v_{i-1~j+1}\right)~=:~\sum\limits_{i+j=n}\left(\del^h_{ij}+(-1)^i\del^v_{ij}\right).
\]
The \ul{total functor} of $\A$ is the following functorial process associated with the total complex:
\[
Tot:\A_0^{\Integer\times\Integer}\ra\A_0^{\Integer},~~\big(C_{\ast\ast}\sr{f_{\ast\ast}}{\ral}C'_{\ast\ast}\big)~\mapsto~ \Big(Tot(C_{\ast\ast})_\ast\sr{Tot(f_{\ast\ast})}{\ral}Tot(C'_{\ast\ast})_\ast\Big),
\]
where ~$Tot(f_{\ast\ast})_n:=\sum_{i+j=n}f_{ij}:Tot(C_{\ast\ast})_n\ra Tot(C'_{\ast\ast})_n$.

When $\A\subset\txt{Sets}$ (e.g., $\A$ $=$ $R$-mod), the differential of the total complex is explicitly given elementwise by the following: For any ~$(c_{ij})_{i+j=n}\in Tot(C_{\ast\ast})_n$,
{\small\begin{align}
&\textstyle \del_n~(c_{ij})_{i+j=n}=\del^{Tot(C_{\ast\ast})}_n~(c_{ij})_{i+j=n}~:=~\sum\limits_{i'+j'=n}\left(\del^h_{i'j'}+(-1)^{i'}\del^v_{i'j'}\right)(c_{ij})_{i+j=n}\nn\\
&\textstyle~~~~\sr{(s)}{=}\Big(\del^h_{i+1,j}c_{i+1,j}+(-1)^i\del^v_{i,j+1}c_{i,j+1}\Big)_{i+j=n}=\Big(\del^h_{i+1,n-i}c_{i+1,n-i}+(-1)^i\del^v_{i,n-i+1}c_{i,n-i+1}\Big)_i\nn\\
&\textstyle~~~~=\Big(\del^h_{ij}c_{ij}+(-1)^{i-1}\del^v_{i-1,j+1}c_{i-1,j+1}\Big)_{i+j=n+1}
=\Big(\del^h_{i~n-i+1}c_{i~n-i+1}+(-1)^{i-1}\del^v_{i-1,n-i}c_{i-1,n-i}\Big)_{i-1}\nn\\
&\textstyle~~~~=\Big(\del^h_{i+1,j-1}c_{i+1,j-1}+(-1)^i\del^v_{ij}c_{ij}\Big)_{i+j=n+1}=\Big(\del^h_{i+1,n-i}c_{i+1,n-i}+(-1)^i\del^v_{i~n-i+1}c_{i~n-i+1}\Big)_{i-1},\nn
\end{align}}where step (s) is obtained using the implicit fact that (for any given entry of $\del_n(c_{ij})_{i+j=n}$) addition in the \index{Formal sum}{\ul{formal sum}} $\sum_{i'+j'=n}\left(\del^h_{i'j'}+(-1)^{i'}\del^v_{i'j'}\right)$ makes sense only for (and is thus determined only by) those differentials in the diagram (\ref{DoubleDiffEq}) with the same destination (i.e., that ``converge'' to the same vertex).
\end{dfn}

\begin{rmk}
The grading factor $(-1)^i$ in the formula for $\del_n^{Tot(C_{\ast\ast})}$ is deliberately chosen to ensure that ~$\del_n^{Tot(C_{\ast\ast})}\circ \del_{n+1}^{Tot(C_{\ast\ast})}=0$. Indeed, for any ~$(c_{ij})_{i+j=n+1}\in Tot(C_{\ast\ast})_{n+1}$,~ we have
{\small\begin{align}
&\textstyle \del^{Tot(C_{\ast\ast})}_n\del^{Tot(C_{\ast\ast})}_{n+1}~(c_{ij})_{i+j=n+1}=\del^{Tot(C_{\ast\ast})}_n\Big(\overbrace{\del^h_{i+1,j}c_{i+1,j}+(-1)^i\del^v_{i,j+1}c_{i,j+1}}^{\rho_{ij}}\Big)_{i+j=n+1}\nn\\
&~~~~=\Big(\del^h_{ij}\rho_{ij}+(-1)^{i-1}\del^v_{i-1,j+1}\rho_{i-1,j+1}\Big)_{i+j=(n+1)+1}\nn\\
&~~~~=\Big(\del^h_{ij}[\del^h_{i+1,j}c_{i+1,j}+(-1)^i\del^v_{i,j+1}c_{i,j+1}]+(-1)^{i-1}\del^v_{i-1,j+1}[\del^h_{ij}c_{i,j+1}+(-1)^{i-1}\del^v_{i-1,j+2}c_{i-1,j+2}]\Big)_{i+j=n+2}\nn\\
&~~~~=\Big(0+(-1)^i[\del^h_{ij}\del^v_{i,j+1}c_{i,j+1}-\del^v_{i-1,j+1}\del^h_{ij}c_{i,j+1}]+0\Big)_{i+j=n}~\sr{(s)}{=}~(0+0+0)_{i+j=n}=0,\nn
\end{align}}
where step (s) follows by the commutativity of the diagram (\ref{DoubleDiffEq}): $\del^h_{i,j-1}\del^v_{ij}=\del^v_{i-1,j}\del^h_{ij}$ (\blue{footnote}\footnote{In practice, when the differentials are being applied on objects/elements that already carry the required indices, we can (for convenience) omit the indices on differentials. Thus, it is sufficient to simply write $\del^h\del^v=\del^v\del^h$, and for example,
\be
\del^{Tot(C_{\ast\ast})}_n~(c_{ij})_{i+j=n}~:=~\big(\del^hc_{i+1,j}+(-1)^i\del^vc_{i,j+1}\big)_{i+j=n}=\big(\del^hc_{i+1,n-i}+(-1)^i\del^vc_{i,n-i+1}\big)_i.\nn
\ee}).
\end{rmk}

\begin{rmk}[\textcolor{blue}{Matrix form of the differential of $Tot(C_{\ast\ast})$}]
As a rule for the matrix form of $\del_n^{Tot}$, maps on the \ul{same column} have the \ul{same domain}, while maps on the \ul{same row} have the \ul{same codomain}. Thus,
{\scriptsize
    \bea
    &&\textstyle \del_n^{Tot(C_{\ast\ast})}=\sum\limits_{i+j=n}~\del^h_{ij}+(-1)^i\del^v_{ij}:=
    \left[
      \begin{array}{cccccc}
        \ddots & \vdots & \vdots &  \vdots & \vdots &  \\\\
        \cdots & \del^h_{i-1,j+1}            & 0             & 0                       & 0 &  \cdots \\\\
        \cdots & (-1)^{i-1}\del^v_{i-1,j+1} & \del^h_{i,j}       & 0                       & 0 &  \cdots \\\\
        \cdots & 0                      & (-1)^i\del^v_{i,j} & \del^h_{i+1,j-1}             & 0 &  \cdots \\\\
        \cdots & 0                      & 0             & (-1)^{i+1}\del^v_{i+1,j-1}  & \del^h_{i+2,j-2} &  \cdots \\\\
         & \vdots & \vdots & \vdots & \vdots & \ddots \\
      \end{array}
    \right].\nn\\
    &&\textstyle \del_n^{Tot(C_{\ast\ast})}:Tot(C_{\ast\ast})_n\ra Tot(C_{\ast\ast})_{n-1},~~
    \left[
      \begin{array}{l}
        \vdots \\
        c_{i-1,j+1} \\
        c_{i,j} \\
        c_{i+1,j-1} \\
        c_{i+2,j-2} \\
        \vdots \\
      \end{array}
    \right]\longmapsto
   \left[
      \begin{array}{l}
        \vdots \\
        \del^h_{i-1,j+1}c_{i-1,j+1}+(-1)^{i-2}\del^v_{i-2,j+2}c_{i-2,j+2} \\
        \del^h_{i,j}c_{i,j}+(-1)^{i-1}\del^v_{i-1,j+1}c_{i-1,j+1} \\
        \del^h_{i+1,j-1}c_{i+1,j-1}+(-1)^i\del^v_{i,j}c_{i,j} \\
        \del^h_{i+2,j-2}c_{i+2,j-2}+(-1)^{i+1}\del^v_{i+1,j-1}c_{i+1,j-1} \\
        \vdots \\
      \end{array}
    \right],\nn
    \eea}
where ~~$Tot(C_{\ast\ast})_n=\bigoplus\limits_{i+j=n}C_{ij}$. Note however that in specific applications (where expressions often can be simplified) conventions can be altered accordingly.
\end{rmk}

\section{Chain homotopy and Homotopy equivalence}

\begin{dfn}[\textcolor{blue}{
\index{Chain! homotopy}{{Chain homotopy}},
\index{Homotopic chain morphisms}{Homotopic chain morphisms},
\index{Chain! homotopy between chain morphisms}{Chain homotopy between chain morphisms},
\index{Chain! homotopy equivalence}{Chain homotopy equivalence},
\index{Homotopy! equivalent chain complexes}{Homotopy equivalent chain complexes},
\index{Chain! homotopy inverse}{Chain homotopy inverse},
\index{Nullhomotopic chain morphism}{Nullhomotopic chain morphism},
\index{Contractible chain complex}{Contractible chain complex}}]

Let $\A$ be an abelian category (or just an additive category) and $(C_\ast,\del_\ast),(C'_\ast,\del_\ast')\in\A_0^\Integer$ chain complexes. A \ul{chain homotopy} $h:C_\ast\ra C_{\ast+1}'$ is a collection of morphisms $\{h_n:C_n\ra C_{n+1}'\}_{n\in\Integer}$.
\bc\adjustbox{scale=0.9}{\bt
\cdots\ar[r] &C_{n+1}\ar[rr,"\del_{n+1}"]&& C_n\ar[ddll,dashed,"h_n"']\ar[dd,"\delta_n(h)"']\ar[rr,"\del_n"] && C_{n-1}\ar[ddll,dashed,"h_{n-1}"']\ar[r,"{\del_{n-1}}"] &\cdots\\
 & && && &\\
\cdots\ar[r]&C_{n+1}'\ar[rr,"\del_{n+1}'"]&& C_n'\ar[rr,"\del_n'"] && C_{n-1}'\ar[r,"{\del_{n-1}'}"] &\cdots
\et}\ec
Note that the morphisms $\delta_n(h):=h_{n-1}\circ\del_n+\del'_{n+1}\circ h_n:C_n\ra C_n'$ define a chain morphism $\delta_\ast(h):C_\ast\ra C_\ast'$ (as can be checked directly).

Two chain morphisms $f_\ast,g_\ast:C_\ast\ra C_\ast'$ are \ul{homotopic}, written $f_\ast\simeq g_\ast$, if there exists a chain homotopy $h:C_\ast\ra C_{\ast+1}'$) such that $f_n-g_n=\delta_n(h):=h_{n-1}\circ\del_n+\del'_{n+1}\circ h_n$, for all $n\in\Integer$ (or briefly, $f_\ast-g_\ast=h\del_\ast+\del_\ast'h$). In this case, we say $h:C_\ast\ra C_{\ast+1}'$ is a \ul{chain homotopy between $f_\ast$ and $g_\ast$}.

\bc\adjustbox{scale=0.9}{\bt
\cdots\ar[r] &C_{n+1}\ar[rr,"\del_{n+1}"]&& C_n\ar[ddll,dashed,"h_n"']\ar[dd,"\substack{f_n-g_n=\delta_n(h)}"description]\ar[rr,"\del_n"] && C_{n-1}\ar[ddll,dashed,near start,"h_{n-1}"]\ar[r,"{\del_{n-1}}"] &\cdots\\
 & && && &\\
\cdots\ar[r]&C_{n+1}'\ar[rr,"\del_{n+1}'"]&& C_n'\ar[rr,"\del_n'"] && C_{n-1}'\ar[r,"{\del_{n-1}'}"] &\cdots
\et}\ec

A chain morphism $f_\ast:C_\ast\ra C_\ast'$ is a \ul{chain homotopy equivalence} (making $C_\ast,C_\ast'$ \ul{homotopy equivalent} chain complexes, written $C_\ast\simeq C_\ast'$) if there exists a chain morphism $g_\ast:C_\ast'\ra C_\ast$ (called a \ul{chain homotopy inverse} of $f_\ast$) such that $g_\ast\circ f_\ast\simeq id_{C_\ast}$ and $f_\ast\circ g_\ast\simeq id_{C_\ast'}$.

A chain morphism $f_\ast:C_\ast\ra C_\ast'$ is \ul{nullhomotopic} if $f_\ast\simeq 0_\ast$ (i.e., $f$ is homotopic to the zero chain morphism). A chain complex $(C_\ast,\del_\ast)$ is \ul{contractible} if $(C_\ast,\del_\ast)\simeq (0_\ast,0_\ast)$ (i.e., it is homotopy equivalent to the zero chain complex).
\end{dfn}

\begin{notation}
We will sometimes abbreviate the word ``homotopy'' as ``hty''.
\end{notation}

\begin{question}
Let $\A$ be an abelian category. Any morphism $A\sr{f}{\ral}A'$ can be viewed as a chain complex $C_f:0\ra A\sr{f}{\ral}A'\ra 0$. Thus (as the category of morphisms $Cat(\Mor\A)$ seen earlier) we have the category ~{\footnotesize $\Mor\A=\Big(\Ob(\Mor\A),\Mor_{\Mor\A}\big(\Ob(\Mor\A),\Ob(\Mor\A)\big),\circ\Big)\subset\A_0^\Integer$}~ whose objects are morphisms in $\A$ and whose morphisms are chain morphisms between (morphisms $A\sr{f}{\ral}A'$ and $B\sr{g}{\ral}B'$ in $\A$ viewed as the) chain complexes $C_f:0\ra A\sr{f}{\ral}A'\ra0$ and $C_g:0\ra B\sr{g}{\ral}B'\ra0$. Does the containment $\Mor\A\subset\A_0^\Integer$ give a subcategory of $\A_0^\Integer$?
\end{question}

\begin{rmk}
Let $\A$ be an abelian category. Then $\A_0^\Integer$ is also an abelian category with respect to the following pointwise operations. Let $A_\ast,B_\ast\in\Ob\A_0^\Integer$, $f_\ast,g_\ast\in \Mor_{\A_0^\Integer}(A_\ast,B_\ast)$, and $n\in\Integer$.
{\flushleft \ul{A zero object exists}}: ~$(0_\ast,\del_\ast):~\cdots\sr{0}{\ral}0\sr{0}{\ral}0\sr{0}{\ral}\cdots$.
{\flushleft \ul{Products exist}}: ~$(A_\ast\oplus B_\ast,\del_\ast^{A\oplus B}):~\bt\cdots\ar[r]&A_{n+1}\oplus B_{n+1}\ar[r,"\del^A_{n+1}\oplus\del^B_{n+1}"]&A_n\oplus B_n\ar[r,"\del^A_n\oplus\del^B_n"]&A_{n-1}\oplus B_{n-1}\ar[r,"\del^A_{n-1}\oplus\del^B_{n-1}"]&\cdots\et$.
{\flushleft \ul{Morphism classes are abelian groups}}: ~$(f_\ast+g_\ast)_n:=f_n+g_n$.
{\flushleft \ul{Kernels and cokernels exist}}: The $n$th component of ~$\ker f_\ast\sr{k_{f_\ast}}{\ral}A_\ast$~ (resp. ~$B_\ast\sr{c_{f_\ast}}{\ral}\coker f_\ast$) is given by ~$\ker f_n\sr{k_{f_n}}{\ral}A_n$~ (resp. ~$B_n\sr{c_{f_n}}{\ral}\coker f_n$). The differential of $\ker f_\ast$ (and that of $\coker f$ similarly) exists automatically due to the universal property of the kernel as in the following commutative diagram:
\bc
\adjustbox{scale=0.8}{\bt
    &                                                & A_n\ar[lld,"f_n"']\ar[rrr,"\del^A_n"] &         &   & A_{n-1}\ar[lld,"f_{n-1}"']\\
B_n\ar[rrr,near start,"\del^B_n"description] &                                                &                    & B_{n-1} &   &  \\
    & \ker f_n\ar[ul,"0"]\ar[uur,near end,"k_{f_n}"description]\ar[rrr,dashed,"\del^{\ker f}_n"]          &                    &         & \ker f_{n-1}\ar[ul,"0"]\ar[uur,"k_{f_{n-1}}"description]  &  \\
    & & & & &\\
    & X_n\ar[uuul,bend left,"0"]\ar[uu,dashed,"h'_n"']
    \ar[ruuuu,bend right,near start,"h_n"description]\ar[rrr,"\del^X_n"]                       &                    &         & X_n\ar[uuul,bend left,"0"]\ar[uu,dashed,"h'_{n-1}"']
                                                                                      \ar[ruuuu,bend right,"h_{n-1}"description]  &
\et}
\ec
{\flushleft \ul{Monomorphisms are kernels and epimorphisms are cokernels}}: These hold by construction.
\end{rmk}

\begin{prp}\label{CompContrPrp}
Let $\A$ be an abelian category. A chain complex $C=(C,d)\in\A_0^\Integer$ is contractible $\iff$ $id_{C}\simeq 0:C\ra C$ (i.e., $id_{C}=hd+dh$ for a homotopy $h:C_\ast\ra C_{\ast+1}$).
\end{prp}
\begin{proof}
($\Ra$): Assume $C\simeq0$, and let $f:C\ra0$ be a homotopy equivalence with homotopy inverse $g:0\ra C$. Then $f\circ g\simeq id_{0}$ and $0=g\circ f\simeq id_{C}$.\\
($\La$): Assume $id_{C}\simeq 0$. Let $f:C\ra0$ and $g=i_{0}:0\hookrightarrow C$. Then $f\circ g=id_{0}$, $g\circ f=0\simeq id_{C}$, and so $f:C\ra0$ is a homotopy equivalence.
\end{proof}

\section{The Homologies of a Chain Complex: Homology Functors}
\begin{dfn}[\textcolor{blue}{
\index{Homologies of a chain complex}{Homologies of a chain complex}}]\label{HlogyDef}
Let $\A$ be an abelian category and $\A_0^\Integer$ the category of chain complexes in $\A$. Given a chain morphism $f_\ast:(C_\ast,\del_\ast)\ra(C_\ast',\del_\ast')$, consider the following decomposition based on the proof of the unitarity theorem (Proposition \ref{UnitaryThm}).
\bea
\label{HomologDec}\adjustbox{scale=0.8}{\bt
\ker\del_n\ar[dd,dashed,"\exists!~k_n(f_\ast)"]\ar[rr,hook,"k_{\del_n}","k_n"']&& C_n\ar[dd,"f_n"]\ar[rrrrrr,bend left,"\del_n"]\ar[rr,two heads,"ck_{\del_n}","ck_n"']&&\ub{\overbrace{\coker k_{\del_n}}^{\coim \del_n}}_{\coker(\ker\del_n)}\ar[dd,dashed,"\exists!~ck_n(f_\ast)"]\ar[rr,hook,two heads,dashed,"u_{\del_n}","u_n"']&&\ub{\overbrace{\ker c_{\del_n}}^{\im\del_n}}_{\ker(\coker\del_n)}\ar[dd,dashed,"\exists!~kc_n(f_\ast)"]\ar[rr,hook,"kc_{\del_n}","kc_n"']&&C_{n-1}\ar[dd,"f_{n-1}"]\ar[rr,two heads,"c_{\del_n}","c_n"']&&\coker\del_n\ar[dd,dashed,"\exists!~c_n(f_\ast)"]\\
 && && && &&\\
\ker\del'_n\ar[rr,hook,"k_{\del'_n}","k'_n"']&& C'_n\ar[rrrrrr,bend right,"\del'_n"]\ar[rr,two heads,"ck_{\del'_n}","ck'_n"']&&\ub{\overbrace{\coker k_{\del'_n}}^{\coim \del'_n}}_{\coker(\ker\del'_n)}\ar[rr,hook,two heads,dashed,"u_{\del'_n}","u'_n"']&&\ub{\overbrace{\ker c_{\del'_n}}^{\im\del'_n}}_{\ker(\coker\del'_n)}\ar[rr,hook,"kc_{\del'_n}","kc'_n"']&&C'_{n-1}\ar[rr,two heads,"c_{\del'_n}","c'_n"']&&\coker \del'_n
\et}
\eea

(\textbf{Note}: The existence of the induced morphisms between kernels/cokernels is immediate from the fact that kernels/cokernels are limits/colimits and are therefore functorial processes.)

From {\small $\del_{n+1}\circ\del_n=0$}, we know that {\small $kc_{k+1}:\im\del_{n+1}\hookrightarrow C_n$} must factor through {\small $k_n:\ker\del_n\hookrightarrow C_n$}, i.e., there is a unique monic morphism (the exactness monomorphism) {\small $m_n:=m_{\del_{n+1}|\del_n}:\im\del_{n+1}\hookrightarrow\ker\del_n$} such that
\bea
(\im\del_{n+1}\hookrightarrow C_n)=(\ker\del_n\hookrightarrow C_n)\circ(\im\del_{n+1}\hookrightarrow\ker\del_n).\nn
\eea
The notion of homology is defined in terms of the cokernel of the above induced (exactness) monomorphism. Specifically, the \ul{$n$th homology} (or \ul{$n$th homology group}) of $(C_\ast,\del_\ast)$ is defined as
\bea
H_n(C_\ast,\del_\ast):={\ker\del_n\over\im\del_{n+1}}:=\coker(\im\del_{n+1}\hookrightarrow\ker\del_n)=\coker(B_n\hookrightarrow Z_n).\nn
\eea
\end{dfn}

\begin{rmk}
The decomposition of the complex $(C_\ast,\del_\ast)$ yields the following two interesting complexes
\bea
&&\bt 0\ar[r]&\overbrace{\ker\del_n}^{Z_n}\ar[r,hook,"k_n"]&C_n\ar[r,two heads,"ck_n"]& \coim\del_n\ar[r,hook,two heads,dashed,"u_n"]&\overbrace{\im\del_n}^{B_{n-1}}\ar[r]&0\et\nn\\
&&\nn\\
&&\bt 0\ar[r]&\ub{\im\del_{n+1}}_{B_n}\ar[r,hook,"m_n"]&\ub{\ker\del_n}_{Z_n}\ar[r,two heads,"cm_n"]& {H_n(C_\ast,\del_\ast)}\ar[r]&0\et\nn
\eea
which can be displayed more explicitly as follows:
\bea\adjustbox{scale=0.9}{\bt
    &                           &                &                  &                          &       0\ar[dr]             &   0                    &          &  \\
    &                           &                &                  &              0\ar[r]     & B_{n-1}\ar[r,hook,"m_{n-1}"]\ar[ur] & Z_{n-1}\ar[dr,hook]\ar[r,two heads,"cm_{n-1}"]  & H_{n-1}={Z_{n-1}\over B_{n-1}}\ar[r]  & 0 \\
\cdots\ar[r] & C_{n+1}\ar[rrr,"\del_{n+1}"]\ar[dr,two heads,"ck_{n+1}"] &                &                  & C_n\ar[rrr,"\del_n"]\ar[ur,two heads,"ck_n"] &                    &                       & C_{n-1}\ar[r]  & \cdots   \\
    &       0\ar[r]             & B_n\ar[r,hook,"m_n"]\ar[dr] & Z_n\ar[ur,hook,"k_n"]\ar[r,two heads,"cm_n"] & H_n={Z_n\over B_n}\ar[r] & 0                  &                       &          &         \\
    &                           &      0\ar[ur]  &   0               &                          &                    &                       &          &
\et}\nn
\eea
where the epimorphisms $ck_n$ in the above diagram are only specified up to unitary equivalence (i.e., up to composition with unitary morphisms, which are isomorphisms since $\A$ is abelian).
\end{rmk}

\begin{dfn}[\textcolor{blue}{
\index{Homology! functors}{Homology functors},
\index{Chains}{Chains},
\index{Boundaries}{Boundaries},
\index{Cycles}{Cycles}}]\label{HlogyDef2}
Based on the above information, including the decomposition in (\ref{HomologDec}), we further define the \ul{$n$th homology functor} of $\A$ as
\bea
H_n:\A_0^\Integer\ra\A,~~\Big((C_\ast,\del_\ast)\sr{f_\ast}{\ral}(C'_\ast,\del'_\ast)\Big)~\mapsto~\Big(H_n(C_\ast,\del_\ast)\sr{H_n(f_\ast)}{\ral}H_n(C'_\ast,\del'_\ast)\Big),\nn
\eea
where
\bit
\item[(1)] $H_n(C_\ast,\del_\ast):=\coker\big(\im\del_{n+1}\sr{m_n}{\hookrightarrow}\ker\del_n\big)$ is the \ul{$n$th homology (group)} of the chain complex $(C_\ast,\del_\ast)$.
\item[(2)] $Z_n(C_\ast,\del_\ast):=\ker\del_n$ is the \ul{$n$th cycle (group)} of the chain complex $(C_\ast,\del_\ast)$.
\item[(3)] $B_n(C_\ast,\del_\ast):=\im\del_{n+1}$ is the \ul{$n$th boundary (group)} of the chain complex $(C_\ast,\del_\ast)$.
\eit
and the morphism $H_n(f_\ast)$ is induced by the functorial property of the cokernel as follows:
\bea\bt
0\ar[r]&\overbrace{\im\del_{n+1}}^{B_n}\ar[d,"kc_n(f_\ast)"]\ar[r,hook,"m_n"]&\overbrace{\ker\del_n}^{Z_n}\ar[d,"k_n(f_\ast)"]\ar[r,two heads,"cm_n"]& \overbrace{H_n(C_\ast,\del_\ast)}^{Z_n/B_n}\ar[d,dashed,"H_n(f_\ast)"]\ar[r]&0 \\
0\ar[r]&\ub{\im\del'_{n+1}}_{B'_n}\ar[r,hook,"m'_n"]&\ub{\ker\del'_n}_{Z'_n}\ar[r,two heads,"cm'_n"]& \ub{H_n(C'_\ast,\del'_\ast)}_{Z'_n/B'_n}\ar[r]&0
\et\nn
\eea
In an abelian category $\A$ whose objects are sets, elements of $C_n$ are called \ul{$n$-chains}, elements of $Z_n(C_\ast,\del_\ast)$ are called \ul{$n$-cycles}, and elements of $B_n(C_\ast,\del_\ast)$ are called \ul{$n$-boundaries}.
\end{dfn}

\begin{lmm}
The homology functor $H_n:\A_0^\Integer\ra\A$ (or $H_\ast:\A_0^\Integer\ra\A_0^\Integer$) in Definition \ref{HlogyDef2} is additive.

(\ul{Byproduct}: The cycle, c-cycle, boundary functors {\small $k_n=Z_n,c_n,kc_n=B_n:\A_0^\Integer\ra\A$} are all additive.)
\end{lmm}
\begin{proof}
Let $f_\ast,g_\ast\in \Mor_{\A_0^\Integer}(A_\ast,A'_\ast)$. Since (by Proposition \ref{LimFunct}) kernels are limits (hence functors), $f_\ast$, $g_\ast$, $f_\ast+g_\ast$ induce respective chain morphisms $\wt{f}_\ast,\wt{g}_\ast,\wt{f_\ast+g_\ast}\in \Mor_{\A_0^\Integer}(\ker \del^A_\ast,\ker\del^{A'}_\ast)$ as follows:
\bea
\adjustbox{scale=1.0}{\bt
\ker\del^A_\ast\ar[dd,dashed,"{\wt{f}_\ast,~\wt{g}_\ast,~\wt{f_\ast+g_\ast}}"]\ar[rrr,hook,"k_\ast"] &&& A_\ast\ar[dd,"{f_\ast,~g_\ast,~f_\ast+g_\ast}"]\ar[rrr,two heads,"ck_\ast"] &&& \coim\del^A_\ast\ar[dd,dashed,"{\ol{f}_\ast,~\ol{g}_\ast,~\ol{f_\ast+g_\ast}}"]\\
 &&& &&& \\
\ker\del^{A'}_\ast\ar[rrr,hook,"k'_\ast"] &&& A'_\ast\ar[rrr,two heads,"ck'_\ast"] &&& \coim\del^{A'}_\ast\\
\et}\nn
\eea
From the left square, $f_\ast\circ k_\ast=k'_\ast\circ\wt{f}_\ast$, ~$g_\ast\circ k_\ast=k'_\ast\circ\wt{g}_\ast$, ~$(f_\ast+g_\ast)\circ k_\ast=k'_\ast\circ\wt{f_\ast+g_\ast}$, and so
\bea
&&k'_\ast\circ\wt{f_\ast+g_\ast}=(f_\ast+g_\ast)\circ k_\ast=f_\ast\circ k_\ast+g_\ast\circ k_\ast=k'_\ast\circ\wt{f}_\ast+k'_\ast\circ\wt{g}_\ast=k'_\ast\circ(\wt{f}_\ast+\wt{g}_\ast),\nn\\
&&~~\Ra~~\wt{f_\ast+g_\ast}=\wt{f}_\ast+\wt{g}_\ast~~~~\txt{(since $k'_\ast$ is monic)}.\nn
\eea
From the right square, the same procedure gives ~$\ol{f_\ast+g_\ast}=\ol{f}_\ast+\ol{g}_\ast$ (because $ck_\ast$ is epic).

Similarly, we have induced maps on homology as follows:
\bea
\adjustbox{scale=1.0}{\bt
\im\del^A_{\ast+1}\ar[dd,dashed,"{\wh{f}_\ast,~\wh{g}_\ast,~\wh{f_\ast+g_\ast}}"]\ar[rrr,hook,"\kappa_\ast"] &&& \ker\del^A_\ast \ar[dd,dashed,"{\wt{f}_\ast,~\wt{g}_\ast,~\wt{f_\ast+g_\ast}}"]\ar[rrr,two heads,"c\kappa_\ast"] &&& H_\ast(A_\ast)\ar[dd,dashed,"{H_\ast(f_\ast),~H_\ast(g_\ast),~H_\ast(f_\ast+g_\ast)}"]\\
 &&& &&& \\
\im\del^{A'}_{\ast+1}\ar[rrr,hook,"{\kappa'_\ast}"] &&& \ker\del^{A'}_\ast\ar[rrr,two heads,"{c\kappa'_\ast}"] &&& {H_\ast(A'_\ast)}\\
\et}\nn
\eea
where, using the previous result $\wt{f_\ast+g_\ast}=\wt{f}_\ast+\wt{g}_\ast$, we likewise see that
\bea
\wh{f_\ast+g_\ast}=\wh{f}_\ast+\wh{g}_\ast~~~~\txt{and}~~~~H_\ast(f_\ast+g_\ast)=H_\ast(f_\ast)+H_\ast(g_\ast).\nn\qedhere
\eea
\end{proof}

\begin{crl}[\textcolor{blue}{\index{Homotopy! invariance of homology}{Homotopy invariance of homology}}]\label{ModHomInvCor}
For any two morphisms of module chain complexes $f,g$ $\in$ $\Mor_{(R\txt{-mod})_0^\Integer}(C,C')$, if $f\simeq g$ (i.e., $f=g+h\circ\del^C+\del^{C'}\circ h$ for a homotopy $h:C_\ast\ra C'_{\ast+1}$), then
\bc
$H_n(f)=H_n(g)$. ~~~~(\blue{footnote}\footnote{We will prove later that any small abelian category $\A$ imbeds (with kernels and cokernels preserved) into some $R$-mod. Therefore this result is actually valid for complexes over any abelian category.})
\ec

\end{crl}
\begin{proof}
Over $R$-mod, the homology functor {\footnotesize $H_n:\big(C\sr{f}{\ral}C'\big)\mapsto \big(H_n(C)\sr{H_n(f)}{\ral}H_n(C')\big)$} is explicitly given by
\[
\textstyle H_n(f):{Z_n\over B_n}={\ker\del_n\over\im\del_{n+1}}\ra {Z'_n\over B'_n}={\ker\del'_n\over\im\del'_{n+1}},~z+B_n\mapsto f(z)+B'_n,~~~~\txt{where}~~~~\del:=\del^C,~~\del':=\del^{C'},\nn
\]
which is well-defined because the chain morphism property $f\del=\del' f$ implies $f(\im\del)\subset\im\del'$.
If $f\simeq g$, then it is clear that $\im(f-g)\subset\im\del'$, and so $H_n(f)-H_n(g)=H_n(f-g)=0$, i.e., $H_n(f)=H_n(g)$.
\end{proof}

\section{Homology and Cohomology: The Cohomologies of a Chain Complex}
\begin{dfn}[\textcolor{blue}{
\index{Stable chain complex}{Stable chain complex},
\index{Homology! functor}{Homology functor},
\index{Homology! of a chain complex}{{Homology of a chain complex}},
\index{Homologically exact complex}{Homologically exact complex},
\index{Quasiisomorphism (qis)}{Quasiisomorphism (qis)},
\index{Quasiisomorphic chain complexes}{Quasiisomorphic chain complexes},
\index{Cohomology of a chain complex}{{Cohomology of a chain complex}},
\index{Cochain complex of a chain complex}{Cochain complex of a chain complex},
\index{Cochains}{Cochains},
\index{Cocycles}{Cocycles},
\index{Coboundaries}{Coboundaries}}]\label{CoHlogyDef}

Let $\A,\B$ be abelian categories and $\A_0^\Integer,\B_0^\Integer$ the categories of chain complexes in $\A,\B$ respectively. A chain complex $C_\ast=(C_\ast,\del_\ast)\in\A_0^\Integer$ is a \ul{stable chain complex} if it has zero differentials (i.e., $\del_\ast=0$, meaning $\del_n=0$ for each $n$), written $C_\ast=(C_\ast,0)$. Denote the (subcategory of) stable chain complexes by $\A_{00}^\Integer\subset\A_0^\Integer\subset\A^\Integer$.

A functor $H_\ast:\A_0^\Integer\ra\B_0^\Integer$ is a \ul{homology functor} if (at least) the following hold.
\begin{enumerate}
\item $H_\ast$ is additive.
\item Homotopic chain morphisms have the same image (and hence homotopy equivalent chain complexes have isomorphic images): i.e., if $f\simeq g$ then $H_\ast(f)=H_\ast(g)$ ~(and hence, if $(C_\ast,\del_\ast)\simeq (C'_\ast,\del'_\ast)$ then $H_\ast(C_\ast,\del_\ast)\cong H_\ast(C'_\ast,\del'_\ast)$).
\item The image $H_\ast(C_\ast,\del_\ast)$, called the \ul{homology of the chain complex} $(C_\ast,\del_\ast)$, is a stable chain complex: i.e., $\im H_\ast\subset\B_{00}^\Integer$.
\item The image $H_\ast(E_\ast,\del_\ast)=0\in\B_0^\Integer$ for any exact complex $(E_\ast,\del_\ast)\in\A_0^\Integer$.

\end{enumerate}
(Note that we have already constructed a concrete homology functor for modules, as can be seen from Definitions \ref{HlogyDef} and \ref{HlogyDef2}, and Corollary \ref{ModHomInvCor}.)

A complex $C_\ast=(C_\ast,\del_\ast)\in\A_0^\Integer$ is \ul{$H_\ast$-exact} (or \ul{homologically exact} with respect to $H_\ast$) if its image $H_\ast(C_\ast,\del_\ast)=0$ (the \ul{zero complex} $0:=0_\ast=(0_\ast,0_\ast)\in\B_{00}^\Integer$). A morphism of complexes $f_\ast\in\Mor\A_0^\Integer$, $f_\ast:(C_\ast,\del_\ast)\ra (C'_\ast,\del'_\ast)$, is an \ul{$H_\ast$-quasiisomorphism}, or an \ul{$H_\ast$-qis}, (making $C_\ast$ and $C'_\ast$ \ul{quiasi-isomorphic}, written $C_\ast\simeq_{\txt{qis}}C'_\ast$) if $H_\ast(f_\ast)$ is an isomorphism in $\B_0^\Integer$. (\blue{footnote}\footnote{As an equivalence of chain complexes, quasi-isomorphism is a generalization of homotopy equivalence.})

Let $A\in\Ob\A$ be a fixed object. Consider the contravariant functor
\bea
\Mor_\A(-,A):\A_0^\Integer\ra Ab_0^{\Integer}\subset Sets_0^{\Integer},~~\Big((C_\ast,\del_\ast)\sr{f_\ast}{\ral}(C'_\ast,\del'_\ast)\Big)\mapsto \Big((C^\ast,\del^\ast)\sr{f^\ast}{\lal}(C'{}^\ast,\del'{}^\ast)\Big)\nn
\eea
where $f^\ast:=\Mor_\A(f_\ast,A):=\circ f_\ast$ (right composition with $f_\ast$), and
\bea
(C^\ast,\del^\ast):=\Mor_\A\big((C_\ast,\del_\ast),A\big),~~~~C^n:=\Mor_\A\big(C_n,A\big),~~~~\del^n:=\Mor_\A\big(\del_{n+1},A\big):=\circ\del_{n+1},\nn
\eea
is the \ul{$A$-cochain complex} of the chain complex $(C_\ast,\del_\ast)$. In diagram form,
\bea
&&\bt
C_\ast=(C_\ast,\del_\ast):~\cdots\ar[r,"\del_{n+2}"]&C_{n+1}\ar[r,"\del_{n+1}"]&C_n\ar[r,"\del_n"]&C_{n-1}\ar[r,"\del_{n-1}"]&\cdots
\et\nn\\
\label{CohomDefEq}&&\bt
C^\ast=(C^\ast,\del^\ast):~\cdots\ar[from=r,"\del^{n+1}"',"\substack{~\\=~\circ\del_{n+2}}"]&C^{n+1}\ar[from=r,"\del^n"',"\substack{~\\=~\circ\del_{n+1}}"]&C^n\ar[from=r,"\del^{n-1}"',"\substack{~\\=~\circ\del_n}"]&C^{n-1}\ar[from=r,"\del^{n-2}"',"\substack{~\\=~\circ\del_{n-1}}"]&\cdots
\et
\eea

If there exists an imbedding $E:\A\hookrightarrow Ab$ (which therefore induces an imbedding $E:\A_0^\Integer\hookrightarrow Ab_0^\Integer$) with respect to which $H_\ast:\A_0^\Integer\ra\B_0^\Integer$ lifts to a unique homology functor $\wt{H}_\ast:Ab_0^\Integer\ra\B_0^\Integer$ in the sense that we have a commutative diagram

\bea
\adjustbox{scale=0.9}{\bt
\A_0^\Integer\ar[d,"{H_\ast}"']\ar[rr,hook,"E"]&& Ab_0^\Integer\ar[dll,dashed,"\exists !~{\wt{H}_\ast}"]\\
\B_0^\Integer
\et}~~~~~~~~H_\ast=\wt{H}_\ast\circ E\nn
\eea
then the image $\wt{H}_\ast(C^\ast,\del^\ast)=\wt{H}_\ast\Big(\Mor_\A\big((C_\ast,\del_\ast),A\big)\Big)=\big(\wt{H}_\ast\circ \Mor_\A(-,A)\big)(C_\ast,\del_\ast)$ is the associated $(A,H_\ast)$-\ul{cohomology} of the chain complex $(C_\ast,\del_\ast)$, and the contravariant functor
\bea
H^\ast:=\wt{H}_\ast\circ \Mor_\C(-,A):\A_0^\Integer\ra\B_0^\Integer\nn
\eea
is the associated $(A,H_\ast)$-\ul{cohomology functor}. (\blue{footnote}\footnote{When $\A=R$-mod (a subcategory of $Ab$) for some ring $R$, then we will simply stay within $R$-mod and set $\wt{H}_\ast:=H_\ast$.})

Based on the above discussion, and especially (\ref{CohomDefEq}), the \ul{$n$th cohomology functor} of the category $\A$ is the contravariant functor $H^n:=\wt{H}_n\circ \Mor_\C(-,A):\A_0^\Integer\ra\A$, which is explicitly given by the following:
\bea
H^n\Big((C_\ast,\del_\ast)\sr{f_\ast}{\ral}(C'_\ast,\del'_\ast)\Big)= \wt{H}_n\Big((C^\ast,\del^\ast)\sr{f^\ast}{\lal}(C'{}^\ast,\del'{}^\ast)\Big)=\wt{H}_n(C^\ast,\del^\ast)\sr{\wt{H}_n(f^\ast)}{\lal}\wt{H}_n(C'{}^\ast,\del'{}^\ast),\nn
\eea
where
\bit
\item[(1)] {\footnotesize $H^n(C_\ast,\del_\ast):=\wt{H}_n(C^\ast,\del^\ast)=\wt{H}_n\Big(\Mor_\A\big((C_\ast,\del_\ast),A\big)\Big)=\coker\left(\im\del^{n-1}\hookrightarrow\ker\del^n\right)={\ker\del^n\over\im\del^{n-1}}={Z^n(C^\ast,\del^\ast)\over B^n(C^\ast,\del^\ast)}$} is the \ul{$n$th cohomology (group)} of the chain complex $(C_\ast,\del_\ast)$.
\item[(2)] {\footnotesize $Z^n(C_\ast,\del_\ast):=\ker\del^n=\ker(\circ\del_{n+1})=\left\{\vphi\in C^n:\del^n(\vphi)=\vphi\circ\del_{n+1}=0\right\}=\left\{\vphi\in C^n:\vphi|_{B_n(C_\ast,\del_\ast)}=0\right\}$} is the \ul{$n$th cocycle (group)} of the chain complex $(C_\ast,\del_\ast)$.
\item[(3)] {\footnotesize $B^n(C_\ast,\del_\ast):=\im\del^{n-1}=\im(\circ\del_n)=\left\{\del^{n-1}(\vphi)=\vphi\circ\del_n:\vphi\in C^{n-1}\right\}=C^{n-1}\circ\del_n\subset Z^n(C_\ast,\del_\ast)$} is the \ul{$n$th coboundary (group)} of the chain complex $(C_\ast,\del_\ast)$.
\eit
Elements of $C^n$ are called \ul{$n$-cochains}, elements of $Z^n(C_\ast,\del_\ast)$ are called \ul{$n$-cocycles}, and elements of $B^n(C_\ast,\del_\ast)$ are called \ul{$n$-coboundaries}.
\end{dfn}

\begin{note*}~ $\left\{\substack{\txt{chain}\\ \txt{isomorphisms}}\right\}$ ~$\subset$~ $\left\{\substack{\txt{homotopy}\\ \txt{equivalences}}\right\}$ ~$\subset$~ $\left\{\substack{\txt{quasi-}\\ \txt{isomorphisms}}\right\}$
\end{note*}

\begin{rmk}[\blue{Homology functor notation}]
In this small section only, we have given $H^n$ (resp. $Z^n$, $B^n$, etc) a meaning that is ``slightly'' different from that of $H_n$ (resp. $Z_n$, $B_n$, etc). Later however, $H^n$ and $H_n$ will most often stand for exactly the same homology functor (after all, cohomology is defined in terms of homology). The appropriate interpretation will be clear from the surrounding context.
\end{rmk} 

%% file: parts/AlgebraCat/AlgebraCatS6.tex
\chapter{Exactness of Functors and some Basic Functors}\label{AlgebraCatS6}

\section{Exactness of functors}
\begin{dfn}[\textcolor{blue}{
\index{Left! exact functor}{Left exact functor},
\index{Left! exact contravariant functor}{Left exact contravariant functor},
\index{Right! exact functor}{Right exact functor},
\index{Right! exact contravariant functor}{Right exact contravariant functor},
\index{Exact! functor}{Exact functor},
\index{Exact! contravariant functor}{Exact contravariant functor}}]
Let $\A,\B$ be abelian categories, $F:\A\ra\B$ an additive functor, and $G:\A\ra\B$ an additive cofunctor.

(Note that these induce an additive functor $F:\A_0^\Integer\ra\B_0^\Integer$ and an additive cofunctor $G:\A_0^\Integer\ra\B_0^\Integer$.)
\begin{enumerate}[leftmargin=0.9cm]
\item $F$ is \ul{left exact} if exactness of {\footnotesize $0\ra A\ra A'\ra A''$} implies exactness of {\footnotesize $0\ra F(A)\ra F(A')\ra F(A'')$}.

$G$ is \ul{left exact} if exactness of {\footnotesize $A\ra A'\ra A''\ra 0$} implies exactness of {\footnotesize $0\ra G(A'')\ra G(A')\ra G(A)$}.

\item $F$ is \ul{right exact} if exactness of {\footnotesize $A\ra A'\ra A''\ra 0$} implies exactness of {\footnotesize $F(A)\ra F(A')\ra F(A'')\ra 0$}.

    $G$ is \ul{right exact} if exactness of {\footnotesize $0\ra A\ra A'\ra A''$} implies  exactness of {\footnotesize $G(A'')\ra G(A')\ra G(A)\ra 0$}.

\item $F$ (resp. $G$) is \ul{exact} if both left exact and right exact.
\end{enumerate}
\end{dfn}

\begin{prp}
Let $\A,\B$ be abelian categories, $F:\A\ra\B$ an additive functor, and $G:\A\ra\B$ an additive cofunctor.
\bit[leftmargin=0.9cm]
\item[(a)] $F$ is exact iff exactness of {\small $0\ra A\ra A'\ra A''\ra 0$} implies exactness of {\small$0\ra F(A)\ra F(A')\ra F(A'')\ra 0$}.

\item[(b)] $G$ is exact iff exactness of {\small $0\ra A\ra A'\ra A''\ra 0$} implies exactness of {\small $0\ra G(A'')\ra G(A')\ra G(A)\ra 0$}.
\item[(c)] $F$ (resp. $G$) is exact iff it preserves exactness of all sequences in $\A_0^\Integer$ (when mapped to $\B_0^\Integer$).
\eit
\end{prp}
\begin{proof}
(a) and (b) are immediate from the definitions. For (c), break (or factor) exact sequences into short exact sequences (see Remarks \ref{SeqFactRmks}).
\end{proof}

\begin{rmk}[\textcolor{blue}{Exact functors commute with homology}]
Let $\A,\B$ be abelian categories and $F:\A\ra\B$ an exact functor. Then $F$ (being both left exact and right exact) commutes with both $\ker$ and $\coker$, i.e.,
\bea
\label{HomolExtEq}F\ker(f)\cong\ker F(f)~~\txt{and}~~F\coker(f)\cong\coker F(f),~~~~\txt{for any}~~f\in\Mor\A,
\eea
and so commutes with homology as well: Recall that
\bea
H_n(C_\ast,\del_\ast):=\coker(\im\del_{n+1}\hookrightarrow\ker\del_n)=\coker\big(\ker(C_n\twoheadrightarrow\coker\del_{n+1})\hookrightarrow\ker\del_n\big).\nn
\eea
That is, given instances of the $n$th homology functor $H_n:\A_0^\Integer\ra\A$ and $H_n:\B_0^\Integer\ra\B$, we have
\[
F\circ H_n\cong H_n\circ F:\A_0^\Integer\ra\A
\]
(where $F:\A_0^\Integer\ra\B_0^\Integer$ is the obvious instance of $F:\A\ra\B$ on complexes).

To see this, let $A_\ast\sr{f_\ast}{\ral}B_\ast$ be a morphism in $\A_0^\Integer$. Then by first applying $H_n$ on $A_\ast\sr{f_\ast}{\ral}B_\ast$ (before applying $F$) we get the following mappings of diagrams:
\[\adjustbox{scale=0.7}{\bt
0      & 0\\
H_n(A_\ast)\ar[u]\ar[r,"H_n(f_\ast)"]           & H_n(B_\ast)\ar[u]\\
Z_n(A_\ast)\ar[u,two heads]\ar[r,"Z_n(f_\ast)"] & Z_n(B_\ast)\ar[u,two heads]\\
B_n(A_\ast)\ar[u,hook]\ar[r,"B_n(f_\ast)"]      & B_n(B_\ast)\ar[u,hook]\\
0\ar[u]                                         & 0\ar[u]
\et}~~~~\sr{F}{\ral}~~~~
\adjustbox{scale=0.7}{\bt
0      && 0\\
F\circ H_n(A_\ast)\ar[u]\ar[rr,"F\circ H_n(f_\ast)"]           && F\circ H_n(B_\ast)\ar[u]\\
F\circ Z_n(A_\ast)\ar[u,two heads]\ar[rr,"F\circ Z_n(f_\ast)"] && F\circ Z_n(B_\ast)\ar[u,two heads]\\
F\circ B_n(A_\ast)\ar[u,hook]\ar[rr,"F\circ B_n(f_\ast)"]      && F\circ B_n(B_\ast)\ar[u,hook]\\
0\ar[u]                                         && 0\ar[u]
\et}~~~~\sr{(\ref{HomolExtEq})}{\cong}~~~~
\adjustbox{scale=0.7}{\bt
0      && 0\\
F\circ H_n(A_\ast)\ar[u]\ar[rr,"F\circ H_n(f_\ast)"]           && F\circ H_n(B_\ast)\ar[u]\\
Z_n\circ F(A_\ast)\ar[u,two heads]\ar[rr,"Z_n\circ F(f_\ast)"] && Z_n\circ F(B_\ast)\ar[u,two heads]\\
B_n\circ F(A_\ast)\ar[u,hook]\ar[rr,"B_n\circ F(f_\ast)"]      && B_n\circ F(B_\ast)\ar[u,hook]\\
0\ar[u]                                         && 0\ar[u]
\et}\nn
\]
On the other hand, by first applying $F$ on $A_\ast\sr{f_\ast}{\ral}B_\ast$, we get the morphism $F(A_\ast)\sr{F(f_\ast)}{\ral}F(B_\ast)$ in $\B_0^\Integer$, on which $H_n$ then gives the following commutative diagram (which is isomorphic to the last diagram above):
\[
\adjustbox{scale=0.7}{\bt
0      && 0\\
H_n\circ F(A_\ast)\ar[u]\ar[rr,"H_n\circ F(f_\ast)"]           && H_n\circ F(B_\ast)\ar[u]\\
Z_n\circ F(A_\ast)\ar[u,two heads]\ar[rr,"Z_n\circ F(f_\ast)"] && Z_n\circ F(B_\ast)\ar[u,two heads]\\
B_n\circ F(A_\ast)\ar[u,hook]\ar[rr,"B_n\circ F(f_\ast)"]      && B_n\circ F(B_\ast)\ar[u,hook]\\
0\ar[u]                                         && 0\ar[u]
\et}~~~~\cong~~~~
\adjustbox{scale=0.7}{\bt
0      && 0\\
F\circ H_n(A_\ast)\ar[u]\ar[rr,"F\circ H_n(f_\ast)"]           && F\circ H_n(B_\ast)\ar[u]\\
Z_n\circ F(A_\ast)\ar[u,two heads]\ar[rr,"Z_n\circ F(f_\ast)"] && Z_n\circ F(B_\ast)\ar[u,two heads]\\
B_n\circ F(A_\ast)\ar[u,hook]\ar[rr,"B_n\circ F(f_\ast)"]      && B_n\circ F(B_\ast)\ar[u,hook]\\
0\ar[u]                                         && 0\ar[u]
\et}\nn
\]
\end{rmk}

\section{Adjoint functors, Tensor functors, and Mor functors}
\begin{dfn}[\textcolor{blue}{
\index{Adjoint relation of categories}{{Adjoint relation of categories}},
\index{Adjoint categories}{Adjoint categories},
\index{Adjoint functors (Adjoint pair)}{Adjoint functors (Adjoint pair)},
\index{Left! adjoint (Right-adjoined functor)}{Left adjoint (Right-adjoined functor)},
\index{Right! adjoint (Left-adjoined functor)}{Right adjoint (Left-adjoined functor)},
\index{Tensor! bifunctor (Tensor product) of a category}{{Tensor bifunctor (Tensor product) of a category}},
\index{Tensor! functor of an object}{Tensor functor associated with an object},
\index{Tensor! bifunctor (Tensor product)}{{Tensor bifunctor (Tensor product)}}
}]
A functor $F:\C\ra\D$ is an \ul{adjoint relation} of categories (making $\C,\D$ \ul{adjoint categories}, with $\D$ an \ul{adjoint} of $\C$, denoted by $\C^\ast$ when unique) if there exists a functor $G:\D\ra\C$ (making $F,G$ \ul{adjoint functors}, with $G$ an \ul{adjoint} of $F$, denoted by $F^\ast$ when unique) such that the following associated bifunctors are isomorphic, i.e., ~$L_\D F\cong R_\C G:\C\times\D\ra\txt{Sets}$,~ where
\bea
&&L_\D F:=\Mor_\D(F(-),-)=\Mor_\D(-,-)\circ(F\times id_\D)~~:~~\C\times\D\ra\txt{Sets},\nn\\
&&R_\C G:=\Mor_\C(-,G(-))=\Mor_\C(-,-)\circ(id_\C\times G)~~:~~\C\times\D\ra\txt{Sets}.\nn
\eea
If $F,G:\C\ra\D$ are adjoint functors, we also say $F,G$ form an \ul{adjoint pair} $(F,G)$, where $F$ is called a \ul{left adjoint} of $G$, and $G$ a \ul{right adjoint} of $F$. In general, a \ul{left adjoint} or \ul{right-adjoined functor} (resp. \ul{right adjoint} or \ul{left-adjoined functor}) is any functor for which a right adjoint (resp. left adjoint) exists.

Let $\C$ be a category. A \ul{tensor bifunctor} (\ul{tensor product}) of $\C$ is an (associative symmetric) bifunctor
\bea
\otimes_\C=-\otimes_\C-:\C\times\C\ra\C,~~(C_1,C_2)\sr{(f,f')}{\ral}(C'_1,C'_2)~~\mapsto~~C_1\otimes_\C C_2\sr{f\otimes f'}{\ral}C'_1\otimes_\C C'_2\nn
\eea
such that for every object $C\in\Ob\C$, the associated functor
\bea
F:=\otimes_\C C:=-\otimes_\C C:\C\ra\C,~~A\sr{f}{\ral}B~~\mapsto~~A\otimes_\C C\sr{f\otimes id_C}{\ral}B\otimes_\C C\nn
\eea
has a unique right adjoint $F^\ast=(\otimes_\C C)^\ast:\C\ra\C$, in the sense that we have a unique isomorphism of functors
\bea
\Mor_\C\big(-\otimes_\C C,-\big)\cong \Mor_\C\big(-,(\otimes_\C C)^\ast(-)\big).\nn
\eea

More generally, let $\A,\B,\C$ be categories. A \ul{tensor bifunctor} (\ul{tensor product}) is a bifunctor of the form
\bea
\otimes=-\otimes-:\A\times\B\ra\C,~~(A,B)\sr{(f,g)}{\ral}(A',B')~~\mapsto~~A\otimes B\sr{f\otimes g}{\ral}A'\otimes B'\nn
\eea
such that for all objects $A\in\Ob\A$ and $B\in\Ob\B$, the associated functors
\bea
&& F_A:=A\otimes :=A\otimes-:\B\ra\C,~~B\sr{g}{\ral}B'~~\mapsto~~A\otimes B\sr{id_A\otimes g}{\ral}A\otimes B',\nn\\
&& F_B:=\otimes B:=-\otimes B:\A\ra\C,~~A\sr{f}{\ral}A'~~\mapsto~~A\otimes B\sr{f\otimes id_B}{\ral}A'\otimes B\nn
\eea
have respective unique right adjoints $F_A^\ast=(A\otimes)^\ast:\C\ra\B$ and $F_B^\ast=(\otimes B)^\ast:\C\ra\A$, in the sense that we have respective unique functor isomorphisms
\[
\Mor_\C\big(A\otimes-,-\big)\cong \Mor_\B\big(-,(A\otimes)^\ast(-)\big)~~~~\txt{and}~~~~\Mor_\C\big(-\otimes B,-\big)\cong \Mor_\A\big(-,(\otimes B)^\ast(-)\big).
\]
\end{dfn}

\begin{rmk}[\blue{Application terminology}]
In applications (e.g., $\C=$ Modules), when we say ``a tensor bifunctor $\otimes_\C:\C \times\C\ra\C$ of $\C$'', we will actually mean a tensor bifunctor system $\otimes_{ij}^k:\C_i\times\C_j\ra\C_k$ for specific/relevant subcategories $\C_i,\C_j,\C_k\subset\C$, where each $\otimes_{ij}^k$ is defined either as a restriction of an existing tensor bifunctor $\otimes_\C:\C\times\C\ra\C$ or directly from the above definition.
\end{rmk}

\begin{dfn}[\textcolor{blue}{
\index{Mor functors of an object}{Mor functors associated with an object},
\index{Mor bifunctor of a category}{{Mor bifunctor of a category}}}]
Let $\C$ be a category and $C\in\Ob\C$. The \ul{mor functor} associated with $C$ is the functor
\bea
\Mor_\C(C,-):\C\ra\txt{Sets},~~(A\sr{f}{\ral}B)
~\longmapsto~
(A_\ast\sr{f_\ast}{\ral} B_\ast):=\Big(
\bt[column sep=large] \Mor_\C(C,A)\ar[r,"{\Mor_\C(C,f)}","~~:=f\circ"'] & \Mor_\C(C,B)\et\Big),\nn
\eea
where $f_\ast(g):=f\circ g$, for $g\in \Mor_\C(B,C)$, is left composition with $f$. Similarly, the \ul{mor cofunctor} associated with $C$ is the cofunctor
\bea
\Mor_\C(-,C):\C\ra\txt{Sets},~~(A\sr{f}{\ral}B)
~\longmapsto~
(A^\ast\sr{f^\ast}{\lal}B^\ast):=
\Big(\bt[column sep=large] \Mor_\C(A,C)\ar[from=r,"{\Mor_\C(f,C)}"',"~~:=\circ f"] & \Mor_\C(B,C)\et\Big),\nn
\eea
where $f^\ast(g):=g\circ f$, for $g\in B^\ast:=\Mor_\C(B,C)$, is right composition with $f$. The \ul{mor bifunctor} of $\C$ is
\bea
&&\Mor_\C(-,-):\C^{op}\times\C\ra\txt{Sets},~~\Big(~\big(A\sr{f^{op}}{\ral}A'~,~B\sr{g}{\ral}B'\big)~\Big)
~=~\Big(\bt(A,B)\ar[r,"{(f^{op},g)}"] & (A',B')\et\Big)\nn\\
&&~~~~~~~~~~~~~\longmapsto~
\Big(\bt A\ast B\ar[r,"f^{op}\ast g"] & A'\ast B'\et\Big):=
\Big(\bt[column sep=large] \Mor_\C(A,B)\ar[r,"{\Mor_\C(f^{op},g)}","~~:=f^\ast\circ g_\ast"'] & \Mor_\C(A',B')\et\Big),\nn
\eea
where $(f^{op}\ast g)(h):=(f^\ast\circ g_\ast)(h)=g\circ h\circ f$, for $h\in A\ast B:=\Mor_\C(A,B)$, is both left and right composition with $g$ and $f$ respectively.
\end{dfn}
\begin{rmk}[\textcolor{blue}{\index{Additivity of! tensor and mor functors}{Additivity of tensor and mor functors}}]
If $\A$ is an additive category, then the mor bifunctor of $\A$ is clearly additive in each argument, and is thus said to be \ul{bi-additive}. The same is true of the tensor bifunctor of $\A$ by Lemma \ref{TensorBiadd} (on page \pageref{TensorBiadd}).
\end{rmk}

\section{Tensor bifunctor and Resolutions in an Additive Category}
\begin{dfn}[\textcolor{blue}{
\index{Bimorphism}{Bimorphism},
\index{Biadditivity}{Biadditivity},
\index{Tensor! bifunctor of an additive category}{{Tensor bifunctor of an additive category}},
\index{Tensor! product of objects}{Tensor product of objects},
\index{Universal! property of tensor}{Universal property of tensor},
\index{Flat object}{Flat object},
\index{Projective! object}{Projective object},
\index{Injective! object}{Injective object},
\index{Projective! resolution (process/functor)}{Projective resolution (process/functor)},
\index{Deleted projective resolution}{Deleted projective resolution},
\index{Injective! resolution (process/functor)}{Injective resolution (process/functor)},
\index{Deleted injective resolution}{Deleted injective resolution},
\index{Free! resolution in $R$-mod}{Free resolution in $R$-mod},
\index{Category with enough projectives}{Additive category with enough projectives},
\index{Category with enough injectives}{Additive category with enough injectives}}]

Let $\A$ be an additive category, $S$ be an associative set, and $f:S\ra\bigsqcup_{A,B\in\Ob\A}\Mor_\A(A,A)\times\Mor_\A(B,B),~s\mapsto f(s)$ a map (with $f(s)=(f_A(s),f_B(s))$ if $f(s)\in\Mor_\A(A,A)\times\Mor_\A(B,B)$). (\blue{footnote}\footnote{When $\A$ is some $R$-mod, then by considering elements $r\in R$ as maps $r:M\ra M,~m\mapsto rm$ of $R$-modules $M={}_RM$, we can then set $S:=R$ and $f(r):=(r,r)$ for each $r\in R$.}). The class of \ul{$(f,S)$-bimorphisms} in $\A$ is a collection of morphisms $Bim\A=\bigsqcup_{A,B,C\in\Ob\A}Bim_\A(A\times B,C)$, for morphism classes $Bim_\A(A\times B,C)\subset Mor_\A(A\times B,C)$ such that the following hold: (\blue{footnote}\footnote{Note that because $\A$ is pre-additive, members of $Bim\A$ (like all morphisms of $\A$) are \ul{biadditive} in the sense that for any $\theta\in Bim_\A(A\times B,C)$ and any morphisms $f_1,f_2\in\Mor_\A(U,A\times B)$, $g_1,g_2\in \Mor_\A(C,V)$, the compositions $U\sr{f_1+f_2}{\ral} A\times B\sr{\theta}{\ral}C$ and $A\times B\sr{\theta}{\ral}C\sr{g_1+g_2}{\ral}V$ (whenever they exist) satisfy
\bea
\theta\circ(f_1+f_2)=\theta\circ f_1+\theta\circ f_2~~~~\txt{and}~~~~(g_1+g_2)\circ\theta=g_1\circ\theta+g_2\circ\theta.\nn
\eea})
\begin{enumerate}[leftmargin=0.9cm]
\item Members of $Bim\A$ are \ul{$(f,S)$-balanced} in the sense for any $\theta:A\times B\ra C$ in $Bim_\A(A\times B,C)$,
\[
\theta\circ(q_Af_A(s)p_A+q_Bp_B)=\theta\circ(q_Ap_A+q_Bf_B(s)p_B),~~~~\txt{for all $s\in S$},
\]
where $q_A:A\hookrightarrow A\times B$, $q_B:B\hookrightarrow A\times B$ are the canonical injections and $p_A:A\times B\twoheadrightarrow A$, $p_B:A\times B\twoheadrightarrow B$ the canonical projections.
\item Members of $Bim\A$ are \ul{bimorphic} in the sense there exists a full imbedding $E:\A\hookrightarrow \C\subset Sets$ such that (i) finite products are preserved and (ii) for each $\theta\in Bim_\A(A\times B,C)$, the image
\[
E(\theta):E(A\times B)\cong E(A)\times E(B)\ra E(C),~(x,y)\mapsto E(\theta)(x,y)
\]
satisfies the following: For every $x\in E(A)$ and $y\in E(B)$,
\begin{align}
&E(\theta)(x,-):E(B)\ra E(C),~y\mapsto E(\theta)(x,y)~~~~\txt{lies in}~~~Hom_\C\big(E(B),E(C)\big),\nn\\
&E(\theta)(-,y):E(A)\ra E(C),~x\mapsto E(\theta)(x,y)~~~~\txt{lies in}~~~Hom_\C\big(E(A),E(C)\big).\nn
\end{align}
\end{enumerate}

When it is helpful, we may also express the above condition for a morphism $\theta\in Hom_\A(A\times B,C)$ to be \ul{bimorphic} in the following form (with an appropriate interpretation if needed (\blue{footnote}\footnote{For example (a situation that is actually sufficiently general for an abelian $\A$ via say the Freyd-Mitchell imbedding (Theorem \ref{FreydMitchIT}) is the following): If $\A\subset\txt{Sets}$, then $\theta\circ(-\times B)\subset \Mor_\A(A,C)$ simply means $\theta(-,b)\in \Mor_\A(A,C)$ for each $b\in B$, and similarly, $\theta\circ(A\times-)\subset \Mor_\A(B,C)$ means $\theta(a,-)\in \Mor_\A(B,C)$ for each $a\in A$.})):
\begin{align}
&\theta_A(B):=\theta\circ(-\times B)\subset \Mor_\A(A,C)~~~~\txt{and}~~~~\theta_B(A):=\theta\circ(A\times-)\subset \Mor_\A(B,C),~~~~\txt{where}~~~~\nn\\
&\theta\circ(-\times B):A\sr{-\times B}{\ral}A\times B\sr{\theta}{\Longrightarrow}C~~~~\txt{and}~~~~\theta\circ(A\times -):B\sr{A\times-}{\ral}A\times B\sr{\theta}{\Longrightarrow}C.\nn
\end{align}
\[
\adjustbox{scale=1.0}{\bt
A\ar[dr,shift left,hook,"q_A"]\ar[from=dr,shift left,two heads,"p_A"] &           & B\ar[dl,shift right,hook,"q_B"']\ar[from=dl,shift right,two heads,"p_B"'] \\
  & A\times B\ar[d,"\theta"] & \\
  &    C\ar[from=uul,bend right,dashed,"{\Mor_\A(A,C)\supset\theta_A(B)}"']\ar[from=uur,bend left,dashed,"{\theta_B(A)\subset \Mor_\A(B,C)}"]      &
\et}
\]

A \ul{tensor bifunctor} of $\A$ is a bifunctor (\blue{footnote}\footnote{This bifunctor will automatically be biadditive, associative, and symmetric: See Lemma \ref{TensorBiadd}.})
\bea
\txt{Tensor}_\A(-,-)=-\otimes_\A-:\A\times\A\ra\A,~~(A,B)\sr{(f,g)}{\ral}(A',B')~~\mapsto~~A\otimes_\A B\sr{f\otimes g}{\ral}A'\otimes_\A B'\nn
\eea
such that any two objects $A,B\in\Ob\A$ satisfy the following:
\bit[leftmargin=0.9cm]
\item[(1)] There exists a bimorphism $\tau_{AB}:A\times B\ra A\otimes_\A B$ (called the \ul{tensor product} of $A$ and $B$) with the following property.
\item[(2)] \ul{Universal property of tensor}: For any bimorphism $\theta:A\times B\ra C$, (i) there exists a unique morphism $\wt{\theta}:A\otimes_\A B\ra C$ satisfying $\theta=\wt{\theta}\circ\tau_{AB}$,
\bea\adjustbox{scale=0.9}{\bt
A\ar[ddr,bend right,dashed]\ar[dr,hook,"q_A"] & & B\ar[dl,hook,"q_B"']\ar[ddl,bend left,dashed,near start] && \\
  & A\times B\ar[d,"\theta"']\ar[rrr,near end,"\tau_{AB}"]&&& A\otimes_\A B\ar[dlll,dashed,"\exists!~\wt{\theta}"]\\
  & C &&&
\et},~~~~~~~~\theta=\wt{\theta}\circ\tau_{AB},\nn
\eea
and, if it is not automatically the case, we also require that (ii) for all objects $A,B,C\in\Ob\A$, the resulting/associated isomorphism of abelian groups
\[
Bim_\A(A\times B,C)\cong \Mor_\A(A\otimes_\A B,C)
\]
is natural (i.e., behaves as a natural transformation) in each of the pairs of objects $(A,B)$, $(A,C)$, and $(B,C)$.
\eit

That is, a tensor bifunctor of $\A$ is a bifunctor
\bea
\txt{Tensor}_\A(-,-)=-\otimes_\A-:\A^2\ra\A,~~(A,B)\sr{(f,g)}{\ral}(A',B')~~\mapsto~~A\otimes_\A B\sr{f\otimes g}{\ral}A'\otimes_\A B\nn
\eea
such that for every object $X\in\Ob\A$, we have the following isomorphisms of bifunctors $A^2\ra Ab^2$, namely,
\begin{align}
& Bim_\A(-\times -,X)\cong \Mor_\A(-\otimes_\A-,X),~~~~Bim_\A(-\times X,-)\cong \Mor_\A(-\otimes_\A X,-),~~~~\txt{and}~~~~\nn\\
&Bim_\A(X\times -,-)\cong \Mor_\A(X\otimes_\A-,-).\nn
\end{align}

An object $A\in\Ob\A$ is called \ul{flat} if the associated functor
\[
-\otimes_\A A:\A\ra\A,~~\big(B\sr{f}{\ral}C\big)\mapsto \big(B\otimes_\A A\sr{f\otimes id_A}{\ral}C\otimes A\big)
\]
maps monic morphisms to monic morphisms.

An object $P\in\Ob\A$ is \ul{projective} if its associated mor functor $\Mor_\A(P,-)$ maps epic morphisms to epic morphisms. Equivalently, given an epic morphism $g\in \Mor_\A(B,C)$, every morphism $h\in \Mor_\A(P,C)$ factors through a morphism $h'\in \Mor_\A(P,B)$ in the form $h=g\circ h'$.
\bea\adjustbox{scale=0.9}{\bt
                      && P\ar[dll,dashed,"\exists~h'"']\ar[d,"h"] & \\
B\ar[rr,two heads,"g"] && C\ar[r]                                     & 0
\et}~~~~~~~~h=g\circ h'.\nn
\eea
An object $I\in\Ob\A$ is \ul{injective} if its associated mor cofunctor $\Mor_\A(-,I)$ maps monic morphisms to epic morphisms. Equivalently, given a monic morphism $f\in \Mor_\A(A,B)$, every morphism $h\in \Mor_\A(A,I)$ factors through a morphism $h'\in \Mor_\A(B,I)$ in the form $h=h'\circ f$.
\bea\adjustbox{scale=0.9}{\bt
0\ar[r] &  A\ar[d,"h"]\ar[rr,hook,"f"]            && B \\
        &  I\ar[from=urr,dashed,"\exists~h'"] &&
\et}~~~~~~~~h=h'\circ f.\nn
\eea

Let $A\in\Ob\A$. A \ul{projective resolution} of $A$ is an exact sequence of the form
\bea\bt[column sep=small]
A_\ast:\cdots\ar[r]& P^A_2\ar[r,"\del_2"]& P^A_1\ar[r,"\del_1"]& P^A_0\ar[r,"\vep_A"]& A\ar[r]& 0
\et~~~~=~~~~
\bt[column sep=small]
P_\ast^A\ar[r,"\vep_A"]& A\ar[r]& 0
\et\nn
\eea
such that each $P_n\in\Ob\A$ is a projective object. The associated chain complex
\bea\bt[column sep=small]
P_\ast^A:~~\cdots\ar[r]& P^A_2\ar[r,"\del_2"]& P^A_1\ar[r,"\del_1"]& P^A_0\ar[r]& 0
\et\nn
\eea
(which is exact at $P_n$ for all $n\neq 0$) is called a \ul{deleted projective resolution} of $A$.

An \ul{injective resolution} of $A$ is an exact sequence of the form
\bea
&&\bt[column sep=small]
A^\ast:~~0\ar[r]& A\ar[r,"\xi_A"]& I_A^0\ar[r,"\del^0"]& I_A^1\ar[r,"\del^1"]& I_A^2\ar[r,"\del^2"]& \cdots
\et~~~~=~~~~
\bt[column sep=small]
0\ar[r] & A\ar[r,"\xi_A"] & I_A^\ast
\et\nn\\
&&\txt{or}\nn\\
&&\bt[column sep=small]
A^\ast:~~\cdots\ar[from=r,"\del^2"']& I_A^2\ar[from=r,"\del^1"']& I_A^1\ar[from=r,"\del^0"']& I_A^0\ar[from=r,"\xi_A"']& A\ar[from=r]& 0
\et~~~~=~~~~
\bt[column sep=small]
I_A^\ast\ar[from=r,"\xi_A"'] & A\ar[from=r] & 0
\et\nn
\eea
such that each $I^n\in\Ob\A$ is an injective object. The associated chain complex
\bea
\bt[column sep=small]
I_A^\ast:~~0\ar[r]& I_A^0\ar[r,"\del^0"]& I_A^1\ar[r,"\del^1"]& I_A^2\ar[r,"\del^2"]& \cdots
\et
~~~~\txt{or}~~~~
\bt[column sep=small]
I_A^\ast:~~\cdots\ar[from=r,"\del^2"']& I_A^2\ar[from=r,"\del^1"']& I_A^1\ar[from=r,"\del^0"']& I_A^0\ar[from=r]& 0
\et\nn
\eea
(which is exact at $I^n$ for all $n\neq 0$) is called a \ul{deleted injective resolution} of $A$.

If $\A\subset\txt{Sets}$ (in which case, we may assume wlog that $\A$ $=$ $R$-mod, for some $R$, by Remarks \ref{FreydMitchRmk}), then a free object in $\A$ is (as verified in the Lemma \ref{FreeIsProj} below) necessarily a projective object: A projective resolution $P_\ast^A\ra A\ra 0$ in which each projective $P_n$ is free is called a \ul{free resolution} of $A$.

We say the category $\A$ has \ul{enough projectives} (resp. \ul{enough injectives}) if for every object $A\in\Ob\A$, there exists a projective object $P\in\Ob\A$ (resp. an injective object $I\in\Ob\A$) and an epic morphism \bt[column sep =small] P\ar[r,two heads]& A\et (a monic morphism \bt[column sep =small] A\ar[r,hook]& I\et).
\end{dfn}

\begin{rmk}[\textcolor{blue}{Existence of injective/projective resolutions}]
It follows immediately using Lemma \ref{CmpDecLmm} below that in an abelian category with enough projectives (resp. enough injectives) every object has a projective resolution (resp. an injective resolution).
\end{rmk}

\begin{lmm}[\textcolor{blue}{\index{Composition-Decomposition Lemma}{Composition-Decomposition Lemma for an abelian category}}]\label{CmpDecLmm}
Let $\A$ be an abelian category.
{\flushleft (1)} In $\A$, any two exact sequences
\bea
\cdots\ral X\sr{x}{\ral}Y\sr{y}{\ral}T\sr{}{\ral}0~~~~\txt{and}~~~~0\ral T\sr{u}{\ral}V\sr{v}{\ral}W\sr{w}{\ral}\cdots\nn
\eea
can be combined (head-to-tail) to get the following exact sequence (with the common end object $T$ removed):
\bea
\cdots\ral X\sr{x}{\ral}Y\sr{uy}{\ral}V\sr{v}{\ral}W\sr{w}{\ral}\cdots\nn
\eea
because monic (resp. epic) invariance of $\ker$ (resp. coker) implies $\ker(uy)\cong\ker y$ and $\coker(uy)\cong\coker u$, and so (i) $\ker(uy)\cong\ker y=\im x$, and
(ii) $\im(uy):=\ker\big(\cod(uy)\ra\coker(uy)\big)\cong\ker\big(\cod(u)\ra\coker(u)\big)=\im u=\ker v$.

{\flushleft (2)} Together with the unitarity theorem (which allows for factoring of a map), the above observation implies that every sequence/chain of morphisms (not necessary a complex) in $\A$
\[
\cdots\ral A^{n-1}\sr{f^{n-1}}{\ral}A^n\sr{f^n}{\ral}A^{n+1}\sr{f^{n+1}}{\ral}\cdots
\]
is a composition (head-to-tail) of short exact sequences.
\end{lmm}

\begin{lmm}\label{FreeIsProj}
Let $\A$ be an additive category. If $\A\subset\Sets$, then every free object $F$ in $\A$ is projective.
\end{lmm}
\begin{proof}
Let $g:B\ra C$ be an epimorphism in $\A$, $h\in\Mor_\A(F,C)$, and $F:=\langle X\rangle\supset X$ for a set $X$. Using the axiom of choice (or otherwise) define a map $\theta:X\ra B,~x\mapsto\theta(x)\in g^{-1}(h(x))$. Then by the definition of the free object $\langle X\rangle$ over the set $X$, a unique morphism $\theta':\langle X\rangle\ra B$ exists such that $\theta'|_X=\theta$.
{\small\bea\adjustbox{scale=0.9}{\bt
X\ar[dd,"\theta"']\ar[rr,hook,"i"] && \langle X\rangle=F\ar[ddll,dashed,bend right=10,"\exists!~\theta'"',"=:\!~h'"]\ar[dd,"h"] & \\
 &&  & \\
B\ar[rr,two heads,"g"] && C\ar[r] & 0
\et}~~~~g\circ h'|_X=g\circ\theta=h|_X~~\sr{(s)}{\Ra}~~g\circ h'=h,\nn
\eea}where step (s) holds by uniqueness in the universal property of the free object $\langle X\rangle$.
\end{proof}

\section{Properties of Tensor Functors of an Additive Category}
\begin{lmm}[\textcolor{blue}{
\index{Biadditivity of tensor}{Biadditivity of tensor},
\index{Associativity of tensor}{Associativity of tensor},
\index{Symmetry of tensor}{Symmetry of tensor}}]\label{TensorBiadd}
Let $\A$ be an additive category with a tensor bifunctor $\otimes_\A:\A^2\ra\A$. Let $A,B,C,C'\in\Ob\A$. Then the following hold:
\bit[leftmargin=0.9cm]
\item[(i)] The associated tensor functors below are additive.
\bea
&&-\otimes_\A B:\A\ra\A,~~C\sr{f}{\ral}C'~~\mapsto~~C\otimes_\A B\sr{f\otimes id_B}{\ral}C'\otimes_\A B\nn\\
&&A\otimes_\A-:\A\ra\A,~~C\sr{f}{\ral}C'~~\mapsto~~A\otimes_\A C\sr{id_A\otimes f}{\ral}A\otimes_\A C'\nn
\eea
\item[(ii)] Associativity: $A\otimes_\A(B\otimes_\A C)\cong (A\otimes_\A B)\otimes_\A C$.
\item[(iii)] Symmetry: $A\otimes_\A B\cong B\otimes_\A A$.
\eit
\end{lmm}
\begin{proof}
(i) Given any morphisms $f,f_1,f_2\in \Mor_\A(A,A')$, $g\in \Mor_\A(B,B')$, we get a unique morphism $f\times g:A\times B\ra A'\times B'$ as in the following commutative diagrams:
\bea\adjustbox{scale=0.8}{\bt
A'\ar[from=dr,"p'"'] &            & B'\ar[from=dl,"q'"] \\
  & A'\times B' &  \\
  & &\\
  & A \times B\ar[uu,dashed,"\exists!","{f\times g}"']  &  \\
A\ar[from=ur,"p"]\ar[uuuu,bend left=15,"f"] &            & B \ar[from=ul,"q"']\ar[uuuu,bend right=15,"g"']\\
\et}~~~~~~~~f\circ p=p'\circ(f\times g),~~g\circ q=q'\circ(f\times g).\nn
\eea
If $f=f_1+f_2$, then $f_1\circ p=p'\circ(f_1\times g)$, $f_2\circ p=p'\circ(f_2\times g)$, and $(f_1+f_2)\circ p=p'\circ\big((f_1+f_2)\times g\big)$, and so
\bea
&&p'\circ\big((f_1+f_2)\times g\big)=(f_1+f_2)\circ p=f_1\circ p+f_2\circ p=p'\circ(f_1\times g)+p'\circ(f_2\times g)\nn\\
&&~~~~=p'\circ(f_1\times g+f_2\times g),~~\Ra~~(f_1+f_2)\times g=f_1\times g+f_2\times g~~\txt{(by uniqueness)}.\nn
\eea
(The above result is also an immediate consequence of the additivity of limits as proved earlier.)

The morphism $f\times g$ in turn induces a unique morphism $f\otimes g:A\otimes B\ra A'\otimes B'$ as in the commutative diagram
\bea\adjustbox{scale=0.9}{\bt
A\times B\ar[d,"{f\times g}"']\ar[rr,"\tau_{AB}"]&& A\otimes_\A B\ar[ddll,dashed,"\exists!~f\otimes g"]\\
A'\times B'\ar[d,"\tau_{A'B'}"'] && \\
A'\otimes B' &&
\et}~~~~~~~~(f\otimes g)\circ\tau_{AB}=\tau_{A'B'}\circ(f\times g).\nn
\eea
In particular, with $f=f_1+f_2$, $B=B'$, $g=id_B$, we get $(f_1+f_2)\otimes id_B=(f_1\otimes id_B)+(f_2\otimes id_B)$, since
\bea
&& \big((f_1+f_2)\otimes id_B\big)\circ\tau_{AB}=\tau_{A'B}\circ\big((f_1+f_2)\times id_B\big)=\tau_{A'B}\circ\big(f_1\times id_B+f_2\times id_B\big)\nn\\
&&~~~~=\tau_{A'B}\circ(f_1\times id_B)+\tau_{A'B}\circ(f_2\times id_B)=(f_1\otimes id_B)\circ\tau_{AB}+(f_2\otimes id_B)\circ\tau_{AB}\nn\\
&&~~~~=\big((f_1\otimes id_B)+(f_2\otimes id_B)\big)\circ\tau_{AB}.\nn
\eea
This shows $-\otimes_\A B$ is additive. The additivity of $A\otimes_\A-$ follows similarly.
{\flushleft(ii)} By the universal property of tensor, the automatic fact $(A\times B)\times C\sr{h}{\cong} A\times (B\times C)$ gives the following commutative diagram, in which $\tau_{12},\tau_{21}$ are respectively the compositions
\begin{align}
&\tau_{12}:(A\times B)\times C\sr{\tau_{AB}\times id_C}{\ral}(A\otimes B)\times C\sr{\tau_{(A\otimes B)C}}{\ral}(A\otimes B)\otimes C,\nn\\
&\tau_{21}:A\times (B\times C)\sr{id_A\times\tau_{BC}}{\ral}A\times(B\otimes C)\sr{\tau_{A(B\otimes C)}}{\ral}A\otimes(B\otimes C),\nn
\end{align}
and $\theta_1,\theta_2$ are (by uniqueness) inverses of each other:
\bea\adjustbox{scale=0.9}{\bt
(A\times B)\times C\ar[d,"h"',"\cong"]\ar[rr,"\tau_{12}"]&& (A\otimes_\A B)\otimes_\A C\ar[ddll,shift right, bend left,dashed,"\exists!~\theta_1"']\ar[from=ddll,bend right,shift right,dashed,"\exists!~\theta_2"']\\
A\times (B\times C)\ar[d,"\tau_{21}"'] && \\
A\otimes_\A (B\otimes_\A C) &&
\et}~~~~~~~~\theta_1\circ\tau_{12}=\tau_{21}\circ h,~~~~\theta_2\circ\tau_{21}=\tau_{12}\circ h^{-1}.\nn
\eea
{\flushleft(iii)} Similarly (commutative diagram below), the automatic fact $A\times B\sr{t}{\cong} B\times A$ implies $A\otimes_\A B\cong B\otimes_\A A$.
\bea\adjustbox{scale=0.9}{\bt
A\times B\ar[d,"\cong","t"']\ar[rr,"\tau_{AB}"] && A\otimes_\A B\ar[ddll,dashed,shift right,"\exists!~h"'] \\
B\times A\ar[d,"\tau_{BA}"'] &&\\
B\otimes_\A A\ar[uurr,dashed,shift right,"\exists!~\wt{h}"'] & &
\et}\nn
\eea
\end{proof}

\section{Left Exactness of Mor Functors}
\begin{lmm}[\textcolor{blue}{\index{Left! exactness of mor functors}{Left exactness of mor functors}}]\label{PaExHomFn}
Let $\A$ be an abelian category and $X\in\Ob\A$. The mor functors $\Mor_\A(X,-)$ and $\Mor_\A(-,X)$ are left exact.

(Note that the existence of an adjoint is not needed in the following proof.)
\end{lmm}
\begin{proof}
Let $0\ra A\sr{f}{\ral}B\sr{g}{\ral}C\ra 0$ be a SES in $\A$. For fixed objects $X,X'\in\Ob\A$, apply the functors $\Mor_\A(X,-)$ and $\Mor_\A(-,X')$ on the SES to get the following sequences in $Ab$:
\bea
&&0\ra \Mor_\A(X,A)\sr{f_\ast}{\ral}\Mor_\A(X,B)\sr{g_\ast}{\ral}\Mor_\A(X,C),\nn\\
&&0\ra \Mor_\A(C,X')\sr{g^\ast}{\ral}\Mor_\A(B,X')\sr{f^\ast}{\ral}\Mor_\A(A,X').\nn
\eea
For any $h\in \Mor_\A(X,A)$ and $h'\in \Mor_\A(C,X')$, we have $f_\ast(h)=f\circ h=0$ and $g^\ast(h')=h'\circ g=0$, which imply $h=0$ and $h'=0$ (since $f$ is monic and $g$ is epic). Thus \ul{$f_\ast$ and $g^\ast$ are monic}. Also, because $g_\ast\circ f_\ast=(f\circ g)_\ast=0_\ast=0$ and $f^\ast\circ g^\ast=(g\circ f)^\ast=0^\ast=0$, we get \ul{$\im f_\ast\subset\ker g_\ast$ and $\im g^\ast\subset\ker f^\ast$}.

Next, by Remarks \ref{SeqFactRmks}, we can expand the given SES as follows:
\bea
\label{PEHFeq1}\adjustbox{scale=0.9}{\bt
  &  &  &  & X\ar[dl,dashed, near start,"\wt{h}"']\ar[d,"h"] &  &  \\
0\ar[r] & A\ar[rrr,bend right=40,"f=k_g\circ\theta\circ u_f"]\ar[r,"u_f","\cong"'] & \im f\ar[r,"\theta","\cong"'] & \ker g\ar[r,hook,"k_g"] & B\ar[r,two heads,"c_f"]\ar[rrr,bend left=40,"g=u_g\circ\theta'\circ c_f"] & \coker f\ar[r,"\theta'","\cong"'] & \coim g\ar[r,"u_g","\cong"'] & C\ar[r] & 0\\
  &  &  &  & X'\ar[from=ur,dashed,"\wt{h}'"]\ar[from=u,"h'"] &  &
\et}
\eea
Let $h\in\ker g_\ast$ (i.e., $g_\ast(h)=g\circ h=0$) and $h'\in\ker f^\ast$ (i.e., $f^\ast(h')=h'\circ f=0$). Then by the universal property of $\ker g$, we have a morphism $\wt{h}:X\ra\ker g$ such that
\bea
&&h=k_g\circ\wt{h}=k_g\circ(\theta\circ u_f)\circ(\theta\circ u_f)^{-1}\circ\wt{h}=f_\ast\left(u_f^{-1}\circ \theta^{-1}\circ\wt{h}\right)\in\im f_\ast,\nn\\
&&~~\Ra~~\ker g_\ast\subset\im f_\ast.\nn
\eea
Similarly, by the universal property of $\coker f$, we have a morphism $\wt{h}':\coker f\ra X'$ such that
\begin{align}
&h'=\wt{h}'\circ c_f=\wt{h}'\circ(u_g\circ\theta')^{-1}\circ(u_g\circ\theta')\circ c_f=g^\ast\left(\wt{h}'\circ(\theta')^{-1}\circ(u_g)^{-1}\right)\in\im g^\ast,\nn\\
&~~\Ra~~\ker f^\ast\subset\im g^\ast.\nn \qedhere
\end{align}
\end{proof}

\begin{crl}[\textcolor{blue}{\index{Projectivity and Injectivity Criteria}{Projectivity and Injectivity Criteria}}]
Let $\A$ be an abelian category and $A,B\in\Ob\A$. Then $A$ is projective $\iff$ ~$\Mor_\A(A,-)$ is exact, or just right exact. Similarly, $B$ is injective $\iff$ $\Mor_\A(-,B)$ exact, or just right exact.
\end{crl}

\section{Exactness and Partial Exactness Criteria for Sequences}

\begin{lmm}[\textcolor{blue}{\index{Exactness! criteria for sequences}{Exactness criteria for sequences}}]\label{ExactCrit}
Let $\A$ be an abelian category. Take any sequence in $\A$
\[
S:~0\ra A\sr{f}{\ral}B\sr{g}{\ral}C\ra 0.
\]
\bit[leftmargin=0.9cm]
\item[(i)] If the sequence ~$\Mor_\A(X,S):~0\ra \Mor_\A(X,A)\sr{f_\ast}{\ral}\Mor_\A(X,B)\sr{g_\ast}{\ral}\Mor_\A(X,C)\ra 0$~ is exact in Ab for all $X\in\Ob\A$, then $S$ is exact in $\A$.
\item[(ii)] If the sequence ~$\Mor_\A(S,X):~0\ra \Mor_\A(C,X')\sr{g^\ast}{\ral}\Mor_\A(B,X')\sr{f^\ast}{\ral}\Mor_\A(A,X')\ra 0$~ is exact in Ab for all $X'\in\Ob\A$, then $S$ is exact in $\A$.
\eit
\end{lmm}
\begin{proof}~
\bit[leftmargin=0.7cm]
\item[(i)] Assume $\Mor_\A(X,S)$ is exact for all $X\in\Ob\A$.
\bit[leftmargin=0.0cm]
\item[] \ul{$f$ is a monomorphism}: This is immediate because $f_\ast=f\circ(-)$ is injective.
\bea\adjustbox{scale=0.9}{\bt
        & X\ar[d,"h"]         &             &    & \\
0\ar[r] & A\ar[r,"f"] & B\ar[r,"g"] & C\ar[r] & 0
\et}~~~~~~f\circ h=f_\ast(h)=0~~\Ra~~h=0.\nn
\eea

\item[] \ul{$\im f\cong\ker g$}: We have ~$\im f_\ast=\ker g_\ast$, i.e.,
\bea\adjustbox{scale=0.9}{\bt
        &          &  X\ar[dl,dashed,"\exists~\wt{h}"']\ar[d,"h"]           &    & \\
0\ar[r] & A\ar[r,"f"] & B\ar[r,"g"] & C\ar[r] & 0
\et}~~~~~~g_\ast(h)=g\circ h=0~~\iff~~h=f\circ\wt{h}=f_\ast(\wt{h}).\nn
\eea
Equivalently, as in the following diagram, ~$g\circ h=0$ $\iff$ $h=f\circ\wt{h}=kc_f\circ(u_f\circ\wt{h})$.
\bea\bt
  &  &    & X\ar[dll,dashed,bend right,"\wt{h}"']\ar[dl,dashed,"u_f\circ\wt{h}"']\ar[d,"h"] &  &  \\
0\ar[r] & A\ar[rr,bend right=40,"f=kc_f\circ u_f"]\ar[r,"u_f","\cong"'] & \im f\ar[r,hook,"kc_f"] & B\ar[r,"g"] & C\ar[r] & 0
\et\nn
\eea

The above relations are simply expressing the universal property for $\ker g$, and so $\im f\cong\ker g$.

\item[] \ul{$g$ is an epimorphism}: If we set $X=C$, then the surjectivity of $g_\ast$ implies $g$ has a right inverse.
\bea\adjustbox{scale=0.9}{\bt
        &             &             &    X=C\ar[dl,dashed,"\exists~\wt{h}"']\ar[d,"h=id_C"] & \\
0\ar[r] & A\ar[r,"f"] & B\ar[r,"g"] & C\ar[r] & 0
\et}~~~~~~~~id_C=h=g\circ\wt{h}=g_\ast(\wt{h}).\nn
\eea
\eit

\item[(ii)] Assume $\Mor_\A(S,X')$ is exact for all $X'\in\Ob\A$.
\bit[leftmargin=0.0cm]
\item[] \ul{$g$ is an epimorphism}: This is immediate because $g^\ast=(-)\circ g$ is injective.
\bea\adjustbox{scale=0.9}{\bt
0\ar[r] & A\ar[r,"f"] & B\ar[r,"g"] & C\ar[r] & 0\\
        &         &             &  X'\ar[from=u,"h"]   & \\
\et}~~~~~~h\circ g=g^\ast(h)=0~~\Ra~~h=0.\nn
\eea

\item[] \ul{$\im f\cong\ker g$}: We have ~$\im g^\ast=\ker f^\ast$.
\bea\adjustbox{scale=0.9}{\bt
0\ar[r] & A\ar[r,"f"] & B\ar[r,"g"] & C\ar[r] & 0\\
        &             &  X'\ar[from=u,"h'"]\ar[from=ur,dashed,"\exists~\wt{h}'"]   &   & \\
\et}~~~~~~f^\ast(h')=h'\circ f=0~~\iff~~h'=\wt{h}'\circ g=g^\ast(\wt{h}').\nn
\eea
Equivalently, as in the following diagram, ~$h'\circ f=0$ $\iff$ $h'=\wt{h}'\circ g=(\wt{h}'\circ u_g)\circ ck_g$.
\bea\adjustbox{scale=0.9}{\bt
0\ar[r] & A\ar[r,"f"] & B\ar[r,two heads,"ck_g"]\ar[rr,bend left=40,"g=u_g\circ ck_g"] & \coim g\ar[r,"u_g","\cong"'] & C\ar[r] & 0\\
     &  & X'\ar[from=ur,dashed,"\wt{h}'\circ u_g"]\ar[from=u,"h'"]\ar[from=urr,dashed,bend left,"\wt{h}'"] &  & &
\et}\nn
\eea
The above relations are simply expressing the universal property for $\coker f$, and so
\bea
\coker f\cong\coim g,~~\Ra~~\im f=\ker\left(B\sr{c_f}{\ra}\coker f\right)\cong \ker\left(B\sr{ck_g}{\ra}\coim g\right)\sr{(s)}{=}\ker g,\nn
\eea
where step (s) holds because $u_g$ is an isomorphism.

\item[] \ul{$f$ is a monomorphism}: If we set $X'=A$, then the surjectivity of $f^\ast$ implies $f$ has a left inverse.
\[\adjustbox{scale=0.9}{\bt
0\ar[r] & A\ar[r,"f"] & B\ar[r,"g"] & C\ar[r] & 0\\
        &  X'=A\ar[from=u,"h=id_A"']\ar[from=ur,dashed,"\exists~\wt{h}"]    &             &     & \\
\et}~~~~~~~~id_A=h=\wt{h}\circ f=f^\ast(\wt{h}).\nn \qedhere
\]
\eit
\eit
\end{proof}

\begin{crl}[\textcolor{blue}{\index{Partial! exactness criteria for sequences}{Partial exactness criteria for sequences}}]\label{ExactCrit2}
Let $\A$ be an abelian category. Take any sequence ~$S:~0\ra A\sr{f}{\ral}B\sr{g}{\ral}C\ra 0$~ in $\A$.
\bit[leftmargin=0.9cm]
\item[(i)] If the sequence ~$\Mor_\A(X,S):~0\ra \Mor_\A(X,A)\sr{f_\ast}{\ral}\Mor_\A(X,B)\sr{g_\ast}{\ral}\Mor_\A(X,C)$~ is exact in Ab for all $X\in\Ob\A$, then ~$0\ra A\sr{f}{\ral}B\sr{g}{\ral}C$~ is exact in $\A$ (i.e., $S$ is left exact).
\item[(ii)] If the sequence ~$\Mor_\A(S,X):~0\ra \Mor_\A(C,X')\sr{g^\ast}{\ral}\Mor_\A(B,X')\sr{f^\ast}{\ral}\Mor_\A(A,X')$~ is exact in Ab for all $X'\in\Ob\A$, then ~$A\sr{f}{\ral}B\sr{g}{\ral}C\ra 0$~ is exact in $\A$ (i.e., $S$ is right exact).
\item[(iii)] If the sequence ~$\Mor_\A(S,X):~\Mor_\A(C,X')\sr{g^\ast}{\ral}\Mor_\A(B,X')\sr{f^\ast}{\ral}\Mor_\A(A,X')$~ is exact in Ab for all $X'\in\Ob\A$, then ~$A\sr{f}{\ral}B\sr{g}{\ral}C$~ is exact in $\A$ (i.e., $S$ is exact at $B$).
\eit
\end{crl}
\begin{proof}
See the proof of Lemma \ref{ExactCrit}, and notice that exactness at $B$ in the sequence $0\ra A\sr{f}{\ral}B\sr{g}{\ral}C\ra 0$ depends neither on $f$ being monic nor on $g$ being epic. This is because in any decomposition of a morphism ~\bt h=kc_h\circ u_h\circ ck_h:U\ar[r,"ck_h"] & \coim h \ar[r,"u_h"] & \im h \ar[r,"kc_h"] & V \et, ~we have
\[
\im h=\im(kc_h\circ u_h\circ ck_h)=\im(kc_h),~~~~\coim h=\coim(kc_h\circ u_h\circ ck_h)=\coim(ck_h). \qedhere
\]
\end{proof}

\begin{crl}[\textcolor{blue}{Right-adjoined $\Ra$ commute with colimit $\varinjlim$, Left-adjoined $\Ra$ commute with limit $\varprojlim$}]
Let $\A,\B$ be abelian categories, $S:\I\ra\A$ a system, $S':\I\ra\B$ a cosystem, and $(F,G)$ = {\footnotesize$\Big(\bt\A\ar[r,shift left,"F"] & \ar[l,shift left,"G"]\B\et\Big)$} an adjoint pair. Then
\bea
\varinjlim F\big(S(i)\big)\cong F\left(\varinjlim S(i)\right)~~~~\txt{and}~~~~\varprojlim G\big(S'(i)\big)\cong G\left(\varprojlim S'(i)\right).\nn
\eea
\end{crl}
\begin{proof}
The result follows from the exactness criteria, since for any objects $X\in\A$, $Y\in\B$, we have
\begin{align}
&\Mor_\B\left(\varinjlim\big(F(S(i))\big),Y\right)\cong \varprojlim \Mor_\B\left(F\big(S(i)\big),Y\right)\cong \varprojlim \Mor_\A\big(S(i),G(Y)\big)\nn\\
&~~~~\cong \Mor_\A\big(\varinjlim S(i),G(Y)\big)\cong \Mor_\B\Big(F\left(\varinjlim S(i)\right),Y\Big).\nn\\
&\Mor_\A\left(X,\varprojlim G(S'(i))\right)\cong \varprojlim \Mor_\A\left(X,G\big(S'(i)\big)\right)\cong \varprojlim \Mor_\B\big(F(X),S'(i)\big)\nn\\
&~~~~\cong \Mor_\B\big(F(X),\varprojlim S'(i)\big)\cong \Mor_\A\Big(X,G\left(\varprojlim S'(i)\right)\Big).\nn \qedhere
\end{align}
\end{proof}

\begin{lmm}[\textcolor{blue}{\index{Partial! exactness of adjoints}{Right-adjoined $\Ra$ right exact, Left-adjoined $\Ra$ left exact})}]\label{PartExactAdj}
Let $\A,\B$ be abelian categories. If $\A\sr{F}{\ral}\B$ and $\A\sr{G}{\lal}\B$ are adjoint functors, then (i) $F$ is right exact (independently of $G$) and (ii) $G$ is left exact (independently of $F$).
\end{lmm}
\begin{proof}
{\flushleft (i)} Consider a SES$:0\ra A\ra A'\ra A''\ra 0$ in $\A$. Then for any object $B\in\Ob\B$, left exactness of $\Mor_\A(-,G(B))$ implies the following exact commutative diagram:
\bea\adjustbox{scale=0.9}{\bt
0\ar[r] & \Mor_\A(A'',G(B))\ar[d,"\cong"]\ar[r] & \Mor_\A(A',G(B))\ar[d,"\cong"]\ar[r] & \Mor_\A(A,G(B))\ar[d,"\cong"] &\\
0\ar[r] & \Mor_\B(F(A''),B)\ar[r] & \Mor_\B(F(A'),B)\ar[r]\ar[r] & \Mor_\B(F(A),B) &
\et}\nn
\eea
It follows by the partial exactness criteria that the sequence $F(A)\ra F(A')\ra F(A'')\ra 0$ is exact.

{\flushleft (ii)} Similarly, consider a SES$:0\ra B\ra B'\ra B''\ra 0$ in $\B$. Then for any object $A\in\Ob\A$, left exactness of $\Mor_\A(F(A),-)$ implies the following exact commutative diagram:
\bea\adjustbox{scale=0.9}{\bt
0\ar[r] & \Mor_\B(F(A),B)\ar[d,"\cong"]\ar[r] & \Mor_\B(F(A),B')\ar[d,"\cong"]\ar[r] & \Mor_\B(F(A),B'')\ar[d,"\cong"] &\\
0\ar[r] & \Mor_\A(A,G(B))\ar[r] & \Mor_\A(A,G(B'))\ar[r] & \Mor_\A(A,G(B'')) &
\et}\nn
\eea
It follows by the partial exactness criteria that the sequence $0\ra G(B)\ra G(B')\ra G(B'')$ is exact.
\end{proof}

\section{Tensor-Mor Adjointness (Tensor-Mor Duality)}
For an abelian category of the form $\A:=R$-mod, using the tensor product for modules from Section \ref{FrMdTePrSec} (page \pageref{FrMdTePrSec}), the following result can be concretely proved (as in Theorems \ref{ModAdjThmI} and \ref{ModAdjThmII}).
\begin{lmm}[\textcolor{blue}{\index{Tensor-Mor adjointness}{Tensor-Mor adjointness}}]
Let $\A$ be an abelian category and $X\in\Ob\A$. The tensor functors $X\otimes_\A -$ and $-\otimes_\A X$ are left adjoints (i.e., they have right adjoints).
\end{lmm}
\begin{proof}[Formal proof]
Fix an object $A\in\Ob\A$. Since $A\otimes_\A X\cong X\otimes_\A A$, it suffices to show $-\otimes_\A X$ is a left adjoint. Recall that directly from the definition of tensor, we have the following natural isomorphisms (i.e., isomorphisms of functors):
\bea
\Mor_\A(-\otimes_\A -,C)\cong Bim_\A(-\times-,C)\sr{\Phi}{\cong} \Mor_\A\big(-,\mor_\A(-,C)\big),~~~~\txt{for any}~~C\in\Ob\A,\nn
\eea
where at step $\Phi$, a bimorphism $\theta:A\times B\ra C$ in $\A$ has a unique \ul{formal} expansion as follows:
{\footnotesize\bea
\left.
  \begin{array}{ll}
   \theta\sr{\Phi}{\longmapsto}\overbrace{\theta\circ(-\times-)}^{\wt{\theta}}: A\sr{\wt{\theta}}{\ral}\overbrace{\theta\circ(A\times-)}^{\theta_B(A)\subset\Mor_\A(B,C)}:B\mathop{\ral}\limits^{A\times-}A\times B\sr{\theta}{\ral}C, \\~\\
    \theta\sr{\Phi}{\longmapsto}\overbrace{\theta\circ(-\times-)}^{\wt{\theta}}:B\sr{\wt{\theta}}{\ral}\overbrace{\theta\circ(-\times B)}^{\theta_A(B)\subset\Mor_\A(A,C)}:A\mathop{\ral}\limits^{-\times B}A\times B\sr{\theta}{\ral}C,
  \end{array}
\right.~~~~
\adjustbox{scale=1.0}{\bt
A\ar[dr,shift left,hook,"q_A"]\ar[from=dr,shift left,two heads,"p_A"] &           & B\ar[dl,shift right,hook,"q_B"']\ar[from=dl,shift right,two heads,"p_B"'] \\
  & A\times B\ar[d,"\theta"] & \\
  &    C\ar[from=uul,bend right,dashed,"{\Mor_\A(A,C)\supset\theta_A(B)}"']\ar[from=uur,bend left,dashed,"{\theta_B(A)\subset \Mor_\A(B,C)}"]      &
\et}\nn
\eea}and so is \ul{formally} a member of $\Mor_\A\big(A,\mor_\A(B,C)\big)$, where $\mor_\A(B,C)$ is the class $\Mor_\A(B,C)$ viewed as a single object in $\A$ (via an appropriate imbedding ~$\A\hookrightarrow Sets$~ if necessary). This shows the functor ~$-\otimes_\A X:\A\ra\A$~ has a \ul{formal} right adjoint ~$\mor_\A(X,-):\A\ra\A$~ in the sense
\[
\Mor_\A(-\otimes_\A X,-)\cong Bim_\A(-\times X,-)\cong \Mor_\A\big(-,\mor_\A(X,-)\big). \qedhere
\]
\end{proof}

\begin{crl}[\textcolor{blue}{Tensor functors (being right-adjoined) commute with the colimit $\varinjlim$}]
Let $\A$ be an abelian category and $S:\I\ra\A$ a system. Then ~$\varinjlim\big(S(i)\otimes_\A X\big)\cong\left(\varinjlim S(i)\right)\otimes_\A X$ ~for any $X\in\Ob\A$.
\end{crl}

\begin{crl}[\textcolor{blue}{\index{Right! exactness of tensor functors}{Tensor functors (being right-adjoined) are right exact}}]
Let $\A$ be an abelian category and $X\in\Ob\A$. The tensor functors $X\otimes_\A -$ and $-\otimes_\A X$ are right exact (as left adjoints by Lemma \ref{PartExactAdj}).
\end{crl}

\begin{crl}[\textcolor{blue}{\index{Flatness criterion}{Flatness criterion}}]
Let $\A$ be an abelian category. Then an object $A\in\Ob\A$ is flat $\iff$ ~$-\otimes_\A A$ (or equivalently, $A\otimes_\A-$) is exact, or just left exact.
\end{crl}

\section{Adjoints and Partial Exactness of the Colimit/Limit Functors}
\begin{lmm}[\textcolor{blue}{\index{Right! adjoint of colimit}{The colimit $\varinjlim$ is right-adjoinded.}, Diagonal functor}]\label{ColimRExact}
\end{lmm}
\begin{proof}
Let $\I,\C$ be categories. Recall that
\bea
\varinjlim:\C^\I\ra\C,~~S\sr{f}{\ral}S'~~\mapsto~~\varinjlim S\sr{\varinjlim f}{\ral}\varinjlim S'.\nn
\eea
We will show that the following (called \index{diagonal functor}{\ul{diagonal functor}}) is a right adjoint of $\varinjlim$. Define
\bea
\Delta:\C\ra\C^\I,~~C\sr{g}{\ral}C'~\mapsto~\Delta(C)\sr{\Delta(g)}{\ral}\Delta(C'),\nn
\eea
where the functor ~$\Delta(C):\I\ra\C$, ~$i\sr{\kappa_{ij}}{\ral}j$ ~$\mapsto$~ $C\sr{id_C}{\ral}C$.~ We need to show
\bea
\Mor_\C\left(\varinjlim(-),-\right)\cong \Mor_{\C^\I}\big(-,\Delta(-)\big).\nn
\eea
It is clear that ~$\varinjlim\circ\Delta=id_\C$. For any system $S\in\C^\I$ and any object $C\in\C$, define a map $\Phi$ as follows:
\bea
\Mor_{\C^\I}\big(S,\Delta(C)\big)\sr{\Phi}{\ral}\Mor_\C\left(\varinjlim S(i),C\right)\cong  \varprojlim \Mor_\C\left(S(i),C\right),\nn
\eea
\bea
\Phi(f)~~:=~~h~~=~~
\adjustbox{scale=0.9}{\bt
S(i)\ar[ddr,bend right,"f_i"']\ar[dr,"q_i"']\ar[rr,"S(\kappa_{ij})"] && S(j)\ar[dl,"q_j"]\ar[ddl,bend left,"f_j"] \\
  &\varinjlim S\ar[d,dashed,"h"]&  \\
 &C&
\et}~~~~\in~~~~\Mor_\C\left(\varinjlim S(i),C\right),\nn
\eea
\bea
~~~~\txt{for}~~~~f:=\{f_i:S(i)\ra C\}~~=~~
\adjustbox{scale=0.9}{\bt
S(i)\ar[d,"f_i"']\ar[rr,"S(\kappa_{ij})"] && S(j)\ar[d,"f_j"] \\
C\ar[rr,"id_C"] && C
\et}~~~~\in~~~~\Mor_{\C^\I}\big(S,\Delta(C)\big).\nn
\eea
By the universal property for the colimit $\varinjlim$, $\Phi$ is bijective.
\end{proof}

\begin{lmm}[\textcolor{blue}{\index{Left! adjoint of limit}{The limit $\varprojlim$ is left-adjoined.}, Diagonal functor}]\label{LimLExact}
\end{lmm}
\begin{proof}
Let $\I,\C$ be categories. Recall that
\bea
\varprojlim:\C^{\I^{op}}\ra\C,~~S\sr{f}{\ral}S'~~\mapsto~~\varprojlim S\sr{\varprojlim f}{\lal}\varprojlim S'.\nn
\eea
We will show that the following (called \index{diagonal functor}{\ul{diagonal functor}}) is a left adjoint of $\varprojlim$. Define
\bea
\Delta:\C\ra\C^{\I^{op}},~~C\sr{g}{\ral}C'~\mapsto~\Delta(C)\sr{\Delta(g)}{\ral}\Delta(C'),\nn
\eea
where the functor ~$\Delta(C):\I\ra\C$, ~$i\sr{\kappa_{ij}}{\ral}j$ ~$\mapsto$~ $C\sr{id_C}{\ral}C$.~ We need to show
\bea
\Mor_\C\left(-,\varprojlim(-)\right)\cong \Mor_{\C^{\I^{op}}}\big(\Delta(-),-\big).\nn
\eea
It is clear that ~$\varprojlim\circ\Delta=id_\C$. For any system $S\in\C^{\I^{op}}$ and any object $C\in\C$, define a map $\Phi$ as follows:
\bea
\Mor_{\C^{\I^{op}}}\big(\Delta(C),S\big)\sr{\Phi}{\ral}\Mor_\C\left(C,\varprojlim S(i)\right)\cong    \varprojlim \Mor_\C\left(C,S(i)\right),\nn
\eea
\bea
\Phi(f)~~:=~~h~~=~~
\adjustbox{scale=0.9}{\bt
S(i)\ar[from=ddr,bend left,"f_i"]\ar[from=dr,"p_i"]\ar[from=rr,"S(\kappa_{ij})"'] && S(j)\ar[from=dl,"p_j"']\ar[from=ddl,bend right,"f_j"'] \\
  &\varprojlim S\ar[from=d,dashed,"h"']&  \\
 &C&
\et}~~~~\in~~~~\Mor_\C\left(C,\varprojlim S(i)\right),\nn
\eea
\bea
~~~~\txt{for}~~~~f:=\{f_i:S(i)\la C\}~~=~~
\adjustbox{scale=0.9}{\bt
S(i)\ar[from=d,"f_i"]\ar[from=rr,"S(\kappa_{ij})"'] && S(j)\ar[from=d,"f_j"'] \\
C\ar[from=rr,"id_C"'] && C
\et}~~~~\in~~~~\Mor_{\C^{\I^{op}}}\big(\Delta(C),S\big).\nn
\eea
By the universal property for the limit $\varprojlim$, $\Phi$ is bijective.
\end{proof}

\begin{crl}[\textcolor{blue}{\index{Partial! exactness of colimits/limits}{Partial exactness of colimits/limits}}]\label{LimPaExFunct}
Let $\I$ be a category and $\A$ an abelian category.
\bit
\item[(i)] The colimit $\varinjlim:\A^\I\ra\A$ is a right exact functor (as a right-adjoined functor, by Lemma \ref{ColimRExact}).
\item[(ii)] The limit $\varprojlim:\A^{\I^{op}}\ra\A$ is a left exact functor (as a left-adjoined functor, by Lemma \ref{LimLExact}).
\eit
\end{crl}

Note that by the above result, the functor $\varinjlim$ (resp. $\varprojlim$) is exact $\iff$ it maps monic (resp. epic) morphisms to monic (resp. epic) morphisms.

\begin{crl}
Let $\A$ be an additive category, $S,S':\I\ra\A$ systems, and $T,T':\I\ra\A$ cosystems.
\bea
(i)~~\varinjlim\coker(S_i\ra S_i')\cong \coker\varinjlim(S_i\ra S_i'),~~~~(ii)~~\varprojlim\ker(T^i\ra T'{}^i)\cong \ker\varprojlim(T^i\ra T'{}^i).\nn
\eea
\end{crl}
\begin{proof}
In additive categories, right (resp. left) exact functors preserve or commute with cokernels (resp. kernels).
\end{proof}

\begin{crl}[\blue{Interchange of colimits/limits}]
Let $\I,\J,\C$ be categories. If $S:\I\times\J\ra\C$ is a bisystem,
\[
\textstyle\varinjlim_i~\varinjlim_j~S(i,j)~=~\varinjlim_j~\varinjlim_i~S(i,j)~~~~\txt{and}~~~~\varprojlim_i~\varprojlim_j~S(i,j)~=~\varprojlim_j~\varprojlim_i~S(i,j).
\]
\end{crl}

\section{Nature of Colimits/Limits in Module and Abelian Categories}
\begin{prp}[\textcolor{blue}{\index{Existence of limits in $R$-mod}{Existence of limits in $R$-mod}: \index{Colimits as cokernels}{Colimits as cokernels}, \index{Limits as kernels}{Limits as kernels}}]\label{ExistLimRmod}
In $R$-mod, consider any system $S:\I\ra R$-mod, written as {\small$\{(S_i,\vphi_{ij})\}:=\big\{S(i),S(\kappa_{ij})~|~i,j\in\I\big\}$}, and any cosystem $T:\I\ra R$-mod, written as {\small$\{(T^i,\psi^{ij})\}:=\big\{T(i),T(\kappa_{ji})~|~i,j\in \I\big\}$}.
\bit[leftmargin=0.9cm]
\item[(i)] The colimit~ $\varinjlim S$ exists, and is given by ~{\footnotesize$\varinjlim S_i\cong\coker\left(\bigoplus\limits_{(i,j)\in \I^2}\dom\vphi_{ij}\sr{\vphi}{\ral}\bigoplus\limits_{i\in \I}S_i\right)$}, ~where  $\vphi$ is the $R$-homomorphism given by
\bea
\vphi|_{\dom\vphi_{ij}}\big(c_i(m)\big):=c_i(m)-c_j\big(\vphi_{ij}(m)\big),~~~~\txt{for all}~~m\in\dom\vphi_{ij}=S_i,~~(i,j)\in\I^2,\nn
\eea
with ~$c_i:S_i\hookrightarrow \bigoplus\limits_{k\in \I}S_k,~m\mapsto c_i(m):=\big(\delta_{ik}m\big)_{k\in \I}$ ~the canonical inclusions.

\item[(ii)] The limit $\varprojlim T$ exists, and is given by ~{\footnotesize$\varprojlim T\cong\ker\left(\prod\limits_{i\in \I}T^i\sr{\psi}{\ral}\prod\limits_{(i,j)\in \I^2}\cod\psi^{ij}\right)$}, ~where  $\psi$ is the $R$-homomorphism given, for all ~$(m^i)_{i\in \I}\in\prod_{i\in \I}T^i$, ~by
\bea
\textstyle\psi\big((m^k)_{k\in \I}\big):=\Big(\big(\psi^{ij}\circ\pi^i-\pi^j\big)\big((m^k)_{k\in\I}\big)\Big)^{\txt{others}=0}_{\psi^{ij}:m^j\in\cod\psi^{ij}}=\big(\psi^{ij}(m^i)-m^j\big)^{\txt{others}=0}_{\psi^{ij}:m^j\in\cod\psi^{ij}}.\nn
\eea
\eit
\end{prp}
\begin{proof}
{\flushleft (i)} Let ~$c_i:S_i\hookrightarrow \bigoplus\limits_{k\in \I}S_k,~m\mapsto c_i(m):=\big(\delta_{ik}m\big)_{k\in \I}$~ be the canonical inclusions. Then a colimit~ $\left\{q_j:S_j\ra\varinjlim S~|~j\in \I\right\}$~ for the system is given by the cokernel of the map $\vphi$ as follows.
\bea
&&\textstyle\varinjlim S:={\bigoplus_{j\in \I}S_j\over W},~~W:=\im\vphi=\left\langle\Big\{c_i(m_i)-c_j\big(\vphi_{ij}(m_i)\big)~\big|~\txt{for all}~~\vphi_{ij},~m_i\in\dom\vphi_{ij}=S_i\Big\}\right\rangle,\nn\\
&&\textstyle q_i:=\pi\circ c_i:S_i\sr{c_i}{\ral}\bigoplus\limits_{j\in \I}S_j\sr{\pi}{\ral}{\bigoplus_{j\in \I}S_j\over W},~~m\mapsto c_i(m)+W.\nn
\eea
Given any $f_i:S_i\ra X$ such that $f_i=f_j\vphi_{ij}$, define $\theta:{\bigoplus_{k\in \I}S_k\over W}\ra X,~(m_j)_{j\in \I}+W\mapsto\sum_jf_j(m_j)$. Then $\theta$ is well defined (\blue{footnote}\footnote{If $(w_i)_{i\in I}\in W$, then $(w_i)_{i\in I}=\sum_{i,j}r_{ij}\big(c_i(m_i)-c_j(\vphi_{ij}(m_i))\big)$ where $r_{ij}\in R$ are $0$ a.e.f., and so $\theta((w_i)_{i\in I})=\sum_{i,j}r_{ij}\big(\theta(c_i(m_i))-\theta(c_j(\vphi_{ij}(m_i)))\big)=\sum_{i,j}r_{ij}\big(f_i(m_i)-f_j(\vphi_{ij}(m_i))\big)=\sum_{i,j}r_{ij}(f_i-f_j\vphi_{ij})(m_i)=0$.}), and $f_i=\theta\circ q_i$.
{\small\bea
\adjustbox{scale=0.9}{\bt
S_i\ar[dddr,bend right,"f_i"']\ar[ddr,"q_i"']\ar[dr,hook,dotted,"c_i"]\ar[rr,"\vphi_{ij}"] && S_j\ar[dl,hook,dotted,"c_j"']\ar[ddl,"q_j"]\ar[dddl,bend left,"f_j"] \\
  & \bigoplus_{k\in \I}S_k\ar[d,two heads,dotted,"\pi"] & \\
  &{\bigoplus_{k\in \I}S_k\over W}\ar[d,dashed,"\theta"]&  \\
 &X&
\et}\nn
\eea}

{\flushleft (ii)} Let {\footnotesize$\pi^i:\prod_{k\in \I}T^k\twoheadrightarrow T^i,~(m^k)_{k\in \I}\mapsto m^i$} be the canonical projections. Then a limit {\footnotesize$\left\{p^i:\varprojlim T\ra T^i~|~i\in \I\right\}$} for the inverse system is given by the kernel of $\psi$ as follows.
\bea
&&\textstyle\varprojlim T:=\ker\psi=\left\{(m^k)_{k\in \I}\in\prod\limits_{k\in \I}T^k:~\psi^{ij}\pi^i\big((m^k)_{k\in \I}\big)=\pi^j\big((m^k)_{k\in \I}\big),~\txt{for all}~\psi^{ij}\right\}\nn\\
&&\textstyle~~~~=\left\{(m^k)_{k\in \I}\in\prod\limits_{k\in \I}T^k:~\psi^{ij}(m^i)=m^j,~\txt{for all}~\psi^{ij},~m^j\in\cod\psi^{ij}\right\}\subset\prod\limits_{k\in \I}T^k,\nn\\
&&\textstyle p^i~=~\pi^i\big|_{\varprojlim T}:~\varprojlim T~\ra~T^i.\nn
\eea
Given any $f^i:X\ra T^i$ such that $f^i=\psi^{ji}f^j$, define ~{\small$\theta:X\ra\varprojlim T,~x\mapsto\left(f^j(x)\right)_{j\in \I}$}. Then $f^i=p^i\circ\theta$.

{\small\bea
\adjustbox{scale=0.9}{\bt
T^i\ar[from=dddr,bend left,"f^i"]\ar[from=ddr,"p^i"]\ar[from=rr,"\psi^{ji}"'] && T^j\ar[from=ddl,"p^j"']\ar[from=dddl,bend right,"f^j"'] \\
  & \prod_{k\in \I}T^k\ar[ul,two heads,dotted,"\pi^i"']\ar[ur,two heads,dotted,"\pi^j"] & \\
  &\varprojlim T\ar[u,hook,dotted]\ar[from=d,dashed,"\theta"']&  \\
 &X&
\et}\nn
\eea}
\end{proof}
It follows from the above result that in $R$-mod, colimits $\varinjlim$ exist \emph{precisely because} (i.e., $\iff$) the corresponding coproducts $\coprod$ (i.e., or direct sums $\bigoplus$) exist, and limits $\varprojlim$ exist precisely because the corresponding products $\prod$ exist. (\blue{footnote}\footnote{The same is true for any abelian category, i.e., in an abelian category, limits (resp. colimits) exist $\iff$ kernels and products (resp. cokernels and coproducts) exist. A more general \ul{partial} result can be found in the literature: In any category, limits exist $\iff$ equalizers and products exist (while the mirror statement ``colimits exist $\iff$ coequalizers and coproducts exist'' might not hold in a general category).})

\begin{crl}[\textcolor{blue}{Construction of limits in an abelian category: Colimits as cokernels, Limits as kernels}]
In an abelian category $\A$, consider any system $S:\I\ra \A$, written as {\small$\{(S_i,\vphi_{ij})\}:=\big\{S(i),S(\kappa_{ij})~|~i,j\in\I\big\}$}, and any cosystem $T:\I\ra \A$, written as {\small$\{(T^i,\psi^{ij})\}:=\big\{T(i),T(\kappa_{ji})~|~i,j\in \I\big\}$}.
\bit[leftmargin=0.9cm]
\item[(i)] If coproducts exists, then the colimit~ $\varinjlim S$ exists, and is given by ~{\footnotesize$\varinjlim S_i\cong\coker\left(\coprod\limits_{(i,j)\in \I^2}\dom\vphi_{ij}\sr{\vphi}{\ral}\coprod\limits_{i\in \I}S_i\right)$}, ~where  $\vphi:=\coprod\limits_{(i,j)\in \I^2}\left(c_i\circ id_{S_i}-c_j\circ\vphi_{ij}\right)$ is a morphism constructed from the associated morphisms
\bea
\textstyle c_i\circ id_{S_i}-c_j\circ\vphi_{ij}:S_i\ra \coprod\limits_{i\in \I}S_i,~~~~\txt{for all}~~(i,j)\in\I^2~~\txt{s.t.}~~\dom\vphi_{ij}=S_i,\nn
\eea
with ~$c_i:S_i\hookrightarrow \coprod\limits_{k\in I}S_k$ ~canonical ``inclusions''.

\item[(ii)] If products exist, then the limit $\varprojlim T$ exists, and is given by ~{\footnotesize$\varprojlim T\cong\ker\left(\prod\limits_{i\in \I}T^i\sr{\psi}{\ral}\prod\limits_{(i,j)\in \I^2}\cod\psi^{ij}\right)$}, ~where  $\psi:=\prod\limits_{(i,j)\in \I^2}\left(id_{T^i}\circ\pi^i-\psi^{ij}\circ\pi^j\right)$ is a morphism constructed from the associated morphisms
\bea
\textstyle id_{T^i}\circ\pi^i-\psi^{ij}\circ\pi^j:T^j\ra \prod\limits_{(i,j)\in \I^2}\cod\psi^{ij},~~~~~~~~\txt{for all}~~(i,j)\in\I^2~~\txt{s.t.}~~\cod\psi^{ij}=T^j,\nn
\eea
with ~$\pi^i:\prod_{k\in I}T^k\twoheadrightarrow T^i$ ~canonical ``projections''.
\eit
\end{crl}

\begin{crl}
For systems in abelian categories, the following hold:
\begin{enumerate}[leftmargin=0.9cm]
\item Right-exact product-preserving functors commute with colimits: If $F$ is a right-exact functor that preserves products, then ~$F\circ\varinjlim\cong\varinjlim\circ F$.
\item Left-exact coproduct-preserving functors commute with limits: If $F$ is a left-exact functor that preserves coproducts, then ~$F\circ\varprojlim\cong\varprojlim\circ F$.
\item Exact product and coproducts-preserving functors commute with colimits and limits: If $F$ is an exact functor that preserves products and coproducts, then ~$F\circ\varinjlim\cong\varinjlim\circ F$ ~and ~$F\circ\varprojlim\cong\varprojlim\circ F$.
\end{enumerate}
(Note: The converses are also trivially true.)
\end{crl}
\begin{proof}
For systems in abelian categories, right-exact functors (resp. left-exact functors) preserve or commute with cokernels (resp. kernels).
\end{proof}

\section{Hom-stability and Exactness of the Limit Functor in $R$-mod}
Here Hom (resp. hom) stands for Mor (resp. mor) because $R$-mod is a category of sets for which \ul{morphisms} are \ul{homomorphisms}. In particular, ``Hom-stability'' means ``Mor-stability'' as used earlier.

The following result has already been established at a general level in Proposition \ref{HomStabLim}. However, the alternative proof here (without assuming Proposition \ref{HomStabLim}) is itself interesting/informative.

\begin{prp}[\textcolor{blue}{\index{Hom-stability of colimits and limits}{Hom-stability of colimits and limits revisited}}]
Let $\I$ be a category. In $R$-mod, consider any system $S:\I\ra R$-mod, written as {\small$\{(S_i,\vphi_{ij})\}:=\big\{S(i),S(\kappa_{ij})~|~i,j\in\I\big\}$}, and any cosystem $T:\I\ra R$-mod, written as {\small$\{(T^i,\psi^{ij})\}:=\big\{T(i),T(\kappa_{ji})~|~i,j\in \I\big\}$}. Then for any $R$-modules $X,Y$,
\bea
(i)~Hom_R\left(\varinjlim S_i,Y\right)~\cong~ \varprojlim~Hom_R(S_i,Y),~~~~(ii)~Hom_R\left(X,\varprojlim T^i\right)~\cong~ \varprojlim~Hom_R(X,T^i).\nn
\eea
\end{prp}
\begin{proof}
{\flushleft (i)} Let ~$S_i\sr{c_i}{\ral}\bigoplus_{k\in {\I}} S_k\sr{\pi_j}{\ral}S_j$~ be the canonical inclusions and projections satisfying
\bea
\textstyle \pi_i\circ c_j=\delta_{ij} id_{S_j},~~~~\sum_{k\in {\I}} c_k\circ \pi_k=id_{\bigoplus_{k\in {\I}} S_k}.\nn
\eea
Observe that the maps
{\small
\bea
&&\textstyle\Phi:Hom_R\Big(\bigoplus\limits_{i\in \I}S_i,Y\Big)\ra\prod\limits_{i\in {\I}}Hom_R(S_i,Y),~f\mapsto (f\circ c_i)_{i\in \I},\nn\\
&&\textstyle\Psi:\prod\limits_{i\in {\I}}Hom_R(S_i,Y)\ra Hom_R\Big(\bigoplus\limits_{i\in \I}S_i,Y\Big),~(f_i)_{i\in \I}\mapsto \sum_{i\in\I} f_i\circ \pi_i\nn
\eea
} satisfy $\Phi\circ\Psi=id_{\txt{dom}\Psi}$ and $\Psi\circ\Phi=id_{\txt{dom}\Phi}$, since
{\footnotesize
\bea
&&\textstyle\Phi\circ\Psi\big((f_i)_{i\in {\I}}\big)=\Phi\left(\sum f_j\circ \pi_j\right)=\big((\sum_j f_j\circ \pi_j)\circ c_i\big)_{i\in {\I}}=\left(f_i\circ \pi_i\circ c_i\right)_{i\in {\I}}=\left(f_i\circ id_{M_i}\right)_{i\in {\I}}=\left(f_i\right)_{i\in {\I}},\nn\\
&&\textstyle\Psi\circ\Phi(f)=\Psi\big((f\circ c_i)_{i\in {\I}}\big)=\sum_i(f\circ c_i)\circ \pi_i=f\circ\sum_i c_i\circ \pi_i=f\circ id_{\bigoplus_i S_i}=f.\nn
\eea
}Recall (Prop. \ref{ExistLimRmod}) that we have an exact sequence $\bigoplus\limits_{(i,j)\in \I^2}\dom\vphi_{ij}\sr{\vphi}{\ral}\bigoplus\limits_{i\in \I}S_i\sr{c}{\ral}\varinjlim S_i:=\coker\vphi\ra 0$. Thus, by the left exactness of $Hom_R(-,Y)$, we get the following commutative diagram with \ul{exact rows}:
\bea
\adjustbox{scale=0.9}{\bt
0\ar[r] & Hom_R\left(\varinjlim S_i,Y\right)\ar[d,dashed,"\exists!"',"g"]\ar[r,"c^\ast"] & Hom_R\Big(\bigoplus\limits_{i\in \I}S_i,Y\Big)\ar[d,"\cong"]\ar[r,"\vphi^\ast"] & Hom_R\Big(\bigoplus\limits_{(i,j)\in \I^2}\dom\vphi_{ij},Y\Big)\ar[d,"\cong"]\\
0\ar[r] & \varprojlim Hom_R(S_i,Y)\ar[r,"c'{}^\ast"] & \prod\limits_{i\in \I} Hom_R\left(S_i,Y\right)\ar[r,"\vphi'{}^\ast"] & \prod\limits_{(i,j)\in \I^2} Hom_R\left(\dom\vphi_{ij},Y\right)
\et}\nn
\eea
It follows (by the commutativity of the diagram) that the induced map of kernels $g$ is an isomorphism.

{\flushleft (ii)} Let $T^i\sr{c_i}{\ral}\prod_{k\in \I}T^k\sr{\pi_j}{\ral}T^j$ be the canonical inclusions and projections satisfying
\bea
\textstyle \pi_i\circ c_j=\delta_{ij} id_{T^j},~~~~\sum_{k\in {\I}} c_k\circ\pi_k=id_{\prod_{k\in \I}T^k},\nn
\eea
where the \emph{sum} has the obvious \emph{interpretation} as a convenient ``\emph{place holder}'' for the \emph{components} of an arbitrary element of $\prod_{k\in \I}T^k$. As in part (i), the maps
{\small
\begin{align}
&\textstyle\Phi:Hom_R\Big(X,\prod\limits_{i\in \I}T^i\Big)\ra\mathop{\prod}\limits_{i\in \I}Hom_R(X,T^i),~f\mapsto (\pi_i\circ f)_{i\in \I},\nn\\
&\textstyle\Psi:\prod\limits_{i\in \I}Hom_R(X,T^i)\ra Hom_R\Big(X,\prod\limits_{i\in \I}T^i\Big),~(f_i)_{i\in \I}\mapsto f:=\sum_i c_i\circ f_i,~~\txt{i.e.,}~~f(x):=\big(f_i(x)\big)_{i\in \I}~\txt{for}~x\in X\nn
\end{align}
} satisfy $\Phi\circ\Psi=id_{\txt{dom}\Psi}$ and $\Psi\circ\Phi=id_{\txt{dom}\Phi}$, since
{\footnotesize
\bea
&&\textstyle\Phi\circ\Psi\big((f_i)_{i\in \I}\big)=\Phi\left(\sum_j c_j\circ f_j\right)=\big(\pi_i\circ(\sum c_j\circ f_j)\big)_{i\in \I}=\left(\pi_i\circ c_i\circ f_i\right)_{i\in \I}=\left(id_{S_i}\circ f_i\right)_{i\in \I}=\left(f_i\right)_{i\in \I},\nn\\
&&\textstyle\Psi\circ\Phi(f)=\Psi\big((\pi_i\circ f)_{i\in \I}\big)=\sum_i c_i\circ(\pi_i\circ f)=\left(\sum_i c_i\circ \pi_i\right)\circ f=id_{\prod_i T^i}\circ f=f.\nn
\eea
}Recall (Prop. \ref{ExistLimRmod}) that we have an exact sequence ~$0\ra\ker\psi:=\varprojlim T^i \sr{k}{\ral} \prod\limits_{i\in \I}T^i\sr{\psi}{\ral}\prod\limits_{(i,j)\in \I^2}\cod\psi^{ij}$. Thus, by the left exactness of $Hom_R(X,-)$, we get the following commutative diagram with \ul{exact rows}:
\bea
\adjustbox{scale=0.9}{\bt
0\ar[r] & Hom_R\left(X,\varprojlim T^i\right)\ar[d,dashed,"\exists!"',"h"]\ar[r,"k_\ast"] & Hom_R\Big(X,\prod\limits_{i\in \I}T^i\Big)\ar[d,"\cong"]\ar[r,"\psi_\ast"] & Hom_R\Big(X,\prod\limits_{(i,j)\in \I^2}\cod\psi^{ij}\Big)\ar[d,"\cong"]\\
0\ar[r] & \varprojlim Hom_R\left(X,T^i\right)\ar[r,"k'_\ast"] & \prod\limits_{i\in \I} Hom_R\left(X,T^i\right)\ar[r,"\psi'_\ast"] & \prod\limits_{(i,j)\in \I^2}Hom_R\left(X,\cod\psi^{ij}\right)
\et}\nn
\eea
It follows (by the commutativity of the diagram) that the induced map of kernels $h$ is an isomorphism.\qedhere
\end{proof}

\begin{thm}[\textcolor{blue}{
\index{Exactness! of the directed colimit in $R$-mod}{\ul{Exactness} of the \ul{directed colimit} in \ul{$R$-mod}}}]\label{DirLimExact}
Let $\I$ be a directed set (viewed as a category). Then the directed colimit ~$\varinjlim:(R\txt{-mod})^\I\ra R\txt{-mod}$ ~is an exact functor.
\end{thm}
\begin{proof} Let ~$S:\I\ra R$-mod ~be a directed system. Recall (Prop. \ref{ExistLimRmod}) that, with the canonical inclusions
\bea
\textstyle c_i:S_i\ra \bigoplus\limits_{k\in \I}S_k,~m\mapsto c_i(m):=\big(\delta_{ik}m\big)_{k\in \I},\nn
\eea
the colimit~ $\left\{q_j:S_j\ra\varinjlim S~|~j\in \I\right\}$~ for the system is given by the following.
\bea
&&\textstyle\varinjlim S:={\bigoplus_{i\in \I}S_i\over W},~~~~W:=\left\langle\Big\{c_i(m)-c_j\big(\vphi_{ij}(m)\big)~\big|~m\in\dom\vphi_{ij}=S_i,~\txt{for all}~~i\leq j,~(i,j)\in\I^2\Big\}\right\rangle,\nn\\
&&\textstyle q_i=\pi\circ c_i:S_i\sr{c_i}{\ral}\bigoplus\limits_{j\in \I}S_j\sr{\pi}{\ral}{\bigoplus_{j\in \I}S_j\over W},~~m\mapsto c_i(m)+W.\nn
\eea
\bit[leftmargin=0.7cm]
\item[(i)] Observe that for any $x\in \varinjlim S$, we have $x=\sum_{k=1}^nc_{i_k}\left(m_{i_k}\right)+W$ for some $m_{i_k}\in S_{i_k}$. Choose an upper bound $j\geq i_1,...,i_n$ in $\I$. Then
\bea
&&\textstyle x=\sum\limits_{k=1}^nc_{i_k}\left(m_{i_k}\right)+W=\sum\limits_{k=1}^n\Big[c_{i_k}\left(m_{i_k}\right)-c_j\big(\vphi_{i_kj}(m_{i_k})\big)+c_j\big(\vphi_{i_kj}(m_{i_k})\big) \Big]+W\nn\\
&&\textstyle~~~~=c_j\big(m_j(x)\big)+W~~\in~~\ol{S}_j:={S_j+W\over W},~~~~\txt{where}~~~~m_j(x):=\sum\limits_{k=1}^n\vphi_{i_kj}(m_{i_k})~\in~S_j.\nn
\eea
That is, every $x\in \varinjlim S$ can be written as ~$x=q_j(m_j)=c_j(m_j)+W$ ~for some $j\in \I$.
\item[(ii)] Let $x\in \varinjlim S$. Then by (i) above, {\small$x=c_t(m_t)+W=q_t(m_t)$} for some $t\in\I$. Thus, if {\small$x=0\in\varinjlim S$}, then
{\small\begin{align}
&\textstyle c_t(m_t)=\sum\limits_{k=1}^n\Big[c_{i_k}(m_{i_k})-c_{j_k}\big(\vphi_{i_kj_k}(m_{i_k})\big)\Big],~~\txt{where wlog}~~j_k>i_k~~\txt{and}~~i_1,...,i_n,~j_1,...,j_n~~\txt{are distinct}.\nn
\end{align}}
We have the following three exhaustive cases for the value of $t$:
\begin{enumerate}
\item \ul{$i_k\neq t\neq j_k$ for all $1\leqslant k\leqslant n$}: Then $c_t(m_t)=0$, i.e., $m_t=0$, and so $\vphi_{tj}(m_t)=0$ for any $j\geq t$.
\item \ul{$t=i_k$ for some $1\leqslant k\leqslant n$}: Then $c_{j_k}\big(\vphi_{tj_k}(m_t)\big)=0$, i.e., $\vphi_{tj_k}(m_t)=0$, where $j_k>t$.
\item \ul{$t=j_k$ for some $1\leqslant k\leqslant n$}: Then $c_{i_k}(m_{i_k})=0$ (i.e., $m_{i_k}=0$) and $c_t(m_t)=-c_t\big(\vphi_{i_kt}(m_{i_k})\big)=-c_t\big(\vphi_{i_kt}(0)\big)=0$ (i.e., $m_t=0$), and so $\vphi_{tj}(m_t)=0$ for any $j\geq t$.
\end{enumerate}
Therefore, if $c_t(m_t)\in W$ (i.e., $x=c_t(m_t)+W=q_t(m_t)=0$) then $\vphi_{tj}(m_t)=0$ for some $j\geq t$. Conversely, if $\vphi_{tj}(m_t)$ for some $j\geq t$, then $c_t(m_t)=c_t(m_t)-0=c_t(m_t)-\vphi_{tj}(m_t)\in W$. Hence,
\bea
\textstyle x=c_t(m_t)+W=q_t(m_t)=0~~\iff~~m_t\in\bigcup_{j\geq t}\ker\big(S_t\sr{\vphi_{tj}}{\ral}S_j\big).\nn
\eea
Equivalently, we have ~$\ker\big(S_i\sr{q_i}{\ral}\varinjlim S\big)=\bigcup_{j\geq i}\ker\big(S_i\sr{\vphi_{ij}}{\ral}S_j\big)$~ for each $i\in \I$.
\eit

{\flushleft Consider} a SES of directed systems in $(R\txt{-mod})^\I$ written simply (with transition morphisms hidden) as
\[
0\ra (S_i)_{i\in\I}\sr{(f_i)}{\ral}(S_i')_{i\in\I}\sr{(f_i')}{\ral}(S_i'')_{i\in\I}\ra 0~~~~\txt{or}~~~~\big(0\ra S_i\sr{f_i}{\ral}S_i'\sr{f_i'}{\ral}S_i''\ra 0\big)_{i\in\I},
\]
and the resulting sequence of limits:
\bea
\adjustbox{scale=0.9}{\bt
 0\ar[r]& S_i\ar[dddd,near start,bend right=45,"\vphi_{ij}"']\ar[rr,"f_i"]\ar[ddr,"q_i"]&& S'_i\ar[dddd,near start,bend right=45,"\vphi'_{ij}"']\ar[rr,"f'_i"]\ar[ddr,"q'_i"] && S''_i\ar[dddd,near start,bend right=45,"\vphi''_{ij}"']\ar[ddr,"q''_i"] \ar[r] & 0 \\
 & &&  &&   &  \\
 & 0\ar[r,crossing over] &\varinjlim S\ar[rr,dashed,crossing over,"{\exists!~f}"]&& \varinjlim S'\ar[rr,dashed, crossing over,"{\exists!~f'}"] && \varinjlim S'' \ar[r,crossing over] & 0 \\
 & &&  && & \\
0\ar[r]& S_j\ar[rr,"f_j"]\ar[uur,"q_j"']&& S'_j\ar[rr,"f'_j"]\ar[uur,"q'_j"'] && S''_j\ar[uur,"q''_j"'] \ar[r] & 0
\et}\nn
\eea
Since the colimit $\varinjlim$ is right exact, to show the sequence of limits is exact, we only need to show $f$ is injective. Let $x\in\varinjlim S$ such that $f(x)=0$. Then for some $t\in\I$,
\begin{align}
&\textstyle 0=f(x)\sr{\txt{(i)}}{=}f(q_t(m_t))=q_t'f_t(m_t)~~\Ra~~f_t(m_t)\in\ker\big(S_t'\sr{q_t'}{\ral}\varinjlim S'\big)\sr{\txt{(ii)}}{=}\bigcup_{j\geq t}\ker\big(S_t'\sr{\vphi'_{tj}}{\ral}S_j'\big),\nn\\
&\textstyle~~~~\Ra~~0=\vphi'_{tj}\big(f_t(m_t)\big)=f_j\big(\vphi_{tj}(m_t)\big)~~~~\txt{for some}~~j\geq t,\nn\\
&~~~~\Ra~~\vphi_{tj}(m_t)=0,~~~~\txt{since $f_j$ is injective},\nn\\
&~~~~\Ra~~0\sr{\txt{(ii)}}{=}q_t(m_t)\sr{\txt{(i)}}{=}x.\nn \qedhere
\end{align}
\end{proof}

\begin{thm}[\textcolor{blue}{\index{Exactness! of the direct sum in $R$-mod}{Exactness of the direct sum in $R$-mod}}]\label{DirSumExact}
Let $\I$ be a set (viewed as a trivial category). Then the direct sum on $\I$-indexings of $R$-modules is an exact functor given by
\[
\textstyle\bigoplus:(R\txt{-mod})^\I\ra R\txt{-mod},~\big\{M_i\sr{f_i}{\ral}N_i\big\}_{i\in\I}~\mapsto~\bigoplus M_i\sr{\bigoplus f_i}{\ral}\bigoplus N_i.
\]
\end{thm}
\begin{proof}
Consider a SES of trivial systems in $(R\txt{-mod})^\I$ written as
\[
0\ra (S_i)_{i\in\I}\sr{(f_i)}{\ral}(S_i')_{i\in\I}\sr{(f_i')}{\ral}(S_i'')_{i\in\I}\ra 0~~~~\txt{or}~~~~\big(0\ra S_i\sr{f_i}{\ral}S_i'\sr{f_i'}{\ral}S_i''\ra 0\big)_{i\in\I},
\]
and the resulting sequence of direct sums:
\bea
\adjustbox{scale=0.9}{\bt
 0\ar[r]& S_i\ar[rr,"f_i"]\ar[dr,"q_i"]&& S'_i\ar[rr,"f'_i"]\ar[dr,"q'_i"] && S''_i\ar[dr,"q''_i"] \ar[r] & 0 \\
 & 0\ar[r,crossing over] &\bigoplus S_i\ar[rr,dashed,crossing over,"{\exists!~f}"]&& \bigoplus S_i'\ar[rr,dashed, crossing over,"{\exists!~f'}"] && \bigoplus S_i'' \ar[r,crossing over] & 0 \\
0\ar[r]& S_j\ar[rr,"f_j"]\ar[ur,"q_j"']&& S'_j\ar[rr,"f'_j"]\ar[ur,"q'_j"'] && S''_j\ar[ur,"q''_j"'] \ar[r] & 0
\et}\nn
\eea
Since the colimit $\varinjlim$ (hence $\bigoplus$) is right exact, to show the sequence of direct sums is exact, we only need to show $f$ is injective. We can check that the definition $f:\bigoplus S_i\ra\bigoplus S'_i$, $(s_i)\mapsto(f_i(s_i))$ satisfies the universal property: Indeed, it is clear that $f q_i=q_i' f_i$ for all $i$, and for any $g:\bigoplus S_i\ra\bigoplus S'_i$ satisfying $g q_i=q_i' f_i$ for all $i$, we have
\bea
(f-g)q_i=0~~~\txt{for all}~~i,~~\Ra~~(f-g)|\im q_i=0~~~\txt{for all}~~i,~~\Ra~~f=g.\nn
\eea
Therefore,
\bea
\textstyle\ker f=\{(s_i)_{i\in\I}~|~(f_i(s_i))_{i\in \I}=0\}=\{(s_i)_{i\in\I}~|~s_i\in\ker f_i\}=\bigoplus_{i\in\I}\ker f_i=\bigoplus_{i\in\I}0=0.\nn\qedhere
\eea
\end{proof}

\begin{crl}[\textcolor{blue}{From the proof of Theorem \ref{DirSumExact}}]
A sequence {\small $0\ra S\ra S'\ra S''\ra 0$} of trivial systems $S,S',S'':\I\ra R\txt{-mod}$ is exact iff the induced sequence {\small $0\ra \bigoplus S(i)\ra \bigoplus S'(i)\ra \bigoplus S''(i)\ra 0$} is exact.
\end{crl}

\begin{crl}[\textcolor{blue}{\index{Exactness! and additivity of Tot on $R$-mod}{Exactness and additivity of Tot on $R$-mod}}]
The total functor
\[
Tot:(R\txt{-mod})_0^{\Integer\times\Integer}\ra(R\txt{-mod})_0^\Integer,~~\big(C_{\ast\ast}\sr{f_{\ast\ast}}{\ral}C'_{\ast\ast}\big)\mapsto~ \big(Tot(C_{\ast\ast})\sr{Tot(f_{\ast\ast})}{\ral}Tot(C'_{\ast\ast})\big)
\]
is an \ul{exact additive} functor (its additivity coming from that of the direct sum $\bigoplus$).

Recall (from Definition \ref{HomCohomDf1}) that the chain complex $Tot(C_{\ast\ast})_\ast=\big(Tot(C_{\ast\ast})_\ast,\del_\ast\big)$ is given formally by
\bea
\textstyle Tot(C_{\ast\ast})_n:=\bigoplus\limits_{i+j=n}C_{ij},~~~~\del_n=\del^{Tot(C_{\ast\ast})}_n:=\sum\limits_{i+j=n}\left(\del^h_{ij}+(-1)^i\del^v_{ij}\right):Tot(C_{\ast\ast})_n\ra Tot(C_{\ast\ast})_{n-1}.\nn
\eea
\end{crl}
\begin{proof}
Let $S:0\ra A_{\ast\ast}\sr{f_{\ast\ast}}{\ral} B_{\ast\ast}\sr{g_{\ast\ast}}{\ral} C_{\ast\ast}\ra 0$ be a SES of bicomplexes of $R$-modules, i.e., $S_{ij}:0\ra A_{ij}\sr{f_{ij}}{\ral} B_{ij}\sr{g_{ij}}{\ral} C_{ij}\ra 0$ is an SES for all $i,j\in\Integer$. Then for each $n\in\Integer$, we get an exact sequence
\bea\bt
S_n:0\ar[r] & \bigoplus_{i+j=n}A_{ij}\ar[r,"\bigoplus f_{ij}"] & \bigoplus_{i+j=n}B_{ij}\ar[r,"\bigoplus g_{ij}"] & \bigoplus_{i+j=n}C_{ij}\ar[r] & 0.\nn
\et\eea
Hence, the sequence of complexes ~{\small $Tot(S):0\ra Tot(A_{\ast\ast})\sr{Tot(f_{\ast\ast})}{\ral} Tot(B_{\ast\ast})\sr{Tot(g_{\ast\ast})}{\ral} Tot(C_{\ast\ast})\ra 0$} ~is exact.
\end{proof}

\section{Free Modules, Tensor Products of Modules and Algebras}\label{FrMdTePrSec}
\begin{dfn}[\textcolor{blue}{\index{Free! object in a category}{Recall: Free object in a category}}]
Let $\C$ be a category and $X$ a set. A free object over $X$ (or generated by $X$) in $\C$ is an object $\langle X\rangle\in\Ob\C$ such that given any imbedding (i.e., mapwise-injective functor) $J:\C\hookrightarrow Sets$, the following hold.
\bit
\item[(i)] There exists an injection $i:X\hookrightarrow J\big(\langle X\rangle\big)$.
\item[(ii)] Universal property: For every map $\theta:X\ra J(C)$, $C\in\Ob\C$, there exists a unique morphism $\theta':J(\langle X\rangle)\ra J(C)$ in Sets (or equivalently, $\theta'':=J^{-1}(\theta'):\langle X\rangle\ra C$ in $\C$) such that $\theta=\theta'\circ i$.
\bea
\adjustbox{scale=0.9}{\bt
X\ar[d,"\theta"]\ar[rr,hook,"i"] && J\big(\langle X\rangle\big)\ar[dll,dashed,"\exists!~\theta'"]\\
J(C)
\et}~~~~~~~~\theta=\theta'|_X~:=~J(\theta')\circ i.\nn
\eea
\eit
\end{dfn}

(Note: If $\C\subset\txt{Sets}$ already, then we simply drop the imbedding $J$, i.e., $J$ becomes the identity functor.)

\begin{dfn}[\textcolor{blue}{Free objects in $R$-mod: \index{Free! module}{Free module}}]
Let $R$ be a ring and $X$ a set. The $R$-module \ul{with basis $X$} (or \ul{freely spanned by $X$}, or \ul{freely generated by $X$}), or equivalently, the \ul{free $R$-module} on $X$ (or spanned by $X$, or generated by $X$), is an $R$-module $R^{\langle X\rangle}$ together with an inclusion $q:X\hookrightarrow R^{\langle X\rangle}$ such that for any map $f:X\ra M$ from $X$ to an $R$-module $M$, there exists a unique $R$-homomorphism $\theta:R^{\langle X\rangle}\ra M$ such that ~$f=\theta|_X:=\theta\circ q$.
\bea
\adjustbox{scale=0.9}{\bt
X\ar[d,"f"]\ar[rr,hook,"q"] && R^{\langle X\rangle}\ar[dll,dashed,"\theta"]\\
M
\et}~~~~f=\theta|_X~(=\theta\circ q).\nn
\eea
{\flushleft\ul{Concrete realization:}} Let $X$ be a set, $R$ a ring, $M$ an $R$-module, and $\bigoplus_{x\in X}R$ the direct sum of an $X$-indexing of $R$-modules $(R_x)_{x\in X}$ with $R_x={}_RR$ for all $x\in X$. Define $q:X\hookrightarrow\bigoplus_{x\in X}R$, $x\mapsto(\delta_{xx'})_{x'\in X}$. Given any map $f:X\ra M$, define $\theta:\bigoplus_{x\in X}R\ra M$, $(r_x)_{x\in X}\mapsto\sum_{x\in X}r_xf(x)$. Then $f=\theta\circ q$, since
\bea
\textstyle(\theta\circ q)(x)=\theta\big(q(x)\big)=\theta\big((\delta_{xx'})_{x'\in X'}\big)=\sum_{x'\in X}\delta_{xx'}f(x)=f(x),~~~~\txt{for all}~~x\in X.\nn
\eea
That is, ~$R^{\langle X\rangle}\cong \bigoplus_{x\in X}R$.
\end{dfn}

Observe that a map $f:X\ra M$ gives a unique collection of $R$-homomorphisms $\{f_x:R\ra M,~1_R\mapsto f(x)\}_{x\in X}$, and such a collection in turn gives a unique map $f:X\ra M$, $x\mapsto f_x(1_R)$. Consequently, we can also directly define a free module using the colimit format as follows.

\begin{dfn}[\textcolor{blue}{Alternative description: Free $R$-module as a colimit $\varinjlim$}]
Let $R$ be a ring and $X$ a set. The $R$-module \ul{with basis $X$} (or \ul{freely spanned by $X$}, or \ul{freely generated by $X$}), or equivalently, the \ul{free $R$-module} on $X$ (or spanned by $X$, or generated by $X$), denoted by $R^{\langle X\rangle}$, is the colimit $\varinjlim S$ of a trivial system of the following form (with $X$ as a trivial category, i.e., $Hom_X(x,x')=\emptyset$ for $x\neq x'$):
\[
S:X\ra R\txt{-mod},~x\mapsto S(x)\cong R.
\]
That is, $R^{\langle X\rangle }$ is given by canonical inclusions $\left\{q_x:S(x)\hookrightarrow\bigoplus S\right\}_{x\in X}$ such that for any collection of $R$-homomorphisms $\{f_x:S(x)\ra M\}_{x\in X}$, there exists a unique $R$-homomorphism $\theta: R^{\langle X\rangle}\ra M$ satisfying $\theta\circ q_x=f_x$ for all $x\in X$.
\bea
\adjustbox{scale=0.8}{\bt
S(x)\cong R\ar[ddr,bend right,"f_x"']\ar[dr,hook,"q_x"] && S(x')\cong R\ar[dl,hook,"q_{x'}"']\ar[ddl,bend left,"f_{x'}"] \\
  &R^{\langle X\rangle}:=\bigoplus S\ar[d,dashed,"\theta"]& \\
 &M&
\et}~~~~f_x=\theta\circ q_x,~~\txt{for all}~~x\in X.\nn
\eea

{\flushleft\ul{Concrete realization  for the alternative description}:} For each $x\in X$, let $S(x):=R$ and define a map $q_x:R\ra \bigoplus_{x\in X}R$, $r\mapsto(r\delta_{xx'})_{x'\in X}$. Given any $R$-homomorphisms $\{f_x:R\ra M\}_{x\in X}$, define a map
\[
\textstyle\theta:\bigoplus_{x\in X}R\ra M,~(r_x)_{x\in X}\mapsto \sum_{x\in X}r_xf_x(1_R).
\]
Then for each $x\in X$, we have ~$\theta\circ q_x=f_x$, ~since
\bea
\textstyle(\theta\circ q_x)(r)=\theta\big(q_x(r)\big)=\theta\big((\delta_{xx'}r)_{x'\in X}\big)=\sum\limits_{x'\in X}\delta_{xx'}rf_{x'}(1_R)=rf_x(1_R)=f_x(r),~~\txt{for all}~~r\in R.\nn
\eea
Hence, with the alternative description, we also get ~$R^{\langle X\rangle}\cong \bigoplus_{x\in X}R$.
\bea
\adjustbox{scale=0.8}{\bt
S(x)\ar[d,"\cong"]\ar[dr,hook] && S(x')\ar[d,"\cong"]\ar[dl,hook,] \\
R\ar[ddr,bend right,"f_x"']\ar[dr,hook,"q_x"'] & R^{\langle X\rangle}:=\bigoplus_{x\in X}S(x)\ar[d,"\cong"]
& R\ar[dl,hook,"q_{x'}"]\ar[ddl,bend left,"f_{x'}"] \\
  &\bigoplus_{x\in X}R\ar[d,dashed,"\theta"]& \\
 &M&
\et}~~~~f_x=\theta\circ q_x,~~\txt{for all}~~x\in X.\nn
\eea

\end{dfn}

\begin{rmks}[\blue{Formal conventions}]
(i) The free $R$-module $R^{\langle X\rangle}\cong\bigoplus_{x\in X}R$ is often written as
\bea
\textstyle\Span_RX:=\left\{\txt{formal finite sums}~\sum_{i=1}^nr_ix_i:~x_i\in X,~n\geq 0\right\},\nn
\eea
where the \ul{formal generator} $x_i\in X$ in the sum actually represents the \ul{unit function} $e_{x_i}:X\ra R$, $x\mapsto\delta_{xx_i}$.

(ii) The free $R$-module $R^{\langle X\rangle}\cong\bigoplus_{x\in X}R$ is also \ul{formally} written as
\bea
\textstyle\sum_{x\in X}Rx:=\left\{\sum_{x\in X}r_xx~|~r_x=0~~\txt{a.e.f.}~~x\in X\right\}=\left\{\sum_{x\in X}r_xx~|~(r_x)_{x\in X}\in\bigoplus_{x\in X}R\right\},\nn
\eea
where a.e.f. means ``for all except finitely many''.
\end{rmks}

\begin{thm}[\textcolor{blue}{Existence of \index{Tensor! products of modules}{tensor products of modules}}]\label{TnsPrdRmod}
Let $R$ be a commutative ring and $R'$ an $R$-algebra. Then there exists a tensor bifunctor
{\small\[
\otimes_{R'}=-\otimes_{R'}-:(\txt{mod-}R')\times(R'\txt{-mod})\ra R\txt{-mod},~~\Big((A,B)\sr{(f,g)}{\ral}(C,D)\Big)~\mapsto \Big(A\otimes_{R'}B\sr{f\otimes g}{\ral}C\otimes_{R'}D\Big).
\]}
\end{thm}
\begin{proof}
Let $A_{R'}\in\txt{mod-}R'\subset R\txt{-mod}$, $_{R'}B\in R'\txt{-mod}\subset R\txt{-mod}$, and $R^{\langle A\times B\rangle}$ the free $R$-module over the \ul{set} $A\times B$. Consider the $R$-submodule $W\subset R^{\langle A\times B\rangle}$ given by
\bea
\label{TensorW}W=W_{R'}=W_{R'}(A,B):=\left\langle
\left\{\substack{(ar',b)-(a,r'b),~~~~~~~~~~~~\\
r(a_1,b)+(a_2,b)-(ra_1+a_2,b),\\
r(a,b_1)+(a,b_2)-(a,rb_1+b_2),~}~\Big|~\substack{r\in R,~r'\in R'\\ a,a_1,a_2\in A,\\ b,b_1,b_2\in B,~}\right\}\right\rangle.
\eea
Let {\small $\ol{A\times B}:={A\times B+W\over W}=\{a\otimes_{R'} b:=(a,b)+W~|~a\in A,b\in B\}$}. Then {\small ${R^{\langle A\times B\rangle}\over W}=\Span_R\big(\ol{A\times B}\big)\neq R^{\langle\ol{A\times B}\rangle}$}, where the inequality means that in general the spanning set $\ol{A\times B}$ is not linearly independent over $R$ (hence not an $R$-basis). Given any $R$-homomorphism $f: A\times B\ra M$, satisfying (for all $r\in R,~~r'\in R',~~a,a_1,a_1\in A,~~b,b_1,b_2\in B$)
\[
f(ar',b)=f(a,r'b),~~f(ra_1+a_2,b)=rf(a_1,b)+f(a_2,b),~~f(a,rb_1+b_2)=rf(a,b_1)+f(a,b_2),
\]
(\blue{footnote}\footnote{Such a map is said to be \index{Balanced map}{\ul{$R'$-balanced}} and \index{Bilinear map}{\ul{$R$-bilinear}}.}) and its unique $R$-linear extension $f_1:R^{\langle A\times B\rangle}\ra M$, consider the $R$-homomorphism $\wt{f}:{R^{\langle A\times B\rangle}\over W}\ra M$ given by ~$\wt{f}(a\otimes_{R'} b):=f_1|_{A\times B}(a,b)=f\big((a,b)\big)$. Then we get $f=f_1\circ q=\wt{f}\circ\tau$ as in the commutative diagram:
\bea
\adjustbox{scale=0.9}{\bt
A\times B\ar[dd,"f"']\ar[rr,bend left,"\tau:=\pi\circ q"]\ar[r,hook,"q"] & R^{\langle A\times B\rangle}\ar[ddl,dashed,"\exists!~f_1"']\ar[r,two heads,"\pi"] & {R^{\langle A\times B\rangle}\over W}\ar[ddll,dashed,"\exists!~\td{f}"]\\
 &&\\
 M
\et}~~~~f=f_1\circ q=\wt{f}\circ\tau,\nn
\eea
which shows the map $\tau:A\times B\ra {R^{\langle A\times B\rangle}\over W}$ has a \ul{universal property} (hence \ul{uniqueness} up to isomorphism).

It follows therefore that in $R$-mod, there exists a tensor product ~$\otimes:(R\txt{-mod})^2\ra R\txt{-mod}$~ such that
{\footnotesize\bea
\otimes_{R'}:=\otimes|_{(\txt{mod-}R')\times(R'\txt{-mod})}:(\txt{mod-}R')\times(R'\txt{-mod})\ra R\txt{-mod},~~\Big((A,B)\sr{(f,g)}{\ral}(C,D)\Big)~\mapsto \Big(A\otimes_{R'}B\sr{f\otimes g}{\ral}C\otimes_{R'}D\Big),\nn
\eea}where ~$A\otimes_{R'} B\cong {R^{\langle A\times B\rangle}\over W}$~ for all $A\in\txt{mod-}R'$ and $B\in R'\txt{-mod}$.

In particular, if $R'=R$ is any commutative ring, then we have a tensor product
{\small\bea
\otimes_R:(R\txt{-mod})^2\ra R\txt{-mod},~~\Big((A,B)\sr{(f,g)}{\ral}(C,D)\Big)~\mapsto \Big(A\otimes_RB\sr{f\otimes g}{\ral}C\otimes_RD\Big),\nn
\eea}where ~$A\otimes_R B\cong {R^{\langle A\times B\rangle}\over W}$~ for all $R$-modules $A,B$.
\end{proof}

\begin{rmk}[\textcolor{blue}{Alternative description: The product $\otimes_{R'}$ as a colimit $\varinjlim$}]
Consider the same setup as in the proof of Theorem \ref{TnsPrdRmod} above, i.e., let $R$ be a commutative ring, $R'$ an $R$-algebra, $A\in \txt{mod-}R'$, and $B\in R'\txt{-mod}$. If $(a,b)\in A\times B$, let ~$R(a,b):=R^{\langle(a,b)\rangle}:=R^{\langle\{(a,b)\}\rangle}\cong R$~ be the free $R$-module generated by the singleton set $\{(a,b)\}$. Then for any $(a,b),(a',b')\in A\times B$, we have the $R$-isomorphism $\vphi:R(a,b)\ra R(a',b'),~(a,b)\mapsto(a',b')$.

Consider the system $S:A\times B\ra R$-mod with (i) objects $S(a,b):=R(a,b)\cong R$ and (ii) transition maps ~{\footnotesize$\vphi_{(a,b)(a',b')}:=
\left\{
  \begin{array}{ll}
    \vphi, &\txt{if}~~(a,b)+W=(a',b')+W, \\
    \emptyset~~\txt{(i.e., no map)}, & \txt{otherwise}.
  \end{array}
\right\}
$}. Then we have the following commutative diagram:
\bc
\adjustbox{scale=0.9}{\bt
\overbrace{S(a,b)}^{R(a,b)}\ar[dddrr,bend right,"f_{(a,b)}"']\ar[ddrr,bend right,"\tau_{(a,b)}"]\ar[drr,hook,"q_{(a,b)}"]\ar[rrrr,bend left=10,"\vphi_{(a,b)(a',b')}"] &&&&
\overbrace{S(a',b')}^{R(a',b')}\ar[dddll,bend left,"f_{(a',b')}"]\ar[ddll,bend left,"\tau_{(a',b')}"']\ar[dll,hook,"q_{(a',b')}"']  \\
 && R^{\langle A\times B\rangle}\ar[d,two heads,"\pi"] &&\\
 &&\ub{{R^{\langle A\times B\rangle}/W}}_{\cong\varinjlim S}\ar[d,dashed,"f"] &&\\
 && M &&
\et}\ec
Given any $R$-homomorphisms $\big\{f_{(a,b)}:S(a,b)\ra M\big\}_{(a,b)\in A\times B}$
satisfying $f_{(a',b')}\circ\vphi_{(a,b)(a',b')}=f_{(a,b)}$, where we note that
{\small\bea
\label{AltTensEq}f_{(a',b')}\circ\vphi_{(a,b)(a',b')}=
\left\{
  \begin{array}{ll}
    f_{(a',b')}\circ\vphi, &\txt{if}~~(a,b)+W=(a',b')+W \\
   \emptyset, & \txt{otherwise}
  \end{array}
\right\},
\eea} define $f:{R^{\langle A\times B\rangle}\over W}\ra M$ to be the $R$-homomorphism given by
\bea
f\big((a,b)+W\big):=f_{(a,b)}(a,b),~~~~\txt{for all}~~(a,b)\in A\times B.\nn
\eea
Then $f$ is well defined because from (\ref{AltTensEq}), $f_{(a',b')}\circ\vphi_{(a,b)(a',b')}=f_{(a,b)}$ implies $f_{(a,b)}(a,b)=(f_{(a',b')}\circ\vphi)(a,b)=f_{(a',b')}(a',b')$ whenever $(a,b)+W=(a',b')+W$. Since we have $f_{(a,b)}=f\circ\tau_{(a,b)}$ for all $(a,b)\in A\times B$, it follows that $\varinjlim S\cong {R^{\langle A\times B\rangle}\over W}$. Hence, because ${R^{\langle A\times B\rangle}\over W}\cong A\otimes_{R'} B$ (by the proof of Theorem \ref{TnsPrdRmod}), we get
\bea
\textstyle \varinjlim S\cong A\otimes_{R'} B.\nn
\eea
\end{rmk}

\begin{rmk}[\blue{\index{Tensor! product of rings/algebras}{Tensor product of rings/algebras}}]
Let $R$ be a commutative ring and $R$-alg the category of $R$-algebras. Then based directly on the construction in the proof of Theorem \ref{TnsPrdRmod}, we can define a tensor bifunctor
\[
\otimes_R:(R\txt{-alg})^2\ra R\txt{-alg},~~(A,B)\sr{(f,g)}{\ral}(C,D)~~\mapsto~~A\otimes_RB\sr{f\otimes g}{\ral}C\otimes_RD
\]
such that for any $R$-algebras $A,B,C$, the following hold:
\begin{enumerate}[leftmargin=0.9cm]
\item There exists an $R$-algebra homomorphism $\tau=\tau_{AB}:A\times B\ra A\otimes_RB$ satisfying the following \ul{$R$-bilinearity} property: For all $r\in R$, $a,a_1,a_2\in A$, $b,b_1,b_2\in B$,
\[
\tau(ar,b)=\tau(a,rb)=r\tau(a,b),~~\tau(a_1+a_2,b)=\tau(a_1,b)+\tau(a_2,b),~~\tau(a,b_1+b_2)=\tau(a,b_1)+\tau(a,b_2).\nn
\]

\item For any $R$-algebra homomorphism $f:A\times B\ra C$ satisfying the \ul{$R$-bilinearity} property above, i.e., for all $r\in R$, $a,a_1,a_2\in A$, $b,b_1,b_2\in B$,
\[
f(ar,b)=f(a,rb)=rf(a,b),~~f(a_1+a_2,b)=f(a_1,b)+f(a_2,b),~~f(a,b_1+b_2)=f(a,b_1)+f(a,b_2),\nn
\]
there exists a unique $R$-algebra homomorphism $\wt{f}:A\otimes_RB\ra C$ such that $\wt{f}\circ\tau=f$.
\end{enumerate}

{\flushleft\ul{Realization}}: The existence of the above tensor product follows steps similar to those in the proof of Theorem \ref{TnsPrdRmod}. Let $Te({}_RR(A\times B))$ be the tensor $R$-algebra of the ring $A\times B$ as an $R$-module, i.e.,
\begin{align}
&\textstyle Te({}_R(A\times B)):=R\langle A\times B\rangle:=\sum_{i\in\Natural}R(A\times B)^i\nn\\
&\textstyle~~~~=\{\sum_{i\in\Natural}r_ic_i:r_i\in R,~c_i\in (A\times B)^i,~r_i=0~\txt{a.e.f.}\},~~~~(A\times B)^0:=R,~~(A\times B)^1=A\times B,\nn
\end{align}
and consider the ideal $J\lhd Te({}_R(A\times B))$ generated as follows:
\bea
J=J(A,B):=\left\langle
\left\{\substack{(ar,b)-(a,rb),~~~~~~~~~~~~\\
r(a_1,b)+(a_2,b)-(ra_1+a_2,b),\\
r(a,b_1)+(a,b_2)-(a,rb_1+b_2),\\
(a,b)(a',b')-(aa',bb'),~~~~~}~\Big|~
\substack{r\in R,~~~~~~~\\
a,a_1,a_2,a'\in A,\\
b,b_1,b_2,b'\in B,~}\right\}
\right\rangle~\lhd ~Te({}_R(A\times B)).\nn
\eea
Given an $R$-algebra homomorphism $f:A\times B\ra C$ satisfying the \ul{$R$-bilinearity} property, it is clear that we get its unique extension to an $R$-algebra homomorphism $f_1:Te({}_R(A\times B))\ra C$ (i.e., such that $f_1|_{A\times B}=f$). We have ${Te({}_R(A\times B))\over J}=\Span_R\ol{A\times B}$, where $\ol{A\times B}:=\{a\otimes b:=(a,b)+J~|~a\in A,b\in B\}$. By construction, in addition to the quotient-induced \ul{$\otimes$-bilinearty} given (for all $r\in R$, $a,a'\in A$, $b,b'\in B$) by $ra\otimes b=a\otimes rb=r(a\otimes b)$, $(a+a')\otimes b=a\otimes b+a'\otimes b$, and $a\otimes(b+b')=a\otimes b+a\otimes b'$, we also have the quotient-induced non-free \ul{$\otimes$-multiplication} rule
\bea
\textstyle (a\otimes b)(a'\otimes b')=aa'\otimes bb',~~~~\txt{for all}~~a\otimes b,a'\otimes b'\in \ol{A\times B}\subset {Te({}_R(A\times B))\over J},\nn
\eea
that is, {\footnotesize $\left(\sum_ir_i~a_i\otimes b_i\right)\left(\sum_ir'_i~a'_i\otimes b'_i\right)=\sum_{i,j}r_ir'_j~a_ia'_j\otimes b_ib'_j$ for all $\sum_ir_i~a_i\otimes b_i$, $\sum_ir'_i~a'_i\otimes b'_i$ in ${Te({}_R(A\times B))\over J}$}. Define {\small $\wt{f}:{Te({}_R(A\times B))\over J}\ra C,~a\otimes b\mapsto f_1(a,b)=f(a,b)$}. Then $\wt{f}\circ\tau=f$, where $\tau$ is given in the diagram:
\bea
\adjustbox{scale=0.9}{\bt
A\times B\ar[dd,"f"']\ar[rr,bend left,"\tau:=\pi\circ q"]\ar[r,hook,"q"] & Te({}_R(A\times B))\ar[ddl,dashed,"\exists!~f_1"']\ar[r,two heads,"\pi"] & {Te({}_R(A\times B))\over J}\ar[ddll,dashed,"\exists!~\td{f}"]\\
 &&\\
 C
\et}~~~~f=f_1\circ q=\wt{f}\circ\tau.\nn
\eea
It follows that we have an $R$-algebra isomorphism ~$A\otimes_R B\cong {Te({}_R(A\times B))\over J}$.
\end{rmk}

\begin{thm}[\blue{\index{Tensor-Hom adjointness I}{Tensor-Hom adjointness I}}]\label{ModAdjThmI}
(\blue{footnote}\footnote{From its proof, it is clear that this result holds for the tensor product of rings/algebras as well.}). Let $R,S$ be rings, $A_R,{}_RB_S,C_S$ (bi)-modules, and $\otimes:=\otimes_R:(\txt{mod-}R)\times(R\txt{-mod})\ra Ab$. Then we have a natural isomorphism of abelian groups given by
\bea
&&\Phi:Hom_S(A_R\otimes{}_RB_S,C_S)\ra Hom_R\big(A_R,Hom_S(_RB_S,C_S)\big),\nn\\
&&~~~~~~~~\left(A\otimes B\sr{f}{\ral}C\right)\longmapsto \left(A\sr{\Phi(f)}{\ral}Hom(B,C)\right),~~~~~~~\Phi(f)(a)(b):=f(a\otimes b).\nn
\eea
\end{thm}
\begin{proof}
$\bullet$ The map $\Phi$ is well defined (with respect to its codomain) because $\Phi(f)$ is indeed and $R$-homomorphism, i.e., $\Phi(f)(ra+a')=r\Phi(f)(a)+\Phi(f)(a')$, as follows: (\blue{footnote}\footnote{This key step works essentially due to the $R$-bilinearity of $\otimes$ (or equivalently, of bimorphisms). Therefore, the isomorphism can also be directly defined for bimorphisms as
\bea
&&\Phi:Bim_S(A_R\times B_S,C_S)\ra Hom_R\big(A_R,Hom_S(_RB_S,C_S)\big),\nn\\
&&~~~~~~~~\Big(A\times B\sr{f}{\ral}C\Big)\longmapsto \Big(A\sr{\Phi(f)}{\ral}Hom(B,C)\Big),~~~~~~~\Phi(f)(a)(b):=f(a,b).\nn
\eea
})
\begin{align}
&\Phi(f)(ra+a')(b)=f((ra+a')\otimes b)=f(ra\otimes b+a'\otimes b)=rf(a\otimes b)+f(a'\otimes b)\nn\\
&~~~~=r\Phi(f)(a)(b)+\Phi(f)(a')(b)=\big(r\Phi(f)(a)+\Phi(f)(a')\big)(b).\nn
\end{align}
Also, by the properties of pointwise function evaluation, we have $\Phi(f+f')=\Phi(f)+\Phi(f')$, since
\begin{align}
&\Phi(f+f')(a)(b)=(f+f')(a\otimes b)=f(a\otimes b)+f'(a\otimes b)=\Phi(f)(a)(b)+\Phi(f')(a)(b)\nn\\
&~~~~=(\Phi(f)(a)+\Phi(f')(a))(b)=(\Phi(f)+\Phi(f'))(b)(a),\nn
\end{align}
and so $\Phi$ is a homomorphism of abelian groups.
{\flushleft $\bullet$} Next, consider the map
\bea
&&\textstyle\Psi:Hom_R\big(A_R,Hom_S(_RB_S,C_S)\big)\ra Hom_S(A_R\otimes{}_RB_S,C_S),\nn\\
&&\textstyle~~~~~~~~\Big(A\sr{g}{\ral}Hom(B,C)\Big)\longmapsto \Big(A\otimes B\sr{\Psi(g)}{\ral}C\Big),~~~~~~~\Psi(g)(\sum_ir_i~a_i\otimes b_i):=\sum_ir_ig(a_i)(b_i),\nn\\
&&\textstyle~~~~\txt{in particular,}~~~~\Psi(g)(a\otimes b):=g(a)(b),~~~~\txt{which suffices by $R$-linearity},\nn
\eea
where it is clear that $\Psi(g)$ is an $R$-homomorphism, and $\Psi$ is a homomorphism of abelian groups since
\begin{align}
&\Psi(g+g')(a\otimes b)=(g+g')(a)(b)=(g(a)+g'(a))(b)=g(a)(b)+g'(a)(b)=\Psi(g)(a\otimes b)+\Psi(g')(a\otimes b)\nn\\
&~~~~=\big(\Psi(g)+\Psi(g')\big)(a\otimes b),~~\Ra~~\Psi(g+g')=\Psi(g)+\Psi(g').\nn
\end{align}

{\flushleft $\bullet$} Then the map~ $\Phi\circ\Psi:Hom_R\big(A_R,Hom_S(_RB_S,C_S)\big)\ra Hom_R\big(A_R,Hom_S(_RB_S,C_S)\big)$,~ that is,
\bea
Hom_R\big(A_R,Hom_S(_RB_S,C_S)\big)\sr{\Psi}{\ral} Hom_S(A_R\otimes{}_RB_S,C_S)\sr{\Phi}{\ral} Hom_R\big(A_R,Hom_S(_RB_S,C_S)\big)\nn
\eea
is given by
\bea
&&\big[(\Phi\circ\Psi)(g)\big](a)(b)=\Phi\big[\Psi(g)\big](a)(b)=\Psi(g)(a\otimes b)=g(a)(b),\nn\\
&&~~\Ra~~\Phi\circ\Psi=\txt{id}_{Hom_R\big(A_R,Hom_S(_RB_S,C_S)\big)},\nn
\eea
and the map~ $\Psi\circ\Phi:Hom_S(A_R\otimes{}_RB_S,C_S)\ra Hom_S(A_R\otimes{}_RB_S,C_S)$,~ that is,
\bea
Hom_S(A_R\otimes{}_RB_S,C_S)\sr{\Phi}{\ral}Hom_R\big(A_R,Hom_S(_RB_S,C_S)\big)\sr{\Psi}{\ral} Hom_S(A_R\otimes{}_RB_S,C_S)\nn
\eea
is given by
\bea
&&\big[(\Psi\circ\Phi)(f)\big](a\otimes b)=\Psi\big[\Phi(f)\big](a\otimes b)=\Phi(f)(a)(b)=f(a\otimes b),\nn\\
&&~~\Ra~~\Psi\circ\Phi=\txt{id}_{Hom_S(A_R\otimes{}_RB_S,C_S)}.\nn
\eea
{\flushleft$ \bullet$} We can verify the naturality of $\Phi$ by directly checking that the following diagram is commutative:
\bea
\bt[column sep=tiny] A\ar[dd,"u"] & B\ar[dd,"v"] & C\ar[dd,"w"] \\ & & \\ A' & B' & C'\et~~\mapsto~~
\bt
 Hom_S(A_R\otimes{}_RB_S,C_S)\ar[from=dd,shift left=5,"{(u\otimes id)^\ast}"]\ar[from=dd,shift right=7,"{(id\otimes v)^\ast}"description]\ar[dd,shift left=15,"w_\ast"]\ar[rr,"\Phi"] && Hom_R\big(A_R,Hom_S(_RB_S,C_S)\big)\ar[from=dd,shift left=10,"{u^\ast}"]\ar[from=dd,shift right=12,"{v^\ast{}_\ast}"description]\ar[dd,shift left=20,"w_{\ast\ast}"] \\
  && \\
 Hom_S(A'_R\otimes{}_RB'_S,C'_S)\ar[rr,"\Phi'"] && Hom_R\big(A'_R,Hom_S(_RB'_S,C'_S)\big)
 \et\nn
\eea
or equivalently,
\bea
\bt[column sep=tiny] A\ar[dd,"u"] & B\ar[dd,"v"] & C\ar[from=dd,"w^{op}"'] \\ & & \\ A' & B' & C'\et~~\mapsto~~
\bt
 Hom_S(A_R\otimes{}_RB_S,C_S)\ar[from=dd,shift left=5,"{(u\otimes id)^\ast}"]\ar[from=dd,shift right=7,"{(id\otimes v)^\ast}"description]\ar[from=dd,shift right=15,"w^{op}_\ast"']\ar[rr,"\Phi"] && Hom_R\big(A_R,Hom_S(_RB_S,C_S)\big)\ar[from=dd,shift left=10,"{u^\ast}"]\ar[from=dd,shift right=12,"{v^\ast{}_\ast}"description]\ar[from=dd,shift right=20,"w^{op}_{\ast\ast}"'] \\
  && \\
 Hom_S(A'_R\otimes{}_RB'_S,C'_S)\ar[rr,"\Phi'"] && Hom_R\big(A'_R,Hom_S(_RB'_S,C'_S)\big)
 \et\nn
\eea
or equivalently,
\bea
\adjustbox{scale=0.9}{\bt
 \Big(~A\otimes B\sr{u\otimes v}{\ral}A'\otimes B'\sr{f'}{\ral}C'\sr{w^{op}}{\ral}C~\Big)\ar[from=ddd,mapsto,shift left=5,"{(u\otimes id)^\ast}"]\ar[from=ddd,mapsto,shift right=7,"{(id\otimes v)^\ast}"description]\ar[from=ddd,mapsto,shift right=15,"w^{op}_\ast"']\ar[rr,mapsto,"\Phi"] && \Big(~A\sr{\Phi\big(w^{op}\circ f'\circ(u\otimes v)\big)}{\ral}Hom_S(B,C)~\Big)\ar[d,equal] \\
  &  &
 \Big(~A\sr{u}{\ral}A'\sr{\Phi'(f')}{\ral}Hom_S(B',C')\sr{w^{op}_\ast\circ v^\ast}{\ral}Hom_S(B,C)~\Big)
 \ar[from=dd,mapsto,shift left=10,"{u^\ast}"]\ar[from=dd,mapsto,shift right=12,"{v^\ast{}_\ast}"description]\ar[from=dd,mapsto,shift right=20,"w^{op}_{\ast\ast}"'] \\
  && \\
 A'\otimes B'\sr{f'}{\ral}C'\ar[rr,mapsto,"\Phi'"] && A'\sr{\Phi'(f')}{\ral}Hom_S(B',C')
 \et}\nn
\eea
\begin{align}
& \Phi\big(w^{op}\circ f'\circ(u\otimes v)\big)(a)(b)=\big(w^{op}\circ f'\circ(u\otimes v)\big)(a\otimes b)=\big(w^{op}\circ f'\big)\big(u(a)\otimes v(b)\big)\nn\\
&~~~~=w^{op}\big(f'(u(a)\otimes v(b))\big)=w^{op}\big(\Phi'(f')(u(a))(v(b))\big)=w^{op}\big([v^\ast\circ\Phi'(f')](u(a))(b)\big)\nn\\
&~~~~=[w^{op}_\ast\circ v^\ast\circ\Phi'(f')](u(a))(b)=[w^{op}_\ast\circ v^\ast\circ u^\ast(\Phi'(f'))](a)(b)\nn\\
&~~~~=[w^{op}_\ast\circ v^\ast\circ\Phi'(f')\circ u](a)(b).\nn \qedhere
\end{align}

\end{proof}

\begin{thm}[\blue{\index{Tensor-Hom adjointness II}{Tensor-Hom adjointness II}}]\label{ModAdjThmII}
Let $R,S$ be rings, $_RA,{}_SB_R,{}_SC$ (bi)-modules, and $\otimes:=\otimes_R:(\txt{mod-}R)\times(R\txt{-mod})\ra Ab$. Then we have a natural isomorphism of abelian groups given by
\bea
&&\Phi:Hom_S({}_SB_R\otimes{}_RA,{}_SC)\ra Hom_R(A_R,Hom_S({}_SB_R,{}_SC)),\nn\\
&&~~~~f\mapsto \Phi(f):a\mapsto\Phi(f)(a):b\mapsto \Phi(f)(a)(b):=f(b\otimes a).\nn
\eea
\end{thm}
\begin{proof}
Consider the new map
\bea
&&\Psi:Hom_R(A_R,Hom_S({}_SB_R,{}_SC))\ra Hom_S({}_SB_R\otimes{}_RA,{}_SC),\nn\\
&&~~~~g\mapsto\Psi(g):(b\otimes a)\mapsto \Psi(g)(b\otimes a):=g(a)(b),\nn
\eea
where the verification of well-definedness for the maps $\Phi$ and $\Psi$ follows the same steps as in the proof of the preceding theorem. Then the map~ $\Phi\circ\Psi:Hom_R(A_R,Hom_S({}_SB_R,{}_SC))\ra Hom_R(A_R,Hom_S({}_SB_R,{}_SC))$~ is given by
\bea
&&\big[(\Phi\circ\Psi)(g)\big](a)(b)=\Phi\big[\Psi(g)\big](a)(b)=\Psi(g)(b\otimes a)=g(a)(b),\nn\\
&&~~\Ra~~\Phi\circ\Psi=\txt{id}_{Hom_R(A_R,Hom_S({}_SB_R,{}_SC))},\nn
\eea
and the map~ $\Psi\circ\Phi:Hom_S({}_SB_R\otimes{}_RA,{}_SC)\ra Hom_S({}_SB_R\otimes{}_RA,{}_SC)$~ is given by
\bea
&&\big[(\Psi\circ\Phi)(f)\big](b\otimes a)=\Psi\big[\Phi(f)\big](a\otimes b)=\Phi(f)(a)(b)=f(b\otimes a),\nn\\
&&~~\Ra~~\Psi\circ\Phi=\txt{id}_{Hom_S({}_SB_R\otimes{}_RA,{}_SC)}.\nn
\eea
Finally, we can again verify the naturality of $\Phi$ as done in the proof of the preceding theorem.
\end{proof}

%% file: parts/AlgebraCat/AlgebraCatS7.tex
\chapter{Imbedding of Categories}\label{AlgebraCatS7}

\begin{dfn}[\textcolor{blue}{Recall: Imbedding of categories, Imbedded (sub)category}]
A functor $F:\C\ra\D$ is called an imbedding (making $F(\C)\subset\D$ an imbedded (sub)category) if it is injective as a map. That is, $F$ is (i) injective on morphisms (i.e., faithful) and also (ii) injective on objects.
\end{dfn}

\section{Pushout and Pullback Theorems}
\begin{dfn}[\textcolor{blue}{\index{Pushout-limit (Fiber coproduct)}{Pushout-limit (Fiber coproduct)}, \index{Pullback-limit (Fiber product)}{Pullback-limit (Fiber product)}}]
Let $\C$ be a category. Let $A\sr{f}{\lal}U\sr{g}{\ral}B$ and $A\sr{h}{\ral}V\sr{k}{\lal}B$ be diagrams in $\C$. Then the colimit $\varinjlim\left(A\sr{f}{\lal}U\sr{g}{\ral}B\right)$ is called a \ul{pushout-limit} (or just a \ul{pushout}), while the limit $\varprojlim\left(A\sr{h}{\ral}V\sr{k}{\lal}B\right)$ is called a \ul{pullback-limit} (or just a \ul{pullback}).

The pushout of $A\sr{f}{\lal}U\sr{g}{\ral}B$ is also called the \ul{fiber coproduct}, \ul{fiber sum}, or \ul{cocartesian square} of $A$ and $B$ (or of $f:U\ra A$ and $g:U\ra B$), written $A\oplus_UB$. The pullback of $A\sr{h}{\ral}V\sr{k}{\lal}B$ is also called the \ul{fiber product} or \ul{cartesian square} of $A$ and $B$ (or of $h:A\ra V$ and $k:B\ra V$), written $A\times_VB$.
\[
\adjustbox{scale=0.7}{\bt U\ar[d,"f"]\ar[r,"g"] & B\ar[d,dashed,"k"] \\ A\ar[r,dashed,"h"] & A\oplus_UB\et}~~~~~~~~
\adjustbox{scale=0.7}{\bt A\times_VB\ar[d,dashed,"f"]\ar[r,dashed,"g"] & B\ar[d,"k"] \\ A\ar[r,"h"] & V \et}
\]

Concerning the latter terminology, observe that in the category Sets, the above pullback is given by
\[
\textstyle A\times_VB\cong\{(a,b)\in A\times B:h(a)=k(b)\}=\bigcup_{v\in V}h^{-1}(v)\times k^{-1}(v)\subset A\times B,
\]
with maps $f:=p_A:A\times_VB\twoheadrightarrow V,~(a,b)\mapsto a$ and $g:=p_B:A\times_VB\twoheadrightarrow V,~(a,b)\mapsto b$. Similarly, in Sets, the above pushout is given by
\[
\textstyle A\oplus_UB\cong{A\sqcup B\over\sim}:={(\{0\}\times A)\cup(\{1\}\times B)\over \sim},~~~~(0,f(u))\sim (1,g(u))~~\txt{for all}~~u\in U,\nn
\]
with maps $h:A\hookrightarrow A\oplus_UB,~a\mapsto[(0,a)]_\sim$ and $k:B\hookrightarrow A\oplus_UB,~b\mapsto[(1,b)]_\sim$.

Pushouts and pullbacks can also be explicitly constructed for groups, rings, modules, and algebras. Our main concern however will be with additive and abelian categories (hence with modules).

\end{dfn}
\begin{rmk}[\textcolor{blue}{\index{Pushout-Pullback symmetries}{Pushout-Pullback symmetries}}]\label{PushPullRmk}
Let $\A$ be an additive category and $A,B\in\Ob\A$. Given diagrams $A\sr{f}{\lal}U\sr{g}{\ral}B$ and $A\sr{h}{\ral}V\sr{k}{\lal}B$, the universal property for the direct sum $A\oplus B$ (and the earlier proof of the fact that $A\oplus B\cong A\times B$ in an additive category) gives induced morphisms $\left[\substack{f\\g}\right]$ and $[h,k]$ as in the following diagram (with only triangles commuting) in which $q_A,q_B$ are the canonical injects and $p_A,p_B$ the canonical projections satisfying
\bea
p_Aq_A=id_A,~~~~p_Bq_B=id_B,~~~~q_Ap_A+q_Bp_B=id_{A\oplus B}.\nn
\eea
\bea
\label{PullPushInd}\adjustbox{scale=0.7}{\bt
U\ar[ddrr,dashed,"{\left[\substack{f\\g}\right]}"description]\ar[dddd,"f"']\ar[rrrr,"g"] && && B\ar[from=ddll,bend left,near start,"p_B"]\ar[ddll,bend left,"q_B"']\ar[dddd,"k"] \\
   && && \\
   && A\oplus B\ar[ddrr,dashed,"{[h,k]}"description] && \\
   && && \\
A\ar[from=uurr,bend right,"p_A"']\ar[uurr,bend right,"q_A"]\ar[rrrr,"h"] && && V
\et}\eea
The induced morphisms satisfy ~$p_A\left[\substack{f\\g}\right]=f$, $p_B\left[\substack{f\\g}\right]=g$, $[h,k]q_A=h$, $[h,k]q_B=k$,~ and so
{\small\bea \label{PullPushEq1}hf+kg=[h,k]q_Ap_A\left[\substack{f\\g}\right]+[h,k]q_Bp_B\left[\substack{f\\g}\right]=[h,k](q_Ap_A+q_Bp_B)\left[\substack{f\\g}\right]=[h,k]id_{A\oplus B}\left[\substack{f\\g}\right]=[h,k]\left[\substack{f\\g}\right].
\eea}
It follows by (\ref{PullPushEq1}) that the diagram
{\scriptsize\bt U\ar[d,"f"]\ar[r,"g"] & B\ar[d,"k"] \\ A\ar[r,"h"] & V \et} commutes if and only if the following sequence is a complex:

\bea
0\ral\bt U\ar[r,"{\left[\substack{f\\g}\right]}"] & A\oplus B\ar[r,"{\left[h,-k\right]}"] & V\et\ral 0~~~~\Big(\txt{equiv.}~~0\ral\bt U\ar[r,"{\left[\substack{f\\-g}\right]}"] & A\oplus B\ar[r,"{\left[h,k\right]}"] & V\et\ral 0\Big)\nn
\eea
Also, whenever cokernels and kernels exist (e.g., when $\A$ is abelian), we observe the following:
\begin{enumerate}[leftmargin=0.7cm]
\item The pushout-limit $\varinjlim\left(A\sr{f}{\lal}U\sr{g}{\ral}B\right)=\coker\left[\substack{f\\-g}\right]$, as shown in the diagram below:
\bc\adjustbox{scale=0.8}{\bt
 U\ar[ddrr,"{\left[\substack{f\\-g}\right]}"description]\ar[dddd,"f"']\ar[rrrr,"g"] && && B\ar[from=ddll,bend left,near start,"p_B"]\ar[ddll,bend left,dotted,"q_B"']\ar[dddd,"\pi q_B"]\ar[dddddr,bend left=10,"\beta"] & \\
  && && & \\
  && A\oplus B\ar[ddrr,two heads,"\pi"] && & \\
  && && & \\
 A\ar[from=uurr,bend right,"p_A"']\ar[uurr,bend right,dotted,"q_A"]\ar[rrrr,"\pi q_A"]\ar[drrrrr,bend right=10,"\al"] && && \coker{\left[\substack{f\\-g}\right]}\ar[dr,dashed,"\gamma"] & \\
 && && & X\\
\et}\ec
where the upper-left outer square commutes because
\bea
\pi q_Af=\pi q_Ap_A\left[\substack{f\\-g}\right]=\pi\big((q_Ap_A+q_Bp_B)-q_Bp_B\big)\left[\substack{f\\-g}\right]=-\pi q_Bp_B\left[\substack{f\\-g}\right]=\pi q_Bg.\nn
\eea
Also, if $\al f=\beta g$, then $\al p_A\left[\substack{f\\-g}\right]=\al f=\beta g=-\beta p_B\left[\substack{f\\-g}\right]$, i.e., $(\al p_A+\beta p_B)\left[\substack{f\\-g}\right]=0$, and so there exists a unique $\gamma:\coker{\left[\substack{f\\-g}\right]}\ra X$ such that $\al p_A+\beta p_B=\gamma\pi$. Therefore,
\[
\al=\al p_Aq_A=(\al p_A+\beta p_B)q_A=\gamma\pi q_A~~~~\txt{and}~~~~\beta=\beta p_Bq_B=(\al p_A+\beta p_B)q_B=\gamma\pi q_B.
\]
It follows by comparing with (\ref{PullPushInd}) that the diagram
{\scriptsize \bt U\ar[d,"f"]\ar[r,"g"] & B\ar[d,"k"] \\ A\ar[r,"h"] & V \et} is a pushout-limit $\iff$ the sequence ~$0\ral\bt U\ar[r,"{\left[\substack{f\\g}\right]}"] & A\oplus B\ar[r,"{\left[h,-k\right]}"] & V\et\ral 0$~ is right-exact.

\item Similarly, the pullback-limit $\varprojlim\left(A\sr{h}{\ral}V\sr{k}{\lal}B\right)=\ker[h,-k]$, as shown in the diagram below:
\bc\adjustbox{scale=0.8}{\bt
 Y\ar[dddddr,bend right=10,"\al"]\ar[rrrrrd,bend left=10,"\beta"]\ar[dr,dashed,"\gamma"] & && && \\
  &\ker[h,-k]\ar[ddrr,hook,"\nu"]\ar[dddd,"p_A\nu"']\ar[rrrr,"p_B\nu"] && && B\ar[from=ddll,bend left,near start,dotted,"p_B"]\ar[ddll,bend left,"q_B"']\ar[dddd,"k"] \\
  & && && \\
  & && A\oplus B\ar[ddrr,"{[h,-k]}"description] && \\
  & && && \\
  & A\ar[from=uurr,bend right,dotted,"p_A"']\ar[uurr,bend right,"q_A"]\ar[rrrr,"h"] && && V
\et}\ec
where the upper-left outer square commutes because
{\footnotesize\[
kp_B\nu=-(-k)p_B\nu=-[h,-k]q_Bp_B\nu=-[h,-k]\big(-q_Ap_A+(q_Ap_A+q_Bp_B)\big)\nu=[h,-k]q_Ap_A\nu=hp_A\nu.\nn
\]}Also, if $h\al=k\beta$, then $[h,-k]q_A\al=h\al=k\beta=-[h,-k]q_B\beta$, i.e., $[h,-k](q_A\al+q_B\beta)=0$, and so there exists a unique $\gamma:Y\ra\ker[h,-k]$ such that $q_A\al+q_B\beta=\nu\gamma$. Therefore,
\[
\al=p_Aq_A\al=p_A(q_A\al+q_B\beta)=p_A\nu\gamma~~~~\txt{and}~~~~\beta=p_Aq_B\beta=p_B(q_A\al+q_B\beta)=p_B\nu\gamma.
\]
It follows by comparing with (\ref{PullPushInd}) that the diagram
{\scriptsize \bt U\ar[d,"f"]\ar[r,"g"] & B\ar[d,"k"] \\ A\ar[r,"h"] & V \et} is a pullback-limit $\iff$ the sequence ~$0\ral\bt U\ar[r,"{\left[\substack{f\\g}\right]}"] & A\oplus B\ar[r,"{\left[h,-k\right]}"] & V\et\ral 0$~ is left-exact.

\item Hence, by comparing with (\ref{PullPushInd}), the diagram
{\scriptsize \bt U\ar[d,"f"]\ar[r,"g"] & B\ar[d,"k"] \\ A\ar[r,"h"] & V \et} is both a pushout-limit and pullback-limit $\iff$ the sequence ~$0\ral\bt U\ar[r,"{\left[\substack{f\\g}\right]}"] & A\oplus B\ar[r,"{\left[h,-k\right]}"] & V\et\ra 0$~ is exact.

\item In (\ref{PullPushInd}), $q_A,q_B$ are monic and $p_A,p_B$ are epic. Indeed, recall (from the proof of $A\times B\cong A\oplus B$ in an additive category) that $p_Aq_A=id_A$, $p_Bq_B=id_B$.

\item \index{Pushout-Pullback symmetry}{\ul{Pushout-Pullback symmetry}}: In (\ref{PullPushInd}), if $f$ or $g$ is monic, then so is $\left[\substack{f\\g}\right]$. Indeed, $q_Af=q_Ap_A\left[\substack{f\\g}\right]$, and so $\left[\substack{f\\g}\right]\al=0$ implies $q_Af\al=0$, which implies $f\al=0$ (since $q_A$ is monic).

Thus, if the diagram {\scriptsize \bt U\ar[d,"f"]\ar[r,"g"] & B\ar[d,dashed,"k"] \\ A\ar[r,dashed,"h"] & P \et}  is a pushout-limit, i.e., $0\ral\bt U\ar[r,"{\left[\substack{f\\g}\right]}"] & A\oplus B\ar[r,"{\left[h,-k\right]}"] & P\et\ral 0$ is right-exact, then
\bit
\item $f$ or $g$ monic implies $0\ral\bt U\ar[r,"{\left[\substack{f\\g}\right]}"] & A\oplus B\ar[r,"{\left[h,-k\right]}"] & P\et\ral 0$ is exact, i.e., the diagram is both a pushout-limit and a pullback-limit.
\eit

\item \index{Pushout-Pullback symmetry}{\ul{Pushout-Pullback symmetry}}: In (\ref{PullPushInd}), if $h$ or $k$ is epic, then so is $[h,k]$. Indeed, $hp_A=[h,k]p_Aq_A$, and so $\al [h,k]=0$ implies $\al hp_A=0$, which implies $\al h=0$ (since $p_A$ is epic).

Thus, if the diagram {\scriptsize \bt P\ar[d,dashed,"f"]\ar[r,dashed,"g"] & B\ar[d,"k"] \\ A\ar[r,"h"] & V \et} is a pullback-limit, i.e., $0\ral\bt P\ar[r,"{\left[\substack{f\\g}\right]}"] & A\oplus B\ar[r,"{\left[h,-k\right]}"] & V\et\ral 0$~ is left-exact, then
\bit
\item $h$ or $k$ epic implies $0\ral\bt P\ar[r,"{\left[\substack{f\\g}\right]}"] & A\oplus B\ar[r,"{\left[h,-k\right]}"] & V\et\ral 0$ is exact, i.e., the diagram is both a pushout-limit and a pullback-limit.
\eit
\end{enumerate}
\end{rmk}

\begin{lmm}[\textcolor{blue}{\index{Pushout and Pullback theorems}{Pushout and Pullback theorems}: \cite[Section 7, p.9]{mitchell1965}, \cite[Section 2.5, p.51]{freyd1964}}]\label{PushPullThm}
Let $\C$ be a category and $A\sr{f}{\lal}C\sr{g}{\ral}B$, $A'\sr{f'}{\ral}C'\sr{g'}{\lal}B'$ diagrams in $\C$. Let the pushout-limit ~$P:=\varinjlim\left(A\sr{f}{\lal}C\sr{g}{\ral}B\right)$ ~be given by the following commutative diagram:
{\footnotesize
\bea\bt
C\ar[d,"f"]\ar[r,"g"] & B\ar[d,dashed,"f_1"] \\
A\ar[r,dashed,"g_1"] & P
\et\nn
\eea}
\bit
\item[(a1)] If $g$ is epic, then so is $g_1$. (Similarly, by symmetry, $f$ epic implies $f_1$.)
\item[(a2)] If $\C$ is abelian, then $g$ is epic (resp. monic) $\iff$ $g_1$ is epic (resp. monic). (Similarly, by symmetry, $f$ is epic (resp. monic) $\iff$ $f_1$ is epic (resp. monic).)
\eit

Similarly, let the pullback-limit ~$P':=\varprojlim\left(A'\sr{f'}{\ral}C'\sr{g'}{\lal}B'\right)$ ~be given by the following commutative diagram:
{\footnotesize
\bea\bt
P'\ar[d,dashed,"g'_1"]\ar[r,dashed,"f'_1"] & B'\ar[d,"g'"] \\
A'\ar[r,"f'"] & C'
\et\nn
\eea}
\bit
\item[(b1)] If $f'$ is monic, then so is $f_1'$. (Similarly, by symmetry, $g'$ monic epic implies $g_1'$ monic.)
\item[(b2)] If $\C$ is abelian, then $f'$ is monic (resp. epic) $\iff$ $f_1'$ is monic (resp. epic). (Similarly, by symmetry, $g'$ is epic (resp. monic) $\iff$ $g_1'$ is epic (resp. monic).)
\eit
\end{lmm}
\begin{proof}
By the ``pushout-pullback symmetry'' discussed in Remark \ref{PushPullRmk}, it suffices in (a2) to prove the epic property for the pushout-limit (since monic property therein will follow from that in (b2)), and suffices in (b2) to prove the monic property for the pullback-limit (since the epic property therein will follow from that in (a2)).

{\flushleft(a1)} Suppose $g$ is epic. If $\al g_1=\beta g_1$, then $\al f_1g=\al g_1f=\beta g_1f=\beta f_1g$ $\Ra$ $\al f_1=\beta f_1$. It follows that $\al=\beta$, by uniqueness in the universal property of the pushout-limit. Hence $g_1$ is also epic.
{\footnotesize
\bea\bt
C\ar[d,"f"]\ar[r,"g"] & B\ar[d,dashed,"f_1"] & \\
A\ar[r,dashed,"g_1"]  & P                    & \\
                      &                      & X\ar[from=ul,near start,bend right,"\al"']\ar[from=ul,bend left,"\beta"']
\ar[from=ull,bend right,"\al g_1=\beta g_1"']\ar[from=uul,bend left,"\beta f_1=\al f_1"]
\et\nn
\eea}

{\flushleft(a2)} Assume $\C$ is abelian. Using the universal property of the limit, we see that if an epic morphism $\bt P\ar[r,two heads,"c_1"] & C_1\et$ is a cokernel of $A\sr{g_1}{\ral}P$ then the composition $\bt c_1f_1:B\ar[r,"f_1"] & P\ar[r,two heads,"c_1"] & C_1\et$ is a cokernel of $C\sr{g}{\ral}B$.
{\footnotesize
\bea\bt
C\ar[d,"f"]\ar[r,"g"] & B\ar[d,dashed,"f_1"] \\
A\ar[r,dashed,"g_1"] & P
\et~~~~~~~~
\bt
C\ar[rrr,bend left,"\txt{Let}~tg=0."]\ar[d,"f"']\ar[r,"g"] & B\ar[d,"f_1"'] && T\ar[from=ll,"t"]\ar[from=dll,bend left=15,dashed,"h"']\ar[from=dl,dotted,near start,"h_1"'] \\
A\ar[urrr,near start,bend right=70,"\exists~0"]\ar[r,"g_1"] & P &  C_1\ar[from=l,two heads,"c_1"]
\et\nn
\eea}Thus, if $g_1$ is epic then so is $g$. Hence, by (a1), $g$ is epic $\iff$ $g_1$ is epic.

{\flushleft(b1)} Suppose $f'$ is monic. If $f_1'\al=f_1'\beta$, then $f'g_1'\al=g'f_1'\al=g'f_1'\beta=f'g_1'\beta$ $\Ra$ $g_1'\al=g_1'\beta$. It follows that $\al=\beta$, by uniqueness in the universal property of the pullback-limit. Hence $f_1'$ is also monic.

{\footnotesize
\bea\bt
 X\ar[dr,bend right,near end,"\al"']\ar[dr,bend left,"\beta"']
 \ar[ddr,bend right,"g_1'\al=g_1'\beta"']\ar[drr,bend left,"f_1'\beta=f_1'\al"] &  & \\
  & P'\ar[d,"g'_1"]\ar[r,"f'_1"] & B'\ar[d,"g'"] \\
  & A'\ar[r,"f'"] & C' \\
\et\nn
\eea}

{\flushleft (b2)} Assume $\C$ is abelian. Using the universal property of the limit, we see that if a monic morphism $\bt K_1\ar[r,hook,"k_1"] & P'\et$ is a kernel of $P'\sr{f_1'}{\ral}B'$ then $\bt K_1\ar[r,hook,"k_1"] & P'\ar[r,"g_1'"] & A'\et$ is a kernel of $A'\sr{f'}{\ral}C'$.
{\footnotesize
\bea\bt
P'\ar[d,dashed,"g'_1"]\ar[r,dashed,"f'_1"] & B'\ar[d,"g'"] \\
A'\ar[r,"f'"] & C'
\et~~~~~~~~
\bt
 & K_1\ar[r,hook,"k_1"] &  P'\ar[d,"g'_1"]\ar[r,"f'_1"] & B'\ar[d,"g'"] \\
 T\ar[urrr,near start,bend left=70,"\exists~0"]\ar[rrr,bend right,"\txt{Let}~f't=0."]\ar[ur,dotted,near end,"h_1"]\ar[urr,dashed,near end,"h"']\ar[rr,"t"] &&  A'\ar[r,"f'"] & C'
\et\nn
\eea}Thus, if $f_1'$ is monic then so is $f'$. Hence, by (b1), $f'$ is monic $\iff$ $f_1'$ is monic.
\end{proof}

\section{Basic Lemmas and Isomorphism Theorems for Abelian Categories}
\begin{lmm}[\blue{\cite[Lemma 2.61, p.54]{freyd1964}}]\label{PsPbPartII1}
Let $\A$ be an abelian category. Consider any commutative diagram with the bottom row left exact:
\[\adjustbox{scale=0.8}{\bt
        & P\ar[d,"\al"]\ar[r,"f'"] & B'\ar[d,"\beta"']\ar[dr,"g'"] &  &  \\
0\ar[r] & A\ar[r,hook,"f"] & B\ar[r,tail,"g"]  & C &
\et}~~~~\txt{or}~~~~~~~~
\adjustbox{scale=0.8}{\bt
0\ar[r] & P\ar[d,"\al"]\ar[r,"f'"] & B'\ar[d,"\beta"]\ar[r,"g'"] & C\ar[r]\ar[d,equal] & 0 \\
0\ar[r] & A\ar[r,hook,"f"] & B\ar[r,tail,"g"]  & C\ar[r] & 0
\et}\]
Then the left square (i.e., containing $P$) is a pullback $\iff$ the top row is exact at $B'$. That is, if $P$ is a pullback and $f=\ker g$, then $f'=\ker(g\beta)$.

Moreover, if $B$ is a pushout and $g=\coker f$, then ~$g\beta=\coker f'$.
\end{lmm}
\begin{proof}
{\flushleft($\Ra$)}: Assume the left square is a pullback. We need to show $\im f'\cong\ker g'$ (or $f'=k_{g'}$).
\[\adjustbox{scale=0.8}{\bt
X\ar[ddr,dashed,near end,bend right=10,"h_1"']\ar[dr,dotted,"h"description]\ar[drr,bend left=10,"x"] & & & & \\
0\ar[r] & P\ar[d,"\al"]\ar[r,near start,"f'"] & B'\ar[d,"\beta"]\ar[r,"g'"] & C\ar[r]\ar[d,equal] & 0 \\
0\ar[r] & A\ar[r,hook,"f"] & B\ar[r,tail,"g"]  & C\ar[r] & 0
\et}\]
If $0=g'x=g\beta x$, then $\beta x:X\ra B$ uniquely factors through $A\cong\ker g$ as $\beta x=fh_1$. By the pullback property, $h_1$ further uniquely factors as $h_1=\al h$, such that $x=f'h$.

{\flushleft($\La$)}: Assume the top row is exact at $B'$, i.e., $\im f'\cong\ker f'$ (or $f'=k_{g'}$). Let $u:X\ra A$ and $v:X\ra B'$ be such that $fu=\beta v$. Then $g'v=g\beta v=gfu=0u=0$, and so $v$ uniquely factors through $P$ as $v=f'h$.
\[\adjustbox{scale=0.8}{\bt
X\ar[ddr,near end,bend right=10,"u"']\ar[dr,dashed,"h"description]\ar[drr,bend left=10,"v"] & & & & \\
0\ar[r] & P\ar[d,"\al"]\ar[r,near start,"f'"] & B'\ar[d,"\beta"]\ar[r,"g'"] & C\ar[r]\ar[d,equal] & 0 \\
0\ar[r] & A\ar[r,hook,"f"] & B\ar[r,tail,"g"]  & C\ar[r] & 0
\et}\]
We also  have $f\al h=\beta f'h=\beta v=fu$, and so $u=\al h$ (since $f$ is monic).

For the moreover part (i.e., assuming $B$ is a pushout and $g=\coker f$), a similar argument gives the result as shown in the following diagram:
\[\adjustbox{scale=0.8}{\bt
        & P\ar[d,"\al"]\ar[r,"f'"] & B'\ar[d,"\beta"'] &  &  \\
0\ar[r] & A\ar[r,hook,"f"] & B\ar[r,tail,"g"]  & C & \\
 & & & & X\ar[from=ull,bend right=10,dashed,"h_1"']\ar[from=ul,dotted,"h"description]\ar[from=uull,bend left,"x"]\ar[from=ulll,bend right,near start,"0"]
\et}~~~~~~\coker f'=g\beta\]
\end{proof}

\begin{crl}[\blue{Proof analog with a pushout instead of a pullback}]
Let $\A$ be an abelian category. Consider any commutative diagram with the top row right exact:
\[\adjustbox{scale=0.8}{\bt
& A\ar[dr,"f'"']\ar[r,"f"] & B\ar[d,"\beta"]\ar[r,tail,two heads,"g"] & C\ar[d,"\gamma"]\ar[r] & 0 \\
 &  & D\ar[r,"g'"]  & P &
\et}~~~~\txt{or}~~~~
\adjustbox{scale=0.8}{\bt
0\ar[r] & A\ar[d,equal]\ar[r,"f"] & B\ar[d,"\beta"]\ar[r,tail,two heads,"g"] & C\ar[d,"\gamma"]\ar[r] & 0 \\
0\ar[r] & A\ar[r,"f'"] & D\ar[r,"g'"]  & P\ar[r] & 0
\et}\]
Then the right square (i.e., containing $P$) is a pushout $\iff$ the bottom row is exact at $D$. That is, if $P$ is a pushout and $g=\coker f$, then $g'=\coker(\beta f)$.

Moreover, if $B$ is a pullback and $f=\ker g$, then $\beta f=\ker g'$.
\end{crl}

\begin{lmm}[\blue{\cite[Lemma 2.63, p.56]{freyd1964}}]\label{PsPbPartII2}
Let $\A$ be an abelian category. Consider the following commutative diagram with the top row exact:
\[\adjustbox{scale=0.8}{\bt
        &                              &                                   & 0\ar[d]                & \\
0\ar[r] & A\ar[d,equal]\ar[r,hook,"f"] & B\ar[d,equal]\ar[r,tail,two heads,"g"] & C\ar[r]\ar[d,"\gamma"] & 0 \\
0\ar[r] & A\ar[r,hook,"f"]             & B\ar[r,"g'"]                      & C'               &
\et}\]
Then the bottom row is exact $\iff$ the right column is exact.
\end{lmm}
\begin{proof}
{\flushleft ($\Ra$):} Assume the right column is exact, i.e., $\gamma$ is a monomorphism. Then $g'=\gamma g$ implies $\ker g'=\ker(\gamma g)\cong\ker g\cong\im f$.

{\flushleft ($\La$):} Assume the bottom row is exact. Consider a pullback $P$ as shown.
\[\adjustbox{scale=0.8}{\bt
        &                              &     & 0\ar[d]                & \\
        &                              & P\ar[dl,dashed,"h"']\ar[d,hook,dashed,"\beta"]\ar[r,dashed,two heads,"\al"] & K\ar[d,hook,"k_\gamma"]\ar[r]                & 0\\
0\ar[r] & A\ar[d,equal]\ar[r,hook,"f"] & B\ar[d,equal]\ar[r,tail,two heads,"g"] & C\ar[r]\ar[d,tail,"\gamma"] & 0 \\
0\ar[r] & A\ar[r,hook,"f"]             & B\ar[r,tail,"g'"]              & C'               &
\et}\]
Because $g'\beta=\gamma g\beta=\gamma k_\gamma\al=0$, there exists $h:P\ra A$ such that $\beta=fh$, and so $k_\gamma\al=g\beta=gfh=0$, i.e., $\al=0$ (since $k_\gamma$ is monic). Hence $K=0$.
\end{proof}

\begin{crl}[\blue{Proof analog with a pushout instead of a pullback}]
Let $\A$ be an abelian category. Consider the following commutative diagram with the bottom row exact:
\[\adjustbox{scale=0.8}{\bt
 & A\ar[d,"\al"]\ar[r,"f"] & B\ar[d,equal]\ar[r,two heads,"g"] & C\ar[d,equal]\ar[r] & 0 \\
0\ar[r] & A'\ar[d]\ar[r,hook,"f'"]             & B\ar[r,tail,two heads,"g"]                      & C\ar[r]              & 0 \\
 & 0 & & &
\et}\]
Then the top row is exact $\iff$ the left column is exact.
\end{crl}

\begin{lmm}[\blue{\cite[Lemma 2.64, p.57]{freyd1964}}]\label{PsPbPartII3}
Let $\A$ be an abelian category. Consider the following commutative diagram with (i) exact columns and (ii) exact middle row:
\[\adjustbox{scale=0.8}{\bt
 & 0\ar[d] & 0\ar[d] & 0\ar[d] & \\
0\ar[r] & A\ar[d,hook,"\al"]\ar[r,"f"] & A'\ar[d,hook,"\al'"]\ar[r,"f'"] & A''\ar[d,hook,"\al''"] & \\
0\ar[r] & B\ar[d,tail,two heads,"\beta"]\ar[r,hook,"g"] & B'\ar[d,tail,"\beta'"]\ar[r,tail,"g'"] & B'' & \\
0\ar[r] & C\ar[d]\ar[r,"h"] & C' &  & \\
  & 0 &  &  &
\et}\]
Then the top row is exact $\iff$ the bottom row is exact.

\end{lmm}
\begin{proof}
$0\ra A\ra A'\ra A''$ is exact $\iff$  $0\ra A\ra A'\ra B''$ is exact (since $\al''$ is monic), $\iff$ $A$ is a pullback (Lemma \ref{PsPbPartII1}), $\iff$ $0\ra A\ra B\ra C'$ is exact (Lemma \ref{PsPbPartII1}), $\iff$ $0\ra A\ra B\ra C$ is exact and $h$ is monic (Lemma \ref{PsPbPartII2}), $\iff$ $h$ is monic (since $0\ra A\ra B\ra C$ is exact by hypotheses).
\end{proof}

\begin{crl}[\blue{Proof analog}]
Let $\A$ be an abelian category. Consider the following commutative diagram with (i) exact columns and (ii) exact middle row:
\[\adjustbox{scale=0.8}{\bt
  &  &   & 0\ar[d] & \\
  &  & A'\ar[d,"\al'"]\ar[r,"f'"] & A''\ar[d,hook,"\al''"]\ar[r] & 0\\
  & B\ar[d,two heads,"\beta"]\ar[r,"g"] & B'\ar[d,tail,two heads,"\beta'"]\ar[r,tail,two heads,"g'"] & B''\ar[d,tail,two heads,"\beta''"]\ar[r] & 0\\
  & C\ar[d]\ar[r,"h"] & C'\ar[d]\ar[r,"h'"] & C''\ar[d]\ar[r] & 0\\
  & 0 & 0 & 0 &
\et}\]
Then the bottom row is exact $\iff$ the top row is exact.
\end{crl}

\begin{lmm}[\blue{\index{Nine (3-by-3) Lemma}{Nine (3-by-3) Lemma}: \cite[Lemma 2.65, p.58]{freyd1964}}]\label{PsPbPartII4}
Let $\A$ be an abelian category. Consider the following commutative diagram with (i) exact columns and (ii) exact middle row:
\[\adjustbox{scale=0.8}{\bt
 & 0\ar[d] & 0\ar[d] & 0\ar[d] & \\
0\ar[r] & A\ar[d,hook,"\al"]\ar[r,"f"] & A'\ar[d,hook,"\al'"]\ar[r,"f'"] & A''\ar[d,hook,"\al''"]\ar[r] & 0 \\
0\ar[r] & B\ar[d,tail,two heads,"\beta"]\ar[r,hook,"g"] & B'\ar[d,tail,two heads,"\beta'"]\ar[r,tail,two heads,"g'"] & B''\ar[d,tail,two heads,"\beta''"]\ar[r] & 0\\
0\ar[r] & C\ar[d]\ar[r,"h"] & C'\ar[d]\ar[r,"h'"] & C''\ar[d]\ar[r] & 0\\
  & 0 & 0 & 0 &
\et}\]
Then top row is exact $\iff$ the bottom row is exact.
\end{lmm}
\begin{proof}
Using Lemma \ref{PsPbPartII3} and its corollary, the top row is exact $\iff$ it is left exact and $f'$ is epic, $\iff$ $h$ is monic and the bottom row is right exact, $\iff$ the bottom row is exact.
\end{proof}

\begin{lmm}[\blue{\index{Split-exactness lemma}{Split-exactness lemma}: \cite[Lemma 2.68, p.59]{freyd1964}}]\label{PsPbPartII5}
Let $\A$ be an abelian category and $f:A\ra B$ a morphism with a left inverse (i.e., there exists a morphism $f':B\ra A$ such that $f'f=id_A$). Then a sequence $0\ra A\sr{f}{\ral}B\sr{g}{\ral}C\ra 0$ is exact $\iff$ there exists a morphism $g':C\ra B$ such that (i) the sequence $0\la A\sr{f'}{\lal}B\sr{g'}{\lal}C\la 0$ is exact, (ii) $gg'=id_C$, and (iii) $ff'+g'g=id_B$.

Moreover, given a sequence {\footnotesize $0\ra A\sr{f}{\ral}B\sr{g}{\ral}C\ra 0$}, the following are equivalent (called \ul{split-exactness}).
\begin{enumerate}[leftmargin=0.9cm]
\item[(1)] The sequence is exact and $f$ is split.
\item[(2)] The sequence is exact and $g$ is split.
\item[(3)] There exist morphisms $f':B\ra A$ and $g':C\ra B$ such that $f'f=id_A$, $gg'=id_C$, $ff'+g'g=id_B$.
\item[(4)] We have an isomorphism of diagrams:
\[\adjustbox{scale=0.7}{\bt
0\ar[r] & A\ar[d,"id_A"]\ar[r,hook,"f"] & B\ar[d,"\al","\cong"']\ar[r,tail,two heads,"g"] & C\ar[d,"id_C"]\ar[r] & 0\\
0\ar[r] & A\ar[r,hook,"q_A"] & A\oplus C\ar[r,tail,two heads,"p_C"] & C\ar[r] & 0\\
\et}\]
\end{enumerate}
\end{lmm}
\begin{proof}
{\flushleft ($\Ra$)} Assume the sequence $0\ra A\sr{f}{\ral}B\sr{g}{\ral}C\ra 0$ is exact. Then using the 3-by-3 lemma, we get an exact commutative diagram in the following form:
\[\adjustbox{scale=0.8}{\bt
 &  & 0\ar[d] & 0\ar[d] & \\
 & 0\ar[d]\ar[r] & C'\ar[d,hook,"k_{f'}"]\ar[r,dashed,"u","\cong"'] & C\ar[d,equal]\ar[r] & 0 \\
0\ar[r] & A\ar[d,equal]\ar[r,hook,"f"] & B\ar[d,tail,two heads,"f'"]\ar[r,tail,two heads,"g"] & C\ar[d]\ar[r] & 0\\
0\ar[r] & A\ar[d]\ar[r,equal] & A\ar[d]\ar[r] & 0 & \\
  & 0 & 0 &  &
\et}\]
Let $g':=k_{f'}u^{-1}$. Then $gg'=gk_{f'}u^{-1}=uu^{-1}=id_C$. We therefore have induced morphisms $\al,\beta$ as follows:
\[\adjustbox{scale=0.7}{\bt
  &   & B\ar[ddl,bend right,"f'"']\ar[d,dashed,"\al"]\ar[ddr,bend left,"g"]  &  \\
  &   & A\times B\ar[dl,"p_A"']\ar[dd,equal]\ar[dr,"p_C"] & \\
  & A\ar[ddr,bend right,"f"']\ar[dr,"q_A"] &   & C\ar[ddl,bend left,"g'"]\ar[dl,"q_C"']\\
  &   & A\oplus B\ar[d,dashed,"\beta"] & \\
0\ar[r] & A\ar[r,"f"'] & B\ar[r,"g"'] & C\ar[r] & 0
\et}\]
By the commutativity of the diagram, we have $f'=p_A\al$, $g=p_C\al$, $f=\beta q_A$, $g'=\beta q_C$, and so
\begin{align}
&ff'+g'g=\beta q_Ap_A\al+\beta q_Cp_C\al=\beta(q_Ap_A+q_Cp_C)\al=\beta id_{A\oplus C}\al=\beta\al,\nn\\
& p_Aq_A=id_A=f'f=p_A\al\beta q_A,~~~~p_Cq_C=id_C=gg'=p_C\al\beta q_C,\nn\\
&~~\sr{\txt{uniqueness}}{\Longrightarrow}~~\al\beta=q_Ap_A+q_Cp_C,~~~~\beta\al=\beta(\al\beta)\al=(\beta\al)^2,\nn\\
&~~\sr{\txt{uniqueness}}{\Longrightarrow}~~\al\beta=id_{A\oplus B}=q_Ap_A+q_Cp_C,~~~~\beta\al=id_B=ff'+g'g.\nn
\end{align}
Explicitly, the isomorphism $\al,\beta$ is given by
\[
\al=q_Af'+q_Cg~~~~\txt{and}~~~~\beta=fp_A+g'p_C.
\]

{\flushleft ($\La$)} Assume (i), (ii), (iii) hold. Then it is clear that $f$ is monic, $g$ is epic, and $gf=0$ (by applying $g$ on $ff'+g'g=id_B$ from the left and then using the fact that $f'$ is epic). Also, given a morphism $x:X\ra B$ such that $gx=0$, if we apply $x$ on $ff'+g'g=id_B$ from the right, we get $ff'x=x$ (i.e., $x$ uniquely factors through $f$, and so $f=k_g$).

The ``moreover'' part (with the stated equivalences) is clear from the above information and symmetry/duality.
\end{proof}

\begin{crl}
Let $\A,\A'$ be additive categories. An additive functor $F:\A\ra\A'$ preserves (i) split-exactness of sequences and hence, in particular, preserves (ii) finite products/coproducts in the sense that $F(A\oplus B)\cong F(A)\oplus F(B)$ for all $A,B\in\Ob\A$.
\end{crl}

\begin{rmk}[\blue{\index{Intersection as a pullback}{Intersection as a pullback} and \index{Union! as a pushout}{Union as a pushout}: \index{Union-Intersection isomorphism}{Union-Intersection isomorphism}}]
Let $\C$ be a category and $A,B\in\Ob\C$. Define $A\cap B$ and $A\cup B$ to be the objects such that the square {\small $A\sr{v'}{\hookleftarrow} A\cap B\sr{u'}{\hookrightarrow} B\sr{v}{\hookrightarrow} A\cup B\sr{u}{\hookleftarrow} A$} is both a pullback and a pushout. If $\C$ is abelian then, as seen earlier, this is equivalent to the assertion that the associated sequence ~{\small $0\ra A\cap B\sr{\left[\substack{v'\\u'}\right]}{\ral}A\oplus B\sr{[u,-v]}{\ral}A\cup B\ra 0$}.~ That is, in an abelian category, ~$A\cup B\cong{A\oplus B\over A\cap B}$.
\end{rmk}

\begin{thm}[\blue{\index{Preimage/Integersection-kernel formulas}{Preimage/Intersection-kernel formulas}}]
Let $\A$ be an abelian category $A,B,C\in\Ob\A$, and $f:A\ra C$.
\bit[leftmargin=0.9cm]
\item[(i)] If $A,B\subset C$, we have an exact sequence ~$0\ra A\cap B\ra A\ra C/B\ra 0$~ in the form
\[
A\cap B:=\varprojlim\big(A\sr{u}{\hookrightarrow}C\sr{v}{\hookleftarrow}B\big)~\cong~\ker\big(A\sr{u}{\hookrightarrow} C\sr{c_v}{\twoheadrightarrow} C/B\big)~\cong~\ker\big(B\sr{v}{\hookrightarrow}C\sr{c_u}{\twoheadrightarrow} C/A\big).
\]
\item[(ii)] If $f:A\ra B$ and $B\subset C$, we have an exact sequence ~$0\ra f^{-1}(B)\ra A\ra C/B\ra 0$~ in the form
\[
f^{-1}(B):=\varprojlim\big(A\sr{f}{\ral}C\sr{v}{\hookleftarrow}B\big)~\cong~\ker\big(A\sr{f}{\ral} C\sr{c_v}{\twoheadrightarrow} C/B\big).
\]
\item[(iii)] $f^{-1}(0)\cong \ker\big(A\sr{f}{\ral} C\sr{id_C}{\ral} C\big)=\ker f$.
\eit
\end{thm}
\begin{proof}
Consider the diagrams below:
\[(i)~~\adjustbox{scale=0.7}{\bt
 & 0\ar[d]        &  0\ar[d]   & 0\ar[d] &     \\
0\ar[r] & A\cap B=P\ar[d,hook,"u'"]\ar[r,hook,"v'"] & A\ar[d,hook,"u"]\ar[r,tail,two heads,"c_{v'}"]   & {A\over A\cap B}\ar[d,hook,dashed,"u''"]\ar[r] &  0   \\
0\ar[r] & B\ar[d,tail,two heads,"c_{u'}"]\ar[r,hook,"v"] & C\ar[d,tail,two heads,"c_u"]\ar[r,tail,two heads,"c_v"] & C/B\ar[r] & 0 \\
0\ar[r] & {B\over A\cap B}\ar[r,hook,dashed,"v''"]\ar[d]     & C/A\ar[d] &  & \\
 &    0     & 0 &  &
\et}~~~~~~~~
(ii)~~\adjustbox{scale=0.7}{\bt
 &         &     & 0\ar[d] &     \\
0\ar[r] & f^{-1}(B)\ar[d,"f'"]\ar[r,hook,"v'"] & A\ar[d,"f"]\ar[r,dashed,"h'"]   & \ker\al\ar[d,hook,"k_\al"]\ar[r] &  0   \\
0\ar[r] & B\ar[d,dashed,"h"]\ar[r,hook,"v"] & C\ar[d,tail,two heads,"c_f"]\ar[r,tail,two heads,"c_v"] & C/B\ar[d,tail,two heads,"\al"]\ar[r] & 0 \\
0\ar[r] & \ker\beta\ar[r,hook,"k_\beta"]\ar[d] & C/A\ar[d]\ar[r,tail,two heads,"\beta"] & P\ar[d]\ar[r] & 0\\
 &    0     & 0 & 0 &
\et}\]
(i) Since $A\cap B$ is a pullback, {Lemma \ref{PsPbPartII1}
} implies
\[
A\cap B=\ker(c_vu)\sr{\txt{symmetry}}{=}\ker(c_uv).
\]
(ii) Similarly, because $f^{-1}(B)$ is a pullback, {Lemma \ref{PsPbPartII1}
} implies
\[
f^{-1}(B)=\ker(c_vf). \qedhere
\]
\end{proof}

\begin{crl}
Let $\A$ be an abelian category and $f:A\ra B$ a morphism and $u:C\hookrightarrow B$ a subobject. If $f$ is \ul{monic}, then ~$f^{-1}(C)\cong A\cap B\cong\im f~\cap~C$.
\end{crl}
\begin{proof}
This is an immediate consequence of part (i) of the preceding theorem:
\[
f^{-1}(C):=\varprojlim(A\sr{f}{\hookrightarrow}B\hookleftarrow C)\cong A\cap C\cong\im f~\cap~C. \qedhere
\]
\end{proof}

\begin{prp}[\blue{Preimage and Kernel properties}]
Let $\C$ be a category. Consider any morphisms $f:A\ra B$ and $g:B\ra C$, and any objects $D,D',D''\in\Ob\C$.
\begin{enumerate}[leftmargin=0.5cm]
\item If $C'\subset C$, then ~~$(gf)^{-1}(C')\cong f^{-1}\big(g^{-1}(C')\big)$.
\item $\ker(gf)\cong f^{-1}(\ker g)$  ~~ (\magenta{Assume $\C$ is abelian.})
\item If $B'\subset B$, then ~~$f^{-1}(B')\cong f^{-1}(\im f\cap B')$. ~~(\magenta{Assume $\C$ is abelian.})
\item If $f$ is monic, then ~$\ker(gf)\cong \im f\cap\ker g$  ~~(\magenta{Assume $\C$ is abelian.})
\item If $B',B''\subset B$, then ~~$f^{-1}(B'\cap B'')\subset f^{-1}(B')\cap f^{-1}(B'')\subset f^{-1}(B')\cup f^{-1}(B'')\subset f^{-1}(B'\cup B'')$.
\item $(D\cap D')\cup(D\cap D'')\subset D\cap(D'\cup D'')$ ~~~~(\magenta{Assume $\C$ is abelian.})
\end{enumerate}
\end{prp}
\begin{proof}
\begin{enumerate}[leftmargin=0.5cm]
\item This is due to the existence of two-way induced morphisms as in the commutative diagram:
\[\adjustbox{scale=0.8}{\bt
 (gf)^{-1}(C')\ar[ddr,hook,near start,bend right=10,"\delta"']\ar[from=dr,bend left=10,dashed,"h"description]\ar[drr,dotted,bend left=10,"h_1"description]\ar[dr,dashed,bend left=10,"h_1'"description]\ar[rrr,"(gf)'"] &  &  & C'\ar[d,equal] \\
 & f^{-1}(g^{-1}(C'))\ar[d,hook,"\gamma"]\ar[r,near start,"f'"] & g^{-1}(C')\ar[d,hook,"\beta"]\ar[r,"g'"] & C'\ar[d,hook,"\al"] \\
 & A\ar[r,"f"] & B\ar[r,"g"] & C
\et}\]
The unique morphism $h$ exists because $(gf)^{-1}(C')$ is a pullback from the $C$ vertex. Next, the unique morphism $h_1$ similarly exists because $g^{-1}(C')$ is a pullback from the $C$ vertex. Finally, due to $h_1$, the unique morphism $h_1'$ then also exists because $f^{-1}(g^{-1}(C'))$ is a pullback from the $B$ vertex.
\item $\ker(gf)\cong(gf)^{-1}(0)\cong f^{-1}(g^{-1}(0))\cong f^{-1}(\ker g)$.

\item Let \bt f=m_fe_f:A\ar[r,two heads,"e_f"] & \im f\ar[r,hook,"m_f"]& B.\et Then
\begin{align}
&f^{-1}(B')=(m_fe_f)^{-1}(B')\cong e_f^{-1}(m_f^{-1}(B'))\cong e_f^{-1}(\im f\cap B')\cong e_f^{-1}(m_f^{-1}(\im f\cap B'))\nn\\
&~~~~\cong(m_fe_f)^{-1}(\im f\cap B')=f^{-1}(\im f\cap B').\nn
\end{align}

\item $\ker(gf)\cong f^{-1}(\ker g)\sr{(s)}{\cong}\im f\cap\ker g$, ~where step (s) holds because $f$ is monic (preceding corollary).

\item The containments are due to the induced maps in the following commutative diagram.
\[\adjustbox{scale=0.8}{\bt
 && f^{-1}(B')\ar[ddd,bend right=60,near start,"f_2"]\ar[dr,hook]\ar[drr,hook,dashed]\ar[drrrr,hook,bend left=10] && && \\
 f^{-1}(B'\cap B'')\ar[ddd,"f_4"]\ar[urr,dashed,hook]\ar[drr,dashed,hook]\ar[r,hook,dotted]\ar[rrrrrr,bend left=25] & f^{-1}(B')\cap f^{-1}(B'')\ar[ur,hook]\ar[dr,hook] &  & f^{-1}(B')\cup f^{-1}(B'')\ar[r,dotted,hook] & f^{-1}(B'\cup B'')\ar[ddd,"f_1"]\ar[rr,hook,near start,"u_1"] &&  A\ar[ddd,"f"] \\
 && f^{-1}(B'')\ar[ddd,bend left=50,near start,"f_3"]\ar[ur,hook]\ar[urr,hook,dashed]\ar[urrrr,hook,bend right=10] && && \\
 && B'\ar[drr,hook] && && \\
B'\cap B''\ar[urr,hook]\ar[drr,hook] &&  && B'\cup B''\ar[rr,hook,"u"]&& B\\
 && B''\ar[urr,hook] && && \\
\et}\]

\item This follows from the preceding result: Consider the inclusion $D\sr{f}{\hookrightarrow}D\cup D'\cup D''$. Then
\[
f^{-1}(D')\cup f^{-1}(D'')\subset f^{-1}(D'\cup D''). \qedhere
\]
\end{enumerate}
\end{proof}

\begin{crl}[\blue{Postimage, Cupimage, Postimage/Cupimage isomorphisms}]
{\flushleft (i)} Let $\C$ be a category and consider morphisms $f:A\ra B$ and $g:B\ra C$. Given a quotient $A\sr{\pi}{\twoheadrightarrow}A'$, define the \ul{$f$-postimage} of $A'$ to be the pushout $f^{+1}(A'):=\varprojlim(A'\sr{\pi}{\twoheadleftarrow}A\sr{f}{\ral}B)$. Then
\[
(gf)^{+1}(A')\cong g^{+1}\big(f^{+1}(A')\big).
\]
In particular, if $\C$ is abelian, then we also get the following:
\[
\coker f\cong f^{+1}(0)~~~~\txt{and}~~~~\coker(gf)\cong g^{+1}(\coker f).
\]

{\flushleft (ii)} Let $\A$ be an abelian category, and consider morphisms $f:A\ra B$ and $g:B\ra C$. Given an inclusion $A\sr{u}{\hookrightarrow}D$, define the \ul{$f$-cuptimage} of $D$ to be the pushout $f^\cup(D):=\varprojlim(D\sr{u}{\hookleftarrow}A\sr{f}{\ral}B)$. Then
\[
(gf)^\cup(D)\cong g^\cup\big(f^\cup(D)\big).
\]
If further given an inclusion $B'\sr{v}{\hookrightarrow}B$, we may also define the \ul{$f$-union} of $D$ and $B'$ as follows:
\[
D\cup_f B':=\big(f|_{f^{-1}(B')}\big)^\cup(D).
\]
In particular, when $f$ is monic (hence $f^{-1}(A')\cong A\cap B'$) then
\[
A\cup_fB':=\big(f|_{f^{-1}(B')}\big)^\cup(A)\cong A\cup B'.
\]
\end{crl}
\begin{proof}
As in the proof of part (1) of the preceding theorem, unique two-way morphisms exist to give the desired isomorphisms as in the diagrams:
\[(i)~~~~\adjustbox{scale=0.7}{\bt
A\ar[d,two heads,"\pi"]\ar[r,"f"] & B\ar[d,two heads,"\pi_1"]\ar[r,"g"] & C\ar[d,two heads,"\pi_2"]\ar[ddr,two heads,bend left,"\pi_3"] & \\
A'\ar[d,equal]\ar[r,"f'"] & f^{+1}(A')\ar[r,"g'"] & g^{+1}(f^{+1}(A')) & \\
A'\ar[rrr,near start,"(gf)^{+1}"] &  &  &  (gf)^{+1}(A')\ar[from=ull,dotted,bend right=10,"h_1"description]\ar[from=ul,dashed,bend right=10,"h_1'"description]\ar[ul,dashed,bend right=10,"h"description]
\et}~~~~
(ii)~~~~\adjustbox{scale=0.7}{\bt
A\ar[d,hook,"u"]\ar[r,"f"] & B\ar[d,hook,"u_1"]\ar[r,"g"] & C\ar[d,hook,"u_2"]\ar[ddr,hook,bend left,"u_3"] & \\
D\ar[d,equal]\ar[r,"f'"] & f^{\cup}(D)\ar[r,"g'"] & g^{\cup}(f^{\cup}(D)) & \\
D\ar[rrr,near start,"(gf)^{\cup}"] &  &  &  (gf)^{\cup}(D)\ar[from=ull,dotted,bend right=10,"h_1"description]\ar[from=ul,dashed,bend right=10,"h_1'"description]\ar[ul,dashed,bend right=10,"h"description]
\et}\]
\end{proof}

\begin{thm}[\blue{\index{Isomorphism theorems}{Isomorphism theorems}: \cite[2.66 and 2.77, pp 58-59]{freyd1964}}]
Let $\A$ be an abelian category and $A,B,C\in\Ob\A$. (i) If $C\subset A\subset B$, then there is an exact sequence ~$0\ra A/C\ra B/C\ra B/A\ra 0$,~ or equivalently, we have an isomorphism
\[\textstyle {B\over A}\cong{B/C\over A/C}.\]
(ii) If $A,B\subset C$, then there is an isomorphism ~${C/A\over B/(A\cap B)}\cong {C/B\over A/(A\cap B)}$. In particular, if $C=A\cup B$, then we have exact sequences $0\ra A\cap B\ra A\ra {A\cup B\over B}\ra 0$ and $0\ra A\cap B\ra B\ra {A\cup B\over A}\ra 0$, or equivalently, we have isomorphisms
\[\textstyle {A\over A\cap B}\cong{A\cup B\over B}~~~~\txt{and}~~~~{B\over A\cap B}\cong{A\cup B\over A}.\]
\end{thm}
\begin{proof}
(i) From the given information, we can obtain the following commutative diagram in which all columns, first row, and middle row are exact. So, the 3-by-3 lemma implies the bottom row is exact.
\[\adjustbox{scale=0.8}{\bt
        & 0\ar[d]& 0\ar[d]   & 0\ar[d]   &   \\
0\ar[r] & C\ar[d,hook,"u"]\ar[r,equal]   & C\ar[d,hook,"vu"]\ar[r]   & 0\ar[d]\ar[r]   & 0 \\
0\ar[r] & A\ar[d,tail,two heads,"c_u"]\ar[r,hook,"v"] & B\ar[d,tail,two heads,"c_{vu}"]\ar[r,tail,two heads,"c_v"]   & B/A\ar[d,equal]\ar[r] & 0 \\
0\ar[r] & A/C\ar[d]\ar[r,dashed] & B/C\ar[d]\ar[r,dashed] & B/A\ar[d]\ar[r] & 0 \\
        &  0   &  0   &   0  &
\et}\]

(ii) Consider the following commutative diagram (obtained directly from the given information). By Lemma \ref{PsPbPartII1}, the induced morphisms $u',v'$ are monic and so the diagram is exact by the 3-by-3 lemma.
\[\adjustbox{scale=0.7}{\bt
        & 0\ar[d]& 0\ar[d]   & 0\ar[d]   &   \\
0\ar[r] & A\cap B=P\ar[d,hook,"\beta"]\ar[r,hook,"\al"] & B\ar[d,hook,"u"]\ar[r,tail,two heads,"c_\al"] & {B\over A\cap B}\ar[d,dashed,hook,"u'"]\ar[r]   & 0 \\
0\ar[r] & A\ar[d,tail,two heads,"c_\beta"]\ar[r,hook,"v"] & C\ar[d,tail,two heads,"c_u"]\ar[r,tail,two heads,"c_v"]   & {C\over A}\ar[d,tail,"c_{u'}"]\ar[r] & 0 \\
0\ar[r] & {A\over A\cap B}\ar[d]\ar[r,dashed,hook,"v'"] & {C\over B}\ar[d]\ar[r,dashed] & {C/A\over B/(A\cap B)}\ar[d]\ar[r] & 0 \\
        &  0   &  0   &   0  &
\et}\]
Because intersection is a pullback and union is a pushout, if we set $C:=A\cup B$ then it follows from Lemma \ref{PsPbPartII1} that the induced morphisms $u',v'$ are isomorphisms.
\end{proof}

\begin{crl}[\blue{\index{Preimage-image formula}{\blue{Preimage-image formula}}}]
Let $\A$ be an abelian category and $f:A\ra B$. If $B'\subset B$, then there exists an exact sequence ~$0\ra f^{-1}(B')\ra A\ra {\im f\over\im f\cap B'}\ra 0$.
\end{crl}
\begin{proof}
Consider the factorization $f=m_fe_f:A\sr{e_f}{\twoheadrightarrow}\im f\sr{m_f}{\hookrightarrow}B$. Then we have the following commutative diagrams:
\[
\adjustbox{scale=0.9}{\bt
0\ar[r] & f^{-1}(B')\ar[d,two heads,"f'"]\ar[r,hook,"u'"] & A\ar[d,"f"] &  & \\
0\ar[r] & B'\ar[r,hook,"u"] & B\ar[r,tail,two heads,"c_u"] & {B\over B'}\ar[r] & 0
\et}~~~~=~~~~
\adjustbox{scale=0.8}{\bt
0\ar[r] & \substack{f^{-1}(B')=\\ e_f^{-1}(m_f^{-1}(B'))}\ar[d,two heads,"e'"]\ar[r,hook,"u'"] & A\ar[d,two heads,"e_f"] &  & \\
0\ar[r] & m_f^{-1}(B')\ar[d,hook,"m'"]\ar[r,hook,"u_1"] & \im f\ar[d,hook,"m_f"]\ar[r,tail,two heads,"c_{u_1}"] & {\im f\over m_f^{-1}(B')}\ar[d,dashed,"\cong"',"h_1"]\ar[r] & 0\\
0\ar[r] & B'\ar[d,equal]\ar[r,hook,"u_2"] & \im f\cup B'\ar[d,hook,"q"]\ar[r,tail,two heads,"c_{u_2}"] & {\im f\cup B'\over B'}\ar[d,dashed,"h_2"]\ar[r] & 0 \\
0\ar[r] & B'\ar[r,hook,"u"] & B\ar[r,tail,two heads,"c_u"] & {B\over B'}\ar[r] & 0
\et}\]
Since $m_f^{-1}(B')$ is a pullback and $\im f\cup B'$ is a pushout, {Lemma \ref{PsPbPartII1}
} implies $h_1$ is an isomorphism. Because we also have $f^{-1}(B')=e_f^{-1}(m_f^{-1}(B'))$, it follows that $\coker u'=c_{u_1}e_f$, and so we get an exact sequence
\[
\textstyle 0\ral f^{-1}(B')\ral A\ral{\im f\over m_f^{-1}(B')}\cong {\im f\cup B'\over B'}\cong{\im f\over\im f\cap B'}\cong{B'\over\im f\cap B'}\ral 0.
\]
\end{proof}

\section{Additive and Abelian Properties of the Category of Functors}
\begin{prp}\label{YonProp}
Let $\I,\A$ be categories, and $\A^\I$ the category of functors $\I\ra\A$. If $\A$ is additive (resp. abelian), then so is $\A^{\I}$.
\end{prp}
\begin{proof}
{\flushleft(1)} \ul{Zero object}: A zero object in $\A^\I$ is given by the zero functor $0:\I\ra\A,~\big(I\sr{f}{\ral}J\big)~\mapsto~\big(0\sr{0}{\ral}0\big)$.
{\flushleft (2)} \ul{Pre-additivity}: Let $F,G\in\Ob\A^\I$. Then morphisms $\eta,\eta'\in \Mor_{\A^\I}(F,G)$ are given by $\eta=\{\eta_I:F(I)\ra G(I)\}_{I\in\Ob\I}$, $\eta'=\{\eta'_I:F(I)\ra G(I)\}_{I\in\Ob\I}$, and so $\Mor_{\A^\I}(F,G)$ is an abelian group with respect to
\bea
\eta+\eta':=\{\eta_I+\eta'_I:F(I)\ra G(I)\}_{I\in\Ob\I},~~~~
\bt I\ar[d,"f"] \\ J\et~~\longmapsto~~
\bt
F(I)\ar[d,"F(f)"] \ar[rrr,"{\eta_I,~\eta'_I,~\eta_I+\eta'_I}"] &&& G(I)\ar[d,"G(f)"] \\
F(J) \ar[rrr,"{\eta_J,~\eta'_J,~\eta_J+\eta'_J}"] &&& G(J)
\et\nn
\eea
Hence, pointwise composition of functors (as defined earlier) is distributive pointwise over the above pointwise addition in $\Mor_{\A^\I}(F,G)$.

{\flushleft (3)} \ul{Finite coproducts/products}: Given functors $F,G:\I\ra\A$, we define $F\oplus G:\I\ra\A$ by
\bea
F\oplus G~:~I\sr{f}{\ral}J~~\mapsto~~F(I)\oplus G(I)\sr{F(f)\oplus G(f)}{\ral}F(J)\oplus G(J),\nn
\eea
where it is clear that $F\oplus G$ satisfies the universal property (pointwise), in the sense that given natural transformations $\al:F\ra H$, $\beta:G\ra H$, there exists a unique natural transformation $\eta:F\oplus G\ra H$ such that the following diagram commutes:
{\scriptsize\bea\bt
F\ar[ddr,bend right,"\al"']\ar[dr,hook,"u"] &   & G\ar[ddl,bend left,"\beta"]\ar[dl,hook,"v"'] \\
  & F\oplus G\ar[d,dashed,"\eta"] &   \\
  & H &
\et:~~I\sr{f}{\ral}J~~~~\mapsto~
\bt
  & F(J)\ar[ddr,bend right,"\al_J"']\ar[dr,hook,"u_J"] &   & G(J)\ar[ddl,bend left,"\beta_J"]\ar[dl,hook,"v_J"'] \\
  & & F(J)\oplus G(J)\ar[d,dashed,"\eta_J"] & \\
  & & H(J) & \\
F(I)\ar[uuur,dashed,"{F(f)}"]\ar[ddr,bend right,"\al_I"']\ar[dr,hook,"u_I"] &   & G(I)\ar[uuur,dashed,near start,"{G(f)}"']\ar[ddl,bend left,"\beta_I"]\ar[dl,hook,near start,"v_I"] &\\
  & F(I)\oplus G(I)\ar[uuur,dashed,"{F(f)\oplus G(f)}"]\ar[d,dashed,"\eta_I"'] &   &\\
  & H(I)\ar[uuur,dashed,near start,"{H(f)}"'] & &
\et\nn
\eea}
where $u,v$ are the canonical inclusions.

{\flushleft (4)} \ul{Kernels and cokernels}: As before, each natural transformation $\eta:F\ra G$ (for functors $F,G:\I\ra \A$) is by definition a collection of morphisms $\left\{F(I)\sr{\eta_I}{\ral}G(I)\right\}_{I\in\Ob\I}$ in $\A$ along with its naturality, which is given componentwise by the following commutative diagram for each morphism $f:I\ra J$ in $\I$:
\bea
\bt I\ar[d,"f"] \\ J\et~~\sr{\eta}{\longmapsto}~~
 \bt
  F(I)\ar[d,"F(f)"] \ar[r,"\eta_I"] & G(I)\ar[d,"G(f)"] \\
  F(J) \ar[r,"\eta_J"] & G(J)
\et\nn
\eea
It thus follows that the kernel and cokernel of $\eta:F\ra G$ can be defined using a diagram (\blue{footnote}\footnote{With universal properties understood as part of the definitions.})
\bea
K\eta C:~\bt
0\ar[r] & \ker\eta \ar[r,"k_\eta"] & F \ar[r,"\eta"] & G\ar[r,"c_\eta"]\ar[r] &\coker\eta \ar[r] & 0
\et\nn
\eea
along with its naturality (i.e., a commutative diagram) given componentwise by the commutative diagram:
\bea
\bt I\ar[d,"f"] \\ J\et~~\sr{K\eta C}{\longmapsto}~~
\bt
0\ar[r] & (\ker\eta)(I)\ar[d,dashed,"(\ker\eta)(f)"] \ar[r,"(k_\eta)_I"] & F(I)\ar[d,"F(f)"] \ar[r,"\eta_I"] & G(I)\ar[d,"G(f)"]\ar[r,"(c_\eta)_I"]\ar[r] &(\coker\eta)(I)\ar[d,dashed,"(\coker\eta)(f)"] \ar[r] & 0 \\
0\ar[r] & (\ker\eta)(J) \ar[r,"(k_\eta)_J"] & F(J) \ar[r,"\eta_J"] & G(J) \ar[r,"(c_\eta)_J"]\ar[r] &(\coker\eta)(J)
  \ar[r] & 0
\et\nn
\eea
where $(\ker\eta)(I):=\ker\eta_I$ and $(\coker\eta)(I):=\coker\eta_I$, while $(\ker\eta)(f)$ and $(\coker\eta)(f)$ are the respective kernel-induced and cokernel-induced morphisms in $\A$ due to the following identities:
\bea
&&\eta_J\circ F(f)\circ(k_\eta)_I=G(f)\circ\eta_I\circ(k_\eta)_I=G(f)\circ 0=0,\nn\\
&&(c_\eta)_J\circ G(f)\circ\eta_I=(c_\eta)_J\circ\eta_J\circ F(f)=0\circ F(f)=0.\nn
\eea
{\flushleft (5)} \ul{Monomorphisms as kernels and Epimorphisms as cokernels}: Since every property of $\A$ (including ``every monic morphism is a kernel'' and ``every epic morphism is a cokernel'') holds for each $\eta_I:F(I)\ra G(I)$, $I\in\Ob\I$, it also holds for $\eta$ by construction.

Hence $\A^{\I}$ is indeed additive (resp. abelian) if $\A$ is additive (resp. abelian).
\end{proof}

In general, given a system $S:\I_S\ra\A^\I$ in $\A^\I$ (essentially a collection of transition morphisms $\{S(\kappa_{ij}):S_i\ra S_j\}_{i,j\in\I_S}$) its naturality is given (for each pair $i,j\in\Ob\I_S$) by the commutative diagram
{\footnotesize
\bea
\bt I\ar[d,"f"] \\ J\et~~\sr{S(\kappa_{ij})}{\longmapsto}~~
\bt
S_i(I)\ar[d,"S_i(f)"]\ar[rr,"S(\kappa_{ij})_I"] && S_j(I)\ar[d,"S_j(f)"]\\
S_i(J)\ar[rr,"S(\kappa_{ij})_J"] && S_j(J)
\et\nn
\eea}
and so the naturality of its \index{Transition kernel of a system}{\ul{transition kernel}} and \index{Transition cokernel of a system}{\ul{transition cokernel}} (\blue{footnote}\footnote{The transition kernel (resp. transition cokernel) of a system consists of the kernels (resp. cokernels) of its transition morphisms.})
\bea
KS(\kappa_{ij})C:\bt
0\ar[r] & \ker S(\kappa_{ij}) \ar[r,"k_{S(\kappa_{ij})}"] & S_i \ar[r,"S(\kappa_{ij})"] & S_j\ar[r,"c_{S(\kappa_{ij})}"]\ar[r] &\coker S(\kappa_{ij}) \ar[r] & 0
\et\nn
\eea
are similarly given by the commutative diagram
{\scriptsize\[
\bt I\ar[dd,"f"] \\~\\ J\et~~\sr{KS(\kappa_{ij})C}{\longmapsto}~~
\bt
0\ar[r] & (\ker{S(\kappa_{ij})})(I)\ar[dd,"(\ker{S(\kappa_{ij})})(f)"] \ar[rr,"(k_{S(\kappa_{ij})})_I"] && S_i(I)\ar[dd,"S_i(f)"] \ar[rr,"{S(\kappa_{ij})}_I"] && S_j(I)\ar[dd,"S_j(f)"]\ar[rr,"(c_{S(\kappa_{ij})})_I"]\ar[rr] &&(\coker{S(\kappa_{ij})})(I)\ar[dd,"(\coker{S(\kappa_{ij})})(f)"] \ar[r] & 0 \\
 & & && && && & \\
0\ar[r] & (\ker{S(\kappa_{ij})})(J) \ar[rr,"(k_{S(\kappa_{ij})})_J"] && S_i(J) \ar[rr,"{S(\kappa_{ij})}_J"] && S_j(J) \ar[rr,"(c_{S(\kappa_{ij})})_J"]\ar[rr] &&(\coker{S(\kappa_{ij})})(J)
  \ar[r] & 0
\et\nn
\]}

\begin{dfn}[\textcolor{blue}{\index{Category of! additive Ab-valued functors}{Category of additive Ab-valued functors $Ab^\A$}, \index{Exactness! in the category of additive Ab-valued functors}{Exactness in $Ab^\A$}, \index{Evaluation functor}{Evaluation functor}}]
Let $\A$ be an abelian category. The \ul{category of additive Ab-valued functors} $F:\A\ra Ab$ is denoted by $Ab^\A$.

A sequence $0\ra F\ra G\ra H\ra 0$ in $Ab^\A$ is \ul{exact} if $0\ra F(A)\ra G(A)\ra H(A)\ra 0$ is \ul{exact} in $Ab$ for all $A\in\Ob\A$ (i.e., the \ul{evaluation functor} $E_A:Ab^\A\ra Ab$, $(F\sr{\eta}{\ral}G)~\mapsto~ \big(F(A)\sr{\eta_A}{\ral}G(A)\big)$ is exact by definition).
\end{dfn}
$\ast$ It follows that if $\A$ is small, then the product ~$E:=\prod_{A\in\A}E_A:Ab^\A\ra Ab$,~ given by
{\scriptsize\[\bt
\big(F\ar[r,"\eta"] & G\big)~~\longmapsto~~\Big(\left(\prod_{A\in\A}E_A\right)(F)\ar[rr,"{\left(\prod_{A\in\A}E_A\right)(\eta)}"] && \left(\prod_{A\in\A}E_A\right)(G)\Big)~~:=~~\Big(\prod_{A\in\A}E_A(F)\ar[rr,"{\prod_{A\in\A}E(\eta)}"] && \prod_{A\in\A}E_A(G)\Big),
\et\nn
\]}
is an exact imbedding.

$\ast$ It follows from Proposition \ref{YonProp} that $Ab^\A$ is an abelian category.

\begin{dfn}[\textcolor{blue}{\index{Intersection of! subobjects}{Intersection of subobjects}}]
Let $\C$ be a category and $C\in\Ob\C$. The \ul{intersection} of a family of subobjects $\{u_i:C_i\hookrightarrow C\}_{i\in I}$ of $C$ is a subobject $u:\bigcap_i C_i\hookrightarrow C$ such that the following hold (\blue{footnote}\footnote{It is clear (from its definition) that the intersection $\bigcap_i C_i$ is the largest common subobject of the subobjects $C_i$.}):
\bit[leftmargin=0.7cm]
\item $u$ factors through each subobject as {\footnotesize $\bt[column sep=small] u=u_iv_i:\bigcap_i C_i\ar[r,"v_i"] & C_i\ar[r,hook,"u_i"] & C\et$}. (Note: As $u$ is monic, so is each $v_i$.)
\item For every subobject $u':C'\hookrightarrow C$ that factors through each subobject as $\bt[column sep=small] u'=u_iv'_i: C'\ar[r,"v'_i"] & C_i\ar[r,hook,"u_i"] & C\et$, there exists a unique morphism $h:C'\ra\bigcap_i C_i$ such that $\bt u'=uh:C'\ar[r,"h"] & \bigcap_i C_i\ar[r,hook,"u"] & C\et$
\eit
{\footnotesize
\bea\bt
C_i\ar[r,hook,"u_i"] & C \\
 & \bigcap_i C_i\ar[ul,hook,"v_i"']\ar[u,hook,"u"] \\
 & C'\ar[uul,bend left,"v_i'"]\ar[u,dashed,"h"]\ar[uu,hook,bend right=50,"u'"']
\et~~=~~
\bt
   & C & \\
C_i\ar[ur,hook,"u_i"] & & C_j\ar[ul,hook,"u_j"'] \\
  & \bigcap^C_iC_i\ar[ul,hook,"v_i"']\ar[ur,hook,"v_j"] & \\
  & X\ar[u,dashed,"x"']\ar[uul,bend left,"x_i"]\ar[uur,bend right,"x_j"'] &
\et~~~~=~~~~\txt{\shortstack[l]{The monic limit of the given\\ system of subobjects of $C$.}}\nn
\eea}
\end{dfn}

\begin{dfn}[\textcolor{blue}{\index{AB5 (Grothendieck) category}{AB5 (Grothendieck) category}: \cite[p.9]{grothdk1957}}]
A category $\C$ is an \ul{AB5 category} (or \ul{Grothendieck category}) if the following hold:
\bit[leftmargin=0.7cm]
\item[(i)] For any family $(C_i)_{i\in I}$ of subobjects of a given object $C$, the coproduct $\coprod_{i\in I} C_i$ exists.
\item[(ii)] For every object $A$, an increasing family of subobjects $A_i\subset A_j\subset A$ for $i\leq j$, $i,j\in I$ (where $I$ is a linearly ordered set), and any subobject $B\subset A$, we have (\blue{footnote}\footnote{By \cite[Propositions 1.4 and 1.6]{mitchell1964}, (\ref{AB5Eq1}) is equivalent to
~$f^{-1}\left(\bigcup A_i\right)=\bigcup f^{-1}(A_i)$ ~for any morphism $f:B\ra A$ and any family of subobjects $A_i\subset A$, as well as to the exactness of (directed) colimits. However, we will not be needing explicit knowledge of these equivalences.})
\bea
\label{AB5Eq1}\textstyle B\cap\left(\bigcup_{i\in I} A_i\right)=\bigcup_{i\in I}(B\cap A_i),
\eea
where for any collection of subobjects $C_i\subset C$ (with $I$ a set), in terms of their coproduct {\footnotesize $\coprod_{i\in I} C_i$}, their \index{Union! of subobjects}{\ul{union}} {\footnotesize $\bigcup_{i\in I}C_i$} is \ul{defined to be} the image of the induced morphism {\footnotesize $g:\coprod_{i\in I}C_i\ra C$} (as in the diagram below). (\blue{footnote}\footnote{Because $u_i:=e_gq_i$ (in the diagram) is monic due to the commutativity of the diagram, this definition of union is equivalent to the statement (from our earlier definition) that {\footnotesize $\bigcup_{i\in I}C_i$} is the smallest subobject of $C$ containing each of the subobjects $C_i$.})
\[\adjustbox{scale=0.7}{\bt
C_i\ar[dddr,bend right=50,hook,"v_i"']\ar[dr,hook,"q_i"] &  & \cdots\ar[dl,hook,"q_j"]\\
     & \coprod C_i\ar[dd,bend right=40,dashed,"g"']\ar[d,two heads,"e_g"] & \\
     & \bigcup C_i\ar[d,hook,"m_g"] &\\
     & C \\
     & \cdots\ar[from=u,"x"] \\
\et}~~=~~
\adjustbox{scale=0.8}{\bt
 & C_i\cap^C C_j\ar[dl,hook,"\al_i"']\ar[dr,hook,"\al_j"]& \\
C_i\ar[dr,hook,"u_i"',"e_gq_i"]\ar[dddr,bend right,"xv_i"'] &&  C_j\ar[dl,hook,"u_j","e_gq_j"']\ar[dddl,bend left,"xv_j"]\\
 & \bigcup^CC_i\ar[dd,dashed,"h"',"xm_g"]& \\
  & & \\
 & X &
\et}~~~~=~~~~
\txt{\shortstack[l]{The monic colimit of the given\\ system of subobjects of $C$.}}
\]
\eit
\end{dfn}

\begin{rmk}[\textcolor{blue}{An aside: \cite[Proposition 1.6]{mitchell1964}}]\label{AB5Equiv1}
For an abelian category, the AB5 condition is equivalent to asserting that the (directed) colimit $\varinjlim$ preserves monic morphisms (and is thus an exact functor). We will not be needing explicit knowledge of this result. Therefore we will not provide the proof here.
\end{rmk}

\begin{prp}[\textcolor{blue}{\cite[Proposition 5.21, p.111]{freyd1964}}]\label{AB5ForLAb}
If $\A$ is an abelian category, then $Ab^\A$ is an abelian AB5 category.
\end{prp}
\begin{proof}
We already know (from Proposition \ref{YonProp}) that $Ab^\A$ is abelian, and so we will prove the AB5 property. Let $F:\A\ra Ab$ be a functor, $I$ a set, and $H,F_i\subset F$ subfunctors (in the sense $H(A),F_i(A)\subset F(A)$ for all $A\in\Ob\A$ along with naturality as in the following diagram).
\bea
\adjustbox{scale=0.9}{\bt A\ar[d,"f"] \\ B\et}
~~\sr{\subset}{\longmapsto}~~
\adjustbox{scale=0.9}{\bt
  H(A),F_i(A)\ar[d,"H(f)"',"F_i(f)"] \ar[r,hook,"\subset_A"] & F(A)\ar[d,"F(f)"] \\
  H(B),F_i(B)\ar[r,hook,"\subset_B"] & F(B)
\et}\nn
\eea
Then with $\bigcup_{i\in I}F_i$ (\blue{footnote}\footnote{To stay in Ab, we need to assume (i) $I$ is a linearly ordered set and (ii) $F_i\subset F_j\subset F$ for $i\leq j$, so that the resulting increasing union of abelian groups is an abelian group.}) given pointwise by $\left(\bigcup_{i\in I}F_i\right)(A):=\bigcup_{i\in I}F_i(A)$, for $A\in\Ob\A$, we have
{\footnotesize\begin{align}
&\textstyle \left(H\cap\left(\bigcup\limits_{i\in I} F_i\right)\right)(A):=H(A)\cap\left(\bigcup\limits_{i\in I} F_i(A)\right)\sr{(s)}{=}\bigcup\limits_{i\in I}(H(A)\cap F_i(A))=:\left(\bigcup\limits_{i\in I}(H\cap F_i)\right)(A),~~~~\txt{for all}~~A\in\Ob\A,\nn\\
&\textstyle~~\Ra~~H\cap\left(\bigcup_{i\in I} F_i\right)=\bigcup_{i\in I}(H\cap F_i),\nn
\end{align}}
where step (s) follows from the familiar fact that $Ab\subset Sets$ is an AB5 category.
\end{proof}

\section{The Yoneda Imbedding (A Left-Exact Full Imbedding)}
In the following discussion, the symbol ``$H$'' originally came from ``Hom'', which we now call ``Mor''.
\begin{lmm}[\textcolor{blue}{
\index{Yoneda's lemma}{Yoneda's lemma},
\index{Representation! functor}{Representation functor}}]\label{YndLmm}
Let $\C$ be a category, $H:\C\ra Sets^\C$ the cofunctor (called \ul{representation functor}) given by
{\footnotesize\bea
H:\C\ra Sets^\C,~~A\sr{f}{\ral}B~~\mapsto~~H^A\sr{H^f}{\lal}H^B~~:=~~\Mor_\C(A,-)\sr{f^\ast:=\circ f}{\lal}\Mor_\C(B,-),\nn
\eea}
and ~$N,E:\C\times Sets^\C\ra Sets$ ~(\blue{footnote}\footnote{Here, the symbols $N,E$ signify the facts that (i) $N$ maps to ``natural transformations'' and (ii) $E$ maps to ``evaluations''.})~ the bifunctors given by
{\footnotesize\bea
&&N:~(A,F)\sr{(f,\eta)}{\ral}(B,G)~\eqv~\left(A\sr{f}{\ral}B,F\sr{\eta}{\ral}G\right)~~\mapsto~~\Mor_{Sets^\C}\big(H^A,F\big)\sr{f^{\ast\ast}\circ \eta_\ast}{\ral}\Mor_{Sets^\C}(H^B,G),\nn\\
&&E:~(A,F)\sr{(f,\eta)}{\ral}(B,G)~\eqv~\left(A\sr{f}{\ral}B,F\sr{\eta}{\ral}G\right)~~\mapsto~~F(A)\sr{\eta(f)}{\ral}G(B),\nn
\eea}
where ~{\small $N(f,\eta):=f^{\ast\ast}\circ\eta_\ast$}~ is given by the diagrams
{\footnotesize\bea
\bt H^A\ar[r,"\al"] & F\et~\mapsto~
\bt H^A\ar[r,"\al"] & F\ar[d,"\eta"] \\
    H^B\ar[u,"f^\ast"] & G \et~=~
\bt H^B\ar[r,"\beta"] & G\et,~~~~
\beta:=\eta\circ\al\circ f^\ast=f^{\ast\ast}\eta_\ast(\al),\nn
\eea}
and ~{\small $E(f,\eta):=\eta(f)$} ~is given by the commutative $\eta$-naturality diagram
{\footnotesize
\bea
\bt A\ar[d,"f"]\\ B \et~~\mapsto~~
\bt
F(A)\ar[drr,dashed,"{\eta(f)}"]\ar[d,"{F(f)}"]\ar[rr,"\eta_A"] && G(A)\ar[d,"{G(f)}"]\\
F(B)\ar[rr,"\eta_B"] && G(B)
\et,~~~~\txt{i.e.,}~~~~\eta(f):=G(f)\circ\eta_A=\eta_B\circ F(f).\nn
\eea}Then we have the following isomorphism $\theta_{(A,F)}$ of sets that is natural in both $A$ and $F$.
\bea
\label{YndLmmEq0}\Mor_{Sets^\C}\big(H^A,F\big)\sr{\theta_{(A,F)}}{\cong} F(A).
\eea
In other words, there exists a \ul{natural isomorphism} ~$\theta:N\ra E$, as in the $\theta$-naturality diagram:
{\footnotesize
\bea
\bt[column sep=tiny]
   A\ar[d,"f"] & F\ar[d,"\eta"]\\
   B & G
\et~~\mapsto~~
\bt
   N(A,F)\ar[d,"{N(f,\eta)}"]\ar[rr,"\theta_{(A,F)}","\cong"'] && E(A,F)\ar[d,"{E(f,\eta)}"]\\
   N(B,G)\ar[rr,"\theta_{(B,G)}","\cong"'] && E(B,G)
\et~~=~~
\bt
\Mor_{Sets^\C}\big(H^A,F\big)\ar[d,"{f^{\ast\ast}\circ\eta_\ast}"]\ar[rr,"\theta_{(A,F)}","\cong"'] && F(A)\ar[d,"{\eta(f):=}"',"{\eta_A\circ F(f)=G(f)\circ\eta_A}"]\\
\Mor_{Sets^\C}(H^B,G)\ar[rr,"\theta_{(B,G)}","\cong"'] && G(B)
\et\nn
\eea}

\end{lmm}
\begin{proof}
Given $\al\in \Mor_{Sets^\A}(H^A,F)$, we define $\theta_{(A,F)}(\al)$ based on the following $\al$-naturality diagram:
{\footnotesize
\bea
\label{YndLmmEq1N}
\bt A\ar[d,"f"]\\ B \et~~\mapsto~~
\bt
H^A(A)=\Mor_\C(A,A)\ar[d,"{f_\ast}"]\ar[rr,"\al_A"] && F(A)\ar[d,"{F(f)}"]\\
H^A(B)=\Mor_\C(A,B)\ar[rr,"\al_B"] && F(B)
\et
\eea}
$\bullet$ So, define ~$\theta=\theta_{(A,F)}:\Mor_{Sets^\C}\big(H^A,F\big)\ra F(A)$ ~by
\bea
\label{YndLmmEq1}\theta(\al):=\al_A(1_A)\in F(A),~~\txt{for}~~\al\in \Mor_{Sets^\C}\big(H^A,F\big).
\eea
Also define a map in the reverse direction ~$\theta'=\theta'_{(A,F)}:F(A)\ra \Mor_{Sets^\C}\big(H^A,F\big)$ ~as follows: For each $x\in F(A)$ we have the naturality diagram for ~$\theta'(x):H^A\ra F$, ~i.e.,
{\footnotesize
\bea
\label{YndLmmEq2N}
\bt A\ar[d,"f"]\\ B \et~~\mapsto~~
\bt
H^A(A)=\Mor_\C(A,A)\ar[d,"{f_\ast}"]\ar[rr,"\theta'(x)_A"] && F(A)\ar[d,"{F(f)}"]\\
H^A(B)=\Mor_\C(A,B)\ar[rr,"\theta'(x)_B"] && F(B)
\et
\eea}
$\bullet$ So, define ~$\theta'=\theta'_{(A,F)}:F(A)\ra \Mor_{Sets^\C}\big(H^A,F\big)$ ~by
\bea
\label{YndLmmEq2}\theta'(x)_B(f):=F(f)(x)\in F(B),~~~~\txt{for}~~x\in F(A),~~~~f\in \Mor_\C(A,B).
\eea
With the above definitions, (\ref{YndLmmEq1}) and (\ref{YndLmmEq2}), we get the identity-compositions
\bea
&&\theta\circ\theta'=1_{F(A)}:F(A)\sr{\theta'}{\ral} \Mor_{Sets^\C}\big(H^A,F\big)\sr{\theta}{\ral}F(A)~~~~\txt{and}~~~~\nn\\
&&\theta'\circ\theta=1_{\Mor_{Sets^\C}\big(H^A,F\big)}:\Mor_{Sets^\C}\big(H^A,F\big)\sr{\theta}{\ral}F(A)\sr{\theta'}{\ral}\Mor_{Sets^\C}\big(H^A,F\big),\nn
\eea
because with $x\in F(A)$ and $f:A\ra B$, we get
\begin{align}
&\left(\theta\circ\theta'\right)(x)_A=\left(\theta\theta'(x)\right)_A\sr{(\ref{YndLmmEq1})}{=}\theta'(x)_A(1_A)\sr{(\ref{YndLmmEq2})}{=}F(1_A)(x)=1_{F(A)}(x)=x~~~~\txt{and}~~~~\nn\\
&\big((\theta'\circ\theta)(\al)\big)_B(f)=\big(\theta'(\theta(\al))\big)_B(f)\sr{(\ref{YndLmmEq1})}{=}\big(\theta'(\al_A(1_A))\big)_B(f)
\sr{(\ref{YndLmmEq2})}{=}F(f)(\al_A(1_A))=\big(F(f)\circ\al_A\big)(1_A)\nn\\
&~~~~\sr{(\ref{YndLmmEq1N})}{=}(\al_B\circ f_\ast)(1_A)=\al_B(f_\ast(1_A))=\al_B(f\circ 1_A)=\al_B(f).\nn
\end{align}

Finally, we can verify the naturality of $\theta=\theta_{(A,F)}$ as in the following diagram:
{\footnotesize
\[
\bt[column sep=tiny]
A\ar[ddd,"f"] & F\ar[ddd,"\eta"] \\
&\\ &\\
B & G
\et~~~~\mapsto~~~~
\bt
H^A\sr{\al}{\ral} F\ar[dd,mapsto,"f^{\ast\ast}\circ(\eta)_\ast"]\ar[rr,mapsto,"\theta_{(A,F)}"] && \al_A(1_A)\ar[d,mapsto,"\eta(f):=G(f)\circ\eta_A=\eta_B\circ F(f)"]  \\
& & \big(\eta_B\circ F(f)\circ\al_A\big)(1_A)\ar[d,equal,"{\txt{naturality of $\al$ from (\ref{YndLmmEq1N})}}"]\\
H^B\sr{\eta\circ\al\circ f^\ast}{\ral}G\ar[rr,mapsto,"\theta_{(B,F)}"] && (\eta\circ\al\circ f^\ast)_B(1_B)=(\eta_B\circ\al_B\circ f_\ast)(1_A)
\et~~~~ \qedhere
\]}

\end{proof}

Note that the hypotheses in the following lemma are the same as those of Lemma \ref{YndLmm}, except that $\C$ is replaced with $\I$ (i.e., the hypotheses of the two lemmas differ only in symbols).

\begin{lmm}[\textcolor{blue}{\index{Yoneda imbedding}{Yoneda imbedding}}]\label{YndImbb}
Let $\I$ be a category $\I$. Then we have a contravariant full imbedding $\I\hookrightarrow Sets^\I$ given by the representation functor
{\small\bea
H:\I\ra Sets^\I,~~\big(I\sr{f}{\ral}J\big)\mapsto \Big(\Mor_{\I}(I,-)\sr{\Mor_{\I}(f,-)}{\lal}\Mor_{\I}(J,-)\Big),\nn
\eea}where the natural transformation ~$H^f:=\Mor_\I(f,-)=\circ f=f^\ast$~ (i.e., right or post composition with $f$).

Moreover, if $\I$ is abelian, then (i) the image $\im H\subset Sets^\I\cap Ab$ and (ii) $H$ is left exact (\blue{footnote}\footnote{That is, if {\scriptsize $0\ra I\sr{f}{\ral}I'\sr{g}{\ral}I''\ra0$} is exact in $\I$, then so is {\scriptsize $0\ra \Mor_\I(I'',-)\sr{\circ g}{\ral}\Mor_\I(I',-)\sr{\circ f}{\ral}\Mor_\A(I,-)$} in $Sets^\I$,~ i.e., so is {\scriptsize $0\ra \Mor_\I(I'',J)\sr{\circ g}{\ral}\Mor_\I(I',J)\sr{\circ f}{\ral}\Mor_\A(I,J)$} in $Sets$ for every object $J\in\Ob\I$.}).
\end{lmm}
\begin{proof}
By Yoneda's lemma (Lemma \ref{YndLmm}) we indeed have a contravariant \ul{mapwise-injective} functor
{\small\bea
H:\I\ra Sets^{\I},~~\big(I\sr{f}{\ral}J\big)\mapsto \Big(H^I\sr{H^f}{\lal}H^J\Big):=\left(\Mor_{\I}(I,-)\sr{f^\ast=\circ f}{\lal}\Mor_{\I}(J,-)\right),\nn
\eea}where the natural transformation ~$H^f:=f^\ast=\circ f$~ (right composition with $f$) is given, for each anticipated morphism $\Big(J\sr{h}{\ral}-\Big)\in \Mor_{\I}(J,-)=H^J=$, by
{\small\bea
f^\ast\Big(J\sr{h}{\ral}-\Big):=\Big(I\sr{f}{\ral}J\sr{h}{\ral}-\Big)=\Big(I\sr{h\circ f}{\ral}-\Big)\in \Mor_{\I}(I,-)=H^I.\nn
\eea}As usual, {\small$H^f=\Big\{(H^I)_K\sr{(H^f)_K}{\lal}(H^J)_K~\big|~K\in\Ob\I\Big\}=\Big\{\Mor_\I(K,I)\sr{f^\ast}{\lal}\Mor_\I(K,J)~\big|~K\in\Ob\I\Big\}$}, and for each morphism $g\in \Mor_\I(K,K')$, we have the naturality commutative diagram
{\small\[
\adjustbox{scale=0.8}{\bt K\ar[d,"g"]\\ K'\et}~~\sr{H^f}{\longmapsto}~~
\adjustbox{scale=0.8}{\bt
{\Mor_\I(I,-)_K}\ar[d,"{\Mor_\I(I,-)^g}"]\ar[from=rr,"(H^f)_K"']&&{\Mor_\I(J,-)_K}\ar[d,"{\Mor_\I(J,-)^g}"]\\
{\Mor_\I(I,-)^{K'}}\ar[from=rr,"(H^f)^{K'}"']&&{\Mor_\I(J,-)^{K'}}
\et}
~~:=~~
\adjustbox{scale=0.8}{\bt
{\Mor_\I(I,K)}\ar[d,"{g_\ast}"]\ar[from=rr,"f^\ast"']&&{\Mor_\I(J,K)}\ar[d,"g_\ast"]\\
{\Mor_\I(I,K')}\ar[from=rr,"f^\ast"']&&{\Mor_\I(J,K)}
\et}~~~~g_\ast\circ f^\ast=f^\ast\circ g_\ast.\nn
\]}
\ul{Naturality} of $H$ is clear because left composition (with $g$) commutes with right composition (with $f$).

Also, the functor $H$ is \ul{full} by Yoneda's lemma, because if we set $F:=H^B$ in the isomorphism from (\ref{YndLmmEq0}), i.e., $\Mor_{Sets^\I}\big(H^A,F\big)\cong F(A)$, we see that on morphisms, $H$ is both injective and surjective, i.e., we have an isomorphism of sets:
\bea
H=\theta^{-1}:\Mor_\A(A,B)\ra \Mor_{Sets^\I}(H^A,H^B).\nn
\eea
Finally, if $\I$ is abelian, then (i) the containment $\im H\subset Ab$ is clear since $H$ is a full imbedding and (ii) \ul{left exactness} of $H$ follows directly from \ul{left exactness of the mor functor} (and $H^{-1}:\im H\ra\I$ is also readily left exact by the \ul{partial exactness criteria}).
\end{proof}

\begin{crl}\label{AbelImbb}
Every additive category fully imbeds (with left exactness preserved) into an abelian category.
\end{crl}
\begin{proof}
Let $\I$ be an additive category. Then for every $I,J\in\Ob\I$, we see that {\small $\Mor_\I(I,J)\cong \Mor_{Sets^\I}\big(H^I,H^J\big)$} is an abelian group (and so we can replace Sets with Ab). Thus, the Yoneda imbedding maps $\I$ into $Ab^\I$, i.e., $H:\I\ra Ab^\I$. Hence, because $Ab^\I$ is abelian by Proposition \ref{YonProp}, $\I$ is fully imbedded (with left exactness preserved) into an abelian category.
\end{proof}

\begin{question}[\textcolor{blue}{See \cite{chih2014} for example}]
Let $\A$ be an additive category. Suppose a left exact full imbedding ~$F:\A\ra\txt{$R$-mod}$~ exists (so, in particular, $F$ maps monomorphisms to monomorphisms). If objects $A,B\in\Ob\A$ are such that $A\subset B$ and $B\subset A$, does it follow that $A\cong B$?
\end{question}

\begin{rmk}[\textcolor{blue}{Improving and using the Yoneda imbedding}]
Even though the representation functor is only left exact (due to the left exactness of Mor), exactness of the imbedding can be improved if we start with the category $\I$ both \ul{small} and \ul{abelian}. Such an improvement will be given in Lemma \ref{RpFnExLmm} (which is also based on Yoneda's lemma). This will then be used to obtain the main imbedding result called the Freyd-Mitchell imbedding (Theorem \ref{FreydMitchIT}).
\end{rmk}

\section{Preliminary Results on Exact Imbedding of Abelian Categories}
\subsection{Relevant properties of abelian categories: Generators and imbedding into $R$-mod}
\begin{dfn}[\textcolor{blue}{\index{Representable! functor}{Representable functor}}]
Let $\C$ be a category. A functor $F:\C\ra Sets$ is \ul{representable} (in the category of functors $Sets^\C$) if $F\cong \Mor_{Sets}(C,-):\C\ra Sets$ for some object $C\in\Ob\C$.
\end{dfn}

\begin{dfn}[\textcolor{blue}{
\index{Generator for a category}{Generator},
\index{Cogenerator for a category}{Cogenerator},
\index{Finitely generated! object}{Finitely generated object},
\index{Small! object}{Small object}}]
Let $\A$ be an abelian category and $A,A_i\in\Ob\A$, $i\in I$. The object $A$ is a \ul{generator for $\A$} (resp. \ul{cogenerator for $\A$}) if $\Mor_\A(A,-):\A\ra Ab$ (resp. $\Mor_\A(-,A):\A\ra Ab$) is an imbedding. The object $A$ is \ul{small} if $\Mor_\A(A,-):\A\ra Ab$ preserves coproducts, i.e., $\Mor_\A\left(A,\bigoplus A_i\right)\cong \bigoplus \Mor_\A(A,A_i)$. (\blue{footnote}\footnote{It is clear that in $R$-mod, $R={}_RR$ is a small projective generator.}).

Given a generator $G\in\Ob\A$, an object $A\in\Ob\A$ is \ul{finitely generated} wrt $G$ if $A$ is a quotient object of an object of the form $G^{\oplus I}:=\bigoplus_{i\in I}G$ for a finite set $I$ (i.e., there exist a finite set $I$ and an epimorphism $e:G^{\oplus I}\ra A$).
\end{dfn}

\begin{dfn}[\textcolor{blue}{\index{Exact! category}{Exact category}}]
A category $\C$ is an \ul{exact category} if every morphism $f:C\ra D$ in $\C$ factors as $f=me:C\sr{e}{\twoheadrightarrow}E\sr{m}{\hookrightarrow}D$, for an epimorphism $e$ and a monomorphism $m$. (\blue{footnote}\footnote{By the isomorphism theorem for abelian categories, every abelian category is an exact category.})
\end{dfn}

\begin{dfn}[\textcolor{blue}{\index{Preservation of a property by a functor}{Preservation} and \index{Reflection of a property by a functor}{Reflection} of Properties by a Functor}]
Let $\C,\D$ be categories, $F:\C\ra\D$ a functor, and $P$ a property of diagrams/systems (e.g., commutativity, morphism type, etc). Then $F$ \ul{preserves} $P$ if for any diagram/system $S\subset\C$, ``if $S$ has property $P$ then so does $F(S)$''. Similarly, $F$ \ul{reflects} $P$ if for any diagram/system $S\subset\C$,`` if $F(S)$ has property $P$ then so does $S$''.
\end{dfn}

\begin{lmm}[\textcolor{blue}{An exact (co)functor that reflects zero objects is an imbedding: \cite[Proposition 7.2]{mitchell1965}}]\label{MitRefLmm}
Let $\A$ be an abelian category, $\E$ an exact additive category, and $T:\A\ra\E$ a functor (resp. $T':\A\ra\E$ a cofunctor). The following are true:
\bit[leftmargin=0.7cm]
\item If $T$ is exact and reflects zero-objects (i.e., $T(A)=0$ implies $A=0$), then $T$ is faithful.
\item If $T'$ is exact and reflects zero-objects (i.e., $T'(A)=0$ implies $A=0$), then $T'$ is faithful.
\eit
\end{lmm}
\begin{proof}
It is enough to consider the functor case only (since the cofunctor case is similar). If $f:A\ra B$ is a morphism, let $f=m_f\circ e_f:A\sr{e_f}{\ral}\im f\sr{m_f}{\ral}B$. Then $T(f)=T(m_f)T(e_f)=0$ implies $T(m_f)=0$ ``or'' $T(e_f)=0$ (since $T$ maps monics to monics ``and'' epics to epics). If $T(e_f)=0$, then $h\circ T(e_f)=0$ implies $h=0$ (whenever the composition exists). In particular, $h:=id_{T(\im f)}=0$, i.e., $T(\im f)=0$, and so $\im f=0$ (since $T$ reflects zero objects). If $T(m_f)=0$, a similar argument again shows $\im f=0$. This in turn implies $m_f=0$, and so $f=0$, i.e., $T$ is injective on morphisms.
\end{proof}
\begin{rmk}[\blue{footnote}\footnote{If $T$ is injective on objects then $T$ reflects zero-objects, but the reverse might not hold in general.}]
In the above proof, suppose $T(A)=T(B)$. Then we have the identity morphism $T(A)\sr{id_{T(A)}=id_{T(B)}}{\ral}T(B)$. This means the morphisms $A\sr{T^{-1}(id_{T(A)})}{\ral}B\sr{T^{-1}(id_{T(B)})}{\ral}A$ and $A\sr{id_A}{\ral}A$ have the same value under $T$, and so $id_A=T^{-1}(id_{T(B)})\circ T^{-1}(id_{T(A)})$, since $T$ is injective on morphisms. It follows by symmetry that $A\cong B$.
\end{rmk}

\begin{question}
In Lemma \ref{MitRefLmm}, if we replace the phrase ``$T$ is exact'' with the phrase ``$T$ maps monics to monics or maps epics to epics'', does this affect the conclusion?
\end{question}

\begin{crl}[\textcolor{blue}{\cite[Proposition 15.3]{mitchell1965}}]\label{MitRefCrl}
Let $\A$ be an abelian category. If $P\in\Ob\A$ is a projective object (resp. $E\in\Ob\A$ is an injective object), then the functor $H^P:=\Mor_\A(P,-):\A\ra Ab$ (resp. $H_E:=\Mor_\A(-,E):\A\ra Ab$) is an imbedding $\iff$ injective on objects (by Lemma \ref{MitRefLmm}) $\iff$ $P$ is a projective generator in $\A$ (resp. $E$ is an injective cogenerator in $\A$).
\end{crl}
\begin{proof}
If a functor $F:\A\ra\E$ is injective on objects, then in particular, it reflects zero objects because if $A\in\Ob\A$, then $F(A)=0=F(0)$ implies $A=0$.
\end{proof}

\begin{lmm}[\textcolor{blue}{\cite[Lemma 2.6]{mitchell1964}}]\label{MitThmLmm1}
Let $\A$ be an abelian category (assumed here to be locally small) and $G\in\Ob\A$. Then $G$ is a generator $\iff$ for each $A\in\Ob\A$, the ``diagonal'' element (\blue{footnote}\footnote{Recall that for a set $S$, ~$\prod_{s\in S} S:=S^S:=\{\txt{maps}~f:S\ra S\}\ni id_S$.}) ~{\small $e_A=(e_h)_{h\in\Mor_\A(G,A)}:=(h)_{h\in\Mor_\A(G,A)}:G^{\oplus \Mor_\A(G,A)}\ra A$}, ~in the set
{\small\bea
\textstyle \Mor_\A\left(G^{\oplus \Mor_\A(G,A)},A\right)=\Mor_\A\left(\bigoplus\limits_{h\in \Mor_\A(G,A)}G~,~A\right)\cong \prod\limits_{h\in \Mor_\A(G,A)}\Mor_\A(G,A),\nn
\eea}is an epimorphism.
\end{lmm}
\begin{proof}
For notational convenience, let $[G,A]:=\Mor_A(G,A)$.
{\flushleft $(\Ra)$}: Suppose $G$ is a generator, i.e., the following associated functor is an imbedding:
\bea
\Mor_\A(G,-):~\A\ra Ab,~~C\sr{f}{\ral} D~\mapsto~~\Mor_\A(G,C)\sr{f\circ }{\ral}\Mor_\A(G,D).\nn
\eea
Given $A\in\Ob\A$, let $f\in \Mor_\A(A,B)$. Then {\small $f\circ e_A=(f\circ h)_{h\in[G,A]}:G^{\oplus \Mor_\A(G,A)}\sr{e_A}{\ral} A\sr{f}{\ral}B$}, and so
\bea
f\circ e_A=(f\circ h)_{h\in[G,A]}=0~~\Ra~~f\circ h=0~~\txt{for all}~~h\in[G,A],~~\Ra~~f\circ =0,~~\sr{\txt{faithfulness}}{\Ra}~~f=0,\nn
\eea
which shows $e_A$ is an epimorphism.

{\flushleft $(\La)$}: Suppose that for each $A\in\Ob\A$, the diagonal element ~$e_A=(e_h)_{h\in[G,A]}:=(h)_{h\in[G,A]}:G^{\oplus \Mor_\A(G,A)}\ra A$ ~in ~$\Mor_\A\left(G^{\oplus \Mor_\A(G,A)},A\right)$ ~is an epimorphism. Consider the functor
\bea
\Mor_\A(G,-):~\A\ra Ab,~~C\sr{f}{\ral} D~\mapsto~~\Mor_\A(G,C)\sr{f\circ }{\ral}\Mor_\A(G,D).\nn
\eea
Let $f\in \Mor_\A(A,B)$. If $\Mor_\A(G,f):=f\circ =0$ (i.e., $f\circ h=0$ for each $h\in[G,A]$) then $f\circ e_A=0$, and so $f=0$ (since $e_A$ is an epimorphism). Therefore $\Mor_\A(G,-)$ is faithful (i.e., injective on morphisms).

If $\Mor_\A(G,A)=0$, then the zero epimorphism ~{\small $e_A=(h)_{h\in[G,A]}=0=e_Aid_A:G^{\oplus \Mor_\A(G,A)}\sr{e_A}{\ral} A\sr{id_A}{\ral}A$} ~implies $id_A=0$, i.e., $A=0$. Hence $\Mor_\A(G,-)$ is injective on objects as well.
\end{proof}

\begin{lmm}[\textcolor{blue}{\cite[Theorem 3.1]{mitchell1964}}]\label{MitThmLmm2}
Let $\A$ be an abelian category. Then $\A$ is equivalent to some $R$-mod $\iff$ $\A$ has a small projective generator $G\in\Ob\A$.

Moreover, if (i) $\A$ has a projective generator $G\in\Ob\A$ and (ii) each object $A\in\Ob\A$ is finitely generated wrt $G$, then $\A$ is equivalent to some $R$-mod (and hence has a small projective generator).
\end{lmm}
\begin{proof}
($\Ra$) Assume $\A$ is equivalent to some $R$-mod through a functor $E:R\txt{-mod}\ra\A$ (recall that as an equivalence of categories, $E$ is full, faithful, and dense). In particular, (1) there exists an object $G\in\Ob\A$ such that $E(R)\cong G$ and (2) for any $R$-modules $M,N$, the map $E:\Mor_R(M,N)\ra\Mor_\A(E(M),E(N))$ is a bijection. We can check that $\Mor_\A(G,-)$ is (i) an imbedding (ii) maps epics to epics, and (iii) preserves coproducts (since $E$ does), and so $G$ is a small projective generator for $\A$.

($\La$) Conversely, assume $\A$ has a small projective generator $G\in\Ob\A$ (or alternatively, $\A$ has a projective generator $G\in\Ob\A$ and each object $A\in\Ob\A$ is finitely generated wrt $G$). We need to construct an equivalence of categories (i.e., a full, faithful, and dense functor) $E:\A\ra R$-mod for some ring $R$. Let $R$ denote the \index{Endomorphism! ring}{\ul{endomorphism ring}} of $G$, defined to be the abelian group $R:=\Mor_\A(G,G)$ as a ring with map composition as multiplication and $1_R:=id_G$ as the identity element. Consider the functor (from $\A$ to the right $R$-modules category, $\txt{mod-}R$):
\bea
E:=\Mor_\A(G,-):\A\ra \txt{mod-}R,~~A\sr{f}{\ral}B~\mapsto~\Mor_\A(G,A)\sr{f_\ast=f\circ}{\ral}\Mor_\A(G,B),\nn
\eea
which is an imbedding because $G$ is a generator, and exact because $G$ is projective, i.e., $E$ is an exact imbedding (hence \ul{faithful}). It remains to show $E$ is \ul{full} and \ul{dense}.

To show $E$ is \ul{full}, we need to show $E$ is surjective on morphisms, i.e., any morphism $M_R\sr{f_R}{\ral}N_R$ in mod-$R$ of the form \bt M_R=\Mor_\A(G,A)\ar[r,"{f_R}"] & N_R=\Mor_\A(G,B)\et actually has the form
\bea
\bt \Mor_\A(G,A)\ar[rr,"{f_\ast=f\circ}"] && \Mor_\A(G,B),~~~~\txt{for some}~~f\in \Mor_\A(A,B).\et\nn
\eea
We consider two cases as follows:
{\flushright \ul{Case 1 (Suppose $A=G$)}:} Since $f_R$ is $R$-linear, for $r\in R=\Mor_\A(G,G)$, we have $f_R(r)=f_R(1_Rr)=f_R(1_R)\circ r$, where $1_R:=id_G$. That is, every morphism $R:=\Mor_\A(G,G)\sr{f_R}{\ral}\Mor_\A(G,B)$ is of the form $f_R(r)=f_\ast(r)=f\circ r$, for some $f:=f_R(1_R)\in \Mor_\A(G,B)$.

{\flushright \ul{Case 2 (General $A$)}:} If $G$ is small, then by Lemma \ref{MitThmLmm1}, for some set $I$ (or alternatively, if every object of $\A$ is finitely generated wrt $G$, then for some finite set $I$), we have a SES (\blue{footnote}\footnote{It is clear that the proof of this step goes through unchanged if we assume $A$ is finitely generated, so that $I$ is a finite set.})
\bea
0\ra K\sr{k_e}{\ral}G^{\oplus I}\sr{e}{\ral} A\ra 0\nn
\eea
For $i\in I$, denote the $i$th inclusion by $\bt[column sep=small] G\ar[r,hook,"u_i"] & G^{\oplus I}\et$ and consider  the associated composition:
\bea
\label{MitIndeq1}\bt f_{R,i}:=f_Re_\ast(u_i)_\ast:~\overbrace{\Mor_\A(G,G)}^R\ar[r,"{(u_i)_\ast}"] & \Mor_\A(G,G^{\oplus I})\ar[r,"e_\ast"] & \Mor_\A(G,A)\ar[r,"f_R"] & \Mor_\A(G,B). \et
\eea
Then, with $f_i:=f_{R,i}(1_R)\in \Mor_\A(G,B)$, we have $f_{R,i}(r)=f_{R,i}(1_R)\circ r=f_i\circ r$ ($r\in R$). Putting the $f_i$'s together, we get the following morphism (a branch of the given SES) in $\A$:
\bea
\label{MitIndeq2}\textstyle\ol{f}:=\bigoplus\limits_{i\in I}f_i\circ p_i:G^{\oplus I}\ra B,~~~~\txt{where $p_i$ is the $i$th canonical projection}~~~~p_i:G^{\oplus I}\twoheadrightarrow G.
\eea
Consider the associated composition:
\bea
\label{MitIndeq3}\bt \ol{f}_R:=f_Re_\ast:~\overbrace{\Mor_\A(G,G^{\oplus I})}^{\sr{\substack{\txt{$G$ small or}\\\txt{$I$ finite~~}}}{\cong} R^{\oplus I}}\ar[r,"e_\ast"] & \Mor_\A(G,A)\ar[r,"f_R"] & \Mor_\A(G,B). \et
\eea
For $i\in I$, denote the $i$th projection by {\small $\bt[column sep=small] G^{\oplus I}\ar[r,two heads,"p_i"] & G\et$ ~( $\sr{\txt{$G$ projective}}{\Longrightarrow}$ $\bt[column sep=small] \Mor_\A(G,G^{\oplus I})\sr{\substack{\txt{$G$ small or}\\\txt{$I$ finite~~}}}{\cong} R^{\oplus I}\ar[r,two heads,"(p_i)_\ast"] & R\et$ )}. ~Then
\[
\textstyle id_{G^{\oplus I}}=\bigoplus\limits_{i\in I} u_i\circ p_i~~\Ra~~id_{R^{\oplus I}}\cong\left(id_{G^{\oplus I}}\right)_\ast\sr{\substack{\txt{$G$ small or}\\\txt{$I$ finite~~}}}{=}\bigoplus\limits_{i\in I} (u_i)_\ast(p_i)_\ast\cong \bigoplus\limits_{i\in I} id_R.\nn
\]
Therefore, with $h=(r_i)_{i\in I}\in R^{\oplus I}$, we have
\bea
\label{MitIndeq4}&&\textstyle \ol{f}_R(h):=f_Re_\ast(h)=f_Re_\ast (id_{G^{\oplus I}})_\ast(h)=f_Re_\ast\left(\bigoplus\limits_{i\in I}(u_i)_\ast(p_i)_\ast\right)(h)=\bigoplus\limits_{i\in I}f_Re_\ast(u_i)_\ast(p_i)_\ast(h)\nn\\
&&\textstyle~~~~=\bigoplus\limits_{i\in I}f_{R,i}(p_i)_\ast(h)=\bigoplus\limits_{i\in I}f_{R,i}\big(p_i\circ h\big)=\bigoplus\limits_{i\in I}f_{R,i}(1_R)\circ(p_i\circ h)=\bigoplus\limits_{i\in I}f_i\circ p_i\circ h\nn\\
&&~~~~=\ol{f}\circ h.
\eea

Now, for any $\al\in \Mor_\A(G,K)$, the composition $\ol{f}\circ k_e:K\sr{k_e}{\ral}G^{\oplus I}\sr{\ol{f}}{\ral}B$ satisfies
\bea
&&\ol{f}\circ k_e\circ\al=\ol{f}_R(k_e\circ\al)=(f_Re_\ast)(k_e\circ\al)=f_R(e\circ k_e\circ\al)=f_R(0)=0,~~\txt{i.e.,}~~\nn\\
&&0=(\ol{f}\circ k_e)_\ast(\al)=\ol{f}\circ k_e\circ\al:G\sr{\al}{\ral}K\sr{k_e}{\ral}G^{\oplus I}\sr{\ol{f}}{\ral}B,~~~~\txt{for all}~~\al\in \Mor_\A(G,K),\nn\\
&&~~\Ra~~0=(\ol{f}\circ k_e)_\ast:\Mor_\A(G,K)\sr{(k_e)_\ast}{\ral}\Mor_\A(G,G^{\oplus I})\sr{(\ol{f})_\ast}{\ral}\Mor_\A(G,B),~~~~\txt{in mod-$R$},\nn\\
&&~~\Ra~~\ol{f}\circ k_e=0~~~~\txt{(since $\Mor_\A(G,-)$ is faithful)}.\nn
\eea
By the universal property of the kernel, we get a morphism $f:A\ra B$ such that ~$\ol{f}=f\circ e$ ~as follows:
\[\adjustbox{scale=0.8}{\bt
0\ar[r] & K\ar[r,hook,"k_e"] & G^{\oplus I}\ar[dr,"\ol{f}"']\ar[r,two heads,"e"] & A\ar[d,dashed,"f"]\ar[r] & 0\\
 &&& B &
\et}\]
Since $G$ is projective, any $h\in \Mor_\A(G,A)$ factors as $h=e\circ h'$ (as in the diagram below):
\[\adjustbox{scale=0.8}{\bt
 &&& G\ar[dl,dashed,"h'"']\ar[d,"h"] &\\
0\ar[r] & K\ar[r,hook,"k_e"] & G^{\oplus I}\ar[dr,"\ol{f}"']\ar[r,two  heads,"e"] & A\ar[d,dashed,"f"]\ar[r] & 0\\
 &&& B &
\et}\]
We can now obtain the desired representation of $f_R$ as follows:
\bea
f_R(h)=f_R(e\circ h')=(f_Re_\ast)(h')=\ol{f}_R(h')=\ol{f}\circ h'=(f\circ e)\circ h'=f\circ h.\nn
\eea
This proves that $E$ is full. It remains to prove $E$ is \ul{dense}.

Let $M=M_R$ be any right $R$-module. Then we know there is an exact sequence of the form
\bea
R^{\oplus I'}\sr{a}{\ral}R^{\oplus I}\sr{b}{\ral}M\ra 0.\nn
\eea
Since $\Mor_\A(G,-)$ commutes with coproducts (as $G$ is small, or alternatively, $I,I'$ are finite sets), we have equalities (up to isomorphism) as follows:
\bea
R^{\oplus I'}=\Mor_\A(G,G^{\oplus I'}),~~~~R^{\oplus I}=\Mor_\A(G,G^{\oplus I}).\nn
\eea
Thus, because $\Mor_\A(G,-)$ is full, we have $a=f_\ast:=\Mor_\A(G,f)$ for some $f\in \Mor_\A(G^{\oplus I'},G^{\oplus I})$. Hence, due to the \ul{exactness} of $\Mor_\A(G,-)$, we get an isomorphism as follows:
\bea\bt
\Mor_\A(G,G^{\oplus I'})\ar[d,equal]\ar[r,"f_\ast"] & \Mor_\A(G,G^{\oplus I})\ar[d,equal]\ar[rr,"b"] && M\ar[from=d,dashed,"\theta","\cong"']\ar[r] & 0\\
\Mor_\A(G,G^{\oplus I'})\ar[r,"f_\ast"] & \Mor_\A(G,G^{\oplus I})\ar[rr,"{(c_f)_\ast=c_{f_\ast}}"] && \Mor_\A(G,\coker f)\ar[r] & 0
\et\nn
\eea
where $b=\theta c_{f_\ast}$ implies the induced map $\theta$ is surjective (as both $b$ and $c_{f_\ast}$ are), and also injective since
\bea
\im f_\ast=\ker b=c_{f_\ast}^{-1}(\ker\theta)=c_{f_\ast}^{-1}(\ker\theta\cap\im c_{f_\ast})~~\Ra~~\ker\theta\sr{(s)}{=}\ker\theta\cap\im c_{f_\ast}=c_{f_\ast}(\im f_\ast)=0,\nn
\eea
where step (s) holds because $c_{f_\ast}$ is surjective.
\end{proof}

\subsection{Relevant results on essential extensions and injective envelopes}

\begin{dfn}[\textcolor{blue}{\index{Essential extension}{Essential extension}, \index{Injective! envelope}{Injective envelope}}]
Let $\A$ be an abelian category and $A,E\in\Ob\A$. A monic morphism $u:A\hookrightarrow E$ is an \ul{essential extension} (of $A$ by $E$) if for every monic morphism $u':A'\hookrightarrow E$,
\bea
&&u'\neq 0~~\Ra~~\im u~\cap~\im u'~\neq~0~~~~\txt{(or ~$A\cap A'\neq 0$)},\nn\\
&&\txt{equivalently},~~\im u~\cap~\im u'~=0~~\Ra~~\im u'=0~(\txt{i.e.,}~u'=0).\nn
\eea
An essential extension $u:A\hookrightarrow E$ (of $A$ by $E$) in which $E$ is an injective object is called an \ul{injective envelope} for $A$.
\end{dfn}

We will not be needing the following alternative description of an essential extension.
\begin{prp}
Let $\A$ be an abelian category. A monic $u:A\hookrightarrow E$ is an \ul{essential extension} (of $A$ by $E$) $\iff$ for every nonzero monomorphism $0\neq u':A'\hookrightarrow E$, there exist nonzero morphisms $0\neq \al:X\ra A$ and $0\neq \beta\ra A'$ (each possibly depending on $u'$) such that ~$u\al=u'\beta:X\ra E$.
\[\adjustbox{scale=0.8}{\bt
  & A\ar[dr,hook,"u"]  & \\
X\ar[ur,dashed,"\al\neq 0"]\ar[dr,dashed,"\beta\neq 0"']\ar[rr,"\al u=\beta u'"] &    & E \\
  & A'\ar[ur,hook,"\forall","u'\neq 0"'] &
\et}\]
Equivalently, $u:A\hookrightarrow E$ is an essential extension (of $A$ by $E$) $\iff$ for every nonzero monomorphism $u':A'\hookrightarrow E$, there exists a nonzero morphism $0\neq h:X\hookrightarrow E$ (possibly depending on $u'$) that factors through both $u$ and $u'$ as $h=\al u=\beta u'$, for some $\al:A\ra E$ and $\beta:A'\ra E$.
\end{prp}
\begin{proof}
By the definition of intersection, if $u:A\hookrightarrow E$ and $u':A'\hookrightarrow E$, then $A\cap A'=0$ $\iff$ for all morphisms $x:X\ra A$ and $x':X\ra A'$, $ux=u'x'$ $\Ra$ $x=0$ and $x'=0$ (i.e., $x\neq 0$ or $x'\neq 0$ $\Ra$ $ux\neq u'x'$). Therefore, $A\cap A'\neq 0$ $\iff$ there exist morphisms $x:X\ra A$ and $x':X\ra A'$ such that (i) $ux=u'x'$ while (ii) $x\neq 0$ or $x'\neq 0$, where the monic properties of $u$ and $u'$ imply that for (i) to be consistent with (ii) we must have both $x\neq 0$ and $x'\neq 0$ in (ii).
\[\adjustbox{scale=0.7}{\bt
 & E & \\
A\ar[ur,hook,"u"] & & A'\ar[ul,hook,"u'"']\\
 & A\cap B\ar[ul,hook,"q_A"']\ar[ur,hook,"q_B"] & \\
 & X\ar[u,dashed,"h"]\ar[uul,bend left,"x"]\ar[uur,bend right,"x'"'] &
\et}\]
That is, $A\cap A'\neq 0$ $\iff$ there exist \ul{nonzero} morphisms $0\neq x:X\ra A$ and $0\neq x':X\ra A'$ such that $ux=u'x':X\ra E$.
\end{proof}

\begin{lmm}[\textcolor{blue}{\cite[Lemma 2.1]{mitchell1964}}]\label{EssExtLmm1}
Let $\A$ be an abelian category. A morphism $u:A\ra B$ is an essential extension iff for every morphism $f:B\ra C$,
\bea
f u:A\sr{u}{\ral}B\sr{f}{\ral}C~~\txt{monic}~~\Ra~~B\sr{f}{\ral}C~~\txt{monic}.\nn
\eea
\end{lmm}
\begin{proof}
($\Ra$): Assume $u:A\ra B$ is an essential extension. Let $f:B\ra C$ be such that $f u:A\sr{u}{\ral}B\sr{f}{\ral}C$ is monic. Then $0=\ker(fu)\cong u^{-1}(\ker f)\cong\im u~\cap~\ker f$, and so $\ker f=0$.

($\La$): Assume $f u:A\sr{u}{\ral}B\sr{f}{\ral}C$ is monic implies $B\sr{f}{\ral}C$ monic. Suppose $u:A\ra B$ is non-essential. Then there is a subobject $0\neq A'\subset B$ such that $A\cap A'=0$. This means the composition $\pi u:A\sr{u}{\ral}B\sr{\pi}{\ral}B/A':=\coker(A'\hookrightarrow B)$ is monic, since
\bea
\ker(\pi u)\sr{h}{\cong}u^{-1}(\ker\pi)=\im u~\cap~\ker\pi= A\cap A'=0.\nn
\eea
This in turn implies $B\sr{\pi}{\ral}B/A'$ is monic (a contradiction).
\end{proof}

\begin{dfn}[\textcolor{blue}{\index{Sublocally small category}{Sublocally small category}, \index{Retract}{Retract}, \index{Retract of an object}{Retract of an object}}]
A category $\C$ is \ul{sublocally small} if for every object $C\in\Ob\C$, the subobjects $\P(C):=\{A\subset C\}$ form a set. A monomorphism $u:C\ra C'$ is a \ul{retract} (making the object $C$ a \ul{retract of the object} $C'$) if it has a left-inverse, i.e., if there exists a morphism $u':C'\ra C$ ``connecting'' $u$ to the identity in the sense
\[
u'u=id_C~~~~~~~~
\adjustbox{scale=0.8}{\bt
0\ar[r] & C\ar[d,hook,"id_C"]\ar[rr,hook,"u"] && C'\ar[dll,dashed,"u'"] \\
  & C&&
\et}\]
\end{dfn}

\begin{lmm}[\textcolor{blue}{\cite[Lemma 2.3]{mitchell1964}}]\label{EssExtLmm01}
Let $\A$ be a sublocally small abelian AB5 category. For any object $C$ in $\A$, the following are equivalent:
\bit
\item[(i)] $C$ is injective.
\item[(ii)] $C$ has no proper essential extension. (\blue{footnote}\footnote{A subobject $u:A\hookrightarrow B$ is called \ul{proper} if $u$ is not an isomorphism.})
\item[(iii)] $C$ is a retract of every object that contains it.
\eit
\end{lmm}
\begin{proof}
{\flushleft\ul{(i)$\Ra$(ii)}}: Assume $C$ is injective. Let $u:C\ra E$ be an essential extension. Then by the injectivity of $C$, there exists a morphism $\theta:E\ra C$ such that $\theta u=id_C$ (i.e., $C$ is a retract of $E$).
\[\adjustbox{scale=0.8}{\bt
0\ar[r] & C\ar[d,hook,"1_C"]\ar[rr,hook,"u"] && E\ar[dll,dashed,"\theta"]\ar[r,tail,two heads] & {E\over C}\ar[r] & 0 \\
  & C&&
\et}\]
From the resulting split-exact sequence, we get $E\cong C\oplus{E\over C}$.
Since $u:C\hookrightarrow E$ is essential, $\im u\cong C$, $\ker\theta\cong{E\over C}$, and $C\cap{E\over C}\cong\im u\cap\ker\theta\cong u^{-1}\big(\ker\theta\big)\cong \ker(\theta u)=0$, we have ${E\over C}=0$ (i.e., $E\cong C$).

{\flushleft \ul{(ii)$\Ra$(iii)}}: Assume $C$ has no proper essential extensions. Suppose a monic morphism $u:C\hookrightarrow A$ is not a retract, i.e., for every morphism $u':A\ra C$, we have
\bea
\bt 1_C\neq u'u:C\ar[r,"u"] & A\ar[r,hook,"u'"] & C.\et\nn
\eea
Let $\I:=\{B\subset A:~B\cap C=0\}$ be the set of subobjects of $A$ intersecting $C$ trivially. If $\{B_i\}\subset\I$ is a chain, then $\left(\bigcup B_i\right)\cap C=\bigcup B_i\cap C=0$, and so $\bigcup B_i\in\I$. This shows $\I$ has a maximal element $A'$ by Zorn's lemma. Since $C$ is not a retract of $A$, {\small $C\oplus A'\subsetneq A$}, and so we get a proper extension
\bea
\textstyle C={C\oplus A'\over A'}\subsetneq{A\over A'}\nn
\eea
of $C$ which by hypotheses is not essential. Thus, there is a nonzero subobject $0\neq \ol{B}\subset A/A'$ such that $C\cap\ol{B}=0$. It follows (by the diagram below, which is exact by Lemma \ref{PsPbPartII3} and its corollary) that
\bea
\ol{B}=B/A',~~~~\txt{for some subobject $B\subset A$ such that $A'\subsetneq B$ (since $\ol{B}\neq0$)}.\nn
\eea
\[\adjustbox{scale=0.8}{\bt
  & 0\ar[d]               & 0\ar[d]   & 0\ar[d] &  \\
0\ar[r] & A'\ar[d,equal]\ar[r,hook,dashed,"h"]  & B:=u'^{-1}(\ol{B})\ar[d,hook,"u"]\ar[r,two heads,"c"] & \ol{B}\ar[d,hook,"\ol{u}"]\ar[r] & 0 \\
0\ar[r] & A'\ar[d]\ar[r,hook,"u'"] & A\ar[d,tail,two heads,"c_u"]\ar[r,tail,two heads,"c_{u'}"] & A/A'\ar[d,tail,two heads,"c_{\ol{u}}"]\ar[r] & 0 \\
0\ar[r] & 0\ar[r]\ar[d]         & A/B\ar[d]\ar[r,dashed,two heads,"h'"] & {A/A'\over\ol{B}}\ar[d]\ar[r] & 0 \\
  & 0               & 0   & 0 &
\et}\]
Since $C\cap{B\over A'}\cong C\cap\ol{B}=0$, the following multiple-pullbacks commutative diagram shows $C\cap B=0$ (which contradicts the maximality of $A'$):
\[\adjustbox{scale=0.7}{\bt
 & \overbrace{A'\cap C}^0\ar[ddd,hook]\ar[rrr,hook]\ar[dr,dashed,hook] &  &  & C\ar[ddddl,hook,bend left=10] & \\
 &  & B\cap C\ar[dd,hook]\ar[urr,hook]\ar[dr,dashed] &  &  & \\
 &  &  & \overbrace{\ol{B}\cap C}^0\ar[d,hook]\ar[uur,hook] &  & \\
0\ar[r] & A'\ar[d,equal]\ar[r,hook] & B\ar[d,hook]\ar[r,tail,two heads] & {B\over A'}\ar[d,dashed,hook]\ar[r] & 0 & \\
0\ar[r] & A'\ar[r,hook] & A~(\txt{pushout})\ar[r,tail,two heads] & {A\over A'}\ar[r] & 0 &
\et}~~~~~~~~\textstyle \ol{B}:={B\over A'}\]

{\flushleft\ul{(iii)$\Ra$(i)}}: Assume $C$ is a retract of every object containing it. Let $m:A\hookrightarrow B$ be a monomorphism, and $h:A\ra C$ any morphism. Then we can take the pushout-limit $P:=\varinjlim\left(C\sr{h}{\lal}A\sr{m}{\ral}B\right)$ to get the following commutative diagram:
\[\adjustbox{scale=0.8}{\bt
0\ar[r] & A\ar[dd,"h"']\ar[rr,hook,"m"] && B\ar[ddll,dashed,"m''h'"']\ar[dd,dashed,"h'"]\\
 & && \\
 & C\ar[rr,hook,dashed,shift left,"m'"] && P\ar[ll,dashed,shift left,"m''"]
\et}~~~~~~~~~~m'h=h'm,~~~~~~m''m'=1_C
\]
Since $m'$ is monic (by the pullback/pushout theorems), i.e., $C\subset P$, we get a morphism $m'':P\ra C$ such that $m''m'=1_C$.
\end{proof}

\begin{lmm}[\textcolor{blue}{\cite[Lemma 2.4]{mitchell1964}}]\label{EssExtLmm2}
Let $\A$ be an abelian category. If $u: A\ra B$ and $v: B\ra C$ are monic, then $vu:A\sr{u}{\ral}B\sr{v}{\ral}C$ is an essential extension if and only if both $u$ and $v$ are essential extensions.
\end{lmm}
\begin{proof}
($\Ra$): Assume $vu:A\sr{u}{\ral}B\sr{v}{\ral}C$ is essential. Suppose $v$ is not essential. Then there is a non-monic $f:C\ra C'$ such that $fv$ is monic. But $fvu$ monic implies $f$ monic (a contraction). Thus $v$ is essential.

Next, suppose $u$ is not essential. Then there exists $0\neq B'\subset B$ such that $B'\cap\im u=0$. Thus, because $B'\subset C$, the composition $A\sr{u}{\ral}B\sr{v}{\ral}C\sr{\pi}{\ral}C/B'$ is monic due to
\bea
\ker(\pi vu)=(vu)^{-1}(\ker\pi)=u^{-1}(v^{-1}(B'))=u^{-1}(\im v~\cap~B')=\im u~\cap~\im v~\cap~B=0.\nn
\eea
But since $vu$ is essential, this implies $\pi$ is monic (a contradiction).

($\La$): Assume $u,v$ are essential. Then for any $f:C\ra C'$,
\[
fvu~~\txt{monic}~~\Ra~~~~fv~~\txt{monic},~~\Ra~~f~~\txt{monic}. \qedhere
\]
\end{proof}

\begin{lmm}[\textcolor{blue}{Uniqueness of the injective envelope (up to a non-unique isomorphism): \cite[Proposition 2.8]{mitchell1964}}]\label{EssExtLmm3}
Let $\A$ be an abelian category, $A\in\Ob\C$, and $u:A\ra E$ and $u':A\ra E'$ injective envelopes for $A$. Then there exists an isomorphism $E\cong E'$.
\end{lmm}
\begin{proof}
Since $E,E'$ are injective, a morphism $\theta:E\ra E'$ exists such that diagram (a) below commutes:
\bea
(a)~~\adjustbox{scale=0.8}{\bt
0\ar[r] & A\ar[d,hook,"u'"]\ar[rr,hook,"u"] && E\ar[dll,dashed,"\theta"] \\
   & E'
\et}~~~~~~~~~~~~(b)~~
\adjustbox{scale=0.8}{\bt
0\ar[r] & E\ar[d,hook,"1_E"]\ar[rr,hook,"\theta"] && E'\ar[dll,dashed,"\theta'"] \\
   & E
\et}
\nn
\eea
By Lemma \ref{EssExtLmm1}, $\theta$ is a monomorphism (which is an essential extension by Lemma \ref{EssExtLmm2}). Hence, by the injectivity of $E$, $\theta$ is an isomorphism.
\end{proof}

\subsection{Abelian AB5 category with a generator: Existence of an injective cogenerator}
\begin{lmm}[\textcolor{blue}{\cite[Lemma 7.12, p.139]{freyd1964}}]\label{EssExtLmm4}
Let $\A$ be an abelian category and $u:F\ra E$ an essential extension in $Ab^\A$. If $F$ preserves monomorphisms, then so does $E$.
\end{lmm}
\begin{proof}
Suppose on the contrary that $E$ does not preserve monomorphisms (even though $F$ does). Then there exists a monic morphism $m:A'\ra A$ and $0\neq x\in E(A')$ such that
\bea
E(m):E(A')\ra E(A),~~x\mapsto E(m)(x)=0~~~~\txt{(i.e., $E(m)$ is not a monomorphism)}.\nn
\eea
Define a subfunctor $E_x\subset E$ by the following (for all $B\in\Ob\A$ and $h\in\Mor\A$):
{\small\begin{align}
&E_x:\A\ra Ab,~~B\sr{h}{\ral}C~~\longmapsto~~E_x(B)\sr{E_x(h)}{\ral}E_x(C),~~~~E_x(h):=E(h)|_{E_x(B)},\nn\\
&E_x(B):=E\big(\Mor_\A(A',B)\big)(x):=\{E(f)(x):f\in\Mor_\A(A',B)\}\nn\\
&~~~~=\left\{y\in E(B):~y=E(f)(x)~\txt{for some}~f:A'\ra B\right\}\subset E(B),\nn
\end{align}}
where it is clear that $E_x(B)\subset E(B)$ is a subgroup, because $E_x(B)+E_x(B)=\{E(f)(x)+E(g)(x):f,g\in E(\Mor_\A(A',B))\}=\{E(f+g)(x):f,g\in E(\Mor_\A(A',B))\}\subset E_x(B)$. We have $E_x\neq 0$, since
\bea
x=id_{A'}(x)\in E_x(A'):=E\big(\Mor_\A(A',A')\big)(x):=\{E(f)(x):f\in\Mor_\A(A',A')\}.\nn
\eea
Since $E$ is essential, $F\cap E_x\neq 0$, and so there exists an object $B\in\Ob\A$ such that
\bea
F(B)\cap E_x(B)\neq 0.\nn
\eea
Let $0\neq y\in F(B)\cap E_x(B)$. Then (by the definition of $E_x$) $y=E(f)(x)$ for some $f\in \Mor_\A(A',B)$. In the following \emph{pushout diagram} (i.e., $P=\varinjlim(B\sr{f}{\lal}A'\sr{m}{\ral}A)$ ), $m'$ is necessarily a monomorphism (by the pushout/pullback theorems, Lemma \ref{PushPullThm}):
\[\adjustbox{scale=0.8}{\bt
A'\ar[d,"f"]\ar[r,hook,"m"] & A\ar[d,dashed,"f'"] \\
B\ar[r,dashed,hook,"m'"] & P
\et}\]
Since $F$ preserves monic morphisms, we have $F(m')(y)\neq 0$, which in turn gives
\bea
&&0\neq E(m')(y)=E(m')\big(E(f)(x)\big)=E(m')E(f)(x)=E(m'f)(x)=E(f'm)(x)\nn\\
&&~~~~=E(f')E(m)(x)=E(f')\big(E(m)(x)\big)=E(f')(0)=0,\nn
\eea
which is a contradiction.
\end{proof}

\begin{lmm}[\textcolor{blue}{\cite[Lemma 2.5]{mitchell1964}}]\label{EssExtLmm5}
Let $\A$ be an abelian AB5 category, $A\in\Ob\A$, and $A'\subset A$. For any family of essential extensions $A_i\subset A'$ of $A'$, the union $\bigcup A_i\subset A'$ is an essential extension of $A'$.
\end{lmm}
\begin{proof}
Let $0\neq B\subset A'$. Then $B\cap\left(\bigcup A_i\right)=\bigcup B\cap A_i\neq 0$, since $B\cap A_i\neq 0$ for each $i$.
\end{proof}

\begin{lmm}[\textcolor{blue}{\cite[Lemma 2.7]{mitchell1964}}]\label{EssExtLmm6}
Let $\A$ be an abelian category, $G\in\Ob\A$ a generator, and $R:=\Mor_\A(G,G)$ the endomorphism ring of $G$. Then the imbedding ~$T:=\Mor_\A(G,-):\A\ra \txt{mod-}R$ ~preserves essential extensions.
\end{lmm}
\begin{proof}
Let $u:A\ra B$ be an essential extension in $\A$. We need to show $u_\ast:\Mor_\A(G,A)\ra \Mor_\A(G,B)$ is an essential extension, i.e., $\Mor_\A(G,A)\cap hR\neq 0$, for all $0\neq h\in \Mor_\A(G,B)$. So, fix $0\neq h\in \Mor_\A(G,B)$. We need to find $r\in R$ such that $0\neq hr\in \Mor_\A(G,A)$, i.e.,
\bea
hr\neq 0~~~~\txt{and}~~~~\im(hr)\subset A.\nn
\eea
Since $u$ is essential, ~$A\cap\im h\neq0$. So, we have a nontrivial preimage $h^{-1}(A)\cong h^{-1}(A\cap\im h)$ given by the following pullback diagram:
{\footnotesize
\bea\bt
G\ar[d,shift left,dashed,"g"]\ar[rr,"h"] && B \\
h^{-1}(A)\ar[u,hook,shift left,"u'"]\ar[rr,"h'"] && A\ar[u,hook,"u"]
\et~~~~~~~~h'\neq 0.\nn
\eea}Since $h'\neq 0$ and $T:=\Mor_\A(G,-)$ is faithful, there exists $g:G\ra h^{-1}(A)$ such that $0\neq (h')_\ast(g)=h'\circ g$.

Set $r:=u'g:G\ra G$. Then $hr=hu'g=uh'g=u_\ast(h'g)\in\im u_\ast\cong \Mor_\A(G,A)$.
\end{proof}

\begin{lmm}[\textcolor{blue}{\cite[Theorem 2.9]{mitchell1964}}]\label{EssExtLmm7}
Let $\A$ be an abelian AB5 category with a generator. Then $\A$ has injective envelopes for each of its objects.
\end{lmm}
\begin{proof}
Let $G\in\Ob\A$ be a generator, $R:=\Mor_\A(G,G)$ the endomorphism ring of $G$, and $T:\A\ra \txt{mod-}R$ the imbedding
\bea
T:A\sr{t}{\ral}B~\mapsto~\Mor_\A(G,A)\sr{t_\ast}{\ral}\Mor_\A(G,B).\nn
\eea
Given $A\in\Ob\A$, let $m:T(A)\hookrightarrow E$ be a monomorphism in mod-$R$, for an injective right $R$-module $E$. Let $\I:=\{(B,u,f)~|~u:A\hookrightarrow B~\txt{an essential extension},~f:T(B)\ra E~~\txt{satisfying}~~fu_\ast=m\}$.
\[\adjustbox{scale=0.9}{\bt
0 \ar[r] & T(A)\ar[rrr,hook,bend left=30,"u'_\ast"]\ar[dd,hook,"m"]\ar[rr,hook,"u_\ast"] && T(B)\ar[ddll,hook,"f"']\ar[r,dashed,hook,"v_\ast"] & T(B')\ar[ddlll,hook,"f'"] \\
 & \\
 & E
\et}~~~~~~~~
\adjustbox{scale=0.9}{\bt
0 \ar[r] & T(A)\ar[rrr,hook,bend left=30,"u_\ast"]\ar[dd,hook,"m"]\ar[rr,hook,"u'_\ast"] && T(B')\ar[ddll,hook,"f'"']\ar[r,dashed,hook,"\wt{v}_\ast"] & T(B)\ar[ddlll,hook,"f"] \\
 & \\
 & E
\et}\]
where $u_\ast$ is an essential extension by Lemma \ref{EssExtLmm6}, and so $f$ is a monomorphism by Lemma \ref{EssExtLmm1}. Order $\I$ such that $(B,u,f)\leq(B',u',f')$ if there exists a morphism $v:B\ra B'$ with $vu=u'$ and $f'v_\ast=f$.

Note that $v$ above is a \ul{monomorphism} by Lemma \ref{EssExtLmm1}, since $u$ is an essential extension. Also, $v$ is \ul{unique} because for any monomorphism $v_1:B\hookrightarrow B'$ satisfying $v_1u=u'$ and $f'(v_1)_\ast=f$, the equality $f'v_\ast=f'(v_1)_\ast$ implies $v_\ast=(v_1)_\ast$ as $f'$ is monic, which means $v=v_1$ since $T$ is faithful. Thus, if we also have $(B',u',f')\leq (B,u,f)$, then we get another unique monomorphism $\wt{v}:B'\ra B$ satisfying $\wt{v}u'=u$ and $f\wt{v}_\ast=f'$, and so the relations $\wt{v}u=u'$ and $\wt{v}u'=u$ along with uniqueness of $v,\wt{v}$ imply $v$ is an \ul{isomorphism}.

Thus, we have an \ul{equivalence relation} on $\I$ given by $(B,u,f)\sim(B',u',f')$ if ``$(B,u,f)\leq(B',u',f')$ and $(B',u',f')\leq(B,u,f)$''.

By Lemma \ref{MitThmLmm1}, we have an epimorphism $e_B:G^{\oplus T(B)}\twoheadrightarrow B$. Also, for any $(B,u,f)\in\I$, we have the cardinality relation $|T(B)|\leq |E|$ due to the monomorphism $f:T(B)\hookrightarrow E$. Thus, we have an epimorphism
\bea
\textstyle\bt e~:~G^{\oplus E}\ar[r,two heads,"e_f"] & G^{\oplus T(B)}\ar[r,two heads,"e_B"] & B\et,~~~~\txt{i.e., $B$ is a quotient object of ~$G^{\oplus E}:=\bigoplus_{x\in E}G$}.\nn
\eea
Also, $\A$ is sublocally small (since $\A$ has a generator), and so the quotient objects of $G^{\oplus E}$ form a set. It follows that the collection of equivalence classes $\I_0:={\I\over\sim}$ forms a set.

Let $\{(B_j,u_j,f_j)\}_{j\in \J}$ be a chain in $\I_0$. Then we get a \ul{directed system} $S_\J:\J\ra\A$ whose objects are $\{A\}\cup\{B_j\}_{j\in \J}$, its transition morphisms are the (equivalence classes of the) \ul{monomorphisms} $v_{jj'}:B_j\hookrightarrow B_{j'}$, for $j\leq j'$, i.e., $(B_j,u_j,f_j)\leq (B_{j'},u_{j'},f_{j'})$, satisfying $v_{jj'}u_j=u_{j'}$ and $f_{j'}(v_{jj'})_\ast=(u_j)_\ast$. By construction, the directed colimit $B_0:=\varinjlim S_\J$ of the system of the directed system $S_\J:\J\ra\A$ is a monic colimit involving the \ul{union} (with respect to $A$) $\bigcup_{j\in\J}^AB_j$ of the $B_j$'s. Specifically, the limit is given by a system of subobjects $q_j:B_j\hookrightarrow B_0$ satisfying appropriate compatibility relations, namely, commutativity of the diagrams below (in which $j\leq j'$):
\[\adjustbox{scale=0.9}{\bt
 &A\ar[dl,hook,"u_j"']\ar[dr,hook,"u_{j'}"]& \\
B_j\ar[ddr,bend right,"x_j"']\ar[dr,hook,"q_j"]\ar[rr,hook,"v_{jj'}"] && B_{j'}\ar[dl,hook,"q_{j'}"']\ar[ddl,bend left,"x_{j'}"] \\
 & B_0:=\varinjlim S_\J\ar[d,dashed,"h"] & \\
 & X &
\et}\]

By Lemma \ref{EssExtLmm5}, the inclusion $u_0:=q_ju_j:A\hookrightarrow B_0$ (independent of $j$) is an essential extension, since each $u_j:A\hookrightarrow B_j$ is essential.
\[\adjustbox{scale=0.9}{\bt
0 \ar[r] & T(A)\ar[rrr,hook,bend left=30,"(u_0)_\ast"]\ar[dd,hook,"m"']\ar[rr,hook,"(u_j)_\ast"] && T(B_j)\ar[ddll,hook,"f_j"']\ar[r,hook,"(q_j)_\ast"] & T(B_0)\ar[ddlll,hook,dashed,"f_0"] \\
 && &&  \\
 & E
\et}~~~~~~~~
\adjustbox{scale=0.9}{\bt
0 \ar[r] & T(A)\ar[rrr,hook,bend left=30,"(u_0)_\ast"]\ar[dd,hook,"m"']\ar[rr,hook,"\bigcup(u_j)_\ast"] && \bigcup T(B_j)\ar[ddll,hook,"\bigcup f_j"']\ar[r,hook,dashed,"\delta"] &  T(B_0)\ar[ddlll,hook,dashed,"f_0"]\\
 &  && &\\
 & E && &
\et}\]
Thus, the induced monomorphism $\bigcup f_j:\bigcup T(B_j)\hookrightarrow E$ on the limit of the directed system ~$T\circ S_\J:\J\ra \txt{mod-}R$~ extends (by the injectivity of $E$) to a monomorphism $f_0:T(B_0)\ra E$. Since (co)limits are functors, and hence preserve commutativity relations, we see that $(B_0,u_0,f_0)\in\I_0$, i.e., the chain $\{(B_j,u_j,f_j)\}_{j\in\J}$ has an upper bound in $\I_0$, and so by Zorn's lemma, $\I_0$ has a maximal element $(B_E,u_E,f_E)$.

Suppose $B_E$ is not injective. Then by Lemma \ref{EssExtLmm01}, $B_E$ has a proper essential extension $v_1:B_E\hookrightarrow B_1$, which by Lemma \ref{EssExtLmm2} gives an essential extension $\bt[] v_1u_E:A\ar[r,hook,"u_E"] & B_E\ar[r,hook,"v_1"] & B_1\et$.
\[\adjustbox{scale=0.9}{\bt
0 \ar[r] & T(A)\ar[rrr,hook,bend left=30,"(v_1u_E)_\ast"]\ar[dd,hook,"m"]\ar[rr,hook,"(u_E)_\ast"] && T(B_E)\ar[ddll,hook,"f_E"']\ar[r,dashed,hook,"(v_1)_\ast"] & T(B_1)\ar[ddlll,hook,dashed,"f_1"] \\
 & \\
 & E
\et}\]
By the injectivity of $E$, we can extend $f_E$ to a monomorphism $f_1:T(B_1)\ra E$. This contradicts the maximality of $(B_E,u_E,f_E)$ in $\I_0$, since $(B_E,u_E,f_E)<(B_1,v_1u_E,f_1)\in\I_0$.
\end{proof}

\begin{lmm}[\textcolor{blue}{\cite[Proposition 2.10]{mitchell1964}}]\label{EssExtLmm8}
Let $\A$ be an abelian category with enough injectives and a generator. Then $\A$ has an injective cogenerator (\blue{footnote}\footnote{Recall that an injective cogenerator for $\A$ is an object $G'\in\A$ such that $\Mor_\A(-,G'):\A\ra Ab$ is an exact imbedding.}).
\end{lmm}
\begin{proof}
Let $G\in\Ob\A$ be a generator. Since $\A$ has enough injectives, let $m:\bigoplus_{V\subset G}{G\over V}\hookrightarrow G'$ be a monomorphism, for an injective object $G'$. We will show that $G'$ is a cogenerator, i.e., $\Mor_\A(-,G'):\A\ra Ab$ is an imbedding. By Lemma \ref{MitRefLmm} (or better, Corollary \ref{MitRefCrl}), it suffices to show that $\Mor_\A(-,G')$ is injective on objects. Fix a nonzero object $A\neq 0$. We need to find a nonzero morphism $0\neq f:A\ra G'$.

Since $G$ is a generator, there exists $0\neq g:G\ra A$. We have
\bea
\textstyle \im g\cong\coim g=\coker(\ker g\ra G)={G\over\ker g},~~~~\ker g\subsetneq G.\nn
\eea
Since $G'$ is injective, the composition {\footnotesize${G\over\ker g}\hookrightarrow\bigoplus_{V\subset G}{G\over V}\hookrightarrow G'$} can be lifted to a morphism $0\neq f:A\ra G'$.
\[\adjustbox{scale=0.8}{\bt
0\ar[r] & {G\over\ker g}\cong\im g\ar[d,hook,"i\neq 0"]\ar[rr,hook,"kc_g"] && A\ar[ddll,dashed,"f\neq 0"] \\
 & \bigoplus_{V\subset G}{G\over V}\ar[d,hook,"m"'] && \\
 & G' &&
\et}~~~~\qedhere
\]
\end{proof}

The following result is an immediate consequence of Lemmas \ref{EssExtLmm7} and \ref{EssExtLmm8}.
\begin{crl}[\textcolor{blue}{\cite[Corollary 2.11]{mitchell1964}}]\label{EssExtLmm9}
An abelian AB5 category with a generator has an injective cogenerator.
\end{crl}

\section{The Category of Left-exact Ab-valued Functors: Abelian and AB5 Properties}

\begin{dfn}[\textcolor{blue}{Category of \index{Category of! left-exact Ab-valued functors}{left-exact} (resp. \index{Category of! right-exact Ab-valued functors}{right-exact}) Ab-valued functors}]
Let $\A$ be an abelian category. The category of left-exact (resp. right-exact) Ab-valued functors is the full subcategory $\L(Ab^\A)\subset Ab^\A$ (resp. $\R(Ab^\A)\subset Ab^\A$) whose objects are left-exact (resp. right-exact) additive functors $F:\A\ra Ab$.
\end{dfn}

In the following result, the abelian property of $\L(Ab^\A)$ is stated in \cite[Theorem 4.2, p.635]{mitchell1964}.
\begin{lmm}\label{AbelLeftEx}
For any abelian category $\A$, the category $\L(Ab^\A)$ is abelian and AB5.
\end{lmm}
\begin{proof}[\magenta{Proof starter}]
To show $\L(Ab^\A)$ is \ul{additive}, we need to show it has a $0$ object, is closed under finite direct sums, and the composition of morphisms is distributive over addition. It is clear that $\L(Ab^\A)$ has a $0$ object and that the composition of morphisms is distributive over addition (as properties inherited from $Ab^\A$). If $F,G\in \L(Ab^\A)$ and $0\ra A\ra B\ra C\ra 0$ is exact in $\A$, then $0\ra F(A)\ra F(B)\ra F(C)$, $0\ra G(A)\ra G(B)\ra G(C)$ are exact in $Ab$, and so $0\ra F(A)\oplus G(A)\ra F(B)\oplus G(B)\ra F(C)\oplus (C)$ is also exact in $Ab$, i.e., $F\oplus G\in \L(Ab^\A)$. Hence $\L(Ab^\A)$ is additive.

To show $\L(Ab^\A)$ is \ul{abelian}, we only need to show $\L(Ab^\A)$ is closed under the kernel and cokernel operations, i.e., if $F,G\in\L(Ab^\A)$, then the kernel and cokernel of $\eta:F\ra G$, defined by a diagram
\[\bt
0\ar[r] & \ker\eta \ar[r,"k_\eta"] & F \ar[r,"\eta"] & G\ar[r,"c_\eta"]\ar[r] &\coker\eta \ar[r] & 0
\et\]
along with its naturality (i.,e., a commutative diagram) given by the commutative diagram (\blue{footnote}\footnote{This diagram was obtained in the proof of Proposition \ref{YonProp}.}):
\[\bt
 A\ar[d,"f"] & 0\ar[r] & (\ker\eta)(A)\ar[d,dashed,"(\ker\eta)(f)"] \ar[r,"(k_\eta)_A"] & F(A)\ar[d,"F(f)"] \ar[r,"\eta_A"] & G(A)\ar[d,"G(f)"]\ar[r,"(c_\eta)_A"]\ar[r] &(\coker\eta)(A)\ar[d,dashed,"(\coker\eta)(f)"] \ar[r] & 0 \\
  B& 0\ar[r] & (\ker\eta)(B) \ar[r,"(k_\eta)_B"] & F(B) \ar[r,"\eta_B"] & G(B) \ar[r,"(c_\eta)_B"]\ar[r] &(\coker\eta)(B)
  \ar[r] & 0
\et\]
are also in $\L(Ab^\A)$. Since $\L(Ab^\A)\subset Ab^\A$ is a full subcategory, we only need to show $\ker\eta$, $\coker\eta$ are left-exact functors. To show the latter, consider the following extended commutative diagram defined/obtained through the same steps as in the proof Proposition \ref{YonProp}.
\[\bt
 0\ar[d] &  & 0\ar[d] &0\ar[d] & 0\ar[d] &0\ar[d] \\
 A\ar[d,"f"] & 0\ar[r] & (\ker\eta)(A)\ar[d,"(\ker\eta)(f)"] \ar[r,"(k_\eta)_A"] & F(A)\ar[d,"F(f)"] \ar[r,"\eta_A"] & G(A)\ar[d,"G(f)"]\ar[r,"(c_\eta)_A"]\ar[r] &(\coker\eta)(A)\ar[d,"(\coker\eta)(f)"] \ar[r] & 0 \\
 B\ar[d,"g"] & 0\ar[r] & (\ker\eta)(B)\ar[d,"(\ker\eta)(g)"] \ar[r,"(k_\eta)_B"] & F(B)\ar[d,"F(g)"] \ar[r,"\eta_B"] & G(B)\ar[d,"G(g)"]\ar[r,"(c_\eta)_B"]\ar[r] &(\coker\eta)(B)\ar[d,"(\coker\eta)(g)"]\ar[r] & 0\\
 C\ar[d] & 0\ar[r] & (\ker\eta)(C) \ar[r,"(k_\eta)_C"] & F(C) \ar[r,"\eta_C"] & G(C) \ar[r,"(c_\eta)_C"]\ar[r] &(\coker\eta)(C)\ar[r] & 0\\
 0 &  & & & &  &
\et\]

Let the sequence $0\ra A\sr{f}{\ral}B\sr{g}{\ral}C\ra 0$ in $\A$ be exact. Then $A\cong\im f\cong\ker g$ (i.e., $f=k_g$). Since $F,G$ are kernel-preserving, we also have $F(A)\cong\ker F(g)$, $G(A)\cong\ker G(g)$, (i.e., $F(f)=k_{F(g)}$, $G(f)=k_{G(g)}$). We need to show the following:
\bit[leftmargin=0.9cm]
\item[(i)] $(\ker\eta)(A)\cong\ker\big((\ker\eta)(g)\big)$, i.e., $(\ker\eta)(f)=k_{(\ker\eta)(g)}$.
\item[(ii)] $(\coker\eta)(A)\cong\ker\big((\coker\eta)(g)\big)$, i.e., $(\coker\eta)(f)=\ker_{(\coker\eta)(g)}$.
\eit
By the commutativity of the diagram, we have (1) $F(f)\circ(k_\eta)_A=(k_\eta)_B\circ(\ker\eta)(f)$, and so by taking kernels, we see that $(\ker\eta)(f)$ is monic. Also, by composing (1) with $F(g)$ on the left, we get
\bea
&&0=F(g)\circ F(f)\circ(k_\eta)_A=F(g)\circ (k_\eta)_B\circ(\ker\eta)(f)=(k_\eta)_C\circ(\ker\eta)(g)\circ(\ker\eta)(f),\nn\\
&&~~\Ra~~(\ker\eta)(g)\circ(\ker\eta)(f)=0,~~\Ra~~\im\big((\ker\eta)(f)\big)\subset\ker\big((\ker\eta)(g)\big).\nn
\eea
Let $b'\in\ker\big((\ker\eta)(g)\big)$, i.e., $(\ker\eta)(g)(b')=0$. Then
\bea
&&0=(k_\eta)_C\circ (\ker\eta)(g)(b')=F(g)\circ(k_\eta)_B(b')~~\Ra~~(k_\eta)_B(b')\in\ker F(g)=\im F(f),\nn\\
&&~~\Ra~~(k_\eta)_B(b')=F(f)(a'),~~\Ra~~0=\eta_B\circ(k_\eta)_B(b')=\eta_B\circ F(f)(a')=G(f)\circ\eta_A(a'),\nn\\
&&~~\Ra~~\eta_A(a')=0,~~a'\in\ker\eta_A=\im(k_\eta)_A,~~\Ra~~a'=(k_\eta)_A(a''),\nn\\
&&~~\Ra~~(k_\eta)_B(b')=F(f)(a')=F(f)\circ(k_\eta)_A(a'')=(k_\eta)_B\circ(\ker\eta)(f)(a''),\nn\\
&&~~\Ra~~b'=(\ker\eta)(f)(a'')\in\im\big((\ker\eta)(f)\big).\nn
\eea
This proves (i) above. \magenta{However, (ii) is not true in general. Consequently, we might need to use the direct definition of the cokernel of a left-exact functor within $\L(Ab^\A)$ where the universal property involves only left-exact functors (instead of all functors as in our naive choice above). We will develop a full proof in the rest of this section.}

\magenta{Although we have not yet completed the argument above, assume temporarily that this has been done.} Then it remains to show $\L(Ab^\A)$ is an \ul{AB5} category, as follows (\blue{footnote}\footnote{Alternative rough argument: This uses the fact that (i) Ab is AB5 and (ii) operations in $\L(Ab^\A)$ are pointwise versions of operations in Ab, and so limits of systems in $\L(Ab^\A)$ are in $\L(Ab^\A)$, since by Remark \ref{AB5Equiv1}, the AB5 condition is equivalent to exactness of the colimit functor $\varinjlim$.}): We know that colimits $\varinjlim$ (resp. limits $\varprojlim$) in Ab = $\Integer$-mod can be expressed in terms of cokernels and coproducts (resp. kernels and products). Since we have seen above that $\L(Ab^\A)$ is closed under kernel and cokernel operations (and the product, as a functor, is left-exact), to show $\L(Ab^\A)$ is closed under colimits and limits, it remains only to argue that $\L(Ab^\A)$ is closed under coproducts. This is indeed the case: Given a SES $0\ra A\ra B\ra C\ra 0$ in $\A$, a family of left-exact functors $F_i:\A\ra Ab$ gives exact sequences in Ab, $0\ra F_i(A)\ra F_i(B)\ra F_i(C)$, whose direct sum
\bea
\textstyle 0\ra \bigoplus F_i(A)\ra \bigoplus F_i(B)\ra \bigoplus F_i(C)\nn
\eea
is also exact in Ab (by the exactness of $\bigoplus$ in Ab = $\Integer$-mod).

Having shown that $\L(Ab^\A)$ is closed under colimits and limits, it now follows that $\L(Ab^\A)$ has an AB5 property induced \ul{pointwise} by the AB5 property of Ab (as done for $Ab^\A$ in Proposition \ref{AB5ForLAb}).
\end{proof}

\begin{lmm}[\textcolor{blue}{Right-exactness of injective objects in $Ab^\A$: \cite[Proposition 7.11, p.138]{freyd1964}}]\label{RightExInj}
Let $\A$ be an additive category. If a functor $T:\A\ra Ab$ is an injective object in $Ab^\A$, then $T$ is right-exact.
\end{lmm}
\begin{proof}
Let $0\ra A\ra B\ra C\ra 0$ be exact in $\A$. Then with $H^A:=\Mor_\A(A,-)$, we get an exact sequence $0\ra H^C\ra H^B\ra H^A$ in $Ab^\A$. Applying $Hom(-,T)$, Yoneda's lemma gives the exact sequence
\[\adjustbox{scale=0.8}{\bt
\Mor_{Ab^\A}(H^A,T)\ar[d,"\cong"]\ar[r] & \Mor_{Ab^\A}(H^B,T)\ar[d,"\cong"]\ar[r] & \Mor_{Ab^\A}(H^C,T)\ar[d,"\cong"]\ar[r] & 0\\
T(A)\ar[r] & T(B)\ar[r] & T(C)\ar[r] & 0
\et}\nn  \qedhere
\]
\end{proof}

\begin{dfn}[\blue{
\index{Injective! functor}{Injective functor},
\index{Monofunctor}{Monofunctor},
\index{Category of! monofunctor}{Category of monofunctor}}]
Let $\A_1$ be an additive category, $\A$ an abelian category, and $\A^{\A_1}$ the category of additive functors $F:\A_1\ra\A$. An additive functor $F:\A_1\ra\A$ is an \ul{injective functor} if it is an injective object in $\A^{\A_1}$. An additive functor $F:\A_1\ra\A$ is a \ul{monofunctor} in $\A^{\A_1}$ if it maps monomorphisms to monomorphisms. Denote by $\M(\A^{\A_1})$ the full subcategory of $\A^{\A_1}$ formed by additive monofunctors. Similarly, as before, denote by $\L(\A^{\A_1})$ and $\R(\A^{\A_1})$ the full subcategory of $\A^{\A_1}$ formed by left-exact and right-exact additive functors respectively.
\end{dfn}

\begin{rmk}[\blue{Properties of the monofunctor subcategory $\M(Ab^\A)$}]
Let $\A$ be an abelian category.
\bit[leftmargin=0.7cm]
\item[(1)] If $F:\A\ra Ab$ is a monofunctor and $G\subset F$ a subfunctor (i.e., $G$ maps a morphism $f:A\ra B$ to $G(f):G(A)\ra G(B)$ such that $G(A)\subset F(A)$, $G(B)\subset F(B)$, and $G(f)=F(f)|_{G(A)}$), then $G:\A\ra Ab$ is also a monofunctor as it is evident from the following associated commutative diagram:
\[\adjustbox{scale=0.9}{\bt A\ar[d,hook,"f"] \\ B\et}~~\mapsto~~
\adjustbox{scale=0.9}{\bt
G(A)\ar[d,"G(f)"]\ar[r,hook,"\subset"] & F(A)\ar[d,hook,"F(f)"] \\
G(B) \ar[r,hook,"\subset"] & F(B)
\et}\]
\item[(2)] Since the product $\prod$ is left exact as a functor, $\M(Ab^\A)$ is closed under products, i.e., a product $\prod_{i\in I}F_i$ of monofunctors is a monofunctor.
\item[(3)] We have seen that if $F\in\Ob Ab^\A$ is a monofunctor, then so is any essential extension $u:F\ra E$. Thus, if $E$ is an injective envelope of $F$, then $E$ is an injective monofunctor, hence an exact functor by Lemma \ref{RightExInj} (i.e., the injective envelope of a monofunctor $F:\A\ra Ab$ is an exact functor).
\eit
\end{rmk}

\begin{crl}[\blue{Properties of the monofunctor subcategory $\M(Ab^\A)$}]
Let $\A$ be an abelian category. Then $\M(Ab^\A)\subset Ab^\A$ (where we know $Ab^\A$ is an abelian AB5 category) is closed under subobjects, products, and essential extensions.
\end{crl}

As in \cite{freyd1964}, we will for convenience introduce the following abstract model suggested by the above result.

\begin{dfn}[\blue{\index{Mono-subcategory}{Mono-subcategory},
\index{Mono-object}{Mono-object}}]
Let $\B$ be an abelian AB5 category. A subcategory $\M\subset\B$ is a \ul{mono-subcategory} if it is closed under subobjects, products, and essential extensions. The objects of $\M$ are called mono-objects in $\B$.
\end{dfn}

\begin{prp}[\blue{\cite[Proposition 7.21, p.141]{freyd1964}}]\label{LeFunAbProp1}
Let $\B$ be an abelian AB5 category, $B\in\Ob\B$, and $\M\subset\B$ a mono-subcategory. Then $B$ has a maximal quotient object $B\twoheadrightarrow M(B)$ with $M(B)\in\M$. (\blue{footnote}\footnote{It is clear that if $B\in\Ob\M$, then $M(B)=B$.})
\end{prp}
\begin{proof}
Let $\F:=(B_i)_{i\in I}$ be an indexing of all quotient objects of $B$ that each lie in $\Ob\M$, where $0\in\F$ and so $\F$ is nonempty.
\[\adjustbox{scale=0.9}{\bt
B_i & & B_j\\
  & \prod_i B_i\ar[ul,"p_i"']\ar[ur,"p_j"] &\\
   &\coim\pi\ar[u,hook]& \\
  & B\ar[uu,dashed, bend left=50,near end,"\pi"]\ar[u,two heads,"ck_\pi"']\ar[uuul,bend left,two heads,"\pi_i"]\ar[uuur,bend right,two heads,"\pi_j"'] &
\et}\]
Consider the coimage $ck_\pi:B\twoheadrightarrow\coim\pi\subset\prod_iB_i\in\Ob\M$. Then $\coim\pi\in\Ob\M$ since $\M$ is closed under products and subobjects. Also, by construction every quotient object of $B$ in $\Ob\M$ factors through $\coim \pi$. Thus, $\coim\pi$ is a maximal quotient object of $B$ in $\M$, with respect to the partial ordering of quotient objects $\pi_1:B\twoheadrightarrow Q_1$ and $\pi_2:B\twoheadrightarrow Q_2$ given by ``$Q_1\geq Q_2$ if $\pi_2$ factors through $\pi_1$ in the sense $\pi_2=p_{12}\circ\pi_1:B\sr{\pi_1}{\ral}Q_1\sr{p_{12}}{\ral}Q_2$ for some $p_{12}:Q_1\ra Q_2$''. (\blue{footnote}\footnote{Since $\pi$ is epic, the factoring here is unique, i.e., $p_{12}$ is unique.}). We can now set $M(B):=\coim\pi$.
\end{proof}

\begin{prp}[\blue{\cite[Proposition 7.22, p.142]{freyd1964}}]\label{LeFunAbProp2}
Let $\B$ be an abelian AB5 category, $B\in\Ob\B$, $\M\subset\B$ a mono-subcategory, and $f:B\ra C$ a morphism such that $C\in\Ob\M$. Then $f:B\ra C$ factors uniquely as $f:B\twoheadrightarrow M(B)\ra C$.
\end{prp}
\begin{proof}
Since $\coim f\subset C$ implies $\coim f\in\Ob\M$, it follows that $B\twoheadrightarrow\coim f$ (hence $f$) uniquely factors through $B\twoheadrightarrow M(B)$.
\end{proof}

\begin{crl}\label{LeFunAbProp3}
Let $\B$ be an abelian AB5 category and $\M\subset\B$ a mono-subcategory. Then any morphism $f:B\ra B'$ in $\B$ induces a unique morphism $M(f):M(B)\ra M(B')$ in $\M$ such that the following diagram commutes:
\[\adjustbox{scale=0.9}{\bt
B\ar[d,"f"]\ar[r,two heads,"e"] & M(B)\ar[d,dashed,"\exists!~M(f)"] \\
B'\ar[r,two heads,"e'"] & M(B')
\et}~~~~~~~~
\adjustbox{scale=0.9}{\bt
B\ar[d,"{f,g,f+g}"]\ar[rrr,two heads,"e"] &&& M(B)\ar[d,dashed,"{M(f),M(g),M(f+g)}"] \\
B'\ar[rrr,two heads,"e'"] &&& M(B')
\et}\]
In particular, we get a functor:
\bea
M:\B\ra\M\subset\B,~~B\sr{f}{\ral}B'~~\mapsto~~M(B)\sr{M(f)}{\ral}M(B'),\nn
\eea
which satisfies $M|_\M=id_\M:\M\ra\M$ (the identity functor) and is also \ul{additive} because
\bea
&&M(f)e=e'f,~~M(g)e=e'g,~~M(f+g)e=e'(f+g)\nn\\
&&~~\Ra~~M(f+g)e=e'(f+g)=e'f+e'g=M(f)e+M(g)e=\big(M(f)+M(g)\big)e,\nn\\
&&~~\Ra~~M(f+g)=M(f)+M(g).\nn
\eea
We also get a natural transformation $e:id_\B\ra M$ with components $\big\{B\sr{e_B}{\twoheadrightarrow}M(B)\big\}_{B\in\Ob\B}$ between the functors ~$id_\B,M:\B\ra\B$,~ where $id_\B:\B\ra\B$ is the identity functor.
\end{crl}

\begin{dfn}[\blue{ \index{Torsion! object}{Torsion object}}]
Let $\B$ be an abelian AB5 category and $\M\subset\B$ a mono-subcategory. An object $T\in\Ob\B$ is \ul{$\M$-torsion} (or torsion wrt $\M$) if the quotient object $M(T)=0$.
\end{dfn}

\begin{prp}[\blue{\cite[Proposition 7.24, p.143]{freyd1964}}]\label{LeFunAbProp4}
Let $\B$ be an abelian AB5 category with injective envelopes (i.e., for all objects), $\M\subset\B$ a mono-subcategory, and $B\in\Ob\B$. Then $\ker(B\twoheadrightarrow M(B))$ is the maximal torsion subobject of $B$.
\end{prp}
\begin{proof}
Let $T\in\Ob\B$ be any torsion object and $f:T\ra B$ any morphism. Then the following diagrams show $ef=0$, and so $\im f\subset\ker e$ (which shows every torsion subobject of $B$ is contained in $\ker e$). Also, $M(k_{e})e_1=ek_{e}=0$ implies $M(k_{e})=0$.
\[\adjustbox{scale=0.9}{\bt
T\ar[d,"f"]\ar[r,two heads,"e_1"] & M(T)=0\ar[d,dashed,"\exists!~M(f)=0"] \\
B\ar[r,two heads,"e"] & M(B)
\et}~~~~~~~~~~~~
\adjustbox{scale=0.9}{\bt
& \ker e\ar[d,hook,"k_{e}"]\ar[r,two heads,"e_1"] & M(\ker e)\ar[d,dashed,"\exists!~M(k_{e})=0"] \\
T\ar[r,"f"]\ar[ur,dashed,"f'"]\ar[rr,bend right,"0"] & B\ar[r,two heads,"e"] & M(B)
\et}\]
Consider any mono-object $B'\in\Ob\M$, any morphism $\al:\ker e\ra B'$, and the injective envelope $u:B'\hookrightarrow E$ of $B'$ (where we know $E\in\Ob\M$ since $B'\in\Ob\M$). Then we get the following commutative diagram:
\[\adjustbox{scale=0.9}{\bt
0\ar[r] & \ker e\ar[d,"\al"]\ar[r,hook,"k_e"] & B\ar[d,dashed,"\beta"]\ar[r,tail,two heads,"e"] & M(B)\ar[dl,"M(\beta)"]\ar[r] & 0 \\
        & B'\ar[r,hook,"u"] & E=M(E) &   &
\et}\]
where the morphism $\beta$ exists because $E$ is injective. $u\al=M(\beta)ek_e=0$ implies $\al=0$. In particular, if $0=\al:\ker e\ra B'$ is set to $0=e_1:\ker e\twoheadrightarrow M(\ker e)$, we see that $M(\ker e)=0$.
\end{proof}

\begin{dfn}[\blue{
\index{Pure subobject}{Pure subobject},
\index{Absolutely pure object}{Absolutely pure object},
\index{Category of! absolutely pure objects}{Category of absolutely pure objects},
\index{Reflection of an object}{Reflection of an object}
}]
Let $\B$ be an abelian AB5 category, $\M\subset\B$ a mono-subcategory, and $M\in\Ob\M$. A subobject $P\subset M$ is a \ul{pure subobject} if the exact sequence $0\ra P\hookrightarrow M\twoheadrightarrow M/P\ra 0$ lies in $\M$ (i.e., $M/P\in\Ob\M$). The object $M$ is an \ul{absolutely pure object} if $M$ is a pure subobject of every object $M'\in\M$ containing $M$ (i.e., whenever $M\subset M'\in\M$, the SES $0\ra M\hookrightarrow M'\twoheadrightarrow M'/M\ra 0$ lies in $\M$).

The full \ul{subcategory of absolutely pure objects} is denoted by $\L\subset\M\subset\B$. Given an object $M\in\M$, a morphism $r:M\ra R$ with $R\in\L$ is a \ul{reflection} of $M$ in $\L$ if every morphism $f:M\ra L$ with $L\in\L$ factors through $r$ as $f=f'\circ r$ for a unique morphism $f':R\ra L$.
\[\adjustbox{scale=0.9}{\bt
M\ar[d,"f"]\ar[rr,"r"] && R\ar[dll,dashed,"\exists!~f'"] \\
L &&
\et}\]
\end{dfn}

\begin{rmk}[\blue{An injective mono-object is absolutely pure}]
Let $\B$ be an abelian AB5 category and $\M\subset\B$ a mono-subcategory. If $E\in\Ob\M$ is injective, then $E$ is absolutely pure. Indeed if $E\subset M$ for $M\in\Ob\M$, then we know $M\cong E\oplus E'$ for a subobject $E'\subset M$, i.e., ${M\over E}\cong E'$ is a mono-object.
\end{rmk}

\begin{lmm}[\blue{\cite[Lemma 7.25, p.144]{freyd1964}}]\label{LeFunAbProp5}
Let $\B$ be an abelian AB5 category with injective envelopes, $\M\subset\B$ a mono-subcategory, and $0\ra M\ra B\ra M'\ra 0$ an exact sequence in $\B$. If $M,M'\in\M$, then $B\in\M$.
\end{lmm}
\begin{proof}
Let $u:M\hookrightarrow E$ be the injective envelope of $M$. Then we get a commutative diagram:
\[\adjustbox{scale=0.9}{\bt
0\ar[r] & M\ar[d,hook,"u"']\ar[r,hook,"f"] & B\ar[dl,dashed,"h"']\ar[d,hook,dotted,"v"]\ar[r,tail,two heads,"g"] & M'\ar[d,dashed]\ar[d,dashed,"\cong"',"w"]\ar[r] & 0 \\
0\ar[r] & E\ar[r,hook,near end,"q"] & P=E\oplus F\ar[r,tail,two heads,"p"] & F\ar[r] & 0
\et}\]
The morphism $h$ is induced because $E$ is injective. Choose a pushout $P$ (which splits as $P=E\oplus F$ since $E$ is injective), in which case $v:=qh$ is monic (as $u$ is). Then by Lemma \ref{PsPbPartII1} $g=\coker f=pv=wg$, and so $w$ is an isomorphism. That is, $B\subset E\oplus M'\in\M$. Hence $B\in\M$.
\end{proof}

\begin{lmm}[\blue{\cite[Lemma 7.26, p.144]{freyd1964}}]\label{LeFunAbProp6}
Let $\B$ be an abelian AB5 category with injective envelopes and $\M\subset\B$ a mono-subcategory. A pure subobject $P\subset A$ of an absolutely pure object $A\in\M$ is absolutely pure.
\end{lmm}
\begin{proof}
Let $P\subset M$ for a mono-object $M\in\M$ and consider a pushout $P'$ of $M\supset P\subset A$ as shown in the diagrams below. By Lemma \ref{PsPbPartII1}, the induced map $h:A/P\ra P'/M$ is an isomorphism (and the following boxed text also gives a direct argument for this):

\fbox{\begin{minipage}{17.5cm}
\footnotesize Since $P'$ is a pullback and $\big(M\sr{0}{\ral}A/P\big)\circ\al=0=c_f\circ f$, there exists $h_1:P'\ra A/P$ as in the diagram. Consider the other induced morphism $h:A/P\ra P'/M$. Then $c_{f'}=hh_1$. So, given $x:X\ra P'$, if $h_1x=0$ then $c_{f'}x=hh_1x=0$, and so there exists a unique $x':X\ra M$ such that $x=f'x'$. This shows $M\cong\ker h_1$, and so we can replace $c_{f'}$ with $h_1$ and still have an exact middle row of the diagram. By symmetry, the same argument shows we also have an exact middle column as shown.
\end{minipage}}

\[\bt
P\ar[d,hook]\ar[r,hook] & A\ar[d,dashed,hook] \\
M\ar[r,dashed,hook]& P'
\et~~~~~~
\adjustbox{scale=0.7}{\bt
       & 0\ar[d]                      & 0\ar[d]                      & 0\ar[d]    & \\
0\ar[r]& P\ar[d,hook,"\al"]\ar[r,hook,"f"]      & A\ar[d,hook,"\al'"]\ar[r,tail,two heads,"{c}_f"] & A/P\ar[d,dashed,"h"]\ar[r]  & 0 \\
0\ar[r]& M\ar[d,tail,two heads,"{c}_\al"]\ar[r,hook,"f'"] & P'\ar[dl,dashed,"h_2"description]\ar[d,tail,two heads,"c_{\al'}"]\ar[ur,dashed,"h_1"description]\ar[r,tail,two heads,"c_{f'}"']   &  P'/M\ar[d]\ar[r] & 0 \\
0\ar[r]& M/P\ar[d]\ar[r,dashed] & P'/A\ar[d]\ar[r]  &  0   &  \\
  X\ar[uur,bend left=20,near end,dashed,"x'"]\ar[uurr,bend right=20,"x"']     & 0                            & 0                            &    &
\et}~~~~=~~~~
\adjustbox{scale=0.7}{\bt
       & 0\ar[d]                      & 0\ar[d]                      & 0\ar[d]    & \\
0\ar[r]& P\ar[d,hook,"\al"]\ar[r,hook,"f"]      & A\ar[d,hook,"\al'"]\ar[r,tail,two heads,"{c}_f"] & A/P\ar[d,equal]\ar[r]  & 0 \\
0\ar[r]& M\ar[d,tail,two heads,"{c}_\al"]\ar[r,hook,"f'"] & P'\ar[d,dashed,tail,two heads,"h_2"]\ar[r,dashed,tail,two heads,"h_1"]   &  A/P\ar[d]\ar[r] & 0 \\
0\ar[r]& M/P\ar[d]\ar[r,equal] & M/P\ar[d]\ar[r]  &  0   &  \\
       & 0                            & 0                            &    &
\et}\]

Since $M,A/P\in\M$ (as $P\subset A$ is pure in $A$) it follows by the preceding lemma that $P'\in\M$. Finally, because $A$ is absolutely pure and $P'\in\M$, it also follows that $M/P\in\M$.
\end{proof}

\begin{lmm}[\blue{A mono-subfunctor of a left-exact functor is a pure $\iff$ left-exact}]\label{LeFunAbProp66}
Let $\A$ be an abelian category, $E\in Ab^\A$ a left-exact functor, and $M\subset E$ a mono-subfunctor. Then  $M$ is pure in $E$ (i.e., in the exact sequence $0\ra M\hookrightarrow E\twoheadrightarrow F:={E/M}\ra 0$, $F$ is a monofunctor) $\iff$ $M$ is left-exact.
\end{lmm}
\begin{proof}
Let $0\ra A\sr{f}{\ral} B\sr{g}{\ral}C$ be exact in $\A$. Then we get a commutative diagram in Ab:
\[\adjustbox{scale=0.8}{\bt
        & 0\ar[d] & 0\ar[d] & 0\ar[d] \\
0\ar[r] & M(A)\ar[d,hook,"\al_A"]\ar[r,hook,"M(f)"] & M(B)\ar[d,hook,"\al_B"]\ar[r,dashed,"M(g)"] & M(C)\ar[d,hook,"\al_C"] \\
0\ar[r] & E(A)\ar[d,tail,two heads,"\beta_A"]\ar[r,hook,"E(f)"] & E(B)\ar[d,tail,two heads,"\beta_B"]\ar[r,tail,"E(g)"] & E(C) \\
0\ar[r] & F(A)\ar[d]\ar[r,"F(f)"] & F(B)\ar[d] & \\
        & 0 & 0 &
\et}\]
The middle row is exact since $E$ is left-exact. It follows from Lemma \ref{PsPbPartII3} that $M$ is left exact (i.e., the top row is exact) $\iff$ $M\subset E$ is pure (i.e., $F$ is a monofunctor).

As an aside, assuming the bottom row is exact, the following boxed text gives a direct proof of the exactness of the first row, i.e., $\im M(f)=\ker M(g)$.

\fbox{\begin{minipage}{17.5cm}
\footnotesize
By the commutativity of the diagram we have the following:
{\flushleft\ul{$\im M(f)\subset\ker M(g)$}:} Indeed, $0=E(g)E(f)\al_A=E(g)\al_BM(f)=\al_CM(g)M(f)$ implies $M(g)M(f)=0$.
{\flushleft\ul{$\ker M(g)\subset\im M(f)$}:} Let $b'\in\ker M(g)$. Since $\al_CM(g)=E(g)\al_B$, we get $\ker M(g)=\al_B^{-1}(\ker E(g))$, and so $\al_B(b')\in\ker E(g)=\im E(f)$. Let $\al_B(b')=E(f)(a')$. Then
\[
0=\beta_B\al_B(b')=\beta_BE(f)(a')=F(f)\beta_A(a'),~~\Ra~~\beta_A(a')=0,~~\Ra~~a'\in\ker\beta_A=\im\al_A.
\]
Let $a'=\al_A(a'')$. Then
\[
\al_BM(f)(a'')=E(f)\al_A(a'')=E(f)(a')=\al_B(b')~~\Ra~~b'=M(f)(a'')\in\im M(f).
\]
\end{minipage}}

~
\end{proof}

\begin{thm}[\blue{\cite[Theorem 7.27, p.145]{freyd1964}}]\label{LeFunAbProp7}
Let $\A$ be an abelian category, and assume knowledge of the fact that $Ab^\A$ has injective envelopes (\blue{footnote}\footnote{We will show later in Lemma \ref{AbelLeftExLmm} that $G:=\bigoplus_{A\in\A}H^A$, with $H^A:=\Mor_\A(A,-)$, is a projective generator for $Ab^\A$. Also, we know from an earlier result that an abelian AB5 category (e.g., $Ab^\A$) with a generator has injective envelopes for its objects.}). Then a monofunctor $F\in Ab^\A$ is absolutely pure (as a mono-object in $Ab^\A$) $\iff$ it is left exact.

(That is, $\L(Ab^\A)\subset\M(Ab^\A)$ consists precisely of absolutely pure objects of $\M(Ab^\A)$.)
\end{thm}
\begin{proof}
Let $E$ be the injective envelope of $F$ in $Ab^\A$. ($\Ra$): Assume $F$ is absolutely pure.  Then $F\subset E$ is left-exact as a pure mono-subfunctor of a left-exact functor. ($\La$): Conversely, assume $F$ is left-exact (hence a monofunctor). Then $F\subset E$ is a pure mono-subfunctor since $E$ is left-exact, and hence absolutely pure since $E$ is absolutely pure (being injective).
\end{proof}

\begin{thm}[\blue{\cite[Theorem 7.28, p.146]{freyd1964}}]\label{LeFunAbProp8}
Let $\B$ be an abelian AB5 category with injective envelopes, $\L\subset\M\subset\B$ the usual subcategories, and $0\ra M\sr{f}{\ral} R\sr{g}{\ral}T\ra 0$ an exact sequence in $\B$. If $M$ is mono, $R$ absolutely pure, and $T$ torsion, then $f:M\ra R$ is a reflection of $M$ in $\L$.
\end{thm}
\begin{proof}
Let $\al:M\ra L$ be any morphism with $L\in\L$. Let $f':L\hookrightarrow E$ be an injective envelope of $L$. Then we get the following commutative diagram (where the \ul{non-unique} $\beta$ is induced by the injectivity of $E$ and the \ul{non-unique} $\gamma$ is induced by the exactness of the rows):
\[\bt
0\ar[r] & M\ar[d,"\al"']\ar[r,hook,"f"] & R\ar[dl,dotted,"\al'"']\ar[d,dashed,"\beta"]\ar[r,tail,two heads,"g"] & T\ar[d,dashed,"\gamma=0"]\ar[r] & 0\\
0\ar[r] & L\ar[r,hook,near end,"f'"] & E\ar[r,tail,two heads,"g'"] & F\ar[r] & 0\\
\et~~~~~~~~
\bt
T\ar[d,"\gamma"]\ar[r,two heads,"e"] & M(T)=0\ar[d,dashed,"\exists!~M(f)=0"] \\
F\ar[r,two heads,"id_F"] & M(F)=F
\et~~~~~~~~\gamma=0 e=0
\]
The object $E$ is mono by the definition of $\M$, and $F$ is mono since $L$ is absolutely pure. Therefore $\gamma=0$ (by the diagram on the right), and so  $g'\beta=0$. Hence, because $L\cong\ker g'$, there exists $\al':R\ra L$ (unique for $\beta$ and $\gamma$ which are however not unique) such that the diagram commutes.

To prove the uniqueness of $\al'$, consider another injective envelope $f_1:L\hookrightarrow E_1$ of $L$.
\[\bt
0\ar[r] & L\ar[r,hook,"f_1"] & E_1\ar[r,tail,two heads,"g_1"] & F_1\ar[r] & 0\\
0\ar[r] & M\ar[dd,"\al"']\ar[u,"\al"]\ar[r,hook,near start,"f"] & R\ar[ddl,dotted,"\al'-\al'_1"description]\ar[ul,dotted,"\al'_1"description]\ar[dd]\ar[u,"\beta_1"']\ar[r,tail,two heads,"g"] & T\ar[dd]\ar[u,"\gamma_1=0"']\ar[r] & 0\\
 &  &  & & \\
0\ar[r] & L\ar[r,hook,near end,"f'"']\ar[from=uurr,dashed,bend left=10,near start,"h"] & E\ar[r,tail,two heads,"g'"'] & F\ar[r] & 0
\et\]
Then $(\al'-\al'_1)f=\al'f-\al'_1f=\al-\al=0:M\ra L$, and so $\al'-\al'_1$ factors through $g:R\ra T$ (via $h:T\ra L$) as $\al'-\al'_1=hg$. Since $T$ is torsion and $L$ mono, we see as before that $h=0$. Hence $\al'=\al'_1$.
\end{proof}

\begin{thm}[\blue{\cite[Theorem 7.29, p.147]{freyd1964}}]\label{LeFunAbProp9}
Let $\B$ be an abelian AB5 category with injective envelopes and $\L\subset\M\subset\B$ the usual subcategories. For every object $M\in\M$, there exists a monomorphism $r:M\hookrightarrow R$ that is a reflection of $M$ in $\L$.
\end{thm}
\begin{proof}
Let $M\in\M$ be a mono-object. Let $u:M\hookrightarrow E$ be an injective envelope. Then we get a commutative diagram (in the order: middle row, right column, bottom row, middle column, top row):

\[\adjustbox{scale=0.8}{\bt
 & 0\ar[d] & 0\ar[d] & 0\ar[d] & \\
0\ar[r] & M\ar[d,equal]\ar[r,hook,dashed,"f"] & R\ar[d,hook,"k_{M(c_u)}"]\ar[r,dashed,"g"] & T\ar[d,hook,"k_e"]\ar[r] & 0\\
0\ar[r] & M\ar[d]\ar[r,hook,"u"] & E\ar[d,tail,two heads,"M(c_u)"]\ar[r,tail,two heads,"c_u"] & F\ar[d,tail,two heads,"e"]\ar[r] & 0\\
0\ar[r] & 0\ar[r] & M(F)\ar[r,equal] & M(F)\ar[r] & 0\\
 & 0\ar[from=u] & 0\ar[from=u] & 0\ar[from=u] & \\
\et}\]
We need to show (i) $R$ is absolutely pure and (ii) the top row is exact (and then conclude using the preceding theorem). For (i), $R$ is absolutely pure as a pure subobject of an absolutely pure object. For (ii), the top row is exact by the 3-by-3 lemma.
\end{proof}

\begin{thm}[\blue{\cite[Theorem 7.31, p.148]{freyd1964}}]\label{LeFunAbProp10}
Let $\B$ be an abelian AB5 category with injective envelopes and $\L\subset\M\subset\B$ the usual subcategories. Then (1) $\L$ is abelian and (2) every object of $\L$ has an injective envelope.
\end{thm}
\begin{proof}
{\flushleft\ul{Zero object}}: The zero object of $\B$ is inherited by $\L$.

{\flushleft\ul{Finite products/coproducts}}: Recall that the reflection $B\twoheadrightarrow M(B)\in\M$ of an object $B\in\B$ gives an additive functor $M:\B\ra\M\subset\B$ satisfying $M|_\M=id_\M$. Similarly, the reflection (from Theorem \ref{LeFunAbProp9}) $M\hookrightarrow R(M)\in\L$ of an object $M\in\M$ gives a functor
\[
R:\M\ra\L\subset\M,~~M\sr{f}{\ral}M'~~\mapsto~~R(M)\sr{R(f)}{\ral}R(M'),
\]
which satisfies $R|_\L=id_\L$ (up to isomorphism of objects), and with the diagrams
\[adjustbox{scale=0.8}{\bt
M\ar[d,"f"]\ar[r,hook,"m"] & R(M)\ar[d,dashed,"\exists!~R(f)"] \\
M'\ar[r,hook,"m'"] & R(M')
\et}~~~~~~~~
\adjustbox{scale=0.8}{\bt
M\ar[d,"{f,g,f+g}"]\ar[rrr,hook,"m"] &&& R(M)\ar[d,dashed,"{R(f),R(g),R(f+g)}"] \\
M'\ar[rrr,hook,"m'"] &&& R(M')
\et}\]
we see that the functor $R$ is additive because
\bea
&&R(f)m=m'f,~~R(g)m=m'g,~~R(f+g)m=m'(f+g)\nn\\
&&~~\Ra~~R(f+g)m=m'(f+g)=m'f+m'g=R(f)m+R(g)m=\big(R(f)+R(g)\big)m,\nn\\
&&~~\Ra~~R(f+g)=R(f)+R(g).\nn
\eea
We also get a natural transformation $m:id_\B\ra M$ with components $\big\{M\sr{m_M}{\hookrightarrow}R(M)\big\}_{M\in\Ob\M}$ between the functors $id_\M,R:\M\ra\M$ with $id_\M:\M\ra\M$ the identity functor.

In particular, since $R(M\oplus M')\cong M(M)\oplus R(M')$ for all $M,M'\in\M$, it follows that
\[
R(L\oplus L')\cong R(L)\oplus R(L')\in\L,~~~~\txt{for all}~~L,L'\in\L.
\]

{\flushleft\ul{Existence of kernels, and monomorphisms as kernels}}: By Lemma \ref{LeFunAbProp6}, $\L$ is closed under kernels, i.e., if $f:L\ra L'$, where both $L,L'\in\L$, then $\ker f\in\L$ as a pure subobject of $L$ since we have the decomposition
\[
0\ra\ker f\sr{k_f}{\hookrightarrow} L\sr{ck_f}{\twoheadrightarrow}\im f\sr{kc_f}{\hookrightarrow}L'.\nn
\]

Also, it is clear that a morphism in $\L$ is a monomorphism in $\L$ $\iff$ a monomorphism in $\B$. Let $f:L_1\hookrightarrow L_2$ be a monomorphism in $\L$. Then we get the exact sequence $0\ra L_1\sr{f}{\hookrightarrow}L_2\sr{c_f}{\twoheadrightarrow}M\ra 0$ in which $M\in\M$ since $L_1$ is absolutely pure. Therefore, we further get the following sequence with $\im f=\ker c_f=\ker(rc_f)$.
\[
0\ra L_1\sr{f}{\hookrightarrow}L_2\sr{c_f}{\twoheadrightarrow}M\sr{r}{\hookrightarrow}R(M).
\]
That is, every monomorphism in $\L$ is the kernel of a morphism in $\L$.

{\flushleft\ul{Existence of cokernels, and epimorphisms as cokernels}}:
Let $f:L_1\ra L_2$ be any morphism in $\L$. Let $L_1\sr{f}{\ral}L_2\sr{g}{\twoheadrightarrow}F\ra 0$ be exact in $\B$ (i.e., $g=\coker f$). We get a commutative diagram (in which $x$ factors through $reg$ whenever $X\in\L$):
\[\adjustbox{scale=0.8}{\bt
  &  & X\ar[from=d,dashed,"h"']\ar[r,two heads,"e_x"] & M(X)\ar[from=d,dashed,"M(h)"']\ar[r,hook,"r_x"] & R(M(X))\ar[from=d,dashed,"R(M(h))"'] \\
L_1\ar[urr,bend left,"0"]\ar[r,"f"] & L_2\ar[ur,"x"]\ar[r,tail,two heads,"g"] & F\ar[r,two heads,"e"] & M(F)\ar[r,hook,"r"] & R(M(F)) 
\et}\]
This shows the $\L$-cokernel of $f$ (i.e., the cokernel of $f:L_1\ra L_2$ in $\L$) exists and is given by
\bea
\label{AbsPurCokerEq}c_f^\L=reg:\bt[column sep=small] L_2\ar[r,two heads,"g"] & F\ar[r,two heads,"e"] & M(F)\ar[r,hook,"r"] & R(M(F)).\et~~~~\txt{($\L$-cokernel formula)}
\eea
It follows that $f:L_1\ra L_2$ is an $\L$-epimorphism (i.e., an epimorphism in $\L$) $\iff$ $0=c^\L_f=reg$, $\iff$ $e=0$, $\iff$ $M(\coker f)=M(F)=0$ (i.e., $\iff$ $\coker f$ is torsion).

Let $f:L_1\twoheadrightarrow_\L L_2$ be an $\L$-epimorphism (i.e., an epimorphism in $\L$). This is the case $\iff$ $T:=\coker f$ is torsion, and so with $M:=\im f$, we have exact sequences
\[
\bt L_1\ar[r,"f"]& L_2\ar[r,tail,two heads,"c_f"]& T\ar[r]& 0,\et~~~~~
\bt 0\ar[r]& M\ar[r,hook,"kc_f"]& L_2\ar[r,tail,two heads,"c_f"]& T\ar[r]& 0.\et
\]
By Theorem \ref{LeFunAbProp8}, $L_2=R(M)$. It follows that $f$ is an $\L$-cokernel as follows: From the diagram
\[\adjustbox{scale=0.8}{\bt
0\ar[r] & L_0\ar[r,hook,"k_f"] & L_1\ar[rr,bend left=50,"f"]\ar[r,tail,two heads,"ck_f"] & M\ar[d,equal]\ar[r,hook,"kc_f"]& L_2\ar[d,equal] & &\\
 & & 0\ar[r]& M\ar[r,hook,"kc_f"]& R(M)\ar[r,tail,two heads,"c_f"]& T\ar[r]& 0
\et}\]
if we set $F:=M=M(M)=M(F)$, $g:=ck_f$, $e:=id_M$, $r:=kc_f$, then from the $\L$-cokernel formula (\ref{AbsPurCokerEq}),
\[
c^\L_{k_f}=reg=f:L_1\sr{g}{\twoheadrightarrow}F\sr{e}{\hookrightarrow} M(F)\sr{r}{\hookrightarrow}M(R(F).
\]
That is, $f$ is the $\L$-cokernel of $k_f:L_0\ra L_1$ (a morphism in $\L$).

{\flushleft\ul{Existence of injective envelopes}}: Since $\L$-monomorphisms are $\B$-monomorphisms, if $L\in\Ob\L$ then a $\B$-injective envelope $u:L\hookrightarrow E$ (which is in $\L$) is injective in $\L$, and so is an $\L$-injective envelope.
\[\adjustbox{scale=0.8}{\bt
0\ar[r] & L_1\ar[d,"h"]\ar[rr,hook,"f"] && L_2\ar[dll,dashed,"\exists~h'"] \\
 & E & \\
\et} \qedhere
\]
\end{proof}

\begin{crl}\label{AbelLeftExCrl}
Let $\A$ be an abelian category, and assume knowledge of the fact that $Ab^\A$ has injective envelopes (\blue{footnote}\footnote{As mentioned before, we will show later in Lemma \ref{AbelLeftExLmm} that $G:=\bigoplus_{A\in\A}H^A$, with $H^A:=\Mor_\A(A,-)$, is a projective generator for $Ab^\A$. Also, we know from an earlier result that an abelian AB5 category (e.g., $Ab^\A$) with a generator has injective envelopes for its objects.}). Then the full subcategory of left-exact additive functors $\L(Ab^\A)$ is an abelian AB5 category with injective envelopes.
\end{crl}

\begin{lmm}\label{AbelLeftExLmm}
Let $\A$ be an abelian small category. The object ~{\small $G:=\bigoplus_{A\in\A}H^A$, $H^A:=\Mor_\A(A,-)$}, is (i) a projective generator in $Ab^\A$ and (ii) a generator in $\L(Ab^\A)$.
\end{lmm}
\begin{proof}
In the following steps, we obtain the desired results at steps (5) and (6).
\begin{enumerate}[leftmargin=0.7cm]
\item $Ab^\A$ is abelian (by Proposition \ref{YonProp}).
\item For $A\in\Ob\A$, the functor $H^A:=\Mor_\A(A,-)$ is a \ul{projective} object in $Ab^\A$:
Indeed, $\eta:F\ra F'$ is an epimorphism in $Ab^\A$ $\iff$ by definition $\eta_A:F(A)\ra F'(A)$ is an epimorphism (which is natural in $A\in\Ob\A$), and so by Yoneda's lemma, we get an epimorphism:
{\footnotesize\bea\bt
\Mor_{Ab^\A}(H^A,F)\ar[d,"\cong"]\ar[r,"\eta_\ast"] & \Mor_{Ab^\A}(H^A,F')\ar[d,"\cong"]\\
F(A)\ar[r,"{\eta_A}"] & F'(A)
\et\nn
\eea}

\item For $A\in\Ob\A$, the functor $H^A:=\Mor_\A(A,-)$ is a \ul{small} (hence \ul{small projective}) object in $Ab^\A$: Indeed, using Yoneda's lemma,
{\footnotesize\bea
&&\textstyle \Mor_{Ab^\A}\Big(H^A,\bigoplus\limits_{i\in I} T_i\Big)\cong \left(\bigoplus\limits_{i\in I} T_i\right)(A)=\bigoplus\limits_{i\in I} T_i(A)\cong\bigoplus\limits_{i\in I} \Mor_{Ab^\A}(H^A,T_i).\nn
\eea}
\item The coproduct $\bigoplus_{A\in\A}H^A$ is a \ul{projective} object in $Ab^\A$, as a \ul{coproduct of projective objects} $\coprod_i P_i$:
\[\adjustbox{scale=0.8}{\bt
 & P_i\ar[ddl,bend right,dotted,"\txt{($fq_i$-induced)}~h_i"']\ar[d,"q_i"] & \\
 & \coprod_iP_i\ar[dl,dashed,"h"']\ar[d,"f"] & \\
B\ar[r,two heads,"g"] & \ar[r] & 0
\et}~~~~~~~~
\adjustbox{scale=0.8}{\bt
P_i\ar[ddr,dotted,bend right,"h_i"']\ar[dr,"q_i"] &  &  P_j\ar[ddl,dotted,bend left,"h_j"]\ar[dl,"q_j"'] \\
 & \coprod_i P_i\ar[d,dashed,"h"] & \\
 & B &
\et}\]
As a side note, a \ul{product of injective objects} $\prod_iI_i$ is injective in a similar way:
\[\adjustbox{scale=0.8}{\bt
0\ar[r] & A\ar[d,"\theta"']\ar[r,hook,"f"] & B\ar[dl,dashed,"h"]\ar[ddl,dotted,bend left,"h_i~\txt{($p_i\theta$-induced)}"] \\
 & \prod_iI_i\ar[d,"p_i"']  & \\
 & I_i &
\et}~~~~~~~~
\adjustbox{scale=0.8}{\bt
I_i\ar[from=ddr,dotted,bend left,"h_i"]\ar[from=dr,"p_i"'] &  &  I_j\ar[from=ddl,dotted,bend right,"h_j"']\ar[from=dl,"p_j"] \\
 & \prod_i I_i\ar[from=d,dashed,"h"'] & \\
 & B &
\et}\]

\item The coproduct $G:=\bigoplus_{A\in\A}H^A$ is a \ul{projective generator} in $Ab^\A$. To check this, it suffices by Corollary \ref{MitRefCrl} to show that $\Mor_{Ab^\A}\left(\bigoplus_{A\in\A}H^A,-\right)$ is injective on objects $T\in Ab^\A$:
{\footnotesize\bea
&&\textstyle \Mor_{Ab^\A}\left(\bigoplus\limits_{A\in\A}H^A,T\right)\cong \prod\limits_{A\in\A}\Mor_{Ab^\A}\left(H^A,T\right)\cong \prod\limits_{A\in\A} T(A)=0\nn\\
&&\nn\\
&&~~\Ra~~T(A)=0~~\txt{for all}~~A\in\Ob\A,~~\Ra~~T=0.\nn
\eea}

\item The coproduct $G:=\bigoplus_{A\in\A}H^A$ is a \ul{generator} in $\L(Ab^\A)$: Indeed, if $0\ra B\ra B'\ra B''\ra 0$ is exact in $\A$, then $0\ra H^A(B)\ra H^A(B')\ra H^A(B'')$ is exact in $Ab$ for each $A\in\Ob\A$, and since the direct sum (as a functor) is exact on $Ab=\Integer$-mod, $0\ra \bigoplus_AH^A(B)\ra \bigoplus_AH^A(B')\ra \bigoplus_AH^A(B'')$ is exact in $Ab$. Hence $G:\A\ra Ab$ is a left-exact functor. Since $G$ is a generator in $Ab^\A$, $G$ is also a generator in $\L(Ab^\A)\subset Ab^\A$ by restriction.

\item The functor $G:=\bigoplus_{A\in\A}H^A:\A\ra Ab$ is an \ul{imbedding}: Indeed, if $0\neq B\in\Ob\A$, then $H^B(B)=\Mor_\A(B,B)\neq 0$, $\Ra$ $G(B)=\bigoplus_{A\in\A}H^A(B)=\bigoplus_{A\in\A}\Mor_\A(A,B)\neq 0$ (i.e., $G$ is injective on objects). Similarly, if $0\neq f:B\ra C$, then $0\neq f_\ast:H^B(B)\ra H^B(C)$, and so $G$ is also injective on morphisms.\qedhere
\end{enumerate}
\end{proof}

\begin{rmk}
In \ul{attempting} to see if $\bigoplus_{A\in\A}H^A$ is a small object in $Ab^\A$, we could consider the path
{\footnotesize
\bea
&&\textstyle \Mor_{Ab^\A}\Big(\bigoplus\limits_{A\in\A}H^A,\bigoplus\limits_{i\in I} T_i\Big)\cong \prod\limits_{A\in\A}\Mor_{Ab^\A}\Big(H^A,\bigoplus\limits_{i\in I} T_i\Big)\cong \prod\limits_{A\in\A}\bigoplus\limits_{i\in I}\Mor_{Ab^\A}(H^A,T_i)\nn\\
&&\textstyle~~~~\cong \prod\limits_{A\in\A}\bigoplus\limits_{i\in I}T_i(A)~\sr{\txt{\textcolor{red}{(s)}}}{\supset}~\bigoplus\limits_{i\in I}\prod\limits_{A\in\A}T_i(A)
\cong\prod\limits_{A\in\A}\Mor_{Ab^\A}(H^A,T_i)\cong \bigoplus\limits_{i\in I}\Mor_{Ab^\A}\Big(\bigoplus\limits_{A\in\A}H^A,T_i\Big),\nn
\eea}where step (s) holds because $\prod$ is always left exact while $\bigoplus$ is exact on $R$-mod, and so the inclusion

{\footnotesize\bea
\textstyle T_i(A)\hookrightarrow\bigoplus\limits_{i\in I}{T_i(A)}\subset \prod\limits_A\prod\limits_i{T_i(A)}\nn
\eea}gives an injection {\footnotesize$\prod\limits_{A\in\A}{T_i(A)}\hookrightarrow\prod\limits_{A\in\A}\bigoplus\limits_{i\in I}{T_i(A)}$},
which in turn gives an injection ~{\footnotesize $\bigoplus\limits_{i\in I}\prod\limits_{A\in\A}{T_i(A)}\hookrightarrow\prod\limits_{A\in\A}\bigoplus\limits_{i\in I}{T_i(A)}$}.~ More directly, we have an induced map as follows:
\[
\adjustbox{scale=0.7}{\bt
\prod_AT_i(A)\ar[ddr,bend right,"q'_i"]\ar[dr,"q_i"] &  & \\
        & \bigoplus_i\prod_AT_i(A)\ar[d,dashed,"c"] &  \\
     & \prod_A\bigoplus_i T_i(A) &
\et}
\]
However, step (s) is not an equality since $\prod$ and $\bigoplus$ in general do not commute on abelian groups.
\end{rmk}

\section{The Freyd-Mitchell Imbedding (An Exact Full Imbedding)}
A possible alternative way of proving this result might be found in \cite[page 16]{swan1968}. The proof of Theorem \ref{FreydMitchIT} (whose preliminary steps we have been developing throughout this chapter) is primarily based on \cite[Corollary 2.2 and Proposition 2.3, p.99]{mitchell1965}, where \cite[Corollary 2.2, p.99]{mitchell1965} is \cite[Lemma 4.1]{mitchell1964}.

\begin{lmm}[\textcolor{blue}{\cite[Lemma 4.3]{mitchell1964}, \cite[Theorem 7.33, p.149]{freyd1964}}]\label{RpFnExLmm}
Let $\A$ be an abelian small category. The representation functor $H:\A\ra\L(Ab^\A)$, $A\sr{f}{\ral}B$ $\mapsto$ $H^A\sr{H^f}{\ral}H^B$ is an exact fully faithful cofunctor.
\end{lmm}
\begin{proof}
Recall from Corollary \ref{AbelLeftExCrl} that $\L(Ab^\A)$ is abelian, and so it is sensible to speak of the ``exactness of $H:\A\ra\L(Ab^\A)$''. By Yoneda's lemma, we know $H$ is fully faithful. Thus, it remains to show $H$ is exact. Let $0\ra A\ra B\ra C\ra 0$ be exact in $\A$. We need to show $0\ra H^C\ra H^B\ra H^A\ra 0$ is exact, which is the case $\iff$ there is an injective cogenerator $G'\in \L(Ab^\A)$ such that the sequence
\bea
0\ra \Mor_{Ab^\A}(H^A,G')\ra \Mor_{Ab^\A}(H^B,G')\ra \Mor_{Ab^\A}(H^C,G')\ra 0\nn
\eea
is exact in Ab (or equivalently, by Yoneda's lemma, the sequence $0\ra G'(A)\ra G'(B)\ra G'(C)\ra 0$ is exact in Ab, i.e., $G':\A\ra Ab$ is an \ul{exact functor}). Thus, we only need to show that $\L(Ab^\A)$ has an \ul{exact injective cogenerator} $G'$.

By Lemma \ref{RightExInj}, if $\L(Ab^\A)$ has an injective cogenerator, then it is automatically exact. Hence, it remains only to show $\L(Ab^\A)$ has an \ul{injective cogenerator}. Since $\L(Ab^\A)$ is AB5 by
Corollary \ref{AbelLeftExCrl}, it suffices by Corollary \ref{EssExtLmm9} to show that $\L(Ab^\A)$ has a generator. By Lemma \ref{AbelLeftExLmm}, $\L(Ab^\A)$ has a generator.
\end{proof}

\begin{thm}[\textcolor{blue}{\index{Freyd-Mitchell imbedding theorem}{Freyd-Mitchell imbedding theorem}: \cite[Theorem 4.4]{mitchell1964}, \cite[Theorem 7.2, p.151]{mitchell1965}}]\label{FreydMitchIT}
If $\A$ is an abelian small category, there exists an exact fully faithful functor $\A\ra\txt{R-mod}$, for some ring $R$.
\end{thm}
\begin{proof}
Let $Ab^\A$ be the category of additive functors $F:\A\ra Ab$, and $\L(Ab^\A)\subset Ab^\A$ the category of left-exact additive functors. By Lemma \ref{RpFnExLmm}, we have a contravariant exact full imbedding $H:\A\ra \L(Ab^\A)$. This in turn gives a covariant exact full imbedding
\bea
H^\ast:=OP\circ H:\A\sr{H}{\ral}\L(Ab^\A)\sr{OP}{\ral}\L(Ab^{\A^{op}}),\nn
\eea
where $\L(Ab^{\A^{op}})$ is now an abelian AB5 category with a \ul{projective generator} (the ``reflection'' $G:=OP(G')$ of the \ul{injective cogenerator} $G'$ from $\L(Ab^\A)$). Since $\A$ is small, by replacing the projective generator $G$ of $\L(Ab^{\A^{op}})$ with $G^{\oplus I}$ for some set $I$, we may assume that every object in ~$\im H^\ast\subset \L(Ab^{\A^{op}})$~ is finitely generated wrt $G^{\oplus I}$, making the category $\Big(\L(Ab^{\A^{op}}),G^{\oplus I}\Big)$ equivalent to some module category (as in Lemma \ref{MitThmLmm2}) via an equivalence of categories
\bea
E:\Big(\L(Ab^{\A^{op}}),G^{\oplus I}\Big)\ra \txt{mod-}R.\nn
\eea
The desired exact full imbedding ~$T:\A\ra\txt{mod-}R$ ~can now be obtained as the following composition:
\bea
T:=E\circ H^\ast=E\circ OP\circ H:\A\sr{H}{\ral}\L(Ab^\A)\sr{OP}{\ral}\Big(\L(Ab^{\A^{op}}),G^{\oplus I}\Big)\sr{E}{\ral} \txt{mod-}R~\cong~R^{op}\txt{-mod},\nn
\eea
where ~$R^{op}:=(R,\cdot_{op})$ ~is the ring $R$ wrt a new product ~$\cdot_{op}:R\times R\ra R$ ~given by
\bea
a\cdot_{op}b:=ba,~~~~\txt{for all}~~~~a,b\in R.\nn
\eea
\end{proof}

\begin{rmks}[\textcolor{blue}{Importance of the Freyd-Mitchell imbedding theorem}]\label{FreydMitchRmk}~
\begin{enumerate}[leftmargin=0.7cm]
\item Let $\C$ be a category whose objects $\Ob\C$ are a class of sets. Then $\C$ is a union of small categories: Indeed, for any \ul{set of objects} $\U\subset\Ob\C$, we have the full subcategory $\C_\U\subset\C$ with objects $\Ob\C_\U:=\U$, and morphisms $\Mor_{\C_\U}(C,D):=\Mor_\C(C,D)$, for all $C,D\in\U$. Therefore, we can write ~$\C=\bigcup\{\C_{\U}~|~\txt{set of objects}~\U\subset\Ob\C\}$.
\item Let $\A$ be an abelian category. The verification of a \emph{practical statement} about systems/diagrams in $\A$ typically involves only a set of objects $\U\subset\Ob\A$. Consequently, we can consider an exact full imbedding $F:\A_\U\hookrightarrow R_\U\txt{-mod}$ for some ring $R_\U$, verify a \emph{practical statement} in $R_\U$-mod, and then return to $\A_\U\subset\A$.
\item Consider a result about systems/diagrams that holds in a general $R$-module category. If (i) the result is non-global in the sense that any instance of its statement involves objects/morphisms from only a small subcategory of $R$-modules (hence does not require $R$-modules to form a proper class of sets as they do) and (ii) the result is preserved by equivalences of categories, then the result also holds in a general abelian category. Indeed if such a result were to fail in a general abelian category, then because of the imbedding the result could not have held in a general module category in the first place.
\end{enumerate}
\end{rmks}

%% file: parts/AlgebraCat/AlgebraCatS8.tex
\chapter{Localization and Localization Functors}\label{AlgebraCatS8}

\section{System Representation Schemes: Motivation for the Localization Operation}
\begin{dfn}[\textcolor{blue}{\index{Naive system representation scheme}{Naive system representation scheme}}]
Given a category $\C$, a naive representation scheme for systems in $\C$ might be a triple $(\I,F,\D)$ consisting of (i) an \ul{auxiliary category} $\I$, (ii) a relatively \ul{``flexible/symmetric/simple'' category} $\D$, and (iii) a \ul{regulated transfer functor} $F:\C^\I\ra\D^\I$. (\blue{footnote}\footnote{For example, a functor $F:\A^\Integer\ra(R\txt{-mod})^\Integer$ for easier understanding of the homology $H:\A_0^\Integer\ra\A_0^\Integer$ of complexes in an abelian category $\A$ in terms of the homology $F(H):(R\txt{-mod})_0^\Integer\ra(R\txt{-mod})_0^\Integer$. An ideal situation happens when $F$ is an exact full imbedding.}). The functor $F$ is regulated by choosing it such that for any relevant system $T\in \C^\I$, the image system $F(T)\in\D^\I$ is in a sense still a good approximation of $T$.
\end{dfn}

In such a scheme, given a system $T:\I\ra\C$, its image $F(T):F(\I)\ra F(\C)$ is relatively simple in the sense that the objects and transition morphisms of the image system $F(T):F(\I)\ra F(\C)$ are easier to describe/understand compared to those of the original system $T:\I\ra\C$. Consequently, if $\tau:T\ra T'$ is a morphism of systems, the $F$-induced morphism of systems $F(\tau):F(T)\ra F(T')$ is also easier to describe.

We may be able to compare two systems $T,T'$ in $\C^\I$ by means of an \emph{equivalence} $F(T)\sim F(T')$ in $\D^\I$. For example, if we manage to discover (while working in $\D^\I$) that $\Mor_{\D^\I}(F(T),F(T'))$ contains an \ul{isomorphism} $\vphi$ in $\D^\I$, then we have also discovered that $\Mor_{\C^\I}(T,T')$ contains an interesting class of morphisms $F^{-1}(\vphi)$, the ``\emph{relative content/size}'' of which provides a measure of ``\emph{difference}'' or ``\emph{disagreement}'' between $T$ and $T'$ in $\C^\I$. In order for such a measurement to be sufficiently stable/robust (hence useful/practical), it is essential to choose the functor $F$ carefully (say in terms of a ``\emph{limit-like}'' universal property).

Note that in the above situation, instead of waiting to discover or stumble upon an isomorphism in $\D^\I$ as we did, we also could have intentionally chosen $F$ in order that a given morphism $s_\vphi\in F^{-1}(\vphi)\subset\Mor\C^\I$ becames the (previously ``discovered'') isomorphism $\vphi\in \Mor_{\D^\I}(F(T),F(T'))\subset\Mor\D^\I$.

Based on the above discussion (especially concerning the existence/creation of a potentially useful \ul{isomorphism} in $\D^\I$), we can ask the following precise question (in which the role of the naive functor $F:\C^\I\ra\D^\I$ above is now played by a precise functor $L:\C^\I\ra\C^\I[S^{-1}]$ ).
\begin{question}[\textcolor{blue}{
\index{Optimized system representation scheme}{Optimized system representation scheme},
\index{Optimized category}{Optimized category}}]\label{CatLocQuest}
Let $\C$ be a category and $\I$ an auxiliary category (for expressing$\slash$parameterizing systems in $\C$ as functors $\I\ra\C$). Given a class of morphisms of systems $S\subset\Mor\C^\I$, what is the ''\ul{smallest}'' (or \ul{optimized}) category $\C^\I[S^{-1}]$, up to isomorphism (or just up to equivalence of categories), with a functor $L:\C^\I\ra \C^\I[S^{-1}]$ such that for each $s\in S$, $L(s)$ is an isomorphism in $\C^\I[S^{-1}]$? (\blue{footnote}\footnote{As mentioned earlier, if $T,T'\in\C^\I$, $s\in \Mor_{\C^\I}(T,T')\cap S$, and the isomorphism $s\in \Mor_{\C^\I[S^{-1}]}(L(T),L(T'))$, then the ``\emph{relative content/size}'' of $L^{-1}(s)\subset\Mor\C^\I$ provides a measure of ``\emph{difference}'' or ``\emph{disagreement}'' between $T$ and $T'$ in $\C^\I$.})
\end{question}

As before, our underlying wish here is for the new category $\C^\I[S^{-1}]$ to be ``richer'' in structure (i.e., structurally ``more interesting'') than $\C^\I$ so as to allow for greater ``flexibility/symmetry/simplicity'' in the description of the original systems $\C^\I\subset\C^\I[S^{-1}]$. For example, the numbers $0$ and $2$ are each more difficult to describe/distinguish in $\Natural$ than in $\Integer$ or $\Rational$ where there is greater freedom of expression (by which we can alternatively represent these numbers as $0=-0=-1+1=4-4=\cdots$ and $2=-1\times -2=-5+7=10-8={6\over 1}\times{1\over 3}={-4\over -2}={3\over 2}+{1\over 2}=\cdots$).

\begin{dfn}[\textcolor{blue}{
\index{Localization}{Localization},
\index{Localizing functor}{Localizing functor}}]
In Question \ref{CatLocQuest}, the optimized category $\C^\I[S^{-1}]$ is known as the \ul{localization} of {\small $\C^\I$ by $S\subset\Mor(\C^\I)$}, and the functor {\small $L:\C^\I\ra\C^\I[S^{-1}]$} is the \ul{localizing functor}.
\end{dfn}

Given a morphism class $S\subset\Mor(\C^\I)$, if any of the collections $\dom S:=\{T\in\C^\I~|~T=\dom s~\txt{for some}~s\in S\}$ and $\cod S:=\{T'\in\C^\I~|~T'=\cod s~\txt{for some}~s\in S\}$ is interesting, then the localizing functor $L$ can be a useful means of comparing a system $T\in\dom S$ against a system $T'\in\cod S$. Moreover, given morphism classes $S_1,S_2\subset\Mor(\C^\I)$, we may also compare (i) $S_1$ and $S_2$, or (ii) the two system collections $\dom S_1$ and $\dom S_2$ (resp. $\cod S_1$ and $\cod S_2$), by comparing the corresponding localizations $\C^\I[S_1^{-1}]$ and $\C^\I[S_2^{-1}]$.

In any specific application of the localization, in order to avoid potential technical obstructions as well as obtain useful results, one may need to impose certain conditions on the morphism class $S\subset\Mor(\C^\I)$. For example, among other things, it is often desirable that $S$ be closed under composition in $\C^\I$. A typical such practical $S$ will be called a \index{Localizing class}{\ul{localizing class}}, even though we will show that the localization exists for any $S$.

\section{Localization of a Category}
\begin{dfn}[\textcolor{blue}{\index{Localization of! a category}{Localization of a category}}]
Let $\C$ be a category and $S\subset\Mor\C$ any class of morphisms. A localization of $\C$ by $S$ is a category $\C[S^{-1}]$ (or $S^{-1}\C$) and a functor $L:\C\ra\C[S^{-1}]$ such that the following are true.
\bit
\item[(a)] For any morphism $s\in S$, its image $L(s)$ is an isomorphism in $\C[S^{-1}]$.
\item[(b)] For any functor $F:\C\ra \D$ mapping elements of $S$ to isomorphisms, there exists a unique functor $F':\C[S^{-1}]\ra \D$ such that $F=F'\circ L$. (I.e., $F$ uniquely factorizes through $L$).
\eit
\begin{figure}[H]
\centering
\begin{tikzcd}
 \C\ar[d,"F"']\ar[rr,"L"] && \C[S^{-1}]\ar[dll,dashed,"F'"]   \\
  \D  &&
\end{tikzcd}
\label{dg14s}
\end{figure}
\end{dfn}
Due to property (b) above, if a localization exists then it is unique up to isomorphism (an equivalence of categories) in an appropriate category of categories. As indicated above, depending on convenience, we will sometimes write the localization $\C[S^{-1}]$ as $S^{-1}\C$, and this will often be the case especially when dealing with certain applications involving (commutative) rings and modules where details about the construction of the localization are not strictly relevant.

\begin{thm}[\textcolor{blue}{\index{Existence of localization}{Existence of the localization of a category}: \cite[pp 144-145]{gelfand-manin2010}}]\label{LocExistThm}
Given a category $\C$ and any class of morphisms $S\subset Mor~\C$, the localization (of $\C$ by $S$) $L:\C\ra \C[S^{-1}]$ exists.
\end{thm}
\begin{proof}
For each $s\in S$, define \ul{any new} arrow $\cod(s)\sr{x_s}{\ral}\dom(s)$, and let $x_S:=\{x_s:s\in S\}$ be the collection of these arrows. Next, consider the category $\C'$ with $Ob~\C':=Ob~\C$, and for $A,B\in Ob~\C'$,
{\small
\begin{align}
&\Mor_{\C'}(A,B):={\left\{\substack{\txt{finite}\\\txt{words}}~~q_0\ast q_1\ast\cdots \ast q_n:A=C_n\sr{q_n}{\ral}\cdots\sr{q_1}{\ral}C_0\sr{q_0}{\ral}C_{-1}=B~\big|~q_i\in x_S\cup \Mor~\C\right\}\over\left\{s\ast x_s=id_{\cod(s)},~x_s\ast s=id_{\dom(s)},~f\ast g=fg~\big|~s\in S,~f,g\in \Mor~\C\right\}}\nn\\
&~~=\left\{[f_0\ast x_{s_1}\ast f_1\ast\cdots\ast x_{s_n}\ast f_n]:A=Y_n\sr{f_n}{\ral}X_n\sr{x_{s_n}}{\ral}\cdots ~\sr{f_1}{\ral}X_1\sr{x_{s_1}}{\ral}Y_0\sr{f_0}{\ral}X_0=B~\big|~f_i\in \Mor~\C,~s_i\in S\right\},\nn
\end{align}}(where $[w]$ denotes the equivalence class of the word $w$) such that (i) the composition of morphisms ~$\circ:\Mor_{\C'}(A,B)\times \Mor_{\C'}(B,C)\ra \Mor_{\C'}(A,C)$~ is given (when possible) by the concatenation of words in the obvious way, i.e.,
{\footnotesize\[
\big[f_0\ast x_{s_1}\ast f_1\ast\cdots\ast x_{s_n}\ast f_n\big]\circ\big[g_0\ast x_{t_1}\ast g_1\ast\cdots\ast x_{t_m}\ast g_m\big]:=\big[f_0\ast x_{s_1}\ast f_1\ast\cdots\ast x_{s_n}\ast f_n\ast g_0\ast x_{t_1}\ast g_1\ast\cdots\ast x_{t_m}\ast g_m\big],\nn\nn
\]}and (ii) the identity morphism in $\Mor_{\C'}(A,A)$ is $id_A^{\C'}:=[id_A]$.
It follows immediately that
\bea
[f_0\ast x_{s_1}\ast f_1\ast\cdots\ast  x_{s_n}\ast f_n]=[f_0]\circ[x_{s_1}]\circ[f_1]\circ\cdots\circ [x_{s_n}]\circ[f_n].\nn
\eea

Consider the map $L:\C\ra\C'$ given by $L(A):=A$ for $A\in Ob~C$ and $L(f):=[f]$ for $f\in Mor~\C$, i.e.,
{\small\bea
L:\C\ra\C',~~\Big(A\sr{f}{\ral}B\Big)\longmapsto \Big(A\sr{[f]}{\ral}B\Big).\nn
\eea}$L$ is a functor, and for any $s\in S$, $L(s)=[s]$ is an isomorphism (with inverse $[s]^{-1}=[x_s]$), since
\bea
L(s)[x_s]=[s][x_s]=[s\ast x_s]=[id_{\cod(s)}]~~~~\txt{and}~~~~[x_s]L(s)=[x_s][s]=[x_s\ast s]=[id_{\dom(s)}].\nn
\eea
\bea\bt
 \C\ar[d,"F"]\ar[rr,"L"] && \C'\ar[dll,dashed,"F'"]   \\
  \D  &&
\et\nn
\eea
Given any functor $F:\C\ra \D$ such that $F(S)$ consists of isomorphisms, define a functor $F':\C'\ra\D$ by
{\small\bea\bt
F': \Big(A\ar[rrrr,"{[f_0][s_1]^{-1}[f_1]\cdots[s_n]^{-1}[f_n]}"] &&&& B\Big)\ar[r,mapsto] &  \Big(F(A) \ar[rrrrrr,"{F(f_0)F(s_1)^{-1}F(f_1)\cdots F(s_n)^{-1}F(f_n)}"] &&&&&& F(B)\Big).
\et\nn
\eea} Then ~{\small$F'L\Big(A\sr{f}{\ral}B\Big)=F'\Big(A\sr{[f]}{\ral}B\Big)=\Big(F(A)\sr{F(f)}{\ral}F(B)\Big)=F\Big(A\sr{f}{\ral}B\Big)$}. ~Hence ~$\C'\cong\C[S^{-1}]$.
\end{proof}

\begin{notation}
Let $\C$ be a category and $S\subset\Mor\C$ any class of morphisms. In the proof of Theorem \ref{LocExistThm}, we have shown that the localization $\C[S^{-1}]$ exists, with $\Ob\C[S^{-1}]=\Ob\C$, and for any $A,B\in\Ob\C[S^{-1}]$,
\bea
\Mor_{\C[S^{-1}]}(A,B)=\left\{\txt{finite compositions}~[f_0][s_1]^{-1}[f_1]\cdots[s_n]^{-1}[f_n]:A\ra B~|~f_i\in\Mor\C,~s_i\in S\right\}.\nn
\eea
For any $f\in\Mor\C$ and $s\in S$, we will denote the equivalence classes $[f]$, $[s]$ by $f$, $s$ respectively. Thus, an arbitrary morphism $[f_0][s_1]^{-1}[f_1]\cdots[s_n]^{-1}[f_n]\in\Mor\C[S^{-1}]$ will be written as $f_0s_1^{-1}f_1\cdots s_n^{-1}f_n$.
\end{notation}

\begin{questions}
Let $\C$ be a category and $S\subset\Mor\C$ a class of morphisms.
\begin{enumerate}
\item If $\C$ is additive, does it follow that $\C[S^{-1}]$ is additive? (Hint: What happens if $0\in S$?)
\item If $\C$ is abelian, does it follow that $\C[S^{-1}]$ is abelian?
\end{enumerate}
\end{questions}

\begin{dfn}[\textcolor{blue}{\index{Derived! category}{Derived category}, \index{Homotopy! category}{Homotopy category}}]\label{DcatHTYcat}
Let $\A$ be an abelian category. The \ul{derived category} $\D\A$ of $\A$ is the localization $D:\A_0^\Integer\ra\D\A$ of $\A_0^\Integer$ by quasiisomorpisms (i.e., $\D\A:=\A_0^\Integer[S_q^{-1}]$ is the smallest category containing $\A_0^\Integer$ such that quasiisomorphisms $S_q\subset\Mor\A_0^\Integer$ become isomorphisms).

The \ul{homotopy category} $\H\A$ of $\A$ is the quotient category ${\A_0^\Integer\over\sim}$, where $\sim$ denotes homotopy of chain morphisms, i.e., for any chain morphisms $f,g\in \Mor\A_0^\Integer$, we have by definition $f\sim g$ iff $f\simeq g$ (where $\simeq$ indeed defines a congruence relation on $\A_0^\Integer$ because if $f\simeq g$ and $f'\simeq g'$ then $f f'\simeq gg'$). (\blue{footnote}\footnote{Note that some quais-iso's in $\A_0^\Integer$ (such as homotopy equivalences) become iso's in the homotopy category $\H\A$, just as in the derived category $\D\A$.})
\end{dfn}

\begin{questions}
Let $\A$ be an abelian category.
\begin{enumerate}
\item Is the homotopy category $\H\A$ (i) additive? (ii) abelian?
\item Is the derived category $\D\A$ (i) additive? (ii) abelian?
\item Does there exist an equivalence of categories ~$E:\H\A\ra\D\A$?
\item Does there exist an equivalence of categories ~$E:\H\A\left[(\txt{quasi-iso's})^{-1}\right]\ra\D\A$?
\end{enumerate}
\end{questions}

Some of these questions will be answered later, inlcuding showing that $\H\A$ is ``triangulated'' (which is a lot more than just additive). However, it is also known in the literature that $\H\A$ is not abelian, but we will not be needing this particular result, and so will not provide a direct proof.

\section{Localization Functors: Localization of Rings and Modules as Categories}
\begin{rmk}[\textcolor{blue}{
\index{Ring! as a category}{A ring as a category},
\index{localizing set in a ring}{Localizing set in a ring},
\index{Localization of! a ring}{Localization of a ring}}]
Given a ring $R$, we can view $R$ as the small \magenta{pre-additive} category $\C(R)$ with one object $R$, and morphisms {\small$Hom_{\C(R)}(R,R)\approx R$} the elements $r\in R$ viewed as the left multiplication maps {\small$\{l_r:R\ra R$, $r'\mapsto rr'\}_{r\in R}$}, with composition of morphisms given by the usual composition of maps and the identity morphism $id_R:=l_{1_R}$. We have the composition and addition formulas
\[
l_{r_1}\circ l_{r_2}=l_{r_1r_2},~~~~l_{r_1}+l_{r_2}=l_{r_1+r_2},
\]
which show that ring homomorphisms $f:R\ra R'$ are precisely \magenta{additive} functors $F_f:\C(R)\ra\C(R')$, i.e., functors $F_f:\C(R)\ra\C(R')$ that are \magenta{additive} in the sense they are homomorphisms of abelian groups $F_f:Hom_{\C(R)}(R,R)\ra Hom_{\C(R')}(R',R')$, or equivalently, $F_f$ satisfies
\[
F_f(l_r)=l_{f(r)}~~~~\txt{for all}~~~~r\in R.
\]

A set $S\subset R$ is called a \ul{localizing set} (or \ul{multiplicative set}: \blue{footnote}\footnote{If adopted, this would be our second notion of a multiplicative set, which is specific to subsets of a ring, and so different from the earlier one.}) if (i) $0\not\in S$, (ii) $1\in S$, and (iii) $s,s'\in S$ $\Ra$ $ss'\in S$, i.e., $S$ is closed under multiplication. With the categorical description of $R$, it is clear (by Theorem \ref{LocExistThm}) that given a localizing set $S\subset R$ (or equivalently, $l_S:=\{l_s:s\in S\}\subset\Mor\C(R)$ ), we get the category localization $L:\C(R)\ra\C(R)[l_S^{-1}]$, in terms of which we define the \ul{localization} $h:R\ra R[S^{-1}]$ of $R$ by $S$ using the following rule: Up to equivalence of categories,
\[
F_h:=L:\C(R)\ra\C(R)[l_S^{-1}]~~~~\txt{and}~~~~\C(R[S^{-1}]):=\C(R)[l_S^{-1}],\nn
\]
where, up to isomorphism of rings, $R[S^{-1}]$ is the ``smallest'' ring containing $R$ such that for every $s\in S$, $h(s)$ is a unit, or equivalently,
\[
L(l_s)=F_h(l_s)=l_{h(s)}\in Hom_{\C(R)[l_S^{-1}]}(R,R)\eqv Hom_{\C(R[S^{-1}])}(R[S^{-1}],R[S^{-1}])
\]
is an isomorphism. Explicitly (in the notation from Theorem \ref{LocExistThm}),
{\footnotesize
\begin{align}
&R[S^{-1}]:=\Mor_{\C(R)[l_S^{-1}]}(R,R)\cong \left\{l_{r_0}l_{s_1}^{-1}l_{r_1}\cdots l_{s_n}^{-1}l_{r_n}~|~r_i\in R,~s_j\in S\right\}\cong \left\{l_{r_0}l_{s_1^{-1}}l_{r_1}\cdots l_{s_n^{-1}}l_{r_n}~|~r_i\in R,~s_j\in S\right\}\nn\\
&~~~~\cong \left\{l_{r_0s_1^{-1}r_1\cdots s_n^{-1}r_n}~|~r_i\in R,~s_j\in S\right\}\cong\left\{r_0s_1^{-1}r_1\cdots s_n^{-1}r_n~|~r_i\in R,~s_j\in S\right\}.\nn
\end{align}}By construction, it is clear that $R[S^{-1}]$ is an \ul{$R$-bimodule}. Also, as an equivalent of the localizing functor $L=F_h:\C(R)\ra\C(R)[l_S^{-1}]\eqv \C(R[S^{-1}])$, the localizing map $h:R\ra R[S^{-1}]$ is a ring homomorphism.
\end{rmk}

In order to connect with relatively familiar terminology, we will now restate the definition of the ring localization.

\begin{dfn}[\textcolor{blue}{\index{Localization of! a ring}{The ring localization} $R[S^{-1}]$: Non-categorical description}]
Let $R$ be a ring and $S\subset R$ a localizing set. The localization of $R$ by $S$ is a ring $R[S^{-1}]$, along with a ring homomorphism $h:R\ra R[S^{-1}]$, such that the following hold:
\bit
\item[(a)] $h(S)\subset U(R[S^{-1}])$, i.e., the elements of $S$ become units in $R[S^{-1}]$.
\item[(b)] Given any ring homomorphism $f:R\ra A$ with $\vphi(S)\subset U(A)$, there exists a unique ring homomorphism $f':R[S^{-1}]\ra A$ such that $f=f'\circ h$.
\eit
\bea
\bt
R\ar[d,"f"]\ar[rr,"h"] && R[S^{-1}]\ar[dll,dashed,"\exists!~f'"]\\
A &&
\et\nn
\eea
\end{dfn}

\begin{rmk}[\textcolor{blue}{\index{Module! as a category}{A module as a category}, \index{Localization of! a module}{Localization of a module}}]
Given an $R$-module $M$, we can view $M$ as the small \magenta{pseudo-preadditive} category $\C_R(M)$ with one object $M$, and morphisms {\small$Hom_{\C_R(M)}(M,M)\approx M\times R$} the elements $(m,r)\in M\times R$ viewed as the maps {\small$\{l_{m,r}:M\ra M,~m_1\mapsto m+rm_1\}_{m\in M,r\in R}$}, with composition of morphisms given by the usual composition of maps and the identity morphisms $id_M:=l_{0,1}$, where $1:=1_R$. We have the composition, factorization, and addition formulas
\[
l_{m_1,r_1}\circ l_{m_2,r_2}=l_{m_1+r_1m_2,~r_1r_2},~~~~l_{m,r}=l_{m,1}\circ l_{0,r},~~~~l_{m_1,r_1}+l_{m_2,r_2}=l_{m_1+m_2,~r_1+r_2},
\]
which show that any two maps $f:M\ra M'$ and $g:R\ra R'$ satisfying
(\blue{footnote}\footnote{However, if there is no $m'\in M'$ such that $Ann_R(m'):=\{r'\in R':r'm'=0\}=\{0\}$, then a functor $F:\C_R(M)\ra\C_{R'}(M')$ that is a homomorphism of abelian groups $F:Hom_{\C_R(M)}(M,M)\ra Hom_{\C_{R'}(M')}(M',M')$ might not take the form $F=F_{f,g}$ for $f,g$ satisfying (1),(2),(3) above. Indeed, for any $(m,r)\in M\times R$,
\begin{align}
&F(l_{m,r})=l_{f(m),g(r)}=l_{f_1(m),g_1(r)}~~\iff~~f(m)+g(r)m'=f_1(m)+g_1(r)m'~~~~\txt{for all}~~m'\in M',\nn\\
&~~~~\iff~~f=f_1~~(\txt{using $m'=0$})~~~~\txt{and}~~~~\big(g(r)-g_1(r)\big)m'=0~~~~\txt{for all}~~m'\in M',\nn
\end{align}
where the functor and abelian group homomorphism properties of $F$ require $f,g$ (as well as $f_1,g_1$) to satisfy (1),(2),(3) above.})
\bit
\item[(1)] $f:M\ra M'$ is a homomorphism of abelian groups,
\item[(2)] $g:R\ra R'$ is a ring homomorphism,
\item[(3)] $f(rm)=g(r)f(m)$ ~for all ~$r\in R$, $m\in M$,
\eit
together give an \magenta{additive} functor $F_{f,g}:\C_R(M)\ra\C_{R'}(M')$, i.e., a functor $F_{f,g}:\C_R(M)\ra\C_{R'}(M')$ that is \magenta{additive} in the sense it is a homomorphism of abelian groups
\[
F_{f,g}:Hom_{\C_R(M)}(M,M)\ra Hom_{\C_{R'}(M')}(M',M'),
\]
or equivalently, $F_{f,g}$ satisfies
\[
F_{f,g}(l_{m,r})=l_{f(m),g(r)},~~~~\txt{for all}~~~~m\in M,~~r\in R.
\]
In particular, if $R':=R[S^{-1}]$ (for a localizing set $S\subset R$) and $g:=h:R\ra R[S^{-1}]$ is the ring localization, then any $f:{}_RM\ra{}_{R[S^{-1}]}M'$ that is $h$-compatible, in the sense the pair $(f,h)$ satisfies (1),(2),(3), must of course be an $R$-homomorphism. In what follows, such an $f$ is a candidate for the module localization $h_M$.

With the categorical description of $M$, it is clear (by Theorem \ref{LocExistThm}) that given a localizing set $S\subset R$ (or equivalently, $l_S:=\{l_{0,s}:s\in S\}\subset\Mor\C_R(M)$ ) and the associated ring localization $h:R\ra R[S^{-1}]$, we get the category localization $L_M:\C_R(M)\ra\C_R(M)[l_S^{-1}]$, in terms of which we define the \ul{localization} $h_M:M\ra M[S^{-1}]$ of $M$ by $S$ using the following rule: Up to equivalence of categories,
\[
F_{h_M,h}:=L_M:\C_R(M)\ra\C_R(M)[l_S^{-1}]~~~~\txt{and}~~~~\C_{R[S^{-1}]}(M[S^{-1}]):=\C_R(M)[l_S^{-1}],
\]
where, up to isomorphism of modules, $M[S^{-1}]$ is the ``smallest'' $R$-module containing $M$ such that for every $s\in S$, $h(s)$ is a unit, or equivalently,
\[
L_M(l_{0,s})=F_{h_M,h}(l_{0,s})=l_{0,h(s)}\in Hom_{\C_R(M)[l_S^{-1}]}(M,M)\eqv Hom_{\C_{R[S^{-1}]}(M[S^{-1}])}(M[S^{-1}],M[S^{-1}])
\]
is an isomorphism. Explicitly (in the notation from Theorem \ref{LocExistThm}),
\begin{align}
&Hom_{\C_R(M)[l_S^{-1}]}(M,M)\cong\big\{l_{m_0,r_0}l_{0,s_1}^{-1}l_{m_1,r_1}\cdots l_{0,s_n}^{-1}l_{m_n,r_n}~|~m_i\in M,~r_i\in R,~s_j\in S\big\}\nn\\
&~~~~\cong\big\{l_{m_0,r_0}l_{0,s_1^{-1}}l_{m_1,r_1}\cdots l_{0,s_n^{-1}}l_{m_n,r_n}~|~m_i\in M,~r_i\in R,~s_j\in S\big\}\nn\\
&~~~~=\big\{l_{g(m_0,r_0,m_1,r_1,s_1^{-1},...,m_n,r_n,s_n^{-1}),~r_0s_1^{-1}r_1,\cdots s_n^{-1}r_n}~|~m_i\in M,~r_i\in R,~s_j\in S\big\}\nn\\
&~~~~=:\big\{l_{x,t}~|~x\in M[S^{-1}],~t\in R[S^{-1}]\big\}~~~~~~~~\big(\approx~~ M[S^{-1}]\times R[S^{-1}]~\big)\nn\\
&~~~~=Hom_{\C_{R[S^{-1}]}(M[S^{-1}])}\big(M[S^{-1}],M[S^{-1}]\big),\nn
\end{align}
where in the last two steps, we have defined $M[S^{-1}]$ using the equivalence ~{\small $\C_{R[S^{-1}]}(M[S^{-1}])\eqv\C_R(M)[l_S^{-1}]$}.

By construction, it is clear that $M[S^{-1}]$ is an \ul{$R[S^{-1}]$-module}. Also, as an equivalent of the localizing functor $L_M=F_{h_M,h}:\C_R(M)\ra\C_R(M)[l_S^{-1}]\eqv\C_{R[S^{-1}]}(M[S^{-1}])$, the localizing map $h_M:M\ra M[S^{-1}]$ is an $R$-homomorphism.
\end{rmk}

Using a specialized tool such as the tensor product of $R$-modules $\otimes_R$, it is possible \ul{under certain conditions} to conveniently describe the localization $M[S^{-1}]$ of the $R$-module $M$ in terms of that of the ring $R$ (as we will do in Lemma \ref{ModLocLmm}). This will \ul{suffice when $R$ is commutative} as proved in Theorem \ref{LocModCmmR} (which uses Theorem \ref{LocCommRing}). In order to connect with relatively familiar terminology, we will now restate the definition of the module localization.

\begin{dfn}[\textcolor{blue}{\index{Localization of! a module}{The module localization} $M[S^{-1}]$: Non-categorical description}]
Let $R$ be a ring, $S\subset R$ a localizing set, and $M$ an $R$-module. The localization of $M$ by $S$ is an $R$-module $M[S^{-1}]$, along with an $R$-homomorphism $h_M:M\ra M[S^{-1}]$, such that the following hold:
\bit
\item[(a)] $M[S^{-1}]$ is an $R[S^{-1}]$-module.
\item[(b)] Given any $R$-homomorphism $f:M\ra X$ with $X$ an $R[S^{-1}]$-module, there exists a unique $R[S^{-1}]$-homomorphism $f':M[S^{-1}]\ra X$ such that $f=f'\circ h_M$.
\eit
\bea\bt
\label{ModLocEqn0}M\ar[d,"f"]\ar[rr,"h_M"] && M[S^{-1}]\ar[dll,dashed,"\exists!~f'"]\\
X &&
\et\eea
\end{dfn}

\begin{lmm}[\textcolor{blue}{Formal structure of the module localization}]\label{ModLocLmm}
Let $R$ be a ring, $S\subset R$ a localizing set, and $h:R\ra R[S^{-1}]$ the localization of $R$ by $S$. Let $M$ be an $R$-module. Then the localization $h_M:M\ra M[S^{-1}]$ of $M$ by $S$ satisfies
\bea
h_M:M\sr{\cong}{\ral}R\otimes_RM\sr{h\otimes id_M}{\ral}R[S^{-1}]\otimes_RM\cong M[S^{-1}],~~m\mapsto 1_{R[S^{-1}]}\otimes m.\nn
\eea
(This is only a formal expression because the actual meaning of the tensor product $\otimes_R:(R\txt{-mod})^2\ra R\txt{-mod}$ is still pending.)
\end{lmm}
\begin{proof}
It is clear that $h_M$ is an R-homomorphism, and that $R[S^{-1}]\otimes_RM$ is a $R[S^{-1}]$-module. Consider the diagram in (\ref{ModLocEqn0}), replacing $M[S^{-1}]$ with $R[S^{-1}]\otimes_RM$, i.e.,
\[\bt
M\ar[d,"f"]\ar[rr,"h_M"] && R[S^{-1}]\otimes_RM\ar[dll,dashed,"\exists!~f'"]\\
X &&
\et\]
Then given an $R$-homomorphism $f:M\ra X$, where $X$ is an $R[S^{-1}]$-module, define $f':R[S^{-1}]\otimes_RM\ra X,~\sum_ic_i\otimes m_i\mapsto \sum_ic_if(m_i)$. Then $f=f'\circ h_M$. Hence $h_M:M\ra R[S^{-1}]\otimes_RM,~m\mapsto 1\otimes m$ is a localization of $M$ by $S$.
\end{proof}

Note that in the above representation $M[S^{-1}]\cong R[S^{-1}]\otimes_RM$ of the module localization, $R$ does not have to be commutative (which is not an issue with the $R$-bimodule $R[S^{-1}]$ ), but the meaning of $\otimes_R$ depends to the category of $R$-modules (left $R$-modules, right $R$-modules, or $R$-bimodules) $M$ being considered. In other words, we need to specify a category of $R$-modules $\C_R\ni M$ (e.g., a category of $R$-bimodules) with a well defined tensor product $\otimes_R:\C_R\times\C_R\ra\C_R$. (See Remark \ref{TnsPrdRmod}).

\begin{dfn}[\textcolor{blue}{\index{Localization functor! (Module)}{Module localization functor}}]
Let $R$ be a ring and $S\subset R$ a localizing set. Then we get a \ul{module localization functor} ~$F_S:=R[S^{-1}]\otimes_R-:\txt{$R$-mod}\ra\txt{$R$-mod}$~ given by
{\small\bea
\Big(M\sr{f}{\ral}M'\Big)\longmapsto\Big(M[S^{-1}]\sr{id\otimes f}{\ral}M'[S^{-1}]\Big),~~~~M[S^{-1}]:=R[S^{-1}]\otimes_RM.\nn
\eea}
\end{dfn}
By the properties of tensor, $F_S$ is a right-exact functor (and becomes an exact functor when $R$ is commutative, as we will soon see).

\begin{dfn}[\textcolor{blue}{\index{Localization functor! (Ring)}{Ring localization functor}}]
For each ring $R$, fix a localizing set $S_R\subset R$ (elements we want to become invertible), and let $\S:=\{S_R\}_{\txt{$R$ a ring}}$ be the collection of these localizing sets. Let $\S$-Rings be the subcategory of Rings with objects all rings, and morphisms
\bea
Hom_{\txt{$\S$-Rings}}\big(R,R'\big)\subset Hom_{\txt{Rings}}(R,R'),\nn
\eea
consisting of those ring homomorphisms $f:R\ra R'$ satisfying $f(S_R)\subset S_{R'}$ (which implies $0\not\in f(S_R)$, and so $f(S_R)\subset R'$ is itself a localizing set). Then we get a \ul{ring localization functor}
\bea
F:\txt{$\S$-Rings}\subset\txt{Rings}\ra\txt{Rings},~~\Big(R\sr{f}{\ral}R'\Big)\mapsto\Big(R[S_R^{-1}]\sr{f_{S_R}}{\ral}R'[S_{R'}^{-1}]\Big),\nn
\eea
where $f_{S_R}$ is the map induced by the universal property of the localization as in the following diagram:
\bc\bt
R\ar[d,"h"]\ar[rr,"f"] && R'\ar[d,"h'"]\ar[rr,"f'"] && R''\ar[d,"h''"]\\
R[S^{-1}]\ar[rr,dashed,"\exists!~f_S"] && R'[S'{}^{-1}]\ar[rr,dashed,"\exists!~f_{S'}"] && R''[S''{}^{-1}]
\et\ec
where $S:=S_R$, $S':=S_{R'}$, and $S'':=S_{R''}$.
\end{dfn}
Note that the condition $f(S_R)\subset S_{R'}$ above is required for consistency because any ring homomorphism (and in particular, the induced ring homomorphism $f_S$) maps invertible elements to invertible elements. The ring localization functor can be generalized to a category localization functor as follows.

\begin{dfn}[\textcolor{blue}{\index{Localization functor! (Category)}{Category localization functor}}]
Let $\C\C$ be a category of categories. For each category $\C\in\Ob\C\C$, fix a class of morphisms $S_\C\subset\Mor\C$ (morphisms we want to become invertible), and let $\S:=\{S_\C\}_{\C\in\Ob\C\C}$ be the collection of these morphism classes. Let $\S$-$\C\C$ be the subcategory of $\C\C$ with objects $\Ob\C\C$ (i.e., the same objects as $\C\C$), and morphisms
\bea
\Mor_{\S\txt{-}\C\C}\big(\C,\C'\big)\subset \Mor_{\C\C}(\C,\C')\nn
\eea
consisting of those functors $F:\C\ra\C'$ satisfying $F(S_\C)\subset S_{\C'}$. Then we get a \ul{category localization functor}
\bea
\L:\S\txt{-}\C\C\subset\C\C\ra\C\C,~~\Big(\C\sr{F}{\ral}\C'\Big)\mapsto\Big(\C[S_\C^{-1}]\sr{F_{S_\C}}{\ral}\C'[S_{\C'}^{-1}]\Big),\nn
\eea
where $F_{S_\C}$ is the map induced by the universal property of the localization as in the following diagram:
\bc\bt
\C\ar[d,"L"]\ar[rr,"F"] && \C'\ar[d,"L'"]\ar[rr,"F'"] && \C''\ar[d,"L''"]\\
\C[S^{-1}]\ar[rr,dashed,"\exists!~F_S"] && \C'[S'{}^{-1}]\ar[rr,dashed,"\exists!~F_{S'}"] && \C''[S''{}^{-1}]
\et\ec
where $S:=S_\C$, $S':=S_{\C'}$, and $S'':=S_{\C''}$.
\end{dfn}

\section{Localization of Commutative Rings}
\begin{thm}[\textcolor{blue}{\index{Localization of! a commutative ring}{The localization of a commutative ring}}]\label{LocCommRing}
Let $R$ be a commutative ring, and $S\subset R$ a localizing set. Then the localization (of $R$ by $S$) ~$h:R\ra R[S^{-1}]$~ satisfies
\bea
h=\pi|_R:R\hookrightarrow R[x_S]\sr{\pi}{\ral} {R[x_S]\over J}\cong R[S^{-1}],~~r\mapsto r+J,~~~~\txt{for indeterminates}~~x_S:=\{x_s~|~s\in S\},\nn
\eea
where $J:=\langle\{sx_s-1:s\in S\}\rangle$ is the ideal of $R[x_S]$ generated by $\{sx_s-1:s\in S\}$. Moreover, its kernel is
\bea
\ker h=\big\{r\in R~|~rs=0~~\txt{for some}~~s\in S\big\}.\nn
\eea
\end{thm}
\begin{proof}
{\flushleft \ul{The localization formula}}: Choose a set of indeterminates $x_S=\{x_s:s\in S\}$ and form the polynomial ring $R[x_S]$. Let
\bea
J:=\big\langle\{sx_s-1|s\in S\}\big\rangle \lhd R[x_S]\nn
\eea
be the ideal of $R[x_S]$ generated by the set $\{sx_s-1|s\in S\}$, and let $\ol{x}_S:=\{x_s+J~|~s\in S\}$. Define
\bea
&&\textstyle R':={R[x_S]\over J}=R[\{x_s+J~|~s\in S\}]=R[\ol{x}_S],\nn\\
&&\textstyle h:=\pi|_R:R\hookrightarrow R[x_S]\sr{\pi}{\ral} {R[x_S]\over J},~r\mapsto r+J.\nn
\eea
We will now show that $(R',h)$ is a localization of $R$ by $S$. Observe that for each $s\in S$,
\bea
h(s)(x_s+J)=(s+J)(x_s+J)=sx_s+J=1+J=1_{R'},~~\Ra~~h(S)\subset U(R').\nn
\eea
Also, with ~$h(S)^{-1}:=\{h(s)^{-1}=x_s+J:s\in S\}=\ol{x}_S$,~ we have
\bea
\textstyle R'={R[x_S]\over J}=R[\ol{x}_S]=R[h(S)^{-1}]=h(R)+\langle h(S)^{-1}\rangle\nn
\eea
Given any ring homomorphism $\theta:R\ra A$ such that $\theta(S)\subset U(A)$, consider its extension
\bea
\textstyle\theta_1:R[x_S]\ra A~~~~\txt{given by}~~~~\theta_1(x_s):=\theta(s)^{-1}~~\txt{and}~~\theta_1\big|_R:=\theta.\nn
\eea
\bea\bt
R~\ar[rrrr,bend left=25,"h"]\ar[dd,"\theta"']\ar[rr,hook,"i"] && R[x_S]\ar[ddll,dashed,"\theta_1"']\ar[rr,two heads,"\pi"] && R':={R[x_S]\over J}\ar[ddllll,dashed,"\theta'"]\\
 && && \\
A
\et\nn
\eea
Because $J\subset \ker \theta_1$, i.e., $\theta_1(J)=0$, we get the following ring homomorphism satisfying ~$\theta'\circ h=\theta$.
\bea
\textstyle\theta':{R[x_S]\over J}\ra A,~~f+J\mapsto \theta_1(f).\nn
\eea
Hence $R'\cong R[S^{-1}]$. It remains to find the kernel $\ker h$.
{\flushleft \ul{The localization kernel}}: Let $K:=\big\{r\in R~|~rs=0~~\txt{for some}~~s\in S\big\}$.
\bit[leftmargin=0.5cm]
\item \ul{$\ker h\supset K$}:~ If $rs=0$ ($s\in S$), then~ $0=h(rs)=h(r)h(s),~\Ra~~h(r)=0~\txt{since}~~h(s)\in U(S^{-1}R)$. Thus,
\bea
\ker h\supset\{r\in R:~rs=0~~\txt{for some}~~s\in S\}.\nn
\eea
\item \ul{$\ker h\subset K$}:~ Let $r\in\ker h$. Then $h(r)=r+J=J$ implies $r\in J$, and so for some $t\in\Natural\backslash\{0\}$,
\bea
\textstyle r=\sum_{i=1}^t(s_ix_{s_i}-1)f_i,~~~~~~\txt{for some}~~f_i\in R[x_S]=\big(R[x_{S'}]\big)[x_{s_1},\cdots,x_{s_t}],~~S':=S\backslash\{s_1,...,s_t\}.\nn
\eea
We have $f_i=\sum_{i_1,...,i_t\in\Natural}a^i_{i_1,...,i_t}x_{s_1}^{i_1}x_{s_2}^{i_2}\cdots x_{s_t}^{i_t}$, where of course only finitely many $a^i_{i_1,...,i_t}\in R[x_{S'}]$ are nonzero, and as usual $\Natural=\{0,1,2,\cdots\}$. To simplify notation, let
\begin{align}
&\textstyle\al:=(i_1,...,i_t)\in\Natural^t,~~|\al|:=i_1+\cdots+i_k,~~x^{|\al|}:=x_{s_1}^{i_1}x_{s_2}^{i_2}\cdots x_{s_t}^{i_t},~~\langle a^i_{\vec{n}},x^{\vec{n}}\rangle:=\sum\limits_{|\al|=n} a_\al^ix^{|\al|},\nn\\
&\textstyle f_i=\sum\limits_\al a_\al^ix^{|\al|}=\sum\limits_{n=0}^N\sum\limits_{|\al|=n} a_\al^ix^{|\al|}=\sum\limits_{n=0}^N\langle a^i_{\vec{n}},x^{\vec{n}}\rangle,~~a^i_{\vec{n}}:=\left(a^i_\al\right)_{|\al|=n},~~x^{\vec{n}}:=\left(x^{|\al|}\right)_{|\al|=n}.\nn
\end{align}
Then with $s:=s_1s_2\cdots s_t\in S$ and $\wh{s}_i:=s_1\cdots s_{i-1}s_{i+1}\cdots s_t\in S$, we see that $s^{N+1}r=0$ as follows:
{\small\begin{align}
&\textstyle r=\sum\limits_{i=1}^t(s_ix_{s_i}-1)\sum\limits_{n=0}^N\langle a^i_{\vec{n}},x^{\vec{n}}\rangle=\sum\limits_{n=0}^N\left\langle \sum\limits_{i=1}^t(s_ix_{s_i}-1)a^i_{\vec{n}},x^{\vec{n}}\right\rangle=\sum\limits_{n=0}^N\left\langle -\sum\limits_{i=1}^ta^i_{\vec{n}}+\sum\limits_{i=1}^ts_ia^i_{\vec{n}}x_{s_i},x^{\vec{n}}\right\rangle\nn\\
&\textstyle~~~~=\sum\limits_{n=0}^N\left[-\left\langle a_{\vec{n}},x^{\vec{n}}\right\rangle+\sum\limits_{i=1}^ts_i\left\langle a^i_{\vec{n}},x^{\vec{n}+1_i}\right\rangle\right],~~~~a_{\vec{n}}:=\sum\limits_{i=1}^ta^i_{\vec{n}},~~~~x^{\vec{n}+1_i}:=x^{\vec{n}}x_{s_i}=\left(x^{|\al|}x_{s_i}\right)_{|\al|=n},\nn
\nn\\
&\textstyle~~~~=-a_{\vec{0}}-\sum\limits_{n=1}^N\left\langle a_{\vec{n}},x^{\vec{n}}\right\rangle+\sum\limits_{n=0}^N\sum\limits_{i=1}^ts_i\left\langle a^i_{\vec{n}},x^{\vec{n}+1_i}\right\rangle=-a_{\vec{0}}-\sum\limits_{n=0}^{N-1}\left\langle a_{\overrightarrow{n+1}},x^{\overrightarrow{n+1}}\right\rangle+\sum\limits_{n=0}^N\sum\limits_{i=1}^ts_i\left\langle a^i_{\vec{n}},x^{\vec{n}+1_i}\right\rangle\nn\\
&\textstyle~~~~=-a_{\vec{0}}+\sum\limits_{n=0}^{N-1}\left(\sum\limits_{i=1}^ts_i\left\langle a^i_{\vec{n}},x^{\vec{n}+1_i}\right\rangle-\left\langle a_{\overrightarrow{n+1}},x^{\overrightarrow{n+1}}\right\rangle\right)+
\sum\limits_{i=1}^ts_i\left\langle a^i_{\vec{N}},x^{\vec{N}+1_i}\right\rangle,\nn
\end{align}}and so upon comparing coefficients (for all $0,1,...,N$ polynomial degrees), we get
{\begin{align}
&\textstyle r=-a_{\vec{0}}=-\sum\limits_{j=1}^ta^j_{\vec{0}},~~~~\left\{s_ia_{\vec{n}}^i=a_{\vec{n}+1_i}=\sum\limits_{j=1}^ta^j_{\vec{n}+1_i}\right\}_{n=0}^{N-1}~~\txt{and}~~s_ia^i_{\vec{N}}=0,~~\txt{for all}~~i=1,...,t,\nn\\
&\textstyle~~\txt{where}~~s_ia^i_{\vec{N}}=0~~\txt{for all}~~i=1,...,t~~\Ra~~sa^i_{\vec{N}}=0~~\txt{for all}~~i=1,...,t,\nn\\
&\textstyle~~\Longrightarrow~~~~s^2a_{\overrightarrow{N-1}}^i=s\wh{s}_is_ia_{\overrightarrow{N-1}}^i=s\wh{s}_ia_{\overrightarrow{N-1}+1_i}=\wh{s}_i\sum\limits_{j=1}^tsa^j_{\overrightarrow{N-1}+1_i}=\wh{s}_i\sum\limits_{j=1}^t0=0,~~\txt{for all}~~i=1,...,t,\nn\\
&\textstyle~~\Longrightarrow~~\txt{(by induction on $k$)}~~s^{k+1}a_{\overrightarrow{N-k}}^i=0,~~\txt{for all}~~i=1,...,t,~~1\leq k\leq N,\nn\\
&\textstyle~~\Longrightarrow~~\txt{(with $k=N$)}~~s^{N+1}a_{\vec{0}}^i=0,~~\txt{for all}~~i=1,...,t,\nn\\
&\textstyle~~\Longrightarrow~~~~s^{N+1}r=-\sum\limits_{j=1}^ts^{N+1}a_{\vec{0}}^j=0,~~\Ra~~\ker h\subset K.\nn\qedhere
\end{align}}

\eit
\end{proof}

It might be possible for the proof of the containment $\ker h\subset K$ above to be presented with less elaborate indexing technology. Nevertheless, dealing with such indexing is itself a useful exercise.

\begin{notation*}[\textcolor{blue}{Fractions convention for the commutative ring localization $R[S^{-1}]$, The $R$-module property of $R[S^{-1}]$, \index{Ring! of quotients}{The ring of quotients $Q(R)$}}]
Let $R$ be a commutative ring and $S\subset R$ a localizing set. Then from Theorem \ref{LocCommRing} the localization $R[S^{-1}]$ is given by
\bea
h=\pi|_R:R\hookrightarrow R[x_S]\sr{\pi}{\ral} {R[x_S]\over J}\cong R[S^{-1}],~~r\mapsto r+J.\nn
\eea
Observe that because (i) $S$ is closed under multiplication and (ii) $h(s)=s+J$ (for each $s\in S$) has inverse $x_s+J$ in $R[S^{-1}]$, it follows that any element of $R[S^{-1}]$ can be written in the form
\bea
\textstyle {r\over s}={r\over 1}{1\over s}~:=~rx_s+J=(r+J)(x_s+J)=h(r)h(s)^{-1},~~\txt{for some}~~r\in R,~~s\in S.\nn
\eea
Consequently, if we introduce an equivalence relation $\sim$ on $S\times R$ by
{\small\bea
\textstyle(s,r)\sim (s',r')~~~\txt{if}~~~{r\over s}={r'\over s'}~~(\txt{i.e., if}~~s'r-sr'\in\ker h,~~\txt{or}~~s''(s'r-sr')=0~\txt{for some}~s''\in S),\nn
\eea}then, in terms of the equivalence classes $[(s,r)]={r\over s}$, we can simply let $R[S^{-1}]$ be the quotient set
\bea
\textstyle R[S^{-1}]:={S\times R\over\sim}:=\left\{{r\over s}~|~r\in R,~s\in S\right\}\nn
\eea
as a ring with addition, multiplication, unity, and zero (in terms of $1=1_R$ and $0=0_R$) given by
\bea
\textstyle{r\over s}+{r'\over s'}:={rs'+r's\over ss'},~~~~{r\over s}{r'\over s'}:={rr'\over ss'},~~~~1_{R[S^{-1}]}:={1\over 1},~~~~0_{R[S^{-1}]}:={0\over 1}.\nn
\eea
In this convention, the \ul{$R$-module} property of $R[S^{-1}]$ is given by the scaling map ~$R\times R[S^{-1}]\ra R[S^{-1}]$,
\bea
\textstyle\left(a,{r\over s}\right)\mapsto a{r\over s}:={a\over 1}{r\over s}={ar\over s}.\nn
\eea
Finally, the \ul{ring of quotients} of $R$, written $Q(R)$, is the localization of $R$ by all nonzero elements of $R$, i.e.,
\bea
\textstyle Q(R):=R[(R\backslash 0)^{-1}]=\left\{{a\over b}~|~a,b\in R,~b\neq 0\right\}.\nn
\eea
\end{notation*}

\begin{rmk}[\textcolor{blue}{Localization of a commutative ring as a colimit of smaller localizations}]
Let $R$ be a commutative ring and $S\subset R$ a localizing set. For each $s\in S$, let ~$\langle s\rangle=\{s^i~|~i\geq 0\}\subset S$~ be the localizing subset of $S$ generated by $s$. Similarly, for each finite set $F=\{s_1,...,s_t\}\subset S$, let ~$\langle F\rangle=\{s_1^{n_1}s_2^{n_2}\cdots s_t^{n_t}~|~n_i\geq 0\}=\langle s_1\rangle\langle s_2\rangle\cdots\langle s_t\rangle\subset S$~ be the localizing subset of $S$ generated by $F$. Then, as shown in the colimit setup diagram below,
\bea
\textstyle R[S^{-1}]\cong\varinjlim_{\txt{finite}~F\subset S}~R\left[\langle F\rangle^{-1}\right],\nn
\eea
where for any localizing subsets $A\subset B\subset S$, we have a ring homomorphism
\[
\textstyle q_{AB}:R[A^{-1}]={R[x_A]\over J_A}\ra R[B^{-1}]={R[x_B]\over J_B},~~rx_a+J_A\mapsto rx_a+J_B,
\]
with $J_A:=\big\langle\{ax_a-1:a\in A\}\big\rangle$ an ideal of $R[x_A]$, and $J_B:=\big\langle\{bx_b-1:b\in B\}\big\rangle$ an ideal of $R[x_B]$. Given ring homomorphisms $f_A:R[A^{-1}]\ra T$ (for each relevant $A\subset S$), we get a unique ring homomorphism $f:R[S^{-1}]\ra T$ (as defined below) satisfying $fq_{AS}=f_A$ (for each relevant $A$).

\bea\bt
{R[A^{-1}]}\ar[ddr,bend right,"f_A"']\ar[dr,"q_{AS}"]\ar[rr,bend left=20,"q_{AB}","A\subset B"'] &&  {R[B^{-1}]}\ar[ddl,bend left,"f_B"]\ar[dl,"q_{BS}"']\\
 & {R[S^{-1}]}\ar[d,dashed,"f"]&\\
 &  T  &
\et~~~~~~~~\shortstack[l]{$A:=\langle F_A\rangle,~~B:=\langle F_B\rangle$, ~~for finite ~~$F_A,F_B\subset S$. \\~\\ $f(rx_s+J_S):=f_A(rx_s+J_A)~~\txt{if}~~s\in A$.}\nn
\eea
\end{rmk}

\begin{dfn}[\textcolor{blue}{Recall: \index{Localization of! a module}{The module localization $M[S^{-1}]$}}]
Let $R$ be a ring, $S\subset R$ a localizing set, and $M$ an $R$-module. The localization of $M$ by $S$ is an $R$-module $M[S^{-1}]$, along with an $R$-homomorphism $h_M:M\ra M[S^{-1}]$, such that the following hold:
\bit
\item[(a)] $M[S^{-1}]$ is an $R[S^{-1}]$-module.
\item[(b)] Given any $R$-homomorphism $f:M\ra X$ with $X$ an $R[S^{-1}]$-module, there exists a unique $R[S^{-1}]$-homomorphism $f':M[S^{-1}]\ra X$ such that $f=f'\circ h_M$.
\eit
\bea
\label{ModLocEqn}\bt
M\ar[d,"f"]\ar[rr,"h_M"] && M[S^{-1}]\ar[dll,dashed,"\exists!~f'"]\\
X &&
\et
\eea
\end{dfn}

\section{Localization of Modules over Commutative Rings}
\begin{thm}[\textcolor{blue}{\index{Localization of! a (commutative R)-module}{The localization of a module over a commutative ring}}]\label{LocModCmmR}
Let $R$ be a commutative ring, $S\subset R$ a localizing set, and $h:R\ra R[S^{-1}]$ the localization of $R$ by $S$. Let $M$ be an $R$-module. Then the localization $h_M:M\ra M[S^{-1}]$ of $M$ by $S$ satisfies
\bea
h_M:M\sr{\cong}{\ral}R\otimes_RM\sr{h\otimes id_M}{\ral}R[S^{-1}]\otimes_RM\cong M[S^{-1}],~~m\mapsto 1_{R[S^{-1}]}\otimes m,\nn
\eea
and its kernel is ~$\ker h_M=\{m\in M:~sm=0~~\txt{for some}~~s\in S\}\cong \ker(h\otimes id_M)$.
\end{thm}
\begin{proof}
{\flushleft\ul{The localization formula}}: It is clear that $h_M$ is an R-homomorphism, and that $R[S^{-1}]\otimes_RM$ is an $R[S^{-1}]$-module. Consider the diagram in (\ref{ModLocEqn}), replacing $M[S^{-1}]$ with $R[S^{-1}]\otimes_RM$.
\[\bt
M\ar[d,"f"]\ar[rr,"h_M"] && R[S^{-1}]\otimes_RM\ar[dll,dashed,"\exists!~f'"]\\
X &&
\et\]
Given an R-homomorphism $f:M\ra X$, where $X$ is an $R[S^{-1}]$-module, let
\[
\textstyle f':R[S^{-1}]\otimes_RM\ra X,~{r\over s}\otimes m\mapsto {r\over s}f(m).
\]
Then $f=f'\circ h_M$. Hence $h_M:M\ra R[S^{-1}]\otimes_RM,~m\mapsto 1\otimes m$ is a localization of $M$ by $S$.

{\flushleft\ul{The localization kernel}}: To find the kernel of $h_M$, let $K:=\{m\in M:sm=0~\txt{for some}~s\in S\}$.

{\flushleft$\bullet$} \ul{$\ker h_M\supset K$}:~ If $sm=0$ ($s\in S$), then~ {\small $h_M(m)={1\over 1}\otimes m={s\over s}\otimes m={1\over s}\otimes sm={1\over s}\otimes 0=0$}, and so $\ker h_M\supset K$.

{\flushleft $\bullet$} \ul{$\ker h_M\subset K$}:~ Let $m\in\ker h_M$. Then with $R[S^{-1}]={R[x_S]\over J}$, where $J:=\big\langle\{ sx_s-1:s\in S\}\big\rangle\subset R[x_S]$,
\begin{align}
&\textstyle J\otimes_RM=0_{M[S^{-1}]}=h_M(m)={1\over 1}\otimes m=(1+J)\otimes m=1\otimes m+J\otimes m~~\Ra~~1\otimes m\in J\otimes_R M,\nn\\
&\textstyle~~\Ra~~1\otimes m=\sum\limits_{i=1}^t(s_ix_{s_i}-1)f_i\otimes m_i~~~~\txt{for some}~~f_i\in R[x_S]=\overbrace{\big(R[x_{S'}]\big)[x_{s_1},\cdots,x_{s_t}]}^{S':=S\backslash\{s_1,...,s_t\}},~~m_i\in M.\nn
\end{align}
Thus, with $s:=s_1s_2\cdots s_t\in S$, via exactly the same procedure as in the proof of Theorem \ref{LocCommRing}, we get
\bea
0=s^l(1\otimes m)=1\otimes s^lm,~~\txt{for some}~~l\geq 1,~~\sr{(\txt{iso})}{\Longrightarrow}~~s^lm=0,~~\Ra~~\ker h_M\subset K,\nn
\eea
where step (iso) holds because $R\otimes_RM\cong M$ via the map $f:R\otimes_RM\ra M$ given by $f(r\otimes m):=rm$ for $r\in R,m\in M$ (with inverse $f^{-1}:M\ra R\otimes_RM$ given by $f^{-1}(m):=1\otimes m$ for $m\in M$).
\end{proof}

\begin{notation*}[\textcolor{blue}{Fractions convention for the commutative module localization $M[S^{-1}]$}]
Let $R$ be a commutative ring, $S\subset R$ a localizing set, and $M$ an $R$-module. Then from Theorem \ref{LocModCmmR} the localization $M[S^{-1}]$ is given by
\bea
h_M:M\sr{\cong}{\ral} R\otimes_RM\sr{h\otimes id}{\ral}R[S^{-1}]\otimes_RM\cong M[S^{-1}],~~m\mapsto 1\otimes m.\nn
\eea
Observe that because we can write $R[S^{-1}]=\left\{{r\over s}~|~r\in R,~s\in S\right\}$, we have the following:
\bit
\item[(i)] ${r\over s}\otimes m={1\over s}\otimes rm$,~ for all ${r\over s}\in R[S^{-1}]$, $m\in M$.
\item[(ii)] ${1\over s}\otimes m+{1\over s'}\otimes m'={1\over ss'}(1\otimes s'm+1\otimes sm')={1\over ss'}\otimes(s'm+sm')$,~ for all $s,s'\in S$, $m,m'\in M$.
\eit
Consequently, with ${m\over s}:={1\over s}\otimes m$, if we introduce an equivalence relation $\sim$ on $S\times M$ by
{\small\bea
\textstyle(s,m)\sim (s',m')~~~\txt{if}~~~{m\over s}={m'\over s'}~~(\txt{i.e., if}~~s'm-sm'\in\ker h_M,~~\txt{or}~~s''(s'm-sm')=0~\txt{for some}~s''\in S),\nn
\eea}then, in terms of the equivalence classes $[(s,m)]={m\over s}$, we can simply let $M[S^{-1}]$ be the quotient set
\bea
\textstyle M[S^{-1}]:={S\times M\over\sim}:=\left\{{m\over s}~|~m\in M,~s\in S\right\}\nn
\eea
as an $R[S^{-1}]$-module with addition, scalar multiplication, and zero (wrt $1=1_R$ and $0=0_M$) given by
\bea
\textstyle{m\over s}+{m'\over s'}:={s'm+sm'\over ss'},~~~~{r\over s}\cdot{m\over t}:={rm\over st},~~~~0_{M[S^{-1}]}:={0\over 1}.\nn
\eea
\end{notation*}

\begin{crl}
Let $R$ be a commutative ring and $h:R\ra S^{-1}R$ a localization. If $F$ is a free (resp. $P$ is a projective) $R$-module then $S^{-1}F$ is a free (resp. $S^{-1}P$ is a projective) $S^{-1}R$-module. (\blue{footnote}\footnote{Here, as usual, $S^{-1}R$ denotes $R[S^{-1}]$, and $S^{-1}M$ denotes $M[S^{-1}]$ for any $R$-module $M$.})
\end{crl}
\begin{proof}
We have {\small $S^{-1}F=S^{-1}R\otimes_RF=S^{-1}R\otimes_R\left(\bigoplus_{i\in I}R\right)\cong \bigoplus_{i\in I}(S^{-1}R\otimes_RR)=\bigoplus_{i\in I}S^{-1}R$}, since the tensor product preserves direct sums. Hence $S^{-1}F$ is free. Also, we have
{\small\[Hom_{S^{-1}R}(S^{-1}P,-)=Hom_{S^{-1}R}(S^{-1}R\otimes_RP,-)\cong Hom_R\big(P,Hom_{S^{-1}R}(S^{-1}R,-)\big)\cong Hom_R\big(P,-\big)\qedhere
\]}
\end{proof}

\begin{thm}[\textcolor{blue}{\index{Flatness of the localization}{Flatness of the commutative ring localization $R[S^{-1}]$}}]
Let $R$ be a commutative ring and $S\subset R$  a localizing set. Then the module-localization functor
\[
S^{-1}R\otimes_R-:R\txt{-mod}\ra R\txt{-mod},~~A\sr{f}{\ral}B~~\mapsto~~S^{-1}A\sr{id\otimes f}{\ral}S^{-1}B~~\big(\txt{where}~~S^{-1}A:=S^{-1}R\otimes_R A\big)
\]
is exact (i.e., $S^{-1}R$ is a flat R-module).
\end{thm}
\begin{proof}
In $R$-mod, let $0\ra A\sr{f}{\ral}B\sr{g}{\ral}C\ra 0$ be exact. Since $S^{-1}R\otimes_R-$ is right exact, it suffices to show $id\otimes f:S^{-1}A\ra S^{-1}B$ is injective. Let ${1\over s}\otimes a\in\ker(id\otimes f)$. Then
\bea
&&\textstyle 0=(id\otimes f)({1\over s}\otimes a)={1\over s}\otimes f(a)={1\over s}\big(1\otimes f(a)\big)={1\over s}h_B\big(f(a)\big)\in\im h_B,\nn\\
&&\textstyle~~\Ra~~f(a)\in\ker h_B=\{b\in B:tb=0~\txt{for some}~t\in S\},\nn\\
&&\textstyle~~\Ra~~0=tf(a)=f(ta),~~\txt{for some}~~t\in S,~~~~\Ra~~~~ta=0,~~~~\txt{(since $f$ is injective)}\nn\\
&&\textstyle~~\Ra~~{1\over s}\otimes a={1\over ts}\otimes ta={1\over ts}\otimes 0=0,\nn
\eea
where ~$h_B:B\ra S^{-1}B,~b\mapsto 1\otimes b$.
\end{proof}

\section{Ideals of the Commutative Ring Localization}
Throughout this section, $R$ will be a commutative ring unless specified otherwise. If $S\subset R$ is a localizing set, recall that we use both $h:R\ra S^{-1}R$ and $h:R\ra R[S^{-1}]$ to denote the localization of $R$ by $S$. Similarly, for an $R$-module $M$, we use both $h_M:M\ra S^{-1}M$ and $h_M:M\ra M[S^{-1}]$ to denote the localization of $M$ by $S$.
\begin{lmm}
Let $R$ be a commutative ring and $h:R\ra S^{-1}R$ a localization. Then
$S^{-1}R=\{h(r)h(s)^{-1}:r\in R,s\in S\}=h(R)h(S)^{-1}$, where $h(S)^{-1}:=\{h(s)^{-1}=1/s:s\in S\}$ and summation is not necessary in the multiplication of subsets $h(R)h(S)^{-1}$.
\end{lmm}
\begin{proof}
Each ${r\over s}\in S^{-1}R$ can be written as ${r\over s}={r\over 1}{1\over s}=h(r)h(s)^{-1}$, and so $S^{-1}R=h(R)h(S)^{-1}$.
\end{proof}

\begin{lmm}
Let $R$ be a commutative ring and $h:R\ra S^{-1}R$ a localization. If $S\subset U(R)$, then $S^{-1}R\cong R$ (i.e., $S^{-1}R$ cannot contain any new elements that are not already found in $R$).
\end{lmm}
\begin{proof}
Let $S\subset U(R)$. Since $h$ is a ring homomorphism, we have $h(s)^{-1}=h(s^{-1})$ for all $s\in S$. Thus, every element ${r\over s}\in S^{-1}R$ can be written as ${r\over s}=h(r)h(s)^{-1}=h(r)h(s^{-1})=h(rs^{-1})$. Hence
\bit
\item[]\hspace{2cm} $S^{-1}R=h(R)h(S)^{-1}=h(R)h(S^{-1})=h(RS^{-1})=h(R)\cong R$. \qedhere
\eit
\end{proof}

\begin{notation}
Recall that the ring homomorphism $h:R\ra S^{-1}R$ makes $S^{-1}R$ an R-module via~
\bea
\textstyle R\times S^{-1}R\ra S^{-1}R,~\left(a,{r\over s}\right)\mapsto a{r\over s}:=h(a){r\over s}={a\over 1}{r\over s}={ar\over s}.\nn
\eea
For an ideal $I\vartriangleleft R$, we will write $S^{-1}I:=I\cdot(S^{-1}R)=h(I)h(S)^{-1}$, because $S^{-1}R=h(R)h(S)^{-1}$ and
\bea
\textstyle I\cdot(S^{-1}R)=\left\{a{r\over s}={ar\over s}~|~a\in I,r\in R,s\in S\right\}=h(I)h(R)h(S)^{-1}=h(IR)h(S^{-1})=h(I)h(S)^{-1}.\nn
\eea
\end{notation}

\begin{lmm}[\blue{Characterization of ideals (and proper ideals) in $S^{-1}R$}]\label{LocIdChrLmm}
Let $R$ be a commutative ring, $h:R\ra S^{-1}R$ a localization, and $I\lhd R$. Then (i) $S^{-1}I\lhd S^{-1}R$, (ii) $S^{-1}I=S^{-1}R$ $\iff$ $I\cap S\neq\emptyset$, and (iii) every ideal $J\lhd S^{-1}R$ has the form $J=S^{-1}I$ for some $I\lhd R$. (\blue{footnote}\footnote{\ul{Agreement with the module localization}: Since ideals of $R$ are R-modules, we have ~$S^{-1}I\cong S^{-1}R\otimes_R I=:I[S^{-1}]$,~ which shows the notation $S^{-1}I$ agrees with that of the module localization $h_I:I\ra S^{-1}I=I[S^{-1}]$.}).
\end{lmm}
\begin{proof}
(i) This is clear from the definition of $S^{-1}I\subset S^{-1}R$. (ii) If $S^{-1}I=S^{-1}R$, then for some $a\in I$ and $s\in S$, we have ${a\over s}={1\over 1}$, i.e., $(a-s)s_1=0$ for some $s_1\in S$, and so $as_1=ss_1\in I\cap S$. Conversely, if $I\cap S\neq\emptyset$, then with $a\in I\cap S$, we have ${1\over 1}={a\over a}\in h(I)h(S)^{-1}=S^{-1}I$, and so $S^{-1}I=S^{-1}R$ as an ideal containing a unit.

(iii) Let $J\lhd S^{-1}R$. Then with the ideal ~$I:=h^{-1}(J)=\{r\in R:h(r)\in J\}\lhd R$~ of $R$ we have
\bea
J=J~S^{-1}R\supset J\cap h(R)~S^{-1}R= h\big(h^{-1}(J)\big)~S^{-1}R=h(I)~S^{-1}R=I\cdot(S^{-1}R)=S^{-1}I.\nn
\eea
Also, for each ${r\over s}={1\over s}{r\over 1}={1\over s}h(r)\in J$, we have $h(r)\in {s\over 1}J\subset J$, and so $r\in h^{-1}(J)=I$. Therefore, $J\subset h(I)h(S)^{-1}=S^{-1}I$. Hence ~$J=S^{-1}I$, where $I:=h^{-1}(J)$.
\end{proof}
From the above proof, we have the useful formula
\bea
\label{InvLocFormEq} S^{-1}h^{-1}(J)=J.
\eea

\begin{crl}
Let $R$ be a commutative ring and $I\lhd R$ an ideal. If $I$ is finitely generated, then so is every localization ~{\small $h_I=h|_R:I\subset R\ra S^{-1}I\subset S^{-1}R$}~ of $I$.
\end{crl}
\begin{proof}
Let $I=\sum_{i=1}^n Rx_i$. Then {\small $S^{-1}I=S^{-1}\sum_{i=1}^n Rx_i\subset \sum_{i=1}^n S^{-1}Rx_i\subset\sum_{i=1}^n S^{-1}I=S^{-1}I$}, and so
\[
\textstyle S^{-1}I=\sum_{i=1}^n S^{-1}Rx_i=\sum_{i=1}^n (S^{-1}R)h(x_i),~~\txt{where}~~h(x_i)={x_i\over 1}. \qedhere
\]
\end{proof}

\begin{crl}
Let $R$ be a commutative ring and $I\lhd R$. Then $S^{-1}{R\over I}\cong{S^{-1}R\over S^{-1}I}$ as $R$-modules.
\end{crl}
\begin{proof}
Let $M$ be an $R$-module. Since $0\ra I\sr{i}{\ral} R\sr{\pi}{\ral} {R\over I}\ra0$ is exact, the right exactness of $(-)\otimes_RM$ gives an exact sequence $I\otimes_RM\sr{i\otimes id}{\ral} R\otimes_RM\sr{\pi\otimes id}{\ral} {R\over I}\otimes_RM\ra0$. Therefore, by the first isomorphism theorem
\[
\textstyle{R\over I}\otimes_RM\cong{R\otimes_RM\over\ker(\pi\otimes id)}={R\otimes_RM\over\txt{im}(i\otimes id)}={R\otimes_RM\over I\otimes_RM}\cong{M\over IM},
\]
where the isomorphism at the last step is given by $\vphi: {R\otimes_RM\over I\otimes_RM}\ra {M\over IM},~~r\otimes m+I\otimes_RM\mapsto rm+IM$ and $\vphi^{-1}: {M\over IM}\ra {R\otimes_RM\over I\otimes_RM},~~m+IM\mapsto 1_R\otimes m+I\otimes_RM$. In particular, if $M:=S^{-1}R$, we get
\[
\textstyle S^{-1}{R\over I}=S^{-1}R\otimes_R{R\over I}\cong {R\over I}\otimes_RS^{-1}R\cong{S^{-1}R\over I(S^{-1}R)}={S^{-1}R\over S^{-1}I}. \qedhere
\]
\end{proof}

\begin{dfn}[\blue{\index{Localization at a prime ideal}{Localization at a prime ideal}}]
Let $R$ be a commutative ring, $P\lhd R$ a prime ideal, and $S_P:=R\backslash P$ the associated localizing set. The localization of $R$ at $P$ is ~$R_P:=S_P^{-1}R=\left\{{r\over s}:r\in R,~s\not\in P\right\}$. Similarly, if $M$ is an $R$-module, the localization of $M$ at $P$ is $M_P:=S_P^{-1}M$. In particular, for any ideal $I\lhd R$, the localization of $I$ at $P$ is $I_P:=S_P^{-1}I$ .
\end{dfn}

\begin{thm}
Let $R$ be a commutative ring. If $P\lhd R$ is prime, then $P_P=S_P^{-1}P$ is the unique maximal ideal of $R_P=S_P^{-1}R$, i.e., $(R_P,P_P)$ is a local ring.
\end{thm}
\begin{proof}
By Lemma \ref{LocIdChrLmm}(ii), any ideal $I\lhd R$ properly containing $P$ (i.e., $P\subsetneq I\lhd R$) intersects $S_P$, and so becomes non-proper in $R_P$. This implies $P_P$ is a maximal ideal in $R_P$. Again by Lemma \ref{LocIdChrLmm}(ii), every other ideal $I\lhd R$ not contained in $P$ intersects $S_P$ and so becomes non-proper in $R_P$. Hence $P_P$ is the only maximal ideal in $R_P$.
\end{proof}

\begin{thm}[\textcolor{blue}{The prime ideals of $S^{-1}R$ correspond to the prime ideals of $R$ that do not intersect $S$}]\label{LocPrIdThm}
Let $R$ be a commutative ring, and $h:R\ra S^{-1}R$ a localization. Then we have a bijective correspondence
\bea
\{\txt{prime ideals}~P\lhd R~\txt{with}~P\cap S=\emptyset\}~~\longleftrightarrow~~\{\txt{prime ideals}~Q\lhd S^{-1}R\}.\nn
\eea
\end{thm}
\begin{proof}
Let $\P:=\{\txt{prime}~P\lhd R~\txt{with}~P\cap S=\emptyset\}$ and $\Q:=\{\txt{prime}~Q\lhd S^{-1}R\}$. Consider the maps
\bea
\al:\P\ra \Q,~~P\mapsto S^{-1}P~~~~\txt{and}~~~~\beta:\Q\ra\P,~Q\mapsto h^{-1}(Q).\nn
\eea
We need to show the maps are well defined.
{\flushleft\ul{$\al$ is well defined}}: If $P\in\P$, then $S^{-1}P$ is prime since for any ${a\over s},{a'\over s'}\in S^{-1}R$,
\bea
&&\textstyle {a\over s}{a'\over s'}={aa'\over ss'}\in S^{-1}P~~\Ra~~h(aa')={aa'\over 1}=ss'{aa'\over ss'}\in S^{-1}P,~~\Ra~~aa'\in h^{-1}(S^{-1}P),~~\txt{where}\nn\\
\label{LocPrIdEq1}&&\textstyle h^{-1}(S^{-1}P)=\left\{r\in R:{r\over 1}\in S^{-1}P\right\}=\left\{r\in R:{r\over 1}={p\over s},~\txt{some}~p\in P,~s\in S\right\}\nn\\
&&\textstyle~~~~=\left\{r\in R:(rs-p)s_1=0,~\txt{some}~p\in P,~s,s_1\in S\right\}=\left\{r\in R:rss_1=ps_1\in P,~s,s_1\in S\right\}\nn\\
&&\textstyle~~~~=\left\{r\in R:r\in P~\txt{or}~ss_1\in P\cap S=\emptyset\right\}=\left\{r\in R:r\in P\right\}=P,\\
&&\textstyle~~\Ra~~a\in P~~\txt{or}~~a'\in P,~~\Ra~~{a\over s}\in h(P)h(S)^{-1}=S^{-1}P~~\txt{or}~~{a'\over s'}\in h(P)h(S)^{-1}=S^{-1}P.\nn
\eea
{\flushleft\ul{$\beta$ is well defined}}: If $Q\vartriangleleft S^{-1}R$ is prime, let $P:=h^{-1}(Q)=\{r\in R:~h(r)\in Q\}$. Then we know (from the proof of Lemma \ref{LocIdChrLmm}) that $Q=S^{-1}P=h(P)h(S)^{-1}$. $P$ is prime because for any $a,b\in R$,
\bea
\textstyle ab\in P~~\Ra~~{ab\over 1}={a\over 1}{b\over 1}\in Q,~~\Ra~~{a\over 1}\in Q~~\txt{or}~~{b\over 1}\in Q,~~\Ra~~a\in P~~\txt{or}~~b\in P.\nn
\eea

It remains to show that $\al\beta=id_{\Q}$ and $\beta\al=id_\P$. For any $Q\in\Q$, we have
\[
\al\beta(Q)=\al(h^{-1}(Q))\sr{(\ref{InvLocFormEq})}{=}S^{-1}h^{-1}(Q)=Q~~~~\txt{(as seen in the proof of Lemma \ref{LocIdChrLmm})}.
\]
Also, for any $P\in\P$, it follows from (\ref{LocPrIdEq1}) that
\[
\beta\al(P)=\beta(S^{-1}P)=h^{-1}(S^{-1}P)=P. \qedhere
\]
\end{proof}

%% file: parts/AlgebraCat/AlgebraCatS9.tex
\chapter{Homology Sequence, Resolvability, and Lifting}\label{AlgebraCatS9}
In this chapter we begin \index{Homological algebra}{\ul{homological algebra}}, which is the systematic/categorical study of the properties of algebraic objects such as groups, rings, modules, algebras using homological/derived functors (which are functors that arise in various ways from homology functors and their generalizations). A primary tool in homological algebra is a functorial process (i.e., a functor) called ``long exact sequence (LES) of homology'', which will be our starting point.
\section{Long Exact Sequence (LES) of Homology and Consequences}
\begin{dfn}[\textcolor{blue}{
\index{Long exact sequence (LES) of homology}{Long exact sequence (LES) of homology},
\index{Induced! LES of homologies}{Induced LES of homologies},
\index{Connecting morphisms}{Connecting morphisms},
\index{LES-category}{LES-category}}]
Let $\A$ be an abelian category and $\txt{SES}(\A_0^\Integer)\subset\A_0^\Integer$ the full subcategory formed by short exact sequences of chain complexes. A \ul{long exact sequence (LES) of homology} for $\A$ is a functor ~$H^\ast:\txt{SES}(\A_0^\Integer)\ra\A_0^\Integer,$~ $S\sr{\al}{\ral}S'~~\mapsto~~H^\ast(S)\sr{H^\ast(\al)}{\ral}H^\ast(S')$~ defined as follows: Given a SES of chain complexes,~ $S:~0\ra A^\ast\sr{f^\ast}{\ra}B^\ast\sr{g^\ast}{\ra}C^\ast\ra0$,~ in $\A_0^\Integer$, as displayed in the diagram

\[\adjustbox{scale=0.9}{\bt
  && \vdots \ar[d,"d_A^{n-1}"] && \vdots \ar[d,"d_B^{n-1}"] && \vdots \ar[d,"d_C^{n-1}"] && \\
0 \ar[rr] && A^n \ar[d,"d_A^n"] \ar[rr,"f^n"] && B^n \ar[d,dashed,"d_B^n"] \ar[rr,dashed,"g^n"] && C^n \ar[d,"d_C^n"] \ar[rr,dashed] &&  0 \\
0 \ar[rr,dashed] && A^{n+1} \ar[d,"d_A^{n+1}"] \ar[rr,dashed,"f^{n+1}"] && B^{n+1} \ar[d,"d_B^{n+1}"] \ar[rr,"g^{n+1}"] && C^{n+1} \ar[d,"d_C^{n+1}"] \ar[rr] &&  0 \\
  && \vdots  && \vdots  && \vdots  && \\
\et}\]
an \ul{induced ``long'' exact sequence (LES) of homologies} of $S$ is an exact sequence of the form
\[\adjustbox{scale=0.9}{\bt
H^\ast(S):~\cdots\ar[r,dashed,"\delta^{n-1}"] &
H^n(A^\ast) \ar[rr,"H^n(f^\ast)"] &&
H^n(B^\ast) \ar[rr,"H^n(g^\ast)"]
\arrow[d, phantom, ""{coordinate, name=Z}] &&
H^n(C^\ast) \arrow[dllll,dashed,
"\delta^n"',
rounded corners,
to path={ -- ([xshift=2ex]\tikztostart.east)
|- (Z) [near end]\tikztonodes
-| ([xshift=-2ex]\tikztotarget.west)
-- (\tikztotarget)}] & \\
& H^{n+1}(A^\ast)  \ar[rr,"H^{n+1}(f^\ast)"] &&
H^{n+1}(B^\ast) \ar[rr,"H^{n+1}(g^\ast)"] &&
H^{n+1}(C^\ast)\ar[r,dashed,"\delta^{n+1}"] & \cdots
\et}\]
where the morphisms $\delta^n$ are called \ul{connecting morphisms}. That is, the induced LES of $S$ an exact sequence of the form
{\footnotesize\bea
H^\ast(S):~\cdots\sr{\delta^{n-1}}{\ral}H^n\left(A^\ast\sr{f^\ast}{\ra}B^\ast\sr{g^\ast}{\ra}C^\ast\right)\sr{\delta^n}{\ral}H^{n+1}\left(A^\ast\sr{f^\ast}{\ra}B^\ast\sr{g^\ast}{\ra}C^\ast\right)\sr{\delta^{n+1}}{\ral}\cdots,\nn
\eea}where
{\footnotesize $H^n\left(A^\ast\sr{f^\ast}{\ra}B^\ast\sr{g^\ast}{\ra}C^\ast\right)=H^n(A^\ast)\sr{H^n(f^\ast)}{\ral}H^n(B^\ast)\sr{H^n(g^\ast)}{\ral}H^n(C^\ast)$}.

The abelian category $\A$ is a \ul{LES-category} if it has a LES of homology (in the obvious sense that every SES of chain complexes in $\A_0^\Integer$ has a natural/functorial induced LES of homologies).
\end{dfn}
\begin{rmk}[\blue{Remembrance aid}]
In the induced LES of homologies above, think of the connecting morphism $\delta^n$ as the three differentials $\left(d_A^n,d_B^n,d_C^n\right)$ collapsed (combined) in some way using $f^n$, $g^n$, $f^{n+1}$, $g^{n+1}$, $H^n$, $H^{n+1}$ into a single morphism
\bea
\delta^n=\delta^n\big(d_A^n,d_B^n,d_C^n\!~|\!~f^\ast,g^\ast,H^\ast\big)\nn
\eea
acting as a forward \ul{head-to-tail connector} of the homology images
\bea
H^n(A^\ast)\sr{H^n(f^\ast)}{\ral}H^n(B^\ast)\sr{H^n(f^\ast)}{\ral}H^n(C^\ast)~~~\txt{and}~~~H^{n+1}(A^\ast)\sr{H^{n+1}(f^\ast)}{\ral}H^{n+1}(B^\ast)\sr{H^{n+1}(f^\ast)}{\ral}H^{n+1}(C^\ast)\nn
\eea
of the two short exact sequences
\bea
0\ra A^n\sr{f^n}{\ral}B^n\sr{f^n}{\ral}C^n\ra 0~~~\txt{and}~~~0\ra A^{n+1}\sr{f^{n+1}}{\ral}B^{n+1}\sr{f^{n+1}}{\ral}C^{n+1}\ra 0.\nn
\eea
\end{rmk}

\begin{thm}[\textcolor{blue}{$R$-mod is a LES-category: \cite[pp 40-44]{gelfand-manin2010}}]\label{RmodLES1}
Given a short exact sequence of complexes of $R$-modules~ $S:~0\ra A^\ast\sr{f^\ast}{\ra}B^\ast\sr{g^\ast}{\ra}C^\ast\ra0$,~ as displayed in the diagram

\[\adjustbox{scale=0.9}{\bt
  && \vdots \ar[dd,"d_A^{n-1}"] && \vdots \ar[dd,"d_B^{n-1}"] && \vdots \ar[dd,"d_C^{n-1}"] && \\
  && && && && \\
0 \ar[rr] && A^n \ar[dd,"d_A^n"] \ar[rr,"f^n"] && B^n \ar[dd,dashed,"d_B^n"] \ar[rr,dashed,"g^n"] && C^n \ar[dd,"d_C^n"] \ar[rr,dashed] &&  0 \\
  && && && && \\
0 \ar[rr,dashed] && A^{n+1} \ar[dd,"d_A^{n+1}"] \ar[rr,dashed,"f^{n+1}"] && B^{n+1} \ar[dd,"d_B^{n+1}"] \ar[rr,"g^{n+1}"] && C^{n+1} \ar[dd,"d_C^{n+1}"] \ar[rr] &&  0 \\
  && && && && \\
  && \vdots  && \vdots  && \vdots  && \\
\et}\]
there exists a \ul{functorial} (i.e., \ul{natural}) induced ``long'' exact sequence (LES) of homologies
\[\adjustbox{scale=0.9}{\bt
H^\ast(S):~\cdots\ar[r,dashed,"\delta^{n-1}"] &
H^n(A^\ast) \ar[rr,"H^n(f^\ast)"] &&
H^n(B^\ast) \ar[rr,"H^n(g^\ast)"]
\arrow[d, phantom, ""{coordinate, name=Z}] &&
H^n(C^\ast) \arrow[dllll,dashed,
"\delta^n"',
rounded corners,
to path={ -- ([xshift=2ex]\tikztostart.east)
|- (Z) [near end]\tikztonodes
-| ([xshift=-2ex]\tikztotarget.west)
-- (\tikztotarget)}] & \\
& H^{n+1}(A^\ast)  \ar[rr,"H^{n+1}(f^\ast)"] &&
H^{n+1}(B^\ast) \ar[rr,"H^{n+1}(g^\ast)"] &&
H^{n+1}(C^\ast)\ar[r,dashed,"\delta^{n+1}"] & \cdots
\et}\]
where for each integer $n$,~ {\small $\delta^n=\delta^n(f^\ast,g^\ast):H^n(C^\ast)={\ker d_C^n\over\im d_C^{n-1}}\ra H^{n+1}(A^\ast)={\ker d_A^{n+1}\over\im d_A^n}$}~ is defined by
\bea
\delta^n(c^n+\im d_C^{n-1}):=a^{n+1}+\im d_A^n,~~~~~~~~a^{n+1}\in (f^{n+1})^{-1}d_B^n(g^n)^{-1}c^n,\nn
\eea
where~ $a^{n+1}\in\ker d_A^{n+1}$~ is related to~ $c^n\in \ker~d_C^n$~ through a ``diagram chase'' as follows:
\bea
f^{n+1}a^{n+1}=d_B^nb^n~~~~\txt{for some}~~~b^n\in B^n~~\txt{such that}~~~~g^n(b^n)=c^n.\nn
\eea
\end{thm}
\begin{proof}
The proof proceeds in four steps. First, we show that $H^\ast(S)$ is a sequence of R-modules by constructing the connecting homomorphisms $\delta^n:H^n(C^\ast)\ra H^{n+1}(A^\ast)$. Next, we show $H^\ast(S)$ is a complex. Then we show $H^\ast(S)$ is exact. Finally, we prove functoriality of the LES in a separate theorem.

{\flushleft\bf \ul{(I) $H^\ast(S)$ is a sequence of R-modules: Existence of connecting homomorphisms}}~
\[\adjustbox{scale=0.9}{\bt
 \cdots\ar[r,"\substack{\delta^{n-1}}"] & ~\substack{H^n(A^\ast)}~\ar[r,"\substack{H^n(f^\ast)}"] & ~\substack{H^n(B^\ast)}~\ar[r,"\substack{H^n(g^\ast)}"] & ~\substack{H^n(C^\ast)}~\ar[r,dashed,"\substack{\delta^n}"] & ~\substack{H^{n+1}(A^\ast)}~\ar[r,"\substack{H^{n+1}(f^\ast)}"] & ~\substack{H^{n+1}(B^\ast)}~\ar[r,"\substack{H^{n+1}(g^\ast)}"] & ~\substack{H^{n+1}(C^\ast)}~\ar[r,dashed,"\substack{\delta^{n+1}}"] & \cdots
\et}\]
\[\adjustbox{scale=0.9}{\bt
   &  \vdots\ar[d,"d_A^{n-1}"] &  \vdots\ar[d,"d_B^{n-1}"] & \vdots\ar[d,"d_C^{n-1}"]  & \\
   & A^n \ar[d,"d_A^n"]\ar[r,"f^n"]  & B^n \ar[d,dashed,"d_B^n"]\ar[r,dashed,"g^n"]  & C^n \ar[d,"d_C^n"]\ar[r,dashed]  & ~ \\
 \ar[r,dashed]  & A^{n+1} \ar[d,"d_A^{n+1}"]\ar[r,dashed,"f^{n+1}"]  & B^{n+1} \ar[d,"d_B^{n+1}"]\ar[r,"g^{n+1}"]  & C^{n+1} \ar[d,"d_C^{n+1}"]  &  \\
   & \vdots & \vdots & \vdots  &
\et}\]

If $c^n\in\ker(d^n_C)$, i.e., $d_C^nc^n=0$, then by the surjectivity of $g^n$, there is $b^n\in B^n$ such that
\bea
&&\label{exist-eq1}g^n(b^n)=c^n,\\
&&~~\Ra~~0=d_C^nc^n=(d_C^ng^n)b^n=(g^{n+1}d_B^n)b^n,~~\Ra~~d_B^nb^n\in\ker(g^{n+1})=\im (f^{n+1}),\nn\\
\label{exist-eq2}&&~~\Ra~~d_B^nb^n=f^{n+1}(a^{n+1}),~~~~\txt{for some}~~a^{n+1}\in A^{n+1}.
\eea
We now define $\delta^n:H^n(C^\ast)\ra H^{n+1}(A^\ast)$ by
\bea
\label{connect-hom}\delta^n:~c^n+\im (d_C^{n-1})~\mapsto~a^{n+1}+\im (d_A^n),
\eea
where $c^n,a^{n+1}$ satisfy (\ref{exist-eq1}), (\ref{exist-eq2}). Note that $\delta^n$ is a well defined $R$-homomorphism, because if we also have $g^n(\wt{b}^n)=c^n$ for some $\wt{b}^n\in B^n$, then
\bea
&&\wt{b}^n-b^n\in\ker g^n=\im f^n,~~\Ra~~\wt{b}^n=b^n+f^n(a^n),~~~~\txt{for some}~~a^n\in A^n,\nn\\
&&~~\Ra~~d_B^n\wt{b}^n=d_B^nb^n+d_B^nf^n(a^n)=f^{n+1}a^{n+1}+f^{n+1}d_A^na^n=f^{n+1}\wt{a}^{n+1},~~~~\wt{a}^{n+1}:=a^{n+1}+d_A^na^n,\nn\\
&&~~\Ra~~\delta^n(c^n+\im (d_C^{n-1}))=a^{n+1}+\im (d_A^n)=\wt{a}^{n+1}+\im (d_A^n),\nn
\eea
where $c^n,\wt{a}^{n+1}$ satisfy (\ref{exist-eq1}), (\ref{exist-eq2}), with $a^{n+1},b^n$ replaced by  $\wt{a}^{n+1},\wt{b}^n$ respectively.

{\flushleft\bf \ul{(II) $H^\ast(S)$ is a complex}}: To show that $H^\ast(S)$ is a complex, we need to check that $H^n(g^\ast)H^n(f^\ast)=0$, $\delta^nH^n(g^\ast)=0$, and $H^{n+1}(f^\ast)\delta^n=0$.
\begin{enumerate}
\item Since $H^n$ is a covariant functor, we indeed have $H^n(g^\ast)H^n(f^\ast)=H^n(g^\ast f^\ast)=H^n(0)=0$.

\item To check that $\delta^nH^n(g^\ast)=0$, let $b^n\in\ker(d_B^n)$, i.e.,
\bea
\label{exist-eq4}d_B^nb^n=0.
\eea
Then,
{\small
\bea
&& \delta^nH^n(g^\ast)\big(b^n+\im d_B^{n-1}\big)=\delta^n\big(g^n(b^n)+\im d_C^{n-1}\big)=\delta^n\big(\wt{c}^n+\im (d_C^{n-1})\big),~~~~\wt{c}^n:=g^n(b^n),\nn\\
&&~~~~=\wt{a}^{n+1}+\im d_A^n,~~\txt{where by (\ref{exist-eq1}),(\ref{exist-eq2}),}~~f^{n+1}(\wt{a}^{n+1})=d_B^n\wt{b}^n,~~~~g^n(\wt{b}^n)=\wt{c}^n=g^n(b^n),\nn\\
&&~~\Ra~~\wt{b}^n-b^n\in\ker g^n=\im f^n,~~\Ra~~\wt{b}^n=b^n+f^na^n~~\txt{for some}~~a^n\in A^n,\nn\\
&&~~\Ra~~f^{n+1}(\wt{a}^{n+1})=d_B^n(b^n+f^na^n)\sr{(\ref{exist-eq4})}{=}d_B^nf^na^n=f^{n+1}d_A^na^n,\nn\\
&&~~\Ra~~\wt{a}^{n+1}=d_A^na^n\in \im d_A^n,~~~~\txt{since $f^{n+1}$ is injective},\nn\\
&&~~\Ra~~\delta^nH^n(g^\ast)=0.\nn
\eea}

\item To check that $H^{n+1}(f^\ast)\delta^n=0$, let $c^n\in\ker(d_C^n)$. Then
\bea
&&H^{n+1}(f^\ast)\delta^n\big(c^n+\im d_C^{n-1}\big)=H^{n+1}(f^\ast)\big(a^{n+1}+\im d_A^n\big)=f^{n+1}(a^{n+1})+\im d_B^n\sr{(\ref{exist-eq2})}{=}d_B^nb^n+\im d_B^n\nn\\
&&~~~~=0+\im d_B^n=0.\nn
\eea
\end{enumerate}

{\flushleft\bf \ul{(III) $H^\ast(S)$ is exact}}: Since $H^\ast(S)$ is a complex, we only need to show that $\ker(H^n(g^\ast))\subset \im(H^n(f^\ast))$, $\ker(\delta^n)\subset \im (H^n(g^\ast))$, and $\ker(H^{n+1}(f^\ast))\subset \im(\delta^n)$.
\begin{enumerate}
\item\ul{$\ker(H^n(g^\ast))\subset \im(H^n(f^\ast))$}:~ Let ~$b^n+\im d_B^{n-1}\in \ker(H^n(g^\ast))$.~ Then,
\bea 0=H^n(g^\ast)\big(b^n+\im d_B^{n-1}\big)~~\Ra~~0=g^n(b^n)+\im d_C^{n-1},~~\Ra~~g^n(b^n)=d_C^{n-1}(c^{n-1}).\nn
\eea
Since $g^\ast$ is surjective, we have $c^{n-1}=g^{n-1}(b^{n-1})$. Thus,
{\small
\bea
&& g^n(b^n)=d_C^{n-1}(c^{n-1})=d_C^{n-1}g^{n-1}(b^{n-1})=g^nd_B^{n-1}(b^{n-1}),\nn\\
&&~~\Ra~~g^n(b^n-d_B^{n-1}b^{n-1})=0,~~\Ra~~b^n-d_B^{n-1}b^{n-1}=f^n(a),~~\Ra~~b^n=f^n(a^n)+d_B^{n-1}b^{n-1},\nn\\
&& ~~\Ra~~b^n+\im d_B^{n-1}=f^n(a^n)+\im d_B^{n-1}= H^n(f^\ast)(a^n+\im d_A^{n-1}).
\eea}

\item\ul{$\ker(\delta^n)\subset \im(H^n(g^\ast))$}:~ Let ~$c^n+\im d_C^{n-1}\in \ker(\delta^n)$.~ Then,
\begin{align}
& 0=\delta^n\big(c^n+\im d_C^{n-1}\big)=a^{n+1}+\im d_A^n,~~\txt{where}~~f^{n+1}(a^{n+1})=d_B^nb^n,~~g^n(b^n)=c^n,\nn\\
&~~\Ra~~a^{n+1}=d_A^na^n,~~~~\txt{for some}~~a^n\in A^n,\nn\\
&~~\Ra~~d_B^nb^n=f^{n+1}(a^{n+1})=f^{n+1}(d_A^na^n)=d_B^nf^n(a^n)\nn\\
&~~\Ra~~\wt{b}^n:=b^n-f^n(a^n)\in\ker(d_B^n),~~~~\txt{i.e.,}~~d_B^n\wt{b}^n=0,\nn\\
&~~\Ra~~c^n=g^n(b^n)=g^n(b^n)-0=g^nb^n-g^nf^n(a^n)=g^n\left(b^n-f^n(a^n)\right)\nn\\
&~~~~~~~~~~~~=g^n(\wt{b}^n),~~~~d_B^n\wt{b}^n=0,\nn\\
&~~\Ra~~c^n+\im d_C^{n-1}=g^n(\wt{b}^n)+\im d_C^{n-1}=H^n(g^\ast)\big(\wt{b}^n+\im d_B^{n-1}\big).
\end{align}
\item\ul{$\ker(H^{n+1}(f^\ast))\subset \im(\delta^n)$}:~ Let~ $a^{n+1}+\im d_A^n\in \ker(H^{n+1}(f^\ast))$.~ Then
\bea
&& 0=H^{n+1}(f^\ast)\big(a^{n+1}+\im d_A^n\big)=f^{n+1}(a^{n+1})+\im d_B^n\nn\\
&&~~\Ra~~f^{n+1}(a^{n+1})=d_B^nb^n,~~~~\txt{for some}~~b^n\in B^n,~~~~~~~~~~\txt{(compare with (\ref{exist-eq2}))},\nn\\
&&~~\Ra~~d_B^nb^n\in \im f^{n+1}=\ker g^{n+1},~~\Ra~~g^{n+1}d_B^nb^n=0,\nn\\
&&~~\Ra~~d_C^ng^n(b^n)=g^{n+1}d_B^nb^n=0,\nn\\
&&~~\Ra~~c^n:= g^n(b^n)\in\ker d_C^n.~~~~~~~~~~~~~~~~~~~~~~~~~~~~~~~~\txt{(compare with (\ref{exist-eq1}))}.\nn
\eea
Therefore, by (\ref{exist-eq2}) and (\ref{exist-eq1}), we get
\bea
a^{n+1}+\im d_A^n=\delta^n\big(c^n+\im d_C^n\big).
\eea
\end{enumerate}
{\flushleft\bf \ul{(IV) Functoriality}}: For functoriality of the LES, see Theorem \ref{RmodLES2} next.
\end{proof}

Notice the switch from covariant to contravariant indices in the following theorem. This should not cause any confusion since we always have covariant-contravariant symmetry in the sense that results proved for covariant systems also automatically hold (up to conventional readjustments) for contravariant systems, and vice versa.

\begin{thm}[\textcolor{blue}{\index{Functoriality of the LES of homology}{Functoriality of the LES of homology in $R$-mod}}]\label{RmodLES2}
The maps in the induced LES of homologies for $R$-mod are funtorial, that is, for any commutative diagram of complexes of $R$-modules
\[\adjustbox{scale=0.9}{\bt
0\ar[r] & A_\ast\ar[d,"\al"]\ar[r,"f"] & B_\ast\ar[d,"\beta"]\ar[r,"g"] & C_\ast\ar[d,"\gamma"]\ar[r] &0\\
0\ar[r] & A'_\ast\ar[r,"f'"] & B'_\ast\ar[r,"g'"] & C'_\ast\ar[r] &0
\et}\]
with exact rows, we get a commutative diagram with the two induced LES's of homologies as rows:
\[\adjustbox{scale=0.9}{\bt
\cdots\ar[r] & H_n(A_\ast)\ar[d,"H_n(\al)"]\ar[r,"H_n(f)"] & H_n(B_\ast)\ar[d,"H_n(\beta)"]\ar[r,"H_n(g)"] & H_n(C_\ast)\ar[d,"H_n(\gamma)"]\ar[r,"\delta_n"] & H_{n-1}(A_\ast)\ar[d,"H_{n-1}(\al)"]\ar[r] &\cdots\\
\cdots\ar[r] & H_n(A'_\ast)\ar[r,"H_n(f')"] & H_n(B'_\ast)\ar[r,"H_n(g')"] & H_n(C'_\ast)\ar[r,"\delta_n'"] & H_{n-1}(A'_\ast)\ar[r] &\cdots\\
\et}\]
where all maps (except $\delta$) are the maps induced on homology by the maps in the first diagram.
\end{thm}
\begin{proof}
(i) The squares with no $\delta$ commute since $H_n$ is a functor., i.e., $H_n(fg)=H_n(f)H_n(g)$ and $f,g$ are chain maps. (ii) To show the squares with $\delta$ also commute, recall that we have commutative diagrams
\[\adjustbox{scale=0.9}{\bt
 && A_n\ar[dl,"d"']\ar[dd,near start,"\al"]\ar[rr,"f_n"] &&  B_n\ar[dl,dashed,"d"']\ar[dd,near start,"\beta"]\ar[rr,dashed,"g_n"] && C_n\ar[dl,"d"']\ar[dd,near start,"\gamma"]\ar[r,dashed] &0  \\
0\ar[r,dashed] & A_{n-1}\ar[dd,near start,"\al"]\ar[rr,dashed,crossing over,near start,"f_{n-1}"] && B_{n-1}\ar[rr,crossing over,near start,"g_{n-1}"] && C_{n-1} && \\
 && A_n'\ar[dl,"d'"']\ar[rr,near start,"f_n'"]            && B_n'\ar[dl,dashed,"d'"']\ar[rr,dashed,near start,"g_n'"]              && C_n'\ar[dl,"d'"'] \ar[r,dashed,] & 0 \\
0\ar[r,dashed,] & A_{n-1}'\ar[rr,dashed,"f_{n-1}'"]            && B_{n-1}'\ar[from=uu,crossing over,near start,"\beta"]\ar[rr,"g_{n-1}'"]              && C_{n-1}'\ar[from=uu,crossing over,near start,"\gamma"] &&
\et}\]
Now, we have $H_{n-1}(\al)\delta(c+\im d)=H_{n-1}(\al)[a+\im d]=\al(a)+\im d'$, where
\bea
\label{S23eq1}f(a)=db,~~c=g(b).
\eea
Similarly, $\delta'H_n(\gamma)(c+\im d)=\delta'[\gamma(c)+\im d']=a'+\im d'$, where
\bea
\label{S23eq2}f'(a')=d'b',~~\gamma(c)=g'(b').
\eea
Since $\gamma g=g'\beta$, we have $\gamma(c)=\gamma g(b)=g'\beta(b)$, and so by (\ref{S23eq2}), $b'\in\beta(b)+\ker g'=\beta(b)+\im f'$, i.e., $b'=\beta(b)+f'(a'')$. Thus, because $\beta$ and $f'$ are chain maps (use at step (s) below),
\bea
&&d'b'=d'\beta(b)+d'f'(a'')\sr{(s)}{=}\beta db+f'd'a''=\beta f(a)+f'd'a''=f'\al(a)+f'd'a''\nn\\
&&~~~~=f'\big(\al(a)+d'a''\big),\nn
\eea
and so by comparing with (\ref{S23eq2}), we see that $f'(a')=f'\big(\al(a)+d'a''\big)$. Since $f'$ is injective, we get $a'=\al(a)+d'a''$, and so ~$\al(a)+\im d'=a'+\im d'$.
\end{proof}

\begin{crl}[\textcolor{blue}{Abelian $\Ra$ LES}]
Every abelian category is a LES-category (by Remarks \ref{FreydMitchRmk}).
\end{crl}

\begin{rmk*}
Based on the properties of abelian categories discussed so far, it is not difficult (Remarks \ref{FreydMitchRmk}) to see that the following corollaries, and similar statements, hold for any LES-category (and not just for $R$-mod). However, for simplicity, we will assume the objects are $R$-modules and so use element-by-element (i.e., element-wise) arguments where applicable.
\end{rmk*}

We have already proved the following lemma (i.e., 3-by-3 lemma) for abelian categories in the previous chapter. There will also be other repetitions, including (i) a review of split-exactness for short exact sequences of modules and (ii) reproving that a coproduct of projective modules is projective while a product of injective modules is injective.

\begin{crl}[\textcolor{blue}{\index{Nine (or three-by-three) lemma}{}\index{Three-by-three (or nine) lemma}{Three-by-three (or nine) lemma}}]\label{3By3Lmm}
Consider any commutative diagram (of $R$-modules) with exact columns:
\[\adjustbox{scale=0.8}{
\begin{tikzcd}[row sep=small,column sep=small]
  & 0 \ar[d] & 0 \ar[d] & 0 \ar[d] &  \\
 0 \ar[r] & A_1 \ar[d] \ar[r] & A_2 \ar[d] \ar[r] & A_3 \ar[d] \ar[r] & 0 \\
 0 \ar[r] & B_1 \ar[d] \ar[r] & B_2 \ar[d] \ar[r] & B_3 \ar[d] \ar[r] & 0 \\
 0 \ar[r] & C_1 \ar[d] \ar[r] & C_2 \ar[d] \ar[r] & C_3 \ar[d] \ar[r] & 0 \\
  & 0 & 0 & 0 &  \\
\end{tikzcd}}
\]

\bit
\item[(a)] If the first and second rows are exact, then the third row is exact.
\item[(b)] If the second and third rows are exact, then the first row is exact.
\eit
\end{crl}

\begin{proof}
Consider the maps in the diagram as follows:
\[\adjustbox{scale=0.8}{
\begin{tikzcd}
0\ar[d] && & 0 \ar[d] & 0 \ar[d] & 0 \ar[d] &  \\
A_\bullet \ar[d,"f"'] && 0 \ar[r] & A_1 \ar[d,"f_1"] \ar[r,"d_1^A"] & A_2 \ar[d,"f_2"] \ar[r,"d_2^A"] & A_3 \ar[d,"f_3"] \ar[r] & 0 \\
B_\bullet \ar[d,"g"'] && 0 \ar[r] & B_1 \ar[d,"g_1"] \ar[r,"d_1^B"] & B_2 \ar[d,"g_2"] \ar[r,"d_2^B"] & B_3 \ar[d,"g_3"] \ar[r] & 0 \\
C_\bullet \ar[d] && 0 \ar[r] & C_1 \ar[d] \ar[r,"d_1^C"] & C_2 \ar[d] \ar[r,"d_2^C"] & C_3 \ar[d] \ar[r] & 0 \\
0 && & 0 & 0 & 0 &  \\
\end{tikzcd}}
\]

(a) Assume the first two rows are exact. The rows form a SES of complexes, since by commutativity,
\bea
d_2^Cd_1^Cg_1=g_3d_2^Bd_1^B=g_30=0,~~~~\Ra~~~~d_2^Cd_1^C=0,~~~~\txt{since $g_1$ is surjective}.\nn
\eea
Therefore, we get the induced LES of homologies
\bea\adjustbox{scale=0.7}{
\begin{tikzcd}
     \ar[d,"\delta_2"] &  \ar[d,"\delta_3"] &  0\ar[d,"\delta_4"]   \\
    H_1(A) \ar[d,"H_1(f)"]  & H_2(A) \ar[d,"H_2(f)"]  & H_3(A) \ar[d,"H_3(f)"]    \\
    H_1(B) \ar[d,"H_1(g)"]  & H_2(B) \ar[d,"H_2(g)"]  & H_3(B) \ar[d,"H_3(g)"]    \\
    H_1(C) \ar[d,"\delta_1"] & H_2(C) \ar[d,"\delta_2"]  & H_3(C) \ar[d,"\delta_3"]    \\
    0 & ~ & ~   \\
\end{tikzcd}}
~~~~\sr{\txt{first two rows exact}}{=}~~~~
\adjustbox{scale=0.7}{
\begin{tikzcd}
    \ar[d,"\delta_2"] &  \ar[d,"\delta_3"] &  0\ar[d,"\delta_4"]   \\
    0\ar[d,"H_1(f)"]  & 0 \ar[d,"H_2(f)"]  & 0 \ar[d,"H_3(f)"]    \\
    0 \ar[d,"H_1(g)"]  & 0 \ar[d,"H_2(g)"]  & 0 \ar[d,"H_3(g)"]    \\
    H_1(C) \ar[d,"\delta_1"] & H_2(C) \ar[d,"\delta_2"]  & H_3(C) \ar[d,"\delta_3"]    \\
    0 & ~ & ~   \\
\end{tikzcd}}\nn
\eea
Hence $H_1(C)=0$, $H_2(C)=0$, $H_3(C)=0$, i.e., the third row is also exact.

(b) By a similar argument, if the 2nd and 3rd rows are exact then the 1st row is also exact.
\end{proof}

\begin{crl}[\textcolor{blue}{\index{Snake lemma}{Snake lemma}}]\label{SnakeLmm}
For any exact commutative diagram (of R-modules)
\[\adjustbox{scale=0.8}{
\begin{tikzcd}
   & A \ar[d,"\al"]\ar[r,"f"] & B \ar[d,"\beta"]\ar[r,"g"] & C \ar[d,"\gamma"]\ar[r]& 0 \\
0 \ar[r] & A' \ar[r,"f'"] & B' \ar[r,"g'"] & C' &
\end{tikzcd}}\]

there exists an exact sequence of R-modules
\bea
\ker\al \ra \ker\beta \ra \ker\gamma \ra \coker\al \ra \coker\beta \ra \coker\gamma.\nn
\eea
\end{crl}
\begin{proof}
We will essentially show that there is an exact diagram of the form
\[\adjustbox{scale=0.9}{
\begin{tikzcd}[row sep=small]
  & & 0\ar[d] & 0\ar[d] & 0\ar[d] & & \\
 0\ar[r] &\ker f\ar[r] & \ker\al\ar[d,hook]\ar[r] & \ker\beta\ar[d,hook]\ar[r] & \ker\gamma\ar[d,hook]\ar[rr,dashed,"\delta"] & & ~ \\
  & & A \ar[d,"\al"]\ar[r,"f"] & B \ar[d,"\beta"]\ar[r,"g"] & C \ar[d,"\gamma"]\ar[r]& 0 &\\
 &0 \ar[r] & A' \ar[d,two heads]\ar[r,"f'"] & B' \ar[d,two heads]\ar[r,"g'"] & C'\ar[d,two heads] & & \\
 \ar[rr,dashed,"\delta"] & & \coker\al\ar[d]\ar[r] & \coker\beta\ar[d]\ar[r] & \coker\gamma\ar[d]\ar[r] & \coker g'\ar[r] & 0\\
  & & 0 & 0 & 0 & & \\
\end{tikzcd}}\]

The maps in the given diagram can be factored as follows.

\[\adjustbox{scale=0.8}{
\begin{tikzcd}[row sep=small]
 0\ar[dr]  &  &  &  &  \\
   & \ol{A}:={A\over\ker f}\ar[dr,"\ol{f}"] &  &  &  \\
   & A \ar[u,"\pi_f"]\ar[d,"\al"]\ar[r,near start,"f"] & B \ar[d,"\beta"]\ar[r,"g"] & C \ar[d,"\gamma"]\ar[r]& 0 \\
0 \ar[r] & A' \ar[r,"f'"] & B' \ar[r,"g'"]\ar[dr,"\ol{g}'"'] & C' &  \\
   &  &  & \ol{C}:=\im g'\ar[u,hook,shift left=1,"q'"']\ar[dr] &  \\
   &  &  &  & 0 \\
\end{tikzcd}}~~~~\Ra~~~~
\adjustbox{scale=0.8}{
\begin{tikzcd}
   & 0\ar[d,hook] & 0\ar[d,hook] & 0\ar[d,hook] &  \\
 0\ar[r]  & \ol{A}\ar[d,"\ol{\al}"]\ar[r,"\ol{f}"] & B \ar[d,"\beta"]\ar[r,"g"] & C \ar[d,"\ol{\gamma}"]\ar[r]& 0 \\
0 \ar[r] & A'\ar[d,two heads]\ar[r,"f'"] & B'\ar[d,two heads]\ar[r,"\ol{g}'"] & \ol{C}'\ar[d,two heads]\ar[r]  & 0 \\
   & 0 & 0 & 0 &  \\
\end{tikzcd}}\]

Therefore, we get the commutative diagram on the right with exact rows, where the shown maps are given by (\blue{footnote}\footnote{Alternatively, we could have considered $\ol{C}'$, $\ol{g}'$, $\ol{\gamma}$, to be the following:
\bea
\textstyle\ol{C}':=\coker f'={B'\over\im f'}={B'\over\ker g'}\cong\im g',~~~~\ol{g}'(b'):=b'+\ker g',~~~~\ol{\gamma}(c):=b'+\ker g'\big|_{g'(b')=\gamma(c)}.\nn
\eea
})
\bea
&&\textstyle \ol{A}:=\coim f={A\over \ker f},~~~~\ol{f}(a+\ker f):=f(a),~~~~\ol{\al}(a+\ker f):=\al(a),\nn\\
&&~~~~g':=q'\circ \ol{g}',~~~~\ol{\gamma}:=\gamma,\nn\\
&& \ol{C}:=\im g'=g'(C),~~~~\ol{g}'(b'):=g'(b'),~~~~\ol{\gamma}(c):=\gamma(c),\nn
\eea
with $\ol{\al}$ and $\ol{\gamma}$ well-defined because $\ker f\subset\ker\al$ and $\im\gamma\subset\im g'$ (by the commutativity and exactness of the given diagram). Thus, the columns form a SES of complexes~ $0\ra A^\ast\ra B^\ast\ra C^\ast\ra 0$. From this, we get the induced LES of homologies
\[\adjustbox{scale=0.8}{
\bt[column sep=small,row sep=tiny]
 0\ar[r]  & H^1(A^\ast)\ar[r,"H(\ol{f})"] & H^1(B^\ast)\ar[r,"H(g)"] & H^1(C^\ast)\ar[r,"\ol{\delta}"]& ~ \\
\ar[r,"\ol{\delta}"] & H^2(A^\ast)\ar[r,"H(f')"] & H^2(B')\ar[r,"H(\ol{g}')"] & H^2(C^\ast)\ar[r]  & 0
\et}~~~~=~~~~\adjustbox{scale=0.8}{
\bt[column sep=small,row sep=tiny]
 0\ar[r]  & \ker\ol{\al}\ar[r] & \ker\beta\ar[r] & \ker\ol{\gamma}\ar[r,"\ol{\delta}"]& ~ \\
\ar[r,"\ol{\delta}"] & \coker\ol{\al}\ar[r] & \coker\beta\ar[r] & \coker\ol{\gamma}\ar[r]  & 0
\et}\]
where
\bea
\textstyle\ker\ol{\al}={\ker\al\over\ker f},~~~~\ker\ol{\gamma}=\ker\gamma,~~~~\coker\ol{\al}=\coker\al,~~~~\coker\ol{\gamma}={\im g'\over\im\gamma}.\nn
\eea
Using the following exact sequences (with the obvious maps),
\bea
&&\textstyle 0\ral\ker f\hookrightarrow\ker\al\sr{\pi_1}{\ral}{\ker \al\over\ker f}=\ker\ol{\al}\ral 0,\nn\\
&&\textstyle 0\ral\coker\ol{\gamma}={\im g'\over\im\gamma}\hookrightarrow\coker\gamma={C'\over\im\gamma}\sr{\pi_2}{\ral}{C'/\im\gamma\over\im g'/\im\gamma}\cong\coker g'\ral 0,\nn
\eea
we get an exact sequence
\[\adjustbox{scale=0.9}{
\bt[column sep=small,row sep=tiny]
0\ar[r] &\ker f\ar[r] & \ker\al\ar[r] & \ker\beta\ar[r] & \ker\gamma\ar[r,"\delta"]& ~ & \\
 & \ar[r,"\delta"] & \coker\al\ar[r] & \coker\beta\ar[r] & \coker\gamma\ar[r] & \coker g'\ar[r]& 0
\et} \qedhere
\]
\end{proof}

\begin{crl}[\textcolor{blue}{\index{Five lemma! short version}{Short five lemma}}]\label{Short5Lmm}
Consider any \ul{commutative diagram} (of $R$-modules) with \ul{exact rows}:
\[\adjustbox{scale=0.9}{\bt
0\ar[r] & A\ar[d,"\al"]\ar[r] & B\ar[d,"\beta"]\ar[r] & C\ar[d,"\gamma"]\ar[r] & 0\\
0\ar[r] & A'\ar[r] & B'\ar[r] & C'\ar[r] & 0
\et}\]
{\flushleft (1)} If any two of the vertical maps are isomorphisms, so is the other.
{\flushleft (2)} If any two of the vertical maps \ul{except $\al,\beta$ when $\al$ is not epic} are monic, so is the other.
{\flushleft (3)} If any two of the vertical maps \ul{except $\beta,\gamma$ when $\gamma$ is not monic} are epic, so is the other.
\end{crl}
\begin{proof}
This is an immediate consequence of the snake lemma, since we have an exact sequence
\bea
0\ra\ker\al\ra\ker\beta\ra\ker\gamma\sr{\delta}{\ra}\coker\al\ra\coker\beta\ra\coker\gamma\ra0,\nn
\eea
and so, for example, if $\al,\gamma$ are isomorphisms, then we get the exact sequence
\[
0\ra 0\ra\ker\beta\ra 0\ra0\ra\coker\beta\ra 0\ra0,
\]
which shows $\beta$ is also an isomorphism.
\end{proof}

\begin{crl}[\textcolor{blue}{\index{Five lemma! full version}{Five lemma}}]\label{Long5Lmm}
Consider any \ul{commutative diagram} (of $R$-modules) of the following form, (i) with \ul{exact rows}, (ii) $f_1$ an \ul{epimorphism}, and (iii) $f_5$ a \ul{monomorphism}:
\begin{figure}[H]
\centering
\adjustbox{scale=0.8}{
\begin{tikzcd}
  &  &  &  &  & 0\ar[d] &  \\
  & {X}_1\ar[d,"{f}_1"]\ar[r] & X_2\ar[d,"f_2"]\ar[r] & X_3\ar[d,"f_3"]\ar[r] & X_4\ar[d,"f_4"]\ar[r] & {X}_5\ar[d,"{f}_5"] &  \\
  & {Y}_1\ar[d]\ar[r] & Y_2\ar[r] & Y_3\ar[r] & Y_4\ar[r] & {Y}_5 & \\
  & 0 &  &  &  &  &
\end{tikzcd}}
\end{figure}
If $f_2$ and $f_4$ are isomorphisms, then so is $f_3$.
\end{crl}
\begin{proof}
By hypotheses, we have the exact commutative diagram (a) below.
\[(a)\adjustbox{scale=0.8}{
\begin{tikzcd}
   &  & 0\ar[d,dashed] &  & 0\ar[d,dashed] & 0\ar[d,dashed] &  \\
  & {X}_1\ar[d,"{f}_1"]\ar[r,"{u}"] & X_2\ar[d,"f_2"]\ar[r,"f"] & X_3\ar[d,"f_3"]\ar[r,"g"] & X_4\ar[d,"f_4"]\ar[r,"{v}"] & {X}_5\ar[d,"{f}_5"] &  \\
  & {Y}_1\ar[d,dashed]\ar[r,"{u}'"] & Y_2\ar[d,dashed]\ar[r,"f'"] & Y_3\ar[r,"g'"] & Y_4\ar[d,dashed]\ar[r,"{v}'"] & {Y}_5 &  \\
   & 0 & 0 &  & 0 &  &
\end{tikzcd}}~~~~(b)~~
\adjustbox{scale=0.8}{
\begin{tikzcd}
   & 0\ar[d,dashed] & 0\ar[d,dashed] & 0\ar[d,dashed] &  \\
 0\ar[r,dashed]  & \coim f\ar[d,"\wt{f}_2"]\ar[r,"\ol{f}"] & X_3 \ar[d,"f_3"]\ar[r,"\pi_g"] & \im g \ar[d,"\wt{f}_4"]\ar[r,dashed]& 0 \\
0 \ar[r,dashed] & \coim f'\ar[d,dashed]\ar[r,"\ol{f}'"] & Y_3\ar[d,dashed]\ar[r,"\pi_{g'}"] & \im g'\ar[d,dashed]\ar[r,dashed] & 0 \\
   & 0 & 0 & 0 &  \\
\end{tikzcd}}\]

Using the first isomorphism theorem, we can factor $f,f',g,g'$ as
 \bea
 && f=\ol{f}\pi_f,~~~~f'=\ol{f}'\pi_{f'},~~~~g=\ol{g}\pi_g,~~~~g'=\ol{g}'\pi_{g'},\nn
 \eea
and get the commutative diagram (b) above with exact rows (and columns being complexes). The maps (which are well-defined as can be checked easily) are given by
\bea
&&\ol{f}(x+\ker f):=f(x),~~~~\ol{f}'(x+\ker f'):=f'(x)\nn\\
&&\wt{f}_2(x+\ker f):=f_2(x)+\ker f',~~~~\wt{f}_4(x):=f_4(x).\nn
\eea
As shown below (using commutativity and hypotheses), $\wt{f}_2$ and $\wt{f}_4$ are isomorphisms, and so by the short five lemma (Lemma \ref{Short5Lmm}), we get the desired result.
{\flushleft\ul{$\wt{f}_2$ is an isomorphism}:} $\wt{f}_2$ is surjective because $f_2$ is surjective, and injective because
\bea
&& x+\ker f\in \ker\wt{f}_2~~\Ra~~f_2(x)\in\ker f'=\im u',\nn\\
&&~~\Ra~~f_2(x)=u'(y)\sr{(s1)}{=}u'(f_1(x_1))=f_2u(x_1),\nn\\
&&~~\sr{(s2)}{\Ra}~~x=u(x_1)\in\im u=\ker f,\nn\\
&&~~\Ra~~x+\ker f=0,\nn
\eea
where step (s1) is due to surjectivity of $f_1$, and step (s2) due to injectivity of $f_2$.

Alternatively, observe that $f_2u=u'f_1$ implies ~$\coim f'=\coker u'=\coker(u'f_1)=\coker(f_2u)=f_2(\coker u)=f_2(\coim f)$.

{\flushleft\ul{$\wt{f}_4$ is an isomorphism}:} Since $f_5 v=v'f_4$, we have ~$\im g=\ker v=\ker(f_5v)=\ker(v'f_4)=f_4^{-1}(\ker v')=f_4^{-1}(\im g')$.

\end{proof}

\section{Object/morphism lifting property (resolvability) of $R$-modules}
\begin{thm}[\textcolor{blue}{\index{Baer criterion}{Baer criterion} for module-injectivity}]
A module $_RE$ is injective $\iff$ every $R$-homomorphism $\al:{}_RI\subset R\ra E$ (for any left ideal $_RI\subset R$) extends to an $R$-homomorphism $\beta:R\ra E$.
\[\adjustbox{scale=0.9}{\bt
0 \ar[r]& I\ar[d,"\al"]\ar[rr,hook] && R\ar[dll,dashed,"\exists~\beta"]\\
   & E &&
\et}\]
\end{thm}
\begin{proof}
{\flushleft($\Ra$)}: (For brevity, $R$-map will stand for $R$-homomorphism). If $_RE$ is injective, it is clear (by the definition of injectivity) that every $R$-map $\al:{}_RI\subset R\ra E$ extends to an $R$-map $\beta:R\ra E$.
{\flushleft($\La$)}: Assume every $R$-map $\al:{}_RI\subset R\ra E$ can be lifted to an $R$-map $\beta:R\ra E$. Since in every exact sequence $0\ra A\sr{f}{\ral} B$, we have $A\cong f(A)\subset B$, it is enough to show that given any R-module $M$, an $R$-map $u:{}_RN\subset M\ra E$ can always be lifted to an $R$-map $v:M\ra E$. So, fix an $R$-module $M$ and consider an $R$-map ${u}:{}_RN\subset M\ra E$. Let the set
\bea
P:=\big\{(X,{v}_X)~|~N\subset{}_RX\subset M,~{v}_X:X\ra E,~{v}_X|_N={u}\big\},~~~~(N,{u})\in P,\nn
\eea
be a poset with ordering given by ~``$(X,{v}_X)\leq (X',{v}_{X'})$ if $X\subset X'$ and ${v}_{X'}|_X={v}_X$''.
\[\adjustbox{scale=0.9}{
\begin{tikzcd}
0 \ar[r] & N \ar[dd, "{u}"'] \ar[rr,hook] && X \ar[ddll,dashed,"{v}_X"' pos=0.5] \ar[rr,hook] && X' \ar[ddllll,dashed,"{v}_{X'}"' pos=0.45] \ar[r,hook] & \cdots \ar[r,hook]& M \ar[ddllllll,dashed,"{v}"]  \\
  & & & & & &\\
 & E & & & &  &
\end{tikzcd}}\]
Then any chain $\{(X_i,{v}_{X_i})\}_{i\in I}$ in $P$ has an upper bound $(\cup X_i,{v}_{\cup X_i})\in{P}$, where ${v}_{\cup X_i}(x):={v}_{X_i}(x)$ if $x\in X_i$. By Zorn's lemma, ${P}$ contains a maximal element $(Y,{v}_Y)$. Next, we will show that $Y=M$.

\[\adjustbox{scale=0.9}{
\begin{tikzcd}
  & I\ar[ddd,bend right=40,near end,"\al"description]\ar[drr,near start,"\times m"description]\ar[rr,hook] &&R\ar[dddll,bend right=30,near end,"\beta"description]  &&   & \\
0 \ar[r] & N \ar[dd, "{u}"'] \ar[rr,hook] && Y \ar[ddll,dashed,"{v}_Y"' pos=0.5] \ar[rr,hook] && Y+Rm \ar[ddllll,dashed,"{v}_{Y+Rm}"' pos=0.3] \ar[r,hook]& M \ar[ddlllll,dashed,"{v}"]  \\
  & & & &  &\\
 & E & & &  &
\end{tikzcd}}\]
Suppose $Y\subsetneq M$. Let $m\in M\backslash Y$. Let $I\lhd R$ be any nonzero ideal such that $Im\subset Y$, e.g., let $I:=\{a\in R:am\in Y\}=\{a\in R:a(m+Y)=Y~\txt{in}~M/Y\}$. Consider the $R$-map $\al:I\sr{\times m}{\ral}Y\sr{v_Y}{\ral}E,~a\mapsto v_Y(am)$, which is well defined since $Im\subset Y$. Then by hypotheses, $\al$ lifts to an $R$-map $\beta:R\ra E$ such that $\beta|_I=\al$. Now consider the pair $(Y+Rm,v_{Y+Rm})$, where $v_{Y+Rm}$ is the well-defined $R$-map
\bea
v_{Y+Rm}:Y+Rm\ra E,~~y+rm\mapsto v_Y(y)+\beta(r),\nn
\eea
where if $y+rm=y'+r'm$, then $y-y'=(r'-r)m\in Y$, and so {\small $v_Y(y)-v_Y(y')=v_Y(y-y')=v_Y((r'-r)m)=\al(r'-r)=\beta(r'-r)=\beta(r')-\beta(r)$}, i.e., $v_{Y+Rm}$ is indeed well-defined since
\[
 y+rm=y'+r'm~~\Ra~~v_{Y+Rm}(y+rm)=v_{Y+Rm}(y'+r'm).
 \]
Then $v_{Y+Rm}|_Y=v_Y$ and $v_{Y+Rm}|_N=v_Y|_N=u$, and so {\small $(Y,v_Y)<(Y+Rm,v_{Y+Rm})\in P$}, which contradicts the maximality of $(Y,v_Y)$.
\end{proof}

\begin{prp}[\textcolor{blue}{An injective module over an ID is divisible}]\label{InjImpDiv}
Let $R$ be an ID (i.e., integral domain). If $_RE$ is an injective module, then $E$ is divisible. (\blue{footnote}\footnote{Recall that a module $_RM$ is divisible if $M\subset rM$ (i.e., $f_r:M\ra M,~m\mapsto rm$ is surjective) for any $0\neq r\in R$.})
\end{prp}
\begin{proof}
    Suppose $E$ is injective. Let $0\neq c\in R$ and $x\in E$. Consider the $R$-map $\al:Rc\ra E,~rc\mapsto rx$, which is well-defined because $R$ is a domain ($rc=r'c$ $\Ra$ $r=r'$, $\Ra$ $rx=r'x$). Then by the injectivity of $E$, there is an $R$-map $\beta:R\ra E$ such that $\beta|_{Rc}=\al$, and so $x=\al(c)=\beta(c)=c\beta(1)$. Hence $E$ is divisible.
\[\adjustbox{scale=0.9}{\bt
0 \ar[r] & Rc \ar[d, "\al"'] \ar[rr,hook] && R \ar[dll,dashed,"\beta"]  \\
 & Rx\subset E &&
\et}\]
\end{proof}

\begin{thm}[\textcolor{blue}{Injectivity criterion for divisible modules over a PID}]\label{InjIffDiv}
Let $R$ be a PID. Then $_RE$ is injective $\iff$ divisible. (Hence, divisible modules over PIDs, such as $_\Integer\Rational$ and $_\Integer(\Rational/\Integer)$, are injective.)
\end{thm}
\begin{proof}
($\Ra$) If $E$ is injective, then $E$ is divisible by Proposition \ref{InjImpDiv}. ($\La$) Conversely, assume $E$ is divisible. Let $\al:I\ra E$, where $0\neq I\lhd R$. We know $I=Rc$ for some $0\neq c\in R$, and that $\al$ is determined completely by $x:=\al(c)\in E$. Since $E$ is divisible, $\al(c)=x=cy$ for some $y\in E$. Setting
\bea
\beta:R\ra E,~r\mapsto ry,~~~~\txt{i.e.,}~~\beta(1):=y,\nn
\eea
we have $\al=\beta|_I$, since $\beta(c)=c\beta(1)=cy=x=\al(c)$. Hence $E$ is injective.
\[\adjustbox{scale=0.9}{\bt
0 \ar[r] & Rc \ar[d, "\al"'] \ar[rr,hook] && R \ar[dll,dashed,"\beta"]  \\
 & Rx\subset E &&
\et}\]
\end{proof}

\begin{lmm}[\textcolor{blue}{$R$-mod has enough projectives and enough injectives}]\label{EnProjInjLmm}
Let $M$ be an $R$-module.
\bit[leftmargin=0.9cm]
\item[(i)] There exists a projective $R$-module $P$ and an epimorphism $g:P\twoheadrightarrow M$.
\item[(ii)] There exists an injective $R$-module $I$ and a monomorphism $f:M\hookrightarrow I$.
\eit
\end{lmm}
\begin{proof}
{\flushleft (i)}: Let $P:=R^{\langle M\rangle}\cong\oplus_{m\in M}R$ be the free (hence projective) $R$-module on $M$ and $g:R^{\langle M\rangle}\ra M$, $(r_m)_{m\in M}\mapsto\sum_{m\in M} r_mm$.

{\flushleft (ii)}: For each $m\in M$, let $D_m=\left\{
                                          \begin{array}{ll}
                                            \Rational, & \hbox{if}~~|m|=\infty \\
                                            {\Rational\over\Integer}, & \hbox{if}~~|m|<\infty
                                          \end{array}
                                        \right\}$, where $|m|$ is the order of $m$ in the group $_\Integer M$, and consider the associated $\Integer$-homomorphism
\[
h_m:\Integer m\ra D_m,~am\mapsto \left\{
                                          \begin{array}{ll}
                                            a, & \hbox{if}~~D_m=\Rational \\
                                            {a\over |m|}+\Integer, & \hbox{if}~~D_m={\Rational\over\Integer}
                                          \end{array}
                                        \right\},
\]
which is a \emph{well-defined} $\Integer$-homomorphism because if $am=a'm$ then $(a-a')m=0$, which means either
\[
a-a'=0~~(\txt{if $|m|=\infty$})~~~~\txt{or}~~~~a-a'=c|m|~~(\txt{if $|m|<\infty$})~~\txt{for some}~~c\in\Integer.
\]
Since $_\Integer D_m$ is injective (as a divisible module over a PID), by the Baer criterion, there exists a $\Integer$-homomorphism $i_m:M\ra D_m$ such that $h_m=i_m|_{\Integer m}$.
\[\adjustbox{scale=0.9}{
\begin{tikzcd}
  0\ar[r] & \Integer m \ar[rr,hook] \arrow[d, "h_m"'] & & M \arrow[dll,dashed,"i_m"]\\
  & D_m &&
\end{tikzcd}}\]
Let $D:=\prod\limits_{m\in M}D_m=\Big\{(q_m)_{m\in M}:q_m\in D_m\Big\}$, where $_\Integer D$ is injective (as a product of injective modules). Then we get the following inclusion:
\[
i:M\hookrightarrow D,~m\mapsto \big(i_m(m')\big)_{m'\in M}.
\]
Since $0\ra{}_\Integer M\sr{i}{\hookrightarrow}{}_\Integer D$ is exact and $Hom_\Integer(R,-)$ is left exact, we get the exact sequence
\[
0\ra Hom_\Integer(R,M)\sr{i_\ast}{\hookrightarrow} Hom_\Integer(R,D).
\]
With the isomorphism $M\cong Hom_R(R,M)$ given by $\phi:Hom_R(R,M)\ra M,~f\mapsto f(1)$, we get
\[
_RM\cong Hom_R(R,M)\subset_R Hom_\Integer (R,M)\sr{i_\ast}{\hookrightarrow} Hom_\Integer(R,D),
\]
where $\subset_R$ denotes ``$R$-submodule'' and $Hom_\Integer(R,X)$ is viewed as an $R$-module via $(r\cdot f)(a):=f(ar)$ for $f\in Hom_\Integer(R,X)$, $a,r\in R$.

If $0\ra A\ra B\ra C\ra 0$ is a SES of R-modules, then by taking ~$Hom_R\Big(-,~Hom_\Integer(R,D)\Big)$~ we get
\[\adjustbox{scale=0.9}{\bt
0\ar[r] & Hom_R\Big(C,~Hom_\Integer(R,D)\Big)\ar[d,"\cong"]\ar[r] & Hom_R\Big(B,~Hom_\Integer(R,D)\Big)\ar[d,"\cong"]\ar[r] & Hom_R\Big(A,~Hom_\Integer(R,D)\Big)\ar[d,"\cong"]\ar[r] & 0\\
0\ar[r] & Hom_\Integer\Big(R\otimes_R C,D\Big)\ar[d,"\cong"]\ar[r] & Hom_\Integer\Big(R\otimes_R B,D\Big)\ar[d,"\cong"]\ar[r] & Hom_\Integer\Big(R\otimes_R A,D\Big)\ar[d,"\cong"]\ar[r] & 0\\
0\ar[r] & Hom_\Integer\Big(C,D\Big)\ar[r] & Hom_\Integer\Big(B,D\Big)\ar[r] & Hom_\Integer\Big(A,D\Big)\ar[r] & 0
\et}\]in which exactness of the last row (by the injectivity of $_\Integer D$) implies exactness of the first row, which implies $Hom_\Integer(R,D)$ is an injective R-module. (\blue{footnote}\footnote{More generally, if $R$ is a commutative ring, $A$ an $R$-algebra, and $E$ an injective $R$-module, then $_AE':=Hom_R(A,E)$ is an injective $A$-module. To prove this, consider an exact sequence of $A$-modules $0\ra M_1\ra M_2\ra M_3\ra 0$, apply $Hom_A(-,Hom_R(A,E))$, and then proceed as before.}).

Finally, the following composition is by construction an injective $R$-homomorphism:
\[
_RM\cong Hom_R(R,M)\subset_R Hom_\Integer (R,M)\sr{i_\ast}{\hookrightarrow} Hom_\Integer(R,D). \qedhere
\]
\end{proof}

\begin{prp}[\blue{Existence of a projective resolution in R-mod}]
Given an R-module $M$, a projective resolution of $M$ exists.
\end{prp}
\begin{proof}
Consider a free (hence projective) $R$-module $F$ and a surjective $R$-homomorphism $\pi:F\twoheadrightarrow M$. Take this as the starting point with $\pi_0:=\pi$, and then repeat/iterative the process to construct an exact sequence as follows:
\begin{figure}[H]
  \centering
\adjustbox{scale=0.9}{ \bt
 \cdots\ar[rr,dashed,"d_3=i_3\circ\pi_3"]\ar[dr,two heads,"\pi_3"] & & F_2\ar[dr,two heads,"\pi_2"]\ar[rr,dashed,"d_2=i_2\circ\pi_2"]  &   & F_1\ar[dr,two heads,"\pi_1"]\ar[rr,dashed,"d_1=i_1\circ\pi_1"]  &   & F_0=F\ar[r,two heads,"\vep=\pi_0"] & M\ar[r] & 0 \\
             &\ker\pi_2\ar[ur,hook,"i_3"]\ar[dr] &  & \ker\pi_1\ar[ur,hook,"i_2"]\ar[dr]  &   & \ker\vep\ar[ur,hook,"i_1"]\ar[dr]  &  &  &  \\
       \cdots\ar[ur] &  & 0\ar[ur] &   & 0\ar[ur] &   & 0 &  &
\et}
\caption{~~Construction of a projective resolution}\label{cpr-construction}
\end{figure}
where $F_n$ is free, with $\pi_n$ sending the basis for $F_n$ to some generators of $\ker\pi_{n-1}$. We have
\[
\im d_{n+1}=\im(i_{n+1}\circ\pi_{n+1})\cong\im(\pi_{n+1})=\ker\pi_n\cong\ker(i_n\circ\pi_n)=\ker d_n. \qedhere
\]
\end{proof}

\begin{prp}[\textcolor{blue}{Existence of an Injective resolution in R-mod}]
Given an R-module $M$, an injective resolution of $M$ exists.
\end{prp}
\begin{proof}
Include $M$ in an injective $R$-module $I^0$ via $M\sr{\xi}{\hookrightarrow}I^0$, and then include $\coker\xi$ (from $I^0\sr{\pi^0}{\twoheadrightarrow}\coker\xi$) into an injective $R$-module $I^1$ via $I^0\sr{i^0}{\hookrightarrow}I^1$. Let $d^0:=i^0\pi^0$. Repeat/iterate the process as follows:

\begin{figure}[H]
  \centering
\adjustbox{scale=0.9}{ \bt
 0\ar[r] & M\ar[r,hook,"\xi"] & I^0\ar[rr,dashed,"d^0:=i^0\circ\pi^0"]\ar[dr,two heads,"\pi^0"] & & I^1\ar[dr,two heads,"\pi^1"]\ar[rr,dashed,"d^1:=i^1\circ\pi^1"]  &   & I^2\ar[dr,two heads,"\pi^2"]\ar[rr,dashed,"d^2:=i^2\circ\pi^2"]  &   & \cdots \\
 & &       & \coker\xi\ar[ur,hook,"i^0"]\ar[dr] &  & \coker i^0\ar[ur,hook,"i^1"]\ar[dr]  &   & \coker i^1\ar[ur,hook,"i^2"]\ar[dr]  &   \\
  & &     0\ar[ur] &  & 0\ar[ur] &   & 0\ar[ur] &   & \cdots
\et}
\caption{~~Construction of an injective resolution}\label{cir-construction}
\end{figure}
\[
\im d^{n-1}=\im(i^{n-1}\pi^{n-1})\cong\im\pi^{n-1}=\coker i^{n-2}=\ker\pi^n\cong\ker(i^n\pi^n)=\ker d^n.\qedhere
\]
\end{proof}

\section{Criteria for Split-SES's, Projectivity, and Injectivity}
\begin{dfn}[\textcolor{blue}{\index{Split! SES of modules}{Split SES of modules}}]
We say a SES of $R$-modules ~$0\ra A\sr{f}{\ral}B\sr{g}{\ral}C\ra 0$ is split if $\im f$ is a direct summand of $B$ (i.e., $B=\im(f)\oplus N=\ker(g)\oplus N$, for some $_RN\subset B$).
\end{dfn}

We will see in Corollary \ref{SplitSeqCrit2} below that a SES $0\ra A\ra B\ra C\ra 0$ is split (or splits) $\iff$ there exists an isomorphism $\theta:B\ra A\oplus C$ giving an exact commutative diagram:
\[\adjustbox{scale=0.9}{\bt
0\ar[r] & A\ar[d,equal]\ar[r,"f"] & B\ar[d,hook,two heads,"\cong"',"\theta"]\ar[r,"g"] & C\ar[d,equal]\ar[r] & 0 \\
0\ar[r] & A\ar[r,"i"] & A\oplus C\ar[r,"p"] & C\ar[r] & 0
\et}\]

\begin{rmk}[\textcolor{magenta}{Caution}]
{\flushleft (i)} A sequence $0\ra A\sr{f}{\ral} B\sr{g}{\ral} C\ra 0$ that is split in the sense $B=\im f~\oplus N$ might not be exact. For example, consider the sequence $0\ra\Integer\sr{f}{\ral}\Integer\oplus(\Integer_2)^\Natural\sr{g}{\ral}(\Integer_2)^\Natural\ra 0$ (in $\Integer$-mod = Ab), where $\Integer_2:={\Integer\over 2\Integer}=\Integer(1+2\Integer)$ and
\bea
f(n):=(n,0,0,\cdots),~~~~g(n,m_1,m_2,\cdots):=(n+2\Integer,m_1,m_2,\cdots).\nn
\eea

{\flushleft (ii)} If a SES $0\ra A\ra B\ra C\ra 0$ splits, then $B=\im(f)\oplus N\cong A\oplus C$, since $A\cong im(f)=\ker g$, and $C=im(g)\cong{B\over\ker g}\cong N$. The converse is false, i.e., a SES $0\ra A\ra B\ra C\ra 0$ such that $B\cong A\oplus C$ might not split. This is because, with an isomorphism $\phi:A\oplus C\ra B$, we only know $B=\phi(A)\oplus\phi(C)$, where $\phi(A)\cong A\cong\im (f)$. But $\phi(A)\cong\im f$ alone $\not\Ra$ $\phi(A)=\im (f)$, (\blue{footnote}\footnote{That is, $\phi(A)\subset B$ and $\im (f)\subset B$ are isomorphic but, in general, can be different submodules of $B$}), and so
\bea
B=\phi(A)\oplus\phi(C)~~\txt{and}~~\phi(A)\cong\im f~~~~\not\Ra~~~~B=\im (f)\oplus N~~\txt{for some}~~N\subset B.\nn
\eea
For example, consider the sequence $0\ra\Integer\sr{f}{\ral}2\Integer\oplus\Rational\sr{g}{\ral}\Rational\ra 0$ (in $\Integer$-mod = Ab), where
\bea
f(n):=(2n,0),~~~~g(n,q):=q.\nn
\eea
\end{rmk}

\begin{prp}[\textcolor{blue}{\index{Split! SES criteria (module)}{Split SES criteria}}]\label{SplitSeqCrit1}
Given a SES (of $R$-modules) ~$0\ra A\sr{f}{\ral}B\sr{g}{\ral}C\ra 0$,~ the following are equivalent.
\begin{enumerate}[leftmargin=0.7cm]
\item[(1)] The SES splits.
\item[(2)] There exists an R-homomorphism $g':C\ra B$ such that $g\circ g'=id_C$. (I.e., $g$ has a right inverse).
\item[(3)] There exists an R-homomorphism $f':B\ra A$ such that $f'\circ f=id_A$. (I.e., $f$ has a left inverse).
\item[(4)] There exist R-homomorphisms $g':C\ra B$, $f':B\ra A$ such that  $g\circ g'=id_C$, $f'\circ f=id_A$.
\item[(5)] There exist R-homomorphisms $g':C\ra B$, $f':B\ra A$ such that ~$0\la A\sr{f'}{\lal}B\sr{g'}{\lal}C\la 0$~ is exact (or just a complex) and $f\circ f'+g'\circ g=id_B$.
\end{enumerate}
\end{prp}

\begin{proof}~
Consider the diagram below (in which $C\cong{B\over\ker g}\cong C'$, and $A\cong\im f=\ker g=K$):
\[\adjustbox{scale=0.9}{\bt
        &          &&     && C'\ar[dll,shift right,"i"']\ar[d,"g\circ i"]  &  \\
0\ar[r] & A\ar[d,"p\circ f"']\ar[rr,"f"]  && \ar[ll,shift left=2,dashed, near end,"f'"]B\ar[dll,shift left,"p"]\ar[rr,"g"']   && \ar[ll,shift right=2,dashed,near start,"g'"']C\ar[r] & 0 \\
        & K        &&     &&   &
\et}~~~~\substack{B=K\oplus C'~~~ \\~\\ g'=i\circ(g\circ i)^{-1} \\~\\ f'=(p\circ f)^{-1}\circ p}
\]
 along with
\[\adjustbox{scale=0.9}{
\begin{tikzcd}
 0 \ar[r] & A \ar[r,shift left=1,"f"] & B \ar[l,dashed,shift left=1,"f'"] \ar[r,shift left=1,"g"]  &  C \ar[l,dashed,shift left=1,"g'"] \ar[r] & 0
\end{tikzcd}}\]
\bit[leftmargin=0.7cm]
\item\ul{(1)$\Ra$(2)}: (1) implies $B=\ker g\oplus C'$. Define $g'$ (automatically satisfying $g\circ g'=id_C$) by
\bea
g':C\ra B=\ker g\oplus C',~c\mapsto (0,c')\in g^{-1}(c),~~\txt{i.e.,}~~g(0,c')=c.\nn
\eea
To show $g'$ is a homomorphism, observe that for $r\in R$, if $g'(r c)=(0,c'')\in g^{-1}(r c)$, then
\bea
g(0,c'')=r c=r g(0,c')=g(0,r c'),~~\sr{(s)}{\Ra}~~c''=r c',~~\Ra~~g'(r c)=r g'(c),\nn
\eea
where step (s) holds because $g|_{\{0\}\oplus C'}$ is an isomorphism. Similarly, if $g'(c_1+c_2)=(0,c'')\in g^{-1}(c_1+c_2)$, then with $g'(c_1)=(0,c_1')\in g^{-1}(c_1)$ and $g'(c_2)=(0,c_2')\in g^{-1}(c_2)$, we have
\bea
&&g(0,c'')=c_1+c_2=g(0,c_1')+g(0,c_2')=g\big((0,c_1')+(0,c_2')\big)=g(0,c_1'+c_2'),\nn\\
&&~~\sr{(s)}{\Ra}~~c''=c_1'+c_2',~~\Ra~~g'(c_1+c_2)=g'(c_1)+g'(c_2).\nn
\eea

\ul{Alternatively}, it suffices to set $g':=i\circ(g\circ i)^{-1}$, where
\bea
i:C'\ra \ker g\oplus C',~~c'\mapsto (0,c').\nn
\eea

\item\ul{(2)$\Ra$(1)$\Ra$(3)}: Given (2), we can show that $B=\ker g\oplus g'(g(B))=f(A)\oplus g'(C)$. Indeed,
\bea
\im (g')\cap\ker g=\{b:b=g'(c),~g(b)=0\}=\{b:b=g'(c),c=gg'(c)=g(b)=0\}=0,\nn
\eea
and any $b\in B$ can be written as
    \bea
    b=b-g'g(b)+g'g(b)\in \ker g+\im g'=\im f+\im g'.\nn
    \eea
Define $f':B\ra A$ by $f'\big(f(a),g'(c)\big)=a$, for all $f(a)\in \im (f)$, $g'(c)\in \im (g')$. Then $f'$ is well-defined because $f$ is injective, and $f'$ is also easily seen to be a homomorphism because $f$ is a homomorphism.

\ul{Alternatively}, it suffices to set $f':=(p\circ f)^{-1}\circ p$, where
    \bea
    p:\ker g\oplus g'(C)\ra \ker g,~~(k,c')\mapsto k.\nn
    \eea

\item\ul{(3)$\Ra$(1)}: If (3) holds, we can check that $B=\ker f'\oplus ff'(B)=\ker f'\oplus f(A)=\ker f'\oplus\ker g$. Indeed,
\bea
\ker f'\cap\im (f)=\{b:f'(b)=0,b=f(a)\}=\{b:f'(b)=0,a=f'f(a)=f'(b)=0\}=0,\nn
\eea
and any $b\in B$ can be written as
    \bea
    b=b-ff'(b)+ff'(b)\in\ker f'+\im f=\ker f'+\ker g.\nn
    \eea

\item\ul{(2) or (3) $\iff$ (4)}: This is clear since we have already shown that $(2)\iff(3)$.

\item\ul{(4)$\Ra$(5)}: Since the equivalences (1) $\iff$ (2) $\iff$ (3) are clear, we only need to check that the two relations $g\circ g'=id_C$, $f'\circ f=id_A$ imply $f\circ f'+g'\circ g=id_B$. Indeed, because $B=\im (f)\oplus \im (g')$,
\bea
&& f'f=id_A,~gg'=id_C~~\Ra~~ff'f=f,~g'gg'=g',~~\Ra~~ff'=id_{\im (f)},~g'g=id_{\im (g')},\nn\\
&&~~\Ra~~ff'+g'g=id_{\im (f)}+id_{\im (g')}\sr{(s1)}{=}id_{\im (f)+\im (g)}=id_B.\nn
\eea
Step (s1) holds because for any $b\in B$, we have $a\in A$, $c\in C$ such that
\bea
b=f(a)+g'(c)=ff'f(a)+g'gg'(c)\sr{(s2)}{=}(ff'+g'g)\big(f(a)+g'(c)\big)=(ff'+g'g)(b),\nn
\eea
where step (s2) holds because $f'g'=0$ and $gf=0$.

\item\ul{(5)$\Ra$(4)}: Observe that applying $g$ on $f\circ f'+g'\circ g=id_B$ we get $0+gg'g=g$, which implies $gg'=id_C$ (by the surjectivity of $g$). Similarly, applying $f'$ on $f\circ f'+g'\circ g=id_B$ we get $f'ff'+0=f'$, which implies $f'f=id_A$ (by the surjectivity of $f'$).
\eit
\end{proof}

\begin{crl}[\textcolor{blue}{\index{Split! SES criterion (module)}{Split SES criterion}}]\label{SplitSeqCrit2}
A sequence $0\ra A\sr{f}{\ral}B\sr{g}{\ral}C\ra 0$ is split-exact (i.e., both exact and split) $\iff$ there exist R-homomorphisms $g':C\ra B$, $f':B\ra A$ such that $g\circ g'=id_C$, $f'\circ f=id_A$, and $f\circ f'+g'\circ g=id_B$.

Consider the canonical inclusions $q_A:A\ra A\oplus C$, $q_C:C\ra A\oplus C$ and projections $p_A:A\oplus C\ra A$, $p_C:A\oplus C\ra C$. Then split-exactness of the SES is equivalent to specifying the isomorphism $\theta:=q_Af'+q_Cg:B\ra A\oplus C$ with inverse $\theta^{-1}:=fp_A+g'p_C:A\oplus C\ra B$, as in the following exact commutative diagram:
\[\adjustbox{scale=0.9}{\bt
0\ar[r] & A\ar[d,equal]\ar[r,"f"] & B\ar[d,hook,two heads,"\cong"',"\theta"]\ar[r,"g"] & C\ar[d,equal]\ar[r] & 0 \\
0\ar[r] & A\ar[r,"i"] & A\oplus C\ar[r,"p"] & C\ar[r] & 0
\et}\]
\end{crl}

\begin{lmm}[\textcolor{blue}{\index{Projectivity criteria (module)}{Projectivity criteria}}]\label{ModProjCrit1}
Let $P$ be an $R$-module. Then the following are equivalent.
\begin{enumerate}
\item[(a)] $P$ is projective.
\item[(b)] Every SES ~$0\ra A\sr{f}{\ral}B\sr{g}{\ral}P\ra 0$~ splits.
\item[(c)] $P$ is a direct summand of a free $R$-module $F$ (i.e., $F=P\oplus P'$ for some ${}_RP'\subset F$).
\end{enumerate}
\end{lmm}
\begin{proof}~
\bit[leftmargin=0.7cm]
\item[]\ul{(a)$\Ra$(b)}: Because $P$ is projective, if $B\sr{g}{\ral}P$ is epic then there is $g':P\ra B$ such that $id_P=gg'$.

\[\adjustbox{scale=0.9}{
\begin{tikzcd}
 && P \ar[dll,dashed,"g'"'] \ar[d, "id_P"]   & \\
B \ar[rr,two heads,"g"]  && P \ar[r] & 0
\end{tikzcd}}\]
\item[]\ul{(b)$\Ra$(c)}: We know the R-module $P$ is a homomorphic image of the free module $F_1:=R^{\langle P\rangle}$ on $P$, and so we have a split exact sequence $0\ra\ker g\ra F_1\sr{g}{\ra} P\ra 0$, where
\bea
\textstyle F_1=\ker g~\oplus~N\cong \ker g\oplus P~~(=F),~~~\txt{since}~~N\cong {F_1\over\ker g}\cong P.\nn
\eea
\item[]\ul{(c)$\Ra$(a)}: Recall from Lemma \ref{FreeIsProj} that a free module is projective. If $F=P\oplus P'$ is free and $0\ra A\ra B\ra C\ra 0$ is exact, then applying $Hom_R(F,-)$, we get the exact sequence
\[\adjustbox{scale=0.8}{\bt
0\ar[r] & Hom_R(F,A)\ar[d,"\cong"]\ar[r] & Hom_R(F,B)\ar[d,"\cong"]\ar[r] & Hom_R(F,C)\ar[d,"\cong"]\ar[r] & 0\\
0\ar[r] & Hom_R(P,A)\oplus Hom_R(P',A)\ar[r] & Hom_R(P,B)\oplus Hom_R(P',B)\ar[r] & Hom_R(P,C)\oplus Hom_R(P',C)\ar[r] & 0
\et}\]
which implies {\footnotesize $0\ra Hom_R(P,A)\ra Hom_R(P,B)\ra Hom_R(P,C)\ra 0$} is exact, and so $P$ is projective. \qedhere
\eit
\end{proof}

\begin{rmk}[\textcolor{blue}{Alternative proof of (c)$\Ra$(a) above}]
If we were not already aware that a free module is projective, we could proceed as follows (an immediate corollary of which is ``Every free module is projective.''):

Let $F=P\oplus P'$ be free. Next, consider an exact sequence $B\sr{g}{\ral}C\ra0$ and any $R$-map $h:P\ra C$. Let $i:P\ra F$, $p:F\ra P$ be the canonical inclusion and projection respectively.
\[\adjustbox{scale=0.9}{
\begin{tikzcd}
 && F=P\oplus P' \ar[dd,shift left=1, "p"] \ar[ddddll,dotted,"\al"'] & \\
 && & \\
 && P \ar[uu,hook,shift left=1,"i"] \ar[dd, "h"] \ar[ddll,dashed,"\beta"']  & \\
 && & \\
B \ar[rr,two heads,"g"]  && C \ar[r] & 0
\end{tikzcd}}
\hspace{1cm}
\adjustbox{scale=0.9}{\begin{tikzcd}
X\ar[dddd,"\al|_X"']\ar[rr,hook] && F=P\oplus P' \ar[dd,shift left=1, "p"] \ar[ddddll,dotted,"\al"'] & \\
 && & \\
 && P \ar[uu,hook,shift left=1,"i"] \ar[dd, "h"] \ar[ddll,dashed,"\beta"']  & \\
 && & \\
B \ar[rr,two heads,"g"]  && C \ar[r] & 0
\end{tikzcd}}\]
Now pick a basis $X=\{x_k\}_{k\in K}$ for $F$. Then because $g$ is surjective, there is $b_k\in B$ such that $g(b_k)=hp(x_k)$, and so we can define a map $\al:F\ra B$ by
    \bea
   \al|_X:X\ra B,~x_k\mapsto b_k\in g^{-1}\big(hp(x_k)\big),~~\txt{i.e.,}~~g(b_k)=hp(x_k),\nn
    \eea
    and consider its restriction to $P$, ~~$\beta:=\al|_P=\al\circ i:P\ra B$. The maps $\al,\beta$ are well-defined due to the axiom of choice and the fact that $\al|_X$ determines $\al$ completely (by the properties of free modules). We can also check as before that $\al,\beta$ are $R$-linear.
    We have $g\beta=h$, because if $\sum_{k\in K} r_kx_k\in P$ then
{\footnotesize\[
\textstyle g\beta p\left(\sum_k r_kx_k\right)=g\al ip\left(\sum_k r_kx_k\right)=g\al\left(\sum_k r_kx_k\right)=\sum_k r_kg\al(x_k)=\sum_k r_kg(b_k)=\sum_k r_khp(x_k)=hp\left(\sum_k r_kx_k\right).\nn
\]}
\end{rmk}

\begin{thm}[\blue{Coproducts of projective modules}]\label{ProjModSumCrit}~
\begin{enumerate}
\item[(a)] If a direct sum $\mathop{\oplus}_{i\in I}M_i$ of $R$-modules is projective, then so is each direct summand $M_i$.
\item[(b)] A direct sum $\mathop{\oplus}_{i\in I}P_i$ of projective $R$-modules is projective.
(\blue{footnote}\footnote{\magenta{Caution}: It is known that a product $\prod_{i\in I}P_i$ of projective $R$-modules need \ul{not} be projective. For \ul{example}, with additional elaboration not of immediate importance to us, one can show that the product of free $\Integer$-modules ~$\prod_{i=1}^\infty \Integer$~ is not a free $\Integer$-module.})
\end{enumerate}
\end{thm}
\begin{proof}
{\flushleft (a)}: Assume $\mathop{\oplus}_{i\in I}M_i$ is projective. Then $(\mathop{\oplus}_{i\in I}M_i)\oplus Q$ is free for some module $Q$. Hence each $M_i$ is projective as a direct summand of a free module.
{\flushleft (b)}: Assume $P_i$ is projective for each $i\in I$. Then for each $i\in I$, $F_i=P_i\oplus N_i$ is free for some module $N_i$. Thus, $\mathop{\oplus}_{i\in I}P_i$ is projective since we have the free module
\[
\textstyle F=\mathop{\oplus}_{i\in I}F_i=\mathop{\oplus}_{i\in I}(P_i\oplus N_i)~\cong~\left(\mathop{\oplus}_{i\in I}P_i\right)\oplus\left(\mathop{\oplus}_{i\in I}N_i\right). \qedhere
\]
\end{proof}

\begin{rmk}[\textcolor{blue}{Alternative of the above proof}]~~
{\flushleft\ul{Alternative proof of (a)}}: Let $g:B\ra C$ be surjective and $h:M_j\ra C$.
\[\adjustbox{scale=0.9}{
\begin{tikzcd}
   && \oplus_iM_i  \ar[ddddll,dotted,"\vphi"'] \ar[dd,shift left=1,"p_j"] &   \\
   &&  & \\
   && M_j \ar[ddll,dashed,"\vphi q_j"'pos=0.2] \ar[uu,hook,shift left=1,"q_j"] \ar[dd,"h"] &  \\
   &&  & \\
 B \ar[rr,two heads,"g"] && C\ar[r] & 0
\end{tikzcd}}\]
Since $\oplus_{i\in I}M_i$ is projective, there is $\vphi:\oplus_{i\in I}M_i\ra B$ such that $hp_j=g\vphi$. So, {\small $h=h~id_{M_j}=hp_jq_j=g~\vphi q_j$}.

{\flushleft\ul{Alternative proof of (b)}}: Let $g:B\ra C$ be surjective and $h:\oplus_{i\in I}P_i\ra C$.
\[\adjustbox{scale=0.9}{
\begin{tikzcd}
   &&  P_j  \ar[ddddll,dotted,"\vphi_j"'] \ar[dd,hook,shift left=1,"q_j"] &   \\
   &&  & \\
   && \oplus_iP_i \ar[ddll,dashed,"\Sigma\vphi_jp_j"pos=0.4] \ar[uu,shift left=1,"p_j"] \ar[dd,"h"] &  \\
   &&  & \\
 B \ar[rr,two heads,"g"] && C\ar[r] & 0
\end{tikzcd}}\]
Since $P_j$ is projective, there is $\vphi_j:P_j\ra B$ such that $hq_j=g\vphi_j$. Applying $p_j$ on the right and summing,
\[
\textstyle h=h~id_{\oplus_iP_i}=\sum_j hq_jp_j=\sum_j g\vphi_jp_j=g\sum_j\vphi_j p_j,
\]
where the expression $\vphi=\sum\vphi_j p_j$ makes sense on $\oplus_{i\in I}P_i$ (as a finite sum) because for all $(m_i)_{i\in I}\in \oplus_{i\in I}P_i$, we have $m_i=0$ except for finitely many $i$. (Alternatively, simply define $\vphi$ by ~$\vphi\big((m_i)_{i\in I}\big):=\sum_{i\in I}\vphi_i(m_i)$.)
\end{rmk}

\begin{note}[\magenta{Caution}]
The analysis in the alternative proof of (b) above does not work for a product $\prod_{i\in I}P_i$ of projective modules because (for the general modules we are considering) the expression $\vphi=\sum\vphi_jp_j$ is well defined (as a finite sum) only on the submodule $\oplus_{i\in I}P_i$ of $\prod_{i\in I}P_i$, but not on all of $\prod_{i\in I}P_i$.
\end{note}

\begin{lmm}[\textcolor{blue}{\index{Injectivity criteria (module)}{Injectivity criteria}}]
Let $E$ be an $R$-module. Then the following are equivalent.
\begin{enumerate}
\item[(a)] $_RE$ is injective.
\item[(b)]  Every SES ~$0\ral E\sr{f}{\ral}B\sr{g}{\ral}C\ral 0$~ splits.
\item[(c)] $E$ is a direct summand of \ul{every} containing module $_RM\supset{}_RE$ (i.e.,, $M=E\oplus X$ for some $_RX\subset M$).
\end{enumerate}
\end{lmm}
\begin{proof}~
\bit[leftmargin=0.7cm]
\item[]\ul{(a)$\Ra$(b)}: Because $E$ is injective, if $E\sr{f}{\ral} B$ is monic then there is $f':B\ra E$ such that $id_E=f'f$.

\[\adjustbox{scale=0.9}{
\begin{tikzcd}
0 \ar[r] & E \ar[rr,hook,"f"] \arrow[d, "id_E"'] & & B \arrow[dll,dashed,"f'"]\\
 & E & &
\end{tikzcd}}\]
\item[]\ul{(b)$\Ra$(c)}: If $_RE\subset{}_RM$ then the sequence $0\ra E\hookrightarrow M\sr{\pi}{\ral}M/E\ra 0$ splits.
\item[]\ul{(c)$\Ra$(a)}: Let $J=E\oplus E'$ be injective (\blue{footnote}\footnote{This is possible by (c) and the fact (from Lemma \ref{EnProjInjLmm}) that every R-module can be embedded in an injective R-module.}). If $0\ra A\ra B\ra C\ra 0$ is exact, then applying $Hom_R(-,J)$, we get the exact sequence
\[\adjustbox{scale=0.8}{\bt
0\ar[r] & Hom_R(C,J)\ar[d,"\cong"]\ar[r] & Hom_R(B,J)\ar[d,"\cong"]\ar[r] & Hom_R(J,A)\ar[d,"\cong"]\ar[r] & 0\\
0\ar[r] & Hom_R(C,E)\oplus Hom_R(C,E')\ar[r] & Hom_R(B,E)\oplus Hom_R(B,E')\ar[r] & Hom_R(A,E)\oplus Hom_R(A,E')\ar[r] & 0
\et}\]
which implies {\footnotesize $0\ra Hom_R(C,E)\ra Hom_R(B,E)\ra Hom_R(A,E)\ra 0$} is exact, and so $E$ is injective.\qedhere

\eit
\end{proof}

\begin{thm}[\blue{Coproducts and products of injective modules}]\label{InjModSumCrit}~
\begin{enumerate}
\item[(a)] If a direct sum $\mathop{\oplus}_{i\in I}M_i$ of $R$-modules is injective, then so is each direct summand $M_i$.
\item[(b)] A product $\prod_{i\in I}E_i$ of injective $R$-modules is injective. (\blue{footnote}{\footnote{\magenta{Caution}: It is known that a direct sum $\mathop{\oplus}_{i\in I} E_i$ of injective $R$-modules need \ul{not} be injective.}})
\end{enumerate}
\end{thm}
\begin{proof}~
\begin{enumerate}[leftmargin=0.7cm]
\item[(a)] Let $f:A\ra B$ be an injective R-homomorphism and $h:A\ra E_j$ an R-homomorphism.
\[\adjustbox{scale=0.9}{
\begin{tikzcd}
 0\ar[r] & A \arrow[dd, "h"'] \ar[rr,hook,"f"]  & & B \arrow[ddll,dashed,"p_j\vphi"'] \arrow[ddddll,dotted,"\vphi"]\\
   & & &\\
   & E_j \ar[dd,hook,shift left=1, "q_j"] & &\\
   & & &\\
   & \oplus_iE_i \ar[uu,shift left=1, "p_j"]& & &
\end{tikzcd}}\]
If $\mathop{\oplus}_{i\in I}E_i$ is injective, then there exists $\vphi:B\ra \mathop{\oplus}_{i\in I}E_i$ such that $q_jh=\vphi f$. This implies
    \bea
    h=id_{E_j}h=p_jq_jh=p_j\vphi~f.\nn
    \eea

\item[(b)] Let $f:A\ra B$ be an injective R-homomorphism and $h:A\ra\prod_{i\in I}E_i$ an R-homomorphism.
    \[\adjustbox{scale=0.9}{
\begin{tikzcd}
 0\ar[r] & A \arrow[dd, "h"'] \ar[rr,hook,"f"]  & & B \arrow[ddll,dashed,"\vphi"'] \arrow[ddddll,dotted,"\vphi_j"]\\
   & & &\\
   & \prod_iE_i \ar[dd, "p_j"'] & &\\
   & & &\\
   & E_j & & &
\end{tikzcd}}\]
Then for each $j\in I$, ~$p_jh:A\ra E_j$~ is an R-homomorphism since for any $r\in R$ and $a\in A$,
\bea
p_jh(ra)=p_j\big(h(ra)_i\big)_{i\in I}=p_j\big(r~h(a)_i\big)_{i\in I}=r~h(a)_j=r~p_j\big(h(a)_i\big)_{i\in I}=r~p_jh(a).\nn
\eea
Since $E_j$ is injective, there exists $\vphi_j:B\ra E_j$ such that $p_jh=\vphi_jf$. Let $\vphi:B\ra \prod_{i\in I}E_i,~b\mapsto \big(\vphi_i(b)\big)_{i\in I}$. Then
~$\vphi f(a)=\big(\vphi_i(f(a))\big)_{i\in I}=\big(\vphi_if(a)\big)_{i\in I}=\big(p_ih(a)\big)_{i\in I}=\big(h(a)_i\big)_{i\in I}=h(a)$.  \qedhere
\end{enumerate}
\end{proof}
\begin{note}[\magenta{Caution}]
The analysis in the proof of (b) above does not work for the direct sum $\oplus_{i\in I}E_i$ of injective modules.
\[\adjustbox{scale=0.9}{
\begin{tikzcd}
 0\ar[r] & A \arrow[dd, "h"'] \ar[rr,hook,"f"]  & & B \arrow[ddll,dashed,"\vphi"'] \arrow[ddddll,dotted,"\vphi_j"]\\
   & & &\\
   & \oplus_iE_i \ar[dd,shift left=1,"p_j"] & &\\
   & & &\\
   & E_j \ar[uu,hook,shift left=1,"q_j"] & & &
\end{tikzcd}}\]
This is because the definition of {\small $\vphi:B\ra\oplus_iE_i,~b\mapsto \big(\vphi_i(b)\big)_{i\in I}$} would require the \emph{additional constraint}
\bea
\label{proj-sum-eq1}\txt{for each}~~b\in B,~~\vphi_i(b)=0~~\txt{except for finitely many}~~i,
\eea
which is automatically satisfied if $b\in \im f$ (\blue{footnote}\footnote{This is due to the relations $p_jh=\vphi_jf$ for all $j$, which we can multiply by $q_j$ on the left and sum to get $h=id_{\oplus_iE_i}h=\sum_jq_jp_jh=\sum q_j\vphi_jf=\vphi f$, and so $\vphi|_{\im f}=\sum q_j\vphi_j$ makes sense on $\im f$ through $\im h$.}), but not necessarily for $b\not\in\im f$. If $f$ is not surjective, i.e., $\im f\neq B$, there might exist no $\{\vphi_j\}_{j\in I}$ that satisfy the \emph{extra requirement} (\ref{proj-sum-eq1}) along with $p_jh=\vphi_jf$ for all $j$.
\end{note}

\begin{question}
Is the following a correct proof of part (a) of Theorem \ref{InjModSumCrit} above? Assume $\mathop{\oplus}_{i\in I}M_i$ is injective. Let $j\in I$. If $_RM_j\subset{}_RM$, then
{\small\[
\textstyle\mathop{\oplus}_{i\in I}M_i\subset M\oplus(\mathop{\oplus}_{i\neq j}M_i),~~\Ra~~M\oplus(\mathop{\oplus}_{i\neq j}M_i)=(\mathop{\oplus}_{i\in I}M_i)\oplus N,~~\Ra~~M\cong M_j\oplus N,~~\mathop{\Ra}\limits^{\textcolor{red}{?}}~~M=M_j\oplus N'.
\]}Hence $M_j$ is injective.
\end{question}

\subsection{Weak criteria for projectivity and injectivity}~\\~
These alternative criteria will be used in the next section to prove fundamental theorems of homological algebra, and later to prove similar results on extendable/extension properties of certain special complexes (including projective and injective resolutions). Once again, the following results hold for any abelian category but for simplicity we will assume wlog that we are in $R$-mod.

Moreover, our proofs for module versions of the results do not strictly require elementwise definitions or arguments, but only require exactness of the abelian category which allows a morphism $f$ to factor as $f=m_fe_f$ with $m_f$ monic and $e_f$ epic. Hence, the proofs generalize directly to proofs for general abelian category versions of the results.

\begin{lmm}[\textcolor{blue}{\index{Projectivity criterion (weak)}{Weak criterion for projectivity}}]
If $P$ is an R-module, the following are equivalent:
\bit
\item[(a)] $P$ is projective.
\item[(b)] Given any exact diagram $B\sr{g}{\ral}C\sr{h}{\ral} D$ and any morphism $f:P\ra C$ satisfying $h f=0$, there exists a morphism $r:P\ra B$ such that $f=g r$.
\[\adjustbox{scale=0.9}{\bt
   && P\ar[dll,dashed,"r"']\ar[d,"f"] & \\
 B\ar[rr,near end,"g"] && C\ar[r,"h"]  & D
\et}~~~~\substack{hg=0~~~~~\\\\ \im g=\ker h}
\]
\eit
\end{lmm}
\begin{proof}
{\flushleft\ul{(a)$\Ra$(b)}}: Assume $P$ is projective. Consider an exact diagram $B\sr{g}{\ral}C\sr{h}{\ral} D$ and a morphism $f:P\ra C$ satisfying $h f=0$.  Then $\im (f)\subset\ker(h)=\im (g)$.
\[\adjustbox{scale=0.9}{\bt
   &&&& P\ar[ddl,dashed,near end,"\ol{f}"']\ar[dddllll,dotted,"r"']\ar[ddd,"f"] & \\
   &&&&   & \substack{hg=0~~~~~\\\\ \im g=\ker h}\\
   &&& {B\over \ker g}\ar[dr,dashed,hook,"\ol{g}"] &   & \\
 B\ar[urrr,dashed,near end,two heads,"\pi"]\ar[rrrr,"g"] &&&& C\ar[r,"h"]  & D
\et}~~~~~~~~~
\adjustbox{scale=0.9}{\bt
   && P\ar[dll,dashed,"r"']\ar[d,"\ol{f}"] &\\
 B\ar[rr,two heads,"\pi"] && {B\over\ker g}\ar[r]  & 0
\et}\]
Thus, along with $\pi(b):=b+\ker g$ and $\ol{g}(b+\ker g):=g(b)$, for $b\in B$, we have a well-defined map
\bea
\textstyle \ol{f}:P\ra {B\over\ker g},~m\mapsto b+\ker g,~~~~g(b)=f(m)~~\txt{(i.e., $b\in g^{-1}(f(m))$)},\nn
\eea
which by construction satisfies ~$f=\ol{g}\ol{f}$. Because $P$ is projective and $\pi$ is epic, there exists $r:P\ra B$ such that $\ol{f}=\pi r$. This implies
\bea
f=\ol{g}\ol{f}=\ol{g}\pi r=gr.\nn
\eea
{\flushleft \ul{(b)$\Ra$(a)}:} If (b) holds, then it holds for $D=0$ (hence $h=0$) in particular, which means $P$ is projective.
\end{proof}

\begin{lmm}[\textcolor{blue}{\index{Injectivity criterion (weak)}{Weak criterion for injectivity}}]
If $E$ is an R-module, the following are equivalent:
 \bit
 \item[(a)] $E$ is injective.
 \item[(b)] Given any exact diagram $U\sr{\al}{\ral}A\sr{f}{\ral}B$, and any morphism $h:A\ra E$ satisfying $h\al=0$, there exists a morphism $\vphi:B\ra E$ such that $h=\vphi f$.
\[\adjustbox{scale=0.9}{\bt
U\ar[r,"\al"] & A\ar[d,"h"]\ar[rr,"f"] && B\ar[dll,dashed,"\vphi"] & \substack{h\al=0~~~~~\\\\ \im\al=\ker f}\\
 & E && &
\et}\]
 \eit
\end{lmm}
\begin{proof}
{\flushleft\ul{(a)$\Ra$(b)}}: Assume $E$ is injective. Consider an exact diagram $U\sr{\al}{\ral}A\sr{f}{\ral}B$, and a morphism $h:A\ra E$ satisfying $h\al=0$. Then $\im \al=\ker f\subset\ker h$.
\[\adjustbox{scale=0.9}{\bt
U\ar[r,"\al"] & A\ar[ddd,"h"']\ar[dr,dashed,two heads,"\pi"']\ar[rrrr,"f"] &&&& B\ar[dddllll,dotted,"\vphi"] & \\
 &  & {A\over\ker f}\ar[ddl,dashed,"\ol{h}"']\ar[urrr,dashed,hook,"\ol{f}"] &&& & \\
\substack{h\al=0~~~~~\\\\ \im\al=\ker f} &  &&&& & \\
 & E &&&& &
\et}~~~~~~~~
\adjustbox{scale=0.9}{\bt
0\ar[r] & {A\over\ker f}\ar[d,"\ol{h}"]\ar[rr,hook,"\ol{f}"] && B\ar[dll,dashed,"\vphi"]\\
 & E &&
\et}\]
Thus, we have the usual maps $\ol{h}:{A\over\ker f}\ra E$ and $\ol{f}:{A\over\ker f}\ra B$ such that $\ol{f}$ is a monomorphism, and
\bea
\textstyle h=\ol{h}\pi,~~~~f=\ol{f}\pi~~~~(~\txt{where}~~\pi(a):=a+\ker f,~~\ol{f}(a+\ker g):=f(a),~~\ol{h}(a+\ker g):=h(a)~).\nn
\eea
Since $E$ is injective and $\ol{f}$ is monic, there exists $\vphi:B\ra E$ such that $\ol{h}=\vphi\ol{f}$, and so
\bea
\ol{h}\pi=\vphi\ol{f}\pi,~~\Ra~~h=\vphi f,~~~~\txt{since $\pi$ is epic}.\nn
\eea
{\flushleft\ul{(b)$\Ra$(a)}}: If (b) holds then it holds for $U=0$ (hence $\al=0$) in particular, which means $E$ is injective.
\end{proof}

\section{Morphism lifting property: Fundamental theorems of homological algebra}
The following two results (named ``Comparison theorem I'' and ``Comparison theorem II'') are also called the \index{Fundamental theorems of homological algebra}{\ul{``fundamental theorems of homological algebra''}}.
\begin{thm}[\textcolor{blue}{\index{Comparison theorem I}{Comparison theorem I}}]\label{CompThmI}
Let $\A$ be an Abelian category (such as R-mod), and let $P_\ast\sr{\vep_X}{\ral}X\ra 0$, $Q_\ast\sr{\vep_Y}{\ral}Y\ra 0$ be projective resolutions of $X,Y\in\Ob\A$. (\blue{footnote}\footnote{Note that $Q_\ast$ can be replaced by any exact sequence whose objects are not necessarily projective.}). Then for any morphism $f:X\ra Y$, we have the following.
\bit[leftmargin=0.7cm]
\item[a)] There exists a morphism of resolutions $f_\ast: P_\ast\ra Q_\ast$ extending $f$ in
the sense that the following diagram commutes.
\begin{figure}[H]
\centering
\adjustbox{scale=0.9}{\begin{tikzcd}
\cdots\ar[r,"d_2^P"] & P_1\ar[dd,dashed,"f_1"]\ar[r,"d_1^P"] & P_0\ar[dd,dashed,"f_0"]\ar[r,"\vep_X"] & X\ar[dd,"f"]\ar[r] & 0\\
  & & & & \\
\cdots\ar[r,"d_2^Q"] & Q_1\ar[r,"d_1^Q"] & Q_0\ar[r,"\vep_Y"] & Y \ar[r] & 0\\
\end{tikzcd}}
\caption{}\label{dg4s}
\end{figure}

\item[b)] Any two such extensions $f_\ast$ and $g_\ast$ are homotopic. (Hence, projective resolutions of the same object are \emph{homotopy equivalent.})
\eit
\end{thm}
\begin{proof}~
\bit[leftmargin=0.7cm]
\item[a)] Since $P_0$ is projective, the morphism $f_0$ exists such that $f\vep_X=\vep_Yf_0$. Similarly, because
    \bea
    (\vep_Yf_0)d_1^P=f\vep_Xd_1^P=f0=0,\nn
    \eea
    the morphism $f_1$ also exists such that $f_0d_1^P=d_1^Qf_1$. The rest follows by induction.

\item[b)] Let $f_\ast,f_\ast':P^\ast\ra Q^\ast$  be two morphisms of resolutions extending $f:X\ra Y$.

\begin{figure}[H]
\centering
\adjustbox{scale=0.9}{\begin{tikzcd}
\cdots\ar[r] & P_1\ar[dd,shift right=1,"f_1"']\ar[dd,shift left=1,"f'_1"]\ar[rr,"d_1^P"] && P_0\ar[dd,shift right=1,"f_0"']\ar[dd,shift left=1,"f'_0"]\ar[rr,"\vep_X"] && X\ar[dd,"f"]\ar[r] & 0\\
  & && && & \\
\cdots\ar[r] & Q_1\ar[rr,"d_1^Q"] && Q_0\ar[rr,"\vep_Y"] && Y \ar[r] & 0\\
\end{tikzcd}}
\caption{}
\end{figure}

The difference $\Delta f_n:=f'_n-f_n$ of the two morphisms gives the diagram with commuting squares:
\begin{figure}[H]
\centering
\adjustbox{scale=0.9}{\begin{tikzcd}
\cdots\ar[r] & P_2\ar[dd,"\Delta f_2"']\ar[rrr,"d_2^P"] &&& P_1\ar[ddlll,dashed,"k_1"']\ar[dd,"\Delta f_1"']\ar[rrr,"d_1^P"] &&& P_0\ar[ddlll,dashed,"k_0"']\ar[dd,"\Delta f_0"']\ar[rr,"\vep_X"] && X\ar[dd,"0=\Delta f"]\ar[r] & 0\\
 & &&& &&& && & \\
\cdots\ar[r] & Q_2\ar[rrr,"d_2^Q"] &&& Q_1\ar[rrr,"d_1^Q"] &&& Q_0\ar[rr,"\vep_Y"] && Y \ar[r] & 0\\
\end{tikzcd}}
\caption{}
\end{figure}
Since $\vep_Y\Delta f_0=0\vep_X=0$, and $P_0$ is projective, the morphism $k_0:P_0\ra Q_1$ exists such that $\Delta f_0=d_1^Qk_0$. Similarly, because
\bea
d_1^Q\Delta f_1=\Delta f_0d_1^P=d_1^Qk_0d_1^P~~\Ra~~d_1^Q(\Delta f_1-k_0d_1^P)=0,\nn
\eea
the morphism $k_1:P_1\ra Q_2$ exists such that $\Delta f_1-k_0d_1^P=d_2^Qk_1$, i.e., such that
\bea
\Delta f_1=k_0d_1^P+d_2^Qk_1.\nn
\eea
By induction, we get morphisms $k_n:P_n\ra Q_{n+1}$ such that $\Delta f_n=k_{n-1}d_n^P+d_{n+1}^Qk_n$. \qedhere
\eit
\end{proof}

\begin{rmk}[\blue{Full imbedding, Homotopy equivalence of resolutions}]~
\begin{enumerate}[leftmargin=0.7cm]
\item In the proof of part (b) of the theorem, only exactness (but not projectivity) of the sequence $Q_\ast\sr{\vep_Y}{\ral}Y\ra 0$ is required. Also, the result of the theorem can be stated as
    \bea
    Hom_\A(X,Y)\cong {Hom_{\A_0^\Integer}(P,Q)\over\sim},\nn
    \eea
where the equivalence relation $\sim$ is given by ``$f_\ast\sim g_\ast$ if $f_\ast\simeq g_\ast$''. Hence we have a full imbedding of $\A$ into the homotopy category $\H\A$,
\[
\A\hookrightarrow\H\A,~~X\sr{f}{\ral}Y~~\mapsto~~P_\ast^X\sr{[f_\ast]}{\ral}P_\ast^Y.
\]

\item With $X=Y$, the theorem says that any two projective resolutions (of the same object $X$),
\bea
P_{\ast 1}\sr{\vep_1}{\ral}X\ra 0,~~~~P_{\ast 2}\sr{\vep_2}{\ral}X\ra 0,\nn
\eea
are \ul{homotopy equivalent}: Indeed, any two extensions $f_\ast:P_{\ast 1}\ra P_{\ast 2}$ and $g_\ast:P_{\ast 2}\ra P_{\ast 1}$, of $id_X$, in turn give yet two more extensions of $id_X$, which are
\bea
f_\ast\circ g_\ast\simeq id_{P_{\ast 2}}:P_{\ast 2}\ra P_{\ast 2},~~~~g_\ast\circ f_\ast\simeq id_{P_{\ast 1}}:P_{\ast 1}\ra P_{\ast 1},\nn
\eea
making $f_\ast,g_\ast$ homotopy inverses of each other.

\item If we replace the projective resolutions by injective resolutions $0\ra X\sr{\vep_X}{\ra}E^\ast$, $0\ra Y\sr{\vep_Y}{\ra}F^\ast$, where
\bea
E^\ast:~~E^1\sr{d_E^1}{\ral}E^2\sr{d_E^2}{\ral}\cdots,~~~~F^\ast:~~F^1\sr{d_F^1}{\ral}F^2\sr{d_F^2}{\ral}\cdots,\nn
\eea
we get the following dual version of the theorem (and of the above results).
\end{enumerate}
\end{rmk}

\begin{thm}[\textcolor{blue}{\index{Comparison theorem II}{Comparison theorem II}}]\label{CompThmII}
Let $\A$ be an Abelian category (such as R-mod), and let $0\ra X\sr{\vep_X}{\ral}E^\ast$, $0\ra Y\sr{\vep_Y}{\ral}F^\ast$ be injective resolutions of $X,Y\in\Ob\A$. (\blue{footnote}\footnote{Note that $E^\ast$ can be replaced by any exact sequence whose objects are not necessarily injective.}). Then for any morphism $f:X\ra Y$, we have the following.
\bit[leftmargin=0.7cm]
\item[a)] There exists a morphism of resolutions $f^\ast: E^\ast\ra F^\ast$ extending $f$ in
the sense that the following diagram commutes.

\begin{figure}[H]
\centering
\adjustbox{scale=0.9}{\begin{tikzcd}
0\ar[r] & X\ar[dd,"f"]\ar[r,"\vep_X"] & E^0\ar[dd,dashed,"f^0"]\ar[r,"d_E^0"] & E^1\ar[dd,dashed,"f^1"]\ar[r,"d_E^1"] & \cdots\\
  & & & & \\
0\ar[r] & Y\ar[r,"{\vep_Y}"] & F^0\ar[r,"{d_F^0}"] & F^1\ar[r,"{d_F^1}"] & \cdots
\end{tikzcd}}
\caption{}
\end{figure}

\item[b)] Any two such extensions $f^\ast$ and $f^\ast{}'$ are homotopic. (Hence, injective resolutions of the same object are \emph{homotopy equivalent})
\eit
\end{thm}

\begin{proof}~
\bit[leftmargin=0.7cm]
\item[(a)] Since $F^0$ is injective, $f^0$ exists such that $\vep_Yf=f^0\vep_X$. Similarly, since
    \bea
    d_F^0f^0\vep_X=d_F^0\vep_Yf=0f=0,\nn
    \eea
    the morphism $f^1$ also exists such that $d_F^0f^0=f^1d_E^0$. The rest follows by induction.
\item[(b)] Let $f^\ast,f^\ast{}':E^\ast\ra F^\ast$  be two morphisms of resolutions extending $f:X\ra Y$.

\begin{figure}[H]
\centering
\adjustbox{scale=0.9}{\begin{tikzcd}
0\ar[r] & X\ar[dd,"f"]\ar[rr,"\vep_X"] && E^0\ar[dd,shift right=1,"f^0"']\ar[dd,shift left=1,"f^0{}'"]\ar[rr,"d_E^0"] && E^1\ar[dd,shift right=1,"f^1"']\ar[dd,shift left=1,"f^1{}'"]\ar[rr,"d_E^1"] && \cdots\\
  & && && && \\
0\ar[r] & Y\ar[rr,"{\vep_Y}"] && F^0\ar[rr,"{d_F^0}"] && F^1\ar[rr,"{d_F^1}"] && \cdots
\end{tikzcd}}
\caption{}
\end{figure}

The difference $\Delta f^n=f^n{}'-f^n$ of the two morphisms gives the diagram with commuting squares:

\begin{figure}[H]
\centering
\adjustbox{scale=0.9}{\begin{tikzcd}
0\ar[r] & X\ar[dd,"\Delta f=0"']\ar[rr,"\vep_X"] && E^0\ar[dd,"\Delta f^0"']\ar[rrr,"d_E^0"] &&& E^1\ar[ddlll,dashed,"k^1"']\ar[dd,"\Delta f^1"']\ar[rrr,"d_E^1"] &&& E^2\ar[ddlll,dashed,"k^2"']\ar[dd,"\Delta f^2"']\ar[r] & \cdots\\
  & && &&& &&& & \\
0\ar[r] & Y\ar[rr,"{\vep_Y}"] && F^0\ar[rrr,"{d_F^0}"] &&& F^1\ar[rrr,"{d_F^1}"] &&& F^2\ar[r] & \cdots
\end{tikzcd}}
\caption{}
\end{figure}

Since $\Delta f^0\vep_X=\vep_Y0=0$, injectivity of $F^0$ implies the morphism $k^1:E^1\ra F^0$ exists such that $\Delta f^0=k^1d_E^0$. Similarly, because
\bea
\Delta f^1d_E^0=d_F^0\Delta R(f)=d_F^0k^1d_E^0~~\Ra~~(\Delta f^1-d_F^0k^1)d_E^0=0,\nn
\eea
the morphism $k^2:E^2\ra F^1$ exists such that $\Delta f^1-d_F^0k^1=k^2d_E^1$, i.e., such that
\bea
\Delta f^1=d_F^0k^1+k^2d_E^1.\nn
\eea
By induction, we get morphisms $k^n:E^n\ra F^{n-1}$ such that $\Delta f^n=d_F^{n-1}k^n+k^{n+1}d_E^n$. \qedhere
\eit
\end{proof}

\begin{rmk}[\blue{Full imbedding, Homotopy equivalence of resolutions}]~
\begin{enumerate}[leftmargin=0.7cm]
\item In the proof of part (b) of the theorem, only exactness (but not injectivity) of the sequence $0\ra X\sr{\vep_X}{\ral}E^\ast$ is required. Also, the result of the theorem can be stated as
    \bea
    Hom_\A(X,Y)\cong {Hom_{\A_0^\Integer}(E,F)\over\sim},\nn
    \eea
where the equivalence relation $\sim$ is given by ``$f^\ast\sim g^\ast$ if $f^\ast\simeq g^\ast$''. Hence we have a full imbedding of $\A$ into the homotopy category $\H\A$,
\[
\A\hookrightarrow\H\A,~~X\sr{f}{\ral}Y~~\mapsto~~I^\ast_X\sr{[f^\ast]}{\ral}I^\ast_Y.
\]

\item With $X=Y$, the theorem says that any two injective resolutions (of the same object $X$),
\bea
0\ra X\sr{\vep_1}{\ral}E_1^\ast,~~~~0\ra X\sr{\vep_2}{\ral}E_2^\ast,\nn
\eea
are \ul{homotopy equivalent}: Indeed, any two extensions, $f^\ast:E^\ast_1\ra E^\ast_2$ and $g^\ast:E^\ast_2\ra E^\ast_1$, of $id_X$, in turn give yet two more extensions of $id_X$, which are
\bea
f^\ast\circ g^\ast\simeq id_{E^\ast_2}:E^\ast_2\ra E^\ast_2,~~~~g^\ast\circ f^\ast\simeq id_{E^\ast_1}:E^\ast_1\ra E^\ast_1,\nn
\eea
making $f^\ast,g^\ast$ homotopy inverses of each other.
\end{enumerate}
\end{rmk}

\section{Short exact sequence (SES) lifting property: Horseshoe lemmas}
As usual, it is apparent that these results hold for any abelian category with enough projectives (resp. enough injectives), but for simplicity, we will assume wlog that we are in $R$-mod.
\begin{lmm}[\textcolor{blue}{\index{Horseshoe lemma I}{Horseshoe lemma I}}]
Consider an exact sequence of $R$-modules ~$0\ra A\sr{f}{\ral}B\sr{g}{\ral}C\ra 0$. Given any projective resolutions $A_\ast:P_A\sr{\vep_A}{\ral}A$ and $C_\ast:P_C\sr{\vep_C}{\ral}C$ of $A$ and $C$, as in the diagram:
\[
\adjustbox{scale=0.8}{\begin{tikzcd}
         & P_A\ar[d,"\vep_A"] &    & P_C\ar[d,"\vep_C"] & \\
 0\ar[r] & A\ar[r,"f"]\ar[d]  & B\ar[r,"g"] & C\ar[d]\ar[r]  & 0 \\
         &         0          &             &        0       &
\end{tikzcd}},\]
there exists a projective resolution $B_\ast:P_B\sr{\vep_B}{\ral}B$ and morphisms of complexes $f_\ast:P_A\ra P_B$, $g_\ast:P_B\ra P_C$ such that we have an exact sequence of projective resolutions  $0\ra A_\ast\sr{f_\ast}{\ral}B_\ast\sr{g_\ast}{\ral}C_\ast\ra 0$, i.e., such that we have an exact commutative diagram of the  form
\[\adjustbox{scale=0.8}{
\begin{tikzcd}
 0\ar[r]      & P_A\ar[d,"\vep_A"]\ar[r,dashed,"f_\ast"] &  P_B\ar[d,dashed,"\vep_B"]\ar[r,dashed,"g_\ast"]  & P_C\ar[d,"\vep_C"]\ar[r] & 0\\
 0\ar[r] & A\ar[r,"f"]\ar[d]  & B\ar[r,"g"]\ar[d] & C\ar[d]\ar[r]  & 0 \\
         &         0          &     0       &        0       &
\end{tikzcd}}
~~~~=~~~~
\adjustbox{scale=0.8}{\begin{tikzcd}
         & \vdots\ar[d] & \vdots\ar[d] & \vdots\ar[d] &  \\
 0\ar[r] & P_1^A\ar[d,"d_1^A"]\ar[r,dashed,"f_1"] & P_1^B\ar[d,dashed,"d_1^B"]\ar[r,dashed,"g_1"] & P_1^C\ar[d,"d_1^C"]\ar[r] &  0 \\
 0\ar[r] & P_0^A\ar[d,"\vep_A"]\ar[r,dashed,"f_0"] & P_0^B\ar[d,dashed,"\vep_B"]\ar[r,dashed,"g_0"] & P_0^C\ar[d,"\vep_C"]\ar[r] &  0 \\
 0\ar[r] & A\ar[d]\ar[r,"f"] & B\ar[d]\ar[r,"g"] & C\ar[d]\ar[r] & 0  \\
         & 0 & 0 & 0 &   \\
\end{tikzcd}}\]
and moreover, the exact sequence of deleted resolutions {\small $0\ra P_A\sr{f_\ast}{\ral}P_B\sr{g_\ast}{\ral}P_C\ra 0$} has \ul{split-exact rows}.
\end{lmm}

\begin{proof}
It suffices to construct an exact diagram of the following form with the desired properties:
\[\adjustbox{scale=0.8}{\begin{tikzcd}
         & \vdots\ar[d] && \vdots\ar[d] && \vdots\ar[d] &  \\
 0\ar[r] & P_2^A\ar[ddrr,dashed,"f_1d_2^A"]\ar[dd,"d_2^A"]\ar[rr,dashed,"f_2"] && P_2^B\ar[dd,dashed,"d_2^B"]\ar[rr,dashed,"g_2"] && P_2^C\ar[ddll,dashed,"h_2"'pos=0.45]\ar[dd,"d_2^C"]\ar[r] &  0 \\
         &   &&  &&  &  \\
 0\ar[r] & P_1^A\ar[ddrr,dashed,"f_0d_1^A"]\ar[dd,"d_1^A"]\ar[rr,dashed,"f_1"] && P_1^B\ar[dd,dashed,"d_1^B"]\ar[rr,dashed,"g_1"] && P_1^C\ar[ddll,dashed,"h_1"'pos=0.45]\ar[dd,"d_1^C"]\ar[r] &  0 \\
         &   &&  &&  &  \\
 0\ar[r] & P_0^A\ar[ddrr,dashed,"f\vep_A"]\ar[dd,"\vep_A"]\ar[rr,dashed,"f_0"] && P_0^B\ar[dd,dashed,"\vep_B"]\ar[rr,dashed,"g_0"] && P_0^C\ar[ddll,dashed,"h_0"'pos=0.45]\ar[dd,"\vep_C"]\ar[r] &  0 \\
         &   &&  &&  &  \\
 0\ar[r] & A\ar[d]\ar[rr,"f"] && B\ar[d]\ar[rr,"g"] && C\ar[d]\ar[r] & 0  \\
         & 0 && 0 && 0 &   \\
\end{tikzcd}}\]

For the first step, we define ~{\small$P_0^B:=P_0^A\oplus P_0^C$,~ $f_0:P_0^A\ra P_0^B,~a\mapsto (a,0)$, ~$g_0:P_0^B\ra P_0^C,~(a,c)\mapsto c$}.\\
Because $P_0^C$ is projective and $g$ is epic, the morphism $h_0$ exists such that $\vep_C=gh_0$. We then define
\bea
\vep_B:P_0^B\ra B,~(a,c)\mapsto f\vep_A(a)+h_0(c),\nn
\eea
which gives the following commutative diagram:

\begin{figure}[H]
\centering
\adjustbox{scale=0.8}{%
\begin{tikzcd}
       &  & 0\ar[d] & 0\ar[d] & 0\ar[d] &  & \\
 &  0\ar[r] & \ker\vep_A\ar[d,hook]\ar[r] & \ker\vep_B\ar[d,hook]\ar[r] & \ker\vep_C\ar[d,hook]\ar[r] &  0\ar[r,dashed,shift left=3,"\delta_0"] & ~\\
 &  0\ar[r] & P_0^A\ar[d,two heads,"\vep_A"]\ar[r,"f_0"] & P_0^B\ar[d,dashed,"\vep_B"]\ar[r,"g_0"] & P_0^C\ar[d,two heads,"\vep_C"]\ar[r] &  0 & \\
 &  0\ar[r] & A\ar[d]\ar[r,"f"] & B\ar[d,dashed]\ar[r,"g"] & C\ar[d]\ar[r] & 0 & \\
\ar[r,dashed,shift right,"\delta_0"]  & ~ & 0 & 0 & 0 &  & \\
\end{tikzcd}}
\caption{}\label{hw3rvdg}
\end{figure}
In Fig \ref{hw3rvdg} the (nonzero) bottom row is exact by hypotheses, the middle row is split exact by construction, and the top row is exact by the snake lemma. Moreover, the snake lemma also implies the sequence of cokernels (of the $\vep$'s) must all be zeros since $\vep_A$ and $\vep_C$ are both epic. Hence the middle column is exact as well.

Next, we repeat the same steps above as follows: Let
\bea
P_1^B:=P_1^A\oplus P_1^C,~~~~f_1:P_1^A\ra P_1^B,~a\mapsto (a,0),~~~~g_1:P_1^B\ra P_1^C,~(a,c)\mapsto c.\nn
\eea
With the morphism $h_1:P_1^C\ra P_0^B$, which exists such that $d_1^C=g_0h_1$, we define
\bea
d_1^B:P_1^B\ra P_0^B,~(a,c)\mapsto f_0d_1^A(a)+h_1(c),\nn
\eea
and so as before, we get the following commutative diagram (with $d^X_n=m_{d^X_n}e_{d^X_n}$):

\begin{figure}[H]
\centering
\adjustbox{scale=0.8}{%
\begin{tikzcd}
       &  & 0\ar[d] & 0\ar[d] & 0\ar[d] &  & \\
 &  0\ar[r] & \ker d_1^A\ar[d,hook]\ar[r] & \ker d_1^B\ar[d,hook]\ar[r] & \ker d_1^C\ar[d,hook]\ar[r] &  0\ar[r,dashed,shift left=3,"\delta_1"] & ~\\
 &  0\ar[r] & P_1^A\ar[d,two heads,"e_{d_1^A}"]\ar[r,"f_1"] & P_1^B\ar[d,dashed,"u"]\ar[r,"g_1"] & P_1^C\ar[d,two heads,"e_{d_1^C}"]\ar[r] &  0 & \\
 &  0\ar[r] & \ker\vep_A\ar[d]\ar[r,"f_0"] & \ker\vep_B\ar[d,dashed]\ar[r,"g_0"] & \ker\vep_C\ar[d]\ar[r] & 0 & \\
\ar[r,dashed,shift right=2,"\delta_1"] & ~ & 0 & 0 & 0 &  & \\
\end{tikzcd}}
\caption{}\label{hw3rvdg2}
\end{figure}
where the morphism $u:P_1^B\ra \ker\vep_B$ is induced because commutativity of the two-step diagram (constructed so far) implies $\vep_Bd_1^Bf_1=f\vep_Ad_1^A=0$, and so $\vep_Bd_1^B=0$.

In Fig \ref{hw3rvdg2} the (nonzero) bottom row is exact by hypotheses, the middle row is split exact by construction, and the top row is exact by the snake lemma. Moreover, the snake lemma also implies the sequence of cokernels (of the $d_1$'s) must all be zeros since $e_{d_1^A}$ and $e_{d_1^C}$ are both epic. Hence the middle column is exact as well.

The rest of the steps follow by induction, and we then glue together the diagrams from all steps (i.e., Figure \ref{hw3rvdg}, Figure \ref{hw3rvdg2}, etc) to get the overall diagram.
\end{proof}

\begin{lmm}[\textcolor{blue}{\index{Horseshoe lemma II}{Horseshoe lemma II}}]
Consider an exact sequence of $R$-modules ~$0\ra A\sr{f}{\ral}B\sr{g}{\ral}C\ra 0$. Given any injective resolutions $A^\ast:A\sr{\vep_A}{\ral}E_A$ and $C^\ast:C\sr{\vep_C}{\ral}E_C$ of $A$ and $C$, as in the diagram
\[\adjustbox{scale=0.9}{\begin{tikzcd}
         & 0\ar[d] &    & 0\ar[d] & \\
 0\ar[r] & A\ar[d,"\vep_A"]\ar[r,"f"]  & B\ar[r,"g"] & C\ar[d,"\vep_C"]\ar[r]  & 0 \\
         &  E_A  &       &  E_C  &
\end{tikzcd}}\]
there exists an injective resolution $B^\ast:E_B\sr{\vep_B}{\ral}B$ and morphisms of complexes $f^\ast:E_A\ra E_B$, $g^\ast:E_B\ra E_C$ such that we have a short exact sequence of injective resolutions $0\ra A^\ast\sr{f^\ast}{\ral}B^\ast\sr{g^\ast}{\ral}C^\ast\ra 0$, i.e., such that we have an exact commutative diagram of the form
\[\adjustbox{scale=0.8}{
\begin{tikzcd}
         &         0 \ar[d]         &     0\ar[d]      &        0\ar[d]       & \\
 0\ar[r] & A\ar[d,"\vep_A"]\ar[r,"f"]  & B\ar[d,dashed,"\vep_B"]\ar[r,"g"] & C\ar[d,"\vep_C"]\ar[r]  & 0 \\
 0\ar[r] & E_A\ar[r,dashed,"f^\ast"] &  E_B\ar[r,dashed,"g^\ast"]  & E_C\ar[r] & 0
\end{tikzcd}}
~~~~=~~~~
\adjustbox{scale=0.9}{\begin{tikzcd}
         & 0\ar[d] & 0\ar[d] & 0\ar[d] &   \\
 0\ar[r] & A\ar[d,"\vep_A"]\ar[r,"f"] & B\ar[d,dashed,"\vep_B"]\ar[r,"g"] & C\ar[d,"\vep_C"]\ar[r] & 0  \\
 0\ar[r] & E_A^0\ar[d,"d_A^0"]\ar[r,dashed,"f^0"] & E_B^0\ar[d,dashed,"d_B^0"]\ar[r,dashed,"g^0"] & E_C^0\ar[d,"d_C^0"]\ar[r] &  0 \\
 0\ar[r] & E_A^1\ar[d]\ar[r,dashed,"f^1"] & E_B^1\ar[d]\ar[r,dashed,"g^1"] & E_C^1\ar[d]\ar[r] &  0 \\
         & \vdots & \vdots & \vdots &
\end{tikzcd}}\]
and moreover, the sequence of deleted injective resolutions {\small $0\ra E_A\sr{f^\ast}{\ral}E_B\sr{g^\ast}{\ral}E_C\ra 0$} has \ul{split exact rows}.
\end{lmm}

\begin{proof}
It suffices to construct an exact diagram of the following form with the desired properties:
\[\adjustbox{scale=0.8}{
\begin{tikzcd}
         & 0\ar[d] && 0\ar[d] && 0\ar[d] &   \\
 0\ar[r] & A\ar[dd,"\vep_A"]\ar[rr,"f"] && B\ar[ddll,dashed,"h"'pos=0.45]\ar[dd,dashed,"\vep_B"]\ar[ddrr,dashed,"\vep_Cg"]\ar[rr,"g"] && C\ar[dd,"\vep_C"]\ar[r] & 0  \\
         &   &&  &&  &  \\
 0\ar[r] & E_A^0\ar[dd,"d_A^0"]\ar[rr,dashed,"f^0"] && E_B^0\ar[ddll,dashed,"h^0"'pos=0.45]\ar[ddrr,dashed,"d_C^0g^0"]\ar[dd,dashed,"d_B^0"]\ar[rr,dashed,"g^0"] && E_C^0\ar[dd,"d_C^0"]\ar[r] &  0 \\
          &   &&  &&  &  \\
 0\ar[r] & E_A^1\ar[dd,"d_A^1"]\ar[rr,dashed,"f^1"] && E_B^1\ar[ddll,dashed,"h^1"'pos=0.45]\ar[ddrr,dashed,"d_C^1g^1"]\ar[dd,dashed,"d_B^1"]\ar[rr,dashed,"g^1"] && E_C^1\ar[dd,"d_C^1"]\ar[r] &  0 \\
         &   &&  &&  &  \\
 0\ar[r] & E_A^2\ar[d]\ar[rr,dashed,"f^2"] && E_B^2\ar[d]\ar[rr,dashed,"g^2"] && E_C^2\ar[d]\ar[r] &  0 \\
         & \vdots && \vdots && \vdots &
\end{tikzcd}}\]

For the first step, we define ~{\small$E_B^0:=E_A^0\oplus E_C^0$,~ $f^0:E_A^0\ra P_B^0,~a\mapsto (a,0)$, ~$g^0:P_B^0\ra P_C^0,~(a,c)\mapsto c$}.\\
Because $E_A^0$ is injective and $f$ is monic, the morphism $h$ exists such that $\vep_A=hf$. We then define
\bea
\vep_B:B\ra E_B^0,~b\mapsto \Big(h(b),\vep_Cg(b)\Big),\nn
\eea
which gives the following commutative diagram:

\begin{figure}[H]
\centering
\adjustbox{scale=0.8}{%
\begin{tikzcd}
 &  & 0\ar[d] & 0\ar[d,dashed] & 0\ar[d] & \ar[r,dashed,shift left=3,"\delta_0"] & ~\\
 &  0\ar[r] & A\ar[d,hook,"\vep_A"]\ar[r,"f"] & B\ar[d,dashed,"\vep_B"]\ar[r,"g"] & C\ar[d,hook,"\vep_C"]\ar[r] & 0 &  \\
 &  0\ar[r] & E_A^0\ar[d,two heads]\ar[r,"f^0"] & E_B^0\ar[d,two heads]\ar[r,"g^0"] & E_C^0\ar[d,two heads,]\ar[r] &  0 & \\
 \ar[r,dashed,shift right=3,"\delta_0"] &  0\ar[r] & \coker \vep_A\ar[d]\ar[r] & \coker \vep_B\ar[d]\ar[r] & \coker \vep_C\ar[d]\ar[r] &  0 & \\
       &  & 0 & 0 & 0 &  & \\
\end{tikzcd}}
\caption{}\label{hw3rvdgg}
\end{figure}
In Figure \ref{hw3rvdgg} the (nonzero) top row is exact by hypotheses, the middle row is split exact by construction, and the bottom row is exact by the snake lemma. Moreover, the snake lemma also implies the sequence of kernels (of the $\vep$'s) must all be zeros since $\vep_A$ and $\vep_C$ are both monic. Hence the middle column is exact as well.

Next, we repeat the same steps above as follows: Let
\bea
E_B^1:=E_A^1\oplus E_C^1,~~~~f^1:E_A^1\ra E_B^1,~a\mapsto (a,0),~~~~g^1:E_B^1\ra E_C^1,~(a,c)\mapsto c.\nn
\eea
With the map $h^0:E_B^0\ra E_A^1$, which exists such that $d_A^0=h^0f^0$, we define
\bea
d_B^0:E_B^0\ra E_B^1,~b\mapsto \Big(h^0(b),d_C^0g^0(b)\Big),\nn
\eea
and so as before, we get the following commutative diagram (with $d_X^n=m_{d_X^n}e_{d_X^n}$):

\begin{figure}[H]
\centering
\adjustbox{scale=0.8}{%
\begin{tikzcd}
 &  & 0\ar[d] & 0\ar[d,dashed] & 0\ar[d] & \ar[r,dashed,shift left=3,"\delta_1"] & ~\\
 &  0\ar[r] & \coker \vep_A\ar[d,hook,"m_{d_A^0}"]\ar[r] & \coker \vep_B\ar[d,dashed,"v"]\ar[r] & \coker \vep_C\ar[d,hook,"m_{d_C^0}"]\ar[r] & 0 &  \\
 &  0\ar[r] & E_A^1\ar[d,two heads]\ar[r,"f^1"] & E_B^1\ar[d,two heads]\ar[r,"g^1"] & E_C^1\ar[d,two heads]\ar[r] &  0 & \\
 \ar[r,dashed,shift right=3,"\delta_1"] &  0\ar[r] & \coker d_A^0\ar[d]\ar[r] & \coker d_B^0\ar[d]\ar[r] & \coker d_C^0\ar[d]\ar[r] &  0 & \\
        &  & 0 & 0 & 0 &  & \\
\end{tikzcd}}
\caption{}\label{hw3rvdgg2}
\end{figure}
where the morphism $v:\coker\vep_B\ra E_B^1$ is induced because commutativity of the two-step diagram (constructed so far) implies $d_B^0\vep_Bf=f_1d_A^0\vep_B=0$, and so $d_B^0\vep_B=0$.

In Figure \ref{hw3rvdgg2} the (nonzero) top row is exact by hypotheses, the middle row is split exact by construction, and the bottom row is exact by the snake lemma. Moreover, the snake lemma also implies the sequence of kernels (of the $d^0$'s) must all be zeros since $m_{d_A^0}$ and $m_{d_C^0}$ are both monic. Hence the middle column is exact as well.

The rest of the steps follow by induction, and we then glue together the diagrams from all steps (i.e., Figure \ref{hw3rvdgg}, Figure \ref{hw3rvdgg2}, etc) to get the overall diagram.
\end{proof}

%% file: parts/AlgebraCat/AlgebraCatS10.tex
\chapter{Classical Derived Functors: Measurement of (partial) exactness}\label{AlgebraCatS10}
\section{The idea of a derived functor: Motivation}
Let $\A,\B$ be abelian categories and $F:\A\ra \B$ an additive functor. Consider an exact sequence in $\A$,
\[
\txt{SES}:~0\ra A\ra A'\ra A''\ra 0.
\]
If $F$ is a \ul{left exact additive} functor (resp. cofunctor) such as a \ul{mor functor}, then its image is the exact sequence
\bea
F(\txt{SES}):~0\ra F(A)\ra F(A')\ra F(A''),~~~~\Big(~\txt{resp.}~~0\ra F(A'')\ra F(A')\ra F(A)~\Big).\nn
\eea
For the same reasons (\blue{footnote}\footnote{Object resolution (or decomposition) in general is done in order to reveal otherwise hidden parameters that can be used for better comparison/contrast of (related) objects/systems. It is instructive to compare this with the reasons behind the category localization process/operation, which is itself simply a type of resolution process (although relatively indirect as compared to the projective/injective resolution operations for a module).}) that we (i) carry out localization and (ii) resolve/extend an object into a projective or injective resolution, we may wish to extend/continue this sequence (and its exactness) to the right until we either reach a stopping point $0$ or continue indefinitely (the length of which would be a measure of the amount of exactness that was destroyed by applying $F$). One realization of such an exact extension of the sequence $F$(SES) to the right is in terms of a sequence of \ul{functorial processes} (obtained using the LES of homology and) called \ul{right derived} functors of $F$.

Similarly, if $G:\A\ra \B$ is a \ul{right exact additive} functor (resp. cofunctor) such as a \ul{tensor functor}, then we get an exact sequence
\bea
G(\txt{SES}):~ G(A)\ra G(A')\ra G(A'')\ra 0,~~~~\Big(~\txt{resp.}~~G(A'')\ra G(A')\ra G(A)\ra 0~\Big)\nn
\eea
an exact continuation of which to the left can be realized in terms of a sequence of \ul{functorial processes} (from the LES of homology) called \ul{left derived} functors of $G$.

As quantities that involve estimates/measures of the amount of \ul{exactness} of the SES that is destroyed by applying $F$, it is but natural for derived functors of $F$ to be defined in terms of \ul{homology} functors (through the LES of homology). To obtain the desired exact continuation of $F$(SES), we proceed in \ul{three steps} as follows (while noting that the process is similar for $G$(SES)):
\begin{enumerate}
\item First, we \ul{vertically resolve} the \ul{horizontally placed} SES (say by means of a horseshoe lemma) into a \ul{split-exact} SES of chain complexes (in the form of \ul{deleted resolutions}) ~$R$(SES) $=$ SE-SES-CC.
\item Next, we apply $F$ on the resolved sequence $R$(SES) $=$ SE-SES-CC to obtain another SES of chain complexes ~$F\big(R(\txt{SES})\big)$ $=$ $F$(SE-SES-CC). (\blue{footnote}\footnote{This SES is still split-exact due to the additivity of $F$.})
\item Finally, we take the induced LES of homologies of the sequence $F\big(R(\txt{SES})\big) = F(\txt{SE-SES-CC})$ to be the desired exact extension of $F(\txt{SES})$, i.e.,
\bea
\label{ExactContEq1}\txt{Exact continuation of}~~F(\txt{SES})~~:=~~\txt{LES}\Big(F\big(R(\txt{SES})\big)\Big)=\txt{LES}\big(F(\txt{SE-SES-CC})\big).
\eea
\end{enumerate}

{\flushleft\ul{Why is $\txt{LES}\big(F\big(R(\txt{SES})\big)\big)$ an exact continuation of $F(\txt{SES})$?}} An explanation goes as follows: A projective (resp. injective) resolution of an object $A\in\Ob\A$, i.e., ~$P^A\sr{\vep}{\ral}A\ra 0$~ or explicitly, $\cdots\sr{d_2}{\ral}P_1\sr{d_1}{\ral}P_0\sr{\vep}{\ral}A\ra 0\ra\cdots$ (where $P^A:\cdots\sr{d_2}{\ral}P_1\sr{d_1}{\ral}P_0\ra 0\ra\cdots$) satisfies
\bea
&&\textstyle H(P^A):\cdots\sr{0}{\ral} \overbrace{H_1(P^A)}^{=0}\ra H_0(P^A)\ra 0,~~\txt{with}~~H_0(P^A)={\ker 0\over\im d_1}={P_0\over\ker\vep}\cong A,\nn\\
&&~~\Ra~~H(P^A)~\cong~ A_\ast:\cdots \ra 0\ra A\ra 0\ra\cdots~~\txt{(i.e., there is a quasi-isomorphism $P^A\ra A_\ast$)}.\nn
\eea
Consequently, we can think of the \ul{homology operation} as a ``\ul{de-resolution operation}'' (i.e., $H$ is ``left-inverse'' to the projective resolution process $P$). Similarly, because (\blue{footnote}\footnote{
The isomorphism $\txt{LES}(R(\txt{SES}))\cong\txt{SES}$ (ignoring $0$'s) implies that the exact sequence of deleted resolutions $R(\txt{SES})$ is a representation (in $\A_0^{\Integer\times\Integer}$) of the $\txt{SES}\in\Ob\A_0^\Integer$, in much the same way that the isomorphism $H(P^A)~\cong~ A$ (ignoring $0$'s) implies that a deleted projective resolution $P^A$ (or a deleted injective resolution $I_A$) is a representation (in $\A_0^\Integer$) of $A\in\Ob\A$.})
\bea
\txt{LES}(R(\txt{SES}))\cong\txt{SES} ~~~~(\txt{ignoring $0$'s}),\nn
\eea
we can also think of the LES operation (which is based on the homology operation) as ``left-inverse'' to our resolution operation $R$, i.e., we may write
\bea
R_l^{-1}:=\txt{LES}.\nn
\eea
We can now restate the original definition in (\ref{ExactContEq1}) as
\bea
\label{ExactContEq2}\txt{Exact continuation of}~~F(\txt{SES})~~:=~~R_l^{-1}(F(R(\txt{SES})))=(R_l^{-1}\circ F\circ R)(\txt{SES}).
\eea
This shows that in $\B_0^\Integer$, the exact sequence $(R_l^{-1}\circ F\circ R)(\txt{SES})$ is a representation of $F(\txt{SES})$ that indeed contains a reasonable measure of the amount of exactness of the SES that is destroyed by applying $F$.

\section{Basic definitions and examples: Torsion and Extension functors}
\begin{dfn}[\textcolor{blue}{
\index{Derived! functor}{Derived functor},
\index{Classical derived functors! Left derived functors}{Left derived functors},
\index{Classical derived functors! Right derived functors}{Right derived functors},
\index{Derived! functors of a functor}{Derived functors of a functor},
\index{Derived! functors of a cofunctor}{Derived functors of a cofunctor}}]
Let $\A,\B$ be abelian categories, $\A\sr{F}{\ral}\B$ a functor (which of course induces a functor $\A_0^\Integer\sr{F}{\ral}\B_0^\Integer$), and
{\small\bea
R:\A\ra \A_0^\Integer,~~A\sr{f}{\ral}A'~\mapsto~R(A)\sr{R(f)}{\ral}R(A')\nn
\eea}a functorial projective/injective deleted-resolution operation on $\A$ (i.e., $R_A:=R(A)$ is a deleted projective or injective resolution of $A$). The \ul{$n$th derived functor} of $F$ (associated with the deleted-resolution operation $R:\A\ra\A_0^\Integer$ and a homology functor $H:\B_0^\Integer\ra\B_0^\Integer$) is the composition
{\small\bea
&&D^nF:=H_n\circ F\circ R:\A\sr{R}{\ral}\A_0^\Integer\sr{F}{\ral}\B_0^\Integer\sr{H_n}{\ral}\B,\nn\\
&& A\sr{f}{\ral}A'~~\mapsto~\bt(D^nF)(A)\ar[rr,"{(D^nF)(f)}"] && (D^nF)(A')\et~:=~\bt H_n(F(R(A)))\ar[rr,"{H_n(F(R(f)))}"] && H_n(F(R(A')))\et.\nn
\eea}

The functors $D^nF:\A\ra\B$ are \ul{left derived} (resp. \ul{right derived}) functors of $F$ if for each object $A\in\Ob\A$, the sequence $D^\ast F(A)=H(F(R(A))):\cdots,~D^{n-1}F(A),~ D^nF(A),~D^{n+1}F(A),~\cdots$ by construction \ul{bounded on the right} (resp. \ul{bounded on the left}) in the sense that the sequence \ul{eventually becomes zeros on the right}, just as the image of a right exact functor does (resp. \ul{eventually becomes zeros on the left}, just as the image of a left exact functor does).

In the following table of various derived functor possibilities/choices (where as usual, we consider \ul{left projective} resolutions and \ul{right injective} resolutions only), we by convention choose the resolution type such that \ul{left-exact} functors/cofunctors always give \ul{right-derived} functors, while \ul{right-exact} functors/cofunctors always give \ul{left-derived} functors:
\begin{table}[H]
  \centering
\scalebox{0.9}{
\begin{tabular}{ll|ll}
   & {\bf Left-derived} $D^nF$ & & {\bf Right-derived} $D^nF$\\
  \hline &&&\\
  1 & {\bf Right-exact} Functor  $F$ ~{\scriptsize plus}~ Projective res. $R$  & ~~~~~1 &  {\bf Left-exact} Functor $F$ ~{\scriptsize plus}~ Injective res. $R$ \\&&&\\
  2 & {\bf Right-exact} Cofunctor $F$ ~{\scriptsize plus}~ Injective res. $R$ & ~~~~~2 & {\bf Left-exact} Cofunctor $F$ ~{\scriptsize plus}~ Projective res $R$.   \\&&&\\
  \hline
\end{tabular}}
\end{table}

{\flushleft More} explicitly, we have \ul{two main cases} as follows: Consider an additive functor of abelian categories
\bea
F:\A\ra\B.~~~~\txt{(This induces a functor ~$F:\A_0^\Integer\ra \B_0^\Integer$.)}\nn
\eea
Let ~$R:\A\ra \A_0^\Integer,~\big(A\sr{f}{\ral}A'\big)\longmapsto\big( R_A\sr{R_f}{\ral}R_{A'}\big)$~ be a functorial deleted-resolution operation, where
{\small\bea
R_A:~\cdots\ra R_A^{-1}\sr{d_A^{-1}}{\ral}R_A^0\sr{d_A^0}{\ral} R_A^1\sr{d_A^1}{\ral}\cdots~~~~\txt{and}~~~~R_{A'}:~\cdots\ra R_{A'}^{-1}\sr{d_{A'}^{-1}}{\ral}R_{A'}^0\sr{d_{A'}^0}{\ral} R_{A'}^1\sr{d_{A'}^1}{\ral}\cdots\nn
\eea}
are deleted resolutions of $A,A'$ in $\A_0^\Integer$, and $R_f:=R(f)$ is a morphism of complexes that extends $f$. (\blue{footnote}\footnote{Recall from the comparison theorems, i.e., the fundamental theorems of homological algebra, that if $\A$ is an abelian category with enough projectives and injective, then $R_f$ is unique up to a homotopy.})

{\flushleft 1. \ul{Functor $F$}}: The $n$th derived functor ~{\footnotesize $D^nF=H_n\circ F\circ R:\A\sr{R}{\ral}\A_0^\Integer\sr{F}{\ral}\B_0^\Integer\sr{H_n}{\ral}\B$}~ of $F$ is given by
{\footnotesize
\bea
&&\bt[column sep=small,row sep=tiny]
\A\ar[r,"R"] & \A_0^\Integer\ar[r,"F"]& \B_0^\Integer\ar[r,"H_n"]&\B\\
\Big(A\sr{f}{\ral}A'\Big)\ar[r,mapsto,"R"] &\Big(R_A\sr{R_f}{\ral}R_{A'}\Big)\ar[r,mapsto,"F"]& \Big(F(R_A)\sr{F(R_f)}{\ral}F(R_{A'})\Big)\ar[r,mapsto,"H_n"]&\Big((D^nF)A\sr{(D^nF)f}{\ral}(D^nF)A'\Big)
\et\nn
\eea}
where ~{\footnotesize $F(R_A)=F\Big(\cdots\ra R_A^{-1}\sr{d_A^{-1}}{\ral}R_A^0\sr{d_A^0}{\ral} R_A^1\sr{d_A^1}{\ral}\cdots\Big)~~=~~\cdots\ra F(R_A^{-1})\sr{F(d_A^{-1})}{\ral}F(R_A^0)\sr{F(d_A^0)}{\ral} F(R_A^1)\sr{F(d_A^1)}{\ral}\cdots$},
{\footnotesize\bea
\textstyle(D^nF)A=H_n\big(F(R_A)\big):={\ker F(d_A^n)\over \im F(d_A^{n-1})}:=\coker\Big(\im F(d_A^{n-1})\hookrightarrow\ker F(d_A^n)\Big),\nn
\eea}
and if $\B\subset\txt{Sets}$, then we explicitly have
\bea
\label{DerFunEq1}(D^nF)f=H_n\big(F(R_f)\big):=F(R_f^n)\big|_{\ker F(d_A^n)}+\im F(d_{A'}^{n-1}).
\eea

{\flushleft 2. \ul{Cofunctor $F$}}: The $n$th derived functor {\footnotesize $D^nF=H_n\circ F\circ R:\A\sr{R}{\ral}\A_0^\Integer\sr{F}{\ral}\B_0^\Integer\sr{H_n}{\ral}\B$} of $F$ is given by
{\footnotesize
\bea
&&\bt[column sep=small,row sep=tiny]
\A\ar[r,"R"] & \A_0^\Integer\ar[r,"F"]& \B_0^\Integer\ar[r,"H_n"]&\B\\
\Big(A\sr{f}{\ral}A'\Big)\ar[r,mapsto,"R"] &\Big(R_A\sr{R_f}{\ral}R_{A'}\Big)\ar[r,mapsto,"F"]& \Big(F(R_A)\sr{F(R_f)}{\lal}F(R_{A'})\Big)\ar[r,mapsto,"H_n"]&\Big((D^nF)A\sr{(D^nF)f}{\lal}(D^nF)A'\Big)
\et\nn
\eea}
where ~{\footnotesize $F(R_A)=F\Big(\cdots\ra R_A^{-1}\sr{d_A^{-1}}{\ral}R_A^0\sr{d_A^0}{\ral} R_A^1\sr{d_A^1}{\ral}\cdots\Big)~~=~~\cdots\la F(R_A^{-1})\sr{F(d_A^{-1})}{\lal}F(R_A^0)\sr{F(d_A^0)}{\lal} F(R_A^1)\sr{F(d_A^1)}{\lal}\cdots$},
{\footnotesize\bea
\textstyle(D^nF)A=H_n\big(F(R_A)\big):={\ker F(d_A^{-(n+1)})\over \im F(d_A^{-n})}:=\coker\Big(\im F(d_A^{-n})\hookrightarrow\ker F(d_A^{-(n+1)})\Big),\nn
\eea}
and if $\B\subset\txt{Sets}$, then we explicitly have
\bea
\label{DerFunEq2}(D^nF)f=H_n\big(F(R_f)\big):=F(R_f^n)\big|_{\ker F(d_{A'}^{-(n+1)})}+\im F(d_A^{-n}).
\eea
\end{dfn}

We note that the functors $D^nF:\A\ra\B$ are well-defined because, by the comparison theorems (Theorems \ref{CompThmI},\ref{CompThmII}), they are independent of the choice of projective or injective resolution.

\begin{dfn}[\textcolor{blue}{Examples of derived functors:
\index{Classical derived functors! Tor functors}{Tor functors},
\index{Classical derived functors! tor functors}{tor functors},
\index{Classical derived functors! Ext functors}{Ext functors},
\index{Classical derived functors! ext functors}{ext functors}}]
Let $\A$ be an abelian category and {\small $A\in\Ob\A$}. Consider a deleted projective (resp. injective) resolution
~{\footnotesize $
P^A_\ast:~\cdots\sr{d_3^A}{\ral} P_2^A\sr{d_2^A}{\ral}P_1^A\sr{d_1^A}{\ral} P_0^A\ra 0~~\big(\txt{resp.}~~I^\ast_A:~0\ra I_A^0\sr{d_A^0}{\ral}I_A^1\sr{d_A^1}{\ral} I_A^2\sr{d_A^2}{\ral}\cdots~\big)\nn
$}~ of $A$.
{\flushleft 1. \ul{Torsion functors (Left-derived functors)}}:
\bit
\item[(i)] Let ~$F=:(-)\otimes_\A X~:~\A\ral\B~(=\A)$,~~ for some object $X\in\Ob\A$. Then
{\footnotesize\bea
\textstyle Tor_n^\A(A,X):=\big(D_nF\big)A=H_n(P_\ast^A\otimes_\A X)={\ker(d_n^A\otimes id_X)\over\txt{im}(d_{n+1}^A\otimes id_X)}:=\coker\left(\im(d_{n+1}^A\otimes id_X)\hookrightarrow\ker(d_n^A\otimes id_X)\right).\nn
\eea}
\item[(ii)] Let ~$F:=X\otimes_\A(-)~:~\A\ral\B~(=\A)$,~~ for some object $X\in\Ob\A$. Then
{\footnotesize\bea
\textstyle tor_n^\A(X,A):=\big(D_nF\big)A=H_n(X\otimes_\A P_\ast^A)={\ker(id_X\otimes d_n^A)\over\txt{im}(id_X\otimes d_{n+1}^A)}:=\coker\left(\im(id_X\otimes d_{n+1}^A)\hookrightarrow\ker(id_X\otimes d_n^A)\right).\nn
\eea}
\eit
{\flushleft 2. \ul{Extension functors (Right-derived functors)}}:
\bit
\item[(i)] Let ~$F:=\Mor_\A(X,-)~:~\A\ral\B~(=Ab)$,~~ for some object $X\in\Ob\A$. Then
{\footnotesize\bea
Ext_\A^n(X,A):=\textstyle\big(D^nF\big)A=H_n\big(\Mor_\A(X,I^\ast_A)\big)={\ker\big(\Mor_\A(X,d_A^n)\big)\over\txt{im}\big(\Mor_\A(X,d_A^{n-1})\big)}.\nn
\eea}
\item[(ii)] Let ~$F:=\Mor_\A(-,X)~:~\A\ral\B~(=Ab)$,~~ for some object $X\in\Ob\A$. Then
{\footnotesize\bea
ext_\A^n(A,X):=\textstyle\big(D^nF\big)A=H_n\big(\Mor_\A(P^A_\ast,X)\big)={\ker\big(\Mor_\A(d_{n+1}^A,X)\big)\over\txt{im}\big(\Mor_\A(d_n^A,X)\big)}.\nn
\eea}
\eit
\end{dfn}

\section{Derived functor LES's: Estimating partial exactness in terms of derived functors}
\begin{prp}[\textcolor{blue}{\index{LES of derived functors}{LES of derived functors}: Exact continuation of $F(\txt{SES})$ using derived functors of $F$}]
Let $\A,\B$ be abelian categories such that (i) $\A$ admits a functorial deleted-resolution operation $R:\A\ra\A_0^\Integer$, ~$A\sr{f}{\ral}A'$ ~$\mapsto$ ~$R_A\sr{R_f}{\ral}R_{A'}$, and (ii) $\B$ is a LES-category.
\begin{enumerate}[leftmargin=0.9cm]
\item If $F:\A\ra\B$ is a functor such that $F\circ R:\A\sr{R}{\ral}\A_0^\Integer\sr{F}{\ral}\B_0^\Integer$ is an exact functor, then for every SES $0\ra A\sr{f}{\ral}A'\sr{f'}{\ral}A''\ra 0$ in $\A$, there exists an induced LES of homologies in $\B$ of the form
\bea
\label{LESeq1}
\adjustbox{scale=0.9}{\bt[row sep=tiny]
\circ\circ\circ\ar[r,"\delta^{n-2}"]& (D^{n-1}F)A\ar[rr,"(D^{n-1}F)f"] && (D^{n-1}F)A'\ar[rr,"(D^{n-1}F)f'"]&& (D^{n-1}F)A''\ar[r,"\delta^{n-1}"] &~ \\
      \ar[r,"\delta^{n-1}"]& (D^nF)A\ar[rr,"(D^nF)f"]        && (D^nF)A'\ar[rr,"(D^nF)f'"]        && (D^nF)A''\ar[r,"\delta^n"]         &~ \\
      \ar[r,"\delta^n"]& (D^{n+1}F)A\ar[rr,"(D^{n+1}F)f"]&& (D^{n+1}F)A'\ar[rr,"(D^{n+1}F)f'"]&& (D^{n+1}F)A''\ar[r,"\delta^{n+1}"] & \circ\circ\circ
\et}\eea
 where ~$(D^nF)(A):=H_n(F(R_A))=(H_n\circ F\circ R)(A)$. (\blue{footnote}\footnote{(i) If the functor $F$ is right-exact and the resolution used is projective, then $D^nF$ is a left-derived functor.\\
 (ii) If the functor $F$ is left-exact and the resolution used is injective, then $D^nF$ is a right-derived functor.})

\item If $F:\A\ra\B$ is a cofunctor such that $F\circ R:\A\sr{R}{\ral}\A_0^\Integer\sr{F}{\ral}\B_0^\Integer$ is an exact cofunctor, then for every SES $0\ra A\sr{f}{\ral}A'\sr{f'}{\ral}A''\ra 0$ in $\A$, there exists an induced LES of homologies in $\B$ of the form
\bea
\label{LESeq2}
\adjustbox{scale=0.9}{\bt[row sep=tiny]
\circ\circ\circ\ar[from=r,"\delta^{n-2}"']& (D^{n-1}F)A\ar[from=rr,"(D^{n-1}F)f"'] && (D^{n-1}F)A'\ar[from=rr,"(D^{n-1}F)f'"']&& (D^{n-1}F)A''\ar[from=r,"\delta^{n-1}"'] &~ \\
      ~\ar[from=r,"\delta^{n-1}"']& (D^nF)A\ar[from=rr,"(D^nF)f"']        && (D^nF)A'\ar[from=rr,"(D^nF)f'"']        && (D^nF)A''\ar[from=r,"\delta^n"']         &~ \\
      ~\ar[from=r,"\delta^n"']& (D^{n+1}F)A\ar[from=rr,"(D^{n+1}F)f"']&& (D^{n+1}F)A'\ar[from=rr,"(D^{n+1}F)f'"']&& (D^{n+1}F)A''\ar[from=r,"\delta^{n+1}"'] & \circ\circ\circ
\et}\eea
which is also written as
\bea
\adjustbox{scale=0.9}{\bt[row sep=tiny]
\circ\circ\circ\ar[r,"\delta^{n+1}"]& (D^{n+1}F)A''\ar[rr,"(D^{n+1}F)f'"] && (D^{n+1}F)A'\ar[rr,"(D^{n+1}F)f"]&& (D^{n+1}F)A\ar[r,"\delta^n"] &~ \\
      \ar[r,"\delta^n"]& (D^nF)A''\ar[rr,"(D^nF)f'"] && (D^nF)A'\ar[rr,"(D^nF)f"]&& (D^nF)A\ar[r,"\delta^{n-1}"] &~ \\
      \ar[r,"\delta^{n-1}"]& (D^{n-1}F)A''\ar[rr,"(D^{n-1}F)f'"] && (D^{n-1}F)A'\ar[rr,"(D^{n-1}F)f"]&& (D^{n-1}F)A\ar[r,"\delta^{n-2}"] &~ \circ\circ\circ
\et}\nn\eea
where ~$(D^nF)(A):=H_n(F(R_A))=(H_n\circ F\circ R)(A)$. (\blue{footnote}\footnote{(a) If the cofunctor $F$ is left-exact and the resolution used is projective, $D^nF$ is a right-derived functor.\\
(b) If the cofunctor $F$ is right-exact cofunctor and the resolution used is injective, $D^nF$ is a left-derived functor.})
\end{enumerate}
\end{prp}
\begin{proof}
By hypotheses, associated with the SES ~$S:~0\ra A\sr{f}{\ral}A'\sr{f'}{\ral}A''\ra 0$~ is a SES of deleted resolutions ~$R S:~0\ra R_A\sr{R_f}{\ral}R_{A'}\sr{R_{f'}}{\ral}R_{A''}\ra 0$, as in the diagram
\bea
\label{GeneralResEq}
\adjustbox{scale=0.9}{\bt
& \ar[d] && \ar[d] && \ar[d] &  \\
0\ar[r]& R_A^{n-1}\ar[d,"d_A^{n-1}"]\ar[rr,"R_f^{n-1}"] && R_{A'}^{n-1}\ar[d,"d_{A'}^{n-1}"]\ar[rr,"R_{f'}^{n-1}"] && R_{A''}^{n-1}\ar[d,"d_{A''}^{n-1}"]\ar[r] & 0\\
0\ar[r]& R_A^n\ar[d,"d_A^n"]\ar[rr,"R_f^n"] && R_{A'}^n\ar[d,"d_{A'}^n"]\ar[rr,"R_{f'}^n"] && R_{A''}^n\ar[d,"d_{A''}^n"]\ar[r] & 0\\
0\ar[r]& R_A^{n+1}\ar[d,"d_A^{n+1}"]\ar[rr,"R_f^{n+1}"] && R_{A'}^{n+1}\ar[d,"d_{A'}^{n+1}"]\ar[rr,"R_{f'}^{n+1}"] && R_{A''}^{n+1}\ar[d,"d_{A''}^{n+1}"]\ar[r] & 0\\
& ~ && ~ && ~ &
\et}\eea
\begin{enumerate}
\item Applying the functor $F$ on the SES of resolutions (\ref{GeneralResEq}) above, we get the exact sequence
\bea
(F\circ R)S:~0\ra F(R_A)\sr{F(R_f)}{\ral}F(R_{A'})\sr{F(R_{f'})}{\ral}F(R_{A''})\ra 0\nn
\eea
as in the diagram
\bc\adjustbox{scale=0.9}{\bt
& \ar[d] && \ar[d] && \ar[d] &  \\
0\ar[r] & F(R_A^{n-1})\ar[d,"F(d_A^{n-1})"]\ar[rr,"F(R_f^{n-1})"] && F(R_{A'}^{n-1})\ar[d,"F(d_{A'}^{n-1})"]\ar[rr,"F(R_{f'}^{n-1})"] && F(R_{A''}^{n-1})\ar[d,"F(d_{A''}^{n-1})"]\ar[r] & 0\\
0\ar[r] & F(R_A^n)\ar[d,"F(d_A^n)"]\ar[rr,"F(R_f^n)"] && F(R_{A'}^n)\ar[d,"F(d_{A'}^n)"]\ar[rr,"F(R_{f'}^n)"] && F(R_{A''}^n)\ar[d,"F(d_{A''}^n)"]\ar[r] & 0\\
0\ar[r] & F(R_A^{n+1})\ar[d,"F(d_A^{n+1})"]\ar[rr,"F(R_f^{n+1})"] && F(R_{A'}^{n+1})\ar[d,"F(d_{A'}^{n+1})"]\ar[rr,"F(R_{f'}^{n+1})"] && F(R_{A''}^{n+1})\ar[d,"F(d_{A''}^{n+1})"]\ar[r] & 0\\
& ~ && ~ && ~ &
\et}\ec
Since $\B$ is a LES category, we get the induced LES of homologies in (\ref{LESeq1}).

\item Applying the cofunctor $F$ on the SES of resolutions (\ref{GeneralResEq}) above, we get the exact sequence
\bea
(F\circ R)S:~0\la F(R_A)\sr{F(R_f)}{\lal}F(R_{A'})\sr{F(R_{f'}^\ast)}{\lal}F(R_{A''})\la 0\nn
\eea
as in the diagram
\bc\adjustbox{scale=0.9}{\bt
& ~ && ~ && ~ &  \\
0 & F(R_A^{n-1})\ar[u]\ar[l] && F(R_{A'}^{n-1})\ar[u]\ar[ll,"F(R_f^{n-1})"'] && F(R_{A''}^n)\ar[u]\ar[ll,"F(R_{f'}^{n-1})"'] & 0\ar[l]\\
0 & F(R_A^n)\ar[l]\ar[u,"F(d_A^{n-1})"'] && F(R_{A'}^n)\ar[ll,"F(R_f^n)"']\ar[u,"F(d_{A'}^{n-1})"'] && F(R_{A''}^n)\ar[ll,"F(R_{f'}^n)"']\ar[u,"F(d_{A''}^{n-1})"'] & 0\ar[l]\\
0 & F(R_A^{n+1})\ar[l]\ar[u,"F(d_A^n)"'] && F(R_{A'}^{n+1})\ar[ll,"F(R_f^{n+1})"']\ar[u,"F(d_{A'}^n)"'] && F(R_{A''}^{n+1})\ar[ll,"F(R_{f'}^{n+1})"']\ar[u,"F(d_{A''}^n)"'] & 0\ar[l]\\
& \ar[u,"F(d_A^{n+1})"'] && \ar[u,"F(d_A{A'}^{n+1})"'] && \ar[u,"F(d_{A''}^{n+1})"'] &
\et}\ec
Since $\B$ is a LES category, we get the induced LES of homologies in (\ref{LESeq2}). \qedhere
\end{enumerate}
\end{proof}

\begin{crl}[\textcolor{blue}{\index{Tor and Ext sequences}{Tor and Ext sequences}}]
Let $\A$ be an abelian category such that the following hold:
\bit[]
\item[(i)] $\A$ has enough projectives, and admits a projective \ul{split-exact} functorial  deleted-resolution operation
\bea
P_\ast:\A\ra\A_0^\Integer,~~A\sr{f}{\ral}A'~~\mapsto~~P_\ast^A\sr{P^f_\ast}{\ral}P_\ast^{A'}.\nn
\eea
\item[(ii)] $\A$ has enough injectives, and admits an injective \ul{split-exact} functorial  deleted-resolution operation
\bea
I^\ast:\A\ra\A_0^\Integer,~~A\sr{f}{\ral}A'~~\mapsto~~I^\ast_A\sr{I_f^\ast}{\ral}I^\ast_{A'}.\nn
\eea
\item[(iii)] $\A$ is a LES-category. (Required by parts (2) and (3) below only.)
\eit
Then we have have the following observations:
\begin{enumerate}[leftmargin=0.9cm]
\item If $F:=-\otimes_\A X:\A\ra\A$, $X\in\Ob\A$ (\blue{footnote}\footnote{It is clear (from the horseshoe lemmas) that
$F\circ P_\ast:\A\sr{P_\ast}{\ral}\A_0^\Integer\sr{F}{\ral}\A_0^\Integer$ is an exact functor.}), then for every SES $0\ra A\sr{f}{\ral}A'\sr{f'}{\ral}A''\ra 0$ in $\A$, there exists an induced LES of homologies in $\A$ of the form

{\footnotesize\[
\cdots\ral Tor_1^\A(A,X)\ra Tor_1^\A(A',X)\ra Tor_1^\A(A'',X)\ra A\otimes_\A X\ra A'\otimes_\A X\ra A''\otimes_\A X\ra 0.
\]}
\item If $F:=X\otimes_\A-:\A\ra\A$, $X\in\Ob\A$ (\blue{footnote}\footnote{It is clear (from the horseshoe lemmas) that $F\circ P_\ast:\A\sr{P_\ast}{\ral}\A_0^\Integer\sr{F}{\ral}\A_0^\Integer$ is an exact functor.}), then for every SES $0\ra A\sr{f}{\ral}A'\sr{f'}{\ral}A''\ra 0$ in $\A$, there exists an induced LES of homologies in $\A$ of the form

{\footnotesize\[
\cdots\ral tor_1^\A(X,A)\ra tor_1^\A(X,A')\ra tor_1^\A(X,A'')\ra X\otimes_\A A\ra X\otimes_\A A'\ra X\otimes_\A A''\ra 0.
\]}
\item If $F:=\Mor_\A(X,-):\A\ra Ab$, $X\in\Ob\A$ (\blue{footnote}\footnote{It is clear (from the horseshoe lemmas) that
$F\circ I^\ast:\A\sr{I^\ast}{\ral}\A_0^\Integer\sr{F}{\ral}(Ab)_0^\Integer$ is an exact functor.}), then for every SES $0\ra A\sr{f}{\ral}A'\sr{f'}{\ral}A''\ra 0$ in $\A$, there exists an induced LES of homologies in $Ab$ of the form

{\footnotesize\[
0\ra \Mor_\A(X,A)\ra \Mor_\A(X,A')\ra \Mor_\A(X,A'')\ra Ext_\A^1(X,A)\ra Ext_\A^1(X,A')\ra Ext_\A^1(X,A'')\ral\cdots
\]}
\item If $F:=\Mor_\A(-,X):\A\ra Ab$, $X\in\Ob\A$ (\blue{footnote}\footnote{It is clear (from the horseshoe lemmas) that $F\circ P_\ast:\A\sr{P_\ast}{\ral}\A_0^\Integer\sr{F}{\ral}(Ab)_0^\Integer$ is an exact cofunctor.}), then for every SES $0\ra A\sr{f}{\ral}A'\sr{f'}{\ral}A''\ra 0$ in $\A$, there exists an induced LES of homologies in $Ab$ of the form

{\footnotesize\[
0\ra \Mor_\A(A'',X)\ra \Mor_\A(A',X)\ra \Mor_\A(A,X)\ra ext_\A^1(A'',X)\ra ext_\A^1(A',X)\ra ext_\A^1(A,X)\ral\cdots\nn
\]}
\end{enumerate}
\end{crl}
\begin{proof}
By hypotheses, we have a SES of deleted projective resolutions with \ul{split-exact rows}:
\bea
0\ra P_\ast^A\sr{P^f_\ast}{\ral} P_\ast^{A'}\sr{P^{f'}_\ast}{\ral} P_\ast^{A''}\ra 0~~~~=~~
\adjustbox{scale=0.8}{\begin{tikzcd}
         & \vdots\ar[d,"d_2^A"] & \vdots\ar[d,"d_2^{A'}"] & \vdots\ar[d,"d_2^{A''}"] &  \\
 0\ar[r] & P_1^A\ar[d,"d_1^A"]\ar[r,"P^f_1"] & P_1^{A'}\ar[d,"d_1^{A'}"]\ar[r,"P^{f'}_1"] & P_1^{A''}\ar[d,"d_1^{A''}"]\ar[r] &  0 \\
 0\ar[r] & P_0^A\ar[d]\ar[r,"P^f_0"] & P_0^{A'}\ar[d]\ar[r,"P^{f'}_0"] & P_0^{A''}\ar[d]\ar[r] &  0 \\
         & 0 & 0 & 0 &
\end{tikzcd}}\label{del-pres-diag}
\eea
~~\hspace{0.5cm}~~and a SES of deleted injective resolutions with \ul{split-exact rows}:
\bea
0\ra I^\ast_A\sr{I_f^\ast}{\ral} I^\ast_{A'}\sr{I_{f'}^\ast}{\ral} I^\ast_{A''}\ra 0~~~~=~~
\adjustbox{scale=0.8}{\begin{tikzcd}
         & 0\ar[d] & 0\ar[d] & 0\ar[d] &   \\
 0\ar[r] & I_A^0\ar[d,"d_A^0"]\ar[r,dashed,"I_f^0"] & I_{A'}^0\ar[d,dashed,"d_{A'}^0"]\ar[r,dashed,"I_{f'}^0"] & I_{A''}^0\ar[d,"d_{A''}^0"]\ar[r] &  0 \\
 0\ar[r] & I_A^1\ar[d,"d_A^1"]\ar[r,dashed,"I_f^1"] & I_{A'}^1\ar[d,"d_{A'}^1"]\ar[r,dashed,"I_{f'}^1"] & I_{A''}^1\ar[d,"d_{A''}^1"]\ar[r] &  0 \\
         & \vdots & \vdots & \vdots &
\end{tikzcd}}\label{del-ires-diag}
\eea

\begin{enumerate}[leftmargin=0.9cm]
\item \ul{Tor sequence}:~ Applying the functor $-\otimes_\A X$ [\blue{footnote}\footnote{Which is additive and thus exact on the split exact rows.}] on (\ref{del-pres-diag}) we get a SES of complexes ~$0\ra P_\ast^A\otimes_\A X\ra P_\ast^{A'}\otimes_\A X\ra P_\ast^{A''}\otimes_\A X\ra 0$,~ as in the diagram
\[\adjustbox{scale=0.8}{\begin{tikzcd}
         & \vdots\ar[dd,"d_2^A\otimes id_X"] && \vdots\ar[dd,"d_2^{A'}\otimes id_X"] && \vdots\ar[dd,"d_2^{A''}\otimes id_X"] &  \\
         & && && & \\
 0\ar[r] & P_1^A\otimes_\A X\ar[dd,"d_1^A\otimes id_X"]\ar[rr,"P^f_1\otimes id_X"] && P_1^{A'}\otimes_\A X\ar[dd,"d_1^{A'}\otimes id_X"]\ar[rr,"P^{f'}_1\otimes id_X"] && P_1^{A''}\otimes_\A X\ar[dd,"d_1^{A''}\otimes id_X"]\ar[r] &  0 \\
         & && && & \\
 0\ar[r] & P_0^A\otimes_\A X\ar[d]\ar[rr,"P^f_0\otimes id_X"] && P_0^{A'}\otimes_\A X\ar[d]\ar[rr,"P^{f'}_0\otimes id_X"] && P_0^{A''}\otimes_\A X\ar[d]\ar[r] &  0 \\
         & 0 && 0 && 0 &
\end{tikzcd}}\]
By the LES property of $\A$, we get the exact sequence
\begin{align}
&\cdots\ra Tor_1^\A(A,X)\ra Tor_1^\A(A',X)\ra Tor_1^\A(A'',X) \nn\\
&~~~~~~~~~~~~\ra Tor_0^\A(A,X)\ra Tor_0^\A(A',X)\ra Tor_0^\A(A'',X)\ra 0,\nn
\end{align}
where the right-exactness of $-\otimes_\A X$ (on a projective resolution) gives
{\small
\begin{align}
&\Big(\cdots\sr{d_1^{A}}{\ral} P_0^{A}\sr{\vep_{A}}{\ral}A\ra 0\Big)\otimes_\A X~~=~~\cdots\sr{d_1^{A}\otimes id_X}{\ral} P_0^{A}\otimes_\A X\sr{\vep_{A}\otimes id_X}{\ral}A\otimes_\A X\ra 0,\nn\\
&~~\Ra~~Tor_0^\A(A,X):=\coker\left(\im(d_1^A\otimes id_X)\hookrightarrow\ker(d_0^A\otimes id_X)\right)\nn\\
&~~~~~~~~~~~~~~~~~~~~=\coker\left(\im(d_1^A\otimes id_X)\hookrightarrow\ker(0\otimes id_X)\right)=\coker\left(\im(d_1^A\otimes id_X)\hookrightarrow \ker(0)\right)\nn\\
&~~~~~~~~~~~~~~~~~~~~=\coker\left(\im(d_1^A\otimes id_X)\hookrightarrow P_0^{A}\otimes_\A X\right)\sr{\txt{exactness}}{\cong}\coker\left(\ker(\vep_{A}\otimes id_X)\hookrightarrow P_0^{A}\otimes_\A X\right)\nn\\
&~~~~~~~~~~~~~~~~~~~~=\im(\vep_{A}\otimes id_X)\cong A\otimes_\A X.\nn
\end{align}}

\item \ul{tor sequence}:~ Applying the functor $X\otimes_\A-$ [\blue{footnote}\footnote{Which is additive and thus exact on the split exact rows.}] on (\ref{del-pres-diag}) we get a SES of complexes ~$0\ra X\otimes_\A P_\ast^A\ra X\otimes_\A P_\ast^{A'}\ra X\otimes_\A P_\ast^{A''}\ra 0$,~ as in the diagram
\[\adjustbox{scale=0.8}{\begin{tikzcd}
         & \vdots\ar[dd,"id_X\otimes d_2^A"] && \vdots\ar[dd,"id_X\otimes d_2^{A'}"] && \vdots\ar[dd,"id_X\otimes d_2^{A''}"] &  \\
         & && && & \\
 0\ar[r] & X\otimes_\A P_1^A\ar[dd,"id_X\otimes d_1^A"]\ar[rr,"id_X\otimes P^f_1"] && X\otimes_\A P_1^{A'}\ar[dd,"id_X\otimes d_1^{A'}"]\ar[rr,"id_X\otimes P^{f'}_1"] && X\otimes_\A P_1^{A''}\ar[dd,"id_X\otimes d_1^{A''}"]\ar[r] &  0 \\
         & && && & \\
 0\ar[r] & X\otimes_\A P_0^A\ar[d]\ar[rr,"id_X\otimes P^f_0"] && X\otimes_\A P_0^{A'}\ar[d]\ar[rr,"id_X\otimes P^{f'}_0"] && X\otimes_\A P_0^{A''}\ar[d]\ar[r] &  0 \\
         & 0 && 0 && 0 &
\end{tikzcd}}\]
By the LES property of $\A$, we get the exact sequence
\begin{align}
&\cdots\ra tor_1^\A(X,A)\ra tor_1^\A(X,A')\ra tor_1^\A(X,A'') \nn\\
&~~~~~~~~~~~~\ra tor_0^\A(X,A)\ra tor_0^\A(X,A')\ra tor_0^\A(X,A'')\ra 0,\nn
\end{align}
where the right-exactness of $X\otimes_\A-$ (on a projective resolution) gives
{\small
\begin{align}
&X\otimes_\A\Big(\cdots\sr{d_1^{A}}{\ral} P_0^{A}\sr{\vep_{A}}{\ral}A\ra 0\Big)~~=~~\cdots\sr{id_X\otimes d_1^{A}}{\ral}X\otimes_\A P_0^{A}\sr{id_X\otimes\vep_{A}}{\ral}X\otimes_\A A\ra 0,\nn\\
&~~\Ra~~tor_0^\A(X,A):=\coker\left(\im(id_X\otimes d_1^A)\hookrightarrow\ker(id_X\otimes d_0^A)\right)\nn\\
&~~~~~~~~~~~~~~~~=\coker\left(\im(id_X\otimes d_1^A)\hookrightarrow\ker(id_X\otimes 0)\right)=\coker\left(\im(id_X\otimes d_1^A)\hookrightarrow\ker(0)\right)\nn\\
&~~~~~~~~~~~~~~~~=\coker\left(\im(id_X\otimes d_1^A)\hookrightarrow X\otimes_\A P_0^{A}\right)\sr{\txt{exactness}}{\cong}\coker\left(\ker(id_X\otimes\vep_{A})\hookrightarrow X\otimes_\A P_0^{A}\right)\nn\\
&~~~~~~~~~~~~~~~~=\im(id_X\otimes\vep_{A})\cong X\otimes_\A A.\nn
\end{align}}

\item \ul{Ext sequence}:~ Applying the functor $\Mor_\A(X,-)$ [\blue{footnote}\footnote{Which is additive and thus exact on the split exact rows.}] on (\ref{del-ires-diag}) we get a SES of complexes
\[
0\ra \Mor_\A(X,I^\ast_A)\ra \Mor_\A(X,I^\ast_{A'})\ra \Mor_\A(X,I^\ast_{A''})\ra 0
\]
as in the diagram
\[\adjustbox{scale=0.8}{\bt
         & 0\ar[d] && 0\ar[d] && 0\ar[d] &   \\
 0\ar[r] & \Mor_\A(X,I_A^0)\ar[dd,"{\Mor_\A(X,d_A^0)}"]\ar[rr,"{\Mor_\A(X,I_f^0)}"] && \Mor_\A(X,I_{A'}^0)\ar[dd,"{\Mor_\A(X,d_{A'}^0)}"]\ar[rr,"{\Mor_\A(X,I_{f'}^0)}"] && \Mor_\A(X,I_{A''}^0)\ar[dd,"{\Mor_\A(X,d_{A''}^0)}"]\ar[r] &  0 \\
          & && && & \\
 0\ar[r] & \Mor_\A(X,I_A^1)\ar[dd,"{\Mor_\A(X,d_A^1)}"]\ar[rr,"{\Mor_\A(X,I_f^1)}"] && \Mor_\A(X,I_{A'}^1)\ar[dd,"{\Mor_\A(X,d_{A'}^1)}"]\ar[rr,"{\Mor_\A(X,I_{f'}^1)}"] && \Mor_R(X,I_{A''}^1)\ar[dd,"{\Mor_\A(X,d_{A''}^1)}"]\ar[r] &  0 \\
         & && && & \\
         & \vdots && \vdots && \vdots &
\et}\]
By the LES property of Ab=$\Integer$-mod, we get the exact sequence
\begin{align}
&0\ra Ext_\A^0(X,A)\ra Ext_\A^0(X,A')\ra Ext_\A^0(X,A'') \nn\\
&~~~~~~~~~~~~\ra Ext_\A^1(X,A)\ra Ext_\A^1(X,A')\ra Ext_\A^1(X,A'')\ral\cdots\nn
\end{align}
where the left-exactness of Hom (on an injective resolution $0\ra A\ra I_A^\ast$) gives
{\footnotesize
\begin{align}
&\Mor_\A\Big(X~,~0\ra A\sr{\vep_A}{\ral}I_A^0\sr{d_A^0}{\ral} \cdots\Big)~~=~~0\ra \Mor_\A(X,A)\sr{\Mor_\A(X,\vep_A)}{\ral}\Mor_\A(X,I_A^0)\sr{\Mor_\A(X,d_A^0)}{\ral} \cdots,\nn\\
&~~\Ra~~Ext_\A^0(X,A):=\coker\left(\im\big(\Mor_\A(X,d_A^{-1})\big)\hookrightarrow\ker\big(\Mor_\A(X,d_A^0)\big)\right)\nn\\
&~~~~~~~~~~~~=\coker\left(\im\big(\Mor_\A(X,0)\big)\hookrightarrow\ker\big(\Mor_\A(X,d_A^0)\big)\right)=\coker\left(0\ra\ker\big(\Mor_\A(X,d_A^0)\big)\right)\nn\\
&~~~~~~~~~~~~=\ker \Mor_\A(X,d_A^0)\cong\im \Mor_\A(X,\vep_A)\cong \Mor_\A(X,A).\nn
\end{align}}

\item \ul{ext sequence}:~ Applying the functor $\Mor_\A(-,X)$ [\blue{footnote}\footnote{Which is additive and thus exact on the split exact rows.}] on (\ref{del-pres-diag}) we get a SES of complexes  $$0\ra \Mor_\A(P_\ast^{A''},X)\ra \Mor_\A(P_\ast^{A'},X)\ra \Mor_\A(P_\ast^A,X)\ra 0$$ as in the diagram

\[\adjustbox{scale=0.8}{\begin{tikzcd}
         & 0\ar[d] && 0\ar[d] && 0\ar[d] &   \\
 0\ar[r] & \Mor_\A(P_0^{A''},X)\ar[dd,"{\Mor_\A(d_1^{A''},X)}"]\ar[rr,"{\Mor_\A(P^{f'}_0,X)}"] && \Mor_\A(P_0^{A'},X)\ar[dd,"{\Mor_\A(d_1^{A'},X)}"]\ar[rr,"{\Mor_\A(P^f_0,X)}"] && \Mor_\A(P_0^A,X)\ar[dd,"{\Mor_\A(d_1^A,X)}"]\ar[r] &  0 \\
          & && && & \\
 0\ar[r] & \Mor_\A(P_1^{A''},X)\ar[dd,"{\Mor_\A(d_2^{A''},X)}"]\ar[rr,"{\Mor_\A(P^{f'}_1,X)}"] && \Mor_\A(P_1^{A'},X)\ar[dd,"{\Mor_\A(d_2^{A'},X)}"]\ar[rr,"{\Mor_\A(P^f_1,X)}"] && \Mor_\A(P_1^A,X)\ar[dd,"{\Mor_\A(d_2^A,X)}"]\ar[r] &  0 \\
         & && && & \\
         & \vdots && \vdots && \vdots &
\end{tikzcd}}\]
By the LES property of Ab=$\Integer$-mod, we get the exact sequence
\begin{align}
&0\ra ext_\A^0(A'',X)\ra ext_\A^0(A',X)\ra ext_\A^0(A,X) \nn\\
&~~~~~~~~~~~~\ra ext_\A^1(A'',X)\ra ext_\A^1(A',X)\ra ext_\A^1(A,X)\ra\cdots,\nn
\end{align}
where the left-exactness of Hom (on a projective resolution $P^A_\ast\ra A\ra 0$) gives
{\footnotesize
\begin{align}
&\Mor_\A\Big(\cdots\sr{d_1^{A}}{\ral} P_0^{A}\sr{\vep_{A}}{\ral}A\ra 0~,~X\Big)~~=~~0\ra \Mor_\A(A,X)\sr{\Mor_\A(\vep_{A},X)}{\ral}\Mor_\A(P_0^{A},X)\sr{\Mor_\A(d_1^{A},X)}{\ral}\cdots,\nn\\
&~~\Ra~~ext_\A^0(A,X):=\coker\left(\im\big(\Mor_\A(d_0^A,X)\big)\hookrightarrow\ker\big(\Mor_\A(d_1^A,X)\big)\right)\nn\\
&~~~~~~~~~~~~~~~~=\coker\left(\im\big(\Mor_\A(0,X)\big)\hookrightarrow\ker\big(\Mor_\A(d_1^A,X)\big)\right)=\coker\left(0\ra\ker\big(\Mor_\A(d_1^A,X)\big)\right)\nn\\
&~~~~~~~~~~~~~~~~=\ker \Mor_\A(d_1^{A},X)\sr{\txt{exactness}}{\cong}\im \Mor_\A(\vep_{A},X)\cong \Mor_\A(A,X).\nn \qedhere
\end{align}}
\end{enumerate}
\end{proof}

\begin{rmk}[\textcolor{blue}{\index{Derived! functor rules}{Derived functor rules}: A summary}] In the following functor correspondences, the 1st \ul{blank argument} uses projective resolutions while the 2nd \ul{blank argument} uses injective resolutions.
\begin{enumerate}
\item ~~$-\otimes_\A X~\longmapsto~Tor_\ast^\A(-,X)$~~~~and~~~~$X\otimes_\A-~\longmapsto~tor_\ast^\A(X,-)$
\item ~~$\Mor_\A(X,-)~\longmapsto~Ext^\ast_\A(X,-)$~~~~and~~~~$\Mor_\A(-,X)~\longmapsto~ext^\ast_\A(-,X)$
\item \ul{Right-exact functors} {\small$-\otimes_\A X$, $X\otimes_\A-$} give \ul{left-derived functors} {\small $Tor_\A(X,-)$, $tor_\A(-,X)$}.
\item \ul{Left-exact functors} {\small $\Mor_\A(X,-)$,~ $\Mor_\A(-,X)$} give \ul{right-derived functors} {\small $Ext_\A(X,-)$, $ext_\A(-,X)$}.
\end{enumerate}
\end{rmk}

\section{Characterization of quasiisomorphism: Mapping cones}
From here on, unless it is specified otherwise, $\A$ is an arbitrary abelian category with a tensor product $\otimes=\otimes_\A:\A\times\A\ra\A$ and a mor ``product''/operator $\Mor=\Mor_\A:\A\times\A\ra Ab$. As before, we will frequently work in $R$-mod for some $R$, but (by the Freyd-Mitchell imbedding) any results that require only a ``small-abelian substructure'' of $R$-mod remain true if $R$-mod is replaced with any abelian category $\A$.
\begin{dfn}[\textcolor{blue}{\index{Translated chain complex}{Translated chain complex},
\index{Negative of a chain complex}{Negative of a chain complex}}]
Let $(A_\ast,d_\ast^{A_\ast})$ be a chain complex in $\A$ (i.e., $A_\ast\in\A_0^\Integer$) and $k\in\Integer$. The \ul{$k$-translation} ~$\big(A_\ast[k],d_\ast^{A_\ast[k]}\big)$~ of $(A_\ast,d_\ast^{A_\ast})$ is the chain complex given by
\bea
A_n[k]:=A_{n+k}~~~~\txt{and}~~~~d_n^{A_\ast[k]}:=(-1)^kd_{n+k}^{A_\ast}.\nn
\eea
The \ul{negative of} the chain complex $A_\ast=(A_\ast,d_\ast^{A_\ast})$ is the associated chain complex $-A_\ast=-(A_\ast,d_\ast^{A_\ast}):=(A_\ast,-d_\ast^{A_\ast})$.
\end{dfn}

Recall from Definition \ref{CoHlogyDef} that a morphism of complexes $f=f_\ast:A_\ast\ra B_\ast$ in $\A$ is a quasiisomorphism if $H(f_\ast)$ is an isomorphism. By viewing $f$ itself as the bicomplex $C_{\ast\ast}:0\ra A_\ast\sr{f}{\ral}B_\ast\ra 0$, we can describe its total complex, which we call ``\emph{the mapping cone of $f$}'' (as in the following definition).
\begin{dfn}[\textcolor{blue}{\index{Mapping cone}{Mapping cone}}]
The \ul{mapping cone} of a chain morphism {\small $f:A_\ast\ra B_\ast$} in $\A$ is defined as
{\small\[
C(f):=-\Tot\left(0\ra A_\ast\sr{f}{\ral}B_\ast\ra 0\right)=
-\Tot\left(\adjustbox{scale=0.9}{\bt
        & \substack{(i=1~\txt{or}~A)\\~\\\vdots}                                               & \substack{(i=0~\txt{or}~B)\\~\\\vdots}                                     &  \\
0\ar[r] & A_{j-1}\ar[u]\ar[u,"d_{j-1}^{A_\ast}"']\ar[r,"f_{j-1}"] & B_{j-1}\ar[u,"d_{j-1}^{B_\ast}"']\ar[u]\ar[r] & 0 \\
0\ar[r] & \overbrace{A_j}^{c_{0j}}\ar[u,"d_j^{A_\ast}"']\ar[r,"f_j"]                & \overbrace{B_j}^{c_{1j}}\ar[u,"d_j^{B_\ast}"']\ar[r]             & 0 \\
0\ar[r] & A_{j+1}\ar[u,"d_{j+1}^{A_\ast}"']\ar[r,"f_{j+1}"]                & B_{j+1}\ar[u,"d_{j+1}^{B_\ast}"']\ar[r]             & 0 \\
        & \vdots\ar[u]                                         & \vdots\ar[u]                          &
\et}\right)=\Tot(C_{\ast\ast}),
\]}
where ~$C_{ij}:=\delta_{i1}A_j\oplus\delta_{i0}B_j$, ~$d^h_{ij}:=-\delta_{i1}f_j$, ~$d^v_{ij}:=-(\delta_{i1}d^{A_\ast}_j+\delta_{i0}d^{B_\ast}_j)$, ~and ( assuming $\A\hookrightarrow\txt{Sets}$ )
\bea
&&d_n^{C(f)}(c_{ij})_{i+j=n}:=\Big(d^h_{i+1,j}c_{i+1,j}+(-1)^id^v_{i,j+1}c_{i,j+1}\Big)_{i+j=n}=\Big(d^hc_{i+1,n-i}+(-1)^id^vc_{i,n-i+1}\Big)_{i\in\{1,0\}}\nn\\
&&~~~~=\Big(d^hc_{2,n-1}-d^vc_{1,n}~,~d^hc_{1,n}+d^vc_{0,n+1}\Big)=-\Big(0-d^{A_\ast}a_n~,~f_na_n+d^{B_\ast}b_{n+1}\Big).\nn
\eea

\end{dfn}
That is (up to sign) the \ul{mapping cone} is $C(f):=A_\ast \oplus B_\ast[1]$, where $C(f)_n:=A_n\oplus B_{n+1}$, and as a diagram,
\bea\bt
C(f):\cdots\ar[r] & C(f)_{n+1}\ar[r,"d_{n+1}"] \ar[r] & C(f)_n\ar[r,"d_n"] & C(f)_{n-1}\ar[r,"d_{n-1}"] &\cdots
\et\nn
\eea
with differential (in Sets-imbedded matrix form) given by
\[
d_n:=d_n^{C(f)}=
-\left[
  \begin{array}{cc}
    -d_n^A & 0 \\
    f_n & d_{n+1}^B\\
  \end{array}
\right]:
\left[
  \begin{array}{l}
    a_n \\
    b_{n+1} \\
  \end{array}
\right]\longmapsto
-\left[
  \begin{array}{l}
    -d^A_na_n \\
    f_na_n+d^B_{n+1}b_{n+1} \\
  \end{array}
\right].
\]

\begin{prp}[\textcolor{blue}{\index{Quasiisomorphism criterion}{Quasiisomorphism criterion on $R$-mod}}]\label{ConeLESprp}
Let $f:A_\ast\ra B_\ast$ be a chain morphism in $R$-mod, i.e., a morphism in $(R\txt{-mod})_0^\Integer$.
{\flushleft (1)} There is a short exact sequence of complexes ~$0\ra B_\ast[1]\sr{i}{\ral}C(f)\sr{p}{\ral}A_\ast\ra 0$,~ where $i,p$ are the obvious inclusion and projection respectively.
{\flushleft (2)} Furthermore, in the associated induced LES of homologies,
\[
\cdots\ra H_{n+1}\big(A_\ast\big)\sr{\delta_{n+1}}{\ral}H_n(B_\ast[1])\ra H_n(C(f))\ra H_n\big(A_\ast\big)\sr{\delta_n}{\ral}H_{n-1}(B_\ast[1])\ra H_{n-1}(C(f))\ra\cdots,
\]
the connecting homomorphism $\delta_{n+1}=-H_{n+1}(f)$, i.e., $\delta$ is induced by $f$, and $H_n\big(B_\ast[1]\big)=H_{n+1}(B_\ast)$.

That is, the sequence takes the form

{\footnotesize\bea
\label{ConeLES1}\bt\cdots\ra H_{n+1}\big(A_\ast\big)\ar[rr,"{-H_{n+1}(f)}"]&& H_{n+1}(B_\ast)\ar[r,"{H_n(i)}"]& H_n(C(f))\ar[r,"{H_n(p)}"]& H_n\big(A_\ast\big)\ar[r,"{-H_n(f)}"]& H_n(B_\ast)\ar[r,"{H_{n-1}(i)}"]& H_{n-1}(C(f))\ra\cdots\et
\eea}
{\flushleft (3)} Hence, $f:A_\ast\ra B_\ast$ is a quasi-iso $\iff$ $H_n(C(f))=0$ ~for all $n$, $\iff$ $C(f)$ is exact.
\end{prp}
\begin{proof}
{\flushleft (1)} The rows are split exact. Thus, we just need to check directly that the diagram commutes.
\[\adjustbox{scale=0.9}{\bt
 & \vdots\ar[d] & \vdots\ar[d] & \vdots\ar[d] & \\
0\ar[r] & B_{n+2}\ar[d,"-d^B_{n+2}"]\ar[r,"i_{n+1}"] & A_{n+1}\oplus B_{n+2}\ar[d,dashed,"d_{n+1}^{C(f)}"]\ar[r,dashed,"p_{n+1}"] & A_{n+1}\ar[d,"d_{n+1}^A"]\ar[r,dashed] & 0\\
0\ar[r,dashed] & B_{n+1}\ar[d,"-d^B_{n+1}"]\ar[r,dashed,"i_n"] & A_n\oplus B_{n+1}\ar[d,"d_n^{C(f)}"]\ar[r,"p_n"] & A_n\ar[d,"d_n^A"]\ar[r] & 0 \\
0\ar[r] & B_n\ar[d,"-d^B_n"]\ar[r,"i_{n-1}"] & A_{n-1}\oplus B_n\ar[d,"d_{n-1}^{C(f)}"]\ar[r,"p_{n-1}"] & A_{n-1}\ar[d,"d_{n-1}^A"]\ar[r] & 0 \\
  & \vdots & \vdots & \vdots &
\et}\]
\begin{align}
& d_n^{C(f)}i_n(b_{n+1})=d_n^{C(f)}(0,b_{n+1})=-\left(-d^{A_\ast}_n(0),f_n(0)+d^{B_\ast}_{n+1}b_{n+1}\right)=-(0,d^{B_\ast}_{n+1}b_{n+1})=-i_{n-1}(d_{n+1}^{B_\ast}b_{n+1}),\nn\\
&~~\Ra~~d_n^{C(f)}i_n=-i_{n-1}d_{n+1}^{B_\ast}.\nn\\
& d_n^{A_\ast}p_n(a_n,b_{n+1})=d_n^{A_\ast}a_n=-p_{n-1}(-d_n^{A_\ast}a_n,f_na_n+d_{n+1}^{B_\ast}b_{n+1})=p_{n-1}d_n^{C(f)}(a_n,b_{n+1})\nn\\
&~~\Ra~~d_n^{A_\ast}p_n=p_{n-1}d_n^{C(f)}.\nn
\end{align}
{\flushleft (2)} By definition, the connecting homomorphism $\delta_{n+1}$ is given (for $a_{n+1}\in\ker d_{n+1}^A$) by
{\footnotesize\[
\delta_{n+1}(a_{n+1}+\im d_{n+2}^A)=b_{n+1}+\im(-d_{n+2}^B),~~\txt{such that}~~i_n(b_{n+1})=d_{n+1}^{C(f)}(a_{n+1},b_{n+2}),~~~~p_{n+1}(a_{n+1},b_{n+2})=a_{n+1},\nn
\]}
where the last condition is clear, and ~$i_n(b_{n+1})=d_{n+1}^{C(f)}(a_{n+1},b_{n+2})$~ implies
{\small\begin{align}
&(0,b_{n+1})=-(-d_{n+1}^Aa_{n+1},f_{n+1}a_{n+1}+d_{n+2}^Bb_{n+2})=-(0,f_{n+1}a_{n+1}+d_{n+2}^Bb_{n+2}),\nn\\
&\Ra~~b_{n+1}=-(f_{n+1}a_{n+1}+d_{n+2}^Bb_{n+2}),\nn\\
&\Ra~~\delta_{n+1}(a_{n+1}+\im d_{n+2}^A)=b_{n+1}+\im(-d_{n+2}^B)=-f_{n+1}a_{n+1}+\im(-d_{n+2}^B)=-H_{n+1}(f)(a_{n+1}+\im d_{n+2}^A).\nn
\end{align}}
Also, we have
\[
\textstyle H_n\big(B_\ast[1]\big)={\ker d_n^{B_\ast[1]}\over\im d_{n+1}^{B_\ast[1]}}={\ker(-d_{n+1}^{B_\ast})\over\im(-d_{n+2}^{B_\ast})}={\ker d_{n+1}^{B_\ast}\over\im d_{n+2}^{B_\ast}}=H_{n+1}(B_\ast).
\]
{\flushleft (3)} From the exact sequence in part (2), $H_n(f)$, for all $n$, is an isomorphism if and only if $H_n(C(f))=0$ for all $n$. Hence $f$ is a quasi-iso if and only if $H_n(C(f))=0$ for all $n$, if and only if $C(f)$ is exact.
\end{proof}

\begin{dfn}[\textcolor{blue}{\index{Bounded! bicomplex}{Bounded bicomplex}}]
Let $\A$ be an abelian category. A bicomplex $C_{\ast\ast}\in\A_0^{\Integer\times\Integer}$,
\[\adjustbox{scale=0.8}{\bt
 ~       &  \ar[from=d]  ~                                                     &&  \ar[from=d]  ~                                               &&   \ar[from=d]  ~                                                    && ~\\
 ~\ar[r] & C_{i+1~j-1}\ar[rr,"\del^h_{i+1~j-1}"] && C_{i~j-1}\ar[rr,"\del^h_{i~j-1}"] && C_{i-1~j-1}\ar[rr,"\del^h_{i-1~j-1}"] && ~\\
 ~       &     ~                                                               &&     ~                                                     &&     ~                                                           && ~\\
 ~\ar[r] & C_{i+1~j}\ar[uu,"\del^v_{i+1~j}"']\ar[rr,"{\del^h_{i+1~j}}"] && C_{i~j}\ar[rr,"\del^h_{i~j}"]\ar[uu,"\del^v_{i~j}"'] && C_{i-1~j}\ar[uu,"{\del^v_{i-1~j}}"']\ar[rr,"\del^h_{i-1~j}"] && ~\\
 ~       &     ~                                                               &&     ~                                                     &&     ~                                                           && ~\\
 ~\ar[r] & C_{i+1~j+1}\ar[uu,"{\del^v_{i+1~j+1}}"']\ar[rr,"\del^h_{i+1~j+1}"] && C_{i~j+1}\ar[uu,"\del^v_{i~j+1}"']\ar[rr,"\del^h_{i~j+1}"]   && C_{i-1~j+1}\ar[uu,"\del^v_{i-1~j+1}"']\ar[rr,"\del^h_{i-1~j+1}"]    && ~\\
 ~       &   \ar[u]  ~                           &&    \ar[u] ~                         &&   \ar[u]  ~                              && ~
\et}\]
is \ul{bounded (above)} if there exists an integer $n$ such that $C_{ij}=0$ for all $i,j\leq n$.
\end{dfn}

\begin{lmm}[\textcolor{blue}{\index{Tot exactness lemma}{Tot exactness lemma}}]\label{TotExactLmm}
Let $C_{\ast\ast}$ be a bicomplex of $R$-modules. If $C_{\ast\ast}$ is (i) bounded and (ii) has exact rows (or exact columns), then the chain complex $\Tot(C_{\ast\ast})$ is exact.
\end{lmm}
\begin{proof}
Recall that for all $(c_{ij})_{i+j=n}\in \Tot(C_{\ast\ast})_n$, we have
\bea
d_n^{\Tot}(c_{ij})_{i+j=n}:=\Big(d^h_{i+1,j}c_{i+1,j}+(-1)^id^v_{i,j+1}c_{i,j+1}\Big)_{i+j=n}=\Big(d^hc_{i+1,n-i}+(-1)^id^vc_{i,n-i+1}\Big)_i.\nn
\eea
Thus, if $(c_{ij})_{i+j=n}$ is a cycle, i.e., $(c_{ij})_{i+j=n}\in\ker d_n^{\Tot}$, then $d_n^{\Tot}(c_{ij})_{i+j=n}=0$ implies
\bea
\label{CycleSyst}d^hc_{i+1,n-i}+(-1)^id^vc_{i,n-i+1}=0,~~~~\txt{for all}~~~~i.
\eea
Since $C_{\ast\ast}$ is bounded, there exists $i_0$ such that $c_{i,n-i+1}=0$ for all $i\leq i_0$. In particular, $c_{i_0,n-i_0+1}=0$, i.e., $i=i_0$ in (\ref{CycleSyst}) $\Ra$ $d^hc_{i_0+1,n-i_0}=0$. Since the rows are exact, there exists $z_{i_0+2,n-i_0}$ such that
\bea
c_{i_0+1,n-i_0}=d^hz_{i_0+2,n-i_0}.\nn
\eea
With $i=i_0+1$ in (\ref{CycleSyst}), we also have
\bea
&&0=d^hc_{i_0+2,n-i_0-1}+(-1)^{i_0+1}d^vc_{i_0+1,n-i_0}=d^hc_{i_0+2,n-i_0-1}+(-1)^{i_0+1}d^vd^hz_{i_0+2,n-i_0}\nn\\
&&~~~~=d^hc_{i_0+2,n-i_0-1}+(-1)^{i_0+1}d^hd^vz_{i_0+2,n-i_0}=d^h\big(c_{i_0+2,n-i_0-1}+(-1)^{i_0+1}d^vz_{i_0+2,n-i_0}\big).\nn
\eea
Again, because the rows are exact, there exists $z_{i_0+3,n-i_0-1}$ such that
\bea
&&c_{i_0+2,n-i_0-1}+(-1)^{i_0+1}d^vz_{i_0+2,n-i_0}=d^hz_{i_0+3,n-i_0-1},\nn\\
&&~~\Ra~~c_{i_0+2,n-i_0-1}=d^hz_{i_0+3,n-i_0-1}-(-1)^{i_0+1}d^vz_{i_0+2,n-i_0}.\nn
\eea
Continuing by induction, we obtain $(z_{ij})_{i+j=n+1}\in \Tot(C_{\ast\ast})_{n+1}$ such that
\bea
&&\textstyle d_{n+1}^{\Tot}(z_{ij})_{i+j=n+1}= (c_{ij})_{i+j=n}\in\im d_{n+1}^{\Tot},\nn\\
&&\textstyle ~~\Ra~~\ker d_n^{\Tot}=\im d_{n+1}^{\Tot},~~\Ra~~H_n(\Tot(C_{\ast\ast}))={\ker d_n^{\Tot}\over\im d_{n+1}^{\Tot}}=0,\nn
\eea
i.e., $\Tot(C_{\ast\ast})$ is an exact complex.
\end{proof}

\begin{lmm}[\textcolor{blue}{\index{Tot quasiisomorphism lemma}{Tot quasiisomorphism lemma}}]\label{TotQuisLmm}
Let $f_{\ast\ast}:A_{\ast\ast}\ra B_{\ast\ast}$ be a morphism of \ul{bounded} bicomplexes of $R$-modules. If the morphisms of horizontal/row chain complexes $f_{\ast j}:A_{\ast j}\ra B_{\ast j}$ (or equivalently,  morphisms of vertical/column chain complexes $f_{i\ast}:A_{i\ast}\ra B_{i\ast}$) are all quasiisomorphisms, then $\Tot(f_{\ast\ast}):\Tot(A_{\ast\ast})\ra \Tot(B_{\ast\ast})$ is a quiasiisomorphism.
\end{lmm}
\begin{proof}
Let $C(f_{\ast j})=C\Big(A_{\ast j}\sr{f_{\ast j}}{\ral}B_{\ast j}\Big)=B_{\ast j}[1]\oplus A_{\ast j}$ be the cone of the $j$th horizontal chain morphism. From Proposition \ref{ConeLESprp}, we have SES's of complexes $0\ra B_{\ast j}[1]\ra C(f_{\ast j})\ra A_{\ast j}\ra 0$ (by which $C(f_{\ast j})$ is exact for each $j$ since each $f_{\ast j}$ is a quasiiso), which together give an SES of bicomplexes
\bea
0\ra B_{\ast\ast}[1]\ra C_h(f_{\ast\ast})\ra A_{\ast\ast}\ra 0~~~~(\txt{$h$ denoting ``horizontal''}),\nn
\eea
where $C_h(f_{\ast\ast})$ is a bounded bicomplex with exact rows $\{C(f_{\ast j})\}_j$.

Taking the Tot of the above sequence, we obtain the following exact sequence (\blue{footnote}\footnote{Recall that Tot is exact on $R$-mod.}):
\bc\bt
0\ar[r]& \Tot(B_{\ast\ast}[1])\ar[d,equal,"{(i)}"]\ar[r]& \Tot\big(C_h(f_{\ast\ast})\big)\ar[d,equal,"{(ii)}"]\ar[r]& \Tot\big(A_{\ast\ast}\big)\ar[r]&0\\
       & \Tot(B_{\ast\ast})[1]  & C\big(\Tot(f_{\ast\ast})\big)  &   &
\et\ec
where (i) follows from translation invariance of Tot, and (ii) follows because Tot is an additive functor (\blue{footnote}\footnote{Recall that additive functors commute with, or preserve, finite direct sums.}), since
\bea
&&\textstyle \Tot\big(C_h(f_{\ast\ast})\big)_n=\bigoplus_{i+j=n}C(f_{\ast j})_i=\bigoplus_{i+j=n}\left(B_{i j}[1]\oplus A_{ij}\right)=\left(\bigoplus_{i+j=n}B_{i j}[1]\right)~\oplus~\left(\bigoplus_{i+j=n}A_{ij}\right)\nn\\
&&~~~~=C(\Tot(f_{\ast\ast})).\nn
\eea
That is, we have the SES
\bc\bt
0\ar[r]& \Tot(B_{\ast\ast})[1]\ar[r]& C\big(\Tot(f_{\ast\ast})\big)\ar[r]& \Tot\big(A_{\ast\ast}\big)\ar[r]&0.\nn
\et\ec
By Lemma \ref{TotExactLmm}, $C\big(\Tot(f_{\ast\ast})\big)$ is exact (since all rows of the bounded bicomplex $C_h(f_{\ast\ast})$ are exact). Hence $\Tot(f_{\ast\ast})$ is a quasiiso by part 3 of Proposition \ref{ConeLESprp}.
\end{proof}

\section{Derived functor isomorphisms (or balancing) and applications}
As before, let $\A$ be an abelian category with a tensor product $\otimes=\otimes_\A:\A\times\A\ra\A$ and a mor ``product''/operator $\Mor=\Mor_\A:\A\times\A\ra Ab$. For concrete cases, we will generally have
\[
\otimes,Hom:\txt{Modules}\times\txt{Modules}\ra\txt{Modules}
\]
 along with the following specific cases (and canonical/natural definitions for ~Tor, tor, Ext, ext~ functors):
\[
\otimes_R:(\txt{mod-}R)\times(R\txt{-mod})\ra Ab~~~~\txt{and}~~~~Hom_R(-,-):R\txt{-mod}\times R\txt{-mod}\ra Ab.
\]

\begin{example}[\textcolor{blue}{\index{Tensor! product of chain complexes}{Tensor product of chain complexes}, \index{Total! torsion}{Total torsion}}]
 Given two chain complexes $A_\ast=\big(\cdots\ra A_1\ra A_0\ra 0\big)$ and $B_\ast=\big(\cdots\ra B_1\ra B_0\ra 0\big)$ in $\A$, we define their \ul{tensor product} $A_\ast\otimes B_\ast$ to be the bicomplex $C_{\ast\ast}=(C_{ij}):=(A_i\otimes B_j)$, i.e.,
\[\adjustbox{scale=0.9}{\bt
        & \vdots\ar[d] & \vdots\ar[d] & \vdots\ar[d] &\vdots\ar[d] &\\
 0 & A_0\otimes B_3 \ar[l]\ar[d,"d_{03}^v"] & A_1\otimes B_3\ar[l,"d_{13}^h"']\ar[d,"d_{13}^v"] & A_2\otimes B_3\ar[l,"d_{23}^h"']\ar[d,"d_{23}^v"] & A_3\otimes B_3\ar[l,"d_{33}^h"']\ar[d,"d_{33}^v"] & \cdots\ar[l]\\
 0 & A_0\otimes B_2 \ar[l]\ar[d,"d_{02}^v"] & A_1\otimes B_2\ar[l,"d_{12}^h"']\ar[d,"d_{12}^v"] & A_2\otimes B_2\ar[l,"d_{22}^h"']\ar[d,"d_{22}^v"] & A_3\otimes B_2\ar[l,"d_{32}^h"']\ar[d,"d_{32}^v"] & \cdots\ar[l]\\
 0 & A_0\otimes B_1\ar[l]\ar[d,"d_{01}^v"] & A_1\otimes B_1\ar[l,"d_{11}^h"']\ar[d,"d_{11}^v"] & A_2\otimes B_1\ar[l,"d_{21}^h"']\ar[d,"d_{21}^v"] & A_3\otimes B_1\ar[l,"d_{31}^h"']\ar[d,"d_{31}^v"] & \cdots\ar[l] \\
 0 & A_0\otimes B_0\ar[l]\ar[d] & A_1\otimes B_0\ar[l,"d_{10}^h"']\ar[d] & A_2\otimes B_0\ar[l,"d_{20}^h"']\ar[d] & A_3\otimes B_0\ar[l,"d_{30}^h"']\ar[d] & \cdots\ar[l] \\
        & 0 & 0 & 0 & 0 &
\et}\]
where ~$d_{ij}^h:=d_i^{A_\ast}\otimes id_{B_j}$ ~and $d_{ij}^v:=id_{A_i}\otimes d_j^{B_\ast}$.

Let $A,B\in\Ob\A$, and consider projective resolutions $P_\ast^A\sr{\vep}{\ral} A\ra 0$ and $P_\ast^B\sr{\eta}{\ral} B\ra 0$, where
\bea
P^A_\ast:~\cdots\ra P^A_2\sr{d^A_2}{\ral} P^A_1\sr{d_1^A}{\ral}P^A_0\ra 0,~~~~P^B_\ast:~\cdots\ra P^B_2\sr{d^B_2}{\ral} P^B_1\sr{d_1^B}{\ral}P^B_0\ra 0.\nn
\eea
The \ul{total torsion} associated with $A,B$ is ~$Tor^{\otimes}(A,B):=H\Big(\Tot(P^A_\ast\otimes P^B_\ast)\Big)$,~ i.e.,
\bea
\textstyle Tor_n^{\otimes}(A,B):=H_n\Big(\Tot(P^A_\ast\otimes P^B_\ast)\Big),~~~~\big(\Tot(P^A_\ast\otimes P^B_\ast)\big)_n:=\bigoplus\limits_{i+j=n}P^A_i\otimes P^B_j,\nn
\eea
with differential ~{\small $d_n^{(A,B)}:=d_n^{\Tot}:=\sum_{i+j=n}\left(d^h_{ij}+(-1)^id^v_{ij}\right)=\sum_{i+j=n}\left(d_i^A\otimes id_{B_j}+(-1)^iid_{A_i}\otimes d_j^B\right)$}.
\end{example}

Although we have used bounded chain complexes to define $A_\ast\otimes B_\ast$ (reflecting intended application to projective/injective resolutions) the definition also applies to any pair of chain complexes $A_\ast,B_\ast$.

\begin{dfn}[\textcolor{blue}{\index{Hom product of chain complexes}{Hom product of chain complexes}, \index{Total! extension}{Total extension}}]
Given two complexes $A_\ast=\big(\cdots\ra A_1\ra A_0\ra 0\big)$ and $B^\ast=\big(0\ra B^0\ra B^1\ra \cdots\big)$ in $\A$,
we define their \ul{mor product} $\Mor(A_\ast,B^\ast)$ to be the bicomplex $C^{\ast\ast}=(C^{ij}):=\big(\Mor(P_i,Q^j)\big)$, i.e.,
\[\adjustbox{scale=0.9}{\bt
        & \vdots & \vdots & \vdots &\vdots &\\
 0\ar[r] & \Mor(A_0,B^3) \ar[r,"d_h^{03}"]\ar[u,"d_v^{03}"'] & \Mor(A_1,B^3)\ar[r,"d_h^{13}"]\ar[u,"d_v^{13}"'] & \Mor(A_2,B^3)\ar[r,"d_h^{23}"]\ar[u,"d_v^{23}"'] & \Mor(A_3,B^3)\ar[r,"d_h^{33}"]\ar[u,"d_v^{33}"'] & \cdots\\
 0\ar[r] & \Mor(A_0,B^2) \ar[r,"d_h^{02}"]\ar[u,"d_v^{02}"'] & \Mor(A_1,B^2)\ar[r,"d_h^{12}"]\ar[u,"d_v^{12}"'] & \Mor(A_2,B^2)\ar[r,"d_h^{22}"]\ar[u,"d_v^{22}"'] & \Mor(A_3, B^2)\ar[r,"d_h^{32}"]\ar[u,"d_v^{32}"'] & \cdots\\
 0\ar[r] & \Mor(A_0,B^1)\ar[r,"d_h^{01}"]\ar[u,"d_v^{01}"'] & \Mor(A_1,B^1)\ar[r,"d_h^{11}"]\ar[u,"d_v^{11}"'] & \Mor(A_2,B^1)\ar[r,"d_h^{21}"]\ar[u,"d_v^{21}"'] & \Mor(A_3, B^1)\ar[r,"d_h^{31}"]\ar[u,"d_v^{31}"'] & \cdots \\
 0\ar[r] & \Mor(A_0,B^0)\ar[r,"d_h^{00}"]\ar[u,"d_v^{00}"'] & \Mor(A_1,B^0)\ar[r,"d_h^{10}"]\ar[u,"d_v^{10}"'] & \Mor(A_2,B^0)\ar[r,"d_h^{20}"]\ar[u,"d_v^{20}"'] & \Mor(A_3, B^0)\ar[r,"d_h^{30}"]\ar[u,"d_v^{30}"'] & \cdots \\
        & 0\ar[u] & 0\ar[u] & 0\ar[u] & 0\ar[u] &\\
\et}\]
where  $d_h^{ij}:=\Mor(d^{A_\ast}_i,B^j)=\big(A_i\sr{(-)}{\ral} B^j\big)\circ d^{A_\ast}_i$ and $d_v^{ij}:=\Mor(A_i,d_{B^\ast}^j)=d_{B^\ast}^j\circ\big(A_i\sr{(-)}{\ral} B^j\big)$.

Let $A,B\in\Ob\A$. Consider a projective resolution $P_\ast^A\sr{\vep}{\ral} A\ra 0$ and an injective resolution $0\ra B\sr{\xi}{\ral} I^\ast_B$, where
\bea
P^A_\ast:~\cdots\ra P^A_2\sr{d^A_2}{\ral} P^A_1\sr{d_1^A}{\ral}P^A_0\ra 0,~~~~I^\ast_B:~ 0\ra I^0_B \sr{d^0_B}{\ral} I^1_B \sr{d^1_B}{\ral} I^2_B\ra \cdots.\nn
\eea
The \ul{total extension} associated with $A,B$ is ~$Ext_{\Mor}(A,B):=H\Big(\Tot\big(\Mor(P^A_\ast,I_B^\ast)\big)\Big)$, i.e.,
\bea
\textstyle Ext_{\Mor}^n(A,B):=H_n\Big(\Tot\big(\Mor(P^A_\ast,I_B^\ast)\big)\Big),~~~~\Big(\Tot\big(\Mor(P^A_\ast,I_B^\ast)\big)\Big)_n:=\bigoplus\limits_{i+j=n}\Mor(P^A_i,I_B^j),\nn
\eea
with differential ~{\small $d_{(A,B)}^n:=d_n^{\Tot}:=\sum_{i+j=n}\left(d^{ij}_h+(-1)^id^{ij}_v\right)=\sum_{i+j=n}\left(\Mor(d^A_i,B^j)+(-1)^i\Mor(A_i,d_B^j)\right)$}.
\end{dfn}

Although we have used bounded chain complexes to define $\Mor(A_\ast,B^\ast)$ (reflecting intended application to projective/injective resolutions) the definition also applies to any pair of chain complexes $A_\ast,B_\ast$.

\begin{thm}[\textcolor{blue}{\index{Balance of Tor}{Balance of Tor}: $\txt{Tor}_R\cong\txt{tor}_R$}]
Let $A_R\in\txt{mod-}R$ and $_RB\in R\txt{-mod}$. For any projective resolutions $P^A_\ast\mathop{\ral}\limits^\vep A\ra 0$ and $P^B_\ast\mathop{\ral}\limits^\eta B\ra 0$,
\bea
Tor_n^R(A,B)~:=~H_n(P^A_\ast\otimes_R B)\cong \ub{H_n\Big(\Tot(P^A_\ast\otimes_R P^B_\ast)\Big)}_{Tor^{\otimes_R}_n(A,B)}\cong H_n(A\otimes_R P^B_\ast)~=:~tor_n^R(A,B).\nn
\eea
\end{thm}
\begin{proof}
 Observe that with $B_\ast:=\big(\cdots\ra 0\ra 0\ra B\ra 0\big)$ and $A_\ast:=\big(\cdots\ra 0\ra 0\ra A\ra 0\big)$, we have
\bea
P^A_\ast\otimes_R B=\Tot(P^A_\ast\otimes_R B_\ast),~~~~A\otimes_R P^B_\ast=\Tot(A_\ast\otimes_R P^B_\ast).\nn
\eea
Also, we have the following quasi-isomorphisms:
\bea
\vep_\ast=(\cdots,0,0,\vep): P^A_\ast\ra A_\ast,~~~~\eta_\ast=(\cdots,0,0,\eta): P^B_\ast\ra B_\ast.\nn
\eea
Thus, we get morphisms of bicomplexes {\footnotesize $\bt P^A_\ast\otimes B_\ast\ar[from=rr,"id_{P^A_\ast}\otimes\eta_\ast"'] && P^A_\ast\otimes P^B_\ast && A_\ast\otimes P^B_\ast\ar[from=ll,"\vep_\ast\otimes id_{P^B_\ast}"]\et$}. Taking Tot,
{\footnotesize\bea
\bt P^A_\ast\otimes B=T(P^A_\ast\otimes B_\ast)\ar[from=rr,"\Tot(id_{P^A_\ast}\otimes\eta_\ast)"'] && \Tot(P^A_\ast\otimes P^B_\ast)&& \Tot(A_\ast\otimes P^B_\ast)=A\otimes P^B_\ast\ar[from=ll,"\Tot(\vep_\ast\otimes id_{P^B_\ast})"]\et .\nn
\eea}
We want to show $\Tot(id_{P^A_\ast}\otimes\eta_\ast)$ is a quasi-isomorphism. By Lemm \ref{TotQuisLmm}, it is enough to show that for each $i$, the chain morphism $id_{P^A_i}\otimes\eta_\ast$ on the $i$th column $P^A_i\otimes P^B_\ast$ of $P^A_\ast\otimes P^B_\ast$ is a quasi-isomorphism.

{\footnotesize
\bea
\left(\adjustbox{scale=0.8}{%
\bt
\vdots\ar[d]  & \vdots\ar[d] & \\
 0\ar[d] & 0\ar[l]\ar[d] & \ar[l]\cdots \\
 P^A_0\otimes B &  P^A_1\otimes B\ar[l] & \ar[l]\cdots
\et}\right)
\sr{id_{P^A_\ast}\otimes\eta_\ast}{\lal}
\left(\adjustbox{scale=0.8}{%
\bt
\vdots\ar[d]  & \vdots\ar[d] & \\
 P^A_0\otimes P^B_1\ar[d] & P^A_1\otimes P^B_1\ar[l]\ar[d] & \ar[l]\cdots \\
 P^A_0\otimes P^B_0 & P^A_1\otimes P^B_0\ar[l] & \ar[l]\cdots
\et}\right)
\sr{\vep_\ast\otimes id_{P^B_\ast}}{\ral}
\left(\adjustbox{scale=0.8}{%
\bt
\vdots\ar[d]  & \vdots\ar[d] & \\
 A\otimes P^B_1\ar[d] & 0\ar[l]\ar[d] & \ar[l]\cdots \\
 A\otimes P^B_0 & 0\ar[l] & \ar[l]\cdots
\et}\right)\nn
\eea}
Indeed, since each $P^A_i$ is projective (hence flat), we know the functor $P^A_i\otimes -$ is exact, where
\bea
P^A_i\otimes -:~~P_\ast^B\sr{\eta_\ast}{\ral}B_\ast~~\longmapsto~~P_i^A\otimes P_\ast^B\sr{id_{P_i^A}\otimes\eta_\ast}{\ral}P_i^A\otimes B_\ast,\nn
\eea
and so preserves quasi-isomorphisms (\blue{footnote}\footnote{
If a functor $F$ is exact then it preserves exactness (of all short exact sequences that make up a complex), and hence preserves (commutes with) the homology of a complex. Also, by Proposition \ref{ConeLESprp}, a module chain morphism $f:A_\ast\ra B_\ast$ is a quasi-iso $\iff$ its mapping cone $C(f)$ is exact. Therefore, the mapping cone of $F(f)$ is also exact. Hence, if $f$ is a quasi-iso, then so is $F(f)$, since ~$H_n(C(f))=0$ ~$\Ra$ ~$0=FH_n(C(f))=H_nF(C(f))=H_n(C(F(f)))$.
}).

By a similar argument (with rows replacing columns), we also conclude that $\Tot(\vep_\ast\otimes id_{P^B_\ast})$ is a quasi-isomorphism.
\end{proof}

\begin{thm}[\textcolor{blue}{\index{Balance of Ext}{Balance of Ext}: $\txt{Ext}_R\cong\txt{ext}_R$}]
Let $A,B$ be $R$-modules. For any projective resolution $P^A_\ast\mathop{\ral}\limits^\vep A\ra 0$ and any injective resolution  $0\ra B\mathop{\ral}\limits^\gamma I_B^\ast$,
\bea
Ext_R^n(A,B)~:=~H_n\Big(Hom_R(A,I_B^\ast)\Big)\cong \ub{H_n\Big(\Tot\big(Hom_R(P^A_\ast,I_B^\ast)\big)\Big)}_{Ext_{Hom_R}^n(A,B)}\cong H_n\Big(Hom_R(P^A_\ast,B)\Big)~=:~ext_R^n(A,B)\nn
\eea
\end{thm}
\begin{proof}
Observe that with $A_\ast=\big(\cdots\ra 0\ra 0\ra A\ra 0\big)$ and $B^\ast=\big(0\ra B\ra 0\ra 0\ra\cdots\big)$, we have
\bea
Hom_R(A,I_B^\ast)=\Tot(Hom_R(A_\ast,I_B^\ast)),~~~~Hom_R(P^A_\ast,B)=\Tot(Hom_R(P^A_\ast,B^\ast)).\nn
\eea
Also, we have the quasi-isomorphisms
\bea
\vep_\ast=(\cdots,0,0,\vep): P^A_\ast\ra A_\ast,~~~~\gamma^\ast=(\gamma,0,0,\cdots): I_B^\ast\ra B^\ast,\nn
\eea
and morphisms of bicomplexes
\bea
Hom_R(A_\ast,I_B^\ast)\sr{\vep^{\ast\ast}}{\ral}Hom_R(P^A_\ast,I_B^\ast)\sr{\gamma^{\ast\ast}}{\lal}Hom_R(P^A_\ast,B^\ast),\nn
\eea
where ~$\vep^{\ast\ast}:=Hom(\vep_\ast,I_B^\ast)$~ and ~$\gamma^{\ast\ast}:=Hom(P^A_\ast,\gamma^\ast)$. Taking Tot, we get
{\small\bea
Hom(A,I_B^\ast)=\Tot(Hom(A_\ast,I_B^\ast))\sr{\Tot(\vep^{\ast\ast})}{\ral}\Tot(Hom(P^A_\ast,I_B^\ast))\sr{\Tot(\gamma^{\ast\ast})}{\lal}\Tot(Hom(P^A_\ast,B^\ast))=Hom(P^A_\ast,B).\nn
\eea}

We want to show $\Tot(\vep^{\ast\ast})$, $\Tot(\gamma^{\ast\ast})$ are quasi-isos.
{\footnotesize
\bea
\left(\adjustbox{scale=0.8}{%
\bt[column sep =small]
\vdots  & \vdots & \\
 Hom(A,I_B^1)\ar[u]\ar[r] & 0\ar[r]\ar[u] & \cdots \\
 Hom(A,I_B^0)\ar[u]\ar[r] & 0\ar[r]\ar[u] & \cdots
\et}\right)
\sr{\vep^{\ast\ast}}{\ral}
\left(\adjustbox{scale=0.8}{%
\bt[column sep =small]
\vdots  & \vdots & \\
 Hom(P^A_0,I_B^1)\ar[u]\ar[r] & Hom(P^A_1,I_B^1)\ar[r]\ar[u] & \cdots \\
 Hom(P^A_0,I_B^0)\ar[u]\ar[r] & Hom(P^A_1,I_B^0)\ar[r]\ar[u] & \cdots
\et}\right)
\sr{\gamma^{\ast\ast}}{\lal}
\left(\adjustbox{scale=0.8}{%
\bt[column sep =small]
\vdots  & \vdots & \\
 0\ar[u]\ar[r] & 0\ar[r]\ar[u] & \cdots \\
 Hom(P^A_0,B)\ar[u]\ar[r] & Hom(P^A_1,B)\ar[r]\ar[u] & \cdots
\et}\right)\nn
\eea}
By Lemma \ref{TotQuisLmm}, it is enough to show that each of the two morphisms $\vep^{\ast\ast}$, $\gamma^{\ast\ast}$ restricts to quasi-isomorphisms on rows (or on columns). This is indeed the case due to the following:
\bit[leftmargin=0.7cm]
\item[(1)] $I_B^j$ injective ~$\Ra$~ $Hom(-,I_B^j)$ exact, ~$\Ra$~ preserves homology and hence quasi-iso also: Recall that
\bea
H_n(C(f))=0~~~~\Ra~~~~0=FH_n(C(f))=H_nF(C(f))=H_n(C(F(f))).\nn
\eea
So, since $P^A_\ast\sr{\vep_\ast}{\ral}A_\ast$ is a quasi-iso, $Hom(P^A_\ast,I_B^j)\sr{\vep^{\ast j}}{\lal}Hom(A,I_B^j)$ is also a quasi-iso (i.e., $\vep^{\ast\ast}$ restricts to a quasi-iso on each row).
\item[(2)] Similarly, $\gamma^{\ast\ast}$ restricts to a quasi-iso on each column, since $B^\ast\sr{\gamma^\ast}{\ral}I_B^\ast$ is a quasi-iso and the exactness of the functor $Hom(P^A_i,-)$ (as $P^A_i$ is projective) gives a quasi-iso ~$Hom(P^A_i,B)\sr{\gamma^{i\ast}}{\ral}Hom(P^A_i,I_B^\ast)$.\qedhere
\eit
\end{proof}

\begin{rmk}[\textcolor{blue}{Alternative derivation of Tor and Ext sequences for $R$-modules}]
Given an object $X\in\Ob\A$ and an exact sequence $\txt{SES}:~0\ra A\ra B\ra C\ra0$ in $\A$, we obtained the tor sequence $tor_n^\A(X,\txt{SES}):=H_n(X\otimes_\A P_\ast^{\txt{SES}})$ by (i) first resolving SES vertically (via horshoe lemmas or otherwise) to obtain a split-exact sequence of deleted projective resolutions $P^{\txt{SES}}:~0\ra P_\ast^A\ra P_\ast^B\ra P_\ast^C\ra0$, (ii) applying $F:=X\otimes_\A-:\A\ra\A$, and then (iii) taking the induced LES of homologies. However (with $\otimes_R:(\txt{mod-}R)\times(R\txt{-mod})\ra Ab$, and modules $X_R$, ${}_RA$, ${}_RB$, ${}_RC$) we now have
\bea
tor_n^R(X,\txt{SES})~:=~H_n(X\otimes_R P_\ast^{\txt{SES}})\cong H_n(P_\ast^X\otimes_R\txt{SES})~=:~Tor_n^R(X,\txt{SES}),\nn
\eea
and so we can obtain the Tor sequences by instead (i) first resolving $X$ into $P_\ast^A$, (ii) forming the bicomplex $P_\ast^A\otimes_R\txt{SES}$ (which is automatically an exact sequence of chain complexes), and then (iii) taking the induced LES of homologies.

A similar argument clearly holds for Ext sequences $Ext_\A^n(X,\txt{SES}):=H_n\big(Hom_\A(X,I^\ast_{\txt{SES}})\big)$, because if $\A$ = $R$-mod then we have
\bea
Ext_R^n(X,\txt{SES})~:=~H_n\big(Hom_R(X,I^\ast_{\txt{SES}})\big)\cong H_n\big(Hom_R(P_\ast^X,\txt{SES})\big)~=:~ext_R^n(X,\txt{SES}),\nn
\eea
where $I_{\txt{SES}}^\ast:0\ra I^\ast_A\ra I^\ast_B\ra I^\ast_C\ra0$ is a split-exact sequence of deleted injective resolutions and the bicomplex $Hom_R(P_\ast^X,\txt{SES})$ is automatically an exact sequence of chain complexes in $Ab$.
\end{rmk}

\begin{crl}
Let $R$ be a commutative ring and $r\in R$. Let $A,B$ be $R$-modules, and let $A\sr{r}{\ral}A$ be multiplication by $r$. Then the following maps induced by $r$ are each also multiplication by $r$.
\bea
&&r\otimes_R id_B:A\otimes_R B\ra A\otimes_R B,~~~~Hom_R(r,B):Hom_R(A,B)\ra Hom_R(A,B),\nn\\
&&~~~~Hom_R(B,r):Hom_R(B,A)\ra Hom_R(B,A),~~~~~~~~\txt{(so far easy to check directly)}\nn\\
&&Tor_n^R(r,B):Tor_n^R(A,B)\ra Tor_n^R(A,B),~~~~Ext_R^n(r,B):Ext_R^n(A,B)\ra Ext_R^n(A,B),\nn\\
&&~~~~Ext_R^n(B,r):Ext_R^n(B,A)\ra Ext_R^n(B,A)~~~~~~\txt{(easy to see using the iso's Tor $\cong$ tor, Ext $\cong$ ext)}.\nn
\eea
\end{crl}
\begin{proof}
(i) Multiplication by $r$ is an $r$-homomorphism since $R$ is commutative. (ii) Multiplication by $r$ commutes with all $R$-homomorphisms (which include multiplications by all elements of $R$, since $R$ is commutative). (iii) With respect to composition of morphisms (including elements of $R$ as multiplication maps), the functors $Hom_R(-,-)$, $-\otimes_R-$ are $R$-bilinear, and likewise the homology functor $H_n(-)$ is $R$-linear in the sense $H_n(rf)=H_n(fr)=rH_n(f)=H_n(f)r$ for any morphism of chain complexes $f$.

Finally, the isomorphisms Tor $\cong$ tor and Ext $\cong$ ext show that (up to isomorphism) the resolved versions of $r:A\ra A$, namely, $P^r:P_\ast^A\ra P_\ast^A$ and $I_r:I^\ast_A\ra I^\ast_A$, are each equivalent to multiplication by $r$.
\end{proof}

\begin{crl}[\textcolor{blue}{Tor$^R$-criterion for Flatness}]
Let $F$ be a right $R$-module. The following are equivalent.
\begin{enumerate}[leftmargin=0.9cm]
\item[(1)] $F$ is flat ~(i.e., $F\otimes_R-$ is exact).
\item[(2)] $Tor_n^R(F,-)=0$, for all $n\geq1$.
\item[(3)] $Tor_1^R(F,-)=0$.
\end{enumerate}
\end{crl}
\begin{proof}
The implication (2)$\Ra$(3) is clear, and so we will prove (1)$\Ra$(2) and (3)$\Ra$(1).
{\flushleft\ul{(1)$\Ra$(2)}}:  $Tor_n^R(F,X)\cong tor_n^R(F,X)=H_n(F\otimes_R P^X_\ast)$. Since $F$ is flat, $F\otimes_R-$ is exact, and so $F\otimes_R P_\ast^X$ is exact except at index $0$ (just as $P_\ast^X$ is). Hence $H_n(F\otimes_R P_\ast^X)=0$ for $n\geq 1$.
{\flushleft\ul{(3)$\Ra$(1)}}: We show $F\otimes_R-$ is exact. For any SES $0\ra A\ra B\ra C\ra 0$, the Tor-LES gives
{\small\bea
\cdots\sr{\delta}{\ral} \ub{Tor_1^R(F,A)}_0\ra \ub{Tor_1^R(F,B)}_0\ra \ub{Tor_1^R(F,C)}_0\sr{\delta}{\ral}F\otimes_RA\ra F\otimes_RB\ra F\otimes_RC\ra 0.\nn\qedhere
\eea}
\end{proof}

\begin{prp}[\textcolor{blue}{Ext$_R$-criterion for Projectivity}]
Let $P$ be an $R$-module. The following are equivalent.
\begin{enumerate}[leftmargin=0.9cm]
\item[(1)] $P$ is projective (i.e., $Hom_R(P,-)$ exact).
\item[(2)] $Ext_R^n(P,-)=0$ for all $n\geq 1$.
\item[(3)] $Ext_R^1(P,-)=0$.
\end{enumerate}
\end{prp}
\begin{proof}
The implication (2)$\Ra$(3) is clear, and so we will prove (1)$\Ra$(2) and (3)$\Ra$(1).
{\flushleft\ul{(1)$\Ra$(2)}}: $Ext_R^n(P,X)=H^n(Hom_R(P,I^\ast_X))$. Since $Hom_R(P,-)$ is exact, $Hom_R(P,I^\ast_X)$ is exact except at index $0$ (just as $I^\ast_X$ is). Hence $H^n(Hom_R(P,I^\ast_X))=0$ for $n\geq 1$.
{\flushleft\ul{(3)$\Ra$(1)}}: We show $Hom_R(P,-)$ is exact. For any SES $0\ra A\ra B\ra C\ra 0$, the Ext-LES gives
{\footnotesize
\bea
0\ra Hom_R(P,A)\ra Hom_R(P,B)\ra Hom_R(P,C)\sr{\delta}{\ral}\ub{Ext_R^1(P,A)}_0\ra \ub{Ext_R^1(P,B)}_0\ra \ub{Ext_R^1(P,C)}_0\sr{\delta}{\ral}\cdots.\nn\qedhere
\eea}
\end{proof}

\begin{prp}[\textcolor{blue}{Ext$_R$-criterion for Injectivity}]
Let $I$ be an $R$-module. The following are equivalent.
\begin{enumerate}[leftmargin=0.9cm]
\item[(1)] $I$ is injective. (i.e., $Hom_R(-,I)$ is exact).
\item[(2)] $Ext_R^n(-,I)=0$ for all $n\geq 1$.
\item[(3)] $Ext_R^1(-,I)=0$.
\end{enumerate}
\end{prp}
\begin{proof}
The implication (2)$\Ra$(3) is clear, and so we will prove (1)$\Ra$(2) and (3)$\Ra$(1).
{\flushleft\ul{(1)$\Ra$(2)}}: $Ext_R^n(X,I)\cong ext_R^n(X,I)=H^n(Hom_R(P_\ast^X,I))$. Since $Hom_R(-,I)$ is exact, $Hom_R(P_\ast^X,I)$ is exact except at index $0$ (just as $P_\ast^X$ is). Hence $H^n(Hom_R(P_\ast^X,I))=0$ for $n\geq 1$.
{\flushleft\ul{(3)$\Ra$(1)}}: We show $Hom_R(-,I)$ is exact. For any SES $0\ra A\ra B\ra C\ra 0$, the Ext-LES gives
{\small
\bea
0\ra Hom_R(C,I)\ra Hom_R(B,I)\ra Hom_R(A,I)\sr{\delta}{\ral}\ub{ext_R^1(C,I)}_0\ra \ub{ext_R^1(B,I)}_0\ra \ub{ext_R^1(A,I)}_0\sr{\delta}{\ral}\cdots\nn\qedhere
\eea}
\end{proof}

\begin{thm}[\textcolor{blue}{Localization-Tor$^R$ symmetry, Localization-Ext$_R$ symmetry}]
Let $R$ be a commutative ring, $S\subset R$ a localizing set, and $A,B$ two $R$-modules,
\begin{enumerate}[leftmargin=1cm]
\item[(1)] $S^{-1}Tor_n^R(A,B)\cong Tor_n^{S^{-1}R}(S^{-1}A,S^{-1}B)$
\item[(2)] If $A$ is \index{Finitely presented module}{\ul{finitely presented}}, in the sense there exists an exact sequence $0\ra R^k\ra R^l\ra A\ra 0$, then
\bea
S^{-1}Ext_R^n(A,B)\cong Ext_{S^{-1}R}^n(S^{-1}A,S^{-1}B).\nn
\eea
This is similar to the same result for Hom, i.e., ~$S^{-1}Hom_R(A,B)\cong Hom_{S^{-1}R}(S^{-1}A,S^{-1}B)$.
\end{enumerate}
\end{thm}
\begin{proof}
Let $R':=S^{-1}R$, which is an R-algebra and a flat R-module. Then for (1), we have
\begin{align}
&S^{-1}Tor_n^R(A,B)\cong S^{-1}R\otimes_RTor_n^R(A,B)=S^{-1}R\otimes_RH_n\big(P^A_\ast\otimes_RB\big)\nn\\
&~~~~\sr{(s1)}{=}H_n\big(S^{-1}R\otimes_R(P^A_\ast\otimes_RB)\big)=H_n\big(R'\otimes_R(P^A_\ast\otimes_RB)\big)=H_n\big((R'\otimes_RP^A_\ast)\otimes_RB\big)\nn\\
&~~~~\sr{(s2)}{\cong} H_n\big(P_\ast^{S^{-1}A}\otimes_RB\big)\sr{(s3)}{\cong} H_n\big((P_\ast^{S^{-1}A}\otimes_{R'}R')\otimes_RB\big)=H_n\big(P_\ast^{S^{-1}A}\otimes_{R'}(R'\otimes_RB)\big)\nn\\
&~~~~\cong H_n\big(P_\ast^{S^{-1}A}\otimes_{R'}S^{-1}B\big)=Tor_n^{S^{-1}R}(S^{-1}A,S^{-1}B),\nn
\end{align}
where step (s1) holds because $S^{-1}R\otimes_R-$ commutes with homology (as an exact functor), and step (s2) holds because $S^{-1}R\otimes_R-$ maps a projective resolution in $R$-mod to a projective resolution in $S^{-1}R$-mod: Recall that $P$ is projective iff $F=P\oplus X$ for a free module $F$, and so
\bea
S^{-1}F=S^{-1}R\otimes(P\oplus X)=S^{-1}P\oplus S^{-1}X.\nn
\eea
Similarly for (2), we have
\begin{align}
&S^{-1}Ext_R^n(A,B)\cong S^{-1}R\otimes_RExt_R^n(A,B)=S^{-1}R\otimes_RH_n\big(Hom_R(A,I^\ast_B)\big)\nn\\
&~~~~\sr{(s1)}{=}H_n\big(S^{-1}R\otimes_RHom_R(A,I^\ast_B)\big)=H_n\big(R'\otimes_RHom_R(A,I^\ast_B)\big)\nn\\
&~~~~\sr{(s)}{\cong} H_n\big(Hom_{S^{-1}R}(S^{-1}A,I^\ast_{S^{-1}B})\big)=Ext_{S^{-1}R}^n(S^{-1}A,S^{-1}B),\nn
\end{align}
where step (s1) is as before and step (s) holds because from the exact sequence $0\ra R^k\ra R^l\ra A\ra 0$ we can obtain two left-exact sequences
\begin{align}
&R'\otimes_R Hom_R\Big(0\ra R^k\ra R^l\ra A\ra 0,I^\ast_B\Big)~~=~~R'\otimes_R\Big(0\ra Hom_R(A,I^\ast_B)\ra I^\ast_{B^l}\ra I^\ast_{B^k}\Big)\nn\\
&~~~~~=~~0\ra R'\otimes_RHom_R(A,I^\ast_B)\ra I^\ast_{(S^{-1}B)^l}\ra I^\ast_{(S^{-1}B)^k}~~~~\txt{and}~~~~\nn\\
&Hom_{S^{-1}R}\Big(R'\otimes_R\big(0\ra R^k\ra R^l\ra A\ra 0\big),I^\ast_{S^{-1}B}\Big)=Hom_{S^{-1}R}\Big(0\ra R'{}^k\ra R'{}^l\ra S^{-1}A\ra 0,I^\ast_{S^{-1}B}\Big)\nn\\
&~~~~=~~0\ra Hom_{S^{-1}R}(S^{-1}A,I^\ast_{S^{-1}B})\ra I^\ast_{(S^{-1}B)^l}\ra I^\ast_{(S^{-1}B)^k},\nn
\end{align}
which give a commutative diagram (with unique kernel-induced maps $\theta,\theta'$, hence the isomorphism):
\[\bt
0\ar[r] & S^{-1}R\otimes_RHom_R(A,I^\ast_B)\ar[from=d,shift left=2,dashed,"\theta'","\cong"']\ar[d,shift left=2.5,dashed,"\theta"]\ar[r,"f"] & I^\ast_{(S^{-1}B)^l}\ar[d,equal,"id"]\ar[r,"g"] & I^\ast_{(S^{-1}B)^k}\ar[d,equal,"id"] \\
0\ar[r] & Hom_{S^{-1}R}(S^{-1}A,I^\ast_{S^{-1}B})\ar[r,"f'"] & I^\ast_{(S^{-1}B)^l}\ar[r,"g'"] & I^\ast_{(S^{-1}B)^k} \\
\et \qedhere \]
\end{proof} 

%% file: parts/AlgebraCat/HomDim.tex
\section{Homological Dimension of Modules and Rings}
In this section, for convenience with counting, the index of the differential in an injective resolution will be the same as that of its codomain (instead of that of its domain as before). This is the convention that is consistent when an injective resolution is viewed as the reflection of a projective resolution.
\subsection{Definitions, Basic results, Projective dimension bounds}
\begin{dfn}[\blue{
\index{Projective-finite module}{Projective-finite module},
\index{Length of a finite projective resolution}{Length of a finite projective resolution},
\index{Projective! dimension of a module}{Projective dimension of a module},
\index{Injective-finite module}{Injective-finite module},
\index{Length of a finite injective resolution}{Length of a finite injective resolution},
\index{Injective! dimension of a module}{Injective dimension of a module}}]
Let $R$ be a ring and $_RM$ an $R$-module. Consider a projective resolution
\bea
M_\ast:\cdots{\ral}P_n\sr{d_n}{\ral}P_{n-1}\sr{d_{n-1}}{\ral}\cdots \sr{d_3}{\ral}P_2\sr{d_2}{\ral}P_1\sr{d_1}{\ral}P_0\sr{\vep}{\ral}M{\ral}0\nn
\eea
(resp. injective resolution {\small $M^\ast:0\ra M\sr{\xi}{\ral}I^0\sr{d^1}{\ral}I^1\sr{d^2}{\ral}\cdots\sr{d^{n-1}}{\ral}I^{n-1}\sr{d^n}{\ral}I^n\ral\cdots$}) and deleted resolution
\bea
PP_\ast^M:~\cdots{\ral}P_n\sr{d_n}{\ral}P_{n-1}\sr{d_{n-1}}{\ral}\cdots \sr{d_3}{\ral}P_2\sr{d_2}{\ral}P_1\sr{d_1}{\ral}P_0{\ral}0\nn
\eea
(resp. $II^\ast_M:0\ra I^0\sr{d^1}{\ral}I^1\sr{d^2}{\ral}\cdots\sr{d^{n-1}}{\ral}I^{n-1}\sr{d^n}{\ral}I^n\ral\cdots$). If there exists $n\geq 1$ such that $P_n':=\ker d_{n-1}$ is projective, making $M$ a \ul{projective-finite module} (defined to be a module with a finite projective resolution), then we know $PP_\ast^M$ is homtotopy equivalent, hence quasi-isomorphic, (in $(R\txt{-mod})_0^\Integer$) to the following finite deleted resolution $P_\ast^M$ of \ul{length} $l(P_\ast^M):=n$.
\bea
P_\ast^M:0{\ral}~P_n'~{\hookrightarrow}~P_{n-1}\sr{d_{n-1}}{\ral}P_{n-2}\ral\cdots{\ral}P_0\sr{\vep}{\ral}M{\ral}0.\nn
\eea
The \ul{projective dimension} $Pd(M)$ of $M$ is the length of the shortest projective resolution of $M$, i.e.,
\[
Pd(M)=Pd_R(M):=\left\{
               \begin{array}{ll}
                 \min\{l(P_\ast^M):P_\ast^M~\txt{a finite proj resolution of $M$}\}, & \txt{if $M$ is projective-finite} \\
                \infty, & \txt{otherwise}
               \end{array}
             \right\}.\nn
\]

Similarly, if there exists $n\geq 1$ such that {\small $I'^n:=\coker d^{n-1}$} is injective, making $M$ an \ul{injective-finite module} (defined to be a module with a finite injective resolution), then we know $II^\ast_M$ is homtotopy equivalent, hence quasi-isomorphic, to the following finite deleted resolution $I^\ast_M$ of \ul{length} $l(I^\ast_M):=n$.
\bea
I^\ast_M:0\ra I^0\sr{d^1}{\ral}I^1\sr{d^2}{\ral}\cdots\sr{d^{n-1}}{\ral} I^{n-1}\sr{d^{n}}{\twoheadrightarrow}I'^n\ra0.\nn
\eea
The \ul{injective dimension} $Id(M)$ of $M$ is the length of the shortest injective resolution of $M$, i.e.,
\[
Id(M)=Id_R(M):=\left\{
               \begin{array}{ll}
                 \min\{l(I^\ast_M):I^\ast_M~\txt{a finite inj resolution of $M$}\}, & \txt{if $M$ is injective-finite} \\
                \infty, & \txt{otherwise}
               \end{array}
             \right\}.\nn
\]
\end{dfn}

\begin{dfn}[\blue{
\index{Left! projective dimension}{Left projective dimension},
\index{Left! injective dimension}{Left injective dimension}}]
Let $R$ be a ring. The \ul{left projective (homological) dimension} $Lpd(R)$ of $R$ is the largest projective dimension of an $R$-module, i.e.,
\bea
\textstyle Lpd(R):=\sup\{Pd(M):M\in\txt{R-mod}\}=\sup_{M\in\txt{R-mod}}\inf_{P_\ast^M}l(P_\ast^M).\nn
\eea
The \ul{left injective (homological) dimension} $Lid(R)$ of $R$ is the largest injective dimension of an $R$-module, i.e.,
\bea
\textstyle Lid(R):=\sup\{Id(M):M\in\txt{R-mod}\}=\sup_{M\in\txt{R-mod}}\inf_{I^\ast_M}l(I^\ast_M).\nn
\eea
\end{dfn}

\begin{rmk}
$Lpd(R)=0$ (resp. $Lid(R)=0$) iff every $R$-module is projective (resp. injective).
\end{rmk}

\begin{lmm}
Let $R$ be a ring. For any $R$-modules $A,B$ with respective projective  and injective resolutions
\bea
&&A_\ast:P_\ast^A\sr{\vep}{\ral}A\ra 0~=~\cdots{\ral}P_n\sr{d_n}{\ral}P_{n-1}\sr{d_{n-1}}{\ral}\cdots \sr{d_1}{\ral}P_0\sr{\vep}{\ral}A{\ral}0,\nn\\
&&B^\ast:0\ra B\sr{\xi}{\ral}I^\ast_B~=~0\ra B\sr{\xi}{\ral} I^0\sr{d^1}{\ral} I^1\sr{d^2}{\ral}\cdots\sr{d^{n-1}}{\ral}I^{n-1}\sr{d^{n}}{\ral}I^n\ral\cdots,\nn
\eea
where $d_n=d_n^A$ and $d^n=d^n_B$, we have the following equalities:
\bea
&& ext^i\big(\ker d_{n-1}^A,B\big)=ext^{n+i}(A,B)~~~~\txt{for all}~~~~i=1,2,\cdots,\nn\\
&& Ext^i\big(A,\coker d_B^{n-1}\big)=Ext^{n+i}(A,B),~~~~\txt{for all}~~~~i=1,2,\cdots.\nn
\eea
\end{lmm}

\begin{proof}
In the projective resolution $A_\ast$, which has the decomposition
\[\adjustbox{scale=0.7}{\bt
        &   &                                    & 0\ar[dr]                            &                        &        0         &      & 0\ar[dr]                            &                              & 0 \\
        &    &                                    &                                     & \ker d_n\ar[dr,hook]\ar[ur] &                  &      &                                     &  \ker d_{n-2}\ar[ur]\ar[dr,hook]  & ~ \\
\cdots \ar[r] &  P_{n+2}\ar[rr]\ar[dr,"d"] &                            & P_{n+1}\ar[ur,"d"]\ar[rr,"d_{n+1}"] &                        & P_n\ar[rr,"d_n"]\ar[dr,"d"] &      & P_{n-1}\ar[ur,"d"]\ar[rr,"d_{n-1}"] &                              & P_{n-2}\ar[r]  & \cdots\\
        &     & \ker d_{n+1}\ar[ur,hook]\ar[ur]\ar[dr]  &                                     &                        &                  &   \ker d_{n-1}\ar[ur,hook]\ar[dr]   &                                     &    & ~ \\
        &   0\ar[ur]  &                                    &     0                               &                        &       0\ar[ur]   &      &   0                                  &   & ~ \\
\et}\]
the subsequence ~$(\ker d_{n-1})_\ast:P^A_{\ast+n}\sr{d_n}{\ral}\ker d_{n-1}~=~\cdots\sr{d_{n+2}}{\ral}P_{n+1}\sr{d_{n+1}}{\ral}P_n\sr{d_n}{\ral}\ker d_{n-1}\ra 0$~ is a projective resolution of $\ker d_{n-1}$. Thus, from the associated complex (i.e., the $Hom_R(-,B)$-image of the deleted resolution of $\ker d_{n-1}$)
\bea
0\ra Hom_R(P_n,B)\sr{Hom_R(d_{n+1},B)}{\ral}Hom_R(P_{n+1},B)\sr{Hom_R(d_{n+2},B)}{\ral}\cdots,\nn
\eea
we find as usual that
{\small
\begin{align}
&\textstyle ext^0(\ker d_{n-1}^A,B)={\ker Hom_R(d_{n+1},B)\over 0}={\im Hom_R(d_n,B)\over 0}\cong Hom_R(\ker d_{n-1}^A,B),\nn\\
&\textstyle ext^1(\ker d_{n-1}^A,B)={\ker Hom_R(d_{n+2},B)\over \im Hom_R(d_{n+1},B)}=H^{n+1}(Hom_R(P_\ast^A,B))=ext^{n+1}(A,B),\nn\\
&\textstyle~~~\vdots\nn\\
&\textstyle ext^i(\ker d_{n-1}^A,B)= ext^{n+i}(A,B),~~~~i=1,2,\cdots\nn
\end{align}}

Similarly, in the injective resolution $B^\ast$, which decomposes as
\[\adjustbox{scale=0.7}{\bt
        &   &                                    & 0\ar[dr]                            &                        &        0         &      & 0\ar[dr]                            &                              & 0 \\
        &    &                                    &                                     & \coker d^{n-1}\ar[dr,hook]\ar[ur] &                  &      &                                     &  \coker d^{n+1}\ar[ur]\ar[dr,hook]  & ~ \\
\cdots \ar[r] &  I^{n-2}\ar[rr]\ar[dr,"d"] &                            & I^{n-1}\ar[ur,"d"]\ar[rr,"d^{n}"] &                        & I^n\ar[rr,"d^{n+1}"]\ar[dr,"d"] &      & I^{n+1}\ar[ur,"d"]\ar[rr,"d^{n+2}"] &                              & I^{n+2}\ar[r]  & \cdots\\
        &     & \coker d^{n-2}\ar[ur,hook]\ar[ur]\ar[dr]  &                                     &                        &                  &   \coker d^{n}\ar[ur,hook]\ar[dr]   &                                     &    & ~ \\
        &   0\ar[ur]  &                                    &     0                               &                        &       0\ar[ur]   &      &   0                                  &   & ~ \\
\et}\]
the subsequence ~$(\coker d^{n-1})^\ast:0\ra \coker d^{n-1}\sr{\ol{d}^n}{\ral}I_B^{\ast+n}~=~0\ra \coker d^{n-1}\sr{\ol{d}^n}{\ral} I^{n}\sr{d^{n+1}}{\ral} I^{n+1}\sr{d^{n+2}}{\ral}\cdots$~ is an injective resolution of ~{\small $\coker d^{n-1}\subset I^n$}. Thus, from the associated complex (i.e., the $Hom_R(A,-)$-image of the deleted resolution of $\coker d^{n-1}$)
\bea
0\ra Hom_R(A,I^{n})\sr{Hom_R(A,d^{n+1})}{\ral} Hom_R(A,I^{n+1})\sr{Hom_R(A,d^{n+2})}{\ral}\cdots\nn
\eea
we find as usual that
{\small
\begin{align}
&\textstyle Ext^0(A,\coker d_B^{n-1})={\ker Hom_R(A,d^{n+1})\over 0}= {\im Hom_R(A,d^n)\over 0}\cong Hom_R\big(A,\coker d_B^{n-1}\big),\nn\\
&\textstyle Ext^1(A,\coker d_B^{n-1})={\ker Hom_R(A,d^{n+2})\over \im Hom_R(A,d^{n+1})}=H^{n+1}(Hom_R(A,I^\ast_B))=Ext^{n+1}(A,B),\nn\\
&\textstyle~~~\vdots\nn\\
&\textstyle Ext^i\big(A,\coker d_B^{n-1}\big)= Ext^{n+i}(A,B),~~~~i=1,2,\cdots\nn \qedhere
\end{align}}
\end{proof}

\begin{thm}\label{ProjDimMnThm}
Let $R$ be a ring and $A,B$ $R$-modules. For any two projective (resp. injective) resolutions,
{\footnotesize\begin{align}
& P_\ast^A:\cdots{\ral}P_n\sr{d_n}{\ral}P_{n-1}\sr{d_{n-1}}{\ral}~\cdots \sr{d_1}{\ral}P_0\ra0~~\big(\txt{resp.}~I^\ast_B:0\ra I^0\sr{d^1}{\ral} I^1\sr{d^2}{\ral}\cdots\sr{d^{n-2}}{\ral}I^{n-2}\sr{d^{n-1}}{\ral}I^{n-1}\sr{d^n}{\ral}\cdots\big),\nn\\
& P'^A_\ast:\cdots{\ral}P'_n\sr{d'_n}{\ral}P'_{n-1}\sr{d'_{n-1}}{\ral}\cdots \sr{d'_1}{\ral}P'_0\ra0~~\big(\txt{resp.}~I'^\ast_B:0\ra I'^0\sr{d'^1}{\ral} I'^1\sr{d'^2}{\ral}\cdots\sr{d'^{n-2}}{\ral}I'^{n-2}\sr{d'^{n-1}}{\ral}I'^{n-1}\sr{d'^n}{\ral}\cdots\big),\nn
\end{align}}$\ker d_{n-1}$ is projective $\iff$ $\ker d'_{n-1}$ is projective (resp. $\coker d^{n-1}$ is injective $\iff$ $\coker d'^{n-1}$ is injective).
\end{thm}
\begin{proof}
By the theorems of homological algebra, which imply $P^A_\ast\simeq P'^A_\ast$ (resp. $I^\ast_B\simeq I'^\ast_B$), we know
{\scriptsize\[ext^n(A,B):=H^n(Hom_R(P^A_\ast,B))=H^n(Hom_R(P'^A_\ast,B))~~\Big(\txt{resp.}~ Ext^n(A,B):=H^n(Hom_R(A,I^\ast_B))=H^n(Hom_R(A,I'^\ast_B)\Big)
\]}is independent of the resolution of $A$ used. Thus, by the preceding lemma,
{\scriptsize\[
ext^i(\ker d_{n-1},B)=ext^{n+i}(A,B)=ext^i(\ker d'_{n-1},B)~~\Big(~\txt{resp.}~Ext^i(A,\coker d^{n-1})=Ext^{n+i}(A,B)=Ext^i(\coker d'^{n-1})~\Big)\nn
\]}for all $i=1,2,\cdots$. Hence, $\ker d_{n-1}$ is projective iff $ext^1(\ker d_{n-1},B)=0$, iff $ext^1(\ker d'_{n-1},B)=0$, iff $\ker d'_{n-1}$ is projective  (resp. $\coker d^{n-1}$ is injective iff $Ext^1(A,\coker d^{n-1})=0$, iff $Ext^1(A,\coker d'^{n-1})=0$, iff $\coker d'^{n-1}$ is injective).
\end{proof}

\begin{crl}\label{HomDimBdCrl}
Let $R$ be a ring and $A,B$ $R$-modules.
(i) $Lpd(R)\leq n$ $\iff$ $ext^{n+1}(A,B)=0$ for all $A,B$ in $R$-mod. (ii) $Lid(R)\leq n$ $\iff$ $Ext^{n+1}(A,B)=0$ for all $A,B$ in $R$-mod.
\end{crl}
\begin{proof}
(i) $Lpd(R)\leq n$ iff $Pd(A)\leq n$ for all $A$ in $R$-mod, iff $\ker d^A_{n-1}$ is projective for every projective resolution of $A$ (for all $A$ in $R$-mod), iff $ext^{n+1}(A,B)=ext^1(\ker d_{n-1}^A,B)=0$ for all $A,B$ in $R$-mod.

(ii) $Lid(R)\leq n$ iff $Id(B)\leq n$ for all $B$ in $R$-mod, iff $\coker d_B^{n-1}$ is injective for every injective resolution of $B$ (for all $B$ in $R$-mod), iff $Ext^{n+1}(A,B)=Ext^1(A,\coker d_B^{n-1})=0$ for all $A,B$ in $R$-mod.
\end{proof}

\begin{thm}
$Lpd(R)=Lid(R)=\sup\left\{Pd\left({R\over I}\right):~_RI\subset R~\txt{a left ideal}\right\}$
\end{thm}
\begin{proof}
By definition, ~$Lpd(R)=\sup_{A\in\txt{R-mod}}Pd(A)$.
\begin{enumerate}[leftmargin=1.2cm]
\item By the isomorphism $ext\cong Ext$ (which ``balances'' left-projective resolutions and right-injective resolutions), the left projective dimension of $R$ equals the left injective dimension of $R$, i.e.,
    \bea
    Lpd(R)=Lid(R).\nn
    \eea
    \item From the preceding corollary, $Lpd(R)\leq n$ $\iff$ $ext^{n+1}_R(A,B)=0$~ for all $A,B$ in R-mod.

\item \ul{Lemma (L1)}: Recall that for any given injective resolution $0\ra B\sr{\vep}{\ral} E^0\sr{d^1}{\ral} E^1\sr{d^2}{\ral}\cdots$ of $B$,
    \bea
    &&Ext_R^{n+1}(A,B)\cong Ext_R^1\big(A,\coker(d^{n-1}_B)\big),~~~~\txt{for all}~~A\in \txt{R-mod}.\nn
    \eea

\item \ul{Lemma (L2)}: By the Baer criterion, an R-module $E$ is injective $\iff$ $ext^1_R(R/I,E)=0$ ($\iff$ $Hom_R(-,E)$ is exact) for every left ideal $_RI\subset R$. This comes from the exact $ext$-sequence associated with the SES ~$0\ra I\hookrightarrow R\ra R/I\ra 0$,~ i.e.,
    \bea
    0\ra Hom_R(R/I,E)\ra Hom_R(R,E)\ra Hom_R(I,E)\ra  ext^1(R/I,E)\ra\cdots\nn
    \eea

\item Thus, the $Ext$-sequence {\footnotesize $\{Ext_R^i(A,B)\}_i$} above zeros off after $n$ {\footnotesize $\sr{(L1)}{\iff}$ $E:=\coker(d^{n-1}_B)$} is injective, $\sr{(L2)}{\iff}$
    \bea
    Ext_R^{n+1}(R/I,B)\sr{(L1)}{=}Ext_R^1\big(R/I,\coker(d^{n-1}_B)\big)\sr{(L2)}{=}0~~~~\txt{for all}~~{}_RI\subset R.\nn
    \eea
    This means~ $Lid(R)\leq n$ iff $Ext^{n+1}_R(R/I,B)=0$ for all $_RI\subset R$ (and all $B$ in R-mod), and so
\bea
\textstyle Lid(R)=\sup_{_RI\subset R}Pd(R/I).\nn
\eea
\item Therefore, in the definition ~$Lpd(R)=\sup_{A\in\txt{R-mod}}Pd(A)$,~ we can replace $A$ by $R/I$ and get
\[
\textstyle Lpd(R)=Lid(R)=\sup_{_RI\subset R}Pd(R/I). \qedhere
\]
\end{enumerate}
\end{proof}

\begin{crl}
Let $R$ be a ring and $_RI\subset R$ a left ideal. Then (i) either $Pd(R/I)=0$ (i.e., $R/I$ is projective), or  $Pd(R/I)=1+Pd(I)$. Hence, (ii) $Lpd(R)\leq 1+\sup_{{}_RI\subset R}Pd(I)$, and (iii) if ~$Lpd(R)\neq 0$,~ then ~$Lpd(R)=1+\sup_{{}_RI\subset R}Pd(I)$.
\end{crl}
\begin{proof}
Consider the exact sequence $0\ra I\hookrightarrow R\sr{\pi}{\ra} R/I\ra 0$, with exact ext-sequence
{\footnotesize\[
\textstyle 0\ra \cdots\ra Hom_R(I,X)\ra ext^1({R\over I},X)\ra \ub{ext^1(R,X)}_0\ra~ext^1(I,X)\ra ext^2({R\over I},X)\ra \ub{ext^2(R,X)}_0\ra ext^2(I,X)\ra\cdots
\]}
This sequence shows that $ext^i(I,X)\cong ext^{i+1}(R/I,X)$ for all $i\geq 1$. Therefore, if $R/I$ is not projective (i.e., $ext^1(R/I,X)\neq 0$), then $Pd(R/I)=1+Pd(I)$. That is, for all $_RI\subset R$, either $Pd(R/I)=0$ (i.e., $ext^1(R/I,X)=0$) or $Pd(R/I)=1+Pd(I)$. Therefore,
\[
\textstyle Lpd(R)=\sup_{_RI\subset R}Pd(R/I)=
\left\{
      \begin{array}{ll}
        0, & \txt{if $R/I$ is projective for all $_RI\subset R$}, \\
        1+\sup_{_RI\subset R}Pd(I), & \txt{otherwise}.
      \end{array}
    \right. \qedhere
\]
\end{proof}

\begin{crl}
Let $R$ be a ring. Then (a) $Lpd(R)=Lid(R)\leq 1$ $\iff$ (b) every submodule of a projective module is projective, $\iff$ (c) every left ideal ${}_RI\subset R$ is projective.

(Equivalently, $Lid(R)=Lpd(R)\leq 1$ $\iff$ every quotient module of an injective module is injective.)
\end{crl}
\begin{proof}
(a)$\Ra$(b): If $Lpd(R)\leq 1$, then for any projective module $_RP$ and a submodule $_RN\subset P$, the exact ext-sequence for the SES $0\ra N\hookrightarrow P\twoheadrightarrow P/N\ra 0$ gives $ext^1(N,X)=0$ as follows:
{\footnotesize\[
0\ra \cdots\ra Hom_R(N,X)\ra ext^1(P/N,X)\ra \ub{ext^1(P,X)}_0\ra ext^1(N,X)\ra \ub{ext^2(P/N,X)}_0\ra \cdots
\]}
(b)$\Ra$(c): If every submodule of a projective is projective, then in particular every ${}_RI\subset R$ is projective. (c)$\Ra$(a): If every $_RI\subset R$ is projective, then ~$Lpd(R)\leq 1+\sup_{{}_RI\subset R}Pd(I)=1+0=1$.

To prove every quotient module of an injective $R$-module is injective, consider an injective module $_RE$, a submodule $_RN\subset E$, and the exact Ext-sequence for the SES $0\ra N\hookrightarrow E\twoheadrightarrow E/N\ra 0$, which gives $Ext^1(X,E/N)=0$ as follows:
{\footnotesize\[
0\ra \cdots\ra Hom_R(X,E/N)\ra Ext^1(X,N)\ra \ub{Ext^1(X,E)}_0\ra Ext^1(X,E/N)\ra \ub{Ext^2(X,N)}_0\ra \cdots \qedhere
\]}
\end{proof}

\begin{crl}
Let $R$ be a PID. Then $Lpd(R)=Lid(R)\leq 1$. Hence, (i) every submodule of a projective $R$-module is projective (\blue{footnote}\footnote{We will see later that (i) every submodule of a free PID-module is free and (ii) every projective PID-module is free.}), and (ii) every quotient module of an injective $R$-module is injective.
\end{crl}
\begin{proof}
With no zerodivisors in $R$, every ideal $_RI=Ra\subset R$ is a free (hence projective) $R$-module.
\end{proof}

\begin{prp}[\blue{Projective dimension (Pd) bounds in SES's}]\label{ProjDimBdPrp}
Let $R$ be a ring and $0\ra A\ra B\ra C\ra 0$ an exact sequence in $R$-mod. Then we have the following bounds:
\begin{enumerate}
\item $Pd(A)<Pd(B)~~~~\Ra~~~~Pd(C)=Pd(B)$
\item $Pd(A)<Pd(C)~~~~\Ra~~~~Pd(B)\leq Pd(C)$
\item $Pd(B)<Pd(A)~~~~\Ra~~~~Pd(C)=Pd(A)+1$
\item $Pd(B)<Pd(C)~~~~\Ra~~~~Pd(A)<Pd(C)$
\item $Pd(C)<Pd(B)~~~~\Ra~~~~Pd(A)=Pd(B)$
\item $Pd(C)<Pd(A)~~~~\Ra~~~~Pd(A)=Pd(B)$
\item $Pd(A)\neq Pd(B)~~~~\Ra~~~~\max\{Pd(A),Pd(B)\}\leq Pd(C)$
\item $Pd(A)\leq \max\{Pd(B),Pd(C)\}$
\item $Pd(B)\leq \max\{Pd(A),Pd(C)\}$
\item $Pd(C)\leq \max\{Pd(A),Pd(B)\}+1$
\end{enumerate}
\end{prp}
\begin{proof}
Recall that the exact ext-sequence for the SES $0\ra A\ra B\ra C\ra 0$ is
{\footnotesize\[
0\ra\cdots\ra ext^n(C,X)\ra ext^n(B,X)\ra ext^n(A,X)
\ra ext^{n+1}(C,X)\ra ext^{n+1}(B,X)\ra ext^{n+1}(A,X)\ra\cdots\nn
\]}Also, by Lemma \ref{HomDimBdCrl} (or its proof), for any $_RM$, $Pd(M)\leq n$ $\iff$ $ext^{n+1}(M,X)=0$ (for all $X\in R$-mod), $\iff$ $ext^{n+i}(M,X)=0$ (for all $X\in R$-mod and $i\geq 1$).

\begin{enumerate}[leftmargin=0.9cm]
\item $Pd(A)<n:=Pd(B)$ $\iff$ $Pd(A)\leq n-1$, $\iff$ $ext^{n-1+i}(A,X)=0$ (all $X$ and $i\geq 1$). So, the ext-sequence gives $Pd(C)=n$ as follows:
{\small\[
0\ra\cdots\ra ext^n(C,X)\ra \ub{ext^n(B,X)}^{\neq 0}\ra 0
\ra ext^{n+1}(C,X)\ra 0\ra 0\ra\cdots\nn
\]}

\item $Pd(A)<n:=Pd(C)=n$ $\iff$ $Pd(A)\leq n-1$, $\iff$ $ext^{n-1+i}(A,X)=0$. So, the ext-sequence gives $ext^{n+1}(B,X)=0$, i.e., $Pd(B)\leq n$, as follows:
{\small\[
0\ra\cdots\ra ext^n(C,X)\ra ext^n(B,X)\ra 0
\ra 0\ra ext^{n+1}(B,X)\ra 0\ra\cdots\nn
\]}

\item $Pd(B)<n:=Pd(A)$ $\iff$ $Pd(B)\leq n-1$, $\iff$ $ext^{n-1+i}(B,X)=0$. So, the ext-sequence gives $ext^n(A,X)\cong ext^{n+1}(C,X)$, i.e., $Pd(C)=Pd(A)+1$, as follows:
{\small\[
0\ra\cdots\ra ext^n(C,X)\ra 0\ra ext^n(A,X)
\ra ext^{n+1}(C,X)\ra 0\ra 0\ra\cdots\nn
\]}

\item $Pd(B)<n:=Pd(C)$ $\iff$ $Pd(B)\leq n-1$, $\iff$ $ext^{n-1+i}(B,X)=0$. So, the ext-sequence gives $Pd(A)<Pd(C)$, as follows:
{\small\[
0\ra\cdots\ra ext^n(C,X)\ra 0\ra ext^n(A,X)
\ra 0\ra 0\ra ext^{n+1}(A,X)\ra\cdots\nn
\]}

\item $Pd(C)<n:=Pd(B)$ $\iff$ $Pd(C)<n-1$, $\iff$ $ext^{n-1+i}(C,X)=0$. So, the $ext$-sequence gives $ext^n(A,X)\cong ext^n(B,X)$, i.e., $Pd(A)=Pd(B)$, as follows:
{\small\[
0\ra\cdots\ra 0\ra ext^n(B,X)\ra ext^n(A,X)
\ra 0\ra 0\ra ext^{n+1}(A,X)\ra\cdots\nn
\]}

\item $Pd(C)<n:=Pd(A)$ $\iff$ $Pd(C)<n-1$, $\iff$ $ext^{n-1+i}(C,X)=0$. So, the $ext$-sequence gives $ext^n(A,X)\cong ext^n(B,X)$, i.e., $Pd(A)=Pd(B)$, as follows:
{\small\[
0\ra\cdots\ra 0\ra ext^n(B,X)\ra ext^n(A,X)
\ra 0\ra ext^{n+1}(B,X)\ra 0\ra\cdots\nn
\]}
\item This is an immediate consequence of the last two results above.

\item If $Pd(B),Pd(C)\leq n$ then $Pd(A)\leq n$, since $ext^{n+1}(A,X)=0$ as follows:
{\small\[
0\ra\cdots\ra ext^n(A,X)\ra \ub{ext^{n+1}(C,X)}_0\ra \ub{ext^{n+1}(B,X)}_0\ra ext^{n+1}(A,X)\ra \ub{ext^{n+2}(C,X)}_0\ra\cdots\nn
\]}
\item If $Pd(A),Pd(C)\leq n$ then $Pd(B)\leq n$, since $ext^{n+1}(B,X)=0$ as follows:
{\small\[
0\ra\cdots\ra ext^n(A,X)\ra \ub{ext^{n+1}(C,X)}_0\ra ext^{n+1}(B,X)\ra \ub{ext^{n+1}(A,X)}_0\ra \ub{ext^{n+2}(C,X)}_0\ra\cdots\nn
\]}

\item If~ $Pd(A),Pd(B)\leq n$ then $Pd(C)\leq n+1$, since $ext^{n+2}(C,X)=0$ as follows:
{\footnotesize\[
0\ra\cdots\ra ext^n(A,X)\ra ext^{n+1}(C,X)\ra \ub{ext^{n+1}(B,X)}_0\ra \ub{ext^{n+1}(A,X)}_0\ra ext^{n+2}(C,X)\ra \ub{ext^{n+2}(B,X)}_0\ra\cdots \qedhere
\]}
\end{enumerate}
\end{proof}

The following is an immediate consequence of the last three bounds above.
\begin{crl}[\blue{Solvability in a SES}]
Let $R$ be a ring and $0\ra A\ra B\ra C\ra 0$ a SES in $R$-mod.
\begin{enumerate}[leftmargin=0.9cm]
\item If any two of $A,B,C$ have finite projective dimension (Pd), then so does the third.
\item Equivalently, if one of $A,B,C$ has infinite Pd, then at least one of the other two also has infinite Pd.
\item Equivalently, if one of $A,B,C$ has finite Pd, then the other two either both have finite Pd or both have infinite Pd.
\end{enumerate}
\end{crl}

\begin{crl}[\blue{\cite[Exercise 10.67, p.955]{rotman2010}, \cite[Exercise 8.9, p.467]{rotman2009}}]
If $\{M_\al\}$ is a family of $R$-modules,
\[\textstyle Pd(\oplus_\al M_\al)=\sup_\al Pd(M_\al).\]
\end{crl}
\begin{proof}
Let $M:=\oplus_\al M_\al$. \ul{Proving $\geqslant$}: Let $P_\ast^{M_\al}$ be a projective resolution of $M_\al$ of minimal length. Then by the exactness of the functor $\oplus$ on $R$-mod, $P_\ast^M:=\oplus_\al P_\ast^{M_\al}$ is a projective resolution of $M$. Thus,
\bea
\textstyle l(P_\ast^M)\geq\max_\al l(P_\ast^{M_\al})~~\Ra~~Pd(M)\geq \sup_\al Pd(M_\al).\nn
\eea
{\flushleft\ul{Proving $\leqslant$}}: Assume wlog the indexing of the modules is in the form $\al\in A\subset\Ordinal\mapsto M_\al$.

(\magenta{Also, we can assume wlog that $A$ is a countable set. Indeed, it suffices to consider only a countable subfamily $\{M_{\al_i}\}_{i\in\Natural}$ given for each $i\in\Natural$ by selecting $M_{\al_i}$ to be any $M_\al$ satisfying ~$Pd(M_\al)=\max\{Pd(M_\beta)\leq i\}_{\beta\in A}$,~ because we then have ~$\sup_iPd(M_{\al_i})=\sup_\al Pd(M_\al)$.})

Let {\footnotesize $d_\al:=\sup_{\al'\leq\al}Pd(M_{\al'})$}, {\footnotesize $L_\al:=\oplus_{\al'\leq\al}M_{\al'}$}, and {\footnotesize $L_{\beta\gamma}:=\oplus_{\beta<\al'\leq\gamma}M_{\al'}$} if $\beta<\gamma$. Then the exact sequences ~{\footnotesize $0\ra L_\beta\ra L_\gamma\ra{L_\gamma\over L_\beta}\cong L_{\beta\gamma}\ra 0$}~ imply
\[
\textstyle Pd(L_\gamma)\leq \max\left\{Pd(L_\beta),Pd\left(L_{\beta\gamma}\right)\right\}.
\]
Let $C:=\{\al\in A:Pd(L_\al)\leq d_\al\}$. Then it is clear that $C\neq\emptyset$, and that $\al\in C$ $\Ra$ $\al^+\in C$ due to the exact sequence $0\ra L_\al\ra L_{\al^+}\ra{L_{\al^+}\over L_\al}\cong M_{\al^+}\ra 0$, which implies
\[
\textstyle Pd(L_{\al^+})\leq \max\left\{Pd(L_\al),Pd\left(M_{\al^+}\right)\right\}\leq \max\left\{d_\al,Pd\left(M_{\al^+}\right)\right\}=d_{\al^+}.
\]
Let $\al\in\Ordinal$ be such that $\al'\in C$ for all $\al'<\al$. If a predecessor $\al^-$ exists, then $\al=(\al^-)^+\in C$.

(\magenta{At this point, the proof is complete since $C$ is countably-inductive and we have already noted it is sufficient to assume $A$ is countable. However, without this observation we can still continue as follows.})

Assume wlog the modules are arranged so that at each limit ordinal $\al\in A$, we have $Pd(L_{\al'\al})\leq d_{\al'}:=\sup_{\al''\leq\al'}Pd(M_{\al''})$ as $\al'\uparrow\al$ (otherwise, if there exists a sequence $\al'_i\uparrow\al$, $i\in\Natural$, such that $d_{\al'_i}<Pd(L_{\al_i'\al})$ for all $i\in\Natural$, then $d_\al=\lim_i d_{\al'_i}\leq \lim\sup_i Pd(L_{\al_i'\al}):=\lim_n\sup_{i\geq n}Pd(L_{\al'_i\al})=0$, a contradiction).

So, if $\al$ is a limit ordinal, then with $\al'<\al$ and $\al'\ra\al$, we get
\bea
Pd(L_\al)\leq\max\big(Pd(L_{\al'}),Pd(L_{\al'\al})\big)\leq \max\big(d_{\al'},Pd(L_{\al'\al})\big)\ra d_\al.\nn
\eea
This shows $C\subset A$ is inductive, and so $C=A$, i.e., $Pd(M)\leq\sup_{\al\in A}M_\al$.
\end{proof}
\begin{crl}
If a SES of $R$-modules {\small $0\ra A\ra B\ra C\ra 0$} is split, then {\small $Pd(B)=\max\{Pd(A),Pd(C)\}$}.
\end{crl}

\begin{prp}
Let $R$ be a ring and $M,N$ $R$-modules. If $Pd_R(M)\leq m$ and $Pd_R(N)\leq n$, then $Pd_R(M\otimes_RN)\leq m+n$.
\end{prp}
\begin{proof}
If $P_\ast^M$ and $P_\ast^N$ are projective resolutions of $M$ and $N$, then the total complex $P_\ast^{M\otimes N}:=Tot(P_\ast^M\otimes_R P_\ast^N)_\ast:=\bigoplus_{i+j=\ast}P_i^M\otimes_RP_j^N$ is a projective resolution of $M\otimes_RN$ of length $\leq n+m$. Recall that (i) a coproduct of free objects is free, (ii) a tensor product of free modules $(\oplus_\al R)\otimes_R(\oplus_\beta R)\cong \oplus_{\al,\beta}R\otimes_R R\cong\oplus_{\al,\beta}R$ is free, and (iii) a projective module is precisely a direct summand of a free module. In particular, if $F_i=P_i\oplus P_i'$ and $F_j=P_j\oplus P_j'$ are free $R$-modules, then so is {\small $F_i\otimes_RF_j=(P_i\oplus P_i')\otimes_R(P_j\oplus P_j')\cong (P_i\otimes_RP_j)\oplus N$}, where {\small $N:=(P_i\otimes_RP_j')\oplus(P_i'\otimes P_j)\oplus(P_i'\otimes_RP_j')$}.
\end{proof}

\begin{thm}[\blue{\index{Ping-pong lemma}{Ping-pong lemma}}]
Let $R$ be a ring and ${}_RI,{}_RJ\subset R$ left ideals. If (i) one of $I,J$ is not projective, and (ii) there exist two exact sequences~ $0\ra I\ra R\ra J\ra 0$ and $0\ra J\ra R\ra I\ra 0$,~ then $Lpd(R)=\infty$. (\blue{footnote}\footnote{For example, $Lpd(\Integer_4)=\infty$ (where $\Integer_n:=\Integer/n\Integer$) due to the exact sequence $0\ra 2\Integer_4\hookrightarrow \Integer_4\sr{\times 2}{\ral}2\Integer_4\ra 0$, which does not split (because $2\Integer_4=2\Integer/4\Integer\cong \Integer/2\Integer=\Integer_2$ while $\Integer_4\not\cong \Integer_2\oplus \Integer_2$). Alternatively, we have the infinite projective resolution
$\cdots\sr{\times 2}{\ral}\Integer_4\sr{\times 2}{\ral}\Integer_4\sr{\times 2}{\ral}\Integer_4\sr{\times 2}{\ral}2\Integer_2\ra 0$.})
\end{thm}
\begin{proof}
 By hypotheses, both $I,J$ are not projective, and so $Pd(I),Pd(J)\geq 1$. This implies $Pd(R)=0<1\leq Pd(I),Pd(J)$, and so $Pd(J)=Pd(I)+1$ and $Pd(I)=Pd(J)+1$. Thus, if $Pd(I)<\infty$ and $Pd(J)<\infty$, then we get a contradiction in the form $0=2$. (Alternatively, observe that using the two exact sequences, we can form a projective resolution of infinite length: $\cdots\ra R\ra R\ra R\ra J\ra 0$.)
\end{proof}

\subsection{Homological dimension of polynomial rings}
\begin{lmm}\label{SyzygyLmm1}
Let $R$ be a ring, and $z\in R$. If $z$ is (i) central in the sense $z\in Z(R)$, (\blue{footnote}\footnote{As usual, $Z(R):=\{z\in R:zr=rz~\txt{for all}~r\in R\}$ denotes the center of $R$.}), (ii) a non-zerodivisor in the sense $zr\neq 0$ for all $r\in R\backslash\{0\}$, and (iii) a non-unit, then $Pd_R\left({R\over Rz}\right)=1$.
\end{lmm}
    \begin{proof}
    Consider the projective resolution $0\ra R\sr{\times z}{\ral}R\ra {R\over Rz}\ra 0$, where $\times z:A\ra R,~r\mapsto zr$ is multiplication by $z$. The sequence does not split, otherwise, we get $R=\im(\times z)\oplus N=Rz\oplus N$ with $N\cong{R\over Rz}$, and so $1=z'z+n$ (for some $z'\in R$) $\Ra$ $z=z(z'z+n)=z'z^2+0$, $\Ra$ $z(1-z'z)=0$, $\Ra$ $z'z=zz'=1$, i.e., $z$ is a unit (a contradiction). Therefore, we conclude that ${R\over Rz}$ is not projective (hence $Pd_R(R/Rz)\neq 0$). Also, the length of the resolution is $1$. Hence $Pd_R(R/Rz)=1$.
    \end{proof}

\begin{lmm}\label{SyzygyLmm2}
Let $R$ be a PID. If $R$ is not a field, then $Lpd(R)=1$.
\end{lmm}
\begin{proof}
Since $R$ is a PID, $Lpd(R)=\sup\left\{Pd\left({R\over I}\right):~_RI\subset R\right\}=\sup\left\{Pd\left({R\over Rz}\right):~z\in R\right\}\leq 1$. Since $R$ is not a field, we have at least one exact sequence $0\ra Rz\hookrightarrow R\ra R/Rz\ra 0$ (for $z\in R$) that does not split (by the proof of Lemma \ref{SyzygyLmm1}), and so $Lpd(R)\neq 0$.
\end{proof}

\begin{crl}[\blue{$Lpd(\Integer)=Lid(\Integer)=1$}]
\end{crl}

\begin{dfn}[\blue{\index{Free! polynomial ring}{Free polynomial ring}}]
Let $R'$ be a ring, $R\subset R'$ a subring, $x,x_1,...,x_n\in R'$ variables (indeterminates), and {\small $X:=\{1,x,x^2,\cdots\}$}. The \ul{free polynomial $R$-ring in one variable} $x$ is the free $R$-module
\[
\textstyle R[x]:=R^{\langle X\rangle}=\sum_{i\in\Natural}Rx^i=\{\sum_if^ix^i:f^i\in R,~f^i=0~\txt{a.e.f.}~i\in\Natural\}\cong\bigoplus_{i\in\Natural}R
\]
with basis $X$ viewed as a polynomial ring in which the basis elements $\{x^i\}_{i\in\Natural}$ commute with elements of $R$, and, addition and multiplication are given by (\blue{footnote}\footnote{When $R$ is commutative, $R[x]$ is equivalent to the usual polynomial ring in one variable.})
{\small\[
\textstyle \sum_i f^ix^i+\sum_i g^ix^i:=\sum_i(f^i+g^i)x^i,~~~~\left(\sum_i f^ix^i\right)\left(\sum_i g^ix^i\right):=\sum_ih^ix^i,~~~~h^i:=\sum_{j+k=i}f^jg^k.
\]}We also inductively define the \ul{free polynomial $R$-ring in $n$-variables} ~$R[x_1,...,x_n]:=\big(R[x_1,...,x_{n-1}]\big)[x_n]$.

If $_RM$ is an $R$-module, let $M[x]:=R[x]\otimes_RM=\{\sum_ix^i\otimes m^i:m^i\in M,~m^i=0~\txt{a.e.f.}~i\in\Natural\}\cong\bigoplus_{i\in\Natural}M$, which can be viewed as a ``free polynomial ring'' with coefficients in $M$.

(It is clear that $M[x]$ is both an $R$-module and an $R[x]$-module, and so has both projective $R$-dimension $Pd_R(M[x])$ and projective $R[x]$-dimension $Pd_{R[x]}(M[x])$, as well as injective $R$-dimension $Id_R(M[x])$ and injective $R[x]$-dimension $Id_{R[x]}(M[x])$.)
\end{dfn}

\begin{lmm}[\blue{\cite[Lemma 8.29, p.468]{rotman2009}}]\label{SyzygyPreLmm1}
For an $R$-module $M$, ~$Pd_R(M)=Pd_R(M[x])=Pd_{R[x]}(M[x])$.
\end{lmm}
\begin{proof}
Let $n\geq 0$. If $Pd_R(M)\leq n$, then there exists a projective $R$-resolution $0\ra P_n\ra\cdots \ra P_0\ra M\ra 0$. Since $R[x]\otimes_R-$ is flat (as $R[x]$ is a free $R$-module), we get a projective $R[x]$-resolution $0\ra R[x]\otimes_RP_n\ra\cdots\ra R[x]\otimes_R P_0\ra M[x]\ra 0$, and so $Pd_{R[x]}(M[x])\leq n$. Therefore, $Pd_R(M)\geq Pd_{R[x]}(M[x])$, since $n\geq 0$ was arbitrary.

Similarly, if $Pd_{R[x]}(M[x])\leq n$, then there exists a projective $R[x]$-resolution $0\ra Q_n\ra\cdots \ra Q_0\ra M[x]\ra 0$, where each $_{R[x]}Q_i$ is projective, hence a direct summand of a free module ${}_{R[x]}F_i={}_{R[x]}Q_i\oplus{}_{R[x]}Q_i'$. Since $F_i$ is also a free $R$-module, it follows that $Q_i$ is also projective as an $R$-module, i.e., $0\ra Q_n\ra\cdots \ra Q_0\ra M[x]\ra 0$ is a projective $R$-resolution of $M[x]$, and so $Pd_R(M[x])\leq n$. That is, $Pd_{R[x]}(M[x])\geq Pd_R(M[x])$. Hence,
\[
\textstyle Pd_R(M)\geq  Pd_{R[x]}(M[x])\geq Pd_R(M[x])=Pd_R(\bigoplus_{i\in\Natural}M)=\sup_iPd_R(M)=Pd_R(M). \qedhere
\]
\end{proof}

\begin{crl}[\blue{\cite[Corollary 8.30, p.468]{rotman2009}}]\label{SyzygyPreCrl1}
If $Lpd(R)=\infty$, then (i) there exists an $R$-module $M$ with $Pd_R(M)=\infty$, and (ii) $Lpd(R[x])=\infty$ as well.
\end{crl}
\begin{proof}
If $Lpd(R)=\sup_{{}_RM} Pd_R(M)=\infty$, pick a sequence of $R$-modules $M_n$ such that $\lim Pd_R(M_n)=\infty$. Then with ${}_RM:=\bigoplus_{n\in\Natural}M_n$, we get $Pd_R(M)=\sup_nPd(M_n)=\infty$. Therefore, $Lpd(R[x])\geq Pd_{R[x]}(M[x])=Pd_R(M)=\infty$.
\end{proof}

\begin{lmm}[\blue{\cite[Lemma 8.31, p.468]{rotman2009}}]\label{SyzygyPreLmm2}
If $N$ is an $R[x]$-module (hence an $R$-module) and $N[x]:=R[x]\otimes_RN$, then there exists an exact sequence ~$0\ra N[x]\ra N[x]\ra N\ra 0$.
\end{lmm}
\begin{proof}
Consider the $R[x]$-homomorphism $\pi:N[x]\ra N$, $\sum x^i\otimes n^i\mapsto \sum x^in^i$. Then we get an exact sequence $0\ra\ker\pi\hookrightarrow N[x]\sr{\pi}{\ral}N\ra 0$. Next, consider the $R[x]$-homomorphism $\phi:N[x]\ra\ker\pi,~\sum x^i\otimes n^i\mapsto\sum x^i(1\otimes x-x\otimes 1)n^i$. Then $\phi$ is injective: Indeed, we can explicitly check that
\[
\textstyle\ker\phi=\left\{\sum_{i=0}^k x^i\otimes n^i:\sum_{i=0}^k x^i(1\otimes x-x\otimes 1)n^i=0\right\}=\{0\},
\] because by expanding the left side of ~$\sum_{i=0}^k x^i(1\otimes x-x\otimes 1)n^i=0$~ we get
\begin{align}
&\textstyle\sum_{i=0}^k x^i(1\otimes x-x\otimes 1)n^i=\sum_{i=0}^k x^i\otimes xn^i-\sum_{i=0}^kx^{i+1}\otimes n^i=\sum_{i=0}^k x^i\otimes xn^i-\sum_{i=1}^{k+1}x^i\otimes n^{i-1}\nn\\
&\textstyle~~~~=1\otimes xn^0+\sum_{i=1}^kx^i\otimes(xn^i-n^{i-1})-x^{k+1}\otimes n^k=0,\nn\\
&~~\Ra~~0=n^k=xn^k-n^{k-1}=xn^{k-1}-n^{k-2}=\cdots=xn^1-n^0,~~\Ra~~0=n^k=n^{k-1}=\cdots=n^0.\nn
\end{align}
Also, $\phi$ is surjective: Indeed, if {\footnotesize $\sum_{i=0}^k x^i\otimes n^i\in\ker\pi$}, then {\footnotesize $\sum x^in^i=n^0+xn^1+x^2n^2+x^3n^3+\cdots+x^kn^k=0$}, which we can rewrite in the for ~$n^0=-xz^0$,~ where
\begin{align}
& z^0:=n^1+xn^2+x^2n^3+\cdots+x^{k-1}n^k=n^1+xz^1,\nn\\
& z^1:=n^2+xn^3+x^2n^4+\cdots+x^{k-2}n^k=n^2+xz^2,\nn\\
& z^2:=n^3+xn^4+x^2n^5+\cdots+x^{k-3}n^k=n^3+xz^3,\nn\\
& \hspace{0.8cm}\vdots\hspace{1.5cm}\vdots\hspace{1.5cm}\vdots \nn\\
& z^{k-1}:=n^{k-1}+xn^k=n^{k-1}+xz^k,\nn\\
& z^k:=n^k,\nn
\end{align}
and so {\footnotesize $\sum x^i(x\otimes 1-1\otimes x)z^i=\sum(x^{i+1}\otimes z^i-x^i\otimes xz^i)=\sum\big(x^{i+1}\otimes z^i-x^i\otimes(z^{i-1}-n^i)\big)=\sum\big(x^{i+1}\otimes z^i-x^i\otimes z^{i-1}+x^i\otimes n^i\big)=\sum x^i\otimes n^i$}, i.e., $\phi|_{\ker\pi}=id_{\ker\pi}$. Since this shows $\phi$ is a $k[x]$-isomorphism, we get the exact sequence ~$0\ra N[x]\sr{\phi}{\ral}N[x]\sr{\pi}{\ral}N\ra 0$.
\end{proof}

\begin{crl}[\blue{\cite[Corollary 8.32, p.469]{rotman2009}}]\label{SyzygyPreCrl2}
Let $R$ be a ring. Then ~$Lpd(R[x])\leq 1+Lpd(R)$.
\end{crl}
\begin{proof}
The case of $Lpd(R)=\infty$ is clear by Corollary \ref{SyzygyPreCrl1}. So assume $n:=Lpd(R)<\infty$. We need to show $Lpd(R[x])\leq 1+n$. Let $N$ be an $R[x]$-module (hence an $R$-module). Then by Lemma \ref{SyzygyPreLmm2}, we have an exact sequence $0\ra N[x]\ra N[x]\ra N\ra 0$, and so using (Lemma \ref{SyzygyPreLmm1} in) Proposition \ref{ProjDimBdPrp}(10), we get ~$Pd_{R[x]}(N)\leq Pd_{R[x]}(N[x])+1=Pd_R(N)+1\leq n+1$.
\end{proof}

\begin{rmk}[\blue{About showing ~$Lpd(R[x])\geq 1+Lpd(R)$.}]
Let $n:=Lpd(R)$, i.e., (i) $ext_R^{n+1}(M',N')=0$ for all ${}_RM',{}_RN'$ and (ii) $ext_R^n(M,N)\neq 0$ for some ${}_RM,{}_RN$ (in particular, we can pick $_RM$ such that $Pd_R(M)=n$). To show $Lpd(R[x])\geq 1+Lpd(R)$, it suffices to show $ext_{R[x]}^{n+1}(A,B)\neq 0$ for some ${}_{R[x]}A,{}_{R[x]}B$.

Note that $R\cong_{\txt{Ring}}{R[x]\over R[x]x}$ (which implies $R$, hence every $R$-module, is an $R[x]$-module), and so we have an exact sequence of $R[x]$-modules ~$0\ra R[x]\sr{x}{\ral}R[x]\ra R\cong_{Ring}{R[x]\over R[x]x}\ra 0$.
Also, by Lemma \ref{SyzygyLmm1},
\[
\textstyle Pd_{R[x]}(R)=Pd_{R[x]}({R[x]\over R[x]x})=1.
\]
\end{rmk}

\begin{lmm}\label{SyzygyPreLmm3}
Let $R$ be a ring and $_RM$ an $R$-module with $n:=Pd_R(M)<\infty$. Then there exists a free $R$-module $F$ such that $Ext_R^n(M,F)\neq 0$.
\end{lmm}
\begin{proof}
It is clear that there exists an $R$-module $X$ such that $Ext_R^n(M,X)\neq 0$, while $Ext_R^{n+1}(M,X)=0$. Consider an exact sequence of $R$-modules $0\ra\ker\pi\ra F\sr{\pi}{\ral} X\ra 0$, where $F$ is free. In the following exact Ext-sequence, we must have $Ext_R^n(M,F)\neq 0$ since $Ext_R^n(M,X)\neq 0$:
{\small\[
0\ra\cdots\ra Ext_R^n(M,\ker\pi)\ra Ext_R^n(M,F)\ra \ub{Ext_R^n(M,X)}_{\neq 0}\ra \ub{Ext_R^{n+1}(M,\ker\pi)}_0\ra\cdots \qedhere
\]}
\end{proof}

\begin{lmm}[\blue{\cite[Proposition 8.33, p.470]{rotman2009}}]\label{ChngRinsHom}
A ring homomorphism $\rho:R\ra R'$ induces an exact additive imbedding of categories ~$F_\rho:R'\txt{-mod}\hookrightarrow R\txt{-mod},~M'\sr{f'}{\ral}N'~\mapsto~M'\sr{f'}{\ral}N'$.~ Moreover, if $\rho$ is surjective, then $F_\rho$ becomes a full imbedding.
\end{lmm}
\begin{proof}
An $R'$-module $M'$ inherits an $R$-action $\ast:R\times M'\ra M',~(r,m')\mapsto r\ast m':=\rho(r)m'$, and every $R'$-homomorphism $f':M'\ra N'$ becomes a unique $R$-homomorphism because
\bea
f'(r\ast m')=f'(\rho(r)m')=\rho(r)f'(m')=r\ast f'(m'),~~\txt{for all}~~r\in R,~m'\in M'.\nn
\eea
By construction, wrt the $R$-action $\ast$, $Hom_R(M',N')$ ``contains'' not only $R'$-linear maps but all $\rho(R')$-linear maps as well, and is therefore strictly larger than $Hom_{R'}(M',N')$ in the sense the map $F_\rho:Hom_{R'}(M',N')\ra Hom_R(M',N')$ is not surjective if $\rho$ is not surjective. If $\rho$ is surjective, then {\small $F_\rho:Hom_{R'}(M',N')\ra Hom_R(M',N')$} is also surjective: Indeed, every $R$-linear map $f:M'\ra N'$ becomes $R'$-linear because, if $r'\in R'$ then $r'=\rho(r)$ for some $r\in R$,
\[
\Ra~~f(r'm')=f(\rho(r)m')=f(r\ast m')=r\ast f(m')=\rho(r)f(m')=r'f(m'),~~~~\txt{for all}~~m'\in M'.\nn
\]
The exactness, additivity, and mapwise-injectivity of $F_\rho$ are clear since $F_\rho$ maps through the identity functor ~$id_{R'\txt{-mod}}:R'\txt{-mod}\ra R'\txt{-mod}$.
\end{proof}

\begin{thm}[\blue{\cite[Theorem 8.34, p.470]{rotman2009}}]\label{SyzygyThm0}
Let $R$ be a ring, and $z\in R\backslash 0$. If $z$ is (i) central in the sense $z\in Z(R)$, (ii) a non-zerodivisor, and (iii) a non-unit, then for every $R$-module $M$ (on which $z$ is \ul{regular} or \ul{non-annihilating} in the sense $z:M\ra M,~m\mapsto zm$ is injective) and every $R/Rz$-module $X^\ast$,
\[
\textstyle ext_{R/Rz}^n(X^\ast,{M/zM})\cong_{Ab} ext_R^{n+1}(X^\ast,M),~~~~\txt{for all}~~n\geq 0.\nn
\]
(By the preceding lemma the surjective ring homomorphism $\pi:R\ra R^\ast:=R^\ast,~r\mapsto r+Rz$ gives an exact full imbedding $F_\pi:R^\ast\txt{-mod}\hookrightarrow R\txt{-mod}$. Thus, $R^\ast$-modules $M^\ast,N^\ast$ are $R$-modules with $Hom_R(M^\ast,N^\ast)\cong_{Ab} Hom_{R^\ast}(M^\ast,N^\ast)$.)
\end{thm}
\begin{proof}
Let $R^\ast:=R/Rz$ and $M^\ast:=M/zM$. The additive functors $\{G^n:=ext^n_{R^\ast}(-,M^\ast):R^\ast\txt{-mod}\ra Ab\}_{n\geq 0}$ are uniquely defined (up to isomorphism) by the following basic facts we proved earlier:
\bit
\item[(a)] associated with any SES $0\ra{}_{R^\ast}A\ra{}_{R^\ast}B\ra{}_{R^\ast}C\ra 0$ is an induced LES $0\ra G^0(C)\ra G^0(B)\ra G^0(A)\sr{\delta^0}{\ral} G^1(C)\ra G^1(B)\ra G^1(A)\sr{\delta^1}{\ral} G^2(X)\ra\cdots$ in which the homomorphisms $\delta^n$ are natural or functorial, in the sense that every diagram
\[\adjustbox{scale=0.7}{\bt
0\ar[r] & A\ar[d]\ar[r] & B\ar[d]\ar[r] & C\ar[d]\ar[r] & 0 \\
0\ar[r] & A'\ar[r] & B'\ar[r] & C'\ar[r] & 0
\et}\]
gives a unique commutative diagram (for each $n\geq 0$) of the following form:
\[\adjustbox{scale=0.7}{\bt
\cdots\ar[r]& G^{n-1}(A)\ar[d,dashed]\ar[r,"\delta^{n-1}"] & G^n(C)\ar[d,dashed]\ar[r] & G^n(B)\ar[d,dashed]\ar[r] & G^n(A)\ar[d,dashed]\ar[r,"\delta^n"] & G^{n+1}(C)\ar[d,dashed]\ar[r]&\cdots \\
\cdots\ar[r]& G^{n-1}(A)\ar[r,"\delta^{n-1}"] & G^n(C')\ar[r] & G^n(B')\ar[r] & G^n(A')\ar[r,"\delta'^n"] & G^{n+1}(C')\ar[r]&\cdots
\et}\]
\item[(b)] for $n=0$, ~$G^0(X^\ast)\cong_{Ab} Hom_{R^\ast}(X^\ast,M^\ast)$~ for every $R^\ast$-module $_RX^\ast$, and
\item[(c)] for $n\geq 1$, ~$G^n(F^\ast)=0$ ~for every free $R^\ast$-module $_{R^\ast}F^\ast$. (``free'' and ``projective'' both work here)
\eit
Indeed if $\{G'^n:R^\ast\txt{-mod}\ra Ab\}_{n\geq 0}$ also satisfy (a),(b),(c) then for any module $_{R^\ast}X^\ast$, an SES $0\ra \ker\pi\ra F^\ast\sr{\pi}{\ral}X^\ast\ra 0$ (where ${}_{R^\ast}F^\ast$ is free) gives the commutative diagram below in which (by construction) unique vertical maps exist in both directions and are thus isomorphisms (since $G^n(F^\ast)=0$ for all $n\geq 1$ and the diagram commutes), i.e., $G^n(X^\ast)\cong_{Ab} G'^n(X^\ast)$ for all $n\geq 0$:
\[\adjustbox{scale=0.8}{\bt[column sep=small]
0\ar[r] &\cdots\ar[d,shift right]\ar[from=d,shift right,"\cong"']\ar[r] &G^0(\ker\pi)\ar[d,shift right]\ar[from=d,shift right,"\cong"']\ar[r] & G^1(X^\ast)\ar[d,shift right,dashed,"a^1"']\ar[from=d,shift right,dashed,"a'^1"']\ar[r] & \overbrace{G^1(F^\ast)}^0\ar[d,shift right]\ar[from=d,shift right,"0"']\ar[r] & G^1(\ker\pi)\ar[d,shift right,dashed,"b^1"']\ar[from=d,shift right,dashed,"b'^1"']\ar[r] & G^2(X^\ast)\ar[d,shift right,dashed,"a^2"']\ar[from=d,shift right,dashed,"a'^2"']\ar[r] & \overbrace{G^2(F^\ast)}^0\ar[d,shift right]\ar[from=d,shift right,"0"']\ar[r] & G^2(\ker\pi)\ar[d,shift right,dashed,"b^2"']=\ar[from=d,shift right,dashed,"b'^2"']\ar[r] &~\\
0\ar[r] &\cdots\ar[r]&G'^0(\ker\pi)\ar[r] & G'^1(X^\ast)\ar[r] & \ub{G'^1(F^\ast)}_0\ar[r] & G'^1(\ker\pi)\ar[r] & G'^2(X^\ast)\ar[r] & \ub{G'^2(F^\ast)}_0\ar[r] & G'^2(\ker\pi)\ar[r] &~
\et}\]
Therefore, it suffices to prove the sequence of functors {\small $\{G'^n:=ext^{n+1}_R(-,M)\}_{n\geq 0}$} satisfies (a),(b),(c).
\bit[leftmargin=1.2cm]
\item[(a)] It is clear that the collection $\{G'^n\}_{n\geq 0}$ satisfies (a) due to the exact ext-sequence for $R^\ast$-mod.
\item[(b)] To show $G'^0:=ext_R^1(-,M)\cong_{Ab} Hom_{R^\ast}(-,M^\ast)$, consider the $R$-exact sequence $0\ra R\sr{z}{\ral}R\ra R/Rz\ra 0$. Applying $-\otimes_RM$, we get $0\ra M\sr{z}{\ral}M\ra{M\over zM}\ra 0$, which is still $R$-exact because $z$ is non-annihilating on $M$. The following $R$-exact ext-sequence (in which the map induced by $z$ is still multiplication by $z$ since $z\in Z(R)$ ) implies ~$Hom_R(X^\ast,{M\over zM})\cong_{Ab} ext_R^1(X^\ast,M)$:
\[\adjustbox{scale=0.8}{
\bt[column sep=small]0\ar[r]& \ub{Hom_R(X^\ast,M)}_0\ar[r,"z","0"']&\ub{Hom_R(X^\ast,M)}_0\ar[r] & Hom_R(X^\ast,{M\over zM})\ar[r] & ext_R^1(X^\ast,M)\ar[r,"z","0"']& ext_R^1(X^\ast,M)\ar[r]&\cdots\et}
\]
Hence ~$Hom_{R^\ast}(X^\ast,{M\over zM})\cong_{Ab} Hom_R(X^\ast,{M\over zM})\cong_{Ab} ext_R^1(X^\ast,M)=:G^0(X^\ast)$,~ where the first isomorphism  holds because Lemma \ref{ChngRinsHom} gives an exact full additive imbedding of categories $F_\pi:R^\ast\txt{-mod}\hookrightarrow R\txt{-mod}$ induced by the surjective ring homomorphism $\pi:R\ra R^\ast,~r\mapsto r+Rz$.

\item[(c)] Next, given a free $R^\ast$-module $F^\ast:=\bigoplus_{i\in I}R^\ast$, let $F\cong_R \bigoplus_{i\in I}R$ be the free $R$-module on the same set/basis as $F^\ast$. Then ${F\over zF}\cong_R{R\over Rz}\otimes_RF=R^\ast\otimes_R(\bigoplus_iR)\cong_R\bigoplus_i R^\ast=F^\ast$. Thus, we have an exact sequence $0\ra F\sr{z}{\ral} F\ra {F\over zF}\cong_R F^\ast\ra 0$, which gives an $R$-exact Ext-sequence (for $n\geq 1$)
{\footnotesize\[
0\ra\cdots\ra \ub{Ext_R^n(F,M)}_0\ra Ext_R^{n+1}(F^\ast,M)\ra \ub{Ext_R^{n+1}(F,M)}_0\ra \ub{Ext_R^{n+1}(F,M)}_0\ra\cdots
\]}which shows ~$G'^n(F^\ast):=ext_R^{n+1}(F^\ast,M)\cong_{Ab} Ext_R^{n+1}(F^\ast,M)=0$ ~for all $n\geq 1$. \qedhere
\eit

\end{proof}

\begin{crl}[\blue{\cite[Corollary 8.35, p.472]{rotman2009}}]\label{SyzygyPreCrl3}
Let $R$ be a ring, and $z\in R\backslash 0$. If $z$ is (i) central in the sense $z\in Z(R)$, (ii) a non-zerodivisor, and (iii) a non-unit, then ~$Lpd(R)\geq Lpd(R/Rz)+1$.

(In particular, ~$Lpd(R[x])\geq Lpd(R[x]/R[x]x)+1=Lpd(R)+1$,~ since $R[x]/R[x]x\cong_{\txt{Ring}} R$.)
\end{crl}
\begin{proof}
Let $R^\ast:=R/Rz$. The case of $Lpd(R)=\infty$ is clear. So assume $Lpd(R)<\infty$. Let $M^\ast$ be an $R^\ast$-module with $n:=Pd_{R^\ast}(M^\ast)<\infty$. Then by Lemma \ref{SyzygyPreLmm3}, there exists a free module $_{R^\ast}F^\ast=\bigoplus_{i\in I}R^\ast$ such that $Ext^n_{R^\ast}(M^\ast,F^\ast)\neq 0$. Let $F\cong_R \bigoplus_{i\in I}R$ be the free $R$-module on the same set/basis as $F^\ast$. Then as before, ${F\over zF}\cong_R{R\over Rz}\otimes_RF=R^\ast\otimes_R(\bigoplus_{i\in I}R)\cong_R\bigoplus R^\ast=F^\ast$. By Theorem \ref{SyzygyThm0}, ~$0\neq Ext_{R^\ast}^n(M^\ast,F^\ast)=Ext_{R/Rz}^n(M^\ast,F/zF)\cong_{Ab} Ext_R^{n+1}(M^\ast,F)$. Therefore,
\[
Pd_R(M^\ast)\geq n+1=Pd_{R^\ast}(M^\ast)+1,~~\Ra~~Lpd(R)\geq Lpd(R^\ast)+1. \qedhere
\]
\end{proof}

\begin{lmm}[\blue{\cite[Theorem 8.36, p.472]{rotman2009}}]\label{SyzygyLmm3}
Let $R$ be a ring. Then~ $Lpd(R[x])=Lpd(R)+1$.
\end{lmm}
\begin{proof}
By Corollary \ref{SyzygyPreCrl1}, the equality holds if $Lpd(R)=\infty$. So assume $Lpd(R)<\infty$. Then $Lpd(R[x])\leq Lpd(R)+1$ by Corollary \ref{SyzygyPreCrl2}, and $Lpd(R[x])\geq Lpd(R)+1$ by Corollary \ref{SyzygyPreCrl3}.
\end{proof}

\begin{thm}[\blue{\index{Hilbert's syzygies theorem}{Hilbert's syzygies theorem}: \cite[Theorem 8.37, p.472]{rotman2009}}]\label{SyzygyThm}
Let $R$ be a ring. Then $Lpd(R[x_1,...,x_n])=Lpd(R)+n$. In particular, if $k$ is a field, then ~$Lpd(k[x_1,...,x_n])=n$.
\end{thm}
\begin{proof}
The result {\footnotesize $Lpd(R[x_1,...,x_n])=Lpd(R)+n$} follows from the recursion relation {\footnotesize $R[x_1,...,x_n]$} $=$ {\footnotesize $(R[x_1,...,x_{n-1}])[x_n]$}, Lemma \ref{SyzygyLmm3}, and induction on $n$. If $k$ is a field, then $Lpd(k)=0$, and so {\footnotesize $Lpd(k[x_1,...,x_n])=Lpd(k)+n=0+n=n$}.
\end{proof}

%% file: parts/AlgebraCat/AlgebraCatS11.tex
\chapter{K\"unneth and Universal Coefficient Theorems}\label{AlgebraCatS11}
The results in this chapter are useful in some computation-oriented applications of homological algebra such as in the description of interconnections between algebraic structures defined over topological spaces in algebraic topology/geometry (to be discussed later). Such applications are not of immediate importance to us, and so we will not give examples here. Our main interest is in giving clean proofs for a few of the major results that underly the said applications.
\section{K\"unneth theorems}
\begin{thm}[\textcolor{blue}{\index{Kunneth's theorem! for modules}{Kunneth's theorem for modules}}]\label{KunnethThm}
Let $F_\ast=(F_\ast,d_\ast)$ be a chain complex of right R-modules such that
\bit
\item[(1)] each $F_n$ is flat (e.g., projective, free) and
\item[(2)] each $B_n:=\im d_{n+1}=d_{n+1}(F_{n+1})\subset F_n$ is flat (\blue{footnote}\footnote{For example, we have already seen (using homological dimension) than when $F_n$ is a projective module over a PID, then every submodule of $F_n$, such as $B_n$ or $Z_n$, is also projective.}).
\eit
Then for each $n$, and each left R-module $M$, there is a short exact sequence
\bea
0\ra H_n(F_\ast)\otimes_R M\ra H_n(F_\ast\otimes_R M)\ra Tor_1^R\big(H_{n-1}(F_\ast),M\big)\ra 0.\nn
\eea
Furthermore, if each $B_n$ is projective then these sequences are split-exact.
\end{thm}
\begin{proof}
{\flushleft \ul{Step 1 (Showing $Z_n:=\ker d_n$ is also flat):}} Recall that the complex $F_\ast$ decomposes as follows:
\[\adjustbox{scale=0.8}{\bt
   &               &                           &                &                  &                          &       0\ar[dr]             &   0                    &          &  \\
   &               &                           &                &                  &              0\ar[r]     & B_{n-1}\ar[r,hook,"m_{n-1}"]\ar[ur] & Z_{n-1}\ar[dr,hook]\ar[r,two heads,"cm_{n-1}"]  & H_{n-1}={Z_{n-1}\over B_{n-1}}\ar[r]  & 0 \\
   & \cdots\ar[r]  & F_{n+1}\ar[rrr,"d_{n+1}"]\ar[dr,two heads,"ck_{n+1}"] &                &                  & F_n\ar[rrr,"d_n"]\ar[ur,two heads,"ck_n"] &                    &                       & F_{n-1}\ar[r]  & \cdots   \\
   &               &       0\ar[r]             & B_n\ar[r,hook,"m_n"]\ar[dr] & Z_n\ar[ur,hook,"k_n"]\ar[r,two heads,"cm_n"] & H_n={Z_n\over B_n}\ar[r] & 0                  &                       &          &         \\
   &               &                           &      0\ar[ur]  &   0               &                          &                    &                       &          &
\et}\]
where the epimorphisms $ck_n$ in the above diagram are unique up to isomorphism (i.e., up to composition with isomorphisms). Consider the SES $0\ra Z_n\sr{i_n}{\ral} F_n\sr{d_n}{\ral}B_{n-1}\ra 0$ (or equivalently, the SES $0\ra Z_\ast\sr{i_\ast}{\ral} F_\ast\sr{d_\ast}{\ral}B_\ast[-1]\ra 0$). Also consider any SES $0\ra A\ra B\ra C\ra 0$ (which we will ``\ul{Tor-tensor}'' with $0\ra Z_n\ra F_n\ra B_{n-1}\ra 0$ in three steps as follows).

First, (i) replace {\small $0\ra Z_n\ra F_n\ra B_{n-1}\ra 0$} with a \ul{split SES} of deleted projective resolutions {\small $0\ra P_\ast^{Z_n}\ra P_\ast^{F_n}\ra P_\ast^{B_{n-1}}\ra 0$}. Then (ii) tensor with {\small $0\ra A\ra B\ra C\ra 0$} to get a tricomplex with \ul{exact columns}:
\[\adjustbox{scale=0.8}{\bt
  & 0\ar[d] & 0\ar[d] & 0\ar[d] & \\
 0\ar[r,dotted] & A\otimes P_\ast^{Z_n}\ar[d]\ar[r] & B\otimes P_\ast^{Z_n}\ar[d]\ar[r] & C\otimes P_\ast^{Z_n}\ar[d]\ar[r,dotted] & 0\\
 0\ar[r] & A\otimes P_\ast^{F_n}\ar[d]\ar[r] & B\otimes P_\ast^{F_n}\ar[d]\ar[r] & C\otimes P_\ast^{F_n}\ar[d]\ar[r] & 0\\
 0\ar[r] & A\otimes P_\ast^{B_{n-1}}\ar[d]\ar[r] & B\otimes P_\ast^{B_{n-1}}\ar[d]\ar[r] & C\otimes P_\ast^{B_{n-1}}\ar[d]\ar[r] & 0\\
  & 0 & 0 & 0 &
\et}\]
Next, (iii) take the induced LES's of homologies for the columns (using naturality/functoriality of the induced LES of homologies) to get the following bicomplex with \ul{exact columns}:
\[\adjustbox{scale=0.8}{\bt
  & \vdots\ar[d] &\vdots\ar[d] & \vdots\ar[d] & \\
  & \overbrace{Tor_1(A,B_{n-1})}^0\ar[d] & \overbrace{Tor_1(B,B_{n-1})}^0\ar[d] & \overbrace{Tor_1(C,B_{n-1})}^0\ar[d] & \\
 0\ar[r,dotted] & A\otimes Z_n\ar[d]\ar[r] & B\otimes Z_n\ar[d]\ar[r] & C\otimes Z_n\ar[d]\ar[r,dotted] & 0\\
 0\ar[r] & A\otimes F_n\ar[d]\ar[r] & B\otimes F_n\ar[d]\ar[r] & C\otimes F_n\ar[d]\ar[r] & 0\\
 0\ar[r] & A\otimes B_{n-1}\ar[d]\ar[r] & B\otimes B_{n-1}\ar[d]\ar[r] & C\otimes B_{n-1}\ar[d]\ar[r] & 0\\
  & 0 & 0 & 0 &
\et}\]
where exactness of the two bottom rows (due to flatness of $F_n$,$B_{n-1}$) implies exactness of the top row (by the 3-by-3 lemma). This shows \ul{$Z_n$ is also flat}.

{\flushleft \ul{Step 2 (Formal expression for the desired SES):}} Since $F_n,B_{n-1}$ are flat, the Tor-LES for the SES $0\ra Z_n\sr{i_n}{\ral} F_n\sr{d_n}{\ral}B_{n-1}\ra 0$ gives the exact sequence
\bea
\cdots\ra \overbrace{Tor^R_1(B_{n-1},M)}^0\ra Z_n\otimes_R M\ra F_n\otimes_RM \ra B_{n-1}\otimes_RM\ra 0,\nn
\eea
and so we get the following SES of chain complexes (in which $id=id_M$ is also written as $1=1_M$):
\bea
\bt 0\ar[r] & Z_\ast\otimes_RM\ar[r,"i_\ast\otimes id"] & F_\ast\otimes_RM\ar[r,"d_\ast\otimes id"] & {(B_\ast\otimes_RM)[-1]}\ar[r] & 0\et\nn
\eea
with differentials induced by those of $F_\ast$. Note that we have
\bea
d^{Z_\ast\otimes M}=0=d^{B_\ast\otimes M}~~~~(\txt{since}~~d^{Z_\ast}=0=d^{B_\ast}).\nn
\eea
The induced LES of homologies for the above SES of chain complexes is
\[\adjustbox{scale=0.8}{\bt[column sep =small]
             &            &  & \substack{\coker\delta_{n+1}\\=\coim f_n\\\cong\im f_n=\ker g_n}\ar[d] &  &  & \\
\cdots\ar[r] & H_{n+1}\big((B_\ast\otimes M)[-1]\big)\ar[d,equal]\ar[r,"\delta_{n+1}"]\ar[r] & H_n(Z_\ast\otimes M)\ar[ur]\ar[dd,equal,"(b)"]\ar[r,"f_n"] & H_n(F_\ast\otimes M)\ar[d]\ar[r,"g_n"] & H_n\big((B_\ast\otimes M)[-1]\big)\ar[dd,equal,"(c)"]\ar[r,"\delta_n"] & H_{n-1}(Z_\ast\otimes M)\ar[dd,equal,"(d)"]\ar[r] & \cdots\\
  & H_n(B_\ast\otimes M)\ar[d,equal,"(a)"] &  & \substack{\coker f_n\\=\coim g_n\\\cong\im g_n=\ker\delta_n}\ar[ur] &  &  & \\
\cdots\ar[r]  & B_n\otimes M\ar[r,"\delta_{n+1}"] & Z_n\otimes M\ar[r,"f_n"] & H_n(F_\ast\otimes M)\ar[r,"g_n"] & B_{n-1}\otimes M\ar[r,"\delta_n"] & Z_{n-1}\otimes M\ar[r] & \cdots\\
\et}\]
where the equalities (a)-(d) hold because the differentials of the complexes involved are all zero. So, from the SES ~$0\ra \im f_n\ra\cod f_n\ra\coker f_n\ra 0$,~ we get the SES
\bea
0\ra \coker \delta_{n+1}\sr{f_n}{\ral} H_n(F_\ast\otimes M)\sr{g_n}{\ral}\ker\delta_n\ra 0,\nn
\eea
where $f_n:=H_n(i_\ast\otimes id_M)\cong H_n(i_\ast)\otimes id_M$ and $g_n:=H_n(d_\ast\otimes id_M)$ are each unique up to isomorphism (i.e., up to composition with isomorphisms).

{\flushleft \ul{Step 3 (Computing the connecting homomorphism $\delta_\ast$ above):}} To get $\delta_\ast$, consider the diagram
\[\adjustbox{scale=0.9}{\bt
         &  \vdots\ar[d]  &  \vdots\ar[d]  &  \vdots\ar[d]  &   \\
 0\ar[r] &  Z_{n+1}\otimes M\ar[d]\ar[r,"i_{n+1}\otimes 1"]  &  F_{n+1}\otimes M\ar[d]\ar[r,"d_{n+1}\otimes 1"]   &  B_n\otimes M\ar[d]\ar[r]  &  0 \\
 0\ar[r] &  Z_n\otimes M\ar[d]\ar[r,"i_n\otimes 1"]  &  F_n\otimes M\ar[d,dashed,"d_n\otimes 1"]\ar[r,dashed,"d_n\otimes 1"]   &  B_{n-1}\ar[d]\otimes M\ar[r,dashed]  &  0 \\
 0\ar[r,dashed] &  Z_{n-1}\otimes M\ar[d]\ar[r,dashed,"i_{n-1}\otimes 1"]  &  F_{n-1}\otimes M\ar[d]\ar[r,"d_{n-1}\otimes 1"]   &  B_{n-2}\otimes M\ar[d]\ar[r]  &  0 \\
         &  \vdots  &  \vdots  &  \vdots  &
\et}\]
Then $\delta_n\big(b_{n-1}\otimes m+\im(d_n\otimes 1)\big)=z_{n-1}\otimes m+\im(d_n\otimes 1)$, where (with $p_n\otimes m\in F_n\otimes M$)
\bea
&&z_{n-1}\otimes m=(i_{n-1}\otimes 1)(z_{n-1}\otimes m)=(d_n\otimes 1)(p_n\otimes m),~~\txt{along with}~~(d_n\otimes 1)(p_n\otimes m)=b_{n-1}\otimes m,\nn\\
&&~~\Ra~~z_{n-1}\otimes m=b_{n-1}\otimes m.\nn
\eea
Therefore, $\delta_n=i_{n-1}\otimes 1$, where $i_{n-1}$ is \ul{now} the inclusion $B_{n-1}\subset Z_{n-1}$, for every $n$.

{\flushleft \ul{Step 4 (Computing the end terms of the SES):}} With $\delta_{n+1}=i_n\otimes 1$, applying the right exact functor $-\otimes_RM$ on the SES ~$0\ra B_n\sr{i_n}{\ral}Z_n\sr{\pi_n}{\ral}{Z_n\over B_n}=H_n\ra 0$,~ we get the exact sequence
\bea\bt
B_n\otimes M\ar[rr,"i_n\otimes id=\delta_{n+1}"] &&Z_n\otimes M \ar[rr,"\pi_n\otimes id"] && {Z_n\over B_n}\otimes M\ar[r] & 0,
\et\nn\eea
and so
\bea
&&\textstyle {Z_n\over B_n}\otimes M=\im(\pi_n\otimes id)\cong \coim(\pi_n\otimes id)={Z_n\otimes M\over\ker(\pi_n\otimes id)}={Z_n\otimes M\over\im(i_n\otimes id)}={Z_n\otimes M\over\im\delta_{n+1}}=\coker\delta_{n+1},\nn\\
&&\textstyle~~\Ra~~\coker\delta_{n+1}\cong H_n(F_\ast)\otimes M.\nn
\eea
Finally, the exact Tor-sequence for the SES
\bea
0\ra B_{n-1}\sr{i_{n-1}}{\ral}Z_{n-1}\ra H_{n-1}\ra 0~~~~\txt{( a flat resolution of $H_{n-1}=H_{n-1}(F_\ast)$ )}\nn
\eea
gives the exact sequence
\bea
&&\cdots \ra \overbrace{Tor_1^R(B_{n-1},M)}^0\ra \overbrace{Tor_1^R(Z_{n-1},M)}^0\ra Tor_1^R(H_{n-1},M)\ra\nn\\
&&~~~~\hspace{3cm}\ra B_{n-1}\otimes M\sr{i_{n-1}\otimes 1}{\ral}Z_{n-1}\otimes M\ra H_{n-1}\otimes M\ra 0,\nn
\eea
which shows ~$\ker \delta_n=\ker(i_{n-1}\otimes 1)\cong  Tor_1^R(H_{n-1},M)$.

{\flushleft \ul{Step 5 (Showing the sequences are split if each $B_n$ is projective):}} See the proof of Theorem \ref{UCTfhom}.
\end{proof}

\begin{thm}[\textcolor{blue}{\index{Kunneth's theorem! for chain complexes}{Kunneth's theorem for chain complexes}}]\label{KunnethThmC}
Let $F_\ast=(F_\ast,d_\ast)$ be a chain complex of right $R$-modules such that
\bit
\item[(1)] each $F_n$ is flat and
\item[(2)] each $B_n:=\im d_{n+1}\subset F_n$ is flat.
\eit
Then for each $n$, and for each chain complex of left $R$-modules $M_\ast$, there exists a short exact sequence
\bea
\textstyle 0\ral\bigoplus\limits_{i+j=n}H_i(F_\ast)\otimes H_j(M_\ast)\ral H_n\big(Tot(F_\ast\otimes M_\ast)\big)\ral \bigoplus\limits_{i+j=n-1}Tor_1\big(H_i(F_\ast),H_j(M_\ast)\big)\ral 0.\nn
\eea

\ul{\textbf{Remark}}: Taking $M_\ast:=(\cdots\ra 0\ra M\ra 0\ra \cdots)$, we get Kunneth's theorem for modules.
\end{thm}
\begin{proof}
We will merely follow the same steps as in the proof of Theorem \ref{KunnethThm}, with appropriate adjustments at some steps.

{\flushleft \ul{Step 1 (Showing $Z_n:=\ker d_n$ is also flat):}} This is the same as in the proof of Theorem \ref{KunnethThm}.

{\flushleft \ul{Step 2 (Formal expression for the desired SES):}} Since $F_n,B_{n-1}$ are flat, the Tor-LES for the SES $0\ra Z_n\ra F_n\sr{d_n}{\ral}B_{n-1}\ra 0$ gives the exact sequence
\bea
\cdots\ra \overbrace{Tor_1^R(B_{n-1},M_\ast)}^0\ra Z_n\otimes_R M_\ast\ra F_n\otimes_RM_\ast \ra B_{n-1}\otimes_RM_\ast\ra 0,\nn
\eea
and so we get the following SES of bicomplexes (in which $1_\ast=1_{M_\ast}$)
\bea
0\ra Z_\ast\otimes_RM_\ast\sr{i_\ast\otimes 1_\ast}{\ral}F_\ast\otimes_RM_\ast\sr{d_\ast\otimes 1_\ast}{\ral}(B_\ast\otimes_RM_\ast)[-1]\ra 0\nn
\eea
with differentials induced by those of $F_\ast$ and $M_\ast$. Taking Tot of the above SES of bicomplexes, we get the following SES (due to the exactness of Tot on $R$-mod):
\bea
0\ra Tot(Z_\ast\otimes_RM_\ast)\sr{Tot(i_\ast\otimes 1_\ast)}{\ral}Tot(F_\ast\otimes_RM_\ast)\sr{Tot(d_\ast\otimes 1_\ast)}{\ral}Tot\big((B_\ast\otimes_RM_\ast)[-1]\big)\ra 0,\nn
\eea
where as usual ~$Tot\big((B_\ast\otimes_RM_\ast)[-1]\big)=Tot\big(B_\ast\otimes_RM_\ast\big)[-1]$.~ Since $d^{Z_\ast}=0=d^{B_\ast}$, we have
\begin{align}
&\textstyle d_n^{Tot(Z_\ast\otimes M_\ast)}=\bigoplus\limits_{i+j=n}\left(d^h_{ij}+(-1)^id^v_{ij}\right)=\bigoplus\limits_{i+j=n}\left(d^{Z_\ast}_i\otimes id_{M_j}+(-1)^iid_{Z_i}\otimes d^{M_\ast}_j\right)=\bigoplus\limits_{i+j=n}(-1)^iid_{Z_i}\otimes d^{M_\ast}_j\nn\\
&\textstyle
d_n^{Tot(B_\ast\otimes M_\ast)}\sr{\txt{similarly}}{=}\bigoplus\limits_{i+j=n}(-1)^iid_{B_i}\otimes d^{M_\ast}_j=\bigoplus\limits_{j}(-1)^{n-j}id_{B_{n-j}}\otimes d^{M_\ast}_j,\nn
\end{align}
using which we can then obtain the following simplifications:
\begin{align}
&\textstyle H_n(Tot(Z_\ast\otimes M_\ast))={\ker d_n^{Tot(Z_\ast\otimes M_\ast)}\over\im d_{n+1}^{Tot(Z_\ast\otimes M_\ast)}}={\ker \bigoplus\limits_{i+j=n}(-1)^iid_{Z_i}\otimes d^{M_\ast}_j\over\im \bigoplus\limits_{i+j=n+1}(-1)^iid_{Z_i}\otimes d^{M_\ast}_j}\sr{(s1)}{=}
{\bigoplus\limits_{i+j=n}(-1)^i\ker(id_{Z_i}\otimes d^{M_\ast}_j)\over \bigoplus\limits_{i+j=n+1}(-1)^i\im(id_{Z_i}\otimes d^{M_\ast}_j)}\nn\\
&\textstyle~~~~\sr{(s2)}{=} {\bigoplus\limits_{i+j=n}Z_i\otimes \ker d^{M_\ast}_j\over \bigoplus\limits_{i+j=n+1}Z_i\otimes \im d^{M_\ast}_j}=
{\bigoplus\limits_iZ_i\otimes \ker d^{M_\ast}_{n-i}\over \bigoplus\limits_iZ_i\otimes \im d^{M_\ast}_{n-i+1}}\cong
\bigoplus\limits_iZ_i\otimes {\ker d^{M_\ast}_{n-i}\over\im d^{M_\ast}_{n-i+1}}=\bigoplus\limits_iZ_i\otimes H_{n-i}(M_\ast)\nn\\
&\textstyle~~~~=\bigoplus\limits_{i+j=n}Z_i\otimes H_j(M_\ast)=Tot(Z_\ast\otimes H_\ast(M_\ast)),\nn\\
&\textstyle H_n(Tot(B_\ast\otimes M_\ast))\sr{\txt{similarly}}{\cong}\bigoplus\limits_{i+j=n}B_i\otimes H_j(M_\ast)=Tot(B_\ast\otimes H_\ast(M_\ast)),\nn
\end{align}
where step (s1) is clear by exactness of the colimit $\bigoplus$ on $R$-mod, and step (s2) is due to flatness of $Z_i$'s.

The induced LES of homologies for the above SES of chain complexes is
\[\adjustbox{scale=0.7}{\bt[column sep =small]
             &            &  & \substack{\coker\delta_{n+1}\\=\coim f_n\\\cong\im f_n=\ker g_n}\ar[d] &  &  & \\
\cdots\ar[r] & H_{n+1}\big(Tot(B_\ast\otimes M_\ast)[-1]\big)\ar[d,equal]\ar[r,"\delta_{n+1}"]\ar[r] & H_n(Tot(Z_\ast\otimes M_\ast))\ar[ur]\ar[dd,equal,"(b)"]\ar[r,"f_n"] & H_n(Tot(F_\ast\otimes M_\ast))\ar[d]\ar[r,"g_n"] & H_n\big(Tot(B_\ast\otimes M_\ast)[-1]\big)\ar[dd,equal,"(c)"]\ar[r,"\delta_n"] & H_{n-1}(Tot(Z_\ast\otimes M_\ast))\ar[dd,equal,"(d)"]\ar[r] & \cdots\\
  & H_n(Tot(B_\ast\otimes M_\ast))\ar[d,equal,"(a)"] &  & \substack{\coker f_n\\=\coim g_n\\\cong\im g_n=\ker\delta_n}\ar[ur] &  &  & \\
\cdots\ar[r]  & Tot(B_\ast\otimes H_\ast(M_\ast))_n\ar[r,"\delta_{n+1}"] & Tot(Z_\ast\otimes H_\ast(M_\ast))_n\ar[r,"f_n"] & H_n(Tot(F_\ast\otimes M_\ast))\ar[r,"g_n"] & Tot(B_\ast\otimes H_\ast(M_\ast))_{n-1}\ar[r,"\delta_n"] & Tot(Z_\ast\otimes H_\ast(M_\ast))_{n-1}\ar[r] & \cdots\\
\et}\]
where the equalities (a)-(d) hold (as shown earlier) because $d^{Z_\ast}=0=d^{B_\ast}$ and Tot is translation-invariant. So, from the SES ~$0\ra \im f_n\ra\cod f_n\ra\coker f_n\ra 0$,~ we get the SES
\bea
0\ra \coker \delta_{n+1}\sr{f_n}{\ral} H_n(Tot(F_\ast\otimes M_\ast))\sr{g_n}{\ral}\ker\delta_n\ra 0.\nn
\eea
where $f_n:=H_n\big(Tot(i_\ast\otimes id_{M_\ast})\big)\cong Tot\left(H_\ast(i_\ast)\otimes id_{H_\ast(M_\ast)}\right)_n$ and $g_n:=H_n\big(Tot(d_\ast\otimes id_{M_\ast})\big)$ are each unique up to isomorphism (i.e., up to composition with isomorphisms).

{\flushleft \ul{Step 3 (Computing the connecting homomorphism $\delta=\delta_\ast$ above):}} To get $\delta$, consider the diagram
\[\adjustbox{scale=0.8}{\bt
         &  \vdots\ar[d]  &&  \vdots\ar[d]  &&  \vdots\ar[d]  &   \\
 0\ar[r] &  Tot(Z_\ast\otimes H_\ast(M_\ast))_{n+1}\ar[d]\ar[rr,"Tot(i_\ast\otimes 1_\ast)_{n+1}"]  &&  Tot(F_\ast\otimes H_\ast(M_\ast))_{n+1}\ar[d]\ar[rr,"Tot(d_\ast\otimes 1_\ast)_{n+1}"]   &&  Tot(B_\ast\otimes H_\ast(M_\ast))_n\ar[d]\ar[r]  &  0 \\
 0\ar[r] &  Tot(Z_\ast\otimes H_\ast(M_\ast))_n\ar[d]\ar[rr,"Tot(i_\ast\otimes 1_\ast)_n"]  &&  Tot(F_\ast\otimes H_\ast(M_\ast))_n\ar[d,dashed,"Tot(d_\ast\otimes 1_\ast)_n"]\ar[rr,dashed,"Tot(d_\ast\otimes 1_\ast)_n"]   &&  Tot(B_\ast\otimes H_\ast(M_\ast))_{n-1}\ar[d]\ar[r,dashed]  &  0 \\
 0\ar[r,dashed] &  Tot(Z_\ast\otimes H_\ast(M_\ast))_{n-1}\ar[d]\ar[rr,dashed,"Tot(i_\ast\otimes 1_\ast)_{n-1}"]  &&  Tot(F_\ast\otimes H_\ast(M_\ast))_{n-1}\ar[d]\ar[rr,"Tot(d_\ast\otimes 1_\ast)_{n-1}"]   && Tot(B_\ast\otimes H_\ast(M_\ast))_{n-2}\ar[d]\ar[r]  &  0 \\
         &  \vdots  &&  \vdots  &&  \vdots  &
\et}\]
where {\footnotesize $Tot(i_\ast\otimes 1_\ast)_n=\bigoplus\limits_{i+j=n}i_i\otimes id_{H_j(M_\ast)}$, $Tot(d_\ast\otimes 1_\ast)_n=\bigoplus\limits_{i+j=n}d_i\otimes id_{H_j(M_\ast)}$, $Tot(F_\ast\otimes H_\ast(M_\ast))_n=\bigoplus\limits_{i+j=n}F_i\otimes H_j(M_\ast)$}, etc.

Then ~$\delta_n\big((b_i\otimes h_j)_{i+j=n-1}+\im Tot(d_\ast\otimes 1_\ast)\big)_n=(z_i\otimes h_j)_{i+j=n-1}+\im Tot(d_\ast\otimes 1_\ast)_n$, \\ where ( with $(p_i\otimes h_j)_{i+j=n}\in Tot(F_\ast\otimes H_\ast(M_\ast))_n$ )
\bea
&&\textstyle (z_i\otimes h_j)_{i+j=n-1}=Tot(i_\ast\otimes 1_\ast)_{n-1}(z_i\otimes h_j)_{i+j=n-1}=Tot(d_\ast\otimes 1_\ast)_n~(p_i\otimes h_j)_{i+j=n},\nn\\
&&\textstyle~~~~~~~~\txt{along with}~~Tot(d_\ast\otimes 1_\ast)_n~(p_i\otimes h_j)_{i+j=n}=(b_i\otimes h_j)_{i+j=n-1},\nn\\
&&\textstyle~~\Ra~~(z_i\otimes h_j)_{i+j=n-1}=(b_i\otimes h_j)_{i+j=n-1}.\nn
\eea
Therefore, {\small $\delta_n=Tot(i_\ast\otimes 1_\ast)_{n-1}=\bigoplus_{i+j=n-1}i_i\otimes id_{H_j(M_\ast)}$}, where $i_i$ is \ul{now} the inclusion $B_i\subset Z_i$, for every $i$.

{\flushleft (iv)} \ul{Step 4 (Computing the end terms of the SES):}
With $\delta_{n+1}=Tot(i_\ast\otimes 1_\ast)_n$, applying the right exact functor $-\otimes_RH_j(M_\ast)$ on the SES ~$0\ra B_i\sr{i_i}{\ral}Z_i\sr{\pi_i}{\ral}{Z_i\over B_i}=H_i(F_\ast)\ra 0$,~ we get the exact sequence
\bea\bt
B_i\otimes H_j(M_\ast)\ar[rr,"i_i\otimes id"] && Z_i\otimes H_j(M_\ast)\ar[rr,"\pi_i\otimes id"] && {Z_i\over B_i}\otimes H_j(M_\ast)\ar[r] & 0,
\et\nn
\eea
whose Tot (\blue{footnote}\footnote{Recall that Tot is exact on $R$-mod.}) gives the exact sequence
\[\adjustbox{scale=0.9}{\bt
Tot(B_\ast\otimes H_\ast(M_\ast))_n\ar[rrr,"Tot(i_\ast\otimes id)_n=\delta_{n+1}"] &&& Tot(Z_\ast\otimes H_\ast(M_\ast))_n\ar[rr,"Tot(\pi\otimes id)"] && Tot\left({Z_\ast\over B_\ast}\otimes H_\ast(M_\ast)\right)_n\ar[r] & 0,
\et}\]
and so
\bea
&&\textstyle\coker\delta_{n+1}={\cod\delta_{n+1}\over\im\delta_{n+1}}={Tot(Z_\ast\otimes H_\ast(M_\ast))_n\over Tot(i_\ast\otimes 1_\ast)_{n+1}\big(Tot(B_\ast\otimes H_\ast(M_\ast))_n\big)}={Tot(Z_\ast\otimes H_\ast(M_\ast))_n\over Tot(B_\ast\otimes H_\ast(M_\ast))_n}\sr{(s1)}{=}Tot\left({Z_\ast\otimes H_\ast(M_\ast)\over B_\ast\otimes H_\ast(M_\ast)}\right)_n\nn\\
&&\textstyle~~~~\sr{(s2)}{=}Tot\left({Z_\ast\over B_\ast}\otimes H_\ast(M_\ast)\right)=Tot\big(H_\ast(F_\ast)\otimes H_\ast(M_\ast)\big),\nn
\eea
where step (s1) is due to the exatness of $\bigoplus$ on $R$-mod, and step (s2) is obtained as follows:
\bea
\textstyle Tot\left({Z_\ast\otimes H_\ast(M_\ast)\over B_\ast\otimes H_\ast(M_\ast)}\right)_n:=\bigoplus\limits_{i+j=n}{Z_i\otimes H_j(M_\ast)\over B_i\otimes H_j(M_\ast)}\sr{(s3)}{\cong} \bigoplus\limits_{i+j=n}{Z_i\over B_i}\otimes H_j(M_\ast)=\bigoplus\limits_{i+j=n}H_i(F_\ast)\otimes H_j(M_\ast),\nn
\eea
with step (s3) obtained as follows:
\bea
\textstyle {Z_i\over B_i}\otimes H_j(M_\ast)=\im(\pi_i\otimes id)\cong \coim(\pi_i\otimes id)={Z_i\otimes H_j(M_\ast)\over\ker(\pi_i\otimes id)}={Z_i\otimes H_j(M_\ast)\over\im(i_i\otimes id)}={Z_i\otimes H_j(M_\ast)\over B_i\otimes H_j(M_\ast)}.\nn
\eea
Finally, the Tor sequence for the SES
\bea
0\ra B_{n-1}\sr{i_{n-1}}{\ral}Z_{n-1}\ra H_{n-1}(F_\ast)\ra 0~~~~\txt{(a flat resolution of $H_{n-1}(F_\ast)$)}\nn
\eea
gives the exact sequence
{\small\bea
&&\cdots \ra \overbrace{Tor_1^R(Tot(Z_\ast,H_\ast(M_\ast))_{n-1})}^0\ra Tor_1^R(Tot(H_\ast(F_\ast),H_\ast(M_\ast))_{n-1})\ra\nn\\
&&~~~~\hspace{1cm}\ra Tot(B_\ast\otimes H_\ast(M_\ast))_{n-1}\sr{Tot(i_\ast\otimes 1_\ast)_{n-1}}{\ral}Tot(Z_\ast\otimes H_\ast(M_\ast))_{n-1}\ra Tot(H_\ast(F_\ast)\otimes H_\ast(M_\ast))_{n-1}\ra 0,\nn
\eea}which shows ~$\ker \delta_n=\ker(Tot(i_\ast\otimes 1_\ast)_{n-1})\cong  Tor_1^R(Tot(H_\ast(F_\ast),H_\ast(M_\ast))_{n-1})$.
\end{proof}

\section{Universal coefficient theorems}
We will write UCT for the phrase ``\emph{universal coefficient theorem}''.
\begin{thm}[\textcolor{blue}{\index{Universal coefficient theorem! for homology}{UCT for homology}}]\label{UCTfhom}
Let $F_\ast=(F_\ast,d_\ast)$ be a chain complex of right $R$-modules such that
\bit
\item[(1)] each $F_n$ is flat (e.g., projective, free) and
\item[(2)] each $B_n:=\im d_{n+1}=d_{n+1}(F_{n+1})\subset F_n$ is projective (\blue{footnote}\footnote{For example, when each $F_n$ is a projective module over a PID, then every submodule, such as $B_n$ or $Z_n$, is also projective.}).
\eit
Then (by K\"unneth's theorem) for each $n$, and for each left $R$-module $M$, there is a short exact sequence
\bea
0\ra H_n(F_\ast)\otimes_R M\ra H_n(F_\ast\otimes_R M)\ra Tor_1^R\big(H_{n-1}(F_\ast),M\big)\ra 0.\nn
\eea
Moreover, these sequences are split. (\blue{footnote}\footnote{\magenta{Caution}: It is known that the splitting process here, $\A_0^\Integer\ra\A_0^\Integer$, is not canonically functorial (in the sense that it might not automatically map a given morphism of SES's $(A,B,C)\sr{(\al,\beta,\gamma)}{\ral}(A',B',C')$ to a split-morphism of split-SES's $(A,A\oplus C,C)\sr{(\al,\beta_1\oplus\beta_2,\gamma)}{\ral}(A',A'\oplus C',C')$).})
\end{thm}
\begin{proof}
We only need to show the above sequences are split. Since $B_{n-1}$ is projective, the sequence
\bc\bt
0\ar[r] &Z_n\ar[r,"i_n"] & F_n\ar[l,dashed,bend left=30,"h_n"]\ar[r,"d_n"] & B_{n-1}\ar[r] & 0~~~~\big(~\txt{with}~h_n\circ i_n=id_{Z_n}~\big)
\et\ec
is split, say via a map $h_n$ as above. Hence, applying the additive functor $-\otimes_R M$, we get the split SES
\bc\bt
0\ar[r] & Z_n\otimes M\ar[r,"i_n\otimes 1"] & F_n\otimes M\ar[l,dashed,bend left=30,"h_n\otimes 1"]\ar[r,"d_n\otimes 1"] & B_{n-1}\otimes M\ar[r] & 0~~~~\big(~\txt{with}~(h_n\otimes 1)\circ(i_n\otimes 1)=id_{Z_n\otimes M}~\big).
\et\ec
Since $(d_n\otimes 1)(Z_n\otimes M)=d_n(Z_n)\otimes M=0$ implies $Z_n\otimes M\subset\ker(d_n\otimes 1)\subset F_n\otimes M$, by restricting the splitting map $h_n\otimes 1$ to $\ker(d_n\otimes 1)$, which gives the split SES
\bc\bt
0\ar[r] & Z_n\otimes M\ar[r,"i_n\otimes 1"] & \ker(d_n\otimes 1)\ar[l,dashed,bend left=30,"h_n\otimes 1"]\ar[r,"d_n\otimes 1"] & {\ker(d_n\otimes 1)\over Z_n\otimes M}\ar[r] & 0~~~~\big(~\txt{with}~(h_n\otimes 1)\circ(i_n\otimes 1)=id_{Z_n\otimes M}~\big),
\et\ec
we see that $Z_n\otimes M$ is a direct summand of $\ker(d_n\otimes 1)$. This in turn implies, via the exact sequences 
\bc\bt[row sep=tiny]
0 \ar[r] & B_n \ar[r,"i_n'"] & Z_n\ar[r,"\pi_n"] &  H_n:={Z_n\over B_n} \ar[r] & 0,\\
 & B_n\otimes M \ar[r,"i_n'\otimes 1"] & Z_n\otimes M\ar[r,"\pi_n\otimes 1"] & H_n\otimes M\ar[r] & 0,
\et\ec
that the quotient module
\bea
\textstyle {Z_n\otimes M\over (i'_n\otimes 1)(B_n\otimes M)}={Z_n\otimes M\over B_n\otimes M}\cong{Z_n\over B_n}\otimes M=H_n(F_\ast)\otimes M\nn
\eea
is a direct summand of the quotient module
\bea
\textstyle {\ker(d_n\otimes 1)\over(i'_n\otimes 1)(B_n\otimes M)}={\ker(d_n\otimes 1)\over B_n\otimes M}={\ker(d_n\otimes 1)\over d_{n+1}F_{n+1}\otimes M}={\ker(d_n\otimes 1)\over\im(d_{n+1}\otimes 1)}=H_n(F_\ast\otimes M).\nn\qedhere
\eea
\end{proof}

We remark that some indices in the following theorem (and its proof) are usually superscripts. This is only necessary in applications where one abbreviates (i.e., one is less explicit with) math expressions and so needs to use upper indices to distinguish cohomological quantities from homological ones. However, because such a change in notation can be needlessly distracting at this advanced stage in our progress, we will be explicit and so not be needing superscript indices.

\begin{thm}[\textcolor{blue}{\index{Universal coefficient theorem! for cohomology}{UCT for Cohomology}}]\label{UCTfcohom}
Let $P_\ast=(P_\ast,d_\ast)$ be a chain complex such that
\bit
\item[(1)] each $P_n$ is projective, (i.e., $Hom(P_n,-)$ is exact for each $n$) and
\item[(2)] each $B_n:=\im d_{n+1}\subset P_n$ is projective. (i.e., $Hom(B_n,-)$ is exact for each $n$)
\eit
Then for each $n$, and for each left $R$-module, there exists a short exact sequence
\bea
0\ra Ext^1(H_{n-1}(P_\ast),M)\ra H_n\big(Hom(P_\ast,M)\big)\ra Hom\big(H_n(P_\ast),M\big)\ra 0\nn
\eea
which is split.  (\blue{footnote}\footnote{\magenta{Caution}: As with the UCT for homology, it is known that the splitting process here, $\A_0^\Integer\ra\A_0^\Integer$, is not canonically functorial.})
\end{thm}
\begin{proof}
This will be similar to the proofs of Theorem \ref{KunnethThm} and Theorem \ref{UCTfhom} (UCT for homology), with $Hom_R(-,M)$ now playing the role of $-\otimes_RM$.

{\flushleft \ul{Step 1 (Showing $Z_n:=\ker d_n$ is also projective):}}

Recall that the complex $P_\ast$ decomposes as follows:
\[\adjustbox{scale=0.8}{\bt
   &               &                           &                &                  &                          &       0\ar[dr]             &   0                    &          &  \\
   &               &                           &                &                  &              0\ar[r]     & B_{n-1}\ar[r,hook,"m_{n-1}"]\ar[ur] & Z_{n-1}\ar[dr,hook]\ar[r,two heads,"cm_{n-1}"]  & H_{n-1}={Z_{n-1}\over B_{n-1}}\ar[r]  & 0 \\
   & \cdots\ar[r]  & P_{n+1}\ar[rrr,"d_{n+1}"]\ar[dr,two heads,"ck_{n+1}"] &                &                  & P_n\ar[rrr,"d_n"]\ar[ur,two heads,"ck_n"] &                    &                       & P_{n-1}\ar[r]  & \cdots   \\
   &               &       0\ar[r]             & B_n\ar[r,hook,"m_n"]\ar[dr] & Z_n\ar[ur,hook,"k_n"]\ar[r,two heads,"cm_n"] & H_n={Z_n\over B_n}\ar[r] & 0                  &                       &          &         \\
   &               &                           &      0\ar[ur]  &   0               &                          &                    &                       &          &
\et}\]
where the epimorphisms $ck_n$ in the above diagram are unique up to isomorphism (i.e., up to composition with isomorphisms). Consider the SES $0\ra Z_n\sr{i_n}{\ral} P_n\sr{d_n}{\ral}B_{n-1}\ra 0$ (or equivalently, the SES $0\ra Z_\ast\sr{i_\ast}{\ral} P_\ast\sr{d_\ast}{\ral}B_\ast[-1]\ra 0$). Since $B_{n-1}$ is projective, the SES splits, and so $Z_n$ is projective as a direct summand of the projective module $P_n$.

{\flushleft \ul{Step 2 (Formal expression for the desired SES):}} Since $P_n,B_{n-1}$ are projective, the ext-LES for the SES $0\ra Z_n\sr{i_n}{\ral} P_n\sr{d_n}{\ral}B_{n-1}\ra 0$ gives the exact sequence
\bea
0\ra Hom(B_{n-1},M)\ra Hom(P_n,M)\ra Hom(Z_n,M)\ra \overbrace{ext_R^1(B_{n-1},M)}^0\ra\cdots,\nn
\eea
and so we get the SES of chain complexes ( with $Hom(-,M)$ written as $h(-,M)$ on maps )
\bea
0\ra Hom(B_\ast,M)[-1] \sr{h(d_\ast,M)}{\ral} Hom(P_\ast,M) \sr{h(i_\ast,M)}{\ral} Hom(Z_\ast,M)\ra 0\nn
\eea
with differentials induced by those of $P_\ast$. Since $d^{Z_\ast}=0=d^{B_\ast}$, we have
\bea
d^{Hom(Z_\ast,M)}=0=d^{Hom(B_\ast,M)}.\nn
\eea
The induced LES of homologies for the above SES of chain complexes is

\[\adjustbox{scale=0.7}{\bt[column sep =small]
             &            &  & \substack{\coker\delta_{n-1}\\=\coim f_n\\\cong\im f_n=\ker g_n}\ar[d] &  &  & \\
\cdots\ar[r] & H_{n-1}(Hom(Z_\ast,M))\ar[dd,equal,"(a)"]\ar[r,"\delta_{n-1}"]\ar[r] & H_n(Hom(B_\ast,M)[-1])\ar[ur]\ar[d,equal]\ar[r,"f_n"] & H_n(Hom(P_\ast,M))\ar[d]\ar[r,"g_n"] & H_n\big(Hom(Z_\ast,M)\big)\ar[dd,equal,"(c)"]\ar[r,"\delta_n"] & H_{n+1}\big(Hom(B_\ast,M)[-1]\big) \ar[dd,equal,"(d)"]\ar[r] & \cdots\\
  &  & H_{n-1}(Hom(B_\ast,M))\ar[d,equal,"(b)"] & \substack{\coker f_n\\=\coim g_n\\\cong\im g_n=\ker\delta_n}\ar[ur] &  &  & \\
\cdots\ar[r]  & Hom(Z_{n-1},M)\ar[r,"\delta_{n-1}"] & Hom(B_{n-1},M)\ar[r,"f_n"] & H_n(Hom(P_\ast,M))\ar[r,"g_n"] & Hom(Z_n,M)\ar[r,"\delta_n"] & Hom(B_n,M)\ar[r] &\cdots \\
\et}\]
where the equalities (a)-(d) hold because the differentials of the complexes involved are all zero. So, from the SES ~$0\ra \im f_n\ra\cod f_n\ra\coker f_n\ra 0$,~ we get the SES
\bea
0\ra \coker \delta_{n-1}\sr{f_n}{\ral} H_n(Hom(P_\ast,M))\sr{g_n}{\ra}\ker\delta_n\ra 0.\nn
\eea
where $f_n:=H_n\big(Hom(d_\ast,M)\big)\cong Ext^1\big(H_n(d_\ast),M\big)$ and $g_n:=H_n\big(Hom(i_\ast,M)\big)$ are each unique up to isomorphism (i.e., up to composition with isomorphisms).

{\flushleft \ul{Step 3 (Computing the connecting homomorphism $\delta=\delta_\ast$ above):}} To get $\delta$, consider the diagram
\[\adjustbox{scale=0.9}{\bt
                &  \vdots\ar[d]                                          &&  \vdots\ar[d]                                                      &&  \vdots\ar[d]                     &    \\
 0\ar[r] &  {Hom(B_{n-2},M)} \ar[d]\ar[rr,"{h(d_{n-1},M)}"]       &&  {Hom(P_{n-1},M)}\ar[d] \ar[rr,"{h(i_{n-1},M)}"]            && {Hom(Z_{n-1},M)} \ar[d]\ar[r]     &  0 \\
 0\ar[r]        &  {Hom(B_{n-1},M)}\ar[d]\ar[rr,"{h(d_n,M)}"]     &&  {Hom(P_n,M)}\ar[d,dashed,"{h(d_{n+1},M)}"]\ar[rr,dashed,"{h(i_n,M)}"]        && {Hom(Z_n,M)} \ar[d]\ar[r,dashed]  &  0 \\
 0\ar[r,dashed]        &  {Hom(B_n,M)}\ar[d]\ar[rr,dashed,"{h(d_{n+1},M)}"]            &&  {Hom(P_{n+1},M)}\ar[d]\ar[rr,"{h(i_{n+1},M)}"]                    &&  {Hom(Z_{n+1},M)} \ar[d]\ar[r]    &  0 \\
                &  \vdots                                                &&  \vdots                                                            &&  \vdots                           &
\et}\]
in which $Hom(-,-)$ is written as $h(-,-)$ for convenience.

Then $\delta_n(h_n^{ZM}+\im h(d_{n-1},M))=h_n^{BM}+\im h(d_{n-1},M)$, where ( with $h_n^{PM}\in Hom(P_n,M)$ )
\bea
&&h_n^{BM}\circ d_{n+1}=:h(d_{n+1},M)(h_n^{BM})=h(d_{n+1},M)(h_n^{PM}):=h_n^{PM}\circ d_{n+1},\nn\\
&&~~~~~~\txt{along with}~~(h_n^{PM}=h_n^{PM}\circ i_n~=:~)~h(i_n,M)(h_n^{PM})=h_n^{ZM},\nn\\
&&~~\Ra~~h_n^{BM}\circ d_{n+1}=h_n^{ZM}\circ d_{n+1},~~\sr{(s)}{\Ra}~~h_n^{BM}=h_n^{ZM},\nn
\eea
where step (s) holds because $P_{n+1}\sr{d_{n+1}}{\ral}B_n$ is surjective. Therefore, $\delta_n=h(i_n,M):=Hom(i_n,M)=(-)\circ i_n$, where $i_n$ is \ul{now} the inclusion $B_n\subset Z_n$, for every $n$.

{\flushleft \ul{Step 4 (Computing the end terms of the SES):}} With $\delta_{n-1}=Hom(i_{n-1},M)=(-)\circ i_{n-1}$, applying the left exact functor $Hom(-,M)$ on the SES ~$0\ra B_n\sr{i_n}{\ral}Z_n\sr{\pi_n}{\ral}{Z_n\over B_n}=H_n\ra 0$,~ we get the exact sequence ~$0\ra Hom\left({Z_n\over B_n},M\right)\sr{h(\pi_n,M)}{\ral}Hom(Z_n,M)\sr{h(i_n,M)}{\ral} Hom(B_n,M)$,~ and so
\bea
\textstyle\ker\delta_n=\ker h(i_n,M)\cong \im h(\pi_n,M)\cong Hom\left({Z_n\over B_n},M\right).\nn
\eea
Finally, the ext-LES for the SES
\bea
0\ra B_{n}\sr{i_{n}}{\ral}Z_{n}\ra H_{n}\ra 0~~~~\txt{( a projective resolution of $H_{n}=H_{n}(P_\ast)$ )}\nn
\eea
gives the exact sequence
\bea
&&~~~~0\ra Hom(H_{n},M)\ra Hom(Z_{n},M) \sr{h(i_{n},M)}{\ral} Hom(B_{n},M)\ra \nn\\
&&~~~~~~~~~~~~\ra ext^1_R(H_n,M)\ra \ub{ext^1_R(Z_n,M)}_0\ra \ub{ext^1_R(B_n,M)}_0\ra\cdots,\nn
\eea
which shows ~$\coker\delta_n=\coker h(i_n,M)\cong ext^1_R(H_n,M)$.

{\flushleft \ul{Step 5 (Showing the sequences are split):}} This follows the same steps as in the proof of Theorem \ref{UCTfhom}, with $Hom_R(-,M)$ playing the role of $-\otimes_R M$. Recall that every additive functor preserves split-exactness of sequences.
\end{proof} 

%% file: parts/AlgebraCat/AlgebraCatS12.tex
\chapter{Grading and Triangulation in Categories}\label{AlgebraCatS12}
The first three short sections below (in which some of the notation and conventions are based on \cite[Section 10.2, pp 616-623]{rotman2009}) are included only for general knowledge/awareness/motivation but they are neither explicitly required nor used in the subsequent discussion on distinguished triangles and triangulation.

\section{Graded Objects and Filtered Objects}
\begin{dfn}[\textcolor{blue}{
\index{Grading! (or Graded object)}{}\index{Graded! object (or Grading)}{Graded object (or Grading) in a category},
\index{Multiple grading}{Multiple grading},
\index{Multiply graded object}{Multiply graded object},
\index{Bigraded! object}{Bigraded object}}]
Let $\I,\C$ be categories. An \ul{$\I$-graded object} (or an \ul{$\I$-grading}) in $\C$ is a functor $S:\I\ra\C$ (i.e., a system in $\C$). (\blue{footnote}\footnote{An $\I$-graded object is an $\I$-indexed object, but of a special type because of the functorial property of the indexing operation, i.e., the grading operation is a functorial indexing operation (i.e., a functor, and not just any map).})

A $\Integer^n$-graded object in $\C$ (i.e., an object $S\in\C^{\Integer^n}$) is said to be \ul{$n$-graded} (or \ul{$n$-multigraded}, or \ul{$n$-multiply graded}). A $1$-graded object (i.e., a system $S\in\C^\Integer$) is said to be \ul{graded}, while a $2$-graded object (i.e., a system $S\in\C^{\Integer^2}$) is said to be \ul{bigraded} (or \ul{doubly graded}).
\end{dfn}

If $\A$ is an abelian category, then $n$-multicomplexes $C\in\A_0^{\Integer^n}$ are examples of $n$-graded objects in $\A$.

\begin{dfn}[\textcolor{blue}{
\index{Filtration}{Filtration},
\index{Filtration of an object}{Filtration of an object},
\index{Filtered! object}{Filtered object},
\index{Finitely filtered object}{Finitely filtered object}}]
Let $\C$ be a category. A graded object $S:\Integer\ra\C$ is called a \ul{filtration} if (i) it is a linear system and (ii) the transition morphisms $s_i:S(i)\ra S(i+1)$ are monomorphisms, i.e., $S$ is a sequence of subobjects
\bea
\cdots\subset S(i-1)\subset S(i)\subset S(i+1)\subset\cdots.\nn
\eea
A filtration $S:\Integer\ra\C$ is a \ul{filtration of an object} $C\in\Ob\C$ (making $C$ a \ul{filtered object}) if
\bea
\cdots\subset S(i-1)\subset S(i)\subset S(i+1)\subset\cdots\subset C.\nn
\eea
\end{dfn}
Henceforth, in any abelian category $\A$, a filtered object $A$ will be assumed (for simplicity) to be a \ul{finitely filtered object} of one of the following decreasing/increasing forms (depending on convenience):
\bea
A=A_0\supset A_1\supset\cdots\supset A_{n+1}=\{0\},~~~~\txt{or}~~~~\{0\}=A^0\subset A^1\subset\cdots\subset A^{n+1}=A.\nn
\eea

\begin{dfn}[\textcolor{blue}{
\index{Graded! morphism}{Graded morphism},
\index{Degree of a graded morphism}{Degree of a graded morphism},
\index{Category of! graded objects}{Category of graded objects},
\index{Bigraded! objects}{Bigraded objects},
\index{Bigraded! morphisms}{Bigraded morphisms},
\index{Category of! bigraded objects}{Category of bigraded objects}}]
Let $\I,\C$ be categories, $g:\I\ra\I$ an isomorphism, and $S,T\in\C^\I$ any $\I$-graded objects. An \ul{$\I$-graded morphism of degree $g$} from $S$ to $T$ (written $f:S\ra T$) is a morphism of systems $\wt{f}:S\ra T\circ g$ (i.e., $\wt{f}_i:S_i\ra T_{g(i)}$ for all $i\in I$). More explicitly,
\bea
\bt i\ar[d,"\kappa_{ij}"]\\ j\et~~\mapsto~~
\bt
S(i)\ar[d,"{S(\kappa_{ij})}"]\ar[r,"\wt{f}_i"] & T\circ g(i)\ar[d,"{T\circ g(\kappa_{ij})}"]\\
S(j)\ar[r,"\wt{f}_j"] & T\circ g(j)
\et,~~~~\wt{f}_jS(\kappa_{ij})=T\circ g(\kappa_{ij})\wt{f}_i.\nn\eea
We write \ul{$d_f:=\deg f:=g$}. If $\deg f=id_\I$, we say $f$ is an \ul{$\I$-graded morphism}. The \ul{category of $\I$-graded objects} in $\C$ is the category $Gr(\C^\I)$ whose (i) objects are $\I$-systems $S\in\Ob\C^\I$, (ii) morphisms $f\in \Mor_{Gr(\C^\I)}(S,T)$ are $\I$-graded morphisms $f:S\ra T$, i.e.,
\bea
\textstyle \Mor_{Gr(\C^\I)}(S,T):=\bigcup\limits_{\txt{isos}~g:\I\ra\I}\Mor_{\C^\I}(S,T\circ g),\nn
\eea
and (iii) the composition of $\I$-graded morphisms $f'\circ f:S\sr{f}{\ral}S'\sr{f'}{\ral}S''$ is given by
\bea
\label{GradedCompRule}\wt{f'\circ f}:=(\wt{f'}\circ d_f)\circ\wt{f}:S\sr{\wt{f}}{\ral}S'\circ d_f\sr{\wt{f'}\circ d_f}{\ral}S''\circ d_f\circ d_{f'},~~\Ra~~~~d_{f'\circ f}=d_f\circ d_{f'}.
\eea

In particular, if $\C$ is a category, $k\in\Integer$, and $S,T\in\C^\Integer$ graded objects, then a \ul{graded morphism of degree $k\in\Integer$} from $S$ to $T$ (written $f:S\ra T$) is a morphism of systems $\wt{f}:S\ra T[k]$ (i.e., $\wt{f}_i:S_i\ra T_{i+k}$ for all $i\in \Integer$, and for the \ul{general case}, see the \blue{footnote}\footnote{More generally, for any system $S\in\C^{\Integer^n}$, we define another system {\tiny $S[k]=S[k_1,...,k_n]\in\C^{\Integer^n}$} (where $k:=(k_1,...,k_n)\in\Integer^n$) by
\bea
S[k]:~i\sr{\kappa_{ij}}{\ral}j~\mapsto~S_i[k]\sr{s_{ij}[k]}{\ral}S_j[k],~~~~\left.
                                                                              \begin{array}{l}
                                                                               S_i[k]:=S_{i+k}=S_{i_1+k_1,...,i_n+k_n},~~ \txt{for all}~~i=(i_1,...,i_n)\in\Integer^n,\\
                                                                               s_{ij}[k]:=s_{i+k~j+k}.
                                                                              \end{array}
                                                                            \right.\nn
\eea
}). We write \ul{$d_f:=\deg f:=k$}. If $\deg f=0$, we say $f$ is a \ul{graded morphism}. The \ul{category of graded objects} in $\C$ is the category $Gr(\C^\Integer)$ whose objects are systems $S\in\Ob\C^\Integer$ and whose morphisms $f\in \Mor_{Gr(\C^\Integer)}(S,T)$ are graded morphisms $f:S\ra T$, i.e.,
\bea
\textstyle \Mor_{Gr(\C^\Integer)}(S,T):=\bigcup\limits_{k\in\Integer}\Mor_{\C^\Integer}(S,T[k]).\nn
\eea

If $\C$ is a category, $k,k'\in\Integer$, and $S,T\in\C^{\Integer^2}$ bigraded objects, then a \ul{bigraded morphism of bidegree $(k,k')\in\Integer^2$} from $S$ to $T$ (written $f:S\ra T$) is a morphism of systems $\wt{f}:S\ra T[k,k']$ (i.e., $\wt{f}_{i,j}:S_{i,j}\ra T_{i+k,j+k'}$ for all $i,j\in \Integer$). We write \ul{$d_f:=\deg f:=(k,k')$}. If $\deg f=0$, we say $f$ is a \ul{graded morphism}. The \ul{category of bigraded objects} in $\C$ is the category $Gr(\C^{\Integer^2})$ whose objects are systems $S\in\Ob\C^{\Integer^2}$ and whose morphisms $f\in \Mor_{Gr(\C^{\Integer^2})}(S,T)$ are bigraded morphisms $f:S\ra T$, i.e.,
\bea
\textstyle \Mor_{Gr(\C^{\Integer^2})}(S,T):=\bigcup\limits_{(k,k')\in\Integer^2}\Mor_{\C^{\Integer^2}}(S,T[k,k']).\nn
\eea

More generally, the \ul{category of $n$-multigraded objects} in $\C$ is the category $Gr(\C^{\Integer^n})$ whose objects are systems {\small $S\in\Ob\C^{\Integer^n}$} and whose morphisms {\small $f\in \Mor_{Gr(\C^{\Integer^n})}(S,T)$} are $n$-multigraded morphisms $f:S\ra T$, i.e.,
\bea
\textstyle \Mor_{Gr(\C^{\Integer^n})}(S,T):=\bigcup\limits_{(k_1,...,k_n)\in\Integer^n}\Mor_{\C^{\Integer^n}}(S,T[k_1,...,k_n]).\nn
\eea
\end{dfn}

If $\A$ is an abelian category and $(C^\ast,d^\ast)\in\A_0^\Integer$ a cocomplex, then the differential $d^\ast:C^\ast\ra C^\ast[1]$ (i.e., $d^n:C^n\ra C^{n+1}$) is a graded morphism $C^\ast\ra C^\ast$ of degree $1$. Similarly, given a chain complex $(C_\ast,d_\ast)\in\A_0^\Integer$, the differential $d_\ast:C_\ast\ra C_\ast[-1]$ is a graded morphism $C_\ast\ra C_\ast$ of degree $-1$. (\blue{footnote}\footnote{Given chain complexes $(C_\ast,d_\ast),(C'_\ast,d'_\ast)\in\A_0^\Integer$, a chain homotopy $h_\ast:C_\ast\ra C'_{\ast+1}$ may also be thought of as a graded ``morphism'' $C_\ast\ra C_\ast$ of degree 1.})

\begin{notation}
Let $\I,\C$ be categories. Given an $\I$-graded morphism $f\in \Mor_{Gr(\C^\I)}(S,T)$, i.e., a morphism $\wt{f}:S\ra T\circ d_f$, the component $\wt{f}_i:S\ra T_{d_f(i)}$ is sometimes written simply as $f_i:S\ra T_{d_f(i)}$. However, to avoid needless confusion in our relatively short discussion we will not do this.

Also, if $\deg f=0$, then it is clear that $\wt{f}=f:S\ra T$ (i.e., it is not strictly necessary to include the \ul{grading-warner} $^\sim$ for zero degree morphisms).
\end{notation}

\begin{lmm}[\textcolor{blue}{\index{Additivity of! the degree}{Additivity of the degree}}]
Let $\C$ be a category and $f,f'$ morphisms in $Gr(\C^{\Integer^n})$. Then
\bea
\deg(f'\circ f)=\deg f'+\deg f.\nn
\eea
\end{lmm}

\begin{dfn}[\textcolor{blue}{\index{Predecessor morphism}{Predecessor morphism} of a graded morphism}]
Let $\C$ be a category and $f\in \Mor_{Gr(\C^{\Integer^n})}(S,T)$ a graded morphism. We define the \ul{predecessor morphism} $f^-\in \Mor_{Gr(\C^{\Integer^n})}(S,T)$ of $f$ by
\[
\wt{f}^-_i:=\wt{f}_{i-\deg f}:S_{i-\deg f}\ra T_i,~~~~\txt{for all}~~~~i\in\Integer^n.\nn
\]
\end{dfn}

It might be possible to avoid using the predecessor morphism in the following definition by explicitly making use of the definition of composition for $\Integer^n$-graded morphisms from the general case given in (\ref{GradedCompRule}).

\begin{dfn}[\textcolor{blue}{\index{Kernel of! a graded morphism}{}\index{Image of! a graded morphism}{}\index{Cokernel of a graded morphism}{}Kernel/Image/Cokernel of a graded morphism,
\index{Exactness! of a sequence of graded morphisms}{Exactness of a sequence of graded morphisms}}]
Let $\A$ be an abelian category and $f\in \Mor_{Gr(\A^{\Integer^n})}(S,T)$. The \ul{kernel} of $f$ is given by $(\ker f)_i:=\ker \wt{f}_i\subset S_i$, where $i=(i_1,...,i_n)$, and $\deg f=(k_1,...,k_n)\in\Integer^n$.
The \ul{image} of $f$ is given by
\bea
(\im f)_i:=\im \wt{f}^-_i=\wt{f}_{i-\deg f}(S_{i-\deg f})\subset T_i,~~~~\txt{since}~~~~\wt{f}_i(S_i)\subset T_{i+\deg f},~~~~\txt{for all}~~i\in\Integer^n.\nn
\eea
Similarly, the \ul{cokernel} of $f$ (through which ~$\im f:=\ker(\cod f\twoheadrightarrow\coker f)$ ) is given by
\bea
(\coker f)_i:=\coker \wt{f}^-_i=\coker \wt{f}_{i-\deg f},~~~~\txt{for all}~~i\in\Integer^n.\nn
\eea
A sequence $S\sr{f}{\ral}S'\sr{f'}{\ral}S''$ is \ul{exact} (at $S'$) if $\im f=\ker f'$, i.e., if $\im \wt{f}_{i-\deg f}=\ker \wt{f'}_i$ for all $i\in\Integer^n$.
\end{dfn}

\begin{dfn}[\textcolor{blue}{\index{Exact! triangle}{Exact triangle} in a LES category}]
Let $\A$ be a LES category. Given a SES of complexes $0\ra A_\ast\sr{f_\ast}{\ral}B_\ast\sr{g_\ast}{\ral}C_\ast\ra 0$ in $\A_0^\Integer$, the induced LES of homologies for the SES, i.e.,
\[\adjustbox{scale=0.9}{
\bt \cdots\ar[r,"\wt{\delta}_{n+1}"] & H_n(A_\ast)\ar[r,"H_n(f_\ast)"] & H_n(B_\ast)\ar[r,"H_n(g_\ast)"] & H_n(C_\ast)\ar[r,"\wt{\delta}_n"] & H_{n-1}(A_\ast)\ar[r,"H_{n-1}(f_\ast)"]\ar[r] & H_{n-1}(B_\ast)\ar[r] &\cdots\et,}\]
is called an \ul{exact triangle} in $\A$, and denoted by ~$H_\ast(A_\ast)\sr{H_\ast(f_\ast)}{\ral}H_\ast(B_\ast)\sr{H_\ast(g_\ast)}{\ral}H_\ast(C_\ast)\sr{\delta_\ast}{\ral}H_\ast(A_\ast)$, ~or
\[\adjustbox{scale=0.8}{\bt
H_\ast(A_\ast)\ar[rr,"H_\ast(f_\ast)"] && H_\ast(B_\ast)\ar[dl,"H_\ast(g_\ast)"]\\
  & H_\ast(C_\ast)\ar[ul,"\delta_\ast"] &
\et}\]
where $H_\ast(f_\ast),H_\ast(g_\ast),\delta_\ast$ are viewed as graded morphisms with degrees
\bea
\deg H_\ast(f_\ast)=\deg H_\ast(g_\ast)=0,~~~~\deg\delta_\ast=-1,\nn
\eea
i.e., the above triangle is a system in the subcategory $Gr(\A_0^\Integer)\subset Gr(\A^\Integer)$ of graded complexes in $\A$.
\end{dfn}

\begin{dfn}[\textcolor{blue}{
\index{Filtered! complex}{Filtered complex},
\index{Factors of a filtered complex}{Factors of a filtered complex},
\index{Grading! of a filtered complex}{Associated grading of a filtered complex}}]
Let $\A$ be an abelian category. A complex $(C_\ast,d_\ast)\in\A_0^\Integer$ (where the differential will be viewed as a morphism in $Gr(\A_0^\Integer)$) is a \ul{filtered complex} if there exists a filtration
\[
\{0\}=C_\ast^0\subset C_\ast^1\subset\cdots C_\ast^{n+1}=C_\ast~~(~\txt{by subcomplexes}~~(C_\ast^k,d_\ast)\in\A_0^\Integer~)~~\txt{such that}~~d_\ast C_\ast^k\subset C_\ast^k,
\]
i.e., $d_iC_i^k\subset C_{i-1}^k$. The \ul{associated graded object} (or \ul{associated grading}) of a filtered complex $(C_\ast,d_\ast)$ is
{\small\begin{align}
\textstyle Gr(C_\ast):=\big(Gr^k(C_\ast)\big)_{k\geq 0}~~\txt{or (\blue{footnote}\footnotemark)}~~\bigoplus\limits_{k\geq 0}Gr^k(C_\ast),~~~~Gr^k(C_\ast):={C_\ast^k\over C_\ast^{k-1}}~~\txt{(called \ul{factors of the filtration})}.\nn
\end{align}}\footnotetext{Recall that a \ul{finite} system $S\in\A^\I$ in an abelian category $\A$ has a direct sum, and the direct sum (by construction) also uniquely determines the \ul{finite} system up to isomorphism.}
\end{dfn}

It is sometimes convenient to separate the filtration index from other indices, say by writing $C_\ast^k$ as
\bea
F^kC_\ast:=C_\ast^k,\nn
\eea
so that the filtration of $C_\ast$ takes the form
\bea
\{0\}=F^0C_\ast\subset F^1C_\ast\subset\cdots\subset F^{n+1}C_\ast=C_\ast,~~~~\txt{with}~~~~d_\ast F^kC_\ast\subset F^kC_\ast,\nn
\eea
and the associated graded object of $C_\ast$ is
{\small\bea
\textstyle Gr(C_\ast):=\big(Gr^k(C_\ast)\big)_{k\geq 0}~~\txt{or}~~\bigoplus\limits_{k\geq 0}Gr^k(C_\ast),~~~~~~~~Gr^k(C_\ast):={F^kC_\ast\over F^{k-1}C_\ast}~~~~\txt{(filtration factors)}.\nn
\eea}

\section{Homology of Graded Objects}
\begin{dfn}[\textcolor{blue}{
\index{Differential bigraded object (DBO)}{Differential bigraded object (DBO)},
\index{Graded! differential}{Graded differential},
\index{Homology! of a DBO}{Homology of a DBO},
\index{Category of! DBO's}{Category of DBO's}}]
Let $\A$ be an abelian category. A \ul{differential bigraded object (DBO)} in $\A$ is a bigraded object $(S,d)=(S_{\ast\ast},d_{\ast\ast})\in Gr(\A^{\Integer^2})$ such that the system transition morphisms $d_{ij}:S_{ij}\ra S_{i'j'}$ form a bigraded morphism $d:S\ra S$ (called \ul{graded differential}) such that $d^2:=d\circ d=0:S\sr{d}{\ral}S\sr{d}{\ral}S$. Recall that by the composition rule (\ref{GradedCompRule}) for bigraded morphisms, if $\deg d_{\ast\ast}:=(a,b)$, i.e., $\wt{d}_{\ast\ast}:S_{\ast\ast}\ra S_{\ast+a,\ast+b}$, then $d^2=0$ precisely means $\wt{d}_{ij}^2:=\wt{d}_{i+a,j+b}\circ \wt{d}_{ij}=0$ for all $i,j$, or explicitly,
\bea\bt
0=\wt{d}_{ij}^2:=\wt{d}_{i+a,j+b}\circ \wt{d}_{ij}: S_{ij}\ar[rr,"\wt{d}_{ij}"] && S_{i+a,j+b} \ar[rr,"\wt{d}_{i+a,j+b}"] && S_{i+2a,j+2b},\nn
\et\nn\eea
which will also be written as {\footnotesize \bt 0=\wt{d}_{\ast\ast}^2:=\wt{d}_{\ast+a,\ast+b}\circ \wt{d}_{\ast\ast}: S_{\ast\ast}\ar[rr,"\wt{d}_{\ast\ast}"] && S_{\ast+a,\ast+b} \ar[rr,"\wt{d}_{\ast+a,\ast+b}"] && S_{\ast+2a,\ast+2b}.\et}

The \ul{homology} of a DBO {\small $(S,d)=(S_{\ast\ast},d_{\ast\ast})$} is the bigraded object {\small $H(S,d)=H_{\ast\ast}(S,d):={\ker d\over\im d}\in Gr(\A^{\Integer^2})$} given componentwise by
\bea
\textstyle H_{ij}(S,d):={\ker\wt{d}_{ij}\over\im\wt{d}^-_{ij}}={\ker\wt{d}_{ij}\over\im \wt{d}_{i-a,j-b}},~~~~\txt{where}~~\deg d:=(a,b).\nn
\eea
The \ul{category of DBO's} in $\A$ will be denoted by $Gr\big(\A_{00}^{\Integer^2}\big)\subset Gr(\A^{\Integer^2})$.

The \ul{homology $H$} (as defined above) is really a \ul{functor}
\bea
H:\A^{\Integer^2}_{00}\ra \A^{\Integer^2}_{00},~~(S,d)\sr{f}{\ral}(S',d')~~\mapsto~~H(S,d)\sr{H(f)}{\ral}H(S',d'),\nn
\eea
where the induced morphism ~$H(f):=f|_{\ker d}+\im d'$, and $H(S,d):={\ker d\over\im d}:=\coker(\ker d\hookrightarrow\im d)$. Equivalently, \ul{the $(i,j)$th homology $H_{ij}$} (as defined above) \ul{is a functor}
\bea
H_{ij}:\A^{\Integer^2}_{00}\ra\A,~~(S,d)\sr{f}{\ral}(S',d')~~\mapsto~~H_{ij}(S,d)\sr{H_{ij}(f)}{\ral}H_{ij}(S',d'),\nn
\eea
where the induced morphism ~$H_{ij}(f):=f|_{\ker \wt{d}_{ij}}+\im \wt{d'}^-_{ij}$, and $H_{ij}(S,d):={\ker \wt{d}_{ij}\over\im \wt{d}^-_{ij}}:=\coker(\ker \wt{d}_{ij}\hookrightarrow\im \wt{d}^-_{ij})$.
\end{dfn}

\begin{rmk*}[\textcolor{blue}{\index{Differential multigraded object (DMO)}{Differential multigraded objects (DMO's)}}]
Let $\A$ be an abelian category. It is clear that the concept of a DBO above generalizes to give a \ul{differential $n$-multigraded object} ($n$-DMO) in $\A$ as an $n$-multigraded object $(S,d)=(S_{\ast_1,...,\ast_n},d_{\ast_1,...,\ast_n})\in Gr(\A^{\Integer^n})$ such that the system transition morphisms $d_i:S_i\ra S_{i'}$ form an $n$-multigraded morphism $d:S\ra S$ (\ul{graded differential}) such that $d^2:=d\circ d=0:S\sr{d}{\ral}S\sr{d}{\ral}S$, where as before, by the composition rule (\ref{GradedCompRule}) for $n$-multigraded morphisms, if $\deg d=a=(a_1,...,a_n)$, i.e., $\wt{d}:S_{\ast_1,...,\ast_n}\ra S_{\ast_1+a_1,...,\ast_n+a_n}$, then $d^2=0$ precisely means $\wt{d}_{i_1\cdots i_n}^2:=\wt{d}_{i_1+a_1,\cdots,i_n+a_n}\circ \wt{d}_{i_1\cdots i_n}=0$ for all $i_1,\cdots,i_n$, or explicitly,
{\small\bea\bt
0=\wt{d}_{i_1,...,i_n}^2:=\wt{d}_{i_1+a_1,...,i_n+a_n}\circ \wt{d}_{i_1,...,i_n}: S_{i_n,...,i_n}\ar[rr,"\wt{d}_{i_1,...,i_n}"] && S_{i_1+a_1,...,i_n+a_n} \ar[rr,"\wt{d}_{i_1+a_1,...,i_n+a_n}"] && S_{i_n+2a_1,...,i_n+2a_n}.\nn
\et\nn\eea}also written {\footnotesize \bt
0=\wt{d}_{\ast_1,...,\ast_n}^2:=\wt{d}_{\ast_1+a_1,...,\ast_n+a_n}\circ \wt{d}_{\ast_1,...,\ast_n}: S_{\ast_n,...,\ast_n}\ar[rr,"\wt{d}_{\ast_1,...,\ast_n}"] && S_{\ast_1+a_1,...,\ast_n+a_n} \ar[rr,"\wt{d}_{\ast_1+a_1,...,\ast_n+a_n}"] && S_{\ast_n+2a_1,...,\ast_n+2a_n}.\nn
\et} The \ul{homology} of the $n$-DMO ~$(S,d)=(S_{\ast_1,...,\ast_n},d_{\ast_1,...,\ast_n})$~ is the $n$-multigraded object
\bea
\textstyle H(S,d)=H_{\ast_1,...,\ast_n}(S,d):={\ker d\over\im d}:=\coker(\ker d\hookrightarrow\im d)~\in~Gr(\A^{\Integer^n}),\nn
\eea
given in component form (with indices $i=(i_1,...,i_n)\in\Integer^n$) by
\bea
\textstyle H_i(S,d):={\ker \wt{d}_i\over\im \wt{d}^-_i}={\ker \wt{d}_i\over\im \wt{d}_{i-\deg d}}={\ker \wt{d}_{i_1,...,i_n}\over\im \wt{d}_{i_1-a_1,...,i_n-a_n}},~~~~\txt{where}~~\deg d:=(a_1,...,a_n).\nn
\eea
The \ul{category of $n$-DMO's} in $\A$ will be denoted by $Gr\big(\A_{0_1,...,0_n}^{\Integer^n}\big)\subset Gr(\A^{\Integer^n})$. Thus, as before, the {homology $H$} (as defined above) is really a {functor}
\bea
H:\A_{0_1,...,0_n}^{\Integer^n}\ra \A_{0_1,...,0_n}^{\Integer^n},~~(S,d)\sr{f}{\ral}(S',d')~~\mapsto~~H(S,d)\sr{H(f)}{\ral}H(S',d'),\nn
\eea
where the induced morphism ~$H(f):=f|_{\ker d}+\im d'$, and $H(S,d):={\ker d\over\im d}:=\coker(\ker d\hookrightarrow\im d)$. Equivalently, for $i=(i_1,...,i_n)\in\Integer^n$, \ul{the $i$th homology $H_i$} (as defined above) \ul{is a functor}
\bea
H_i:\A_{0_1,...,0_n}^{\Integer^n}\ra\A,~~(S,d)\sr{f}{\ral}(S',d')~~\mapsto~~H_i(S,d)\sr{H_i(f)}{\ral}H_i(S',d'),\nn
\eea
where the induced morphism ~$H_i(f):=f|_{\ker \wt{d}_i}+\im \wt{d'}^-_i$, and $H_i(S,d):={\ker \wt{d}_i\over\im \wt{d}^-_i}:=\coker(\ker \wt{d}_i\hookrightarrow\im \wt{d}^-_i)$.
\end{rmk*}

\section{Exact Couples and Spectral Sequences}
\begin{dfn}[\textcolor{blue}{\index{Exact! couple}{Exact couple in an abelian category}}]
Let $\A$ be an abelian category. An \ul{exact couple} in $\A$ is an exact sequence in the category $Gr(\A^{\Integer^2})$ of bigraded objects in $\A$ of the form
\[\adjustbox{scale=0.8}{\bt
A_{\ast\ast}\ar[rr,"\al_{\ast\ast}"] && A_{\ast\ast}\ar[dl,"\beta_{\ast\ast}"]\\
  & B_{\ast\ast}\ar[ul,"\gamma_{\ast\ast}"] &
\et},~~~~
\im\al_{\ast\ast}=\ker\beta_{\ast\ast},~~
\im\beta_{\ast\ast}=\ker\gamma_{\ast\ast},~~
\im\gamma_{\ast\ast}=\ker\al_{\ast\ast}.
\]
We also say the $5$-tuple $(A,B,\al,\beta,\gamma)$ is an exact couple in $\A$. Using more familiar notation, we can also say the sequence ~$A\sr{\al}{\ral}A\sr{\beta}{\ral}B\sr{\gamma}{\ral}A$~ is an exact couple in $\A$. Componentwise, this means
\[
\bt A_{ij}\ar[r,"\wt{\al}_{ij}"] & A_{(i,j)+d_\al}\ar[rr,"\wt{\beta}_{(i,j)+d_\al}"] && B_{(i,j)+d_\al+d_\beta}\ar[rr,"\wt{\gamma}_{(i,j)+d_\al+d_\beta}"] && A_{(i,j)+d_\al+d_\beta+d_\gamma}\et,~~~~~~\txt{for each pair}~~i,j.
\]
\end{dfn}
In general, a \ul{triangle} ~$A\sr{f}{\ral}B\sr{g}{\ral}C\sr{h}{\ral}A$~ in $Gr(\A^{\Integer^2})$ has the componentwise form
\[
\bt A_{ij}\ar[r,"\wt{f}_{ij}"] & B_{(i,j)+d_f}\ar[rr,"\wt{g}_{(i,j)+d_f}"] && C_{(i,j)+d_f+d_g}\ar[rr,"\wt{h}_{(i,j)+d_f+d_g}"] && A_{(i,j)+d_f+d_g+d_h}\et,~~~~~~\txt{for each pair}~~i,j.
\]

\begin{prp}[\textcolor{blue}{\cite[Proposition 10.8, p.618]{rotman2009}}]\label{FiltExCoup1}
Let $\A$ be a LES category and $C_\ast\in\A_0^\Integer$ a complex in $\A$. Then every filtration $(F^kC_\ast)_{k\in\Integer}$ of $C_\ast$ determines an exact couple
\[\adjustbox{scale=0.8}{\bt
A\ar[rr,"{\al,{(1,-1)}}"] && A\ar[dl,"{\beta,{(0,0)}}"]\\
  & B\ar[ul,"{\gamma,{(-1,0)}}"] &
\et}\]
where the pair of numbers accompanying each bigraded morphism is its degree (bidegree).
\end{prp}
\begin{proof}
Let $C_\ast^k:=F^kC_\ast$. We have a SES of complexes $0\ra C_\ast^{k-1}\sr{f_\ast}{\ral}C_\ast^k\sr{g_\ast}{\ral}{C_\ast^k\over C_\ast^{k-1}}\ra 0$, with the associated induced LES of homologies
\bea\bt
\cdots\ar[r] & H_n(C_\ast^{k-1})\ar[r,"{H_n(f_\ast)}"] & H_n(C_\ast^k)\ar[r,"{H_n(g_\ast)}"] & H_n\left(C_\ast^k/C_\ast^{k-1}\right)\ar[r,"\delta_n"] & H_{n-1}(C_\ast^{k-1})\ar[r,"{H_{n-1}(f_\ast)}"] & \cdots
\et\nn\eea
The desired exact couple is obtained by letting $i:=k$, $j:=n-i=n-k$ (so $n=i+j$), and then defining
\[
A_{ij}:=H_{i+j}(C_\ast^{i-1}),~~~B_{ij}:=H_{i+j}(C_\ast^i/C_\ast^{i-1}),~~~\wt{\al}_{ij}:=H_{i+j}(f_\ast),~~~\wt{\beta}_{ij}:=H_{i+j}(g_\ast),~~~~\wt{\gamma}_{ij}:=\delta_{i+j}. \qedhere
\]
\end{proof}

\begin{prp}[\textcolor{blue}{\cite[Proposition 10.9, p.620]{rotman2009}}]\label{FiltExCoup2}
In $R$-mod, given an exact couple $A\sr{\al}{\ral}A\sr{\beta}{\ral}B\sr{\gamma}{\ral}A$, we get the following:
\bit
\item[(i)] A differential graded object ~$(B,d^1)\in Gr(R\txt{-mod}^\Integer_0)$, ~where ~$d^1:=\beta\gamma:B\sr{\gamma}{\ral} A\sr{\beta}{\ral}B$.
\item[(ii)] Another exact couple (called \ul{derived couple}) $A^2\sr{\al^2}{\ral}A^2\sr{\beta^2}{\ral}B^2\sr{\gamma^2}{\ral}A^2$ such that
\[
\deg\al^2=\deg\al,~~~~\deg\beta^2=\deg\beta-\deg\al,~~~~\deg\gamma^2=\deg\gamma.
\]
\eit
\end{prp}
\begin{proof}
Let $\deg\al=(i_\al,j_\al)$, $\deg\beta=(i_\beta,j_\beta)$, and $\deg\gamma=(i_\gamma,j_\gamma)$.
{\flushleft (i)} It is clear that the morphism $d^1:=\beta\gamma:B\ra B$ makes $(B,d^1)$ a DBO, since $d^1 d^1=(\beta\gamma)(\beta\gamma)=\beta(\gamma\beta)\gamma=\beta 0\gamma=0$, by the exactness of the given sequence. The degree of $d^1$ is of course given by
\[
\deg d^1=\deg(\beta\gamma)=\deg\beta+\deg\gamma=(i_\beta+i_\gamma,j_\beta+j_\gamma)=:(i_{d^1},j_{d^1}).
\]

{\flushleft (ii)} Using the given information, we can construct the following diagram, in which $A^2:=\im\al$ and $B^2:=H(B,d^1)$ are given componentwise by
\[
A^2_{ij}:=(\im\al)_{ij}=\im\wt{\al}^-_{ij}\subset A_{ij}~~~~\txt{and}~~~~B^2_{ij}:=H_{ij}(B,d^1):={\ker \wt{d}^1_{ij}\over\im \wt{d}^{1-}_{ij}}.
\]

\[\adjustbox{scale=0.9}{\bt
   &               &                           &                 &   \overbrace{\im\al}^{A^2}\ar[d,hook]              &                          &       0\ar[dr]             &   0                    &          &  \\
   &               &              A\ar[urr,two heads]\ar[drrr,dashed,"\beta"']\ar[rr,"\al"]             &                &   A\ar[dr,"\beta"]               &              0\ar[r]     & \im d^1\ar[r,hook,""]\ar[ur] & \ker d^1\ar[dr,hook]\ar[r,two heads,""]  & \overbrace{H(B,d^1)}^{B^2}\ar[r]  & 0 \\
   & \cdots\ar[r]  & B\ar[u,dashed,"\gamma"]\ar[rrr,near start,"d^1"]\ar[dr,two heads,"d^1"] &                &                  & B\ar[rrr,"d^1"]\ar[ur,two heads,"d^1"] &                    &                       & B\ar[r]  & \cdots   \\
   &               &       0\ar[r]             & \im d^1\ar[r,hook,""]\ar[dr] & \ker d^1\ar[ur,hook,""]\ar[r,two heads,""] & H(B,d^1)\ar[r] & 0                  &                       &          &         \\
   &               &                           &      0\ar[ur]  &   0               &                          &                    &                       &          &
\et}\]
\ul{Define} $\al^2:=\al|_{A^2}:A^2\ra A^2,~\al(a)\mapsto\al^2(a)$. Then with the inclusion $q_{A^2}:A^2\hookrightarrow A$, we have $\deg\al^2=\deg(\al|_{A^2})=\deg(\al q_{A^2})=\deg\al+\deg q_{A^2}=\deg\al+(0,0)=\deg\al$. Therefore, explicitly,
\[
\wt{\al}_{ij}^2:=\wt{\al}_{ij}|_{A_{ij}^2}:A_{ij}^2\ra A_{i+i_{\al},j+j_{\al}}^2,~~\wt{\al}^-_{ij}(a)\mapsto \wt{\al}_{ij}\wt{\al}^-_{ij}(a),~~~~\txt{for}~~a\in A_{i-i_\al,j-j_\al}.\nn
\]
Since $d^1\beta=\beta\gamma\beta =\beta 0=0$, we \ul{define} $\beta^2:A^2\ra B^2$ by $\beta^2(\al(a)):=\beta a+\im d^1$, for $a\in A^2=\im\al$, i.e.,
\[
\wt{\beta}_{ij}^2:A_{ij}^2\ra B_{i+i_{\beta^2},j+j_{\beta^2}}^2,~~\wt{\al}^-_{ij}(a)\mapsto \wt{\beta}_{ij}(a)+\im\wt{d}^{1-}_{ij},
\]
where by construction, we have ~$\deg\beta^2=\deg\beta-\deg\al=(i_\beta-i_\al,j_\beta-j_\al)=:(i_{\beta^2},j_{\beta^2})$.

Finally, if $z\in\ker d^1$, then $0=d^1z=\beta\gamma z$ implies $\gamma z\in\ker\beta=\im\al$, and so we \ul{define} $\gamma^2:B^2\ra A^2$ by $\gamma^2(z+\im d^1):=\gamma z$ (where it is clear that $\deg\gamma^2=\deg\gamma$), i.e.,
\[
\wt{\gamma}_{ij}^2:B_{ij}^2\ra A_{i+i_\gamma,j+j_\gamma}^2,~~z+\im\wt{d}^{1-}_{ij}\mapsto \wt{\gamma}_{ij}(z).
\]
It now remains to show the resulting sequence ~$A^2\sr{\al^2}{\ral}A^2\sr{\beta^2}{\ral}B^2\sr{\gamma^2}{\ral}A^2$~ is exact. It is clear that $\al^2\gamma^2=\beta^2\al^2=\gamma^2\beta^2=0$. Also, we have the following:
\bit[leftmargin=0.2cm]
\item[]\ul{$\ker\beta^2\subset\im\al^2$}: If $0=\beta^2(\al(a)):=\beta(a)+\im d^1$, then $\beta(a)\in\im d^1$, and so $\beta(a)=d^1b=\beta\gamma(b)$, i.e., {\small $a-\gamma(b)\in\ker\beta=\im\al$}. Let $a=\gamma(b)+\al(a')$. Then {\small $\al(a)=\al\gamma(b)+\al\al(a')=\al|_{A^2}(\al(a'))=\al^2(\al(a'))\in\im\al^2$}.
\item[]\ul{$\ker\gamma^2\subset\im\beta^2$}: If $0=\gamma^2(z+\im d^1):=\gamma(z)$, then $z\in\ker\gamma=\im\beta$, and so $z=\beta(a)$, i.e., $z+\im d^1=\beta(a)+\im d^1=\beta^2(\al(a))\in\im\beta^2$.
\item[]\ul{$\ker\al^2\subset\im\gamma^2$}: If $0=\al^2(\al(a)):=\al|_{A^2}(\al(a))$, then $\al(a)\in\ker\al=\im\gamma$, and so $\al(a)=\gamma(b)$. Since $d^1b=\beta\gamma(b)=\beta\al(a)=0$, we have $\gamma(b)=\gamma^2(b+\im d^1)$, i.e., $\al(a)=\gamma(b)=\gamma^2(b+\im d^1)\in \im\gamma^2$. \qedhere
\eit
\end{proof}

\begin{dfn}[\textcolor{blue}{
\index{Multispectral sequence}{Multispectral sequence},
\index{Spectral sequence}{Spectral sequence},
\index{Bispectral sequence}{Bispectral sequence}}]
Let $\A$ be an abelian category. An \ul{$n$-(multi)spectral sequence} in $\A$ is a sequence $(E^r,d^r)_{r\geq 1}$ of differential $n$-multigraded objects
{\small\bea
(E^r,d^r)=(E_{\ast_1,...,\ast_n}^r,d_{\ast_1,...,\ast_n}^r)\in Gr\big(\A_{0_1,...,0_n}^{\Integer^n}\big)\nn
\eea} such that the following hold:
\begin{enumerate}[leftmargin=0.9cm]
\item $E^{r+1} = H(E^r,d^r)$, for all $r\geq 1$. (I.e., it is \ul{a sequence of successive homologies}). In diagram form,
\bea\bt
(E^1,d^1)\ar[r,"H"] & (E^2,d^2)\ar[r,"H"] & (E^3,d^3)\ar[r,"H"] & \cdots
\et,~~~~E^r=(\circ H)^{r-1}(E^1,d^1).\nn\eea
\item The dependence of ~$\deg d^r=\big(a_1(r),...,a_n(r)\big)\in\Integer^n$~ on $r$ is determined by the purpose/application of the $n$-spectral sequence.
\item Additional conditions may apply on the objects $E^r$, depending on their intended application.
\end{enumerate}

A \ul{spectral sequence} (\blue{footnote}\footnote{Note that the phrase ``\emph{spectral sequence}'' usually instead stands for ``\emph{bispectral sequence}''.}) is a $1$-spectral sequence (i.e., a sequence $(E^r,d^r)_{r\geq 1}$  of differential graded objects $(E^r,d^r)=(E_{\ast}^r,d_{\ast}^r)\in Gr(\A_0^\Integer)$ satisfying 1,2,3 above).

A \ul{bispectral sequence} (\blue{footnote}\footnote{A ``\emph{bispectral sequence}'' is usually simply called a ``\emph{spectral sequence}''.}) is a $2$-spectral sequence (i.e., a sequence $(E^r,d^r)_{r\geq 1}$  of differential bigraded objects $(E^r,d^r)=(E_{\ast\ast}^r,d_{\ast\ast}^r)\in Gr\big(\A_{00}^{\Integer^2}\big)$  satisfying 1,2,3 above).
\end{dfn}

\begin{rmk}[\blue{\cite[Theorem 10.11, p.622]{rotman2009}}]
By Propositions \ref{FiltExCoup1} and \ref{FiltExCoup2}, every filtration of a complex of $R$-modules gives a bispectral sequence as follows (where $A$ and $B$ are as defined in Proposition \ref{FiltExCoup1}):
\bea
&&\adjustbox{scale=0.8}{\bt
A\ar[rr,"{\al,{(1,-1)}}"] && A\ar[dl,"{\beta,{(0,0)}}"]\\
  & (B,d^1)\ar[ul,"{\gamma,{(-1,0)}}"] &
\et}~~~~\adjustbox{scale=0.8}{\bt
A^2\ar[rr,"{\al^2,{(1,-1)}}"] && A^2\ar[dl,"{\beta^2,{(-1,1)}}"]\\
  & (B^2,d^2)\ar[ul,"{\gamma^2,{(-1,0)}}"] &
\et}~~~~\cdots~~~~\adjustbox{scale=0.8}{\bt
A^{r+1}\ar[rr,"{\al^{r+1},{(1,-1)}}"] && A^{r+1}\ar[dl,"{\beta^{r+1},{(-r,r)}}"]\\
  & (B^{r+1},d^{r+1})\ar[ul,"{\gamma^{r+1},{(-1,0)}}"] &
\et}\nn\\
&&\textstyle A^{r+1}:=\im\al^r,~~~~B^{r+1}=H(B^r,d^r)={\ker d^r\over\im d^r},~~~~d^{r+1}=\beta^r\gamma^r,~~~~\deg d^{r+1}=(-r,r-1),\nn\\
&&\textstyle \al^{r+1}(\al^r(a)):=\al^r|_{A^{r+1}}(\al^r(a)),~~~~\beta^{r+1}(\al^r(a)):=\beta^r(a)+\im d^r,~~~~\gamma^{r+1}(z+\im d^r):=\gamma^r(z).\nn
\eea
\end{rmk}

\section{Distinguished Triangles: Mapping Cones and Mapping Cylinders}
This section is of central importance in the subsequent discussion on the triangulation of derived categories. Also, seemingly, the only relevant thing this section has in common with the previous sections of this chapter is the notion of \ul{exact triangles} (e.g., exact couples) in a category.

Throughout this short section, we will let $\A$ be an abelian category. As usual, most (if not all) of the concepts (along with results) for $R$-mod directly carry over to any $\A$ via the Freyd-Mitchell imbedding.

\begin{dfn}[\textcolor{blue}{Recall: \index{Translated chain complex}{Translated chain complex}, \index{Translation! functor}{Translation functor}}]
If $(A_\ast,d_\ast^{A_\ast})$ is a chain complex in $\A$ (i.e., $A_\ast\in\A_0^\Integer$), then for any $k\in\Integer$, its \ul{$k$-translation} ~$\big(A_\ast[k],d_\ast^{A_\ast[k]}\big)$~ is the chain complex given by
\bea
A_n[k]:=A_{n+k},~~~~d_n^{A_\ast[k]}:=(-1)^kd_{n+k}^{A_\ast}.\nn
\eea
The \ul{$k$-translation functor} is the functor ~$T_k:\A_0^\Integer\ra \A_0^\Integer$ ~$A_\ast\sr{f_\ast}{\ral}B_\ast~\mapsto~A_\ast[k]\sr{f_\ast[k]}{\ral}B_\ast[k]$, ~where
\bea
f_n[k]:=(-1)^kf_{n+k}:A_n[k]\ra B_n[k].\nn
\eea
\end{dfn}

\begin{dfn}[\blue{\index{Negative of a chain complex}{Recall: Negative of a chain complex}}]
Let $C_\ast=(C_\ast,d_\ast)\in\A_0^\Integer$. Then the negative of $C_\ast$ is the associated chain complex $-C_\ast=-(C_\ast,d_\ast):=(C_\ast,-d_\ast)$.
\end{dfn}

In the following recall of the definition of the mapping cone, we will for convenience (i) forward-shift the indices of $A_\ast,B_\ast$ (i.e., the vertical or $j$-index in the definition of $C_{\ast\ast}$) by one unit in the sense $j\mapsto j-1$, (ii) switch the positions of $A_\ast,B_\ast$, and (iii) change the sign of the mapping cone complex by $-1$. Consequently, we will also repeat both the statement and proof of Proposition \ref{ConeLESprp} in terms of the new convention.

\begin{dfn}[\textcolor{blue}{Recall: \index{Mapping cone}{Mapping cone}}]
Let $f_\ast:A_\ast\ra B_\ast$ be a chain morphism in $\A$. The \ul{mapping cone} of $f_\ast$ is defined as
\[
C(f_\ast):=Tot\left(0\ra A_\ast\sr{f_\ast}{\ral}B_\ast\ra 0\right)=
Tot\left(\adjustbox{scale=0.8}{\bt
        & \substack{(i=1~\txt{or}~A)\\~\\\vdots}                                               & \substack{(i=0~\txt{or}~B)\\~\\\vdots}                                     &  \\
0\ar[r] & \overbrace{A_{j-1}}^{c_{0j}}\ar[u]\ar[u,"d_{j-1}^{A_\ast}"']\ar[r,"f_{j-1}"] & \overbrace{B_{j-1}}^{c_{1j}}\ar[u,"d_{j-1}^{B_\ast}"']\ar[u]\ar[r] & 0 \\
0\ar[r] & A_j\ar[u,"d_j^{A_\ast}"']\ar[r,"f_j"]                & \overbrace{B_j}^{c_{1j}}\ar[u,"d_j^{B_\ast}"']\ar[r]             & 0 \\
0\ar[r] & A_{j+1}\ar[u,"d_{j+1}^{A_\ast}"']\ar[r,"f_{j+1}"]                & B_{j+1}\ar[u,"d_{j+1}^{B_\ast}"']\ar[r]             & 0 \\
        & \vdots\ar[u]                                         & \vdots\ar[u]                          &
\et}\right)=Tot(C_{\ast\ast}),
\]
where ~$C_{ij}:=\delta_{i1}A_{j-1}\oplus\delta_{i0}B_{j-1}$, ~$d^h_{ij}:=\delta_{i1}f_{j-1}$, ~$d^v_{ij}:=\delta_{i1}d^{A_\ast}_{j-1}+\delta_{i0}d^{B_\ast}_{j-1}$, ~and ( assuming $\A\hookrightarrow\txt{Sets}$ )
\bea
&&d_n^{C(f_\ast)}(c_{ij})_{i+j=n}:=\Big(d^h_{i+1,j}c_{i+1,j}+(-1)^id^v_{i,j+1}c_{i,j+1}\Big)_{i+j=n}=\Big(d^hc_{i+1,n-i}+(-1)^id^vc_{i,n-i+1}\Big)_{i\in\{0,1\}}\nn\\
&&~~~~=\Big(d^hc_{1,n}+d^vc_{0,n+1}~,~d^hc_{2,n-1}-d^vc_{1,n}\Big)=\Big(f_{n-1}a_{n-1}+d^{B_\ast}b_n~,~0-d^{A_\ast}a_{n-1}\Big).\nn
\eea
\end{dfn}
That is, the \ul{mapping cone} is $C(f_\ast):=B_\ast\oplus A_\ast[-1]$, where $C(f_\ast)_n:=B_n\oplus A_{n-1}$. As a diagram,
\bea\bt
C(f_\ast):\cdots\ar[r] & C(f_\ast)_{n+1}\ar[r,"d_{n+1}"] \ar[r] & C(f_\ast)_n\ar[r,"d_n"] & C(f_\ast)_{n-1}\ar[r,"d_{n-1}"] &\cdots
\et\nn
\eea
with the differential (in Sets-imbedded matrix form) given by
\bea
d_n^{C(f_\ast)}=
\left[
  \begin{array}{cc}
    d_n^B & f_{n-1} \\
    0 & -d_{n-1}^A\\
  \end{array}
\right]:
\left[
  \begin{array}{l}
    b_n \\
    a_{n-1} \\
  \end{array}
\right]\longmapsto
\left[
  \begin{array}{l}
    d_n^Bb_n + f_{n-1}a_{n-1} \\
    -d_{n-1}^Aa_{n-1} \\
  \end{array}
\right].\nn
\eea

\begin{prp}[\textcolor{blue}{Recall: \index{Quasiisomorphism criterion}{Quasiisomorphism criterion}}]\label{ConeLES_prp}
Let $f_\ast:A_\ast\ra B_\ast$ be a chain morphism in $R$-mod.
\begin{enumerate}[leftmargin=0.8cm]
\item There is a SES of chain complexes ~$0\ra B_\ast\sr{i}{\ral}C(f_\ast)\sr{p}{\ral}A_\ast[-1]\ra 0$,~ where $i,p$ are the obvious inclusion and projection respectively.
\item Furthermore, in the associated induced LES of homologies,
\bea
\cdots\ra H_{n+1}\big(A_\ast[-1]\big)\sr{\delta_{n+1}}{\ral}H_n(B_\ast)\ra H_n(C(f_\ast))\ra H_n\big(A_\ast[-1]\big)\sr{\delta_n}{\ral} \cdots\nn
\eea
the connecting homomorphism $\delta_n=H_n(f_\ast)$, i.e., $\delta$ is induced by $f_\ast$, and $H_{n+1}\big(A_\ast[-1]\big)=H_n(A_\ast)$.

That is, the sequence takes the form
{\small
\bea
\label{ConeLES1}
\bt\cdots\ar[r]& H_{n+1}(C(f))\ar[r,"{H_{n+1}(p)}"]& H_n(A_\ast)\ar[r,"{H_{n+1}(f)}"]& H_n(B_\ast)\ar[r,"{H_n(i)}"]& H_n(C(f_\ast))\ar[r,"H_n(p)"]& H_{n-1}(A_\ast)\ar[r,"{H_n(f_\ast)}"]& \cdots\et
\eea}
\item Hence, $f_\ast:A_\ast\ra B_\ast$ is a quasi-iso $\iff$ $H_n(C(f_\ast))=0$ for all $n$, $\iff$ $C(f_\ast)$ is exact.
\end{enumerate}
\end{prp}
\begin{proof}
This is the same as the proof of Proposition \ref{ConeLESprp}
{\flushleft (1)} The rows are split exact. Thus, we just need to check directly that the diagram commutes.
\[\adjustbox{scale=0.9}{\bt
 & \vdots\ar[d] & \vdots\ar[d] & \vdots\ar[d] & \\
0\ar[r] & B_{n+1}\ar[d,"d^B_{n+1}"]\ar[r,"i_{n+1}"] & B_{n+1}\oplus A_n\ar[d,dashed,"d_{n+1}^{C(f)}"]\ar[r,dashed,"p_{n+1}"] & A_n\ar[d,"-d_n^A"]\ar[r,dashed] & 0\\
0\ar[r,dashed] & B_n\ar[d,"d^B_n"]\ar[r,dashed,"i_n"] & B_n\oplus A_{n-1}\ar[d,"d_n^{C(f)}"]\ar[r,"p_n"] & A_{n-1}\ar[d,"-d_{n-1}^A"]\ar[r] & 0 \\
0\ar[r] & B_{n-1}\ar[d,"d^B_{n-1}"]\ar[r,"i_{n-1}"] & B_{n-1}\oplus A_{n-2}\ar[d,"d_{n-1}^{C(f)}"]\ar[r,"p_{n-1}"] & A_{n-2}\ar[d,"-d_{n-2}^A"]\ar[r] & 0 \\
  & \vdots & \vdots & \vdots &
\et}\]

\begin{align}
& d_{n+1}^{C(f)}i_{n+1}(b_{n+1})=d_{n+1}^{C(f)}(b_{n+1},0)=\big(d_{n+1}^{B_\ast}b_{n+1}+f_{n}0~,~-d_n^{A_\ast}0\big)
=\big(d^{B_\ast}b_{n+1},0\big)=i_nd_{n+1}^{B_\ast}b_{n+1},\nn\\
&~~\Ra~~d_{n+1}^{C(f)}i_{n+1}=i_nd_{n+1}^{B_\ast},\nn\\
& -d_n^Ap_{n+1}(b_{n+1},a_n)=-d_n^Aa_n=p_n(d_{n+1}^Bb_{n+1} + f_{n}a_{n},-d_{n}^Aa_{n})
=p_nd_{n+1}^{C(f)}(b_{n+1},a_n),\nn\\
&~~\Ra~~-d_n^Ap_{n+1}=p_nd_{n+1}^{C(f)}.\nn
\end{align}

{\flushleft (2)} By definition, the connecting homomorphism $\delta_{n+1}$ is given (for $a_{n}\in\ker(-d_{n}^A)$) by
{\footnotesize\[
\delta_{n}(a_{n}+\im(-d_{n+1}^A))=b_{n}+\im d_{n+1}^B,~~\txt{such that}~~i_n(b_n)=d_{n+1}^{C(f)}(b_{n+1},a_n),~~~~p_{n+1}(b_{n+1},a_n)=a_n,\nn
\]}
where the last condition is clear, and ~$i_n(b_n)=d_{n+1}^{C(f)}(b_{n+1},a_n)$~ implies

{\small\begin{align}
&(b_{n},0)=(d_{n+1}^Bb_{n+1} + f_{n}a_{n},-d_{n}^Aa_{n})=(d_{n+1}^Bb_{n+1} + f_{n}a_{n},0),\nn\\
&\Ra~~b_{n}=d_{n+1}^Bb_{n+1} + f_{n}a_{n},\nn\\
&\Ra~~\delta_{n}(a_{n}+\im(-d_{n+1}^A))=b_{n}+\im d_{n+1}^B
=f_{n}a_{n}+\im d_{n+1}^B=H_n(f)(a_{n}+\im(-d_{n+1}^A)).\nn
\end{align}}

Also, we have
\[
\textstyle H_{n+1}\big(A_\ast[-1]\big)={\ker d_{n+1}^{A_\ast[-1]}\over\im d_{n+2}^{A_\ast[-1]}}={\ker(- d_{n}^{A_\ast})\over\im(-d_{n+1}^{A_\ast})}={\ker d_{n}^{A_\ast}\over\im d_{n+1}^{A_\ast}}=H_n(A_\ast).
\]
{\flushleft (3)} From the exact sequence in part (2), $H_n(f)$, for all $n$, is an isomorphism if and only if $H_n(C(f))=0$ for all $n$. Hence $f$ is a quasi-iso if and only if $H_n(C(f))=0$ for all $n$, if and only if $C(f)$ is exact.
\end{proof}

\begin{dfn}[\textcolor{blue}{\index{Mapping cylinder}{Mapping cylinder}}]
Let $f_\ast:A_\ast\ra B_\ast$ be a chain morphism in $R$-mod, and $p_\ast^A:C(f_\ast)[1]\ra A_\ast$ the obvious projection. The \ul{mapping cylinder} of $f_\ast$ is
\bea
\label{mapping_cyl_eqn1}Cyl(f_\ast):= C\Big(C(f_\ast)[1]\sr{p_\ast^A}{\ral}A_\ast\Big)=Tot\Big(0\ra C(f_\ast)[1]\sr{p_\ast^A}{\ral}A_\ast\ra 0\Big).
\eea
\end{dfn}
That is, ~$Cyl(f_\ast)_n:=A_n\oplus C(f_\ast)[1]_{n-1}=A_n\oplus C(f_\ast)_n=A_n\oplus B_n\oplus A_{n-1}$, ~and so
\bea
Cyl(f_\ast):=A_\ast\oplus C(f_\ast)=A_\ast\oplus B_\ast\oplus A_\ast[-1],\nn
\eea
and the differential (in Sets-embedded matrix form) is given by
{\footnotesize\begin{align}
d_n^{Cyl(f_\ast)}:=\left[
               \begin{array}{cc}
                 d_n^A & p^A_{n-1} \\
                 0 & -d_{n-1}^{C(f_\ast)[1]} \\
               \end{array}
             \right]=
             \left[
               \begin{array}{cc}
                 d_n^A & p^A_{n-1} \\
                 0 & d_n^{C(f_\ast)} \\
               \end{array}
             \right]
=\left[
         \begin{array}{ccc}
           d_n^A & 0 & id_{A_{n-1}} \\
           0 & d_n^B & f_{n-1} \\
           0 & 0 & -d_{n-1}^A \\
         \end{array}
       \right]~:~\left[
        \begin{array}{c}
          a_n \\
          b_n \\
          a_{n-1} \\
        \end{array}
      \right]\longmapsto
      \left[
        \begin{array}{l}
          d_n^Aa_n+a_{n-1} \\
          d_n^Bb_n+f_{n-1}(a_{n-1}) \\
          -d_{n-1}^Aa_{n-1} \\
        \end{array}
      \right].\nn
\end{align}}

\begin{rmk}[\textcolor{blue}{Equivalent description of the mapping cylinder}]\label{EquivConvRmk}
Observe that, by swapping positions of $A_\ast$ and $B_\ast$ (hence also swapping 1st,2nd rows, and swapping 1st,2nd columns, in the matrix of $d_n^{Cyl(f_\ast)}$), we can rewrite the cylinder and its differential as follows:
{\small
\begin{align}
\label{mapping_cyl_eqn2}& Cyl(f_\ast)=A_\ast\oplus B_\ast\oplus A_\ast[-1]\sr{\txt{rearrange}}{\ral}B_\ast\oplus A_\ast\oplus A_\ast[-1]=B_\ast\oplus C(1_A)=C\left(C(1_A)[1]\sr{[0,f_\ast]}{\ral}B_\ast\right),\\
&\nn\\
&d_n^{Cyl(f_\ast)}=\left[
         \begin{array}{ccc}
           d_n^A & 0 & id_{A_{n-1}} \\
           0 & d_n^B & f_{n-1} \\
           0 & 0 & -d_{n-1}^A \\
         \end{array}
       \right]\sr{\txt{rearrange}}{\ral}
\left[
         \begin{array}{ccc}
           d_n^B & 0 & f_{n-1} \\
           0 & d_n^A & id_{A_{n-1}} \\
           0 & 0 & -d_{n-1}^A \\
         \end{array}
       \right]=
\left[
         \begin{array}{cc}
           d_n^B & [0,f_{n-1}] \\
           0 &  d_n^{C(1_A)} \\
         \end{array}
       \right]=
\left[
         \begin{array}{cc}
           d_n^B & [0,f_{n-1}] \\
           0 &  -d_n^{C(1_A)[1]} \\
         \end{array}
       \right],\nn
\end{align}}where
~$d_n^{C(1_A)}:=
     \left[
         \begin{array}{cc}
           d_n^A & id_{A_{n-1}} \\
           0 & -d_{n-1}^A \\
         \end{array}
       \right]$~ and ~$[0,f_{n-1}](a_n,a_{n-1}):=f_{n-1}a_{n-1}$,~ for all $(a_n,a_{n-1})\in C(1_A)[1]_{n-1}$.

Henceforth, an \ul{equivalent convention for the mapping cylinder} is
{\small\begin{align}
& Cyl(f_\ast):=C\left(C(1_A)[1]\sr{[0,f_\ast]}{\ral}B_\ast\right)=B_\ast\oplus C(1_A)=B_\ast\oplus A_\ast\oplus A_\ast[-1],\nn\\
&\nn\\
&d_n^{Cyl(f_\ast)}:=
\left[
         \begin{array}{cc}
           d_n^B & [0,f_{n-1}] \\
           0 &  d_n^{C(1_A)} \\
         \end{array}
       \right]
=\left[
         \begin{array}{ccc}
           d_n^B & 0 & f_{n-1} \\
           0 & d_n^A & id_{A_{n-1}} \\
           0 & 0 & -d_{n-1}^A \\
         \end{array}
       \right]~:~\left[
        \begin{array}{c}
          b_n \\
          a_n \\
          a_{n-1} \\
        \end{array}
      \right]\longmapsto
      \left[
        \begin{array}{l}
          d_n^Bb_n+f_{n-1}(a_{n-1}) \\
          d_n^Aa_n+a_{n-1} \\
          -d_{n-1}^Aa_{n-1} \\
        \end{array}
      \right].\nn
\end{align}}where
~$d_n^{C(1_A)}:=
     \left[
         \begin{array}{cc}
           d_n^A & id_{A_{n-1}} \\
           0 & -d_{n-1}^A \\
         \end{array}
       \right]$.
\end{rmk}

\begin{lmm}[\textcolor{blue}{Quasiisomorphism/homotopy-equivalence between Codomain and Cylinder}]\label{CylQisLmm}
Let $f_\ast:A_\ast\ra B_\ast$ be a chain morphism in $R$-mod. (1) There is a SES of chain complexes ~$0\ra B_\ast\sr{i}{\ral}Cyl(f_\ast)\sr{p}{\ral}C(1_A)\ra 0$. Moreover, (2) $i$ is a quasi-iso, and in fact, (3) $i$ is a homotopy equivalence.
\end{lmm}
\begin{proof}
(1) With the \ul{equivalent convention} in Remark \ref{EquivConvRmk}, i.e.,
\[
C(f_\ast)=B\oplus A_\ast[-1],~~~~Cyl(f_\ast)=C([0,f_\ast])=B\oplus C(1_A)=B\oplus A\oplus A_\ast[-1],
\]
it follows (by Proposition \ref{ConeLES_prp}) that we have a SES
\bc\adjustbox{scale=0.9}{\bt[row sep=tiny]
0\ar[r] & B_\ast\ar[r,"i"] & Cyl(f_\ast)\ar[d,equal]\ar[r,"p"] & C(1_A)\ar[r] & 0\\
 & & C\Big(C(1_A)[1]\sr{{[0,f_\ast]}}{\ral}B_\ast\Big) & &
\et}\ec
{\flushleft (2)} The above SES in turn gives the induced LES of homologies
{\small\bea\bt[row sep=tiny]
\cdots\ar[r] & H_n\big(C(1_A)[1]\big)\ar[r,"{H_n({[0,f_\ast]})}"] & H_n(B_\ast)\ar[r,"{H_n(i)}"] & H_n(Cyl(f_\ast))\ar[r,"{H_n(p)}"] & H_{n-1}\big(C(1_A)[1]\big)\\
 & \ar[r,"{H_{n-1}({[0,f_\ast]})}"] & H_{n-1}(B_\ast)\ar[r,"{H_{n-1}(i)}"] & H_{n-1}(Cyl(f_\ast))\ar[r,"{H_{n-1}(p)}"] & H_{n-2}\big(C(1_A)[1]\big)\ar[r] & \cdots
 \et\nn
\eea}Note that $C(1_A)$ is exact, since $1_A$ is an iso (hence a quasi-iso). (\blue{footnote}\footnote{By the definition of $C(1_A)$, we also have a SES~~ $0\ra A_\ast\sr{i'}{\ral}C(1_A)\sr{p'}{\ral}A_\ast[-1]\ra 0$. We can also obtain exactness of $C(1_A)$ directly as follows:
$C(1_A)=A\oplus A_\ast[-1]$ is exact because
\bea
&&\left[
  \begin{array}{cc}
    d_n^A & 1_{A_{n-1}} \\
    0 & -d_{n-1}^A \\
  \end{array}
\right]\left[
  \begin{array}{c}
    a_n  \\
    a_{n-1}\\
  \end{array}
\right]=\left[
  \begin{array}{c}
    d_n^Aa_n+a_{n-1}  \\
    -d_{n-1}^Aa_{n-1}\\
  \end{array}
\right]=
\left[
  \begin{array}{c}
    0  \\
    0\\
  \end{array}
\right]
~~\Ra~~
\left[
  \begin{array}{c}
    a_n  \\
    a_{n-1}\\
  \end{array}
\right]=
\left[
  \begin{array}{c}
    a_n  \\
    -d_n^Aa_n\\
  \end{array}
\right]=
\left[
  \begin{array}{cc}
    d_{n+1}^A & 1_{A_n} \\
    0 & -d_n^A \\
  \end{array}
\right]\left[
  \begin{array}{c}
    0  \\
    a_n\\
  \end{array}
\right],\nn\\\nn\\
&&~~\Ra~~\ker d_n^{C(1_A)}\subset\im d_{n+1}^{C(1_A)}.\nn
\eea
}). Thus, $H_n\big(C(1_A)[1]\big)=0$, and so the LES above becomes
\bea
\cdots\ra 0\ra H_n(B_\ast)\sr{H_n(i)}{\ral} H_n(Cyl(f_\ast))\ra 0\ra H_{n-1}(B_\ast)\sr{H_{n-1}(i)}{\ral}H_{n-1}(Cyl(f_\ast))\ra 0\ra \cdots,\nn
\eea
which shows ~$B_\ast\sr{i}{\ral}Cyl(f_\ast)$~ is a quasi-iso.

{\flushleft (3)} To show ~$B_\ast\sr{i}{\ral}Cyl(f_\ast)=B_\ast\oplus A_\ast\oplus A_\ast[-1],~b_n\mapsto (b_n,0,0)$~ is a homotopy equivalence, define
\bea
Cyl(f_\ast)\sr{j}{\ral}B_\ast,~~(b_n,a_n,a_{n-1})\mapsto b_n-f_n(a_n)=b_n+f_n(-a_n),\nn
\eea
Then we can check that $i,j$ are chain morphisms, $ji=1_{B_\ast}$, and $ij\simeq 1_{Cyl(f_\ast)}$ via the homotopy
\bea
s_n=\left[
      \begin{array}{ccc}
        0 & 0 & 0 \\
        0 & 0 & 0 \\
        0 & -1 & 0 \\
      \end{array}
    \right]: Cyl(f_\ast)_n\ra Cyl(f_\ast)_{n+1},~~(b_n,a_n,a_{n-1})\mapsto (0,0,-a_n).\nn
\eea
Explicitly, $j_ni_n=id_{B_n}$ and {\small $i_nj_n=id_{Cyl(f_\ast)_n}+d_{n+1}^{Cyl(f_\ast)}s_n+s_{n-1}d_n^{Cyl(f_\ast)}$}, where
{\scriptsize $d_n^{Cyl(f_\ast)}=\left[
         \begin{array}{ccc}
           d_n^B & 0 & f_{n-1} \\
           0 & d_n^A & id_{A_{n-1}} \\
           0 & 0 & -d_{n-1}^A \\
         \end{array}
       \right]$}.
\end{proof}

\begin{lmm}[\textcolor{blue}{\index{Triangle diagram}{Triangle diagram}}]\label{TrngDiagLmm}
Let $f_\ast:A_\ast\ra B_\ast$ be a chain morphism in $R$-mod. There exist \ul{SES's} of chain complexes, $0\ra B_\ast\ra C(f_\ast)\ra A_\ast[-1]\ra 0$ and $0\ra A_\ast\ra Cyl(f_\ast)\ra C(f_\ast)\ra 0$, that fit into a \ul{commutative diagram} of the following form (\blue{footnote}\footnote{If the lower-right square is included, then the whole diagram commutes up to homotopy because $p_{1,3}i=i_1$ while $i,j$ are homotopy inverses of each other, and so $i_1j\simeq p_{1,3}$.}):
\bea
\label{triangle-diag}\adjustbox{scale=0.9}{\bt
         & 0\ar[r]  & B_\ast\ar[d,"i"]\ar[r,"i_1"] & C(f_\ast)\ar[d,equal]\ar[r,"{p}_2"] & A_\ast[-1]\ar[r] & 0 \\
 0\ar[r] & A_\ast\ar[d,"-id_{A_\ast}"']\ar[r,"i_2"] & Cyl(f_\ast)\ar[d,shift right,"j"']\ar[r,"{p}_{1,3}"] & C(f_\ast)\ar[d,equal,dashed,"\txt{up to hty}"]\ar[r] & 0 \\
  & A_\ast\ar[r,"f_\ast"] & B_\ast\ar[u,shift right,dashed,"i"']\ar[r,dashed,"i_1","\txt{up to hty}"'] & C(f_\ast) &
\et}
\eea
such that ~$ji=id_{B_\ast}$~ and ~$ij\simeq id_{Cyl(f_\ast)}$ (i.e., $i,j$ are \ul{homotopy inverses}). (\blue{footnote}\footnote{Note that the lower row is an underlying element of the LES in (\ref{ConeLES1}).})

Moreover, we have the following \ul{nullhomotopic composition} (\blue{footnote}\footnote{This nullhomotopic composition shows that the natural choice for a distinguished triangle (next def.) in the homotopy category $\H\A:={\A_0^\Integer\over\simeq}$ is the following sequence obtained from (\ref{triangle-diag}): ~$A_\ast\sr{f_\ast}{\ral}B_\ast\sr{i_1}{\ral}C(f_\ast)\sr{{p}_2}{\ral}A_\ast[-1]$, ~because in $\A_0^\Integer$ (resp. in $\H\A$), the above sequence is quasisomorphic (resp. isomorphic) to the sequence ~$A_\ast\sr{i_2}{\ral}Cyl(f_\ast)\sr{{p}_{1,3}}{\ral}C(f_\ast)\sr{{p}_2}{\ral}A_\ast[-1]$.
}):
\bea
(i_1f)_n=d^{C(f_\ast)}_{n+1}h_n+h_{n-1}d_n^A,~~~~\txt{where}~~h_n:A_n\ra C_{n+1}(f_\ast),~~h_n(a_n):=(0,a_n).\nn
\eea
\end{lmm}
\begin{proof}
The top SES is Proposition \ref{ConeLES_prp}(1). Next, consider the equivalent convention (Remark \ref{EquivConvRmk})
\bea
C(f_\ast)=B\oplus A_\ast[-1],~~~~Cyl(f_\ast)=C([0,f_\ast])=B\oplus C(1_A)=B\oplus A\oplus A_\ast[-1].\nn
\eea
Then the middle SES is also clear by Proposition \ref{ConeLES_prp}(1), the definition of $Cyl(f_\ast)$ in (\ref{mapping_cyl_eqn1}) or (\ref{mapping_cyl_eqn2}), and the fact that $C(f_\ast)[1][-1]=C(f_\ast)$. Also, in the equivalent convention of Remark \ref{EquivConvRmk}, the maps shown (easily verified to be \ul{chain morphisms}) are the obvious ones:
\bea
&&i_1(b_n):=(b_n,0),~~i_2(a_n):=(0,a_n,0),~~{p}_{1,3}(b_n,a_n,a_{n-1}):=(b_n,a_{n-1}),~~{p}_2(b_n,a_{n-1}):=a_{n-1}\nn\\
&& i(b_n):=(b_n,0,0),~~~~j(b_n,a_n,a_{n-1}):=b_n-f_n(a_n).\nn
\eea
It now remains to directly verify (with the above definitions) that the diagram commutes.

The properties of $i,j$ (i.e., $ji=id_{B_\ast}$ and $ij\simeq id_{Cyl(f_\ast)}$) come from the proof of Lemma \ref{CylQisLmm}(3), while the property $i_1f\simeq 0$ can be verified directly using the given homotopy {\small $h_n:A_n\ra C_{n+1}(f_\ast),~a_n\mapsto (0,a_n)$}.
\end{proof}

\begin{dfn}[\textcolor{blue}{
\index{Triangle in a category}{Triangle in a category},
\index{Morphism of! triangles}{Morphism of triangles},
\index{Isomorphism of! triangles}{Isomorphism of triangles},
\index{Category of! triangles}{Category of triangles},
\index{Distinguished triangle in $(R\txt{-mod})_0^\Integer$}{Distinguished triangle in $(R\txt{-mod})_0^\Integer$}: \cite[Definition 4, pp 156-157]{gelfand-manin2010}}]~
\bit[leftmargin=0.8cm]
\item A \ul{triangle} in $\A_0^\Integer$ is a sequence of the form ~$\Delta:~A_\ast\sr{u}{\ral}B_\ast\sr{v}{\ral}C_\ast\sr{w}{\ral}A_\ast[-1]$.
\item A \ul{morphism of triangles} $\eta:\Delta\ra \Delta'$ is a commutative diagram of the form
\bc\adjustbox{scale=0.9}{\bt
\Delta\ar[d,"\eta"] & A_\ast\ar[d,"f"]\ar[r,"u"] & B_\ast\ar[d,"g"]\ar[r,"v"] & C_\ast\ar[d,"h"]\ar[r,"w"] & A_\ast[-1]\ar[d,"{f[-1]}"]\\
\Delta' & A'_\ast\ar[r,"u'"]       & B'_\ast\ar[r,"v'"]       & C'_\ast\ar[r,"w'"]       & A'_\ast[-1]
\et}\ec
The morphism of triangles $\eta:\Delta\ra \Delta'$ is an \ul{isomorphism of triangles} (making $\Delta,\Delta'$ \ul{isomorphic triangles}, written $\Delta\cong \Delta'$) if $f,g,h$ are isomorphisms.

\item The \ul{category of triangles} in $\A_0^\Integer$ is the category $\Delta(\A_0^\Integer)$ whose objects are triangles and whose morphisms are morphisms of triangles.

\item In $(R\txt{-mod})_0^\Integer$, a triangle ~$\Delta:~U_\ast\sr{u}{\ral}V_\ast\sr{v}{\ral}W_\ast\sr{w}{\ral}U_\ast[-1]$ ~is a \ul{distinguished triangle} if there exists a chain morphism $A_\ast\sr{f_\ast}{\ral}B_\ast$ such that ${\Delta}\cong \Delta_{f_\ast}:~A_\ast\sr{i_2}{\ral}Cyl(f_\ast)\sr{{p}_{1,3}}{\ral}C(f_\ast)\sr{{p}_2}{\ral}A_\ast[-1]$.
\eit
\end{dfn}
\begin{note}
In the commutative diagram (\ref{triangle-diag}), since ${p}_2{p}_{1,3}\simeq 0$ and ${p}_{1,3} i_2\simeq0$, it follows that in any distinguished triangle ~$\Delta:~A_\ast\sr{u}{\ral}B_\ast\sr{v}{\ral}C_\ast\sr{w}{\ral}A_\ast[-1]$~ in $(R\txt{-mod})_0^\Integer$, we have
\[
wv\simeq 0~~~~\txt{and}~~~~vu\simeq 0.
\]

\end{note}
\begin{rmk}[\textcolor{blue}{\index{Induced! LES of homologies of a distinguished triangle}{Induced LES of homologies of a distinguished triangle}}]\label{ModLesForDT}
A \ul{distinguished triangle} in $(R\txt{-mod})_0^\Integer$ is essentially a generalization of the cyclic sequence contained in (\ref{triangle-diag}) of the form
\bea
\Delta_f:~A_\ast\ral Cyl(f_\ast)\ral C(f_\ast)\ral A_\ast[-1].\nn
\eea
By the preceding results leading to the following associated long exact sequence from (\ref{ConeLES1}):
{\footnotesize\begin{align}
&\cdots\ral H_{n+1}(A_\ast)\ral H_{n+1}(B_\ast)\ral H_{n+1}(C(f_\ast))\ral\nn\\
&~~~~~~~~\ral H_n(A_\ast)~\ral~H_n(B_\ast)~\ral~H_n(C(f_\ast))\ral\nn\\
&~~~~~~~~~~~~~\ral H_{n-1}(A_\ast)\ral H_{n-1}(B_\ast)\ral H_{n-1}(C(f_\ast))\ral \cdots,\nn
\end{align}}the distinguished triangle $\Delta_f$ has an \ul{induced LES of homologies} in the compact form
\[\adjustbox{scale=0.9}{\bt
H_n(A_\ast)\ar[d,equal]\ar[r] & H_n(Cyl(f_\ast))\ar[d,"\cong"]\ar[r] & H_n(C(f_\ast))\ar[d,equal]\ar[r] & H_n(A_\ast[-1])\ar[d,equal]\\
H_n(A_\ast)\ar[r] & H_n(B_\ast)\ar[r] & H_n(C(f_\ast))\ar[r] & H_{n-1}(A_\ast)
\et}\]
Hence, every distinguished triangle in $(R\txt{-mod})_0^\Integer$ has an associated induced LES of homologies.
\end{rmk}

\begin{prp}[\textcolor{blue}{\cite[Proposition 5, p.157]{gelfand-manin2010}}]
In $R$-mod, every SES of complexes $0\ra A_\ast\sr{f_\ast}{\ral}B_\ast\sr{g_\ast}{\ral} C_\ast\ra 0$ is quasiisomorphic to a SES of the form $0\ra A_\ast\sr{{p}_2}{\ral}Cyl(f_\ast)\sr{{p}_{1,3}}{\ral} C(f_\ast)\ra 0$ (i.e., the middle row in (\ref{triangle-diag}) for the given monomorphism $f_\ast:A_n\ra B_\ast$ ). (\blue{footnote}\footnote{It follows that, in $(R\txt{-mod})_0^\Integer$ or $\H(R\txt{-mod})$, the induced LES of homologies for the SES $0\ra A_\ast\sr{f_\ast}{\ral}B_\ast\sr{g_\ast}{\ral} C_\ast\ra 0$ is precisely the induced LES of homologies for the associated distinguished triangle $A_\ast\sr{f_\ast}{\ral}B_\ast\sr{i}{\ral}C(f_\ast)\sr{p}{\ral}A_\ast[-1]$.})
\end{prp}
\begin{proof}
Consider the following diagram (and assume the equivalent notation of Remark \ref{EquivConvRmk}):
\[\adjustbox{scale=0.9}{\bt
0\ar[r] & A_\ast\ar[r,"f_\ast"] & B_\ast\ar[r,"g_\ast"] & C_\ast\ar[r] & 0\\
0\ar[r] & A_\ast\ar[u,"-id_A"']\ar[r,"i_2"] & Cyl(f_\ast)\ar[u,"j"']\ar[r,"{p}_{1,3}"] & C(f_\ast)\ar[u,"\theta"']\ar[r] & 0
\et}\]
where ~$i_2(a_n):=(0,a_n,0)$, ~$j(b_n,a_n,a_{n-1}):=b_n-f(a_n)$, ~${p}_{1,3}(b_n,a_n,a_{n-1}):=(b_n,a_{n-1})$,~ and
\bea
\theta(b_n,a_{n-1}):=g\big(b_n-f(a_n)\big)=g(b_n),~~~~\txt{(easily verified to be a \ul{chain morphism})}.
\eea
By construction, the above diagram commutes, and $-id_A$ and $j$ are quasi-iso's. Hence, by applying the LES of homology (and using its functoriality) we see that $H_n(\theta)$ is an iso, i.e., $\theta$ is a quasi-isomorphism.
\end{proof}

\section{Contractible Complexes of Modules (a digression)}
Throughout this section $\A=R\txt{-mod}$ for some ring $R$. Observe that a sequence in $\A$ of the form $0\ra A\ra A\oplus C\ra C\ra 0$ is the direct sum of the two sequences
\bea
0\ra A\ra A\ra 0\ra 0,~~~~0\ra 0\ra C\ra C\ra 0.\nn
\eea
\begin{prp}
Let $X_\ast$ be a complex of $R$-modules. The following are equivalent.
\bit
\item[(i)] $X_\ast$ is contractible.
\item[(ii)] $1_{X_\ast}\simeq 0_{X_\ast}$ (i.e., $1_{X_\ast}=dh+hd$ for a homotopy $h:X_\ast\ra X_{\ast+1}$)
\item[(iii)] $X_\ast$ is exact and \ul{totally split} (in the sense for each $n$, $Z_n=B_n$ and {\small $0\ra Z_n\sr{i}{\ral}X_n\sr{\del_n}{\ral}B_{n+1}\ra 0$} is split).
\item[(iv)] $X_\ast$ is isomorphic to a direct sum of shifts of complexes of the form $0\ra A\ra A\ra 0$.
\eit
\end{prp}
\begin{proof}
{\flushleft \ul{(i)$\iff$(ii)}}: This is Proposition \ref{CompContrPrp}. \ul{(ii)$\Ra$(iii)}: This is clear.

{\flushleft \ul{(iii)$\Ra$(ii)}}: Recall that if $X_\ast$ is exact, we have $B_{n+1}:=\im d_{n+1}=\ker d_n=:Z_n$ for all $n$.
\[\adjustbox{scale=0.7}{\bt
        &   &                                    & 0\ar[dr]                            &                        &        0         &      & 0\ar[dr]                            &                              & 0 \\
        &    &                                    &                                     & Z_n\ar[dr,hook]\ar[ur] &                  &      &                                     &  Z_{n-2}\ar[ur]\ar[dr,hook]  & ~ \\
\cdots \ar[r] &  X_{n+2}\ar[rr]\ar[dr,"d"] &                            & X_{n+1}\ar[ur,"d"]\ar[rr,"d_{n+1}"] &                        & X_n\ar[rr,"d_n"]\ar[dr,"d"] &      & X_{n-1}\ar[ur,"d"]\ar[rr,"d_{n-1}"] &                              & X_{n-2}\ar[r]  & \cdots\\
        &     & Z_{n+1}\ar[ur,hook]\ar[ur]\ar[dr]  &                                     &                        &                  &   Z_{n-1}\ar[ur,hook]\ar[dr]   &                                     &    & ~ \\
        &   0\ar[ur]  &                                    &     0                               &                        &       0\ar[ur]   &      &   0                                  &   & ~ \\
\et}\]
Thus, if $X_\ast$ is also totally split, then $X_\ast$ is isomorphic to the sequence
\bc\bt
 \cdots\ar[r] & Z_{n+1}\oplus Z_n\ar[dd,"1_{X_{n+1}}"']\ar[r,"d_{n+1}"] & Z_n\oplus Z_{n-1}\ar[dd,"1_{X_n}"']\ar[ddl,dashed,"h_n"']\ar[r,"d_n"] & Z_{n-1}\oplus Z_{n-2}\ar[dd,"1_{X_{n-1}}"]\ar[ddl,dashed,"h_{n-1}"']\ar[r,"d_{n+1}"] & \cdots \\
 &   &   &  &  \\
 \cdots\ar[r] & Z_{n+1}\oplus Z_n\ar[r,"d_{n+1}"] & Z_n\oplus Z_{n-1}\ar[r,"d_n"] & Z_{n-1}\oplus Z_{n-2}\ar[r,"d_{n+1}"] & \cdots
\et\ec
where by construction $d_{n+1}|_{0\oplus Z_n}:0\oplus Z_n\ra Z_n\oplus 0$ is an isomorphism for each $n$. Therefore, with $d_n(z_n,z_{n-1}):=(d_nz_{n-1},0)$ and $h_n(z_n,z_{n-1}):=(0,d_{n+1}^{-1}z_n)$, we get
\bea
1_{X_n}=d_{n+1}h_n+h_{n-1}d_n,~~~~\Ra~~~~1_{X_\ast}=hd+dh.\nn
\eea

{\flushleft \ul{(iii)$\iff$(iv)}}: This is also clear.
\end{proof}

\begin{crl}\label{CCMcrl3}
If $P_\ast:\cdots\ra P_2\ra P_1\ra P_0\ra 0$ is an \ul{exact} complex of projective $R$-modules, then $P_\ast$ is \ul{totally split} (hence \ul{contractible}).
\end{crl}

\begin{proof}
As usual, let $Z_n:=\ker \del_n^P$ and $B_n:=\im\del_{n+1}^P$. The exact sequence $0 \ra Z_1 \ra P_1\sr{\del_1}{\ral} P_0=B_0\ra 0$ splits because $P_0$ is projective. By induction on $n$, assume the exact sequence $0 \ra Z_n \ra P_n \sr{\del_n}{\ral} B_{n-1} \ra 0$ splits. Then $P_n = Z_n\oplus \wt{B}_{n-1}\cong Z_n\oplus B_{n-1}$ (where $\wt{B}_{n-1}\cong B_{n-1}$). Therefore $Z_n$ is projective (as a direct summand of a projective module). Moreover, since the exactness of $P_\ast$ implies $Z_n = B_n$ (which is projective), the exact sequence $0 \ra Z_{n+1} \ra P_{n+1} \sr{\del_{n+1}}{\ral} B_n = Z_n \ra 0$ splits as well.
\end{proof}
\begin{prp}
If $X_\ast:\cdots\ra X_1\ra X_0\ra 0$ and $Y_\ast:\cdots\ra Y_1\ra Y_0\ra 0$ are complexes of projective $R$-modules and $X_\ast\sr{f}{\ral}Y_\ast$ is a quasi-isomorphism, then $f$ is a homotopy equivalence.
\end{prp}
\begin{proof} The mapping cone $C(f):\cdots\ra Y_2\oplus X_1\sr{\del_2}{\ral}Y_1\oplus X_0\sr{\del_1}{\ral}Y_0\oplus 0\ra 0$ is a complex of projective modules \ul{bounded} in the same way as $X_\ast,Y_\ast$. Also, since $f$ is a quasi-isomorphism we know $C(f)$ is exact. Therefore, by Corollary \ref{CCMcrl3}, $C(f)$ is totally-split/contractible (with $1_{C(f)} = h\del+\del h$ for a homotopy $h:C(f)_\ast\ra C(f)_{\ast+1}$). It follows that in each degree, $h$ is a map from one direct sum to another, and so we can rewrite $h$ in matrix form as
\bea
h=\left[
  \begin{array}{cc}
    s & u \\
    t & v \\
  \end{array}
\right]:C(f)_\ast=Y_\ast\oplus X_\ast[-1]\ra C(f)_\ast[1]=Y_\ast[1]\oplus X_\ast,\nn
\eea
where $s : Y\ra Y[1]$, $u : X[-1] \ra Y[1]$, $t : Y \ra X$, and $v : X[-1] \ra X$. Thus $1_{C(f)} = h\del + \del h$ becomes
\bea
&&\left[
  \begin{array}{cc}
    1_Y & 0 \\
    0 & 1_X \\
  \end{array}
\right]=
\left[
  \begin{array}{cc}
    s & u \\
    t & v \\
  \end{array}
\right]\left[
  \begin{array}{cc}
    \del_Y & f \\
    0 & -\del_X \\
  \end{array}
\right]+\left[
  \begin{array}{cc}
    \del_Y & f \\
    0 & -\del_X \\
  \end{array}
\right]\left[
  \begin{array}{cc}
    s & u \\
    t & v \\
  \end{array}
\right]\nn\\
&&~~~~=\left[
  \begin{array}{cc}
    s\del+\del s+ft & -u\del+\del u+sf+fv \\
    t\del-\del t & -v\del-\del v+tf \\
  \end{array}
\right].\nn
\eea
By equating matrix entries, we have $0 = t\del + \del t$, and so $t$ is a chain morphism. Similarly, $1_Y = s\del+\del s+ft$ (i.e., $ft\simeq 1_Y$) and $1_X=-v\del-\del v+tf$ (i.e., $tf\simeq 1_X$). Hence $f$ is a homotopy equivalence.
\end{proof}

\begin{crl}[\blue{of the proof above}]
Let $X_\ast,Y_\ast$ be complexes of $R$-modules and $f_\ast:X_\ast\ra Y_\ast$ a quasi-isomorphism (i.e., $C(f)$ is exact). (i) If $C(f)$ is contractible then $f$ is a homotopy equivalence. (ii) If $f$ is an isomorphism, then $C(f)$ is contractible.
\end{crl}
\begin{proof}
Consider the situation in the proof of the preceding proposition. Then $f$ has a homotopy inverse $g:Y\ra X$, with $fg=1_Y-s\del -\del s$ and $gf=1_X+v\del+\del v$, $\iff$ the equation $-u\del+\del u+sf+fv=0$ (with $u$ arbitrary) can be solved. In particular, if $s=0$ and $v=0$, then we can set $u=0$.
\end{proof}

\section{Triangulated Categories}
Let $\A$ be an abelian category and $S_q:=$ quasiisomorphisms in $\A_0^\Integer$. Briefly, a triangulated category $\T$ is an additive, but possibly non-abelian, category such as the homotopy category ~$\H\A:={\A_0^\Integer\over\simeq}$~ or the derived category ~$\D\A:=\A_0^\Integer[S_q^{-1}]$~ (additivity to be proved later) in which the role originally played in $\A_0^\Integer$ by SES's, as basis for homological computations, is now played by DT's (distinguished triangles). The LES of homology in $\A_0^\Integer$ (a useful computational tool) exists in $\T$, and counterparts of \ul{derived functors} called \ul{homological functors} (Definition \ref{HomcalFun}) also exist in $\T$.

\subsection{Triangulation axioms}
\begin{dfn}[\textcolor{blue}{
\index{Translation! (Shift) functor}{Translation (Shift) functors} of an additive category,
\index{Triangle}{Triangle},
\index{Morphism of! triangles}{Morphism of triangles},
\index{Isomorphism of! triangles}{Isomorphism of triangles}}]
Let $\T$ be an additive category. Fix an additive autofunctor (called a \ul{shift functor} or \ul{translation functor} of $\T$):
\bea
T:\T\ra\T,~~\Big(X\sr{f}{\ral}Y\Big)\mapsto \Big(X[1]\sr{f[1]}{\ral}Y[1]\Big):=\Big(T(X)\sr{T(f)}{\ral}T(Y)\Big),\nn
\eea
and for $n\in\Integer$, consider the associated \ul{$n$-shift functor} or \ul{$n$-translation functor}
\bea
T^n:\T\ra\T,~\Big(X\sr{f}{\ral}Y\Big)\mapsto \Big(X[n]\sr{f[n]}{\ral}Y[n]\Big):=\Big(T^n(X)\sr{T^n(f)}{\ral}T^n(Y)\Big).\nn
\eea
For the case where $\T=\A_0^\Integer$, if a translation functor is not specified, we will assume it is given by
{\tiny\begin{align}
T^n:\A_0^\Integer\ra\A_0^\Integer,~\Big((X_\ast,d_\ast^X)\sr{f_\ast}{\ral}(Y_\ast,d_\ast^Y)\Big)\mapsto \Big((X_\ast[n],d_\ast^{X[n]})\sr{f_\ast[n]}{\ral}(Y_\ast[n],d_\ast^{Y[n]})\Big):=\Big((X_{\ast+n},(-1)^nd_{\ast+n}^X)\sr{(-1)^nf_{\ast+n}}{\ral}(Y_{\ast+n},(-1)^nd_{\ast+n}^Y)\Big).\nn
\end{align}}

A \ul{triangle} in $\T$ is a triple of morphisms $(u,v,w)$ given by a diagram of the form
\[
(u,v,w):~X\sr{u}{\ral}Y\sr{v}{\ral}Z\sr{w}{\ral}X[-1]~~~~~~~~\txt{or}~~~~
\adjustbox{scale=0.7}{\bt
  & Z\ar[ddl,"w"',"{[-1]}"] & \\
  &  & \\
 X\ar[rr,"u"] &  & Y\ar[uul,"v"']
\et}\]
A \ul{morphism of triangles} $(u,v,w)\sr{(f,g,h)}{\ral}(u',v',w')$ in $\T$ is a commutative diagram of the form
\bc\adjustbox{scale=0.9}{\bt
 X\ar[d,"f"]\ar[r,"u"] & Y\ar[d,"g"]\ar[r,"v"] & Z\ar[d,"h"]\ar[r,"w"] & X[-1]\ar[d,"{f[-1]}"]  \\
 X'\ar[r,"u'"] & Y'\ar[r,"v'"] & Z'\ar[r,"w'"] & X'[-1]
\et}\ec
A morphism of triangles $(u,v,w)\sr{(f,g,h)}{\ral}(u',v',w')$ is an \ul{isomorphism of triangles} if $f,g,h$ are isomorphisms.
\end{dfn}

\begin{dfn}[\textcolor{blue}{\index{Triangulated category}{Triangulated category}, \index{Distinguished triangles}{Distinguished triangles}, \index{Triangulation axioms}{Triangulation axioms}}]
A \ul{triangulated category} $\T=(\T,T,\Delta)$ is an additive category $\T$ together with a fixed translation functor $T:\T\ra\T$
and a fixed collection of triangles $\Delta$ (named \ul{distinguished triangles}, or DT's) in $\T$ satisfying the following defining axioms (called \ul{triangulation axioms}):
\bit[leftmargin=1.5cm]
\item[(TR1)] \index{Triangulation axiom! TR1: Construction axiom}{\ul{Construction axiom}}: DT's arise from morphisms and isomorphic images of other DT's:
   \bit
   \item[(a)] $X\sr{1_X}{\ral}X\ral 0\ral X[-1]$~ is a DT (distinguished triangle).
   \item[(b)] Any triangle isomorphic to a DT is also a DT.
   \item[(c)] Any morphism $X\sr{u}{\ral}Y$ can be completed to a DT, say,
       \bea
       \big(u,\al(u),\beta(u)\big):~X\sr{u}{\ral}Y\sr{\al(u)}{\ral}E(u)\sr{\beta(u)}{\ral}X[-1].\nn
       \eea
   \eit
\item[(TR2)] \index{Triangulation axiom! TR2: Shift symmetry axiom}{\ul{Shift symmetry axiom}}: If ~{\small $(u,v,w):~X\sr{u}{\ral}Y\sr{v}{\ral}Z\sr{w}{\ral}X[-1]$}~ is a DT, then so are the following:
 \bit
 \item[(a)] Its \ul{left shift} ~(removing $X$ and using using $Y[-1]$ to repair the right end)
    \bea
    (v,w,-u[-1]):~Y\sr{v}{\ral}Z\sr{w}{\ral}X[-1]\sr{-u[-1]}{\ral}Y[-1]\nn
     \eea
 \item[(b)] Its \ul{right shift} ~(removing $X[-1]$ and using $Z[1]$ to repair the left end)
     \bea
     (-w[1],u,v):~Z[1]\sr{-w[1]}{\ral}X\sr{u}{\ral}Y\sr{v}{\ral}Z\nn
     \eea
 \eit

\item[(TR3)] \index{Triangulation axiom! TR3: Stability axiom}{\ul{Stability axiom}}: Given DT's $(u,v,w)$, $(u',v',w')$ and morphisms $f,g$ such that the \ul{left square} in the following diagram \ul{commutes}
\bc\adjustbox{scale=0.9}{\bt
 X\ar[d,"f"]\ar[r,"u"] & Y\ar[d,"g"]\ar[r,"v"] & Z\ar[d,dashed,"h"]\ar[r,"w"] & X[-1]\ar[d,"{f[-1]}"]  \\
 X'\ar[r,"u'"] & Y'\ar[r,"v'"] & Z'\ar[r,"w'"] & X'[-1]
\et}\ec
there exists a morphism $h$ completing the diagram to a \ul{morphism of triangles}. (\blue{footnote}\footnote{(i) Similarly, by the shift-symmetry axioms, if $f,h$ (resp. $g,h$) above form the obvious partial commutative diagram without $g$ (resp. without $f$) then a $g$ (resp. an $f$) will exist completing the diagram. (ii) Also, it follows from the axioms so far that in any DT, $X\sr{u}{\ral}Y\sr{v}{\ral}Z\sr{w}{\ral}X[-1]$, we have $vu=wv=u[-1]w=0$.})

\item[(TR4)] \index{Triangulation axiom! TR4: Octahedral axiom}{\ul{Octahedral axiom}}: Given morphisms\adjustbox{scale=0.7}{%
\bt X\ar[rr,bend left,"vu"]\ar[r,"u"] & Y\ar[r,"v"] & Z\et
}and completing DT's for $u,v,vu$ (by TR1 c)
{\small\begin{align}
&\big(u,\al(u),\beta(u)\big):~X\sr{u}{\ral}Y\sr{\al(u)}{\ral}E(u)\sr{\beta(u)}{\ral}X[-1]\nn\\
&\big(v,\al(v),\beta(v)\big):~Y\sr{v}{\ral}Z\sr{\al(v)}{\ral}E(v)\sr{\beta(v)}{\ral}Y[-1]\nn\\
&\big(vu,\al(vu),\beta(vu)\big):~X\sr{vu}{\ral}Z\sr{\al(vu)}{\ral}E(vu)\sr{\beta(vu)}{\ral}X[-1]\nn
\end{align}}the following associated sequence (which exists by TR3: \blue{footnote}\footnote{By TR3, the morphisms $f,g$ always exist. What is new, and thus an axiom, is the requirement that $f,g$ can be chosen so that the indicated resulting sequence is a DT.}) can be chosen to be a DT.
{\small\[
\bt[column sep=small]\Big(f,g,\al(u)[-1]\beta(v)\Big):~E(u)\ar[r,dashed,"f"]& E(vu)\ar[r,dashed,"g"]& E(v)\ar[rrrr,"{\al(u)[-1]\beta(v)}"]&&&& E(u)[-1]\et
\]}
In commutative-diagram form, this can be expressed in several ways, including the following:

\bea
\adjustbox{scale=0.7}{%
\bt
X\ar[d,equal]\ar[r,"u"] & Y\ar[d,"v"]\ar[r,"\al(u)"] & E(u)\ar[d,dashed,"f"]\ar[r,"\beta(u)"] & X[-1]\ar[d,equal] \\
X\ar[d,"u"]\ar[r,"vu"] & Z\ar[d,equal] \ar[r,"\al(vu)"] & E(vu)\ar[d,dashed,"g"]\ar[r,"\beta(vu)"] & X[-1]\ar[d,"{u[-1]}"] \\
Y\ar[dd,"0"]\ar[r,"v"] & Z\ar[dd,"0"] \ar[r,"\al(v)"] & E(v)\ar[dd,"{\al(u)[-1]\beta(v)}"]\ar[r,"\beta(v)"] & Y[-1]\ar[dd,"0"] \\
 & & & \\
X[-1]\ar[r,"{u[-1]}"] & Y[-1] \ar[r,"{\al(u)[-1]}"] & E(u)[-1]\ar[r,"{\beta(u)[-1]}"] & X[-2] \\
\et}~~~~
\adjustbox{scale=0.7}{%
\bt
E(u)\ar[ddddrr,dashed, bend right=40,"f"']\ar[ddr,"\beta(u)","{[-1]}"']  &  & Y\ar[ll,"\al(u)"']\ar[ddr,"v"] &  & E(v)\ar[ll,"\beta(v)"',"{[-1]}"]\ar[llll,bend right=40,"{\al(u)[-1]\beta(v)}"',"{[-1]}"] \\
  &  & \curvearrowright &  &  \\
  & X\ar[uur,"u"]\ar[rr,"vu"] &  & Z\ar[uur,"\al(v)"]\ar[ddl,"\al(vu)"'] &  \\
  &  &  &  &  \\
  &  & E(vu)\ar[uul,"\beta(vu)"',"{[-1]}"]\ar[uuuurr,dashed, bend right=40,"g"'] &  &
\et}\nn
\eea

\bea
\adjustbox{scale=0.7}{%
\bt
X\ar[d,equal]\ar[r,"u"] & Y\ar[d,"v"]\ar[r,"\al(u)"] & E(u)\ar[d,dashed,"f"]\ar[r,"\beta(u)"] & X[-1]\ar[d,equal] \\
X\ar[d,"u"]\ar[r,"vu"] & Z\ar[d,equal] \ar[r,"\al(vu)"] & E(vu)\ar[d,dashed,"g"]\ar[r,"\beta(vu)"] & X[-1]\ar[d,"{u[-1]}"] \\
Y\ar[d,"v"]\ar[r,"v"] & Z\ar[d,"\al(v)"] \ar[r,"\al(v)"] & E(v)\ar[d,"\beta(v)"]\ar[r,"\beta(v)"] & Y[-1]\ar[d,"{v[-1]}"] \\
Z\ar[dd,"0"]\ar[r,"\al(v)"] & E(v)\ar[dd,"\beta(v)"] \ar[r,"\beta(v)"] & Y[-1]\ar[dd,"{\al(u)[-1]}"]\ar[r,"{v[-1]}"] & Z[-1]\ar[dd,"0"] \\
 & & & \\
X[-1]\ar[r,"{u[-1]}"] & Y[-1] \ar[r,"{\al(u)[-1]}"] & E(u)[-1]\ar[r,"{\beta(u)[-1]}"] & X[-2] \\
\et}~~~~
\adjustbox{scale=0.7}{%
\bt
   &  & E(vu)\ar[dddddll,dotted,Rightarrow,blue,"\beta(vu)"']\ar[dddddrr,dashed,Rightarrow,blue,"g"] &   &  \\
   &  &  &   &  \\
   &  &  &   &  \\
   &  &  &   &  \\
   & E(u)\ar[dl,dotted,Rightarrow,blue,"\beta(u)"]\ar[uuuur,dashed,"f"'] &  &   &  \\
 X\ar[dddddrr,"u"']\ar[drrr,"vu"]  &  &  &   & E(v)\ar[ulll,"{\al(u)[-1]\beta(v)}"]\ar[dddddll,dotted,Rightarrow,blue,"\beta(v)"] \\
   &  &  & Z\ar[uuuuuul,dotted,"\al(vu)"description]\ar[ur,dotted,"\al(v)"]  &  \\
   &  &  &   &  \\
   &  &  &   &  \\
   &  &  &   &  \\
   &  & Y\ar[uuuuuul,dotted,"\al(u)"description]\ar[uuuur,"v"] &   &  \\
\et}\nn
\eea
\eit
\end{dfn}

\subsection{Properties of triangulated categories}
All results below (i.e., Lemma \ref{HomPropLmm}, Corollary \ref{HomPropCrl}, Theorem \ref{FiveLmmTC}, and Corollary \ref{uniqueDT-prop}) do not require the fourth axiom TR4.

\begin{lmm}[\blue{A distinguished triangle is a complex}]\label{HomPropCrl}
Let $\T$ be a triangulated category and ~$X\sr{u}{\ral}Y\sr{v}{\ral}Z\sr{w}{\ral}X[-1]$~ a DT in $\T$. Then ~$vu=wv=u[-1]w=0$.
\end{lmm}
\begin{proof}
By applying the stability axiom TR3 to $X\sr{id}{\ra} X\ra 0\ra X[-1]$ and the given distinguished triangle $X\ra Y\ra Z\ra X[-1]$, we get the following commutative diagram (which shows $vu=v\al=0$):
\bc\adjustbox{scale=0.9}{\bt
 X\ar[d,"id"]\ar[r,"1"] & X\ar[d,dashed,"\al"]\ar[r,"0"] & 0\ar[d,"0"]\ar[r,"0"] & X[-1]\ar[d,"{id[-1]}"]  \\
 X\ar[r,"u"] & Y\ar[r,"v"] & Z\ar[r,"w"] & X[-1]
\et}\ec
Therefore, by applying the shift-symmetry axioms TR2, we get $wu=0$ and $u[-1]w=0$ as well.
\end{proof}

\begin{dfn}[\blue{
\index{Periodic complex of a DT}{Periodic complex of a DT},
\index{Morphism of! triangulated categories}{Morphism of triangulated categories}}]
Let $\T$ be a triangulated category and {\footnotesize $\Delta:~X\sr{u}{\ral}Y\sr{v}{\ral}Z\sr{w}{\ral}X[-1]$}. The \ul{periodic complex} of $\Delta$ is the following sequence of DT's obtained by repeatedly applying TR2.
\bea
P(\Delta):~\cdots\ral Z[1]\sr{w[1]}{\ral}X\sr{u}{\ral}Y\sr{v}{\ral}Z\sr{w}{\ral}X[-1]\sr{u[-1]}{\ral}Y[-1]\sr{v[-1]}{\ral}Z[-1]\sr{w[-1]}{\ral}X[-2]\sr{u[-2]}{\ral}\cdots\nn
\eea
The image of $\Delta$ (resp., of the periodic complex $P(\Delta)$) under an additive functor might not be a DT (resp., might not be as periodic as $P(\Delta)$). An additive functor (between triangulated categories) that (i) commutes with translation and (ii) maps DT's to DT's is known as a \ul{morphism of triangulated categories}. Also, there can be additive functors under which P(DT)'s become exact. These are called \ul{homological functors} (see the next definition).
\end{dfn}

\begin{dfn}[\textcolor{blue}{\index{Homological functor}{Homological functor}}]\label{HomcalFun}
Let $\T$ be a triangulated category (with translation $T:\T\ra\T$) and $\A$ an abelian category. A functor functor $F:\T\ra\A$ is a \ul{homological functor} if (i) it is additive and (ii) for any DT, {\small $X\sr{u}{\ral}Y\sr{v}{\ral}Z\sr{w}{\ral}X[-1]$}, the following sequence is exact:
\bea
\cdots\ral F(Z[1])\sr{F(w[1])}{\ral}F(X)\sr{F(u)}{\ral}F(Y)\sr{F(v)}{\ral}F(Z)\sr{F(w)}{\ral}F(X[-1])\sr{F(u[-1])}{\ral}F(Y[-1])\sr{F(v[-1])}{\ral}\cdots\nn
\eea
(A \ul{homological cofunctor}, i.e., \ul{contravariant homological functor}, is similarly defined except with arrows reversed in the above sequence.)
\end{dfn}
\begin{rmk}[\blue{Homology functors are homological functors}]
Let $\T\subset\H\A:={\A_0^\Integer\over\simeq}$ be a triangulated subcategory of the homotopy category $\H\A$ (\blue{footnote}\footnote{We will show later that $\H(R\txt{-mod})$ (and hence $\H\A$ by imbedding) is triangulated.}). Then by the induced LES of homologies for a distinguished triangle in $\H\A$ (from Remark \ref{ModLesForDT}), the \ul{zeroth homology} functor $H_0=H^0:\T\ra\A$ is a homological functor. Consequently, the \ul{$i$th homology functor} $H_i:=H_0\circ T^{-i}:\T\ra\A$ (or $H^i:=H^0\circ T^i:\T\ra\A$) is also a homological functor.
\end{rmk}

\begin{lmm}[\textcolor{blue}{Mor functors are homological functors}]\label{HomPropLmm}
Let $\T$ be a triangulated category, $U\in\Ob\T$. Then $\Mor_\T(U,-):\T\ra Ab$ is a homological functor (resp. $\Mor_\T(-,U):\T\ra Ab$ is a homological cofunctor).
\end{lmm}
\begin{proof}
We will explain $\Mor_\T(U,-)$ only, because the proof for $\Mor_\T(-,U)$ is similar. It is enough to check that given any distinguished triangle $X\sr{u}{\ral}Y\sr{v}{\ral}Z\sr{w}{\ral}X[-1]$, the sequence
\bea
\Mor_\T(U,X)\sr{u_\ast}{\ral}\Mor_\T(U,Y)\sr{v_\ast}{\ral}\Mor_\T(U,Z)\nn
\eea
is exact, because exactness of the rest of the long sequence follows from shifting (TR2).
(\blue{footnote}\footnote{Recall that ~$u_\ast:=u\circ(-)~=~(-)\sr{(-)}{\ral}\dom u\sr{u}{\ral}\cod u$.})

\bit[leftmargin=0.7cm]
\item \ul{$\im u_\ast\subset\ker v_\ast$}:~ Recall that $vu=0$. Since $\Mor_\T(U,-)$ is an additive functor, we get $v_\ast u_\ast=(vu)_\ast=0$. So, $\im u_\ast\subset\ker v_\ast$. (Recall that $\Mor_\T(U,W)$ is an abelian group, i.e., $\Integer$-module for every $W\in\T$).
\item \ul{$\im u_\ast\supset\ker v_\ast$}:~ Let $f\in\ker v_\ast$. Then by the definition of $v_\ast$, $0=v_\ast(f)=vf$. By TR1, $U\sr{id}{\ral} U\sr{0}{\ral} 0\sr{0}{\ral} U[-1]$ is a distinguished triangle. So, we have the commutative diagram
\bc\adjustbox{scale=0.9}{\bt
 U\ar[d,dashed,"g"]\ar[r,"1"] & U\ar[d,"f"]\ar[r,"0"] & 0\ar[d,"0"]\ar[r,"0"] & U[-1]\ar[d,dashed,"{g[-1]}"]  \\
 X\ar[r,"u"] & Y\ar[r,"v"] & Z\ar[r,"w"] & X[-1]
\et}\ec
Using TR2 (to rotate it), TR3 (to complete it with a morphism $U\sr{g}{\ral}X$), and then TR2 again (to rotate back), we get a morphism $U\sr{g}{\ral}X$ with the first square commuting:~ $f=f1=ug=u_\ast(g)$, ~$\Ra$ ~$f\in im(u_\ast)$.
\eit
Hence ~$\im u_\ast=\ker v_\ast$, ~which proves $\Mor_\T(U,-)$ is homological.

Similarly, to show $\Mor_\T(-,U)$ is homological, it suffices to show the following sequence is exact.
\bea
\bt \Mor_\T(X,U)\ar[from=r,"u^\ast"'] & \Mor_\T(Y,U)\ar[from=r,"v^\ast"'] & \Mor_\T(Z,U).\et\nn
\eea
As before, $vu=0$ implies $u^\ast v^\ast=(vu)^\ast=0$, and so $\im v^\ast\subset\ker u^\ast$. On the other hand, if $f\in\ker u^\ast$, i.e.,
\bea
0=u^\ast(f)~=~fu~=~X\sr{u}{\ral}Y\sr{f}{\ral}U,\nn
\eea
then by TR1 (now with TR2 also), $0\sr{0}{\ral}U\sr{id}{\ral} U\sr{0}{\ral}0[-1]$ is a DT, giving a commutative diagram
\[\adjustbox{scale=0.9}{\bt
 X\ar[r,"u"] & Y\ar[r,"v"] & Z\ar[r,"w"] & X[-1]\\
  0\ar[from=u,"0"]\ar[r,"0"] & U\ar[from=u,"f"]\ar[r,"1"] & U\ar[from=u,dashed,"g"]\ar[r,"0"] & 0[-1]\ar[from=u,"{0[-1]}"]
\et}~~~~\Ra~~~~f=1f=gv=v^\ast(g)\in\im v^\ast.
\]
\end{proof}

\begin{thm}[\textcolor{blue}{\index{Five Lemma for triangulated categories}{Five Lemma for triangulated categories (or for distinguished triangles)}}]\label{FiveLmmTC}
Let $\T$ be a triangulated category and consider the following diagram from the statement of TR3 for $\T$.
\[\adjustbox{scale=0.9}{\bt
 X\ar[d,"f"]\ar[r,"u"] & Y\ar[d,"g"]\ar[r,"v"] & Z\ar[d,dashed,"h"]\ar[r,"w"] & X[-1]\ar[d,"{f[-1]}"]  \\
 X'\ar[r,"u'"] & Y'\ar[r,"v'"] & Z'\ar[r,"w'"] & X'[-1]
\et}~~~~~~~~
\txt{If $f,g$ are isomorphisms, then so is $h$.}\]

(By shifting, a more general result is immediate: If any two of $f,g,h$ are isos, then so is the other.)
\end{thm}
\begin{proof}
Applying the \ul{homological} functor $\Mor_\T(Z',-)$ to the given diagram in TR3 (with $f,g$ isos), we get a commutative diagram
\[\adjustbox{scale=0.8}{\bt
 \Mor_\T(Z',X)\ar[d,"\cong"',"f_\ast"]\ar[r,"u_\ast"] & \Mor_\T(Z',Y)\ar[d,"\cong"',"g_\ast"]\ar[r,"v_\ast"] & \Mor_\T(Z',Z)\ar[d,dashed,"h_\ast"]\ar[r,"w_\ast"] & \Mor_\T(Z',X[-1])\ar[d,"\cong"',"{f[-1]_\ast}"]\ar[r,"{u[-1]_\ast}"] & \Mor_\T(Z',Y[-1])\ar[d,"\cong"',"{g[-1]_\ast}"]  \\
 \Mor_\T(Z',X')\ar[r,"u_\ast'"] & \Mor_\T(Z',Y')\ar[r,"v_\ast'"] & \Mor_\T(Z',Z')\ar[r,"w_\ast'"] & \Mor_\T(Z',X'[-1])\ar[r,"{u'[-1]_\ast}"] & \Mor_\T(Z',Y'[-1])
\et}\]
whose rows are exact by Lemma \ref{HomPropLmm}. Since $f,g$ are isos and the translation $T:=[1]:\T\ra\T$ is a functor, $f[-1]$, $g[-1]$ are also isos. Thus, $h_\ast$ is an iso by the five-lemma (for $\Integer$-modules), since $f_\ast,g_\ast,f[-1]_\ast,g[-1]_\ast$ are isos in $\Integer$-mod. By the surjectivity of $h_\ast$, there exists $Z'\sr{k}{\ral}Z$ such that $hk=h_\ast(k)=1_{Z'}:Z'\sr{k}{\ral}Z\sr{h}{\ral}Z'$ in $\Mor_\T(Z',Z')$ (i.e., $k$ is a \ul{right inverse} of $h$).

Similarly, by applying the \ul{homological} functor $\Mor_\T(-,Z)$ to the given diagram in TR3 (with $f,g$ isos), we get a morphism $Z'\sr{l}{\ral}Z$ such that {\small $lh=1_Z:Z\sr{h}{\ral}Z'\sr{l}{\ral}Z$} (a \ul{left inverse} of $h$).
\[\adjustbox{scale=0.8}{\bt
Hom(X,Z) & Hom(Y,Z)\ar[l,"u^\ast"'] & Hom(Z,Z)\ar[l,"v^\ast"'] & Hom(X[-1],Z)\ar[l,"w^\ast"'] & Hom(Y[-1],Z)\ar[l,"{u[-1]^\ast}"'] \\
Hom(X',Z)\ar[u,"\cong","f^\ast"'] & Hom(Y',Z)\ar[l,"u'{}^\ast"']\ar[u,"\cong","g^\ast"'] & Hom(Z',Z)\ar[l,"v'{}^\ast"']\ar[u,dashed,"h^\ast"'] & Hom(X'[-1],Z)\ar[l,"w'{}^\ast"']\ar[u,"\cong","{f[-1]{}^\ast}"'] & Hom(Y'[-1],Z)\ar[l,"{u'[-1]^\ast}"']\ar[u,"\cong","{g[-1]{}^\ast}"']
\et}\]
Hence $h$ is an iso.
\end{proof}

\begin{crl}\label{uniqueDT-prop}
Let $\T$ be a triangulated category. A distinguished triangle in $\T$,
\[
(u,v,w):~X\sr{u}{\ral}Y\sr{v}{\ral}Z\sr{w}{\ral}X[-1],
\]
is uniquely determined (up to isomorphism) by any one of its transition morphisms $u,v,w$.

(So, in TR1(c), the extension of a morphism to a distinguished triangle is unique up to isomorphism.)
\end{crl}
\begin{proof}
Observe that the existing morphisms $h,h'$ in the following diagram are mutually inverse isomorphisms (by Theorem \ref{FiveLmmTC}, \blue{footnote}\footnote{Note that this result can also be obtained directly from the triangulation axioms, without using the five lemma (Theorem \ref{FiveLmmTC}).}).
\bc\adjustbox{scale=0.9}{\bt
 X\ar[d,equal]\ar[r,"u"] & Y\ar[d,equal]\ar[r,"v"] & Z\ar[d,dashed,shift right,"h"']\ar[r,"w"] & X[-1]\ar[d,equal]  \\
 X\ar[r,"u"] & Y\ar[r,"\al(u)"] & E(u)\ar[u,dashed,shift right,"h'"']\ar[r,"\beta(u)"] & X[-1]
\et}\ec
Thus, the triangle is determined by the first morphism $u$ (up to isomorphism).
By shift-symmetry TR2, the triangle is similarly determined (up to isomorphism) by $v$, and by $w$.
\end{proof}

%% file: parts/AlgebraCat/AlgebraCatS13.tex
\chapter{Derived Categories and their Triangulation}\label{AlgebraCatS13}
\section{Derived Category. Localization of Homotopy Category. Localizing Class Construction}
For this chapter it is essential for the reader to review the earlier chapter on the localization of a category. Throughout, $\A$ will denote an abelian category, but by the Freyd-Mitchell imbedding, we will also occasionally assume we are in $R$-mod (for some ring $R$) and therefore use element-wise arguments in our proofs. Therefore, wherever it is seems necessary, the reader can simply assume $\A=R$-mod.
\begin{dfn}[\textcolor{blue}{Recall: \index{Derived! category}{Derived category}, \index{Homotopy! category}{Homotopy category}}]
Let $\A$ be an abelian category, and $S_q:=\{\txt{quasi-iso's in $\A_0^\Integer$}\}$. The \ul{derived category} $\D\A$ of $\A$ is the localization $L_q:\A_0^\Integer\ra\D\A$ of $\A_0^\Integer$ by quasiisomorpisms (i.e., $\D\A:=\A_0^\Integer\left[S_q^{-1}\right]$ is the smallest category containing $\A_0^\Integer$ such that quasiisomorphisms are isomorphisms). The \ul{homotopy category} $\H\A$ of $\A$ is the quotient category ${\A_0^\Integer\over\simeq}$, where $\simeq$ denotes homotopy of chain morphisms. (\blue{footnote}\footnote{Note that some quais-iso's in $\A_0^\Integer$ (such as homotopy equivalences) become iso's in the homotopy category $\H\A$, just as in the derived category $\D\A$.})
\end{dfn}
Recall (from the definition of a quotient category) that $\H\A$ is the category with the \ul{same objects as $\A_0^\Integer$} and \ul{morphisms of $\H\A$} being equivalence classes of homotopic chain morphisms, i.e.,
{\small\begin{align}
\textstyle \Mor_{\H\A}(C,C'):={\Mor_{\A_0^\Integer}(C,C')\over\simeq}:=\left\{[f]_\simeq~\big|~f\in \Mor_{\A_0^\Integer}(C,C')\right\},~~~~\txt{for any chain complexes $C,C'\in\A_0^\Integer$}.\nn
\end{align}}

\begin{lmm}[\textcolor{blue}{\cite[Lemma 3, p.159]{gelfand-manin2010}}]
Let $\A$ be an abelian category and $L=L_q:\A_0^\Integer\ra\D\A$ its derived category. Then homotopic morphisms in $\A_0^\Integer$ become equal in $\D\A$.
\end{lmm}
\begin{proof}
Let $\al\in\Mor\A_0^\Integer$ be a homotopy equivalence (hence a quasi-iso), i.e., there exists $\beta:\cod \al\ra\dom \al$ such that $\al\beta\simeq id_{\cod \al}$, $\beta\al\simeq id_{\dom \al}$ ($\iff$ $[\al][\beta]=id_{\cod \al}$, $[\beta][\al]=id_{\dom \al}$ in $\H\A$). Then it is clear that in the derived category $\D\A$, the latter relations become $L(\al)L(\beta)=id_{\cod \al}$, $L(\beta)L(\al)=id_{\dom \al}$ (i.e., $\al$ is an isomorphism in $\D\A$).

Let $f,g\in\Mor\A_0^\Integer$ be homotopic (where we assume $\A=R$-mod by imbedding), i.e., $f=g+d^Bh+hd^A:A\ra B$, for a homotopy $h_\ast:A_\ast\ra B_{\ast+1}$. Then from before, we have the diagram (a) below in $\A_0^\Integer$ which commutes except for the left square which only commutes up to homotopy (where the new map $c(h)$ is explicitly given by $c(h)(b_n,a_{n-1}):=(b_n+h_{n-1}(a_{n-1}),a_{n-1})$ and, as can be verified directly, it is a chain morphism). By applying the induced LES of homologies on diagram (a), we see that $c(h)$ is a quasi-isomorphism.
\[(a)~~\adjustbox{scale=0.9}{\bt
A\ar[d,equal]\ar[r,"f"] & B\ar[d,equal]\ar[r,"i_1^f"]& C(f)\ar[d,"c(h)"]\ar[r,"{p}_2^f"] & A[-1]\ar[d,equal] \\
A\ar[r,"g"] & B\ar[r,"i_1^g"]& C(g)\ar[r,"{p}_2^g"] & A[-1] \\
\et}~~~~~~~~(b)~~
\adjustbox{scale=0.9}{\bt
  & B\ar[d,shift right,dotted,"i_f"']\ar[from=d,shift right,"j_f"'] &  &  \\
A\ar[d,equal]\ar[ur,bend left,"-f"]\ar[r,"i_2^f"] & Cyl(f)\ar[d,"cyl(h)"]\ar[r,"{p}_{1,3}^f"]& C(f)\ar[d,"c(h)"]\ar[r,"{p}_2^f"] & A[-1]\ar[d,equal] \\
A\ar[dr,bend right,"-g"']\ar[r,"i_2^g"] & Cyl(g)\ar[r,"{p}_{1,3}^g"]& C(g)\ar[r,"{p}_2^g"] & A[-1] \\
  & B\ar[from=u,shift right,"j_g"']\ar[u,shift right,dotted,"i_g"']  &  &
\et}\]
Similarly, we have the diagram (b) above which commutes in $\A_0^\Integer$ wrt all solid (i.e., non-dotted) morphisms and commutes up to homotopy wrt the dotted morphisms (where the map $cyl(h)$ is explicitly given by $cyl(h)(b_n,a_n,a_{n-1}):=\big(b_n+h_{n-1}(a_{n-1}),a_n,a_{n-1}\big)$ and, as can be checked using our equivalent convention for the mapping cylinder, it is a chain morphism). By again applying the induced LES of homologies, this time on (b), we see that $cyl(h)$ is also a quasi-isomorphism. Since (i) all vertical morphisms in diagram (b) are quasi-isomorphisms and (ii) $j_g\circ cyl(f)\circ i_f=id_B$, by applying the localizing functor $L$ on diagram (b), we see that $L(f)=L(g)$, i.e., $f=g$ in $\D\A$, as follows:
{\small\begin{align}
-L(g)&=L(-g)=L(j_g\circ i_2^g)=L\big(j_g\circ cyl(h)\circ i_2^f\big)=L\big(j_g\circ cyl(h)\big)\circ L(i_2^f)
=L\big(j_g\circ cyl(h)\big)\circ L(i_f\circ j_f\circ i_2^f)\nn\\
&~~~~=L\big(j_g\circ cyl(h)\circ i_f\big)\circ L(j_f\circ i_2^f)=L\big(j_g\circ cyl(h)\circ i_f\big)\circ L(-f)=L(-f)=-L(f).\nn \qedhere
\end{align}}
\end{proof}

\begin{prp}[\blue{Derived category via homotopy category: \cite[Proposition 2, p.159]{gelfand-manin2010}}]
Let $\A$ be an abelian category, $L=L_q:\A_0^\Integer\ra\D\A$ its derived category, $\H\A:={\A_0^\Integer\over\simeq}$ the homotopy category, and $\wt{L}_q:\H\A\ra\H\A[S_q^{-1}]$ the localization of $\H\A$ by quasi-isomorphisms. Then we can also obtain the derived category through the homotopy category as follows:
\bea
\bt L_q=\wt{L}_q\circ(/\simeq):\A_0^\Integer\ar[r,"/\simeq"] & \H\A:={\A_0^\Integer\over\simeq}\ar[r,"\wt{L}_q"] & \H\A[S_q^{-1}]\cong\D\A.\et\nn
\eea
\end{prp}
\begin{proof}
A quasi-isomorphism $q'$ in $\H\A$ is naturally defined to be the homotopy class $q':=[q]_\simeq$ of a quasi-isomorphism $q$ in $\A_0^\Integer$. Therefore, the composition {\footnotesize $\wt{L}_q\circ(/\simeq):\A_0^\Integer\sr{/\simeq}{\ral}\H\A\sr{\wt{L}_q}{\ral}\H\A(S_q^{-1})$} takes all quasi-iso's to isomorphisms, and so there exists a unique functor $F:\D\A\ra\H\A[S_q^{-1}]$ such that $F\circ L_q=\wt{L}_q\circ(/\simeq)$ as follows:
\[\adjustbox{scale=0.9}{\bt \A_0^\Integer\ar[d,"/\simeq"']\ar[rr,"L_q"] && \D\A=\A_0^\Integer[S_q^{-1}]\ar[ddll,dashed,"F"',"\cong"]\\
\H\A\ar[d,"\wt{L}_q"'] & & \\
\H\A[S_q^{-1}] & &
\et}\]
By construction, $F$ is bijective on objects. Since the functor ~$/\simeq$~ is surjective on morphisms, so is $F$, i.e.,
\[
\substack{f_0s_1^{-1}f_1\cdots s_n^{-1}f_n~\sim~g_0t_1^{-1}g_1\cdots t_n^{-1}g_n\\~\\ \txt{in}~\D\A=\A_0^\Integer[S_q^{-1}]}~~~~\Ra~~~~\substack{[f_0]_\simeq[s_1]_\simeq^{-1}[f_1]_\simeq\cdots[s_n]_\simeq^{-1}[f_n]_\simeq~\sim~ [g_0]_\simeq[t_1]_\simeq^{-1}[g_1]_\simeq\cdots[t_n]_\simeq^{-1}[g_n]_\simeq\\~\\~\txt{in}~~\H\A[S_q^{-1}].}
\]
As homotopic chain morphisms become equal in $\D\A$ (lemma), $F$ is also injective on morphisms, because
\[
\substack{[f_0]_\simeq[s_1]_\simeq^{-1}[f_1]_\simeq\cdots[s_n]_\simeq^{-1}[f_n]_\simeq~\sim~ [g_0]_\simeq[t_1]_\simeq^{-1}[g_1]_\simeq\cdots[t_n]_\simeq^{-1}[g_n]_\simeq\\~\\~\txt{in}~~\H\A[S_q^{-1}]}~~~~\Ra~~~~\substack{f_0s_1^{-1}f_1\cdots s_n^{-1}f_n~\sim~g_0t_1^{-1}g_1\cdots t_n^{-1}g_n\\~\\~\txt{in}~\D\A=\A_0^\Integer[S_q^{-1}].}
\]
\end{proof}

\begin{notation}[\textcolor{blue}{Recall}]
Let $\C$ be a category and $S\subset\Mor\C$ any class of morphisms. In the proof of Theorem \ref{LocExistThm}, we have shown that the localization $\C[S^{-1}]$ exists, with $\Ob\C[S^{-1}]=\Ob\C$, and for any $A,B\in\Ob\C[S^{-1}]$,
\bea
\Mor_{\C[S^{-1}]}(A,B)=\left\{\txt{finite compositions}~[f_0][s_1]^{-1}[f_1]\cdots[s_n]^{-1}[f_n]:A\ra B~|~f_i\in\Mor\C,~s_i\in S\right\}.\nn
\eea
For any $f\in\Mor\C$ and $s\in S$, we will denote the equivalence classes $[f]$, $[s]$ by $f$, $s$ respectively. Thus, an arbitrary morphism $[f_0][s_1]^{-1}[f_1]\cdots[s_n]^{-1}[f_n]\in\Mor\C[S^{-1}]$ will be written as $f_0s_1^{-1}f_1\cdots s_n^{-1}f_n$.
\end{notation}

In order to \ul{simplify} the above expression for morphisms in $\C[S^{-1}]$, we will impose on $S$ additional conditions, which cause no loss of generality in the sense that they will be shown (in Theorem \ref{QuisLocLmm}) to be sufficient for our intended purpose, namely, the construction of the derived category
\bea
\D\A:=\A_0^\Integer\left[S_q^{-1}\right]:=\A_0^\Integer\left[\{\txt{quasi-iso's}\}^{-1}\right].\nn
\eea
The morphism class $S\subset\Mor\C$ with certain additional constraints will be called a \ul{localizing class} of morphisms as follows:

\begin{dfn}[\textcolor{blue}{\index{Localizing class}{Localizing class} of morphisms}]\label{loc-class-dfn}
Let $\C$ be a category. A class of morphisms $S\subset\Mor\C$ is \ul{localizing} (or \ul{multiplicative}: \blue{footnote}\footnote{If adopted, this would be our second notion of a multiplicative class, which is specific to subclasses of morphisms of a category, and so different from the earlier one.}) if it satisfies the following \ul{three} conditions:
\bit
\item[(a)] \ul{$S$ is a ``composable'' class (``Composability'' conditions)}
\bit
\item[(1)] $id_C$ (or $1_C$) $\in S$, for every object $C\in\Ob\C$.
\item[(2)] $s,t\in S$ ~$\Ra$~ $st\in S$.
\eit
\item[(b)] \ul{$S$ is extendable (Extension conditions)}: (\blue{footnote}\footnote{Except lacking universal properties, these resemble the limits known as ``pullbacks'' and ``pushouts''.})~
\bit
\item[(1)] For every $f\in\Mor \C$, $s\in S$ with the same codomain (i.e., $\cod f=\cod s$), there exist $f'\in \Mor\C$, $s'\in S$ such that $fs'=sf'$, i.e., the diagram
{\footnotesize
\bt
C' \ar[d,dashed,"s'"]\ar[r,dashed,"f'"] & D'\ar[d,"s"] \\
C\ar[r,"f"] & D
\et} commutes.

\item[(2)] Equivalently: For every $f\in\Mor\C$, $s\in S$ with the same domain (i.e., $\dom f=\dom s$), there exist $f'\in\Mor\C$, $s'\in S$ such that $s'f=f's$, i.e., the diagram
{\footnotesize
\bt
C \ar[d,"s"]\ar[r,"f"] & D\ar[d,dashed,"s'"] \\
C'\ar[r,dashed,"f'"] & D'
\et} commutes.
\eit
\item[(c)] \ul{$S$ is equalizer-balanced (Balance of left and right equalizers)}:  For any two morphisms $f,g\in\Mor\C$,
\bea
sf=sg~~\txt{for some}~~s\in S~~~~\iff~~~~fs'=gs'~~\txt{for some}~~s'\in S.\nn
\eea
That is,~~~~
$\bt
 C\ar[rrr,bend left,"sf=sg"]\ar[rr,shift left=1,"f"]\ar[rr,shift right=1,"g"'] && D\ar[r,dashed,"s"] & D'
\et$ ~~~~$\iff$~~~~
$\bt
  C'\ar[r,dashed,"s'"]\ar[rrr,bend left,"fs'=gs'"] & C \ar[rr,shift left=1,"f"]\ar[rr,shift right=1,"g"'] && D
\et$
\eit
\end{dfn}

\begin{rmk}[\textcolor{blue}{Expression of morphisms in $\C[S^{-1}]$ when $S$ is localizing:
\index{Left! roof}{Left roof},
\index{Right! roof}{Right roof},
\index{Equivalent left roofs}{Equivalent left roofs}}]~
Let $\C$ be a category and $S\subset\Mor\C$ a localizing class. Let $f\in\Mor\C$, $s\in S$. Consider the extension property of $S$ in Definition \ref{loc-class-dfn}(b):

\bit
\item[b(1):] Since $fs'=sf'$ $\Ra$ $s^{-1}f=f's'^{-1}$ in $\C[S^{-1}]$, we can \ul{move the denomenators} $\{s^{-1}:s\in S\}$ in morphisms of $\C[S^{-1}]$ \ul{to the right} and so (up to equivalences) write each morphism in $\C[S^{-1}]$ as $fs^{-1}$, for $f\in \Mor \C$, $s\in S$. We write the morphism {\small $(s,f):=fs^{-1}:\bullet\sr{s^{-1}}{\ral}\bullet\sr{f}{\ral}\bullet$} as a diagram
\adjustbox{scale=0.7}{%
\bt
 & \bullet\ar[dl,"s"']\ar[dr,"f"] & \\
 \bullet && \bullet
\et} called a \ul{left roof}.

\item[b(2):] Since $s'f=f's$ $\Ra$ $fs^{-1}=s'^{-1}f'$ in $\C[S^{-1}]$, we can \ul{move the denomenators} $\{s^{-1}:s\in S\}$ in morphisms of $\C[S^{-1}]$ \ul{to the left} and so (up to equivalences) write each morphism in $\C[S^{-1}]$ as $s^{-1}f$, for $f\in \Mor \C$, $s\in S$. We write the morphism {\small $(f,s):=s^{-1}f:\bullet\sr{f}{\ral}\bullet\sr{s^{-1}}{\ral}\bullet$} as a diagram
\adjustbox{scale=0.7}{%
\bt
 & \bullet & \\
 \bullet\ar[ur,"f"] && \bullet\ar[ul,"s"']
\et} called a \ul{right roof}.
\eit

{\flushleft \ul{Observe}} that for any left roofs $(s,f),(t,g)\in \Mor_{C[S^{-1}]}(X,Y)$,
as in the diagram

\bc\adjustbox{scale=0.7}{\bt
  & & X'''\ar[dl,dashed,"f_s"']\ar[dr,dashed,"g_t"] & & \\
  & X'\ar[dl,"s"']\ar[drrr,"f"'pos=0.6] &  & X''\ar[dr,"g"]\ar[dlll,"t"pos=0.6] & \\
 X & & & & Y
\et}\ec
{\flushleft by} ``\ul{pulling back}'' the \ul{converging pair of morphisms} $s,t$ using the extension property of $S$, we know that there exist morphisms $f_s,g_t\in\Mor\C$
such that $sf_s=tg_t$ (but in general, we might have $gg_t\neq ff_s$, i.e., the diagram might not commute). This shows that our simplification of the expression for morphisms in $\C[S^{-1}]$ has possibly introduced \ul{new redundancy}, which can arise as follows: If we can choose $f_s,g_t$ such that (i) $gg_t=ff_s$, meaning the diagram commutes, and (ii) $sf_s\in S$, i.e., $(sf_s,gg_t)$ is a left roof, then both $f_s,g_t$ are necessarily invertible in $\C[S^{-1}]$ since $sf_s=tg_t$, (\blue{footnote}\footnote{Recall that if both $s$ and $sf_s$ (resp. $t$ and $tg_t$) are invertible, then so is $f_s$ (resp. $g_t$). Note however that, in general, $f_s,g_t$ are elements of $\Mor\C$ that do not have to lie in $S$ (except in special situations that might arise in practice).}), and so as morphisms in $\C[S^{-1}]$,
\bea
fs^{-1}=ff_s(sf_s)^{-1}=gg_t(tg_t)^{-1}=gt^{-1},\nn
\eea
i.e., we have the following relations of left roofs:
\bea
\bt(s,f)\ar[rr,equal,"f_s^{-1}~\txt{exists}"] && (sf_s,ff_s)\ar[rr,equal,"\txt{commutativity}"] && (tg_t,gg_t)\ar[rr,equal,"g_t^{-1}~\txt{exists}"] && (t,g).
\et\nn
\eea

\bea
\adjustbox{scale=0.7}{%
\bt
  & X'''\ar[d,dashed,"f_s"'] & \\
  & X'\ar[dl,"s"]\ar[dr,"f"'] & \\
X &    & Y
\et}=
\adjustbox{scale=0.7}{%
\bt
  & X'''\ar[ddl,dashed,"sf_s"']\ar[ddr,dashed,"ff_s"] & \\
  & X'\ar[dl,"s"]\ar[dr,"f"'] & \\
X &    & Y
\et}~~~~~~~~
\adjustbox{scale=0.7}{%
\bt
  & X'''\ar[ddl,dashed,"tg_t"']\ar[ddr,dashed,"gg_t"] & \\
  & X''\ar[dl,"t"]\ar[dr,"g"'] & \\
X &    & Y
\et}=
\adjustbox{scale=0.7}{%
\bt
  & X'''\ar[d,dashed,"g_t"']& \\
  & X''\ar[dl,"t"]\ar[dr,"g"'] & \\
X &    & Y
\et}\nn
\eea
Consequently, if (i),(ii) above hold, then we will say that $(s,f)$ and $(t,g)$, as diagrams in $\C\subset \C[S^{-1}]$, are \ul{equivalent left roofs} , written $(s,f)\sim(t,g)$.

Moreover, if we can show that $\sim$ is indeed an equivalence relation, then we can conclude that $\sim$ completely encodes all new redundancy due to the extension property of $S$.
\bc\fbox{\begin{minipage}{17.5cm}
\footnotesize
\blue{To see why, suppose $\sim$ is an equivalence relation. Then, the condition ``\ul{For every `pullback' of $\bullet\sr{s}{\ral}\bullet\sr{t}{\lal}\bullet$, at least one of (i),(ii) fails.}'' implies $[(s,f)]_\sim\cap[(t,g)]_\sim=\emptyset$, and so $(s,f)\neq(t,g)$ in $\C[S^{-1}]$.}
\end{minipage}}\ec
For this reason, in order to complete our characterization of morphisms of $\C[S^{-1}]$ in terms of left roofs, it will suffice to simply show (as done done in Proposition \ref{LocClassPrp}b) that $\sim$ is an equivalence relation on $\Mor_{\C[S^{-1}]}(X,Y)$.
\end{rmk}
\begin{summary}
Based on the above remark, we have the following characterization/properties of the localization by a localizing class of morphisms, which states that morphisms of $\C[S^{-1}]$ are given by
{\small\bea
\textstyle \Mor_{\C[S^{-1}]}(X,Y)={\big\{\txt{left roofs}~(s,f)=fs^{-1}:X\sr{s^{-1}}{\ral}Z\sr{f}{\ral}Y\big\}\over \sim}={\bigcup_{Z\in\Ob\C}~\Mor_\C(Z,X)\cap S~\times~\Mor_\C(Z,Y)\over\sim},\nn
\eea}where of course, ~$\Ob\C[S^{-1}]=Ob~\C$ ~by definition.
\end{summary}

\begin{prp}[\textcolor{blue}{Character/properties of the localization by a localizing class: \cite[Lemma 8, p.148]{gelfand-manin2010}}]\label{LocClassPrp}
Let $\C$ be a category and $S\subset\Mor\C$ a localizing class. Then $\C[S^{-1}]$ can be defined as follows.
\bit[leftmargin=0.9cm]
\item[(a)] $\Ob \C[S^{-1}]:=\Ob\C$
\item[(b)]  A morphism $X\ra Y$ in $\C[S^{-1}]$ is given by an equivalence class $[(s,f)]_\sim$ of left roofs, where
\bea
(s,f):=fs^{-1}=~
\adjustbox{scale=0.7}{%
\bt
 & X'\ar[dl,"s"']\ar[dr,"f"] & \\
 X && Y
\et}\nn
\eea
and $(s,f)\sim (t,g)$ if there exist morphisms $r,h\in\Mor\C$ giving a \ul{commutative diagram}
\bc
\adjustbox{scale=0.7}{%
\bt
  & & X'''\ar[dl,dashed,"r"']\ar[dr,dashed,"h"] & & \\
  & X'\ar[dl,"s"']\ar[drrr,"f"'pos=0.6] &  & X''\ar[dr,"g"]\ar[dlll,"t"pos=0.6] & \\
 X & & & & Y
\et} with ~$sr\in S$, ~i.e., ~
\adjustbox{scale=0.7}{%
\bt
 & X'''\ar[dl,"sr"']\ar[dr,"gh"] & \\
 X && Y
\et} ~is a \ul{left roof}.
\ec
(\ul{Note}: It is clear that the morphisms $r,h\in\Mor\C$ in the above diagram are necessarily iso's in $\C[S^{-1}]$, and so in $\C[S^{-1}]$ we have $(s,f)=(sr,fr)=(sr,gh)=(th,gh)=(t,g)$.)

\item[(c)] The identity morphism $(id_X,id_X)$ is a left roof, for every $X\in\Ob\C=\Ob\C[S^{-1}]$.
\item[(d)] The composition $(t,g)\circ(s,f)$ of two left roofs $(s,f),(t,g)$ is given by the following left roof $(st',gf')$:
\bea
(t,g)\circ(s,f)~:=~
\adjustbox{scale=0.7}{%
\bt
  & & X''\ar[dl,dashed,"t'"']\ar[dr,dashed,"f'"] & & \\
  & X'\ar[dl,"s"']\ar[dr,"f"] &  & Y'\ar[dl,"t"']\ar[dr,"g"] & \\
 X & & Y & & Z
\et}~~~~\sim~~~~
\adjustbox{scale=0.7}{%
\bt
 & X''\ar[dl,"st'"']\ar[dr,"gf'"] & \\
 X && Z
\et}\nn
\eea
where the square is based on the extension property in Definition \ref{loc-class-dfn}b (i.e., $t^{-1}f=f't'^{-1}$ in $\C[S^{-1}]$).
\bc\fbox{\begin{minipage}{17.5cm}
\footnotesize
\blue{As a formal motivation for this definition of composition, observe that by the commutativity of the diagram $tf'=ft'$, which implies $t^{-1}f=f't'^{-1}$, and so
{\scriptsize $(t,g)(s,f)=(gt^{-1})(fs^{-1})=g(t^{-1}f)s^{-1}=g(f't'^{-1})s^{-1}=gf'(st')^{-1}=(st',gf')$}.}
\end{minipage}}\ec

\item[(e)] The localization functor~ $L:\C\ra \C[S^{-1}]$ is given by
{\small\bea
\Big(X\sr{f}{\ral}Y\Big)\mapsto \Big(L(X)\sr{L(f)}{\ral}L(Y)\Big):=\Big(X\sr{[f]}{\ral}Y\Big)=\Big(X\sr{(1_X,f)}{\ral}Y\Big),\nn
\eea}
where $[f]=(1_X,f)$ is of course also written as $f 1_X^{-1}$, or simply as $(1,f)$, $f1^{-1}$, $f$, etc.
\eit
\end{prp}
\begin{proof}
It is clearly sufficient to simply verify that (i) the definitions in each part are consistent with (or satisfy) the defining properties of the localization and (ii) all parts (a)-(e) together exhaust the defining properties of the localization. In part (e) in particular, we verify that the given functor satisfies the universal property of the localization.
{\flushleft (a)} This is just part of the original definition of the localization.
{\flushleft (b)} By construction of the localization, it is clear (by the extension property of $S$) that every morphism can be written as a left roof, and every left roof is a morphism. Thus, we only need to show the given relation $\sim$ is an equivalence relation. \ul{Symmetry}, i.e., $(s,f)\sim(t,g)$ $\Ra$ $(t,g)\sim(s,f)$, holds by construction, and \ul{reflexivity} $(s,f)\sim(s,f)$ follows from the commutative diagram
\bc
\adjustbox{scale=0.7}{%
\bt
  & & X'\ar[dl,dashed,"id_{X'}"']\ar[dr,dashed,"id_{X'}"] & & \\
  & X'\ar[dl,"s"']\ar[drrr,"f"'pos=0.6] &  & X'\ar[dr,"f"]\ar[dlll,"s"pos=0.6] & \\
 X & & & & Y
\et}
\ec
To prove \ul{transitivity}, suppose $(s,f)\sim({s'},{f'})\sim({s''},{f''})$ are roofs. Then we have a diagram

\bea
\adjustbox{scale=0.8}{%
\bt[column sep=small]
    &    &    & W\ar[dl,dotted,"S\ni s_0"']\ar[dr,dotted,"f_1"]   &     &   Z'''\ar[ll,bend right=10,dashed,"s'_1\in S"']   &  \\
    &    & Z'\ar[dl,dashed,"r"']\ar[rd,dashed,"h"] &     & Z''\ar[dl,dashed,"r'"']\ar[dr,dashed,"h'"] &      &  \\
    & X'\ar[dl,"s"']\ar[drrrrr,"f"'pos=0.6] &    & X''\ar[dlll,"{s'}"'pos=0.6]\ar[drrr,"{f'}"pos=0.6] &     & X'''\ar[dlllll,"{s''}"pos=0.6]\ar[dr,"{f''}"] &  \\
 X  &    &    &     &     &      & Y
\et}~~~~
\fbox{\begin{minipage}{9.5cm}
\ul{Left roof relations}: Recall that as morphisms in $\C[S^{-1}]$,
{\footnotesize
\begin{align}
&(s,f)=(sr,fr)=(sr,f'h)=(s'r',f'h)=(s',f'),\nn\\
&(s',f')=(s'r',f'r')=(s'r',f''h')=(s''h',f''h')=(s'',f'').\nn
\end{align}}
\end{minipage}}
\nn\eea
in which the \ul{lower part from $Z',Z''$ downwards is commutative}, with left roofs $(sr,f'h),(s'r',f''h')$. So, by the localizing property, there exist $s_0\in S$, $f_1\in\Mor\C$ giving a commutative diagram as follows:
\bea
\adjustbox{scale=0.7}{%
\bt
  & W\ar[dl,dashed,"s_0\in S"']\ar[dr,dashed,"f_1"] &  \\
 Z'\ar[dr,"sr"'] &  & Z''\ar[dl,"{s'}r'\in S"]   \\
  & X &   \\
\et}~~~~~~~~srs_0=s'r'f_1.\nn
\eea
By the commutativity in the above two diagrams, we have ${s'}r'f_1=srs_0={s'}hs_0$, and so by the localizing property, there exists $s'_1\in S$ such that $(r'f_1)s'_1=(hs_0)s'_1$. Thus, we get a commutative diagram
\bea
\adjustbox{scale=0.7}{%
\bt
  & & Z'''\ar[dl,dashed,"rs_0s'_1"']\ar[dr,dashed,"h'f_1s'_1"] & & \\
  & X'\ar[dl,"s"']\ar[drrr,"f"'pos=0.6] &  & X'''\ar[dlll,"{s''}"pos=0.6]\ar[dr,"{f''}"] & \\
 X & & & & Y
\et}\nn
\eea
which is a left roof since $s(rs_0s'_1)=(sr)(s_0s'_1)\in S$. Hence, ~$(s,f)\sim (s'',f'')$.

{\flushleft (c)} This is clear because $id_X\in S$ for every $X\in\Ob\C=\Ob\C[S^{-1}]$.

{\flushleft (d)} Here, we need to show that composition (as given) is well defined. That is, if ~$(s,f)\sim(\ol{s},\ol{f})$ and $(t,g)\sim(\ol{t},\ol{g})$, then $(t,g)\circ(s,f)\sim(\ol{t},\ol{g})\circ(\ol{s},\ol{f})$. It suffices to show that ~$(s,f)\sim(\ol{s},\ol{f})~\Ra~(t,g)\circ(s,f)\sim(t,g)\circ(\ol{s},\ol{f})$ and $(s,f)\circ(u,h)\sim(\ol{s},\ol{f})\circ(u,h)$ ~for any left roofs $(t,g),(u,h)$. (\blue{footnote}\footnote{If ~$(s,f)\sim(\ol{s},\ol{f})~\Ra~(t,g)\circ(s,f)\sim(t,g)\circ(\ol{s},\ol{f})$ and $(s,f)\circ(u,h)\sim(\ol{s},\ol{f})\circ(u,h)$ ~for any left roofs $(t,g),(u,h)$, then
\bea
(s,f)\sim(\ol{s},\ol{f}),~(t,g)\sim(\ol{t},\ol{g})~~~~\Ra~~~~
(t,g)\circ(s,f)\sim (t,g)\circ(\ol{s},\ol{f})\sim (\ol{t},\ol{g})\circ(\ol{s},\ol{f}).\nn
\eea}). So, assume $(s,f)\sim(\ol{s},\ol{f})$.

$\bullet$ To show $(t,g)\circ(s,f)\sim(t,g)\circ(\ol{s},\ol{f})$, consider the diagram
\bc
\adjustbox{scale=0.7}{%
\bt
&&  & X_4\ar[ddll,dashed,"S\ni x_4"']\ar[ddr,dashed,"x_4'"]  &   & X_5\ar[ll,dotted,"x_5\in S"']& \\
      &     &        &     &      &     & \\
 & X_3\ar[ddl,dotted,"x_3"']\ar[ddr,dotted,"x_3'\in S"]& &  &  Y_3\ar[ddll,dotted,"S\ni y_3"']\ar[ddr,dotted,"y_3'"]  &  & \\
      &     &        &     &      &     & \\
  X_2\ar[dddr,dashed,"S\ni x_2"']\ar[dddrrrr,dashed,"x_2'"'pos=0.2] &  &  Y_2\ar[dddl,dashed,"y_2"'pos=0.2]\ar[ddd,dashed,"y_2'"pos=0.2] &  &  & Z_2\ar[dddlll,dashed,"S\ni z_2"'pos=0.2]\ar[dddl,dashed,"z_2'"]   & \\
      &     &        &     &      &     & \\
      &     &        &     &      &     & \\
  & X_1\ar[ddl,"s"']\ar[ddrr,"f"'] &  Y_1\ar[ddll,"\ol{s}"]\ar[ddr,"\ol{f}"] &  & Z_1\ar[ddl,"t"']\ar[ddr,"g"]  &   & \\
      &     &        &     &      &     & \\
  X  & &  & Y  &  & Z  &
\et}
\ec
in which the lower part from $X_2,Y_2,Z_2$ downwards is commutative, with the following:
\bit
\item[(i)] The left roof $XY_2Y$ implements the equivalence ~~$(s,f)\sim(\ol{s},\ol{f})$.
\item[(ii)] The left roof $XX_2Z$ denotes the composition ~~$(t,g)\circ (s,f)$.
\item[(iii)] The left roof $XZ_2Z$ denotes the composition ~~$(t,g)\circ(\ol{s},\ol{f})$.
\eit
Using the localizing properties, the morphisms in the diagram (except $x_5$) exist such that the primitive (or undivided) squares commute, along with the lower part from $X_2,Y_2,Z_2$ downwards. Thus (with the $X_5$-link not yet active),
\bea
&& (tz'_2)y'_3x'_4=\ol{f}(z_2y'_3)x'_4=(\ol{f}y'_2)(y_3x'_4)=f(y_2x'_3)x_4=(fx_2)x_3x_4=tx'_2x_3x_4,\nn\\
&&~~\Ra~~\txt{the morphism ~$x_5\in S$~ exists such that}~~z'_2y'_3x'_4x_5=x'_2x_3x_4x_5,\nn\\
&&~~\Ra~~gz'_2y'_3x'_4x_5=gx'_2x_3x_4x_5.\nn
\eea
Since it is also clear that {\small $sx_2x_3x_4=\ol{s}z_2y'_3x'_4$}, so {\small $sx_2x_3x_4x_5=\ol{s}z_2y'_3x'_4x_5$}, we get the commutative diagram
\bc
\adjustbox{scale=0.7}{%
\bt
  & & X_5\ar[dl,dashed,"x_3x_4x_5"']\ar[dr,dashed,"y'_3x'_4x_5"] & & \\
  & X_2\ar[dl,"sx_2"']\ar[drrr,"gx'_2"'pos=0.6] &  & Z_2\ar[dlll,near end,"\ol{s}z_2"]\ar[dr,"gz'_2"pos=0.6] & \\
 X & & & & Z
\et}
\ec
forming a left roof such that
\bea
(s,f)\circ(t,g)~~~~\sim~~~~
\adjustbox{scale=0.7}{%
\bt[column sep =small]
 & X_5\ar[dl,"sx_2x_3x_4x_5"']\ar[dr,"gz_2'y_3'x_4'x_5"] & \\
X &  & Z
\et}
~~~~\sim~~~~(\ol{s},\ol{f})\circ(t,g).\nn
\eea

$\bullet$ Similarly, to show $(s,f)\circ(u,h)\sim(\ol{s},\ol{f})\circ(u,h)$, we consider the diagram
\bc
\adjustbox{scale=0.7}{%
\bt
      &     &  \substack{X_4~(\txt{via}~u),\\\txt{then via}~sy_2}\ar[ddl,dashed,"S\ni x_4"']\ar[ddrr,dashed,"x_4'"]   &     &      & X_5\ar[lll,dotted,"x_5\in S"'] &  \\
      &     &        &     &      &     &  \\
      & X_3~(\txt{via}~s)\ar[ddl,dotted,"S\ni x_3"']\ar[ddrr,dotted,"x_3'"] &        &     & Y_3~(\txt{via}~\ol{s})\ar[ddl,dotted,"y_3"']\ar[ddr,dotted,"y_3'\in S"]   &     &  \\
      &     &        &     &      &     &  \\
  X_2\ar[dddr,dashed,"S\ni x_2"']\ar[dddrrr,dashed,near start,"x_2'"] &     &        & Y_2\ar[ddd,dashed,near start,"y_2"']\ar[dddr,dashed,near start,"y_2'"] &      & Z_2\ar[dddllll,dashed,near start,"z_2\in S"]\ar[dddl,dashed,"z_2'"] &  \\
      &     &        &     &      &     &  \\
      &     &        &     &      &     &  \\
      & X_1\ar[ddl,"u"']\ar[ddr,"h"] &        & Y_1\ar[ddl,"s"']\ar[ddrr,"f"'] & Z_1\ar[ddll,"\ol{s}"]\ar[ddr,"\ol{f}"]  &     &  \\
      &     &        &     &      &     &  \\
  X   &     &   Y    &     &      &  Z  &
\et}
\ec
in which the lower part from $X_2,Y_2,Z_2$ downwards is commutative, with the following:
\bit
\item[(i)] The left roof $YY_2Z$ implements the equivalence ~~$(s,f)\sim(\ol{s},\ol{f})$.
\item[(ii)] The left roof $XX_2Z$ denotes the composition ~~$(s,f)\circ (u,h)$.
\item[(iii)] The left roof $XZ_2Z$ denotes the composition ~~$(\ol{s},\ol{f})\circ (u,h)$.
\eit
Using the localizing properties, the morphisms in the diagram (except $x_5$) exist such that the primitive (or undivided) squares commute, along with the lower part from $X_2,Y_2,Z_2$ downwards. Note that in this case, $X_4$ is obtained via $u$ in the sense that $X_3\sr{x_4}{\lal}X_4\sr{x_4'}{\ral}Y_3$ is obtained by ``pulling back'' the diagram $\bt X_3 & X\ar[l,"ux_2x_3"']\ar[r,"uz_2y_3'"] & Y_3\et$, and then ensuring commutativity of the upper simple square by equalization through the following equality (with the $X_5$-link not yet active):
\bea
&&(sy_2)y_3x_4'=(\ol{s}z_2'y_3')x_4'=(hz_2y_3')x_4'=hx_2x_3x_4=sx_2'x_3x_4=(sy_2)x_3'x_4,\nn\\
&&~~\Ra~~~~y_3x_4'x_5=x_3'x_4x_5~~~~\txt{for some}~~x_5\in S.\nn
\eea
Thus, we get the commutative diagram
\bc
\adjustbox{scale=0.7}{%
\bt
  & & X_5\ar[dl,dashed,"x_3x_4x_5"']\ar[dr,dashed,"y'_3x'_4x_5"] & & \\
  & X_2\ar[dl,"ux_2"']\ar[drrr,"fx'_2"'pos=0.6] &  & Z_2\ar[dlll,near end,"u z_2"]\ar[dr,"\ol{f}z'_2"pos=0.6] & \\
 X & & & & Z
\et}
\ec
forming a left roof such that
\bea
(u,h)\circ(s,f)~~~~\sim~~~~
\adjustbox{scale=0.7}{%
\bt[column sep =small]
 & X_5\ar[dl,"ux_2x_3x_4x_5"']\ar[dr,"\ol{f}z_2'y_3'x_4'x_5"] & \\
X &  & Z
\et}
~~~~\sim~~~~(u,h)\circ(\ol{s},\ol{f}).\nn
\eea

{\flushleft (e)} Parts (a) to (d) define the \ul{category} $\R$ whose \ul{objects} are the same as those of $\C$, and whose \ul{morphisms} are left roofs. Define a functor~ $L:\C\ra \R$~ by
\bea
\Big(X\sr{f}{\ral}Y\Big)\mapsto \Big(L(X)\sr{L(f)}{\ral}L(Y)\Big):=\Big(X\sr{(1_X,f)}{\ral}Y\Big).\nn
\eea
We will show $\R$ satisfies the universal property for $\C[S^{-1}]$. Consider any functor $\C\sr{F}{\ral}\D$ such that $F(s)$ is an iso in $\C$ for $s\in S$. We want to show there exists a unique functor $\R\sr{G}{\ral}\C$ such that $GL=F$.
\bc
\adjustbox{scale=0.7}{%
\bt
\C\ar[d,"F"']\ar[rr,"L"] & & \R \\
\D\ar[from=urr,dashed,"G"] &  &
\et}\ec
Define $G$ as follows: $\bt G:\R\ra\D,~\Big(X\ar[r,"{(s,f)}"] & Y\Big)\ar[r,mapsto] & \Big(X\ar[rr,"{F(f)F(s)^{-1}}"] && Y\Big).\et$  Then $G$ is \ul{well defined} (\blue{footnote}\footnote{In the sense that $(s,f)\sim(t,g)$ $\Ra$ $G(s,f)=G(t,g)$.}) because for any left roofs $(s,f)\sim(t,g)$ and any $r,h\in\Mor\C$ giving a commutative diagram
\bc
\adjustbox{scale=0.7}{%
\bt
  & & X'''\ar[dl,dashed,"r"']\ar[dr,dashed,"h"] & & \\
  & X'\ar[dl,"s"']\ar[drrr,"f"'pos=0.6] &  & X''\ar[dr,"g"]\ar[dlll,"t"pos=0.6] & \\
 X & & & & Y
\et} with ~$sr\in S$, ~i.e., ~
\adjustbox{scale=0.7}{%
\bt
 & X'''\ar[dl,"sr"']\ar[dr,"gh"] & \\
 X && Y
\et} ~is a \ul{left roof},
\ec
where $sr\in S$ and $th\in S$ imply $F(r)^{-1}$ and $F(h)^{-1}$ exist respectively, we have
\bea
&&G(s,f)=F(f)F(s)^{-1}=F(f)F(r)(F(s)F(r))^{-1}=F(fr)F(sr)^{-1}=F(fr)F(th)^{-1}=F(gh)F(th)^{-1}\nn\\
&&~~~~=F(g)F(h)(F(t)F(h))^{-1}=F(g)F(t)^{-1}=G(t,g).\nn
\eea
It is also clear that $GL=F$, and $G(id_X)=id_{G(X)}$ for any $X\in\Ob\C$. It remains to check that $G$ preserves composition (and so is a functor). For any left roofs $(s,f),(t,g)\in\Mor\C[S^{-1}]$, using the composition rule
\bea
(t,g)\circ(s,f)~:=~
\adjustbox{scale=0.7}{%
\bt
  & & X''\ar[dl,dashed,"t'"']\ar[dr,dashed,"f'"] & & \\
  & X'\ar[dl,"s"']\ar[dr,"f"] &  & Y'\ar[dl,"t"']\ar[dr,"g"] & \\
 X & & Y & & Z
\et}~~~~\sim~~~~
\adjustbox{scale=0.7}{%
\bt
 & X''\ar[dl,"st'"']\ar[dr,"gf'"] & \\
 X && Z
\et}\nn
\eea
we see that $G$ preserves composition as follows:
\begin{align}
&G((t,g)(s,f))=G(st',gf')=F(gf')F(st')^{-1}=F(g)F(f')F(t')^{-1}F(s)^{-1}\sr{(m)}{=}F(g)F(t)^{-1}F(f)F(s)^{-1}\nn\\
&~~~~=G(t,g)F(s,f),\nn
\end{align}
where step (m) holds because $ft'=tf'$ implies $F(f')F(t')^{-1}=F(t)^{-1}F(f)$.

Hence (up to equivalence of categories) we have $\R=\C[S^{-1}]$.
\end{proof}

\begin{rmk*}[\textcolor{blue}{Factorization and Inversion of left roofs}]
Let $\C$ be a category, $S\subset\Mor\C$ a localizing class, $s,t\in S$ and $f,g\in\Mor\C$.
\begin{enumerate}[leftmargin=0.8cm]
\item $(s,f)\sim(1,f)(s,1)$ ~by the commutative diagram
\bea
(1,f)\circ(s,1):=\adjustbox{scale=0.7}{%
\bt
  & & X'\ar[dl,dashed,"1"']\ar[dr,dashed,"1"] & & \\
  & X'\ar[dl,"s"']\ar[dr,"1"] &  & X'\ar[dl,"1"']\ar[dr,"f"] & \\
 X & & X' & & Y
\et}~~~~\sim~~~~
\adjustbox{scale=0.7}{%
\bt
 & X'\ar[dl,"s"']\ar[dr,"f"] & \\
 X && Y
\et}~~=~~(s,f).\nn
\eea
\item $(1,s)^{-1}\sim(s,1)$, ~hence ~$(s,f)\sim(1,f)(s,1)\sim(1,f)(1,s)^{-1}$, ~by the commutative diagram
\bea
(s,1)\circ(1,s):=
\adjustbox{scale=0.7}{%
\bt
  & & X\ar[dl,dashed,"1"']\ar[dr,dashed,"1"] & & \\
  & X\ar[dl,"1"']\ar[dr,"s"] &  & X\ar[dl,"s"']\ar[dr,"1"] & \\
 X & & Y & & X
\et}~~~~\sim~~~~
\adjustbox{scale=0.7}{%
\bt
 & X\ar[dl,"1"']\ar[dr,"1"] & \\
 X && X
\et}~~=(1,1),\nn
\eea
where, by (2) above, $(1,1)\sim(s,s)=(1,s)\circ(s,1)$ via the commutative diagram
\bea
\adjustbox{scale=0.7}{%
\bt
  & &X\ar[dl,dashed,"1"']\ar[dr,dashed,"s"] & & \\
  & X\ar[dl,"s"']\ar[drrr,"s"'pos=0.6] &  & Y\ar[dlll,near end,"1"]\ar[dr,"1"pos=0.6] & \\
 Y & & & & Y
\et}~~=~~(s,s)=(1,s)\circ(s,1)~~=~~
\adjustbox{scale=0.7}{%
\bt
  & & X\ar[dl,dashed,"1"']\ar[dr,dashed,"1"] & & \\
  & X\ar[dl,"s"']\ar[dr,"1"] &  & X\ar[dl,"1"']\ar[dr,"s"] & \\
 Y & & Y & & Y
\et}\nn
\eea

\item $(1,g)\circ(s,f)\sim(s,gf)$ ~by the commutative diagram
\bea
(1,g)\circ(s,f):=\adjustbox{scale=0.7}{%
\bt
  & & X'\ar[dl,dashed,"1"']\ar[dr,dashed,"f"] & & \\
  & X'\ar[dl,"s"']\ar[dr,"f"] &  & Z\ar[dl,"1"']\ar[dr,"g"] & \\
 X & & Z & & Y
\et}~~~~\sim~~~~
\adjustbox{scale=0.7}{%
\bt
 & X'\ar[dl,"s"']\ar[dr,"gf"] & \\
 X && Y
\et}~~=~(s,gf)\nn
\eea
\item $(t,g)\circ(1,f)\sim(t',gf')$, ~where $tf'=ft'$, ~by the commutative diagram
\bea
(t,g)\circ(1,f):=\adjustbox{scale=0.7}{%
\bt
  & & X'\ar[dl,dashed,"t'"']\ar[dr,dashed,"f'"] & & \\
  & X\ar[dl,"1"']\ar[dr,"f"] &  & Y\ar[dl,"t"']\ar[dr,"g"] & \\
 X & & Z & & Y
\et}~~~~\sim~~~~
\adjustbox{scale=0.7}{%
\bt
 & X'\ar[dl,"t'"']\ar[dr,"gf'"] & \\
 X && Y
\et}~~=(t',gf')\nn
\eea
\end{enumerate}
\end{rmk*}

\section{Triangulation of Homotopy Category. Additivity of Derived Category}
\begin{lmm}[\textcolor{blue}{Triangulation of $\H\A$}]
Let $\A=R$-mod (for some ring $R$). The homotopy category $\H\A:={\A_0^\Integer\over\simeq}$ is triangulated.
\end{lmm}
\begin{proof}
We define \ul{translation} in $\H\A$ to be as in $\A_0^\Integer$, i.e., for any $n\in\Integer$, $T^n:\H\A\ra\H\A$ is given by
{\footnotesize\begin{align}
\Big((X_\ast,d_\ast^X)\sr{f_\ast}{\ral}(Y_\ast,d_\ast^Y)\Big)\mapsto \Big((X_\ast[n],d_\ast^{X[n]})\sr{f_\ast[n]}{\ral}(Y_\ast[n],d_\ast^{Y[n]})\Big):=\Big((X_{\ast+n},(-1)^nd_{\ast+n}^X)\sr{(-1)^nf_{\ast+n}}{\ral}(Y_{\ast+n},(-1)^nd_{\ast+n}^Y)\Big).\nn
\end{align}}
Based on our earlier discussion, we will say a triangle $(u,v,w):X_\ast\sr{u}{\ral}Y_\ast\sr{v}{\ral}Z_\ast\sr{w}{\ral}X_\ast[-1]$ in $\H\A$ is a \ul{distinguished triangle} if it is isomorphic to a triangle of the form $(\mu,i^{\mu},p^{\mu}):A_\ast\sr{\mu}{\ral}B_\ast\sr{i^{\mu}}{\ral}C(\mu)\sr{p^{\mu}}{\ral}A_\ast[-1]$, where $p^{\mu}i^{\mu}=i^{\mu}\mu=0$, $C(\mu):=B_\ast\oplus A_\ast[-1]$ is the cone of $\mu$, $i^{\mu}$ the canonical inclusion, and $p^{\mu}$ the canonical projection. Since it is clear that $\H\A$ is \ul{additive}, it remains to verify the triangulation axioms.

\bit[leftmargin=1.2cm]
\item[\ul{TR1}:] \ul{Construction axiom}: DT's arise from morphisms and isomorphic images of other DT's:
   \bit
   \item[(a)] In the DT ~$X\sr{1_X}{\ral}X\sr{i}{\ral}C(1_X)\sr{p}{\ral}X[-1]$,~ the mapping cone $C(1_X)$ is exact and $1_X$ is an isomorphism in $\A_0^\Integer$, and so $C(1_X)$ is contractible (i.e., homotopy equivalent to 0) in $\A_0^\Integer$, i.e., $C(1_X)$ is isomorphic to $0$ in $\H\A$. Hence ~$X\sr{1_X}{\ral}X\sr{}{\ral}0\sr{}{\ral}X[-1]$~ is a DT in $\H\A$.
   \item[(b)] It is clear that \emph{any triangle isomorphic to a DT is also a DT}, by the definition of a DT in $\H\A$ (since a composition of isomorphisms is an isomorphism).
   \item[(c)] Clearly, any morphism $X\sr{u}{\ral}Y$ can be completed into the distinguished triangle
       \bea
       X\sr{u}{\ral}Y\sr{i}{\ral}C(u)\sr{p}{\ral}X[-1].\nn
       \eea
   \eit

\item[\ul{TR2}:] \ul{Shift symmetry axiom}: The right-shift of {\small $S:~X\sr{u}{\ral}Y\sr{i}{\ral}C(u)\sr{p}{\ral}X[-1]$} is isomorphic to a distinguished triangle (through the mapping cylinder of $u$):
\bea
\adjustbox{scale=0.9}{\bt
R(S)\ar[r,draw=none,":"description]& C(u)[1]\ar[d,equal]\ar[r,"{-p[1]}"] & X\ar[d,equal]\ar[r,"u"] & Y\ar[d,"\cong"]\ar[r,"i"] & C(u)\ar[d,equal]\\
 & C(u)[1]\ar[d,equal]\ar[r,"{-p[1]}"] & X\ar[d,equal]\ar[r,"\ol{u}"] & Cyl(u)\ar[d,equal,"\txt{def.}"]\ar[r,"\ol{i}"] & C(u)\ar[d,equal]\\
 & C(u)[1]\ar[r,"{-p[1]}"] & X\ar[r,"\ol{u}"] & C(-p[1])\ar[r,"\ol{i}"] & C(u)
\et}\nn
\eea
If we now further right-shift, we get a DT that is isomorphic to the left-shift of $S$ as follows:
\bea
\adjustbox{scale=0.9}{\bt
R^2(S)\ar[d,"T^{-1}"]\ar[r,draw=none,":"description]& Y[1]\ar[d,"\cong"']\ar[r,"{-i[1]}"] & C(u)[1]\ar[d,"\cong"]\ar[r,"{-p[1]}"] & X\ar[d,"\cong"]\ar[r,"u"] & Y\ar[d,"\cong"]\\
L(S)\ar[r,draw=none,":"description]& Y\ar[r,"i"] & C(u)\ar[r,"p"] & X[-1]\ar[r,"{-u[-1]}"] & Y[-1]
\et}\nn
\eea

\item[\ul{TR3}:] \ul{Stability axiom}: Given~~
\adjustbox{scale=0.9}{\bt
 X\ar[d,"f"]\ar[r,"u"] & Y\ar[d,"g"]\ar[r,"i"] & C(u)\ar[d,dashed,"h"]\ar[r,"p"] & X[-1]\ar[d,"{f[-1]}"]  \\
 X'\ar[r,"u'"] & Y'\ar[r,"i'"] & C(u')\ar[r,"p'"] & X'[-1]
\et},~~ define $h:=g\oplus f[-1]$.
\bit
\item[(i)] $h$ is a chain morphism (morphism of complexes): With $h_n=-g_n\oplus f_{n-1}$ and the differentials
\bea
\left.
  \begin{array}{l}
   d_n^{C(u)}(y_n,x_{n-1}):=(d_n^Yy_n + u_{n-1}x_{n-1},-d_{n-1}^Xx_{n-1}), \\
   d_n^{C(u')}(y_n',x_{n-1}'):=(d_n^{Y'}y'_n + u'_{n-1}x'_{n-1},-d_{n-1}^{X'}x'_{n-1})
  \end{array}
\right\},~~\txt{we see that}~~h_{n-1}d_n^{C(u)}=d_n^{C(u')}h_n.\nn
\eea

\item[(ii)] It is also clear that the squares commute.
\eit

\item[\ul{TR4}:] \ul{Octahedral axiom}: Consider morphisms~ $A\sr{u}{\ral}B$ and $B\sr{v}{\ral}C$. Working up to isomorphism, we will assume the distinguished triangles extending $u,v,vu$ are of the standard form with cones. We have the diagram
\bc
\adjustbox{scale=0.8}{%
\bt
\ub{C(u)}_{B\oplus A[-1]}\ar[ddddrr,dashed, bend right=40,"f"']\ar[ddr,"p^u","{[-1]}"']  &  & B\ar[ll,"i^u"']\ar[ddr,"v"] &  & \ub{C(v)}_{C\oplus B[-1]}\ar[ll,"p^v"',"{[-1]}"]\ar[llll,bend right=40,"{i^u[-1]p^v}"',"{[-1]}"] \\
  &  & \curvearrowright &  &  \\
  & A\ar[uur,"u"]\ar[rr,"vu"] &  & C\ar[uur,"i^v"]\ar[ddl,"i^{vu}"'] &  \\
  &  &  &  &  \\
  &  & \ub{C(vu)}_{C\oplus A[-1]}\ar[uul,"p^{vu}"',"{[-1]}"]\ar[uuuurr,dashed, bend right=40,"g"'] &  &
\et}\ec
where ~$f:=v\oplus 1_{A[-1]}$ ~and ~$g:=1_C\oplus u[-1]$, ~and (as seen before) these are \ul{chain morphisms}.
\bit
\item[(i)] It is clear that the two \ul{new triangles created commute}, i.e., $p^{vu}f=p^u$ and $gi^{vu}=i^u$.
\item[(ii)] To see that {\small \bt[column sep=small] C(u)\ar[r,"{f}"]& C(vu)\ar[r,"{g}"]&C(v)\ar[rr,"{i^u[-1]p^v}"]&&C(u)[-1]\et} is a distinguished triangle, recall that
\bea
&&C_n(v)=C_n\oplus B_{n-1},~~~~C_n(f)=C_n(vu)\oplus C_{n-1}(u)=(C_n\oplus A_{n-1})\oplus(B_{n-1}\oplus A_{n-2}),\nn\\
&& d_n^{C(v)}=
\left[
  \begin{array}{cc}
    d_n^C & v_{n-1} \\
    0 & -d_{n-1}^B\\
  \end{array}
\right]=\left[
  \begin{array}{cc}
    d_n^C & -v_n[-1] \\
    0 & d_n^{B[-1]}\\
  \end{array}
\right]:
\left[
  \begin{array}{l}
    c_n \\
    b_{n-1} \\
  \end{array}
\right]\longmapsto
\left[
  \begin{array}{l}
    dc_n + v(b_{n-1}) \\
    -db_{n-1} \\
  \end{array}
\right],\nn\\
&& d_n^{C(f)}=
\left[
  \begin{array}{cc}
    d_n^{C(vu)} & f_{n-1} \\
    0 & -d_{n-1}^{C(u)}\\
  \end{array}
\right]=
\left[
  \begin{array}{cccc}
    d_n^C & (vu)_{n-1} & v_{n-1}     & 0 \\
    0 & -d_{n-1}^A     & 0           & 1_{A_{n-2}} \\
    0 &  0             & -d_{n-1}^B  & -u_{n-2} \\
    0 &  0             & 0           & d_{n-2}^A \\
  \end{array}
\right]~:~\left[
  \begin{array}{c}
    c_n \\
    a_{n-1} \\
    b_{n-1} \\
    a_{n-2} \\
  \end{array}
\right]\mapsto\nn\\
&&~~~~~~\mapsto~~\big[dc_n+vu(a_{n-1})+v(b_{n-1})~,~-da_{n-1}+a_{n-2}~,~-db_{n-1}-u(a_{n-2})~,~da_{n-2}\big]^T.\nn
\eea
We have a \ul{commutative diagram} in $\H\A$ (with $\al i^f=g$ and $p^f\beta=i^u[-1]p^v$ in $\A_0^\Integer$):
\bea\adjustbox{scale=0.9}{\bt
\label{HtyCatTR4}
C(u)\ar[r,"f"] & C(vu)\ar[r,"i^f"] & C(f)\ar[d,dashed, shift right,"\al"']\ar[from=d,dashed,shift right,"\beta"']\ar[rr,"p^f"] && C(u)[-1] \\
C(u)\ar[from=u,equal]\ar[r,"f"] & C(vu)\ar[from=u,equal]\ar[r,"g"] & C(v)\ar[rr,"{i^u[-1]p^v}"] && C(u)[-1]\ar[from=u,equal]  \\
\et}\eea
in which $\al:C(f)\ra C(v),~(c_n,a_{n-1},b_{n-1},a_{n-2})\mapsto \big(c_n,b_{n-1}+u(a_{n-1})\big)$, a \ul{chain morphism}, is a \ul{homotopy equivalence} with homotopy inverse the following natural inclusion:
\bea
\beta:C(v)\ra C(f),~~(c_n,b_{n-1})\mapsto(c_n,0,b_{n-1},0),~~~~\txt{(a \ul{chain morphism})}.\nn
\eea
We can see this directly as follows: With the maps ~$h_n:C_n(f)\ra C_{n+1}(f)$~ given by
\[
~~~~(c_n,a_{n-1},b_{n-1},a_{n-2})\mapsto (0,0,0,-a_{n-1}),~~~~\txt{i.e.},~~~~
h_n:=
\txt{\footnotesize$\left[
  \begin{array}{cccc}
    0 & 0  & 0 & 0 \\
    0 & 0 & 0 & 0 \\
    0 & 0  & 0 & 0 \\
    0 & -1  & 0 & 0 \\
  \end{array}
\right]$},~~~~\txt{for all}~~n,
\]
we can check that
\bea
&&\al\beta(c_n,b_{n-1})=\al(c_n,0,b_{n-1},0)=(c_n,b_{n-1}),~~\Ra~~\boxed{~\al_n\beta_n=1_{C_n(v)}~},\nn\\
&&(\beta\al-1_{C(f)})(c_n,a_{n-1},b_{n-1},a_{n-2})=(0,-a_{n-1},u(a_{n-1}),-a_{n-2}),\nn\\
&&~~\Ra~~\boxed{~\beta_n\al_n=1_{C_n(f)}+h_{n-1}d_n^{C(f)}+d_{n+1}^{C(f)}h_n~}.\nn
\eea

\eit
\eit
\end{proof}

The proof of the following result relies heavily on the triangulation of the homotopy category above.

\begin{thm}[\textcolor{blue}{Localizing property of quasiisomorphisms: \cite[Theorem 4, p.160]{gelfand-manin2010}}]\label{QuisLocLmm}
Let $\A$ $=$ $R$-mod for some ring $R$. The class $S_q:=\{\txt{quasi-iso's}\}\subset\Mor\A_0^\Integer$ is localizing in $\H\A:={\A_0^\Integer\over\simeq}$.
\end{thm}
\begin{proof}
We need to verify parts (a)-(c) of Definition \ref{loc-class-dfn}, as follows:
{\flushleft (a)} \ul{Composability property}: It is clear that $S_q$ contains the identity morphisms and is closed under composition.
{\flushleft(b)} \ul{Extension property}: Let $f\in \Mor_{\A_0^\Integer}(C_\ast,D_\ast)$ be any chain morphism, and $s\in S_q\cap \Mor_{\A_0^\Integer}(D_\ast',D_\ast)$ a quasi-iso. Since the category $\A_0^\Integer$ is abelian, morphisms $\{g:C_\ast'\ra D_\ast',~t:C_\ast'\ra C_\ast\}$ exist as the colimit (i.e., pullback-limit)
{\footnotesize
\bea
\varprojlim\{D_\ast'\sr{s}{\ral}D_\ast\sr{f}{\lal}C_\ast\},~~~~~~~~
\bt
C_\ast' \ar[d,dashed,"g"]\ar[r,dashed,"t"] & C_\ast\ar[d,"f"] \\
D_\ast'\ar[r,"s"] & D_\ast
\et\nn
\eea}
To ensure that $t\in S_q$ (with commutativity up to homotopy), we can instead obtain the diagram
\bea
\label{QuisLocEq1}\adjustbox{scale=0.9}{\bt
\overbrace{C(if)[1]}^{C_\ast'}\ar[d,dashed,"f'"]\ar[r,dashed,"s'"] & C_\ast\ar[d,"f"]\ar[r,"if"] & C(s)\ar[d,equal]\ar[r,"q"] & C(if)\ar[d,"{f'[-1]}"]\\
D_\ast'\ar[r,"s"] & D_\ast\ar[r,"i"] & C(s)\ar[r,"p"] & D_\ast'[-1]
\et},~~~~~~q(d_{n+1},d_n'):=(d_{n+1},d_n',0),
\eea
in which the lower sequence is (a distinguished triangle) based on Lemma \ref{TrngDiagLmm} and the upper sequence is introduced directly by observation, and the new \ul{chain morphisms} $s',f'$ are defined as follows:
\bea
&& s'\big(c_{n+1}(s),c_n\big)=s'\big(d_{n+1},d_n',c_n\big):=-c_n,~~~~ f'\big(c_{n+1}(s),c_n\big)=f'\big(d_{n+1},d_n',c_n\big):=d_n',\nn\\
&&~~\Ra~~(fs'-sf')(d_{n+1},d_n',c_n)=-f(c_n)-s(d_n').\nn
\eea
Recall that if $h\in \Mor_{\A_0^\Integer}(A_\ast,B_\ast)$, then $C(h)=B_\ast\oplus A_\ast[-1],~C_n(h)=B_n\oplus A_{n-1}$, and
\bea
d_n^{C(h_\ast)}=
\left[
  \begin{array}{cc}
    d_n^B & h_{n-1} \\
    0 & -d_{n-1}^A\\
  \end{array}
\right]=\left[
  \begin{array}{cc}
    d_n^B & -h_n[-1] \\
    0 & d_n^{A[-1]}\\
  \end{array}
\right]:
\left[
  \begin{array}{l}
    b_n \\
    a_{n-1} \\
  \end{array}
\right]\longmapsto
\left[
  \begin{array}{l}
    d_n^Bb_n + h_{n-1}a_{n-1} \\
    -d_{n-1}^Aa_{n-1} \\
  \end{array}
\right].\nn
\eea
Therefore, for $C(if)=C(s)\oplus C_\ast[-1]=D_\ast\oplus D_\ast'[-1]\oplus C_\ast[-1]$,
\bea
&&d_n^{C(if)}=\left[
  \begin{array}{cc}
    d_n^{C(s)} & (if)_{n-1} \\
    0 & -d_{n-1}^C\\
  \end{array}
\right]=\left[
  \begin{array}{ccc}
    d_n^D & s_{n-1}       & f_{n-1} \\
     0    & -d_{n-1}^{D'} & 0 \\
     0    & 0             &-d_{n-1}^C\\
  \end{array}
\right]\nn\\
&&~~\Ra~~
d_n^{C(if)[1]}=-d_{n+1}^{C(if)}:\left[
  \begin{array}{c}
    d_{n+1} \\
     d_n'  \\
     c_n \\
  \end{array}
\right]\mapsto
-\left[
  \begin{array}{c}
    d_{n+1}^Dd_{n+1}+s(d_n')+f(c_n) \\
     -d_n^{D'}d_n'  \\
     -d_n^Cc_n \\
  \end{array}
\right].\nn
\eea
Let $\rho_n:C(if)[1]_n\ra D_{n+1}$ be given by $\rho_n(d_{n+1},d_n',c_n):=d_{n+1}$. Then
\bea
&&\adjustbox{scale=0.9}{\bt
\cdots\ar[r] & C(if)[1]_{n+1}\ar[dd]\ar[r,"d_{n+1}^{C(if)[1]}"] & C(if)[1]_n\ar[ddl,"\rho_n"']\ar[dd]\ar[r,"d_n^{C(if)[1]}"] & C(if)[1]_{n-1}\ar[ddl,"\rho_{n-1}"']\ar[dd]\ar[r,"d_{n-1}^{C(if)[1]}"] &\cdots\\
 & & & & \\
\cdots\ar[r] & D_{n+1}\ar[r,"d_{n+1}^D"]                       & D_n\ar[r,"d_n^D"]                       & D_{n-1}\ar[r,"d_{n-1}^D"]                 &\cdots\\
\et}\nn\\
&&f_ns_n'-s_nf_n'=\rho_{n-1}d_n^{C(if)[1]}+d_{n+1}^D\rho_n,~~\Ra~~fs'\simeq sf'.\nn
\eea
Now, the upper sequence in (\ref{QuisLocEq1}) is (equivalent to) the following distinguished triangle:
\bea
C_\ast\sr{if}{\ral}C(s)\ral C(if)\sr{s'[-1]}{\ral}C_\ast[-1],\nn
\eea
in which $C(s)$ is exact (since $s$ is a quasi-iso). It thus follows that $s'$ is also a quasi-iso.

This proves part (1) of the extension property. Part (2) also follows by a dual/symmetric argument. That is, given $B\sr{f}{\lal}A\sr{s}{\ral}A'$, we can consider a diagram of the following form:
\[\bt
C(s)[1]\ar[d,equal]\ar[r,"{-p[1]}"]& A\ar[d,"f"]\ar[r,"s"]& A'\ar[d,dashed,"f'"]\ar[r,"i"]& C(s)\ar[d,equal] \\
C(s)[1]\ar[r,"{-fp[1]}"]& B\ar[r,dashed,"s'"]& C(-fp[1])\ar[r,"p"]& C(s)
\et\]
where the top row is a distinguished triangle, the bottom row is obtained by direct observation, and (as before) we need to define $f',s'$ and then construct a homotopy $\rho:A_\ast\ra C(-fp[1])_{\ast+1}$ such that $f's-s'f=d^{C(-fp[1])}\rho+\rho d^A$.

{\flushleft (c)} \ul{Equalizer-balance property}: Since $\H\A$ is additive, for morphisms $f,g\in\Mor\A_0^\Integer$, $s\in S_q$, the condition $[s]_\simeq[f]_\simeq=[s]_\simeq[g]_\simeq$ (i.e., $sf\simeq sg$) is equivalent to $[s]_\simeq([f]_\simeq-[g]_\simeq)=0$ (i.e., $s(f-g)\simeq 0$). Thus it suffices to show the following:
    \bit
    \item[$\bullet$] Given $f\in \Mor_{\A_0^\Integer}(C,D)$, there exists $s\in S_q\cap \Mor_{\A_0^\Integer}(D,D')$ with $sf\simeq 0$ $\iff$ there exists $s'\in S_q\cap \Mor_{\A_0^\Integer}(C',C)$ with $fs'\simeq 0$.
    \eit
{\flushleft ($\Ra$)} Assume we are given $f,s$ satisfying $sf\simeq 0$. Consider $s$ as part of a distinguished triangle $(u,s,w):\bt Z\ar[r,"u"] & D\ar[r,"s"] & D'\ar[r,"w"] & Z[-1]\et$ in $\A_0^\Integer$, the induced LES of homologies for which implies $H_n(Z)=0$ for all $n$ (since $s$ is a quasi-iso). Since $\Mor_{\A_0^\Integer}(C,-)$ is a homological functor, we have the exact sequence
\bea
\bt Hom(C,Z)\ar[r,"u_\ast"] & Hom(C,D)\ar[r,"s_\ast"] & Hom(C,D')\ar[r,"w_\ast"] & Hom(C,Z[-1])\et.\nn
\eea
Thus, because $0\simeq sf=s_\ast(f)$ (i.e., $f\in\ker s_\ast=\im u_\ast$ up to homotopy), we have $f\simeq u_\ast(g)$, i.e.,
\bea
f\simeq ug,~~~~\txt{for some}~~~~g:C\ra Z.\nn
\eea
Consider $g$ as part of a distinguished triangle $(s',g,w'):\bt[column sep=small] C'\ar[r,"s'"] & C\ar[r,"g"] & Z\ar[r,"w'"] & C'[-1]\et$ in $\A_0^\Integer$.
\bea
\adjustbox{scale=0.9}{\bt
                       & Z\ar[r,"u"]                  & D\ar[r,"s"]                & D'\ar[r,"w"]   & Z[-1]   \\
 C'[-1]\ar[from=r,dashed,"w'"'] & Z\ar[u,equal]\ar[from=r,dashed,"g"'] & C\ar[u,"f"']\ar[from=r,dashed,"s'"'] & C'                   &
\et}\nn
\eea
Then ~$gs'\simeq0~~\Ra~~fs'\simeq ugs'\simeq0$. Since $H_n(Z)=0$ for all $n$, the induced LES of homologies for the sequence $(s',g,w')$ implies $s'$ is a quasi-iso.

($\La$): Assume we are given $f,s'$ satisfying $fs'\simeq 0$. Then by symmetry in the above case, we also get $s$ such that $sf\simeq 0$.
\end{proof}

The above result means that quasiisomorphisms in $\A_0^\Integer$ are ``localizing up to homotopy'', in the sense that the required commutativity properties in Definition \ref{loc-class-dfn} hold up to homotopy as follows:

\begin{dfn}[\textcolor{blue}{\index{Homotopy-localizing class}{Homotopy-localizing class of morphisms}}]\label{loc-class-dfn2}
Let $\A$ be an abelian category. A class of chain morphisms $S\subset\Mor\A_0^\Integer$ is \ul{homotopy-localizing} (or \ul{localizing up to homotopy}) if it satisfies the following \ul{three} conditions:
\bit[leftmargin=0.9cm]
\item[(a)] \ul{$S$ is a ``composable'' class (``composability'' conditions)}
\bit
\item[(1)] $id_C$ (or $1_C$) $\in S$, for every object $C\in\Ob\C$.
\item[(2)] $s,t\in S$ ~$\Ra$~ $st\in S$.
\eit
\item[(b)] \ul{$S$ is extendable (Extension conditions)}:
\bit
\item[(1)] For every $f\in\Mor\A_0^\Integer$, $s\in S$ with the same codomain (i.e., $\cod f=\cod s$), there exist $f'\in \Mor\A_0^\Integer$, $s'\in S$ such that $fs'\simeq sf'$, i.e., the diagram
{\footnotesize
\bt
C' \ar[d,dashed,"s'"]\ar[r,dashed,"f'"] & D'\ar[d,"s"] \\
C\ar[r,"f"] & D
\et} commutes up to homotopy.

\item[(2)] (Equivalently) for every $f\in\Mor\A_0^\Integer$, $s\in S$ with the same domain (i.e., $\dom f=\dom s$), there exist $f'\in\Mor\A_0^\Integer$, $s'\in S$ such that $s'f\simeq f's$, i.e., the diagram
{\footnotesize
\bt
C \ar[d,"s"]\ar[r,"f"] & D\ar[d,dashed,"s'"] \\
C'\ar[r,dashed,"f'"] & D'
\et} commutes up to homotopy.
\eit
\item[(c)] \ul{$S$ is equalizer-balanced (Balance of left and right equalizers)}:  For any two morphisms $f,g\in\Mor\A_0^\Integer$,
\bea
sf\simeq sg~~\txt{for some}~~s\in S~~~~\iff~~~~fs'\simeq gs'~~\txt{for some}~~s'\in S.\nn
\eea
That is,~~~~
$\bt
 C\ar[rrr,bend left,"sf\simeq sg"]\ar[rr,shift left=1,"f"]\ar[rr,shift right=1,"g"'] && D\ar[r,dashed,"s"] & D'
\et$ ~~~~$\iff$~~~~
$\bt
  C'\ar[r,dashed,"s'"]\ar[rrr,bend left,"fs'\simeq gs'"] & C \ar[rr,shift left=1,"f"]\ar[rr,shift right=1,"g"'] && D
\et$
\eit
\end{dfn}

\begin{note}
In $\A_0^\Integer$, by replacing equality of maps $f=g$ with (i.e. equality up to) homotopy of maps $f\simeq g$, we see that certain results involving a ``\emph{localizing class}'' in $\H\A$ are also valid for a ``\ul{homotopy-localizing class}'' in $\A_0^\Integer$, with ``\emph{commutative diagram}'' replaced by ``\emph{\ul{homotopy-commutative diagram}}'' (defined to be a diagram that is commutative up to homotopy). In other words, some results proved for localizing classes in $\H\A$ can (and will) be directly applied to homotopy-localizing classes in $\A_0^\Integer$.
\end{note}

\begin{thm}[\textcolor{blue}{Additivity of $\D\A$}]\label{HtyCatAdd}
Let $\A=R$-mod (for some $R$). The derived category $\D\A$ is additive.
\end{thm}
\begin{proof}
We need to show (1) a \ul{zero object} $0\in\Ob\D\A$ exists (which is clear), (2) $\Mor_{\D\A}(X,Y)$ is an abelian group such that composition is distributive over addition, and (3) any two objects $X,Y\in\Ob\D\A$ have a \ul{product object} $X\times Y\in\Ob\D\A$. We only need to prove (2) and (3).

{\flushleft (2)}: Since quasi-iso's are localizing, every morphism in $\D\A$ is a roof. Given any morphisms $(s_1,f_1),(s_2,f_2)\in \Mor_{\D\A}(X,Y)$ represented by roofs $X\sr{s_1}{\lal}U_1\sr{f_1}{\ral}Y$ and $X\sr{s_2}{\lal}U_2\sr{f_2}{\ral}Y$ respectively, by extending the open triangle $U_1\sr{s_1}{\ral}X\sr{s_2}{\lal}U_2$ we have the following \ul{partly commutative} diagram:
\bea
\adjustbox{scale=0.7}{%
\bt
  & & U\ar[dl,dashed,"S_q\ni s_2'"']\ar[dr,dashed,"s_1'\in S_q"] & & \\
  & U_1\ar[dl,"s_1"']\ar[drrr,"f_1"'pos=0.6] &  & U_2\ar[dr,"f_2"]\ar[dlll,"s_2"pos=0.6] & \\
 X & & & & Y
\et}~~~~s:=s_1s_2'=s_2s_1'~~~~\txt{(a \ul{common denominator} for both morphisms)},\nn
\eea
which shows that $(s_1,f_1)\sim(s_1s_2',f_1s_2')=(s,f_1s_2')$ and $(s_2,f_2)\sim(s_2s_1',f_2s_1')=(s,f_2s_1')$. We \ul{define the sum} of these morphisms as ~$(s_1,f_1)+(s_2,f_2):\sim(s,f)$, ~where ~$f:=f_1s_2'+f_2s_1'$,~ i.e.,
\bea
\Big(X\sr{s_1}{\lal}U_1\sr{f_1}{\ral}Y\Big)+\Big(X\sr{s_2}{\lal}U_2\sr{f_2}{\ral}Y\Big)~:\sim~\Big(X\sr{s}{\lal}U\sr{f}{\ral}Y\Big),\nn
\eea
and the symbol $:\sim$ means ``$:=$ up to $\sim$-equivalence'' (and so will be written simply as $:=$ or $=$ once consistency with related operations has been verified). In seemingly more familiar notation, we can also express the sum above as
\bea
f_1s_1^{-1}+f_2s_2^{-1}~:\sim~(f_1s'_2+f_2s'_1)(s_1s_2')^{-1}=(f_1s'_2+f_2s'_1)(s_2s_1')^{-1}.\nn
\eea
To show addition is \ul{well-defined}, we need to check that if $(t_1,g_1)\sim (s_1,f_1)$ and $(t_2,g_2)\sim (s_2,f_2)$, then
\bea
(t_1,g_1)+(t_2,g_2)~\sim~(s_1,f_1)+(s_2,f_2).\nn
\eea
For $(t_1,g_1):~X\sr{t_1}{\lal}V_1\sr{g_1}{\ral}Y$ and $(t_2,g_2):~X\sr{t_2}{\lal}V_2\sr{g_2}{\ral}Y$, we also have a \ul{partly commutative} diagram:
\bea
\adjustbox{scale=0.7}{%
\bt
  & & V\ar[dl,dashed,"S_q\ni t_2'"']\ar[dr,dashed,"t_1'\in S_q"] & & \\
  & V_1\ar[dl,"t_1"']\ar[drrr,"g_1"'pos=0.6] &  & V_2\ar[dr,"g_2"]\ar[dlll,"t_2"pos=0.6] & \\
 X & & & & Y
\et}~~~~t:=t_1t_2'=t_2t_1'~~~~\txt{(a \ul{common denominator} for both morphisms)},\nn
\eea
which shows that $(t_1,g_1)\sim(t_1t_2',g_1t_2')=(t,g_1t_2')$ and $(t_2,g_2)\sim(t_2t_1',g_2t_1')=(t,g_2t_1')$. Hence, by the transitivity of $\sim$ (as expressed by the \ul{diagram below}), we have
\bea
(t_1,g_1)+(t_2,g_2)~:\sim~(t,g_1t_2'+g_2t_1')~\sim~(s,f_1s_2'+f_2s_1')~\sim:~(s_1,f_1)+(s_2,f_2).\nn
\eea
The \ul{following diagram} is a combination of the $(s_1,f_1)+(s_2,f_2)$ diagram and the $(t_1,g_1)+(t_2,g_2)$ diagram (in which only the \ul{purely quasi-iso's squares} commute, i.e., \ul{those without $f$'s or $g$'s}) with the top part obtained using the extension properties of quasi-iso's:

\bea\adjustbox{scale=0.7}{%
\bt
  &     &          &     & U'\ar[ddl,dotted,"t_2'''"']\ar[ddr,dotted,"s_1'''"] &     &          &   W\ar[lll,dotted,"w\in S_q"']   &   \\
  &     &          &     &    &     &          &      &   \\
  &     &          &  C'\ar[ddl,dotted,"\al_1'"']\ar[ddr,dotted,"s_1''"]   &    &  C''\ar[ddl,dotted,"t_2''"']\ar[ddr,dotted,"\beta_2'"]   &          &      &   \\
  &     &          &     &    &     &          &      &   \\
  &     & \ub{U}_+\ar[ddl,dashed,"s_2'"']\ar[ddr,dashed,"s_1'"] &     & C\ar[ddl,dotted,"\al_1"']\ar[ddr,dotted,"\beta_2"]  &     & \ub{V}_+\ar[ddl,dashed,"t_2'"']\ar[ddr,dashed,"t_1'"] &      &   \\
  &     &          &     &     &     &          &      &   \\
  & U_1\ar[ddl,"s_1"']\ar[ddrrrrrrr,near end, bend right=5,"f_1"'] &          & U_2\ar[ddlll,"s_2"']\ar[ddrrrrr,"f_2"'] &    & V_1\ar[ddlllll,"t_1"]\ar[ddrrr,"g_1"] &          &  V_2\ar[ddlllllll,near end, bend left=5,"t_2"]\ar[ddr,"g_2"] &   \\
  &     &          &     &    &     &          &      &   \\
X &     &          &     &    &     &          &      & Y
\et}\nn
\eea
In the above diagram (with the $W$-link not yet active),
\bea
&&t_1t_2'\beta_2's_1'''=t_1\beta_2t_2''s_1'''=t_1\beta_2s_1''t_2''' ~~\Ra~~t_2'\beta_2's_1'''w=\beta_2s_1''t_2'''w,~~\txt{for some}~~w\in S_q,\nn\\
&&~~\Ra~~g_1t_2'\beta_2's_1'''w=g_1\beta_2s_1''t_2'''w,\nn\\
&&~~\txt{and}~~~~s_2\al_1t_2''s_1'''=s_2s_1'\al_1't_2'''~~\Ra~~s_2\al_1t_2''s_1'''w=s_2s_1'\al_1't_2'''w,\nn
\eea
and so we get the following left roof expressing the equivalence $(s_1,f_1)+(s_2,f_2)\sim(t_1,g_1)+(t_2,g_2)$:
\bea
\adjustbox{scale=0.7}{%
\bt
  &                                     &W\ar[dl,dashed,"t_2'''w"']\ar[dr,dashed,"s_1'''w"] &                                    &  \\
  & C'\ar[dl,"s_2s_1'\al_1'"']\ar[drrr,"g_1\beta_2s_1''"'pos=0.6] &                                       & C''\ar[dlll,"s_2\al_1t_2''"pos=0.6]\ar[dr,"g_1t_2'\beta_2'"] &  \\
X &                                     &                                       &                                    & Y
\et}~~\sim~~
\adjustbox{scale=0.7}{%
\bt
   & W\ar[dl,"s_2s_1'\al_1't_2'''w"']\ar[dr,"g_1t_2'\beta_2's_1'''w"] &  \\
 X &                          & Y
\et}\nn
\eea
This addition operation clearly makes $G:=\Mor_{\D\A}(X,Y)$ an abelian group with $0_G:=(1,0)\sim(s,0)$, i.e., the class of the roof $(s,0)$ for any quasi-iso $s:X\ra X'$.

It remains to check that \ul{composition is distributive over addition} of morphisms: \ul{Distributing from the right}, by first composing before adding, we get from the diagram below:
\bea
fs^{-1}h\ld^{-1}+gt^{-1}h\ld^{-1}\sim fh'(\ld s')^{-1}+gh''(\ld t'')^{-1}\sim
\big(fh'(\ld t'')'''+gh''(\ld s')'''\big)\big((\ld s')(\ld t'')'''\big)^{-1}\nn
\eea
(At the ``$+$ vertex'' in the following diagram, only the purely quasi-iso's part is commutative.)
\bea
\label{HtyCatAddEq1}\adjustbox{scale=0.7}{%
\bt
  &                &                                   &                              &                                       &                   &    &&  \\
  &                &                                   & \ub{X''''}_+\ar[dl,dotted,"(\ld t'')'''"']
                                                           \ar[dr,dotted,"(\ld s')'''"]   &                                       &                   &    &&  \\
  &                & \ub{X''}_\circ\ar[dl,dashed,"s'"']
                        \ar[dr,dashed,near start,"h'"] &                              & \ub{X'''}_\circ\ar[dlll,dashed,near start,"t''"']
                                                                                         \ar[dr,dashed,"h''"]                 &                   &    &&  \\
  & X'\ar[dl,"\ld"']
      \ar[dr,"h"]  &                                   &Y'\ar[dl,"s"']
                                                         \ar[drrr,"f"']                &                                      & Y''\ar[dlll,"t"]
                                                                                                                                 \ar[dr,"g"]      &    &&  \\
X &                & Y                                 &                               &                                      &                   &  Z &&
\et}
\eea
and, by first adding before composing, we get from the diagram below:
\bea
&&(fs^{-1}+gt^{-1})h\ld^{-1}\sim (ft'+gs')(st')^{-1}h\ld^{-1}\sim (ft'+gs')h''(\ld(st')'')^{-1}\nn\\
&&~~~~\sim (ft'h''+gs'h'')(\ld(st')'')^{-1}.\nn
\eea
(At the ``$+$ vertex'' in the following diagram, only the purely quasi-iso's square is commutative.)
\bea
\label{HtyCatAddEq2}\adjustbox{scale=0.7}{%
\bt
  &                &    &                                          &                                &                 &    &&  \\
  &                &    & \ub{Z''''}_\circ\ar[ddll,dotted,"(st')''"']
                                        \ar[dr,dotted,"h''"]       &                                &                 &    &&  \\
  &                &    &                                          &\ub{Y'''}_+\ar[dl,dashed,"t'"']
                                                                               \ar[dr,dashed,"s'"]  &                 &    &&  \\
  & X'\ar[dl,"\ld"']
      \ar[dr,"h"]  &    &Y'\ar[dl,"s"']
                          \ar[drrr,"f"']                           &                                & Y''\ar[dlll,"t"]
                                                                                                       \ar[dr,"g"]    &    &&  \\
X &                & Y  &                                          &                                &                 &  Z &&
\et}
\eea
To show the above two left roofs (\ref{HtyCatAddEq1}) and (\ref{HtyCatAddEq2}) are equivalent, we need to be able to construct from them a \ul{commutative diagram} of the following form giving a left roof:
\bc
\adjustbox{scale=0.7}{%
\bt
  & &U\ar[dl,dashed,"r"']\ar[dr,dashed,"l"] & & \\
  & X''''\ar[dl,"a"']\ar[drrr,"\al"'pos=0.6] &  & Z''''\ar[dlll,"b"pos=0.6]\ar[dr,"\beta"] & \\
 X & & & & Z
\et} with ~$ar\in S$, ~i.e., ~
\adjustbox{scale=0.7}{%
\bt
 & U\ar[dl,"ar"']\ar[dr,"\beta l"] & \\
 X && Z
\et} ~is a \ul{left roof}.
\ec
 We can combine the two diagram (\ref{HtyCatAddEq1}) and (\ref{HtyCatAddEq2}) into a single diagram as follows (where at each ``$+$ vertex'', only the purely quasi-iso's part/square is commutative.):
\bea
\label{HtyCatAddEq3}\adjustbox{scale=0.7}{%
\bt
               &                 &  V\ar[dl,dotted,"S_q\ni x"']\ar[dr,dotted,"S_q\ni y"']               &                                                  &      U\ar[ll,dotted,"z\in S_q"']                  &                 &    &&  \\
               &\ub{X''''}_+\ar[dl,dotted,"(\ld t'')'''"']
                          \ar[dr,dotted, near start,"(\ld s')'''"']         &                  & \ub{Z''''}_\circ\ar[ddll,dotted,near start,bend right,"(s\bar{t}')''"']
                                                                                          \ar[dr,dotted,"\bar{h}''"]       &                                &                 &    &&  \\
\ub{X''}_\circ\ar[dr,dashed,"s'"']
 \ar[drrr,dashed,near end,"h'"']   &                 &  \ub{X'''}_\circ\ar[dl,dashed,near end,"t''"]
                                                  \ar[drrr,dashed,near end,"h''"']               &                             &\ub{Y'''}_+\ar[dl,dashed,near start,"\bar{t}'"']
                                                                                                                                \ar[dr,dashed,"\bar{s}'"]   &                 &    &&  \\
               & X'\ar[dl,"\ld"']
                   \ar[dr,"h"]   &                  &Y'\ar[dl,"s"']
                                                       \ar[drrr,"f"']                                                &                                & Y''\ar[dlll,"t"]
                                                                                                                                                          \ar[dr,"g"]    &    &&  \\
X              &                 & Y                &                                                                &                                &                 &  Z &&
\et}
\eea
Observe that (with the $U$-link not yet active),
\bea
&&s\bar{t}'\bar{h}''y=h(s\bar{t}')''y=ht''(\ld s')'''x=hs'(\ld t'')'''x=sh'(\ld t'')'''x,\nn\\
&&~~~~\Ra~~\bar{t}'\bar{h}''yz=h'(\ld t'')'''xz~~~~\txt{for some}~~~~z\in S_q,\nn\\
&&~~\Ra~~f\bar{t}'\bar{h}''yz=fh'(\ld t'')'''xz,\nn
\eea
and so we get the following left roof:
\bc
\adjustbox{scale=0.7}{%
\bt
  & &U\ar[dl,dashed,"xz"']\ar[dr,dashed,"yz"] & & \\
  & X''''\ar[dl,"\ld s'(\ld t'')'''"']\ar[drrr,"fh'(\ld t'')'''"'pos=0.6] &  & Z''''\ar[dlll,"\ld(s\bar{t}')''"pos=0.6]\ar[dr,"f\bar{t}'\bar{h}''"] & \\
 X & & & & Z
\et} with ~$\ld s'(\ld t'')'''xz\in S$, ~i.e., ~
\adjustbox{scale=0.7}{%
\bt
 & U\ar[dl,"\ld s'(\ld t'')'''xz"']\ar[dr,"f\bar{t}'\bar{h}''yz"] & \\
 X && Z
\et} ~is a \ul{left roof},
\ec
which proves that ~$fs^{-1}h\ld^{-1}+gt^{-1}h\ld^{-1}\sim (fs^{-1}+gt^{-1})h\ld^{-1}$.

Similarly, \ul{distributing from the left}, we get
\bea
&& h\ld^{-1}(fs^{-1}+gt^{-1})\sim h\ld^{-1}(ft'+gs')(st')^{-1}
\sim h(ft'+gs')''(st'\ld'')^{-1}~~~~\txt{(i.e., add before composing)},\nn\\
&&h\ld^{-1}fs^{-1}+h\ld^{-1}gt^{-1}\sim hf'(s\ld')^{-1}+hg''(t\ld'')^{-1}
\sim \big(hf'(t\ld')'''+hg''(s\ld'')'''\big)\big((s\ld'')(t\ld')'''\big)^{-1}\nn\\
&&~~~~\sim h\big(f'(t\ld')'''+g''(s\ld'')'''\big)\big((s\ld'')(t\ld')'''\big)^{-1}~~~~\txt{(i.e., compose before adding)},\nn
\eea
where, by the obvious \ul{left-right symmetry} of the situation, an analogous argument (as above) should show that ~$ h\ld^{-1}(fs^{-1}+gt^{-1})\sim h\ld^{-1}fs^{-1}+h\ld^{-1}gt^{-1}$.

{\flushleft (3)}: Let $X,Y\in\Ob\H\A=\Ob\D\A$ ($:=\Ob\A_0^\Integer$), and recall that in $\H\A$, $X\times Y\cong X\oplus Y$. Given any morphisms $f:X\ra Z$, $g:Y\ra Z$ in $\H\A$, there exists a morphism $h_{f,g}:X\oplus Y\ra Z$ in $\H\A$ such that the following diagram commutes:
\bea
\adjustbox{scale=0.7}{%
\bt
X\ar[ddr,bend right,"f"']\ar[dr,"i_X"] &                            &  Y\ar[ddl,bend left,"g"]\ar[dl,"i_Y"'] \\
                                       & X\oplus Y\ar[d,dashed,"h_{f,g}"] &  \\
                                       & Z                          &
\et}\nn
\eea
Similarly, in $\D\A$, it suffices to check that given morphisms $(s,f):X\sr{s}{\lal}X'\sr{f}{\ral} Z$ and $(t,g):Y\sr{t}{\lal}Y'\sr{g}{\ral} Z$ in $\D\A$, there exists a morphism $(\ld,h):X\oplus Y\sr{\ld}{\lal}(X\oplus Y)'\sr{h}{\ral} Z$ in $\D\A$ such that the following diagram commutes:
\bea
\adjustbox{scale=0.7}{%
\bt
X\ar[ddr,bend right,"{(s,f)}"']\ar[dr,"{(1,i_X)}"] &                            &  Y\ar[ddl,bend left,"{(t,g)}"]\ar[dl,"{(1,i_Y)}"'] \\
                                       & X\oplus Y\ar[d,dashed,"{(\ld,h)}"] &  \\
                                       & Z                          &
\et}\nn
\eea
Set ~$(\ld, h):=(s\oplus t,h_{f,g}):X\oplus Y\sr{s\oplus t}{\lal}X'\oplus Y'\sr{h_{f,g}}{\ral} Z$. Then we get the following:
\bea
(s\oplus t,h_{f,g})\circ(1,i_X):~
\adjustbox{scale=0.7}{%
\bt
  & & X'\ar[dl,dashed,"s"']\ar[dr,dashed,"i_{X'}"] & & \\
  & X\ar[dl,"1"']\ar[dr,"i_X"] &  & X'\oplus Y'\ar[dl,"s\oplus t"']\ar[dr,"h_{f,g}"] & \\
 X & & X\oplus Y & & Z
\et}~~~~\sim~~~~
\adjustbox{scale=0.7}{%
\bt
 & X'\ar[dl,"s"']\ar[dr,"f"] & \\
 X && Z
\et}\nn
\eea
and
\bea
(s\oplus t,h_{f,g})\circ(1,i_Y):~
\adjustbox{scale=0.7}{%
\bt
  & & Y'\ar[dl,dashed,"t"']\ar[dr,dashed,"i_{Y'}"] & & \\
  & Y\ar[dl,"1"']\ar[dr,"i_Y"] &  & X'\oplus Y'\ar[dl,"s\oplus t"']\ar[dr,"h_{f,g}"] & \\
 Y & & X\oplus Y & & Z
\et}~~~~\sim~~~~
\adjustbox{scale=0.7}{%
\bt
 & Y'\ar[dl,"t"']\ar[dr,"g"] & \\
 Y && Z
\et}\nn
\eea
\end{proof}

\begin{crl}[\textcolor{blue}{of the proof of Theorem \ref{HtyCatAdd}}]
Let $\A$ be an additive category. If $S\subset\Mor\A$ is a localizing class, then $\A[S^{-1}]$ is an additive category.
\end{crl}

\section{Derived Category of a Subcategory}
\begin{prp}[\textcolor{blue}{\cite[Proposition 10, p.151]{gelfand-manin2010}}]\label{SubLocPrp}
Let $\C$ be a category, and $S\subset\Mor\C$ a localizing class. Let $\F\subset\C$ be a full subcategory and $S_\F:=S\cap\Mor\F$. (\blue{footnote}\footnote{Since $\F\subset\C$ is full, we have $\Mor\F:=\Mor_\F(\Ob\F,\Ob\F)=\Mor_\C(\Ob\F,\Ob\F)$.}). Then $\F[S_\F^{-1}]\subset\C[S^{-1}]$ is a \ul{full} \ul{subcategory} (\blue{footnote}\footnote{That is, the inclusion functor ~{\scriptsize $I:\F[S_\F^{-1}]\ra\C[S^{-1}],~\Big(U\sr{s^{-1}}{\ral}W\sr{f}{\ral}V\Big)~\mapsto~\Big(U\sr{s^{-1}}{\ral}W\sr{f}{\ral}V\Big)$} ~is both \ul{faithful} and \ul{full}.}) if one of the following is true:
\bit[leftmargin=0.9cm]
\item[i.] (a) $S_\F\subset\Mor\F$ is localizing, and (b) for every {\small $s\in S\cap \Mor_\C(\Ob\C,\Ob\F)$} there exists {\small $f\in \Mor_\C(\Ob\F,\Ob\C)$} such that ~{\small $sf\in S\cap \Mor_\C(\Ob\F,\Ob\F)\sr{\txt{fullness}}{=}S\cap \Mor_\F(\Ob\F,\Ob\F)=S_\F$}, ~i.e.
\bea
\adjustbox{scale=0.9}{\bt\Ob\F\ar[rr,dashed,"\exists~f"]\ar[rrrr,dashed,bend left,"sf\in S"] && \Ob\C \ar[rr,"\forall~s\in S"] && \Ob\F.\et}\nn
\eea
\item[ii.] (a) $S_\F\subset\Mor\F$ is localizing, and (b) for every {\small $s\in S\cap \Mor_\C(\Ob\F,\Ob\C)$} there exists {\small $f\in \Mor_\C(\Ob\C,\Ob\F)$} such that ~{\small $sf\in S\cap \Mor_\C(\Ob\F,\Ob\F)\sr{\txt{fullness}}{=}S\cap \Mor_\F(\Ob\F,\Ob\F)=S_\F$}, ~i.e.
\bea
\adjustbox{scale=0.9}{\bt\Ob\F\ar[from=rr,dashed,"\exists~f"']\ar[from=rrrr,dashed,bend right,"sf\in S"'] && \Ob\C \ar[from=rr,"\forall~s\in S"'] && \Ob\F.\et}\nn
\eea
\eit
\end{prp}
\begin{proof}
It suffices to verify i., which is essential in the left roof characterization (while noting that the verification of ii., which is essential in the right roof characterization, is similar or even symmetric/dual).

Since $\F\subset\C$ is full, we have $\Mor\F:=\Mor_\F(\Ob\F,\Ob\F)=\Mor_\C(\Ob\F,\Ob\F)$, and so
{\small\begin{align}
&\textstyle \Mor_{\F[S_\F^{-1}]}(U,V)={\Big\{\txt{left roofs}~(s,f)=fs^{-1}:U\sr{s^{-1}}{\ral}W\sr{f}{\ral}V\Big\}\over \sim_\F}={\bigcup_{W\in\Ob\F} \Mor_\F(W,U)\times \Mor_\F(W,V)\over\sim_\F},~~~~\txt{for all}~~U,V\in\Ob\F,\nn\\
&\textstyle~~~~\sr{\txt{fullness}}{=}{\bigcup_{W\in\Ob\F} \Mor_\C(W,U)\times \Mor_\C(W,V)\over\sim_\F}\sr{(1)}{\subset}{\bigcup_{W\in\Ob\F} \Mor_\C(W,U)\times \Mor_\C(W,V)\over\sim_\C}
\sr{(2)}{\subset}{\bigcup_{W\in\Ob\F} \Mor_\C(W,U)\times \Mor_\C(W,V)\over\sim_\F},\nn\\
&~~\Ra~~\Mor_{\F[S_\F^{-1}]}(U,V)\sr{(1),(2)}{=}\Mor_{\C[S^{-1}]}(U,V),\nn
\end{align}}where step (1) holds because $(s,f)\sim_\F(t,g)$ $\Ra$ $(s,f)\sim_\C(t,g)$, since $\Ob\F\subset\Ob\C$, i.e.,
{\small
\bea
(s,f)\sim_\F(t,g)~\Ra~(s,f)\sim_\C(t,g):
\adjustbox{scale=0.7}{%
\bt
  & & U'''\ar[dl,dashed,"r"']\ar[dr,dashed,"h"] & & \\
  & U'\ar[dl,"s"']\ar[drrr,"f"'pos=0.6] &  & U''\ar[dr,"g"]\ar[dlll,"t"pos=0.6] & \\
 U & & & & V
\et},~~~~U'''\in\Ob\F~~\Ra~~U'''\in\Ob\C,\nn
\eea
and step (2), or injectivity of the inclusion $(1)$, is due to the following:
\bea
(s,f)\sim_\C(t,g)~\Ra~(s,f)\sim_\F(t,g):
\adjustbox{scale=0.7}{%
\bt
  & U'''\ar[r,dashed,"\exists~f'''"] & X\ar[dl,dashed,"r"']\ar[dr,dashed,"h"] & (sr)f'''\in S & \\
  & U'\ar[dl,"s"']\ar[drrr,"f"'pos=0.6] &  & U''\ar[dr,"g"]\ar[dlll,"t"pos=0.6] & \\
 U & & & & V
\et}=
\adjustbox{scale=0.7}{%
\bt
  & & U'''\ar[dl,dashed,"rf'''"']\ar[dr,dashed,"hf'''"] & & \\
  & U'\ar[dl,"s"']\ar[drrr,"f"'pos=0.6] &  & U''\ar[dr,"g"]\ar[dlll,"t"pos=0.6] & \\
 U & & & & V
\et}
\nn
\eea}
The fullness of the subcategory ~$\F[S_\F^{-1}]\subset\C[S^{-1}]$ ~is clear from the equality $(1),(2)$.
\end{proof}

\begin{dfn}[\textcolor{blue}{
\index{Category of! left (right-bounded) complexes}{Category of left (right-bounded) complexes},
\index{Category of! right (left-bounded) complexes}{Category of right (left-bounded) complexes},
\index{Category of! bounded complexes}{Category of bounded complexes},
\index{Left! derived category}{Left derived category},
\index{Right! derived category}{Right derived category},
\index{Bounded! derived category}{Bounded derived category}}]
Let $\A$ be an abelian category. The \ul{category of left complexes} of $\A$ is the full subcategory $\A_{0-}^\Integer\subset\A_0^\Integer$ of all left (or right-bounded) complexes:
\bea
(C,d):\cdots\sr{d^{i_0-3}}{\ral}C^{i_0-2}\sr{d^{i_0-2}}{\ral}C^{i_0-1}\sr{d^{i_0-1}}{\ral}C^{i_0}\ra 0\ra 0\ra\cdots\nn
\eea
The \ul{category of right complexes} of $\A$ is the full subcategory $\A_{0+}^\Integer\subset\A_0^\Integer$ of all right (or left-bounded) complexes:
\bea
(C,d):\cdots\ra 0\ra 0\ra C^{i_0}\sr{d^{i_0}}{\ral}C^{i_0+1}\sr{d^{i_0+1}}{\ral}C^{i_0+2}\sr{d^{i_0+2}}{\ral}\cdots,\nn
\eea
The \ul{category of bounded complexes} of $\A$ is the full subcategory $\A_{0b}^\Integer:=\A_{0-}^\Integer\cap \A_{0+}^\Integer\subset\A_0^\Integer$ of all bounded complexes:
\bea
(C,d):\cdots\ra 0\ra 0\ra C^{i_0}\sr{d^{i_0}}{\ral}C^{i_0+1}\sr{d^{i_0+1}}{\ral}C^{i_0+2}\sr{d^{i_0+2}}{\ral}\cdots\sr{d^{j_0-2}}{\ral}C^{j_0-1}\sr{d^{j_0-1}}{\ral}C^{j_0}\ra 0\ra 0\ra\cdots\nn
\eea
The \ul{left derived category} of $\A$ is $\D^-\A:=\A_{0-}^\Integer[\{\txt{quasi-iso's}\}^{-1}]$, with localizing functor $L_\A^-:\A_{0-}^\Integer\ra\D^-\A$. The \ul{right derived category} of $\A$ is $\D^+\A:=\A_{0+}^\Integer[\{\txt{quasi-iso's}\}^{-1}]$, with localizing functor $L_\A^+:\A_{0+}^\Integer\ra\D^+\A$. The \ul{bounded derived category} of $\A$ is $\D^b\A:=\A_{0b}^\Integer[\{\txt{quasi-iso's}\}^{-1}]$, with localizing functor $L_\A^b:\A_{0b}^\Integer\ra\D^b\A$.
\end{dfn}

\begin{dfn}[\textcolor{blue}{
\index{Category of! left-exact complexes}{Category of left-exact complexes},
\index{Category of! right-exact complexes}{Category of right-exact complexes},
\index{Category of! bounded-exact complexes}{Category of bounded-exact complexes},
\index{Left-exact derived category}{Left-exact derived category},
\index{Right-exact derived category}{Right-exact derived category},
\index{Bounded-exact derived category}{Bounded-exact derived category}}]
Let $\A$ be an abelian category. The \ul{category of left-exact complexes} of $\A$ is the full subcategory $LE(\A_{0}^\Integer)\subset\A_0^\Integer$ of all left-exact complexes, defined to be complexes $C:\cdots\ra C^{i-1}\ra C^i\ra C^{i+1}\ra\cdots$ in $\A_0^\Integer$ with
\bea
H_i(C)=0~~~~\txt{for all}~~i\leq i_0,~~~~\txt{for some}~~i_0\in\Integer.\nn
\eea
The \ul{category of right-exact complexes} of $\A$ is the full subcategory $RE(\A_{0}^\Integer)\subset\A_0^\Integer$ of all \ul{right-exact complexes}, defined to be complexes $C:\cdots\ra C^{i-1}\ra C^i\ra C^{i+1}\ra\cdots$ in $\A_0^\Integer$ with
\bea
H_i(C)=0~~~~\txt{for all}~~i\geq i_0,~~~~\txt{for some}~~i_0\in\Integer.\nn
\eea

The \ul{category of bounded-exact complexes} of $\A$ is the full subcategory $BE(\A_0^\Integer):=LE(\A_0^\Integer)\cap RE(\A_0^\Integer)\subset\A_0^\Integer$ of all \ul{bounded-exact complexes}, defined to be complexes $C:\cdots\ra C^{i-1}\ra C^i\ra C^{i+1}\ra\cdots$ in $\A_0^\Integer$ with
\bea
H_i(C)=0~~~~\txt{for all}~~i\leq i_0,~i\geq j_0,~~~~\txt{for some}~~i_0\leq j_0,~~i_0,j_0\in\Integer.\nn
\eea
The \ul{left-exact derived category} of $\A$ is $LE(\D\A):=LE(\A_0^\Integer)[\{\txt{quasi-iso's}\}^{-1}]$. The \ul{right-exact derived category} of $\A$ is $RE(\D\A):=RE(\A_0^\Integer)[\{\txt{quasi-iso's}\}^{-1}]$. The \ul{bounded-exact derived category} of $\A$ is $BE(\D\A):=BE(\A_0^\Integer)[\{\txt{quasi-iso's}\}^{-1}]$.
\end{dfn}

\begin{crl}
Let $\A$ be an abelian category. Then the following equalities hold:
\bea
\D^-\A=RE(\D\A),~~~~\D^+\A=LE(\D\A),~~~~\D^b\A=BE(\D\A).\nn
\eea
\end{crl}
\begin{proof}
Consider the following diagram (in which the equalities hold because every exact (sub)complex is by definition quasi-isomorphic to the zero (sub)complex):
\[
\adjustbox{scale=0.7}{%
\bt
    &   \D^-\A \ar[r,equal]   & {RE(\D\A)}\ar[dr] & \\
  \D^b\A\ar[ur]\ar[dr]\ar[r,equal]  &  {BE(\D\A)} \ar[ur]\ar[dr] &   & \D\A \\
    &  \D^+\A \ar[r,equal]  & {LE(\D\A)}\ar[ur] &
\et}\]
We will show that the arrows in the above diagram are fully faithful imbeddings by verifying the conditions of Proposition \ref{SubLocPrp} (with ``\emph{localizing}'' replaced by ``\emph{localizing up to homotopy}'') for the given subcategories.

Let $S_q:=\{\text{quasi-iso's}\}$, $S_q^\pm:=S_q\cap\Mor\A_{0\pm}^\Integer$, and $S_q^b:=S_q\cap\Mor\A_{0b}^\Integer$.
To show $S_q^\ld$ (for each $\ld=+,-,b$) is \emph{localizing up to homotopy} we only need to replace the existing morphisms (and domains/codomains) in \emph{homotopy-commutative} diagrams such as
\bea
\adjustbox{scale=0.9}{\bt
C' \ar[drr,draw=none,"f_\ld s'\simeq s_\ld f'"description]\ar[d,dashed,"s'"']\ar[rr,dashed,"f'"] && D'\ar[d,"s_\ld"] \\
C_\ld\ar[rr,"f_\ld"'] && D_\ld
\et}~~~~\txt{and}~~~~
\adjustbox{scale=0.9}{\bt
 C_\ld\ar[rrr,bend left,"sf_\ld\simeq sg_\ld"]\ar[rr,shift left=1,"f_\ld"]\ar[rr,shift right=1,"g_\ld"'] && D_\ld\ar[r,dashed,"s"] & D'
\et}~\iff~
\adjustbox{scale=0.9}{\bt
  C'\ar[r,dashed,"s'"]\ar[rrr,bend left,"f_\ld s'\simeq g_\ld s'"] & C_\ld \ar[rr,shift left=1,"f_\ld"]\ar[rr,shift right=1,"g_\ld"'] && D_\ld
\et}\nn
\eea
with their appropriate ``\emph{systematic-restrictions}'' to morphisms (hence homotopy-commutative diagrams) in the underlying category $\A_{0\ld}^\Integer\subset \A_0^\Integer$. Observe that such ``restrictions'' can be obtained simply by replacing appropriate structure morphisms/differentials (within chain complexes as $\Integer$-systems in $\A$) and their domains/codomains with $0$, since such an operation (being a componentwise/degreewise operation) does not spoil commutativity up to homotopy, even when the given morphisms $s_\ld,f_\ld,g_\ld$ are assumed to be arbitrary members of $\Mor\A_0^\Integer$.

It remains to verify the property
\bea
\adjustbox{scale=0.9}{\bt
\Ob\A_{0\ld}^\Integer\ar[rr,dashed,"\exists~f"]\ar[rrrr,dashed,bend left,"sf\in S_q"] && \Ob\A_0^\Integer \ar[rr,"\forall~s\in S_q"] && \Ob\A_{0\ld}^\Integer,
\et}\nn
\eea
which indeed holds because given a quasi-isomorphism $C\sr{s}{\ral}D_\ld$, where $D_\ld\in\Ob\A_{0\ld}^\Integer$, we have
\bea
\adjustbox{scale=0.9}{\bt
C_\ld\ar[rr,hook,dashed,"\exists~i"]\ar[rrrr,dashed,bend left,"si\in S_q"] && C \ar[rr,"\forall~s\in S_q"] && D_\ld,
\et}\nn
\eea
where $C_\ld:=C|_{\A_{0\ld}^\Integer}$ is the ``restriction'' of $C$ obtained by replacing appropriate differentials in $C$ and their domains/codomains with $0$ (i.e., the ``restriction'' of a quasi-isomorphism is a quasi-isomorphism).
\end{proof} 

%% file: parts/AlgebraCat/AlgebraCatS14.tex
\section{Triangulation of Derived Categories}
Let $\T=(\T,T)$ be a triangulated category with translation functor (i.e., a fixed additive autofunctor)
\bea
T:\T\ra\T,~X\sr{f}{\ral}Y~\mapsto~X[1]\sr{f[1]}{\ral}Y[1].\nn
\eea
If $S\subset\Mor\T$ is a localizing class, we define \ul{($T$-induced) translation} in $\T[S^{-1}]$ by
\bea
T:\T[S^{-1}]\ra \T[S^{-1}],~~(s,f):X\sr{s}{\lal}Z\sr{f}{\ral}Y~~\mapsto~~\big(s[1],f[1]\big):X[1]\sr{s[1]}{\lal}Z[1]\sr{f[1]}{\ral}Y[1].\nn
\eea

\begin{note}
Observe that if the localizing class $S$ satisfies ``$s\in S\iff s[1]\in S$'', then the translation $T:\T[S^{-1}]\ra \T[S^{-1}]$ above automatically \ul{preserves equivalence/composition of roofs} (recall that $T$, being a functor $\T\ra\T$, preserves commutativity of diagrams in $\T$).
\end{note}

\begin{dfn}[\textcolor{blue}{\index{Triangulation-compatible localizing class}{Triangulation-compatible localizing class}}]
Let $\T$ be a triangulated category. A localizing class $S\subset\Mor\T$ is \ul{triangulation-compatible} (or \ul{compatible with the triangulation} on $\T$) if it has the following properties:
\bit[leftmargin=0.8cm]
\item[(i)] Translation-compatibility: $s\in S$ $\iff$ $s[1]\in S$ (i.e., $s\in S$ $\Ra$ $s[1],s[-1]\in S$).
\item[(ii)] TR3-compatibility: Given DT's $(u,v,w)$, $(u',v',w')$ and morphisms $s,t\in S$ such that the \ul{left square} in the following diagram \ul{commutes}
\bc\adjustbox{scale=0.9}{\bt
 X\ar[d,"s"]\ar[r,"u"] & Y\ar[d,"t"]\ar[r,"v"] & Z\ar[d,dashed,"r"]\ar[r,"w"] & X[-1]\ar[d,"{s[-1]}"]  \\
 X'\ar[r,"u'"] & Y'\ar[r,"v'"] & Z'\ar[r,"w'"] & X'[-1]
\et}\ec
there exists a morphism $r\in S$ completing the diagram to a \ul{morphism of triangles}.
\eit
\end{dfn}

\begin{thm}[\textcolor{blue}{\index{Triangulation of a localization}{Triangulation of a localization}: \cite[Theorem 2, p.251]{gelfand-manin2010}}]\label{TDCthm}
Let $\T=(\T,T)$ be a triangulated category, $S\subset\Mor\T$ a localizing class, and ~$L:\T\ra\T[S^{-1}],~X\sr{f}{\ral}Y~\mapsto~X\sr{(1,f)}{\ral}Y$~ the localizing functor. Declare a triangle in $\T[S^{-1}]$,
\bea
Q:\bt
A\ar[r,"{(s,u)}"] & B\ar[r,"{(t,v)}"] & C\ar[r,"{(r,w)}"] & A[-1],
\et\nn
\eea
\ul{distinguished} if there exists a distinguished triangle
$Q_0:\bt[column sep=small]
X\ar[r,"{\al}"] & Y\ar[r,"\beta"] & Z\ar[r,"\gamma"] & X[-1]
\et$ in $\T$ such that
\bea
Q\cong L(Q_0):
\bt
X\ar[r,"{(1,\al)}"] & Y\ar[r,"{(1,\beta)}"] & Z\ar[r,"{(1,\gamma)}"] & X[-1].
\et\nn
\eea
If $S$ is triangulation-compatible, then the above make $\T[S^{-1}]=\big(\T[S^{-1}],T\big)$ a triangulated category.
\end{thm}
\begin{proof}
For convenience, all morphisms in $\T[S^{-1}]$ will be explicitly expressed in terms of morphisms $f,s_f$ of $\T$ as $(s_f,f)$, and any morphism $f$ of $\T$ will also be written as $(1,f)$ if necessary, i.e., $f=(1,f)$. Also, it is clear that a left roof $X\sr{(s,f)}{\ral}Y$ or
\adjustbox{scale=0.8}{%
\bt[column sep=small] X\ar[from=r,"s"'] & Z\ar[r,"f"] & Y\et
} (as a diagram in $\T$) is equivalent to a commutative triangle in $\T[S^{-1}]$ of the form
\adjustbox{scale=0.7}{%
\bt
& Z\ar[dl,"{(1,s)}"']\ar[dr,"{(1,f)}"] &  \\
X\ar[rr,"{(s,f)}"] &  & Y
\et
} or
\adjustbox{scale=0.7}{%
\bt
& Z\ar[dl,"s"']\ar[dr,"f"] &  \\
X\ar[rr,"{(s,f)}"] &  & Y
\et
}
 or
\adjustbox{scale=0.7}{%
\bt
& Z\ar[dl,"s"']\ar[dr,"f"] &  \\
X\ar[rr,"fs^{-1}"] &  & Y
\et
}, since we have the \ul{factorization} $(s,f)=(1,f)(s,1)$, where $(s,1)$ is invertible, with $(s,1)(1,s)= (1,1)\sim(s,s)=(1,s)(s,1)$, i.e., $(s,f)=(1,f)(s,1)=(1,f)(1,s)^{-1}=fs^{-1}$.

Since $\T$ is additive and $S$ is localizing, we know from a previous result that $\T[S^{-1}]$ is additive. It therefore remains to verify the triangulation axioms.

\bit[leftmargin=1.2cm]
\item[(TR1)] \ul{Construction axiom}: DT's arise from morphisms and isomorphic images of other DT's:
   \bit[leftmargin=0.1cm]
   \item[(a)] $X\sr{(1,1)}{\ral}X\ral 0\ral X[-1]$ ~equals ~$L\big(X\sr{1}{\ral}X\ral 0\ral X[-1]\big)$, ~and so is a DT.
   \item[(b)] It is clear that \emph{any triangle isomorphic to a DT is also a DT}, by the definition of a DT in $\T[S^{-1}]$ (since a composition of isomorphisms is an isomorphism).

   \item[(c)]Given any morphism $X\sr{(s,u)}{\ral}Y$, i.e., $X\sr{s}{\lal}Z\sr{u}{\ral}Y$, in $\T[S^{-1}]$ we can complete $Z\sr{u}{\ral}Y$ in $\T$ into a DT (the middle part below)
       \bea
       \left(X\sr{s}{\lal}Z\right)~~Z\sr{u}{\ral}Y\sr{\al^u}{\ral}E(u)\sr{\beta^u}{\ral}Z[-1]~~\left(Z[-1]\sr{s[-1]}{\ral}X[-1]\right),\nn
       \eea
       and so (through the obvious appending) get the following associated triangle in $\T[S^{-1}]$;
    \[\adjustbox{scale=0.9}{\bt
       X\ar[r,"{(s,u)}"] & Y\ar[r,"{(1,\al^u)}"] & E(u)\ar[rr,"{(1,s[-1]\beta^u)}"] && X[-1]\\
       Z\ar[u,"{(1,s)}"']\ar[r,"{(1,u)}"] & Y\ar[u,equal]\ar[r,"{(1,\al^u)}"] & E(u)\ar[u,equal]\ar[rr,"{(1,\beta^u)}"] && Z[-1]\ar[u,"{(1,s[-1])}"']
       \et}~~~~\txt{(An isomorphism of triangles in $\T[S^{-1}]$.)}
       \]
Hence, the morphism $X\sr{(s,u)}{\ral}Y$ can be completed into a DT in $\T[S^{-1}]$, namely, the top row above in our conventional form for such a completion:
$\adjustbox{scale=0.8}{\bt X\ar[r,"{(s,u)}"] & Y\ar[r,"{\al^{(s,u)}}"] & E((s,u))\ar[r,"{\beta^{(s,u)}}"]& X[-1]\et}$,~ where $\al^{(s,u)}:=(s_\al,f_\al):=(1,\al^u)$, $E((s,u)):=E(u)$, and $\beta^{(s,u)}:=(s_\beta,f_\beta):=(1,s[-1]\beta^u)$.
   \eit

\item[(TR2)] \ul{Shift symmetry axiom}: This is clear by the definition of both translation and DT's in $\T[S^{-1}]$, which are in terms of the \ul{functors} $T:\T\ra\T$ and $L:\T\ra\T[S^{-1}]$ respectively.

\item[(TR3)] \ul{Stability axiom}: It is enough to consider only DT's in $\T[S^{-1}]$ that are images of DT's from $\T$. Let $X\sr{u}{\ral}Y\sr{v}{\ral}Z\sr{w}{\ral}X[-1]$ and $X'\sr{u'}{\ral}Y'\sr{v'}{\ral}Z'\sr{w'}{\ral}X'[-1]$ be DT's in $\T$. In $\T[S^{-1}]$, given two morphisms $(s,f),(t,g)$ as in the following diagram with a \ul{commuting left square},
\bea
\adjustbox{scale=0.8}{%
\label{TDCthmEq0}\bt
 X\ar[dd,"{(s,f)}"]\ar[rr,"{(1,u)}"] && Y\ar[dd,"{(t,g)}"]\ar[rr,"{(1,v)}"] && Z\ar[dd,dashed,"{(\ld,h)}"]\ar[rr,"{(1,w)}"] && X[-1]\ar[dd,"{(s[-1],f[-1])}"]  \\
  && && && && \\
 X'\ar[rr,"{(1,u')}"] && Y'\ar[rr,"{(1,v')}"] && Z'\ar[rr,"{(1,w')}"] && X'[-1]
\et}\eea
we need to construct a morphism $(\ld,h)$ completing the diagram to a \ul{morphism of triangles}. The above commutative diagram (with $(\ld,h)$ yet to be constructed) can be \ul{factored} into the following diagram (with $\T$-morphisms in $\T[S^{-1}]$ simplified) in which \ul{the left pyramid commutes}:
\bea
\label{TDCthmEq1}
\adjustbox{scale=0.8}{%
\bt
 & X\ar[dd,near start,"{(s,f)}"]\ar[rr,"u"] && Y\ar[dd,near start,"{(t,g)}"]\ar[rr,"v"] && Z\ar[dd,dashed,near start,"{(\ld,h)}"]\ar[rr,"w"] && X[-1]\ar[dd,near start,"{(s[-1],f[-1])}"]  \\
X''\ar[ur,"s"]\ar[dr,"f"]\ar[rr,dotted,near start,"u''"] & & Y''\ar[ur,"t"]\ar[dr,"g"]\ar[rr,dotted,near start,"v''"] & &Z''\ar[ur,dashed,"\ld"]\ar[dr,dashed,"h"]\ar[rr,dotted,near start,"w''"]& &X''[-1]\ar[ur,"{s[-1]}"]\ar[dr,"{f[-1]}"]& && \\
 & X'\ar[rr,"u'"] && Y'\ar[rr,"v'"] && Z'\ar[rr,"w'"] && X'[-1]
\et}\eea
where the diagram is obtained (beginning from the left) as follows:
   \bit[leftmargin=0.1cm]
    \item[(i)] First, with a \ul{starting} representative $(s_o,f_o)\sim(s,f)$ for the morphism $(s,f)$, the $\T$-morphism $u''$ is easily seen to exist (\blue{footnote}\footnote{The morphism $u''$ exists up to equivalence wrt the left roof $(s,f)$ in the sense that to obtain $u''$, we need to choose an appropriate representative for the equivalence classes of the left roof $(s,f)$.}) using the extension properties of $S$ applied to $X_o''\sr{us_o}{\ral}Y\sr{t}{\lal}Y''$,
\bea
\label{TDCthmEq2}
\adjustbox{scale=0.9}{%
\bt
 && &   & X\ar[ddd,bend left=10,"{(s_o,f_o)}"pos=0.65]\ar[from=dl,"s_o"]\ar[rr,"u"]&  & Y\ar[ddd,"{(t,g)}"] \\
 && &  X_o''\ar[ddr,near end,bend left=30,"f_o"']\ar[urrr,"us_o"description] &  & Y''\ar[ddr,"g"]\ar[ur,"t"']  &  \\
X''\ar[rr,dotted,"r\in S"] && X_1''\ar[drr,bend right=10,"f_1:=f_ot'"']\ar[uurr,bend left=40,"s_1:=s_ot'"]\ar[ur,dotted,"t'"]\ar[urrr,dotted,"{u_1'':=(us_o)'}"'pos=0.1] & && & \\
 && &  & X'\ar[rr,"u'"] & & Y'
\et}\eea
 making the top square $(X,Y,X_o'',Y'')$ commute in $\T$, and the bottom square $(X',Y',X_o'',Y'')$ commute in $\T[S^{-1}]$ as follows:
\bea
&&u'(s_o,f_o)=(t,g)u~~\Ra~~u'f_os_o^{-1}=gt^{-1}u=gt^{-1}us_os_o^{-1}=g(us_o)'t'^{-1}s_o^{-1},\nn\\
&&~~\Ra~~u'f_ot'=g(us_o)',~~\Ra~~u'f_1=gu_1'':X_1''\ra Y'~~~~\txt{in}~~\T[S^{-1}],\nn
\eea
but not yet in $\T$. To make the latter square commute in $\T$, observe that the equality $u'f_1=gu_1''$ in $\T[S^{-1}]$ is equality in $\T$ up to equivalence of left roofs via a commutative diagram:
\bea
\label{TDCthmEq3}
\adjustbox{scale=0.7}{%
\bt
& & X''\ar[dl,dashed,"r"']\ar[dr,dashed,"r\in S"] & & \\
& X_1''\ar[dl,"1"']\ar[drrr,"u'f_1"'pos=0.6] &  & X_1''\ar[dlll,"1"pos=0.6]\ar[dr,"gu_1''"] & \\
 X_1'' & & & & Y'
\et} \sim
\adjustbox{scale=0.7}{%
\bt
 & X''\ar[dl,"r"']\ar[dr,"gu_1''r"] & \\
 X_1'' && Y'
\et}=
\adjustbox{scale=0.7}{%
\bt
 & X''\ar[dl,"r"']\ar[dr,"u'f_1r"] & \\
 X_1'' && Y'
\et}
\eea
That is, there exists $r\in S$ such that ~$u'fr=gu''r:X''\ra Y'$. We can now set
\bea
s:=s_1r=s_ot'r,~~~~f:=f_1r=f_ot'r,~~~~u'':=u_1''r=(us_o)'r.\nn
\eea

\item[(ii)] Now complete $X''\sr{u''}{\ral}Y''$ into a DT in $\T$ given by $(u'',v'',w'')$ in the diagram (\ref{TDCthmEq1}). By TR3 for $\T$, the morphism $h\in\Mor\T$ exists such that $(f,g,h)$ is a morphism of triangles in $\T$. Similarly, by TR3 for $\T$ and compatibility of $S$ with the triangulation in $\T$, the morphism $\ld\in S$ exists such that $(s,t,\ld)$ is a morphism of triangles in $\T$.
    \eit

\item[(TR4)] \ul{Octahedral axiom}: Consider morphisms\adjustbox{scale=0.7}{%
\bt X\ar[rrrr,bend left,"{(s_v,v)(s_u,u)}"]\ar[rr,"{(s_u,u)}"] && Y\ar[rr,"{(s_v,v)}"] && Z\et
}. Then $(s_v,v)(s_u,u)=(s_us_v',vu')=:(s_{vu'},vu')$, where $s_vu'=us_v'$ and $s_{vu'}:=s_us_v'$ are as in the commutative diagram
\bea
(s_v,v)\circ(s_u,u)~:=~
\adjustbox{scale=0.7}{%
\bt
  & & X''\ar[dl,dashed,"s_v'"']\ar[dr,dashed,"u'"] & & \\
  & X'\ar[dl,"s_u"']\ar[dr,"u"] &  & Y'\ar[dl,"s_v"']\ar[dr,"v"] & \\
 X & & Y & & Z
\et}~~~~\sim~~~~
\adjustbox{scale=0.7}{%
\bt
 & X''\ar[dl,"s_{vu'}:=s_us_v'"']\ar[dr,"vu'"] & \\
 X && Z
\et}\nn
\eea
Using completing DT's for these three morphisms, namely,
$\adjustbox{scale=0.8}{\bt X\ar[r,"{(s_u,u)}"] & Y\ar[r,"{\al^{(s_u,u)}}"] & E((s_u,u))\ar[r,"{\beta^{(s_u,u)}}"]& X[-1]\et}$, $\adjustbox{scale=0.8}{\bt Y\ar[r,"{(s_v,v)}"] & Z\ar[r,"{\al^{(s_v,v)}}"] & E((s_v,v))\ar[r,"{\beta^{(s_v,v)}}"]& Y[-1]\et}$, and $\adjustbox{scale=0.7}{\bt X\ar[rr,"{(s_{vu'},vu')}"] && Z\ar[rr,"{\al^{(s_{vu'},vu')}}"] && E((s_{vu'},vu'))\ar[rr,"{\beta^{(s_{vu'},vu')}}"]&& X[-1]\et}$, which have \ul{simplified forms} as in the proof of TR1(c) above, we need to construct a DT in $\T[S^{-1}]$ of the form:
\[\adjustbox{scale=0.9}{\bt
E((s_u,u))\ar[r,dashed,"{(s_f,f)}"] & E((s_v,v)(s_u,u))\ar[r,dashed,"{(s_g,g)}"] & E((s_v,v))\ar[rrr,"{\al^{(s_u,u)}[-1]\beta^{(s_v,v)}}"] &&& E((s_u,u))[-1].
\et}\]
TR1(c) gives completing DT's for $(s_u,u)$, $(s_v,v)(s_u,u)$, $(s_v,v)$ respectively as follows:
\bea
&&\adjustbox{scale=0.8}{\bt
    X\ar[r,"{(s_u,u)}"] & Y\ar[r,"\al^u"] & E(u)\ar[rr,"{s_u[-1]\beta^u}"] && X[-1]\\
    X'\ar[u,"s_u"']\ar[r,"u"] & Y\ar[u,equal]\ar[r,"\al^u"] & E(u)\ar[u,equal]\ar[rr,"\beta^u"] && X'[-1]\ar[u,"{s_u[-1]}"']
    \et}~~~~\txt{(an iso of triangles in $\T[S^{-1}]$)}\nn\\
&&\adjustbox{scale=0.8}{\bt
    X\ar[r,"{(s_{vu'},vu')}"] & Z\ar[r,"\al^{vu'}"] & E(vu')\ar[rr,"{s_{vu'}[-1]\beta^{vu'}}"] && X[-1]\\
    X''\ar[u,"{s_{vu'}}"']\ar[r,"vu'"] & Z\ar[u,equal]\ar[r,"\al^{vu'}"] & E(vu')\ar[u,equal]\ar[rr,"\beta^{vu'}"] && X''[-1]\ar[u,"{s_{vu'}[-1]}"']
    \et}~~~~\txt{(an iso of triangles in $\T[S^{-1}]$)}\nn\\
&&\adjustbox{scale=0.8}{\bt
    Y\ar[r,"{(s_v,v)}"] & Z\ar[r,"\al^v"] & E(v)\ar[rr,"{s_v[-1]\beta^v}"] && Y[-1]\\
    Y'\ar[u,"s_v"']\ar[r,"v"] & Z\ar[u,equal]\ar[r,"\al^v"] & E(v)\ar[u,equal]\ar[rr,"\beta^v"] && Y'[-1]\ar[u,"{s_v[-1]}"']
    \et}~~~~\txt{(an iso of triangles in $\T[S^{-1}]$)}\nn
\eea
These in turn give an isomorphism of commutative diagrams in $\T[S^{-1}]$ as follows:
\[\adjustbox{scale=0.7}{\bt
X'\ar[from=d,"s_v'"]\ar[r,"u"] & Y\ar[from=d,"s_v"]\ar[r,"\al^u"] & E(u)\ar[from=d,dashed,"s_f"]\ar[r,"\beta^u"] & X'[-1]\ar[from=d,"{s_v'[-1]}"] \\
X''\ar[d,equal]\ar[r,"u'"] & Y'\ar[d,"v"]\ar[r,"\al^{u'}"] & E(u')\ar[d,dashed,"f"]\ar[r,"\beta^{u'}"] & X''[-1]\ar[d,equal] \\
X''\ar[d,"u'"]\ar[r,"vu'"] & Z\ar[d,equal] \ar[r,"\al^{vu'}"] & E(vu')\ar[d,dashed,"g"]\ar[r,"\beta^{vu'}"] & X''[-1]\ar[d,"{u'[-1]}"] \\
Y'\ar[dd,"0"]\ar[r,"v"] & Z\ar[dd,"0"] \ar[r,"\al^v"] & E(v)\ar[dd,"{\al^u[-1]\beta^v}"]\ar[r,"\beta^v"] & Y'[-1]\ar[dd,"0"] \\
 & & & \\
X'[-1]\ar[r,"{u[-1]}"] & Y[-1] \ar[r,"{\al^u[-1]}"] & E(u)[-1]\ar[r,"{\beta^u[-1]}"] & X'[-2] \\
&&&\\
\ar[rrr,draw=none,"\txt{\large Setup in $\T$}"] &~&~&~
\et}~~\cong~~
\adjustbox{scale=0.7}{\bt
&& & && \\
&& & && \\
X\ar[d,equal]\ar[rr,"{(s_u,u)}"] && Y\ar[d,"{(s_v,v)}"]\ar[r,"\al^u"] & E(u)\ar[d,dashed,"{(s_f,f)}"]\ar[rr,"{s_u[-1]\beta^u}"] && X[-1]\ar[d,equal] \\
X\ar[d,"{(s_u,u)}"]\ar[rr,"{(s_{vu'},vu')}"] && Z\ar[d,equal] \ar[r,"\al^{vu'}"] & E(vu')\ar[d,dashed,"{(1,g)}"]\ar[rr,"{s_{vu'}[-1]\beta^{vu'}}"] && X[-1]\ar[d,"{(s_u,u)[-1]}"] \\
Y\ar[dd,"0"]\ar[rr,"{(s_v,v)}"] && Z\ar[dd,"0"] \ar[r,"\al^v"] & E(v)\ar[dd,"{\al^u[-1]\beta^v}"]\ar[rr,"{s_v[-1]\beta^v}"] && Y[-1]\ar[dd,"0"] \\
 && & && \\
X[-1]\ar[rr,"{(s_u,u)[-1]}"] && Y[-1] \ar[r,"{\al^u[-1]}"] & E(u)[-1]\ar[rr,"{s_u[-2]\beta^u[-1]}"] && X[-2] \\
 && & && \\
\ar[rrrrr,draw=none,"\txt{\large Setup in $\T[S^{-1}]$}"] &~&~&~&~&~
\et}\]
Omitted horizontal links in the following associated commutative diagram are equalities.

\[\adjustbox{scale=0.7}{\bt
                                           &                                                &                                                                  & X'[-1]\ar[ddddrrrrrr,dotted,"{s_u[-1]}"description]
&&&       &       &         &       \\
                                           &                                                &                                                                  &
&&&       &       &         &       \\
                                           &                                                & E(u)\ar[uur,"\beta^u"]                                           & X''[-1]\ar[dd,equal]\ar[uu,"{s_v'[-1]}"'] \ar[ddrrrrrr,dotted,"{s_{vu'}[-1]}"description]
&&&       &       &         &       \\
                                           &                                                &                                                                  &
&&&       &       &         &       \\
                                           & Y\ar[uur,"\al^u"]                              & E(u')\ar[uu,dashed,"s_f"]\ar[dd,dashed,"f"]\ar[uur,"\beta^{u'}"]\ar[ddrrrrrr,dashed,"s_f"description] & X''[-1]\ar[dd,"{u'[-1]}"]\ar[ddrrrrrr,dotted,"{s_{vu'}[-1]}"description]
&&&       &       &         & X[-1]\ar[dd,equal] \\
                                           &                                                &                                                                  &
&&&       &       &         &       \\
 X'\ar[uur,"u"]\ar[ddddrrrrrr,dotted,near end,"s_u"description]                            & Y'\ar[uu,"s_v"]\ar[dd,"v"']\ar[uur,"\al^{u'}"]\ar[ddrrrrrr,dotted,near end,"s_v"description] & E(vu')\ar[dd,dashed,"g"]\ar[uur,"\beta^{vu'}"description]  & Y'[-1]\ar[dd,"0"]\ar[ddrrrrrr,dotted,"{s_v[-1]}"description]
&&&       &       & E(u)\ar[dd,dashed,"{(s_f,f)}"description]\ar[uur,"{s_u[-1]\beta^u}"description]    & X[-1]\ar[dd,"{(s_u,u)[-1]}"] \\
                                           &                                                &                                                                  &
&&&       &       &         &       \\
 X''\ar[uu,"s_v'"]\ar[uur,"u'"]\ar[ddrrrrrr,dotted,near end,"s_{vu'}"description]  & Z\ar[dd,equal]\ar[uur,near end,"\al^{vu'}"description]  & E(v)\ar[dd,near end,"{\al^u[-1]\beta^v}"description]\ar[uur,"\beta^v"description] & X'[-2]\ar[ddrrrrrr,dotted,near start,"{s_u[-2]}"description]
&&&       & Y\ar[dd,near start,"{(s_v,v)}"description]\ar[uur,"\al^u"description]     & E(vu')\ar[dd,dashed,"{(s_g,g)}"description]\ar[uur,shift right,"{s_{vu'}[-1]\beta^{vu'}}"description]  & Y[-1]\ar[dd,"0"] \\
                                           &                                                &                                                                  &
&&&       &       &         &       \\
 X''\ar[uu,equal]\ar[dd,"u'"]\ar[uur,"vu'"]\ar[ddrrrrrr,dotted,"s_{vu'}"description]& Z\ar[dd,"0"]\ar[uur,"\al^v"description]   & E(u)[-1]\ar[uur,near end,"{\beta^u[-1]}"description]                                &
&&& X\ar[uur,near start,"{(s_u,u)}"description]     & Z\ar[dd,equal]\ar[uur,near end,"\al^{vu'}"description]     & E(v)\ar[dd,"{\al^u[-1]\beta^v}"description]\ar[uur,"{s_v[-1]\beta^v}"description]    & X[-2] \\
                                           &                                                &                                                                  &
&&&       &       &         &       \\
 Y'\ar[dd,"0"]\ar[uur,"v"]\ar[ddrrrrrr,dotted,"s_v"description]  & Y[-1]\ar[uur,"{\al^u[-1]}"'] &                                                                  &
&&& X\ar[uu,equal]\ar[dd,"{(s_u,u)}"description]\ar[uur,"{(s_{vu'},vu')}"description]     & Z\ar[dd,"0"]\ar[uur,"\al^v"]  & E(u)[-1]\ar[uur,"{s_u[-2]\beta^u[-1]}"']&       \\
                                           &                                                &                                                                  &
&&&       &       &         &       \\
 X'[-1]\ar[uur,"{u[-1]}"']\ar[ddrrrrrr,dotted,"{s_u[-1]}"description]                &                                                &                                                                  &
&&& Y\ar[dd,"0"']\ar[uur,"{(s_v,v)}"description]     & Y[-1]\ar[uur,"{\al^u[-1]}"description] &         &       \\
                                           &                                                &                                                                  &
&&&       &       &         &       \\
                                           &                                                &                                                                  &
&&& X[-1]\ar[uur,"{(s_u,u)[-1]}"description] &       &         &
\et}\]

Hence, we get induced morphisms $(s_f,f)$ and $(s_g,g)=(1,g)$ that give the desired DT as follows:
\[\adjustbox{scale=0.8}{\bt
E((s_u,u))\ar[r,dashed,"{(s_f,f)}"]\ar[d,equal,"\txt{construction}"] & E((s_{vu'},vu'))\ar[r,dashed,"{(s_g,g)}"]\ar[d,equal,"\txt{construction}"] & E((s_v,v))\ar[rrr,"{\al^{(s_u,u)}[-1]\beta^{(s_v,v)}}"]\ar[d,equal,"\txt{construction}"] &&& E((s_u,u))[-1]\ar[d,equal,"\txt{construction}"] \\
E(u)\ar[r,dashed,"{(s_f,f)}"] & E(vu')\ar[r,dashed,"{(1,g)}"] & E(v)\ar[rrr,"{(1,\al^u)[-1](1,\beta^v)}"] &&& E(u)[-1] \\
E(u')\ar[u,"{(1,s_f)}"]\ar[r,dashed,"{(1,f)}"] & E(vu')\ar[u,equal]\ar[r,dashed,"{(1,g)}"] & E(v)\ar[u,equal]\ar[rrr,"{(1,\al^u)[-1](1,\beta^v)}"] &&& E(u)[-1]\ar[u,equal]
\et} \qedhere
\]
\eit
\end{proof}

\begin{crl}
Let $\A:=R$-mod. The derived categories $\D\A$, $\D^-\A$, $\D^+\A$, $\D^b\A$ are triangulated.
\end{crl}
\begin{proof}
We only need to verify triangulation-compatibility of quasi-iso's $S=S_q$ in the homotopy category $\H\A$. (i) It is clear that $s\in S_q$ iff $s[1]\in S_q$ (\blue{footnote}\footnote{Recall that the components of $s[k]$ are given by $s[k]_n:=(-1)^ks_{n+k}$ for each $n\in\Integer$.}). (ii) The TR3-compatibility of $S_q$ follows from the five lemma applied to the $H_n$-image (\blue{footnote}\footnote{As usual, $H_n$ (or $H^n$) is the $n$th homology functor.}) of the diagram below (whose rows are DT's and the left square commutes), because in the morphism of triangles $(s,t,r)$ that exists due to TR3 for $\H\A$, if $s,t\in S_q$ then so is $r$.
\bc\adjustbox{scale=0.9}{\bt
 X\ar[d,"s"]\ar[r,"u"] & Y\ar[d,"t"]\ar[r,"v"] & Z\ar[d,dashed,"r"]\ar[r,"w"] & X[-1]\ar[d,"{s[-1]}"]  \\
 X'\ar[r,"u'"] & Y'\ar[r,"v'"] & Z'\ar[r,"w'"] & X'[-1]
\et}\ec
\end{proof}

\section{Derived Functors from Derived Categories}
Our discussion in this section will be mainly motivational. For a further reading on this section, see for example \cite[Section III.6, p.185]{gelfand-manin2010} and \cite[Section 10.5, p.390]{weibel1994}.

Given an abelian category $\A$, recall that the LES of homology was introduced as a functor $LES:SES(\A_0^\Integer)\ra \A_0^\Integer$, where $SES(\A_0^\Integer)\subset\A_0^\Integer$ is the full subcategory with short exact chain complexes as its objects. Moreover, in the homotopy category $\H\A$, the LES of homology (now based on the triangulation of $\H\A$ instead of some abelian structure as in the case of $\A_0^\Integer$) becomes a functor $LES:DT(\H\A)\ra\H\A$, where $DT(\H\A)\subset\H\A$ is the full subcategory with distinguished triangles as its objects.

\begin{rmk}[\textcolor{blue}{The derived functor as an extension of a functor to the derived categories}]\label{DerFunctRmk1}
Let $\A,\B$ be abelian categories. Consider a \ul{left exact} additive functor $F:\A\ra\B$. Consider a SES $0\ra A\sr{f}{\ral} B\sr{g}{\ral} C\ra 0$ in $\A$. We \ul{observe} that the obvious \ul{extension $F:\A_0^\Integer\ra\B_0^\Integer$ of $F$ to complexes is ``resolution-exact''} in the sense that it maps the split-exact SES of deleted projective resolutions $0\ra P_\ast^A\sr{f_\ast}{\ral} P_\ast^B\sr{g_\ast}{\ral} P_\ast^C\ra 0$ (equivalently, the associated distinguished triangle $P_\ast^A\sr{f_\ast}{\ral} P_\ast^B\sr{i}{\ral}C(f_\ast)\sr{p}{\ral}P_\ast^A[-1]$ in the homotopy category) to another SES of complexes $0\ra F(P_\ast^A)\sr{F(f_\ast)}{\ral} F(P_\ast^B)\sr{F(g_\ast)}{\ral} F(P_\ast^C)\ra 0$ (equivalently, another distinguished triangle
\bea
\adjustbox{scale=0.8}{%
\bt[row sep=tiny]
F(P_\ast^A)\ar[r,"{F(f_\ast)}"] &  F(P_\ast^B)\ar[r,"{F(i)}"] & F(C(f_\ast))\ar[d,equal]\ar[r,"{F(p)}"] & F(P_\ast^A[-1])\ar[d,equal]\\
 &  & C(F(f_\ast)) & F(P_\ast^A)[-1]
\et}\nn\eea
in the homotopy category). We can then obtain the LES of homologies
\bea\bt
\cdots\ar[r] & R^iF(A)\ar[r] & R^iF(B)\ar[r] & R^iF(C)\ar[r] & R^{i+1}F(A)\ar[r] & \cdots
\et\nn
\eea
(viewed as either obtained from the SES of complexes or, equivalently, from the distinguished triangle) as usual in terms of the right-derived functors of $F$, i.e.,
\bea
R^iF:\A\ra\B,~~A\sr{f}{\ral}A'~~\mapsto~~\bt R^iF(A)\ar[r,"{R^iF(f)}"]& R^iF(A')\et,~~~~\txt{with}~~~~R^iF(A):=H^i(F(P_\ast^A)).\nn
\eea
The above \ul{observation} shows that in general, we can obtain the right-derived functors of $F$ through an \ul{extension} of $F$ to a functor $RF:\H\A\ra\H\B$ (also to be called \ul{the right-derived functor of $F$}) that is \ul{exact on distinguished triangles}. Using $RF$, we can then define the usual right-derived functors of $F$ as
\bea
R^iF:\A\ra\B,~~A\sr{f}{\ral}A'~~\mapsto~~\bt R^iF(A)\ar[r,"{R^iF(f)}"]& R^iF(A')\et,~~~~\txt{with}~~~~R^iF(A):=H^i(RF(A_\ast[0])),\nn
\eea
where $A_\ast[j]:\cdots\ra 0\ra A\ra 0\ra\cdots$ is the complex with $A_i[j]:=A\delta_{ij}$. This is sensible/consistent because in the homotopy category $\H\A$, we have $A_\ast[0]\cong P_\ast^A$ (recall that an exact complex of projective modules of the form $\cdots\ra P_{i_0+2}\ra P_{i_0+1}\ra P_{i_0}\ra 0$ is contractible, i.e., homotopy equivalent to the zero complex).

A similar story as above goes for the \ul{left-derived} functor $LF:\H\A\ra\H\B$ of a given \ul{right-exact} additive functor $F:\A\ra\B$. Even more interestingly, the same reasoning can be applied to the derived category, i.e., given a  left-exact (resp. right-exact) additive functor $F:\A\ra\B$, the right-derived (resp. left-derived) functor of $F$ may be introduced as a functor $RF:\D\A\ra\D\B$ (resp. $LF:\D\A\ra\D\B$)  extending $F$ in such a way that it is (i) exact on distinguished triangles (\blue{footnote}\footnote{This is so that we can obtain, as before, the induced LES of homologies associated with any SES of complexes in $\A_0^\Integer\subset\D\A$ (which also include SES's in $\A\subset\A_0^\Integer\subset\D\A$ as a subcategory). We may also interpret this particular property of the derived functor to mean that ``\ul{the derived functor restores exactness of $F$}''.}) and it is also (ii) as efficient as possible in representing $F$ in the sense that it has a universal property (\blue{footnote}\footnote{This universal property for the derived functor is assigned in accordance with the universal properties of the derived categories as localizations.}).
\end{rmk}

\begin{rmk}[\textcolor{blue}{\index{Construction of derived functors}{Construction of derived functors}}]\label{DerFunctRmk2}
Let $\A,\B$ be abelian categories and $F:\A\ra\B$ a left-exact additive functor. One notices in the above discussion that in order to successfully obtain the desired right-derived functors $R^iF$, one needs to apply $F$ on (i) special complexes $\O_{0F}^\Integer\subset\A_0^\Integer$ involving (ii) special types of objects $O_F\subset\Ob\A$ (in this case, left-complexes of projective objects), and $F:\A\ra\B$, $F:\A_0^\Integer\ra\B_0^\Integer$, or $\H(F):\H\A\ra\H\B$ needs to be ``\ul{well adapted to}'' (i.e., behave in certain specific ways on) these complexes and constituent objects for the whole process to succeed.

We will thus think of the right-derived functor $RF:\D\A\ra\D\B$ as arising from the universal property of the localization in the following way: With $O_F\subset\Ob\A$ a class of objects, $\O_F\subset\A$ the full subcategory whose objects are $O_F$, $\H\O_F\subset\H\A$ the full subcategory whose objects are complexes $\O_{0F}^\Integer\subset\A_0^\Integer$ of the objects $O_F$, we get the following diagram (in which $S_{q,F}:=S_q\cap\Mor\H\O_F$) which commutes except for the fully-dashed part which only commutes up to right-equalization by the localizing functor $L_{\O_F}$ in the sense
\[
\D(F)\al L_{\O_F}=\beta L_{\O_F}~~~~(\iff~\D(F)L_\A=L_\B\H(F)~).
\]
\[\adjustbox{scale=0.7}{%
\bt
F(\O_F)\subset\B\ar[d] & \O_F\subset\A\ar[l,"F"']\ar[d] &&&&\\
 \H F(\O_F)\ar[d,hook] & \H\O_F\ar[l,"\H(F)"']\ar[d,hook]\ar[rrrr,"L_{\O_F}"] &&&& \H\O_F[S_{q,F}^{-1}]\ar[ddllll,dashed,bend right=10,"\exists !~\al"']\ar[ddlllll,dashed,bend left=45,"\exists !~\beta"']\\
 \H\B\ar[d,"L_\B"'] & \H\A\ar[l,"\H(F)"']\ar[d,"L_\A"'] &&&& \\
 \D\B & \D\A\ar[l,dashed,"\D(F)"'] &&&&
\et},~~
\adjustbox{scale=0.7}{\bt[row sep=tiny, ampersand replacement=\&]
\left.
   \begin{array}{l}
    RF:=\lim\limits_G\big\{G:\D\A\ra\D\B~|~G\al=\beta,~G|_{\im(\al L_{\O_F})}=\D(F)\big\}\\
    ~~~~=\lim\limits_G\big\{G:\D\A\ra\D\B~|~G\al=\beta,~GL_\A=\D(F)L_\A=L_\B\H(F)\big\}\\
~~~~=\lim\limits_G\big\{G:\D\A\ra\D\B~|~G\al=\beta,~GL_\A=L_\B\H(F)\big\}\\
~~~~=\lim\limits_G\big\{G:\D\A\ra\D\B~|~\eta_G:=(\al,\beta):L_\B\H(F)\ra GL_\A\big\},
   \end{array}
\right.
\et}
\]
where the existence of $\al$ is due to the fact that $L_\A:\H\O_F\ra\D\A$ maps quasi-iso's to iso's, and we have \ul{assumed} (for the existence of $\beta$) that $L_\B\circ\H(F):\H\O_F\ra\D\B$ maps quasi-iso's to iso's (a \ul{requirement} that places a restriction on the class of objects $O_F$, as well as on the existence of the right-derived functor $RF$).
\end{rmk}

Consequently, for a general construction of derived functors, we introduce the following definition of a class of objects that is adapted to a given half-exact additive functor.

\begin{dfn}[\textcolor{blue}{\index{Adapted class of objects}{Adapted class of objects} for a half-exact additive functor}]
Let $\A,\B$ be abelian categories. Let $F:\A\ra\B$ be a left-exact additive functor. A class of objects $O\subset\Ob\A$ is adapted to $F$ if the following hold (with $\O\subset\A$ the full subcategory whose objects are $O$):
\begin{enumerate}
\item \ul{$O$ is additive}: $A,B\in O$ $\Ra$ $A\oplus B\in O$.
\item \ul{$F$ is exact on $\O\subset\A$}: If $E^\ast\in\O_{0+}^\Integer$ is an exact right-complex, then so is $F(E^\ast)$.
\item \ul{$O$ is large enough}: Every object $A\in\Ob\A$ is a subobject of an object from $O$, i.e., there exists a monomorphism $m:A\hookrightarrow A'$ for some $A'\in O$.
\end{enumerate}

Similarly, if $F:\A\ra\B$ is a right-exact additive functor, a class of objects $O\subset\Ob\A$ is adapted to $F$ if the following hold (with $\O\subset\A$ the full subcategory whose objects are $O$):
\begin{enumerate}
\item \ul{$O$ is additive}: $A,B\in O$ $\Ra$ $A\oplus B\in O$.
\item \ul{$F$ is exact on $\O\subset\A$}: If $E_\ast\in\O_{0-}^\Integer$ is an exact left-complex, then so is $F(E_\ast)$.
\item \ul{$O$ is large enough}: Every object $A\in\Ob\A$ is a quotient object of an object from $O$, i.e., there exists an epimorphism $e:A'\twoheadrightarrow A$ for some $A'\in O$.
\end{enumerate}
\end{dfn}
Since our discussion here is mainly motivational, a reader who is interested in details on the role of adapted classes on the existence of derived functors can see \cite[Section III.6, p.185]{gelfand-manin2010}.

\begin{dfn}[\textcolor{blue}{\index{Exact! functor between triangulated categories}{Exact functor between triangulated categories}: \index{Morphism of! triangulated categories}{Morphism of triangulated categories}}]
Let $\T,\T'$ be triangulated categories. A functor $F:\T\ra\T'$ is \ul{exact} if it is (i) additive, (ii) commutes with translation, and (iii) maps distinguished triangles to distinguished triangles. (\blue{footnote}\footnote{The requirement for the exact functor to commute with translation can be relaxed if necessary.})
\end{dfn}

Let $\A,\B$ be abelian categories and $L_\A:\H\A\ra\D\A$, $L_\B:\H\B\ra\D\B$ their derived categories. If $F:\A\ra\B$ is an additive functor, then (by the universal property of the localization) as in the diagram
\bea
\adjustbox{scale=0.7}{%
\bt
\A\ar[d,"F"]\ar[r,hook] & \H\A\ar[d,"\H(F)"]\ar[rr,"L_\A"]  &\ar[from=d,dashed,"\eta_F"']& \D\A\ar[d,dashed,"\exists!~\D(F)"]     \\
\B\ar[r,hook] & \H\B\ar[rr,"L_\B"']    &~& \D\B
\et}\nn
\eea
we get induced functors $\H(F):\H\A\ra\H\B$ and $\D(F):\D\A\ra\D\B$ such that the diagram commutes, i.e.,
\bea
\D(F)\circ L_\A=L_\B\circ\H(F).\nn
\eea
This \ul{equality} equivalently expresses a \ul{natural isomorphism} (i.e., an isomorphism of functors)
\bea
\eta_F:L_\B\circ\H(F)\ra\D(F)\circ L_\A~~~~\txt{(also briefly expressed as ~$\eta_F:L_\B\ra L_\A$)}.\nn
\eea
Based on this setup and the motivational discussions in Remarks \ref{DerFunctRmk1} and \ref{DerFunctRmk2}, we introduce the following definition for derived functors.

\begin{dfn}[\textcolor{blue}{\index{Right-derived functor}{Right-derived functor} of a left-exact functor, \index{Left-derived functor}{Left-derived functor} of a right-exact functor}]
Let $\A,\B$ be abelian categories and $F:\A\ra\B$ a left-exact additive functor.

Let $\H^-\A:={\A_{0-}^\Integer\over\simeq}$ (resp. $\H^+\A:={\A_{0+}^\Integer\over\simeq}$) and $\H^-\B:={\B_{0-}^\Integer\over\simeq}$ (resp. $\H^+\B:={\B_{0+}^\Integer\over\simeq}$) be the homotopy categories associated with left-complexes (resp. right-complexes). Similarly, let $L_\A^-:\A_{0-}^\Integer\ra\D^-\A$ (resp. $L_\A^+:\A_{0+}^\Integer\ra\D^+\A$) and $L_\B^-:\B_{0-}^\Integer\ra\D^-\B$ (resp. $L_\B^+:\B_{0+}^\Integer\ra\D^+\B$) be the associated left-derived (resp. right-derived) categories. As usual, $\H^\pm(F):\H^\pm\A\ra\H^\pm\B$ will denote the extensions of $F$ to the homotopy categories.

The \ul{right-derived functor} of $F$ is a pair $(RF,\vep_F)$ consisting of an \ul{exact functor} ~$RF:\D^+\A\ra\D^+\B$ ~and a \ul{morphism of functors} ~$\vep_F:L^+_\B\circ\H^+(F)\ra RF\circ L^+_\A$~ (which we represent by the following diagram)
\bea
~~~~~~~~
\adjustbox{scale=0.7}{%
\bt
\A\ar[d,"F"]\ar[r,hook] & \H^+\A\ar[d,"\H^+(F)"]\ar[rr,"L^+_\A"]  &\ar[from=d,"\vep_F"']& \D^+\A\ar[d,"RF"]     \\
\B\ar[r,hook] & \H^+\B\ar[rr,"L^+_\B"']    &~& \D^+\B
\et}~~~~\txt{(not commutative)}\nn
\eea
with the following universal property:
\bit
\item \ul{Universal property}: For any pair $(G,\vep)$ consisting of an exact functor ~$G:\D^+\A\ra\D^+\B$ ~and a morphism of functors ~$\vep:L^+_\B\circ\H^+(F)\ra G\circ L^+_\A$
\bea
\adjustbox{scale=0.7}{%
\bt
\A\ar[d,"F"]\ar[r,hook] & \H^+\A\ar[d,"\H^+(F)"]\ar[rr,"L^+_\A"]  &\ar[from=d,"\vep"']& \D^+\A\ar[d,"G"]     \\
\B\ar[r,hook] & \H^+\B\ar[rr,"L^+_\B"']    &~& \D^+\B
\et}~~~~\txt{(not commutative)}\nn
\eea
there exists a unique morphism of functors $\eta:RF\ra G$ such that the following diagram commutes
\bea
\adjustbox{scale=0.7}{%
\bt
L^+_\B\circ\H^+(F)\ar[d,"\vep"]\ar[rr,"\vep_F"] && RF\circ L^+_\A\ar[dll,dashed,"\eta\circ L^+_\A"]\\
G\circ L^+_\A &&
\et}\nn
\eea
\eit

Given a right-exact additive functor $F:\A\ra\B$, its \ul{left-derived functor} $(LF,\xi_F)$ consists of an exact functor  $LF:\D^-\A\ra\D^-\B$, and a morphism of functors $\xi_F:L^-_\B\circ\H^-(F)\ra LF\circ L^-_\A$ satisfying a universal property similar (in an obvious way) to the one above.
\end{dfn}

\begin{dfn}[\textcolor{blue}{\index{Classical! left-derived functor}{Classical left-derived functor}, \index{Classical! right-derived functor}{Classical right-derived functor}}]
Let $\A,\B$ be abelian categories. If $F:\A\ra\B$ is a left-exact additive functor with right-derived functor $RF:\D\A\ra\D\B$, the $i$th \ul{classical right-derived functor} of $F$ is
\bea
R^iF:=H^i\circ RF:\D\A\ra\B,~~~~\txt{i.e.,}~~~~RF_{\txt{classical}}:=H^\ast\circ RF:\D\A\ra\D\B,\nn
\eea
where $H^i:=H^0\circ T^i$ is the $i$th homology functor (with $H^0$ the $0$th homology functor and $T:\D\A\ra\D\A$ the translation functor). Similarly, if $G:\A\ra\B$ is a right-exact additive functor with left-derived functor $LG:\D\A\ra\D\B$, the $i$th \ul{classical left-derived functor} of $F$ is
\bea
L^iG:=H^i\circ LG:\D\A\ra\B,~~~~\txt{i.e.,}~~~~LG_{\txt{classical}}:=H^\ast\circ LG:\D\A\ra\D\B.\nn
\eea
\end{dfn}

Recall that for every distinguished triangle $A^\ast\ra B^\ast\ra C^\ast\ra A^\ast[1]$ in $\D\A$, we have the LES in $\B$
\bea
\cdots\ra H^i(A^\ast)\ra H^i(B^\ast)\ra H^i(C^\ast)\ra H^{i+1}(A^\ast)\ra\cdots\nn
\eea
Since derived functors preserve distinguished triangles (and $T$), it follows that we also have LES's in $\B$
\bea
&&\cdots\ra R^iF(A^\ast)\ra R^iF(B^\ast)\ra R^iF(C^\ast)\ra R^{i+1}F(A^\ast)\ra\cdots\nn\\
&&\cdots\ra L^iG(A^\ast)\ra L^iG(B^\ast)\ra L^iG(C^\ast)\ra L^{i+1}G(A^\ast)\ra\cdots\nn
\eea

\begin{dfn}[\textcolor{blue}{\index{Classical! Ext functor}{Classical Ext functor}, \index{Classical! Tor functor}{Classical Tor functor}}]
Let $\A$ be an abelian category, $Y\in\Ob\A$, and $R\Mor_\A(-,Y):\D\A\ra\D(\txt{Ab})$ the right-derived functor of $\Mor_\A(-,Y):\A\ra Ab$. We define
\bea
\Ext_\A^i(-,Y):=R^i\Mor_\A(-,Y):=H^i\circ  R\Mor_\A(-,Y):\A\ra Ab.\nn
\eea
Similarly, if $\A=R$-mod, let $M=M_R$ be a right $R$-module, and $M\mathop{\otimes}\limits^L-:\D(R\txt{-mod})\ra \D(Ab)$ the left-derived functor of $M\otimes_R-:R\txt{-mod}\ra Ab$. We define
\bea
\Tor^R_{-i}(-,M):=M\mathop{\otimes}\limits^{L^i}-:=H^i\circ(M\mathop{\otimes}\limits^L-):R\txt{-mod}\ra Ab.\nn
\eea
\end{dfn}

\begin{questions}
Consider the situation in the above definition. (i) What is the adapted class of objects for $F:=\Mor_\A(-,Y)$? (ii) What is the adapted class of objects for $G:=M\otimes_R-$?
\end{questions}

\section{Special Chain Morphisms and Complexes. The Classical Ext Functor}
\begin{dfn}[\textcolor{blue}{Alternative definition of \index{Ext bifunctor}{Ext bifunctor}}]
Let $\A$ be an abelian category, $A\in\Ob\A$, $n\in\Integer$, and $A_\ast[n]\in\A_0^\Integer$ the complex with $A_i[n]:=A\delta_{in}$, i.e.,
\bea
A_\ast[n]~:~\cdots\ra0\ra\ub{A}_{A_n[n]}\ra 0\ra \cdots.\nn
\eea
We define a bifunctor ~$\Ext^n:\A^{op}\times\A\ra Ab$ ~by
\bea
\Ext_\A^n(X,Y):=\Mor_{\D\A}(X_\ast[0],Y_\ast[n]).\nn
\eea
\end{dfn}
As we will show in Corollary \ref{NullHtyCrit3}, this new definition of Ext is equivalent to the classical definition of Ext in terms of the homology functor.

\begin{dfn}[\textcolor{blue}{
\index{Category of! projective objects}{Category of projective objects},
\index{Category of! injective objects}{Category of injective objects},
\index{Category of! (left) complexes of projective objects}{Category of (left) complexes of projective objects},
\index{Category of! (right) complexes of injective objects}{Category of (right) complexes of injective objects}}]
Let $\A$ be an abelian category, $\P\A\subset\A$ (resp. $\I\A\subset\A$) the full subcategory formed by projective (resp. injective) objects of $\A$, and $\P\A_0^\Integer\subset\A_0^\Integer$ (resp. $\I\A_0^\Integer\subset\A_0^\Integer$) be the full subcategory whose objects are complexes of projective (resp. injective) objects of $\A$. As before, let $\P\A_{0-}^\Integer\subset\P\A_0^\Integer$ (resp. $\I\A_{0+}^\Integer\subset\I\A_0^\Integer$) be the full subcategory whose objects are left-complexes of projective (resp. right-complexes of injective) objects of $\A$.
\end{dfn}

\begin{lmm}[\textcolor{blue}{\cite[p.180]{gelfand-manin2010}}]\label{NullHtyCrit}
Let $\A$ be an abelian category, $E_\ast,E^\ast\in\A_0^\Integer$ exact complexes, $P_\ast\in\P\A_{0-}^\Integer$, and $I^\ast\in\I\A_{0+}^\Integer$. Then the following hold:
\bit[leftmargin=0.8cm]
\item[(i)] Any chain morphism $f:P_\ast\ra E_\ast$ is homotopic to $0$, i.e., $f\simeq 0$.
\item[(ii)] Any chain morphism $f:E^\ast\ra I^\ast$ is homotopic to $0$, i.e., $f\simeq 0$.
\eit
\end{lmm}
\begin{proof}~
\bit[leftmargin=0.7cm]
\item[(i)] In the following diagram, since $P_0$ is projective and $d_0^Ef_0=0$, the morphism $h_0$ exists such that $f_0=d_1^Eh_0$. Thus, $d_1^Ef_1=f_0d_1^P=d_1^Eh_0d_1^P$, and so $d_1^E(f_1-h_0d_1^P)=0$.
\bea\adjustbox{scale=0.7}{%
\bt
\cdots\ar[rrrr] &&&& P_2\ar[ddllll,dashed,near start,bend left=15,"h_2"']\ar[dd,"f_2~(~-h_1d_2^P~)"]\ar[rrrr,"d_2^P"] &&&& P_1\ar[ddllll,dashed,near start,bend left=15,"h_1"']\ar[dd,"f_1~(~-h_0d_1^P~)"]\ar[rrrr,"d_1^P"] &&&& P_0\ar[ddllll,dashed,near start,bend left=15,"h_0"']\ar[dd,"f_0"]\ar[r] & 0\ar[dd]\ar[r] & 0\ar[dd]\ar[r]&\cdots\\
  &&&& &&&& &&&& && && & \\
\cdots\ar[rrrr,near end,"d_3^E"] &&&& E_2\ar[rrrr,near end,"d_2^E"] &&&& E_1\ar[rrrr,near end,"d_1^E"] &&&& E_0\ar[r,"d_0^E"] & E_{-1}\ar[r,"d_{-1}^E"]& E_{-2}\ar[r]&\cdots
\et}\nn
\eea
Since $P_1$ is projective, the morphism $h_1$ exists such that $f_1-h_0d_1^P=d_2^Eh_1$, i.e.,
\bea
f_1=h_0d_1^P+d_2^Eh_1.\nn
\eea
Continuing this way, by induction, we obtain morphisms $h_n:P_n\ra E_{n+1}$ such that
\bea
f_n=h_{n-1}d_n^P+d_{n+1}^Eh_n,~~~~h_n=0~~\txt{for}~~n\leq -1.\nn
\eea

\item[(ii)] In the following diagram, since $I^0$ is injective and $f^0d_E^{-1}=0$, the morphism $h^1$ exists such that $f^0=h^1d_E^0$. Thus, $f^1d_E^0=d_I^0f^0=d_I^0h^1d_E^0$, and so $(f^1-d_I^0h^1)d_E^0=0$.
\bea
\adjustbox{scale=0.7}{%
\bt
\cdots\ar[r]& E^{-2}\ar[dd]\ar[r,"d_E^{-2}"]& E^{-1}\ar[dd]\ar[r,"d_E^{-1}"] & E^0\ar[dd,"f^0"]\ar[rrrr,"d^0_E"] &&&& E^1\ar[ddllll,dashed,near start,bend left=15,"h^1"']\ar[dd,"f^1~(~-d_I^0h^1~)"]\ar[rrrr,"d^1_E"] &&&& E^2\ar[ddllll,dashed,near start,bend left=15,"h^2"']\ar[dd,"f^2~(~-d_I^1h^2~)"]\ar[rrrr,"d^2_E"] &&&& \ar[ddllll,dashed,near start,bend left=20,"h^3"']\cdots\\
 && &&&& &&&& &&&& \\
\cdots\ar[r] & 0\ar[r]& 0\ar[r] & I^0\ar[rrrr,near end,"d^0_I"] &&&& I^1\ar[rrrr,near end,"d^1_I"] &&&& I^2\ar[rrrr,near end,"d^2_I"] &&&& \cdots
\et}\nn
\eea
Since $I^1$ is injective, the morphism $h^2$ exists such that $f^1-d_I^0h^1=h^2d_E^1$, i.e.,
\bea
f^1=d_I^0h^1+h^2d_E^1.\nn
\eea
Continuing this way, by induction, we obtain morphisms $h^n:E^n\ra I^{n-1}$ such that
\bea
f^n=d_I^{n-1}h^n+h^{n+1}d_E^n,~~~~h^n=0~~\txt{for}~~n\leq 0.\nn
\eea
\eit
\end{proof}

\begin{crl}[\textcolor{blue}{\cite[p.180]{gelfand-manin2010}}]\label{NullHtyCrit1}
Let $\A$ be an abelian category, $A_\ast\in\A_0^\Integer$ (or just $A_\ast\in\A_{0-}^\Integer$), $P_\ast\in\P\A_{0-}^\Integer$, $I^\ast\in\I\A_{0+}^\Integer$, and $A^\ast\in\A_0^\Integer$ (or just $A^\ast\in\A_{0+}^\Integer$).
\bit[leftmargin=0.8cm]
\item[(i)] For any quasi-iso $A_\ast\sr{s}{\ral}P_\ast$, there exists a chain morphism $P_\ast\sr{t}{\ral}A_\ast$ such that $st\simeq 1_{P_\ast}$.
\item[(ii)] For any quasi-iso $I^\ast\sr{s}{\ral}A^\ast$, there exists a chain morphism $A^\ast\sr{t}{\ral}I^\ast$ such that $ts\simeq 1_{I^\ast}$.
\item[(iii)] Any quasi-iso $I^\ast\sr{s}{\ral}P_{-\ast}$ is a homotopy equivalence. (\blue{footnote}\footnote{By our conventions so far, lower indices such as on $P_\ast$ are decreasing (to the right) while upper indices such as on $I^\ast$ are increasing (to the right). Thus, the indices on $P_{-\ast}$ are increasing to the right, and so match those on $I^\ast$.})
\eit
\end{crl}
\begin{proof}~
\bit[leftmargin=0.7cm]
\item[(i)] In $\H\A$, consider the distinguished triangle $A_\ast\sr{s}{\ral}P_\ast\sr{i^s}{\ral}C_\ast(s)\sr{p^s}{\ral}A_\ast[-1]$, where $C_\ast(s):=P_\ast\oplus A_\ast[-1]$ with differential
~{\footnotesize $d_n^{C_\ast(s)}=
\left[
  \begin{array}{cc}
    d_n^P & s_{n-1} \\
    0 & -d_{n-1}^A\\
  \end{array}
\right]:
\left[
  \begin{array}{l}
    p_n \\
    a_{n-1} \\
  \end{array}
\right]\longmapsto
\left[
  \begin{array}{l}
    d_n^Pp_n + s_{n-1}a_{n-1} \\
    -d_{n-1}^Aa_{n-1} \\
  \end{array}
\right]$}
~is exact since $s$ is a quasi-iso. By Lemma \ref{NullHtyCrit}(i), $i^s\simeq 0$ via a homotopy {\footnotesize $h_n=[h'_n,h''_n]^T:P_n\ra C_{n+1}(s)=P_{n+1}\oplus A_n$} such that

\bea
&&i^s_n=h_{n-1}d_n^P+d_{n+1}^{C(s)}h_n,~~~~h_n=0~~\txt{for}~~n\leq -1,\nn\\
&&~~\iff~~
[id_{P_n},0]=[h'_{n-1},h''_{n-1}]^Td_n^P+
\left[
  \begin{array}{cc}
    d_{n+1}^P & s_n \\
    0 & -d_n^A\\
  \end{array}
\right][h'_n,h''_n]^T,\nn\\
&&~~\iff~~id_{P_n}=h'_{n-1}d_n^P+ d_{n+1}^Ph'_n + s_nh''_n~~~~\txt{and}~~~~0=h''_{n-1}d_n^P-d_n^Ah''_n.\nn
\eea
Hence, with $t:=h'':P_\ast\ra A_\ast$, we see that $st\simeq 1_{P_\ast}$ via the homotopy $h'_n:P_n\ra P_{n+1}$.

\item[(ii)] In $\H\A$, consider the distinguished triangle $I^\ast\sr{s}{\ral}A^\ast\sr{i_s}{\ral}C^\ast(s)\sr{p_s}{\ral}I^\ast[1]$, where $C^\ast(s):=A^\ast\oplus I^\ast[1]$ with differential
~{\footnotesize $d^n_{C^\ast(s)}=
\left[
  \begin{array}{cc}
    d^n_A & s^{n+1} \\
    0 & -d^{n+1}_I\\
  \end{array}
\right]:
\left[
  \begin{array}{l}
    a^n \\
    i^{n+1} \\
  \end{array}
\right]\longmapsto
\left[
  \begin{array}{l}
    d^n_Aa^n + s^{n+1}i^{n+1} \\
    -d^{n+1}_Ii^{n+1} \\
  \end{array}
\right]$}
~is exact since $s$ is a quasi-iso. By Lemma \ref{NullHtyCrit}, $p_s\simeq 0$ via a homotopy $h^n=[h_1^n,h_2^n]:C^n(s)=A^n\oplus I^{n+1}\ra I^{n-1}[1]=I^n$ such that
\bea
&&p_s^n=d_{I[1]}^{n-1}h^n+h^{n+1}d_{C^\ast(s)}^n,~~~~h^n=0~~\txt{for}~~n\leq 0,\nn\\
&&~~\iff~~
[0,id_{I^n[1]}]=d_{I[1]}^{n-1}[h_1^n,h_2^n]+[h_1^{n+1},h_2^{n+1}]
\left[
  \begin{array}{cc}
    d^n_A & s^{n+1} \\
    0 & -d^{n+1}_I\\
  \end{array}
\right],\nn\\
&&~~\iff~~0=d_{I[1]}^{n-1}h_1^n+h_1^{n+1}d_A^n~~~~\txt{and}~~~~id_{I^n[1]}=d_{I[1]}^{n-1}h_2^n+h_1^{n+1}s^{n+1}-h_2^{n+1}d^{n+1}_I.\nn
\eea
Hence, with $t:=h_1:A^\ast\ra I^\ast$, we see that $ts\simeq 1_{I^\ast}$ via the homotopy $h_2^n:I^n\ra I^{n-1}$.

\item[(iii)] By (i) and (ii), we have $t,t':P_{-\ast}\ra I^\ast$ such that $st\simeq 1_{P_{-\ast}}$ and $t's\simeq 1_{I^\ast}$, and so $t\simeq (t's)t=t'(st)\simeq t'$ is a homotopy inverse of $s$.
\qedhere
\eit
\end{proof}

\begin{crl}\label{NullHtyCrit2}
Let $\A$ be an abelian category. Then for any $P_\ast\in\P\A_{0-}^\Integer$, $I^\ast\in\I\A_{0+}^\Integer$, and $X^\ast,Y_\ast\in\A_0^\Integer$,
\bea
\Mor_{\H\A}(P_\ast,Y_\ast)\cong \Mor_{\D\A}(P_\ast,Y_\ast)~~~~\txt{and}~~~~\Mor_{\H\A}(X^\ast,I^\ast)\cong \Mor_{\D\A}(X^\ast,I^\ast).\nn
\eea
\end{crl}
\begin{proof}
Consider the map ~$L:\Mor_{\H\A}(P_\ast,Y_\ast)\ra \Mor_{\D\A}(P_\ast,Y_\ast),~f\mapsto(1,f)$ ~induced by the localization. We need to show $L$ is bijective. Let $(s,f)\in \Mor_{\D\A}(P_\ast,Y_\ast)$. By Corollary \ref{NullHtyCrit1}(i), there exists a morphism $t:P_\ast\ra Z_\ast$ giving equivalences of left roofs as follows:
{\scriptsize\bea\bt
 & Z_\ast\ar[dl,"s"']\ar[dr,"f"] & \\
P_\ast & & Y_\ast
\et~~~~\sim~~~~
\bt
 & P_\ast\ar[d,dashed,"t"]& \\
 & Z_\ast\ar[dl,"s"']\ar[dr,"f"] & \\
P_\ast & & Y_\ast
\et~~~~\sim~~~~
\bt
 & P_\ast\ar[dl,"1"']\ar[dr,"ft"] & \\
P_\ast & & Y_\ast
\et\nn
\eea}
which shows $L$ is surjective. Since the localization functor is additive, to show $L$ is injective, it suffices to show $(1,f)\sim(1,0)$ implies $f=0$. This is indeed the case by the following commutative diagram:
\bea
\adjustbox{scale=0.7}{%
\bt
& & Z'_\ast\ar[dl,dashed,"s"']\ar[dr,dashed,"s"] & & \\
& P_\ast\ar[dl,"1"']\ar[drrr,"f"'pos=0.6] &  & P_\ast\ar[dlll,"1"pos=0.6]\ar[dr,"0"] & \\
 P_\ast & & & & Y_\ast
\et}~~\sim~~
\adjustbox{scale=0.7}{%
\bt
 & P_\ast\ar[dl,"s"']\ar[dr,"0=fs"] & \\
P_\ast & & Y_\ast
\et}~~\Ra~~fs=0~~\txt{for all}~~s,~~\Ra~~f=0.\nn
\eea

The isomorphism $\Mor_{\H\A}(X^\ast,I^\ast)\cong \Mor_{\D\A}(X^\ast,I^\ast)$ holds similarly, except that one uses left-right symmetry to conclude that the description of morphisms in $\D\A$ in terms of left roofs is equivalent to their description in terms of right roofs, and so (with Corollary \ref{NullHtyCrit1}(ii)) proving bijectivity of
\bea
L:\Mor_{\H\A}(X^\ast,I^\ast)\ra \Mor_{\D\A}(X^\ast,I^\ast),~f\mapsto(1,f)\nn
\eea
is equivalent to proving bijectivity of the map induced by the localization in terms of right roofs, i.e.,
\bea
L':\Mor_{\H\A}(X^\ast,I^\ast)\ra \Mor_{\D\A}(X^\ast,I^\ast),~f\mapsto(f,1).\nn
\eea
\end{proof}

\begin{crl}\label{NullHtyCrit3}
Let $\A$ be an abelian category with enough projectives. For any objects $X,Y\in\Ob\A$,
\bea
\textstyle\Ext_\A^i(X,Y)\cong H^i\big(\Mor_\A(P^X_\ast,Y)\big):={\ker \Mor_\A\big(d^{P_\ast^X}_{i+1}\big)\over\im \Mor_\A \big(d^{P_\ast^X}_i\big)},\nn
\eea
where $P_\ast^X$ is a deleted projective resolution of $X$.
\end{crl}
\begin{proof}
Given $X,Y\in\Ob\A$, consider a deleted projective resolution $P^X_\ast\sr{\simeq}{\ral} X_\ast[0]$ of $X$ (where $\simeq$ denotes the obvious hty equivalence in $\A_0^\Integer$, i.e., iso in $\H\A$), and display any chain morphism $f_\ast=(\cdots,0,f_i,0,\cdots)\in \Mor_{\H\A}(P_\ast^X,Y_\ast[i])$ and any homotopy $h_\ast:P_\ast^X\ra Y_{\ast+1}[i]$ as in the following diagram:
\bc
\adjustbox{scale=0.8}{%
\bt
 &P_\ast^X\ar[dl,"h_\ast"']\ar[d,"f_\ast"] &\cdots\ar[r] & P^X_{i+1}\ar[d,"0"]\ar[rr,"d_{i+1}"] &&P_i\ar[dll,bend left=15,near start,"h_i=0"']\ar[d,"f_i"]\ar[rr,"d_i"] && P^X_{i-1}\ar[dll,bend left=15,near start,"h_{i-1}"']\ar[d,"0"]\ar[r] &\cdots\ar[r] & P^X_0\ar[r] & 0 \\
Y_{\ast+1}[i] & Y_\ast[i]  & \cdots\ar[r] & 0\ar[rr,near end,"0"] && Y\ar[rr,near end,"0"]  && 0\ar[r]  & \cdots  & &
\et}\ec
Then, by Corollary \ref{NullHtyCrit2} at step (a) below, and using the above diagram at step (b) below, we get

{\small\bea
&&\textstyle Ext_\A^i(X,Y):= \Mor_{\D\A}(X_\ast[0],Y_\ast[i])= \Mor_{\D\A}(P^X_\ast,Y_\ast[i])\sr{(a)}{\cong} \Mor_{\H\A}(P^X_\ast,Y_\ast[i]):= {\Mor_{\A_0^\Integer}(P^X_\ast,Y_\ast[i])\over\simeq}\nn\\
&&\textstyle~~~~\sr{(b)}{=}{\left\{\txt{morphisms ~$P^X_i\sr{f_i}{\ral} Y$ ~such that the left square in the diagram commutes, i.e., such that ~$\Mor_\A(d_{i+1},Y)(f_i):=f_id_{i+1}=0$}\right\}\over \txt{$f_i=0$~ if ~$f_i\simeq 0$,~ i.e., if ~$f_i=h_{i-1}d_i+0h_i=\Mor_\A(d_i,Y)(h_{i-1})$}}\nn\\
&&\textstyle~~~~~~={\ker \Mor_\A(d_{i+1},Y)\over\im \Mor_\A(d_i,Y)}~=:~H^i\big(\Mor_\A(P^X_\ast,Y)\big).\nn
\eea}
\end{proof}
Similarly, if $\A$ is an abelian category with enough injectives, then for any objects $X,Y\in\Ob\A$,
\bea
\textstyle\Ext_\A^i(X,Y)\cong H^i\big(\Mor_\A(X,I^\ast_Y)\big):={\ker \Mor_\A\big(d_{I^\ast_Y}^i\big)\over\im \Mor_\A \big(d_{I^\ast_Y}^{i-1}\big)},\nn
\eea
where $I^\ast_Y$ is a deleted injective resolution of $Y$. The above two characterizations of $Ext$ give an alternative proof of the earlier encountered Ext-ext isomorphism (or balance)
\[
H^i\big(\Mor_\A(P^X_\ast,Y)\big)\cong H^i\big(\Mor_\A(X,I^\ast_Y)\big).\nn
\]

%% file: parts/AlgebraM/RepsTheoryI.tex
\chapter{Representation Theory I: Groups, Symmetry, and Object Representation}\label{RepsTheoryI}
For initial ideas about representation theory, see for example \cite{fulton-harris}.

Given a ring $R$, a (left) $R$-module is really a set, with an abelian group structure, on which the ring $R$ ``\emph{acts}'' (from the left). In other words, an R-module $M$ ``\emph{represents}'' (or carries a ``\emph{representation}'' of) $R$ through the ``\emph{action}'' of $R$ on $M$. Thus, the \emph{theory of modules} can be called the ``\emph{representation theory of rings}''.

Similarly, thinking of a group $G=(G,\cdot,e)=(G,\cdot,1)$ as a ring $R=(R,\cdot,+,0,1)$ in which (i) addition $+$ has been dropped, (ii) $0$ has been removed, and (iii) every other element is invertible, an analog of the $R$-module is called a $G$-set, accordingly defined as a set on which $G$ acts (where of course, unlike for the $R$-module, an abelian group structure is not required for the $G$-set). Therefore, the \emph{theory of $G$-sets} or the ``\emph{representation theory of groups}'' (which is easier in one way due to the invertibility of all group elements, while more difficult in another way due to the presence of only a single operation) is similar in spirit to the \emph{theory of modules}, i.e., the ``\emph{representation theory of rings}''.

The main purpose of this chapter (and others related to it) is to discuss some basic concepts of representation theory with a categorical perspective in mind. After using the relatively simple structure of groups to introduce required basic tools (such as the concepts of ``\emph{symmetry}'', ``\emph{action}'', and ``\emph{invariance}'' (\blue{footnote}\footnote{Although we will not explicitly use ``invariance'' terminology, it is important to note that a homomorphism of $G$-sets or $G$-homomorphism (still to be defined) is also called a ``$G$-invariant map''. Likewise, an $R$-homomorphism (i.e., an ``$R$-linear map'') is partly an ``$R$-invariant map'' (or ``$R$-scale invariant map''). For objects that are both $G$-sets and $R$-modules simultaneously (call them $RG$-modules), their homomorphisms are not only both $G$-invariant and $R$-linear but also $G$-linear.})) a systematic/categorical discussion of representation theory simultaneously involving groups, rings, modules, and algebras will be given in section \ref{SystRepsThry} (starting on page \pageref{SystRepsThry}).

An idea that is of central importance in representation theory is that of an \emph{irreducible} (or \emph{elementary}) object in a given category. Under very special/limited circumstances, irreducible objects are also called ``\emph{simple}'' or ``\emph{indecomposable}'' objects (usually named with respect/accordance to some preexisting process/procedure/operation, such as limits and colimits, in an underlying category).

\begin{dfn}[\blue{\index{Irreducible! (elementary) object in a category}{Irreducible (elementary) object in a category}}]
Let $\C$ be a category . An object $C\in\Ob\C$ is \ul{irreducible} (\ul{elementary}) if it is nontrivially-minimal with respect to inclusion $\subset$ in $\C$, i.e., (i) $C$ is not a special object such as an initial object and (ii) $C$ has no proper subobjects in the sense that if $C'\subset C$ is a subobject of $C$, then $C\cong C'$.
\end{dfn}

Our end goal is to describe a given object (of interest) in terms of irreducible subobjects. However, the object reduction process typically includes intermediate steps involving objects of intermediary complexity (e.g., ``\emph{cyclic}'' and ``\emph{finitely generated}'' objects). In general therefore, we will consider a stepwise object reduction approach to representation theory where the idea is to first express a given object in terms of intermediate objects (e.g., cyclic and finitely generated objects) and then, if necessary (\blue{footnote}\footnote{There exist very special/accidental situations where some cyclic objects are automatically irreducible.}), eventually in terms of irreducible objects.

Even though there is no explicit reference due to extensive reformulation, the next three sections (cyclic groups, finite products, symmetry groups) are partly inspired by \cite{gallian2013}.

\section{Cyclic Groups}
Let $\C\subset Sets$ be a subcategory of Sets. A \index{Cyclic! object}{\ul{cyclic object}} (resp. \index{Finitely generated! object}{\ul{finitely generated object}}) in $\C$ is an object $C=\langle c\rangle$ (resp. $C=\langle c_1,...,c_n\rangle$) that is generated in a precise way by a single element $c\in C$ (resp. finite number of elements $c_1,...,c_n\in C$). Cyclic (resp. finitely generated) groups, rings, modules, and algebras were defined in the first part of our notes.

In typical subcategory of Sets, cyclic objects are among the most basic/fundamental types of objects, in the sense that (i) they approximately resemble minimal/irreducible objects of the category, with respect to inclusion, (ii) every object of the subcategory is a relatively simple combination (e.g., union, sum, etc) of cyclic subobjects, (iii) many properties of a given object can be easily characterized in terms of those of its cyclic subobjects. Consequently, in any given subcategory of Sets (e.g., Groups, Rings, Modules, Algebras, etc) a thorough first understanding of cyclic objects of the category is often essential.

Recall that in a group $G$, (i) the \ul{order of $G$} is the cardinality $|G|$, and (ii) the \ul{order} $|g|$ of an element $g\in G$ is defined to be the order (cardinality) of the cyclic subgroup $\langle g\rangle\leq G$, i.e., $|g|:=|\langle g\rangle|$. Also, recall that for two elements $d,m$ in a commutative ring, we write $d|m$ to mean ``$d$ divides $m$'' (i.e., $d$ is a \ul{divisor} of $m$, or $m$ is a \ul{multiple} of $d$). Furthermore, we will use the following new piece of notation:
\begin{notation}
Let $R$ be a commutative ring and $d,m\in R$. If $d|m$ (i.e., $d$ divides $m$), then in $m=qd$ we will denote the element $q\in R$ by ${m\over d}$ (which should not be mistaken for an arbitrary element of the ring of quotients $Q(R):=\{{a\over b}:a,b\in R,b\neq 0\}$).
\end{notation}

\subsection{Subgroups of the additive group of integers}
\begin{thm}[\blue{\index{Minimal! generator}{Minimal generator} of a subgroup of the PID $\Integer$}]\label{spigt}
Let $S$ be a subgroup of $\Integer=(\Integer,+)$. Then either $S=\{0\}$, or $S=\Integer a$, where $a$ is the least positive integer in $S$.
\end{thm}
\begin{proof}
We know $0\in S$. If $S$ contains no other elements, then $S=\{0\}$. Next, assume $S\neq \{0\}$. Let $n\in S\backslash\{0\}$. Since $-n,n\in S$, and one of $n,-n$ must be positive, it follows $S$ contains a positive integer. Let $a$ be the smallest positive integer in $S$. It is clear that $\Integer a\subset S$, since $\Integer a$ is the smallest subgroup of $(\Integer,+)$ containing $a$. Let $s\in S$. By the division algorithm, $s=qa+r$ for integers $q,r$ with $0\leq r<a$. (\blue{footnote}\footnote{That is, $\Integer=\Integer a+\Integer\cap[0,a)$.}). Since $r=s-q a~\in~S$, we see that $r=0$, since $a$ is the smallest positive integer in $S$. Hence $s=qa$, and so $S\subset \Integer a$.
\end{proof}

\begin{lmm}
Let $R$ be a ring and $a_1,...,a_n\in R$. If $\sum_i a_i$ generates $R$, then so does $\{a_1,...,a_n\}$.
\end{lmm}
\begin{proof}
If $\sum_i a_i$ generates $R$, then $R=R\sum_i a_i\subset \sum_i Ra_i\subset R$, and so $R=\sum_i Ra_i$.
\end{proof}

\begin{dfn}[\blue{\index{Prime! number}{Prime number}}]
A \ul{prime number} (or just a \ul{prime}) $p\in\Integer$ is a positive prime integer (equivalently, an integer greater than $1$ whose only divisors are $1$ and itself).
\end{dfn}
\begin{lmm}[\blue{\index{Euclid lemma}{Euclid lemma}}]
Let {\small $p\in\Integer$} be a prime and {\small $a_1,...,a_n\in\Integer$}. If {\small $p|a_1a_2...a_n$}, then $p|a_i$ for some $i$ (i.e., if $p$ divides $a_1a_2\cdots a_n$, then $p$ divides some $a_i$).
\end{lmm}
\begin{proof}
Suppose for each $i$, $p$ does not divide $a_i$ (which implies $gcd(p,a_i)=1$, or $\Integer p+\Integer a_i=\Integer$, since $p$ is prime). Then $\Integer=\Integer p+\Integer a_1=\Integer p+(\Integer p+\Integer a_2)a_1=\Integer p+\Integer pa_1+\Integer a_1a_2=\Integer p+\Integer a_1a_2=\cdots=\Integer p+\Integer a_1a_2\cdots a_n$, i.e., $gcd(p,a_1a_2\cdots a_n)=1$, and so $p$ does not divide $a_1a_2\cdots a_n$.
\end{proof}

\begin{thm}[\blue{\index{Gcd-Lcm balance}{Gcd-Lcm balance}}]
For any positive integers $a,b\geq 1$, we have $\txt{gcd}(a,b)~\txt{lcm}(a,b)=ab$.
\end{thm}
\begin{proof}
Let $d:=gcd(a,b)$ and $l:=lcm(a,b)$. Observe that ${ab\over d}=a{b\over d}={a\over d}b$ is a common multiple of $a$ and $b$, and so ${ab\over d}\geq l$. (\blue{footnote}\footnote{The same steps show that for any positive integers $a_1,...,a_n\geq 1$, we have $gcd(a_1,...,a_n)lcm(a_1,...,a_n)\leq a_1a_2\cdots a_n$.}). On the other hand, ${ab\over l}$ is a common divisor of $a$ and $b$ since $a/\big({ab\over l}\big)={l\over b}$ and $b/\big({ab\over l}\big)={l\over a}$, and so ${ab\over l}\leq d$.
\end{proof}

\subsection{Properties of cyclic groups}~\\~
Recall that a group $G$ is cyclic if {\footnotesize $G=\langle g\rangle:=\{g^k:k\in\Integer\}$} for some $g\in G$. We have the group epimorphism
\bea
h_g:(\Integer,+)\ra\langle g\rangle,~k\mapsto g^k,\nn
\eea
where $\Integer=(\Integer,+)$ denotes the additive group of integers. In the ensuing discussion, we will see that $h_g$ is an isomorphism (i.e., $\langle g\rangle\cong(\Integer,+)$ ) $\iff$ $|g|=\infty$. Meanwhile, if $|g|<\infty$, then $\ker h_g=\{k\in\Integer:g^k=e\}=\Integer|g|$, in which case the first isomorphism gives $\langle g\rangle\cong(\Integer_{|g|},+):=\big({\Integer\over \Integer|g|},+\big)$.

\begin{thm}[\blue{Multiplication of Cyclic Subgroups}]
In any cyclic group $G=\langle a\rangle$, the subgroup $\langle a^m\rangle \langle a^n\rangle$ is generated by $a^{\txt{gcd}(m,n)}$, i.e., $\langle a^m\rangle \langle a^n\rangle=\langle a^{\txt{gcd}(m,n)}\rangle$.
\end{thm}
\begin{proof}
$\langle a^m\rangle\langle a^n\rangle=\{a^{im}:i\in\Integer\}\{a^{in}:i\in\Integer\}=\{a^{im}a^{jn}=a^{im+jn}:i,j\in\Integer\}=\{a^k:k\in \Integer m+\Integer n=\Integer~gcd(m,n)\}=\{a^{i~gcd(m,n)}:i\in\Integer\}=\langle a^{gcd(m,n)}\rangle$.
\end{proof}

\begin{thm}[\blue{Intersection of Cyclic Subgroups}]
In any cyclic group $G=\langle a\rangle$, the subgroup $\langle a^m\rangle\cap\langle a^n\rangle$ is generated by $a^{\txt{lcm}(m,n)}$, i.e., $\langle a^m\rangle\cap\langle a^n\rangle=\langle a^{\txt{lcm}(m,n)}\rangle$.
\end{thm}
\begin{proof}
Notice $x\in \langle a^m\rangle\cap\langle a^n\rangle$ if and only if $x=a^k$ where $k=im=jn$ for some $i,j\in\Integer$. Therefore,
\bea
\langle a^m\rangle\cap\langle a^n\rangle=\{a^k:k\in\Integer m\cap\Integer n\}=\{a^k:k\in\Integer~\txt{lcm}(m,n)\}=\{a^{i~\txt{lcm(m,n)}}:i\in\Integer\}=\langle a^{\txt{lcm}(m,n)}\rangle.\nn
\eea
\end{proof}

\begin{thm}[\blue{Criterion for $a^i=a^j$}]
Let $G$ be a group and $g\in G$. Then for any $i,j\in\Integer$,
\[g^i=g^j\iff
\left\{
  \begin{array}{ll}
    i=j, & \txt{if ~$|g|=\infty$}, \\
    |g|~~\txt{divides}~~i-j, & \txt{if ~$|g|<\infty$.}
  \end{array}
\right.
\]
\end{thm}
\begin{proof}
Assume wlog that $i\geq j$. ($\Ra$): Observe that $g^i=g^j~~\iff~~g^{i-j}=e$. If $|g|=\infty$, then $i=j$, otherwise if $i>j$, then $g^{i-j}=e$ implies
$\{e\}\neq\langle g^{i-j}\rangle=\langle e\rangle=\{e\}$ (a contradiction). If $|g|<\infty$, then $\langle g\rangle=\{e,g,g^2,\cdots,g^{|g|-1}\}$, and so by the division algorithm for $\Integer$ (a Euclidean domain, in which $i-j=q|g|+r$ for a quotient $q$ and remainder $r$ satisfying $0\leq r<|g|$ ), $g^{i-j}=e$ implies $g^r=g^{i-j}=e$, which implies $r=0$, i.e., $|g|$ divides $i-j$.

($\La$): Conversely, if $i=j$ or $|g|$ divides $i-j$, it is clear that $g^i=g^j$.
\end{proof}

\begin{crl}
Let $G$ be a group and $g\in G$ such that $|g|<\infty$. If $g^k=e$, then $|g|$ divides $k$.
\end{crl}

Having obtained a criterion for $g^i=g^j$, we similarly would like to obtain a criterion for $\langle g^i\rangle=\langle g^j\rangle$ (which turns out, for the case $|g|<\infty$, to be the same as the criterion for $|g^i|=|g^j|$). Also, we would like be able to express the order $|g^k|:=|\langle g^k\rangle|$ in terms of the order $|g|:=|\langle g\rangle|$. (\blue{footnote}\footnote{Note that because of the cardinal equivalence $\Integer k\approx\Integer$ for all $k\neq 0$, if $|g|=\infty$ (so $\langle g\rangle\approx\langle g^k\rangle$ for all $k\neq 0$), then $|g^k|=\infty$ for all $k\neq 0$ (even when $k$ satisfies $\langle g^k\rangle\neq\langle g\rangle$). Consequently, meanwhile $\langle g^i\rangle=\langle g^j\rangle$ always implies $|g^i|=|g^j|$, if $|g|=\infty$ then we can have $|g^i|=|g^j|$ with $\langle g^i\rangle\neq\langle g^j\rangle$.}).

Even though the rest of the discussion below mostly involves \ul{finite} cyclic groups (important for studying/characterizing properties of finite groups in general) each of the results can be generalized simply by appending a statement for the case where at least one of the cyclic groups in question has \ul{infinite} order.

\begin{thm}
Let $G$ be a group and $g\in G$ such that $|g|<\infty$. Then for any positive integer $k\geq 1$,
\bea
\textstyle\langle g^k\rangle=\langle g^{\txt{gcd}(|g|,k)}\rangle~~~~\txt{and}~~~~|g^k|={|g|\over \txt{gcd}(|g|,k)}~~\left(~={\txt{lcm}(|g|,k)\over k}~\right).\nn
\eea
\end{thm}
\begin{proof}
Let $d:=\txt{gcd}(|g|,k)$. Observe that
\bea
\langle g^k\rangle=\{e\}\langle g^k\rangle=\langle e\rangle\langle g^k\rangle=\langle g^{|g|}\rangle\langle g^k\rangle=\langle g^{gcd(|g|,k)}\rangle=\langle g^d\rangle.\nn
\eea
Also, observe that for any divisor $d_n$ of $n:=|g|$, we have
\bea
&&\textstyle |g^{d_n}|=\min_{q\geq 1}\{q:~(g^{d_n})^q=e\}=\min_{q\geq 1}\{q:~qd_n=rn,~\txt{for the smallest possible integer}~r\geq 1\}\nn\\
&&\textstyle~~~~=\min_{q\geq 1}\left\{q:~qd_n=rn,~r=1\right\}=\min_{q\geq 1}\left\{q:~qd_n=n\right\}={n\over d_n},\nn
\eea
and so ~$|g^k|=|\langle g^k\rangle|=|\langle g^d\rangle|=|g^d|={n\over d}={|g|\over \txt{gcd}(|g|,k)}$.
\end{proof}

\begin{crl}[\blue{Element order divides (finite cyclic) group order}]\label{GroupOrdDiv1}
In a finite cyclic group, the order of an element divides the order of the group.
\end{crl}
\begin{proof}
Let $G$ be a finite cyclic group and consider a generator $c\in G$. Then $\langle c\rangle=G$, and so $|G|=|c|$. Any element $g\in G$ has the form $g=c^k$ for some integer $k\geq 0$. The result holds trivially if $k=0$. So assume $k\geq 1$. Then by the preceding theorem, $|g|=|c^k|={|c|\over\txt{gcd}(|c|,k)}={|G|\over\txt{gcd}(|G|,k)}$,
which by construction divides the group order $|G|=\txt{gcd}(|G|,k)~|g|$.
\end{proof}

\begin{crl}[\blue{Criterion for $\langle g^i\rangle=\langle g^j\rangle$ and $|g^i|=|g^j|$}]
Let $g\in G$ with $|g|<\infty$. Then
\bea
\langle g^i\rangle=\langle g^j\rangle~~\iff~~\txt{gcd}(|g|,i)=\txt{gcd}(|g|,j)~~\iff~~|g^i|=|g^j|.\nn
\eea
\end{crl}
\begin{proof}
It is clear that $|g^i|=|g^j|$ $\iff$ ${|g|\over \txt{gcd}(|g|,i)}=|g^i|=|g^j|={|g|\over \txt{gcd}(|g|,j)}$ $\iff$ $\txt{gcd}(|g|,i)=\txt{gcd}(|g|,j)$. It is also clear that ~$\txt{gcd}(|g|,i)=\txt{gcd}(|g|,j)~~\Ra~~\langle g^i\rangle=\langle g^{\txt{gcd}(|g|,i)}\rangle=\langle g^{\txt{gcd}(|g|,j)}\rangle=\langle g^j\rangle$,~ and conversely,
\[
\langle g^i\rangle=\langle g^j\rangle~~\Ra~~|g^i|=|g^j|,~~\Ra~~\txt{gcd}(|g|,i)=\txt{gcd}(|g|,j). \qedhere
\]
\end{proof}

\begin{dfn}[\blue{The ring of \index{Integers modulo $n$}{integers modulo $n$}}]
Let $n\in\Natural$. Then $\Integer_n:={\Integer\over\Integer n}=\{\ol{k}:=k+\Integer n:k\in\Integer\}=\left\{\ol{0},\ol{1},\cdots,\ol{n-1}\right\}$ as a ring under \ul{addition modulo $n$} and \ul{multiplication modulo $n$}  (i.e., addition and multiplication up to the remainder wrt $n$):
\bea
&&\ol{a}+\ol{b}:=\ol{a+b}=\ol{q_{n,a+b}n+r_{n,a+b}}=\ol{r_{n,a+b}},~~~~\txt{for}~~\ol{a},\ol{b}\in\Integer_n,\nn\\
&&\ol{a}\ol{b}:=\ol{ab}=\ol{q_{n,ab}n+r_{n,ab}}=\ol{r_{n,ab}},~~~~\txt{for}~~\ol{a},\ol{b}\in\Integer_n.\nn
\eea
The additive identity in $\Integer_n$ is $\ol{0}=\Integer n$, the additive inverse of $\ol{k}\in\Integer_n$ is $-\ol{k}=\ol{n-k}$, and the unity or multiplicative identity in $\Integer_n$ is $\ol{1}:=1+\Integer n$.
\end{dfn}
\begin{convention}[\blue{\index{Modular arithmetic}{Modular arithmetic}}]
For convenience, we often write {\small $\ol{k}=k+\Integer n\in\Integer_n=\left\{\ol{0},\ol{1},\cdots,\ol{n-1}\right\}$} simply as {\small $k\in\Integer_n=\{0,1,\cdots,n-1\}$} (i.e., by treating the \ul{quotient-absorption term} $\Integer n$ as invisible, we write $k+\Integer n$ simply as $k$), with addition and multiplication modulo $n$ understood, in which case,
\bea
a+b:=r_{n,a+b},~~ab:=r_{n,ab},~~~~\txt{for}~~a,b\in\Integer_n,\nn
\eea where $r_{n,a+b},r_{n,ab}$ are the respective remainders after dividing $a+b,ab\in\Integer$ by $n$. Sometimes (when it seems essential), we also explicitly define
\bea
k~\txt{mod}~n~:=~k+\Integer n,\nn
\eea
and then carry out integer addition and multiplication as usual, but only up to modulo $n$.
\end{convention}

\begin{dfn}[\blue{The \index{Group! of units modulo $n$}{group of units modulo $n$}}]
For temporary convenience, let $U(n):=U(\Integer_n)$, where
\bea
&&U(\Integer_n):=\{u+\Integer n\in\Integer_n:uu'+\Integer n=1+\Integer n~\txt{for some}~u'\in\Integer\}\nn\\
&&~~~~=
\{u+\Integer n\in\Integer_n:uu'+nu''=1~\txt{for some}~u',u''\in\Integer\}\nn\\
&&~~~~=\{u+\Integer n\in\Integer_n:gcd(n,u)=1\}\nn
\eea
is the group of units of the ring $(\Integer_n,+,\cdot)$.
\end{dfn}

\begin{crl}[\blue{Generators of a finite cyclic group}]\label{GenFinCycGp}
Let $g\in G$ with $|g|<\infty$. Then
\bea
\langle g\rangle=\langle g^j\rangle~~\iff~~|g|=|g^j|~~\iff~~j\in U(|g|).\nn
\eea
\end{crl}
Thus, if $G$ is cyclic, then given one generator $c$ (i.e., $G=\langle c\rangle$) of order $|c|$, all generators $c^j$ (i.e., $G=\langle c^j\rangle$) of $G$ are given by ~$gen(G)=\{c^j:~j\in \mathbb{Z},~\txt{gcd}(|c|,j)=1\}=\{c^j:~j\in U(|c|)\}$.

\begin{crl}[\blue{Generators of $(\mathbb{Z}_n,+)$}]
An integer $k\in \mathbb{Z}_n$ is a generator of $(\mathbb{Z}_n,+)$ $\iff$ $\txt{gcd}(n,k)=1$,
\bea
~~\iff~~\txt{Gen}(\Integer_n,+)=U(n).\nn
\eea
\end{crl}
\begin{proof}
Observe that for any cyclic group $G=\langle g\rangle$, we have $(\Integer_{|g|},+)\cong\langle g\rangle$, i.e., $(\Integer_{|G|},+)\cong G$, through the group homomorphism $h_g:(\Integer,+)\ra\langle g\rangle,~k\mapsto g^k$ with kernel ~$\ker h_g=\{k\in\Integer:g^k=e\}=\Integer|g|$.
\end{proof}

\subsection{Subgroups of a cyclic group}~\\~
Here, we see how many subgroups a finite cyclic group has, and how to find these subgroups.
\begin{thm}[\blue{\index{Fundamental theorem of! cyclic groups}{Fundamental theorem of cyclic groups}}]\label{FunThCycGps}
If $G=\langle c\rangle$ is a cyclic group, the following hold.
\bit[leftmargin=0.7cm]
\item[(1)] Every subgroup of $G$ is cyclic.
\item[(2)] If $|G|<\infty$, the order of every subgroup divides $|G|$.
\item[(3)] For any positive divisor $d\geq 1$ of $|G|=|c|$, the only subgroup of $G$ of order $d$ is $\langle c^{|c|\over d}\rangle$.
\item[(4)] Every subgroup of $G$ is of the form $\langle c^{|c|\over d}\rangle$ for a positive divisor $d\geq 1$ of $|G|=|c|$. That is, if $H\leq G$ is a subgroup, then (i) $|H|$ divides $|G|$ and (ii) $H=\langle c^{|G|\over|H|}\rangle$.
\eit
\end{thm}
\begin{proof}
{\flushleft(1)} Let $H\subset G=\langle c\rangle=\{c^k:k\in\Integer\}$ be a subgroup. Observe that $H=\{g^k:k\in H_0\}$, where $H_0:=\{k\in\Integer:c^k\in H\}\subset\Integer$. It is easy to see that $H_0$ is a subgroup of $(\Integer,+)$. Therefore, $H_0=\Integer m$, where $m\geq 1$ is the smallest positive integer in $H_0$. Hence, ~$H=\{c^k:k\in H_0=\Integer m\}=\langle c^m\rangle$.

{\flushleft(2)} Let $H=\langle c^m\rangle$, where $m$ is as defined as above, i.e., $m:=\min\{k\in \mathbb{Z}:~k\geq 1,~c^k\in H\}$. Then by Corollary \ref{GroupOrdDiv1}, $|H|=|c^m|$ divides $|G|$.

{\flushleft(3)} Let $d\geq 1$ be any positive divisor of $|G|=|c|$. Then $\langle c^{|c|\over d}\rangle$ is certainly a subgroup of $G$ of order ${|c|\over\txt{gcd}(|c|,|c|/d)}=d$. Let $H=\langle c^m\rangle\subset G$ (with $m$ as defined above) be any subgroup of order $d\geq 1$. Then
\bea
\textstyle d=|H|=|c^m|=|c^{\txt{gcd}(|c|,m)}|={|c|\over \txt{gcd}(|c|,m)}\sr{(s)}{=}{|c|\over m},~~\Ra~~m={|c|\over d},\nn
\eea
where step (s) holds because $m$ divides $|c|$ (since $c^{|c|}=e\in H$).
{\flushleft(4)} This follows immediately from (2) and (3) above.
\end{proof}

\begin{crl}[\blue{Subgroups of $\mathbb{Z}_n$}]
For each positive divisor $d=d_n$ of $n$, $\langle{n\over d}\rangle$ is the unique subgroup of $\mathbb{Z}_n$ of order $d$. Moreover, these are the only subgroups of $\mathbb{Z}_n$ (since the order of every element divides $n=|\Integer_n|$).
\end{crl}

\begin{rmk}
Let $d,n\in\Integer$. If $d|n$ (i.e., $d$ is a divisor of $n$), then the subgroup $\langle{n\over d}\rangle\subset\Integer_n$ satisfies
\bea
\textstyle\langle{n\over d}\rangle:=\{{n\over d}(k+\Integer n):k+\Integer n\in\Integer_n\}=\{{n\over d}k+\Integer n:0\leq k<d\}\cong\Integer_d.\nn
\eea
Hence, we have a bijective correspondence of the form
\[
\textstyle\{\txt{subgroups of}~\Integer_n\}\longleftrightarrow\{\langle{n\over d}\rangle\leq\Integer_n:d|n\}\longleftrightarrow\{\Integer_d:d|n\}.
\]
\end{rmk}

\begin{thm}[\blue{Number of elements of each order in a finite cyclic group}]
Let $G=\langle c\rangle$ be a cyclic group of order $|G|<\infty$. For any positive divisor $d\geq 1$ of $|G|$, the number of elements of order $d$ is
\bea
|\{g\in G:|g|=d\}|=\phi(d):=|U(d)|.\nn
\eea
\end{thm}
\begin{proof}
Let $d$ be a positive divisor of $n:=|G|$. The number of elements of order $d$ is precisely the number of generators of the unique subgroup $\langle c^{n\over d}\rangle$ of order $d$ (where the uniqueness is due to Theorem \ref{FunThCycGps}). From Corollary \ref{GenFinCycGp} the number of generators of $\langle c^{n\over d}\rangle$ is $|U(|c^{n\over d}|)|=|U(d)|$.
\end{proof}

\begin{crl}[\blue{Number of elements of order $d$ in a finite group}]\label{NoWithCertOrd}
In a finite group, the number of elements of order $d$ is a multiple of~ $\phi(d):=|U(d)|$.
\end{crl}
\begin{proof}
Let $G$ be a finite group, i.e., $|G|<\infty$, and $g\in G$ with order $|g|=d$. Then by the preceding theorem, the cyclic subgroup $\langle g\rangle\subset G$ has a total of $\phi(d):=|U(d)|$ elements of order $d$, i.e., with $[g]_d:=\{a\in\langle g\rangle:|a|=d\}\subset\langle g\rangle$, we have $|[g]_d|=\phi(d):=|U(d)|$. Moreover, for any $g,g'\in G$ of order $d$, either
\bea
[g]_d=[g']_d~~~~\txt{or}~~~~[g]_d\cap[g']_d=\emptyset,\nn
\eea
since any member of $[g]_d\subset\langle g\rangle$ can generate the whole subgroup $\langle g\rangle\subset G$. Hence, the total number of elements of $G$ of order $d$ must be a multiple of $\phi(d):=|U(d)|$.
\end{proof}

\begin{lmm}[\blue{Arithmetic of integer ideals}]
Let $a,b\in \Integer$ be positive integers. Then the following hold:
\bit[leftmargin=0cm]
\item[] (i) $\Integer a+\Integer b=\Integer gcd(a,b)$, (ii) $\Integer a\cap\Integer b=\Integer lcm(a,b)$, (iii) $(\Integer a)(\Integer b)=\Integer ab$, (iv)  If $a|b$, then ${\Integer a\over\Integer b}\cong{\Integer\over\Integer(b/a)}$.
\eit
In particular, if $a,b$ are relatively prime, then ~$\Integer a\cap\Integer b=\Integer ab=(\Integer a)(\Integer b)$.
\end{lmm}
\begin{proof}
We already know (i),(ii),(iii). So we will prove (iv). Consider the group epimorphism
\bea
\phi:\Integer\ra \Integer a/\Integer b,~~k\mapsto ka+\Integer b.\nn
\eea
Then by the 1st isomorphism theorem, ~${\Integer\over\ker\phi}\cong\im\phi={\Integer a\over\Integer b}$,~ where
{\small\[
\textstyle\ker\phi=\{k\in\Integer: ka+\Integer b=\Integer b\}=\{k\in\Integer: ka\in\Integer b\}=\left\{k\in\Integer: k\in\Integer {b\over a}\right\}=\Integer {b\over a}. \qedhere
\]}
\end{proof}

\begin{lmm}[\blue{\index{Gcd-Lcm remainder theorem}{Gcd-Lcm remainder theorem}}]\label{GcdLcmRemThm}
For any positive integers $a,b\geq 1$, ~${\Integer a+\Integer b\over\Integer a\cap\Integer b}\cong{\Integer a+\Integer b\over\Integer a}\times{\Integer a+\Integer b\over\Integer b}$. In particular, if $a,b$ are relatively prime, i.e., $\txt{gcd}(a,b)=1$, then ~${\Integer\over\Integer ab}\cong{\Integer\over\Integer a}\times{\Integer\over\Integer b}$.
\end{lmm}
\begin{proof}
Let $d:=gcd(a,b)$. Then $\Integer d=\Integer a+\Integer b$. Consider the group homomorphism
\bea
\textstyle \psi:\Integer d\ra {\Integer d\over\Integer a}\times{\Integer d\over\Integer b},~kd\mapsto(kd+\Integer a,kd+\Integer b).\nn
\eea
Then $\ker\psi=\Integer a\cap\Integer b=\Integer l$ (where $l:=lcm(a,b)$). Since $\txt{gcd}(a,b)=d$, we have $xa+yb=d$ for some integers $x,y\in\Integer$. For any $(u+\Integer a,v+\Integer b)\in{\Integer d\over\Integer a}\times {\Integer d\over\Integer b}$, let $kd:=x{v\over d}a+y{u\over d}b$. Then
{\small $\psi(kd)=\big(kd+\Integer a,kd+\Integer b\big)=\big(x{v\over d}a+y{u\over d}b+\Integer a,x{v\over d}a+y{u\over d}b+\Integer b\big)=\big(y{u\over d}b+\Integer a,x{v\over d}a+\Integer b\big)=\big(u+\Integer a,v+\Integer b\big)$}. That is, $\im\psi={\Integer d\over\Integer a}\times {\Integer d\over\Integer b}$. Hence by the 1st isomorphism theorem, {\small${\Integer d\over\Integer a\cap\Integer b}={\Integer d\over\ker\psi}\cong\im\psi={\Integer d\over\Integer a}\times{\Integer d\over\Integer b}$}.
\end{proof}

\begin{thm}
Let $p\geq 0$ be a prime and $n\geq 1$ a positive integer. (As before, let $\phi(k):=|U(k)|:=|U(\Integer_k)|$ for any positive integer $k\geq 1$). Then (i) $\phi(p^n)=p^n-p^{n-1}$. Also, (ii) if $m$ and $n$ are relatively prime, then $\phi(mn)=\phi(m)\phi(n)$.
\end{thm}
\begin{proof}
{\flushleft(i)} Since $p$ is a prime, an element of $\Integer_{p^n}$ is not relatively prime to $p^n$ (i.e., does not lie in $U(p^n):=\{k:gcd(p^n,k)=1\}\subset\Integer_{p^n}$) $\iff$ it has $p$ as a factor, and so the set of such elements is $\{kp:~1\leq k\leq p^{n-1}\}$. Therefore, ~$\phi(p^n)=|U(p^n)|=|\Integer_{p^n}|-|\{kp:~1\leq k\leq p^{n-1}\}|=p^n-p^{n-1}$.
{\flushleft (ii)} From Lemma \ref{GcdLcmRemThm}, if $m$ and $n$ are relatively prime, then {\footnotesize $U(mn)\cong U(m)\times U(n)$}, (\blue{footnote}\footnote{The group isomorphism ~$\mathbb{Z}_m\times \mathbb{Z}_n\cong \mathbb{Z}_{mn}$~ upgrades to a ring isomorphism $\mathbb{Z}_m\times\mathbb{Z}_n\cong \mathbb{Z}_{mn}$. Moreover, an isomorphism of rings $R\cong R'$ restricts to an isomorphism of groups ~$U(R)\cong U(R')$.}) and thus
\[
\phi(mn)=|U(mn)|=|U(m)\times U(n)|=|U(m)||U(n)|=\phi(m)\phi(n). \qedhere
\]
\end{proof}

\section{Review of Finite Products}
Recall that groups (as well as rings, modules, and algebras) can be combined in various ways to form larger ones. One method involves taking products of groups (resp. rings, modules, algebras). When the process is viewed in reverse, it also reveals how larger groups (resp. rings, modules, algebras) can be decomposed into products of smaller ones.

\subsection{Definitions and well-definedness}
\begin{dfn}[\blue{\index{Finite! products}{Finite product of groups}}]
The product of a collection of groups $G_1,...,G_n$ is the set of $n$-tuples $G_1\times\cdots\times G_n=\{(g_1,...,g_n),~~g_i\in G_i\}$ as a group
 with \ul{multiplication}, \ul{inverse}, and \ul{identity element} given by the following componentwise operations:
{\small\bea
(g_1,...,g_n)(g_1',...,g_n')~:=~(g_1g_1',...,g_ng_n'),~~~~(g_1,...,g_n)^{-1}:=(g_1^{-1},...,g_n^{-1}),~~~~e_{G_1\times\cdots\times G_n}:=(e_{G_1},\cdots,e_{G_n}).\nn
\eea}
\end{dfn}

\begin{thm}[\blue{The product of groups is well defined}]
Given groups $G_1,...,G_n$, their product $G_1\times\cdots\times G_n$ is indeed a group with respect to the given structure.
\end{thm}
\begin{proof}
Associativity of multiplication follows from associativity of the component-group multiplications, i.e., for any $(a_1,...,a_n),(b_1,...,b_n),(c_1,c_2,...,c_n)\in G_1\times\cdots\times G_n$,
{\small \bea
 &&\big((a_1,...,a_n)(b_1,...,b_n)\big)(c_1,...,c_n)=(a_1b_1,...,a_nb_n)(c_1,...,c_n)=\big((a_1b_1)c_1,...,(a_nb_n)c_n\big)\nn\\
 &&~~~~=\big(a_1(b_1c_1),...,a_n(b_nc_n)\big)=(a_1,...,a_n)(b_1c_1,...,b_nc_n)=(a_1,...,a_n)\big((b_1,...,b_n)(c_1,...,c_n)\big).\nn
 \eea}With $e_i:=e_{G_i}$, ~$(e_1,...,e_n)$ is an identity in $G_1\times\cdots\times G_n$ since for any $(g_1,...,g_n)\in G_1\times\cdots\times G_n$,
{\small\bea
(g_1,...,g_n)(e_1,...,e_n)=(g_1e_1,...,g_ne_n)=(e_1g_1,...,e_ng_n)=(e_1,...,e_n)(g_1,...,g_n)=(g_1,....,g_n).\nn
\eea}Similarly, any $(g_1,...,g_n)\in G_1\times G_2\times\cdots\times G_n$ has an inverse $(g_1^{-1},g_2^{-1},...,g_n^{-1})$ since
{\small\[
(g_1,...,g_n)(g_1^{-1},...,g_n^{-1})=(g_1g_1^{-1},...,g_ng_n^{-1})=(g_1^{-1}g_1,...,g_n^{-1}g_n)=(g_1^{-1},...,g_n^{-1})(g_1,...,g_n)=(e_1,...,e_n). \qedhere
\]}
\end{proof}

For example, if $R$ is a ring, then $R^n=R\times R^{n-1}$ is an abelian group under componentwise addition:
\bea
(a_1,...,a_n)+(b_1,....,b_n)=(a_1+b_1,...,a_n+b_n).\nn
\eea
It is easy to see that ~$|G_1\times G_2\times\cdots\times G_n|=|G_1|~|G_2|~\cdots~|G_n|$, especially when the groups are finite.

\begin{rmk}[\textcolor{blue}{Finite products of rings, modules, algebras}]
The above procedure can be repeated for rings, modules, and algebras. In particular, given a commutative ring $R$ and $R$-algebras $A_1,...,A_n$, their product is the set of $n$-tuples $A_1\times\cdots\times A_n=\{(a_1,...,a_n),~~a_i\in A_i\}$ as an $R$-algebra
 with \ul{addition}, \ul{multiplication}, \ul{scalar multiplication}, \ul{zero}, and \ul{identity element} given by
{\small\bea
&&(a_1,...,a_n)+(a_1',...,a_n'):=(a_1+a_1',...,a_n+a_n'),~~~~(a_1,...,a_n)(a_1',...,a_n'):=(a_1a_1',...,a_na_n')\nn\\
&&r(a_1,...,a_n):=(ra_1,\cdots,ra_n),~~~~0_{A_1\times\cdots\times A_n}:=(0_{A_1},\cdots,0_{A_n}),~~~~1_{A_1\times\cdots\times A_n}:=(1_{A_1},\cdots,1_{A_n}).\nn
\eea}
\end{rmk}

Given our interest in cyclic groups, it is natural to ask a question such as ``When is $\mathbb{Z}_m\times\mathbb{Z}_n\cong \mathbb{Z}_{mn}$?'', or more generally, ``When is a product of groups cyclic?''.

\subsection{Cyclic nature of finite products of groups}~\\~
Since every subgroup of a cyclic group is cyclic, it follows that $G_1\times\cdots\times G_n$ cannot be cyclic if some $G_i$ fails to be cyclic. Therefore, in theorems characterizing cyclic nature of products of groups, it will be enough to consider only products of cyclic groups.
\begin{dfn}[\blue{\index{Order! types}{Order types} of elements}]
Let $G_1,...,G_n$ be groups. The \ul{order type} of a group element $(g_1,...,g_n)\in G_1\times...\times G_n$ is the tuple of component-group orders $(|g_1|,...,|g_n|)$.
\end{dfn}
\begin{thm}[\blue{Order of a product-group element}]\label{external-product-order}
Let $G_1,...,G_n$ be groups. The order of a group element $(g_1,...,g_n)\in G_1\times...\times G_n$ is the lcm of its order type, i.e., ~$|(g_1,...,g_n)|=\txt{lcm}(|g_1|,...,|g_n|)$.
\end{thm}
\begin{proof}
For any integer $k\geq 1$, ~$(g_1,...,g_n)^k=(g_1{}^k,...,g_n{}^k)$, and so if $(g_1,...,g_n)^k=(g_1{}^k,...,g_n{}^k)=(e_1,...,e_n)$, then each of the orders $|g_i|$ must divide $k$, i.e., $k$ must be a common multiple of the orders $|g_1|$,...,$|g_n|$. Hence, ~$|(g_1,...,g_n)|:=\inf\{k\in\Integer:k\geq 1,~(g_1,...,g_n)^k=(e_1,...,e_n)\}=\txt{lcm}(|g_1|,...,|g_n|)$.
\end{proof}

\begin{thm}[\blue{Criterion for $G\times G'$ to be cyclic}]
Let $G$ and $G'$ be finite cyclic groups. Then $G\times G'$ is cyclic if and only if $|G|$ and $|G'|$ are relatively prime, i.e., $\txt{gcd}(|G|,|G'|)=1$.
\end{thm}
\begin{proof}
($\Ra$) Assume $G\times G'$ is cyclic. Then $G\times G'=\langle(g,g')\rangle$ for some $(g,g')\in G\times G'$, and so
\[
\textstyle |G||G'|=|G\times G'|=|(g,g')|=\txt{lcm}(|g|,|g'|)\leq\txt{lcm}(|G|,|G'|)={|G||G'|\over\txt{gcd}(|G|,|G'|)}.
\]
That is, $\txt{gcd}(|G|,|G'|)\leq 1$. Hence $\txt{gcd}(|G|,|G'|)=1$.

($\La$) Conversely, assume $\txt{gcd}(|G|,|G'|)=1$. Let $G=\langle a\rangle$ and $G'=\langle a'\rangle$. Then
{\small\[
|G\times G'|=|G||G'|=\txt{gcd}(|G|,|G'|)\txt{lcm}(|G|,|G'|)=\txt{lcm}(|G|,|G'|)=\txt{lcm}(|a|,|a'|)=|(a,a')|=|\langle(a,a')\rangle|.\nn
\]}Hence $(a,a')$ generates $G\times G'$ (i.e., $G\times G'$ is cyclic) since $G\times G'$ is finite.
\end{proof}

\begin{crl}[\blue{Criterion for {\small $G_1\times...\times G_n$} to be cyclic}]
Let $G_1,...,G_n$ be finite cyclic groups. Then {\small $G_1\times\cdots\times G_n$} is cyclic if and only if $|G_i|$ and $|G_j|$ are relatively prime (i.e., {\small $\txt{gcd}(|G_i|,|G_j|)=1$}) for all $i\neq j$.
\end{crl}
\begin{proof}
{\flushleft ($\Ra$)} Assume $G_1\times...\times G_n$ is cyclic. Then $G_i\times G_j$ (for each $i\neq j$) is cyclic as a subgroup.

{\flushleft ($\La$)} Assume $\txt{gcd}(|G_i|,|G_j|)=1$ for each $i\neq j$. We will use induction on $n$. By hypotheses $G_1$ is cyclic. Let $1\leq k<n$. If $G_1\times\cdots\times G_k$ is cyclic, then because (\blue{footnote}\footnote{For any $a,b,c\in\Integer$, if $\gcd(a,b)=\gcd(a,c)=\gcd(b,c)=1$, then $\gcd(ab,c)=\gcd(a,bc)=1$, since $\Integer=\Integer a+\Integer c=\Integer(b+c)a+\Integer c=\Integer ba+\Integer ca+\Integer c=\Integer ab+\Integer c$ and by symmetry we also have $\Integer=\Integer a+\Integer bc$.})
{\small\begin{align}
1\leq\gcd(|G_1\times\cdots\times G_k|,|G_{k+1}|)=\gcd(\txt{lcm}(|G_1|,\cdots, |G_k|),|G_{k+1}|)\leq\gcd\big(|G_1||G_2|\cdots|G_k|,|G_{k+1}|\big)=1,\nn
\end{align}}
it follows that $G_1\times\cdots\times G_{k+1}=(G_1\times\cdots\times G_k)\times G_{k+1}$ is also cyclic.

Hence $G_1\times...\times G_n$ is cyclic by induction.
\end{proof}

\begin{crl}[\blue{Criterion for $\mathbb{Z}_{n_1}\times \mathbb{Z}_{n_2}\times...\times \mathbb{Z}_{n_k}\cong \mathbb{Z}_{n_1n_2\cdots n_k}$}]
We have ~{\small $\mathbb{Z}_{n_1}\times \mathbb{Z}_{n_2}\times...\times \mathbb{Z}_{n_k}\cong\mathbb{Z}_{n_1n_2\cdots n_k}$} $\iff$ $n_i$ and $n_j$ are relatively prime (i.e., $\txt{gcd}(n_i,n_j)=1$) for each $i\neq j$.
\end{crl}

\begin{crl}[\blue{Criterion for $U(n_1)\times...\times U(n_k)\cong U(n_1n_2\cdots n_k)$}]
We have ~{\small $U(n_1n_2...n_k)\cong U(n_1)\times U(n_2)\times...\times U(n_k)$} $\iff$ $n_i$ and $n_j$ are relatively prime (i.e., $\txt{gcd}(n_i,n_j)=1$) for each $i\neq j$. (\blue{footnote}\footnote{Recall that ~$U(n):=U(\Integer_n)=\{k+\Integer n\in\Integer_n:gcd(k,n)=1\}$.})
\end{crl}
\begin{proof}
The group isomorphism ~$\mathbb{Z}_{n_1}\times \mathbb{Z}_{n_2}\times...\times \mathbb{Z}_{n_k}\cong \mathbb{Z}_{n_1n_2\cdots n_k}$~ upgrades to a ring isomorphism. Moreover, an isomorphism of rings $R\cong R'$ restricts to an isomorphism of groups ~$U(R)\cong U(R')$.
\end{proof}

\begin{crl}
Let $G_1,...,G_n$, $H_1,...,H_n$ be groups and $(g_1,...,g_n)\in G_1\times\cdots\times G_n$ an element of finite order (i.e., $|(g_1,...,g_n)|<\infty$). The following hold:
\bit
\item[(a)] $\langle(g_1,\cdots,g_n)\rangle=\langle g_1\rangle\times\cdots\times\langle g_n\rangle$ $\iff$ $|g_i|$ and $|g_j|$ are relatively prime for each $i\neq j$.
\item[(b)] $G_1\times\cdots\times G_n$ is abelian $\iff$ each $G_i$ is abelian.
\item[(c)] If $G_i\cong H_i$ for each $i$, then $G_1\times\cdots\times G_n\cong H_1\times\cdots\times H_n$.
\eit
\end{crl}
\begin{proof}
(a) follows from the previous results. (b) follows from the definitions.
For (c), if $f_i:G_i\ra H_i$ are group isomorphisms, then so is {\small $f:G_1\times\cdots\times G_n\ra H_1\times\cdots\times H_n,~(g_1,...,g_n)\mapsto(f_1(g_1),\cdots,f_n(g_n))$}.
\end{proof}

\section{Symmetry Groups (a technical digression)}
\begin{rmk}[\magenta{Caution}]
The word ``cyclic permutation'' (as used below) is a notion of ``\ul{cyclic group element}'', which is different from our notion of a cyclic object in a subcategory of Sets (and instead closer to a notion of ``\ul{cyclic morphism}'' in such a category).
\end{rmk}
\subsection{Introduction and cyclic finite permutations}
\begin{dfn}[\blue{
\index{Symmetry (Permutation)}{Symmetry (Permutation)},
\index{Symmetry group}{Symmetry group},
\index{Finite! permutations}{Finite permutations}}]
Let $X$ be a set. Recall that a bijection of the form $f:X\ra X$ is called a \ul{symmetry} (or \ul{permutation}) of $X$. The \ul{symmetry group} of $X$ is the set
\bea
Sym(X):=Perm(X):=\{\txt{all bijections ~$f:X\ra X$}\}\nn
\eea
as a group with multiplication the map composition $\circ: (f,f')\mapsto f\circ f'$ and identity element the identity map $id_X:X\ra X$. The symmetry (permutation) group of $n$ objects $o_1,o_2,...,o_n$ is denoted by
\[S_n:=Sym(\mathbb{X}_n):=Perm(\mathbb{X}_n),~~~~\mathbb{X}_n:=\{o_1,o_2,...,o_n\}\approx\{1,2,...,n\}\approx n:=\{0,1,\cdots,n-1\}.\]
We will refer to the elements of $S_n$ as \ul{$n$-permutations} (or \ul{finite permutations} of degree $n$). The identity element $id_{\mathbb{X}_n}$ (or $1_{\mathbb{X}_n}$) in $S_n$ will be denoted by $\ep\in S_n$.
\end{dfn}
\begin{notation*}[\blue{\index{Two-row notation}{Two-row notation}}]
For each $\al\in S_n$, we give $\al$ a convenient \ul{two-row notation} as follows:
{\footnotesize\[
\al=\left[
      \begin{array}{cccc}
        1 & 2 & \cdots & n \\
        \al(1) & \al(2) & \cdots & \al(n) \\
      \end{array}
    \right],~~~~
\al^{-1}=\left[
      \begin{array}{cccc}
        \al(1) & \al(2) & \cdots & \al(n) \\
        1 & 2 & \cdots & n \\
      \end{array}
    \right].\nn
\]}If $\al,\beta\in S_n$, the product (composition) of $\al$ and $\beta$ is given by
{\footnotesize\[
\al\beta=\left[
      \begin{array}{cccc}
        1 & 2 & ... & n\\
        \al(1) & \al(2) & ... & \al(n) \\
      \end{array}
    \right]\left[
      \begin{array}{cccc}
        1 & 2 & ... & n\\
        \beta(1) & \beta(2) & ... & \beta(n) \\
      \end{array}
    \right]
:=\left[
      \begin{array}{cccc}
        1 & 2 & ... & n\\
        \al\big(\beta(1)\big) & \al\big(\beta(2)\big) & ... & \al\big(\beta(n)\big) \\
      \end{array}
    \right].\nn
\]}
\end{notation*}

\begin{dfn}[\blue{
\index{Cyclic! permutation (or Cycle)}{Cyclic permutation (or Cycle)},
\index{Length of a cycle}{Cycle length},
\index{Transposition}{Transposition}}]
For any distinct $x_1,...,x_l\in\mathbb{X}_n$, consider the permutation $(x_1,...,x_l):\mathbb{X}_n\ra\mathbb{X}_n$ defined by
\bea
(x_1,...,x_l):x_i\mapsto x_{i+1\!\!\mod l},~~~~\txt{and}~~~~(x_1,...,x_l):x\mapsto x~~\txt{if}~~x\not\in \{x_1,...,x_l\},\nn
\eea
i.e., $(x_1,...,x_l)$ moves the elements $x_1,...,x_l$ in a cycle, while fixing any element of $\mathbb{X}_n$ not in $\{x_1,...,x_l\}$. The permutation $(x_1,...,x_l)$ is called a \ul{cyclic permutation} (or a \ul{cycle}) of \ul{length} $l$, or simply an \ul{$l$-cycle}. A $2$-cycle is called a \ul{transposition}.
\end{dfn}

Note that in the cycle $(x_1,...,x_l)$, since $(x_1,...,x_l)(x_i)=x_{i+1\!\!\mod l}$, we have
\bea
(x_1,...,x_j,x_{j+1},...,x_l)\neq (x_1,...,x_{j+1},x_j,...,x_l),\nn
\eea
i.e., the ordering of the components $x_1,...,x_l$ is important, while it is clear (by definition) that
\bea
(x_1,...,x_l)=(x_l,x_1,...,x_{l-1})=(x_{l-1}x_l,x_1,...,x_{l-2})=\cdots=(x_3,...,x_l,x_1,x_2)=(x_2,...,x_l,x_1),\nn
\eea
i.e., cyclic shifting or translation of the components does not change the permutation.

Notice that any $1$-cycle $(x)$, $x\in\mathbb{X}_n$, is equivalent to the identity permutation $\ep$, which as usual we also denote simply by $1$.

\begin{lmm}[\blue{2-cycle decomposition of a cycle, Inverse of a cycle}]\label{TwoCycDecLmm}
(i) Every cycle is a product of $2$-cycles. (ii) The inverse of the cycle $(x_1x_2...x_r)$ is $(x_rx_{r-1}...x_1)$.
\end{lmm}
\begin{proof}
(i) A basic observation is that each $l$-cycle can be written as a product of $l-1$ transpositions in two different ways as follows (as can be directly verified using the definition of the cycle):
\begin{align}
(x_1,...,x_l)&=(x_1,x_l)(x_1,x_{l-1})\cdots (x_1,x_2)\nn\\
\label{TransDecEq}&=(x_1,x_2)(x_2,x_3)(x_3,x_4)\cdots (x_{l-1},x_l).
\end{align}
(ii) The inverse of a cycle is therefore given by
\begin{align}
&(x_1,...,x_l)^{-1}=[(x_1,x_2)(x_2,x_3)(x_3,x_4)\cdots (x_{l-1},x_l)]^{-1}=(x_{l-1},x_l)^{-1}(x_{l-2},x_{l-1})^{-1}\cdots (x_2,x_3)^{-1}(x_1,x_2)^{-1}\nn\\
&~~~~=(x_l,x_{l-1})(x_{l-1},x_{l-2})\cdots (x_3,x_2)(x_2,x_1)=(x_l,x_{l-1},...,x_1).\nn \qedhere
\end{align}
\end{proof}

\begin{rmk}[\blue{The order of a cycle}]
Since a cycle moves its components (as elements of $\mathbb{X}_n$) in a cycle while fixing all other elements of $\mathbb{X}_n$, an $l$-cycle has order $l$ since it returns any one of its components to its original position after exactly $l$ steps (and thus requires a minimum of $l$, and only $l$, steps to simultaneously return every element of $\mathbb{X}_n$ to itself). That is,
\[
|(x_1,...,x_l)|=l.\nn
\]
\end{rmk}

\begin{rmk}[\blue{The conjugate of a cycle}]
Using the definition of the $l$-cycle, we can also directly verify that given any permutation $\al$ and an $l$-cycle $(x_1,...,x_l)$, the permutation $\al(x_1,...,x_l)\al^{-1}$, which is clearly of order $l$ (\blue{footnote}\footnote{Recall that if $G=[G,\cdot,e]$ is a group and $\al,g\in G$, then $(\al g\al^{-1})^k=\al g^k\al^{-1}$, and so $(\al g\al^{-1})^k=e$ iff $g^k=e$. Alternatively, observe that the map $f_\al:\langle g\rangle\ra \langle \al g\al^{-1}\rangle,~h\mapsto\al h\al^{-1}$ is a bijection.}), is also an $l$-cycle given by
\[
\al(x_1,...,x_l)\al^{-1}=\big(\al(x_1),...,\al(x_l)\big).
\]
\end{rmk}

\begin{lmm}[\blue{Disjoint cycles commute}]
If two cyclic permutations $\al=(x_1,...,x_r)$ and $\beta=(y_1,...,y_s)$ are disjoint, in the sense they have no common entries, then $\al\beta=\beta\al$.
\end{lmm}
\begin{proof}
Consider our standard finite set $\mathbb{X}_n$ (on which $\al$ and $\beta$ act) in the form
\[
\mathbb{X}_n=\{x_1,...,x_r,y_1,...,y_s,z_1,...,z_t\},
\]
where $z_1,...,z_t$ are the elements fixed under both $\al$ and $\beta$. Since $\al$ and $\beta$ are disjoint, each element $x\in \mathbb{X}_n$ is fixed by at least one of $\al$ and $\beta$. (\blue{footnote}\footnote{Let $x\in\mathbb{X}_n$. It is clear that if $\al$ fixes $x$, then $\al$ also fixes $\beta^k(x)$ for all $k\in\Integer$, and similarly, if $\beta$ fixes $x$, then $\beta$ also fixes $\al^k(x)$ for all $k\in\Integer$.}). Hence it suffices to consider the effect of $\al$ and $\beta$ on $x\in \mathbb{X}_n$ under the following two exhaustive cases:
\bit[leftmargin=0.7cm]
\item If $x$ is fixed by $\al$, then $(\al\beta)(x)=\al(\beta(x))=\beta(x)=\beta(\al(x))=(\beta\al)(x)$.
\item
If $x$ is fixed by $\beta$, then $(\al\beta)(x)=\al(\beta(x))=\al(x)=\beta(\al(x))=(\beta\al)(x)$. \qedhere
\eit
\end{proof}

\subsection{Cycle decomposition and parity of finite permutations}
\begin{thm}[\blue{Disjoint-cycle decomposition}]\label{DisjCycDec}
Every permutation $\al\in S_n$ can be uniquely written as the product of a finite number of disjoint cycles.
\end{thm}
\begin{proof}
Recall that $S_n:=Sym(\mathbb{X}_n)$, where $\mathbb{X}_n\approx\{1,2,...,n\}$. Using $\al$, define an equivalence relation $\sim_\al$ on ${\mathbb{X}_n}$ as follows: For $x,y\in {\mathbb{X}_n}$,
\bea
x\sim_\al y~~~~\txt{if}~~~~\al^k(x)=y~~~~\txt{for some}~~~~k\in\Integer.\nn
\eea
The equivalence class $[x]_\al:=\left\{y\in \mathbb{X}_n:~x\sim_\al y\right\}$ is precisely the \index{Orbit}{\ul{orbit}} $\langle\al\rangle(x):=\{\al^k(x)~|~k\in\Integer\}$ of $x$ under $\al$ (i.e., under the subgroup $\Integer_{|\al|}\cong\langle\al\rangle\leq S_n$, where $|\al|$ is the order of $\al$ in $S_n$). These orbits are of course disjoint because for any $x,y\in\mathbb{X}_n$, if $\langle\al\rangle(x)\cap\langle\al\rangle(y)\neq\emptyset$, i.e., $\al^i(x)=\al^j(y)$ for some $i,j$ ($i\geq j$ wlog), then $y=\al^{i-j}(x)\in\langle\al\rangle(x)$, and so $\langle\al\rangle(x)=\langle\al\rangle(y)$.

Now, for a fixed $x\in {\mathbb{X}_n}$, let
\[
l_\al(x):=\min\left\{k\geq 1:~\al^k(x)=x\right\}=|\langle\al\rangle(x)|~~\leq~~|\al|.
\]
Then for each $0\leq k\leq l_\al(x)-1$, we have $\al^{-k}(x)=\al^{-k}\al^{l_\al(x)+1}(x)=\al^{l_\al(x)-k+1}(x)$, and so $l_\al(x)$ is the order of the subgroup $\Integer_{l_\al(x)}\cong\{\al^k:k\in\Integer_{l_\al(x)}\}\leq\langle\al\rangle\cong\Integer_{|\al|}$.

By applying $\al$ repeatedly on the fixed $x\in {\mathbb{X}_n}$ we get the finite sequence
\bea
x_1:=x,~~~~x_2:=\al(x),~~~~x_3:=\al^2(x),~~~~\cdots,~~~~x_{l_\al(x)}:=\al^{l_\al(x)-1}(x),\nn
\eea
using which we can then write
\bea
[x]_\al:=\left\{y\in \mathbb{X}_n:~x\sim_\al y\right\}=\langle\al\rangle(x)=\left\{x_1,...,x_{l_\al(x)}\right\}.\nn
\eea

Since $\sim_\al$ partitions $\mathbb{X}_n$, it follows that $\al$ (being a pointwise operation) can be expressed as the ``\emph{product}'' (or composition) of a finite number of disjoint (hence commuting) cycles as follows:
\begin{align}
\al&=\big([x]_\al\big)\big([y]_\al\big)\big([z]_\al\big)\cdots~=~\al|_{\langle\al\rangle(x)}~\al|_{\langle\al\rangle(y)}~\al|_{\langle\al\rangle(z)}~\cdots,~~~~\txt{for some}~~~~x,y,z,\cdots\in \mathbb{X}_n,\nn\\
\label{CycleDecEq}&=(x_1,...,x_{l_\al(x)})(y_1,...,y_{l_\al(y)})(z_1,...,z_{l_\al(z)})\cdots
\end{align}
That is, with the partition {\footnotesize $X:=[x]_\al=\langle\al\rangle(x)=\{x_1,...,x_{l_\al(x)}\}$, $Y:=[y]_\al=\langle\al\rangle(y)=\{y_1,...,y_{l_\al(y)}\}$ , $Z:=[z]_\al=\langle\al\rangle(z)=\{z_1,...,z_{l_\al(z)}\}$,} $\cdots$ of $\mathbb{X}_n$, the map $\al:\mathbb{X}_n\ra\mathbb{X}_n$ splits up into parallel/disjoint submaps
\[
(x_1,...,x_{l_\al(x)})\eqv\al|_X:X\ra X,~~~(y_1,...,y_{l_\al(y)})\eqv\al|_Y:Y\ra Y,~~~(z_1,...,z_{l_\al(z)})\eqv\al|_Z:Z\ra Z,~~~\cdots. \qedhere
\]
\end{proof}

\begin{crl}[\blue{$2$-cycles decomposition}]
Every permutation in $S_n,~n>1,$ is a product of $2$-cycles.
\end{crl}
\begin{proof}
This follows immediately from Lemma \ref{TwoCycDecLmm} and Theorem \ref{DisjCycDec}.
\end{proof}

\begin{dfn}[\blue{\index{Cycle type}{Cycle type}}]
Let $\al\in S_n$ be a permutation with disjoint-cycle decomposition $\al=(x_1,{\cdots},x_r)(y_1,{\cdots},y_s){\cdots}(z_1,{\cdots},z_t)$. Then the \ul{cycle type} of $\al$ is the tuple of cycle lengths $(r,s,...,t)$.
\end{dfn}

\begin{crl}[\blue{\index{}{Order of a permutation}}]
The order of a permutation $\al\in S_n$ is the lcm of its cycle type. That is, if the disjoint-cycle decomposition of $\al$ is $\al=(x_1,...,x_{l_\al(x)})(y_1,...,y_{l_\al(y)})(z_1,...,z_{l_\al(z)})\cdots$, then
\[
|\al|=|(x_1,...,x_{l_\al(x)})(y_1,...,y_{l_\al(y)})(z_1,...,z_{l_\al(z)})\cdots|=\txt{lcm}\big(l_\al(x),l_\al(y),l_\al(z),...\big),\nn
\]
\end{crl}
\begin{proof}
Recall that the order of a cycle is equal to its length. If $\al=ab$, where $a=(x_1,...,x_m)$ and $b=(y_1,...,y_n)$ are disjoint cycles, then $\al^k=(ab)^k=a^kb^k=\ep$ $\iff$ $a^k=\ep$ and $b^k=\ep$, $\iff$ $m$ divides $k$ and $n$ divides $k$ (i.e., $k$ is a common multiple of $m$ and $n$). Hence, the order of $\al$ (i.e., the smallest possible value of $k$ such that $\al^k=\ep$) is precisely the lcm of $m$ and $n$. The general case follows by induction.
\end{proof}

\begin{dfn}[\blue{
\index{Even permutation}{Even permutation},
\index{Odd permutation}{Odd permutation},
\index{Sign (parity) of a permutation}{Sign (parity) of a permutation}}]
A permutation $\al\in S_n$ is \ul{even} (resp. \ul{odd}) if it can be written as the product of an even (resp. odd) number of $2$-cycles.

It will be shown in Theorem \ref{EvenOddPerm} that every permutation $\al\in S_n$ is either even or odd, i.e., it cannot be both even and odd. In particular, an $l$-cycle $(x_1,...,x_l)$ is odd (even) iff $l$ is even (odd), because
\bea
\txt{number of transpositions in $(x_1,...,x_l)$}~~=~~l-1.\nn
\eea
Since the \ul{determinant} (\blue{footnote}\footnote{The determinant of a matrix is defined in Lemma \ref{MatrixDetLmm} on page \pageref{MatrixDetLmm}.}) of a transposition $\mathbb{X}_n\ra\mathbb{X}_n$, which defines a unique linear transformation $Span_\Rational(\mathbb{X}_n)\ra Span_\Rational(\mathbb{X}_n)$, is $-1$, it follows from (\ref{TransDecEq}) and (\ref{CycleDecEq}) that
\bea
\det\al=(-1)^{\txt{number of transpositions in $\al$}}.\nn
\eea
This number $\in\{+1,-1\}$ is called the \index{Sign (parity) of a permutation}{\ul{sign} or \ul{parity}} of $\al$, written ~$sign(\al):=\det\al$.
\end{dfn}

\begin{lmm}\label{odd-even-lemma}
Consider the identity $\ep\in S_n$. If $\ep=\beta_1\beta_2...\beta_r$ for $2$-cycles $\beta_1,...,\beta_r$, then $r$ is even.
\end{lmm}
\begin{proof}
Consider the determinant (\blue{footnote}\footnote{The determinant of a matrix is defined in Lemma \ref{MatrixDetLmm} on page \pageref{MatrixDetLmm}.}) on the permutations $S_n:=Sym(\mathbb{X}_n)=\{\txt{bijetions}~~\mathbb{X}_n\ra \mathbb{X}_n\}$ which clearly define unique linear transformations $Span_\Rational(\mathbb{X}_n)\ra Span_\Rational(\mathbb{X}_n)$. Then every transposition has determinant equal to $-1$. Therefore $r$ is even, since
\[
1=\det(\ep)=\det(\beta_1\beta_2...\beta_r)=\det(\beta_1)\det(\beta_2)...\det(\beta_r)=(-1)^r. \qedhere
\]
\end{proof}

\begin{thm}[\blue{Always even or always odd}]\label{EvenOddPerm}
If a permutation $\al$ can be written as a product of an even (odd) number of $2$-cycles, then every other decomposition of $\al$ into a product of $2$-cycles must have an even (odd) number of $2$-cycles. That is, if
\bea
\al=\beta_1\beta_2...\beta_r~~~~\txt{and}~~~~\al=\gamma_1\gamma_2...\gamma_s,\nn
\eea
where the $\beta$'s and $\gamma$'s are $2$-cycles, then $r$ and $s$ are both even or both odd.
\end{thm}
\begin{proof}
Since $\beta_1\beta_2...\beta_r=\gamma_1\gamma_2...\gamma_s$ $\iff$ $\ep=\beta_r\beta_{r-1}...\beta_1\gamma_1\gamma_2...\gamma_s$, it follows from Lemma \ref{odd-even-lemma} that $r+s$ must be even.  Hence $r$ and $s$ are both even or both odd.
\end{proof}

\begin{thm}[\blue{Even permutations form a group}]
The set $A_n\subset S_n$ of even permutations is a subgroup.
\end{thm}
\begin{proof}
$A_n$ is nonempty since the identity permutation can be written as $\ep=(x,y)(x,y)$ for any $x,y\in \mathbb{X}_n$. Since $S_n$ is a finite group, we only need to check that $A_n$ is closed under the composition of permutations. This closedness follows from the fact that the sum of two even numbers is an even number, i.e., if we multiply (i.e., compose) any two even permutations, we get an even permutation.
\end{proof}

\begin{dfn}[\blue{\index{Alternating Group}{Alternating group}}]
The alternating $n$-group (or \ul{alternating group} of degree $n$), denoted by $A_n$, is the subgroup of even permutations in $S_n$.
\end{dfn}

\begin{thm}
For $n>1$, the order of $A_n$ is ${n!\over 2}$.
\end{thm}
\begin{proof}
Fix any $x,y\in \mathbb{X}_n$. Let $\phi:A_n\ra \{\txt{odd $n$-permutations}\},~\al\mapsto (x,y)\al$, which is a surjection because for any odd $n$-permutation $\rho$, $\rho=(x,y)(x,y)\rho=\phi\big((x,y)\rho\big)$, and an injection because for any $\al,\beta\in A_n$, $\phi(\al)=\phi(\beta)$ $\Ra$ $(x,y)\al=(x,y)\beta$, $\Ra$ $\al=\beta$. Since $\phi$ is a bijection, the number of even permutations equals the number of odd permutations in $S_n$, and since each permutation is either odd or even, we have $|A_n|={|S_n|\over 2}={n!\over 2}$.
\end{proof}

\subsection{Multiplication of finite permutations in cycle notation}~\\~
As before, by examining the relationship between the upper and lower rows of a permutation (in the two-row notation) one observes that the permutation consists of disjoint finite length cycles, where a cycle of length $k$, or a $k$-cycle, $(x_1,x_2,...,x_k)\in S_n:=Sym(\mathbb{X}_n)$ is based on $k$ elements $x_1,...,x_k\in\mathbb{X}_n:=\{1,2,\cdots,n\}$ that inter-map repeatedly/endlessly among themselves as one zigzags between the upper and lower rows. For example
{\footnotesize\bea
\al=\left[
      \begin{array}{cccccc}
        1 & 2 & 3 & 4 & 5 & 6\\
        2 & 1 & 4 & 6 & 5 & 3
      \end{array}
    \right]=
\left[\bt[row sep=small, column sep=small]
1\ar[d] & 2\ar[d] & 3\ar[d] & 4\ar[d] & 5 & 6\ar[d] \\
2\ar[ur] & 1 & 4\ar[ur] & 6\ar[urr] & 5 & 3
\et\right]=(1,2)(3,4,6)(5).\nn
\eea}
\begin{rmks*}[\blue{Writing permutations: Useful conventions and basic facts}]~
\bit[leftmargin=0.5cm]
\item If the permutation involves only single-digit numbers, commas can be omitted without confusion, e.g.,
\bea
\al:=(1,2)(3,4,6)(5)=(12)(346)(5).\nn
\eea
\item Recall that the order of elements is important in a cycle but the elements can be cyclically shifted or translated (i.e., uniformly shifted/translated modulo the cycle length), e.g.,
\bea
(346)=(463)=(634)\neq (436).\nn
\eea
However, for convenience, we begin each cycle with its smallest element.
\item We can omit $1$-cycles (since a 1-cycle equals the identity permutation $\ep$) with the understanding that any element absent is fixed (i.e., forms a $1$-cycle), e.g.,
\bea
\al:=(12)(346)(5)=(12)(346).\nn
\eea
In this case, since the identity permutation $\ep$ fixes all elements, we can (as an exception to the rule) write it as any of the $1$-cycles, e.g, $\ep=(1)$, or $\ep=(2)$, etc, or simply as $1$ or generically as $\ep$.
\item If cycles are disjoint (i.e., have no elements in common) their product can be written in any order, e.g.,
\bea
\al:=(12)(346)=(346)(12).\nn
\eea
\eit
\end{rmks*}

\begin{rmks*}[\blue{Multiplication of permutations: \index{Reduction process}{The reduction process}}]~
We can directly multiply permutations $\al,\beta,\gamma,\cdots\in S_n$ in cycle form (through a \ul{reduction-to-disjoint-cycles process}) by noting the following:
\begin{enumerate}[leftmargin=0.7cm]
\item We are applying the product of permutations $P:=\al\beta\gamma\cdots$ on one element $x\in \{1,2,\cdots,n\}$ at a time, i.e., $P$ is a pointwise mapping. Begin with the smallest element $x_1:=1\in\{1,2,\cdots,n\}$, and repeatedly apply $P$, say $l_1$ times, until a cycle is formed as follows.
    \bea
    C_1=(x_{11},x_{12},...,x_{1l_1})=\big(P^{l_1}x_1\big)\eqv\left(x_1=x_{11}\sr{P}{\ral}x_{12}\sr{P}{\ral}\cdots \sr{P}{\ral}x_{1l_1}\sr{P}{\ral}x_{11}\right)\nn
     \eea
Then take the next smallest element $x_2\in \{1,2,\cdots,n\}$ not already in $C_1$ (i.e., $x_2:=$ smallest element of $\{1,2,\cdots,n\}\backslash C_1$) and repeat the process until another cycle is formed as follows:
\bea
C_2=(x_{21},x_{22},...,x_{2l_2})=\big(P^{l_2}x_2\big)\eqv \left(x_2=x_{21}\sr{P}{\ral}x_{22}\sr{P}{\ral}\cdots \sr{P}{\ral}x_{2l_2}\sr{P}{\ral}x_{21}\right).\nn
\eea
At the $k$th step, we take $x_k:=$ smallest element of $\{1,2,....,n\}\backslash(C_1\cup\cdots\cup C_{k-1})$ and, beginning with $x_k$, repeatedly apply $P$ until a cycle is formed as follows:
\bea
C_k=(x_{k1},x_{k2},...,x_{kl_k})=\big(P^{l_k}x_k\big)\eqv \left(x_k=x_{k1}\sr{P}{\ral}x_{k2}\sr{P}{\ral}\cdots \sr{P}{\ral}x_{kl_k}\sr{P}{\ral}x_{k1}\right).\nn
\eea
Continue until $\{1,2,\cdots,n\}$ is exhausted at the final step $k_f$. The result of the product is
\bea
P:=\al\beta\gamma\cdots=C_1C_2\cdots C_{k_f}=\big(P^{l_1}x_1\big)\big(P^{l_2}x_2\big)\cdots\big(P^{l_{k_f}}x_{k_f}\big).\nn
\eea
\item A cycle fixes any element of $\{1,2,\cdots,n\}$ not in it.
\item For any two adjacent (i.e., nearby or ``next-to-each-other'') elements in a cycle, the permutation sends the first (or left) element to the second (or right) element.
\item Composition of $m$ permutations proceeds from right to left (by exactly $m$ mapping steps before the final result on the left is written down, and then used to restart the composition again from the right) and continues repeatedly until a cycle is formed.
\item By convention, we begin each new (disjoint) cycle of the product with the smallest element not already contained in the already formed previous cycles.
\item Note that for clarity (or when necessary), we may include $1$-cycles when taking a product in cycle notation.
\end{enumerate}
For example, (i) if $\al:=(13)(27)(456)(8)$ and $\beta:=(1237)(648)(5)$, then
\bea
\al\beta=(13)(27)(456)(8)~(1237)(648)(5)=(1732)(48)(56),\nn
\eea
where $(1732)\eqv(\al\beta)^4:1\sr{\al\beta}{\ral}7\sr{\al\beta}{\ral}3\sr{\al\beta}{\ral}2\sr{\al\beta}{\ral}1$,~ $(48)\eqv(\al\beta)^2:4\sr{\al\beta}{\ral}8\sr{\al\beta}{\ral}4$,~ and
$(56)\eqv(\al\beta)^2:5\sr{\al\beta}{\ral}6\sr{\al\beta}{\ral}5$. Similarly, (ii) if $\al:=(12)(3)(45)$ and $\beta:=(153)(24)$, then
\bea
\al\beta=(12)(3)(45)~(153)(24)=(14)(253),\nn
\eea
where $(14)\eqv(\al\beta)^2:1\sr{\al\beta}{\ral}4\sr{\al\beta}{\ral}1$~ and~
$(253)\eqv(\al\beta)^3:2\sr{\al\beta}{\ral}5\sr{\al\beta}{\ral}3\sr{\al\beta}{\ral}2$.
\end{rmks*}

\section{Group Representation Theory: Group Actions and Size Estimates}
\subsection{Cayley and Lagrange theorems}~\\~
Cayley's theorem will show that any (structured) group can be fully described by viewing (or ``representing'') it as a group of symmetries of some (structured) set. Meanwhile, Lagrange's theorem underlies some estimates for the relative (i.e., reasonably restricted) number of ways of describing a given group.
\begin{thm}[\blue{\index{Cayley's theorem}{Cayley's theorem}}]
Every group $G$ is isomorphic to a group of permutations $G'$.
\end{thm}
\begin{proof}
Let $G=(G,\cdot,e)$. For some set $X$, we want to construct a subgroup $G'\subset Sym(X):=\{\txt{bijections}~f:X\ra X\}$ such that $G\cong G'$. Define $X:=G$. Next, for each $g\in G$, define $T_g:G\ra G,~x\mapsto gx$. Then each $T_g$ is bijective (hence a permutation of $X:=G$). Now consider the set $G'=\{T_g:g\in G\}$ as a multiplicative set with multiplication $G'\times G'\ra G',~(T_g,T_h)\mapsto T_gT_h:=T_g\circ T_h$ (i.e., composition of maps), which is well-defined because $T_gT_h=T_{gh}\in G'$. This is clear since $(T_gT_h)(x)=T_g(T_h(x))=T_g(hx)=ghx=T_{gh}(x)$, for all $x\in G$. Thus, $G'$ is a group with multiplication $T_gT_h:=T_{gh}$, identity element $T_e$, and inverse $T_g{}^{-1}:=T_{g^{-1}}$ for each $T_g\in G'$.

The map ~$\phi:G\ra G',~g\mapsto T_g$ is a group isomorphism (seen as follows): (\blue{footnote}\footnote{The map $\phi:G\ra G'$ (or equivalently, the group $G'$) is known as the \index{Regular! representation}{\ul{left-regular representation}} of $G$.})
\bit[leftmargin=0.8cm]
\item $\phi$ is surjective: This is clear by construction.
\item $\phi$ is injective: For any $g,h\in G$, $\phi(g)=\phi(h)$ ~$\Ra$~ $T_g=T_h$, ~$\Ra$~ $gx=hx$ ~for any $x\in G$, $\Ra$ $g=h$.
\item $\phi$ preserves multiplication: For any $g,h\in G$, $\phi(gh)=T_{gh}=T_gT_h=\phi(g)\phi(h)$. \qedhere
\eit
\end{proof}

\begin{lmm}
Let $G$ be a group and $H\subset G$ a subgroup. Then $|gH|=|H|$ for all $g\in G$.
\end{lmm}
\begin{proof}
The map ~$H\ra gH:~h\mapsto gh$~ is a bijection. Hence $H\approx gH$.
\end{proof}

\begin{thm}[\blue{\index{Lagrange's theorem}{Lagrange's theorem}}]
Let $G$ be a group, $H\subset G$ a subgroup, and $G/H:=\{gH:g\in G\}$ the set of left cosets. Then $|G|=|H|~|G/H|$. (\blue{footnote}\footnote{Recall that the cardinality $|G/H|$ is called the \ul{index} of $H$ in $G$.}).

In particular, if $G$ is finite, then the order of every subgroup (hence of every element) of $G$ divides the order of $G$.
\end{thm}
\begin{proof}
Recall that distinct cosets are disjoint (\blue{footnote}\footnote{For any $a,b\in G$, ~$a\in bH$ iff $b^{-1}a\in H$, iff $b^{-1}aH=H$, iff $aH=bH$. Thus, if $a\in gH\cap g'H$, then $gH=aH=g'H$.}). If $G$ is finite, then there are only finitely many distinct left cosets $a_1H,...,a_nH$ for some $a_1,...,a_n\in G$, where $n=|G/H|$, and so
\[
|G|=|a_1H\cup...\cup a_nH|=|a_1H|+...+|a_nH|=|H|+...+|H|=nH=|G/H|~|H|.
\]If $G$ is not necessarily finite, let $(g_tH)_{t\in T}$ be an injective indexing of the distinct cosets. Then the axiom of choice gives a map $f:T\ra G=\bigsqcup_{t\in T}g_tH,~t\mapsto f(t)\in g_tH$. Hence we get a bijection (\blue{footnote}\footnote{Observe that $\phi$ consists of disjoint/parallel bijective submaps $\phi|_{\{t\}\times H}:\{t\}\times H\ra f(t)H=g_tH$. Alternatively, for injectivity in particular, if $f(t)h=f(t')h'$ then $g_tH=f(t)hH=f(t')h'H=g_{t'}H$, and so $t=t'$.})
\[
\phi:T\times H\ra G,~(t,h)\mapsto f(t)h,
\]which implies ~$|G|=|T\times H|=|T||H|=|G/H||H|$.
\end{proof}

\begin{crl}
{\flushleft (i)} If $G$ is a finite group, the order of any $g\in G$ divides the order of $G$.
{\flushleft (ii)} If $\phi:G\ra G'$ is a group homomorphism, then ~$|G|=|\ker\phi||\im\phi|$.
\end{crl}
\begin{proof}
{\flushleft(i)} For any $g\in G$, ~$|G|=|G/\langle g\rangle||\langle g\rangle|=[G/\langle g\rangle||g|$,
and so $|g|$ divides $|G|$.
{\flushleft (ii)} By the 1st isomorphism theorem, $\im\phi\cong G/\ker\phi$, and so by Lagrange's Theorem,
\[
|G|=|\ker\phi||G/\ker\phi|=|\ker\phi||\im\phi|. \qedhere
\]
\end{proof}

\begin{crl}[\textcolor{blue}{Prime order groups}]
A group of prime order is cyclic.
\end{crl}
\begin{proof}
Let $G$ be a group such that $|G|$ is a prime. Then the only possible subgroup orders are $1=|\{e\}|$ and $|G|$. Because $1$ is not a prime, we know $G\neq\{e\}$. Let $g\in G\backslash\{e\}$. Since the order $|g|$ divides $|G|$, we get $|\langle g\rangle|=|g|=|G|$. Hence $G=\langle g\rangle$, because $G$ is finite.
\end{proof}

\begin{crl}[\textcolor{blue}{Order-idempotence}]\label{order-power-corollary}
If $G=(G,\cdot,e)$ is a finite group, then $g^{|G|}=e$ for all $g\in G$.
\end{crl}
\begin{proof}
Since $|g|$ divides $|G|$, we have $g^{|G|}=g^{|g|q}=(g^{|g|})^q=e^q=e$.
\end{proof}

\begin{crl}[\textcolor{blue}{\index{Fermat's little theorem}{Fermat's little theorem}}]
For every integer $k\geq 0$ and every prime $p\geq 1$,
\[
k^p+\Integer p=k+\Integer p.
\]
\end{crl}
\begin{proof}
By the division algorithm, $k=pq+r$, $0\leq r<p$, and so $k+\Integer p=r+\Integer p$. If $r=0$, the result is trivially true. So, let {\small $r\in\{1,2,...,p-1\}=U(p):=U(\Integer_p)$}. Then by Corollary \ref{order-power-corollary} at step (s) below,
\[
k^p+\Integer p=(pq+r)^p+\Integer p=r^p+\Integer p=r^{p-1}r+\Integer p=r^{|U(p)|}r+\Integer p\sr{(s)}{=} r+\Integer p=k+\Integer p. \qedhere
\]
\end{proof}

\begin{lmm}[\blue{Simultaneous translation and disjointness}]
Let $G$ be a group, $g\in G$ an element, and $A,B\subset G$ subsets. Then ~$A\cap B=\emptyset$ $\iff$ $(gA)\cap(gB)=\emptyset$ $\iff$ $(Ag)\cap(Bg)=\emptyset$.
\end{lmm}
\begin{proof}
For any $a\in A$ and $b\in B$, we have ~$a=b$ $\iff$ $ga=gb$ $\iff$ $ag=bg$.
\end{proof}

\begin{crl}
Let $G$ be a group and $H,K\leq G$ subgroups. If $H\cap K$ is finite, then ~$|HK|={|H||K|\over |H\cap K|}$.
\end{crl}
\begin{proof}
Since $U=H\cap K$ is a subgroup of both $H$ and $K$, we observe that (with $I,J$ possibly infinite)
\begin{align}
&\textstyle HK=\Big(\bigsqcup\limits_{i\in I}h_i(H\cap K)\Big)K=\bigcup\limits_ih_i(H\cap K)K=\bigcup\limits_ih_iK\sr{(s)}{=}\bigsqcup\limits_ih_iK=\bigsqcup\limits_i\Big(h_i\bigsqcup\limits_{j\in J}(H\cap K)k_j\Big)\nn\\
&\textstyle~~~~=\bigsqcup\limits_i\Big(\bigsqcup\limits_jh_i(H\cap K)k_j\Big)=\bigsqcup\limits_{i,j}h_i(H\cap K)k_j,\nn\\
&\textstyle~~\Ra~~|HK|=|I||K|=|H||J|=|I|~|H\cap K|~|J|={|H|\over |H\cap K|}~|H\cap K|{|K|\over |H\cap K|}={|H||K|\over |H\cap K|},\nn
\end{align}
where (i) $|I|$ is the index of $H\cap K$ in $H$, (ii) $|J|$ is the index of $H\cap K$ in $K$, (iii) step (s) holds because for all $h,h'\in H$, we have $\big(h(H\cap K)\big)\cap \big(h'(H\cap K)\big)\neq\emptyset$ iff $h(H\cap K)=h'(H\cap K)$ iff $h^{-1}h'\in H\cap K\subset K$ iff $hK=h'K$ iff $hK\cap h'K\neq\emptyset$, and finally, (iv) it is clear that the bicosets $\{h(H\cap K)k:h\in H,k\in K\}$ have the same cardinality, i.e., $|h(H\cap K)k|=|h(H\cap K)|=|(H\cap K)k|=|H\cap K|$ for all $h\in H,k\in K$.
\end{proof}

Using knowledge of the product of groups, the above result can be improved by removing the finiteness of $H\cap K$ as follows.

\begin{prp}
Let $G$ be a group and $H,K\leq G$ subgroups. Then ~{\small $|H||K|=|HK||H\cap K|$}.
\end{prp}
\begin{proof}
The fibers of the map ~$f:H\times K\ra HK,~(h,k)\mapsto hk$~ are given by
\bea
&&f^{-1}(hk):=\{(h',k')\in H\times K~|~hk=h'k'\}=\{(h',k')\in H\times K~|~h^{-1}h'=kk'^{-1}\in H\cap K\}\nn\\
&&~~~~=\{(h',k')\in H\times K~|~h'=hu,~k'=u^{-1}k,~u:=h^{-1}h'=kk'^{-1}\in H\cap K\}\nn\\
&&~~~~=\{(hu,u^{-1}k)~|~u\in H\cap K\},~~~~\txt{for each}~~~~hk\in HK.\nn
\eea
Also, (i) the fibers partition $H\times K$, i.e., ~$H\times K=\bigsqcup_{hk\in HK}f^{-1}(hk)$,~ (ii) the fibers each have the same cardinality $|H\cap K|$, and (iii) $f$ is surjective. Consider an injective indexing ${H\times K\over\sim}=\{[(h_t,k_t)]:=f^{-1}(h_tk_t)\}_{t\in T}$ of the fibers. The axiom of choice gives a map $\al:T\ra H\times K=\bigsqcup_{t\in T}[(h_t,k_t)],~t\mapsto\al(t)=(h_\al(t),k_\al(t))\in[(h_t,k_t)]$. We get a bijection of sets
\[
\beta:T\times(H\cap K)\ra H\times K,~(t,u)\mapsto (h_\al(t)u,u^{-1}k_\al(t)).\nn
\]
Hence ~$|H\times K|=|T||H\cap K|=|HK||H\cap K|$.
\end{proof}

\subsection{Group symmetry operations}~\\~
The notion of a $k$-space (i.e., a vector space over a field $k$) generalizes in more than one direction, one of which is the notion of an $R$-module (i.e., a module over a ring $R$). Another generalization of a $k$-space is a \emph{G-set} (i.e., a set with an \emph{action} of a group $G$ on it), which thus automatically inherits some of the concepts for a $k$-space.

Accordingly, given a group $G$, we have the category of G-sets, whose objects are $G$-sets and whose morphisms are homomorphisms of G-sets. These concepts and more will be defined in this section. Note that the $R$-module retains the abelian group structure of the $k$-space, but the $G$-set does not (in general). On the other hand, the G-set retains the invertibility structure of the field $k$ of a $k$-space, but the $R$-module does not (except when $R$ happens to be a division ring).

Objects that have both an $R$-module structure and a $G$-set structure also exist. We will see an example of such an object known as an $RG$-module, where $RG$ denotes a ring (called \emph{group $R$-ring} or \emph{group $R$-algebra}) that is also a free $R$-module containing $G$ as a basis.

\begin{dfn}[\blue{
\index{Group! action}{Group action},
\index{G-set}{G-set},
\index{Orbit}{Orbit},
\index{Transitive! action}{Transitive action},
\index{Transitive! G-set}{Transitive G-set},
\index{Similar elements}{Similar elements},
\index{Statilizer}{Statilizer}}]
Let $G=(G,\cdot,e)=(G,\cdot,1)$ be a group and $S$ a set. A \ul{left $G$-action} or \ul{left action of $G$} on $S$ (making $S$ a left G-set) is a map ~$G\times S\ra S$, $(g,s)\mapsto gs$~ such that (i) $1s=s$ and (ii) $g(g's)=(gg')s$, for all $s\in S$, $g,g'\in G$.  A \ul{right $G$-action} or \ul{right action of $G$} on $S$ is similarly a map
\bea
S\times G\ra S,~~(s,g)\mapsto sg,~~\txt{with}~~s1=s,~~(sg)g'=s(gg')~~\txt{for all}~~s\in S,~g,g'\in G.\nn
\eea

If $s\in S$, the \ul{orbit} of $s$ in $S$ is $\sigma_s:=\{gs:g\in G\}=Gs$. The group $G$ acts \ul{transitively} on the G-set $S$ (making $S$ a \ul{transitive G-set}) if $S$ is a single orbit under $G$, i.e., if $S=Gs$ for some $s\in S$. Two elements $s,s'\in S$ are \ul{similar} (written $s\sim s'$) if $gs=g's'$ for some $g,g'\in G$ (which is the case iff $s',s\in\sigma_s=\sigma_{s'}$, i.e., $s,s'$ share the same orbit). The \ul{stabilizer} of an element $s\in S$ is the subgroup ~$G_s:=\{g\in G:gs=s\}\leq G$, i.e., all elements of $G$ that fix $s$.
\end{dfn}
Henceforth, ``$G$-action'' (resp. ``$G$-set'') will mean ``left $G$-action'' (resp. `` left $G$-set'') unless it is stated otherwise. Also, because of remark (i) below, a left-right distinction is not even necessary (although it is necessary in the case of a ring $R$ and its $R$-modules or $R$-module actions).

\begin{rmks}
Let $S$ be a left $G$-set.
{\flushleft(i)} If $S$ is a left G-set, then $S$ is also a right G-set with respect to the operation
\bea
S\times G\ra S,~(s,g)\mapsto sg:=g^{-1}s.\nn
\eea
{\flushleft(ii)} The power set $\P(S)$ is a left $G$-set under the operation $G\times\P(S)\ra\P(S),~(g,A)\mapsto gA:=\{gs:s\in A\}$. The orbit of $A\subset S$ in $\P(S)$ is $\sigma_A:=\{gA:g\in G\}\subset\P(S)$.
{\flushleft(iii)} For any $s\in S$, the orbit $\sigma_s=Gs\subset S$ is a G-set, with inherited G-action
\bea
G\times \sigma_s\ra \sigma_s,~(g,g's)\mapsto g(g's)=(gg')s.\nn
\eea
This is because $g\sigma_s=gGs=Gs=\sigma_s$ for all $g\in G$. In this case, the orbit of any $gs\in \sigma_s$ is
\bea
\sigma_{gs}=Ggs=Gs=\sigma_s~~~~\txt{(i.e., $G$ acts transitively on an orbit)}\nn
\eea
{\flushleft(iv)} A G-set is a disjoint union of orbits (\blue{footnote}\footnote{This is because similarity $\sim$ is an equivalence relation on $S$, and the equivalence class of $s\in S$ is $[s]_\sim=\sigma_s$.}): I.e., there is $T\subset S$ such that {\small $S=\bigsqcup_{t\in T}\sigma_t=\bigsqcup_{t\in T}Gt$}.
{\flushleft(v)} The stabilizer $G_s\subset G$ (of any $s\in S$) is indeed a subgroup: To see this, observe that for all $a,b\in G_s$,
\bea
as=s=bs~~\Ra~~ab^{-1}s=as=s,~~\Ra~~ab^{-1}\in G_s.\nn
\eea
\end{rmks}

\begin{prp}[\blue{G-sets are ``permutation representations'' of $G$}]
Let $S$ be a G-set and $m:G\times S\ra S$ the action of $G$ on $S$. Then for each $g\in G$, the map $m_g:S\ra S,~s\mapsto gs$ is a permutation of $S$. That is, $G\subset Sym(S):=\{\txt{bijections}~f:S\ra S\}$.

Hence, we get a unique group homomorphism $m_\ast:G\ra Sym(S),~g\mapsto m_g$ (called a \ul{permutation representation} of $G$), which in turn uniquely defines the given action as ~$m:G\times S\ra S,~(g,s)\mapsto m_\ast(g)(s)=m_g(s)=gs$.
\end{prp}
\begin{proof}
Let $g\in G$. The map $m_g$ is surjective because any $s\in S$ can be written as $s=gg^{-1}s=m_g(g^{-1}s)$. Also, $m_g$ is injective because for any $s,s'\in S$, $gs=gs'$ $\Ra$ $s=s'$.
\end{proof}

\begin{thm}[\blue{\index{Conjugate-stabilizer theorem}{Conjugate-stabilizer theorem}}]\label{conjugate-stabilizers-theorem}
Let $S$ be a G-set and $s,s'\in G$. If $s\sim s'$ (i.e., $s,s'$ are in the same $G$-orbit), then $G_{s'}=gG_sg^{-1}$ for some $g\in G$.
\end{thm}
\begin{proof}
Since $s\sim s'$, we have $s'=gs$ for some $g\in G$. Therefore,
\begin{align}
&G_{s'}=G_{gs}=\{g'\in G:~g'gs=gs\}=\{g'\in G:~g^{-1}g'gs=s\}=\{g'\in G:~g^{-1}g'g\in G_s\}\nn\\
&~~~~=\{g'\in G:~g'\in gG_sg^{-1}\}=gG_sg^{-1}. \qedhere
\end{align}
\end{proof}

\begin{dfn}[\blue{
\index{Homomorphism of G-sets}{Homomorphism of G-sets},
\index{Endomorphism! of a G-set}{Endomorphism of a G-set},
\index{Monomorphism! of G-sets}{Monomorphism of G-sets},
\index{Epimorphism of G-sets}{Epimorphism of G-sets},
\index{Isomorphism of! G-sets}{Isomorphism of G-sets}}]
Given G-sets $S,S'$, a map $\vphi:S\ra S'$ is a \ul{homomorphism} of G-sets (\ul{$G$-homomorphism} or \ul{$G$-invariant map}) if ~$\vphi(gs)=g\vphi(s)$~ for all $s\in S,g\in G$.

An \ul{endomorphism} of the G-set $S$ is a homomorphism of $G$-sets $\vphi:S\ra S$. A \ul{monomorphism} (resp. \ul{epimorphism}) of G-sets is an injective (resp. surjective) homomorphism of G-sets. An \ul{isomorphism} of G-sets is a bijective homomorphism of G-sets.
\end{dfn}

\begin{thm}[\blue{\index{Isomorphism theorem for G-sets}{Isomorphism theorem for G-sets}: $G/G_s\cong Gs$}]\label{G-set-iso-thm}
Let $S$ be a G-set, $s\in S$, and $G/G_s:=\{gG_s:g\in G\}$. Then (i) $G/G_s$ is a transitive G-set wrt the action $G\times G/G_s\ra G/G_s$, $(g,g's)\mapsto gg's$, and (ii) the map $\vphi:G\ra\sigma_s$, $g\mapsto gs$ is a homomorphism of G-sets that uniquely factors as
\bea
\vphi=\ol{\vphi}\circ\pi:G\sr{\pi}{\ral}G/G_s\sr{\ol{\vphi}}{\ral}\sigma_s,~~~~\Ra~~G/G_s\cong\sigma_s,\nn
\eea
where the canonical map $\pi:G\ra G/G_s,~g\ra gG_s$ is a homomorphism of G-sets, and the map $\ol{\vphi}:G/G_s\ra \sigma_s$, $gG_s\mapsto gs$ is an isomorphism of G-sets.
\begin{center}
\adjustbox{scale=0.9}{\begin{tikzcd}
G \ar[dr,"\pi"] \ar[rr,"\vphi=\ol{\vphi}\pi"] & & \sigma_s \\
 &  G/G_s \ar[ur,"\ol{\vphi}"]&
\end{tikzcd}}
\end{center}
\end{thm}
\begin{proof}
{\flushleft (i)} It is clear that $G/G_s$ is a G-set. Also, for any $gG_s,g'G_s\in G/G_s$, we have $g'G_s=(g'g^{-1})gG_s$, and so $G/G_s$ is a transitive G-set.
{\flushleft (ii)} It is clear that $\vphi=\ol{\vphi}\circ\pi$, and all three maps $\vphi,\ol{\vphi},\pi$ are well-defined homomorphisms of G-sets. Also, if $\vphi=\ol{\vphi}\circ\pi=\phi\circ\pi$ (for any $\phi:G/G_s\ra\sigma_s$), then $\phi=\ol{\vphi}$ since $\pi$ is surjective (i.e., the factorization is unique). Finally, $\ol{\vphi}$ is a bijective homomorphism (hence an isomorphism) of G-sets, since $gG_s=g'G_s$ $\iff$ $gs=g's$ for all $g,g'\in G$, and $\ol{\vphi}$ is clearly surjective.
\end{proof}

\begin{crl}[\blue{\index{Orbit-Stabilizer theorem}{Orbit-Stabilizer theorem}}]\label{OrbitStabThm}
Let $G$ be a group and $S$ a G-set. Then for any $s\in S$,
\bea
|G|=|G/G_s||G_s|=|\sigma_s||G_s|.\nn
\eea
\end{crl}
\begin{proof}
This follows from Lagrange's theorem and the isomorphism of G-sets $G/G_s\cong\sigma_s$.
\end{proof}

\begin{crl}[\blue{\index{Counting theorem I}{Counting theorem I}}]\label{CountThmI}
Let $G$ be a group and $S$ a G-set. Then
\bea
\textstyle |S|=\left|\bigsqcup_{t\in T}\sigma_t\right|\geq\sum_{t\in T}|\sigma_t|=\sum_{t\in T}|G/G_t|,~~~~\txt{for some}~~T\subset S,\nn
\eea
where $\sum_{t\in T}|\sigma_t|:=\sup_{\txt{finite}~F\subset T}\sum_{t\in F}|\sigma_t|$. Furthermore, if (i) $T$ is finite or (ii) $S\subset G$ and the orbits are right cosets $\{\sigma_t=Ht\}_{t\in T}$ of some subgroup $H\leq G$ or (iii) $\sum_{t\in T}|\sigma_t|=\infty$, then
\[
\textstyle |S|=\sum_{t\in T}|\sigma_t|=\sum_{t\in T}|G/G_t|~~~~~~(\txt{which holds if $G$ is finite}). \qedhere
\]
\end{crl}
\begin{proof}
(i) If $T$ is finite, equality holds by the definition of cardinality. (ii) If $T$ is finite, equality holds by (i). So, assume $|T|=\infty$. The axiom of choice gives a map $f:T\ra\bigsqcup_{t\in T}\sigma_t=S,~t\mapsto f(t)\in\sigma_t=Ht$, which in turn gives a bijection $\phi:H\times T\ra S,~(h,t)\mapsto hf(t)$. (iii) Since
\[
\textstyle|S|=\left|\bigsqcup_{t\in T}\sigma_t\right|\geq\sup_{\txt{finite}~F\subset T}\left|\bigsqcup_{t\in F}\sigma_t\right|=\sup_{\txt{finite}~F\subset T}\sum_{t\in F}|\sigma_t|=\sum_{t\in T}|\sigma_t|,
\]
equality holds if $\sum_{t\in T}|\sigma_t|=\infty$.

Finally, assume {\small $|T|=\infty$} and {\small $\sum_{t\in T}|\sigma_t|<\infty$}. In this case, the inequality still holds by the proof of (iii) above.
\end{proof}

\subsection{Group symmetry operations on subsets of a G-set}~\\~
Recall that if $S$ is a G-set, then the power set $\P(S)$ is also a G-set, and the action of $G$ on $\P(S)$ is
\bea
G\times\P(S)\ra\P(S),~(g,U)\mapsto gU:=\{gu:~u\in U\}.\nn
\eea
Two subsets $U,U'\in\P(S)$ of $S$ are \ul{similar} in $\P(S)$, written $U\sim U'$, if $gU=g'U'$ for some $g,g'\in G$ (or equivalently, $U=gU'$ for some $g\in G$). (\blue{footnote}\footnote{As usual, $U=gU'$ means $U\subset gU'$ and $gU'\subset U$.}). As before, the \ul{orbit} of $U\subset S$ in $\P(S)$ is given by
\bea
\sigma_U:=\{gU:~g\in G\}=\{V\in\P(S):~V=gU~\txt{for some}~g\in G\},\nn
\eea
Similarly, the \ul{stabilizer} of $U\in\P(S)$ in $G$ is the subgroup of $G$ consisting of permutations of $U$ as follows:
\bea
G_U:=\{g\in G:~gU=U\}=G\cap Sym(U),~~~~Sym(U):=\{\txt{bijections ~$f:U\ra U$}\}.\nn
\eea

\begin{prp}[\blue{\index{G-subset criterion}{G-subset criterion}}]\label{SubGsetCrit}
Let $S$ be a G-set and $U\subset S$. Then $U$ is a $G$-set (i.e., $G_U=G$) $\iff$ $U$ is a union of orbits in $S$.
\end{prp}
\begin{proof}
{\flushleft ($\Ra$)} Assume $G_U=G$. Since $G_U=\{g\in G:gU=U\}$, we have $gU\subset U$ for all $g\in G$, and so $U$ is a G-set with G-action $G\times U\ra U,~(g,u)\mapsto gu$ (i.e., the restriction of the action $G\times S\ra S$ to $G\times U$). Hence $U$ is a union of orbits in $S$.

{\flushleft ($\La$)} Assume $U$ is a union of orbits in $S$, i.e., $U=\bigcup_{u\in V}\sigma_u$ for some $V\subset U\cap S=U$. Then for all $g\in G$,
\bea
\textstyle gU=g\bigcup_{u\in V}\sigma_u=\bigcup_{u\in V}g\sigma_u=\bigcup_{u\in V}\sigma_u=U.\nn
\eea
Hence ~$G_U=\{g\in G:~gU=U\}=G$.
\end{proof}

\begin{prp}[\blue{\index{Counting theorem II}{Counting theorem II}}]\label{left-G-subset-prop}
Let $G$ be a group and $U\subset G$. Then (i) $G$ is a $G$-set under left multiplication $G\times G\ra G$, $(g,a)\mapsto ga$ in $G$, and (ii) $|U|=|T||G_U|$ for some subset $T\subset G$, where $G_U:=\{g\in G:gU=U\}$ is the stabilizer of $U$ wrt the $G$-set $\P(G)$.
\end{prp}
\begin{proof}
Part (i) is clear, and so we will prove part (ii). Since the stabilizer $G_U$ of $U$ is a group, it follows that $G$ (hence $U$) is a $G_U$-set with action
\bea
G_U\times G\ra G,~(g,s)\mapsto gs,~~\Ra~~G_U\times U\ra U,~(g,u)\mapsto gu.\nn
\eea
Hence $U$ is a union of $G_U$-orbits in $G$ (which are right cosets $\{\sigma_g^{G_U}=G_Ug:g\in U\}$ of $G_U$ in $G$), i.e.,
\bea
\textstyle U=\bigcup_{g\in U}\sigma^{G_U}_g=\bigcup_{g\in U}G_Ug.\nn
\eea
Hence, with a set of orbit (or right coset) representatives $T\subset U$, Theorem \ref{CountThmI} gives
\[
\textstyle |U|=\sum_{g\in T}|\sigma^{G_U}_g|=\sum_{g\in T}|G_Ug|=\sum_{g\in T}|G_U|=|T||G_U|. \qedhere
\]
\end{proof}

\begin{dfn}[\blue{\index{Conjugation action}{Conjugation action},
\index{Conjugacy classes}{Conjugacy classes},
\index{Centralizer subgroups}{Centralizer subgroups},
\index{Conjugate subgroups}{Conjugate subgroups},
\index{Normalizer of a subgroup}{Normalizer of a subgroup}}]
The conjugation (or \ul{conjugation action}) of a group $G$ on itself is the map $G\times G\ra G,~(g,a)\mapsto gag^{-1}$, with orbits (called \ul{conjugacy classes}) $[g]_{cc}:=\sigma_g=\{aga^{-1}:a\in G\}$ and stabilizers (called \ul{centralizer subgroups} or just \ul{centralizers})
\bea
C(g):=G_g=\{a\in G:aga^{-1}=g\}=\{a\in G:ag=ga\}.~(\blue{footnote}\footnotemark).\nn
\eea
\footnotetext{By the isomorphism of G-sets $G/G_g\cong\sigma_g$, the orbit-stabilizer theorem implies ~$|G|=\big|[g]_{cc}\big|~\big|C(g)\big|$.}This induces the conjugation action of $G$ on its subgroups ~$Sub(G):=\big\{H\in\P(G):H\leq G\big\}$~ given by
\bea
G\times Sub(G)\ra Sub(G),~~(g,H)\mapsto gHg^{-1}.\nn
\eea
The conjugation-orbit~ $\sigma_H=\{gHg^{-1}:g\in G\}$~ of a subgroup $H$ in $Sub(G)$ is the set of \ul{conjugate subgroups} to $H$. (\blue{footnote}\footnote{Observe that $\sigma_H=\{H\}$ $\iff$ $G_H=G$ $\iff$ $H\vartriangleleft G$ (i.e., $\iff$ $H\subset G$ is a normal subgroup).}). The conjugation-stabilizer $G_H$ of a subgroup $H\in Sub(G)$ is called the \ul{normalizer} of $H$, and is denoted by $N(H)$, i.e., ~$N(H):= G_H=\{g\in G:~gHg^{-1}=H\}$.
\end{dfn}

\begin{thm}[\blue{Properties of the normalizer}]
Let $G$ be a group and $H\leq G$ a subgroup. The normalizer subgroup $N(H)\leq G$ has the following properties.
\begin{enumerate}[leftmargin=0.9cm]
\item[(1)]   $H\vartriangleleft N(H)$, i.e., $H$ is a normal subgroup of $N(H)$.
\item[(2)] $N(H)$ is the largest subgroup of $G$ containing $H$ as a normal subgroup.
\item[(3)] $H\vartriangleleft G$ if and only if $N(H)=G$.
\end{enumerate}
\end{thm}
\begin{proof}
{\flushleft(1)} It is clear that $H\leq N(H)$. For all $g\in N(H)$, $gHg^{-1}=H$ and so $H\vartriangleleft N(H)$.
{\flushleft(2)} Let $H\vartriangleleft K\leq G$, and let $k\in K$. Then $kHk^{-1}=H$ implies $k\in N(H)$. Hence $K\leq N(H)$.
{\flushleft(3)} This is clear by the definition of a normal subgroup.
\end{proof}

\begin{prp}[\blue{Conjugate subgroup count}]
Let $G$ be a group and $H\leq G$  a subgroup. Then
\bea
\big|\{\txt{conjugate subgroups to $H\leq G$}\}\big|=|G/N(H)|.\nn
\eea
\end{prp}
\begin{proof}
Consider the conjugation action of $G$ on its subgroups Sub(G). Then we have the isomorphism of G-sets $\sigma_H\cong G/G_H$, where $\sigma_H=\{gHg^{-1}:g\in G\}$ consists of the conjugate subgroups to $H$, and $G_H=N(H)$. Hence the result follows.
\end{proof}

\subsection{Number of orbits in a G-set and Burnside's theorem}~\\~
We have seen how to determine the size of a $G$-orbit and of a G-set, but not how many orbits there are in an independent way (i.e., without using knowledge of sizes for the orbits and the G-set).

\begin{dfn}[\textcolor{blue}{\index{Pointwise invariant domain}{Pointwise (point) invariant domain} of a group element}]
Let $S$ be a G-set and $g\in G$. The \ul{point-invariant domain} (or \ul{pointwise-invariant domain}) of $g$ in $S$ is the set $S_g=\{s\in S:~gs=s\}$.
\end{dfn}

\begin{thm}[\textcolor{blue}{\index{Burnside's theorem}{Burnside's theorem}}]
Let $S$ be a G-set and $T\subset S$ a full set of orbit representatives (i.e., $S=\bigsqcup_{t\in T}\sigma_t$ as a disjoint union of $G$-orbits). For each $g\in G$, let $S_g:=\{s\in S:gs=s\}$ be the point-invariant domain of $g$ in $S$. Then, with suitable summing, the number $|T|$ of $G$-orbits in $S$ satisfies
\[
\textstyle |T||G|=\sum_{s\in S}|G_s|=\sum_{g\in G}|S_g|.
\]
In particular, if $G$ is finite then ~$|T|={1\over|G|}\sum_{s\in S}|G_s|={1\over|G|}\sum_{g\in G}|S_g|$.
\end{thm}

\begin{proof}
Consider the set $A:=\{(g,s)\in G\times S:~gs=s\}=\bigcup_{g\in G}A^g=\bigcup_{s\in S}A_s$, where $A^g:=\{(g,s): s\in S:~gs=s\}$ and $A_s:=\{(g,s): g\in G:~gs=s\}$. Since the above unions are disjoint unions, we have ~$\sum_{g\in G}|A^g|=|A|=\sum_{s\in S}|A_s|$~ under suitable conditions (e.g., cardinality restrictions and appropriate specifications/realizations for the sums involved), and so
{\footnotesize\begin{align}
\textstyle\sum\limits_{g\in G}|S_g|
=\sum\limits_{s\in S}|G_s|
=\sum\limits_{s\in \bigsqcup_{t\in T\subset S}\sigma_t}|G_s|
=\sum\limits_{t\in T\subset S}\sum\limits_{s\in\sigma_t}|G_s|
\sr{(a)}{=}\sum\limits_{t\in T}\sum\limits_{s\in\sigma_t}|G_t|
=\sum\limits_{t\in T}|\sigma_t||G_t|
\sr{(b)}{=}\sum\limits_{t\in T}|G|
=|T||G|,\nn
\end{align}}where step (a) is Theorem \ref{conjugate-stabilizers-theorem} ($s\in\sigma_t$ iff $\sigma_s=\sigma_t$, and so $G_s=gG_tg^{-1}$ for some $g\in G$), and step (b) is Corollary \ref{OrbitStabThm} (orbit-stabilizer theorem).
\end{proof}

\begin{dfn}[\blue{
\index{Inner! automorphisms}{Inner automorphisms},
\index{Inner! automorphism group}{Inner automorphism group}}]
Let $G$ be a group. Recall the conjugation action $G\times G\ra G,~(g,a)\mapsto gag^{-1}$. The associated permutations of $G$, i.e., $\phi_g:G\ra G,~a\mapsto gag^{-1}$ for each $g\in G$, are called the \ul{inner automorphisms} of $G$, and the \ul{inner automorphism group} of $G$ is the permutation subgroup ~$Inn(G):=\{\phi_g:g\in G\}\subset Sym(G):=\{\txt{bijections $f:G\ra G$}\}$.
\end{dfn}

\begin{lmm}[\blue{Character of inner automorphism group}]
For any group $G$, we have~ $Inn(G)\cong{G\over Z(G)}$,~ where $Z(G)$ is the center of $G$.
\end{lmm}
\begin{proof}
Define the surjective map ~$f:G\ra Inn(G),~g\mapsto\phi_g$. Then ~$f(gg')=\phi_{gg'}=\phi_g\phi_{g'}=f(g)f(g')$.~ Also, ~$\ker f=\{g\in G:\phi_g=e_{Inn(G)}\}=\{g\in G:gag^{-1}=a~\txt{for all}~a\in G\}=Z(G)$. Hence, by the 1st isomorphism theorem, ~${G\over Z(G)}={G\over\ker f}\cong\im f=Inn(G)$.
\end{proof}

\begin{thm}[\textcolor{blue}{Number of conjugacy classes}]
Let $G$ be a group. If $g\in G$, let $[g]_{cc}:=\{aga^{-1}:a\in G\}$ be the conjugacy class at $g$, and $C(g):=\{a\in G:ag=ga\}$ the centralizer subgroup at $g$. Then the number $Ncc(G)$ of conjugacy classes of $G$ satisfies ~{\footnotesize $Ncc(G)~|G/Z(G)|=\sum_{aZ(G)\in G/Z(G)}|C(a)|$}.~ In particular, if $G/Z(G)$ is finite, then  ~{\footnotesize $Ncc(G)={1\over |G/Z(G)|}\sum_{aZ(G)\in G/Z(G)}|C(a)|$}.
\end{thm}
\begin{proof}
Observe that $S:=G$ is a $Inn(G)$-set, where $\txt{Inn}(G)=\{\phi_g:g\in G\}$ acts as
\bea
Inn(G)\times G\ra G,~(\phi_g,s)\mapsto \phi_g(s)=gsg^{-1}.\nn
\eea
The $Inn(G)$-orbits are the conjugacy classes $[s]_{cc}:=\sigma_s=\{gsg^{-1}:g\in G\}$, and for each $\phi_g\in Inn(G)$, the point-invariant domain of $\phi_g$ in $S:=G$ is $S_{\phi_g}=\{s\in G:gsg^{-1}=s\}=C(g)$, i.e., the point-invariant domains of elements of Inn(G) are the centralizer subgroups of $G$. Hence, by Burnsides' theorem (and the fact that $Inn(G)\cong G/Z(G)$) the number of Inn(G)-orbits (hence the number $Ncc(G)$ of conjugacy classes of $G$) satisfies ~{\small $ Ncc(G)~|Inn(G)|=\sum_{\phi_g\in Inn(G)}|S_{\phi_g}|$}~ or ~{\small $Ncc(G)~|G/Z(G)|=\sum_{gZ(G)\in G/Z(G)}|C(g)|$}.
\end{proof}

\vspace{0.2cm}
\section{Systematic Representation Theory}\label{SystRepsThry}
We will now begin a general discussion of representation theory that simultaneously involves groups, rings, modules, and algebras.
\subsection{Representations of objects in categories}
\begin{dfn}[\textcolor{blue}{
\index{Endomorphism! class}{Endomorphism class},
\index{Representation}{Representation},
\index{Representable! object}{Representable object},
\index{Representable! system}{Representable system}}]
Let $\C,\D$ be categories and $C,C'\in\Ob\C$ objects of $\C$. The \ul{endomorphism class} of $C$ is $End_\C(C):=\Mor_\C(C,C)$. (\blue{footnote}\footnote{The endomorphism class $End_\C(C)$ is a \ul{multiplicative class} with respect to composition of morphisms.}). By the construction of a category, $C=C'$ $\iff$ $End_\C(C)=End_\C(C')$, i.e., objects in $\C$ are uniquely determined by their endomorphisms collections, and so $C$ is uniquely determined by the subcategory $\C(C)\subset\C$ with a single object $C$ and morphisms {\footnotesize $\Mor_{\C(C)}(C,C):=\Mor_\C(C,C)$}. A \ul{$\D$-representation} of $C$ (or a \ul{representation of $C$ in $\D$}) is an appropriate nontrivial functor of the form (\blue{footnote}\footnote{Here, if we view $End_\C(C)$ as an object in the category of multiplicative classes, then the phrases ``a $\D$-representation of $C$'' and ``a $\D$-representation of $End_\C(C)$'' should mean the same thing after the phrase ``a $\D$-representation of a multiplicative class'' has been suitably defined.})
\bea
F:\E\subset\C(C)\ra\D,~~C\sr{f}{\ral}C~~\longmapsto~~F(C)\sr{F(f)}{\ral}F(C).\nn
\eea
If such a $\D$-representation of $C$ exists, then we say $C$ is \ul{$\D$-representable} (or \ul{representable in $\D$}). An object of $\C$ is \ul{representable} if it is $Sets$-representable (i.e., representable in Sets).

A system $S:\I\ra\C$ is a \ul{$\D$-representable system} if it is $\D$-representable as an object in the category of functors $\C^\I$. (\blue{footnote}\footnote{The notion of a representable system here (although related) is independent of the earlier notion of a representable functor.})
\end{dfn}

We will now consider explicit constructions of important special cases of the above generic situation. Parts of the following definition are recalls from earlier definitions.

\begin{dfn}[\blue{
\index{Homomorphism! of multiplicative sets}{Homomorphism of multiplicative sets},
\index{Category of! multiplicative sets}{Category of multiplicative sets},
\index{Left$\slash$Right-cancelling element}{Left$\slash$Right cancelling element},
\index{Left$\slash$Right-representable multiplicative set}{Left$\slash$Right-representable multiplicative set}}]
Let $S=(S,\cdot)$ and $T=(T,\cdot)$ be multiplicative sets. A map $\rho:S\ra T$ is a \ul{homomorphism of multiplicative sets} if $\rho(ss')=\rho(s)\rho(s')$ for all $s,s'\in S$. The \ul{category of multiplicative sets} (\ul{MultSets}) is the category with multiplicative sets as objects and homomorphisms of multiplicative sets as morphisms.

An element $l\in S$ is \ul{left-cancelling} (making $S=(S,l)$ \ul{right-representable}) if for any $s,s'\in S$, $ls=ls'$ $\Ra$ $s=s'$. An element $r\in S$ is \ul{right-cancelling} (making $S=(S,r)$ \ul{left-representable}) if for any $s,s'\in S$, $sr=s'r$ $\Ra$ $s=s'$.

\end{dfn}
Observe that by construction, the elements $s\in S$ of a right-representable (resp. left-representable)  multiplicative set $S$ can be uniquely viewed as left-multiplication maps $s:S\ra S$, $x\mapsto sx$ (resp. right-multiplication maps $s:S\ra S$, $x\mapsto xs$). It is clear that a group (or any identity set) is a representable multiplicative set, and that a ring/algebra is representable wrt multiplication in the ring/algebra.

\begin{dfn}[\blue{
\index{Representable! multiplicative set}{Representable multiplicative set as a category},
\index{Representation}{Representation},
\index{Linear! representation}{Linear representation}}]
Let $\D$ be a locally small category. It is clear that for any object $D\in\Ob\D$, the set $End_\D(D):=\Mor_\D\big(D,D\big)$ is a multiplicative set with respect to composition of morphisms.

Let $S$ be a left-representable multiplicative set. By construction, $S$ can be viewed as the category $\C(S)$ with a single object $S$ and morphisms being the elements $s\in S$ as left-multiplication maps $s:S\ra S$, $x\mapsto sx$. A \ul{representation} of $S$ in $\D$ (or a $\D$-representation of $S$) is a functor (\blue{footnote}\footnote{If there is a category $\I$ and an object $I\in\Ob\I$ such that $S=End_\I(I)$, then a $\D$-representation of $S$ (as a multiplicative set) is equivalently a $\D$-representation of $I$ (as an object in $\I$).})
\bea
F:\C(S)\ra\D,~~S\sr{s}{\ral}S~~\mapsto~~F(S)\sr{F(s)}{\ral}F(S),\nn
\eea
that is, we have a homomorphism of multiplicative sets (\blue{footnote}\footnote{That is, given a left-representable multiplicative set $S$ and a locally small category $\D$, a \ul{representation} of $S$ in $\D$ (or a $\D$-representation of $S$) is a homomorphism of multiplicative sets ~$\rho:S\ra End_\D(D)$~ for an object $D\in\Ob\D$.}) ~$F:S\ra End_\D(D)$,~
for some object $D:=F(S)\in\Ob\D$. If $\D$ is the category of $k$-spaces (or $k$-modules for a commutative ring $k$), the representation is a \ul{$k$-linear representation} of $S$.
\end{dfn}

\begin{dfn}[\blue{
\index{Automorphism group}{Automorphism group},
\index{Group! as a category}{Group as a category},
\index{Group! representation}{Group representation},
\index{Permutation (symmetry) group}{Permutation (symmetry) group},
\index{Permutation representation}{Permutation representation},
\index{General linear group}{General linear group},
\index{Linear! group representation}{Linear group representation}}]~\\~
Let $\D$ be a locally small category. For any object $D\in\Ob\D$, the \ul{automorphism group} of $D$ is the set
\bea
Aut_\D(D):=\left\{\txt{isomorphisms in}~\Mor_\D\big(D,D\big)\right\}=\left\{\txt{invertibles in}~End_\D(D):=\Mor_\D\big(D,D\big)\right\}\nn
\eea
as a group with multiplication given by composition of morphisms and identity element $e:=id_D$.

Let $G$ be a group. By construction, $G$ can be viewed as the category $\C(G)$ with a single object $G$ and morphisms being the group elements $g\in G$ as left multiplication maps $g:G\ra G$, $a\mapsto ga$. A \ul{representation} of $G$ in $\D$ (or a $\D$-representation of $G$) is a functor
\bea
F:\C(G)\ra\D,~~G\sr{g}{\ral}G~~\mapsto~~F(G)\sr{F(g)}{\ral}F(G),\nn
\eea
that is, i.e., we have a group homomorphism (\blue{footnote}\footnote{That is, given a group $G$ and a locally small category $\D$, a \ul{representation} of $G$ in $\D$ (or a $\D$-representation of $G$) is a group homomorphism ~$\rho:G\ra Aut_\D(D)$~ for an object $D\in\Ob\D$.})
\bea
F:G\ra Aut_\D\big(D\big),\nn
\eea
for some object $D:=F(G)\in\Ob\D$. If $\D=Sets$, then $Aut_\D(D)=Sym(D):=\{\txt{bijections}~D\ra D\}$ is called the \ul{permutation (symmetry) group} of $D$ and, the representation is called a \ul{permutation representation} of $G$. If $\D$ is the category of $k$-spaces (or $k$-modules for a commutative ring $k$), then $GL_k(D):=Aut_k(D)$ is called the \ul{general linear group} of $D$, and the representation is called a \ul{$k$-linear representation} of $G$.
\end{dfn}

\begin{dfn}[\blue{
\index{Endomorphism! ring}{Endomorphism ring},
\index{Ring! as a category}{Ring as a category},
\index{Ring! representation}{Ring representation},
\index{Linear! ring representation}{Linear ring representation}}]
Let $\A$ be a \magenta{pre-additive} category. For any object $A\in\Ob\A$, the \ul{endomorphism ring} of $A$ is the abelian group ~$End_\A(A):=\Mor_\A\big(A,A\big)$~ as a ring with multiplication given by composition of morphisms and unity $1_{End(A)}:=id_A$.

Let $R$ be a ring. By construction, $R$ can be viewed as the small \magenta{pre-additive} category $\C(R)$ with a single object $R$ and morphisms being the ring elements $r\in R$ as left multiplication maps $r:R\ra R$, $x\mapsto rx$. A \ul{representation} of $R$ in $\A$ (or an $\A$-representation of $R$) is an \magenta{additive} functor
\bea
F:\C(R)\ra\A,~~R\sr{r}{\ral}R~~\mapsto~~F(R)\sr{F(r)}{\ral}F(R),\nn
\eea
that is, we have a ring homomorphism (\blue{footnote}\footnote{That is, given a ring $R$ and a pre-additive category $\A$, a \ul{representation} of $R$ in $\A$ (or an $\A$-representation of $R$) is a ring homomorphism ~$\rho:R\ra End_\A(A)$~ for an object $A\in\Ob\A$.}) ~$F:R\ra End_\A\big(A\big)$,~ for some object $A:=F(R)\in\Ob\C$.

If $\A$ is the category of $k$-spaces, the representation is a \ul{$k$-linear representation} of $R$.
\end{dfn}

\begin{dfn}[\blue{
\index{Ring-Preadditive category}{Pre-additivity category},
\index{Ring-additive (Linear) functor}{Ring-additive (Linear) functor}}]
Let $R$ be a ring and $\B$ a category. $\B$ is \ul{(left) pre-additive over $R$} (or \ul{(left) $R$-preadditve}) if the following hold:
\bit[leftmargin=0.9cm]
\item[(i)] A zero object $0\in\Ob\B$ exists.
\item[(ii)] For any objects $B,B'\in\Ob\B$, the class $\Mor_\B(B,B')$ is an $R$-module over which composition is $R$-distributive in the sense that given {\small $g,h\in\Mor_\B(B,B')$, $f':B_1\ra B$, $f':B'\ra B'_1$, $r\in R$}, we have
{\small\[
f'\circ(g+h)=f'\circ g+f'\circ h,~~~~(g+h)\circ f=g\circ f+h\circ f,~~~~f'\circ (rg)=r(f'\circ g),~~~~(rg)\circ f=r(g\circ f).
\]}
\eit
In particular, a category $\A$ is \ul{pre-additive} $\iff$ it is $\Integer$-preadditive.

(\ul{Right $R$-additivity} is similarly defined in the obvious way.)

Let $\B,\D$ be $R$-additive categories. A functor $F:\B\ra\D$ is an \ul{$R$-additive functor} (or \ul{$R$-linear functor}) if for all objects $B,B'\in\Ob\B$, $F:\Mor_\B(B,B')\ra \Mor_\D(F(B),F(B'))$ is an $R$-homomorphism, i.e.,
\[
F(f+g)=F(f)+F(g)~~~\txt{and}~~~F(rf)=rF(f),~~~~\txt{for all}~~~~f,g\in\Mor_\B(B,B'),~r\in R.
\]
\end{dfn}

\begin{dfn}[\blue{
\index{Endomorphism! algebra}{Endomorphism algebra},
\index{Algebra! as a category}{Algebra as a category},
\index{Algebra! representation}{Algebra representation},
\index{Linear! algebra representation}{Linear algebra representation}}]
Let $R$ be a commutative ring and $\B$ an \magenta{$R$-preadditive} category. For any object $B\in\Ob\B$, the \ul{endomorphism $R$-algebra} of $B$ is the $R$-module ~$End_\B(B):=\Mor_\B\big(B,B\big)$~ as a ring with multiplication given by composition of morphisms and unity $1_{End_\B(B)}:=id_B$.

Let $A$ be an $R$-algebra. By construction, $A$ can be viewed as the \magenta{$R$-preadditive} category $\C(A)$ with a single object $A$ and morphisms being the algebra elements $a\in A$ as left multiplication maps $a:A\ra A$, $x\mapsto ax$. A \ul{representation} of $A$ in $\B$ (or a $\B$-representation of $A$) is an \magenta{$R$-additive} functor
\bea
F:\C(A)\ra\B,~~A\sr{a}{\ral}A~~\mapsto~~F(A)\sr{F(a)}{\ral}F(A),\nn
\eea
that is, we have an $R$-algebra homomorphism (\blue{footnote}\footnote{That is, given an $R$-algebra $A$ and an $R$-preadditive category $\B$, a \ul{representation} of $A$ in $\B$ (or a $\B$-representation of $A$) is an $R$-algebra homomorphism ~$\rho:A\ra End_\B(B)$~ for an object $B\in\Ob\B$.}) ~$F:A\ra End_\B\big(B\big)$,~ for some object $B:=F(A)\in\Ob\B$.

If $\B$ is the category of $k$-spaces (or $k$-modules for a commutative ring $k$), the representation is a \ul{$k$-linear representation} of $A$.
\end{dfn}

\begin{rmk}[\blue{
\index{Lie! algebra representations}{Lie algebra representations},
\index{Lie! algebra homomorphism}{Lie algebra homomorphism},
\index{Universal! enveloping algebra}{Universal enveloping algebra}}]
Let $R$ be a commutative ring.

Although we will not be concerned with representations of non-associative algebras, they are included in the discussion on representations of multiplicative sets in general. For Lie $R$-algebras in particular, every associative $R$-algebra $A$ is itself a Lie $R$-algebra $Lie(A)$ with respect to the operation {\small $A\times A\ra A,~(a,b)\mapsto a\ast b:=ab-ba$}. Thus, if $L$ is a Lie $R$-algebra (\magenta{with multiplication denoted by {\small $L\times L\ra L,~(a,b)\mapsto a\cdot b$}}) and $f:L\ra Lie(A)$ a \ul{Lie $R$-algebra homomorphism} (\magenta{i.e., {\small $f(ra+b)=rf(a)+f(b)$, $f(a\cdot b)=f(a)\ast f(b)$} for all $a,b\in L$, $r\in R$, and if necessary, $f(1)=1$}), then a representation $\rho:A\ra End_\B(B)$ of $A$ in an $R$-preadditive category $\B$ gives a representation {\small $\rho\circ f:L\sr{f}{\ral}Lie(A)\sr{\rho}{\ral} Lie(End_\B(B))$} of $L$ in $\B$.

This is sufficient because we can always choose $A$ to be the ``universal enveloping algebra'' $U(L)$ of $L$ defined as follows: Let $(a,b)\mapsto a\cdot b$ denote multiplication of $a,b\in L$ as before. With the tensor algebra
\[
\textstyle Te(L):=R\langle L\rangle\cong R^{\langle X\rangle}=\bigoplus_{x\in X}R=\bigoplus_{n\in\Natural}R^{\langle L^n\rangle},~~~~\txt{where}~~~~X:=\bigcup_{n\in\Natural}L^n,
\]
and the ideal $I:=Te(L)Q\lhd Te(L)$ generated by elements of the form $Q:=\{ab-ba-a\cdot b~|~a,b\in L\}\subset Te(L)$, the \ul{universal enveloping algebra} of $L$ is the associative $R$-algebra
\bea
\textstyle U(L):={Te(L)\over I}={R\langle L\rangle\over R\langle L\rangle Q}=R\langle\{a+I:a\in L\}\rangle,~~~~\txt{where}~~~~ab-ba+I=a\cdot b+I~~\txt{for all}~~a,b\in L.\nn
\eea
The map $f:L\hookrightarrow Lie(U(L)),~a\mapsto a+I$ is a Lie $R$-algebra monomorphism and, based on $f$, $U(L)$ has the following universal property: Given a Lie $R$-algebra homomorphism $h:L\ra Lie(A_1)$ with $A_1$ an associative $R$-algebra, we get the unique (Lie) $R$-algebra homomorphism $h_1:Lie(U(L))\ra Lie(A_1),~a+I\mapsto h(a)$ satisfying $h_1\circ f=h$.
\[\adjustbox{scale=0.9}{\bt
L\ar[d,"h"]\ar[rr,hook,"f"] && Lie(U(L))\ar[dll,dashed,"h_1"]\\
Lie(A_1)
\et}\]
\end{rmk}

\subsection{Action-representation correspondences}
\begin{thm}[\blue{Group actions = Permutation representations}]\label{GpActRepThm}
Let $G$ be a group and $S$ a set. We have a bijection
\bea
{Act}(G\times S,S):=\left\{
  \begin{array}{l}
    \txt{left actions} \\
    m:G\times S\ra S
  \end{array}
\right\} \longleftrightarrow
\left\{
  \begin{array}{l}
    \txt{group homomorphisms} \\
    h:G\ra Sym(S)
  \end{array}
\right\}=:Hom_{Grp}\big(G,Sym(S)\big),\nn
\eea
where~ $Sym(S):=\{\txt{bijections}~f:S\ra S\}$ is the group of permutations of $S$.
\end{thm}
\begin{proof}
Consider the following maps:
\bea
&&\phi:{Act}(G\times S,S)\ra {Hom_{Grp}}\big(G,Sym(S)\big),~~G\times S\sr{m}{\ral}S~~\mapsto~~G\sr{\phi(m)}{\ral}Sym(S),\nn\\
&&\psi:{Hom_{Grp}}\big(G,Sym(S)\big)\ra{Act}(G\times S,S),~~G\sr{h}{\ral}Sym(S)~~\mapsto~~G\times S\sr{\psi(h)}{\ral}S,\nn
\eea
where the maps $\phi(m)$ and $\psi(h)$ are defined as follows: For all $g\in G$, $s\in S$,
\bea
\phi(m):g\mapsto m_g:=m(g,-)~~~~\txt{and}~~~~\psi(h):(g,s)\mapsto h(g)s.\nn
\eea
Then $\phi,\psi$ are (a) well-defined and (b) inverses of each other (as shown below):
{\flushleft(a)}
For any $m\in {Act}(G\times S,S)$, ~$\phi(m)$ is a \ul{group homomorphism} (\blue{footnote}\footnote{We have already shown that for all $g\in G$, the map $m_g:=m(g,-):S\ra S,~s\mapsto gs$ is a permutation of $S$.}): Indeed, for all $g,g'\in G$,
\bea
\phi(m)(gg')=m_{gg'}=m(gg',-)=m_g\circ m_{g'}=\phi(m)(g)\circ\phi(m)(g').\nn
\eea
For any $h\in {Hom_{Grp}}\big(G,Sym(S)\big)$, ~$\psi(h)$ is an \ul{action} of $G$ on $S$: Indeed, for all $g,g'\in G$ and $s\in S$,
\bea
&& \psi(h)(gg',s)=h(gg')s=\big(h(g)\circ h(g')\big)(s)=h(g)\big(h(g')s\big)=\psi(h)\big(g,\psi(h)(g',s)\big),\nn\\
&& \psi(h)(1_G,s)=h(1_G)s=1_{Sym(S)}s=s.\nn
\eea

{\flushleft (b)}
The compositions $\psi\circ\phi$ and $\phi\circ\psi$ satisfy the following:
{\small\bea
&& \psi\circ\phi: {Act}(G\times S,S)\ra {Act}(G\times S,S),~~G\times S\sr{m}{\ral}S~~\sr{\phi}{\mapsto}~~G\sr{\phi(m)}{\ral}Sym(S)~~\sr{\psi}{\mapsto}~~G\times S\sr{\psi(\phi(m))}{\ral}S,\nn\\
&&\psi(\phi(m))(g,s)=\phi(m)(g)s=m_gs=m(g,s),~~\Ra~~~\psi\circ\phi=1_{{Act}(G\times S,S)},\nn
\eea}and similarly,
{\footnotesize\bea
&& \phi\circ\psi: {Act}(G,Sym(S))\ra {Act}(G,Sym(S)),~~G\sr{h}{\ral}Sym(S)~~\sr{\psi}{\mapsto}~~G\times S\sr{\psi(h)}{\ral}S~~\sr{\phi}{\mapsto}~~G\sr{\phi(\psi(h))}{\ral}Sym(S),\nn\\
&&\phi(\psi(h))(g)=\psi(h)_g=\psi(h)(g,-)=h(g),~~\Ra~~\phi\circ\psi=1_{{Act}(G,Sym(S))}.\nn
\eea}
\end{proof}

\begin{dfn}[\blue{\index{Linear! group action}{Linear group action}}]
Let $R$ be a commutative ring and $M$ an $R$-module that is a G-set with action $G\times M\ra M,~(g,m)\mapsto gm$. The $G$-action is $R$-linear if for all $r\in R$, $g\in G$, $m,m'\in M$, we have (i) $g~rm=r gm$ and (ii) $g(m+m')=gm+gm'$.
\end{dfn}

The following is a corollary of the proof of Theorem \ref{GpActRepThm}.
\begin{crl}[\blue{Linear group actions = Linear group representations}]\label{LinGpActRepThm}
Let $R$ be a commutative ring and $M$ an $R$-module. We have a bijection
{\small\begin{align}
{Act_R}(G\times M,M):=\left\{
  \begin{array}{l}
    \txt{$R$-linear left actions} \\
    m:G\times M\ra M
  \end{array}
\right\} \longleftrightarrow
\left\{
  \begin{array}{l}
    \txt{group homomorphisms} \\
    h:G\ra GL_R(M)
  \end{array}
\right\}=:Hom_{Grp}\big(G,GL_R(M)\big).\nn
\end{align}}where {\footnotesize $GL_R(M):=Aut_R(M)=U\big(End_R(M)\big)$} is the group of units of the endomorphism ring {\footnotesize $End_R(M):=Hom_R(M,M)$}.
\end{crl}
\begin{proof}
In the proof of Theorem \ref{GpActRepThm}, replace $Act(...)$ with $Act_R(...)$, $S$ with $M$, $Sym(S)$ with $GL_R(M)$. As before, consider the following maps:
\bea
&&\phi:{Act}_R(G\times M,M)\ra {Hom_{Grp}}\big(G,GL_R(M)\big),~~G\times M\sr{m}{\ral}M~~\mapsto~~G\sr{\phi(m)}{\ral}GL_R(M),\nn\\
&&\psi:{Hom_{Grp}}\big(G,GL_R(M)\big)\ra{Act}_R(G\times M,M),~~G\sr{h}{\ral}GL_R(M)~~\mapsto~~G\times M\sr{\psi(h)}{\ral}M,\nn
\eea
where the maps $\phi(m)$ and $\psi(h)$ are defined as follows: For all $g\in G$, $s\in M$,
\bea
\phi(m):g\mapsto m_g:=m(g,-)~~~~\txt{and}~~~~\psi(h):(g,s)\mapsto h(g)s.\nn
\eea
Then $\phi,\psi$ are (a) well-defined and (b) inverses of each other (as explained below):
{\flushleft(a)}
For any $m\in {Act}_R(G\times M,M)$, ~$\phi(m)$ is a \ul{group homomorphism} (shown exactly as before), where for each $g\in G$,
$\phi(m)(g)=m_g=m(g,-):M\ra M$ is an $R$-isomorphism with inverse $\phi(m)(g)^{-1}:=\phi(m)(g^{-1})$.

For any $h\in {Hom_{Grp}}\big(G,GL_R(M)\big)$, ~$\psi(h)$ is an \ul{$R$-linear operation} of $G$ on $M$ (shown to be a $G$-action exactly as before, and $R$-linear because for any $g\in G$, $h(g):M\ra M$ is an $R$-isomorphism).
{\flushleft (b)}
Exactly as before, we have ~$\psi\circ\phi=1_{{Act}_R(G\times M,M)}$~ and ~$\phi\circ\psi=1_{{Act}_R(G,GL_R(M))}$.
\end{proof}

\begin{rmk}[\blue{Left-right symmetry of group representations}]
Since every left group action $G\times S\ra S,~(g,s)\mapsto gs$ is equivalent to a right group action $S\times G\ra S,~(s,g)\mapsto sg:=g^{-1}s$, it is not necessary to distinguish between left-representations and right-representations of a group. This is however not the case for representations of a ring/algebra.
\end{rmk}

\begin{dfn}[\blue{
\index{Faithful! action}{Faithful action},
\index{Faithful! representation}{Faithful representation}}]
Let $G$ be a group, $\C\subset Sets$ a subcategory, and $C\in\Ob\C$ an object in $\C$. An action ~$m:G\times C\ra C,~(g,c)\ra g\cdot c$~ is a \ul{faithful action} if the corresponding representation ~$m_\ast:G\ra Aut_\C(C),~g\mapsto g\cdot(-)$~ is injective. An injective representation is also called a \ul{faithful representation}. The same terminology applies as well to representations of other types of objects, for example, faithful ring (resp. algebra) actions will correspond to faithful/injective ring (resp. algebra) homomorphisms.
\end{dfn}

\begin{dfn}[\blue{
\index{Action! of a ring}{Action of a ring} on an abelian group,
\index{Module}{{Module}}}]
Let $R$ be a ring and $M$ an abelian group. A \ul{left action} of $R$ on $M$ is a map ~$R\times M\ra M,~(r,m)\mapsto rm$~ with the following properties (which precisely mean that $M$ is a \ul{left $R$-module}): For all $r,r'\in R$, $m,m'\in M$,
\bit[leftmargin=0.7cm]
\item (i) $1_Rm=m$, (ii) $r(r'm)=(rr')m$, (iii) $(r+r')m=rm+r'm$ and $r(m+m')=rm+rm'$.
\eit
\end{dfn}

The following is another corollary of the proof of Theorem \ref{GpActRepThm}.

\begin{crl}[\blue{Modules = Ring representations}]\label{RingActRepThm}
Let $R$ be a ring, and $M$ an abelian group. We have a bijection
{\footnotesize\bea
Act(R\times M,M):=\left\{
  \begin{array}{l}
    \txt{$R$-module actions} \\
    m:R\times M\ra M
  \end{array}
\right\} \longleftrightarrow
\left\{
  \begin{array}{l}
    \txt{$R$-homomorphisms} \\
    h:A\ra End_\Integer(M)
  \end{array}
\right\}=:Hom_{Rings}\big(R,End_\Integer(M)\big),\nn
\eea}where ~$End_\Integer(M):=\Mor_\Integer(M,M)$~ as a ring.
\end{crl}
\begin{proof}
In the proof of Theorem \ref{GpActRepThm}, replace $G$ with $R$, $S$ with $M$, $Sym(S)$ with $End_\Integer(M)$, and $Hom_{Grp}(...)$ with $Hom_{Rings}(...)$. Consider the following maps:
\bea
&&\phi:{Act}(R\times M,M)\ra {Hom_{Rings}}\big(R,End_\Integer(M)\big),~~R\times M\sr{m}{\ral}M~~\mapsto~~R\sr{\phi(m)}{\ral}End_\Integer(M),\nn\\
&&\psi:{Hom_{Rings}}\big(R,End_\Integer(M)\big)\ra{Act}(R\times M,M),~~R\sr{h}{\ral}End_\Integer(M)~~\mapsto~~R\times M\sr{\psi(h)}{\ral}M,\nn
\eea
where the maps $\phi(m)$ and $\psi(h)$ are defined as follows: For all $r\in R$, $x\in M$,
\bea
\phi(m):r\mapsto m_r:=m(r,-)~~~~\txt{and}~~~~\psi(h):(r,x)\mapsto h(r)x.\nn
\eea
Then $\phi,\psi$ are (a) well-defined and (b) inverses of each other (as shown/explained below):
{\flushleft(a)}
For any $m\in {Act}(R\times M,M)$, ~$\phi(m)$ is a \ul{ring homomorphism} (\blue{footnote}\footnote{It is clear that for all $r\in R$, the map ~$m_r:=m(r,-):M\ra M,~x\mapsto rx$~ is a $\Integer$-homomorphism.}): Indeed, for all $r,r'\in R$,
\bea
&&\phi(m)(rr')=m_{rr'}=m(rr',-)=m_r\circ m_{r'}=\phi(m)(r)\circ\phi(m)(r'),\nn\\
&&\phi(m)(r+r')=m_{r+r'}=m(r+r',-)=m(r,-)+m(r',-)=m_r+m_{r'}=\phi(m)(r)+\phi(m)(r'),\nn\\
&&\phi(m)(1_R)=m_{1_R}=m(1_R,-)=id_M=1_{End_\Integer(M)}.\nn
\eea
For any $h\in {Hom_{Rings}}\big(R,End_\Integer(M)\big)$, ~$\psi(h)$ is an \ul{action} of $R$ on $M$: Indeed, for all $r,r'\in R$ and $x,x'\in M$,
{\small\bea
&& \psi(h)(rr',x)=h(rr')x=\big(h(r)\circ h(r')\big)(x)=h(r)\big(h(r')x\big)=\psi(h)\big(r,\psi(h)(r',x)\big),\nn\\
&& \psi(h)(r+r',x)=h(r+r')x=\big(h(r)+h(r')\big)x=h(r)x+ h(r')x=\psi(h)(r,x)+\psi(h)(r',x),\nn\\
&&\psi(h)(r,x+x')=h(r)(x+x')=h(r)x+h(r)x'=\psi(h)(r,x)+\psi(h)(r,x'),\nn\\
&& \psi(h)(1_R,x)=h(1_R)x=1_{End_\Integer(M)}x=x.\nn
\eea}
{\flushleft (b)} Exactly as before, we have ~$\psi\circ\phi=1_{{Act}(R\times M,M)}$~ and ~$\phi\circ\psi=1_{{Act}(R,End_\Integer(M))}$.
\end{proof}

\begin{dfn}[\blue{
\index{Action! of an algebra}{Action of an algebra} on a module,
\index{Algebra! module}{{Algebra module}}}]
Let $R$ be a commutative ring, $A$ an $R$-algebra, and $M$ an $R$-module. A \ul{left action} of $A$ on $M$ is a map ~$A\times M\ra M,~(a,x)\mapsto ax$~ with the following properties (which precisely mean that $M$ is a \ul{left $A$-module}): For all $r\in R$, $a,a'\in A$, $x,x'\in M$,
\bit[leftmargin=0.7cm]
\item (i) $1_Ax=x$, (ii) $a(a'x)=(aa')x$, (iii) $(a+a')x=ax+a'x$ and $a(x+x')=ax+ax'$, (iv) $arx=rax$.
\eit
\end{dfn}

Since an arbitrary ring $R$ is precisely a $\Integer$-algebra, a left action of $R$ on an abelian group (i.e., a $\Integer$-module) $M$, making $M$ a (left) $R$-module, is also included in the last definition above. That is, a module is precisely a $\Integer$-algebra module. Consequently, (i) any results that hold for algebra modules automatically hold for modules, and (ii) results that hold for modules might extend to algebra modules if appropriate modifications can be made.

The following is one more corollary of the proof of Theorem \ref{GpActRepThm}.

\begin{crl}[\blue{Algebra modules = Algebra representations}]\label{AlgActRepThm}
Let $R$ be a commutative ring, $A$ an $R$-algebra, and $M$ an $R$-module. We have a bijection
{\footnotesize\bea
Act(A\times M,M):=\left\{
  \begin{array}{l}
    \txt{$A$-module actions} \\
    m:A\times M\ra M
  \end{array}
\right\} \longleftrightarrow
\left\{
  \begin{array}{l}
    \txt{$R$-algebra homomorphisms} \\
    h:A\ra End_R(M)
  \end{array}
\right\}=:Hom_{R\txt{-Alg}}\big(A,End_R(M)\big),\nn
\eea}where ~$End_R(M):=Hom_R(M,M)$ as a ring.
\end{crl}
\begin{proof}
In the proof of Theorem \ref{GpActRepThm}, replace $G$ with $A$, $S$ with $M$, $Sym(S)$ with $End_R(M)$, and $Hom_{Grp}(...)$ with $Hom_{R\txt{-Alg}}(...)$. Consider the following maps:
\bea
&&\phi:{Act}(A\times M,M)\ra {Hom_{R\txt{-Alg}}}\big(A,End_R(M)\big),~~A\times M\sr{m}{\ral}M~~\mapsto~~A\sr{\phi(m)}{\ral}End_R(M),\nn\\
&&\psi:{Hom_{R\txt{-Alg}}}\big(A,End_R(M)\big)\ra{Act}(A\times M,M),~~A\sr{h}{\ral}End_R(M)~~\mapsto~~A\times M\sr{\psi(h)}{\ral}M,\nn
\eea
where the maps $\phi(m)$ and $\psi(h)$ are defined as follows: For all $a\in A$, $x\in M$,
\bea
\phi(m):a\mapsto m_a:=m(a,-)~~~~\txt{and}~~~~\psi(h):(a,x)\mapsto h(a)x.\nn
\eea
Then $\phi,\psi$ are (a) well-defined and (b) inverses of each other (as shown/explained below):
{\flushleft(a)}
For any $m\in {Act}(A\times M,M)$, ~$\phi(m)$ is an \ul{$R$-algebra homomorphism} (\blue{footnote}\footnote{It is clear that for all $a\in A$, the map ~$m_a:=m(a,-):M\ra M,~x\mapsto ax$~ is an $R$-homomorphism.}): Indeed, for all $a,a'\in A$, $r\in R$,
\bea
&&\phi(m)(aa')=m_{aa'}=m(aa',-)=m_a\circ m_{a'}=\phi(m)(a)\circ\phi(m)(a'),\nn\\
&&\phi(m)(a+a')=m_{a+a'}=m(a+a',-)=m(a,-)+m(a',-)=m_a+m_{a'}=\phi(m)(a)+\phi(m)(a'),\nn\\
&&\phi(m)(ra)=m_{ra}=m(ra,-)=rm(a,-)=rm_a=r\phi(m)(a), \nn\\
&&\phi(m)(1_A)=m_{1_A}=m(1_A,-)=id_M=1_{End_R(M)}.\nn
\eea
For any $h\in {Hom_{R\txt{-Alg}}}\big(A,End_R(M)\big)$, ~$\psi(h)$ is an \ul{action} of $A$ on $M$: Indeed, for all $a,a'\in A$ and $x,x'\in M$ and $r\in R$,
{\small\bea
&& \psi(h)(aa',x)=h(aa')x=\big(h(a)\circ h(a')\big)(x)=h(a)\big(h(a')x\big)=\psi(h)\big(a,\psi(h)(a',x)\big),\nn\\
&& \psi(h)(a+a',x)=h(a+a')x=\big(h(a)+h(a')\big)x=h(a)x+ h(a')x=\psi(h)(a,x)+\psi(h)(a',x),\nn\\
&&\psi(h)(a,x+x')=h(a)(x+x')=h(a)x+h(a)x'=\psi(h)(a,x)+\psi(h)(a,x'),\nn\\
&&\psi(h)(a,rx)=h(a)(rx)=rh(a)(rx)=r\psi(h)(a,x),\nn\\
&& \psi(h)(1_A,x)=h(1_A)x=1_{End_R(M)}x=x.\nn
\eea}
{\flushleft (b)} Exactly as before, we have ~$\psi\circ\phi=1_{{Act}(A\times M,M)}$~ and ~$\phi\circ\psi=1_{{Act}(A,End_R(M))}$.
\end{proof}

\subsection{The group algebra and linear group representations}
\begin{dfn}[\blue{\index{Group! algebra}{Group algebra}}]
Let $R$ be a commutative ring and $G$ a group. The $R$-algebra of $G$ (denoted by $RG$ or $R[G]$ ) is the free $R$-module {\small $R[G]:=R^{\langle G\rangle}=\sum_{g\in G}Rg=\{\sum_{g\in G}r_gg:r_g\in R,~r_g=0~\txt{a.e.f.}~g\in G\}$} with basis $G$, as a ring with \ul{addition} and \ul{multiplication} given by
\bea
&&\textstyle\sum\al_gg+\sum\beta_gg:=\sum(\al_g+\beta_g)g,~~\left(\sum\al_gg\right)\left(\sum\beta_gg\right):=\sum\gamma_g(\al,\beta)g,~~\gamma_g(\al,\beta):=\sum_{hh'=g}\al_h\beta_{h'},\nn
\eea
\ul{zero} and \ul{unity} given by ~$0_{R[G]}:=0_Re_G$~ and ~$1_{R[G]}:=1_Re_G$, and \ul{scalar multiplication} given by
\bea
\textstyle R\times R[G]\ra R[G],~~\left(r,\sum\al_gg\right)\ra r\cdot \sum\al_gg:=\sum r\al_g~g.\nn
\eea
We have $G\cong 1_RG\subset U(R[G])$ and $R\cong Re_G\subset Z(R[G])$.
\end{dfn}

\begin{thm}[\blue{Linear group representations = Group algebra representations}]\label{FT-LinGp-Reps}
Let $R$ be a commutative ring, $M$ an $R$-module, and $G$ a group. Then we have a bijection
{\footnotesize \bea
\bt
\left\{
  \begin{array}{l}
    \txt{group homomorphisms} \\
    h:G\ra GL_R(M)
  \end{array}
\right\}
\ar[rr,leftrightarrow] &&
\left\{
  \begin{array}{l}
    \txt{ring homomorphisms} \\
    h:R[G]\ra End_R(M)
  \end{array}
\right\} \ar[d,leftrightarrow,"\txt{Cor. \ref{LinGpActRepThm}}"]\\
\left\{
  \begin{array}{l}
    \txt{$R$-linear left actions} \\
    m:G\times M\ra M
  \end{array}
\right\}\ar[u,leftrightarrow,"\txt{Cor. \ref{AlgActRepThm}}"']\ar[rr,dotted,leftrightarrow]
 &&
\left\{
  \begin{array}{l}
    \txt{$R[G]$-module actions} \\
    m:R[G]\times M\ra M
  \end{array}
\right\},\et\nn
\eea}
where ~{\small $GL_R(M):=Aut_R(M)=U\big(End_R(M)\big)$}~ are the units in ~{\small $End_R(M):=Hom_R(M,M)$}~ as a ring. Hence, we have a bijection
~{\footnotesize
$\left\{
  \begin{array}{l}
    \txt{$R$-linear group actions} \\
    m:G\times M\ra M
  \end{array}
\right\} \longleftrightarrow
\left\{
  \begin{array}{l}
    \txt{$R[G]$-module actions} \\
    m:R[G]\times M\ra M
  \end{array}
\right\}.$}
\end{thm}
\begin{proof}
Let $h:G\ra GL_R(M):=Aut_R(M)=U(End_R(M))$ be a group homomorphism. Then we know there exists a \ul{unique} $R$-homomorphism $f:R[G]\ra End_R(M)$ such that $f|_G=h$.
\[\adjustbox{scale=0.8}{\bt
G\ar[d,"h"]\ar[rr,hook] && R[G]:=R^{\langle G\rangle}\ar[ddll,dashed,"f"] \\
Aut_R(M)\ar[d,hook] && \\
End_R(M)
\et}\]To see that $f$ is a ring homomorphism, observe that for any $\sum_g \al_g g$ and $\sum_g\beta_gg$ in $R[G]$, we have
{\small\begin{align}
&\textstyle f\left(\sum_g \al_g g~\sum_g\beta_gg\right)=f\left(\sum_{gg'}\al_g\beta_{g'}gg'\right)=\sum_{gg'} \al_g\beta_{g'}f(gg')=\sum_{gg'}\al_g\beta_{g'} h(gg')=\sum_{gg'}\al_g\beta_{g'} h(g)h(g')\nn\\
&\textstyle~~~~=\sum_{gg'}\al_g\beta_{g'} f(g)f(g')=\left(\sum_g\al_gf(g)\right)\left(\sum_g\beta_gf(g)\right)=f\left(\sum_g\al_gg\right)f\left(\sum_g\beta_gg\right).\nn
\end{align}}
This extension process gives a \ul{well-defined} map {\small $\phi:Hom_{Grp}(G,Aut_R(M))\ra Hom_{Rings}(R[G],End_R(M))$, $h\mapsto\phi(h):=f$}. Moreover, $\phi$ is injective because if $\phi(h)=\phi(h')$, then $h=\phi(h)|_G=\phi(h')|_G=h'$.

On the other hand, any ring homomorphism $f:R[G]\ra End_R(M)$ restricts to a group homomorphism $f|_G:G\ra Aut_R(M)$, since a ring homomorphism maps units to units, i.e., if $a,b\in R[G]$ satisfy $ab=ba=1_{R[G]}$, then applying $f$ on both sides, we get $f(a)f(b)=f(b)f(a)=1_{End_R(M)}$. The restriction $f|_G$ extends, as before, to a unique ring homomorphism $\phi(f|_G):R[G]\ra End_R(M)$ which must therefore be $f$ (i.e., $f|_G$ uniquely extend back to $f$), and so the map $\phi$ above is surjective.
\end{proof}

\subsection{Analysis of irreducible representations}
\begin{dfn}[\blue{In Groups:
\index{Representation}{Representation},
\index{Representation! space}{Representation space},
\index{Subrepresentation}{Subrepresentation},
\index{Irreducible! group representation}{Irreducible representation}}]
Let $G$ be a group and $S$ a set. By the equivalence of actions and representations established in Theorem \ref{GpActRepThm}, an action $m:G\times S\ra S,~(g,s)\mapsto gs:=m(g,s)$ of $G$ on $S$ makes $S=(S,m)$ a \ul{representation} (or \ul{representation space}) of $G$.

If $U\subset G$, then $U$ is a \ul{subrepresentation} of $G$ in $S$ if $gU=\{gs:s\in U\}\subset U$ for all $g\in G$, i.e., the restriction $m|_{G\times U}:G\times U\ra U$ gives a well-defined action of $G$ on $U$. We have already seen (in Proposition \ref{SubGsetCrit}) that $U\subset S$ is a subrepresentation if and only if $U$ is a union orbits in $S$.

A representation $(S,m)$ is \ul{irreducible} if it has no proper subrepresentation (i.e., it is \ul{transitive}, or equivalently, the representation space $S$ is a single $G$-orbit). If $S=(S,m)$ is a representation of $G$, then examples of irreducible subrepresentations of $G$ in $S$ are the $G$-orbits $\sigma_s:=\{gs:g\in G\}$, for all $s\in S$.
\end{dfn}
At this point, besides other concerns such as counting of (irreducible) representations, one should expect that a further study of representations of a given group $G$ might involve expressing a given representation of $G$ in terms of irreducible representations, say using the product or otherwise. Also, for the linear representations of $G$ in particular, it follows from Theorem \ref{FT-LinGp-Reps} that a study of the representations of algebras (to be discussed later) should suffice. (\blue{footnote}\footnote{Given a group $G$, except perhaps for special cases such abelian groups (i.e., $\Integer$-modules), one should in general expect more representations of $G$ than just the linear representations that can be characterized as in Theorem \ref{FT-LinGp-Reps}. This follows from the fact that in general $R$-mod is a proper subcategory of Sets, and so has fewer objects.})

\begin{question}[\blue{See below}]
Let $G$ be a group. How many irreducible representations does $G$ have?
\end{question}
To answer this question, observe that irreducible representations correspond one-to-one with $G$-orbits (i.e., transitive $G$-sets). Also, we have seen that as a G-set, every $G$-orbit $\sigma\cong{G\over G_\sigma}$ for some subgroup $G_\sigma\leq G$. Hence, the number of irreducible representations of $G$ is the same as the number of non-isomorphic G-sets of the form ${G\over H}$ for subgroups $H\leq G$.

Also, the $G$-set $G/H$ (i.e., the representation $m_\ast:G\ra Sym(G/H),~g\mapsto m_g:=m(g,-)$ of $G$ corresponding to the subgroup $H\leq G$) is faithful/injective $\iff$  for any $g\in G$, if $g(aH)=aH$ for all $a\in G$ then $g=e$, i.e., for any $g\in G$,
\[
\{a^{-1}ga:a\in G\}\subset H~~\Ra~~g=e,~~~\txt{or equivalently},~~~~g\neq e~~\Ra~~\{a^{-1}ga:a\in G\}\not\subset H.
\]This shows that for any nontrivial normal subgroup $\{e\}\neq N\lhd G$, $G/N$ is not a faithful $G$-set. In particular, if $G$ is abelian, then $G/H$ is a faithful $G$-set $\iff$ $H=\{e\}$, i.e., the only such faithful $G$-set is $G$ itself.

Now, given subgroups $H,K\leq G$, when are ${G\over H}\cong{G\over K}$ as G-sets? That is, when is a map $f:G/H\ra G/K$, $gH\mapsto\wt{f}(g)K$ an isomorphism of $G$-sets? The condition is as follows.

\begin{thm}
Let $G$ be a group and $H,K\leq G$ subgroups. Then a $G$-isomorphism $G/H\cong G/K$ exists $\iff$ $H$ and $K$ are conjugate subgroups of $G$ (i.e., there exists $x\in G$ such that $H=xKx^{-1}$).
\end{thm}
\begin{proof}
{\flushleft ($\Ra$)} Let $f:G/H\ra G/K$, $gH\mapsto\wt{f}(g)K$ be an isomorphism of $G$-sets. Then (i) well-definedness and injectivity mean $gH=g'H$ iff $\wt{f}(g)K=\wt{f}(g')K$ for all $g,g'\in G$, meanwhile, (ii) the morphism property $f(gg'H)=gf(g'H)$ for all $g,g'\in G$ implies $\wt{f}(g)K=g\wt{f}(e)K$ for all $g\in G$. Let $x_f:=\wt{f}(e)$. Then $h\in H$ iff $hH=H$, iff $f(hH)=f(H)$, iff $\wt{f}(h)K=\wt{f}(e)K$, iff $h\wt{f}(e)K=\wt{f}(e)K$, iff $x_f^{-1}hx_f\in K$. Therefore, $x_f^{-1}Hx_f\subset K$. The inverse {\small $f^{-1}:G/K\ra G/H$, $gK\mapsto\wt{f^{-1}}(g)H$} gives expressions
\bea
&&f^{-1}f(gH)=gH~\txt{and}~~ff^{-1}(gK)=gK~~~~\txt{for all}~~g\in G,\nn\\
&&~~\txt{or equivalently,}~~\wt{f^{-1}}(\wt{f}(e))H=H~~\txt{and}~~\wt{f}(\wt{f^{-1}}(e))K=K,\nn\\
&&~~\txt{or equivalently,}~~\wt{f}(e)\wt{f^{-1}}(e)H=H~~\txt{and}~~\wt{f^{-1}}(e)\wt{f}(e)K=K,\nn\\
&&~~\txt{or equivalently,}~~\wt{f^{-1}}(e)H=\wt{f}(e)^{-1}H~~\txt{and}~~\wt{f}(e)K=\wt{f^{-1}}(e)^{-1}K,\nn\\
&&~~\txt{or equivalently,}~~x_{f^{-1}}H=x_f^{-1}H~~\txt{and}~~x_fK=x_{f^{-1}}^{-1}K,~~~~\txt{where}~~~~x_f:=\wt{f}(e),~~x_{f^{-1}}:=\wt{f^{-1}}(e),\nn
\eea
and so the reverse situation $x_{f^{-1}}^{-1}Kx_{f^{-1}}\subset H$, or equivalently $x_{f^{-1}}^{-1}Kx_{f^{-1}}H\subset H$, implies $x_fKx_f^{-1}\subset H$. Hence $H=x_fKx_f^{-1}$.

{\flushleft ($\Ra$)} Assume there exists $x\in G$ such that $H=xKx^{-1}$. The map ~{\small $f:G/H\ra G/K,~gH\mapsto gxK$}~ is well-defined and injective since $gH=g'H$ iff $gxKx^{-1}=g'xKx^{-1}$, iff $gxK=gxK$. It is also clear that $f$ is a homomorphism of G-sets. Similarly, the map $F:G/K\ra G/H,~gK\mapsto gx^{-1}H$ is a well-defined injective homomorphism of G-sets, and satisfies $F\circ f=id_{G/H}$ and $f\circ F=id_{G/K}$. Hence ${G\over H}\cong{G\over K}$ as $G$-sets.
\end{proof}

\begin{crl}[\blue{Irreducible group representations}]
Let $G$ be a group and $Sub(G):=\{H:H\leq G\}$ the set of subgroups of $G$. (i) The set Irred(G) of irreducible representations of $G$ is cardinality equivalent to the set $CcSub(G):=\big\{[H]_{cc}:=\{gHg^{-1}:g\in G\}~|~H\in Sub(G)\big\}$ of conjugacy classes of subgroups of $G$. (\blue{footnote}\footnote{Observe that the conjugacy class $\{g\langle a\rangle g^{-1}:g\in G\}=\{\langle gag^{-1}\rangle:g\in G\}\subset Sub(G)$ of a cyclic subgroup $\langle a\rangle\leq G$ involves elements of the conjugacy class $\{gag^{-1}:g\in G\}\subset G$ of $a$. However, there is no 1-1 correspondence even when $G$ is abelian (in which case we can have $\langle a\rangle=\langle b\rangle$ while $a\neq b$).}) (ii) The cardinality ~$|Irred(G)|=|CcSub(G)|$~ satisfies
\bea
\textstyle |CcSub(G)|\cdot|G/Z(G)|=\sum_{gZ(G)\in G/Z(G)}|\{H\in Sub(G):gHg^{-1}=H\}|.\nn
\eea
In particular, if $G/Z(G)$ is finite, then
\bea
\textstyle |CcSub(G)|={1\over|G/Z(G)|}\sum_{gZ(G)\in G/Z(G)}|\{H\in Sub(G):gHg^{-1}=H\}|.\nn
\eea
Furthermore, if $G$ is abelian, then this cardinality is the number of subgroups of $G$.
\end{crl}
\begin{proof}
(i) The cardinal equivalence of the said collections is clear from the preceding theorem. (ii) To describe (and hence count the number of) such conjugacy classes, consider the action
\bea
Inn(G)\times Sub(G)\ra Sub(G),~(\phi_g,H)\mapsto \phi_g(H):=gHg^{-1}.\nn
\eea
The Inn(G)-orbits {\footnotesize $\sigma_H:=\{gHg^{-1}:\phi_g\in Inn(G)\}=\{gHg^{-1}:g\in G\}=[H]_{cc}$, for $H\in Sub(G)$}, in Sub(G) are the conjugacy classes of subgroups of $G$. The point-invariant domains are $Sub(G)_{\phi_g}:=\{H\in Sub(G):gHg^{-1}=H\}$, for $\phi_g\in Inn(G)$. Therefore, by Burnside's theorem, the number of orbits $|CcSub(G)|$ in Sub(G), and hence the number of irreducible representations of $G$, satisfies
{\footnotesize\[
\textstyle |CcSub(G)|\cdot|Inn(G)|=\sum_{\phi_g\in Inn(G)}|Sub(G)_{\phi_g}|~~~~\txt{or}~~~~|CcSub(G)|\cdot|G/Z(G)|=\sum_{gZ(G)\in G/Z(G)}|Sub(G)_{\phi_g}|.\qedhere
\]}
\end{proof}
To count the number of faithful irreducible representations of $G$, we simply need to replace $Sub(G)$ in the above theorem with ~$Sub_\ast(G):=\big\{H\leq G~|~\{a^{-1}ga:a\in G\}\not\subset H~\txt{for all}~e\neq g\in G\big\}$.

\begin{question}
Let $G$ be a group and $R$ a commutative ring. (i) Up to isomorphism, what are the irreducible $R$-linear representations of $G$? (Recall from Theorem \ref{FT-LinGp-Reps} that these are the same as the ``irreducible modules'' of the ring $RG$). (ii) How can we estimate the number of such irreducible representations of $G$?
\end{question}

\begin{question}[\blue{\index{Irreducible! algebra representation}{Irreducible representations of an algebra}}]
Let $R$ be a commutative ring and $A$ an $R$-algebra. How do we define an irreducible representation (i.e., an irreducible module) of $A$?
\end{question}
To answer the latter question, let $A_r$ denote $A$ viewed as a ring (i.e., $A$ with its $R$-module structure ignored). Then, along with the zero $A_r$-module, \ul{simple $A_r$-modules} (i.e., $A_r$-modules with no nontrivial proper submodules) are precisely the \ul{irreducible $A_r$-modules} (i.e., minimal subobjects in $A_r\txt{-mod}$). Since $A\txt{-mod}:=A_r\txt{-mod}\cap R\txt{-mod}$, it follows that the irreducible $A$-modules are precisely those \ul{irreducible $A_r$-modules} that are also \ul{irreducible $R$-modules}. That is, the problem of finding irreducible representations of an algebra reduces to that of finding irreducible representations of rings (i.e., \ul{simple modules} over rings).

\begin{thm}
Let $R$ be a ring. We have a bijection of the following form:
\[
\textstyle \L:={\{\txt{maximal left ideals ~${}_RL\subset R$}\}\over\sim}
\longleftrightarrow
\big\{\txt{non-isomorphic simple $R$-modules $S$}\big\}=:\S,
\]
where we define $L\sim L'$ if ${R\over L}\cong{R\over L'}$. (We will show in Corollary \ref{MaxSpecEquiv} that if $R$ is commutative, then $\sim$ is not needed.)
\end{thm}
\begin{proof}
Consider the map $\phi:\L\ra\S,~L\mapsto R/L$. (\blue{footnote}\footnote{In the map $\phi:L\mapsto R/L$, $L$ is really the equivalence class $[L]_\sim:=\{\txt{maximal}~{}_RL'\subset R:L'\sim L\}$, and similarly, $R/L$ is really the equivalence class $[R/L]_{\cong}:=\{\txt{simple}~{}_RS\in R\txt{-mod}:S\cong R/L\}$.}). If $_RL\subset R$ is a maximal left ideal, then $S:={R\over L}=R(1+L)$ is a simple $R$-module because for any $L\neq r+L\in R/L$ (i.e., $r\not\in L$) we have $R(r+L)={Rr+L\over L}={R\over L}$, where $Rr+L=R$ due to the maximality of the left ideal $L\subset R$. This shows the map $\phi$ is well-defined. It is also clear that $\phi$ is injective (up to $R$-isomorphism) because if ${}_RL,{}_RL'\subset R$ are maximal left ideals then, $R/L\cong R/L'$ (as objects in $R$-mod) iff $L\sim L'$. It remains to show $\phi$ is surjective (up to $R$-isomorphism).

Given a simple $R$-module $S$, if we pick any $0\neq m\in S$ and define the map $f:R\ra S,~r\mapsto rm$, then by the 1st isomorphism theorem, ${R\over\ker f}\cong f(R)=Rm=S$ (where the surjectivity of $f$ follows because $S$ is simple). Moreover, because $S$ is simple, $\ker f=\{r\in R:rm=0\}\subset R$ is a maximal left ideal.
\end{proof}

Therefore, with the definition ~$\txt{MaxLeftSpec}(R):=\{\txt{maximal left ideals}~~{}_RL\subset R\}$,~ finding irreducible representations of $R$ is equivalent to finding equivalence classes of maximal left ideals:
\bea
\textstyle [L]_{R\txt{-iso}}:=\left\{L'\in \txt{MaxLeftSpec}(R):{R\over L}\cong_R{R\over L'}\right\},~~~~L\in \txt{MaxLeftSpec}(R).\nn
\eea
An immediate hopeful consequence of this result is that, if $R$ is a ring then, an $R$-module $M$ might be completely expressible in terms of the $R$-modules $\left\{{R\over L}~|~\txt{for maximal left ideals}~L\in \txt{MaxLeftSpec(R)}\right\}$, say using the product, coproduct, etc. More generally, if $R$ is a commutative ring and $A$ an $R$-algebra then, a natural means of study/classification of $A$-modules (i.e., representations of $A$) is concerned with \ul{semisimplicity}/\ul{decomposability} of $A$-modules, which will be discussed in chapter \ref{RepsTheoryIII}.

\begin{rmk}
With $L\sim L'$ if ${R\over L}\cong_R{R\over L'}$ (for left ideals $_RL,_RL'\subset R$), we similarly have a bijection
\[
\textstyle\L:={\{\txt{left ideals ~${}_RL\subset R$}\}\over\sim}
\longleftrightarrow \big\{\txt{non-isomorphic cyclic $R$-modules $Rm$}\big\}=:\C\nn
\]
via the maps {\footnotesize $\phi:\L\ra\C,~L\mapsto R/L$} and {\footnotesize $\psi:\C\ra\L,~Rm\cong{R\over K_m}\mapsto K_m:=\ker(R\ra Rm)=\{r\in R:rm=0\}$}, where
\bea
\textstyle\phi\psi(Rm)=\phi(K_m)={R\over K_m}\cong Rm,~~~~\psi\phi(L)=\psi({R\over L})=\psi\big(R(1+L)\big)=K_{1+L}=L,\nn
\eea
that is (up to $R$-isomorphism) we have ~$\phi\psi=id_\C$~ and ~$\psi\phi=id_\L$. (\blue{footnote}\footnote{In the maps $\phi(L)$ and $\psi(Rm)$, $L$ is really the equivalence class $[L]_\sim:=\{{}_RL'\subset R:L'\sim L\}$, and similarly, $Rm$ is really the equivalence class $[Rm]_{\cong}:=\{Rm'\in R\txt{-mod}:Rm'\cong Rm\}$.})
\end{rmk}

\begin{question}
Let $R$ be a ring and {\small ${}_RI,{}_RJ\subset R$} \ul{left ideals}. (i) When is $f$ below an $R$-isomorphism?
\bea
f:R/I\ra R/J,~r+I\mapsto \wt{f}(r)+J\nn
\eea
(ii) If the left ideals $I,J$ are maximal, how does the answer change? (iii) If $R$ is commutative, how does the answer change?
\end{question}
The answer to the above question is given by the following result.

\begin{thm}
Let $R$ be a ring and ${}_RI,{}_RJ\subsetneq R$ proper left ideals. (i) An $R$-isomorphism ${R\over I}\cong{R\over J}$ exists $\iff$ there exist elements $x,y\in R$ satisfying the following:
\bit
\item[(a)] $xy+I=1+I$ and $yx+J=1+J$ ~~(required for \ul{mutual inversion} of such an iso).
\item[(b)] $Ix\subset J$ and $Jy\subset I$ ~~(required by \ul{well-definedness} of such an iso).
\eit
Moreover, (a) implies $x\not\in J$, $y\not\in I$, $Rxy+I=Ry+I=R$, and $Ryx+J=Rx+J=R$, while (b) implies $Ixy\subset Jy\subset I$ and $Jyx\subset Ix\subset J$.

(ii) If $I$ and $J$ are maximal then, because ${R\over I},{R\over J}$ are simple $R$-modules, every nontrivial (hence surjective) $R$-homomorphism $f:{R\over I}\ra {R\over J}$ is an $R$-isomorphism (since $\ker f=0$ is the only possible choice for $\ker f$), i.e., a nontrivial $R$-homomorphism ${R\over I}\ra{R\over J}$ exists $\iff$ there exist $x,y\in R$ satisfying (a) and (b).

(iii) If $R$ is commutative, then for any $x,y\in R$ satisfying (a) and (b), we have ~$x,y,xy\in (R\backslash I)\cap(R\backslash J)=R\backslash(I\cup J)$ and moreover, we have corresponding mutually inverse units $x+I,y+I\in U({R\over I})$ and $x+J,y+J\in U({R\over J})$.
\end{thm}
\begin{proof}
It suffices to prove (i) only, since (ii) and (iii) follow immediately from (i).
{\flushleft($\Ra$)}: If elements $x,y\in R$ exist with the given properties, the following are well-defined mutually inverse $R$-homomorphisms (i.e., $f g=id_{R/J}$ and $gf=id_{R/I}$):
\bea
f:R/I\ra R/J,~r+I\mapsto rx+J,~~~~g:R/J\ra R/I,~r+J\mapsto ry+I.\nn
\eea
{\flushleft ($\La$)}: Conversely, given any well-defined mutually inverse $R$-homomorphisms
\bea
f:R/I\ra R/J,~r+I\mapsto \wt{f}(r)+J,~~~~g:R/J\ra R/I,~r+J\mapsto \wt{g}(r)+I,\nn
\eea
we have $\wt{f}(r)=r\wt{f}(1)$ and $\wt{g}(r)=r\wt{g}(1)$, and so with $x:=\wt{f}(1)$ and $y:=\wt{g}(1)$ we get
\bit
\item[(a)] $xy+I=1+I$ and $yx+J=1+J$ ~~(since $gf(1+I)=1+I$ and $fg(1+J)=1+J$),
\item[(b)] $Ix\subset J$ and $Jy\subset I$ ~~(since $r\in I$ $\Ra$ $rx\in J$, and $r\in J$ $\Ra$ $ry\in I$, for all $r\in R$). \qedhere
\eit
\end{proof}

\begin{crl}
Let $R$ be a PID and $I,J\lhd R$ distinct maximal ideals. Then there is \ul{no} $R$-isomorphism ${R\over I}\cong{R\over J}$. Hence, all non-isomorphic simple $R$-modules are of the form ${R\over Rp}$ for primes $p\in R$. (\blue{footnote}\footnote{That is, the number of irreducible representations of a PID is equal to its number of maximal ideals (or number of primes).}). In other words, if $R$ is a PID, then we have a bijection
\[
\big\{\txt{maximal ideals of $R$}\big\}\longleftrightarrow\big\{\txt{simple $R$-modules}\big\}.
\]
\end{crl}
\begin{proof}
We know $I=Rp$, $J=Rq$ for primes $p\neq q$. Suppose an $R$-isomorphism ${R\over Rp}\cong{R\over Rq}$ exists. Let $x,y\in R$ such that $Rpx\subset Rq$ and $Rqy\subset Rp$, i.e., $q|px$ and $p|qy$. Since $p,q$ are primes, it follows that $q|x$ and $p|y$. This contradicts the relations $Ry+I=R$ and $Rx+J=R$ which say $gcd(y,p)=1=gcd(x,q)$.
\end{proof}
An immediate hopeful consequence of the above result is that if $R$ is a PID, then an $R$-module $M$ might be completely expressible in terms of the $R$-modules $\{R\}\cup\left\{{R\over Rp}~|~\txt{for primes}~p\in R\right\}$, say using the product, coproduct, etc, where (as for any integral domain) $Rr\cong_RR$ for all $r\in R\backslash\{0\}$. This complete expressibility will be realized for \ul{finitely generated} modules over a PID in chapter \ref{RepsTheoryII}.

\begin{crl}\label{MaxSpecEquiv}
Let $R$ be a commutative ring and $I,J\lhd R$ prime ideals not containing each other (i.e., $I\not\subset J$ and $J\not\subset I$). Then ${R\over I}\not\cong_R{R\over J}$. In particular, for any commutative ring $R$, we have a bijection
\[
\big\{\txt{maximal ideals of $R$}\big\}\longleftrightarrow\big\{\txt{simple $R$-modules}\big\}.
\]
\end{crl}
\begin{proof}
Suppose an $R$-isomorphism ${R\over I}\cong{R\over J}$ exists. Let $x,y\in R$ such that $I(Rx)=Ix\subset J$ and $J(Ry)=Jy\subset I$. By hypotheses, this implies $Rx\subset J$ and $Ry\subset I$ (a contradiction of the relations $Ry+I=R$ and $Rx+J=R$).
\end{proof}

{\flushleft\hrulefill}

\begin{exercise}
Based on the discussion of this chapter, consider writing a \emph{fully technical} essay (say in the form of a typical section of this chapter) on what is known in the mathematics literature as \index{Sylow Theorems}{``Sylow Theorems''}.
\end{exercise}

%% file: parts/AlgebraM/RepsTheoryII.tex
\chapter{Representation Theory II: Finitely Generated Modules}\label{RepsTheoryII}
In this chapter we begin discussing representations of a ring (i.e., modules over a ring). In view of the stepwise reduction approach to representation theory mentioned earlier, the basic idea is to first express a given module in terms of modules of intermediate complexity, and then eventually in terms of irreducible (i.e., simple) modules. In particular, our main goal in this chapter is to characterize finitely generated modules (modules that are finite sums of cyclic modules) over a PID. (\blue{footnote}\footnote{Depending on the ring $R$ considered, some cyclic $R$-modules can become automatically irreducible (i.e., simple) $R$-modules.})

\section{Noetherian and Finitely Generated (FG) Modules}
\subsection{Recall of concepts and preliminary remarks}

\begin{dfn}[\blue{\small
\index{Minimal! (or simple) ideal}{Minimal (or simple) ideal},
\index{Maximal! ideal}{Maximal ideal},
\index{Prime! ideal}{Prime ideal}}]
Let $R$ be a ring. A nonzero proper ideal $0\neq I\vartriangleleft R\neq I$ is a \ul{minimal ideal} or \ul{simple ideal} if it cannot properly contain a nonzero ideal in the sense that if $I\supset J\vartriangleleft R$, then $J=0$ or $J=I$ (i.e., $RaR=I$ for all $0\neq a\in I$). A proper ideal $I\vartriangleleft R\neq I$ is a \ul{maximal ideal} if it cannot be properly contained by a proper ideal in the sense that if $I\subset J\vartriangleleft R$, then $J=I$ or $J=R$ (i.e., $RaR+I=R$ for all $a\not\in I$). A proper ideal $P\vartriangleleft R\neq P$ is a \ul{prime ideal} if for any $a,b\in R$,~ $aRb\subset P$ $\Ra$ $a\in P$ or $b\in P$ (or equivalently, for any $I,J\lhd R$,~ $IJ\subset P$ $\Ra$ $I\subset P$ or $J\subset P$).
\end{dfn}
In particular, if $R$ is commutative then, a proper ideal $P\vartriangleleft R\neq P$ is prime iff for any $a,b\in R$,~ $ab\in P$ $\Ra$ $a\in P$ or $b\in P$.

\begin{dfn}[\blue{\small
\index{Minimal! (or simple) submodule}{Minimal (or simple) submodule},
\index{Maximal! submodule}{Maximal submodule},
\index{Minimal! left ideal}{Minimal left ideal},
\index{Maximal! left ideal}{Maximal left ideal},
\index{Minimal! right ideal}{Minimal right ideal},
\index{Maximal! right ideal}{Maximal right ideal}}]
Let $M={}_RM$ be an $R$-module. A nonzero proper submodule $0\neq{}_RN\subset{}_RM\neq N$ is a \ul{minimal submodule} or \ul{simple submodule} if it cannot properly contain a nonzero submodule in the sense that if~ $_RN\supset{}_RX\subset M$, then $X=0$ or $X=N$ (i.e., $Rn=N$ for all $0\neq n\in N$).
A proper submodule $_RN\subset{}_RM\neq N$ is a \ul{maximal submodule} if it cannot be properly contained by a proper submodule in the sense that if~ $_RN\subset{}_RX\subset M$, then $X=N$ or $X=M$ (i.e., $Rm+N=M$ for all $m\not\in N$).

A minimal (resp. maximal) left ideal ${}_RI\subset R$ is a minimal (resp. maximal) submodule of $_RR$. Similarly, a minimal (resp. maximal) right ideal $I_R\subset R$ is a minimal (resp. maximal) submodule of $R_R$. Accordingly, a minimal (resp. maximal) ideal $I={}_RI_R\subset R$ is a minimal (resp. maximal) submodule of ${}_RR_R$.
\end{dfn}

\begin{dfn}[\blue{\small
\index{Simple! ring}{Simple ring},
\index{Simple! module}{Simple module}}]
A ring $R$ is a \ul{simple ring} if $0$ and $R$ are the only ideals of $R$ (i.e., $RaR=R$ for all $a\in R\backslash 0$). A nonzero module $_RS\neq 0$ is a \ul{simple module} if $0$ and $S$ are the only submodules of $S$ (i.e., $Rs=S$ for all $s\in S\backslash 0$).
\end{dfn}
By the correspondence theorem for rings, an ideal $I\vartriangleleft R$ is maximal iff $R/I$ is a simple ring. Similarly, by the correspondence theorem for modules, a submodule $_RN\subset M$ is maximal iff $M/N$ is a simple module. In particular, a left ideal $_RI\subset R$ is maximal iff $R/I$ is a simple $R$-module.

\begin{dfn}[\blue{\small
\index{Noetherian module}{Noetherian module},
\index{Left! noetherian ring}{Left noetherian ring},
\index{Right! noetherian ring}{Right noetherian ring},
\index{Noetherian ring}{Noetherian ring}}]~
\begin{enumerate}
\item A module $_RM$ is a \ul{noetherian module} if given any ascending chain of submodules $M_1\subset M_2\subset M_3\subset\cdots$, there exists $n\geq 1$ such that $M_n=M_{n+1}=M_{n+2}=\cdots$. That is, every ascending chain of submodules terminates or stabilizes.
\item A ring $R$ is a \ul{left noetherian ring} if $_RR$ (resp. \ul{right noetherian ring} if $R_R$) is noetherian.
\item A ring is a \ul{noetherian ring} if it is both left noetherian and right noetherian.
\end{enumerate}
\end{dfn}

\begin{dfn}[\blue{\index{Finitely generated! module}{Finitely generated (FG) module}}]
A module $_RM$ is a \ul{finitely generated module} (or \ul{FG module}) if $M=Rm_1+\cdots+Rm_n$ for some $m_1,...,m_n\in M$ and $n\geq 1$. That is, $_RM$ is finitely generated iff there exists an exact sequence of $R$-modules ~$0\ra K\ra R^n\ra M\ra 0$,~ where $R^n:=R\oplus R^{n-1}$.
\end{dfn}

\subsection{Results on noetherian modules}
\begin{thm}[\blue{Noetherian characterization I}]\label{NoethChrI}
A module $_RM$ is noetherian iff every nonempty collection of submodules has a maximal element with respect to inclusion.
\end{thm}
\begin{proof}
Assume $M$ is noetherian and let $\S$ be a nonempty collection of submodules of $M$. Then every nonempty chain, with respect to inclusion, in $\S$ has an upper (resp. a lower) bound since $M$ is noetherian, and so $\S$ has a maximal element by Zorn's lemma.

Conversely, if every nonempty collection of submodules has a maximal element, then in particular, any ascending chain of submodules terminates, and so $M$ is noetherian.
\end{proof}

\begin{lmm}
Let $M$ be an $R$-module and ${}_RA,{}_RB,{}_RN\subset{}_RM$ submodules. If $A\subset B$, then
\bea
A\cap N=B\cap N~~\txt{and}~~A+N=B+N~~~~\Ra~~~~A=B.\nn
\eea
\end{lmm}
\begin{proof}
Let $C:=A\cap N=B\cap N$. Then by the module isomorphism theorems, we have
\bea
\textstyle {A\over C}={A\over A\cap N}\cong{A+N\over N}={B+N\over N}\cong {B\over B\cap N}={B\over C},\nn
\eea
where by construction the isomorphism is explicitly given by {\small $f:A/C\ra B/C,~a+C\ra a+C$} (which is clearly well-defined and injective). Since $f$ is surjective, $f(A/C)=B/C$ implies $A/C=B/C$, and so $A=B$.
\end{proof}

\begin{thm}[\textcolor{blue}{Noetherian characterization II}]\label{NoethChrII}
Let $N\subset M$ be $R$-modules. $M$ is noetherian $\iff$ $N$ and $M/N$ are both noetherian.
\end{thm}
\begin{proof}
{\flushleft ($\Ra$)}
Assume $M$ is noetherian. Then every increasing sequence of submodules of $N$ is an increasing sequence of submodules of $M$, and so terminates. Hence $N$ is noetherian. Similarly, every increasing sequence of submodules $M_1/N\subset M_2/N\subset\cdots$ of $M/N$ gives (by the correspondence theorem) an increasing sequence of submodules $M_1\subset M_2\subset\cdots$ of $M$, and so terminates. Hence $M/N$ is noetherian.
{\flushleft ($\La$)}
Assume $N$ and $M/N$ are both noetherian. Let $A_1\subset A_2\subset\cdots$ be a chain of submodules of $M$. Then the following two sequences of submodules of $N$ and $M/N$ (respectively) each stabilize:
\bea
\label{NoethChIIIEq}&&\textstyle (a)~~A_1\cap N\subset A_2\cap N\subset\cdots,~~~~~~(b)~~{A_1+N\over N}\subset {A_2+N\over N}\subset\cdots
\eea
Thus for some $k$, ${A_k+N\over N}={A_{k+1}+N\over N}=\cdots$ (i.e., $A_k+N=A_{k+1}+N=\cdots$) as the stabilization point for (\ref{NoethChIIIEq})(b). Without loss of generality assume the stabilization point for (\ref{NoethChIIIEq})(a) is also $A_k\cap N=A_{k+1}\cap N=\cdots$ (otherwise, take the higher of the two). That is, $A_k+N=A_{k+1}+N=\cdots$ and $A_k\cap N=A_{k+1}\cap N=\cdots$, and since $A_k\subset A_{k+1}\subset\cdots$, the preceding lemma implies $A_k=A_{k+1}=\cdots$, i.e., the chain $A_1\subset A_2\subset\cdots$ stabilizes. Hence $M$ is noetherian.
\end{proof}

\begin{crl}\label{NthCh2Crl1}
Let {\small $0\ra A\ra B\ra C\ra 0$} be an exact sequence of $R$-modules. Then $B$ is noetherian iff $A,C$ are both noetherian.
\end{crl}

\begin{crl}\label{NthCh2Crl2}
Let $M,M'$ be R-modules. Then $M\oplus M'$ is noetherian iff $M,M'$ are both noetherian.
\end{crl}
\begin{proof}
Observe that we have an exact sequence $0\ra M\ra M\oplus M'\ra M'\ra 0$.
\end{proof}

\begin{crl}\label{NthCh2Crl3}
Let {\small $\{M_i\}_{i\in I}$} be $R$-modules. (i) If $\bigoplus_{i\in I} M_i$ is noetherian, so is $M_i$ for every $i\in I$. (ii) If $I=\{1,...,n\}$ (i.e., a finite set) then {\small $M_1\oplus \cdots\oplus M_n$} is noetherian iff each $M_i$ is noetherian.
\end{crl}
\begin{proof}
(i) For each $j\in I$, we have the exact sequence $0\ra M_j\ra\bigoplus_{i\in I}M_i\ra\bigoplus_{i\neq j}M_i\ra 0$.
(ii) This follows by induction on $n$, due to associativity of the direct sum, i.e.,
\bit
\item[]\hspace{2cm} $M_1\oplus M_2\oplus\cdots\oplus M_n=(M_1\oplus M_2\oplus\cdots\oplus M_{n-1})\oplus M_n$.\qedhere
\eit
\end{proof}

\begin{dfn}[\blue{\small
\index{Decomposable module}{Decomposable module},
\index{Indecomposable module}{Indecomposable module}}]
Let $M$ be an $R$-module. $M$ is \ul{decomposable} if $M$ a direct sum of two nonzero submodules (i.e., $M\cong N_1\oplus N_2$ for nonzero submodules $0\neq N_1,N_2\subset M$), otherwise $M$ is \ul{indecomposable} (i.e., for any nonzero submodules $0\neq N_1,N_2\subset M$ we have $M\not\cong N_1\oplus N_2$, or equivalently, for any submodules $N_1,N_2\subset M$, if $M\cong N_1\oplus N_2$ then $N_1=0$ or $N_2=0$).
\end{dfn}

\begin{crl}
If an $R$-module $M$ is noetherian, then $M=N_1\oplus\cdots\oplus N_r$ for indecomposable submodules $N_1,...,N_r\subset M$.
\end{crl}
\begin{proof}
Assume $M$ is noetherian (and that $M$ is decomposable, otherwise there is nothing to prove). Let $M=N_1\oplus M_1$, where $M_1$ is a maximal direct summand of $M$ wrt $N_1$ (so $N_1$ is indecomposable, otherwise $M_1$ would not be a maximal direct summand). This can be repeated since $M_1\cong {M\over N_1}$ is again noetherian (Theorem \ref{NoethChrII}). At the $i$th step, we have $M_i=N_{i+1}\oplus M_{i+1}$, where $M_{i+1}$ is a maximal direct summand of $M_i$ wrt $N_{i+1}$ (so $N_{i+1}$ is indecomposable), and so we get
\bea
M=N_1\oplus M_1=N_1\oplus (N_2\oplus M_2)=N_1\oplus N_2\oplus\cdots.\nn
\eea
By the noetherian property, the sequence $N_1\subsetneq N_1\oplus N_2\subsetneq N_1\oplus N_2\oplus N_3\subsetneq\cdots$ terminates at some integer $r-1$ with $N_{r-1}$ indecomposable. Set $N_r:=M_{r-1}$ in $M_r=N_{r-1}\oplus M_{r-1}$.
\end{proof}

\begin{lmm} Let $M$ be an $R$-module and $f\in End_R(M):=Hom_R(M,M)$. If $M$ is Noetherian and $f$ is surjective, then $f$ is an isomorphism.
\end{lmm}
\begin{proof}
Let $f^{(n)}:=f\circ f^{(n-1)}=f^{(n-1)}\circ f$ be the $n$-fold composition of $f$, where $f^{(0)}:=id_M$. Since $f$ is surjective, so is $f^{(n)}$. Also, ~{\small $\ker f^{(n+1)}=\ker (f\circ f^{(n)})=(f^{(n)})^{-1}(\ker f)\supset\ker f^{(n)}$},~ and so
\bea
&&\ker f^{(n)}\subset\ker f^{(n+1)},~~\sr{(s1)}{\Ra}~~\ker f^{(N)}=\ker f^{(N+1)}=\ker f^{(N+2)}=\cdots~~\txt{for some}~~N,\nn\\
&&~~\Ra~~f\left(\ker f^{(N)}\right)=f\left(\ker f^{(N+1)}\right),~~\sr{(s2)}{\Ra}~~\ker f^{(N-1)}=\ker f^{(N)},~~\cdots,\nn\\
&&~~\sr{\txt{induction}}{\Ra}~~0=\ker f^{(0)}=\ker f^{(1)}=\ker f,\nn
\eea
where step (s1) is the noetherian property and step (s2) the surjectivity of $f$. (\blue{footnote}\footnote{Because $f$ is surjective, $f\big(f^{-1}(m)\big)=m$ for all $m\in M$.}).
\end{proof}

\subsection{Results on FG modules}
\begin{thm}[\textcolor{blue}{Noetherian-FG characterization}]\label{NoethFGChr}
A module $M$ is noetherian iff every submodule of $M$ is finitely generated.
\end{thm}
\begin{proof}
{\flushleft ($\Ra$)} Assume $M$ is noetherian, and let $N\subset M$ be a submodule. If $N$ is not finitely generated, then for any $x_1\in N$, some $x_2\in N\backslash Rx_1$, some $x_3\in N\backslash(Rx_1+Rx_2)$, some $x_4\in N\backslash(Rx_1+Rx_2+Rx_3)$, and so on. We get the infinite sequence $Rx_1\subsetneq Rx_1+Rx_2\subsetneq \cdots$ (a contradiction).

{\flushleft ($\La$)} Conversely, assume every submodule of $M$ is finitely generated. Let $M_1\subset M_2\subset\cdots$ be a chain of submodules of $M$. Then we know $N:=\bigcup_{i=1}^\infty M_i\subset M$ is a submodule. Since $N=Rx_1+\cdots+Rx_k$ is finitely generated, we have $x_1,\cdots, x_k\in M_t$ for some $t\geq 1$, and so $N\subset M_t$. Hence $M_t=M_{t+1}=M_{t+2}=\cdots$.
\end{proof}

\begin{crl}[\blue{A PID is noetherian}]
\end{crl}

\begin{thm}[\textcolor{blue}{A finitely generated module over a left noetherian ring is noetherian}]
Let $M$ be an $R$-module. If $R$ is left noetherian and $M$ is finitely generated, then $M$ is noetherian. (\blue{footnote}\footnote{Here, the ring need only be left-noetherian, since we assume the modules are left modules.})
 \end{thm}
\begin{proof}
Recall that $M$ is FG iff we have an exact sequence of $R$-modules ~$0\ra K\ra R^n\ra M\ra 0$.
\end{proof}

\begin{crl}[\blue{A finitely generated module over a PID is noetherian}]\label{FgOvPidNoe}
\end{crl}

\begin{thm}[\textcolor{blue}{Maximal submodule for a nonzero finitely generated module}]
Every nonzero finitely generated module has a maximal submodule.
\end{thm}
\begin{proof}
Let $M$ be a nonzero finitely generated module. As seen earlier under ``Linear Algebra of Modules'', because the union of an increasing sequence of proper submodules of $M$ is itself a proper submodule of $M$, the set of proper submodules of $M$ is a regular poset under inclusion (in the sense every chain has an upper bound), and so has a maximal element by Zorn's lemma.
\end{proof}

\begin{crl}[\textcolor{blue}{Containing maximal submodule in a finitely generated module}]
If $M$ is a finitely generated module, then every proper submodule $W\subsetneq M$ is contained in a maximal submodule.
\end{crl}
\begin{proof}
Since $M=Rm_1+\cdots+Rm_n$ is finitely generated, so is ${M\over W}={Rm_1+\cdots+Rm_n\over W}={Rm_1+W\over W}+\cdots+{Rm_n+W\over W}=R(m_1+W)+\cdots+R(m_n+W)$. Thus $M/W$ contains a maximal submodule (preceding theorem), say $W'$. By the correspondence theorem, $W'=N/W$ where $N\subset M$ is a submodule such that $W\subset N$.

To see that $N\subset M$ is a maximal submodule, observe that by the 3rd isomorphism theorem for modules, we have an isomorphism of simple $R$-modules ~${M\over N}\cong{M/W\over N/W}={M/W\over W'}$.
\end{proof}

\section{Some Intermediate Results}
\subsection{Sectioning (or right-splitting) of a module homomorphism}
\begin{dfn}[\textcolor{blue}{Recall: \index{Internal direct sum}{Internal direct sum}}]
Recall that if $M$ is an $R$-module and $A,B\subset M$ submodules, then the map $A\oplus B\ra M,~(a,b)\mapsto a+b$ is an $R$-isomorphism $\iff$ $A+B=M$ and $A\cap B=\{0\}$, in which case we say ``$M$ is the \ul{internal direct sum} of $A$ and $B$'', written as ``$M=A\oplus B$''.
\end{dfn}

\begin{dfn}[\textcolor{blue}{\index{Section of an R-homomorphism}{Section of an R-homomorphism}}]
Let $g:M\ra N$ be an R-homomorphism. A section of $g$ is an R-homomorphism $s:N\ra M$ such that $g\circ s=id_N:N\ra N,~y\mapsto y$ (i.e., a right-inverse of $g$). Equivalently, a section of $g$ is an R-homomorphism of the form $s:N\ra M,~y\mapsto s(y)\in g^{-1}(y)$.
\bea\bt
        &              & s(N)\ar[d,hook]\ar[from=rr,,"s"',"\cong"]             && N\ar[dll,dashed,near end,bend left=10,"s"']\ar[d,"id_N"] &\\
0\ar[r] & \ker g\ar[r] & M\ar[rr,near end,"g"] && N\ar[r] & 0
\et\nn
\eea
\end{dfn}

\begin{rmk}
Let $g:M\ra N$ be an R-homomorphism. If a section $s:N\ra M$ of $g$ exists, then (i) it makes $g$ surjective, (ii) it is injective, and so (iii)  $N\cong s(N)$, where the isomorphism is given by the maps $s:N\ra s(N)$ and $g|_{s(N)}:s(N)\ra N$, which satisfy $g|_{s(N)}\circ s=id_N$ and $s\circ g|_{s(N)}=id_{s(N)}$.
\end{rmk}

\begin{lmm}[\blue{Existence of a section of an R-homomorphism}]
Let $g:M\ra N$ be an R-homomorphism. A section of $g$ exists $\iff$ the sequence ~$0\ra \ker g\ra M\sr{g}{\ral}N\ra 0$ is split-exact. (\blue{footnote}\footnote{Recall that an epimorphism with a right-inverse (section) is also called a \ul{split epimorphism}.}). Equivalently, a section $s:N\ra M$ exists $\iff$ (i) $g$ is surjective and (ii) $M=\ker g\oplus N'$ for some submodule $N'\subset M$. In this case, $g|_{N'}:N'\ra N$ is an isomorphism.
\end{lmm}
\begin{proof}
{\flushleft ($\Ra$)} Assume a section $s:N\ra M$ exists. Then it is clear that $g$ is surjective. For $m\in M$,
\bea
m=(m+sg(m))-sg(m)\in\ker g+s(N),~~\Ra~~M=\ker g+s(N),\nn
\eea
where $m+sg(m)\in\ker g$, since $g(m+sg(m))=g(m)+gsg(m)=g(m)-g(m)=0$. Also,
{\small\begin{align}
\ker g\cap s(N)=\{m\in M:g(m)=0,~m=s(n),~\txt{some}~n\in N\}=0,~~\Ra~~M=\ker g\oplus s(N).\nn
\end{align}}
{\flushleft ($\La$)} Conversely, assume $g$ is surjective and $M=\ker g\oplus N'$ for some submodule $N'\subset M$. Then the exact sequence $0\ra\ker g\ra M\sr{g}{\ral}N\ra 0$ is split, and so we get the map  $s:N\ra M,~n\mapsto s(n)\in g^{-1}(n)\cap N'$ (which is well-defined because $g|_{N'}:N'\ra N$ is bijective, and also clearly satisfies $g\circ s=id_N$), i.e.,
\bea
s:N\ra M,~g(m)\mapsto sg(m)\in g^{-1}(g(m))\cap N'=(m+\ker g)\cap N'.\nn
\eea
Finally, $s$ is an $R$-homomorphism because (for $r\in R,n\in N$) we have $gs(rn)=rn=rgs(n)=g(rs(n))$, and so $s(rn)=rs(n)$ since the restriction $g|_{N'}$ is injective.
\end{proof}

\begin{crl}\label{section-existence-lmm}
Let $P$ be a projective $R$-module (e.g., a free $R$-module) and $g:M\ra P$ a surjective R-homomorphism. Then a section $s:P\ra M$ of $g$ exists. Hence ~$M=\ker g\oplus s(P)\cong \ker g\oplus P$.
\end{crl}

\subsection{Torsion in modules and decomposition}
\begin{dfn}[\textcolor{blue}{
\index{Annihilator left ideal}{Annihilator left ideal},
\index{Annihilator submodule}{Annihilator submodule},
\index{Torsion! (or zero-division) element}{Torsion (or zero-division) element},
\index{Torsion! subset}{Torsion subset},
\index{Torsion! module}{Torsion module},
\index{Torsion-free module}{Torsion-free module}}]
Let $R$ be any ring, $M$ an $R$-module, and $A\subset M$ a subset. The \ul{$R$-annihilator} of $A$ (i.e., annihilator of $A$ in $R$) is the left ideal (\blue{footnote}\footnote{If $N\subset M$ is a submodule, then it is clear that $Ann_R(N)\subset R$ is an ideal.})
\[
Ann_R(A):=\big\{r\in R:rA=\{0\}\big\}=\{r\in R:ra=0~\txt{for all}~a\in A\}\subset R.
\]
Similarly, if $B\subset R$, the \ul{$M$-annihilator} of $B$ (i.e., annihilator of $B$ in $M$) is the subgroup (\blue{footnote}\footnote{If $B\subset R$ is a right ideal or $B\subset Z(R)$, then it is clear that $Ann_M(B)\subset M$ is a submodule.})
\[
Ann_M(B):=\big\{m\in M:Bm=\{0\}\big\}=\{m\in M:bm=0~\txt{for all}~b\in B\}\subset M.
\]
An element $m\in M$ is a \ul{torsion element} of $M$ if $Ann_R(m)\neq\{0\}$.  (\blue{footnote}\footnote{That is, $m\in M$ is \ul{torsion} if $m$ is killed by a nonzero scalar, in the sense that $rm=0$ for some $r\in R\backslash 0$. Also, observe that if $R$ is commutative, then \ul{torsion elements} in an R-module are analogous to \ul{zerodivisors} in $R$. In particular, if $R$ is commutative then $T_R(R)=\{\txt{zerodivisors in $R$}\}$, and so $R$ is an ID $\iff$ the module $_RR$ is torsion-free.}). The \ul{torsion subset} of $M$ is
\[
T(M):=T_R(M):=\{\txt{torsion elements of $M$}\}.
\]
The module $M$ is a \ul{torsion module} if $T_R(M)=M$ (resp. a \ul{torsion-free module} if $T_R(M)=0$).
\end{dfn}
When the $R$-module $M$ is being considered over more than one ring (e.g., when the fact that $M={}_\Integer M$ is important), the word ``torsion'' as used above needs to be preceded by the name of the ring under consideration e.g., ``$R$-torsion'', ``$\Integer$-torsion'', etc. Consequently, the subscripts $R$ and $M$ in the expressions $Ann_R(A)$, $Ann_M(B)$, and $T_R(M)$ are always important. Nevertheless, for convenience, we will often write $T_R(M)$ simply as $T(M)$ when the ring $R$ under consideration is well understood.

If $V$ is a $k$-vector space, then $T_k(V)=0$, i.e., vector spaces are torsion-free (because every nonzero element $a\in k$ is a unit, and so if $v\in V$ then $av=0$ $\iff$ $v=a^{-1}av=0$). On the other hand, if the $R$-module $M={}_\Integer M$ is viewed as an abelian group, then ~$T_\Integer(M)=\{\txt{elements of finite order in $M$}\}$. Also, the 1st isomorphism theorem implies cyclic $R$-submodules of $M$ satisfy the following: For any $m\in M$,
\bea
\textstyle{R\over Ann_R(m)}={R\over\ker(R\ra Rm)}\cong\im(R\ra Rm)=Rm.\nn
\eea

\begin{prp}[\textcolor{blue}{
\index{Torsion! submodule}{Torsion submodule},
\index{Torsion-free quotient module}{Torsion-free quotient module}}]\label{TorFreeIso1}
Let $R$ be an ID and $M$ an $R$-module. Then (i) $T(M)$ is a submodule of $M$, (ii) $T(M)$ is torsion, and (iii) $M/T(M)$ is torsion-free.
\end{prp}
\begin{proof}
{\flushleft (i)} It is clear that $0\in T(M)$. Let $m,m'\in T(M)$. Then $rm=0,~r'm'=0$ for some $r,r'\in R\backslash 0$, and so $rr'(m-m')=r'(rm)+r(r'm)=r'0+r0=0$. Since $rr'\neq 0$ (as $R$ is an ID) we have $m-m'\in T(M)$. Also, $r(Rm)=R(rm)=R0=0$, i.e., $Rm\in T(M)$. Hence $T(M)\subset M$ is a submodule.
{\flushleft (ii)} For any submodule $N\subset M$, it is clear that $T(N)=N\cap T(M)$, and so in particular,
\bea
T(T(M))=T(M)\cap T(M)=T(M).\nn
\eea
{\flushleft (iii)} Let $m+T(M)\in T(M/T(M))$. Then $rm+T(M)=T(M)$, i.e., $rm\in T(M)$, for some $r\in R\backslash 0$. This means $srm=0$ for some $s\in R\backslash 0$. Because $R$ is an ID, $sr\neq 0$, and so $m\in T(M)$. Hence $m+T(M)=T(M)$.
\end{proof}

\begin{rmk}\label{TorFreeRmk1}
 Let $R$ be an ID, $M$ an $R$-module, and $m,m'\not\in T(M)$. Then we have the following.
 \bit[leftmargin=0.7cm]
\item[(1)] \ul{$Rm$ is a free module with basis $\{m\}$}: For any $r\in R$, we have $rm=0$ iff $r=0$ (since $m\not\in T(M)$).\\
\ul{If $Rm\cap Rm'=\{0\}$, then $Rm+Rm'$ is a free module with basis $\{m,m'\}$.}
\item[(2)] \ul{If $M/T(M)$ is free then $M\cong T(M)\oplus M/T(M)$}: Recall that we have the SES~ $0\ra T(M)\ra M\ra M/T(M)\ra 0$, which is split if $M/T(M)$ is free (hence projective).
\item[(3)] \ul{If (i) $R$ is a PID and (ii) $M$ if FG (hence $M/T(M)$ is FG), then $M/T(M)$ is free}: Because, we will see in Theorem \ref{TorFreeIso2} that a FG torsion-free module over a PID is free.
\item[(4)] \ul{If $M$ is free then $M$ is torsion-free}: To see why, let $B$ be a basis for $M$. Then $x\in M\backslash 0$ can be written as $x=r_1b_1+\cdots+r_nb_n$ for some distinct $b_i\in B$ and $r_i\in R\backslash 0$. Therefore, for any $a\in R$, if $ax=ar_1b_1+\cdots+ar_nb_n=0$, then $ar_i=0$ for all $i$, and so $a=0$ (since $R$ is an ID).
\eit
\end{rmk}

\begin{thm}[\textcolor{blue}{(Over an ID) ``FG''  $\Ra$ ``Torsion Quotient Module''}]\label{TorIso1}
Let $R$ be an ID and $M$ a FG (finitely generated) R-module. Then there exists a free submodule $F\subset M$ such that $M/F$ is torsion.
\end{thm}
\begin{proof}
Let $M=Rx_1+\cdots+Rx_t$. If $M$ is torsion, then we can set $F:=0$. So, assume $M$ is not torsion, i.e., $T(M)\subsetneq M$. Then by the preceding remark, we have the \emph{nonempty} set $P:=\{RB:\txt{linearly independent}~B\subset\{x_1,...,x_t\}\}$ which has a maximal element (wrt inclusion) by Zorn's lemma.

Thus, we can choose a maximal subset $B\subset\{x_1,...,x_t\}$ that forms a basis for a free submodule $F:=R^{\langle B\rangle}\subset M$, i.e., $B$ is a maximal linearly independent subset of $\{x_1,...,x_t\}$. Without loss of generality, take this subset to be $B:=\{x_1,...,x_n\}$, where $n\leq t$. Then $F=Rx_1+\cdots+Rx_n$.

Let $1\leq i\leq n$ (i.e., $x_i\in B$). Since $F=0_{M/F}\in T(M/F)$, it is clear that {\small $x_i+F=F\in T(M/F)$}. Next, let $n<i\leq t$. Then $B\cup\{x_i\}=\{x_1,...,x_n,x_i\}$ is a dependent set (by the choice of $B$), and so for some $r_1,...,r_n\in R$ not all zero, and $r\in R\backslash 0$, we have $\sum_{j=1}^nr_jx_j+rx_i=0$. This implies $rx_i=-\sum_{j=1}^nr_jx_j\in F$, and so $r(x_i+F)=rx_i+F=F$, i.e., $x_i+F\in T(M/F)$. Hence $M/F=T(M/F)$.
\end{proof}

\begin{dfn}[\textcolor{blue}{Recall: \index{Rank! of a free IBN-module}{Rank of a free IBN-module}}]
Let $R$ be an IBN ring (e.g., a commutative ring), and $F$ a free R-module. The \ul{rank} of $F$ is the cardinality of any basis for $F$. (\blue{footnote}\footnote{Recall that $R$ is an IBN ring iff for any integers $m,n\geq 1$, $R^m\cong R^n$ $\Ra$ $m=n$.}). In particular, if $F$ is a finitely generated free R-module, the rank of $F$ is the number of elements in any basis for $F$.
\end{dfn}

\begin{lmm}[\textcolor{blue}{(Over a PID) ``FG Free'' is ``Hereditary''}]
Let $R$ be a PID, and $F$ a finitely generated (FG) free R-module. Then every submodule $_RN\subset F$ is also a finitely generated (FG) free R-module with $\txt{rank}(N)\leq \txt{rank}(F)$.
\end{lmm}
\begin{proof}
We know $N$ is FG by Corollary \ref{FgOvPidNoe} (but this also emerges from this proof). Let $n:=\rank(F)$. Without loss of generality, set $F:=R^n$. We will proceed by induction on $n$. If $n=1$, i.e., $F=R$, then submodules are either $0$ or $Ra,~a\neq 0$ (which are FG and Free on $\emptyset$ or one generator). So assume $n>1$, and consider the R-homomorphism ~$\pi_n:R^n\ra R,~(a_1,...,a_n)\mapsto a_n$, whose kernel is
\bea
\ker\pi_n=R^{n-1}\oplus 0~\subset~R^n.\nn
\eea
Consider the restriction $g:=\pi_n|_N:N\ra \pi(N)$, which is a surjective R-homomorphism with kernel
\bea
\ker g=\ker\pi_n|_N=N\cap\ker\pi_n=N\cap(R^{n-1}\oplus 0)\subset R^n,\nn
\eea
which by the induction hypothesis is FG and Free with $\txt{rank}(\ker g)\leq n-1$. Also, $\pi_n(N)\subset R$ is FG and Free on $\emptyset$ or one generator.

Since $\pi(N)$ is free, $g$ has a section $s:\pi_n(N)\ra N$, and so we have the direct sum decomposition
\bea
N=\ker g\oplus s(\pi_n(N))=\big(N\cap(R^{n-1}\oplus 0)\big)\oplus s(\pi_n(N)),\nn
\eea
where $\ker g$ is FG and Free with $\txt{rank}(\ker g)\leq n-1$ and $s(\pi_n(N))$ is FG and Free with $\txt{rank}(s(\pi_n(N)))\leq 1$. Hence $N$ is FG and Free with $\txt{rank}(N)\leq(n-1)+1=n$.
\end{proof}

\begin{thm}[\textcolor{blue}{(Over a PID) ``Free'' is ``Hereditary'': \cite[Theorem 8.9, pp 639-640]{rotman2010}}]\label{FrModPidHered}
Let $R$ be a PID, and $F$ a free R-module. Every submodule $_RN\subset F$ is also a free R-module with $\txt{rank}(N)\leq\txt{rank}(F)$.
\end{thm}
\begin{proof}
If $\rank(F)<\infty$ (i.e., $F$ is FG) the result holds by the preceding lemma. So assume $F$ is not FG, i.e., $\rank(F)=\infty$. By the earlier rank theorems, if $N$ is free, then we necessarily have $\rank(N)\leq\rank(F)$. Thus, it remains only to show $N$ is free. Let $\{f_\al\}_{\al\in A}$ be an infinite basis for $F$.

If $N$ is FG then $N=Rx_1+\cdots+Rx_n$ (where $x_j=\sum_{i=1}^{k_j}r_{ij}f_{\al_{ij}}$, $j=1,...,n$) is free as a submodule of the FG free module $\sum_{j=1}^n\sum_{i=1}^{k_j}Rf_{\al_{ij}}$. So, further assume $N$ is not FG.

By the axiom of choice, assume wlog that the index set $A$ in the basis $\{f_\al\}_{\al\in A}$ for $F$ is well-ordered (i.e., it is linearly ordered and every nonempty subset has a least element). For each $\al\in A$, let
\bea
\textstyle F'_\al:=\sum\limits_{\beta<\al}Rf_\beta,~~F_\al:=\sum\limits_{\beta\leq\al}Rf_\beta=F'_\al\oplus Rf_\al,~~N'_\al:=N\cap F'_\al,~~N_\al:=N\cap F_\al,\nn
\eea
in terms of which we have $F=\bigcup_{\al\in A}F_\al$ and $N=\bigcup_{\al\in A}N_\al$. Since {\small $N'_\al=N\cap F'_\al=N\cap F_\al\cap F'_\al=N_\al\cap F'_\al$},
\bea
\textstyle {N_\al\over N'_\al}={N_\al\over N_\al\cap F'_\al}\cong{N_\al+F'_\al\over F'_\al}\subset{F_\al\over F'_\al}\cong Rf_\al\cong R,\nn
\eea
which is projective and so shows the exact sequence $0\ra N'_\al\ra N_\al\sr{\pi_\al}{\ral} N_\al/N'_\al\ra 0$ is split, with $N_\al=N'_\al\oplus C$, where $C\cong N_\al/N'_\al\subset Rf_\al\cong R$ is free on $\emptyset$ or on one generator. Therefore, either
\bea
N_\al=N'_\al=N'_\al\oplus R0~~~~\txt{or}~~~~N_\al=N'_\al\oplus Rn_\al\cong N'_\al\oplus R~~\txt{for some nonzero}~~n_\al\in N_\al=N\cap F_\al.\nn
\eea

Since $A$ is well-ordered, for any $f\in F$, there exists a least index $\al(f)\in A$ such that $f\in F_\al$ for all $\al\geq\al(f)$, i.e., $\al(f):=\min\{\al\in A~|~f\in F_\al\}$, and so $f\not\in F_\al$ for all $\al<\al(f)$. That is, if $f\in F_\al$, then $\al(f)\leq\al$. In particular, if $f'\in F'_\al\subsetneq F_\al=F'_\al\oplus Rf_\al$ then $\al(f')<\al$. (We can also express this briefly as ~$f\not\in F_{\al<\al(f)}$,~ or equivalently, as ~$f\in F_{\al\geq\al(f)}$.)

Consider the index set $A^\ast:=\{\al\in A:~N_\al\cong N'_\al\oplus Rn_\al\}$ and the submodule ~$N^\ast:=\sum_{\al\in A^\ast}Rn_\al\subset N$, ~where each $n_\al\in N_\al=N\cap F_\al$ is as selected above such that $Rn_\al\cong R$.

Suppose $N^\ast\subsetneq N$. Let ${\gamma}:=\min\{\al(n)~|~n\in N\backslash N^\ast\}=\al(n_{\min})$, for some $n_{\min}\in N\backslash N^\ast$. Pick any
\[
\wt{n}\in (N\backslash N^\ast)\cap F_{{\gamma}}=(N\cap F_{{\gamma}})\backslash N^\ast=N_{{\gamma}}\backslash N^\ast,~~~~\txt{i.e., ~$\al(\wt{n})\leq{\gamma}$~ and ~$\wt{n}\in N\backslash N^\ast$},
\]
which means $\al(\wt{n})={\gamma}$ by the minimality of ${\gamma}$ in the set $\{\al(n)~|~n\in N\backslash N^\ast\}$. Since $\wt{n}\in N'_{{\gamma}}\subset F'_{{\gamma}}$ would imply $\al(\wt{n})<{\gamma}$ (a contradiction), we see that $\wt{n}\not\in N'_{{\gamma}}$, i.e., $\wt{n}\in N_\gamma\backslash (N^\ast\cap N'_\gamma)$. This means we must have $N_\gamma=N'_\gamma\oplus Rn_\gamma$ (i.e., $N_\gamma\neq N'_\gamma$ since $n_\gamma\neq 0$), and so
\bea
&&\wt{n}\in N_\gamma=N'_\gamma\oplus Rn_\gamma~~\Ra~~\wt{n}=\wt{n}'+\wt{r}n_\gamma~~\txt{for some}~~\wt{n}'\in N'_\gamma,~~\wt{r}\in R\backslash 0.\nn
\eea
where $\wt{n}'=\wt{n}-\wt{r}n_\gamma\in N'_\gamma\backslash N^\ast\subset(N\backslash N^\ast)\cap F_\gamma$ (since $\wt{n}'\in N^\ast$ would imply $\wt{n}=\wt{n}'+\wt{r}n_\gamma\in N^\ast+Rn_{\gamma}=N^\ast$). But because $N'_\gamma\subset F'_\gamma\subsetneq F_\gamma$, this implies $\al(\wt{n}')<\gamma$, which  contradicts the minimality of $\gamma$ in the set $\{\al(n)~|~n\in N\backslash N^\ast\}$. Hence $N^\ast=N$.

It remains to show that $\{n_\al\}_{\al\in A}\subset N^\ast=N$ is linearly independent (hence a basis for $N$). Consider any linear equation ~$LE:r_1n_{\al_1}+\cdots+r_tn_{\al_t}=0$~ in $N^\ast$,~ where wlog $\al_1<\al_2<\cdots<\al_t$, and so by construction $LE$ is a linear equation in $N_{\al_t}=N'_{\al_t}\oplus Rn_{\al_t}\subset N^\ast$. If $r_t\neq 0$, then $r_tn_{\al_t}=-r_1n_{\al_1}-\cdots-r_{t-1}n_{\al_{t-1}}\in N'_{\al_t}\cap Rn_{\al_t}=0$ (a contradiction). Hence, by induction on $t$, we get ~$r_1=r_2=\cdots=r_t=0$.
\end{proof}

\begin{crl}[\textcolor{blue}{(Over a PID) ``Projective'' $\iff$ ``Free''}]
Let $R$ be a PID. Then an $R$-module is projective $\iff$ free.
\end{crl}
\begin{proof}
An $R$-module is projective iff a direction summand (hence iff a submodule) of a free $R$-module.
\end{proof}

\begin{thm}[\textcolor{blue}{(Over a PID) ``FG Torsion-Free'' $\iff$ ``FG Free''}]\label{TorFreeIso2}
Let $R$ be a PID and $M$ an $R$-module. If $M$ is FG and torsion-free, then $M$ is free. (The converse is also true by Remark \ref{TorFreeRmk1}(4).)
\end{thm}
\begin{proof}
Since $M$ is FG, by Theorem \ref{TorIso1} there exits a free submodule $F\subset M$ such that
\bea
\rank(F)<\infty~~~~\txt{and}~~~~T(M/F)=M/F.\nn
\eea
Let $M=Rx_1+\cdots+Rx_n$. For each $i$, since $x_i+F\in M/F=T(M/F)$, there exists $r_i\in R\backslash 0$ such that
\bea
r_i(x_i+F)=r_ix_i+F=F,~~~~\txt{i.e.,}~~~~r_ix_i\in F.\nn
\eea
Let $r:=r_1r_2\cdots r_n$. Then $rx_i=(\prod_{j\neq i}r_j)r_ix_i\in F$ for all $i$, and so ~$rM\subset F$. Since $rM\subset F$ and $F$ is FG and free, by Theorem \ref{FrModPidHered}, $rM\subset F\subset M$ is also free with finite rank. We get the surjective $R$-linear map
\bea
\al:M\ra rM,~m\mapsto rm,\nn
\eea
which is injective (since $M$ is torsion-free), and hence an isomorphism, i.e.,
\[
M\cong rM\subset F\subset M. \qedhere
\]
\end{proof}

\begin{thm}[\textcolor{blue}{(Over a PID) ``FG'' is ``Decomposable into Torsion-plus-Free''}]\label{TorFreeDec}
Let $R$ be a PID.
\bit[leftmargin=0.8cm]
\item[(1)] If $_RM$ is FG, then $M=T(M)\oplus F$, where $F\cong M/T(M)$ is FG and free. (\blue{footnote}\footnote{See also Remark \ref{TorFreeRmk1} parts (2) and (3).})
\item[(2)] Let $M=T(M)\oplus F$ and $N=T(N)\oplus G$ be FG $R$-modules, where $F$ and $G$ are each FG and free of finite rank. If ${}_RM\cong{}_RN$, then $T(M)\cong T(N)$ and $\txt{rank}(F)=\txt{rank}(G)$. (\blue{footnote}\footnote{The converse is also trivially true. Thus if $M,N$ are FG modules over a PID then $M\cong N$ if and only if $T(M)\cong T(N)$ and $M/T(M)\cong N/T(N)$.})
\eit
\end{thm}
\begin{proof}
{\flushleft (1)} Since $M$ is finitely generated, so is $M/T(M)$. Since $M/T(M)$ is also torsion-free, $M/T(M)$ is free of finite rank (by Theorem \ref{TorFreeIso2}). Thus the exact sequence $0\ra T(M)\ra M\sr{\pi}{\ral} M/T(M)\ra 0$ is split, with a map $s:M/T(M)\ra M$ such that $M=\ker \pi\oplus s(M/T(M))=T(M)\oplus F$, where $F=s(M/T(M))\cong M/T(M)$ is free of finite rank.
{\flushleft (2)} Let $\al:M\ra N$ be an isomorphism. If $m\in T(M)$, then there exists $r\in R\backslash 0$ such that $rm=0$, which implies $r\al(m)=\al(rm)=\al(0)=0$, which implies $\al(m)\in T(N)$, and so $\al(T(M))\subset T(N)$. Similarly, $\al^{-1}(T(N))\subset T(M)$, which implies $T(N)\subset \al(T(M))$. That is, $\al(T(M))=T(N)$. So $T(M)\cong T(N)$ via
\bea
\al|_{T(M)}:T(M)\ra T(N).\nn
\eea
Consider the R-homomorphism ~$\phi=\pi\circ\al:M\sr{\al}{\ral}N\sr{\pi}{\ral} N/T(N)$.~ Then
\bea
\ker\phi=\phi^{-1}(0_{N/T(N)})=\al^{-1}(\pi^{-1}(0_{N/T(N)}))=\al^{-1}(T(N))=T(M),\nn
\eea
and so by the 1st isomorphism theorem, $M/F(M)\cong N/T(N)$. Hence
\bea
F\cong M/T(M)\cong N/T(N)\cong G,~~~~\Ra~~~~\txt{rank}(F)=\txt{rank}(G).\nn \hspace{2cm}\qedhere
\eea
\end{proof}

\begin{dfn}[\textcolor{blue}{
\index{Dual module}{Dual module},
\index{Dual basis of a basis}{Dual basis of a basis}}]
Let $R$ be a ring and $M$ an $R$-module. The \ul{dual module} of $M$ is the abelian group $M^\ast:=\txt{Hom}_R(M,R)$ viewed as an $R$-module with scalar multiplication $R\times M^\ast\ra M^\ast,~(r,f)\mapsto rf$ given by $(rf)(m):=r(f(m))$ for all $m\in M$.

Now suppose $M$ is free with a basis $B:=\{y_i\}_{i\in I}$, in which case
\begin{align}
&\textstyle M=\sum_{i\in I}Ry_i:=\big\{\sum_{i\in I}r_iy_i:~r_i\in I,~r_i=0~\txt{a.e.f.}~i\in I\big\}\cong\bigoplus_{i\in I}R,~~~~\txt{and}\nn\\
&\textstyle M^\ast:=Hom_R(M,R)\cong Hom_R\big(\bigoplus _{i\in I}R_i,R\big)\cong \prod_{i\in I}Hom_R(R,R)\cong\prod_{i\in I}R.\nn
\end{align}
The \ul{dual basis} of $B$ is the linearly independent set $B^\ast:=\{y^\ast_i\}_{i\in I}\subset M^\ast$ given by
\[
\textstyle y^\ast_i:M\ra R,~~\sum_{j\in I}r_jy_j\mapsto r_i,~~~~\txt{i.e.,}~~~~y_i^\ast(y_j)=\delta_{ij},
\]
which is well-defined because $\{y_i\}_{i\in I}$ is a linearly independent set in $M$. (\blue{footnote}\footnote{The dual basis is a basis for $M^\ast$ only when $|B|$ is finite (in which case $M^\ast\cong R^{|B|}$ is a finite product).})
\end{dfn}
\begin{rmk}\label{TorFreeRmk2}
Let $R$ be a ring and $M$ a free $R$-module with a basis $B:=\{y_i\}_{i\in I}$. Consider the dual module $M^\ast$ and the dual basis $B^\ast:=\{y_i^\ast\}_{i\in I}\subset M^\ast$ of $B$. We have the following observations:
\bit[leftmargin=0.9cm]
\item[(i)] For any $m=\sum_{j=I}r_jy_j\in M$, applying $y^\ast_i$ on both sides gives $r_i=y^\ast_i(m)$, and so the $I$-expansion of $m$ becomes ~$m=\sum_{j\in I}y^\ast_j(m)y_j$.
\item[(ii)] Any $f\in M^\ast$ is uniquely determined by the values $\{f(y_i)\}_{i\in I}\in\prod_{i\in I}R$ in the form
\[
\textstyle f:M\ra R,~~\sum_i r_iy_i\mapsto\sum_i r_if(y_i),
\]
and so can be uniquely expressed in terms of the dual basis as $f=\sum_{i\in I} f(y_i)y^\ast_i$. (\blue{footnote}\footnote{The sum makes since because it is a finite sum on every element of $M$.}).
 That is, if $M$ is free with a basis $B=\{y_i\}_{i\in I}$, then with the dual basis $B^\ast=\{y^\ast_i\}_{i\in I}$ of $B$, we have
\[
\textstyle M\cong\bigoplus_{i\in I}R\cong\Span_R(B^\ast)\subset M^\ast\cong\prod_{i\in I}R.
\] In particular, if $n:=|B|<\infty$, then $M\cong R^n\cong M^\ast$ (where $R^n\cong R^{n'}$ does not imply $n=n'$, unless $R$ is an IBN ring).
\eit
\end{rmk}

\subsection{Coprime ideals, division, and remainder theorems}
\begin{dfn}[\textcolor{blue}{
\index{Relatively! prime (or coprime) ideals}{Relatively prime (or coprime) ideals}}]
Let $R$ be a ring and $I,J\vartriangleleft R$ ideals. Then $I$ and $J$ are \ul{relatively prime} (or \ul{coprime}) if $I+J=R$.
\end{dfn}
Recall that if $R$ is a PID, then for any $a,b,d\in R$, we have $\txt{gcd}(a,b)=d$ iff $Ra+Rb=Rd$, and so $Ra$ and $Rb$ are coprime iff $gcd(a,b)=1$.
\begin{notation}[\blue{\index{Modulo convention in a ring}{Modulo convention in a ring}}]
Let $R$ be a ring, $I\vartriangleleft R$ an ideal, and $a,b\in R$ elements. Then ``$a+I=b+I$'' (i.e., ``$a-b\in I$'') is sometimes denoted by ``$a=b~(\modf I)$''.
\end{notation}
We will not make explicit use of the modulo notation but it is used in the literature.

\begin{lmm}\label{CoprimeIdProd}
Let $R$ be a ring. If $I_1,...,I_n\vartriangleleft R$ are coprime in pairs (i.e., $I_i+I_j=R$ for all $i\neq j$), then for each $k=1,...,n$, the ideals $I_k$ and $I_1\cap\cdots\cap I_{k-1}\cap I_{k+1}\cap\cdots\cap I_n$ are coprime.
\end{lmm}
\begin{proof}
By symmetry (which will be evident in the proof of the chosen case), it suffices to show that $I_1$ and $I_2\cap I_3\cap\cdots\cap I_n$ are coprime. Since $I_i$ is an ideal for each $i$, we have
\[
I_1I_2\cdots I_n=(I_1I_2\cdots I_{i-1})I_i(I_{i+1}\cdots I_n)\subset RI_iR\subset I_i,~~\Ra~~I_1I_2\cdots I_n\subset I_1\cap I_2\cap\cdots\cap I_n.
\]
Since $(I_1+I_i)(I_1+I_j)\subset I_1I_1+I_1I_j+I_iI_1+I_iI_j\subset I_1+I_iI_j$, it follow by induction that
\bea
&&R=R^{n-1}=(I_1+I_2)(I_1+I_3)\cdots(I_1+I_n)~\subset~I_1+I_2I_3\cdots I_n\subset~I_1+I_2\cap I_3\cap\cdots\cap I_n\subset R,\nn
\eea
which implies ~$R=I_1+I_2I_3\cdots I_n=I_1+I_2\cap I_3\cap\cdots\cap I_n$.
\end{proof}

\begin{lmm}[\textcolor{blue}{Unit remainder wrt an ideal}]
Let $R$ be a ring and $I_1,...,I_n\vartriangleleft R$ coprime in pairs (i.e., $I_i+I_j=R$ if $i\neq j$). Then for each $k=1,...,n$ there exists an element $x_k\in R$ such that
\[
x_k+I_k=1+I_k~~~~\txt{and}~~~~x_k+I_j=I_j~~\txt{for all}~~j\neq k.
\]
\end{lmm}
\begin{proof}
By symmetry, it suffices to consider $k=1$ only. By the preceding lemma, ~$R=I_1+I_2\cap I_3\cap\cdots\cap I_n$.~ Therefore, $1=s_1+x_1$ for some $s_1\in I_1$ and some $x_1\in I_2\cap I_3\cap\cdots\cap I_n$. Hence $x_1+I_1=1+I_1$ and $x_1+I_j=I_j$ for all $j\neq 1$.
\end{proof}

\begin{thm}[\textcolor{blue}{\index{Chinese remainder theorem}{Chinese remainder theorem}}]\label{ChinRemThm}
Let $R$ be a ring and $I_1,...,I_n\vartriangleleft R$ coprime in pairs. Then for any $r_1,...,r_n\in R$, there exists $r\in R$ such that $r+I_i=r_i+I_i$ for all $i=1,\cdots,n$.
\end{thm}
\begin{proof}
For each $k=1,\cdots, n$, the preceding lemma implies there is an element $x_k\in R$ such that $x_k+I_k= 1+I_k$ and $x_k+I_j=I_j$ for all $j\neq k$.  Let $r:=r_1x_1+r_2x_2+\cdots+r_nx_n$. Then for each $k=1,\cdots,n$,
\[
r+I_k=r_1x_1+r_2x_2+\cdots+r_nx_n+I_k = r_10+\cdots+r_k1+\cdots+r_n0+I_k=r_k+I_k. \qedhere
\]
\end{proof}

\begin{crl}\label{ChinRemCrl1}
Let $R$ be a ring and $I_1,...,I_n\vartriangleleft R$ coprime in pairs (i.e., $I_i+I_j=R$ for all $i\neq j$). Then {\small ${R\over I_1\cap I_2\cap\cdots\cap I_n}\cong{R\over I_1}\times{R\over I_2}\times\cdots\times{R\over I_n}$} (as $R$-modules, and also as rings).
\end{crl}
\begin{proof}
The map {\small $f:R\ra {R\over I_1}\times{R\over I_2}\times\cdots\times{R\over I_n}$, $r\mapsto(r+I_1,...,r+I_n)$} is surjective by the Chinese remainder theorem, and {\small $\ker f=I_1\cap I_2\cap\cdots\cap I_n$}. Hence the result follows from the 1st isomorphism theorem (for modules, and also for rings).
\end{proof}

\begin{lmm}
Let $R$ be a ring and $I_1,\cdots,I_n\vartriangleleft R$ coprime in pairs (i.e., $I_i+I_j=R$ for all $i\neq j$). Then ~$I_1\cap I_2\cap\cdots\cap I_n=\sum_{i_1\neq i_2\neq\cdots\neq i_n}I_{i_1}I_{i_2}\cdots I_{i_n}$.~ That is,
\bea
I_1\cap I_2=I_1I_2+I_2I_1,~~I_1\cap I_2\cap I_3=I_1I_2I_3+I_3I_1I_2+I_2I_3I_1+I_2I_1I_3+I_3I_2I_1+I_1I_3I_2,~~\cdots\nn
\eea
In particular, if $R$ is commutative, then ~{\small $I_1\cap I_2\cap\cdots\cap I_n=I_1I_2\cdots I_n$}.
\end{lmm}
\begin{proof}
If $I,J\vartriangleleft R$ are such that $I+J=R$, then $I\cap J=IJ+JI$ because
\[
I\cap J=(I\cap J)R=(I\cap J)(I+J)\subset (I\cap J)I+(I\cap J)J\subset JI+IJ\subset I\cap J+I\cap J=I\cap J.
\]
The general case thus follows by induction on $n$, since (by the proof of Lemma \ref{CoprimeIdProd}) for all ~$1\leq s\leq t\leq n$,
{\small\[
I_s+I_1\cap I_2\cap\cdots\cap I_{s-1}\cap I_{s+1}\cap\cdots\cap I_t=I_s+I_1I_2\cdots I_{s-1}I_{s+1}\cdots I_t=R. \qedhere
\]}
\end{proof}

\begin{crl}\label{ChinRemCrl2}
If $R$ is a PID and $Ra_1,...,Ra_n\subset R$ are coprime in pairs, then
\bea
\textstyle Ra_1\cap Ra_2\cap\cdots\cap Ra_n=Ra_1Ra_2\cdots Ra_n=Ra_1a_2\cdots a_n.\nn
\eea
\end{crl}

\begin{crl}\label{ChinRemCrl3}
Let $R$ be a PID (hence a UFD). Let $d\in R\backslash 0$ with factorization $d=p_1^{k_1}p_2^{k_2}\cdots p_r^{k_r}$, where $p_1,...,p_r$ are distinct irreducible (prime) elements and  $k_i\geq 1$. Then ~${R\over Rd}\cong {R\over Rp_1^{k_1}}\times {R\over Rp_2^{k_2}}\times\cdots\times {R\over Rp_r^{k_r}}$ (as $R$-modules, and also as rings).
\end{crl}
\begin{proof}
In Corollaries \ref{ChinRemCrl1} and \ref{ChinRemCrl2}, we can set $I_1:=Rp_1^{k_1}$, $I_2:=Rp_2^{k_2}$, $...$, $I_r:=Rp_r^{k_r}$ to obtain
\bea
\textstyle {R\over Rd}={R\over Rp_1^{k_1}Rp_2^{k_2}\cdots Rp_r^{k_r}}={R\over Rp_1^{k_1}\cap Rp_2^{k_2}\cap\cdots \cap Rp_r^{k_r}}\cong {R\over Rp_1^{k_1}}\times {R\over Rp_2^{k_2}}\times\cdots\times {R\over Rp_r^{k_r}}.\nn
\eea
\end{proof}

\begin{rmk}
If $R$ is a PID, then for any $a,b\in R\backslash 0$, ~$a{R\over Rb}:=\{ar+Rb:~r\in R\}={Ra+Rb\over Rb}={R~\!\txt{gcd}(a,b)\over Rb}$. Thus if $a,b$ are coprime then $a{R\over Rb}={R\over Rb}$. On the other hand, if $a|b$ then $a=\txt{gcd}(a,b)$, and so $a{R\over Rb}={Ra\over Rb}$.
\end{rmk}

\begin{lmm}\label{DivRuleLmm}
Let $R$ be an ID and $c,d\in R\backslash 0$. Let $M,N$ be $R$-modules and $m\in M\backslash 0$ such that $Ann_R(m)=0$. Then (i) ${Rm\over Rdm}\cong {R\over Rd}$, (ii) $Ann_{R\over Rdc}(c)\cong {R\over Rc}$, and (iii) $Ann_{M\oplus N}(c)= Ann_M(c)\oplus Ann_N(c)$.
\end{lmm}
\begin{proof}
{\flushleft (i)} The map $\al:R\ra {Rm\over Rdm},~r\mapsto rm+Rdm$ is a surjective R-homomorphism with kernel
\bea
\ker\al=\{r\in R:rm+Rdm=Rdm\}=\{r\in R:r+Rd=Rd\}=Rd,\nn
\eea
and so the result follows by the 1st isomorphism theorem.
{\flushleft (ii)} $Ann_{R\over Rdc}(c)=\left\{r+Rdc\in {R\over Rdc}:cr+Rdc=Rdc\right\}=\left\{r+Rdc\in {R\over Rdc}:r+Rd=Rd\right\}={Rd\over Rdc}\cong {R\over Rc}$.
{\flushleft (iii)} $(m,n)\in Ann_{M\oplus N}(c)$ $\iff$ $c(m,n)=(cm,cn)=0$, $\iff$ $cm=0$ and $cn=0$, $\iff$ $m\in Ann_M(c)$ and $n\in Ann_N(c)$, $\iff$ $(m,n)\in Ann_M(c)\oplus Ann_N(c)$.
\end{proof}

\begin{lmm}\label{DirSumQtLmm}
If $M_1,...,M_n$ are $R$-modules, $N_i\subset M_i$ submodules, $M:=M_1\oplus\cdots \oplus M_n$, and $N:=N_1\oplus\cdots\oplus N_n\subset M$, then ${M\over N}={M_1\oplus\cdots \oplus M_n\over N_1\oplus\cdots \oplus N_n}\cong {M_1\over N_1}\oplus\cdots\oplus{M_n\over N_n}$.
\end{lmm}
\begin{proof}
The map {\footnotesize $\phi: M=M_1\oplus\cdots \oplus M_n\ra {M_1\over N_1}\oplus\cdots\oplus{M_n\over N_n}$, $(m_1,...,m_n)\mapsto (m_1+N_1,\cdots,m_n+N_n)$} is a surjective R-homomorphism with kernel $N_1\oplus\cdots \oplus N_n$. So the result follows from the 1st isomorphism theorem.
\end{proof}

\section{Characterization of FG Modules over a PID}
\subsection{Submodule structure of a FG free module over a PID}
\begin{lmm}\label{FreeScalHomLmm}
Let $R$ be a PID, $S$ a free $R$-module, and $0\neq N\subset S$ a nonzero submodule. Then there exists an $R$-epimorphism $\ld:S\twoheadrightarrow R$ such that $\ld(N)$ is maximal in the set of (left) ideals ~$S^\ast(N):=\big\{f(N):f\in S^\ast:=Hom_R(S,R)\big\}$.
\end{lmm}
\begin{proof}
By Theorem \ref{NoethChrI}, ${}_RR$ is noetherian (recall that $_RR$ is FG over a PID, which is a noetherian ring). Thus, we can choose $\ld\in S^\ast:=\txt{Hom}_R(S,R)$ such that its image $\ld(N)=Rd\subset R$ (for some $d\in R$) is maximal in the set of $S^\ast|_N$ images
\[
S^\ast(N):=\{f(N):f\in S^\ast\}\subset\txt{LeftSpec}(R):=\{\txt{left ideals}~{}_RL\subset R\}.\nn
\]
Let $Y=\{y_j\}_{j\in J}$ be a basis for $S$, i.e., $S=\{\sum_{j\in J}r_jy_j:r_j\in R,~r_j=0~\txt{a.e.f.}~j\in J\}$.
{\flushleft \ul{(1) $d\neq 0$}}: For any $n=\sum_{j\in J}r_jy_j=\sum_{j\in J}y_j^\ast(n)y_j~\in~N\backslash 0$, there is $j'\in J$ such that $r_{j'}=y^\ast_{j'}(n)\neq 0$. Since $0\neq y^\ast_{j'}\in S^\ast$, the maximality of $\ld(N)$ in $S^\ast(N)$ implies $\ld(N)\neq 0$, otherwise we would have $\ld(N)=0~\subsetneq~y^\ast_{j'}(N)$, which is a contradiction. Hence $d\neq 0$ in $\ld(N)=Rd$.
{\flushleft \ul{(2) For any $\gamma\in S^\ast$, ~$\gamma\left(\ld|_N^{-1}(d)\right)=\gamma\big(\ld^{-1}(d)\cap N\big)\subset \ld(N)=Rd$}}: Suppose on the contrary that for some $\gamma\in S^\ast$, we have $\gamma\left(\ld^{-1}(d)\cap N\right)\not\subset Rd$, i.e., there is $n\in \ld^{-1}(d)\cap N$ such that $w:=\gamma(n)\not\in Rd$ (equivalently, $d\nmid w$). Then for some $r,s\in R$, we have
\begin{align}
\ld(N)&=Rd~\!\subsetneq\!~Rd+Rw=R\!~\txt{gcd}(d,w)=R(rw+sd)=R\big(r\gamma(n)+s\ld(n)\big)\nn\\
      &=R(r\gamma+s\ld)(n)=(r\gamma+s\ld)(Rn)\subset (r\gamma+s\ld)(N),\nn
\end{align}
which is a contradiction since $r\gamma+s\ld\in S^\ast$ and $\ld(N)$ is maximal in $S^\ast(N)$.

{\flushleft \ul{(3) $\ld:S\ra R$ is surjective}}: Let {\small $n=\sum_{j\in J}r_jy_j=\sum_{j\in J}y^\ast_j(n)y_j~\in~\ld^{-1}(d)\cap N$}. Then by (2), {\small $y^\ast_j(n)\in y_j^\ast\big(\ld^{-1}(d)\cap N\big)\subset Rd$} for each $j\in J$, that is, for each $j\in J$, we have $y^\ast_j(n)=r_j'd$ for some $r_j'\in R$, and so

\begin{align}
&\textstyle n=\sum_{j\in J}y^\ast_j(n)y_j=dn',~~n':=\sum_{j\in J}r_j'y_j\in S,~~~~\Ra~~~~d=\ld(n)=\ld(dn')=d\ld(n'),\nn\\
&\textstyle~~\Ra~~\ld(n')=1~\in~\ld(S),~~\txt{since}~~d\neq 0,~~~~\Ra~~~~\ld(S)=R~~\txt{(i.e., $\ld$ is surjective)},\nn
\end{align}
where $n'\in S$ holds because $R$ is an ID, and so ~$0=y^\ast_j(n)=r_j'd$~ $\iff$ ~$r_j'=0$.
\end{proof}

\begin{thm}[\textcolor{blue}{Submodule basis in a FG free module over a PID}]
Let $R$ be a PID, ${}_RM$ a FG \ul{free} module of rank $k$, and $N\subset M$ a submodule. Then there exists a \ul{basis} $x_1,\cdots,x_k$ for $M$ and \ul{nonzero} scalars $d_1,\cdots,d_t\in R\backslash\{0\}$, $1\leq t\leq k$, such that (i) $d_j|d_{j+1}$ for all $1\leq j<t$, and (ii) $d_1x_1,\cdots,d_tx_t$ is a basis for $N$. Furthermore, (iii) $d_1,\cdots,d_t$ are uniquely determined up to associates.
\end{thm}
\begin{proof}
We will proceed by induction on $k$. If $N=0$, the result holds trivially. So assume $N\neq 0$. Observe that given a basis $x_1,...,x_k$ for $M$, the existence of nonzero scalars $d_j$ satisfying $d_jx_j\in N$ is guaranteed if $M/N$ is torsion, i.e., if $T(M/N)=M/N$, because if $d_j\in R\backslash 0$ then
\bea
d_jx_j\in N~~\iff~~d_jx_j+N=N,~~\iff~~d_j(x_j+N)=N,~~\iff~~x_j+N\in T(M/N).\nn
\eea
Therefore we will first extract a submodule $S\subset M$ such that (1) $N\subset S$ and (2) $S/N$ is torsion.

{\flushleft \ul{Step I (Restricting to a submodule $S$, $N\subset S\subset M$, such that $S/N$ is torsion):}} By the correspondence theorem, there is a submodule $S$, $N\subset S\subset M$, such that $T(M/N)=S/N$. By Proposition \ref{TorFreeIso1}, $S/N$ is torsion since
\bea
T(S/N)=T(T(M/N))=T(M/N)=S/N.\nn
\eea
By the 3rd isomorphism theorem, ${M\over S} \cong {M/N\over S/N}={M/N\over T(M/N)}$, and so $M/S$ is torsion-free, hence FG and free. So by Corollary \ref{section-existence-lmm}, the SES $0\ra S\ra M\sr{\pi}{\ral} M/S\ra 0$ is split, with a section $j:M/S\ra M$, such that
\bea
M=\ker\pi\oplus j(M/S)=S\oplus F,~~~~\txt{where $F:=j(M/S)$ is FG and free}.\nn
\eea
We also know $S$ is FG and free as a submodule of $M$ (a FG free module over a PID). If $\{x_1,...,x_{k'}\}$ is a basis for $S$ and $\{z_1,...,z_l\}$ a basis for $F$, then it is clear $\{x_1,\cdots,x_{k'},z_1,\cdots, z_l\}$ is a basis for $M=S\oplus F$.

{\flushleft \ul{Step II (Initial induction step):}} By Lemma \ref{FreeScalHomLmm} we have an epimorphism $\ld:S\ra R$ such that $\ld(N)=Rd$ (for some $d\in R$) is maximal in $S^\ast(N):=\big\{f(N):f\in S^\ast:=Hom_R(S,R)\big\}$.

Since $\ld:S\ra R$ is surjective and $R$ is free of rank 1, a section $j:R\ra S$ of $\ld$ exists. In particular, for any $m\in \ld^{-1}(1)$, the map $j_m:R\ra S,~r\mapsto rm$ is a section of $\ld$ since $\ld j_m(r)=\ld(rm)=r\ld(m)=r1=r$ for all $r\in R$. It follows therefore that, with $m\in\ld^{-1}(1)$ and $S_{k'-1}:=\ker\ld$, we have
\bea
S=\ker\ld\oplus j_m(R)=\ker\ld\oplus Rm=S_{k'-1}\oplus Rm.\nn
\eea
For the restriction {\small $\ld|_N:N\ra \ld(N)=Rd$}, since $Rd$ is free on one generator, a section of $P$ also exists. In particular, a section of {\small $\ld|_N$ is $j_m|_{Rd}:Rd\ra N,~rd\mapsto rdm$}. So, with {\small $N_{{k'}-1}:=N\cap\ker\ld=N\cap S_{{k'}-1}$}, we get
\bea
N=\ker\ld|_N\oplus j_m|_{Rd}(Rd)=(N\cap\ker\ld)\oplus Rdm=(N\cap S_{{k'}-1})\oplus Rdm=N_{{k'}-1}\oplus Rdm.\nn
\eea

{\flushleft \ul{Step III (Final induction step):}} Observe that $S_{{k'}-1}:=\ker\ld$ is free on ${k'}-1$ generators. So by the induction hypotheses, we can find a basis $x_2,...,x_{k'}$ for $S_{{k'}-1}:=\ker\ld$ and $d_2,...,d_t\in R\backslash 0$, $1\leq t\leq k'$, such that (i) $d_j|d_{j+1}$ for all $2\leq j<t$ and (ii) $d_2x_1,\cdots, d_tx_t$ is a basis for $N_{{k'}-1}:=N\cap S_{{k'}-1}=N\cap\ker\ld$.

Let $x_1:=m\in\ld^{-1}(1)$ and $d_1:=d$. Then it is clear (from the construction in Step II above) that $d_1x_1,d_2x_2,\cdots,d_tx_t$ is a basis for $N$ and $x_1,x_2,\cdots, x_{k'}$ is a basis for $S$. It remains to show that $d_1|d_2$.

{\flushleft \ul{Step IV (Checking divisibility $d_1|d_2$):}} The map $\gamma:S\ra R,~\sum_{i=1}^{k'}r_ix_i\mapsto r_1+r_2$ lies in the dual module of $S$, i.e., $\gamma\in S^\ast=\txt{Hom}_R(S,R)$, and so by the maximality of $\ld(N)=Rd$ in $S^\ast(N)$,
\bea
&&\ld(N)=Rd=Rd_1\subset Rd_1+Rd_2=\gamma(Rd_1x_1+Rd_2x_2)\subset\gamma(N)~~\Ra~~\gamma(N)=\ld(N),\nn\\
&&~~\Ra~~Rd_1+Rd_2=Rd_1,~~\Ra~~d_1|d_2.\nn
\eea

{\flushleft \ul{Step V (Checking uniqueness):}} This will be done more automatically in Theorem \ref{ST-FG-PID}.
\end{proof}

\subsection{Structure theorem of FG modules over a PID}
\begin{rmk}[\blue{${R\over I}$-module correspondence}]\label{RoverImod}
Let $R$ be a ring, $I\lhd R$ an ideal, and $M$ an $R$-module. If $IM=0$, then the map
${R\over I}\times M\ra M,~(r+I,m)\mapsto rm$ is well-defined (\blue{footnote}\footnote{If $(r+I,m)=(r'+I,m')$, then $r+I=r'+I$ (i.e., $r-r'\in I$) and $m=m'$, and so $rm=rm'=(r-r')m'+r'm'=r'm'$.}), thereby making $M$ an ${R\over I}$-module. Moreover, because the original $R$-scalar multiplication $R\times M\ra M,~(r,m)\mapsto rm$ has exactly the same effect on $M$ as the ${R\over I}$-scalar multiplication, we have an automatic bijective correspondence
\bea
\textstyle \left\{\txt{$R$-submodules ~${}_RN\subset M$}\right\}~~\sr{id}{\longleftrightarrow}~~\left\{\txt{${R\over I}$-submodules ${}_{R\over I}N\subset M$}\right\}.\nn
\eea
\end{rmk}

\begin{crl}
Let $M$ be an $R$-module and $I\lhd R$ an ideal. Then we have a bijective correspondence
\bea
\textstyle \left\{\txt{$R$-submodules ${}_RN\subset{M\over IM}$}\right\}~~\sr{id}{\longleftrightarrow}~~\left\{\txt{${R\over I}$-submodules ${}_{R\over I}N\subset {M\over IM}$}\right\}.\nn
\eea
\end{crl}

\begin{thm}[\textcolor{blue}{\index{Structure theorem of! FG modules over a PID}{Structure theorem of FG modules over a PID}}]\label{ST-FG-PID}
Let $R$ be a PID and ${}_RM$ a FG module. Then there exist (1) unique scalars $d_1,...,d_t\in R\backslash\{0\}$ satisfying $d_i|d_{i+1}$ for all $1\leq i<t$ and (2) a unique free R-module of finite rank $F_0\cong R^k$ such that (\blue{footnote}\footnote{The uniqueness of scalars is only up to associates. Moreover, in this theorem, ``$d_1,...,d_t\in R\backslash\{0\}$'' is equivalent to ``$d_1,...,d_t\in R\backslash(\{0\}\cup U(R))$'' because any $d_i\in U(R)$ contributes 0 to the direct sum decomposition of $M$.})
\bea
\textstyle M\cong {R\over Rd_1}\oplus {R\over Rd_2}\oplus\cdots\oplus {R\over Rd_t}\oplus F_0.\nn
\eea
\end{thm}
\begin{proof}
{\flushleft\ul{Existence}}: Since $M$ is FG, let $M=Rx_1+\cdots+Rx_n$. Then with the free R-module $F=Rx_1\oplus\cdots\oplus Rx_n$ over $\{x_1,...,x_n\}$, we have the surjective R-homomorphism $\pi:F\ra M,~(r_1x_1,...,r_nx_n)\ra\sum r_ix_i$. Let $N:=\ker\pi$. Then by Theorem \ref{FrModPidHered}, we know $N$ is a free R-module of finite rank $\leq n$. By the preceding theorem, assume wlog that $\{x_1,...,x_n\}$ was chosen (\blue{footnote}\footnote{The map $\pi$ is now only unique up to composition with a change of basis isomorphism.}) such that we have a basis $\{d_1x_1,...,d_tx_t\}$, $1\leq t\leq n$, for $N$ with the scalars $d_i\in R\backslash\{0\}$ satisfying $d_i|d_{i+1}$ for all $1\leq i<t$. By the 1st isomorphism theorem, Lemma \ref{DirSumQtLmm}, and Lemma \ref{DivRuleLmm} (plus Remark \ref{TorFreeRmk1}(4), i.e., a free ID-module has no torsion),
\begin{align}
&\textstyle M\cong{F\over N}={Rx_1\oplus Rx_2\oplus\cdots\oplus Rx_t\oplus Rx_{t+1}\oplus\cdots\oplus Rx_n\over Rd_1x_1\oplus Rd_2x_2\oplus\cdots\oplus Rd_tx_t\oplus 0\oplus\cdots\oplus 0}\cong
{Rx_1\over Rd_1x_1}\oplus{Rx_2\over Rd_2x_2}\oplus\cdots\oplus{Rx_t\over Rd_tx_t}\oplus {Rx_{t+1}\over 0}\oplus\cdots\oplus {Rx_n\over 0}\nn\\
&\textstyle~~~~\cong {R\over Rd_1}\oplus{R\over Rd_2}\oplus\cdots\oplus{R\over Rd_t}\oplus F_0,~~~~F_0:=Rx_{t+1}\oplus\cdots\oplus Rx_n\cong R^{\langle n-t\rangle}.\nn
\end{align}

{\flushleft\ul{Uniqueness}}: Assume ~$D:={R\over Rd_1}\oplus{R\over Rd_2}\oplus\cdots\oplus{R\over Rd_t}\oplus F_0~\cong~{R\over Rw_1}\oplus{R\over Rw_2}\oplus\cdots\oplus{R\over Rw_s}\oplus G_0=:W$,~ where $d_i,w_j\in R\backslash\{0\}$, $d_i|d_{i+1}$, $w_j|w_{j+1}$, and $F_0,G_0$ are free of finite rank. We need to show (a) $t=s$, (b) $Rd_i=Rw_i$ for all $i=1,...,t$, and (c) $\txt{rank}(F_0)=\txt{rank}(G_0)$.

Recall that if $R$ is a PID and ${}_RD\cong{}_RW$, then as shown before in Theorem \ref{TorFreeDec}, (1) {\small $T(D)\cong T(W)$}, (2) {\small ${D\over T(D)}\cong {W\over T(W)}$}, and clearly, (3) if $r\in R$, then {\small $Ann_D(r)\cong Ann_W(r)$} as R-modules, where by Lemma \ref{DivRuleLmm},
\bea
&&\textstyle Ann_D(r):=\{x\in D:rx=0\}=\bigoplus_{i=1}^t Ann_{R\over Rd_i}(r),\nn\\
&&\textstyle Ann_W(r):=\{x\in W:rx=0\}=\bigoplus_{i=1}^s Ann_{R\over Rw_i}(r).\nn
\eea

Since $T(D)={R\over Rd_1}\oplus{R\over Rd_2}\oplus\cdots\oplus{R\over Rd_t}$ and $T(W)={R\over Rw_1}\oplus{R\over Rw_2}\oplus\cdots\oplus{R\over Rw_s}$, we get $F_0\cong {D\over T(D)}\cong{W\over T(W)}\cong G_0$, and so $\txt{rank}(F_0)=\txt{rank}(G_0)$. Because $T(D)\cong T(W)$, we can now assume $F_0=G_0=0$ in what follows.

We now use induction on the number of prime factors (counting multiplicity) of the product
$$d_1d_2\cdots d_tw_1w_2\cdots w_s.$$ Pick a prime $p\in R$ such that $p|d_1$. This implies $p|d_i$ for all $i$. Write $d_i=pd_i'$. Then by Lemma \ref{DivRuleLmm},
\bea
\textstyle\left({R\over Rp}\right)^t\cong\bigoplus_{i=1}^tAnn_{R\over Rd_i}(p)= Ann_D(p)\cong Ann_W(p)=\bigoplus_{i=1}^sAnn_{R\over Rw_i}(p)\cong\left({R\over Rp}\right)^u,\nn
\eea
where $u$ is the number of $w_i$'s such that $p|w_i$. By Remark \ref{RoverImod}, we can view $\left({R\over Rp}\right)^t$ and $\left({R\over Rp}\right)^u$ as ${R\over Rp}$-modules (i.e., vector spaces), where we note ${R\over Rp}$ is a field (as a simple commutative ring). So, by vector space dimension, $t=u\leq s$. By symmetry, $s\leq t$, and hence $t=u=s$. Therefore,
\bit
\item[] ${R\over Rd_1}\oplus{R\over Rd_2}\oplus\cdots\oplus{R\over Rd_t}=T(D)\cong T(W)={R\over Rw_1}\oplus{R\over Rw_2}\oplus\cdots\oplus{R\over Rw_t}$, ~which implies
\item[] ${R\over Rd_1'}\oplus\cdots\oplus{R\over Rd_t'}\cong pT(D)\cong pT(W)\cong {R\over Rw_1'}\oplus\cdots\oplus{R\over Rw_t'}$,
\eit
where $d_i=pd_i'$ and $w_j=pw_j'$. Hence, by induction on the number of prime factors, we can conclude that $Rd_i=Rw_i$ for all $i$.
\end{proof}

\begin{dfn}[\textcolor{blue}{\index{Invariant! divisors}{Invariant divisors}, \index{Invariant! factors}{Invariant factors}}]
In the above theorem, the we will call the ideals $Rd_1,\cdots, Rd_t$ (resp. the modules ${R\over Rd_1},\cdots,{R\over Rd_t}$) the \ul{invariant divisors} (resp. \ul{invariant factors}) of $M$.
\end{dfn}

\subsection{Fundamental theorem of FG modules over a PID}
\begin{thm}[\textcolor{blue}{\index{Fundamental theorem of! FG modules over a PID}{Fundamental theorem of FG modules over a PID}}]\label{FT-FG-PID}
Let $R$ be a PID and ${}_RM$ a FG module. Then there exist (1) primes $p_1,...,p_s\in R$ (not necessarily distinct), and (2) a free R-module $F_0$ of unique finite rank such that
\bea
\textstyle M\cong {R\over Rp_1^{k_1}}\oplus{R\over Rp_2^{k_2}}\oplus\cdots\oplus{R\over Rp_s^{k_s}}\oplus F_0,~~~~\txt{for some}~~k_i\geq1.\nn
\eea
Furthermore, the pairing $\{(p_i,k_i):1\leq i\leq s\}$ is unique up to the ordering of the decomposition factors.
\end{thm}
\begin{proof}
In Theorem \ref{ST-FG-PID}, consider the prime factorizations $d_i=p_{1i}^{k_{1i}}p_{2i}^{k_{2i}}\cdots p_{r_ii}^{k_{r_ii}}$, $i=1,...,t$. Then, in the decomposition of $M$, use Corollary \ref{ChinRemCrl3} to expand each summand ${R\over Rd_i}$ so as to obtain the decomposition in the desired form.
\end{proof}

\begin{dfn}[\textcolor{blue}{\index{Elementary divisors}{Elementary divisors}, \index{Elementary factors}{Elementary factors}}]
In the above theorem, we will call the ideals $Rp_1^{k_1},\cdots, Rp_s^{k_s}$ (resp. the modules {\small ${R\over Rp_1^{k_1}},\cdots,{R\over Rp_s^{k_s}}$}) the \ul{elementary divisors} (resp. \ul{elementary factors}) of $M$.
\end{dfn}

\begin{rmk}
It is clear by construction that Theorems \ref{ST-FG-PID} and \ref{FT-FG-PID} are equivalent since they follow from each other. In particular, in the process of obtaining the invariant factors ${R\over Rd_j}$ (satisfying $d_j|d_{j+1}$) by rewriting the (already given) elementary factors ${R\over Rp_i^{k_i}}$, one should find that the resulting number of invariant factors $t=|\{d_1,...,d_t\}|$ is precisely the highest prime power, i.e., $t=|\{d_1,...,d_t\}|=\max\{k_1,...,k_s\}$.
\end{rmk}

\begin{crl}[\textcolor{blue}{\index{Fundamental theorem of! FG abelian groups}{Fundamental theorem of FG abelian groups}}]
Let $A$ be a FG abelian group. Then
\bea
\textstyle A\cong \Integer_{d_1}\oplus\Integer_{d_2}\oplus\cdots\oplus\Integer_{d_t}\oplus\Integer^n,~~~~\Integer_d:={\Integer\over\Integer d},\nn
\eea
for unique integers $d_i\geq 1$ satisfying $d_i|d_{i+1}$, and a unique integer $n\geq 0$. Equivalently, we have a unique expression (up to reordering)
\bea
A\cong\Integer_{p_1^{k_1}}\oplus\Integer_{p_2^{k_2}}\oplus\cdots\oplus\Integer_{p_s^{k_s}}\oplus\Integer^n,\nn
\eea
for primes $p_1,...,p_s$ (not necessarily distinct), integers $k_i\geq 1$, and a unique integer $n\geq 0$.
\end{crl}

\begin{crl}[\textcolor{blue}{\index{Fundamental theorem of! finite abelian groups}{Fundamental theorem of finite abelian groups}}]
Let $A$ be a finite abelian group. Then
\bea
\textstyle A\cong \Integer_{d_1}\oplus\Integer_{d_2}\oplus\cdots\oplus\Integer_{d_t},~~~~\Integer_d:={\Integer\over\Integer d},\nn
\eea
for unique integers $d_i\geq 1$ satisfying $d_i|d_{i+1}$. Equivalently, we have a unique expression (up to reordering)
\bea
A\cong \Integer_{p_1^{k_1}}\oplus\Integer_{p_2^{k_2}}\oplus\cdots\oplus\Integer_{p_s^{k_s}},\nn
\eea
for primes $p_1,...,p_s$ (not necessarily distinct) and integers $k_i\geq 1$.
\end{crl}

\section{Canonical Forms of Linear Transformations}
In this section $k$ will denote a field unless specified otherwise or its alternative role is clear from the context. Some of the concepts and results might still hold even when $k$ is replaced with any commutative ring.
\begin{dfn}[\blue{
\index{Matrix! of a linear transformation}{Matrix of a linear transformation},
\index{Coordinate expansion}{Coordinate expansion},
\index{Coordinates}{Coordinates},
\index{Coordinate representation}{Coordinate representation},
\index{Matrix! transformation}{Matrix transformation},
\index{Change of basis}{Change of basis}}]
Let $V,W$ be (finite-dimensional) $k$-spaces with respective bases $\beta=\{\beta_1,...,\beta_n\}$, $\gamma=\{\gamma_1,...,\gamma_m\}$, and respective dual bases $\beta^\ast=\{\beta_1^\ast,...,\beta_n^\ast\}$, $\gamma^\ast=\{\gamma_1^\ast,...,\gamma_m^\ast\}$). Then we know that any $v\in V$ can be written as follows:
{\footnotesize\[
\textstyle v=\sum_i \beta_i^\ast(v)\beta_i=\beta^\ast(v)^T\beta,~~\txt{for a unique}~~\beta^\ast(v):=\big(\beta_1^\ast(v),\cdots,\beta_n^\ast(v)\big)=
\left[
     \begin{array}{c}
       \beta_1^\ast(v) \\
       \beta_2^\ast(v) \\
       \vdots \\
       \beta_n^\ast(v) \\
     \end{array}
   \right]\in k^n:=k^{n\times 1}:=Mat_{n\times 1}(k),
\]} where $\beta^\ast(v)^T:=
\left[
  \begin{array}{cccc}
    \beta_1^\ast(v) & \beta_2^\ast(v) & \cdots & \beta_n^\ast(v)\\
  \end{array}
\right]\in k^{1\times n}:=Mat_{1\times n}(k)
$ and, for convenience, we similarly write
{\footnotesize $\beta:=(\beta_1,...,\beta_n):
=\left[
     \begin{array}{c}
       \beta_1 \\
       \beta_2 \\
       \vdots \\
       \beta_n \\
     \end{array}
   \right]\in V^n:=V^{n\times 1}:=Mat_{n\times 1}(V)
$} and
{\footnotesize $\beta^T:=
\left[
  \begin{array}{cccc}
    \beta_1 & \beta_2 & \cdots & \beta_n\\
  \end{array}
\right]\in V^{1\times n}:=Mat_{1\times n}(V)
$}. In the expansion $v=\sum_i \beta_i^\ast(v)\beta_i$ (called \ul{$\beta$-coordinate expansion} of $v$), the coefficients $\beta_i^\ast(v)$ are called \ul{$\beta$-coordinates} of $v$. Hence, we have $k$-linear isomorphisms (called \ul{coordinate representations} wrt the bases)
\bea
&&\textstyle\tau_\beta=\tau_\beta^V:V\ra k^n,~v\mapsto\beta^\ast(v),~~~~\tau_\beta^{-1}:k^n\ra v,~x\mapsto x^T\beta=\sum_i x_i\beta_i,\nn\\
&&\textstyle \tau_\gamma=\tau_\gamma^W:W\ra k^m,~w\mapsto\gamma^\ast(w),~~~~\tau_\gamma^{-1}:k^m\ra W,~y\mapsto y^T\gamma=\sum_jy_j\gamma_j.\nn
\eea
It is clear that a $k$-linear transformation ${\tau}:V\ra W$ is uniquely given by linear relations of the form
{\footnotesize\[
\textstyle {\tau}(\beta_j)=\sum_{i=1}^m\tau_{ij}\gamma_i=
\tau_{1j}\gamma_1+...+\tau_{mj}\gamma_m=
\left[
  \begin{array}{ccc}
    \gamma_1 & \cdots & \gamma_m \\
  \end{array}
\right]\cdot_{op}
\left[
  \begin{array}{c}
    \tau_{1j} \\
    \tau_{2j} \\
    \vdots \\
    \tau_{mj} \\
  \end{array}
\right]=
\gamma^T\cdot_{op}
\left[
  \begin{array}{c}
    \tau_{1j} \\
    \tau_{2j} \\
    \vdots \\
    \tau_{mj} \\
  \end{array}
\right],~~~~j=1,...,n.\nn
\]}
These can be written as a matrix relation
{\scriptsize\bea
&&\label{MatRepEq}{\tau}(\beta)^T=
\left[
  \begin{array}{cccc}
    {\tau}(\beta_1) & {\tau}(\beta_2) & \cdots & {\tau}(\beta_n) \\
  \end{array}
\right] =
\left[
  \begin{array}{cccc}
    \gamma_1 & \gamma_2 & \cdots & \gamma_m \\
  \end{array}
\right]\cdot_{op} \left[
          \begin{array}{cccc}
            \tau_{11} & \tau_{12} & \cdots & \tau_{1n} \\
            \tau_{21} & \tau_{22} & \cdots & \tau_{2n}\\
            \vdots & \vdots & \ddots & \vdots \\
            \tau_{m1} & \tau_{m2} & \cdots & \tau_{mn} \\
          \end{array}
        \right]=\gamma^T\cdot_{op}[\tau]_{\beta\gamma},
\eea}
or equivalently (after taking the transpose and using $(x^T)^T=x$, $(xy)^T=y^Tx^T$) we have
\bea
\label{MatRepEq2}\tau(\beta)=[\tau]_{\beta\gamma}^T\gamma,
\eea
where $[{\tau}]_{\beta\gamma}:=[\tau_{ij}]\in Mat_{m\times n}(k)$ is the \ul{matrix} of ${\tau}$ with respect to $\beta$ and $\gamma$. When $\gamma=\beta$, we will write $[{\tau}]_{\beta}:=[{\tau}]_{\beta\beta}$. Given $v\in V$, we have
\bea
\tau(v)=\tau\big(\beta^\ast(v)^T\beta\big)=\beta^\ast(v)^T\tau(\beta)=\beta^\ast(v)^T[\tau]_{\beta\gamma}^T\gamma=\big([\tau]_{\beta\gamma}\beta^\ast(v)\big)^T\gamma=:\big(\tau_{\beta\gamma}(\beta^\ast(v))\big)^T\gamma,\nn
\eea
where $\tau_{\beta\gamma}(\beta^\ast(v)):=[\tau]_{\beta\gamma}\beta^\ast(v)$. Thus, we can equivalently view $\tau:V\ra W$ as the \ul{matrix transformation}
\bea
\tau_{\beta\gamma}:k^n\ra k^m,~~x\mapsto [\tau]_{\beta\gamma}x.\nn
\eea
In terms of the isomorphisms $\tau_\beta:V\ra k^n$, $\tau_\gamma:W\ra k^m$, and the matrix transformation $\tau_{\beta\gamma}:k^n\ra k^m$,
\begin{align}
~~~\tau&=\tau_\gamma^{-1}\circ\tau_{\beta\gamma}\circ\tau_\beta:(V,\beta)\sr{\tau_\beta}{\ral}k^n\sr{\tau_{\beta\gamma}}{\ral}k^m\sr{\tau_\gamma^{-1}}{\ral}(W,\gamma),\nn\\
~~~\Ra~~~\tau_{\beta\gamma}&=\tau_\gamma\circ\tau\circ\tau_\beta^{-1}:k^n\sr{\tau_\beta^{-1}}{\ral}(V,\beta)\sr{\tau}{\ral}(W,\gamma)\sr{\tau_\gamma}{\ral} k^m,~~x\mapsto\gamma^\ast\big(\tau(x^T\beta)\big)=:[\tau]_{\beta\gamma}x.\nn
\end{align}
That is, the \ul{matrix} of $\tau$ is given by $[\tau]_{\beta\gamma}:=[\tau_{ij}]:=\big[\gamma^\ast_i\big(\tau(\beta_j)\big)\big]\in Mat_{m\times n}(k)$, where the entries are
\bea
\label{MatRepEq3}\tau_{ij}:=\gamma^\ast_i\big(\tau(\beta_j)\big),~~~~i=1,...,m,~~j=1,...,n.
\eea
When $W=V$ and $\tau=\mathds{1}:=id_V$, then $\mathds{1}_{\beta\gamma}$ (hence $[\mathds{1}]_{\beta\gamma}$) is called a \ul{change of basis transformation}:
\[
\mathds{1}_{\beta\gamma}=\tau_\gamma\circ\tau_\beta^{-1}:k^n\ra (V,\beta)\sr{\mathds{1}}{\ral}(V,\gamma)\ra k^n,~~x\mapsto\gamma^\ast(x^T\beta)=:[\mathds{1}]_{\beta\gamma}x,
\]
where the \ul{change of basis matrix} ~$[\mathds{1}]_{\beta\gamma}=[\mathds{1}_{ij}]$~ has entries ~$\mathds{1}_{ij}:=\gamma^\ast_i(\beta_j)$, ~$i,j=1,...,n$. Observe that in any basis for $\beta$ for $V$, we have $\mathds{1}_{\beta\beta}=id_{k^n}$ (which denote by $[\mathds{1}]$ or just by $\mathds{1}$), since $[\mathds{1}]_{\beta\beta}=[\delta_{ij}]$.
\end{dfn}

\begin{rmk}
Let $V,V',V''$ be $k$-spaces with respective bases $\beta,\beta',\beta''$. (i) We can directly see that
\begin{align}
&[\tau_1+\tau_2]_{\beta\beta'}=[\tau_1]_{\beta\beta'}+[\tau_2]_{\beta\beta'},~~~~\txt{for all}~~~~\tau_1,\tau_2\in Hom_k(V,V'),\nn\\
&[\tau''\circ\tau']_{\beta\beta''}=[\tau'']_{\beta\beta'}[\tau']_{\beta'\beta''},~~~~\txt{for all}~~~~\tau\in Hom_k(V,V'),~\tau'\in Hom_k(V',V'').\nn
\end{align}
(ii) For convenience, we will define ~$[\tau]_{\beta}:=[\tau]_{\beta\beta}$.
\end{rmk}

\begin{dfn}[\blue{
\index{Monic! polynomial}{Monic polynomial},
\index{Conjugate matrices}{Conjugate matrices}}]
Let $R$ be a commutative ring. A polynomial $f=f^0+f^1x+\cdots+f^{d_f}x^{d_f}\in R[x]$ is \ul{monic} if the leading coefficient $f^{d_f}=1$. Matrices $A,B\in Mat_n(R):=Mat_{n\times n}(R)$ are \ul{conjugate}, written $A\sim B$, if there exists an invertible matrix $P\in GL_n(R):=U(Mat_{n\times n}(R))\cong Aut_R(R^n)$ such that $A=PBP^{-1}$.
\end{dfn}

\subsection{Polynomial decomposition of a finite-dimensional vector space}
\begin{dfn}[\blue{\index{Vector space as a polynomial-module}{Vector space as a polynomial-module}}]
Let $k$ be a field, $V$ a $k$-space, and $\tau\in End_k(V):=Hom_k(V,V)$ as a ring. Then $V$ is a $k[x]$-module via the map
\bea
m_\tau:k[x]\times V\ra V,~\big(f,v\big)\mapsto fv:=f(\tau)(v),~~~~f(\tau):=f^0+f^1\tau+f^2\tau^2+\cdots+f^{\deg(f)}\tau^{\deg(f)}.\nn
\eea
We will denote $V$ with this $k[x]$-module structure by $_{k[x]}V(\tau)$, or simply by $_{k[x]}V$ when the choice/identity of the endomorphism $\tau$ is unimportant (where we note that in general ${}_{k[x]}V(\tau)\not\cong{}_{k[x]}V(\tau')$ if $\tau\not\cong\tau'$ in the category of linear transformations). Henceforth, we will assume that every $k$-space $V$ has a $k[x]$-module structure as defined above, i.e., $V={}_kV={}_{k[x]}V(\tau)$.
\end{dfn}

\begin{rmk}[\blue{Polynomial-module as a vector space}]
Let $M$ be a $k[x]$-module with scalar multiplication $\al:k[x]\times M\ra M$. Then because $k\subset k[x]$ is a subring, $M$ is a $k$-space with scalar multiplication $\al|_{k\times M}:k\times M\subset k[x]\times M\ra M$. Also, it is clear that any $k[x]$-homomorphism $h:{}_{k[x]}M\ra{}_{k[x]}M'$ gives a $k$-linear transformation $h:{}_kM\ra{}_kM'$.
\end{rmk}

\begin{lmm}\label{VecPolTor}
Let $V$ be a $k$-space and $\tau\in End_k(V)$. If $V$ is finite-dimensional, then the module $_{k[x]}V(\tau)$ is torsion (i.e., for each $v\in V$, there exists $f\in k[x]\backslash 0$ such that $f(\tau)(v)=0$).
\end{lmm}
\begin{proof}
Let $M:={}_{k[x]}V(\tau)$. Suppose $T_{k[x]}(M)\neq M$. Let $v\in M\backslash T_{k[x]}(M)$, and consider the map
\bea
\phi:k[x]\ra M,~f\mapsto fv:=f(\tau)(v).\nn
\eea
Then $\phi$ is injective. Since $k[x]$ is infinite-dimensional as a k-space, i.e., $\dim_k(k[x])=\infty$, $M$ contains an infinite-dimensional k-subspace $\phi(k[x])=k[x]v$. This implies $\dim_kV=\dim_kM\geq\dim_k\phi\big(k[x]v\big)=\infty$ (a contradiction).
\end{proof}

\begin{lmm}[\blue{Polynomial decomposition of a vector space}]\label{RatDecLmm}
Let $V$ be a finite-dimensional $k$-space and $\tau\in End_k(V)$. Then there exist nonzero polynomials $f_1,...,f_t\in k[x]\backslash 0$, satisfying $f_i|f_{i+1}$, such that
\bea
\textstyle V={}_kV={}_{k[x]}V(\tau)\cong{k[x]\over k[x]f_1}\oplus\cdots\oplus{k[x]\over k[x]f_t}~~~~\txt{(as $k$-spaces and as $k[x]$-modules)}.\nn
\eea
\end{lmm}
\begin{proof}
Let $B=\{v_1,...,v_n\}$ be a basis for $V$. Then $V=kv_1+\cdots+kv_n$, which implies $V=k[x]V=k[x]v_1+\cdots+k[x]v_n$, i.e., $_{k[x]}V(\tau)$ is finitely generated (\blue{footnote}\footnote{Here, the choice of $\tau\in End_k(V)$ for the $k[x]$-module $_{k[x]}V={}_{k[x]}V(\tau)$ is unimportant.}). Since $k[x]$ is a PID and $_{k[x]}V(\tau)$ is both FG and torsion (by Lemma \ref{VecPolTor}), we get (by Theorem \ref{ST-FG-PID}) an isomorphism of $k[x]$-modules
\bea
\textstyle{}_{k[x]}V(\tau)\cong{k[x]\over k[x]f_1}\oplus\cdots\oplus{k[x]\over k[x]f_t},~~~~\txt{for some}~~~~f_i\in k[x]\backslash\{0\}~~\txt{satisfying}~~f_i|f_{i+1},\nn
\eea
which of course implies we also have an isomorphism of $k$-spaces
\[
\textstyle{}_kV\cong{k[x]\over k[x]f_1}\oplus\cdots\oplus{k[x]\over k[x]f_t}. \qedhere
\]
\end{proof}

\subsection{Invariant subspace decomposition and polynomial-module equivalences}
\begin{dfn}[\textcolor{blue}{\index{Invariant! subspace}{Invariant subspace}}]
Let $V$ be a $k$-space, $\tau\in\txt{End}_k(V):=\txt{Hom}_k(V,V)$, and ${}_kW\subset V$ a $k$-subspace. Then $W$ is \ul{$\tau$-invariant} if $\tau(W)\subset W$.
\begin{note*}
The $k$-subspace ${}_kW\subset V$ is a submodule of $_{k[x]}V(\tau)$ $\iff$ $k[x]W\subset W$, $\iff$ $xW=\tau(W)\subset W$, $\iff$ $W$ is $\tau$-invariant.
\end{note*}
\end{dfn}

\begin{dfn}[\textcolor{blue}{\index{Direct sum of matrices (or linear transformations)}{Direct sum of matrices (or linear transformations)}}]
Given matrices $A\in Mat_{n}(k)$ and $A'\in Mat_{n'}(k)$, their direct sum is the matrix
{\footnotesize $A\oplus A':=\left[
           \begin{array}{c|c}
             A & 0 \\
             \hline
             0 & A' \\
           \end{array}
         \right]\in Mat_{n+n'}(k)$}. We will denote the collection of such matrices by ~$Mat_{n}(k)\oplus Mat_{n'}(k)\subset Mat_{n+n'}(k)$.
\end{dfn}

\begin{prp}
Let $V$ be a finite-dimensional $k$-space, ${\tau}\in\txt{End}_k(V)$, and $W,W'\subset V$ ${\tau}$-invariant subspaces. If (i) $V=W\oplus W'$ and (ii) $\gamma=\{\gamma_1,...,\gamma_l\}$ and $\gamma'=\{\gamma'_1,...,\gamma'_{l'}\}$ are bases for $W$ and $W'$ respectively, then $\beta:=\gamma\cup \gamma'$ is a basis for $V$, and the matrix of $\tau$ wrt $\beta$ is given by
\bea
[{\tau}]_\beta=\big[{\tau}|_{W}\big]_{\gamma}\oplus\big[{\tau}|_{W'}\big]_{\gamma'}=A\oplus A',
      ~~~~\txt{where}~~~~A=\big[{\tau}|_{W}\big]_{\gamma},~~A'=\big[{\tau}|_{W'}\big]_{\gamma'}.\nn
\eea

\end{prp}
\begin{proof}
Indeed $\beta:=\gamma\cup \gamma'$ is a basis for $V$ since by the associativity of the direct sum $\oplus$,
\bea
&&V=W\oplus W'=(k\gamma_1\oplus k\gamma_2\oplus\cdots\oplus k\gamma_l)\oplus(k\gamma'_1\oplus k\gamma'_2\oplus\cdots\oplus k\gamma'_{l'})\nn\\
&&~~~~=k\gamma_1\oplus k\gamma_2\oplus\cdots\oplus k\gamma_l\oplus k\gamma'_1\oplus k\gamma'_2\oplus\cdots\oplus k\gamma'_{l'}.\nn
\eea
Also, we know that
{\footnotesize\bea
&&\left[
    \begin{array}{cccc}
      {\tau}(\gamma_1) & {\tau}(\gamma_2) & \cdots & {\tau}(\gamma_l) \\
    \end{array}
  \right]=
  \left[
    \begin{array}{cccc}
     \gamma_1 & \gamma_2 & \cdots & \gamma_l \\
    \end{array}
  \right]\cdot_{op} \big[{\tau}|_{W}\big]_{\gamma}
=\left[
    \begin{array}{cccc}
     \gamma_1 & \gamma_2 & \cdots & \gamma_l \\
    \end{array}
  \right]\cdot_{op} A,\nn\\
&&\left[
    \begin{array}{cccc}
      {\tau}(\gamma'_1) & {\tau}(\gamma'_2) & \cdots & {\tau}(\gamma'_{l'}) \\
    \end{array}
  \right]=
  \left[
    \begin{array}{cccc}
     \gamma'_1 & \gamma'_2 & \cdots & \gamma'_{l'} \\
    \end{array}
  \right]\cdot_{op} \big[{\tau}|_{W'}\big]_{\gamma'}
=\left[
    \begin{array}{cccc}
     \gamma'_1 & \gamma'_2 & \cdots & \gamma'_{l'} \\
    \end{array}
  \right]\cdot_{op} A'.\nn
\eea}Therefore,
~{\footnotesize $\left[
           \begin{array}{c|c}
             A & 0 \\
             \hline
             0 & A' \\
           \end{array}
         \right] = [{\tau}]_\beta$},~ since
{\footnotesize\[
\left[
    \begin{array}{cccccc}
      \gamma_1 & \cdots & \gamma_l & \gamma'_1 & \cdots & \gamma'_{l'} \\
    \end{array}
  \right]\cdot_{op}\left[
           \begin{array}{c|c}
             A & 0 \\
             \hline
             0 & A' \\
           \end{array}
         \right]
  = \left[
    \begin{array}{cccccc}
      {\tau}(\gamma_1) & \cdots & {\tau}(\gamma_l) & {\tau}(\gamma'_1) & \cdots & {\tau}(\gamma'_{l'}) \\
    \end{array}
  \right]. \qedhere
\]}
\end{proof}
By induction, we get the following:
\begin{crl}\label{InvSubspDec}
If $V=V_1\oplus V_2\oplus\cdots\oplus V_t$ is a direct sum of $k$-spaces, and for some ${\tau}\in\txt{End}_k(V)$ each $V_i$ is ${\tau}$-invariant and has a basis $\beta_i$, then with $A_i:=\big[{\tau}|_{V_i}\big]_{\beta_i}$ and $\beta:=\beta_1\cup \beta_2\cup\cdots\cup \beta_t$,
{\scriptsize\bea
    [{\tau}]_\beta=\big[{\tau}|_{V_1}\big]_{\beta_1}\oplus\big[{\tau}|_{V_2}\big]_{\beta_2}\oplus\cdots\oplus\big[{\tau}|_{V_t}\big]_{\beta_t}=A_1\oplus A_2\oplus\cdots\oplus A_t=\left[
            \begin{array}{c|c|c|c}
              A_1 & 0 & \cdots & 0 \\
              \hline
              0 & A_2 & \cdots & 0 \\
              \hline
              \vdots & \vdots & \ddots & \vdots \\
              \hline
              0 & 0 & \cdots & A_t \\
            \end{array}
          \right].\nn
\eea}
\end{crl}

\begin{prp}[\textcolor{blue}{Equivalence of $k[x]$-module structures for a $k$-space}]\label{EquivKmodStr}
Let $V$ be a finite-dimensional $k$-space of dimension $n:=\dim_kV<\infty$. For $\tau,{\tau'}\in \txt{End}_k(V)$, let $V_\tau:={}_{k[x]}V(\tau)$ and $V_{\tau'}:={}_{k[x]}\!V({\tau'})$ be the $k[x]$-modules on which $x$ acts as $\tau$ and as ${\tau'}$ respectively. Then the following are equivalent.
\begin{enumerate}[leftmargin=0.9cm]
\item[(1)] $V_\tau\cong V_{\tau'}$ as $k[x]$-modules.
\item[(2)] There is a $k$-isomorphism $\phi:V\ra V$ such that ${\tau'}=\phi \tau\phi^{-1}$.
\item[(3)] $[\tau]_\beta\sim[{\tau'}]_\beta$, for any basis $\beta=\{\beta_1,...,\beta_n\}$ for $V$.
\end{enumerate}
\end{prp}
\begin{proof}
{\flushleft \ul{(1)$\Ra$(2)}:} Let $\phi:V_\tau\ra V_{\tau'}$ be a $k[x]$-isomorphism. Then $\phi$ is automatically a $k$-isomorphism since $k\subset k[x]$. For every $v\in V$, $\phi(xv)=x\phi(v)$ $\Ra$ $\phi(\tau(v))={\tau'}(\phi(v))$, and so
\bea
\phi\tau={\tau'}\phi,~~\Ra~~{\tau'}=\phi \tau\phi^{-1}.\nn
\eea
{\flushleft \ul{(2)$\Ra$(3)}:} Let $P:=[\phi]_\beta$, for any basis $\beta$ for $V$. Then ${\tau'}=\phi \tau\phi^{-1}$ implies $[\tau]_\beta\sim [{\tau'}]_\beta$, since
\bea
[{\tau'}]_\beta=[\phi \tau\phi^{-1}]_\beta=[\phi]_\beta[\tau]_\beta[\phi]_\beta^{-1}=P [\tau]_\beta P^{-1}.\nn
\eea
{\flushleft \ul{(3)$\Ra$(1)}:} Assume that for any basis $\beta$ for $V$,~ $[{\tau'}]_\beta=P[\tau]_\beta P^{-1}$ for some $P\in GL_n(k)$. Let $\phi:V\ra V$ be a $k$-linear map such that $[\phi]_\beta=P$. Since $P$ is invertible, $\phi$ is a $k$-isomorphism. Therefore,
    \bea
    && [{\tau'}]_\beta=P[\tau]_\beta P^{-1}=[\phi]_\beta [\tau]_\beta[\phi]_\beta^{-1}=[\phi \tau\phi^{-1}]_\beta~~~~\txt{for all}~~\beta,~~\Ra~~\phi \tau\phi^{-1} = {\tau'},\nn\\
    &&~~\Ra~~\phi \tau={\tau'}\phi\in \txt{End}_k(V),~~\Ra~~\phi \tau^l={\tau'}^l\phi,~~~~\txt{for all}~~l\geq 0.\nn
    \eea
    To show $\phi$ is $k[x]$-linear, observe that if $f=\sum_{i=1}^mf^ix^i\in k[x]$, then for $v\in V_\tau$,
\bea
\textstyle\phi(fv)=\phi\left(\sum_{i=1}^mf^i\tau^i(v)\right)=\sum_{i=1}^mf^i\phi \tau^i(v)=\sum_{i=1}^mf^i{\tau'}^i\phi(v)=f~\phi(v),~~~~\txt{where}~~\phi(v)\in V_{\tau'},\nn
\eea
and so ~$\phi:V_\tau\ra V_{\tau'}$ (i.e., $\phi:V\ra V$ is $k[x]$-linear).
~~~~$\adjustbox{scale=0.7}{\bt
(V,\beta)\ar[d,"{[\phi]_\beta}"]\ar[rr,"{[\tau]_\beta}"] && (V,\beta)\ar[d,"{[\phi]_\beta}"] \\
(V,\beta)\ar[rr,"{[{\tau'}]_\beta}"] && (V,\beta)
\et}~~~~\phi \tau=\phi x=x\phi={\tau'}\phi$
\end{proof}

The following restatement of Proposition \ref{EquivKmodStr} is important only for a more direct alternative description of the effect of a change of basis, since finite-dimensional $k$-spaces of the same dimension are isomorphic.
\begin{crl}[\textcolor{blue}{Equivalence of $k[x]$-module structures for $k$-spaces}]\label{EquivKmodStr2}
Let $V,V'$ be finite-dimensional $k$-spaces of dimension $n<\infty$. For $\tau\in\txt{End}_k(V)$ and $\tau'\in End_k(V')$, let $V_\tau:={}_{k[x]}V(\tau)$ and $V'_{\tau'}:={}_{k[x]}V'(\tau')$ be the $k[x]$-modules on which $x$ acts as $\tau$ and as $\tau'$ respectively. Then the following are equivalent.
\begin{enumerate}[leftmargin=0.9cm]
\item[(1)] $V_\tau\cong V'_{\tau'}$ as $k[x]$-modules.
\item[(2)] There is a $k$-isomorphism $\phi:V\ra V'$ such that $\tau'=\phi \tau\phi^{-1}$.
\item[(3)] $[\tau]_{\beta}\sim[\tau']_{\beta'}$, for any bases $\beta=\{\beta_1,...,\beta_n\}$ and $\beta'=\{\beta'_1,...,\beta'_n\}$ for $V$ and $V'$ respectively.
\end{enumerate}
\end{crl}
\begin{proof}
{\flushleft \ul{(1)$\Ra$(2)}:} Let $\phi:V_\tau\ra V'_{\tau'}$ be a $k[x]$-isomorphism. Then $\phi$ is automatically a $k$-isomorphism since $k\subset k[x]$. For every $v\in V$, $\phi(xv)=x\phi(v)$ $\Ra$ $\phi(\tau(v))={\tau'}(\phi(v))$, and so
\bea
\phi\tau={\tau'}\phi,~~\Ra~~{\tau'}=\phi \tau\phi^{-1}.\nn
\eea
{\flushleft \ul{(2)$\Ra$(3)}:} Let $\beta$ be a basis for $V$ and $\beta'$ a basis for $V'$. Then in terms of the isomorphisms $\tau_\beta:V\ra k^n$, $\tau'_{\beta'}:V'\ra k^n$, we have $\tau=\tau_\beta^{-1}\tau_{\beta\beta}\tau_\beta$, $\tau'=\tau'^{-1}_{\beta'}\tau'_{\beta'\beta'}\tau'_{\beta'}$ and $\phi=\tau'^{-1}_{\beta'}\phi_{\beta\beta'}\tau_\beta$. Let $P:=\phi_{\beta\beta'}=[\phi]_{\beta\beta'}$. Then ${\tau'}=\phi \tau\phi^{-1}$ implies $\tau'_{\beta'\beta'}=\phi_{\beta\beta'}\tau_{\beta\beta}\phi_{\beta\beta'}^{-1}=P\tau_{\beta\beta}P^{-1}$, and so $[\tau]_\beta\sim [{\tau'}]_{\beta'}$.

{\flushleft \ul{(3)$\Ra$(1)}:} Assume that for any bases $\beta$ for $V$ and $\beta'$ for $V'$,~ $[{\tau'}]_{\beta'}=P[\tau]_\beta P^{-1}$ for some $P\in GL_n(k)$. Let $\phi:V\ra V'$ be a $k$-linear map such that $[\phi]_{\beta\beta'}:=\phi_{\beta\beta'}=P$. Since $P$ is invertible, $\phi$ is a $k$-isomorphism. Therefore,
    \bea
    && [{\tau'}]_{\beta'}=P[\tau]_\beta P^{-1}=\phi_{\beta\beta'}\tau_{\beta\beta}\phi_{\beta\beta'}^{-1}=[\phi \tau\phi^{-1}]_{\beta'}~~~~\txt{for all}~~\beta,\beta',~~~~\Ra~~~~\phi \tau\phi^{-1} = {\tau'},\nn\\
    &&~~\Ra~~\phi \tau={\tau'}\phi\in Hom_k(V,V'),~~\Ra~~\phi \tau^l={\tau'}^l\phi,~~~~\txt{for all}~~l\geq 0.\nn
    \eea
    To show $\phi$ is $k[x]$-linear, observe that if $f=\sum_{i=1}^mf^ix^i\in k[x]$, then for $v\in V_\tau$,
\bea
\textstyle\phi(fv)=\phi\left(\sum_{i=1}^mf^i\tau^i(v)\right)=\sum_{i=1}^mf^i\phi \tau^i(v)=\sum_{i=1}^mf^i{\tau'}^i\phi(v)=f~\phi(v),~~~~\txt{where}~~\phi(v)\in V'_{\tau'},\nn
\eea
and so ~$\phi:V_\tau\ra V'_{\tau'}$ (i.e., $\phi:V\ra V'$ is $k[x]$-linear).
~~~~$\adjustbox{scale=0.7}{\bt
(V,\beta)\ar[d,"{[\phi]_{\beta\beta'}}"]\ar[rr,"{[\tau]_\beta}"] && (V,\beta)\ar[d,"{[\phi]_{\beta\beta'}}"] \\
(V',\beta')\ar[rr,"{[\tau']_{\beta'}}"] && (V',\beta')
\et}~~~~\phi \tau=\phi x=x\phi=\tau'\phi$
\end{proof}

\subsection{Rational canonical form (RCF)}
\begin{dfn}[\blue{\index{Companion matrix of a polynomial}{Companion matrix of a polynomial}}]\label{CompMatDfn}
Let $k$ be a field and $f=f^0+f^1x+\cdots+f^{d_f-1}x^{d_f-1}+x^{d_f}\in k[x]$ a monic polynomial. The \ul{companion matrix} of $f$ is
{\scriptsize\bea
C(f):=\left[
        \begin{array}{ccccccc}
          0 &  &  & &  &  & -f^0 \\
          1 & 0 &  & &  &  & -f^1 \\
           & 1 & 0 & &  &  & -f^2 \\
           &  & \ddots & \ddots &  &  & \vdots \\
           &  &  & 1 & 0 &  & -f^{d_f-3} \\
           &  &  &  & 1 & 0 & -f^{d_f-2} \\
           &  &  &  &  & 1 & -f^{d_f-1} \\
        \end{array}
      \right]:=\left[
        \begin{array}{ccccccc}
          0 & 0 & 0 &\cdots & 0 & 0 & -f^0 \\
          1 & 0 & 0 &\cdots & 0 & 0 & -f^1 \\
          0 & 1 & 0 &\cdots & 0 & 0 & -f^2 \\
          \vdots & \vdots & \vdots & \ddots & \vdots & \vdots & \vdots \\
          0 & 0 & 0 & \cdots & 0 & 0 & -f^{d_f-3} \\
          0 & 0 & 0 & \cdots & 1 & 0 & -f^{d_f-2} \\
          0 & 0 & 0 & \cdots & 0 & 1 & -f^{d_f-1} \\
        \end{array}
      \right]\in Mat_{d_f}(k).\nn
\eea}
\end{dfn}

\begin{lmm}\label{CompMatLmm}
Let $k$ be a field, $f=f^0+f^1x+\cdots+f^{d_f-1}x^{d_f-1}+x^{d_f}\in k[x]$ a monic polynomial, and $I=I_f:=k[x]f\lhd k[x]$. Then with $\beta_i:=x^{i-1}+I=x^{i-1}(1+I)$, the set
\bea
\beta_f:=\{1+I,x+I,\cdots, x^{{d_f}-1}+I\}=\{\beta_1,\beta_2,...,\beta_{d_f}\}\nn
\eea
is a $k$-basis for~ $_{k[x]}W:={k[x]\over I}={k[x]\over k[x]f}$,~ and the matrix of the $k$-linear transformation $x:W\ra W,~w\mapsto xw$ with respect to $\beta_f$ is the companion matrix of $f$, i.e., ~$[x]_{\beta_f}=C(f)$.
\end{lmm}
\begin{proof}
By the polynomial division algorithm (i.e., any $g\in k[x]$ can be written as $g=qx^{d_f}+r$ for $q,r\in k[x]$ with $r=0$ or $0\leq\deg r<\deg x^{d_f}=d_f$), it is clear that ${}_{k[x]}W:={k[x]\over I}$ has the $k$-basis $\beta_f$. Since $f\in I$ (used in the last equation below), the matrix of the $k$-linear transformation $x:W\ra W$, $w\mapsto xw$ (where $_kW:={k[x]\over I}$) is given by
\[\left.
  \begin{array}{l}
    x(\beta_1)=x(1+I)=x+I=\beta_2, \\
    x(\beta_2)=x(x+I)=x^2+I=\beta_3, \\
     \hspace{2cm}\vdots~~~~\vdots~~~~\vdots \\
    x(\beta_{d_f-1})=x(x^{{d_f}-2}+I)=x^{d_f-1}+I=\beta_{d_f},\\
    x(\beta_{d_f})=x(x^{{d_f}-1}+I)=x^{d_f}+I=-f^0\beta_1-f^1\beta_2-\cdots-f^{d_f-1}\beta_{d_f}.
  \end{array}
\right.\]
That is,
{\scriptsize
$\left[
        \begin{array}{cccc}
         x(\beta_1)  & x(\beta_2) & \cdots & x(\beta_n) \\
        \end{array}
      \right] =
      \left[
        \begin{array}{cccc}
         \beta_1  & \beta_2 & \cdots & \beta_n \\
        \end{array}
      \right]\cdot_{op}
      \left[
        \begin{array}{ccccccc}
          0 & 0 & 0 &\cdots & 0 & 0 & -f^0 \\
          1 & 0 & 0 &\cdots & 0 & 0 & -f^1 \\
          0 & 1 & 0 &\cdots & 0 & 0 & -f^2 \\
          \vdots & \vdots & \vdots & \ddots & \vdots & \vdots & \vdots \\
          0 & 0 & 0 & \cdots & 0 & 0 & -f^{d_f-3} \\
          0 & 0 & 0 & \cdots & 1 & 0 & -f^{d_f-2} \\
          0 & 0 & 0 & \cdots & 0 & 1 & -f^{d_f-1} \\
        \end{array}
      \right]$,
} and so the matrix of $x$ with respect to $\beta_f$ satisfies ~$[x]_{\beta_f}=C(f)$.
\end{proof}

\begin{thm}[\textcolor{blue}{\index{Rational! canonical form (RCF)}{Rational canonical form (RCF)} of a matrix}]\label{RCF-Theorem}
Let $k$ be a field and $A\in Mat_n(k)$ a matrix. Then there exist \ul{unique} monic polynomials $f_1,...,f_t\in k[x]$ such that $f_i|f_{i+1}$ and
{\scriptsize\bea
A\sim C(f_1)\oplus C(f_2)\oplus\cdots\oplus C(f_t)=\left[
        \begin{array}{c|c|c|c}
          C(f_1) & 0 & \cdots & 0 \\
          \hline
         0  & C(f_2) & \cdots & 0 \\
           \hline
          \vdots & \vdots & \ddots & \vdots \\
           \hline
          0 & 0 & \cdots & C(f_t) \\
        \end{array}
      \right].\nn
\eea}This result is called the \ul{rational canonical form} of $A$.
\end{thm}
\begin{proof}
Let $V$ be an $n$-dimensional $k$-space, ${\tau}:V\ra V$ a $k$-linear transformation, and $\beta$ a $k$-basis for $V$ such that $[{\tau}]_\beta=A$. We know $V$ is a $k[x]$-module wrt the scalar multiplication
\bea
k[x]\times V\ra V,~\big(f,v\big)\mapsto fv:=f({\tau})(v).\nn
\eea
By Lemma \ref{RatDecLmm}, we have unique monic polynomials $f_1,...,f_t\in k[x]\backslash 0$, satisfying $f_i|f_{i+1}$, such that
{\small\bea
\textstyle V={}_{k[x]}V(\tau)\sr{\phi}{\cong} {k[x]\over k[x]f_1}\oplus {k[x]\over k[x]f_2}\oplus\cdots\oplus {k[x]\over k[x]f_t}=W_1\oplus W_2\oplus\cdots\oplus W_t=:V',~~~~\txt{where}~~~W_i:={k[x]\over k[x]f_i}.\nn
\eea}
By Corollary \ref{InvSubspDec}, the matrix of the action of $x:V'\ra V'$ (equivalently, of $\tau$) in the basis $\beta':=\beta_{f_1}\cup\beta_{f_2}\cup\cdots\cup\beta_{f_t}$, based on Lemma \ref{CompMatLmm}, is
~{\footnotesize $[x]_{\beta'}=\big[x|_{W_1}\big]_{\beta_{f_1}}\oplus \big[x|_{W_2}\big]_{\beta_{f_2}}\oplus\cdots\oplus\big[x|_{W_t}\big]_{\beta_{f_t}}=C(f_1)\oplus C(f_2)\oplus\cdots\oplus C(f_t)$.}~ Hence, by the properties (from Corollary \ref{EquivKmodStr2}(1) with $\tau=\tau'=x$) of the isomorphism $\phi:V\ra V'$, we get
\[
A=[\tau]_\beta=[x]_\beta\sim [x]_{\beta'}=C(f_1)\oplus C(f_2)\oplus\cdots\oplus C(f_t).~~~~~~~~
\adjustbox{scale=0.7}{\bt
(V,\beta)\ar[d,"{[\phi]_{\beta\beta'}}"]\ar[rr,"{[x]_\beta}"] && (V,\beta)\ar[d,"{[\phi]_{\beta\beta'}}"] \\
(V',\beta')\ar[rr,"{[x]_{\beta'}}"] && (V',\beta')
\et} \qedhere
\]
\end{proof}

\subsection{Jordan canonical form (JCF), diagonalization, and determinant}
\begin{dfn}[\textcolor{blue}{\index{Root (or zero) of a polynomial}{Root (or zero) of a polynomial}}]
Let $k$ be a field, $a\in k$, and $f\in k[x]$. Then $a$ is a \ul{root} (or \ul{zero}) of $f$ if $f(a)=0$.
\end{dfn}

\begin{lmm}[\textcolor{blue}{Polynomial root-factorization}]
Let $k$ be a field, $a\in k$, and $f\in k[x]$. Then $f(a)=0$ (i.e., $a$ is a root of $f$ in $k$) $\iff$ $(x-a)|f$, i.e., $\iff$ $f\in k[x](x-a)$.
\end{lmm}
\begin{proof}
($\Ra$): Assume $f(a)=0$. By the division algorithm (for $k[x]$ as a ED),
\bea
f=q(x-a)+r,~~\txt{for some}~~q,r\in k[x]~~\txt{with}~~r=0~~\txt{or}~~0\leq\deg r<\deg (x-a)=1,\nn
\eea
and so $\deg r=0$, i.e., $r\in k$. Since $f(a)=0$, we get $r=f(a)=0$, and so $(x-a)|f$.

($\La$): Conversely, if $(x-a)|f$, then $f=q(x-a)$ for some $q\in k[x]$, and so $f(a)=0$.
\end{proof}

\begin{dfn}[\textcolor{blue}{\index{Algebraically closed field}{Algebraically closed field}}]
A field $k$ is \ul{algebraically closed} if every monic polynomial $f\in k[x]\backslash 0$ has a root in $k$, i.e., $f(a)=0$ for some $a\in k$ (equivalently, $f$ can be written as $f=(x-a_1)(x-a_2)\cdots(x-a_l)$ for some $a_i\in k$ and $l=\deg f$). (\blue{footnote}\footnote{The word ``monic'' is included only for convenience of expression. A field $k$ is \ul{algebraically closed} iff every polynomial $f\in k[x]\backslash 0$ has a root in $k$.})
\end{dfn}

\begin{thm}[\textcolor{blue}{\index{Fundamental theorem of! algebra}{Fundamental theorem of algebra}}]
The field $\Complex$ is algebraically closed.
\end{thm}
\begin{proof}
Let $f\in k[x]$ be monic (hence nonconstant). Since $|f(z)|\ra\infty$ as $|z|\ra\infty$, $\inf|f|:=\inf_{z\in\Complex}|f(z)|=\inf_{z\in D}|f(z)|$ for some disk $D=\{z\in\Complex:|z|\leq r_D\}$ of radius $r_D>0$. (\magenta{Fast-forward to Corollaries \ref{ClosedImLmm} and \ref{HeinBorelCrl}: The function $|f|:\Complex\ra\Real,~z\mapsto|f(z)|$ is continuous and so maps the compact set $D\subset\Complex$ to a compact set $|f(D)|\subset\Real$, which is closed because $(\Real,|~|)$ is Hausdorff, and so $|f(D)|$ contains its infimum}). Thus we can pick some $a\in D$ such that $|f(a)|=\inf|f|$. It is clear that we can rewrite $f$ in the form
\[
f(x)=f(a)+c_{i_a}(a)(x-a)^{i_a}+c_{i_a+1}(a)(x-a)^{i_a+1}+\cdots+c_{d_f}(a)(x-a)^{d_f},~~~~c_{i_a}(a)\neq 0,~~i_a\geq 1.\nn
\]
Suppose $f$ has no root in $\Complex$, i.e., $f(a)\neq 0$. Let $g(x):=f(a)+c_{i_a}(a)(x-a)^{i_a}$. Then
\begin{align}
&\textstyle|f(z)|\leq|g(z)+f(z)-g(z)|\leq |g(z)|+|f(z)-g(z)|=|g(z)|+|z-a|^{i_a+1}{|f(z)-g(z)|\over|z-a|^{i_a+1}},\nn\\
&~~\Ra~~|f(a+re^{i\theta})|\leq |f(a)+r^{i_a}e^{i\theta i_a}c_{i_a}(a)|+r^{i_a+1}C=\big||f(a)|+r^{i_a}e^{i\theta' i_a}|c_{i_a}(a)|\big|+r^{i_a+1}C,\nn
\end{align}
where $0<r<<|D|$, the bound $C>0$ is independent of $r$, and $\theta':=\theta+\arg c_{i_a}(a)-\arg f(a)$ is given by
\[
c_{i_a}(a)=|c_{i_a}(a)|e^{i\arg c_{i_a}(a)}~~~~\txt{and}~~~~f(a)=|f(a)|e^{i\arg f(a)}.
\]
If we choose $r$ such that $|f(a)|>r^{i_a}|c_{i_a}(a)|$, then by setting $\theta=\theta_0:=\pi+\arg f(a)-\arg c_{i_a}(a)$, we get
\[
|f(a+re^{i\theta_0})|\leq |f(a)|-r^{i_a}|c_{i_a}(a)|+r^{i_a+1}C.
\]
Hence, with a sufficiently small $r$, we get $|f(a+re^{i\theta_0})|<|f(a)|$, which is a contradiction.
\end{proof}

\begin{dfn}[\textcolor{blue}{\index{Jordan block}{Jordan block}}]\label{JordBlockDfn}
Let $k$ be a field and $a\in k$. A Jordan block of \ul{size} $l$ and \ul{constant} $a$ is
{\scriptsize\[
J_l(a):=\left[
         \begin{array}{cccccc}
          a & 1 &   &       &   &   \\
            & a & 1 &       &   &   \\
            &   & a &  1     &   &   \\
            &   &   &\ddots & \ddots  &   \\
            &   &   &       & a & 1 \\
            &   &   &       &   & a \\
         \end{array}
       \right]:=\left[
         \begin{array}{cccccc}
           a & 1 & 0 &\cdots& 0 & 0 \\
           0 & a & 1 & \cdots& 0 & 0 \\
           0 & 0 & a & \cdots& 0 & 0 \\
           \vdots & \vdots &\vdots &\ddots & \vdots & \vdots \\
           0 & 0 & 0 & \cdots & a & 1 \\
           0 & 0 & 0 & \cdots & 0 & a \\
         \end{array}
       \right]\in Mat_l(k).\nn
\]}
\end{dfn}

\begin{lmm}\label{JordBlockLmm}
Let $k$ be a field, $a\in k$, and $I=I_{a,l}:=k[x](x-a)^l\lhd k[x]$, where $l\geq 1$. Then the set
\bea
\beta_{a,l}:=\{\beta_1,...,\beta_l\},~~~~\beta_i:=(x-a)^{l-i}+I=(x-a)^{l-i}(1+I),\nn
\eea
is a k-basis for~ $_{k[x]}W:={k[x]\over I}={k[x]\over k[x](x-a)^l}$,~ and the matrix of the $k$-linear transformation $x:W\ra W,~w\mapsto xw$ with respect to $\beta_{a,l}$ is the Jordan block of length $l$ and constant $a$, i.e., ~$[x]_{\beta_{a,l}}=J_l(a)$.
\end{lmm}
\begin{proof}
It is clear that $(x-a)^{l-1},\cdots,(x-a)^1,(x-a)^0\in k[x]$ are linearly independent over $k$, and so span an $l$-dimensional subspace of $k[x]$ as a $k$-space. Therefore, by the division algorithm (i.e., any $g\in k[x]$ can be written as $g=q(x-a)^l+r$ for some $q,r\in k[x]$ with $r=0$ or $0\leq\deg r<\deg(x-a)^l=l$), $\beta_{a,l}$ is a $k$-basis for $W:={k[x]\over k[x](x-a)^l}={k[x]\over I}$, where $I:=k[x](x-a)^l$. With $\beta_i=(x-a)^{l-i}+I$, we have
{\small\bea
&&x(\beta_i)=x(x-a)^{l-i}+I=(x-a)(x-a)^{l-i}+a(x-a)^{l-i}+I=(x-a)^{l-(i-1)}+a(x-a)^{l-i}+I\nn\\
&&~~~~=[(x-a)^{l-(i-1)}+I]+a[(x-a)^{l-i}+I]
=\left\{
                    \begin{array}{ll}
                      a\beta_1, & i=1 \\
                      \beta_{i-1}+a\beta_i, & 2\leq i\leq l
                    \end{array}
                  \right.\nn
\eea}That is, the matrix of $x:W\ra W,~w\mapsto xw$ with respect to $\beta_{a,l}$ satisfies ~$[x]_{\beta_{a,l}}=J_l(a)$,~ since
{\scriptsize
\[
\left[
           \begin{array}{cccc}
             x(\beta_1) & x(\beta_2) & \cdots & x(\beta_l) \\
           \end{array}
         \right]=
         \left[
           \begin{array}{cccc}
             \beta_1 & \beta_2 & \cdots & \beta_l \\
           \end{array}
         \right]\cdot_{op}
        \left[
         \begin{array}{cccccc}
           a & 1 & 0 &\cdots& 0 & 0 \\
           0 & a & 1 & \cdots& 0 & 0 \\
           0 & 0 & a & \cdots& 0 & 0 \\
           \vdots & \vdots &\vdots &\ddots & \vdots & \vdots \\
           0 & 0 & 0 & \cdots & a & 1 \\
           0 & 0 & 0 & \cdots & 0 & a \\
         \end{array}
       \right].\nn \hspace{2cm}\qedhere
\]}
\end{proof}

\begin{thm}[\textcolor{blue}{\index{Jordan canonical form (JCF)}{Jordan canonical form (JCF)} of a matrix}]\label{JCF-Theorem}
Let $k$ be an algebraically closed field and $A\in Mat_n(k)$. There exist Jordan blocks $J_{l_1}(a_1),\cdots,J_{l_t}(a_t)$ (not necessarily distinct) such that we have a unique expression (up to reordering)
{\scriptsize\bea
A\sim J_{l_1}(a_1)\oplus J_{l_2}(a_2)\oplus\cdots\oplus J_{l_t}(a_t)
=\left[
        \begin{array}{c|c|c|c}
          J_{l_1}(a_1) & 0 & \cdots & 0 \\
          \hline
         0  & J_{l_2}(a_2) & \cdots & 0 \\
           \hline
          \vdots & \vdots & \ddots & \vdots \\
           \hline
          0 & 0 & \cdots & J_{l_2}(a_t) \\
        \end{array}
      \right].\nn
\eea}This result is called the \ul{Jordan canonical form} (JCF) of $A$.
\end{thm}
\begin{proof}
Since $k$ is algebraically closed, the only irreducible (prime) elements of $k[x]$ are $\{x-a:~a\in k\}$. Let $V$ be an $n$-dimensional $k$-space, ${\tau}:V\ra V$ a $k$-linear transformation, and $\beta$ a $k$-basis for $V$ such that $[{\tau}]_\beta=A$. Then $V$ is a $k[x]$-module with scalar multiplication $k[x]\times V\ra V,~(f,v)\mapsto f(\tau)(v)$. Since $_{k[x]}V$ is a FG torsion $k[x]$-module (by Lemma \ref{VecPolTor}), it follows from Theorem \ref{FT-FG-PID} that
{\footnotesize\bea
\textstyle V={}_{k[x]}V(\tau)\sr{\phi}{\cong} {R\over R(x-a_1)^{l_1}}\oplus {R\over R(x-a_2)^{l_2}}\oplus\cdots\oplus {R\over R(x-a_t)^{l_t}}=W_1\oplus W_2\oplus\cdots\oplus W_t=:V',~~~\txt{where}~~~W_i:={R\over R(x-a_i)^{l_i}},\nn
\eea}for some $a_i\in k$, $l_i\geq 1$ (where the isomorphism is both of $k[x]$-modules and of $k$-spaces).

By Corollary \ref{InvSubspDec}, the matrix of the action of $x:V'\ra V'$ (equivalently, of $\tau$) in the basis $\beta':=\beta_{a_1,l_1}\cup\beta_{a_2,l_2}\cup\cdots\cup\beta_{a_t,l_t}$, based on Lemma \ref{JordBlockLmm}, is
~{\footnotesize $[x]_{\beta'}=\big[x|_{W_1}\big]_{\beta_{a_1,l_1}}\oplus \big[x|_{W_2}\big]_{\beta_{a_2,l_2}}\oplus\cdots\oplus\big[x|_{W_t}\big]_{\beta_{a_t,l_t}}=J_{l_1}(a_1)\oplus J_{l_2}(a_2)\oplus\cdots\oplus J_{l_t}(a_t)$.}~ Thus by the properties (from Corollary \ref{EquivKmodStr2}(1) with $\tau=\tau'=x$) of the isomorphism $\phi:V\ra V'$, we get
\[
A=[\tau]_\beta=[x]_\beta\sim [x]_{\beta'}=J_{l_1}(a_1)\oplus J_{l_2}(a_2)\oplus\cdots\oplus J_{l_t}(a_t).
~~~~~~~~
\adjustbox{scale=0.7}{\bt
(V,\beta)\ar[d,"{[\phi]_{\beta\beta'}}"]\ar[rr,"{[x]_\beta}"] && (V,\beta)\ar[d,"{[\phi]_{\beta\beta'}}"] \\
(V',\beta')\ar[rr,"{[x]_{\beta'}}"] && (V',\beta')
\et} \qedhere
\]
\end{proof}

\begin{dfn}[\blue{\index{Diagonal matrix}{Diagonal matrix}}]
Let $k$ be a field. A \ul{diagonal matrix} over $k$ is a matrix of the form
\[
a_1\oplus a_2\oplus\cdots\oplus a_n
=\left[
   \begin{array}{cccc}
     a_1 & 0 & \cdots & 0 \\
     0 & a_2 & \cdots & 0 \\
     \vdots & \vdots & \ddots & \vdots \\
     0 & 0 & \cdots & a_n \\
   \end{array}
 \right]\in Mat_n(k)
\]
for some $a_i\in Mat_1(k)=k$.
\end{dfn}

\begin{dfn}[\blue{
\index{Diagonalizable matrix}{Diagonalizable matrix},
\index{Eigenvalue}{Eigenvalue},
\index{Diagonalizable linear transformation}{Diagonalizable linear transformation},
\index{Eigenvector}{Eigenvector}}]
Let $k$ be a field, $A\in Mat_n(k)$, $V$ a finite-dimensional $k$-space, and $\tau\in End_k(V):=Hom_k(V,V)$. Then the matrix $A$ is \ul{diagonalizable} if it is conjugate to a diagonal matrix, i.e., if
\[
A\sim\Ld,~~~\txt{for some diagonal matrix}~~~\Ld=\ld_1\oplus \ld_2\oplus\cdots\oplus\ld_n\in Mat_n(k),
\]where the scalars $\ld_i\in k$ are called \ul{eigenvalues} of $A$. The linear transformation $\tau$ is \ul{diagonalizable} if in some basis (and hence in any basis) $\beta$ for $V$ the matrix $[\tau]_\beta$ is diagonalizable. (\blue{footnote}\footnote{It is possible to replace $k$ with any commutative ring $R$, replace $V$ with any free $R$-module $F$ of finite rank, and consider the notion of a matrix representation $[\tau]_\beta$ (in a basis $\beta$ for $F$) of an $R$-endomorphism $\tau\in End_R(F):=Hom_R(F,F)$.})

A vector $v\in V\backslash 0$ is an \ul{eigenvector} of $\tau$ if $\tau(v)=\ld v$ for a scalar $\ld$ (called an \ul{eigenvalue} of $\tau$ for the eigenvector $v$). (\blue{footnote}\footnote{Equivalently, a scalar $\ld\in k$ is an \ul{eigenvalue} of $\tau$ if $\tau(v)=\ld v$ for a nonzero vector $v\in V\backslash 0$ (called an \ul{eigenvector} of $\tau$ with an eigenvalue $\ld$).})
\end{dfn}

The following is a corollary of Theorem \ref{JCF-Theorem}.
\begin{crl}[\blue{Diagonalizability over an algebraically closed field}]\label{JCF-Thm-Cor1}
Let $k$ be an algebraically closed field and $A\in Mat_n(k)$. Then $A$ is diagonalizable $\iff$ (1) in the JCF ~$A\sim J_{l_1}(a_1)\oplus J_{l_2}(a_2)\oplus\cdots\oplus J_{l_t}(a_t)$,~ we have $l_1=l_2=\cdots=l_t$, and (2) the scalars $J_{l_i}(a_i)=J_1(a_i)=a_i$ are \ul{unique eigenvalues} of $A$.
\end{crl}

From here onward given a linear transformation $\tau:V\ra W$ we will write $\tau(v)$ simply as $\tau v$ (for $v\in V$).

\begin{dfn}[\textcolor{blue}{
\index{Generalized eigenvector}{Generalized eigenvector},
\index{Generalized eigenvalue}{Generalized eigenvalue},
\index{Exponent of a generalized eigenvector}{Exponent of a generalized eigenvector}}]
Let $V$ be a $k$-space, $\tau\in End_k(V)$, and $\mathds{1}:=id_V\in End_k(V)$ the identity transformation. A \ul{generalized eigenvector} (\ul{gen-eigenvector}) of $\tau$ is a nonzero vector $v\in V\backslash 0$ such that
\bea
(\tau-\ld \mathds{1})^nv=0~~~~\txt{for some $\ld\in k$ and some positive integer $n\geq1$},\nn
\eea
where the scalar $\ld\in k$ is a \ul{generalized eigenvalue} (\ul{gen-eigenvalue}) of $\tau$ for the gen-eigenvector $v$. The \ul{$\ld$-exponent} of the gen-eigenvector $v$ is the smallest positive integer $d=d_\ld$ such that $(\tau-\ld \mathds{1})^dv=0$.
\end{dfn}
We note that an eigenvector of $\tau$ is a gen-eigenvector of $\tau$ of exponent $1$.

\begin{prp}
Let $V$ be a $k$-space and $\tau\in End_k(V)$. Let $v$ be a gen-eigenvector of $\tau$, with a gen-eigenvalue $\ld$ and exponent $d$. For $1\leq j\leq d$, set {\small $w_j:=(\tau-\ld \mathds{1})^{d-j}v$}. Let {\small $\beta:=(w_1,...,w_d)$} and $W:=\Span \beta$. Then (i) $W$ is a $\tau$-invariant subspace of $V$ and (ii) $\beta$ is a basis for $W$.
\end{prp}

\begin{proof}
{\flushleft (i)} By direct computation, we have (for each $1\leq j\leq d$)
{\small\begin{align}
&\tau w_j=\tau(\tau-\ld \mathds{1})^{d-j}v=(\tau-\ld\mathds{1})(\tau-\ld \mathds{1})^{d-j}v+\ld\mathds{1}(\tau-\ld \mathds{1})^{d-j}v=(\tau-\ld \mathds{1})^{d-(j-1)}v+\ld (\tau-\ld \mathds{1})^{d-j}v\nn\\
\label{JordProp1Eq0}&~~~~=w_{j-1}+\ld w_j=
\left\{
                          \begin{array}{ll}
                            \ld w_1, & \hbox{if}~~j=1 \\
                            \ld w_j+w_{j-1}, & \hbox{if}~~2\leq j\leq d
                          \end{array}
                        \right\}~\in~W,~~~~\Ra~~\tau(W)\subset W.
\end{align}}
{\flushleft (ii)} It is clear that $w_j\neq 0$ for all $j=1,...,d$ by the definition of the exponent $d$. For $c_1,...,c_d\in k$, if
\bea
\label{JordProp1Eq1} c_1w_1+...+c_dw_d=0,
\eea
then by applying appropriate powers $(\tau-\ld \mathds{1})^r$ with $1\leq r\leq d-1$ on (\ref{JordProp1Eq1}), we obtain $c_1=c_2=...=c_d=0$, and thus $\beta$ is linearly independent, hence a basis for $W$.
\end{proof}

\begin{crl}[\blue{Gen-eigenvalues are eigenvalues}]
Let $V$ be a $k$-space and $\tau\in End_k(V)$. Let $v$ be a gen-eigenvector of $\tau$, with a gen-eigenvalue $\ld$ and exponent $d$. Then $\ld$ is an eigenvalue of $\tau$.
\end{crl}
\begin{proof}
From (\ref{JordProp1Eq0}), the vector {\small $w_{d-1}:=(\tau-\ld \mathds{1})^{d-1}v$} is an eigenvector of $\tau$ with an eigenvalue $\ld$.
\end{proof}

\begin{dfn}[\textcolor{blue}{Recall: \index{Jordan block}{Jordan block}}]
Let $V$ be a $k$-space and $\tau\in End_k(V)$. Let $v$ be a generalized eigenvector of $\tau$, with a generalized eigenvalue $\ld$ and exponent $d$. The action of $\tau$ on $\beta=(w_1,...,w_d)$, $w_j:=(\tau-\ld \mathds{1})^{d-j}v$, as shown in equation (\ref{JordProp1Eq0}) produces a nonzero matrix $J_d(\ld)\in Mat_d(k)$ called a \ul{Jordan block} of $\tau$ of size $d$ corresponding to $\ld$, i.e., $\tau(\beta)^T=\beta^T\cdot_{op}[\tau]_\beta=\beta^T\cdot_{op}J_d(\ld)$, where
{\scriptsize\[
[\tau]_\beta:=J_d(\ld):=\left[
        \begin{array}{ccccccc}
       \ld & 1   &     &        &      &   \\
           & \ld &  1  &        &      &   \\
           &     & \ld &  1    &      &   \\
           &     &     & \ddots &   \ddots   &   \\
           &     &     &        &  \ld & 1 \\
           &     &     &        &      & \ld \\
        \end{array}
      \right].\nn
\]}
\end{dfn}

\begin{lmm}[\blue{
\index{Levi-Civita symbol}{Levi-Civita symbol},
\index{Sign of a permutation}{Sign of a permutation},
\index{Determinant}{Determinant},
\index{Adjugate matrix}{Adjugate matrix},
\index{Minors}{Minors},
\index{Determinant-adjugate formula}{Determinant-adjugate formula},
\index{Inverse! matrix over a field}{Inverse matrix over a field}}]\label{MatrixDetLmm}
Let $R$ be a commutative ring and $A=[A_{ij}]\in Mat_n(R)$. Consider the function (called \ul{Levi-Civita symbol}) $\vep:\{1,2,\cdots,n\}^n\ra\Integer$ that gives the collection of integers $\{\vep_{i_1i_2\cdots i_n}:1\leq i_1,...,i_n\leq n\}$ defined as follows: $\vep_{12\cdots n}:=1$ and for each adjacent pair of indices $i_k,i_{k+1}$ ($k=1,...,n$), {\small $\vep_{i_1i_2\cdots i_ki_{k+1}\cdots i_n}=-\vep_{i_1i_2\cdots i_{k+1}i_k \cdots i_n}$}. Equivalently,
\bea
\textstyle \vep_{i_1i_2\cdots i_n}:=\sum_{\sigma\in S_n}\txt{sign}(\sigma)\delta_{\sigma(1)i_1}\delta_{\sigma(2)i_2}\cdots\delta_{\sigma(n)i_n}=\sum_{\sigma\in S_n}(-1)^{N(\sigma)}\delta_{\sigma(1)i_1}\delta_{\sigma(2)i_2}\cdots\delta_{\sigma(n)i_n},\nn
\eea
where $\txt{sign}(\sigma)=1$ (resp. $-1$) if $\sigma$ is an even (resp. odd) permutation, and $N(\sigma)$ is the number of transpositions required to return $\sigma$ to the identity permutation in $S_n$.

Then the \ul{determinant} of $A$ is the signed sum of products of entries of $A$ given by
\bea
\textstyle\det A:=\sum \vep_{i_1i_2...i_n}A_{1i_1}A_{2i_2}...A_{ni_n}={1\over n!}\sum\vep_{i_1i_2...i_n}\vep_{j_1j_2...j_n}A_{i_1j_1}A_{i_2j_2}...A_{i_nj_n},\nn
\eea
where the sum is over the induces $i_1,...,i_n$ and/or $j_1,...,j_n$. Consider the matrix ~$adj(A)\in Mat_n(R)$~ (called \ul{adjugate matrix} of $A$) with entries given by
{\small\begin{align}
&\textstyle (adj(A))_{ts}=(adj(A)^T)_{st}:=\sum\vep_{i_1i_2\cdots i_{t-1}si_{t+1}\cdots i_n}A_{1i_1}A_{2i_2}\cdots A_{t-1,i_{t-1}}A_{t+1,i_{t+1}}\cdots A_{ni_n}=(-1)^{s+t}\det M_{st}(A),\nn
\end{align}}where the contents of the expression are described as follows:
\bit[leftmargin=0.7cm]
\item the sum is over the indices $i_1,i_2,\cdots,i_{t-1},i_{i+1},\cdots,i_n$ whenever they make sense,
\item $M_{st}(A):=[A_{ij}]_{i\neq s,j\neq t}\in Mat_{n-1}(k)$ is the submatrix of $A$ formed by deleting both the row and the column containing (i.e., intersecting at) the entry $A_{st}$ of $A$, and
\item the (determinant of the) matrix $M_{st}(A)$ is called the $(s,t)^{th}$ \ul{minor} of $A$.
\eit
Then by direct observation we get the following equality (call it the \ul{determinant-adjugate formula}):
{\small\begin{align}
\textstyle \big(A~adj(A)^T\big)_{st}=\sum_j A_{sj}~adj(A)^T_{jt}=
\sum \vep_{i_1i_2...i_{t-1}ji_{t+1}...,i_n}A_{1i_1}A_{2i_2}...A_{t-1,i_{t-1}}A_{sj}A_{t+1,i_{t+1}}...A_{ni_n}=\det(A)\delta_{st}.\nn
\end{align}}
That is, ~$A\!~adj(A)^T=\det(A)\mathds{1}_n=adj(A)^TA$,~ where ~$\mathds{1}_n:=[\delta_{ij}]\in Mat_n(k)$~ is the identity matrix.

It follows therefore that for any field $k$, a matrix $A\in Mat_n(k)$ is invertible $\iff$ $\det A\neq 0$, in which case it is clear that the \ul{inverse matrix} of $A$ is given by the equality ~$A^{-1}={1\over\det A}adj(A)^T$.
\end{lmm}

\begin{thm}[\blue{
\index{Finding eigenvalues}{Finding eigenvalues},
\index{Characteristic! polynomial}{Characteristic polynomial},
\index{Characteristic! equation}{Characteristic equation}}]
Let $V$ be a $k$-space and $\tau\in End_k(V)$. Then a scalar $\ld\in k$ is an eigenvalue of $\tau$ $\iff$ it is a root of the polynomial
\bea
p_\tau(x):=\det(\tau-x\mathds{1}):=\det\big([\tau]_\beta-x[\mathds{1}]\big)\in k[x],~~~~\txt{for any basis $\beta$ for $V$},\nn
\eea
which is called the \ul{characteristic polynomial} of $\tau$ (and of any matrix $[\tau]_\beta$ of $\tau$ in any basis $\beta$ for $V$). The \ul{characteristic equation} of $\tau$ (and of any matrix $[\tau]_\beta$ of $\tau$ in any basis $\beta$ for $V$) is the equation
\[ p_\tau(x)=0.\]
\end{thm}
\begin{proof}
If $\ld\in k$ is an eigenvalue of $\tau$, then there is a nonzero vector $v\in V\backslash 0$ such that $(\tau-\ld\mathds{1})v=0$. Since $v\neq 0$, $\tau-\ld\mathds{1}$ (hence any matrix $[\tau-\ld\mathds{1}]_\beta$) is not invertible, and so $\det[\tau-\ld\mathds{1}]_\beta=0$, since a matrix over a field is invertible iff its determinant is nonzero (Lemma \ref{MatrixDetLmm}). Because any two bases $\beta,\beta'$ for $V$ are expressible in terms of each other as $\beta'=P^T\beta$ for an invertible matrix $P:=[id_V]_{\beta\beta'}\in GL(k)$, the respective matrix representations $\tau(\beta)=[\tau]_\beta^T\beta$ and $\tau(\beta')=[\tau]_{\beta'}^T\beta'$ imply that $[\tau]_\beta\sim [\tau]_{\beta'}$, since
{\small\[
P^T\tau(\beta)=\tau(P^T\beta)=\tau(\beta')=[\tau]_{\beta'}^T\beta'=[\tau]_{\beta'}^TP^T\beta=P^T\big(P[\tau]_{\beta'}P^{-1}\big)^T\beta~~~\Ra~~~[\tau]_{\beta'}=P^{-1}[\tau]_{\beta}P.\nn
\]}Therefore ~$\det[\tau-\ld\mathds{1}]_\beta=\det([\tau]_\beta-\ld[\mathds{1}])$~ is independent of the choice of basis $\beta$ for $V$.
\end{proof}

\begin{crl}
Let $k$ be a field and $f\in k[x]$ a monic polynomial. The roots of $f$ are the eigenvalues of its companion matrix, because
\[
p_{C(f)}(x):=\det\big(C(f)-x\mathds{1}\big)=-f(x).
\]
\end{crl}

\begin{crl}[\textcolor{blue}{Distinct eigenvalue diagonalizability}]
Let $k$ be an algebraically closed field, $V$ a k-space of dimension $n$, and $\tau\in End_k(V)\cong Mat_n(k)$. If the characteristic polynomial of $\tau$ has $n$ distinct roots (i.e., $n$ distinct zeros), then the matrix of $\tau$ is diagonalizable.
\end{crl}
\begin{proof}
Let $\tau\in End_k(V)$ be such that $p_\tau$ has $n$ distinct roots. Then in the expansion
\[
\textstyle p_\tau(x)=p_{[\tau]}(x)=p_{JCF([\tau])}(x)=\prod_{i=1}^tp_{J_{l_i}(a_i)}(x),~~~~t\leq n=l_1+l_2+\cdots+l_t,
\]
we must have $t=n$ since a Jordan block cannot give more than one eigenvalue. So, all Jordan blocks are of size 1. Hence, the matrix of $\tau$ is diagonalizable by Corollary \ref{JCF-Thm-Cor1}.
\end{proof}

\begin{thm}[\blue{\index{Cayley-Hamilton theorem}{Cayley-Hamilton theorem}}]\label{CayHamThm}
Let $k$ be a field, $A\in Mat_n(k)$, and ~$p(x):=\det(x\mathds{1}-A)=x^n+c_{n-1}x^{n-1}+...+c_1x+c_0$.~ Then ~$p(A)=A^n+c_{n-1}A^{n-1}+...+c_1A+c_0\mathds{1}=0$~ in $Mat_n(k)$. That is, a matrix satisfies its characteristic equation. (Equivalently, if $V$ is a $k$-space and $\tau\in End_k(V)$, then $p_\tau(\tau)=0$, i.e., a linear transformation satisfies its characteristic equation.)
\end{thm}
\begin{proof}
Recall from Lemma \ref{MatrixDetLmm} that a square matrix $Q=[q_{ij}]\in Mat_n(R)$ over any commutative ring $R$ satisfies {\small $\det Q~\mathds{1}=Q~\txt{adj}(Q)^T$}, where {\small $\txt{adj}(Q):=[(-1)^{i+j}\det M_{ij}(Q)]$}, with $M_{ij}(Q)$ the submatrix of $Q$ formed by deleting the $i$th row and the $j$th column of $Q$. In particular, with {\small $Q:=x\mathds{1}-A\in Mat_n(k[x])$},
\bea
\label{CayleyHamEqn1}p(x)~\mathds{1}=\det(x\mathds{1}-A)~\mathds{1}=(x\mathds{1}-A)~\big(\txt{adj}(x\mathds{1}-A)\big)^T,
\eea
where all entries of~ $\txt{adj}(x\mathds{1}-A)$~ are polynomials of degree $\leq~n-1$. Let $c_n:=1$ (just for convenience) so that $p(x)=c_nx^n+c_{n-1}x^{n-1}+...+c_1x+c_0$. For some matrices $B_i\in Mat_n(k)$, we have
\bea
\txt{adj}(x\mathds{1}-A)=B_{n-1}x^{n-1}+...+B_1x+B_0,\nn
\eea
and so (\ref{CayleyHamEqn1}) takes the expanded form
\bea
&&(c_nx^n+c_{n-1}x^{n-1}+...+c_1x+c_0)\mathds{1}=(x\mathds{1}-A)(B_{n-1}x^{n-1}+...+B_1x+B_0)\nn\\
&&~~~~=B_{n-1}x^n+(B_{n-2}-AB_{n-1})x^{n-1}+...+(B_1-AB_2)x^2+(B_0-AB_1)x-AB_0.\nn
\eea
Equating the two coefficients of $x^i$ from both sides (for each $i=0,1,...,n$) we obtain Table \ref{EqCoeffTab1} below.

\begin{center}
\adjustbox{scale=0.9}{
\begin{minipage}{5cm}
\begin{table}[H]
  \centering
\begin{tabular}{l|l}
  Var. & Equation of coefficients \\
\hline
  $x^n$ & $c_n\mathds{1}=B_{n-1}$ \\
  $x^{n-1}$ & $c_{n-1}\mathds{1}=B_{n-2}-AB_{n-1}$ \\
  $x^{n-2}$ & $c_{n-2}\mathds{1}=B_{n-3}-AB_{n-2}$ \\
  $\vdots$ & \hspace{1.2cm} $\vdots$ \\
  $x$ & $c_1\mathds{1}=B_0-AB_1$ \\
  $x^0$ & $c_0\mathds{1}=-AB_0$ \\
  \hline
\end{tabular}
  \caption{}\label{EqCoeffTab1}
\end{table}
\end{minipage}}
\hspace{1.5cm}
\adjustbox{scale=0.9}{
\begin{minipage}{12cm}
\begin{table}[H]
  \centering
\begin{tabular}{l|l}
  Var. & Multiplied equation of coefficients \\
\hline
  $x^n$ & $c_nA^n=A^nB_{n-1}$ \\
  $x^{n-1}$ & $c_{n-1}A^{n-1}=A^{n-1}(B_{n-2}-AB_{n-1})=A^{n-1}B_{n-2}-A^nB_{n-1}$ \\
  $x^{n-2}$ & $c_{n-2}A^{n-2}=A^{n-2}(B_{n-3}-AB_{n-2})=A^{n-2}B_{n-3}-A^{n-1}B_{n-2}$ \\
  $\vdots$ & \hspace{1.2cm} $\vdots$ \\
  $x$ & $c_1A=A(B_0-AB_1)=AB_0-A^2B_1$ \\
  $x^0$ & $c_0\mathds{1}=-AB_0$ \\
  \hline
\end{tabular}
  \caption{}\label{EqCoeffTab2}
\end{table}
\end{minipage}}
\end{center}

In Table \ref{EqCoeffTab1}, left-multiplying each equation of coefficients by the corresponding power of $A$, we get Table \ref{EqCoeffTab2} above. Adding the equations in Table \ref{EqCoeffTab2} together, the RHS vanishes, and so we get ~$p(A)=0$.
\end{proof}

\begin{dfn}[\blue{\index{Minimal! polynomial}{Minimal polynomial}}]
Let $V$ be a $k$-space and $\tau\in End_k(V)$. The \ul{minimal polynomial} of $\tau$ is the polynomial $m_\tau\in k[x]$ of smallest degree such that $m_\tau(\tau)=0$. Equivalently (by the division algorithm for $k[x]$), $m_\tau\in k[x]$ is the PID generator of the kernel of the ring homomorphism {\small $h_\tau:k[x]\ra End_k(V),~f\mapsto f(\tau)$}, i.e., $\ker h_\tau=k[x]m_\tau$.
\end{dfn}

\begin{crl}
Let $V$ be a $k$-space and $\tau\in End_k(V)$. Then $m_\tau|p_\tau$, i.e., the minimal polynomial $m_\tau$ divides the characteristic polynomial $p_\tau$.
\end{crl}
\begin{proof}
Since $m_\tau|p_\tau$ $\iff$ $p_\tau\in k[x]m_\tau$, we only need to show that $p_\tau(\tau)=0$, i.e., $p_\tau\in\ker h_\tau=k[x]m_\tau$ $:=$ kernel of the ring homomorphism $h_\tau:k[x]\ra End_k(V),~f\mapsto f(\tau)$. Hence, the conclusion follows from Theorem \ref{CayHamThm}.
\end{proof}

\subsection{Rank of a linear transformation}
\begin{dfn}[\textcolor{blue}{
\index{Rank}{Rank},
\index{Nullity}{Nullity},
\index{Dimension formula (or rank-nullity theorem)}{Dimension formula (or rank-nullity theorem)}}]
Let $V,W$ be $k$-spaces (hence free $k$-modules). Then we know that a $k$-linear transformation $\tau:V\ra W$ is split, and so $V\cong \ker \tau\oplus \im \tau$. The \ul{nullity} of $\tau$, \ul{rank} of $\tau$, and the \ul{dimension formula} for $\tau$ are given (respectively) by
\[
\txt{nullity}_k(\tau):=\dim_k(\ker \tau),~~~~\rank_k(\tau):=\dim_k(\txt{im}\tau),~~~~\txt{and}~~~~\dim_k(\ker \tau)+\dim_k(\txt{im}\tau)=\dim_k V.
\]
\end{dfn}

\begin{dfn}[\blue{
\index{Associated linear transformation}{Associated linear transformation},
\index{Rank! of a matrix}{Rank of a matrix},
\index{Column rank}{Column rank},
\index{Row rank}{Row rank}}]
Let $k$ be a field, $A=[a_{ij}]\in Mat_{m\times n}(k)$, and consider the \ul{associated linear transformation} $A:k^n\ra k^m,(c_i)\mapsto \big(\sum_ja_{ij}c_j\big)\in\txt{ColSpace}(A):=k(a_{i1})+ k(a_{i2})+\cdots + k(a_{in})\subset k^m$. The \ul{rank} (or \ul{column rank}) of $A$ is
\bea
\rank(A):=\dim_k(\im A)=\dim_k~\txt{ColSpace}(A),\nn
\eea
i.e., the number of linearly independent columns of $A$. The \ul{row rank} of $A$ is the rank of the transpose $A^T$.
\end{dfn}

\begin{prp}
Let $V,W$ be finite-dimensional $k$-spaces with $n:=\dim V$ and $m:=\dim W$, $\tau\in Hom_k(V,W)$ a linear transformation, and $r\geq 1$ an integer. Then $\rank \tau=r$ $\iff$ there exists a basis $\al$ for $V$ and a basis $\beta$ for $W$ such that the matrix $[\tau]_{\al\beta}\in Mat_{m\times n}(k)$ of $\tau$ takes the form
{\footnotesize\[\left[
           \begin{array}{c|c}
             0 & 0 \\
             \hline
             0 & A \\
           \end{array}
         \right]=\left[
           \begin{array}{c|c}
             0_{(m-r)\times(n-r)} & 0_{(m-r)\times r} \\
             \hline
             0_{r\times(n-r)} & A_{r\times r} \\
           \end{array}
         \right],~~\txt{for some invertible matrix $A=A_{r\times r}\in GL_r(k)$}.
\]}
\end{prp}
\begin{proof}
Let $K:=\ker \tau$ and $I:=\im \tau$. Then we know $V=K\oplus K'$ (for a subspace $K'\subset V$) and $W=I'\oplus I$ (for a subspace $I'\subset W$). Let $\al_K$, $\al_{K'}$, $\beta_{I'}$, $\beta_I$ be bases for $K$, $K'$, $I'$, $I$ respectively. Then $\al:=\al_K\cup\al_{K'}$ is a basis for $V$ and $\beta:=\beta_{I'}\cup\beta_I$ a basis for $W$. By construction, $\tau|_K:K\ra 0\subset W$ and $\tau|_{K'}:K'\ra I\subset W$. Therefore, using the definition of the matrix representation of $\tau$ (in terms of the coordinate representation isomorphisms $\tau_\al:V\ra k^n$ and $\tau_\beta:W\ra k^m$)
\begin{align}
\tau_{\al\beta}&=\tau_\beta\circ\tau\circ\tau_\al^{-1}:k^n\sr{\tau_\al^{-1}}{\ral}(V,\al)\sr{\tau}{\ral}(W,\beta)\sr{\tau_\beta}{\ral} k^m,~~x\mapsto\beta^\ast\big(\tau(x^T\al)\big)=:[\tau]_{\al\beta}x,\nn
\end{align}
we get
{\footnotesize\bea
[\tau]_{\al\beta}=\left[
           \begin{array}{c|c}
             [\tau|_K]_{\al_K\beta_{I'}} & [\tau|_K]_{\al_K\beta_I} \\
             \hline
             [\tau|_{K'}]_{\al_{K'}\beta_{I'}} & [\tau|_{K'}]_{\al_{K'}\beta_I} \\
           \end{array}
         \right]
=\left[
           \begin{array}{c|c}
             0_{(m-r)\times(n-r)} & 0_{(m-r)\times r} \\
             \hline
             0_{r\times(n-r)} & A_{r\times r} \\
           \end{array}
         \right],~~~~A_{r\times r}:=[\tau|_{K'}]_{\al_{K'}\beta_I},\nn
\eea}where it is clear (since $\tau|_{K'}:K'\ra I$ is an isomorphism) that $A_{r\times r}$ is invertible.
\end{proof}

\begin{crl}
Let $k$ be a field, $A\in Mat_{m\times n}(k)$, and $r\geq 1$ an integer. Then $\rank A=r$ $\iff$ there exist invertible matrices $U\in GL_m(k)$ and $U'\in GL_n(k)$ such that
~{\footnotesize $UAU'=\left[
           \begin{array}{c|c}
             0_{(m-r)\times(n-r)} & 0_{(m-r)\times r} \\
             \hline
             0_{r\times(n-r)} & \mathds{1}_{r\times r} \\
           \end{array}
         \right]$}.
\end{crl}
\begin{proof}
Consider the associated $k$-linear transformation $A:k^n\ra k^m$. Choose bases $\al$ for $V:=k^n$ and $\beta$ for $W:=k^m$ as in the proof of the preceding proposition. Then the matrix representation of $A$ is
{\scriptsize\[
[A]_{\al\beta}=
\left[
           \begin{array}{c|c}
             0_{(m-r)\times(n-r)} & 0_{(m-r)\times r} \\
             \hline
             0_{r\times(n-r)} & A_{r\times r} \\
           \end{array}
         \right]
=\left[
           \begin{array}{c|c}
             \mathds{1}_{(m-r)\times(m-r)} & 0_{(m-r)\times r} \\
             \hline
             0_{r\times(m-r)} & A_{r\times r} \\
           \end{array}
         \right]\left[
           \begin{array}{c|c}
             0_{(m-r)\times(n-r)} & 0_{(m-r)\times r} \\
             \hline
             0_{r\times(n-r)} & \mathds{1}_{r\times r} \\
           \end{array}
         \right]\left[
           \begin{array}{c|c}
             \mathds{1}_{(n-r)\times(n-r)} & 0_{(n-r)\times r} \\
             \hline
             0_{r\times(n-r)} & \mathds{1}_{r\times r} \\
           \end{array}
         \right].\nn
\]}Hence, using the change of base matrices, we get the desired result. (Multiplication by an invertible matrix preserves rank, since a $k$-isomorphism preserves $k$-independence of vectors.)
\end{proof}

\begin{crl}[\blue{\index{Rank! balance}{Rank balance}}]
Let $k$ be a field. The rank of a matrix $A\in Mat_{m\times n}(k)$ equals that of its transpose $A^T\in Mat_{n\times m}(k)$. Equivalently, the row rank and column rank of a matrix are equal, i.e.,
\bea
\rank A=\rank A^T\leq\min(m,n).\nn
\eea
\end{crl}

\begin{crl}[\blue{Rank and conjugate transpose}]
Let $A=[a_{ij}]\in Mat_{m\times n}(\Complex)$, $\ol{A}:=[\ol{a}_{ij}]\in Mat_{m\times n}(\Complex)$ the \ul{complex conjugate matrix} of $A$, and $A^\ast:=\ol{A}^T\in Mat_{n\times m}(\Complex)$ the \ul{conjugate transpose} of $A$. Then
\bea
\rank A=\rank A^T=\rank\ol{A}=\rank A^\ast=\rank A^\ast A=\rank AA^\ast.\nn
\eea
In particular, for a real matrix $A=[a_{ij}]\in Mat_{m\times n}(\Real)\subset Mat_{m\times n}(\Complex)$, we have
\bea
\rank A=\rank A^T=\rank A^T A=\rank AA^T.\nn
\eea
\end{crl}
\begin{proof}
It is clear that $\rank A=\rank\ol{A}$ because the complex conjugation map $\Complex^k\ra\Complex^k,~x=(x_i)\mapsto\ol{x}:=(\ol{x}_i)$ preserves linear independence of vectors. Also, the preceding corollary implies $\rank A=\rank A^T$ and $\rank\ol{A}=\rank\ol{A}^T=\rank A^\ast$. That is, ~$\rank A=\rank A^T=\rank\ol{A}=\rank A^\ast$.

Next, observe that for any $x\in \Complex^n$, $Ax=0$ $\Ra$ $A^\ast Ax=0$, $\Ra$ $(Ax)^\ast Ax=x^\ast A^\ast Ax=0$, $\Ra$ $Ax=0$. That is, for any $x\in \Complex^n$, $Ax=0$ iff $A^\ast Ax=0$ (and similarly, for any $x\in \Complex^m$, $A^\ast x=0$ iff $AA^\ast x=0$), and so
\bea
\label{KerEqual}\ker A=\ker A^\ast A,~~~~~~~~\ker A^\ast=\ker AA^\ast.
\eea
Since $A:\Complex^n\ra \Complex^m$, $A^\ast A:\Complex^n\ra \Complex^n$, $A^\ast:\Complex^m\ra \Complex^n$, and $AA^\ast:\Complex^m\ra \Complex^m$, the dimension formula gives
{\small\begin{align}
&n=\dim_k(\ker A)+\dim_k(\im A)=\dim_k(\ker A^\ast A)+\dim_k(\im A^\ast A),~~\sr{(\ref{KerEqual})}{\Longrightarrow}~~\rank A=\rank A^\ast A,\nn\\
&m=\dim_k(\ker A^\ast)+\dim_k(\im A^\ast)=\dim_k(\ker AA^\ast)+\dim_k(\im AA^\ast),~~\sr{(\ref{KerEqual})}{\Longrightarrow}~~\rank A^\ast=\rank AA^\ast.\nn \hspace{1cm}\qedhere
\end{align}}
\end{proof}

{\flushleft\hrulefill}

\begin{exercise}
Based on the discussion of this chapter, consider writing a \emph{fully technical} essay (say in the form of a typical section of this chapter) on what is known in the mathematics literature as \index{Galois Theory}{``Galois Theory''}.
\end{exercise}

%% file: parts/AlgebraM/RepsTheoryIII.tex
\chapter{Representation Theory III: Semisimple Modules and Rings}\label{RepsTheoryIII}
In this chapter we continue discussing representations of a ring (i.e., modules over a ring). Once more, the end goal is to express a given module in terms of irreducible (i.e., simple) modules. We will consider the possibility of directly expressing a module as a sum of its simple submodules. The main goal in this chapter is to characterize semisimple modules (modules that can be expressed as a sum of their simple submodules).

\section{Noetherian, Artinian, and Finitely Generated (FG) Modules}
\subsection{Recall of concepts and preliminary remarks (almost as in previous chapter)}
\begin{dfn}[\blue{\small
\index{Minimal! (or simple) ideal}{Minimal (or simple) ideal},
\index{Maximal! ideal}{Maximal ideal},
\index{Prime! ideal}{Prime ideal}}]
Let $R$ be a ring. A nonzero proper ideal $0\neq I\vartriangleleft R\neq I$ is a \ul{minimal ideal} or \ul{simple ideal} if it cannot properly contain a nonzero ideal in the sense that if $I\supset J\vartriangleleft R$, then $J=0$ or $J=I$ (i.e., $RaR=I$ for all $0\neq a\in I$). A proper ideal $I\vartriangleleft R\neq I$ is a \ul{maximal ideal} if it cannot be properly contained by a proper ideal in the sense that if $I\subset J\vartriangleleft R$, then $J=I$ or $J=R$ (i.e., $RaR+I=R$ for all $a\not\in I$). A proper ideal $P\vartriangleleft R\neq P$ is a \ul{prime ideal} if for any $a,b\in R$,~ $aRb\subset P$ $\Ra$ $a\in P$ or $b\in P$ (or equivalently, for any $I,J\lhd R$,~ $IJ\subset P$ $\Ra$ $I\subset P$ or $J\subset P$).
\end{dfn}
In particular, if $R$ is commutative then, a proper ideal $P\vartriangleleft R\neq P$ is prime iff for any $a,b\in R$,~ $ab\in P$ $\Ra$ $a\in P$ or $b\in P$.

\begin{dfn}[\blue{\small
\index{Minimal! (or simple) submodule}{Minimal (or simple) submodule},
\index{Maximal! submodule}{Maximal submodule},
\index{Minimal! left ideal}{Minimal left ideal},
\index{Maximal! left ideal}{Maximal left ideal},
\index{Minimal! right ideal}{Minimal right ideal},
\index{Maximal! right ideal}{Maximal right ideal},
\index{Minimal! ideal}{Minimal ideal},
\index{Maximal! ideal}{Maximal ideal}}]
Let $M$ be an $R$-module. A nonzero proper submodule $0\neq{}_RN\subset{}_RM\neq N$ is a \ul{minimal submodule} or \ul{simple submodule} if it cannot properly contain a nonzero submodule in the sense that if~ $_RN\supset{}_RX\subset M$, then $X=0$ or $X=N$ (i.e., $Rn=N$ for all $0\neq n\in N$).
A proper submodule $_RN\subset{}_RM\neq N$ is a \ul{maximal submodule} if it cannot be properly contained by a proper submodule in the sense that if~ $_RN\subset{}_RX\subset M$, then $X=N$ or $X=M$ (i.e., $Rm+N=M$ for all $m\not\in N$).

A minimal (resp. maximal) left ideal ${}_RI\subset R$ is a minimal (resp. maximal) submodule of $_RR$. Similarly, a minimal (resp. maximal) right ideal $I_R\subset R$ is a minimal (resp. maximal) submodule of $R_R$. Accordingly, a minimal (resp. maximal) ideal $I={}_RI_R\subset R$ is a minimal (resp. maximal) submodule of ${}_RR_R$.
\end{dfn}

\begin{dfn}[\blue{\small
\index{Simple! ring}{Simple ring}, 
\index{Simple! module}{Simple module}}]
A ring $R$ is a \ul{simple ring} if $0$ and $R$ are the only ideals of $R$ (i.e., $RaR=R$ for all $a\in R\backslash 0$). A nonzero module $_RS\neq 0$ is a \ul{simple module} if $0$ and $S$ are the only submodules of $S$ (i.e., $Rs=S$ for all $s\in S\backslash 0$).
\end{dfn}
By the correspondence theorem for rings, an ideal $I\vartriangleleft R$ is maximal iff $R/I$ is a simple ring. Similarly, by the correspondence theorem for modules, a submodule $N\subset M$ is maximal iff $M/N$ is a simple module. In particular, a left ideal $_RI\subset R$ is maximal iff $R/I$ is a simple $R$-module.

\begin{dfn}[\blue{\small
\index{Noetherian module}{Noetherian module},
\index{Left! noetherian ring}{Left noetherian ring},
\index{Right! noetherian ring}{Right noetherian ring},
\index{Noetherian ring}{Noetherian ring}}]~
\begin{enumerate}
\item A module $_RM$ is a \ul{noetherian module} if given any ascending chain of submodules $M_1\subset M_2\subset M_3\subset\cdots$, there exists $n\geq 1$ such that $M_n=M_{n+1}=M_{n+2}=\cdots$. That is, every ascending chain of submodules terminates or stabilizes.
\item A ring $R$ is a \ul{left noetherian ring} if $_RR$ (resp. \ul{right noetherian ring} if $R_R$) is noetherian.
\item A ring is a \ul{noetherian ring} if it is both left noetherian and right noetherian.
\end{enumerate}
\end{dfn}

\begin{dfn}[\blue{\small
\index{Artinian module}{Artinian module},
\index{Left! artinian ring}{Left artinian ring},
\index{Right! artinian ring}{Right artinian ring},
\index{Artinian ring}{Artinian ring}}]~
\begin{enumerate}
\item A module $_RM$ is an \ul{artinian module} if given any descending chain of submodules $M_1\supset M_2\supset M_3\supset\cdots$, there exists $n\geq 1$ such that $M_n=M_{n+1}=M_{n+2}=\cdots$. That is, every descending chain of submodules terminates or stabilizes.
\item A ring $R$ is a \ul{left artinian ring} if ${}_RR$ (resp. \ul{right artinian ring} if the module $R_R$) is artinian.
\item A ring is an \ul{artinian ring} if it is both left artinian and right aritian.
\end{enumerate}
\end{dfn}

Recall that every simple module is clearly both artinian and noetherian.

\begin{dfn}[\blue{\index{Finitely generated! module}{Finitely generated (FG) module}}]
A module $_RM$ is a \ul{finitely generated module} (or \ul{FG module}) if $M=Rm_1+\cdots+Rm_n$ for some $m_1,...,m_n\in M$ and $n\geq 1$. That is, $_RM$ is finitely generated iff there exists an exact sequence of $R$-modules ~$0\ra K\ra R^n\ra M\ra 0$,~ where $R^n:=R\oplus R^{n-1}$.
\end{dfn}

\subsection{Results on noetherian and artinian modules (almost as in previous chapter)}
\begin{thm}[\blue{Noetherian and Artinian characterization I}]\label{NoeArtChrI}
A module $_RM$ is noetherian (resp. artinian) iff every nonempty collection of submodules has a maximal (resp. minimal) element with respect to inclusion.
\end{thm}
\begin{proof}
Assume $M$ is noetherian (resp. artinian) and let $\S$ be a nonempty collection of submodules of $M$. Then every nonempty chain, with respect to inclusion, in $\S$ has an upper (resp. a lower) bound since $M$ is noetherian (resp. artinian), and so $\S$ has a maximal (resp. minimal) element by Zorn's lemma.

Conversely, if every nonempty collection of submodules has a maximal (resp. minimal) element, then in particular, any ascending (resp. descending) chain of submodules terminates, and so $M$ is noetherian (resp. artinian).
\end{proof}

\begin{lmm}
Let $M$ be an $R$-module and ${}_RA,{}_RB,{}_RN\subset{}_RM$ submodules. If $A\subset B$, then
\bea
A\cap N=B\cap N~~\txt{and}~~A+N=B+N~~~~\Ra~~~~A=B.\nn
\eea
\end{lmm}
\begin{proof}
Let $C:=A\cap N=B\cap N$. Then by the isomorphism theorems, we have
\bea
\textstyle {A\over C}={A\over A\cap N}\cong{A+N\over N}={B+N\over N}\cong {B\over B\cap N}={B\over C},\nn
\eea
where by construction the isomorphism is explicitly given by $f:A/C\ra B/C,~a+C\ra a+C$ (which is clearly well-defined and injective). Since $f$ is surjective, $f(A/C)=B/C$ implies $A/C=B/C$, and so $A=B$.
\end{proof}

\begin{thm}[\textcolor{blue}{Noetherian and Artinian characterization II}]\label{NoeArtChrII}
Let $N\subset M$ be $R$-modules. $M$ is noetherian (resp. artinian) $\iff$ $N$ and $M/N$ are both noetherian (resp. artinian).
\end{thm}
\begin{proof}
The artinian case will follow by symmetry of the following argument for the noetherian case (simply by replacing ``noetherian'' with ``artinian'', ``increasing'' with ``decreasing'', and ``$\subset$'' with ``$\supset$'').
{\flushleft ($\Ra$)}
Assume $M$ is \magenta{noetherian}. Then every \magenta{increasing} sequence of submodules of $N$ is an \magenta{increasing} sequence of submodules of $M$, and so terminates. Hence $N$ is \magenta{noetherian}. Similarly, every \magenta{increasing} sequence of submodules $M_1/N\magenta{\subset} M_2/N\magenta{\subset}\cdots$ of $M/N$ gives (by the correspondence theorem) an \magenta{increasing} sequence of submodules $M_1\magenta{\subset} M_2\magenta{\subset}\cdots$ of $M$, and so terminates. Hence $M/N$ is \magenta{noetherian}.
{\flushleft ($\La$)}
Assume $N$ and $M/N$ are both \magenta{noetherian}. Let $A_1\magenta{\subset} A_2\magenta{\subset}\cdots$ be a chain of submodules of $M$. Then the following two sequences of submodules of $N$ and $M/N$ (respectively) each stabilize:
\bea
\label{NoeArtChIIIEq}&&\textstyle (a)~~A_1\cap N\magenta{\subset} A_2\cap N\magenta{\subset}\cdots,~~~~~~(b)~~{A_1+N\over N}\magenta{\subset} {A_2+N\over N}\magenta{\subset}\cdots
\eea
Thus for some $k$, ${A_k+N\over N}={A_{k+1}+N\over N}=\cdots$ (i.e., $A_k+N=A_{k+1}+N=\cdots$) as the stabilization point for (\ref{NoeArtChIIIEq})(b). Without loss of generality assume the stabilization point for (\ref{NoeArtChIIIEq})(a) is also $A_k\cap N=A_{k+1}\cap N=\cdots$ (otherwise, take the higher of the two). That is, $A_k+N=A_{k+1}+N=\cdots$ and $A_k\cap N=A_{k+1}\cap N=\cdots$, and since $A_k\magenta{\subset} A_{k+1}\magenta{\subset}\cdots$, the preceding lemma implies $A_k=A_{k+1}=\cdots$, i.e., the chain $A_1\magenta{\subset} A_2\magenta{\subset}\cdots$ stabilizes. Hence $M$ is \magenta{noetherian}.
\end{proof}

\begin{crl}\label{NthArtCh2Crl1}
Let {\small $0\ra A\ra B\ra C\ra 0$} be an exact sequence of $R$-modules. Then $B$ is noetherian (resp. artinian) iff $A,C$ are both noetherian (resp. artinian).
\end{crl}

\begin{crl}\label{NthArtCh2Crl2}
Let $M,M'$ be R-modules. Then $M\oplus M'$ is noetherian (resp. artinian) iff $M,M'$ are both noetherian (resp. artinian).
\end{crl}
\begin{proof}
Observe that we have an exact sequence $0\ra M\ra M\oplus M'\ra M'\ra 0$.
\end{proof}

\begin{crl}\label{NthArtCh2Crl3}
Let {\small $\{M_i\}_{i\in I}$} be $R$-modules. (i) If $\bigoplus_{i\in I} M_i$ is noetherian (resp. artinian), so is $M_i$ for every $i\in I$. (ii) if $I=\{1,...,n\}$ (i.e., a finite set) then {\small $M_1\oplus \cdots\oplus M_n$} is noetherian (resp. artinian) iff each $M_i$ is noetherian (resp. artinian).
\end{crl}
\begin{proof}
(i) For each $j\in I$, we have the exact sequence $0\ra M_j\ra\bigoplus_{i\in I}M_i\ra\bigoplus_{i\neq j}M_i\ra 0$.
(ii) This follows by induction on $n$, due to associativity of the direct sum, i.e.,
\bit
\item[]\hspace{2cm} $M_1\oplus M_2\oplus\cdots\oplus M_n=(M_1\oplus M_2\oplus\cdots\oplus M_{n-1})\oplus M_n$.\qedhere
\eit
\end{proof}

\begin{lmm} Let $M$ be an $R$-module and $f\in End_R(M):=Hom_R(M,M)$. Then the following hold:
\bit[leftmargin=0.8cm]
\item[(1)] If $M$ is Noetherian and $f$ is surjective, then $f$ is an isomorphism.
\item[(2)] If $M$ is Artinian and $f$ is injective, then $f$ is an isomorphism.
\eit
\end{lmm}
\begin{proof}
Let $f^{(n)}:=f\circ f^{(n-1)}=f^{(n-1)}\circ f$ be the $n$-fold composition of $f$, where $f^{(0)}:=id_M$.
{\flushleft(1)} Since $f$ is surjective, so is $f^{(n)}$. Also, ~{\small $\ker f^{(n+1)}=\ker (f\circ f^{(n)})=(f^{(n)})^{-1}(\ker f)\supset\ker f^{(n)}$},~ and so
\bea
&&\ker f^{(n)}\subset\ker f^{(n+1)},~~\sr{(s1)}{\Ra}~~\ker f^{(N)}=\ker f^{(N+1)}=\ker f^{(N+2)}=\cdots~~\txt{for some}~~N,\nn\\
&&~~\Ra~~f\left(\ker f^{(N)}\right)=f\left(\ker f^{(N+1)}\right),~~\sr{(s2)}{\Ra}~~\ker f^{(N-1)}=\ker f^{(N)},~~\cdots,\nn\\
&&~~\sr{\txt{induction}}{\Ra}~~0=\ker f^{(0)}=\ker f^{(1)}=\ker f,\nn
\eea
where step (s1) is the noetherian property and step (s2) the surjectivity of $f$. (\blue{footnote}\footnote{Because $f$ is surjective, $f\big(f^{-1}(m)\big)=m$ for all $m\in M$.}).
{\flushleft (2)} Since $f$ is injective, so is $f^{(n)}$. Also, ~{\small $\im f^{(n+1)}=\im (f^{(n)}\circ f)=f^{(n)}(\im f)\subset\im f^{(n)}$},~ and so
\bea
&&\txt{im}~f^{(n)}\supset\txt{im}~f^{(n+1)},~~\sr{(s1)}{\Ra}~~\txt{im}~f^{(N)}=\txt{im}~f^{(N+1)}=\txt{im}~f^{(N+2)}=\cdots~~\txt{for some}~~N,\nn\\
&&~~\Ra~~f^{-1}\left(\txt{im}~f^{(N)}\right)=f^{-1}\left(\txt{im}~f^{(N+1)}\right),~~\sr{(s2)}{\Ra}~~\txt{im}~f^{(N-1)}=\txt{im}~f^{(N)},~~\cdots,\nn\\
&&~~\sr{\txt{induction}}{\Ra}~~M=\im f^{(0)}=\im f^{(1)}=\txt{im}~f,\nn
\eea
where step (s1) is the artinian property and step (s2) the injectivity of $f$. (\blue{footnote}\footnote{Because $f$ is injective, $f^{-1}\big(f(m)\big)=m$ for all $m\in M$.})
\end{proof}

\begin{thm}[\textcolor{blue}{A finitely generated module over a left noetherian (resp. left artinian) ring is noetherian (resp. artinian)}]\label{NoeArtFGthm}
Let $M$ be an $R$-module. If $R$ is left noetherian (resp. left artinian) and $M$ is finitely generated, then $M$ is noetherian (resp. artinian). (\blue{footnote}\footnote{Here, the ring need only be left-noetherian (resp. left-artinian), since we assume the modules are left modules.})
 \end{thm}
\begin{proof}
Recall that $M$ is FG iff we have an exact sequence of $R$-modules ~$0\ra K\ra R^n\ra M\ra 0$.
\end{proof}

\section{Semisimple Modules}
\subsection{Decomposability of modules}
\begin{dfn}[\blue{\small
\index{Decomposable module}{Decomposable module},
\index{Indecomposable module}{Indecomposable module},
\index{Completely decomposable (Completely reducible) module}{Completely decomposable (Completely reducible) module}}]
Let $R$ be a ring. An $R$-module $M$ is \ul{decomposable} if it is an internal direct sum, or simply ``a direct sum'', of two nonzero submodules (i.e., $M=N_1\oplus N_2$ for nonzero submodules $0\neq N_1,N_2\subset M$), otherwise $M$ is \ul{indecomposable} (i.e., for any nonzero submodules $0\neq N_1,N_2\subset M$ we have $M\neq N_1\oplus N_2$, or equivalently, for any submodules $N_1,N_2\subset M$, if $M=N_1\oplus N_2$ then $N_1=0$ or $N_2=0$).

An $R$-module $M$ is \ul{completely decomposable} (or \ul{completely reducible}) if every submodule of $M$ is a direct summand of $M$ (i.e., for every submodule $_RN\subset M$, there exists a submodule $_RN'\subset M$ such that $M=N\oplus N'$), equivalently, if every SES of $R$-modules $0\ra A\ra M\ra B\ra 0$ is split.
\end{dfn}

\begin{lmm}\label{IndecSumLmm}
If an $R$-module $M$ is noetherian or artinian, then $M=N_1\oplus\cdots\oplus N_r$ for indecomposable submodules $N_1,...,N_r\subset M$.
\end{lmm}
\begin{proof}
{\flushleft $\bullet$} \ul{Assume $M$ is noetherian} (and that $M$ is decomposable, otherwise there is nothing to prove). Let $M=N_1\oplus M_1$, where $M_1$ is a maximal direct summand of $M$ wrt $N_1$ (so $N_1$ is indecomposable, otherwise $M_1$ would not be a maximal direct summand). This can be repeated since $M_1\cong {M\over N_1}$ is again noetherian (Theorem \ref{NoeArtChrII}). At the $i$th step, we have $M_i=N_{i+1}\oplus M_{i+1}$, where $M_{i+1}$ is a maximal direct summand of $M_i$ wrt $N_{i+1}$ (so $N_{i+1}$ is indecomposable), and so we get
\bea
M=N_1\oplus M_1=N_1\oplus (N_2\oplus M_2)=N_1\oplus N_2\oplus\cdots.\nn
\eea
By the noetherian property, the sequence $N_1\subsetneq N_1\oplus N_2\subsetneq N_1\oplus N_2\oplus N_3\subsetneq\cdots$ terminates at some integer $r-1$ with $N_{r-1}$ indecomposable. Set $N_r:=M_{r-1}$ in $M_r=N_{r-1}\oplus M_{r-1}$.
{\flushleft $\bullet$} \ul{Assume $M$ is artinian} (and that $M$ is decomposable, otherwise there is nothing to prove). Let $M=N_1\oplus M_1$, where $N_1$ is a minimal direct summand of $M$ wrt $M_1$ (so $N_1$ is indecomposable). This can be repeated since $M_1\cong {M\over N_1}$ is again artinian (Theorem \ref{NoeArtChrII}). At the $i$th step, we have $M_i=N_{i+1}\oplus M_{i+1}$, where $N_{i+1}$ is a minimal direct summand of $M_i$ wrt $M_{i+1}$ (so $N_{i+1}$ is indecomposable), and so we get
\bea
M=N_1\oplus M_1=N_1\oplus (N_2\oplus M_2)=N_1\oplus N_2\oplus\cdots.\nn
\eea
By the artinian property the sequence $M_1\supsetneq M_2\supsetneq M_3\supsetneq\cdots $ terminates at some integer $r-1$ such that $M_{r-1}\neq0$ (which must be indecomposable). Set $N_r:=M_{r-1}$ in $M_r=N_{r-1}\oplus M_{r-1}$.
\end{proof}

\subsection{Review of the product, coproduct, and internal direct sum for modules}
\begin{dfn}[\blue{Recall:
\index{Direct product (Product)}{Direct product (Product)},
\index{Direct sum (Coproduct)}{Direct sum (Coproduct)}}]
 Let $\F:=\{M_i:~i\in I\}$ be a family of R-modules. The \ul{direct product} of $\F$ is the set
{\small
\bea
\textstyle\prod_{i\in I}M_i:=\big\{(m_i)=(m_i)_{i\in I}:~m_i\in M_i\big\}=\left\{\txt{maps}~f:I\ra\bigcup_{j\in I} M_j,~f(i)\in M_i\right\},\nn
\eea}
as an R-module with addition and scalar multiplication given componentwise by
\bea
(m_i)+(n_i):=(m_i+n_i),~~~~~r\cdot(m_i):=(rm_i).\nn
\eea
The \ul{direct sum} $\coprod_{i\in I}M_i$ of $\F$ is the submodule of $\prod_{i\in I}M_i$ given by (\blue{footnote}\footnote{This definition of the coproduct as a subobject of the product works only because $R$-mod is an abelian category (for which finite products are equivalent to finite coproducts). The modules $\prod_{i\in I}M_i$ and $\coprod_{i\in I}M_i$ are equal when $I$ is finite.})
{\small
\bea
\textstyle\coprod_{i\in I}M_i:=\big\{(m_i)=(m_i)_{i\in I}:~m_i\in M_i,~m_i=0~~a.e.f.\big\}=\left\{\txt{finite maps}~f:I\ra\bigcup_{j\in I} M_j,~f(i)\in M_i\right\}.\nn
\eea}
\end{dfn}

\begin{dfn}[\blue{\index{Internal direct sum}{Recall: Internal direct sum}}]
Let $M$ be an $R$-module and {\small $\F:=\{M_i:~i\in I\}$} a family of submodules of $M$. If
  \bit
  \item[(a)] $M=\sum_{i\in I}M_i~:=~\left\{\sum_{i\in I}m_i:~m_i\in M_i,~m_i=0~~a.e.f.\right\}$, and
  \item[(b)] every $m\in M$ can be written as a \emph{unique sum} ~$m=\sum_{i\in I}m_i$, $m_i\in M_i$, $m_i=0$ a.e.f. (\blue{footnote}\footnote{That is, if ~$\sum_{i\in I}m_i=0$, where $m_i\in M_i$ and $m_i=0$ a.e.f., then $m_i=0$ for all $i\in I$.}),
  \eit
then we write $M=\bigoplus_{i\in I}M_i$, and call it an \ul{internal direct sum} of $\F$. It is clear that $M=\bigoplus M_i$ $\iff$ the map ~{\small $\coprod_{i\in I}M_i\ra M,~(m_i)\mapsto \sum_i m_i$}~ is an isomorphism of $R$-modules.
\end{dfn}

Because of the above isomorphism, the module direct sum and module internal direct sum do not need to be distinguished. In particular, we often write the expression $\bigoplus_{i\in I}M_i$ for all direct sums of modules, which is to be understood as $\coprod_{i\in I}M_i$ in the appropriate context.

\begin{lmm}
Let $M$ be an $R$-module and $\F:=\{M_i:i\in I\}$ a family of submodules of $M$. Then $M=\bigoplus M_i$ $\iff$ (a) $M=\sum_{i\in I} M_i$ and (b) $M_i\cap \sum_{j\neq i}M_j=0$ for each $i\in I$.
\end{lmm}
\begin{proof}
Recall from its definition that $M=\bigoplus M_i$ $\iff$ the map {\small $\coprod_{i\in I}M_i\ra M,~(m_i)\mapsto \sum_i m_i$} is an $R$-isomorphism, which (as we can check directly) is the case $\iff$ (a) and (b) both hold, where (a) gives surjectivity and (b) gives injectivity. Surjectivity is clear. For injectivity, observe that if $(m_i),(m'_i)\in\coprod M_i$ satisfy $\sum_i m_i=\sum_i m'_i$, then for each $i\in I$, we have
\[
\textstyle m_i-m'_i=-\sum_{j\neq i}(m_j-m'_j)~\in~M_i\cap\sum_{j\neq i}M_j. \qedhere
\]
\end{proof}

\begin{rmk}\label{ModSmSpCrRmk}
Given a family of R-modules $\F=\{M_i:~i\in I\}$, Let $N:=\coprod M_i$ and $N_i:=\big\{(m_j)_{j\in I}\in N:m_j=0~\txt{for all}~j\neq i\big\}=\{n_i:=(m_i\delta_{ij})_{j\in I}\in N~|~m_i\in N_i\}$ for each $i\in I$. Then $N=\sum N_i$ is the internal direct sum $\bigoplus N_i$ of the family of submodules $\G:=\{N_i:~i\in I\}$ of $N$. This follows from observing that the map ~$\coprod N_i\ra N,~(n_i)\mapsto\sum_i n_i$~ is an isomorphism of $R$-modules.
\end{rmk}

\begin{lmm}[\blue{Recall}]
Let $\{M_i\}_{i\in I}$ be $R$-modules and $N_i\subset M_i$ $R$-submodules. Then ${\bigoplus M_i\over \bigoplus N_i}\cong\bigoplus{M_i\over N_i}$ via the $R$-homomorphism $f:\bigoplus M_i\ra \bigoplus{M_i\over N_i},~(m_i)\mapsto(m_i+N_i)$ and the 1st isomorphism theorem.
\end{lmm}

\begin{lmm}\label{ComDecSimSub}
Let $M$ be an $R$-module and $N\subset M$ a submodule. If $M$ is completely decomposable (CD), then (i) $N$ is CD and (ii) $M$ has a simple submodule.
\end{lmm}
\begin{proof}
Assume $M$ is CD. (i) Pick a submodule $A\subset N$. Then $M=N\oplus N'=A\oplus A'$ for submodules $N',A'\subset M$. Let $A'_N:=A'\cap N$. For any $n=a+a'\in N$ (where $a\in A$, $a'\in A'$) we have $a'=n-a\in N\cap A'=A'_N$, i.e., $N=A+A'_N$. Since $A\cap A'_N=A\cap A'\cap N=\{0\}$, we get $N=A\oplus A'_N$.
{\flushleft (ii)} Let $0\neq m\in M$. Then the submodule $Rm\subset M$ is CD, where we know $Rm\cong{R\over I}$ with $I:=Ann_R(m)$. Let $I\subset{}_RJ\subset R$ for a maximal left ideal $J$. Then $Rm=Jm\oplus S$ (by CD), where ${}_RS\subset Rx$ satisfies
\[
\textstyle S\cong{Jm\oplus S\over Jm}={Rm\over Jm}\sr{f}{\cong}{R/I\over J/I}\cong{R\over J}~~\txt{(a simple $R$-module)}.
\]
The isomorphism $f$ above can be obtained using the map $g:{R\over I}\ra {Rm\over Jm},~r+I\mapsto rm+Jm$ with kernel
\begin{align}
&\textstyle\ker g=\{r+I\in R/I:rm\in Jm\}=\{r+I\in R/I:rm=jm,~j\in J\}\nn\\
&\textstyle~~~~=\{r+I\in R/I:r-j\in I,~\txt{i.e.,}~r\in j+I\subset J,~j\in J\}={J+I\over I}={J\over I}.\nn \qedhere
\end{align}
\end{proof}

\subsection{Semisimplicity of modules}
\begin{dfn}[\blue{
\index{Socle of a module}{Socle of a module},
\index{Semisimple module}{Semisimple module},
\index{Left! semisimple ring}{Left semisimple ring},
\index{Right! semisimple ring}{Right semisimple ring},
\index{Semisimple ring}{Semisimple ring}}]
Let $R$ be a ring and $M$ an $R$-module. The \ul{socle} of $M$ is the sum of all simple submodules of $M$, i.e., $Soc(M):=\sum\{S:\txt{for simple}~{}_RS\subset M\}$. The module $M$ is a \ul{semisimple module} if $Soc(M)=M$. The ring $R$ is a \ul{left semisimple ring} (resp. a \ul{right semisimple ring}) if the module $_RR$ (resp. $R_R$) is semisimple. The ring $R$ is a \ul{semisimple ring} if it is both left semisimple and right semisimple.
\end{dfn}

We will see in Corollary \ref{ArtWeddCrl1} that a ring is left semisimple iff it is right semisimple iff  it is semisimple. The following theorem establishes that a module is semisimple iff it is completely decomposable.

\begin{thm}[\blue{\index{Module! semi-simplicity criterion}{Module semi-simplicity criterion}}]\label{ModSmSpCrThm}
Given any module $_RM$, the following are equivalent.
\begin{enumerate}
\item[(a)] $M$ is a direct sum of simple modules, i.e.,~ $M=\bigoplus_{i\in I}S_i$,~ where $_RS_i$ are simple modules.
\item[(b)] $M$ is a sum of simple submodules, i.e.,~ $M=\sum_{i\in I}S_i$,~ where $_RS_i\subset M$ are simple submodules.
\item[(c)] $M$ is completely decomposable (i.e., every SES of $R$-modules $0\ra A\ra M\ra B\ra 0$ is split).
\end{enumerate}
\end{thm}
\begin{proof}
{\flushleft \ul{(a)$\Ra$(b)}}: Assume $M=\bigoplus_{i\in I}S_i$, where each $_RS_i$ is a simple module. Then as in Remark \ref{ModSmSpCrRmk}, we have $M=\sum_{i\in I}M_i$, where {\small $M_i:=\big\{(s_j)_{j\in I}\in M:s_j=0~\txt{for all}~j\neq i\big\}$} are simple submodules of $M$ because for each $i\in I$, $M_i\cong S_i$ via the $R$-isomorphism
\bea
f_i:M_i\ra S_i,~(s_j)_{j\in I}\mapsto s_i.\nn
\eea

{\flushleft \ul{(b)$\Ra$(c)}}: Assume $M=\sum_{i\in I}S_i$,~ where $_RS_i\subset M$ are simple submodules. Pick any $_RN\subset M$. Let the collection of submodules $\C:=\{{}_RX\subset M:~N\cap X=0\}\ni 0$ be ordered by inclusion. Then any chain $\{X_\ld\}_{\ld\in\Ld}\subset\C$ has an upper bound $\bigcup_{\ld\in\Ld}X_\ld$, since $\bigcup_{\ld\in\Ld}X_\ld\subset M$ is a submodule satisfying
\bea
\textstyle N\cap\left(\bigcup_\ld X_\ld\right)=\bigcup_\ld(N\cap X_\ld)=\bigcup_\ld 0=0,~~\Ra~~\bigcup_\ld X_\ld\in\C,\nn
\eea
and so by Zorn's lemma, $\C$ contains a maximal element $Y$. Because $M=\sum_{i\in I}S_i$ is a sum of simple modules, we have $N+Y=M$. Otherwise, if $N+Y\neq M$, then $S_j\not\subset N+Y$ for some $j\in I$, and so
{\small\begin{align}
&S_j\cap(N+Y)=0,~~\Ra~~N\cap(Y+S_j)=\Big\{n=y+s:~s\in S_j\cap(N-Y)=S_j\cap(N+Y)=0\Big\}=N\cap Y=0,\nn\\
&~~\Ra~~Y\subsetneq Y+S_j\in \C~~(\txt{a contradiction}).\nn
\end{align}}
{\flushleft \ul{(c)$\Ra$(a)}}: Assume $M$ is completely decomposable. Let $S_I:=\{S_i:i\in I\}$ be the family of simple submodules of $M$. By Lemma \ref{ComDecSimSub}, $S_I$ is nonempty. Consider the collection of ``internal direct-summing'' index subsets
\bea
\textstyle\J:=\left\{J\subset I~|~\sum S_I=\sum S_J=\bigoplus S_J\right\},~~\txt{where}~~\sum S_J:=\sum_{j\in J}S_j~~\txt{and}~~\bigoplus S_J:=\bigoplus_{j\in J}S_j,\nn
\eea
ordered by inclusion. For any $i,i'\in I$, if $S_i\cap S_{i'}\neq0$ then $S_i=S_{i'}$, i.e., either $S_i\cap S_{i'}=0$ or $S_i=S_{i'}$. If we let ``$i\sim i'$ if $S_i=S_{i'}$'', we get an equivalence relation on $I$, and so $\J$ is nonempty (because with an injective indexing $\{[j]\}_{j\in J\subset I}$ of the $\sim$-equivalence classes in $I$, the axiom of choice gives a map $f:J\ra I=\bigcup_{j\in J}[j]$, $j\mapsto f(j)\in[j]$). Any chain $\{J_\ld\}_{\ld\in \Ld}\subset\J$ has an upper bound $\cup_\ld J_\ld$, since for any $\ld\in\Ld$,
\begin{align}
&\textstyle\sum S_I=\sum S_{J_\ld}=\bigoplus S_{J_\ld}\subset\sum S_{\cup_{\ld'} J_{\ld'}}=\sum_{\ld'}\sum S_{J_{\ld}'}\sr{(a)}{=}\bigcup_{\ld'}\sum S_{J_{\ld}'}=\bigcup_{\ld'}\sum S_I=\sum S_I\nn\\
&\textstyle~~\Ra~~\sum S_I=\sum S_{\cup_\ld J_\ld}\sr{(b)}{=}\bigoplus S_{\cup_\ld J_\ld},\nn
\end{align}
where step (a) holds because $\{\sum S_{J_\ld}\}_{\ld\in\Ld}\subset\P(M)$ is a chain wrt inclusion (since $\{J_\ld\}_{\ld\in\Ld}\subset\J$ is a chain wrt inclusion), and step (b) holds because by step (a) every $\sum_{j\in\cup_\ld J_\ld} s_j\in\sum S_{\cup_\ld J_\ld}=\bigcup_{\ld}\sum S_{J_\ld}$ lies in some $\sum S_{J_\ld}$, and so $\sum_{j\in\cup_\ld J_\ld} s_j=0$ $\iff$ $s_j=0$ for all $j\in\cup_\ld J_\ld$. So, by Zorn's lemma, $\J$ contains a maximal element $J^\ast$.

Because $M$ is completely decomposable, $M=\sum S_{J^\ast}$. Otherwise, if $M\neq\sum S_{J^\ast}$, then by complete decomposability (CD), $M=N\oplus N'$, where $N:=\sum S_{J^\ast}$ and $0\neq{}_RN'\subset M$. Let $x\in N'\backslash 0$. Then $M\sr{\txt{CD}}{=}(N\oplus Rx)\oplus Y$ for some $_RY\subset M$. Since $Rx$ is finitely generated, $Rx$ has a maximal submodule $_R\mu\subset Rx$, and so we get a simple module ~${Rx\over {\mu}}\cong{(Rx\oplus(N\oplus Y))/(N\oplus Y)\over ({\mu}\oplus(N\oplus Y))/(N\oplus Y)}\cong{N\oplus Rx\oplus Y\over N\oplus {\mu}\oplus Y}={M\over N\oplus {\mu}\oplus Y}$. Therefore,
\[
\textstyle M\sr{\txt{CD}}{=}(N\oplus {\mu}\oplus Y)\oplus S_{i_x}=(N\oplus S_{i_x})\oplus({\mu}\oplus Y),~~\txt{for some}~~i_x\in I~~\txt{such that}~~~S_{i_x}\cong {M\over N\oplus {\mu}\oplus Y}\cong{Rx\over\mu}.
\]
This contradicts the maximality of $J^\ast$ in $\J$, because with $J_x:=J^\ast\cup\{i_x\}\subset I$, we have
\[
\textstyle \sum S_I=\sum S_I+S_{i_x}=\sum S_{J^\ast}+S_{i_x}=\left(\bigoplus S_{J^\ast}\right)+ S_{i_x}=\bigoplus S_{J_x}=\sum S_{J_x},
\]
which implies $J^\ast\subsetneq J^\ast\cup\{i_x\}=J_x\in \J$.
\end{proof}

\begin{lmm}[\blue{Simple image}]
Let $\phi:S\ra M$ be an R-homomorphism. If $_RS$ is simple, then either $\phi(S)=0$ or $\phi(S)\cong S$.
\end{lmm}
\begin{proof}
By the first isomorphism theorem, $\phi(S)\cong{S\over\ker\phi}$, where either $\ker\phi=S$ or $\ker\phi=0$.
\end{proof}

\begin{crl}[\blue{\index{Schur's lemma}{Schur's lemma}}]\label{SchurLmm}
If $_RS$ is simple, then $End_R(S):=Hom_R(S,S)$ is a division ring. (Recall that a ring is called divisible if every nonzero element is a unit.)
\end{crl}

\begin{crl}[\blue{Semisimplicity image}]
Let $\phi:M\ra N$ be an R-homomorphism. If $_RM$ is semisimple, then either $\phi(M)=0$ or $\phi(M)$ is semisimple.
\end{crl}
\begin{proof}
Let $M=\sum_{i\in I}S_i$, where $_RS_i\subset M$ are simple submodules. Then $\phi(M)=\sum_{i\in I}\phi(S_i)$, where for each $i$, either $\phi(S_i)=0$ or $\phi(S_i)\cong S_i$.
\end{proof}

\begin{crl}[\blue{Semisimplicity of submodules and quotient modules}]
If $_RM$ is semisimple, then so is (i) every \ul{nonzero submodule} and (ii) every \ul{nonzero quotient module} of $M$.
\end{crl}
\begin{proof}
Let $M=\sum_iS_i$, where ${}_RS_i\subset M$ are simple. (i) If $_RA\subsetneq M$, let $\pi:M\ra {M\over A},~m\mapsto m+A$. Then the \ul{nonzero quotient module} ${M\over A}=\pi(M)$ is semisimple as a nontrivial homomorphic image of a semisimple module. (ii) Next, consider a \ul{nonzero submodule} $0\neq{}_RB\subset M$. Since $M$ is completely decomposable, $M=B\oplus B'$ for some $_RB'\subset M$, and so $B\cong{B\oplus B'\over B'}={M\over B'}$ is semisimple as a nonzero quotient module. (\blue{footnote}\footnote{As a quicker alternative to (ii), recall that a submodule of a completely decomposable (CD) module is itself a CD module.})
\end{proof}

\begin{caution}
What is wrong with the following claim? Let $R$ be a PID (so, an $R$-module is projective $\iff$ free). Then every projective (i.e., free) $R$-module is completely decomposable, hence semisimple.
\end{caution}

\begin{questions}
Let $R$ be a ring. If $Lpd(R)\leq 1$ (Lpd the homological dimension), then we know (with a true converse as well) that every submodule of a projective $R$-module is projective. (i) Does it follow from this fact alone that a projective $R$-module $P$ cannot have nonzero torsion (i.e., that if $Hom_R(P,-)$ is exact then $T_R(P)=0$)? (ii) Does the answer to (i) change if $R$ is a PID?
\end{questions}

\begin{rmk}[\blue{Split sequences and semisimplicity}]
Consider an exact sequence of $R$-modules
\[
0\ra A\ra B\ra C\ra 0.
\]
It is clear that if $B$ is semisimple, then so are $A$ and $C$. Conversely, if (i) both $A$ and $C$ are semisimple and (ii) the sequence is split (e.g., when $A$ is injective or $C$ is projective), then it is also clear that $B$ is semisimple.
\end{rmk}

\section{Semisimple Rings}
\subsection{Relevant basic concepts}
\begin{dfn}[\blue{\index{Opposite! ring}{Recall: Opposite ring}}]
Let $R$ be a ring. The \ul{opposite ring} $R^{op}:=(R,\ast)$ of $R$ is the same set $R$ as a ring wrt a new multiplication (in terms of multiplication $R\times R\ra R,~(r,s)\mapsto rs$ in $R$) given by
\bea
\ast:R\times R\ra R,~(r,s)\mapsto r\ast s:=sr.\nn
\eea
Recall that $R^{op}$-modules ${}_{R^{op}}M$ (resp. left ideals $_{R^{op}}I\subset R^{op}$) are precisely right $R$-modules $M_R$ (resp. right ideals $I_R\subset R$). Similarly, $R$-modules ${}_RM$ (resp. left ideals $_RI\subset R$) are precisely right $R^{op}$-modules $M_{R^{op}}$ (resp. right ideals $I_{R^{op}}\subset R^{op}$).
\end{dfn}
\begin{dfn}[\textcolor{blue}{Recall:
\index{Annihilator left ideal}{Annihilator left ideal},
\index{Annihilator submodule}{Annihilator submodule},
\index{Faithful! module}{Faithful module}}]
Let $R$ be ring, $M$ an $R$-module, and $A\subset M$ a subset. The \ul{$R$-annihilator} of $A$ (i.e., annihilator of $A$ in $R$) is the left ideal (\blue{footnote}\footnote{If $N\subset M$ is a submodule, then it is clear that $Ann_R(N)=\bigcap_{n\in N}Ann_R(n)\subset R$ is an ideal (hence an intersection of left ideals that gives an ideal).})
\bea
Ann_R(A):=\big\{r\in R:rA=\{0\}\big\}=\{r\in R:ra=0~\txt{for all}~a\in A\}\subset R.\nn
\eea
Similarly, if $B\subset R$, the \ul{$M$-annihilator} of $B$ (i.e., annihilator of $B$ in $M$) is the subgroup (\blue{footnote}\footnote{If $B\subset R$ is a right ideal or $B\subset Z(R)$, then it is clear that $Ann_M(B)\subset M$ is a submodule.})
\bea
Ann_M(B):=\big\{m\in M:Bm=\{0\}\big\}=\{m\in M:bm=0~\txt{for all}~b\in B\}\subset M.\nn
\eea
The module $M$ is \ul{faithful} if $Ann_R(M)=0$. Equivalently, $_RM$ is faithful iff the ring homomorphism
\bit
\item[] $\phi:R\ra End_\Integer(M)=\Mor_\Integer(M,M),~r\mapsto \phi_r:m\mapsto rm$ ~~is injective.
\eit
\end{dfn}

\begin{dfn}[\blue{
\index{(Left) primitive ideal}{(Left) primitive ideal},
\index{Right-primitive ideal}{Right-primitive ideal},
\index{(Left) primitive ring}{(Left) primitive ring},
\index{Right-primitive ring}{Right-primitive ring},
\index{Jacobson radical}{{Jacobson radical}},
\index{Semi-primitive (or Jacobson semisimple) ring}{Semi-primitive (or Jacobson semisimple) ring},
\index{Idempotent subset}{Idempotent subset},
\index{Idempotent system}{Idempotent system},
\index{Nilpotent subset}{Nilpotent subset},
\index{Nilradical}{{Nilradical}},
\index{Nil subset}{Nil subset}}]\label{PrmtvIdDfn}
Let $R$ be a ring. An ideal $P\lhd R$ is a \ul{(left) primitive ideal} if it is the annihilator of a simple $R$-module (i.e., $P=Ann_R({}_RS)$ for some simple module $_RS$). An ideal $P\lhd R$ is a \ul{right-primitive ideal} if $P\lhd R^{op}$ is primitive. The ring $R$ is a \ul{(left) primitive ring} if the zero ideal $0\lhd R$ is primitive (i.e., there exists a faithful simple module $_RS$ in the sense $Ann_R(S)=0$). The ring $R$ is a \ul{right-primitive ring} if $R^{op}$ is primitive. The \ul{Jacobson radical} $J(R)$ of $R$ is the intersection of all maximal left ideals of $R$, i.e.,
\bea
\textstyle J(R):=\bigcap\{\txt{maximal left ideals}~{}_RM\subset R\}.~~~~(\txt{\blue{footnote}}\footnotemark).\nn
\eea
\footnotetext{If $R$ is primitive, then $J(R)=0$. Indeed, given any simple module $_RS$ and any $m\in S\backslash 0$, $Ann_R(m)\subset R$ is a maximal left ideal because $S\cong{R\over Ann_R(m)}$ via the map $R\ra S,~r\mapsto rm$, and so $S$ satisfies ~$Ann_R(S)=\bigcap\{Ann_R(m):m\in S\backslash 0\}\supset J(R)$.}
The ring $R$ is a \ul{semi-primitive ring} (or \ul{Jacobson-semisimple ring}) if $J(R)=0$. (\blue{footnote}\footnote{It is clear that if $R$ is primitive (which implies $J(R)=0$), then $R$ is semi-primitive.})

A subset $E\subset R$ is \ul{idempotent} if $EE=E$. A collection of elements $\{e_1,...,e_n\}\subset R$ is an \ul{idempotent system} if $\sum_i e_i=1$ and $e_ie_j=\delta_{ij}e_i$. (\blue{footnote}\footnote{For any $a\in R$, $R=R1=R(a+1-a)=Ra+R(1-a)$. Thus, if $e\in R$ is an idempotent element, i.e., $e^2=e$, then $R=Re+R(1-e)$, where $Re\cap R(1-e)=\{re=r'(1-e):r,r'\in R\}=\{0\}$, and so $R=Re\oplus R(1-e)$. More generally, for any idempotent system $\{e_1,...,e_n\}\in R$, it follows by induction on $n$ that $R=Re_1\oplus Re_2\oplus\cdots\oplus Re_n$.
}). A subset $S\in R$ is \ul{nilpotent} if $S^n=\{0\}$ for some $n\geq 1$ (where $S^n:=SS^{n-1}$). The \ul{nilradical} of $R$ is the collection of all nilpotent elements of $R$, i.e.,
\bea
Nilrad(R):=\{a\in R:~a~\txt{is nilpotent}\}=\{a\in R:~a^n=0~\txt{for some}~~n\geq 1\}. ~~(\txt{\blue{footnote}}\footnotemark)\nn
\eea
\footnotetext{$Nilrad(R)\subset R$ is not an ideal in general, but it is easy to see that it is an ideal when $R$ is commutative.}
A subset $N\subset R$ is \ul{nil} if every element of $N$ is nilpotent, i.e., $N\subset Nilrad(R)$.
\end{dfn}

Every nilpotent subset $N\subset R$ of a ring is nil (because $N^n=\{0\}$ $\iff$ $x_1x_2\cdots x_n=0$ for all $x_1,...,x_n\in N$) but a nil subset need not be nilpotent since in the ring $R:={k[x_1,x_2,x_3,\cdots]\over\langle x_1,x_2^2,x_3^3,\cdots\rangle}$ (where $k$ is a field), for example, the nilradical contains the nilpotent elements $x_i+\langle x_1,x_2^2,x_3^3,\cdots\rangle$ for all $i\geq 1$, and so cannot be nilpotent.

\begin{lmm}[\blue{A maximal ideal is primitive}]
(\blue{footnote}\footnote{A maximal \emph{left} ideal $_RM\subset R$ is not primitive in general since a primitive ideal is necessarily two-sided. However, we always have {\tiny $Ann_R(R/M)=\{a\in R:a(R/M)={aR+M\over M}=0\}=\{a\in R:aR\subset M\}\subset M$} (while {\tiny $M(R/M)={MR+M\over M}\neq 0$} in general).})
\end{lmm}
\begin{proof}
Let $M\vartriangleleft R$ be a maximal ideal. Then $Ann_R(R/M)=\left\{a\in R:a(R/M)={aR+M\over M}=0\right\}=\{a\in R:aR\subset M\}\subset M$, and since $MR=M\subset M$, we also have $M\subset Ann_R(R/M)$. Hence $M=Ann_R(R/M)$. Since ${R\over M}$ is a simple $R$-module, it follows that $M$ is a primitive ideal.
\end{proof}

\begin{lmm}[\blue{A primitive ideal is prime}]
\end{lmm}
\begin{proof}
  Let $P\lhd R$ be primitive, i.e., $P=\txt{Ann}_R(S)$, where $_RS$ is simple. Let $I,J\lhd R$ be such that $IJ\subset P$.  If $IJ\subset P$ but $I\not\subset P$ and $J\not\subset P$, then $IS=S=JS$ (since $_RS$ is simple) and so
\bit
\item[] $0=(IJ)S=I(JS)=IS=S$ ~~(a contradiction). \qedhere
\eit
 \end{proof}
\begin{rmk}
We therefore have inclusions ~$\{\txt{Maximal ideals}\}\subset\{\txt{Primitive ideals}\}\subset\{\txt{Prime ideals}\}$.
\end{rmk}

\begin{lmm}\label{JacRadAltLmm}
Let $R$ be a ring. A primitive ideal $P\lhd R$ is an intersection of maximal left ideals. Hence
\bit
\item[] $J(R):=\bigcap\{\txt{maximal left ideals}\}=\bigcap\{\txt{primitive ideals}\}=\bigcap\{Ann_R(S):\txt{for simple}~{}_RS\}$.
\eit
\end{lmm}
\begin{proof}
Let $J_1:=\bigcap\{\txt{primitive ideals}\}$. For any simple module $_RS$ and any $s\in S\backslash 0$, $Ann_R(s)\subset R$ is a maximal left ideal because $S\cong{R\over Ann_R(s)}$ via the map $R\ra S,~r\mapsto rs$. Therefore, every primitive ideal $P=Ann_R(_RS)=\bigcap\{M_s:=Ann_R(s)~|~s\in S\backslash 0\}$ is an intersection of maximal left ideals, and so
$J(R)\subset \bigcap_PP=J_1$. Meanwhile, for any maximal left ideal $_RM\subset R$, we have
\[
\textstyle P_M:=\txt{Ann}_R\left(R\over M\right)=\left\{a\in R:a(R/M)={aR+M\over M}=0\right\}=\{a\in R:aR\subset M\}\subset M,
\]
and so $J_1=\bigcap_PP\subset\bigcap_MP_M\subset\bigcap_M M=J(R)$. Hence $J_1=J(R)$.
\end{proof}

\subsection{The Jacobson radical in terms of units}
\begin{dfn}[\blue{
\index{Set! of left units}{Set of left units},
\index{Set! of right units}{Set of right units},
\index{Set! of units}{Set of units}}]
Let $R$ be a ring. We define
\bit
\item $U_l(R):=\{y\in R:~Ry=R\}=\{\txt{left units: elements with left inverses}\}$
\item $U_r(R):=\{y\in R:~yR=R\}=\{\txt{right units: elements with right inverses}\}$
\item $U(R):=\{y\in R:~Ry=R=yR\}=U_l(R)\cap U_r(R)=\{\txt{units: elements with inverses}\}$
\eit
\end{dfn}

\begin{thm}\label{JacRadThm}
Let $R$ be a ring. Each of the following equals~ $J(R):=\bigcap\{{}_RM\subset R:M~\txt{maximal}\}$.
\bit
\item[(1)] $J_1:=\bigcap\{P\lhd R:P~\txt{primitive}\}=\bigcap\{Ann_R({}_RS):{}_RS~\txt{simple}\}$
\item[(2)] $J_2:=\{x\in R:~1+Rx\subset U_l(R)\}$

\item[(3)] $J_3:=\{x\in R:~1+Rx\subset U(R)\}$

\item[(4)] $J_4:=\{x\in R:~1+xR\subset U(R)\}$

\item[(5)] $J_5:=\{x\in R:~1+xR\subset U_r(R)\}$
\item[(6)] $J_6:=\bigcap\{P\lhd R:P~\txt{right primitive}\}=\bigcap\{Ann_R(S_R):S_R~\txt{simple}\}$
\item[(7)] $J_r(R):=\bigcap\{M_R\subset R:M~\txt{maximal}\}$
\eit
\end{thm}
 \begin{proof}
{\flushleft (1)}
It follows from Lemma \ref{JacRadAltLmm} that $J_1=J(R)$.
\bit
\item \ul{About (6)}: With ``right'' replacing ``left'', the same argument as for $J_1=J(R)$ shows $J_6=J_r(R)$.
\eit
{\flushleft (2)} Let $x\in J(R)$, and suppose $x\not\in J_2$. Then (for some $r\in R$ and maximal $_RM\subset R$)
\bea
&& R(1+rx)\subsetneq R~\Ra~R(1+rx)\subset M~\Ra~1+rx,rx\in M~\Ra~1=(1+rx)-rx\in M\nn\\
&&~\Ra~M=R~~\txt{(a contradiction)}.\nn
\eea

Let $x\in J_2$, and suppose $x\not\in J(R)$. Then (for some maximal $_RM\subset R$, $r\in R$, and $m\in M$)
{\small\bea
 x\not\in M~\Ra~Rx+M=R~\Ra~rx+m=1~\Ra~m=1-rx\in U_l(R)~\Ra~R=Rm\subset M~~\txt{(a contradiction)}.\nn
 \eea}
\bit
\item \ul{About (5)}: With ``right'' replacing ``left'', the same argument as for $J_2=J(R)$ shows $J_5=J_r(R)$.
\eit

{\flushleft (3)}
 Let $x\in J(R)$, and suppose $x\not\in J_3$. Then (for some $r,r'\in R$ and maximal ${}_RM,{}_RM'\subset R$)
 \bea
 &&R(1+rx)\subsetneq R~~\txt{or}~~(1+r'x)R\subsetneq R~~~~\Ra~~~~R(1+rx)\subset M~~\txt{or}~~(1+r'x)R\subset M',\nn\\
 &&~~\Ra~~1+rx,~rx\in M~~\txt{or}~~1+r'x,~r'x\in M',\nn\\
 &&~~\Ra~~1=(1+rx)-rx\in M~~\txt{or}~~1=(1+r'x)-r'x\in M'~~\txt{( a contradiction)}.\nn
 \eea

 Let $x\in J_3$. Then ~$1+Rx\subset U(R)\subset U_l(R)$~ implies ~$x\in J_2=J(R)$~ by part (2).
\bit
\item \ul{About (4)}: With ``right'' replacing ``left'', the same argument as for $J_3=J(R)$ shows $J_4=J_r(R)$.
\eit

{\flushleft (4)}  Observe (from above) that both $J_3$ and $J_4$ are ideals. Therefore $J_4=J_3=J(R)$, since
\bit
\item[] \hspace{2cm} $1+J_3R=1+J_3=1+RJ_3\subset U(R)~~\Ra~~J_3\subset J_4$,
\item[] \hspace{2cm} $1+RJ_4=1+J_4=1+J_4R\subset U(R)~~\Ra~~J_4\subset J_3$. \qedhere
\eit
 \end{proof}

\begin{crl}
$J(R)=\{x\in R:1+RxR \subset U(R)\}$.
\end{crl}
\begin{proof}
Let $J':=\{x\in R:~1+RxR \subset U(R)\}$. Then $J'\subset J_3\cap J_4=J(R)$, since
 \bea
 1+J'R,1+RJ'\subset 1+RJ'R\subset U(R).\nn
 \eea
Similarly, ~$1+RJ(R)R=1+J(R)\subset U(R)$~ implies ~$J(R)\subset J'$.
\end{proof}

\begin{rmks}
{\flushleft (a)} Observe that $J(R)$ is the unique largest ideal $I\lhd R$ with the property that $1+I\in U(R)$. That is, if $I\lhd R$, then $1+I~\subset~U(R)$ iff $I\subset J(R)$. This is because $I=RI=IR=RIR$, and so $1+I\subset U(R)$ iff $1+RI\subset U(R)$, iff $1+IR\subset U(R)$, iff $1+RIR\subset U(R)$. Hence
\bea
\textstyle J(R)=\bigcup\{I\lhd R:1+I\subset U(R)\}.\nn
\eea
{\flushleft (b)} Let $x\in R$. If we simply have $1+x\in U(R)$, it is possible that $x\not\in J(R)$. (\blue{footnote}\footnote{Such as when $R$ is a field. Indeed $J(R)$ contains no units, while if $k$ is a field and $x\not\in\{0,-1\}$ in $k$, then $1+x,x\in U(k)$ and so $x\not\in J(k)=\{0\}$.}). To be certain that $x\in J(R)$, we must have $1+xr\in U(R)$ for all $r\in R$, or $1+rx\in U(R)$ for all $r\in R$, or $1+rxr'\in U(R)$ for all $r,r'\in R$.
\end{rmks}

\begin{crl}
If $R$ is a simple ring, then $J(R)=0$.
\end{crl}
\begin{proof}
Since every primitive ideal $P\lhd R$ is prime, it equals $0$ (as $R$ is simple). Hence $J(R)=0$.
\end{proof}

\begin{crl}\label{JacRadNilCrl}
If $_RI\subset R$ is a nil left ideal (or $I_R\subset R$ is a nil right ideal), then $I\subset J(R)$.
\end{crl}
\begin{proof}
Let ${}_RI\subset R$ be a nil left ideal and $x\in I$. Then $x^n=0$ implies {\small $1=1-x^n=(1-x)(1+x+\cdots+x^{n-1})=(1+x+\cdots+x^{n-1})(1-x)$}, and so $1-x\in U(R)$. That is, $1+RI=1+I=1-I\subset U(R)$, and so $I\subset J(R)$.
\end{proof}

\subsection{Module correspondences and Nakayama's lemma}
\begin{prp}[\blue{Recall: ${R\over I}$-module correspondence}]\label{RoverImod-RC}
Let $M$ be an $R$-module and $I\lhd R$ an ideal. Then $M/IM$ is both an $R$-module and an $R/I$-module, and we have a bijective correspondence (\blue{footnote}\footnote{Moreover, since the functor $-\otimes_RM$ is right-exact, the SES $0\ra I\hookrightarrow R\sr{\pi}{\ral} R/I\ra 0$ gives an exact sequence $I\otimes_RM\ra R\otimes_RM\sr{\pi\otimes id_M}{\ral}{R\over I}\otimes_R M\ra 0$, and so
~${R\over I}\otimes_RM\cong {R\otimes_RM\over\ker(\pi\otimes id_M)}={R\otimes_RM\over I\otimes_RM}\cong{M\over IM}$ (as $R$-modules and as $R/I$-modules).})
\bea
\textstyle \left\{\txt{$R$-submodules ${}_RN\subset {M\over IM}$}\right\}~~\sr{id}{\longleftrightarrow}~~\left\{\txt{${R\over I}$-submodules ${}_{R\over I}N\subset {M\over IM}$}\right\}.\nn
\eea
\end{prp}
\begin{proof}
Let $X$ be any $R$-module such that $IX=0$. Then the map
${R\over I}\times X\ra X,~(r+I,x)\mapsto rx$ is well-defined (\blue{footnote}\footnote{If $(r+I,x)=(r'+I,x')$, then $r+I=r'+I$ (i.e., $r-r'\in I$) and $x=x'$, and so $rx=rx'=(r-r')x'+r'x'=0+r'x'$.}), thereby making $X$ an ${R\over I}$-module. Moreover, because the original $R$-scalar multiplication $R\times X\ra X,~(r,x)\mapsto rx$ has exactly the same effect on $X$ as the ${R\over I}$-scalar multiplication, we have an automatic bijective correspondence
\bea
\textstyle \left\{\txt{$R$-submodules ${}_RN\subset X$}\right\}~~\sr{id}{\longleftrightarrow}~~\left\{\txt{${R\over I}$-submodules ${}_{R\over I}N\subset X$}\right\}.\nn
\eea
In particular, with $X:=M/IM$, we get the desired result.
\end{proof}

\begin{crl}
Let $R$ be a ring and $J:=J(R)$. Then $R$ and $R/J$ have the same simple modules, up to a bijective correspondence
~$\left\{\txt{simple modules}~{}_RS\right\}~~\sr{id}{\longleftrightarrow}~~\left\{\txt{simple modules}~{}_{R/J}S\right\}$.
\end{crl}
\begin{proof}
($\Ra$): If $_RS$ is simple, then $JS=0$ since $J\subset P:=\txt{Ann}_R(S)$. This shows $S={S\over JS}$ is an ${R\over J}$-module. To see that $S$ is simple as an ${R\over J}$-module, observe that (i) the correspondence theorem for modules and Proposition \ref{RoverImod-RC} imply that ${}_RM\subset R$ is a maximal left ideal $\iff$ $_{R/J}M/J\subset R/J$ is a maximal left ideal, and (ii) $_RS$ is simple $\iff$ there is a maximal left ideal ${}_RM\subset R$ such that $S\cong{R\over M}$, meanwhile $S\cong{R\over M}\cong {R/J\over M/J}$ (as $R$-modules and as $R/J$-modules by Proposition \ref{RoverImod-RC}).
{\flushleft ($\La$)}: Conversely, an $R/J$-module $_{R\over J}S$ is simple iff there is a maximal left ideal $_{R/J}M/J\subset R/J$ such that $S\cong{R/J\over M/J}$ (as $R/J$-modules). Since ${R/J\over M/J}$ is an $R$-module, so is $S$. Moreover, the $R$-module isomorphism ${R/J\over M/J}\cong{R\over M}$ and Proposition \ref{RoverImod-RC} imply that $S\cong{R/J\over M/J}\cong{R\over M}$ as $R$-modules. Hence $_RS$ is simple.
\end{proof}

\begin{crl}[\blue{Quotient-invariance of $J(R)$}]
Let $I\lhd R$ be an ideal. If $I\subset J(R)$, then $J\left({R\over I}\right)={J(R)\over I}$. In particular, if $I\lhd R$ is a nil ideal, then $I\subset J(R)$ by Corollary \ref{JacRadNilCrl}, and so $J(R/I)=J(R)/I$.
\end{crl}
\begin{proof}
By Proposition \ref{RoverImod-RC} at step (a) below and the correspondence theorem at step (b) below,
{\footnotesize\bea
\textstyle J\left({R\over I}\right):=\bigcap\limits_{\txt{maximal}~{}_{R/I}M/I\subset R/I}{M\over I}~\sr{(a)}{=}~\bigcap\limits_{\txt{maximal}~{}_RM/I\subset R/I}{M\over I}~\sr{(b)}{=}~\bigcap\limits_{\txt{maximal}~I\subset{}_RM\subset R}{M\over I}~\sr{(c)}{=}~{\mathop{\bigcap}_{\txt{maximal}~I\subset{}_RM\subset R}~M\over I}={J(R)\over I},\nn
\eea}where step (c) holds because using surjectivity of the natural map $\pi:R\ra R/I,~r\mapsto r+I$, we see that
\bit
\item[] $\bigcap_M{M\over I}=\pi\pi^{-1}\left(\bigcap_M{M\over I}\right)=\pi\left(\bigcap_M\pi^{-1}\left({M\over I}\right)\right)=\pi\left(\bigcap_M(M+I)\right)=\pi\left(\bigcap_MM\right)={\bigcap_MM\over I}$. \qedhere
\eit
\end{proof}

\begin{thm}[\blue{\index{Nakayama's lemma}{Nakayama's lemma}}]
Let $_RM$ be finitely generated (e.g., noetherian $_RM$) and $J:=J(R)$. If $JM=M$, then $M=0$. (That is, for every nonzero FG module ${}_RM\neq 0$, we have $J(R)M\subsetneq M$.)
\end{thm}
\begin{proof}
Suppose $M\neq 0$, and let $M=Rm_1+\cdots+Rm_t$, where $t$ is minimal. Then
\bea
&&Rm_1+\cdots+Rm_t=M=JM=Jm_1+\cdots+Jm_t~~\Ra~~m_t=a_1m_1+a_2m_2+\cdots+a_tm_t,~~a_i\in J,\nn\\
&&~~\Ra~~m_t=(1-a_t)^{-1}(a_1m_1+\cdots+a_{t-1}m_{t-1})\in Rm_1+\cdots+Rm_{t-1},\nn
\eea
which is a contradiction since $t$ was minimal.
\end{proof}

\begin{crl}[\blue{$J(R)$ kills all simple $R$-modules}]
Let $R$ be a ring and $_RS$ simple. Then $J(R)S=0$.
\end{crl}
\begin{proof}
Suppose $J(R)S\neq 0$. Then $J(R)S=RJ(R)S=S$, and so $S=0$ (a contradiction).
\end{proof}

\begin{crl}
For any ring $R$, the additive functor ~$F:={R\over J(R)}\otimes_R-:R\txt{-mod}\ra R\txt{-mod}$~ takes nonzero FG modules to nonzero FG modules (i.e., $F$ reflects zero-modules in the sense that for any $R$-module $M$, $F(M)=0$ implies $M=0$).
\end{crl}

\begin{crl}[\blue{Nilpotence of $J(R)$}]\label{JacRadNilp}
If $R$ is left (or right) artinian, then $J:=J(R)$ is nilpotent.
\end{crl}
\begin{proof}
It is clear (by the artinian property) that there is $n$ such that {\small $J\supset J^2\supset\cdots\supset J^n=J^{n+k}$} for all $k\geq 1$. Suppose $J^n\neq 0$, and let {\small $\C:=\{{}_RI\subset R:J^nI\neq 0\}\ni R$}. Then by the artinian property, $\C$ contains a minimal element $I_0$ with respect to inclusion $\subset$. Since $J^nI_0\neq 0$, there exists $x_0\in I_0$ such that $J^nx_0\neq 0$. This implies $J^nx_0=I_0=Rx_0$ (since {\small $J^nRx_0\supset J^nJ^nx_0=J^{n+n}x_0=J^nx_0\neq 0$} implies $Rx_0,J^nx_0\in\C$, meanwhile $I_0\supset Rx_0,J^nx_0$ is minimal in $\C$). That is, $J^nx_0=Rx_0\in\C$ is finitely generated, and so
\bit
\item[] \hspace{2cm} $JJ^nx_0=J^{n+1}x_0=J^nx_0~~\sr{\txt{Nakayama}}{\Longrightarrow}~~J^nx_0=0$ (a contradiction). \qedhere
\eit
\end{proof}

\subsection{Semisimplicity of rings and Maschke's theorem}
\begin{thm}[\blue{\index{Ring! semisimplicity criterion}{Ring semisimplicity criterion}}]\label{RngSmSpCrThm}
Let $R$ be a ring. The following are equivalent.
\bit
\item[(1)] $_RR$ is artinian and $J(R)=0$.~~~~ (Equivalently, ${}_{R/J(R)}R/J(R)$ is artinian)
\item[(2)] $_RR$ is semisimple.~~~~ (Equivalently, ${}_{R/J(R)}R/J(R)$ is semisimple)
\item[(3)] Every R-module is semisimple.~~~~ (Equivalently, every $R/J(R)$-module is semisimple)
\eit
\end{thm}
\begin{proof}
{\flushleft (1)$\Ra$(2)}: Assume $R$ is left artinian (i.e., $_RR$ is artinian). Pick any maximal left ideal $_RM_1\subset R$. If $M_1\neq J(R)$ choose another maximal left ideal $_RM_2\subset R$. If $M_1\cap M_2\neq J(R)$, choose another maximal left ideal $_RM_3\subset R$, and continue similarly.  Since $R$ is artinian, we get a finite sequence
    \bea
   \textstyle M_1\supsetneq M_1\cap M_2\supsetneq\cdots\supsetneq\bigcap_{i=1}^tM_i=J(R).\nn
    \eea
    Thus, because $\bigoplus_{i=1}^t{R\over M_i}$ is semisimple and the map {\small $\phi:R\ra\bigoplus_{i=1}^t{R\over M_i},~r\mapsto (r+M_1,...,r+M_t)$} has kernel {\small $\ker\phi=\bigcap_{i=1}^tM_i=J(R)$}, we get by the 1st isomorphism theorem (and the preceding corollaries) that
    \bea
    \textstyle {R\over J(R)}={R\over\ker\phi}\cong\phi(R)\subset\bigoplus_{i=1}^t{R\over M_i}~~~~\txt{(as $R$-modules and as $R/J(R)$-modules)}\nn
    \eea
    is semisimple as a nonzero submodule of a semisimple module.

{\flushleft (2)$\Ra$(3)}: Assume $_RR$ is semisimple. Let $R=\bigoplus_{i\in I}L_i$, where $_RL_i\subset R$ are minimal (simple) left ideals. Then for any R-module $_RM$, we have ~{\small $M=\sum_{m\in M}Rm=\sum_{m\in M}\sum_{i=1}^tL_im~\sr{(s)}{=}~\txt{Soc}~M$}~ where step (s) holds because the R-homomorphism $\phi:L_i\ra L_im,~l\mapsto lm$ (with $\ker\phi$ either $0$ or $L_i$) implies $L_im$ is either $0$ or a simple submodule of $M$. Hence $M$ is semisimple.

{\flushleft (3)$\Ra$(1)}: Assume every R-module is semisimple. Then, in particular, $_RR$ is semisimple. So, let {\small $R=\bigoplus_{i\in I}L_i$}, where $_RL_i\subset R$ are minimal left ideals. Then {\small $1=\sum_{i\in I}l_i$} (where $l_i=0$ a.e.f.) implies $1=l_{i_1}+l_{i_2}+\cdots+l_{i_t}$ for some finite $t$. Thus, {\small $R=Rl_{i_1}\oplus Rl_{i_2}\oplus\cdots\oplus Rl_{i_t}=L_{i_1}\oplus L_{i_2}\oplus\cdots\oplus L_{i_t}$}, which implies $R$ is left artinian (by induction on $t$) as a finite direct sum of simple (hence left artinian) R-modules.

Moreover, since {\small $M_j:=L_{i_1}\oplus\cdots\oplus L_{i_{j-1}}\oplus 0\oplus L_{i_{j+1}}\oplus\cdots\oplus L_{i_t}\cong {R\over L_{i_j}}$} are maximal left ideals, we have the descending chain {\small $R\supset M_1\supset M_1\cap M_2\supset M_1\cap M_2\cap M_3\supset\cdots\supset \bigcap_{j=1}^n M_j\supset\cdots$} that must terminate with a nilpotent $J(R)$, since $_RR$ is artinian, and so {\small $J(R)=\bigcap_{j=1}^t M_j=0$}.
\end{proof}

\begin{rmk}
In the proof of (3)$\Ra$(1) above, we have actually obtained a \ul{composition series} (\blue{footnote}\footnote{To be introduced later, a \ul{composition series} of an $R$-module $_RM$ is any (finite) descending chain of submodules $M:=M_0\supset M_1\supset M_2\supset\cdots\supset M_t=0$ such that the quotients (called \ul{factors} of the series) $S_i:={M_i\over M_{i+1}}$ are simple $R$-modules.})
    \bea
    R=C_0\supset C_1\supset C_2\supset\cdots\supset C_{t+1}=0,~~~~C_j:=L_{i_j}\oplus L_{i_{j+1}}\oplus\cdots\oplus L_{i_t},\nn
    \eea
    for $R$. We shall see later that a ring $R$ has a composition series if and only if it is left (or right) artinian.
\end{rmk}

\begin{crl}[\blue{Semisimplicity of a simple ring}]
Let $R$ be a simple ring (in particular, $J(R)=0$). Then $R$ is semisimple if and only if $R$ is left artinian.
\end{crl}

\begin{crl}[\textcolor{blue}{\index{Maschke's theorem}{Maschke's theorem}}]
Let $G$ be a finite group, $k$ a field, and $R:=kG:=\sum_{g\in G}kg$ the $k$-algebra of $G$. Then $_RR$ is semisimple $\iff$ $char(k)$ does not divide $|G|$ in $\Integer$.

(In particular, if $char(k)=0$, then $kG$ is left-semisimple for every finite group $G$.)
\end{crl}
\begin{proof}
Let $a:=\sum_{g\in G}g1_k\in R\backslash 0$. Then $ga=a$ for all $g\in G$. This implies (i) $a^2=|G|a$ and (ii) $ra=ar\in ka$ for all $r\in R$, i.e., $Ra=ka\subset Z(R)$. Hence, $ka=Ra\lhd R$ is a 1-dimensional ideal.
{\flushleft ($\Ra$)}: Assume $char(k)$ divides $|G|$ in $\Integer$. (\blue{footnote}\footnote{Recall that $char(k)$ is either $0$ or a positive prime $p\in\Integer$.}). Then $char(k)=p=|1_k|_+$ (the additive order of $1_k$) for a positive prime $p\in\Integer$, and so $|G|1_k=cp1_k=c|1_k|_+1_k=0$ in $R:=kG$. Thus $a^2=|G|a=|G|1_ka=0$, and so the nilpotent ideal $0\neq Ra=ka\subset J(R)$, which implies $J(R)\neq 0$. Hence $R$ is not semisimple.
{\flushleft ($\La$)}: Assume $char(k)$ does not divide $|G|$ in $\Integer$. Then the additive order $|1_k|_+$ does not divide $|G|$ in $\Integer$, and so $|G|1_k\neq 0$ in $R:=kG$. Let $M$ be an $R$-module and $_RN\subset{}_RM$ an $R$-submodule. Then we know $_kN\subset{}_kM$ is a $k$-subspace, and that $_kM={}_kN\oplus{}_kN'$ for some $k$-subspace $_kN'\subset{}_kM$. Let $\pi:{}_kM\ra {}_kN$ be the obvious $k$-linear projection onto $_kN$ in $_kM$ satisfying $\pi|_N=id_N$ and $\pi(M)=N$ (i.e., the inclusion $i_N:N\hookrightarrow M$ is a $k$-section of $\pi$). Define
\bea
\textstyle \phi:=\txt{average}_G(\pi):M\ra N,~m\mapsto {1\over |G|}\sum_{g\in G} g^{-1}\pi(gm).\nn
\eea
Then $\phi$ is both $k$-linear (by the $k$-linearity of $\pi$) and $G$-linear: Indeed, for $m\in M$, $g\in G$, we have
{\small
\bea
\textstyle \phi(gm)={1\over |G|}\sum_{h\in G} h^{-1}\pi(hgm)={1\over |G|}\sum_{h\in G} g~(hg)^{-1}\pi(hgm)=g~\phi(m).\nn
\eea
}Therefore, $\phi$ is an R-homomorphism satisfying $\phi|_N=id_N$ and $\phi(M)=N$ (i.e., the inclusion $i_N:N\hookrightarrow M$ is an $R$-section of $\phi$), and so ${}_RM=\im i_N\oplus_R\ker\phi=N\oplus_R\ker\phi$.

Thus ${}_RM$ is completely decomposable (hence semisimple). Since every $R$-module is semisimple, it follows from Theorem \ref{RngSmSpCrThm} that $_RR$ is semisimple.
\end{proof}

\section{Basic Results on Idempotents and Matrix Rings}
The following result shows that a simple ring (i.e., a ring $R$ whose only ideals are $0$ and $R$) can have a large number of maximal left ideals (which may or may not be isomorphic as $R$-modules).
\begin{prp}
Let $k$ be an infinite field. The simple ring $R:=Mat_n(k)$ has infinitely many maximal left ideals.
\end{prp}
\begin{proof}
For each $\al\in k$, let
~{\scriptsize
$M_\al:=\left\{
\left[
  \begin{array}{cccccc}
    \ast & \ast & \cdots &\ast & a_1 & \al a_1 \\
     \ast & \ast & \cdots &\ast & a_2 & \al a_2 \\
     \vdots & \vdots & \vdots &\vdots & \vdots & \vdots \\
     \ast & \ast & \cdots &\ast & a_n & \al a_n \\
  \end{array}
\right]:~\ast,\ast,\cdots,\ast,a_i\in k
\right\}$}. ~These are maximal left ideals of $R$, since
{\scriptsize
$S_\al:={R\over M_\al}=\left\{
\left[
  \begin{array}{ccccc}
    0 &0 & \cdots &0 & a_1 \\
   0 &0 & \cdots &0 & a_2 \\
   \vdots & \vdots & \vdots & \vdots & \vdots \\
   0 & 0 & \cdots &0 & a_n \\
  \end{array}
\right]+M_\al:~a_i\in k
\right\}\cong k^n$} is a simple R-module. Also, if
{\scriptsize $M_\al\cap M_\beta
\neq L_{n-2}:=\left\{
\left[
  \begin{array}{cccccc}
    \ast & \ast & \cdots &\ast & 0 & 0 \\
     \ast & \ast & \cdots &\ast & 0 & 0 \\
     \vdots & \vdots & \vdots &\vdots & \vdots & \vdots \\
     \ast & \ast & \cdots &\ast & 0 & 0 \\
  \end{array}
\right]:\ast,\ast,\cdots,\ast\in k
\right\}$}, then with $m\in (M_\al\cap M_\beta)\backslash L_{n-2}$,
{\scriptsize
\[
m=\left[
  \begin{array}{cccccc}
    \ast & \ast & \cdots &\ast & a_1 & \al a_1 \\
     \ast & \ast & \cdots &\ast & a_2 & \al a_2 \\
     \vdots & \vdots & \vdots &\vdots & \vdots & \vdots \\
     \ast & \ast & \cdots &\ast & a_n & \al a_n \\
  \end{array}
\right]=\left[
  \begin{array}{cccccc}
    \ast' & \ast' & \cdots &\ast' & a_1' & \beta a_1' \\
     \ast' & \ast' & \cdots &\ast' & a_2' & \beta a_2' \\
     \vdots & \vdots & \vdots &\vdots & \vdots & \vdots \\
     \ast' & \ast' & \cdots &\ast' & a_n' & \beta a_n' \\
  \end{array}
\right]~~~~\Ra~~~~\al=\beta.
\]}Hence, if $\al\neq\beta$, then $M_\al\cap M_\beta=L_{n-2}$, and so $M_\al\neq M_\beta$.
\end{proof}

\subsection{Ideals of a homogeneous matrix ring}
\begin{dfn}[\blue{
\index{Matrix! group}{Matrix group},
\index{Matrix! ring}{Matrix ring},
\index{Unit matrices}{Unit matrices},
\index{Unit idempotent system}{Unit idempotent system},
\index{Matrices over a subset of a ring}{Matrices over a subset of a ring}}]
Let $R$ be a ring and $M_{m\times n}(R):=\{[a_{ij}]:a_{ij}\in R\}$ the additive \ul{group of $m\times n$ matrices} (which is a \ul{matrix ring} when $m=n$). The \ul{unit matrices} in $M_{m\times n}(R)$ are
\bea
e_{ij}:=[\delta_{ri}\delta_{sj}],~~~~1\leq i\leq m,~~1\leq j\leq n,~~\txt{where if $m=n$ we have}~~e_{ij}e_{kl}=\delta_{jk}e_{il}.\nn
\eea
In the \ul{matrix ring} $M_n(R):=M_{n\times n}(R)=\sum_{ij}Re_{ij}$, the set $\{e_{11},e_{22},\cdots,e_{nn}\}$ is an idempotent system of matrices (\blue{footnote}\footnote{The collection $\{e_{11},e_{22},\cdots,e_{nn}\}$ is an idempotent system in $M_n(R)$, since $\sum_i e_{ii}=\mathds{1}_n$ and $e_{ii}e_{jj}=\delta_{ij}e_{ii}.$}) (call it \ul{unit idempotent system}) in terms of which
\bea
\textstyle M_n(R)=\sum_{ij}Re_{ij}=\sum_{j=1}^n M_n(R)e_{jj}=\bigoplus_{j=1}^nL_j(R)\nn
\eea
is a direct sum of left ideals $L_j(R):=M_n(R)e_{jj}\cong R^n$ that correspond to the columns of $M_n(R)$. Given a subset $A\subset R$, the set of \ul{matrices over $A$} (or the set of \ul{$A$-matrices}) is
\bea
\textstyle M_n(A):=\sum_{ij}Ae_{ij}=\{[a_{ij}]:a_{ij}\in A\}\subset M_n(R).\nn
\eea
\end{dfn}
\begin{lmm}
Let $R$ be a ring. (i) For any ideal $I\lhd R$, we get an ideal $M_n(I)\lhd M_n(R)$. (ii) Every ideal of $M_n(R)$ is of the form $M_n(I)$ for an ideal $I\lhd R$, i.e., $J\lhd M_n(R)$ $\iff$ $J=M_n(I)$ for some $I\lhd R$.
\end{lmm}
\begin{proof}
Let $S:=M_n(R)$. Given an ideal $I\lhd R$, let $J:=M_n(I)$.
{\flushleft (i)} We need to show $J\lhd S$. Observe that $(J,+)$ is a subgroup of $(S,+)$ since $J-J=\sum Ie_{ij}-\sum Ie_{ij}=\sum(I-I)e_{ij}\subset\sum Ie_{ij}=J$. Also, $SJ=\left(\sum Re_{ij}\right)(\sum Ie_{kl})=\sum RIe_{ij}e_{kl}=\sum I\delta_{jk}e_{il}=\sum Ie_{il}=J$, i.e., $SJ\subset J$, and similarly, $JS\subset J$.
{\flushleft (ii)} By part (i) if $J=M_n(I)$ for some $I\lhd R$, then $J\lhd S$. Conversely, let $J\lhd S$. It is clear that as a subset of
{\scriptsize $M_n(R)=\sum Re_{ij}=\left[
          \begin{array}{cccc}
            R & R & \cdots & R \\
            R & R & \cdots & R \\
            \vdots & \vdots & \ddots & \vdots \\
            R & R & \cdots & R \\
          \end{array}
        \right]$},
$J$ has the general form
{\scriptsize $J=\sum I_{ij}e_{ij}=\left[
          \begin{array}{cccc}
            I_{11} & I_{12} & \cdots & I_{1n} \\
            I_{21} & I_{22} & \cdots & I_{2n} \\
            \vdots & \vdots & \ddots & \vdots \\
            I_{n1} & I_{n2} & \cdots & I_{nn} \\
          \end{array}
        \right]$} for some $I_{ij}\subset R$. We will use the fact that $J$ is an ideal to show that $I_{ij}$ are all equal to some ideal $I\lhd R$. With $I_{ij}-I_{ij}:=\{u-v:u,v\in I_{ij}\}$, since $(J,+)\subset\big(S,+\big)$ is a subgroup,
{\small\[
\textstyle J-J=\sum I_{ij}e_{ij}-\sum I_{ij}e_{ij}=\sum (I_{ij}-I_{ij})e_{ij}\subset J=\sum I_{ij}e_{ij}~~\Ra~~I_{ij}-I_{ij}\subset I_{ij},~~\txt{for each pair}~(i,j),\nn
\]}and so each $(I_{ij},+)$ is a subgroup of $(R,+)$. Also, {\small $SJ=\left(\sum R_{ik}e_{ik}\right)\left(\sum I_{k'j}e_{k'j}\right)=\sum R_{ik}I_{k'j}e_{ik}e_{k'j}=\sum R_{ik}I_{k'j}\delta_{kk'}e_{ij}=\sum R_{ik}I_{kj}e_{ij}\subset J=\sum I_{ij}e_{ij}$} implies {\footnotesize $\sum_{k=1}^n R_{ik}I_{kj}=\sum_{k=1}^n RI_{kj}\subset I_{ij}$} for all $i,j$, i.e.,
\bea
\label{ex2-eq1}&&\textstyle SJ\subset J~~\Ra~~\txt{for all $i,j$},~~\sum_{k=1}^n R_{ik}I_{kj}=\sum_{k=1}^n RI_{kj}\subset I_{ij},~~\sr{(s)}{\Longrightarrow}~~I_{kj}\subset I_{ij}~~\txt{for all}~~k,\\
\label{ex2-eq2}&&\textstyle JS\subset J~~\Ra~~\txt{for all $i,j$},~~\sum_{k=1}^n I_{ik}R_{kj}=\sum_{k=1}^n I_{ik}R\subset I_{ij},~~\sr{(s)}{\Longrightarrow}~~I_{ik}\subset I_{ij}~~\txt{for all}~~k,
\eea
where step (s) holds because $R$ has a unity. Since (\ref{ex2-eq1}) implies the $I$'s on each column of $J$ are equal and (\ref{ex2-eq2}) implies the $I$'s on each row of $J$ are equal, it follows that the subsets $I_{ij}\subset R$ are all equal to some ideal $I\lhd R$. Hence
{\scriptsize $J=\sum Ie_{ij}=M_n(I)
  =\left[
          \begin{array}{cccc}
            I & I & \cdots & I \\
            I & I & \cdots & I \\
            \vdots & \vdots & \ddots & \vdots \\
            I & I & \cdots & I \\
          \end{array}
        \right]$}.
\end{proof}

\subsection{Idempotent projections: Jacobson radicals of matrix rings and product rings}
\begin{lmm}\label{JacIdLmm1}
Let $R$ be a ring and $e\in R\backslash 0$ a nonzero idempotent. Then {\small $R_e:=eRe=\{ere:r\in R\}$} is a ring with unity $e$.
\end{lmm}
\begin{proof}
With the addition and multiplication inherited from $R$, it suffices to check the following three things:
(i) It is clear that $(eRe,+)\subset (R,+)$ is a subgroup since $eRe-eRe=e(Re-Re)=e(R-R)e\subset eRe$. (ii) Also, $(eRe)(eRe)=eReRe\subset eRe$. (iii) Finally, for each $r\in R$, ~$e(ere)=(ere)e=ere$. Hence $eRe$ is a ring with unity $e$.
\end{proof}

\begin{lmm}\label{JacIdLmm2}
Let $R$ be a ring, $M$ an $R$-module, and $e\in R\backslash 0$ a nonzero idempotent. The following hold.
 \bit[leftmargin=0.9cm]
 \item[(a)] $eM$ is an $eRe$-module.
 \item[(b)] If $S$ is a simple $R$-module and $eS\neq 0$, then $eS$ is a simple $eRe$-module.
 \item[(c)] If $M$ is a semisimple $R$-module and $eM\neq 0$, then $eM$ is a semisimple $eRe$-module.
 \eit
 \end{lmm}
\begin{proof}
{\flushleft (a)} With the scalar multiplication $eRe\times eM\ra eM$ inherited from $_RM$ as a restriction of the scalar multiplication $R\times M\ra M$, it is enough to show that $(eM,+)\subset(M,+)$ is a subgroup, $(eRe)eM\subset eM$, and $e(em)=em$ for all $m\in M$ (which are immediate).
{\flushleft (b)} For any $0\neq es\in eS\subset S$, we have~
$(eRe)es=e(Res)=eS$, since $_RS$ is simple (and so $Res=S$).
{\flushleft (c)} Let $M=\sum_{i\in I}S_i$, where $_RS_i\subset M$ are simple. Then $eM=\sum_{i\in I}eS_i$ is a sum of simple $eRe$-submodules.
\end{proof}

\begin{thm}\label{JacIdThm1}
Let $R$ be a ring and $e\in R\backslash 0$ a nonzero idempotent. Then $J(eRe)=eJ(R)e$.
\end{thm}
\begin{proof}
{\flushleft \ul{Showing $eJ(R)e\subset J(eRe)$}}: Let $x\in eJ(R)e\subset J(R)$, and recall that $x=ex=xe=exe$. Then for any $r\in R$, $1+rx\in U(R)$, which implies $1=(1+rx)(1+y_r)=(1+y_r)(1+rx)$ for some $y_r\in R$. Multiplying on both sides by $e$, we get ~$e=(e+erx)(e+ey_re)=(e+ey_re)(e+erx)$ for all $r\in R$,
\bea
~~\Ra~~e+erx=e+erex\in U(eRe)~~\txt{for all}~~ere\in eRe,~~\Ra~~x\in J(eRe).\nn
\eea
{\flushleft\ul{Showing $J(eRe)\subset eJ(R)e$}}: It suffices to show $J(eRe)\subset J(R)$, because
\[
J(eRe)\subset J(R)~~\Ra~~J(eRe)=eJ(eRe)e\subset eJ(R)e,
\]and so also suffices to show $J(eRe)\subset Ann_R(S)$ for all simple $R$-modules $S$. So, let $_RS$ be simple. If $eS=0$, then $(eRe)S=0$ implies $J(eRe)\subset eRe\subset Ann_R(S)$. If $eS\neq 0$, then $eS$ is a simple $eRe$-module by Lemma \ref{JacIdLmm2}, and so $0=J(eRe)eS=J(eRe)S$, $\Ra$ $J(eRe)\subset Ann_R(S)$.
\end{proof}

\begin{crl}\label{JacRadMatCrl}
Let $R$ be a ring and $M_n(R)$ the matrix ring. Then $J\big(M_n(R)\big)=M_n\big(J(R)\big)$.
\end{crl}
\begin{proof}
Recall that $J\big(M_n(R)\big)=M_n(I)$ for some ideal $I\lhd R$. Therefore, using Theorem \ref{JacIdThm1},
\bea
Ie_{ii}=e_{ii}M_n(I)e_{ii}=e_{ii}J(M_n(R))e_{ii}=J(e_{ii}M_n(R)e_{ii})=J(Re_{ii})=J(R)e_{ii},~~\txt{for all}~~i=1,...,n,\nn
\eea
and so $I=J(R)$, since $\sum_{i=1}^n e_{ii}=\mathds{1}_n$.
\end{proof}

\begin{crl}\label{JacRadProdCrl}
Let $R_1,...,R_n$ be rings. Then $J(R_1\times R_2\times\cdots\times R_n)=J(R_1)\times J(R_2)\times\cdots\times J(R_n)$.
\end{crl}
\begin{proof}~
By induction on $n$, it suffices to prove that for two rings $R,S$, we have $J(R\times S)=J(R)\times J(S)$. Let $R':=R\times S$, and consider the idempotent system $e_1:=(1,0)$, $e_2:=(0,1)$ in $R'$. Then $R\times S=R'=R'e_1+R'e_2=Re_1+Se_2=e_1R'e_1+e_2R'e_2$. Therefore,
\begin{align}
&J(R\times S)=J(R\times S)(e_1+e_2)=J(R\times S)e_1+J(R\times S)e_2=e_1J(R\times S)e_1+e_2J(R\times S)e_2\nn\\
&~~~~=J(e_1(R\times S)e_1)+J(e_2(R\times S)e_2)=J(Re_1)+J(Se_2)=J(R)e_1+J(S)e_2=J(R)\times J(S).\nn \qedhere
\end{align}

\end{proof}

\begin{lmm}\label{JacIdLmm3}
Let $R$ be a ring. A proper ideal $P\lhd R$ is prime $\iff$ for any ideal $I,J\lhd R$,
\bea
IJ\subset P\subset I\cap J~~\Ra~~I=P~~\txt{or}~~J=P.\nn
\eea
\end{lmm}
\begin{proof}
{\flushleft ($\Ra$)}: Assume $P$ is prime. Then for any ideals $I,J\lhd R$, the chain of containments $IJ\subset P\subset I\cap J$ implies $I\subset P\subset I$ or $J\subset P\subset J$, and so $I=P$ or $J=P$. ($\La$): Conversely, assume that for any ideal $I,J\lhd R$, the containments $IJ\subset P\subset I\cap J$ imply $I=P$ of $J=P$. Then for any ideals $A,B\lhd R$,
\bit
\item[] \hspace{1cm} $AB\subset P$ ~~$\Ra$~~ $(A+P)(B+P)\subset AB+P=P\subset(A+P)\cap(B+P)$,
\item[] \hspace{1cm} ~~~~$\Ra$~~ ~$A+P=P$ ~or~ $B+P=P$, ~~$\Ra$~~ $A\subset P$ ~or~ $B\subset P$. \qedhere
\eit
\end{proof}

\begin{rmk}
If $P\lhd R$ is prime, then the preceding result implies for all $A,B,C\lhd R$, the following hold:
\bit
\item $ABC\subset P$ ~~$\Ra$~~ $AB\subset P$ ~or~ $C\subset P$ ~~$\Ra$~~ $A\subset P$ ~or~ $B\subset P$ ~or~ $C\subset P$.
\item $ABC\subset P\subset A\cap B\cap C$ ~~$\Ra$~~ $AB\subset P\subset A\cap B$ ~or~ $C=P$, ~~$\Ra$~~ $A=P$ or $B=P$ or $C=P$.
\eit
Since these two conditions also each imply that $P$ is prime, it follows by induction that $P\lhd R$ is prime $\iff$ for all $I_1,...,I_n\lhd R$,~ $I_1I_2\cdots I_n\subset P$ $\Ra$ $I_i\subset P$ for some $i$, $\iff$ for all $I_1,...,I_n\lhd R$,~ $I_1I_2\cdots I_n\subset P\subset I_1\cap I_2\cap\cdots\cap I_n$ $\Ra$ $I_i=P$ for some $i$.
\end{rmk}

\begin{thm}\label{JacIdThm2}
Let $R$ be a ring, $e\in R\backslash 0$ a nonzero idempotent, and $\I(R),\J(eRe)$ the sets of ideals of $R,eRe$ respectively. Then, with the observation that $eIe=I\cap eRe$ (\blue{footnote}\footnote{This is because ~$eIe\subset I\cap eRe=e(I\cap eRe)e\subset e(I\cap R)e=eIe$.}), the maps
\begin{align}
\phi:\I(R)\ra \J(eRe),~I\mapsto eIe=I\cap eRe~~~~\txt{and}~~~~\psi:\J(eRe)\ra \I(R),~J\mapsto RJR\nn
\end{align}
have the following properties: (1) $\phi\circ\psi=id_{\J(eRe)}$. Thus $\phi$ is surjective (\blue{footnote}\footnote{By the surjectivity of $\phi$, every ideal of $eRe$ is of the form $eIe=I\cap eRe$ for some ideal $I$ of $R$. In particular, $R$ has at least as many ideals as $eRe$.}) and $\psi$ is injective. \\
(2) $\phi$ and $\psi$ induce a bijection between $\{\txt{prime ideals of $R$ not containing $e$}\}$ and $\{\txt{prime ideals of $eRe$}\}$.
\end{thm}
\begin{proof}
{\flushleft (1)} $\phi\psi(J)=\phi(RJR)=e(RJR)e=(eRe)J(eRe)=J$, since $J=eJ=Je=eJe$.
{\flushleft (2)} Let $\P(R):=\{\txt{prime ideals $e\not\in P\lhd R$}\}$ and $\Q(eRe):=\{\txt{prime ideals $Q\lhd eRe$}\}$, where $e\not\in P$ $\Ra$ $ePe\subsetneq eRe$. If $P\in \P(R)$, then $\phi(P)=ePe$ is a prime ideal of $eRe$, since
\bea
&&(eIe)(eI'e)\subset ePe~~\Ra~~R(eIe)(eI'e)R=(ReIeR)(ReI'eR)\subset RePeR\subset P,\nn\\
&&~~\Ra~~ReIeR\subset P~~\txt{or}~~ReI'eR\subset P,~~\Ra~~eIe=eReIeRe\subset ePe~~\txt{or}~~eI'e=eReI'eRe\subset ePe.\nn
\eea
This shows that ~$\phi|_{\P(R)}:\P(R)\ra \Q(eRe),~P\ra ePe$. It \ul{remains} to show that (i) $\psi|_{\Q(eRe)}:\Q(eRe)\ra \P(R),~Q\mapsto RQR$ and (ii) $\phi|_{\P(R)}$ is injective (i.e., if $ePe=eP'e$, then $P=P'$).

So, let $Q\lhd eRe$ be prime. We will show $P:=\psi(Q)=RQR\lhd R$ is the unique prime ideal of $R$ such that $\phi(P)=Q$ (i.e, with $Prime(R)$ denoting the set of prime ideals of $R$, we have ~$\phi^{-1}(Q)\cap Prime(R)=\{P\}$). Consider, as a poset under inclusion, the nonempty set
\bea
\E:=\phi^{-1}(Q)=\{I\lhd R:~eIe=Q\}~~\ni~~RQR.\nn
\eea
Every nonempty chain $\{I_\ld\}_{\ld\in\Ld}$ in $\E$ has an upper bound $I^u:=\bigcup_\ld I_\ld$ in $\E$, since $a,b\in I^u$ $\Ra$ $a,b\in I_\ld$ for some $\ld$, $\Ra$ $a-b\in I_\ld\subset I^u$ and $Ra,aR\subset I_\ld\subset I^u$ (making $I^u\lhd R$ an ideal), and
\[
eI^ue=\textstyle e(\bigcup_\ld I_\ld)e=\bigcup_\ld eI_\ld e=\bigcup_\ld Q=Q,~~\Ra~~I^u\in\E.\nn
\]
By Zorn's lemma, $\E$ contains a maximal element $P\lhd R$. \ul{Suppose $P$ is not prime}. Then by Lemma \ref{JacIdLmm3} there exist ideals $I,I'\lhd R$ such that $II'\subset P\subset I\cap I'$ but $I\neq P$ and $I'\neq P$ (i.e., $P\subsetneq I$ and $P\subsetneq I'$). Since $P$ is maximal in $\E$, it follows that $eIe\neq Q=ePe$ and $eI'e\neq Q=ePe$. But because $Q$ is prime,
\bea
&&II'\subset P\subset I\cap I'~~\Ra~~(eIe)(eI'e)\subset eII'e\subset ePe=Q\subset e(I\cap I')e\subset (eIe)\cap(eI'e),\nn\\
&&~~\Ra~~eIe=Q~~\txt{or}~~eI'e=Q~~(\txt{a \ul{contradiction}}).\nn
\eea
Finally, we need to show $P$ is unique. Suppose $P'\lhd R$ is a prime ideal such that $eP'e=ePe$. Then
\bea
&&(ReR)P'(ReR)=ReP'eR=RePeR\subset P~~\Ra~~ReR\subset P~~\txt{or}~~P'\subset P,\nn\\
&&~~\Ra~~P'\subset P~~(\txt{because $e\not\in P$ $\Ra$ $ReR\not\subset P$}).\nn
\eea
By symmetry, $P\subset P'$ as well, and so $P'=P$.
\end{proof}

\begin{question}
Let $R$ be a commutative ring and $e\in R\backslash 0$ a nonzero idempotent. What is the relation between the prime ideals of the localization $R[\{1,e\}^{-1}]$ and those of the ring ~$Re=eRe$?
\end{question}

\section{The Artin-Wedderburn Theorem}
\subsection{Endomorphism ring of a direct sum of modules: Matrix representation}
\begin{rmk}[\blue{Evaluation module of the endomorphism ring}] Let $R$ be a ring and $_RM$ an $R$-module. Recall that the endomorphism ring of $M$ is $End_R(M):=Hom_R(M,M)$. $M$ is naturally an $End_R(M)$-module with scalar multiplication the evaluation map ~$End_R(M)\times M\ra M,~(f,m)\mapsto fm:=f(m)$.
\end{rmk}

\begin{dfn}[\blue{
\index{Nonhomogeneous matrix ring}{Nonhomogeneous matrix ring},
\index{Channel matrix}{Channel matrix},
\index{Channel maps}{Channel maps}}]
Let $R$ be a ring, $\{M_i:i=1,...,n\}$ a finite set of $R$-modules, $H_{ij}:=Hom_R(M_j,M_i)$ the associated hom sets, and
\bea
\textstyle H:=\big\{[f_{ij}]_{ij}:~f_{ij}\in H_{ij}\big\}=\sum\limits_{i,j}H_{ij}e_{ij}=\left[
                                                                                  \begin{array}{ccc}
                                                                                   H_{11}  & \cdots & H_{1n}\\
                                                                                   \vdots  & \ddots & \vdots \\
                                                                                   H_{n1}  & \cdots & H_{nn} \\
                                                                                  \end{array}
                                                                                \right],~~~~e_{ij}e_{kl}=\delta_{jk}e_{il},\nn
\eea
the set of $n\times n$ matrices with $(i,j)$th entry $f_{ij}\in H_{ij}$. Then $H$ is a ring with respect to matrix addition and matrix multiplication given respectively by
\[
\textstyle [f_{ij}]+[g_{ij}]:=[f_{ij}+g_{ij}]~~\txt{and}~~[f_{ij}][g_{ij}]:=\left[\sum_lf_{il}\circ g_{lj}\right],~~~\txt{where}~~~f_{il}\circ g_{lj}:M_j\sr{g_{lj}}{\ral}M_l\sr{f_{il}}{\ral} M_i.
\]

Let ${}_RM:=\bigoplus_{i=1}^n M_i$ be the direct sum and \bt[column sep=small] M_j\ar[r,hook,"q_j"] & M\ar[r,two heads,"p_i"]& M_i\et the canonical injections and projections. If $f\in End_R(M):=Hom_R(M,M)$, the \ul{channel matrix} of $f$ in $H$ is the matrix $[f_{ij}]\in H$ with $(i,j)$th entry  given by the composition $f_{ij}:=p_ifq_j:M_j\sr{q_j}{\ral}M\sr{f}{\ral}M\sr{p_i}{\ral} M_i$ (the $(i,j)$th \ul{channel map} of $f$).

We clearly have $p_iid_Mq_j=p_iq_j=\delta_{ij}~id_{M_j}$, which upon left-multiplication with $q_i$ and right-multiplication with $q_j$ gives idempotents $\eta_i:=q_ip_i\in End_R(M)$ satisfying $\eta_i\eta_j=\delta_{ij}\eta_j$ or $(\eta_i-\delta_{ij}id_M)\eta_j=0$. This implies $(\sum_i\eta_i-id_M)\eta_j=0$ or $\sum_i\eta_i=id_M$ (by properties of the direct sum), i.e., $\{\eta_1,...,\eta_n\}$ is an idempotent system in $End_R(M)$. Hence we have the basic identities
{\small\begin{align}
\label{BasicChIds}\textstyle p_iid_Mq_j=p_iq_j=\delta_{ij}~id_{M_j}~\txt{(diagonal channel matrix)}~~\txt{and}~~\sum_{i=1}^nq_ip_i=id_M~~(\txt{completeness}).
\end{align}}
\end{dfn}

\begin{lmm}\label{ArtWeddLmm}
Let $R$ be a ring, $\{M_i:i=1,...,n\}$ $R$-modules, {\small $H_{ij}:=Hom_R(M_j,M_i)$, $H:=\sum_{i,j}H_{ij}e_{ij}$}, {\small $M:=\bigoplus_iM_i$}, and $End_R(M):=Hom_R(M,M)$. We have an isomorphism of rings $End_R(M)\cong H$, i.e.,
\[
\textstyle Hom_R\left(\bigoplus_iM_i,\bigoplus_iM_i\right)\cong_{\txt{Ring}}\sum_{i,j}Hom_R(M_i,M_j)e_{ij},~~~e_{ij}\in Mat_{n\times n}(Hom_R),~~~e_{ij}e_{kl}=\delta_{jk}e_{il}.\nn
\]
\end{lmm}
\begin{proof}
Consider the maps {\small $\phi:End_R(M)\ra H,~f\mapsto [p_ifq_j]$} and {\small $\psi:H\ra End_R(M),~[f_{ij}]\mapsto\sum_{i,j}q_if_{ij}p_j$}. It is clear by matrix addition that $\phi(f+g)=\phi(f)+\phi(g)$ and $\psi([f_{ij}]+[g_{ij}])=\psi([f_{ij}])+\psi([g_{ij}])$. Also,
{\small\begin{align}
&\textstyle\phi(fg)=[p_ifgq_j]=[p_ifid_Mgq_j]=[\sum_kp_ifq_kp_kgq_j]=[p_ifq_j][p_igq_j]=\phi(f)\phi(g),~~~~\txt{for all}~~f,g\in End_R(M),\nn\\
&\textstyle\psi([f_{ij}][g_{ij}])=\psi([\sum_kf_{ik}g_{kj}])=\sum_{i,j}q_i(\sum_kf_{ik}g_{kj})p_j=\sum_{i,j,k}q_if_{ik}g_{kj}p_j=\sum_{i,j,k,k'}q_if_{ik}p_kq_{k'}g_{k'j}p_j\nn\\
&\textstyle~~~~=(\sum_{i,k}q_if_{ik}p_k)(\sum_{k',j}q_{k'}g_{k'j}p_j)=\psi([f_{ij}])\psi([g_{ij}]),~~~~\txt{for all}~~[f_{ij}],[g_{ij}]\in H,\nn\\
&\textstyle \phi(1_{End_R(M)})=\phi(id_M)=[p_iid_Mq_j]=[\delta_{ij}id_{M_j}]=1_H,\nn\\
&\textstyle \psi(1_H)=\psi([\delta_{ij}id_{M_j}])=\sum_{i,j} q_i\delta_{ij}p_j=\sum_iq_ip_i=1_{End_R(M)},\nn
\end{align}}and so $\phi$ and $\psi$ are ring homomorphisms. Finally, $\psi\circ\phi=id_{End_R(M)}$ and $\phi\circ\psi=id_H$, since
{\small\begin{align}
&\textstyle\psi\phi(f)=\psi\left([p_ifq_j]\right)=\sum_{i,j}q_i(p_ifq_j)p_j=\sum_{i,j}q_ip_ifq_jp_j=id_Mf~id_M=f,~~~~\txt{for all}~~f\in End_R(M),\nn\\
&\textstyle\phi\psi\big([f_{ij}]\big)=\phi\big(\sum_{i,j} q_if_{ij}p_j\big)=\big[p_i\big(\sum_{i',j'} q_{i'}f_{i'j'}p_{j'}\big)q_j\big]=\big[\sum_{i',j'} p_iq_{i'}f_{i'j'}p_{j'}q_j\big]\nn\\
&\textstyle~~~~=\big[\sum_{i',j'} \delta_{ii'}id_{M_{i'}}f_{i'j'}\delta_{j'j}id_{M_j}\big]=[id_{M_i}f_{ij}id_{M_j}]=[f_{ij}],~~~~\txt{for all}~~[f_{ij}]\in H.\nn \qedhere
\end{align}}
\end{proof}

\subsection{Matrix Representation of Semisimple Rings}
\begin{dfn}[\blue{Recall: \index{Division! (or Divisible) ring}{Division ring}}]
A ring $D$ is a division ring if $U(D)=D\backslash 0$ (\blue{footnote}\footnote{That is, every nonzero element of $D$ is a unit.}).
\end{dfn}

\begin{dfn}[\blue{\index{Centralizer in a ring}{Centralizer in a ring}}]
Let $R$ be a ring and $X\subset R$ a subset. The centralizer of $X$ in $R$ is the set ~$C(X):=\{c\in R:~cx=xc~\txt{for all}~x\in X\}$.
\end{dfn}

\begin{thm}[\blue{\index{Artin-Wedderburn theorem}{Artin-Wedderburn theorem}}]\label{ArtWeddThm}
Let $R$ be a ring. The following are equivalent.
\begin{enumerate}[leftmargin=0.7cm]
\item[(1)] {\small $R\cong\prod_{i=1}^tM_{n_i}(D_i)$}~ for division rings $D_i$. (I.e., $R$ is a finite product of matrix rings over division rings.)
\item[(2)] $R$ is left-artinian and $J(R)=0$. (I.e., $_RR$ is artinian and $J(R)=0$.)
\item[(3)] $R$ is left-semisimple. (I.e., $_RR$ is semisimple.)
\end{enumerate}
\end{thm}
\begin{proof}
Since (2)$\iff$(3) was proved in Theorem \ref{RngSmSpCrThm}, it is enough to prove (1)$\Ra$(2) and (3)$\Ra$(1).
{\flushleft \ul{(1)$\Ra$(2)}}: Assume $R\cong \prod_{i=1}^tM_{n_i}(D_i)$, where $D_i$ are division rings. Observe that for each $i$,
{\small\[
\textstyle R_i:=M_{n_i}(D_i)=\sum_{k,k'}D_ie^{(i)}_{kk'}=\sum_{k=1}^{n_i}R_ie^{(i)}_{kk}=\bigoplus_{k=1}^{n_i}R_ie^{(i)}_{kk},
\]}
where $\{e^{(i)}_{kk}:k=1,...,n_i\}$ is the unit idempotent system in $R_i$, is left-semisimple (hence left-artinian) as a sum of simple (i.e., minimal) left ideals $L_k:=R_ie^{(i)}_{kk}\subset R_i$. So, by induction on $t$, $R$ is left-artinian. Moreover, by Corollaries \ref{JacRadMatCrl} and \ref{JacRadProdCrl}, we also have $J(R)\cong \prod_{i=1}^tM_{n_i}(J(D_i))=\prod_{i=1}^tM_{n_i}(0)=0$.

{\flushleft \ul{(3)$\Ra$(1)}}: Assume $_RR$ is semisimple. Let $R=\bigoplus_{i=1}^tL_i^{n_i}=\bigoplus_{\al\in A}L_\al$, where $L_i$ are non-isomorphic simple left ideals, {\small $A:=\{(i,k_i):1\leq k_i\leq n_i,~1\leq i\leq t\}=\bigsqcup_{i=1}^t~\{i\}\times(n_i+1)$}, and $L_{(i,k_i)}:=L_i$ for all $1\leq k_i\leq n_i$. Let {\small $End_R(R):=Hom_R(R,R)$}, {\small $H:=\sum_{\al,\beta} H_{\al\beta}e_{\al\beta}$} with $H_{\al\beta}:=Hom_R(L_\al,L_\beta)$ as in Lemma \ref{ArtWeddLmm}, and $D_i:=End_R(L_i):=Hom_R(L_i,L_i)$. Then using Lemma \ref{ArtWeddLmm} and Schur's lemma (Lemma \ref{SchurLmm}), we get
{\small\begin{align}
&\textstyle R\cong End_R(R)\cong H=\sum_{\al,\beta}H_{\al\beta}e_{\al\beta}=\sum_{(i,k_i)(j,k'_j)}H_{(i,k_i)(j,k'_j)}e_{(i,k_i)(j,k'_j)}=\sum_{i,j}\sum_{k_i,k'_j} Hom_R(L_i,L_j)e_{(i,k_i)(j,k'_j)}\nn\\
&\textstyle~~~~=\sum_{i,j}\sum_{k_i,k'_j}\delta_{ij}D_ie_{(i,k_i)(j,k'_j)}=\sum_i\sum_{k_i,k'_i}D_ie_{(i,k_i)(i,k'_i)},~~~~\txt{where}~~~~e_{(i,k_i)(i,k'_i)}:=e_{k_ik'_i}^{(i)}E_{ii},\nn\\
&\textstyle~~~~=\sum_i\sum_{k_i,k'_i}D_ie_{k_ik'_i}^{(i)}E_{ii},~~~~~~\txt{with $\{E_{ii}\}_{i=1}^t\subset H$ an idempotent system, i.e.,}~~E_{ii}E_{jj}=\delta_{ij}E_{jj},~~\sum_i E_{ii}=1_H,\nn\\
&\textstyle~~~~=\sum_iM_{n_i}(D_i)E_{ii}\cong\prod_iM_{n_i}(D_i),~~~~~~\txt{where}~~~~M_{n_i}(D_i):=\sum_{k_i,k'_i}D_ie_{k_ik'_i}^{(i)}.\nn \qedhere
\end{align}}
\end{proof}

\begin{crl}[\blue{\cite[Corollary 7.45, p.553]{rotman2010}}]\label{ArtWeddCrl1}
``Left'' can be replaced by ``Right'' in every part of Theorem \ref{ArtWeddThm}. Hence, a ring is left-semisimple $\iff$ right-semisimple $\iff$ semisimple.
\end{crl}
\begin{proof}
The matrix rings have left-right symmetry (equivalently, recall that $J(R)\lhd R$ is a two-sided ideal that is determined by maximal left ideals and maximal right ideals in exactly the same way).
\end{proof}

\begin{crl}\label{ArtWeddCrl2}
Let $R$ be a ring. The following are equivalent.
\begin{enumerate}
\item[(1)] $R$ is a finite product of matrix rings over division rings.
\item[(2)] $R$ is left-artinian and $J(R)=0$. (I.e., $_RR$ is artinian and $J(R)=0$.)
\item[(3)] $R$ is left-semisimple. (I.e., $_RR$ is semisimple.)
\item[(4)] Every $R$-module is semisimple.
\item[(5)] Every $R$-module is completely decomposable.
\item[(6)] Every SES of $R$-modules $0\ra A\ra B\ra C\ra 0$ is split.
\item[(7)] Every $R$-module is projective.
\item[(8)] Every $R$-module is injective.
\item[(9)] $Lpd(R)=Lid(R)=0$. (Recall that $Lpd$ and $Lid$ are the homological dimensions of $R$.)
\end{enumerate}
\end{crl}

\begin{crl}[\blue{Artin-Wedderburn theorem for simple rings}]\label{ArtWeddCrl3}
If $R$ is a simple ring (so, $J(R)=0$), then the following are equivalent.
\bit
\item[(1)] $R\cong M_n(D)$ for a division ring $D$ (and $n\geq1$).
\item[(2)] $R$ is left-artinian. (I.e., $_RR$ is artinian.)
\item[(3)] $R$ is left-semisimple. (I.e., $_RR$ is semisimple.)
\eit
\end{crl}
\begin{proof}
A finite product of matrix rings $\prod_{i=1}^tM_{n_i}(D_i)$, for division rings $D_i$, is simple $\iff$ $t=1$.
\end{proof}

\section{Noetherian-Artinian Interaction}
\subsection{Composition series and the Jordan-Holder theorem}
\begin{dfn}[\blue{
\index{Series of a module}{Series of a module},
\index{Length of a series}{Length of a series},
\index{Refinement of a series}{Refinement of a series},
\index{Composition! factor}{Composition factor},
\index{Multiplicity of a composition factor}{Multiplicity of a composition factor},
\index{Factor-equivalent series}{Factor-equivalent series},
\index{Composition! series}{Composition series},
\index{Module! of finite length}{Module of finite length},
\index{Equivalent composition series}{Equivalent composition series}}]
Let $R$ be a ring and $M$ an $R$-module. A \ul{series} of $M$ (or \ul{series} for $M$) is a finite descending chain of submodules $\C:M=M_0\supset M_1\supset M_2\supset\cdots\supset M_t=0$ (also briefly written ~$\C=\{M_i\}_{i=0}^t$~ or ~$\C=\{M_i:0\leq i\leq t\}$). The integer $t$ is the \ul{length} of the series. A \ul{refinement} of the series $C$ is any series $\C':M=N_0\supset N_1\supset N_2\supset\cdots\supset N_s=0$ containing $\C$, i.e., such that $\C\subset \C'$ in the sense $\{M_0,...,M_t\}\subset\{N_0,...,N_s\}$.

Each of the quotient modules ${M_{i-1}\over M_i}$ is a \ul{composition factor} of the series $\C$. The \ul{multiplicity} of a given composition factor is the number of times that composition factor occurs up to isomorphism (i.e., the number of composition factors of the series that are isomorphic to the given composition factor). Two series $\C_1:~M=M_0\supset M_1\supset M_2\supset\cdots\supset M_t=0$ and $\C_2:~M=N_0\supset N_1\supset N_2\supset\cdots\supset N_s=0$ of $M$ are \ul{factor-equivalent series} if (i) $t=s$ and (ii) there is a permutation $\sigma\in S_t$ of the indices such that ${M_{i-1}\over M_i}\cong {N_{\sigma(i)-1}\over N_{\sigma(i)}}$ for all $i=1,...,t$.

A series $M=M_0\supset M_1\supset M_2\supset\cdots\supset M_t=0$ of $M$ is a \ul{composition series} (making $M$ a \ul{finite length module} or \ul{module of finite length}) if each of its composition factors is a simple R-module (i.e., $M_i\subset M_{i-1}$ is a maximal submodule of $M_{i-1}$ for each $1\leq i\leq t$). Two composition series are \ul{equivalent composition series} if they are factor-equivalent series.
\end{dfn}

\begin{lmm}\label{SchrRefLmm}
Let $M$ be an $R$-module and $A,B,A'\subset M$ submodules. Then we have the following.
\bit
\item[] (1) ${A\over A\cap B}\cong {A+B\over B}$ and ${B\over A\cap B}\cong{A+B\over A}$. (2) If $A'\subset A$, then $A\cap(A'+B)=A'+A\cap B$.
\eit
\end{lmm}
\begin{proof}
(1) is the 2nd isomorphism theorem. For (2), observe that
\[
A\cap A',A\cap B\subset A\cap(A'+B)~~\Ra~~A\cap A'+A\cap B=A'+A\cap B\subset A\cap(A'+B),
\]
and we also have $A\cap(A'+B)\subset A'+A\cap B$, since
\[
a=a'+b\in A\cap(A'+B)~~\Ra~~b=a-a'\in A\cap B,~~\Ra~~a=a'+b\in A'+A\cap B. \qedhere
\]
\end{proof}

\begin{thm}[\textcolor{blue}{\index{Schreier refinement theorem}{Schreier refinement theorem}}]\label{SchrRefThm}
Let $R$ be a ring and $M$ an $R$-module. Any two series
\[
\C:M=M_0\supset M_1\supset\cdots\supset M_t=0~~~~\txt{and}~~~~\D:M=N_0\supset N_1\supset\cdots\supset N_s=0
\] of $M$ can be refined to obtain two factor-equivalent series of $M$.
\end{thm}
\begin{proof}
For each $0\leq i\leq t$, define ~$M_{i,j}:=(M_i+N_j)\cap M_{i-1}=M_i+N_j\cap M_{i-1}$,~ for all ~$0\leq j\leq s$,~ i.e.,
\bea
\cdots\supset M_{i-2}\supset M_{i-1}=\ub{M_{i,0}\supset M_{i,1}\supset\cdots\supset M_{i,s-1}\supset M_{i,s}}_{\txt{Refinement of $M_{i-1}\supset M_i$ using all of $\D$}}=M_i\supset M_{i+1}\supset\cdots,\nn
\eea
which gives a refinement~ $\C_1:=\{M_{i,j}:0\leq i\leq t,~0\leq j\leq s\}$~ of $\C$.

Similarly, for each $0\leq j\leq s$, define ~$N_{j,i}:=(N_j+M_i)\cap N_{j-1}=N_j+M_i\cap N_{j-1}$,~ for all ~$0\leq i\leq t$,~ i.e.,
\bea
\cdots\supset N_{j-2}\supset N_{j-1}=\ub{N_{j,0}\supset N_{j,1}\supset\cdots\supset N_{j,t-1}\supset N_{j,t}}_{\txt{Refinement of $N_{j-1}\supset N_j$ using all of $\C$}}=N_j\supset N_{j+1}\supset\cdots,\nn
\eea
which gives a refinement~ $\D_1:=\{N_{j,i}:0\leq i\leq t,~0\leq j\leq s\}$~ of $\D$.

Using Lemma \ref{SchrRefLmm} we see (as detailed below) that ~${M_{i,j-1}\over M_{i,j}}\cong{M_{i-1}\cap N_{j-1}\over M_{i-1}\cap(M_i+N_j)\cap N_{j-1}}\cong {N_{j,i-1}\over N_{j,i}}$~ for all $i,j$.
\[
\textstyle {M_{i,j-1}\over M_{i,j}}={M_{i,j}+M_{i-1}\cap N_{j-1}\over M_{i,j}}\cong {(M_{i-1}\cap N_{j-1})\over M_{i,j}\cap(M_{i-1}\cap N_{j-1})}={M_{i-1}\cap N_{j-1}\over M_{i-1}\cap(M_i+N_j)\cap N_{j-1}}\sr{\txt{symmetry}}{\cong}{N_{j,i-1}\over N_{j,i}}.
\]
This shows the refinements $\C_1$ and $\D_1$ are factor-equivalent series.
\end{proof}

\begin{thm}[\textcolor{blue}{\index{Jordan-Holder theorem}{Jordan-Holder theorem}: Uniqueness of composition series}]\label{JordHoldThm}
Let $M$ be an $R$-module and $\C:M=M_0\supset M_1\supset\cdots\supset M_t=0$ a composition series of $M$. Then given any series of unequal terms $\D:M=N_0\supsetneq N_1\supsetneq\cdots\supsetneq N_s=0$, the following hold: (i) We can refine $\D$ to obtain a composition series of $M$. (ii) $s\leq t$. (iii) The simple factors occurring in $\C$ are unique up to isomorphism and multiplicity.

(That is, a composition series of a module is the finest possible series of the module.)
\end{thm}
\begin{proof}
Since $\C$ is a composition series, every refinement of $\C$ has the same number of nonzero factors, because if $M_{i-1}\supset M_{i,j}\supset M_i$, then ${M_{i-1}\over M_{i,j}}\cong {M_{i-1}/M_i\over M_{i,j}/M_i}$ where, as a submodule of a simple module, either $M_{i,j}/M_i=0$ or $M_{i,j}/M_i=M_{i-1}/M_i$. Thus refinement cannot change the length of $\C$.

(i) Also, $\C$ and $\D$ have factor-equivalent refinements (by Theorem \ref{SchrRefThm}), and so the factor-equivalent refinement for $\D$ will have the same number of nonzero factors as $\C$. (ii) It follows therefore that $s\leq t$. (iii) By the factor-equivalence of the refinements, the resulting composition factors are unique up to isomorphism and permutation (hence unique up to multiplicity or number of occurrences of a given factor).
\end{proof}

\begin{dfn}[\blue{
\index{Module of finite length}{Module of finite length},
\index{Length of a module}{Length of a module}}]
Let $R$ be a ring. An $R$-module $M$ has \ul{finite length} (and is therefore called a \ul{module of finite length}) if it has a composition series, in which case the \ul{$R$-length} of $M$, written ~$\txt{length}(M):=\txt{length}_R(M)$,~ is the length of any composition series for $M$.
\end{dfn}

\begin{thm}[\textcolor{blue}{Module composition series: Existence I}]\label{ComSeExThm1}
A module $_RM$ has a composition series if and only if it is both noetherian and artinian.
\end{thm}
\begin{proof}
{\flushleft ($\Ra$)}: Assume $M$ has a composition series $M=M_0\supset M_1\supset\cdots\supset M_t=0$. Then by induction on $t$, and the fact that a simple module is both artinian and noetherian (along with the fact that in a SES $0\ra A\ra B\ra C\ra 0$, $B$ is artinian, or noetherian, iff both $A,B$ are), we see that $M$ is both noetherian and artinian. ($\La$): Conversely, assume $M$ is both artinian and noetherian. Then by the noetherian property (which ensures every submodule is FG and so guarantees the existence of a maximal submodule within each submodule) we have a chain of submodules $M=M_0\supset M_1\supset M_2\supset\cdots$, where $M_i$ is a maximal submodule of $M_{i-1}$. By the artinian property, this chain terminates, and so we get a composition series $M=M_0\supset M_1\supset M_2\supset\cdots\supset M_t=0$. (\blue{footnote}\footnote{That is, the noetherian property allows us to progressively choose maximal submodules within maximal submodules, while the artinian property ensures the process terminates, resulting in a composition series for $M$.})
\end{proof}

\begin{crl}\label{ComSeExCrl}
A module $_RM$ has a composition series if and only if some submodule $_RN\subset{}_RM$ and the quotient module ${M\over N}$ each have a composition series. (Equivalently, in a SES ~$0\ra A\ra B\ra C\ra 0$, ~$B$ has a composition series iff $A,B$ each have a composition series).
\end{crl}
\begin{proof}
$M$ has a composition series $\iff$ $M$ is both noetherian and artinian $\iff$ $N$ and $M/N$ are each both noetherian and artinian $\iff$ $N$ and $M/N$ each have a composition series.
\end{proof}

\begin{thm}[\textcolor{blue}{Module composition series: Existence II}]\label{ComSeExThm2}
Let $_RM$ be a semisimple module. Consider any direct sum expansion $M=\bigoplus_{i\in I}S_i$, for simple modules $_RS_i$. Then $M$ has a composition series (equivalently, the index set $I$ is finite) $\iff$ $M$ is artinian, $\iff$ $M$ is noetherian.

(That is, a semisimple module has a composition series iff it is noetherian, iff it is artinian.)
\end{thm}
\begin{proof}
If $I$ is finite, then it is clear by induction on $t:=|I|$ that $M$ is artinian (resp, noetherian). Conversely, if $I$ is not finite, then we have a strictly descending (resp, strictly ascending) chain of submodules that does not stabilize, and so $M$ is not artinian (resp, not noetherian).
\end{proof}

\subsection{The Hopkins-Levitzki theorems: Noetherian-Artinian interaction}
\begin{thm}[\textcolor{blue}{\index{Hopkins-Levitzki theorem I}{Hopkins-Levitzki theorem I}: For rings, left-artinian $\Ra$ left-noetherian}]\label{HopLevThm1}
Let $R$ be a ring. If $R$ is left-artinian (resp. right-artinian), then $R$ is also left-noetherian (resp. right-noetherian).
\end{thm}
\begin{proof}
Assume $_RR$ is artinian. Let $J:=J(R)$. Then $_{R/J}R/J$ is a semisimple artinian module, and so has a composition series. Also, because $J$ is nilpotent (by Corollary \ref{JacRadNilp}), we have the series
\bea
\label{J-series}{}_RR=J^0\supset J^1\supset\cdots\supset J^n=0.
\eea
Since $J{J^i\over J^{i+1}}=0$ and ${}_R({J^i\over J^{i+1}})$ is artinian, we get a semisimple module ${}_{R/J}({J^i\over J^{i+1}})$ that is artinian and so has a composition series (by Theorem \ref{ComSeExThm2}). (\blue{footnote}\footnote{Again by Theorem \ref{ComSeExThm2}, each ${}_{R/J}({J^i\over J^{i+1}})$ is noetherian.}). Recall that $R$ and $R/J$ have the same simple modules. Thus, using the composition series for the modules $_{R/J}({J^i\over J^{i+1}})$ we can refine the series (\ref{J-series}) to obtain a composition series for $_RR$ (making $_RR$ noetherian by Theorem \ref{ComSeExThm1}) as follows:
\[
\textstyle{}_{J/R}({J^i\over J^{i+1}})=S_{i1}\oplus S_{i2}\oplus\cdots\oplus S_{it_i}={J^i_1\over J^{i+1}}\oplus{J^i_2\over J^{i+1}}\oplus\cdots \oplus{J^i_{t_i}\over J^{i+1}}={J^i_1+\cdots+J^i_{t_i}\over J^{i+1}},~~J^i_1\cap J^i_2\cap\cdots\cap J^i_{t_i}=J^{i+1},
\]
and so with $M^i_r:=\sum_{k=r}^{t_i}J_k^i$, for $1\leq r\leq t_i$, we get a refinement of the $i$th step $J^i\supset J^{i+1}$ of (\ref{J-series}) as follows:
\bea
J^i=M^i_1\supset M^i_2\supset\cdots \supset M^i_{t_i}=J^i_{t_i}\supset J^{i+1},\nn
\eea
where by construction the composition factors are simple, since
{\small\[
\textstyle{M^i_r\over M^i_{r+1}}={\sum_{k=r}^{t_i}J^i_k\over \sum_{k=r+1}^{t_i}J^i_k}\cong {\left(\sum_{k=r}^{t_i}J^i_k\right)/J^{i+1}\over \left(\sum_{k=r+1}^{t_i}J^i_k\right)/J^{i+1}}={\bigoplus_{k=r}^{t_i}J^i_k/J^{i+1}\over \bigoplus_{k=r+1}^{t_i}J^i_k/J^{i+1}}\cong{J^i_r\over J^{i+1}}\cong {}_{R/J}S_{ir}~\eqv~{}_RS_{ir}. \qedhere
\]}
\end{proof}

\begin{rmks}
(i) Theorem \ref{HopLevThm1} says $R$ has a left-composition series $\iff$ $R$ is left-artinian.

(ii) A left-semisimple ring (with unity) is both left-noetherian and left-artinian, and so $_RR$ has a left-composition series.

(iii) However a ring that has a left-composition series might not be left-semisimple. For example, if $R$ is left-artinian and $J(R)\neq 0$, then $R$ is not left-semisimple, but has a left-composition series by Theorem \ref{HopLevThm1}. It follows therefore that \ul{a left-noetherian ring might not be left-artinian}.

(iv) The following result shows that if a ring $R$ is already left-artinian (or right-artinian), then its noetherian and artinian properties become equivalent.
\end{rmks}

\begin{thm}[\textcolor{blue}{\index{Hopkins-Levitzki theorem II}{Hopkins-Levitzki theorem II}: For right-artinian rings, left-artinian $\Lra$ left-noetherian}]\label{HopLevThm2}
Let $R$ be a left-artinian (resp. right-artinian) ring. Then $R$ is right-artinian $\iff$ right-noetherian (resp. left-artinian $\iff$ left-noetherian).
\end{thm}
(The following proof mirrors the proof of Theorem \ref{HopLevThm1} accordingly.)
\begin{proof}
($\Ra$) If $R$ is right-artinian, it is clear by symmetry in Theorem \ref{HopLevThm1} that $R$ is right-noetherian. ($\La$) Conversely, assume $R_R$ is noetherian. Let $J:=J(R)$. Since $_RR$ is artinian, $_{R/J}R/J$ is semisimple, and so $R/J_{R/J}$ is semisimple by Corollary \ref{ArtWeddCrl1}. Also, because $J$ is nilpotent (by Corollary \ref{JacRadNilp}), we have the series
\bea
\label{J-series2}R_R=J^0\supset J^1\supset\cdots\supset J^n=0.
\eea
Since ${J^i\over J^{i+1}}J=0$ and $({J^i\over J^{i+1}})_R$ is noetherian, we get a semisimple module $({J^i\over J^{i+1}})_{R/J}$ that is noetherian and so has a composition series (by left-right symmetry in Theorem \ref{ComSeExThm2}). (\blue{footnote}\footnote{Again by left-right symmetry in Theorem \ref{ComSeExThm2}, each $({J^i\over J^{i+1}})_{R/J}$ is artinian.}). Recall that $R$ and $R/J$ have the same simple modules. Thus, using the composition series for the modules $({J^i\over J^{i+1}})_{R/J}$ we can refine the series (\ref{J-series2}) to obtain a composition series for $R_R$ (making $R_R$ artinian by Theorem \ref{ComSeExThm1}) as follows:
\bea
&&\textstyle({J^i\over J^{i+1}})_{J/R}=S_{i1}\oplus S_{i2}\oplus\cdots\oplus S_{it_i}={J^i_1\over J^{i+1}}\oplus{J^i_2\over J^{i+1}}\oplus\cdots \oplus{J^i_{t_i}\over J^{i+1}}={J^i_1+\cdots+J^i_{t_i}\over J^{i+1}},~~J^i_1\cap J^i_2\cap\cdots\cap J^i_{t_i}=J^{i+1},\nn
\eea
and so with $M^i_r:=\sum_{k=r}^{t_i}J_k^i$, for $1\leq r\leq t_i$, we get a refinement of the $i$th step $J^i\supset J^{i+1}$ of (\ref{J-series2}) as follows:
\bea
J^i=M^i_1\supset M^i_2\supset\cdots \supset M^i_{t_i}=J^i_{t_i}\supset J^{i+1},\nn
\eea
where by construction the composition factors are simple, since
{\small\[
\textstyle{M^i_r\over M^i_{r+1}}={\sum_{k=r}^{t_i}J^i_k\over \sum_{k=r+1}^{t_i}J^i_k}\cong {\left(\sum_{k=r}^{t_i}J^i_k\right)/J^{i+1}\over \left(\sum_{k=r+1}^{t_i}J^i_k\right)/J^{i+1}}={\bigoplus_{k=r}^{t_i}J^i_k/J^{i+1}\over \bigoplus_{k=r+1}^{t_i}J^i_k/J^{i+1}}\cong{J^i_r\over J^{i+1}}\cong S_{ir}{~}_{R/J}~\eqv~S_{ir}{~}_R. \qedhere
\]}
\end{proof}

\begin{dfn}[\blue{Semi-primary ring, Semi-primary module}]
A ring $R$ is a \ul{semi-primary ring} if (i) ${R\over J(R)}$ is semisimple and (ii) $J(R)$ is nilpotent. A \ul{semi-primary module} is a module over a semi-primary ring.
\end{dfn}

\begin{thm}[\textcolor{blue}{\index{Hopkins-Levitzki theorem III}{Hopkins-Levitzki theorem III}: For semi-primary modules, artinian $\Lra$ noetherian}]\label{HopLevThm3}
Let $R$ be a semi-primary ring and $M$ an $R$-module. Then the following are equivalent.
\bit
\item[(1)] $M$ is artinian. ~~(2) $M$ is noetherian. ~~(3) $M$ has a composition series.
\eit
\end{thm}
(Once again, the following proof mirrors the proofs of Theorems \ref{HopLevThm1} and \ref{HopLevThm2} accordingly.)
\begin{proof}
``(3) $\Ra$ (1) and (2)'' are clear by Theorem \ref{ComSeExThm1}. We will prove ``(1) or (2) $\Ra$ (3)''. Assume $M$ is artinian or noetherian. Let $J:=J(R)$. Since $J$ is nilpotent, we have a series
\bea
\label{JM-series}{}_RM\supset JM\supset J^2M\supset\cdots\supset J^nM=0.
\eea
Since $J{J^iM\over J^{i+1}M}=0$ and $_{R/J}R/J$ is semisimple, we get a semisimple module $_{R/J}({J^iM\over J^{i+1}M})$. By Theorem \ref{ComSeExThm2}, since $_RM$ (hence $_{R/J}({J^iM\over J^{i+1}M})$, in addition to being semisimple) is noetherian or artinian, $_{R/J}({J^iM\over J^{i+1}M})$ has a composition series. Recall that $R$ and $R/J$ have the same simple modules. Thus, using the composition series for the modules $_{R/J}({J^iM\over J^{i+1}M})$ we can refine the series (\ref{JM-series}) to obtain a composition series for $_RM$ as follows:
{\footnotesize\begin{align}
\textstyle{}_{R/J}({J^iM\over J^{i+1}M})=S_{i1}\oplus S_{i2}\oplus\cdots\oplus S_{it_i}={N_{i1}\over J^{i+1}M}\oplus{N_{i2}\over J^{i+1}M}\oplus\cdots \oplus{N_{it_i}\over J^{i+1}M}={N_{i1}+\cdots+N_{it_i}\over J^{i+1}M},~~N_{i1}\cap N_{i2}\cap\cdots\cap N_{it_i}=J^{i+1}M,\nn
\end{align}}
and so with $M^i_r:=\sum_{k=r}^{t_i}N_{ik}$, for $1\leq r\leq t_i$, we get a refinement of the $i$th step $J^iM\supset J^{i+1}M$ of (\ref{JM-series}) as follows:
\bea
J^iM=M^i_1\supset M^i_2\supset\cdots \supset M^i_{t_i}=N_{it_i}\supset J^{i+1}M,\nn
\eea
where by construction the composition factors are simple, since
{\small\[
\textstyle{M^i_r\over M^i_{r+1}}={\sum_{k=r}^{t_i}N_{ik}\over \sum_{k=r+1}^{t_i}N_{ik}}\cong {\left(\sum_{k=r}^{t_i}N_{ik}\right)/J^{i+1}M\over \left(\sum_{k=r+1}^{t_i}N_{ik}\right)/J^{i+1}M}={\bigoplus_{k=r}^{t_i}N_{ik}/J^{i+1}M\over \bigoplus_{k=r+1}^{t_i}N_{ik}/J^{i+1}M}\cong{N_{ir}\over J^{i+1}M}\cong {}_{R/J}S_{ir}~\eqv~{}_RS_{ir}. \qedhere
\]}
\end{proof}

\subsection{Homological dimension of a ring with a composition series}~\\~
The following remark is simply an explicit reconstruction of Corollary \ref{ComSeExCrl} for a special situation, but the construction holds in general: I.e., if $R$ is a ring and $0\ra A\ra B\ra C\ra 0$ an exact sequence of $R$-modules, then $B$ has a composition series $\iff$ $A,C$ both have composition series, and moreover,
\[
\txt{length}_R(B)=\txt{length}_R(A)+\txt{length}_R(C).
\]
\begin{rmk}
Let $R$ be a ring, {\small $_RI\subset R$} a left ideal, and {\footnotesize $_RR=J_0\supset J_1\supset\cdots\supset J_t=0$} a composition series.

(i) The left ideal $I$ has a composition series
\[
{}_RI=J_0\cap I\supset J_1\cap I\supset\cdots\supset J_t\cap I=0.\nn
\]
Indeed ${J_{i-1}\cap I\over J_i\cap I}$ is simple, since for any $a_{i-1}\in (J_{i-1}\cap I)\backslash(J_i\cap I)$ we have
\[
\textstyle R(a_{i-1}+J_i\cap I)={Ra_{i-1}+J_i\cap I\over J_i\cap I}\sr{(s1)}{=}{(Ra_{i-1}+J_i)\cap I\over J_i\cap I}\sr{(s2)}{=}{J_{i-1}\cap I\over J_i\cap I},
\]
where step (s1) holds by Lemma \ref{SchrRefLmm} and step (s2) holds because $J_i$ is maximal in $J_{i-1}$ and $a_{i-1}\not\in J_i$ (otherwise $a_{i-1}\in J_i\cap I$, since we already have $a_{i-1}\in I$) .

(ii) The cyclic $R$-module ${R\over I}=R(1+I)$ has a composition series
\[
\textstyle {}_RR(1+I)\supset J_1(1+I)\supset\cdots \supset J_t(1+I)=0,~~~~\txt{where}~~J_i(1+I)={J_i+I\over I}.\nn
\]Indeed, {\small ${J_{i-1}(1+I)\over J_i(1+I)}={(J_{i-1}+I)/I\over(J_i+I)/I}\cong{J_{i-1}+I\over J_i+I}$} is simple, since for any $a_{i-1}(1+I)\in J_{i-1}(1+I)\backslash J_i(1+I)$, we have
\[
\textstyle R\big(a_{i-1}(1+I)+J_i(1+I)\big)={Ra_{i-1}(1+I)+J_i(1+I)\over J_i(1+I)}\sr{(s)}{=}{(Ra_{i-1}+J_i)(1+I)\over J_i(1+I)}={J_{i-1}(1+I)\over J_i(1+I)},
\]
where as before step (s) holds because $J_i$ is maximal in $J_{i-1}$ and $a_{i-1}\not\in J_i$ (otherwise $a_{i-1}(1+I)\in J_i(1+I)$).

(iii) Moreover, because of the exact sequence $0\ra I\ra R\ra R/I\ra 0$, we have
\[
\txt{length}(R)=\txt{length}(I)+\txt{length}(R/I).
\]Indeed given a composition series $_R{R\over I}={A_0\over I}\supset {A_1\over I}\supset\cdots\supset {A_t\over I}=0$ for $R/I$ and a composition series $_RI=B_0\supset B_1\supset\cdots\supset B_s=0$ for $I$, we get a composition series for $R$ given by
\[
{}_RR=A_0\supset A_1\supset\cdots\supset A_{t-1}\supset B_0\supset B_1\supset\cdots\supset B_s=0.
\]
\end{rmk}

\begin{thm}
Let $R$ be a ring. If ${}_RR$ has a composition series, then {\small $Lpd(R)=\sup\{Pd(S):{}_RS~\txt{simple}\}$}.
\end{thm}
\begin{proof}
If {\footnotesize $\sup\{Pd(S):{}_RS~\txt{simple}\}=\infty$}, the result is trivially true. So assume {\footnotesize $\sup\{Pd(S):{}_RS~\txt{simple}\}<\infty$}. Let {\footnotesize $n:=\sup\{Pd(S):{}_RS~\txt{simple}\}$}. Recall that {\footnotesize $Lpd(R)=\sup\{Pd(R/I):~{}_RI\subset R~\txt{a left ideal}\}$}. Moreover, if $R$ has a composition series, so do all cyclic R-modules (since they are all of the form ${R\over I}=R(1+I)$ for $_RI\subset R$ (\blue{footnote}\footnote{Recall that if $_RM=Rm$ is a cyclic $R$-module, then with $f:R\ra Rm,~r\mapsto rm$, the first isomorphism theorem gives $M=Rm\cong R/\ker f=R/Ann_R(m)$.})). It thus suffices to show that for every R-module $_RM$ with a composition series (which include the quotient modules $R/I$), we have $Pd(M)\leq n:=\sup\{Pd(S):{}_RS~\txt{simple}\}$.

We will proceed by induction on the length $r$ of the composition series for $M$. If $_RM$ has a composition series of length 1, then $M$ is simple, i.e., $M\cong {R\over L}=R(1+L)$ for some maximal $_RL\subset R$, and so it is clear that $Pd(M)\leq n$. Assume the result holds for all R-modules with composition series of length $<r$. Let $_RM$ have a composition series $M=M_0\supsetneq M_1\supsetneq\cdots\supsetneq M_r=0$. Then we know $M_1,...,M_r$ each have composition series of length $<r$, and so by the induction hypothesis, $Pd(M_i)\leq n$ for $i=1,...,r$. Finally, because of the exact sequence $0\ra M_1\ra M\ra S:={M\over M_1}\ra 0$, we also have $Pd(M)\leq \max\{Pd(M_1),Pd(S)\}\leq n$.
\end{proof}

\begin{crl}
Let $R$ be a ring. If $_RR$ has a composition series and $Lpd(R)\neq 0$ (\blue{footnote}\footnote{That is, $_RR$ is not semisimple (or equivalently in this case, $J(R)\neq 0$).}), then
\bea
Lpd(R)=1+\sup\{Pd(M):{}_RM\subset R~\txt{maximal}\}.\nn
\eea
\end{crl}
\begin{proof}
$Lpd(R)=\sup\{Pd(S):{}_RS~\txt{simple}\}=\sup\{Pd(R/M):{}_RM\subset R~\txt{maximal}\}\sr{(s)}{=}1+\sup\{Pd(M):{}_RM\subset R~\txt{maximal}\}$,
where step (s) holds because ~$Lpd(R)\neq 0$.
\end{proof}

{\flushleft\hrulefill}

\begin{exercise}
Based on the discussion of these notes, consider writing a \emph{fully technical} essay (say in the form of a typical chapter of these notes) on what is known in the mathematics literature as \index{Commutative! Algebra}{``Commutative Algebra''}. If you consider topology to be an essential part of the discussion, then you can of course refer to the discussion on topology from the later chapters of these notes.
\end{exercise}

%% file: parts/AlgebraM/GeomAnaI.tex
\chapter{Geometry and Analysis I: Topology and Application Instances}\label{GeomAnaI}
By \ul{geometry} in general, we will mean the study/classification of objects (in any given category, although we will eventually focus on categories of sets) based on any notions of ``shape'' or ``topology'' that can be successfully defined on them. (\blue{footnote}\footnote{The word ``geometry'' need not be practically sensible in this general abstract or categorical sense, because not every category can be ``exactly'' realized as a category of sets.}). By \ul{analysis} in general, we will mean various tools/methods (both qualitative and quantitative) that can be successfully employed in such a study (i.e., in geometry).

An isomorphism in the category of topological spaces (to be defined below) is called a homeomorphism. More traditionally, \ul{topology} (resp. \ul{geometry}) is the study of homeomorphism-invariant (resp. specific homeomorphism-invariant) properties of spaces using various categorical objects/structures that can be fruitfully associated with those spaces. Such a study can be used to distinguish/classify certain spaces that might be of interest.

For a traditional introduction to general topology, metric spaces, and some applications, see for example \cite[Appendices A,B,C, pp 128-284]{akoforth2020} and the references therein.

\begin{note}[\textcolor{blue}{A new convention for countability}]
Recall that a set $A$ is countable iff it has the same cardinality as $\Natural$, in which case we write $A\approx \Natural$ or $|A|=|\Natural|$. However, for brevity, ``countable'' from this chapter onward will mean ``at most countable'' (i.e., ``finite or countable'') unless it is specified otherwise. That is, by default, a set will now be called ``countable'' if $|A|\leq|\Natural|$.
\end{note}

\section{Topological Systems and Geometry of Spaces}\label{TopSyGeoSec}
\subsection{Basic definitions and preliminary remarks}
Topological systems here (as well as measurable systems later) are meant to be examples of $Sets$-representable systems in a category. The notion of a topological system (as defined below) is related to the notion of a sheaf, to be defined in section \ref{PreshShSec}. More specifically, a sheaf resembles a pseudo-topological system. In order to quickly establish a smooth connection between our new categorical definitions and those of familiar topology, we will begin with an unusually long definition containing footnotes that are likewise unusually long/elaborate.
\begin{dfn*}[\textcolor{blue}{
\index{Subsystem}{Subsystem},
\index{Finite! subsystem}{Finite subsystem},
\index{Sub(co)limit}{Sub(co)limit},
\index{Finite! sub(co)limit}{Finite sub(co)limit},
\index{Monic! (co)limit-closed system}{Monic (co)limit-closed system},
\index{Finite! monic (co)limit-closed system}{Finite monic (co)limit-closed system},
\index{Topology! (Topological system) in a category}{Topology (Topological system) in a category},
\index{Pseudo-topology (Pseudo-topological system)}{Pseudo-topology (Pseudo-topological system)},
\index{Generated! Topology}{Generated topology},
\index{Generator of! (subbase of) a topology}{Generator (Subbase) of a topology},
\index{Base for a topology}{Base for a topology},
\index{Topology! on a set}{Topology on a set},
\index{Topological! space}{Topological space},
\index{Open! sets}{Open sets},
\index{Closed sets}{Closed sets},
\index{Neighborhood (nbd)}{Neighborhood (nbd)},
\index{Open! neighborhood (open nbd)}{Open neighborhood (open nbd)},
\index{Cover}{Cover},
\index{Open! cover}{Open cover},
\index{Compact space}{Compact space},
\index{Connected! space}{Connected space},
\index{Continuous map}{Continuous map},
\index{Open! map}{Open map},
\index{Closed map}{Closed map},
\index{Category of! topological spaces (Top)}{Category of topological spaces (Top)},
\index{Homeomorphism}{Homeomorphism},
\index{Homeomorphic spaces}{Homeomorphic spaces},
\index{Imbedding (in Top)}{Imbedding},
\index{Geometry of a space}{Geometry of a space},
\index{Morphism of! geometries (Geometric invariant)}{Morphism of geometries (Geometric invariant)},
\index{Geometric! property/invariant}{Geometric property/invariant},
\index{Topological! property/invariant}{Topological property/invariant},
\index{Localized! space}{Localized space},
\index{Localized! geometry}{Localized geometry},
\index{Localized! geometric invariant}{Localized geometric invariant},
\index{Classical! geometric invariant}{Classical geometric invariant}}]\label{TopGeomDfn1}

Let $\I,\C$ be categories and $S:\I\ra\C$ a system. If $\J\subset\I$ is a subcategory, then $S|_{\J}:\J\ra\C$ is a \ul{subsystem} of $S$. (\blue{footnote}\footnote{Note that this (map restriction)-based notion of a subsystem (although related) is different from the notion of a subfunctor introduced earlier in a fixed category of functors $\D^\C$. For greater generality (if necessary), we can also say that given systems $S:\I\ra\C$, $S':\I'\ra\C$, $S$ is a subsystem of $S'$, written $S\subset S'$, if (i) $\I\subset\I'$ is a subcategory and (ii) every object of $S$ is a subobject of an object of $S'$, i.e., for each $i\in\Ob\I$, $S_i\subset S_{i'}'$ is a subobject for some $i'\in\Ob\I'$.}). A subsystem $S|_{\J}:\J\subset\I\ra\C$ is a \ul{finite subsystem} if $\Ob\J$ is a finite set. A \ul{sublimit} (resp. \ul{subcolimit}) of $S$ is the limit (resp. colimit) of a subsystem of $S$. A \ul{finite sub(co)limit} of $S$ is the (co)limit of a finite subsystem.

A system $S:\I\ra\C$ is \ul{monic limit-closed} (resp. \ul{monic colimit-closed}) if every monic sublimit $\varprojlim_{mo} S|_{\J}$ (resp. monic subcolimit $\varinjlim_{mo} S|_{\J}$) is an object of $S$, i.e., $\varprojlim_{mo} S|_{\J}=S(i_\J)$ (resp. $\varinjlim_{mo} S|_{\J}=S(i_\J)$) for some $i_\J\in\Ob\I$. Similarly, a system $S:\I\ra\C$ is \ul{finite monic limit-closed} if every finite monic sublimit $\varprojlim_{mo} S|_{\J}$ (resp. \ul{finite monic colimit-closed} if every finite monic subcolimit $\varinjlim_{mo} S|_{\J}$) is an object of $S$.

A system $T:\I\ra\C$ is a \ul{topology} (or \ul{topological system}) in $\C$ if $T$ is both (i) finite monic limit-closed and (ii) monic colimit-closed. (\blue{footnote}\footnote{Equivalently (up to a contravariant equivalence of categories $\C\ra\C^{op}$), a system $T:\I\ra\C$ is a \ul{topology} (or \ul{topological system}) in $\C$ if $T$ is both (i) finite epic colimit-closed and (ii) epic limit-closed.}). By removing ``monic'', a system $T:\I\ra\C$ is a \ul{pseudo-topology} (or \ul{pseudo-topological system}) in $\C$ if $T$ is both (i) finite limit-closed and (ii) colimit-closed. (\blue{footnote}\footnote{Equivalently (up to a contravariant equivalence of categories $\C\ra\C^{op}$), a system $T:\I\ra\C$ is a \ul{pseudo-topology} (or \ul{pseudo-topological system}) in $\C$ if $T$ is both (i) finite colimit-closed and (ii) limit-closed.}).

If $T:\I\ra\C$ is a topology, $\I$ is directed, and $\{T(i)\}_{i\in\Ob\I}$ has a unique maximal object $M=T(i)$ for some $i$, then we may call $T$ a \ul{topology on the object} $M\in\Ob\C$. Given a system $S:\I\ra\C$, the \ul{topology generated} by $S$ is the smallest topology $\langle S\rangle:\J\ra\C$ containing $S$ as a subsystem, i.e., such that $\I\subset\J$ is a subcategory and $S=\langle S\rangle|_\I$. Given a topology $T:\J\ra\C$, a system $S:\I\ra\C$ is a \ul{generator (or subbase)} of $T$ if $T=\langle S\rangle$ (where it is also understood that $\I\subset \J$ is a subcategory). Given a topology $T:\J\ra\C$, a system $S:\I\ra\C$ is a \ul{base} for $T$ if (i) $T=\langle S\rangle$ and (ii) every object of $T$ is a monic subcolimit of $S$ (i.e., a monic colimit of objects of $S$), i.e., for each $j\in\Ob\J$ we have $T_j=\varinjlim_{mo} S_{\al_j(i)}$ for some map $\al_j:\I\ra\J$.

Let $X$ be a set and $\X=\Big(\P(X),Hom_\X\big(\P(X),\P(X)\big),\circ\Big)\subset Sets$ the full subcategory whose objects $\Ob\X:=\P(X)$ are the subsets of $X$, and $Hom_\X(A,B):=Hom_{Sets}(A,B)$ maps from $A$ to $B$, for any subsets $A,B\subset X$. (\blue{footnote}\footnote{If necessary, we may also consider $X$ to be an object in any category $\C$ with an initial object $I\in\Ob\C$ (replacing $\emptyset$ in $Sets$), and let $\X=\Big(\P(X),Hom_\X\big(\P(X),\P(X)\big),\circ\Big)\subset \C$ be the full subcategory whose objects $\Ob\X:=\P(X)$ are the subobjects of $X$, and $Hom_\X(A,B):=Hom_\C(A,B)$, for any subobjects $A,B\subset X$.}). Then a topology $S:\I\ra\X$ such that both $\emptyset$ and $X$ are objects of $S$ is called a \ul{topology on $X$}, making $X=(X,S)$ a \ul{topological space} (or just a \ul{space}) with the sets $\{S(i)\subset X\}_{i\in\Ob\I}$  called \ul{open sets} in (or \ul{open subsets} of) the space $X$ and their complements $\{S(i)^c:=X-S(i)\}_{i\in\Ob\I}$ called \ul{closed sets} in (or \ul{closed subsets of}) the space $X$. (\blue{footnote}\footnote{Traditionally, we specify a topology for a space $X$ in terms of its \emph{open sets}, namely, a collection of subsets $\O\subset\P(X)$ such that (i) $\emptyset,X\in\O$, (ii) $O,O'\in\O$ $\Ra$ $O\cap O'\in\O$, and (iii) $O_\al\in\O$ $\Ra$ $\bigcup_\al O_\al\in\O$. We then define \emph{closed sets} in the space $X$ to be the complements $C:=O^c$ of the open sets $O$. However, by DeMorgan's laws, the closed sets of a topological space also completely/uniquely determine the topology (which was given in terms of open sets). Consequently, we can also specify a topology for a space $X$ in terms of its \emph{closed sets}, namely, a collection of subsets $\C\subset\P(X)$ such that (i) $\emptyset,X\in\C$, (ii) $C,C'\in\C$ $\Ra$ $C\cup C'\in\C$, and (iii) $C_\al\in\C$ $\Ra$ $\bigcap_\al C_\al\in\C$. We then define \emph{open sets} in the space $X$ to be the complements $O:=C^c$ of the closed sets $C$.}). As already indicated above, given a collection of sets $\C\subset\P(X)$, the topology $\langle\C\rangle$ on $X$ generated by $\C$ (making $\C$ a \ul{subbase} for $\langle\C\rangle$) is the smallest topology on $X$ containing $\C$. Given a topology $\T$ on $X$ and a collection of sets $\B\subset\P(X)$, $\B$ is a \ul{base} for $\T$ if (i) $\B$ generate $\T$, i.e., $\T=\langle\B\rangle$ and (ii) every member of $\T$ is a union of members of $\B$.

Let $X$ be a space. Given sets $A,N\subset X$, the set $N$ is a \ul{neighborhood} (nbd) of $A$ if there is an open set $O\subset X$ such that $A\subset O\subset N$. A neighborhood is called an \ul{open neighborhood} if it is itself an open set. It is often enough to assume that neighborhoods are open, and so we will do this, i.e. (unless stated otherwise) ``neighborhood'' will mean ``open neighborhood''. A family of sets $\{A_\al\}_{\al\in \Gamma}\subset\P(X)$ is a \ul{cover} of $X$ if $X\subset \bigcup_\al A_\al$. A cover consisting of open sets is called an \ul{open cover}. A space $X$ is \ul{compact} if every open cover of $X$ has a finite subcover (where a cover is finite if it is finite as a set, and a subcover of a cover is a subset of the cover as a set). A space $X$ is \ul{connected} if it cannot be written as a disjoint union of two nonempty open sets.

Given spaces $(X,S_X:\I\ra\X)$ and $(Y,S_Y:\J\ra\Y)$, and a map $f:X\ra Y$, consider any process
\bea
P_f:\Y\ra\X,~B\sr{v}{\ral}B'~\mapsto~f^{-1}(B)\sr{P_f(v)}{\ral}f^{-1}(B').\nn
\eea
We say the map $f:X\ra Y$ is a \ul{continuous map} if the process $P_f:\Y\ra\X$ takes objects of $S_Y$ to objects of $S_X$, i.e., for each $j\in\Ob\J$, ~$P_f(S_Y(j)):=f^{-1}(S_Y(j))=S_X(i)$ for some $i\in\Ob\I$. In other words, a map of spaces $f:X\ra Y$ is a \ul{continuous map} if for every open set $B\subset Y$, the primage set $f^{-1}(B)\subset X$ is open. A map of spaces $f:X\ra Y$ is an \ul{open map} (resp. a \ul{closed map}) if for every open (resp. closed) set $A\subset X$, the image set $f(A)\subset Y$ is open (resp. closed).

The \ul{category of spaces}, denoted by $\Top$, is the category whose objects are spaces and its morphisms $Hom_{\Top}(X,Y)$ are continuous maps $f:X\ra Y$. An isomorphism $X\sr{f}{\ral}Y$ in $\Top$ is called a \ul{homeomorphism} (making $X,Y$ \ul{homeomorphic spaces}). A monomorphism (or equivalently, a subobject) $f:X\hookrightarrow Y$ (i.e., an injective continuous map) in the category Top is called an \ul{imbedding} of $X$ into $Y$ if it is a homeomorphism onto its image (i.e., an imbedding is an injective continuous map such that $f:X\ra f(X)$ is a homeomorphism). (\magenta{footnote}\footnote{Thus, if an imbedding $f:X\ra Y$ exists, then we can view $X$ as a subspace of $Y$ (see Definition \ref{TopSubSpDfn}). Note that in Top, an imbedding (subspace) is a monomorphism (subobject) but a monomorphism (subobject) is not an imbedding (subspace). \magenta{In order for a subobject (in the usual sense) to be automatically homeomorphic to a subspace, we need to restrict to the subcategory $CatTop\subset Top$ whose objects are the same as those of Top (i.e., all topological spaces) but whose morphisms are \ul{categorically continuous maps} (or \ul{cat-continuous maps}), where a map between topological spaces $f:X\ra Y$ is cat-continuous if $B\subset Y$ is open $\iff$ $f^{-1}(B)\subset X$ is open (i.e., $f$ is continuous and preimages $f^{-1}(B)\subset X$ of non-open sets $B\subset Y$ are non-open sets). A surjective cat-continuous map is called a \ul{quotient map} (see Definition \ref{QuoTopDef}). A subobject in CatTop (i.e., an injective cat-continuous map) $f:X\hookrightarrow Y$ is an injective open continuous map and hence an imbedding: Indeed if $A\subset X$ is open, then so is $f(A)\subset Y$, since $f^{-1}(f(A))=A\subset X$ is open.\\
Also, observe that a surjective open continuous map $f:X\ra Y$ is cat-continuous (hence a quotient map: Definition \ref{QuoTopDef}), since for any $B\subset Y$ we have ~$B\cap f(X)=f(f^{-1}(B\cap f(X)))=f(f^{-1}(B))$.
}})

Let $X$ be a space, and $\H_X\subset\Top$ the subcategory whose objects $\Ob\H_X:=\{Y\in\Ob\Top:~Y\cong X\}$ are spaces that are homeomorphic to $X$, and whose morphisms
{\small\bea
\textstyle \Mor\H_X:=\bigcup\limits_{Y,Y'\in\Ob\H_X}Hom_{\H_X}(Y,Y'),~~~~ Hom_{\H_X}(Y,Y'):=\big\{\txt{isomorphisms in}~Hom_{\Top}(Y,Y')\big\},\nn
\eea}are homeomorphisms of $X$ in $\Top$. A \ul{geometry of $X$} is a system in $\H_X$, i.e., a functor $S:\J\ra\H_X$, (\blue{footnote}\footnote{We may also consider a geometry of $X$ to be any system $S:\J\ra\Top$ with a subsystem in $\H_X$ (i.e., $S|_{\J'}:\J'\ra\H_X$ for some subcategory $\J'\subset\J$).}), for a category $\J$ (the detailed structure of which depends on the intended application). Let $Geom(X)$ denote the category whose objects are the geometries $\bigcup_\J\H_X^\J\subset\bigcup_{\J}\Top^\J$ of $X$, and whose morphisms (i.e., \ul{morphisms of geometries}) $Hom_{Geom(X)}(S,S')$ are morphisms of systems $\eta:S\ra S'$ (when defined, in the sense $\dom S=\dom S'$). A \ul{geometric invariant} of $X$ is a morphism of geometries of $X$.

For example, if $\J:=2:=(\{0,1\},\{0\ra 1\},\circ)$ as a poset, then examples of geometries of $X$ are homeomorphisms $X\sr{h}{\ral}Y$ and morphisms between them (i.e., geometric invariants of $X$) $Hom_{Geom(X)}(h,h')\subset \Top^2$ are commutative diagrams of the form
\adjustbox{scale=0.7}{%
\bt
X\ar[d,"\al"]\ar[rr,"h"] && Y\ar[d,"\beta"]\\
X\ar[rr,"h'"] && Y'
\et}.

In familiar language, a \ul{geometric property/invariant} (resp. \ul{topological property/invariant}) of a space is a property/statement about the space that is unchanged under specific homeomorphisms (resp. all homeomorphisms) of $X$.
Examples of geometries and their geometric invariants can be obtained through the localization process.

A space $X$ is a \ul{localized space} if $\H_X$ is equivalent to a traditional localization in the sense that there exists a subcategory $\K\subset\Top$ and a localizing class of morphism (i.e., continuous maps) $Q\subset\Mor\K$ such that we have an equivalence of categories $\H_X\cong_w\K[Q^{-1}]$. A \ul{localized geometry} (resp. \ul{localized geometric invariant}) is a geometry (resp. geometric invariant) of a localized space.

Let $\I,\C$ be categories. A \ul{classical $(\I,\C)$-geometric invariant} (or a classical geometric invariant with grading in $\I$ and coefficients in $\C$) is a functor $S:\J\subset\Top\ra Gr(\C^\I)$, i.e., a system in $Gr(\C^\I)$ of the form $S:\J\ra Gr(\C^\I)$, for some subcategory $\J\subset\Top$. The category of ($\I,\C$)-geometric invariants is the category $GeomInv(\I,\C)$ whose objects are all functors of the type $\bigcup_{\J\subset\Top}Gr(\C^\I)^\J\subset \bigcup_{\J\subset\Top}Gr(\C^{\I\times\J})$, and whose morphisms $Hom_{GeomInv(\I,\C)}(S,S')$ are morphisms of systems $\eta:S\ra S'$ (when defined, in the sense $\dom S=\dom S'$).

Given a space $X$ (and categories $\I,\C$ as above), the above situation becomes local to $X$ if we consider only systems $S:\J\subset\Top\ra Gr(\C^\I)$ for which $\J\subset \H_X\subset \Top$. In this case, we speak of a ``\ul{classical $(\I,\C)$-geometric invariant of $X$}'' (or a ``classical geometric invariant of $X$ with grading in $\I$ and coefficients in $\C$''). Accordingly, the category $GeomInv(\I,\C)$ becomes the category $GeomInv(X;\I,\C)$ of $(\I,\C)$-geometric invariants of $X$, defined to be the category whose objects are all functors of the type $\bigcup_{\J\subset\H_X}Gr(\C^\I)^\J\subset \bigcup_{\J\subset\H_X}Gr(\C^{\I\times\J})$, and whose morphisms $Hom_{GeomInv(X;\I,\C)}(S,S')$ are morphisms of systems $\eta:S\ra S'$.
\end{dfn*}

For a typical application in algebraic topology (section \ref{AlgTopSec}), $\I\subset\Integer^n$ for some integer $n\geq1$, and $\C\subset\A$ for an abelian category such as $R$-mod for a ring $R$. An example of a classical $(\Integer,Ab)$-geometric invariant is the singular chain complex $C:Top\ra Ab^\Integer$ (see Definition \ref{SingCCdef1}). For a typical application in algebraic geometry (section \ref{AlgGeomSec}), $\I=\I_X$ is the poset-category of open subsets of some topological space $X$ ordered by inclusion $\subset$, and $\C\subset Rings$ is a subcategory of Rings.

\begin{lmm}[\textcolor{blue}{\index{Open! set criterion}{Open set criterion}}]\label{OpenCrit}
Let $X$ be a space. A set $A\subset X$ is open $\iff$ for every point $x\in A$, there is an open set $O$ such that $x\in O\subset A$.
\end{lmm}
\begin{proof}
If $A$ is open, then for every $x\in A$, we have $x\in O:=A\subset A$. Conversely, if for every $x\in A$, we have $x\in O_x\subset A$ for some open set $O_x$, then $A=\bigcup_{x\in A}O_x$ is open as a union of open sets.
\end{proof}

\begin{rmk}[\textcolor{blue}{\index{Base from a subbase}{Base from a subbase}}]
Let $X$ be a set and $\C\subset\P(X)$. Then it is clear (from the definitions) that the topology $\langle\C\rangle$ on $X$ generated by $\C$ has the base
\[
\B=\big\{\txt{finite intersections of sets from $\C\cup\{\emptyset,X\}$}\big\}.
\]
That is, sets in $\langle\C\rangle$ are unions of sets from $\B$.
\end{rmk}

\begin{lmm}[\blue{\index{Base criterion}{Base criterion}}]
Let $X=(X,\T)$ be a space. A collection of open sets $\B\subset\T$ is a base for $\T$ $\iff$ for any point $x_0\in X$ and any open neighborhood $U\in\T$ of $x_0$, we have $x_0\in B\subset U$ for some $B\in\B$.
\end{lmm}
\begin{proof}
If $\B$ is a base for $\T$ and $U$ is an open neighborhood of $x_0\in X$, then with $U=\bigcup_\al B_\al$ for $B_\al\in\B$, we have $x_0\in B_\al\subset U$ for some $\al$. Conversely, if every open neighborhood $U$ of every point $x_0\in X$ satisfies $x_0\in B_{x_0,U}\subset U$ for some $B_{x_0,U}\in\B$, then every open set $O\in\T$ satisfies $O=\bigcup_{x\in O}B_{x,O}$, and so $\B$ is a base for $\T$.
\end{proof}

\begin{dfn}[\textcolor{blue}{\index{Local! topology}{Local topology}, \index{Neighborhood! base}{Neighborhood base}}] Let $X=(X,\T)$ be a space and $A\subset X$ a subset. All open sets containing $A$ together with the empty set, $\T_A=\{\emptyset\}\cup\{O\in\T:A\subset O\}$, form a topology on $X$ called the \ul{local topology} at $A$. Any base $\B_A$ for $\T_A$ is called a \ul{neighborhood base} (at $A\subset X$). Note that $\T=\bigcup_{A\subset X}\T_A$, and $\B=\bigcup_{A\subset X}\B_A$ is a base for $\T$.
\end{dfn}

\begin{dfn}[\textcolor{blue}{\index{First countable space}{First countable space}, \index{Second countable space}{Second countable space}}] Let $X=(X,\T)$ be a space. $X$ is \ul{first countable} if every point $x_0\in X$ has a countable neighborhood base. $X$ is \ul{second countable} if the topology $\T$ has a countable base.
\end{dfn}

\begin{examples}[\textcolor{blue}{Discrete topology, \index{Discrete space}{Discrete space}, Indiscrete topology, \index{Indiscrete space}{Indiscrete space}}]
Let $X=(X,\T)$ be a space. The topology $\T$ is \ul{discrete} (making {\small $X$} a discrete space) if every point $x_0\in X$ is an open set $\{x_0\}\subset X$ (equivalently, every subset of $X$ is an open set, i.e., {\small $\T=\P(X)$} ). The topology $\T$ is \ul{indiscrete} (making {\small $X$} an indiscrete space) if $\T=\{\emptyset,X\}$.
\end{examples}

\begin{rmk}[\textcolor{blue}{\index{Order! topology}{Order topology}, \index{Open! intervals}{Open intervals}}]
Let $L=(L,\leq)$ be any linearly ordered set, such as $\Real$. Unless specified otherwise, we will assume $L$ is a space with the \ul{order topology}, defined as the topology with base (whose elements are called \ul{open intervals}) given by the following collection:
\bea
\B:=\big\{(a,b):=\{x\in L:a<x<b\}\subset L~|~\txt{for all}~~a,b\in L\big\}\cup\{L\}.\nn
\eea
\end{rmk}

\begin{rmk}[\textcolor{blue}{\index{Induced! topologies}{Topologies induced} via a map}]
Let $(X,\T_X),(Y,\T_Y)$ be spaces, and $f:X\ra Y$ a map. Then through $f$, (i) $\T_Y$ induces a topology on $X$ given by
\bea
f^{-1}(\T_Y):=\{f^{-1}(B):B\in\T_Y\},~~~~\txt{(where $f$ is a continuous map iff ~$f^{-1}(\T_Y)\subset\T_X$)}\nn
\eea
and (ii) $\T_X$ induces a topology on $Y$ given by
\bea
f^\ast(\T_X):=\{B\subset Y: f^{-1}(B)\in\T_X\},\nn
\eea
where $f$ is a quotient map (Definition \ref{QuoTopDef}) iff (i) $f$ is surjective and (ii) $f^\ast(\T_X)=\T_Y$.

Note that if $f$ is surjective, then ~$f^\ast(\T_X)=f\left(\T_X\cap f^{-1}(\P(Y))\right)$,~ where
\[
f\left(\T_X\cap f^{-1}(\P(Y))\right):=\{f(f^{-1}(B)):B\subset Y,~f^{-1}(B)\in\T_X\}=\{B\cap f(X):B\subset Y,~f^{-1}(B)\in\T_X\}.
\]
\end{rmk}

\subsection{Preliminary concepts/results and the protosheaf concept}
\begin{dfn}[\textcolor{blue}{
\index{Subspace (relative) topology}{Subspace (relative) topology},
\index{Subspace}{Subspace},
\index{Relatively! open set}{Relatively open set},
\index{Relatively! closed set}{Relatively closed set},
\index{Space-subspace property}{Space-subspace property}}]\label{TopSubSpDfn}
Let $X=(X,\T)$ be a space. Every subset $A\subset X$ inherits the topology $A\cap\T :=\{A\cap O~|~O\in\T\}$ on $A$, called \ul{subspace topology} on $A$, which makes $A=(A,A\cap\T)$ a \ul{subspace} of $X$. The subspace-complements $\{O^{c_A}:=A-(A\cap O)=A\cap O^c~|~O\in\T\}$ are the \ul{closed sets} of the subspace topology, where $O^c=X-O$ is a closed set in $X$. (\blue{footnote}\footnote{If $X$ is a space, a \emph{subobject} of $X$ in the category Top is an \emph{injective continuous map} $f:Z\ra X$, but $Z$ does not have to be (homeomorphic to) a \emph{subspace} of $X$, since $f$ does not have to be an \emph{imbedding}.}). The subspace topology is also called the \ul{relative topology}, and accordingly, open (resp. closed) sets in the subspace topology are said to be \ul{relatively open} (resp. \ul{relatively closed}).

Let $P$ be a property of a topology (in general) and $X$ a space. We say the space $X$ has property $P$ (or $X$ is a \ul{$P$-space}) if the topology of $X$ has property $P$. (Irrespective of whether $X$ is a $P$-space or not) we also say a subset $A\subset X$ has the property $P$ (or $A$ is a \ul{$P$-subset} or \ul{$P$-subspace} of $X$) if $A$ has the property $P$ as a subspace of $X$ (i.e., if the space $(A,A\cap\T)$ is a $P$-space). In particular, a set $A\subset X$ is compact (resp. connected) in $X$ if it is compact (resp. connected) as a subspace.
\end{dfn}

\begin{lmm}
\begin{enumerate*}[leftmargin=0.7cm]
\item A relatively open subset of an open set is itself open (as an intersection of two open sets).
\item Similarly, a relatively closed subset of a closed set is itself closed (as an intersection of closed sets).
\item A continues map takes a compact (resp. connected) set to a compact (resp. connected) set.
\end{enumerate*}
\end{lmm}
\begin{proof}
Since (1) and (2) are clear, we will prove (3). Let $f:X\ra Y$ be continuous. If $A\subset X$ is compact and $f(A)\subset \bigcup_\al B_\al$ (i.e., $A\subset f^{-1}(\bigcup_\al B_\al)=\bigcup_\al f^{-1}(B_\al)$) for open sets $B_\al\subset Y$, then a finite subcover $\{f^{-1}(B_{\al_i})\}_{i=1}^n$ of $\{f^{-1}(B_\al)\}$ gives a finite subcover $\{B_{\al_i}\}_{i=1}^n$ of $f(A)$. Also, if $A\subset X$ is connected and $f(A)=B_1\cup B_2$ in $Y$, then with $A_j:=f^{-1}(B_j)$, i.e., $f(A_j)=f(f^{-1}(B_j))=B_j$, we have
\[
A=A\cap f^{-1}(f(A))=(A\cap f^{-1}(B_1))\cup(A\cap f^{-1}(B_2))=A_1\cup A_2,
\]
and so if $B_j=f(A)\cap O_j$, i.e., $A_j=A\cap f^{-1}(O_j)$, for open sets $O_j\subset Y$, then either $B_1=\emptyset$ or $B_2=\emptyset$.
\end{proof}

\begin{dfn}[\textcolor{blue}{\index{Quotient! map}{Quotient map}, \index{Quotient! topology}{Quotient topology}, \index{Quotient! space}{Quotient space}}]\label{QuoTopDef}
A map of spaces $f:X\ra Y$ is called a \ul{quotient map} (making $Y$ a \ul{quotient space} of $X$, and the topology of $Y$ a \ul{quotient topology}) if (i) it is surjective and (ii) $B\subset Y$ is open $\iff$ $f^{-1}(B)$ is open in $X$. (\blue{footnote}\footnote{If $X$ is a space, a \emph{quotient object} of $X$ in the category Top is a \emph{surjective continuous map} $f:X\ra Y$, but $Y$ does not have to be a \emph{quotient space} of $X$, since $f$ does not have to be a \emph{quotient map}.}).

(Note that a quotient map is continuous, but does not have to be an open map).
\end{dfn}

\begin{examples}
(1) Given a space $(X,\T)$ and any equivalence relation $\sim$ on $X$ (as a set), the map $q:X\ra X_\sim$, $x\mapsto x_\sim$ onto the set of equivalence classes {\small$X_\sim={X\over\sim}=\left\{x_\sim:=\{x'\in X:x\sim x'\}~|~x\in X\right\}$} is a quotient map if $X_\sim$ is given the topology {\small$\T_\sim:=\big\{B\subset X_\sim:q^{-1}(B)\in\T\big\}$}.\\
(2) Given a space $X$, any map of sets $f:X\ra A$ induces an equivalence relation $\sim_f$ on $X$, where $x\sim_f x'$ iff $f(x)=f(x')$, with equivalence classes
\bit
\item[] {\small$X_{\sim_f}={X\over\sim_f}=\left\{x_{\sim_f}:=f^{-1}(f(x))~|~x\in X\right\}=\{f^{-1}(a):a\in f(X)\}$}.
\eit
(3) Given a space $X$ and any subspace $A\subset X$, we can define an equivalence relation $\sim$ on $X$ as follows: $x\sim x'$ iff $x=x'$ or $x,x'\in A$. The space $X_\sim$ is denoted by ${X\over A}$ or $X/A$.
\end{examples}

\begin{prp}[\textcolor{blue}{\index{Correspondence theorem for! topological spaces}{Correspondence theorem}}]
Let $X$ be a space and $A\subset X$. Then we have a bijective correspondence as follows (\blue{footnote}\footnote{By its proof, this is really a correspondence of subsets that does not depend on topology.}):
{\small\[
\P(X\pm A):=\left\{
  \begin{array}{l}
    \txt{subspaces of $X$ containing} \\
    \txt{$A$ or disjoint from $A$}
  \end{array}
\right\}\longleftrightarrow \big\{\txt{subspaces of $X/A$}\big\}=:\P(X/A).
\]}
\end{prp}
\begin{proof}
Let $q:X\ra X/A$ be the quotient map. Consider the clearly well-defined map
\[
\phi:\P(X\pm A)\subset\P(X)\ra \P(X/A),~~U\mapsto q(U)
\]
from subspaces of $X$ that exclude or contain $A$ to subspaces of $X/A$, and the reverse map
\[
\psi:\P(X/A)\ra \P(X\pm A)\subset\P(X),~~V\mapsto q^{-1}(V).
\]
If $V\subset X/A$, then either (i) $q(A)\in V$, in which case $A\subset q^{-1}(V)$, or (ii) $q(A)\not\in V$, in which case $q^{-1}(V)\cap A=\emptyset$. This shows $\psi$ is also well-defined, since $q$ is surjective and
\[
q^{-1}(V)=q^{-1}(V')~~~\Ra~~~V\cap q(X)=q(q^{-1}(V))=q(q^{-1}(V'))=V'\cap q(X).
\]
We have $\phi\psi=id_{\P(X/A)}$, since $\phi\psi(V)=qq^{-1}(V)=V\cap q(X)=V$. Also, $\psi\phi=id_{\P(X\pm A)}$, because for any $U\in \P(X\pm A)$,
\begin{align}
&\psi\phi(U)=q^{-1}(q(U))=q^{-1}\big(q(U-A)\cup q(U\cap A)\big)=q^{-1}\big(q(U-A)\big)\cup q^{-1}\big(q(U\cap A)\big)\nn\\
&~~~~\sr{(s1)}{=}(U-A)\cup q^{-1}\big(q(U\cap A)\big)\sr{(s2)}{=}U,\nn
\end{align}
where step (s1) holds because $q|_{U-A}$ is injective and step (s2) holds because either $A\subset U$ or $U\cap A=\emptyset$.
\end{proof}

\begin{thm}[\textcolor{blue}{\index{Universal! property of quotient maps}{Universal property of quotient maps}}]\label{UPQM}
Let $X$ be a space and $\sim$ an equivalence relation on $X$. For any continuous map $f:X\ra Y$ such that $x\sim x'$ $\Ra$ $f(x)=f(x')$, there exists a unique continuous map $\ol{f}:X_\sim\ra Y$ such that the following diagram commutes (i.e., $f$ factors through $q$ as $f=\ol{f}\circ q$).
\bc\bt
X\ar[d,"f"']\ar[rr,"q"]&& X_\sim\ar[dll,dashed,"\ol{f}"]\\
Y
\et\ec
\end{thm}
\begin{proof}
Let $X_\sim=\{x_\sim:x\in X\}$, where $x_\sim=\{x'\in X:x\sim x'\}$. \ul{\emph{Existence}}: The hypotheses imply $\ol{f}:X_\sim\ra Y$, $x_\sim\mapsto f(x)$ is well defined, since
\bea
x_\sim=x'_\sim~~\iff~~x\sim x'~~\Ra~~\ol{f}(x_\sim)=f(x)=f(x')=\ol{f}(x'_\sim),\nn
\eea
and moreover, $\ol{f}\circ q(x)=\ol{f}(x_\sim)=f(x)$, and so $\ol{f}\circ q=f$. \ul{\emph{Uniqueness}}: If $g:{X\over\sim}\ra Y$ is any map (continuous or not) such that $g\circ q=f$, then because $q$ is surjective, $f=g\circ q=\ol{f}\circ q$ ~$\Ra$~ $g=\ol{f}$. \ul{\emph{Continuity}}: If $B\subset Y$ is open, then because $f^{-1}(B)=(\ol{f}\circ q)^{-1}(B)=q^{-1}\left(\ol{f}^{-1}(B)\right)$, $f$ is continuous, and $q$ is a quotient map, we see that $\ol{f}^{-1}(B)$ is open. Hence, $\ol{f}$ is continuous.
\end{proof}

\begin{crl}[\textcolor{blue}{\index{Continuous map! injective companion of}{Injective companion of a continuous map}}]
Every continuous map of spaces $f:X\ra Y$ induces a unique injective continuous map $\ol{f}:X_{\sim_f}\ra Y$ such that $f=\ol{f}\circ q:X\sr{q}{\ral}X_{\sim_f}\sr{\ol{f}}{\ral}Y$, where $q:X\ra X_{\sim_f}$ is the quotient map, with ~$X_{\sim_f}:=\left\{f^{-1}(f(y)):y\in f(X)\right\}$.
\end{crl}

\begin{crl}[\textcolor{blue}{Continuous maps on a quotient space}]
Let $X$ be a space, $\sim$ an equivalence relation on $X$, and $Y:=X_\sim$ the associated quotient space of $X$. Then continuous maps $\C(Y,Z)$ are precisely continuous maps $\C(X,Z)$ that are constant on equivalence classes $x_\sim:=\{x'\in X~|~x'\sim x\}$ for all $x\in X$.
\end{crl}
\begin{proof}
Let $q:X\ra Y$ be the quotient map. By the universal property of quotient maps, any continuous map $f:X\ra Z$ that is constant on equivalence classes (i.e., $x\sim x'$ $\Ra$ $f(x)=f(x')$~) gives a unique continuous map $\wt{f}:Y\ra Z$ such that $f=\ol{f}\circ q:X\sr{q}{\ral}Y\sr{\ol{f}}{\ral}Z$. Conversely, any continuous map $h:Y\ra Z$ also gives the unique continuous map $H=h\circ q:X\sr{q}{\ral}Y\sr{h}{\ral}Z$, which is by construction constant on equivalence classes (i.e., $x\sim x'$ $\Ra$ $H(x)=H(x')$).
\end{proof}

\begin{dfn}[\textcolor{blue}{\index{Product topology}{Product topology}, \index{Product space}{Product space}, \index{Open! rectangles}{Open rectangles}}]\label{ProdTop}
Let $\{(X_i,\T_i)\}_{{i\in I}}$ be a family of spaces (where $I$ is a set), and consider the cartesian product {\footnotesize
\bea
\textstyle X:=\prod_{{i\in I}}X_{i}=\{(x_{i})_{{i\in I}}:x_{i}\in X_{i}\}=\left\{I\sr{x}{\ral}\bigcup_i X_{i}~|~x({i})\in X_{i}\right\}.\nn
\eea}
The \ul{product topology} $\T$ on $X$ (making $X=(X,\T)$ a \ul{product space}) is the topology with base given by the following collection (called \ul{open rectangles}):
{\small\begin{align}
\label{ProdTopBasEq}\textstyle\B:=\left\{{\textstyle\prod_{{i\in I}}}U_{i}~\big|~U_{i}\in\T_{i},~U_{i}=X_{i}~\txt{a.e.f.}~i\right\}=\left\{U_F:=\prod_{{i}\in F}U_{i}\times\prod_{i\not\in F}X_{i}~\big|~\txt{finite}~F\subset I,~ U_{i}\in\T_{i}\right\},
\end{align}}
where ~{\small $\prod_{{i\in I}}U_{i}:=\big\{(x_{i})_{{i\in I}}:x_{i}\in U_{i}\big\}$}~ and ~{\small $U_F=\left\{x\in X:x_{i}\in U_{i}\in\T_{i},~\txt{for}~{i}\in F\right\}$}. (\blue{footnote}\footnote{This product space is indeed the product $\prod_i X_i$ in the category Top, i.e., it satisfies the desired universal property (by Corollary \ref{SpaceProdThm}).})
\end{dfn}

\begin{notation}[\textcolor{blue}{Product space}]
Let $\{X_i\}_{i\in I}$ be a family of sets. Whenever the cartesian product $\prod_{{i\in I}}X_{i}$ is a space, we will always assume it is given the product topology, unless it is stated otherwise.
\end{notation}

\begin{lmm}[\textcolor{blue}{\cite[Theorem 19.6, p 117]{munkres}}]\label{ComFunCont}
Let $\{X_{i}\}_{i\in I}$ be a family of spaces, $X:=\prod_iX_i$ the product space, $f:Z\ra X$ a map, and $p_i:X\ra X_i,~x\mapsto x_i$ the $i$th projection. Then the map $f$ is continuous $\iff$ each of the maps $f_{i}:=p_i\circ f:Z\ra X_{i}$ is continuous.
\end{lmm}
\begin{proof}
Observe that $f(z)=\left(f_{i}(z)\right)_{{i\in I}}$ for each $z\in Z$. ($\Ra$): If $f$ is continuous, then because the maps $p_{i}$ are continuous due to $p_{i}^{-1}(U_{i})=U_{i}\times\prod_{j\neq{i}}X_j$, we see that each $f_{i}=p_{i}\circ f$ is continuous. ($\La$): Conversely, if each $f_{i}$ is continuous, then so is $f$, since for each finite set $F\subset I$,
\begin{align}
&\textstyle f^{-1}\left(\prod_{{i}\in F}U_{i}\times\prod_{i\not\in F}X_{i}\right)=\big\{z\in Z:f_{i}(z)\in U_{i}~\txt{for}~{i}\in F,~f_{i}(z)\in X_{i}~\txt{for}~{i}\not\in F\big\}\nn\\
&\textstyle~~~~=\bigcap_{{i}\in F}f_{i}^{-1}(U_{i})~\cap~\bigcap_{{i}\not\in F}f_{i}^{-1}(X_{i})=\bigcap_{{i}\in F}f_{i}^{-1}(U_{i})\nn
\end{align}
is open as a finite intersection of open sets.
\end{proof}

\begin{crl}[\textcolor{blue}{Product of spaces in Top: \index{Universal! property of the product space }{Universal property of the product space}}]\label{SpaceProdThm}
Let $\{X_{i}\}_{i\in I}$ be a family of spaces. The product space $\prod_i X_{i}$ is the product in the category Top.
\end{crl}
\begin{proof}
Let {\small $p_i:\prod_iX_i\ra X_i,~x\mapsto x_i$} be the projections (which we know to be continuous). Given continuous maps $f_{i}:Z\ra X_{i}$, define a map {\small $f:Z\ra\prod_i X_{i}$} by {\small $f(z):=(f_{i}(z))_{i\in I}$}. Then $p_{i}\circ f=f_{i}$, for all ${i}$.
\bea\adjustbox{scale=0.8}{%
\bt
X_{i} & & \\
 & \prod_i X_{i}\ar[ul,two heads,"p_{i}"'] & \\
 & Z\ar[uul,bend left,"f_{i}"]\ar[u,dashed,"f"'] &
\et}~~~~p_{i}\circ f=f_{i}.\nn
\eea
Since each $f_i$ is continuous, so is $f$ (by the preceding lemma).
\end{proof}

\begin{dfn}[\textcolor{blue}{
\index{Disjoint! union space}{Disjoint union space},
\index{Disjoint! union topology}{Disjoint union topology}}]
Let {\small $(X_{i},\T_{i})_{{i}\in I}$} be a family of spaces. The \ul{disjoint union space} {\small $\coprod X_{i}$} of the spaces {\small $X_{i}$} is the disjoint union of sets
\[
\textstyle\coprod_{i\in I}X_i:=\bigcup_{i\in I}X_i\times\{i\}=\bigcup_{i\in I}\big\{(x_i,i):x_i\in X_i\big\}
\]
(along with inclusions map $q_i:X_i\hookrightarrow\coprod_{i\in I}X_i,~x\mapsto(x,i)$, for each $i$, i.e., $X_i\subset\coprod X_i$ for each $i$) as a space with the topology (called \ul{disjoint union topology}) whose base $\B$ is as follows (\blue{footnote}\footnote{This union space is indeed the coproduct $\coprod_i X_{i}$ in the category Top, i.e., it satisfies the desired universal property (by Corollary \ref{SpaceCoprodThm}).}):
\bea
\textstyle\B:=\left\{\coprod_{i\in I}U_i:=\bigcup_{i\in I}(U_i\times\{i\})~\big|~U_i\in\T_i\right\}.\nn
\eea
By construction, the powerset {\small $\P(\coprod_i X_{i}):=\left\{\coprod_i A_{i}:=\bigcup_i(A_i\times\{i\}):A_{i}\subset X_{i}\right\}$} admits componentwise union and intersection rules given by
\bea
\label{DistUnEq2}\textstyle
(\coprod_i A_{i})\cup(\coprod_i A_{i}'):=\coprod_i(A_{i}\cup A_{i}'),~~~~(\coprod_i A_{i})\cap(\coprod_i A_{i}'):=\coprod_i(A_{i}\cap A_{i}').
\eea
\end{dfn}

\begin{crl}[\textcolor{blue}{Coproduct (Direct sum) of spaces in Top: \index{Universal! property of the disjoint union space}{Universal property of the disjoint union space}}]\label{SpaceCoprodThm}
Let $\{X_i\}_{i\in I}$ be spaces. The disjoint union space $\coprod_i X_i$ is the coproduct in the category Top.
\end{crl}
\begin{proof}
Given continuous maps $f_i:X_i\ra Y$, define a map $f:\coprod_i X_i\ra Y$ by $f((x,i)):=f_i(x)$. With the inclusion $q_i:X_i\hookrightarrow\coprod_i X_i$, $x\mapsto (x,i)$, it is clear that $f\circ q_i:=f|_{X_i}=f_i$.
\bea\adjustbox{scale=0.8}{%
\bt
X_i\ar[ddr,bend right,"f_i"']\ar[dr,hook,"q_i"] & & \\
 & \coprod_i X_i\ar[d,dashed,"f"] & \\
 & Y &
\et}~~~~f_i=f\circ q_i.\nn
\eea
Each inclusion $q_i$ is continuous because {\small $q_i^{-1}(\coprod_iU_i)=\{x\in X_i:(x,i)\in\coprod_i U_i\}=\{x\in X_i:x\in U_i\}=U_i$}.
Also, $f$ is continuous, since for any $B\subset Y$,
\[
\textstyle f^{-1}(B)=\left\{(x,i)\in\coprod_i X_i~|~f((x,i))\in B\right\}=\left\{x\in\coprod_i X_i~|~f_i(x)\in B\right\}=\coprod_i f_i^{-1}(B). \qedhere
\]
\end{proof}

\begin{dfn}[\textcolor{blue}{
\index{Local! homeomorphism}{Local homeomorphism},
\index{Sheet}{Sheet},
\index{Stalk (fiber)}{Stalk (fiber)},
\index{Protosheaf (sheaf prototype)}{Protosheaf (sheaf prototype)},
\index{Sheaf space (total space)}{Sheaf space (total space)},
\index{Projection}{Projection},
\index{Base space}{Base space},
\index{Protosheaf morphism}{Protosheaf morphism},
\index{Category of! protosheafs}{Category of protosheafs},
\index{Covering space}{Covering space}}]
A continuous map $f:X\ra Y$ is a \ul{local homeomorphism} (\blue{footnote}\footnote{In other words, $f$ is ``locally an open imbedding''. One obvious generalization is to consider $f$ to be a ``\ul{locally open map}'' (i.e., each $x\in X$ has a neighborhood $O_x$ such that $f|_{O_x}:O_x\subset X\ra Y$ is an open map.)}) if every point $x\in X$ has an open neighborhood $U\subset X$ (called a \ul{sheet} of $f$ at $x$) such that $f|_U:U\subset X\ra Y$ is an open imbedding (equivalently, (i) $f|_U:U\ra f(U)$ is a homeomorphism and (ii) $f(U)\subset Y$ is open). For each $y\in Y$, the \ul{stalk} of $f$ at $y$ is the fiber $E_y:=f^{-1}(y)\subset X$ over $y$ as a subspace of $X$.

A surjective local homeomorphism $p:E\ra X$ is called an \ul{$X$-protosheaf} (or \ul{protosheaf} over $X$, or \ul{sheaf prototype} over $X$). The $X$-protosheaf $p:E\ra X$ is also written as $(E,p,X)$ or $(E,p)$, where $E$ is called the \ul{sheaf space} (or \ul{total space}), $p$ the \ul{projection}, $X$ the \ul{base space}, and for each $x\in X$, the subspace $E_x:=p^{-1}(x)\subset E$ (i.e., the fiber over $x$) the \ul{stalk} at $x$. Given $X$-protosheafs $(E,p,X)$ and $(E',p',X)$ over $X$, a \ul{protosheaf morphism} $\vphi:(E,p,X)\ra(E',p',X)$ is a continuous map $\vphi:E\ra E'$ such that $p'\circ\vphi=p$ (i.e., the restriction of $\vphi$ on each fiber $E_x$, $x\in X$, gives a continuous map ~$\vphi|_{E_x}:E_x\ra E_x'$).
\bea\bt
E\ar[d,"p"]\ar[rr,"\vphi"] && E'\ar[d,"p'"] \\
X\ar[rr,"id_X"] && X
\et~~~~~~p'\circ \vphi=id_X\circ p=p.\nn
\eea
The \ul{category of $X$-protosheafs}, denoted by \ul{$X$-prosh}, is the category whose objects are $X$-protosheafs $(E,p)$ and whose morphisms are protosheaf morphisms
\[
Hom_X\big((E,p),(E',p')\big):=\{\txt{continuous maps $f:E\ra E'$ with $p'f=p$}\}.
\]

An $X$-protosheaf $(E,p,X)$ is called a \ul{covering space} of $X$ if each point $x\in X$ has an open neighborhood $x\in U\subset X$ such that $p^{-1}(U)=\bigsqcup_{i\in I}O_i$ is a disjoint union of open subsets $O_i=O_{i,x,U}\subset E$ on which $p$ restricts to a homeomorphism $p|_{O_i}:O_i\ra U$ for each $i\in I=I_{x,U}$.
\end{dfn}

\begin{lmm}\label{LocHomeoOpen}
Let $f:X\ra Y$ be a local homeomorphism. Then the following hold:
\bit[leftmargin=1cm]
\item[(a)] The sheets of $f$ form a base for the topology of $X$. (\blue{footnote}\footnote{Each sheet $f:O\ra f(O)$ (i.e., an open imbedding $O\subset X\hookrightarrow Y$) represents an open continuous section $s_O:f(O)\subset Y\ra O\subset X$. Thus, the images of these special sections (which are open imbeddings) form a base for the topology of $X$.})
\item[(b)] $f$ is an open map. (Hence, a bijective local homeomorphism is a homeomorphism.) (\blue{footnote}\footnote{This remains true if $f$ is only locally open (i.e., a locally open map is an open map). Also, in general, if $U\subset\bigcup_{x\in U}O_x$,
\bea
\textstyle f(U)=f\left(\bigcup_{x\in U} O_x\cap U\right)=\bigcup_{x\in U} f(O_x\cap U)\subset \bigcup_{x\in U} f(O_x)\cap f(U)=f\left(\bigcup_{x\in U}O_x\right)\cap f(U)=f(U)\cap f(U)=f(U).\nn
\eea
})
\item[(c)] Each stalk of $f$ is a discrete space.
\eit
\end{lmm}
\begin{proof}
{\flushleft (a)}: Let $U\subset X$ be open, and $x\in U$. Let $O_x$ be a sheet of $f$ at $x$. Then $x\in O_x\cap U\subset U$.

{\flushleft (b):} Let $U\subset X$ be open. For each $x\in U$, let $O_x$ be a sheet of $f$ at $x$, i.e., $f|_{O_x}:O_x\ra f(O_x)$ is a homeomorphism (where $f(O_x)$ is open). Then $f(U)=f\left(\bigcup_{x\in U}O_x\cap U\right)=\bigcup_{x\in U}f|_{O_x}(O_x\cap U)$ is open as a union of open sets (since each $f|_{O_x}$ is a homeomorphism, hence an open map). This proves $f$ is an open map. (Hence, if $f$ is bijective then $f^{-1}:Y\ra X$ is continuous, since $f=(f^{-1})^{-1}$ is open.)

{\flushleft (c)}: Let $y\in Y$ and consider the stalk $E_y:=f^{-1}(y)\subset X$. If $E_y=\emptyset$ the result is clear. So, assume $E_y\neq\emptyset$. For each $x\in E_y$, let $O_x$ be a sheet of $f$ at $x$. Then $E_y\cap O_x=\{x\}$ since $f|_{O_x}$ is injective (hence $\{x\}\subset E_y$ is open in $E_y$).
\end{proof}

\begin{lmm}[\blue{\index{Topological! gluing lemma}{Topological gluing lemma}}]\label{TopGluLmm}
(i) Maps that agree on a cover agree everywhere. (ii) Continuous maps on elements of an open cover are restrictions of a unique continuous map, provided they agree on the intersections of their domains.
\end{lmm}
\begin{proof}
{\flushleft (i)} Let $X$ be a space and $\{A_i\}_{i\in I}$ subsets of $X$. Let $f,g:\bigcup_i A_i\subset X\ra Y$ be maps such that $f|_{A_i}=g|_{A_i}$ for all $i\in I$. Then for each $x\in\bigcup_i A_i$, $x\in A_{i_x}$ for some $i_x\in I$, and so $f(x)=f|_{A_{i_x}}(x)=g|_{A_{i_x}}(x)=g(x)$. Hence $f=g$ on $\bigcup_i A_i$.
{\flushleft (ii)} Let $\{O_i\}_{i\in I}$ be open subsets of $X$. Let $f_i:O_i\subset X\ra Y$ be continuous maps such that $f_i|_{O_i\cap O_j}=f_j|_{O_i\cap O_j}$ for all $i,j\in X$. Then the map ~$f:\bigcup_i O_i\subset X\ra Y,~x\mapsto f_i(x)$ if $x\in O_i$~ clearly satisfies $f|_{O_i}=f_i$. Moreover, $f$ is continuous: For any open set $U\subset Y$, the preimage $f^{-1}(U)=(\bigcup_i O_i)\cap f^{-1}(U)=\bigcup_i \big(O_i\cap f^{-1}(U)\big)=\bigcup_i f_i^{-1}(U)$ is open in $X$, since each $f_i^{-1}(U_i)$ is open in $O_i$ (which is open in $X$) and so open in $X$. Also, $f$ is unique by part (i) above.
\end{proof}

\begin{crl}
Let $f:X\ra Y$ be a continuous map, $\{U_i\}_{i\in I}\subset \P(Y)$ open subsets of $Y$, and $\{s_i:U_i\ra X\}_{i\in I}$ continuous sections of $f$ (i.e., each $s_i$ is continuous and $f\circ s=id_{U_i}$). If $s_i|_{U_i\cap U_j}=s_j|_{U_i\cap U_j}$ for all $i,j\in I$, then there exists a unique continuous section $s:\bigcup_i U_i\subset Y\ra X$ such that $s|_{U_i}=s_i$.
\end{crl}
\begin{proof}
The existence and uniqueness of $s$ follow from the preceding result. Also, $s$ is a section of $f$ since each $x\in \bigcup_i U_i$ lies in some $U_{i_x}$, and so $f\circ s(x)=f\circ s_{i_x}(x)=x$, i.e., $f\circ s=id_{\bigcup_i U_i}$.
\end{proof}

\begin{rmk}[\blue{An aside: Nondestructive gluing}]
Let $X,Y$ be spaces, $\{X_i\}_{i\in I}$ a cover of $X$, and $f_i:X_i\ra Y$ maps that agree at intersections in the sense $f_i|_{X_i\cap X_j}=f_j|_{X_i\cap X_j}$ for all $i,j\in I$. We have seen that (i) if the cover $\{X_i\}$ is ``\emph{continuous}'' in the sense it is \emph{open} and (ii) the maps $\{f_i\}$ are \emph{continuous} then we can glue them to get a \emph{unique continuous} map $f:X\ra Y$ such that $f|_{X_i}=f_i$ for each $i\in I$. This nondestructive gluing was made possible by the strong dependence between continuity and openness. Similarly, (i) if the cover is ``\emph{constant}'' in the sense $X_i=X$ for all $i\in I$ and (ii) the maps $\{f_i\}$ are \emph{constant} then it is also clear that we can again glue them to get a \emph{unique constant} map $f:X\ra Y$ such that $f|_{X_i}=f_i$ for all $i\in I$.

We can thus also ask the following question (assuming ``differentiability'' for a map $X\ra Y$ is at least as elaborate as for the case of a function $\Real\ra\Real$): If the maps $\{f_i\}$ are \emph{differentiable}, what property should the cover $\{X_i\}$ have (making it a ``\emph{differentiable cover}'') so that we can glue the maps to get a \emph{unique differentiable} map $f:X\ra Y$ satisfying $f|_{X_i}=f_i$ for all $i\in I$? Openness of the cover is not enough, because constant maps for example are differentiable but gluing constant maps with distinct values on a finite disjoint open cover clearly gives a map that is neither constant nor differentiable. For the case where $X=Y=\Real$ and $I$ is finite, a ``differentiable'' cover is impossible as can be seen using constant functions (with distinct values) on such a cover.
\end{rmk}

\section{Algebraic Topology}\label{AlgTopSec} 
Algebraic topology involves a geometric study/classification of topological spaces using special modules/rings defined over (i.e., assigned in some way to) these spaces. In algebraic topology, common methods used to describe certain geometric invariants of a space are known as (co)homologies of the space. An example of (co)homology of a space is singular (co)homology. In such a study, geometric invariants of a space $X$ are obtained as module/ring-valued functions of \ul{geometric objects in $X$}, defined to be continuous maps (from certain basic spaces) to $X$ called \ul{simplices in $X$} (e.g., points, lines/paths, planes/surfaces, etc).

Throughout this section the real \ul{unit interval} will be denoted by $I:=[0,1]\subset\Real$, unless stated otherwise. If $X$ is a space or an abelian group, the \ul{identity map} $id_X:X\ra X,~x\mapsto x$ will as usual also be denoted by $1_X$. Meanwhile, as usual, we will denote \ul{homeomorphism} of spaces $X,Y$ by $X\cong Y$, and also denote \ul{isomorphism} of groups $G,H$ by $G\cong H$.

\subsection{The homotopy types and homotopy groups of spaces}
{\flushleft Let} $f:X\ra Y$ be any (continuous) map. Given subsets $A\subset X$, $B\subset Y$ such that $f(A)\subset B$, we will write $f:(X,A)\ra (Y,B)$, and call $f$ a \index{Map! of pairs}{\ul{map of pairs}} of spaces $(X,A)$ and $(Y,B)$.
\begin{dfn}[\textcolor{blue}{
\index{Path!}{Path},
\index{Loop}{Loop},
\index{Reverse path}{Reverse path},
\index{Reverse map}{Reverse map},
\index{Constant map}{Constant map},
\index{Constant path}{Constant path}}]
Let $X$ be a space. A continuous map of the form $\gamma:I\ra X$ is called a \ul{path} in $X$ from $\gamma(0)$ to $\gamma(1)$. A path $\gamma:I\ra X$ is a \ul{loop} if $\gamma(0)=\gamma(1)$. Given a path $\gamma:I\ra X$, the \ul{reverse path} of $\gamma$ is the path $\ol{\gamma}:I\ra X$, $\ol{\gamma}(t):=\gamma(1-t)$. Similarly, given a continuous map $f:I^n\ra X$, we define (up to homeomorphism of $I^n$) its \ul{reverse map} $\ol{f}:I^n\ra X$ by $\ol{f}(t_1,...,t_n):=f(1-t_1,t_2,\cdots,t_n)$.

For spaces $X,Y$ and a point $y_0\in Y$, the \ul{constant map} at $y_0$, i.e., $X\ra Y$, $x\mapsto y_0$, will be denoted by $c_{y_0}:X\ra Y$ (with the understanding that $c_{y_0}$ depends on $X$). In particular, for any $x\in X$, the \ul{constant path} at $x$, i.e., $I\ra X$, $t\mapsto x$, will similarly be denoted by $c_x:I\ra X$.
\end{dfn}

\begin{dfn}[\textcolor{blue}{
\index{Path! components}{Path components} of a space,
\index{Path! connected space}{Path connected space}}]
Let $X$ be a space. Define a relation $\sim_p$ on $X$ by $x\sim_p y$ if there exists a path $\gamma:[0,1]\ra X$ between $x$ and $y$, i.e., such that $\gamma(0)=x$ and $\gamma(1)=y$. Then $\sim_p$ is an equivalence relation, and for each $x\in X$, the equivalence class $[x]_p:=\{y\in X:~y\sim_p x\}\subset X$ is called a \ul{path component} of $X$. The space $X$ is \ul{path connected} if $X=[x]_p$ for some $x\in X$ (i.e., any two points of $X$ are connected by a path).
\end{dfn}

\begin{dfn}[\textcolor{blue}{
\index{Homotopy! topological}{Homotopy},
\index{Invariant! homotopy}{Invariant homotopy}}]
Let $X,Y$ be spaces.
A family of maps $\{f_t:X\ra Y\}_{t\in I}$ is called a \ul{homotopy} (or continuous family of continuous maps) if the joint map ~$F:X\times I\ra Y$, $(x,t)\mapsto f_t(x)$~ is continuous. The homotopy is \ul{invariant} if $F(x,t)$ is constant in $t$, i.e., if the map $F:X\times I\ra Y$ is equivalent to a single map $f:X\ra Y$.
\end{dfn}
The joint map $F$ is itself called the homotopy, because any continuous map $F:X\times I\ra Y,~(x,t)\mapsto F(x,t)$ gives a homotopy $\big\{f_t:X\ra Y,~x\mapsto f_t(x):=F(x,t)\big\}_{t\in I}$; recall that a multivariate continuous map $h:X_1\times\cdots\times X_n\ra Y,~(x_1,...,x_n)\mapsto h(x_1,...,x_n)$ is continuous in each variable $x_i$ separately. Also, a homotopy $F:X\times I\ra Y$ is equivalently a continuous family of paths $\gamma_X$ $:=$ $\left\{\gamma_x:I\ra Y\right\}_{x\in X}$, where for each $x\in X$, $\gamma_x(t):=F(x,t)=f_t(x)$ is a path in $Y$ from $F(x,0)=f_0(x)$ to $F(x,1)=f_1(x)$. Thus, the homotopy is a \emph{continuous deformation} in $Y$ of the image $F(X,0)=f_0(X)\subset Y$, through the images $F(X,t)=f_t(X)\subset Y$ into the image $F(X,1)=f_1(X)\subset Y$, of $X$. Therefore,
\[\adjustbox{scale=0.9}{%
\bt[row sep=tiny]
 X\times I\ar[r,"F"] & Y \ar[r,draw=none,"\eqv"description] & \Big\{I_x\ar[r,"\gamma_x"] & Y\Big\}_{x\in X} \ar[r,draw=none,":="description] & \Big\{\{x\}\times I \ar[r,"{F(x,~)}"] & Y\Big\}_{x\in X}, & \gamma_X:I\ra\P(Y),~t\mapsto F(X,t),\\
 X\times I\ar[r,"F"] & Y \ar[r,draw=none,"\eqv"description] & \Big\{X_t\ar[r,"f_t"] & Y\Big\}_{t\in I}\ar[r,draw=none,":="description] & \Big\{X\times\{t\} \ar[r,"{F(~,t)}"] & Y\Big\}_{t\in I}, & f_I:X\ra \P(Y),~x\mapsto F(x,I).
\et}\]

\begin{dfn}[\textcolor{blue}{\index{Homotopy! relative to a set}{Homotopy relative to a set}}]
Let $X$ be a space and $A\subset X$. A homotopy $F:X\times I\ra Y$ is called a \ul{homotopy relative} to $A$ if the \ul{subhomotopy} $F|_{A\times I}:A\times I\ra Y$ is invariant, in the sense ~$F|_{A\times\{t\}}(-,t)=F|_{A\times\{t'\}}(-,t')$~ for all $t,t'\in I$. That is, the homotopy leaves the image $F(A,0)=F(A,t)=F(A,1)\subset Y$ of $A$ undeformed.
\end{dfn}

\begin{dfn}[\textcolor{blue}{\index{Homotopic maps}{Homotopic maps}}] Continuous maps $f,g:X\ra Y$ are \ul{homotopic} (written $f\simeq g$) if there is a homotopy $\{f_t:X\ra Y\}_{t\in I}$ such that $f_0=f$, $f_1=g$ (i.e., such that the continuous map $F(x,t)=f_t(x)$ satisfies ~$F|_{X\times\{0\}}=f$,~ $F|_{X\times\{1\}}=g$). That is, in $Y$, the image $f(X)$ can be continuously deformed into the image $g(X)$.
\end{dfn}

\begin{dfn}[\textcolor{blue}{\index{Nullhomotopic map}{Nullhomotopic map}}] A continuous map $f:X\ra Y$ is \ul{nullhomotopic} if it is homotopic to a constant map $c_{y_0}:X\ra Y,~x\mapsto y_0$ (for some $y_0\in Y$), i.e., $f\simeq c_{y_0}$. That is, the image $f(X)\subset Y$ can be continuously deformed into a point $y_0\in Y$.
\end{dfn}

\begin{dfn}[\textcolor{blue}{
\index{Homotopy! equivalence}{Homotopy equivalence},
\index{Homotopy! inverse}{Homotopy inverse},
\index{Homotopy! equivalent spaces}{Homotopy equivalent spaces},
\index{Homotopy! type}{Homotopy type of a space}}]
A continuous map $f:X\ra Y$ is a \ul{homotopy equivalence} if there exists a continuous map $g:Y\ra X$ (called \ul{homotopy inverse} of $f$) such that $f\circ g\simeq 1_Y$ and $g\circ f\simeq 1_X$. If a homotopy equivalence $X\ra Y$ exists, we say $X$ and $Y$ are \ul{homotopy equivalent} (or of the \ul{same homotopy type}), written $X\simeq Y$.
\end{dfn}

\begin{dfn}[\textcolor{blue}{\index{Contractible space}{Contractible space}}]
A space $X$ is \ul{contractible} if $X\simeq\{x_0\}$ for a point $x_0\in X$, or equivalently (by the following remark), if the identity $1_X:X\ra X$ is nullhomotopic.
\end{dfn}

\begin{rmk}
Suppose $X\simeq\{x_0\}$. Let $f:X\ra\{x_0\}$ be a hty equivalence with hty inverse $g:\{x_0\}\ra X$, i.e., $fg\simeq id_{\{x_0\}}$ and $gf\simeq id_X$. Since both $f$ and $g$ are constant maps, we have $id_{\{x_0\}}=fg=gf|_{\{x_0\}}=c_{x_0}|_{\{x_0\}}$ (which is completely trivial) and $id_X\simeq gf=c_{x_0}$. Hence $X\simeq\{x_0\}$ $\iff$ $id_X\simeq c_{x_0}$ (i.e., $\iff$ there is a homotopy $H:X\times I\ra X$ such that $H(x,0)=x_0$ and $H(x,1)=x$).
\end{rmk}

\begin{crl}[\textcolor{blue}{Contractibility criterion}]  Let $X$ be a space. (i) $X$ is contractible $\iff$ (ii) every continuous map $f:Z\ra X$ is nullhomotopic $\iff$ (iii) every continuous map $f:X\ra Z$ is nullhomotopic.
\end{crl}
\begin{proof}
{\flushleft \ul{(i)$\iff$(ii)}}: For (i)$\Ra$(ii), let $X\simeq\{x_0\}$ via a hty $H:X\times I\ra X$ such that $H|_{X\times\{0\}}=c_{x_0}$ and $H|_{X\times\{1\}}=id_X$. Given a continuous map $f:Z\ra X$, we get a homotopy $H_f:=H\circ(f\times id_I):Z\times I\sr{f\times id_I}{\ral}X\times I\sr{H}{\ral}X$ between $f:Z\ra X$ and $c_{x_0}:Z\ra X,~z\mapsto x_0$. ((ii)$\La$(i) is immediate.)
{\flushleft \ul{(i)$\iff$(iii)}}: For (i)$\Ra$(iii), let $X\simeq\{x_0\}$ via a hty $H:X\times I\ra X$ such that $H|_{X\times\{0\}}=c_{x_0}$ and $H|_{X\times\{1\}}=id_X$. Given a continuous map $f:X\ra Z$, we get a homotopy $H_f:=f\circ H:X\times I\sr{H}{\ral}X\sr{f}{\ral}Z$ between $f:X\ra Z$ and $c_{f(x_0)}:X\ra Z,~x\mapsto f(x_0)$. ((iii)$\La$(i) is immediate.)
\end{proof}

\begin{crl}\label{NullHomExtCrl}
Let $X$ be a space and $f:A\subset X\ra Z$ a continuous map. If (i) $X$ is contractible and (ii) $f$ has a continuous extension $F:X\ra Z$ ($f=F|_A$), then $f$ is nullhomotopic.
\end{crl}
\begin{proof}
Assume $X$ is contractible. Let $F:X\ra Z$ be a continuous extension of $f$. Then $F\simeq c_{z_0}$ for a point $z_0\in Z$ via a homotopy $H:X\times I\ra Z$, and so $f\simeq c_{z_0}$ via the homotopy $H|_{A\times I}:A\times I\ra Z$.
\end{proof}

\begin{dfn}[\textcolor{blue}{\index{Gluing! of continuous maps}{Gluing of compatible continuous maps on $[0,1]^n$}}]
Let $X$ be a space. Given continuous maps $f,g:[0,1]\ra X$ such that $f(1)=g(0)$, we can glue them to obtain the continuous map
~{\small
\[
f\cdot g:[0,1]\ra X,~~f\cdot g(t)=
\left\{
  \begin{array}{ll}
    f(2t), & t\in[0,1/2] \\
    g(2t-1), & t\in[1/2,1]
  \end{array}
\right\},\]}
where in particular,
~{\small
\[
f\cdot\ol{f}:[0,1]\ra X,~~f\cdot \ol{f}(t)=
\left\{
  \begin{array}{ll}
    f(2t), & t\in[0,1/2] \\
    \ol{f}(2t-1)=f\big(2(1-t)\big), & t\in[1/2,1]
  \end{array}
\right\}.\]}
Similarly, if continuous maps $f,g:I^n\ra X$ are such that for the two ``opposite'' faces $A=\{0\}\times I^{n-1}\subset\del I^n$ and $B=\{1\}\times I^{n-1}\subset\del I^n$ of $I^n$ we have (up to homeomorphism of $I^n$)
\begin{align}
f(1,t_2,\cdots,t_n)=g(0,t_2,\cdots,t_n),~~\txt{for all}~~(1,t_2,\cdots,t_n)\in B,~~(0,t_2,\cdots,t_n)\in A,\nn
\end{align}
then we can glue them together (along the coordinate $t_1$) to obtain the continuous map
\bit
\item[] $f\cdot g(t_1,...,t_n):=\left\{
  \begin{array}{ll}
    f(2t_1,t_2,\cdots,t_n), & t_1\in[0,1/2] \\
   g(2t_1-1,t_2,\cdots,t_n), & t_1\in[1/2,1]
  \end{array}
\right\}$.
\eit
We can also consider gluing such maps based on a single common value only. For example, given continuous maps $f,g:I^n\ra X$ such that $f(1,...,1)=g(0,...,0)$, we can glue them to get the continuous map
\[
f\star g(t_1,...,t_n):=\left\{
  \begin{array}{ll}
    f(2t_1,...,2t_n), & t_1,...,t_n\in[0,1/2] \\
   g(2t_1-1,...,2t_n-1), & t_1,...,t_n\in[1/2,1]
  \end{array}
\right\}.
\]
\end{dfn}
\begin{rmk}[\textcolor{blue}{\index{Eckmann-Hilton identity}{Eckmann-Hilton identity}}]\label{ProdHtyRmk}
Let $X$ be a space. Given maps ~$f,g:(I^2,\del I^2)\ra (X,x_0)$,~ define maps ~$f\cdot g,f\ast g:(I^2,\del I^2)\ra (X,x_0)$~ by
\begin{align}
f\cdot g(t,s):=\left\{
  \begin{array}{ll}
    f(2t,s), & t\in[0,1/2] \\
   g(2t-1,s), & t\in[1/2,1]
  \end{array}
\right\},~~~~f\ast g(t,s):=\left\{
  \begin{array}{ll}
    f(t,2s), & s\in[0,1/2] \\
   g(t,2s-1), & s\in[1/2,1]
  \end{array}
\right\}.\nn
\end{align}
Then for any maps $f_1,f_2,g_1,g_2:(I^2,\del I^2)\ra (X,x_0)$, we see (by direct verification at step (s) below) that
{\footnotesize\begin{align}
\label{EckHiltId}&\Big((f_1\ast f_2)\cdot(g_1\ast g_2)\Big)(t,s)=\left\{
  \begin{array}{ll}
    f_1\ast f_2(2t,s), & t\in[0,1/2] \\
   g_1\ast g_2(2t-1,s), & t\in[1/2,1]
  \end{array}
\right\}\sr{(s)}{=}\left\{
  \begin{array}{ll}
    f_1\cdot g_1(t,2s), & s\in[0,1/2] \\
   f_2\cdot g_2(t,2s-1), & s\in[1/2,1]
  \end{array}
\right\}\nn\\
&~~=\Big((f_1\cdot g_1)\ast(f_2\cdot g_2)\Big)(t,s).
\end{align}}It is clear that this remains true if we consider any maps $f,g:(I^n,\del I^n)\ra (X,x_0)$ and suppress/fix all but any two of the $n$ arguments in the definitions of $f\cdot g,f\ast g:(I^n,\del I^n)\ra (X,x_0)$.
\end{rmk}

\begin{dfn}[\textcolor{blue}{
\index{Homotopy! groups}{Homotopy groups} of a space,
\index{Homotopy! class}{Homotopy class} of a map}]
Let $X$ be a space. Fix $x_0\in X$, and an integer $n\geq 0$. The $n$th \ul{homotopy group} of $X$ at $x_0$ is the set of maps
{\small\begin{align}
&\textstyle\pi_n(X,x_0):=\big\{\txt{continuous maps}~f:(I^n,\del I^n)\ra (X,x_0),~\txt{up to homotopy relative to $\del I^n$}\big\}\nn\\
&\textstyle~~~~={\big\{\txt{continuous maps}~f:(I^n,\del I^n)\ra (X,x_0)\big\}\over \txt{$f\sim g$ if $f\simeq_\del g$ in the sense ``$f\simeq g$ relative to $\del I^n$''}}=\Big\{[f]~\big|~\txt{for continuous maps}~f:(I^n,\del I^n)\ra (X,x_0)\Big\}\nn
\end{align}}as a multiplicative set with multiplication given by
\[
[f][g]:=[f\cdot g],~~\txt{for all}~~[f],[g]~\in~\pi_n(X,x_0).
\]
The set of homotopic maps $[f]:=\{g:g\simeq_\del f\}$ is called the \ul{homotopy class} of $f$.
\end{dfn}

\begin{rmks}
\begin{enumerate}[leftmargin=0.8cm]
\item By the universal property of quotient maps, any continuous map $f:(I^n,\del I^n)\ra (X,x_0)$ corresponds to a unique continuous map ${f'}:\big(I^n/\del I^n,\{\del I^n\}\big)\mathop{\cong}\limits^{\vphi}\big(S^n,\vphi(\{\del I^n\})\big)\ra (X,x_0)$.
\item $\pi_n(X,x_0)$ is a \emph{\ul{group}} for every $n\geq 1$, with \emph{\ul{identity}} $e:=[c_{x_0}]$ and \emph{\ul{inverse}} $[f]^{-1}:=[\ol{f}]$, where
\[
c_{x_0}\simeq_\del f\cdot\ol{f}~~~~\txt{via}~~~~
H:I^2\times[0,1]\ra X,~~H\big((x,y),t\big):=
\left\{
  \begin{array}{ll}
    f(2(2-t)x,y), & x\in[0,{t\over 2}] \\
   f(2t(1-x),y), & x\in[{t\over 2},1]
  \end{array}
\right\}.
\]
\item $\pi_n(X,x_0)$ is an \emph{\ul{abelian group}} for every $n\geq 2$. This is because the equality (\ref{EckHiltId}) in Remark \ref{ProdHtyRmk} shows (by setting $f_2\simeq c_{x_0}$ and $g_1\simeq c_{x_0}$, and/or otherwise) that homotopy groups defined using the product $\cdot$ are isomorphic to homotopy groups defined using the product $\ast$. Hence, by setting $f_1\simeq g_2\simeq c_{x_0}$ in (\ref{EckHiltId}), we see that the homotopy groups are abelian for $n\geq 2$. Equivalently, we have $f\cdot g\simeq g\cdot f$.
\end{enumerate}
\end{rmks}

\begin{dfn}[\textcolor{blue}{\index{Homotopy! $n$-chains}{Homotopy $n$-chains}}]
Let $X$ be a space. The \ul{homotopy $n$-chain} of $X$ is
\[
\C_n(X,x_0):=\{\txt{continuous maps}~f:(I^n,\del I^n)\ra (X,x_0)\}.
\]
Using $\C_n(X,x_0)$, we can write {\small $\pi_n(X,x_0)={\C_n(X,x_0)\over\simeq_\del}=\big\{[f]:f\in\C_n(X,x_0)\big\}$}, where as before, $f\simeq_\del g$ is homotopy relative to $\del I^n$ (i.e., $f\simeq_\del g$ iff $f\simeq g$ relative to $\del I^n$).
\end{dfn}

\begin{dfn}[\textcolor{blue}{\index{Induced! homomorphisms}{Induced homomorphisms}}]
Given a continuous map $f:(X,x_0)\ra (Y,y_0)$, we get homomorphisms ~{\small $\C_\#(X,x_0)\sr{f_\#}{\ral}\C_\#(Y,y_0),~u\mapsto f\circ u$},
~{\small $\pi_\ast(X,x_0)\sr{f_\ast}{\ral}\pi_\ast(Y,x_0),~[u]\mapsto [f_\#(u)]=[f\circ u]$}.
\end{dfn}
For any continuous maps $X\sr{g}{\ral} Y\sr{f}{\ral}Z$, the induced maps satisfy the relations
\begin{align}
(f\circ g)_\#=f_\#\circ g_\#,~~(1_X)_\#=1_{\C_\#(X,x_0)},~~~~(f\circ g)_\ast=f_\ast\circ g_\ast,~~(1_X)_\ast=1_{\pi_\ast(X,x_0)},\nn
\end{align}
and so define functors
\[
\#:\Top\ra \txt{Associative Sets}~~~~\txt{and}~~~~\ast:\Top\ra \txt{Groups}.
\]

\subsection{The singular homology groups of spaces}
\begin{dfn}[\textcolor{blue}{
\index{Convex hull}{Convex hull},
\index{Standard! unit vectors}{Standard unit vectors},
\index{Standard! $n$-simplex}{Standard $n$-simplex},
\index{Vertices of the standard $n$-simplex}{Vertices of the standard $n$-simplex},
\index{Orientation}{Orientation}}]
Given points $v_0,...,v_n\in\Real^{n+1}$, their \ul{convex hull} is the set {\small$\Conv(v_0,...,v_n):=\big\{{\textstyle\sum}_{i=0}^n\al_iv_i:{\textstyle\sum}_{i=0}^n\al_i=1,~\al_i\geq0\big\}$}.
Let $e_i:=(\delta_{ij})_{j=0}^n\in\Real^{n+1}$ be the \ul{standard unit vectors}. The \ul{standard $n$-simplex}
\bea
\Delta^n=[e_0,...,e_n]:=(\Conv(e_0,\cdots,e_n),O)\nn
\eea
with \ul{vertices} $e_0,...,e_n$ is the space $\Conv(e_0,\cdots,e_n)\subset\Real^{n+1}$ along with an \ul{orientation} $O\in\{\pm1\}$ based on a linear ordering of the vertices $e_0\prec e_1\prec\cdots\prec e_n$, such that
{\small\begin{align}
&[e_0,...,e_i,e_{i+1},...,e_n]:=\Big(\Conv(e_0,...,e_i,e_{i+1},...,e_n),O\Big)=\Big(\Conv(e_0,...,e_{i+1},e_i,...,e_n),-O\Big)\nn\\
&~~~~=:-[e_0,...,e_{i+1},e_i,...,e_n],~~~~\txt{for all}~~0\leq i<n.\nn
\end{align}}Thus, as a set, the standard $n$-simplex in $\Real^{n+1}$ is given by
{\small\begin{align}
\Delta^n=\Conv(e_0,...,e_n):=\big\{{\textstyle\sum}_{i=0}^nt_ie_i:{\textstyle\sum}_{i=0}^nt_i=1,~t_i\in[0,1]\big\}=\left\{(t_0,...,t_n)\in[0,1]^{n+1}:{\textstyle\sum}_{i=0}^nt_i=1\right\}.\nn
\end{align}}
\end{dfn}

\begin{dfn}[\blue{
\index{Singular! Simplex}{Singular $n$-simplex in a space},
\index{Set! of singular simplices}{Set of singular simplices},
\index{Simplex! faces of}{Faces of a simplex},
\index{Simplex! boundary faces of}{Boundary faces of a simplex},
\index{Simplex! vertices of}{Vertices},
\index{Simplex! edges of}{Edges}}]
Let $X$ be a space. A \ul{singular $n$-simplex} (an $n$-dimensional simplex) in $X$ is a continuous map $\sigma=\sigma^n:\Delta^n\ra X$. The \ul{set of singular $n$-simplices} in $X$ is denoted by
\[
S_n(X):=\C(\Delta^n,X)=\{\txt{singular $n$-simplices in $X$}\}.
\]
It is clear that in $\Delta^n=[e_0,...,e_n]$, each subset of vertices $\{e_{i_0},...,e_{i_k}\}\subset\{e_0,...,e_n\}$ determines a $k$-dimensional \ul{subsimplex} $[e_{i_0},...,e_{i_k}]:=\big(\Conv(e_{i_0},...,e_{i_k}),O_{e_{i_0},...,e_{i_k}}\big)\subset\Delta^n$, with orientation $O_{e_{i_0},...,e_{i_k}}\in\{\pm1\}$ inherited from that of $\Delta^n$ via the ordering $e_0\prec e_1\prec\cdots\prec e_n$.

Let $X$ be a space and $\sigma^n\in S_n(X)$. A \ul{$k$-face} of $\sigma^n$ is a $k$-simplex of the following form (where $\cong$ is the obvious homeomorphism): For any $\left\{e_{i_0},...,e_{i_k}\right\}\subset\{e_0,...,e_n\}:=\{\txt{the vertices of $\Delta^n$}\}$,
{\footnotesize\bea
\sigma^n|_{[e_{i_0},...,e_{i_k}]}:=\sigma^n\circ F_{\{e_{i_0},...,e_{i_k}\}}^n:\bt\Delta^k\ar[rr,bend left,"{F_{\{e_{i_0},...,e_{i_k}\}}^n}"]\ar[r,hook,two heads,"\cong"] & {[e_{i_0},...,e_{i_k}]}\ar[r,hook] & \Delta^n\ar[r,"\sigma^n"] & X.\et\nn
\eea}In particular, with ~$\Delta^{n-1}_i:=[e_0,...,\wh{e}_i,\cdots,e_n]:=[e_0,....,e_{i-1},e_{i+1},....,e_n]\subset[e_0,...,e_n]$ ~the \ul{$i$th boundary face} of an $n$-simplex $\sigma^n\in S_n(X)$ is the $(n-1)$-simplex
{\small\bea
\sigma^n|_{\Delta^{n-1}_i}:=\sigma^n\circ F^n_i:\bt\Delta^{n-1}\ar[rr,bend left,"F_i^n"]\ar[r,hook,two heads,"\cong"] & {\Delta^{n-1}_i}\ar[r,hook] & \Delta^n\ar[r,"\sigma^n"] & X.\et\nn
\eea}The \ul{$0$-faces} and \ul{$1$-faces} of a simplex are called its \ul{vertices} and \ul{edges} respectively.
\end{dfn}

\begin{dfn}[\blue{
\index{Singular! chain complex/group}{Singular chain complex/group},
\index{Boundary map (singular)}{Boundary map},
\index{Singular! homology}{Singular homology},
\index{Singular! cochain complex/group}{Singular cochain complex/group},
\index{Coboundary map (singular)}{Coboundary map},
\index{Singular! cohomology}{Singular cohomology}}]\label{SingCCdef1}
Let $X$ be a space, $R$ a commutative ring, and $M$ an $R$-module. The \ul{singular chain complex} of $X$ with coefficients in $M$ is the left chain complex
(\blue{footnote}\footnote{By construction, we have a functor ~$C_\#(-;M):\Top\ra(R\txt{-mod})_0^\Integer$, ~$X\sr{f}{\ral}Y$ $\mapsto$ $C_\#(X;M)\sr{f_\#\otimes id_M}{\ral}C_\#(Y;M)$, ~with $f_\#:C_\#(X)\ra C_\#(Y)$ ~given by ~$f_n(\sigma):=f_n\circ\sigma$~ for $\sigma\in S_n(X)$.})
\bea
\label{SingChComWC}
&&\bt[row sep = tiny,column sep = small]
C_\#(X;M):~~\cdots\ar[r,"\del^M_3"]&C_2(X;M)\ar[r,"\del^M_2"] & C_1(X;M)\ar[r,"\del^M_1"]&C_0(X;M)\ar[r]& 0\ar[r] &0\ar[r] & \cdots,
\et\\
&&C_n(X;M):=C_n(X)\otimes_\Integer M\cong C_n(X;{}_RR)\otimes_R M,~~~~\del^M_n:=\del_n\otimes id_M,\nn
\eea
where ~$C_n(X):=\Integer^{\langle S_n(X)\rangle}=\{\txt{singular $n$-chains in $X$}\}$~ (called the \ul{$n$th singular chain group}  of $X$) is the free abelian group on the set of simplices $S_n(X)$, i.e.,
\begin{align}
\textstyle C_n(X):=\FreeSpan_\Integer S_n(X)=\sum_{\sigma\in S_n(X)}\Integer\sigma=\big\{\sum_{\sigma\in S_n(X)}n_\sigma\sigma:~n_\sigma\in\Integer,~n_\sigma=0~\txt{a.e.f.}\big\},\nn
\end{align}
and the \ul{boundary map} $C_n(X)\sr{\del_n}{\ral} C_{n-1}(X)$ is given by ~$\del_n\sum_\sigma n_\sigma\sigma:=\sum n_\sigma\del_n\sigma$, where
\begin{align}
\del_n\sigma:={\textstyle\sum\limits_{i=0}^n}(-1)^i\sigma|_{[e_0,...,\wh{e}_i,...,e_n]}:={\textstyle\sum\limits_{i=0}^n}(-1)^i\sigma\circ F_i^n,~~~~\txt{for}~~n\geq 0,~~~~~~\txt{for all}~~\sigma\in S_n(X).\nn
\end{align}
The $n$th \ul{singular homology group} of $X$ with coefficients in $M$ is
{\small\bea
\textstyle H_n(X;M):=H_n\big(C_\#(X;M)\big)={Z_n\big(C_\#(X;M)\big)\over B_n\big(C_\#(X;M)\big)}=:{Z_n(X;M)\over B_n(X;M)}=:{\ker\del_n^M\over\im\del_{n+1}^M}.\nn
\eea}The \ul{singular cochain complex} of $X$ with coefficients in $M$ is the right complex
(\blue{footnote}\footnote{By construction, we have a functor ~$C^\#(-;M):\Top\ra(R\txt{-mod})_0^\Integer$, ~$X\sr{f}{\ral}Y$ $\mapsto$ $C^\#(Y;M)\sr{Hom_R(f_\#,M)}{\ral}C^\#(X;M)$, ~with $f_\#:C_\#(X)\ra C_\#(Y)$ ~given by ~$f_n(\sigma):=f_n\circ\sigma$~ for $\sigma\in S_n(X)$.})
{\footnotesize\bea
\label{SingCChComWC}\bt[row sep = tiny,column sep = small]
C^\#(X;M):~~\cdots \ar[r]&0\ar[r]&0\ar[r]&C^0(X;M)\ar[r,"{d}_M^0"]& C^1(X;M)\ar[r,"{d}_M^1"]&C^2(X;M)\ar[r,"{d}_M^2"]&C^3(X;M)\ar[r,"{d}_M^3"]&\cdots
\et
\eea}
where  the \ul{$n$th singular cochain group} of $X$ (with coefficients in $M$) is
\[
C^n(X;M):=Hom_\Integer\big(C_n(X),M\big)=\{\txt{singular $n$-cochains in $X$}\}\approx \txt{Maps}\big(S_n(X),M\big)~~~~~~~~\txt{(set maps)}\nn
\]
and the \ul{coboundary map} ~${d}_M^n:=Hom_\Integer(\del_{n+1},M):C^n(X;M)\ra C^{n+1}(X;M)$ ~is given by
\bea
\vphi:C_n(X)\ra M~~~~\longmapsto~~~~\vphi\circ\del_{n+1}:C_{n+1}(X)\ra M.\nn
\eea
The $n$th \ul{singular cohomology group} of $X$ with coefficients in $M$ is
\bea
\textstyle H^n(X;M):=H_n\big(C^\#(X;M),{d}_M^\#\big)={\ker({d}_M^n)\over\im({d}_M^{n-1})}=\txt{the $n$th homology group of $\big(C^\#(X;M),{d}_M^\#\big)$}.\nn
\eea
\end{dfn}

\begin{lmm}
Let $X$ be a space and $(C(X),\del)$ its singular chain complex. The boundary map {\small $\del_n:C_n(X)\ra C_{n-1}(X)$} satisfies {\small $\del_{n+1}\circ\del_n=0$}.
\end{lmm}
\begin{proof}
Recall that for sets $A\subset B$, the characteristic function of $A$ with respect to $B$ is the map
\bea
\chi_A:B\ra 2:=\{0,1\}\subset\Real,~b\mapsto
\left\{
  \begin{array}{ll}
    1, & b\in A \\
    0, & b\not\in A
  \end{array}
\right\}.\nn
\eea
It suffices by linearity to verify the desired property on $\sigma\in S_n(X)$.
\begin{align}
&\textstyle(\del_n\circ\del_{n+1})(\sigma)=\del_n\sum\limits_{i=0}^{n+1}(-1)^i\sigma|_{[e_0,...,\wh{e}_i,...,e_{n+1}]}=\sum\limits_{i=0}^{n+1}(-1)^i\del_n\sigma|_{[e_0,...,\wh{e}_i,...,e_{n+1}]}\nn\\
&\textstyle~~~~=\sum\limits_{i=0}^{n+1}(-1)^i\bigg(\sum\limits_{j=0}^{i-1}(-1)^j\sigma|_{[e_0,...,\wh{e}_j,...,\wh{e}_i,...,e_{n+1}]}+\sum\limits_{j=i+1}^{n+1}(-1)^{j-1}\sigma|_{[e_0,...,\wh{e}_i,...,\wh{e}_j,...,e_{n+1}]}\bigg)\nn\\
&\textstyle~~~~\sr{(a)}{=}\sum\limits_{j=0}^{n+1}(-1)^j\sum\limits_{i=0}^{j-1}(-1)^i\sigma|_{[e_0,...,\wh{e}_i,...,\wh{e}_j,...,e_{n+1}]}+\sum\limits_{i=0}^{n+1}(-1)^i\sum\limits_{j=i+1}^{n+1}(-1)^{j-1}\sigma|_{[e_0,...,\wh{e}_i,...,\wh{e}_j,...,e_{n+1}]}\nn\\
&\textstyle~~~~=\sum\limits_{i,j=0}^{n+1}\Big[\chi_{\{0,...,j-1\}}(i)-\chi_{\{i+1,...,n+1\}}(j)\Big](-1)^{i+j}\sigma|_{[e_0,...,\wh{e}_i,...,\wh{e}_j,...,e_{n+1}]}\sr{(b)}{=}0,\nn
\end{align}
where at (a) we swap $i$ and $j$ in the first term, and (b) holds because ~$\chi_{\{i+1,...,n+1\}}(j)=\chi_{\{0,...,j-1\}}(i)$.
\end{proof}

\begin{exercise}
Based on the K\"unneth theorems or universal coefficient theorems, how are $H_n(X;M)$, $H_n(X;R)$, $H_n(X):=H_n(X;\Integer)$ related, and how are $H^n(X;M)$, $H^n(X;R)$, $H^n(X):=H^n(X;\Integer)$ related?
\end{exercise}

For convenience of expression, we often also define the following.
\begin{rmk*}[\textcolor{blue}{\index{Augmented singular (co)chain complex}{Augmented singular (co)chain complex}, \index{Reduced singular (co)homology}{Reduced singular (co)homology}}]
Let $X$ be a space, $R$ a commutative ring, and $M$ an $R$-module. The \ul{augmented singular chain complex} of $X$, with coefficients in $M$, is the left chain complex
\bea
\label{SingChComAugWC}
&&\bt[row sep = tiny,column sep = small]
\wt{C}_\#(X;M):~~\cdots\ar[r,"\del^M_3"]&C_2(X;M)\ar[r,"\del^M_2"] & C_1(X;M)\ar[r,"\del^M_1"]&C_0(X;M)\ar[r,two heads,"\vep^M"]& \overbrace{\Integer\otimes_\Integer M}^{\cong M}\ar[r] &0\ar[r] & \cdots,
\et\\
&&C_\#(X;M):=C_\#(X)\otimes_\Integer M,~~~~\del^M_\#:=\del_\#\otimes id_M~~\txt{(as before)},~~~~\vep^M:=\vep\otimes id_M,\nn
\eea
where the map $\vep:C_0(X)\ra\Integer$ is given by ~$\vep\left(\sum_\sigma n_\sigma\sigma\right):=\sum n_\sigma$, with $\sigma\in S_0(X)$. In particular, with $M:=\Integer$, the \ul{augmented singular chain complex} of $X$ is
\bea
\label{SingChComAug}\bt[row sep = tiny, column sep = small]
\wt{C}_\#(X):~~\cdots\ar[r,"\del_3"]&C_2(X)\ar[r,"\del_2"] & C_1(X)\ar[r,"\del_1"]&C_0(X)\ar[r,two heads,"\vep"]&\Integer\ar[r]& 0\ar[r]& \cdots.
\et
\eea

The $n$th \ul{reduced singular homology group} of $X$, with coefficients in $M$, is the homology
{\small\bea
\textstyle \wt{H}_n(X;M):=H_n\big(\wt{C}_\#(X;M)\big)={Z_n\big(\wt{C}_\#(X;M)\big)\over B_n\big(\wt{C}_\#(X;M)\big)}=:{\wt{Z}_n(X;M)\over \wt{B}_n(X;M)}=:{\ker\wt{\del}_n^M\over\im\wt{\del}_{n+1}^M},\nn
\eea}where we have an exact sequence
{\footnotesize $\bt[column sep=small]0\ar[r] & \wt{H}_0(X;M)={\ker\vep^M\over\im\del_1^M}\ar[r,hook] & H_0(X;M)={C_0(X;M)\over\im\del_1^M}\ar[r,two heads] & {C_0(X;M)\over\ker\vep^M}\cong M\ar[r] & 0\et$} and

\bea
\wt{H}_n(X;M)=H_n(X;M),~~\txt{for all}~~n\neq 0.\nn
\eea
In particular, the $n$th \ul{reduced singular homology group} of $X$ is the homology ~$\wt{H}_n(X):=H_n\big(\wt{C}_\#(X)\big)\subset H_n(X)$, where the exact sequence {\footnotesize $\bt[column sep=small]0\ar[r] & \wt{H}_0(X)\ar[r,hook] & H_0(X)\ar[r,two heads] & \Integer\ar[r] & 0\et$} is split (as $_\Integer\Integer$ is projective), and so
\bea
H_0(X)\cong\wt{H}_0(X)\oplus\Integer~~~~\txt{and}~~~~H_n(X)=\wt{H}_n(X)~~\txt{for all}~~n\neq 0.\nn
\eea
The \ul{augment singular cochain complex} and \ul{reduced singular cohomology groups} of $X$, with coefficients in an $R$-module $M$, can be defined from these in the obvious way. That is, with {\footnotesize $\wt{C}^n(X;M):=Hom_\Integer\big(\wt{C}_n(X),M\big)$ and $\wt{{d}}_M^n:=Hom_\Integer(\wt{\del}_{n+1},M):\wt{C}^n(X;M)\ra \wt{C}^{n+1}(X;M)$},
{\footnotesize\bea
\label{SingCChComAugWC}\bt[row sep = tiny,column sep = small]
\wt{C}^\#(X;M):~~\cdots \ar[r]& 0\ar[r]&M\ar[r,hook,"\vep'"]&C^0(X;M)\ar[r,"{d}_M^0"]& C^1(X;M)\ar[r,"{d}_M^1"]&C^2(X;M)\ar[r,"{d}_M^2"]&C^3(X;M)\ar[r,"{d}_M^3"]&\cdots
\et
\eea}where for $m\in M$, ~$\vep'(m):C_0(X)\ra M$~ is given by ~$\vep'(m)(s):=\vep(s)m$~ for all $s\in C_0(X)$, and
\bea
\textstyle \wt{H}^n(X;M):=H_n\big(\wt{C}^\#(X;M),\wt{{d}}_M^\#\big)={\ker(\wt{{d}}_M^n)\over\im(\wt{{d}}_M^{n-1})},\nn
\eea
where we have an exact sequence {\footnotesize \bt[column sep=small]0\ra M\ar[r,hook,"\vep'"] & {H^0(X;M)=\ker{d}_M^0}\ar[r,two heads] & {\wt{H}^0(X;M)={\ker d_M^0\over\im\vep'}={H^0(X;M)\over\im\vep'}}\ra0\et} and
\bea
H^n(X;M)=\wt{H}^n(X;M)~~~~\txt{for all}~~n\neq 0.\nn
\eea
In particular, if $M=\Integer$, the exact sequence {\footnotesize \bt[column sep=small]0\ra \Integer\ar[r,hook] & H^0(X)\ar[r,two heads] & \wt{H}^0(X)\ra 0\et} is split (\blue{footnote}\footnote{According to the UCT, $H^0(X)\cong Hom_\Integer(H_0(X),\Integer)$ and $\wt{H}^0(X)\cong Hom_\Integer(\wt{H}_0(X),\Integer)$, and so the given sequence is the $Hom_\Integer$-image of the split-exact sequence ~$0\ra \wt{H}_0(X)\hookrightarrow H_0(X)\twoheadrightarrow \Integer\ra 0$.}), and so
\bea
H^0(X)\cong\wt{H}^0(X)\oplus\Integer,~~~~H^n(X)=\wt{H}^n(X)~~\txt{for all}~~n\neq 0,~~~~\txt{where}~~H^\ast(X):=H^\ast(X;\Integer).\nn
\eea
\end{rmk*}

\begin{exercise}
Is it true that {\small $H_0(X;M)\cong\wt{H}_0(X;M)\oplus M$}? Is it true that {\small $H^0(X;M)\cong\wt{H}^0(X;M)\oplus M$}?
\end{exercise}

\begin{rmk}[\textcolor{blue}{Interpretation of the homology of a space $X$: \index{Abelianization}{Abelianization},
\index{Commutator subgroup}{Commutator subgroup},
\index{Commutator}{Commutator}}]
Let $X$ be a space and $(C(X),\del)$ the singular chain complex of $X$. Then based on Corollary \ref{NullHomExtCrl}, we have the following:
{\small\begin{align}
&\textstyle H_0(X)={\ker\del_0\over\im\del_1}={C_0(X)\over\im\del_1}\cong{\FreeSpan_\Integer\{\txt{points in $X$}\}\over\FreeSpan_\Integer\{\gamma(1)-\gamma(0)~:~\txt{$\gamma$ a path in $X$}\}}
={\sum_{x\in X}\Integer x\over\sum_{x\sim_py}\Integer(x-y)}\nn\\
&\textstyle~~~~={\textstyle\sum}_{x\in X}\Integer\ol{x},~~~~\txt{with}~~~~\ol{x}:=x+{\textstyle\sum}\big\{\Integer(x-y):x\sim_py~\txt{in $X$}\big\},\nn\\
&\textstyle~~~~\cong \FreeSpan_\Integer\{\txt{path components of $X$}\},\nn
\end{align}}
where $x\sim_py$ iff there is a path between $x$ and $y$, iff ~$\Integer \ol{x}=\Integer\ol{y}$. For step (s) below, see \cite[p.166]{hatcher2001}.
\begin{align}
&\textstyle H_1(X)={\ker\del_1\over\im\del_2}\cong{\FreeSpan_\Integer\{\gamma:~\txt{$\gamma$ a loop in $X$}\}\over\FreeSpan_\Integer\{\gamma:~\txt{$\gamma$ a nullhomotopic loop in $X$}\}}\sr{(s)}{\cong}\txt{abelianization of $\pi_1(X)$, if $X$ is path-connected},\nn\\
&\textstyle H_n(X)={\ker\del_n\over\im\del_{n+1}}\cong{\FreeSpan_\Integer\{f:S^n\ra X~|~\txt{$f$ a continuous map}\}\over\FreeSpan_\Integer\{f:S^n\ra X~|~\txt{$f$ a nullhomotopic continuous map}\}}\nn\\
&\textstyle~~~~~~~~\cong\FreeSpan_\Integer\{\txt{certain $n$-dimensional holes in $X$}\},\nn
\end{align}
where ~$S^0:=\{x\in\Real:x^2=1\}=\{\pm1\}$ ~~(i.e., a two-point set), and the \ul{abelianization} of a group $G$ is the abelian quotient group $G_{ab}:=G/[G,G]$, where
\begin{align}
[G,G]:=\big\langle [a,b]:a,b\in G\big\rangle=\big\langle aba^{-1}b^{-1}:a,b\in G\big\rangle\nn
\end{align}
is the normal subgroup of $G$ (called the \ul{commutator subgroup} of $G$) generated by elements of the form $[a,b]:=(ab)(ba)^{-1}=aba^{-1}b^{-1}$ called \ul{commutators}. The subgroup $[G,G]\subset G$ is normal because the commutators satisfy ~$g[a,b]g^{-1}=[gag^{-1},gbg^{-1}]$~ for all $a,b,g\in G$.
\end{rmk}

\begin{dfn}[\textcolor{blue}{\index{Induced! homomorphisms}{Induced homomorphisms}}]
Given a continuous map $f:X\ra Y$, we get homomorphisms $f_\#:C_\#(X)\ra C_\#(Y)$, $f_\#(\sigma):=f\circ\sigma$ (extended by linearity) that give a chain map $f_\#:C_\#(X)\ra C_\#(Y)$ because $f_\#\circ\del^{C_\#(X)}=\del^{C_\#(Y)}\circ f_\#$, and thus also give homomorphisms $f_\ast:H_\ast(X)\ra H_\ast(Y)$, $f_\ast([z]):=[f_\#(z)]=[f\circ z]$. The induced maps satisfy the immediate relations
\begin{align}
(f\circ g)_\#=f_\#\circ g_\#,~~(1_X)_\#=1_{C_\#(X)},~~~~(f\circ g)_\ast=f_\ast\circ g_\ast,~~(1_X)_\ast=1_{H_\ast(X)},\nn
\end{align}
and so define functors
\[
\#:\Top\ra\txt{Ab}~\txt{or}~(Ab)_0^\Integer~~~~\txt{and}~~~~\ast:\Top\ra \txt{Ab}~\txt{or}~(Ab)_0^\Integer.
\]
\end{dfn}

\begin{lmm}[\textcolor{blue}{\cite[Theorem 2.10, p.111]{hatcher2001}}]
Let $f,g:X\ra Y$ be continuous maps. If $f\simeq g$, then $f_\#\simeq g_\#:C_\#(X)\ra C_\#(Y)$, and so $f_\ast=g_\ast:H_\ast(X)\ra H_\ast(Y)$.
\end{lmm}
\begin{proof}
Let $\Delta^n=[e_0,...,e_n]\subset\Real^{n+1}$ be the standard $n$-simplex. Let $v_i:=e_i\times\{0\}$ and $w_i:=e_i\times\{1\}$. Then $[v_0,...,v_n]$ (resp. $[w_0,...,w_n]$) can been called/considered the lower (resp. upper) face of the figure $\Delta^n\times I$. Moreover, for each $0\leq i\leq n$, the $(n+1)$-simplex $[v_0,...,v_i,w_i,...,w_n]\subset\Delta^n\times I$ has $[v_0,...,v_i,\wh{w}_i,...,w_n]:=[v_0,...,v_i,w_{i+1},...,w_n]$ as ``lower face'' and $[v_0,...,\wh{v}_i,w_i,...,w_n]:=[v_0,...,v_{i-1},w_i,...,w_n]$ as ``upper face'', and $\Delta^n\times I=\bigcup_{i=0}^n [v_0,...,v_i,w_i,...,w_n]$. Let $F:X\times I\ra Y$ be a homotopy from $f$ to $g$. If $\sigma\in S_n(X)$, then based on the above simplex-decomposition of the solid cylinder $\Delta^n\times I$, and hence of its image $F\circ(\sigma\times 1_I):\Delta^n\times I\sr{\sigma\times 1_I}{\ral}X\times I\sr{F}{\ral}Y$, we can define a map $h:C_n(X)\ra C_{n+1}(Y)$ by ~$h(\sigma):=\sum_{i=0}^n(-1)^iF\circ(\sigma\times 1_I)|_{[v_0,...,v_i,w_i,...,w_n]}$. We have (after some cancellation at the last step below)
{\small\begin{align}
&\textstyle\del h(\sigma)=\sum\limits_{i\geq j}(-1)^{i+j}F\circ(\sigma\times 1_I)|_{[v_0,...,\wh{v}_j,...,v_i,w_i,...,w_n]}+\sum\limits_{i\leq j}(-1)^{i+j+1}F\circ(\sigma\times 1_I)|_{[v_0,...,v_i,w_i,...,\wh{w}_j,...,w_n]}\nn\\
&\textstyle~~=\sum_{i=0}^n(-1)^{2i}F\circ(\sigma\times 1_I)|_{[v_0,...,\wh{v}_i,w_i,...,w_n]}+\sum_{i=0}^n(-1)^{2i+1}F\circ(\sigma\times 1_I)|_{[v_0,...,v_i,\wh{w}_i,...,w_n]}\nn\\
&\textstyle~~~+\sum\limits_{i>j}(-1)^{i+j}F\circ(\sigma\times 1_I)|_{[v_0,...,\wh{v}_j,...,v_i,w_i,...,w_n]}+\sum\limits_{i<j}(-1)^{i+j+1}F\circ(\sigma\times 1_I)|_{[v_0,...,v_i,w_i,...,\wh{w}_j,...,w_n]}\nn\\
&\textstyle~~=F\circ(\sigma\times 1_I)|_{[w_0,...,w_n]}-F\circ(\sigma\times 1_I)|_{[v_0,...,v_n]}-h(\del\sigma)=g\circ\sigma-f\circ\sigma-h(\del\sigma),\nn
\end{align}}where at the last step, by inspecting the special case $n=2$ and then using induction, we see that
\[
\textstyle\sum_{i=1}^n(-1)^{2i}F\circ(\sigma\times 1_I)|_{[v_0,...,\wh{v}_i,w_i,...,w_n]}+\sum_{i=0}^{n-1}(-1)^{2i+1}F\circ(\sigma\times 1_I)|_{[v_0,...,v_i,\wh{w}_i,...,w_n]}=0.
\]
Hence, ~$g_\#-f_\#=\del\circ h+h\circ\del$.
\end{proof}

\begin{crl}\label{HtyInvCrl1}
Let $X,Y$ be spaces. If $X\simeq Y$, then $C_\#(X)\simeq C_\#(Y)$, and so $H_\ast(X)\cong H_\ast(Y)$.
\end{crl}

In order to relate homologies of spaces to those of their subspaces, we often make the following definitions (in which the role of the module $M$ in Definition \ref{SingCCdef1} is suppressed, i.e., $M$ is replaced with $\Integer$ simply for convenience of expression).

\begin{dfn}[\textcolor{blue}{
\index{Pairs of spaces}{Pairs of spaces},
\index{Map! of pairs}{Map of pairs},
\index{Homotopy! topological}{Homotopy},
\index{Homotopy! of pairs}{Homotopy of pairs},
\index{Homotopic maps of pairs}{Homotopic maps of pairs},
\index{Homotopy! equivalent pairs}{Homotopy equivalent pairs},
\index{Homotopy! equivalence of pairs}{Homotopy equivalence of pairs}}]
A \ul{pair of spaces} is an ordered pair $(X,A)$ in which $X$ is a space and $A\subset X$. A \ul{map of pairs} $f:(X,A)\ra(Y,B)$ is a map $f:X\ra Y$ such that $f(A)\subset B$ (i.e., $f|_A:A\ra B$).

Two maps of pairs $f,g:(X,A)\ra (Y,B)$ are \ul{homotopic as maps of pairs} (written $f\simeq_{mp} g$) if they are homotopic through a \ul{homotopy of pairs} $H:(X,A)\times I\ra(Y,B)$, defined to be a \ul{homotopy} (i.e., continuous map) $H:X\times I\ra Y$ such that $H|_{X\times\{t\}}:(X,A)\ra(Y,B)$ is a map of pairs for each $t\in I$.

Two pairs $(X,A)$ and $(Y,B)$ are \ul{homotopy equivalent}, written $(X,A)\simeq(Y,B)$, if there exists a \ul{homotopy equivalence of pairs} $f:(X,A)\ra(Y,B)$, defined to be homotopy equivalence $f:X\ra Y$ with homotopy inverse a map of pairs $g:(Y,B)\ra(X,A)$ such that $f\circ g\simeq_{mp}1_Y$ and $g\circ f\simeq_{mp} 1_X$.
\end{dfn}

\begin{dfn}[\blue{
\index{Relative singular chain complex}{Relative singular chain complex},
\index{Relative singular homology}{Relative singular homology},
\index{Local! singular homology}{Local singular homology}}]\label{SingCCdef2}
Let $(X,A)$ be a pair of spaces. Since the simplices satisfy $S_n(A)\subset S_n(X)$ for each $n$, the inclusion $i:A\hookrightarrow X$ induces (by its injectivity) an inclusion $i_\#:C_\#(A)\hookrightarrow C_\#(X)$, i.e., $C_\#(A)\subset C_\#(X)$. Therefore, we have an exact sequence of singular chain complexes $0\ra C_\#(A)\hookrightarrow C_\#(X)\sr{\pi}{\twoheadrightarrow} C_\#(X,A)\ra 0$, where
\[
C_\#(X,A):=C_\#(X)/C_\#(A)
\]
is called the \ul{relative singular chain complex} of $(X,A)$. The \ul{$n$th relative singular homology group} of $(X,A)$ is defined as
\[
H_n(X,A):=H_n\big(C_\#(X,A)\big).
\]
The homology group $H_n(X,A)$ is equivalently generated/spanned by classes $c:=[z]\in H_n(X)$, of $n$-cycles $z\in Z_n(X)\subset C_\#(X)$, that are not contained in $A$ in the sense that
\[
\textstyle z=\sum_\sigma n_\sigma\sigma=\sum_{i=1}^m n_i\sigma_i~~\txt{for}~~\sigma_i\in S_n(X)~~\txt{such that}~~\big(\bigcup_i\sigma_i\big)\cap(X-A)\neq\emptyset.
\]
We also define the \ul{$n$th local singular homology group} of $(X,A)$, or of $X$ at $A\subset X$, by
\[
H_n(X|A):=H_n\big(C_\#(X|A)\big),~~~~C_\#(X|A):=C_\#(X,X-A).
\]
The homology group $H_n(X|A)$ is equivalently generated/spanned by classes $c:=[z]\in H_n(X)$, of $n$-cycles $z\in Z_n(X)\subset C_\#(X)$, that intersect $A$ in the sense that
\[
\textstyle z=\sum_\sigma n_\sigma\sigma=\sum_{i=1}^m n_i\sigma_i~~\txt{for}~~\sigma_i\in S_n(X)~~\txt{such that}~~\big(\bigcup_i\sigma_i\big)\cap A\neq\emptyset.
\]
\end{dfn}

\begin{rmk}
Observe that associated with any triple of spaces $(X,A,B)$, i.e., $X\supset A\supset B$, is an exact sequence of relative singular chain complexes given by
\bea
\textstyle 0\ra C_\#(A,B)\hookrightarrow C_\#(X,B)\sr{\pi}{\twoheadrightarrow} {C_\#(X,B)\over C_\#(A,B)}\cong C_\#(X,A)\ra 0.\nn
\eea
Similarly, further associated with the triple of spaces $(X,X-B,X-A)$, i.e., $X\supset X-B\supset X-A$, is an exact sequence of local singular chain complexes given by
\[
\textstyle 0\ra C_\#(X-B|A)\hookrightarrow C_\#(X|A)\sr{\pi}{\twoheadrightarrow} {C_\#(X|A)\over C_\#(X-B|A)}\cong C_\#(X|B)\ra 0.
\]
\end{rmk}

\section{Presheafs and Sheafs}\label{PreshShSec}
The familiar notion of a sheaf here is related (as will be indicated in Proposition \ref{TopSysShPrp}) to the notion of a topological system from section \ref{TopSyGeoSec}.

\subsection{Basic definitions and the protosheaf representation}
\begin{dfn}[\textcolor{blue}{
\index{TopIndex-category}{TopIndex-category of a space},
\index{Presheaf on a space}{Presheaf on a space},
\index{Sections}{Sections},
\index{Restriction! morphisms}{Restriction morphisms},
\index{Sheaf on a space}{Sheaf on a space},
\index{Gluing! condition}{Gluing condition},
\index{Compatible family of sections}{Compatible family of sections}}]
Let $\C$ be a category and $X$ a space. The \ul{topological indexing category} (\ul{TopIndex-category}) of $X$, written $\I_X$, is the poset of open subsets of $X$ ordered by inclusion, i.e., $\I_X$ is the category with objects ~$\Ob\I_X:=\{\txt{open sets}~U\subset X~\txt{ordered such that}~U\leq V~\txt{iff}~U\subset V\}$ ~and morphisms
\[Hom_{\I_X}(U,V):=
\left\{
  \begin{array}{ll}
    \!\!\!\!\bt\big\{U\ar[r,hook,"{i_{UV}}"]& V\big\}\et, & \txt{if}~U\leq V, \\
    \emptyset, & \txt{otherwise}.
  \end{array}
\right.\]
A \ul{presheaf} on $X$ with values in $\C$ (equiv., in $\C$ over $X$) is an $\I_X$-cosystem in $\C$ (i.e., a cofunctor)
\[
S:\I_X\ra\C,~~\bt U\ar[r,hook,"{i_{UV}}"]& V\et~~\mapsto~~\bt S(V)\ar[r,"{S(i_{UV})}"] & S(U)\et.
\]
If $\C\subset Sets$ is a category whose objects are sets, the elements $s\in S(U)$ are called \ul{sections} of $S$ over $U$ (in the sense that wlog, wrt common applications, they may be viewed as maps $s:U\ra T_U$, for some target $T_U$ (\blue{footnote}\footnote{Say via an imbedding or otherwise. The ``\emph{sections}'' terminology (along with ``\emph{restriction morphisms}'' terminology) is used mainly because in practice, $\C$ is usually a subcategory of Top (a subcategory of Sets) and there exists a surjective continuous map $p_S:E_S\ra X$ such that for each open set $U\subset X$, $S(U)$ is the set of continuous local sections $s:U\subset X\ra E_S$ of $p_S$ (i.e., $p_S\circ s=id_U$), in which case, $S(i_{UV})(s):=s|_U$ is the restriction of a map in the usual sense. In the associated triple $(E_S,p_S,X)$, $E_S$ is called the \ul{total space} (or \ul{sheaf space}), $p_S$ the \ul{projection}, and $X$ the \ul{base space}. Depending on the intended application, the projection $p_S$ is often given other properties besides continuity and surjectivity, including further structuring/categorizing of the fibers $E_x:=p_S^{-1}(x)$ for all $x\in X$. The main point (as shown in Theorem \ref{ProtShRepThm}) is that for certain practical purposes, given $(E_S,p_S,X)$ only, one can reconstruct $S:\I_X\ra\C$ uniquely, and likewise, given $S:\I_X\ra\C$ only (along with appropriate refinements), one can reconstruct $(E_S,p_S,X)$ uniquely, and so $(E_S,p_S,X)$ and $S:\I_X\ra\C$ become different but equivalent ways of specifying the same underlying structure, namely, a \ul{(pre)sheaf}. In particular, in the process of recovering or reconstructing $(E_S,p_S,X)$ from $S:\I_X\ra\C$, for each $x\in X$, the \ul{fiber} $E_x$ is related in some way to (and so might be obtainable from) an object $S_x$ in $\C$ called the \ul{stalk} of $S:\I_X\ra \C$ at $x$.
})). (\blue{footnote}\footnote{In the case where $\C=$ Rings, the ring operations are pointwise addition and multiplication ~$+,\cdot:S(U)\times S(U)\ra S(U)$ given respectively by ~$(ss')(u):=s(u)s'(u)$~ and ~$(s+s')(u):=s(u)+s'(u)$,~ for $u\in U$.}). The system transition morphisms
\bea
S(i_{UV}):S(V)\ra S(U),~~V\sr{s}{\ral}T_V~~\mapsto~~U\sr{s|_U}{\ral}T_U\subset T_V,~~~\txt{for}~~~U\subset V,\nn
\eea
are called \ul{restriction morphisms}.

If $\C\subset Set$ is a category whose objects are sets,  then a presheaf $S:\I_X\ra\C$ is called a \ul{sheaf} if it satisfies the following condition (called ``\ul{gluing condition}''):
\bit[leftmargin=0.6cm]
\item\ul{Gluing condition}: For any open set $U\subset X$ and any open cover $\{U_\al\}_{\al\in A}$ of $U$ (i.e., {\small $U\subset\bigcup U_\al$}), if a family of sections {\small $\{s_\al\in S(U_\al)\}_{\al\in A}$} agree at intersections (hence called \ul{compatible}) in the sense that they satisfy
\bea
s_\al|_{U_\al\cap U_\beta}:=S(i_{U_\al\cap U_\beta,U_\al})(s_\al)=S(i_{U_\al\cap U_\beta,U_\beta})(s_\beta)=:s_\beta|_{U_\al\cap U_\beta},~~\txt{for each pair $\al,\beta\in A$},\nn
\eea
then there exists a unique section $s\in S(U)$ such that $s|_{U_\al}:=S(i_{U_\al U})(s)=s_\al$ for all $\al$.

Alternatively, for any collection of sections $\S\subset\big\{s\in S(U)~|~U\subset X~\txt{open}\big\}$ that agree at intersections (hence called \ul{compatible}) in the sense that (with $U_s:=\dom s$ and $U_0:=\bigcup_{s\in\S} U_s$)
\bea
s|_{U_s\cap U_{s'}}:=S(i_{U_s\cap U_{s'},U_s})(s)=S(i_{U_s\cap U_{s'},U_{s'}})(s')=:s'|_{U_s\cap U_{s'}},~~\txt{for each pair $s,s'\in\S$},\nn
\eea
there exists a unique section $s_0\in S(U_0)$ such that $s_0|_{U_s}:=S(i_{U_sU_0})(s_0)=s$ for all $s\in\S$.
\eit
\end{dfn}

\begin{notation}[\blue{
\index{Membership of a restriction}{Membership of a restriction},
\index{Domain! of a (pre)sheaf section}{Domain of a (pre)sheaf section}}]
Let $S:\I_X\ra\C\subset Sets$ be a presheaf. If $s\in S(U)$ and $U'\subset U$, then for convenience (with the \ul{membership of a restriction}) $s\in S(U')$ will simply mean $s|_{U'}\in S(U')$. We also define the \ul{domain of a (pre)sheaf section} $s\in S(U)$ to be $\dom s:=U$ (even though we do not know whether $s$ is a map). Consequently, $s\in S(U)$ means the same thing as $U\subset\dom s$, i.e., $s\in S(U)$ $\iff$ $U\subset\dom s$.
\end{notation}

\begin{dfn}[\blue{
\index{Local! TopIndex-category}{Local TopIndex-category of a space},
\index{Local! section collection}{Local section collection},
\index{Stalk of a (pre)sheaf}{Stalk of a (pre)sheaf},
\index{Germs! of a (pre)sheaf}{Germs of a (pre)sheaf}}]
Let $X$ be a space and $\I_X$ its TopIndex-category. If $A\subset X$, the \ul{local TopIndex-category} of $X$ at $A$, written $\I_{X,A}\subset\I_X$, is the full subcategory whose objects are open sets $U\subset X$ containing $A$. If $x\in X$ the \ul{local TopIndex-category} of $X$ at $x$ is $\I_{X,x}:=\I_{X,\{x\}}$. Given a (pre)sheaf $S:\I_X\ra S$, $A\subset X$, and $x\in X$, the \ul{local section collection} of $S$ at $A$ is $S_A(X):=\bigsqcup_{A\subset U\subset X}S(U)$ and the \ul{local section collection} of $S$ at $x$ is $S_x(X):=S_{\{x\}}(X)$.

Given a (pre)sheaf ~$S:\I_X\ra\C$,  the \ul{stalk} of $S$ at $x$ is the colimit (\blue{footnote}\footnote{We can more generally define the stalk $S_A$ of $S$ over a set $A\subset X$ in the same way, i.e.,
\bea
\textstyle S_A:=\varinjlim S|_{\I_{X,A}}\sr{\substack{\txt{if}~\C\subset Sets}}{\cong}{\big\{(U,s)~|~U\supset A,~s\in S(U)\big\}\over\sim}=\big\{[(U,s)]~|~U\supset A,~s\in S(U)\big\},\nn
\eea
where $(U,s)\sim(U',s')$ iff $s|_V=s'|_V$ for some neighborhood $A\subset V\subset U\cap U'$.
})
\bea
\textstyle S_x:=\varinjlim S|_{\I_{X,x}}\sr{\substack{\txt{if}~\C\subset Sets}}{\cong}{\big\{(U,s)~|~x\in U,~s\in S(U)\big\}\over\sim}=\big\{[(U,s)]~|~x\in U,~s\in S(U)\big\},\nn
\eea
where $(U,s)\sim(U',s')$ iff $s|_V=s'|_V$ for some open neighborhood $x\in V\subset U\cap U'$. (See the diagram below.)
\[
\adjustbox{scale=0.8}{%
\bt
S(V)\ar[ddr,bend right,"f_V"']\ar[dr,hook,"q_V"]\ar[rr,"{S(i_{UV})}"] && S(U)\ar[dl,hook,"q_U"']\ar[ddl,bend left,"f_U"]\\
 & S_x\ar[d,dashed,"f"] & \\
 & Q &
\et}~~~~
\begin{minipage}{12cm}\footnotesize
Compatibility: ~$q_V(s)=q_U(s|_U)$~ and ~$f_V(s)=f_U(s|_U)$,~ for ~$V\supset U\ni x$,\\~\\
ensures the following maps are well defined.\\~\\
$q_U(s):=[(U,s)]$, ~~$f([(U,s)]):=f_U(s)$,~~ ~$\Ra$~ ~~$f_V=fq_V$,~~ ~~$f_U=fq_U$.
\end{minipage}
\]
The elements of the stalk $S_x$ (i.e., the equivalence classes) are called \ul{germs} of $S$ at $x$. By construction, it is clear that a \emph{germ} of $S$ at $x$ (i.e., an element of $S_x$) can also be viewed as consisting of sections that are in \emph{nontrivial pairwise-agreement} near $x$ (i.e., in a nonempty neighborhood of $x$), in the sense that
\[
\textstyle S_x:=\varinjlim S|_{\I_{X,x}}\cong{S_x(X)\over\sim}=\left\{[s]:s\in S_x(U)\right\},
\]
where $s\sim s'$ iff there exists a nbd $x\in U\subset\dom s\cap\dom s'$ such that $s|_U=s'|_U$. (\blue{footnote}\footnote{As before, we can more generally define the stalk $S_A$ of $S$ over a set $A\subset X$ in the same way, i.e.,
\bea
\textstyle S_A:=\varinjlim S|_{\I_{X,A}}\cong{S_A(X)\over\sim}=\big\{[s]:s\in S_A(X)\big\},\nn
\eea
where $s\sim s'$ iff $s|_V=s'|_V$ for some neighborhood $A\subset V\subset\dom s\cap\dom s'$.
}). The latter representation of the stalk $S_x$ has the following slightly more compact (but possibly not any more useful) explicit look:
\[
\adjustbox{scale=0.8}{%
\bt
S(V)\ar[ddr,bend right,"f_V"']\ar[dr,hook,"q_V"]\ar[rr,"{S(i_{UV})}"] && S(U)\ar[dl,hook,"q_U"']\ar[ddl,bend left,"f_U"]\\
 & S_x\ar[d,dashed,"f"] & \\
 & Q &
\et}~~~~
\begin{minipage}{12cm}\footnotesize
Compatibility: ~$q_V(s)=q_U(s|_U)$~ and ~$f_V(s)=f_U(s|_U)$,~ for ~$V\supset U\ni x$,\\~\\
ensures the following maps are well defined.\\~\\
$q_V(s):=[s]$, ~~$f([s]):=f_{\dom(s)}(s)$,~~ ~$\Ra$~ ~~$f_V=fq_V$,~~ ~~$f_U=fq_U$.
\end{minipage}
\]
\end{dfn}
\begin{dfn}[\blue{
\index{Morphism of! (pre)sheafs ((pre)sheaf-morphism)}{Morphism of (pre)sheafs ((pre)sheaf-morphism)},
\index{Category of! (pre)sheafs}{Category of (pre)sheafs}}]
Let $X$ be a space. Given (pre)sheafs $S,S':\I_X\ra\C$, a \ul{morphism of (pre)sheafs} (or \ul{(pre)sheaf-morphism}) $\eta:S\ra S'$ is a morphism of systems. More explicitly, for all open sets $U\subset V\subset X$, we have a commutative diagram:
\[
\adjustbox{scale=0.9}{\bt U\ar[d,hook,"i_{UV}"]\\ V\et}~~~~
\adjustbox{scale=0.9}{\bt
S(U)\ar[rr,"\eta_U"] && S'(U) \\
S(V)\ar[u,"{|_U:=S(i_{UV})}"']\ar[rr,"\eta_V"] && S'(V)\ar[u,"{|_U:=S'(i_{UV})}"']
\et}~~~~~~~~\eta_V(s)|_U=\eta_U(s|_U),~~\txt{for}~~s\in S(V).
\]
The \ul{category of presheafs} (resp. \ul{category of sheafs}) on $X$ with values in $\C$ is the category $\C^{\I_X}$ (resp. $\C^{\I_X,Sh}$) whose objects are presheafs (resp. sheafs) on $X$ and whose morphisms are morphisms of presheafs (resp. sheafs).
\end{dfn}

We will not be making explicit use of the following observation but it gives a potentially useful alternative way of describing/explaining/motivating the gluing condition for a sheaf.

\begin{rmk}[\blue{The sheaf gluing condition as a colimit}]
Let $X$ be a space, $S:\I_X\ra\C$ a sheaf, and $\{s_\al\in S(U_\al)\}_{\al\in A}$ a compatible family of sections of $S$. With $U:=\bigcup_\al U_\al$, let SEC$(U)$ be the category with objects given by all sections $\{(V,t)~|~V\subset U,~t\in S(V)\}$ and morphisms $f:(V,t)\ra(V',t')$ given by commutative diagrams
\[\adjustbox{scale=0.8}{\bt
V\ar[d,"t"]\ar[r,"f"] & V'\ar[d,"t'"] \\
T_U\ar[r,"id"] & T_U
\et}~~~~t'f=t
\]
The gluing condition (which is meant to mimic/retain our ability to define a continuous map piecewise using an open cover) can be seen as arising through a colimit in SEC$(U)$. Specifically, the gluing condition suggests that the compatible family of sections $\{s_\al\in S(U_\al)\}_{\al\in A}$ defines a section given by the colimit
\[
\textstyle(U,s):=\varinjlim_{\al\in A} (U_\al,s_\al),~~~~\txt{with}~~s\in S(U),
\]
which is as usual explicitly described as follows:
\[
\begin{minipage}{9cm}\footnotesize
Given any maps $f_\al:(U_\al,s_\al)\ra (V,t)$ in diagram from
\[\bt
U_\al\ar[d,"s_\al"]\ar[r,"f_\al"] & V\ar[d,"t"] \\
T_U\ar[r,"id"] & T_U
\et~~~~tf_\al=s_\al
\]
that satisfy $f_\al|_{U_\al\cap U_\beta}=f_\beta|_{U_\al\cap U_\beta}$, there exists a unique map
$f:(U,s)\ra (V,t)$ in diagram form
\[\bt
U\ar[d,"s"]\ar[r,"f"] & V\ar[d,"t"] \\
T_U\ar[r,"id"] & T_U
\et~~~~tf=s
\]
that makes the main diagram on the right commutative.
\end{minipage}
~~~~
\adjustbox{scale=0.8}{%
\bt
 (U_\al\cap U_\beta,s_\al)\ar[d,hook]\ar[rr,equal] && (U_\al\cap U_\beta,s_\beta)\ar[d,hook]\\
 (U_\al,s_\al)\ar[ddr,bend right,"f_\al"']\ar[dr,hook,"q_\al"] && (U_\beta,s_\beta)\ar[dl,hook,"q_\beta"']\ar[ddl,bend left,"f_\beta"]  \\
 & (U,s)\ar[d,dashed,"\exists!~f"] & \\
 & (V,t) &
\et}
\]
\end{rmk}

\begin{rmk}[\blue{Isolating uniqueness in the gluing condition}]\label{ShfGluUniq}
Let $S:\I_X\ra\C$ be a sheaf. The gluing condition implies that for any open set $U\subset X$ and open subsets $\{U_i\}_{i\in I}$ of $X$ that cover $U$, i.e., $U\subset\bigcup_i U_i$, if two sections $s,s'\in S(U)$ satisfy $s|_{U_i\cap U}=s'|_{U_i\cap U}$ for all $i\in I$, then $s=s'$. Indeed, with $s_i:=s|_{U_i\cap U}\in S(U_i\cap U)$ and $s_i':=s'|_{U_i\cap U}\in S(U_i\cap U)$, we have
\[
\left.
  \begin{array}{l}
   (1)~ s_i|_{U_i\cap U_j}=s|_{U_i\cap U_j}=s_j|_{U_i\cap U_j}~~\txt{for all}~~i,j\in I \\
   (2)~ s'_i|_{U_i\cap U_j}=s'|_{U_i\cap U_j}=s'_j|_{U_i\cap U_j}~~\txt{for all}~~i,j\in I \\
   (3)~ s_i=s_i'~~\txt{for all}~~i\in I
  \end{array}
\right.
\]
Hence, the collection $\{s_i\}_{i\in I}=\{s'_i\}_{i\in I}$ glues into a unique section $s=s'\in S(U)$.
\end{rmk}

\begin{lmm}[\blue{\cite[Lemma 5.69, p.284]{rotman2009}}]\label{MorOfShfUniq}
Let $S,S':\I_X\ra\C\subset Sets$ be presheafs and $f,g:S\ra S'$ morphisms of presheafs. If (i) $S'$ is a sheaf and (ii) $f,g$ agree on every stalk of $S$ in the sense that the limit morphisms $f_x:=\varinjlim f|_{\I_{X,x}}:S_x\ra S'_x$ and $g_x:=\varinjlim g|_{\I_{X,x}}:S_x\ra S'_x$ are equal, i.e., $f_x=g_x:S_x\ra S'_x$, for all $x\in X$, then $f=g$.
\end{lmm}
\begin{proof}
Let $U\subset X$ be open. We need to show $f_U=g_U:S(U)\ra S'(U)$. By construction/functoriality of the colimit, we have commutative diagrams (in which $q_U(s):=[(U,s)]$ and $q'_U(s'):=[(U,s')]$):
\[\bt
S(U)\ar[d,hook,"q_U"]\ar[drr,bend left=10,"q'_Uf_U"description]\ar[rr,"f_U"] && S'(U)\ar[d,hook,"q'_U"] \\
S_x\ar[rr,dashed,"f_x"] && S_x'
\et~~~~~~~~
\bt
S(U)\ar[d,hook,"q_U"]\ar[drr,bend left=10,"q'_Ug_U"description]\ar[rr,"g_U"] && S'(U)\ar[d,hook,"q'_U"] \\
S_x\ar[rr,dashed,"g_x"] && S_x'
\et~~~~~~~~\txt{for each}~~~~x\in U,
\]
and so the morphisms $f_x,g_x$ are given by
\[
f_x:S_x\ra S'_x,~[(O_x,s_x)]\mapsto[(O_x,f_{O_x}(s_x))],~~~~g_x:S_x\ra S'_x,~[(O_x,s_x)]\mapsto[(O_x,g_{O_x}(s_x))].\nn
\]
Recall that we have $[(O_x,f_{O_x}(s_x))]=[(O_x,g_{O_x}(s_x))]$ iff there exists an open nbd $x\in W_x\subset O_x$ such that
\begin{align}
\label{StalkEqualEq1}f_{W_x}(s_x|_{W_x})\sr{(a)}{=}f_{O_x}(s_x)|_{W_x}=g_{O_x}(s_x)|_{W_x}\sr{(a)}{=}g_{W_x}(s_x|_{W_x}),
\end{align}
where the equalities $(a)$ hold because $f,g$ are natural transformations.

Given $s\in S(U)$, $x\in U$, and any open nbd $x\in U_x\subset U$, let $s_x:=s|_{U_x}$ (where it is clear that $s_x|_{U_x\cap U_{x'}}=s_{x'}|_{U_x\cap U_{x'}}$). Then by hypotheses $[(U_x,f_{U_x}(s_x))]=[(U_x,g_{U_x}(s_x))]$ for all $x\in U$. Therefore, with $x\in W_x\subset U_x$ as before, the equalities $f_{W_x}(s_x|_{W_x})\sr{(a)}{=}f_U(s)|_{W_x}$ and $g_{W_x}(s_x|_{W_x})\sr{(a)}{=}g_U(s)|_{W_x}$ imply the following: with $t_x:=s_x|_{W_x}$ and $t_{x'}:=s_{x'}|_{W_{x'}}$,
\bea
f_{W_x}(t_x)|_{W_x\cap W_{x'}}\sr{(a)}{=}f_{W_{x'}}(t_{x'})|_{W_x\cap W_{x'}}\sr{(\ref{StalkEqualEq1})}{=}g_{W_{x'}}(t_{x'})|_{W_x\cap W_{x'}}\sr{(a)}{=}g_{W_x}(t_x)|_{W_x\cap W_{x'}},~~\txt{for all}~~x,x'\in U,\nn
\eea
where the equalities (a) again hold because $f,g$ are natural transformations. Since $\{W_x\}_{x\in U}$ is an open cover of $U$, the gluing condition for the sheaf $S'$ implies the collections $\{f_{W_x}(s_x|_{W_x})\}_{x\in U}=\{f_U(s)|_{W_x}\}_{x\in U}$ and $\{g_{W_x}(s_x|_{W_x})\}_{x\in U}=\{g_U(s)|_{W_x}\}_{x\in U}$ determine the same unique section in $S'(U)$, and so $f_U(s)=g_U(s)$ for any $s\in S(U)$, i.e., $f_U=g_U$.
\end{proof}

The following result is related to \cite[Theorem 5.68, p.281]{rotman2009}.

\begin{thm}[\textcolor{blue}{\index{Protosheaf representation theorem}{Protosheaf (sheaf prototype) representation theorem}}]\label{ProtShRepThm}
Let $X$ be a space and $S:\I_X\ra\C\subset Sets$ a presheaf. We have the following: There exists
\bit[leftmargin=0.9cm]
\item[(i)] a unique protosheaf ~$p_S:E_S\ra X$~ such that each fiber ~$E_{S,x}:=p_S^{-1}(x)\cong S_x$~ in $\C$. (\blue{footnote}\footnote{Recall that a protosheaf is defined to be a surjective local homeomorphism (and hence a surjective open continuous map).})
\item[(ii)] a presheaf morphism $\phi:S\ra S^\ast$, where $S^\ast$ is the sheaf of continuous sections of {\footnotesize $p_S:E_S\ra X$} given by
\[
S^\ast:\I_X\ra\C,~~U\sr{i_{UV}}{\hookrightarrow}V~~\mapsto~~\Gamma(V,E_S)\sr{|_U}{\ral}\Gamma(U,E_S):=\{\txt{continuous}~~t:U\ra E_S~|~p_S\circ t=id_U\}.
\]
\eit

Furthermore, if $S$ is a sheaf, then we also have the following:
\bit[leftmargin=0.9cm]
\item[(iii)] The morphism $\phi:S\ra S^\ast$ is an isomorphism of sheafs: For each open set $U\subset X$, we have a bijection
\[
S(U)\longleftrightarrow S^\ast(U):=\Gamma(U,E_S).
\]
\item[(iv)] For any sheaf-morphism $\eta:S\ra S'$ (where $S,S':\I_X\ra\C$), there exists a unique continuous map $\wt{\eta}:E_S\ra E_{S'}$ such that the following diagram commutes (i.e., $p_{S'}\circ\wt{\eta}=id_X\circ p_S=p_S$):
\[\adjustbox{scale=0.9}{\bt
E_S\ar[d,"p_S"]\ar[rr,"\wt{\eta}"] && E_{S'}\ar[d,"p_{S'}"] \\
X\ar[rr,"id_X"] && X
\et}~~~~~~p_{S'}\circ\wt{\eta}=id_X\circ p_S=p_S,
\]
i.e., the restriction of $\wt{\eta}$ on each fiber $E_{S,x}$, $x\in X$, gives a morphism $\wt{\eta}|_{E_{S,x}}:E_{S,x}\ra E_{S',x}$ in $\C\cap Top$.

The converse is true as well: Given sheafs $S,S':\I_X\ra\C$, any continuous map $\wt{\eta}:E_S\ra E_{S'}$ (such that the above diagram commutes) induces a unique sheaf-morphism $\eta:S\ra S'$.
\eit
\end{thm}
\begin{proof}
(i) Let $E=E_S:=\bigsqcup_{x\in X}S_x$ as a space with the topology generated by subsets of the form
\begin{align}
\textstyle\langle U,s\rangle:&=\textstyle\bigsqcup_{x\in U}\{q_O(s):x\in O,~s\in S(O)\},~~~~\txt{for open sets $U\subset X$ and sections $s\in S(U)$},\nn\\
~~~~&=\textstyle\bigsqcup_{x\in U}\{[(O,s)]:x\in O,~s\in S(O)\}=\bigsqcup_{x\in U}S_{x,s},\nn
\end{align}
where $S_{x,s}:=\{[(O,s)]:x\in O,~s\in S(O)\}=\langle U,s\rangle\cap S_x$ = singleton (i.e., one element set). (\blue{footnote}\footnote{It might be helpful to notice that $S_{x,s}$ consists of a single element (which makes $\langle U,s\rangle$ a sheet of the space $E$), even though this is not strictly necessary since we might desire the following more general situation (of a preferred topology $\T_{S_x}$ on $S_x$ other than the discrete topology) especially if $S$ is a presheaf that is not necessarily a sheaf: Let $F$ be a space that gives every fiber $S_x$ its topology through a homeomorphism $h_x:F\ra S_x$. Let $\U\subset\P(X)$ be an open cover of $X$ and $\V\subset\P(F)$ an open cover of $F$. Then we could consider (sub)base sets for the topology of $E$, say, of the form $\langle U,V\rangle=\bigsqcup_{x\in U}\langle U,V\rangle\cap S_x:=\bigsqcup_{x\in U}S_{x,V}\cong U\times V$, for $U\in\U$ and $V\in\V$, such that $S_{x,V}:=h_x(V)=\{h_x(v)=[(O_v,s_v)]\in S_x:v\in V\}=\bigsqcup_{v\in V}S_{x,s_v}=\bigsqcup_{v\in V}\langle U,s_v\rangle\cap S_x$.}). That is, an open subset of $E$ is a union of finite intersections of the form
\bea
&&\textstyle\langle U_1,s_1\rangle\cap \langle U_2,s_2\rangle\cap\cdots\cap \langle U_n,s_n\rangle=\bigsqcup_{x\in U_1\cap\cdots\cap U_n}S_{x,s_1}\cap\cdots\cap S_{x,s_n},~~~~\txt{for}~~n\geq 1,\nn\\
&&~~~~=\langle U_1\ast U_2\ast\cdots\ast U_n,s\rangle,~~~~\txt{for some}~~s\sim s_1\sim s_2\sim\cdots\sim s_n~~\txt{wrt}~~U_1\ast U_2\ast\cdots\ast U_n,\nn
\eea
where the set $V:=U_1\ast U_2\ast\cdots\ast U_n:=\{x\in U_1\cap\cdots\cap U_n:S_{x,s_1}=S_{x,s_2}=\cdots=S_{x,s_n}\}$ is open in $X$ due to the following: Let $x\in V$. Then there is an open set $V_x\ni x$ such that $s_1|_{V_x}=s_2|_{V_x}=\cdots=s_n|_{V_x}$ (and vice versa if needed). If $x'\in V_x\cap U_1\cap\cdots\cap U_n$, then it is clear (by the openness of the intersection) that we also have an open set $x'\in V_{x,x'}\subset V_x$ such that $s_1|_{V_{x,x'}}=s_2|_{V_{x,x'}}=\cdots=s_n|_{V_{x,x'}}$ (i.e., $S_{x',s_1}=S_{x',s_2}=\cdots=S_{x',s_n})$, and so $x'\in V$ as well. Since $x\in V_x\cap U_1\cap\cdots\cap U_n\subset V$ for all $x\in V$, it follows that $V\subset X$ is open.

Therefore, the sets $\langle U,s\rangle$ form a base for the topology of $E$. Define the surjective map
\[
p=p_S:E\ra X,~[(U,s)]\mapsto x~~\txt{if}~~[(U,s)]\in S_x,~~\txt{i.e.,}~~p|_{S_x}:=\txt{constant}=x~~\txt{for all}~~x\in X.
\]
For any open set $U\subset X$ and any section $s\in S(U)$, we have
\bea
&&\textstyle p^{-1}(U)=E\cap p^{-1}(U)=\big(\bigcup_\al\langle U_\al,s_\al\rangle\big)\cap p^{-1}(U)=\left(\bigcup_\al\bigsqcup_{x\in U_\al}S_{x,s_\al}\right)\cap p^{-1}(U)\nn\\
&&\textstyle ~~~~=\bigcup_\al\bigsqcup_{x\in U_\al}S_{x,s_\al}\cap p^{-1}(U)=\bigcup_\al\bigsqcup_{x\in U_\al\cap U}S_{x,s_\al}=\bigcup_\al\langle U_\al\cap U,s_\al\rangle,\nn
\eea
and so $p$ is continuous. Moreover, for each $e\in E$, if $e\in \langle U_e,s_e\rangle$, then $p(\langle U_e,s_e\rangle)=U_e$ shows $p$ is a locally open map (hence an open map). Also, since each $S_{x,s}$ is a singleton, we get an open imbedding
\[
p|_{\langle U_e,s_e\rangle}:\langle U_e,s_e\rangle\subset E\ra X,\nn
\]
and so $p:E\ra X$ is a protosheaf.

(ii) We can define a map $\phi_U:S(U)\ra\Gamma(U,E)$ as follows: Using the axiom of choice or otherwise, any section $s\in \bigcup_{U\in\T_X} S(U)$ can be associated the map ~$\phi_{\dom s}(s):=\wt{s}\in \Gamma\big(\dom(s),E\big)$ given by
\[
\wt{s}:\dom s\subset X\ra E,~x\mapsto \wt{s}(x)\in\langle U_x,s\rangle\cap S_x=S_{x,s},~~\txt{for an open nbd}~~x\in U_x\subset X,
\]which is continuous because for any base set $\langle U,s'\rangle\subset E$,
\begin{align}
&\wt{s}^{-1}(\langle U,s'\rangle)=\{x'\in\dom s:\wt{s}(x')\in \langle U,s'\rangle\}=\{x'\in\dom s:\wt{s}(x')\in \langle U,s'\rangle\cap\langle U_{x'},s\rangle\cap S_{x'}\}\nn\\
&~~~~=\{x'\in(\dom s)\cap U:S_{x',s}=S_{x',s'}\}=\langle(\dom s)\ast U,s\rangle=\langle(\dom s)\ast U,s'\rangle.\nn
\end{align}

The naturality diagram below holds by construction (from the definition of $\phi:S\ra S^\ast$): Given a section $s\in S(V)$ and $U\subset V$, it is clear that $\wt{s|_U}=\wt{s}|_U$.
\[\adjustbox{scale=0.9}{\bt U\ar[d,hook,"i_{UV}"]\\ V\et}~~~~
\adjustbox{scale=0.9}{\bt
S(U)\ar[rr,"\phi_U"] && \Gamma(U,E) \\
S(V)\ar[u,"{|_U:=S(i_{UV})}"']\ar[rr,"\phi_V"] && \Gamma(V,E)\ar[u,"{|_U}"']
\et}~~~~~~~~\phi_V(s)|_U=\phi_U(s|_U),~~\txt{for}~~s\in S(V).
\]

(iii) Now assume $S$ is a sheaf. To prove (iii), we will construct a bijection as follows:

{\flushleft \ul{Defining an injective map $\phi_U:S(U)\ra\Gamma(U,E),~s\mapsto\wt{s}$}}: By Remark \ref{ShfGluUniq}, if $s,s'\in S(U)$, then $s=s'$ $\iff$ there exist open nbds $\{U_x\ni x:x\in U\subset X\}_{x\in U}$ of points of $U$ such that $(U_x,s)\sim (U_x,s')$ in $S_x$ for all $x\in U$. Consequently, (using the axiom of choice or otherwise), any section $s\in \bigcup_{U\in\T_X} S(U)$ can be uniquely represented by the continuous section
\[
\wt{s}:\dom s\subset X\ra E,~x\mapsto \wt{s}(x)\in\langle U_x,s\rangle\cap S_x=S_{x,s},~~\txt{for an open nbd}~~x\in U_x\subset X,
\]
because for all $s_1,s_2\in S(U)$, we have {\small $\wt{s_1}=\wt{s_2}$ $\iff$ $\wt{s_1}(x)=\wt{s_2}(x)$} for all $x\in U$, $\iff$ {\small $(W_x,s_1)\sim(W_x,s_2)$} for all $x\in U$ (and open nbds $x\in W_x\subset U\subset\bigcup_{x\in U} W_x$), $\iff$ $s_1=s_2$.

{\flushleft \ul{Surjectivity of $\phi_U$}}: Consider any open set $U\subset X$ and any continuous section $t:U\ra E,~x\mapsto t(x)\in S_x$. Then $t(x)\in\langle U_x,s_x\rangle$ for some open neighborhood $x\in U_x$ (\blue{footnote}\footnote{Without loss of generality, we assume $U_x\subset U$, otherwise we can take the intersection $U_x\cap U$ at every step where the intersection is relevant.}) and a section $s_x\in S(U_x)$. That is,
\begin{align}
&t:U\ra E,~x\mapsto t(x)\in \langle U_x,s_x\rangle\cap S_x=S_{x,s_x},~~s_x\in S(U_x),\nn
\end{align}
where for each $x\in U$, all such $s_x\in S(U_x)$ are $\sim$-equivalent through $t(x)$. Our goal is to glue the collection $\{s_x\}_{x\in U}$ into a unique section $s\in S(U)$ such that $s|_{U_x}=s_x$ for each $x\in U$ (and hence $t=\wt{s}$).

Since {\small $\langle U_x,s_x\rangle\cap\langle U_{x'},s_{x'}\rangle=\langle U_x\ast U_{x'},s_{xx'}\rangle\neq\emptyset$ $\iff$ $s_x\sim s_{x'}$ (wrt $U_x\ast U_{x'}$), the resulting equivalence relation on both $U$ and $\{\langle U_x,s_x\rangle\}_{x\in U}$} gives disjoint equivalence classes {\small $[\langle U_x,s_x\rangle]=\{\langle U_y,s_y\rangle:y\in[x]\}$}, for $x\in U$. Thus, for each $x\in U$, we have a cover of $U_{[x]}:=\bigcup_{y\in [x]}U_y$ by open sets {\small $\{W_y:y\in W_y\subset U_x\ast U_y\}_{y\in[x]}$} such that {\small $s_y|_{W_y\cap W_{y'}}=s_{y'}|_{W_y\cap W_{y'}}$} for all $y,y'\in [x]$. By the gluing condition, the elements of each equivalence class $[\langle U_x,s_x\rangle]$ glue into a unique section $s_{[x]}\in S(U_{[x]})$ such that $s_{[x]}|_{U_y}=s_y$ for each $y\in[x]$, and so
\bea
t|_{U_{[x]}}=\wt{s_{[x]}},~~~~\txt{for each}~~x\in U.\nn
\eea
Since $U=\bigcup_{x\in U}U_{[x]}$ and the restrictions $\{t|_{U_{[x]}}\}_{x\in U}$ of $t$ (and hence the corresponding sections $\{s_{[x]}\}_{x\in U}$) agree on intersections of their domains, the sections $\{s_{[x]}\}_{x\in U}$ further glue into a unique section $s\in S(U)$ such that $s|_{U_{[x]}}=s_{[x]}$  (and so $s|_{U_x}=s|_{U_{[x]}}|_{U_x}=s_{[x]}|_{U_x}=s_x$) for each $x\in U$. Hence $t=\wt{s}$.

{\flushleft (iv)} \ul{Transfer of natural transformations}: Given sheafs $S,S':\I_X\ra\C$, consider the associated protosheafs $p:E\ra X$ and $p':E'\ra X$, where $E:=E_S=\bigsqcup_{x\in X}S_x$ and $E':=E_{S'}=\bigsqcup_{x\in X}S'_x$ (and similarly, $E_x:=E_{S,x}=S_x$ and $E'_x:=E_{S',x}=S'_x$). Let $U\subset V\subset X$ be open sets and $i_{UV}:U\hookrightarrow V$ the inclusion. Given a natural transformation $\eta:S\ra S'$, we get the following commutative diagram:
\[\adjustbox{scale=0.8}{\bt
      &&&&                             S(V)\ar[dd,"\eta_V"]\ar[ddddllll,"\phi_V"']\ar[rrrr,"S(i_{UV})"]  &&&& S(U)\ar[dd,"\eta_U"] \\ 
      &&&&                                   &&&&       \\
      &&&&                             S'(V)\ar[ddddllll,near start,"\phi'_V"']\ar[rrrr,near start,"S'(i_{UV})"] &&&& S'(U)\ar[ddddllll,"\phi'_U"'] \\
      &&&&                                   &&&&       \\
\Gamma(V,E)\ar[dd,"\eta^\ast_V=\phi'_V\eta_V\phi_V^{-1}"']\ar[rrrr,crossing over,near start,"|_{U}"]    &&&& \Gamma(U,E)\ar[dd,"\eta^\ast_U=\phi'_U\eta_U\phi_U^{-1}"']\ar[from=uuuurrrr,crossing over,near end,"\phi_U"']          &&&&       \\
                 &&&&                        &&&&       \\
\Gamma(V,E')\ar[rrrr,"|_{U}"] &&&& \Gamma(U,E')       &&&&
\et}~~~~\substack{\eta_US(i_{UV})=S'(i_{UV})\eta_V\\~\\~or~\\~\\ \eta_U(s|_U)=\eta_V(s)|_U}
\]
($\La$): The above diagram (along with Lemma \ref{MorOfShfUniq}) shows that a map $\wt{\eta}:E\ra E'$ (which gives a map $\eta^\ast_U:\Gamma(U,E)\ra\Gamma(U,E'),~t\mapsto\wt{\eta}\circ t$) that preserves fibers (in the sense its restrictions on the fibers give morphisms $\wt{\eta}|_{E_x}:E_x\ra E'_x$ for each $x\in X$) induces a unique natural transformation $\eta:S\ra S'$ through the bijections $\phi_U:S(U)\ra\Gamma(U,E)$ and $\phi'_U:S'(U)\ra\Gamma(U,E')$. (\blue{footnote}\footnote{The fiber preserving property of $\wt{\eta}$ is necessary so that the colimit functor $\varinjlim(-)|_{\I_{X,x}}$ which takes $\eta:S\ra S'$ to $\eta_x:=\wt{\eta}|_{E_x}:S_x\ra S'_x$ is well-defined.}).

($\Ra$): Conversely, we have the following: Given a natural transformation $\eta:S\ra S'$, define a map $\wt{\eta}:E\ra E'$ (where {\footnotesize $E:=E_S=\bigsqcup_{x\in X}S_x$} and {\footnotesize $E':=E_{S'}=\bigsqcup_{x\in X}S'_x$}) by the following restrictions (which correspond to morphisms in $\C$):
\[
\wt{\eta}|_{S_x}:=\eta_x:S_x\ra S'_x,~[(U,s)]\mapsto[(U,\eta_U(s))],~~~~\txt{for all}~~x\in X.
\]
Then the following diagram commutes: If $e\in E$, then with $x:=p(e)$ and $e=[(U_e,s_e)]\in S_x$, we have
\[
p'\circ\wt{\eta}(e)=p'\circ\wt{\eta}([(U_e,s_e)])=p'([(U_e,\eta_{U_e}(s_e))])=x=p(e).
\]
\[\adjustbox{scale=0.9}{\bt
E\ar[d,"p"]\ar[rr,"\wt{\eta}"] && E'\ar[d,"p'"] \\
X\ar[rr,"id_X"] && X
\et}~~~~~~p'\circ\wt{\eta}=p.
\]
With a base open set ~$B:=\langle U,s'\rangle=\bigsqcup_{x\in U}S_{x,s'}\subset E'$,~ we have
{\small\begin{align}
&\textstyle \wt{\eta}^{-1}(B)=\bigsqcup_{x\in U}\wt{\eta}^{-1}(S_{x,s'})=\bigsqcup_{x\in U}\wt{\eta}|_{S_x}^{-1}(S_{x,s'})=\bigsqcup_{x\in U}\{[(W,s)]\in S_x:[(W,\eta_W(s))]\in S_{x,s'}\}\nn\\
&\textstyle~~~~=\bigsqcup_{x\in U}\{[(W,s)]\in S_x:\eta_W(s)|_C=s'|_C,~C\subset W\cap\dom s'=W\cap U\}\nn\\
&\textstyle~~~~=\bigsqcup_{x\in U}\{[(W,s)]\in S_x:\eta_C(s|_C)=s'|_C,~C\subset W\cap U\}
=\bigsqcup_{x\in U}\{[(C,s|_C)]\in S_x:s|_C\in\eta_C^{-1}(s'|_C),~C\subset U\}\nn\\
&\textstyle~~~~\sr{(a)}{=}\bigsqcup_{x\in U}\bigcup_{C\subset U}\bigcup_{s\in S(C)\cap\eta_C^{-1}(s'|_C)}S_{x,s}
=\bigsqcup_{x\in U}\bigcup_{C\subset U}\bigcup_{s\in S(C)\cap\eta_C^{-1}(s'|_C)}\langle C,s\rangle\cap S_x\nn\\
&\textstyle~~~~\sr{(b)}{=}(\bigsqcup_{x\in U}S_x)~\cap~\bigcup_{C\subset U}\bigcup_{s\in S(C)\cap\eta_C^{-1}(s'|_C)}\langle C,s\rangle=\bigcup_{C\subset U}\bigcup_{s\in S(C)\cap\eta_C^{-1}(s'|_C)}\langle C,s\rangle,\nn
\end{align}}where step (a) holds because $S_{x,s}=\langle C,s\rangle\cap S_x$ (with $x\in C$) is a singleton and step (b) holds because by construction $\langle C,s\rangle\cap S_x=\bigsqcup_{c\in C}S_{c,s}\cap S_x=\emptyset$ if $x\not\in C$, and so in the union we can use all $C\subset U$ (instead of only $x\in C\subset U$). This shows $\wt{\eta}:E\ra E'$ is continuous. Moreover, $\wt{\eta}$ is unique by construction (due to Lemma \ref{MorOfShfUniq}).
\end{proof}

\subsection{Sheafification, pseudo-topologies, and gluing of (pre)sheafs}
\begin{dfn}[\blue{Sheafification of a presheaf}]
Let $X$ be a space and $S:\I_X\ra\C\subset Sets$ a presheaf. The \ul{sheafification} of $S$ is (i) a presheaf-morphism $\psi:S\ra S^\ast$ (with $S^\ast$ a sheaf), such that (ii) for any presheaf-morphism $\eta:S\ra S'$ (with $S'$ a sheaf), there exists a unique sheaf-morphism $\eta_1:S^\ast\ra S'$ satisfying $\eta_1\circ\psi=\eta$.
\[\adjustbox{scale=0.9}{\bt
S\ar[d,"\eta"']\ar[rr,"\psi"] && S^\ast\ar[dll,dashed,"\eta_1"]\\
S' &&
\et}\]
\end{dfn}
The following result is related to \cite[Theorem 5.70, pp 284-285]{rotman2009}.

\begin{crl}
Let $S:\I_X\ra\C\subset Sets$ be a presheaf and $p_S:E_S\ra X$ the associated protosheaf from Theorem \ref{ProtShRepThm}(i). Then up to isomorphism of sheafs, the sheafification of $S$ is given by the presheaf-morphism $\phi:S\ra S^\ast$ from Theorem \ref{ProtShRepThm}(ii), where $S^\ast$ is the sheaf of continuous sections of $p_S$, i.e.,
\bea
S^\ast=\Gamma(-,E_S):\I_X\ra\C,~~U\sr{i_{UV}}{\hookrightarrow}V~~\mapsto~~\Gamma(V,E_S)\sr{|_U}{\ral}\Gamma(U,E_S).\nn
\eea
\end{crl}
\begin{proof}
Let $\eta:S\ra S'$ be a presheaf-morphism (with $S'$ a sheaf). As before, define $\phi:S\ra S^\ast$ and $\phi':S'\ra S'^\ast$ as follows: Given $s\in S(U)$ and $s'\in S'(U)$, define $\phi_U(s)\in\Gamma(U,E_S)$ and $\phi'_U(s)\in\Gamma(U,E_{S'})$ to be the following continuous sections (as obtained in the proof of Theorem \ref{ProtShRepThm}(ii)):
\begin{align}
&\phi_U(s):=\wt{s}:U\subset X\ra E,~x\mapsto \wt{s}(x)\in\langle U_x,s\rangle\cap S_x=S_{x,s},~~\txt{for an open neighborhood}~~x\in U_x\subset X,\nn\\
&\phi'_U(s'):=\wt{s'}:U\subset X\ra E,~x\mapsto \wt{s'}(x)\in\langle U_x,s'\rangle\cap S'_x=S'_{x,s'},~~\txt{for an open neighborhood}~~x\in U_x\subset X.\nn
\end{align}
The presheaf-morphism $\eta:S\ra S'$ induces a unique continuous map $\wt{\eta}:E_S\ra E_{S'}$ (as in Theorem \ref{ProtShRepThm}(iv)), which in turn induces a sheaf-morphism $\eta^\ast:S^\ast\ra S'^\ast$ through {\footnotesize $\eta^\ast_U:\Gamma(U,E_S)\ra\Gamma(U,E_{S'}),~t\mapsto\wt{\eta}\circ t$}.

Since $S'$ is a sheaf, $\phi'$ is a bijection (as shown in the proof of Theorem \ref{ProtShRepThm}(iii)).
\[\adjustbox{scale=0.9}{\bt
S\ar[dd,"\eta"']\ar[rrrr,"\phi"] &&&& S^\ast\ar[ddllll,dashed,"\eta_1:=\phi'^{-1}\eta^\ast"']\ar[dd,"\eta^\ast"]\\
 &&&& \\
S'\ar[rrrr,"\phi'"] &&&& S'^\ast
\et}\]
Hence, we get the desired sheaf-morphism ~$\eta_1:=\phi'^{-1}\eta^\ast:S^\ast\ra S'$.
\end{proof}

\begin{crl}[\blue{\cite[Proposition 1.1, p.64]{hartshorne1977}}]\label{MorOfShfUniq2}
Let $S,S':\I_X\ra\C\subset Sets$ be sheafs (over a space $X$). A sheaf-morphism $\eta:S\ra S'$ is an isomorphism if and only if for each $x\in X$, the map induced on stalks ~$\eta_x:=\varinjlim_{U\ni x}\eta_U=\wt{\eta}|_{S_x}:S_x\ra S'_x$~ is an isomorphism.
\end{crl}
\begin{proof}
This follows from the construction in the proof of Theorem \ref{ProtShRepThm}(iv). ($\Ra$): If $\eta$ is an isomorphism, it is clear that each $\eta_x$ is an isomorphism. ($\La$): Conversely, assume each $\eta_x$ is an isomorphism. Then the associated protosheaf-morphism $\wt{\eta}:E_S\ra E_{S'}$ is an isomorphism (otherwise, if some $e'\in E_{S'}$ is not in $\im\wt{\eta}=\bigcup_x\im\eta_x$ then $\wt{\eta}$ cannot be surjective on all fibers, and if $\wt{\eta}(e_1)=\wt{\eta}(e_2)$ for some distinct $e_1,e_2\in E_S$ then $\wt{\eta}$ cannot preserve fibers). Hence $\eta$ is an isomorphism, because by uniqueness, the sheaf-morphism $\wt{\eta}_\phi=\phi'^{-1}\wt{\eta}\phi:S\ra S'$ induced by $\wt{\eta}$ (through the sheaf-isomorphisms $\phi:S\ra S^\ast$ and $\phi':S'\ra S'^\ast$) is the same as $\eta$, up to isomorphism.
\end{proof}

The strength of the following proposition comes from Theorem \ref{ProtShRepThm}, in which if $S:\I_X\ra \C$ is a sheaf then its equivalent given by the sheaf $S^\ast:=\Gamma(-,E_S):\I_X\ra Sets$ of sections of $p_S:E_S\ra X$ features the space $E_S$ with some of its subsets (namely, the fibers of $p_S$) as objects in $\C$. One way to test the said potence would be to compare $E_S$ with the topological system $\langle S\rangle:\J\supset\I_X\ra\C$ generated by $S$.

\begin{prp}[\textcolor{blue}{Pseudo-topologies from sheafs}]\label{TopSysShPrp}
Let $X$ be a space and $S:\I_X\ra\C\subset Sets$ a presheaf. (i) With some extendability assumption, $S$ is finite colimit-closed. (ii) If $S$ has an extension $S':\I\supset\I_X\ra\C$ that maps colimits to limits, then $S$ is \ul{potentially} a sheaf. (\blue{footnote}\footnote{Here $S'$ is a cofunctor that maps colimits to limits. An example of a cofunctor with the desired property (due to hom-stability of limits/colimits) is $S':=Hom_R(-,E)$ for a ring $R$ and an $R$-module $E$, where $E$ is viewed as an object in Top (i.e., $E$ is a topological space).}). (iii) If $S$ is a sheaf, the gluing condition makes $S$ limit-closed (and so moves $S$ closer to a pseudo-topology, because of part (i) above).
\end{prp}
\begin{proof}
By Theorem \ref{ProtShRepThm}, if necessary, we can assume sections of a sheaf are continuous maps.\\
(i) By the restriction of maps, we have the cardinal inequality $S(U)\cup S(V)\lesssim S(U\cap V)$. Therefore, $S$ is \ul{finite colimit-closed} if we \ul{assume} a morphism $h:(S(U)\cup S(V))|_{U\cap V}\subset S(U\cap V)\ra Z$ on restrictions of sections uniquely extends to a morphism $h:S(U\cap V)\ra Z$. (That is, with the above \ul{assumption}, $S(U\cap V)\cong\varinjlim\{S(U),S(V)\}$ in $\C$.)
\[\adjustbox{scale=0.9}{\bt
 & S(U)\cap S(V) & \\
 S(U)\ar[from=ur,hook]\ar[dddr,bend right,"f_U"']\ar[ddr,bend right=10,"q_U"']\ar[dr,dotted,hook] &&  S(V)\ar[from=ul,hook]\ar[dl,dotted,hook]\ar[ddl,bend left=10,"q_V"]\ar[dddl,bend left,"f_V"] \\
 & S(U)\cup S(V)\ar[d,dotted,near start,"|_{U\cap V}"] & \\
 &  S(U\cap V)\ar[d,dashed,"h"] & \\
 & Z &
\et}~~~~~~~~~h(s|_{U\cap V}):=
\left\{
  \begin{array}{ll}
    f_U(s),& \txt{if}~~s\in S(U) \\
    f_V(s),& \txt{if}~~s\in S(V)
  \end{array}
\right.
\]
This proves (i). To prove (ii) and (iii), we begin with the following setup:

Let $\{U_\al\}_{\al\in A}$ be a collection of open sets in $X$. Observe that the topological gluing condition (for continuous maps $f_\al:U_\al\ra T_1$ on the cover $\{U_\al\}_{\al\in A}$ of $U:=\bigcup_\al U_\al$) resembles the expression of union-closedness of some system in Top (that contains $\I_U\subset\I_X$ as a subsystem), as in diagram (a) below (which expresses the universal property of union as a monic colimit): (\blue{footnote}\footnote{Since the targets of sections of a sheaf are themselves spaces, the sheaf gluing condition can also be obtained from diagram (a) by setting $f_\al:=s_\al\in S(U_\al)$ and $f:=s\in S(\bigcup U_\al)$. I.e., the sheaf gluing condition can be seen as a special case of topological gluing.})
\bea
(a)~~~
\adjustbox{scale=0.9}{%
\bt
 & U_\al\cap U_\beta\ar[dl,hook,"i_\al"']\ar[dr,hook,"i_\beta"] & \\
U_\al\ar[ddr,bend right,"f_\al"']\ar[dr,hook,"q_\al"] & & U_\beta\ar[ddl,bend left,"f_\beta"]\ar[dl,hook,"q_\beta"']\\
 & \bigcup U_\al\ar[d,dashed,"\exists!~f"] & \\
 & T_1 &
\et}\hspace{2cm}
(b)~~~
\adjustbox{scale=0.9}{%
\bt
 & S(U_\al\cap U_\beta)\ar[from=dl,"S(i_\al)"]\ar[from=dr,"S(i_\beta)"'] & \\
S(U_\al)\ar[from=ddr,bend left,"S'(f_\al)"description]\ar[from=dr,"S(q_\al)"'] & & S(U_\beta)\ar[from=ddl,bend right,"S'(f_\beta)"description]\ar[from=dl,"S(q_\beta)"]\\
 & S\left(\bigcup U_\al\right)\ar[from=d,dashed,"\exists!~ S'(f)"'] & \\
 & S'(T_1) & \\
 & T\ar[u,dashed,"\ol{m}"]\ar[uu,dashed,crossing over,bend left=60,"m"]\ar[uuul,bend left=40,"m_\al"] \ar[uuur,bend right=40,"m_\beta"'] &
\et}\nn
\eea
{\flushleft (ii)} Assume $S$ is extended by a presheaf $S':\I\supset\I_X\ra\C$ (i.e., $S'|_{\I_X}=S$) that maps colimits to limits. Applying the cofunctor $S'$ on diagram (a) above we obtain diagram (b), which by hypotheses defines a limit in $\C$. To show the gluing condition is \ul{potentially} realizable, consider compatible sections $s_\al\in S(U_\al)$ for all $\al$, i.e., sections that agree at intersections in the sense $s_\al|_{U_\al\cap U_\beta}=s_\beta|_{U_\al\cap U_\beta}$ for all $\al,\beta$. Fix $t_0\in T$ (see diagram (b)). \ul{Assume} that for each $\al$, we can choose a morphism $m_\al:T\ra S(U_\al)$ such that $m_\al(t_0)=s_\al$, and the collection $\{m_\al\}_{\al\in A}$ satisfies the compatibility condition in diagram (b).
(\blue{footnote}\footnote{By the hypotheses, it might be possible to express the morphisms $m_\al$ in the form $m_\al=S'(f_\al)$ for continuous maps $f_\al$ as follows: Consider compatible continuous maps $f_\al:U_\al\ra T_1$ for all $\al$ (i.e., $f_\al|_{U_\al\cap U_\beta}=f_\beta|_{U_\al\cap U_\beta}$). Topological gluing gives a unique continuous map $f:\bigcup_\al U_\al\ra T_1$ such that $f|_{U_\al}=f_\al$ for all $\al$. Using the axiom of choice, define a map $\ol{m}: T\ra S'(T_1),~t\mapsto\ol{m}(t)\in S'(f)^{-1}(m(t))$. \magenta{In general, neither the uniqueness nor the $\C$-morphism property of $\ol{m}$ is guaranteed.} It is clear that $m=S'(f)\circ\ol{m}$ and $S'(f_\al)\circ\ol{m}=m_\al$ for all $\al$. Now \ul{assume} $\ol{m}$ is a morphism. Because of the factorability of the $m_\al$'s above, wlog we can set $T:=S'(T_1)$, $m_\al:=S'(f_\al)$, and $m:=S'(f)$ or $\ol{m}=id$. This means (by the hypotheses) the limit {\tiny $\varprojlim_\al S(U_\al)\cong S(\bigcup_\al U_\al)$} in $\C$ is completely determined by the associated colimit {\tiny $\varinjlim_\al U_\al\cong\bigcup_\al U_\al$} in Top.}).
Since $S$ maps colimits to limits, the limit {\small $\varprojlim S(U_\al)\cong S(\bigcup_\al U_\al)$} exists. So, there is a unique morphism {\small $m:T\ra S(\bigcup_\al U_\al)$} such that $S(q_\al)m=m_\al$ for each $\al$. We get a unique section {\small $s:=m(t_0)\in S(\bigcup_\al U_\al)$} with
\begin{align}
s|_{U_\al}=m(t_0)|_{U_\al}=S(q_\al)m(t_0)=m_\al(t_0)=s_\al,~~~~\txt{for all}~~~~\al.\nn
\end{align}

{\flushleft (iii)} Assume $S$ is a sheaf. Given morphisms $m_\al:T\ra S(U_\al),~t\mapsto m_\al(t)$ satisfying compatibility as in diagram (b) above, define a map $f:T\ra S(\bigcup_\al U_\al)$ as follows: For $t\in T$, the collection $\{m_\al(t)\}_\al$ satisfies $m_\al(t)|_{U_\al\cap U_\beta}=m_\beta(t)|_{U_\al\cap U_\beta}$ for all $\al,\beta$, and so gives a unique section $m(t)\in S(\bigcup_\al U_\al)$ such that $S(q_\al)m(t)=m(t)|_{U_\al}=m_\al(t)$ for all $\al$. Therefore, $S(q_\al)m=m_\al$ for all $\al$, and so {\small $S(\bigcup_\al U_\al)\cong\varprojlim_\al S(U_\al)$}. Hence (assuming the map $f$ above is automatically a morphism in $\C$), $S$ is limit-closed.
\end{proof}

\begin{question}
Does a sheaf ~$S:\I_X\ra\C\subset Sets$~ necessarily map colimits to limits?
\end{question}

\begin{dfn}[\blue{
\index{Restriction! (pre)sheaf}{Restriction (pre)sheaf},
\index{Restriction! presheaf support (RPS) category}{Restriction presheaf support (RPS) category}}]
Let $X$ be a space, $\C$ a category, and $S:\I_X\ra\C$ a (pre)sheaf. For any open subset $U\subset X$, the \ul{restriction (pre)sheaf} over $U$ is the (pre)sheaf
\[
S|_U:=S|_{\I_U}:\I_U\subset\I_X\ra\C,~~\bt V\ar[r,hook,"{i_{VV'}}"]& V'\et~~\mapsto~~\bt S(V')\ar[r,"{S(i_{VV'})}"] & S(V)\et,\nn
\]
where $\I_U\subset\I_X$ is the full subcategory formed by all open subsets $V\subset U$ of the subspace $U\subset X$. The \ul{Restriction presheaf support (RPS) category} (associated with $X$ and $\C$) is the category of categories RPS$(X,\C)$ with objects $\{\I_U:\txt{open}~U\subset X\}\cup\{\C\}$ and morphisms cofunctors between the objects.
\end{dfn}

It is worth noting (so as to avoid possible confusion) that $\I_U\neq \I_{X,U}$, where $\I_{X,U}\subset\I_X$ was defined to be the full subcategory formed by all open sets containing $U$. That is, as defined earlier, $\I_U$ is the TopIndex-category of the subspace $U\subset X$, while $\I_{X,U}$ is the local TopIndex-category of $X$ at $U$.

\begin{lmm}[\blue{Gluing lemma: \cite[Proposition 5.76, pp 291-292]{rotman2009}}]\label{SheafGluLmm}
Let $X$ be a space, $\C\subset Sets$ a category, $\{U_i\}_{i\in I}$ an open cover of $X$, $\{S_i:\I_{U_i}\ra\C\}_{i\in I}$ sheafs, and $\theta_{ij}:S_i|_{U_i\cap U_j}\ra S_j|_{U_i\cap U_j}$ sheaf-isomorphisms satisfying the following: For all $i,j,k\in I$, with $U_{ij}:=U_i\cap U_j$ and $U_{ijk}:=U_i\cap U_j\cap U_k$,
\bit
\item[(i)] $\theta_{ii}=1_{S_i}$ ~(not strictly essential, except perhaps for convenience in specific cases).
\item[(ii)] the restrictions to $U_{ijk}$ satisfy ~$\bt \theta_{ik}=\theta_{jk}\theta_{ij}:S_i|_{U_{ijk}}\ar[r,"{\theta_{ij}|}"]& S_j|_{U_{ijk}}\ar[r,"{\theta_{jk}|}"]& S_k|_{U_{ijk}}\et$.
\eit
\[\adjustbox{scale=0.9}{\bt
S_i|_{U_{ij}}\ar[d,"|_{U_{ijk}}"]\ar[rr,"\theta_{ij}"] && S_j|_{U_{ij}}\ar[d,"|_{U_{ijk}}"] \\
S_i|_{U_{ijk}}\ar[rr,"\theta_{ij}|"] && S_j|_{U_{ijk}} \\
 &S_k|_{U_{ijk}}\ar[from=ul,"\theta_{ik}|"']\ar[from=ur,"\theta_{jk}|"] &
\et}\]

Then there exists a unique sheaf $S:\I_X\ra\C$ and sheaf-isomorphisms $\eta_i:S_i\ra S|_{U_i}$ such that
\[
\bt\eta_j^{-1}\eta_i=\theta_{ij}:S_i|_{U_{ij}}\ar[r,"{\eta_i|}"]&S|_{U_{ij}}\ar[r,"{\eta_j^{-1}|}"]&S_j|_{U_{ij}}\et~~\txt{over}~~U_{ij}~~\txt{for all}~~i,j\in I.\nn
\]
\[\adjustbox{scale=0.9}{\bt
S_i|_{U_{ij}}\ar[rr,"\theta_{ij}"] && S_j|_{U_{ij}} \\
 &S|_{U_{ij}}\ar[from=ul,"\eta_{i}|"']\ar[from=ur,"\eta_{j}|"] &
\et}\]

\end{lmm}
\begin{proof}
{\flushleft \ul{Existence}}: Define {\footnotesize $S:\I_X\ra\C, \bt U\ar[r,hook,"{i_{UV}}"]&V\et \mapsto \bt S(V)\ar[r,"{S(i_{UV})}"]&S(U)\et$} by ~{\footnotesize $S(U):=\varinjlim_{i\in I}S_i(U\cap U_i)$} and {\small $S(i_{UV}):=\varinjlim_{i\in I}S_i(i_{U\cap U_i,V\cap U_i})$},~ i.e., ~{\small $S|_O:=\varinjlim_{i\in I}S_i|_{U_i\cap O}$} for any open set $O\subset X$.
\[\adjustbox{scale=0.8}{\bt
 S_i|_{U_{ij}}(U\cap U_{ij})\ar[ddr,bend left=10,"{q_i|}"]\ar[rr,"\theta_{ij}"] && S_j|_{U_{ij}}(U\cap U_{ij})\ar[ddl,bend right=10,"{q_j|}"']\\
 S_i(U\cap U_i)\ar[u,"|_{U_{ij}}"']\ar[dddr,bend right,"f_i"']\ar[dr,"q_i"] &&  S_j(U\cap U_j)\ar[u,"|_{U_{ij}}"']\ar[dl,"q_j"']\ar[dddl,bend left,"f_j"] \\
 & S(U):=\varinjlim_{i\in I}S_i(U\cap U_i)\ar[dd,dashed,"f"] & \\
 &  & \\
 & Z &
\et}~~~~f_i(s_i|_{U_{ij}})=f_j\theta_{ij}(s_i|_{U_{ij}}).
\]
\ul{Showing $S$ is a presheaf}: Observe that given an inclusion $\bt U\ar[r,hook,"{i_{UV}}"]&V\et$ of open subsets of $X$, we have
\[
\textstyle \bt\Big(U\cap U_i\ar[rr,hook,"{i_{U\cap U_i,V\cap U_i}}"]&&V\cap U_i\Big)\et~~\sr{S_i}{\longmapsto}~~\bt\Big(S_i(V\cap U_i)\ar[rr,"{S_i(i_{U\cap U_i,V\cap U_i})}"]&& S_i(U\cap U_i)\Big)\et,~~~~\txt{for each}~~i\in I,
\]
and so taking the colimit $\varinjlim_i$ (as a functor $\txt{RPS}(X,\C)^I\ra\txt{RPS}(X,\C)$, where RPS$(X,\C)$ is the associated restriction presheaf support category) we get
\[
\textstyle  \bt\Big(U\ar[r,hook,"{i_{UV}}"]&V\Big)\et~~\sr{S}{\longmapsto}~~\bt\Big(S(V)\ar[r,"{S(i_{UV})}"]&S(U)\Big)\et,~~~~\txt{where}~~~~S(i_{UV}):=\varinjlim_{i\in I}S_i(i_{U\cap U_i,V\cap U_i}).\nn
\]
\ul{Finding the presheaf-isomorphisms}: Let $\eta_i:=\varinjlim_{k\in I}\theta_{ik}:\varinjlim_{k\in I}\left(S_i|_{U_{ik}}\ra S_k|_{U_{ik}}\right)=\left(S_i\ra S|_{U_i}\right)$, where for any fixed $i\in I$, if we set
\[
U':=U_i,~~~~U'_k:=U'\cap U_k=U_i\cap U_k=U_{ik}~~~~\txt{and}~~~~U'_{ik}:=U'_i\cap U'_k=U_{ik}~~~~\txt{for all}~~~~k\in I,
\]
then $\varinjlim_{k\in I}S_i|_{U_{ik}}=S_i$ as the \ul{unique (up to iso)} colimit $S'\cong S_i$ of the presheaf-system
\[
\textstyle \{S'_k|_{U'_{ik}}=S_i|_{U_{ik}}:\I_{U'_{ik}}=\I_{U_{ik}}\ra\C\}_{k\in I}~~\txt{on}~~U'_i=\bigcup_kU'_{ik}\subset U'
\]
based on the constant-cover presheaf-system $\{S'_k=S_i:\I_{U'_k}=\I_{U_i}\ra\C\}_{k\in I}$ on $U'=\bigcup_kU'_k$ ~(see \blue{footnote}\footnote{This is in the same way that the presheaf-system $\{S_k|_{U_{ik}}:\I_{U_{ik}}\ra\C\}_{k\in I}$ on $U_i=\bigcup_kU_{ik}\subset U$ is based on the (non constant-cover) presheaf-system $\{S_k:\I_{U_k}\ra\C\}_{k\in I}$ on $U=\bigcup_kU_k$.}).

Then taking $\varinjlim_k$ in (ii) {\footnotesize $\bt \theta_{ik}=\theta_{jk}\theta_{ij}:S_i|_{U_{ijk}}\ar[r,"{\theta_{ij}|}"]&S_j|_{U_{ijk}}\ar[r,"{\theta_{jk}|}"]&S_k|_{U_{ijk}},\et$} we get (by the functoriality of $\varinjlim$)
\[
\textstyle \bt \eta_i=\eta_j\theta_{ij}:S_i|_{U_{ij}}\ar[r,"{\theta_{ij}|}"]&S_j|_{U_{ij}}\ar[r,"{\eta_j|}"]&S|_{U_{ij}}.\et
\]
\ul{Showing $S:\I_X\ra\C$ is a sheaf}: Let $U\subset X$ be open, $U\subset \bigcup_{\al\in A} U_\al$ for open sets $U_\al\subset U$, and $s_\al\in S(U_\al)$ such that $s_\al|_{U_\al\cap U_\beta}=s_\beta|_{U_\al\cap U_\beta}$ for all $\al,\beta$. Fix $i\in I$. Let $s_{\al,i}:=\eta_i^{-1}(s_\al|_{U_\al\cap U_i})\in S_i(U_\al\cap U_i)$ for all $\al\in A$, which (because $\eta_i^{-1}$ is a natural transformation) are compatible, i.e.,
\[
s_{\al,i}|_{U_\al\cap U_\beta\cap U_i}=s_{\beta,i}|_{U_\al\cap U_\beta\cap U_i}~~~~\txt{for all}~~~~\al,\beta.
\]
Then $\{s_{\al,i}\}_{\al\in A}$ glue together into a unique section $s_i\in S_i(U\cap U_i)$ such that $s_i|_{U_\al\cap U_i}=s_{\al,i}|_{U_\al\cap U_i}$ for all $\al\in A$. Consider the colimit morphisms $q_i:S_i(U\cap U_i)\ra S(U)=\varinjlim_iS_i(U\cap U_i)$. Then, depending on the structure of the category $\C$, the collections $\wt{s}:=\{q_i(s_i)\}_{i\in I}\subset S(U)$ and $\wt{s}_\al:=\{q_i(s_{i,\al})|_{U_\al}\}_{i\in I}\subset S(U_\al)$, which by construction (including naturality of $q_i$ and $\eta_i^{-1}$) satisfy
\[
\wt{s}|_{U_\al}:=\{q_i(s_i)|_{U_\al}\}_{i\in I}=\{q_i(s_{i,\al})|_{U_\al}\}_{i\in I}=\{q_i\eta_i^{-1}(s_\al|_{U_\al\cap U_i})|_{U_\al}\}_{i\in I}=\wt{s}_\al|_{U_\al}~~~~\txt{for all}~~\al\in A,
\]
will represent/describe a unique section $s\in S(U)$ satisfying
\[
s|_{U_\al}=s_\al|_{U_\al}~~~~\txt{for all}~~\al\in A.
\]

{\flushleft \ul{Uniqueness (up to isomorphism)}}: Consider any sheafs $S,S':\I_X\ra\C$ with isomorphisms $\eta_i:S_i\ra S|_{U_i}$ and $\eta'_i:S_i\ra S'|_{U_i}$, for all $i\in I$, satisfying (ii). Then we get isomorphisms
\[
\bt \eta_{ij}:=\eta'_j\theta_{ij}\eta_i^{-1}:S|_{U_{ij}}\ar[r,"{\eta_i^{-1}|}"]&S_i|_{U_{ij}}\ar[r,"{\theta_{ij}}"]&S_j|_{U_{ij}}\ar[r,"{\eta'_j|}"]&S'|_{U_{ij}}\et~~\txt{for all}~~i,j\in I,
\]
which give isomorphisms $\eta_{ii}:S|_{U_i}\ra S'|_{U_i}$ for all $i\in I$.~~~~
\bt[column sep=small]
 &S|_{U_i}\ar[dl,bend right=10,"\eta_i"']\ar[d,"q_{\eta_i}"]\ar[rr,"\eta_{ii}"] && S'|_{U_i}\ar[d,"q_{\eta'_i}"']\ar[dr,bend left=10,"\eta'_i"]& \\
S\ar[r,draw=none,"\cong"description]& \varinjlim_iS|_{U_i}\ar[rr,"\varinjlim_i\eta_{ii}"] && \varinjlim_iS'|_{U_i} \ar[r,draw=none,"\cong"description]& S'
\et
\end{proof}

It is worth noting that (unlike with the presheaf) the gluing condition defining the sheaf allows us to easily understand certain things through a description of the sheaf in terms of its stalks only. For example, we have already seen in Lemma \ref{MorOfShfUniq} that for morphisms of presheafs to agree, they need only agree on stalks provided the target presheaf is a sheaf. Similarly, by Corollary \ref{MorOfShfUniq2}, a morphism of sheafs is an isomorphism if and only if the morphisms it induces on stalks are isomorphisms.

For further details on (pre)sheafs and their applications, see for example \cite[Section 5.4, pp 273-303]{rotman2009}, \cite[Section 6.3, pp 377-403]{rotman2009}, \cite[Chapter II.1, pp 60-69]{hartshorne1977}, and \cite[Chapter 3, p.36]{grothdk1957}.

\section{Classical Manifolds and the Tangent Space}
In this section, $\mathbb{F}$ will denote a classical field such as $\Real$ or $\Complex\cong\Real^2$ (as spaces).
\begin{dfn}[\textcolor{blue}{
\index{Manifold!}{Manifold},
\index{Manifold! without boundary}{Manifold without boundary},
\index{Manifold! internal point of}{Internal point},
\index{Manifold! boundary point of}{Boundary point},
\index{Lie! group}{Lie group},
\index{Manifold! boundary of}{Boundary of a manifold},
\index{Manifold! $n$-manifold}{$n$-Manifold},
\index{Manifold! Curve ($1$-manifold)}{Curve},
\index{Manifold! Surface ($2$-manifold)}{Surface}}]
Let $X,Y$ be spaces (\blue{footnote}\footnote{The topology (homeo) class of each of the spaces $X,Y$ is often narrowed down, depending on the intended application.}) and $\C$ a subcollection of closed subsets of $X$. $Y$ is an \ul{$(X,\C)$-manifold} (or a manifold over $(X,\C)$) if each point of $Y$ has a neighborhood that is either (i) homeomorphic to an \ul{open} subset of $X$ or (ii) homeomorphic to a \ul{closed} subset of $X$ in $\C$. (\blue{footnote}\footnote{Additional conditions are often applied to the space $Y$, depending on the intended application.}). If $Y$ is an $(X,\emptyset)$-manifold, we say $Y$ is an \ul{$X$-manifold without boundary}.

If $Y$ is an $(X,\C)$-manifold, a point $y\in Y$ is an \ul{$(X,\C)$-internal point} of $Y$ if it has a neighborhood that is homeomorphic to an open subset of $X$, otherwise $y$ is an \ul{$(X,\C)$-boundary point} of $Y$. The \ul{$(X,\C)$-boundary} of $Y$ is the collection $\del_{(X,\C)} Y:=$ $\{$\txt{all $(X,\C)$-boundary points of $Y$}$\}$.

An $(X,\C)$-manifold $G$ is called a \ul{Lie group} if (i) $G$ is a group and (ii) both group multiplication $G\times G\ra G,~(g,g')\mapsto gg'$ and group inverse $G\ra G,~g\mapsto g^{-1}$ are continuous maps.

An $\big(\Real^n,\{\Real^{n-1}\times[0,\infty)\}\big)$-manifold is called an \ul{$n$-manifold} (or an \ul{$n$-dimensional manifold}). A $1$-manifold is called a \ul{curve}, and a $2$-manifold is called a \ul{surface}.
\end{dfn}

Henceforth, unless it is specified otherwise, ``\emph{manifold}'' we will mean ``\emph{$n$-manifold without boundary (for some $n\geq 1$)}'', even though the techniques used in the discussion (along with appropriate modifications) can also be applied to manifolds with boundary.

\begin{dfn}[\textcolor{blue}{
\index{Coordinate charts}{Coordinate chart},
\index{Atlas}{Atlas},
\index{Classical! manifold}{Classical manifold},
\index{Geometric linear approximation}{Geometric linear approximation},
\index{Tangent space}{Tangent space},
\index{Tangent vectors}{Tangent vectors}}]
Let $\mathbb{F}$ be a classical field such as $\Real$ or $\Complex$ (where as a space, $\Complex\cong\Real^2$), $n\geq 1$ and integer, and $\M$ a topological space. An \ul{$\mathbb{F}^n$-coordinate chart} (or just a \ul{chart}) on $\M$ is a pair $(\vphi,U)$ consisting of an open set $U\subset\M$ and an open imbedding $\vphi:U\ra \mathbb{F}^n$. (\blue{footnote}\footnote{A coordinate chart is also called a ``coordinate system'' or a ``coordinate parametrization''.}). An \ul{$\mathbb{F}^n$-atlas} $\A$ on $\M$ (making $\M=(\M,\A)$ an \ul{$\mathbb{F}^n$-manifold} or an \ul{$n$-dimensional manifold over $\mathbb{F}$}) is a collection of charts $\A=\{(\vphi_\al,U_\al)\}_{\al\in A}$ that covers $\M$ in the sense $\{U_\al\}_{\al\in A}$ is an open cover of $\M$.

Given a point $p\in\M$, the \ul{geometric linear approximation} (i.e., description or expression up to homeomorphism using curves/lines) of a ``small'' open neighborhood of $p$ is called the \ul{tangent space} $T_p\M$ to $\M$ at $p$. (\blue{footnote}\footnote{Technically, this definition for the tangent space is not yet completely precise but it is sufficient (or suitable) for this brief opening discussion. For brevity, a technically precise definition will be stated only inside the proof of Lemma \ref{TanSpLmm}.}). $T_p\M$ typically consists of copies of $\mathbb{F}^n$ (as an $\mathbb{F}$-vector space), in the form
{\small\begin{align}
&\textstyle\bigoplus_{\al\in A_p}\mathbb{F}^n\cong T_p\M\sr{(s)}{\subset}\left\{\txt{$\mathbb{F}$-linear functions}~~v:\mathbb{F}^U\ra \mathbb{F}~|~U\subset \M,~U\ni p\right\},\nn
\end{align}}where the motivation for the containment at step (s) will be explained in the proof of Lemma \ref{TanSpLmm}. The elements of $T_p\M$ are called \ul{tangent vectors} to $\M$ at $p$.
\end{dfn}

\begin{dfn}
[\textcolor{blue}{
\index{Transition maps}{Transition maps},
\index{Differentiable! map}{Differentiable map},
\index{Set! of differentiable maps}{Set of differentiable maps},
\index{Derivative}{Derivative},
\index{Differentiable! atlas}{Differentiable atlas},
\index{Differentiable! manifold}{Differentiable manifold},
\index{Compatible charts}{Compatible charts},
\index{Differentiable! structure}{Differentiable structure},
\index{Derivations}{Derivations}}]~\\~
A map between manifolds $f:(\M,\{(\vphi_\al,U_\al)\})\ra(\M',\{(\vphi'_{\al'},U'_{\al'})\})$ is a \ul{differentiable map}, written $f\in D(\M,\M')$, if the associated maps (to be called \ul{$f$-transition maps})
\bea
f_{\al\al'}:=\vphi'_{\al'}f\vphi_\al^{-1}:\vphi_\al(U_\al)\subset \mathbb{F}^n\ra \vphi'_{\al'}(U'_{\al'})\subset \mathbb{F}^n\nn
\eea
are \ul{differentiable}, i.e., for all $a\in \dom f_{\al\al'}$ and all $b\in \mathbb{F}^n$, the limit ~$\del_bf_{\al\al'}(a):=\lim\limits_{t\in \mathbb{F},~t\ra 0} {f_{\al\al'}(a+tb)-f_{\al\al'}(a)\over t}$~ (called the \ul{derivative} of $f_{\al\al'}$ at $a$ along $b$)~ exists or makes sense. (\blue{footnote}\footnote{Equivalently, $f_{\al\al'}$ is differentiable if for every $a\in\dom f_{\al\al'}$, there exist (i) a unique $\mathbb{F}$-linear map $\del f_{\al\al'}(a)\in Hom_{\mathbb{F}}(\mathbb{F}^n,\mathbb{F}^n)$ and (ii) a continuous map $F(a):N(0)\subset \mathbb{F}^n\ra Hom_{\mathbb{F}}(\mathbb{F}^n,\mathbb{F}^n)$ (for a neighborhood $N_a(0)\subset \mathbb{F}^n$ of $0\in \mathbb{F}^n$) such that
\[
f_{\al\al'}(a+h)=f_{\al\al'}(a)+\del f_{\al\al'}(a)h+F(a)(h)h,~~\txt{for all}~~h\in N_a(0),~~\txt{and}~~\lim_{h\ra 0}F(a)(h)=0,
\]
where by convention, the limit ~$\lim\limits_{h\ra 0}F(a)(h)$~ makes sense $\iff$ the partial limit ~$\lim\limits_{t\ra 0,~t\in \mathbb{F}}F(a)(tb)$~ makes sense for each $b\in \mathbb{F}^n$, ~$Hom_{\mathbb{F}}(\mathbb{F}^n,\mathbb{F}^n)\subset(\mathbb{F}^n)^{\mathbb{F}^n}=\prod_{h\in \mathbb{F}^n}\mathbb{F}^n$~ as a product-top subspace, and we often write ~$\del_hf_{\al\al'}(a):=\del f_{\al\al'}(a)h$. For the precise meaning of the notion of topological limits, see the discussions in the later chapters of these notes.
}).

Let $\M=(\M,\A)$ be a manifold. The atlas $\A=\{(\vphi_\al,U_\al)\}_{\al\in A}$ is called a \ul{differential atlas} (making $\M$ a \ul{differentiable manifold}) if the identity map $id_\M:(\M,\{(\vphi_\al,U_\al)\})\ra(\M,\{(\vphi_{\beta},U_{\beta})\})$ is differentiable, i.e., the following associated \ul{transition maps} are differentiable: (\blue{footnote}\footnote{With further assumptions and a suitable standard atlas on $\mathbb{F}^n$, if each of the coordinate charts $\vphi_\al:U_\al\ra \mathbb{F}^n$ is differentiable as a map between manifolds, then we may be able to conclude that the atlas $\{(\vphi_\al,U_\al)\}_{\al\in A}$ is differentiable.})
\bea
f_{\al\beta}:=\vphi_\beta\vphi_\al^{-1}:\vphi_\al(U_\al)\subset \mathbb{F}^n\ra \vphi_\beta(U_\beta)\subset \mathbb{F}^n.\nn
\eea

Given charts $\vphi:U\subset\M\ra \mathbb{F}^n$ and $\vphi':U'\subset\M\ra \mathbb{F}^n$, call them \ul{compatible charts} (written $\vphi\sim\vphi'$) if the associated transition maps $\vphi\circ\vphi'^{-1}:\vphi'(U')\ra\vphi(U)$ and $\vphi'\circ\vphi^{-1}:\vphi(U)\ra\vphi'(U')$ are both differentiable. This is clearly an equivalence relation. Each of the resulting equivalence classes is therefore a maximal differentiable atlas on $\M$, which is called a \ul{differentiable structure} on $\M$ (i.e., a \ul{differentiable structure} on $\M$ is defined to be a maximal differentiable atlas on $\M$).

It is clear that a differentiable manifold $\M=(\M,\A)$ has a unique differentiable structure, namely, the maximal differentiable atlas $\langle\A\rangle$ containing $\A$.

Consider a differentiable manifold $\M=(\M,\{(\vphi_\al,U_\al)\})$. Given $p\in \M$, let $D_p(\M):=\bigsqcup_{p\in U\subset\M}D(U,\mathbb{F})$. Then we might (as in the proof of Lemma \ref{TanSpLmm} below) be able to use a bijection
\begin{align}
\textstyle \M\longleftrightarrow &\left\{\txt{evaluations}~~e_p:D_p(\M)\ra \mathbb{F},~f\mapsto f(p)\right\}_{p\in\M}\nn\\
&\textstyle\subset Hom_{\mathbb{F}}(D_p(\M),\mathbb{F})\subset \mathbb{F}^{D_p(\M)}=\prod_{f\in D_p(\M)}\mathbb{F}~~~~\txt{(as subsets)}\nn
\end{align}
to show that for each $p\in\M$, we have $T_p\M\cong \mathbb{F}^n$ and the $\mathbb{F}$-vectors $v\in T_p\M$ satisfy
\bea
v(fg)=v(f)g(p)+f(p)v(g),~~~\txt{for all differentiable functions}~~~f,g\in D(U,\mathbb{F}),~U\ni p,\nn
\eea
and are then called \ul{derivations at $p$}.
\end{dfn}

\begin{lmm}[\textcolor{blue}{The tangent space of a differentiable manifold}]\label{TanSpLmm}
Let $\M$ be a differentiable $\mathbb{F}^n$-manifold, and ~$S:\I_\M\ra \mathbb{F}\txt{-Algebras}\subset Rings$ ~the presheaf of functions given by
\[
S(U):=D(U,\mathbb{F}):=\left\{\txt{differentiable functions}~f:U\ra \mathbb{F}\right\},~~~~\txt{for any open set $U\subset \M$}
\] (with restriction given by the usual restriction of maps). Suppose, with {\footnotesize $S_p(\M):=\bigsqcup_{p\in U\subset\M}S(U)$}, the map
\[
\textstyle\M\ra\big\{\txt{evaluations}~~e_p:S_p(\M)\ra \mathbb{F},~f\mapsto f(p)\big\}_{p\in\M},~p\mapsto e_p
\] is a bijection. Then for each $p\in\M$, we have an $\mathbb{F}$-linear isomorphism of the following form:
\bea
\bt[column sep=small] T_p\M\ar[r,hook,two heads]&\left(S_p/S_p^{st}\right)^\ast:=Hom_{\mathbb{F}}(S_p/S_p^{st},\mathbb{F})\cong \mathbb{F}^n\et,\nn
\eea
where $S_p$ is the stalk of $S$ at $p$, (\blue{footnote}\footnote{Recall that given $p\in\M$ and sections $f,g\in S_p(\M)\:=\bigsqcup_{p\in U\subset\M}S(U)$, we write $f\sim g$ if $f|_V=g|_V$ for some open set $p\in V\subset\dom f\cap\dom g$, and define $S_p:=\varinjlim S|_{\I_{\M,p}}\cong{S_p(\M)\over\sim}=\left\{[f]:f\in S_p(\M)\right\}$.}), and $S_p^{st}\subset S_p$ is the ideal (consisting of \index{Germs! stationary}{\ul{stationary germs}})
{\footnotesize\bea
\textstyle S_p^{st}:=\big\{[f]\in S_p:(f\circ\gamma)'(0)=0~\txt{for every curve}~\gamma:(-1,1)\ra\M~\txt{with}~\gamma(0)=p\big\},\nn
\eea}where ~$(f\circ\gamma)'(0):={df(\gamma(t))\over dt}\big|_{t=0}$,~ with ~${df(\gamma(t))\over dt}:=\lim_{\vep\ra 0}{f(\gamma(t+\vep))-f(\gamma(t))\over\vep}$.
\end{lmm}
\begin{proof}
View $S_p(\M):=\bigsqcup_{p\in U\subset\M}S(U)$ as a set with operations of addition $+:S_p(\M)\times S_p(\M)\ra S_p(\M),~(f,g)\mapsto f+g\in S(\dom f\cap\dom g)$ and multiplication $S_p(\M)\times S_p(\M)\ra S_p(\M),~(f,g)\mapsto fg\in S(\dom f\cap\dom g)$ given respectively by
\[
(f+g)(q):=f(q)+g(q),~~~~(fg)(q):=f(q)g(q),~~~~\txt{for all}~~~~q\in \dom f\cap\dom g,
\]
as well as scalar multiplication $\mathbb{F}\times S_p(\M)\ra S_p(\M),~(\ld,g)\mapsto \ld g\in S(\dom f)$ given by
\[
(\ld f)(q):=\ld f(q),~~~~\txt{for all}~~~~q\in\dom f.
\]
Then by the hypotheses, for any $p,q\in\M$, we have
\bea
&&\label{TanSpEq1}p=q\iff f(p)=f(q)~~\txt{for all}~~f\in S_p(\M)=S_q(\M),\nn\\
&&~~~~\iff f(p)=f(q)~~\txt{for all}~~[f]\in S_p=S_q,
\eea
and so we can view each point $p\in\M$ as the evaluation function $p:S_p(\M)\ra \mathbb{F}$, $f\mapsto f(p)$. Likewise if $p\in\M$, then for any curve $\gamma:(-1,1)\ra\M$ with $\gamma(0)=p$, we can view each image point $\gamma(t)\in\M$, $t\in(-1,1)$, as the map
\bea
\label{TanSpEq2}\gamma(t):S_{\gamma(t)}(\M)\ra \mathbb{F},~~f\mapsto f(\gamma(t))=(f\circ\gamma)(t),
\eea
which we can ``differentiate'' to define the \index{Curve differential}{\ul{curve differential}} {\small $\dot{\gamma}:(-1,1)\ra Hom_{\mathbb{F}}(S_{\gamma(t)}(\M),\mathbb{F})$} given by
\bea
\label{TanSpEq3}\textstyle t\in(-1,1)~~\longmapsto~~\dot{\gamma}_t:S_{\gamma(t)}(\M)\ra \mathbb{F},~~f\mapsto(f\circ\gamma)'(t):={df(\gamma(t))\over dt}.
\eea
Using the definition of the derivative, it is clear that for each $t\in(-1,1)$ the \ul{linear approximation} to the curve $\gamma$ near the point $\gamma(t)\in\M$ is fully determined by the associated $\mathbb{F}$-linear map $\dot{\gamma}_t$, in the sense for some $\delta_t>0$ and all $\vep\in(0,\delta_t)$, the following relation holds (pointwise at each {\small $f\in S_{\gamma(t)}(\M)\cap S_{\gamma(t+\vep)}(\M)$}):
\[
\gamma(t+\vep)=\gamma(t)+\vep~\dot{\gamma}_t+\vep^2\Gamma(t,\vep),~~\txt{for some}~~\Gamma(t,\vep)\in S_{\gamma(t)}(\M)\cap S_{\gamma(t+\vep)}(\M)~~\txt{with}~~\lim_{\vep\ra 0}\Gamma(t,\vep)=0.\nn
\]
In particular, the linear approximation to $\gamma$ near $p=\gamma(0)\in\M$ is fully determined by $\dot{\gamma}_0$ (which is thus a \ul{representative of the tangent vector} to $\M$ at $p$). Consequently, we can now describe the tangent space at $p\in\M$ as follows: With $\Gamma_p:=\{\txt{curves}~\gamma:(-1,1)\ra\M,~\gamma(0)=p\}$ and $S_p^{st}:=\big\{[f]\in S_p:(f\circ\gamma)'(0)=0~\txt{for all}~\gamma\in\Gamma_p\big\}$, the tangent space at $p$ is given by the following \ul{definition} and $\mathbb{F}$-linear correspondence:
{\small\begin{align}
&\textstyle T_p\M:={\Gamma_p\over\sim_p}=\big\{[\gamma]_p:\gamma\in\Gamma_p\big\},~~~~\txt{where}~~~~\gamma_1\sim_p\gamma_2~~\txt{iff}~~(f\circ\gamma_1)'(0)=(f\circ\gamma_2)'(0)~~\txt{for all}~~f\in S_p(\M),\nn\\
&\textstyle~~~~\cong \left\{v_{[\gamma]_p}=\dot{\gamma}_0:S_p(\M)\ra \mathbb{F},~f\mapsto(f\circ\gamma)'(0),~\gamma\in\Gamma_p\right\},~~~~\txt{where}~~~~(f\circ\gamma)'(0):={df(\gamma(t))\over dt}\big|_{t=0},\nn\\
&\textstyle~~~~\cong \left\{v_{[\gamma]_p}=\dot{\gamma}_0:S_p\ra \mathbb{F},~[f]\mapsto(f\circ\gamma)'(0),~\gamma\in\Gamma_p\right\}\nn\\
&\textstyle~~~~\cong \left\{v_{[\gamma]_p}=\dot{\gamma}_0:S_p/S_p^{st}\ra \mathbb{F},~[f]+S_p^{st}\mapsto(f\circ\gamma)'(0),~\gamma\in\Gamma_p\right\}\nn\\
&~~~~\cong \mathbb{F}^n,~~~~\txt{(proved below)}\nn
\end{align}}where the proof of the last isomorphism is given by the rest of the discussion below. We can easily check directly using the above definition of a tangent vector based on (\ref{TanSpEq3}) that $v\in T_p\M$ is a \index{Derivation}{\ul{derivation}} (i.e., $v$ is linear and satisfies the \index{Leibnitz rule}{\ul{Leibnitz rule}}) with the following basic properties: For any $\al,\beta\in \mathbb{F},~~f,g\in S_p(\M)$,
\bea
\label{tangent-der1}&& v(\al f+\beta g)=\al v(f)+\beta v(g),~~~~~~~~\txt{(linearity)}\\
\label{tangent-der2}&& v(fg)=v(f)~g(p)+f(p)~v(g).~~~~~~~~\txt{(Leibnitz rule)}
\eea
Having realized the tangent space $T_p\M$ as a space of derivations on differentiable functions $f\in S_p(\M)$ that represent \index{Germs! non-stationary}{\ul{non-stationary germs}} $[f]+S_p^{st}\in S_p/S_p^{st}$ at $p$, let us see what a tangent vector $v\in T_p\M$ looks like in a coordinate chart. Using a coordinate chart $\vphi:U_p\ra \mathbb{F}^n,~p\mapsto x:=\vphi(p)=\big(\vphi_1(p),\cdots,\vphi_n(p)\big)=:(x_1,...,x_n)$ in an open neighborhood $U_p$ of $p$, we have
{\small
\begin{align}
&\textstyle v(f)={df(\gamma(t))\over dt}\big|_{t=0}={df\big(\vphi^{-1}\circ\vphi(\gamma(t))\big)\over dt}\big|_{t=0}={d\big(f\circ\vphi^{-1}\big)\big(\vphi(\gamma(t))\big)\over dt}\big|_{t=0}\\
\label{function-pullback}&\textstyle~~~~={df_\vphi(\gamma_\vphi(t))\over dt}\big|_{t=0},~~~~f_\vphi:=f\circ\vphi^{-1},~~\gamma_\vphi:=\vphi\circ \gamma:(-1,1)\ra \mathbb{F}^n,~\gamma_\vphi(0)=\vphi(p)=x,\\
&\textstyle~~~~\sr{(s)}{=}\sum_{i=1}^n{d\gamma_{\vphi_i}(t)\over dt}\big|_{t=0}{\del f_\vphi(\gamma_\vphi(t))\over \del \gamma_{\vphi_i}(t)}\big|_{t=0}=\sum_{i=1}^n{d\gamma_{\vphi_i}(t)\over dt}\big|_{t=0}{\del f_\vphi(\gamma_\vphi(0))\over \del \gamma_{\vphi_i}(0)}=\sum_{i=1}^n{d\gamma_{\vphi_i}(t)\over dt}\big|_{t=0}{\del f_\vphi(\vphi(p))\over \del \vphi_i(p)}\nn\\
\label{coordinate-vector}&\textstyle~~~~=\sum_{i=1}^n{d\gamma_{\vphi_i}(t)\over dt}\big|_{t=0}{\del f_\vphi(x)\over \del x_i}=\sum_{i=1}^nv_i\del_{x_i}(f),~~~~~~v_i:={d\gamma_{\vphi_i}(t)\over dt}\big|_{t=0},~~\del_{x_i}(f):={\del f(p)\over \del x_i}={\del f_\vphi(x)\over \del x_i},
\end{align}
}where the quantity ${\del f_{\vphi}(x)\over\del x_i}:=\lim\limits_{\vep\ra 0}{f_\vphi(x+\vep e_i)-f_\vphi(x)\over\vep}$ (i.e., the derivative of $f_\vphi$ at $x$ along the standard unit vector $e_i:=(\delta_{ij})_{j=1}^n$) is called the \index{Partial! derivative}{\ul{partial derivative}} of $f_\vphi$ at $x$ in the variable $x_i$, and the easily verified property of the derivative used at step (s), i.e., ${d f_\vphi(\gamma_\vphi(t))\over dt}=\sum_{i=1}^n{d\gamma_{\vphi_i}(t)\over dt}{\del f_{\vphi}(\gamma_\vphi(t))\over\del\gamma_{\vphi_i}(t)}$, is called the \index{Chain! rule}{\ul{chain rule}}. That is,
\[
\textstyle v= \sum_{i=1}^nv_i\del_{x_i},~~~~v_i=v_i^\vphi~~\txt{(i.e., the $v_i$'s depend on the chart)},
\]
where~ the coefficients $v_i\in \mathbb{F}$, $\del_{x_i}\in T_p\M$ as we can verify using the definition of $\del_{x_i}(f)\in \mathbb{F}$ in (\ref{coordinate-vector}), and the set $\{\del_{x_1},\cdots,\del_{x_n}\}\subset T_p\M$ is linearly independent over $\mathbb{F}$. Hence, in the coordinate chart $\vphi$, we have a natural basis $\del_x=\{\del_{x_i}\}_{i=1}^n$~ for $T_p\M$ as an $\mathbb{F}$-space, i.e., $T_p\M\cong \mathbb{F}^n$.
\end{proof}

Lemma \ref{TanSpLmm} shows that if $\M=(\M,\A)$ is an $\mathbb{F}^n$-manifold, then under certain \ul{regularity conditions} (including differentiability) on the atlas $\A=\{\vphi_\al,U_\al\}_{\al\in A}$, we can \ul{equivalently} view $\M$ as a pair $(\M,S)$ consisting of the topological space $\M$ and the following presheaf of rings (in which {\small $S(U):=\left\{\txt{differentiable functions}~f:U\ra \mathbb{F}\right\}$})
\bea
S:\I_\M\ra\txt{Rings},~~U\sr{i_{UV}}{\hookrightarrow}V~~\mapsto~~S(V)\sr{S(i_{UV})}{\ral}S(U),\nn
\eea
with the property that (i) for all $p,q\in\M$, we have $p=q$ iff $f(p)=f(q)$ for all $[f]\in S_p=S_q$, and/or (ii) for each $p\in\M$, as $\mathbb{F}$-spaces, ~$Hom_{\mathbb{F}}(S_p/J_p,\mathbb{F})\cong \mathbb{F}^n$~ for some \ul{special} ideal $J_p\lhd S_p$. Moreover, using the fact that each open set $U\subset U_\al\subset\M$ (for some $\al$) can be viewed as an open subset $U\subset \mathbb{F}^n$, we might be able to more efficiently describe $\M$ (at least locally) using only the polynomial ring $\mathbb{F}[x_1,...,x_n]\subset D(\mathbb{F}^n,\mathbb{F})$ instead of the full collection of functions on all open subsets of $\M$.

For a further discussion of how ``ringed spaces'' (Definition \ref{RingSpDfn}) such as $(\M,S)$ above arise from classical geometry, see for example \cite[pp 93-100]{gelfand-manin2010}, \cite[Section 5.4, pp 273-303]{rotman2009}, and \cite[Section 6.3, pp 377-403]{rotman2009}.

\section{Schemes in Algebraic Geometry: A Brief Touch}\label{AlgGeomSec}
In view of the earlier encountered relationship between pseudo-topological systems and sheafs, this section considers/realizes rings as sources of topological systems and geometry. For further details on this section, see for example \cite[Chapter II.2, pp 69-82]{hartshorne1977} and \cite[2002/2003 version, pp 74-91]{gathmann2002}.

Algebraic geometry involves the study of geometric objects/structures based on commutative rings $R$ such that (in the spirit of representation theory) each commutative ring $R$ is viewed as a set of maps (functions) on some underlying topological space (called the \emph{spectrum} of $R$) with a natural topology called \emph{Zariski topology}. The following basic concepts introduced for an arbitrary commutative ring $R$ are motivated (as evident in the literature on the subject) by extensive studies involving the special case where $R=k[x_1,...,x_n]$, i.e., a commutative $k$-algebra or polynomial ring over a given field $k$.

\begin{dfn}[\textcolor{blue}{
\index{Spectrum of a commutative ring}{Spectrum of a commutative ring},
\index{Zero! map}{Zero map},
\index{Zariski topology}{Zariski topology},
\index{Distinguished open sets}{Distinguished open sets},
\index{Radical of an ideal}{Radical of an ideal},
\index{Nilradical of a commutative ring}{Nilradical of a commutative ring},
\index{Local! ring at a prime ideal}{Local ring at a prime ideal}}]
Let $R$ be a commutative ring. The \ul{spectrum} of $R$ is the set of prime ideals ~$\Spec R:=\{\txt{prime ideals}~p\lhd R\}$ ~of $R$. (\blue{footnote}\footnote{The spectrum of a commutative ring is often specified in the the form $(\Spec R,S)$ as a topological space $\Spec R$ along with a sheaf of rings $S$ on $\Spec R$, but for us at this early stage the spectrum is merely a set.}). View $R$ as a set of maps
{\footnotesize
\bea
\textstyle R\subset\prod_{p\in\Spec R}Q(R/p):=\left\{\txt{maps}~f:\Spec R\ra\bigcup_{p\in\Spec R}Q(R/p),~f(p)\in Q(R/p)\right\},\nn
\eea} where $Q(R/p):=(R/p)[(R/p-\{0\})^{-1}]=\{{r+p\over a+p}:r,a\in R,~a\not\in p\}$ is the ring of quotients of the ring $R/p$. That is, $r\in R$ is a map as follows (along with pointwise addition and multiplication in $R$):
{\small\bea
\textstyle r:\Spec R\ra\bigcup\limits_{p\in\Spec R}Q(R/p),~p\mapsto r(p):={r+p\over 1+p}\in Q(R/p),~~~~\ub{(r+r')(p):=r(p)+r'(p)}_{\txt{pointwise addition}},~~\ub{(rr')(p):=r(p)r'(p)}_{\txt{pointwise multiplication}}.\nn
\eea}The \ul{zero map} of $R$ is the map $Z:\P(R)\ra\P(\Spec R)$ given (for $A\subset R$) by
\bea
&&\textstyle Z(A):=\{p\in\Spec R:r(p)=0~\txt{for all}~r\in A\}=\{p\in\Spec R:~p\supset A\}=\{p\in\Spec R:p\supset {I(A)}\}\nn\\
&&\textstyle ~~~~=Z(I(A))\cong\Spec{R\over{I(A)}},~~~~\txt{where ${I(A)}:=RA\lhd R$ is the ideal of $R$ generated by $A$.}\nn
\eea
Observe that because $Z(A)=Z({I(A)})$, with $\I(R):=\{\txt{ideals}~I\lhd R\}\subset\P(R)$, it suffices to define the \ul{zero map} of $R$ by its restriction:
\bea
\textstyle Z|_{\I(R)}:\I(R)\subset\P(R)\ra\P(\Spec R),~~I\mapsto Z(I):=\{p\in\Spec R:~p\supset I\}\cong\Spec{R\over I}.\nn
\eea
By Lemma \ref{ZarTopLmm}, the zero map $Z$ defines the closed sets of a topology on $\Spec R$ called \ul{Zariski topology}. Henceforth, unless stated otherwise, we will assume $\Spec R$ is a topological space with the Zarisky topology. The \ul{distinguished open subset} of $\Spec R$ associated with an element $r\in R$ is the set
\bea
(\Spec R)_r:=(\Spec R)\backslash Z(r)=\{p\in \Spec R~|~r\not\in p\}=\{p\in \Spec R~|~p\subset R\backslash r\}.\nn
\eea

For any ideal $I\in\I(R)$, the \ul{radical} of $I$ (or \ul{nilradical} of the ring $R/I$) is the ideal of $R$ given by
\bea
\sqrt{I}:=\{r\in R:r^n\in I~\txt{for some}~n\geq 1\}=\{r\in R:(r+I)^n=0~\txt{in}~R/I~\txt{for some}~n\geq 1\}.\nn
\eea
The \ul{local ring at a prime ideal} $p\in\Spec R$ is the localization of $R$ by the localizing set $R-p$, i.e.,
\bea
\textstyle R_p:=R[(R-p)^{-1}]=\left\{{a\over b}:a,b\in R,~b\not\in p\right\},\nn
\eea
where it was proved earlier that {\footnotesize $R_p=(R_p,\mfm_p)$, $\mfm_p:=p[(R-p)^{-1}]=\{{a\over b}:a\in p,b\not\in p\}$}, is indeed a local ring.
\end{dfn}

\begin{lmm}[\blue{Zariski topology}]\label{ZarTopLmm}
Let $R$ be a commutative ring with zero map $Z:\I(R)\ra\P(\Spec R)$. For any ideals ~$I,J,I_i\in\I_R$,~ we have
\bit
\item[] (a) $Z(I)\cup Z(J)=Z(IJ)$, ~and~ (b) $\bigcap_i Z(I_i)=Z\big(\sum_i I_i\big)$.
\eit
\end{lmm}
\begin{proof}
(a): Let $p\in\Spec R$. If $I\subset p$ or $J\subset p$, then {\small $IJ\subset pJ\subset p$} or {\small $IJ\subset Ip\subset p$}, and so {\small $Z(I)\cup Z(J)\subset Z(IJ)$}. Also, if $IJ\subset p$, then $I\subset p$ or $J\subset p$ (since $p$ is prime), and so {\small $Z(IJ)\subset Z(I)\cup Z(J)$}.\\
(b): $p\supset I_i$~ for all $i$ $\iff$ $p\supset\sum I_i$,~ because $\sum I_i$ is the smallest ideal containing all of the ideals $I_i$.
\end{proof}

Lemmas \ref{MaxPrimeLmm} and \ref{ZarTopLmm0} below present some properties of the zero-map $Z$ that are useful for further studies, but which are not of immediate relevance to us.
\begin{lmm}\label{MaxPrimeLmm}
Let $R$ be a commutative ring, $S\subset R$ a localizing set, and $J$ a maximal element of the set $\I(R)\cap(R\backslash S):=\{I\lhd R~|~I\cap S=\emptyset\}$ of ideals of $R$ that are disjoint from $S$. Then $J$ is a prime ideal.
\end{lmm}
\begin{proof}
By the maximality of $J$, for any $r\in R$, we have $(J+Rr)\cap S\neq\emptyset$ $\iff$ $r\not\in J$. Let $a,b\in R$ be such that $a,b\not\in J$. Then by the maximality of $J$, we have $(J+Ra)\cap S\neq\emptyset$ and $(J+Rb)\cap S\neq\emptyset$. Since $S$ is closed under multiplication, we also have $\big((J+Ra)(J+Rb)\big)\cap S\neq\emptyset$, which implies $(J+Rab)\cap S\neq\emptyset$ as well (since $(J+Ra)(J+Rb)\subset J+Rab$). Hence, $ab\not\in J$ (by the maximality of $J$).
\end{proof}

\begin{lmm}\label{ZarTopLmm0}
Let $R$ be a commutative ring with zero map $Z:\I(R)\subset\P(R)\ra\P(\Spec R)$. For any ideals $I,J\in\I(R)$, we have
\bit
\item[](i) $Z(\sqrt{I})=Z(I)$,~ (ii) $Z(I)=Z(J)\iff\sqrt{I}=\sqrt{J}$, ~(iii)~$Z(I)\subset Z(J)$ $\iff$ $\sqrt{I}\supset\sqrt{J}$.
\eit
\end{lmm}
\begin{proof}
{\flushleft (i)} Since $Z$ is clearly inclusion reversing, $Z(\sqrt{I})\subset Z(I)$ and $I\subset\sqrt{I}$. Also,
 \bea
&&{p}\in Z(I)~\Ra~I\subset{p}~~\Ra~~\sqrt{I}\subset\sqrt{{p}}={p},~~\Ra~~p\in Z(\sqrt{I}),~~\Ra~~Z(I)\subset Z(\sqrt{I}).\nn
\eea
{\flushleft (ii)} If $Z(I)=Z(J)$, then ~$\sqrt{I}~\sr{(s)}{=}~\bigcap\limits_{I\subset{p}\in Spec(R)}{p}=\bigcap\limits_{{p}\in Z(I)}{p}=\bigcap\limits_{{p}\in Z(J)}{p}\sr{(s1)}{=}\bigcap\limits_{J\subset{p}\in Spec(R)}{p}=\sqrt{J}$,~ where step (s) is due to the following: If $I\subset{p}$, then
\bea
\textstyle r\in\sqrt{I}~~\Ra~~r^t\in I\subset{p},~~\Ra~~r\in{p},~~\Ra~~\sqrt{I}\subset p,~~\Ra~~ \sqrt{I}\subset\bigcap\limits_{I\subset{p}\in Spec(R)}{p}.\nn
\eea
On the other hand, let $r\in R$ such that $r\not\in\sqrt{I}$. Then $r^t\not\in I$ for all $t$, and so the localizing set $S:=\{1,r,r^2,\cdots\}$ satisfies $I\cap S=\emptyset$ (and $\sqrt{I}\cap S=\emptyset$). By Zorn's lemma, the nonempty set of ideals $\{J\lhd R~|~I\subset J,~J\cap S=\emptyset\}\ni I$, as a poset under inclusion, has a maximal element $q$. By Lemma \ref{MaxPrimeLmm}, $q$ is prime. That is, $R\backslash S$ contains a prime ideal $q$ containing $I$. Hence, $r\not\in \bigcap\limits_{I\subset{p}\in Spec(R)}{p}$.
{\flushleft (iii)} As shown in the proof of part (ii) above, $\sqrt{I}$ = $\bigcap\{p\in Z(I)\}$ = the intersection of all prime ideals containing $I$. Hence, $Z(I)\subset Z(J)$ $\iff$ $\sqrt{I}\supset\sqrt{J}$.
\end{proof}

\begin{rmk}[\blue{Intersection up to representatives of equivalence classes, Locally constant selection map}]
Let $R$ be a commutative ring and $q,q'\in\Spec R$. Then we have $R_q\cap R_{q'}=\emptyset$ if $q=q'$. This is because in $R_q={R\times(R\backslash q)\over\sim_q}$ and $R_{q'}={R\times(R\backslash q')\over\sim_{q'}}$ the equivalence relations $\sim_q$ and $\sim_{q'}$ (and hence the localizing maps $h_q:R\ra R_q$ and $h_{q'}:R\ra R_{q'}$) are different if $q\neq q'$. Nevertheless, we will introduce a notion of ``intersection'' between $R_q$ and $R_{q'}$, even if $q\neq q'$.

Given $a,b\in R$ with $b\not\in q$, we will write ${a\over b}|_q:=[(a,b)]_{\sim q}\in R_q$ to mean ``${a\over b}$ as an element of $R_q$''. Let $A\subset\Spec R$. The phrase ``$a,b\in R$ with $b\not\in q$ for all $q\in A$'' will be also be written as ``$a,b\in R$ with ${a\over b}|_q\in R_q$ for all $q\in A$'', or equivalently, as ``$a,b\in R$ with ${a\over b}\in\bigcap^{\txt{rep}}_{q\in A}R_q$'' (\ul{intersection up to representatives of equivalence classes}).

Consider an open set $U\subset\Spec R$ and a \ul{selection map} $f:U\ra\bigcup_{q\in U}R_q,~p\mapsto f(p)\in R_p$. We will say $f$ is \ul{locally constant} if for each $p\in U$, there exists an open set $p\in V\subset U$ such that $f|_V=\txt{constant}$ in the sense there exist $a,b\in R$, with $b\not\in q$ for all $q\in V$ (i.e., with ${a\over b}|_q\in R_q$ for each $q\in V$), such that
\[
\textstyle f(q)={a\over b}|_q:=[(a,b)]_{\sim q}\in R_q:={R\times(R\backslash q)\over\sim_q}~~~~\txt{for all}~~~~q\in V.
\]
For convenience, we will restate this as follows: $f$ is \ul{locally constant} if for each $p\in U$, there exists an open set $p\in V\subset U$ such that (i) `` $f|_V={a\over b}$ for some ${a\over b}\in\bigcap_{q\in V}^{\txt{res}}R_q$'', or equivalently, (ii) ``$f|_V=\txt{constant}={a\over b}$'' (for some $a,b\in R$, $b\not\in q\in V$) or equivalently, (ii) ``$f(q)={a\over b}$ for all $q\in V$'' (for some $a,b\in R$, $b\not\in q\in V$).
\end{rmk}

\begin{dfn}[\textcolor{blue}{
\index{Ringed space}{Ringed space},
\index{Structure sheaf of a ringed space}{Structure sheaf of a ringed space},
\index{Morphism of! ringed spaces}{Morphism of ringed spaces},
\index{Locally! ringed space}{Locally ringed space},
\index{Morphism of! locally ringed spaces}{Morphism of locally ringed spaces},
\index{Affine scheme}{Affine scheme},
\index{Locally! constant section}{Locally constant section},
\index{Scheme}{Scheme},
\index{Morphism of! schemes}{Morphism of schemes},
\index{Category of! schemes}{Category of schemes}}]\label{RingSpDfn}
A \ul{ringed space} $(X,S)$ is a topological space $X$ together with a sheaf of rings $S:\I_X\ra\txt{Rings}$ called the \ul{structure sheaf} of $X=(X,S)$. A \ul{morphism of ringed spaces} $f:(X,S)\ra (X',S')$ is a continuous map $f:X\ra X'$ that induces a morphism of (pre)sheafs $f^\ast:S'\ra S\circ f^{-1}$, i.e., there exists a family of ring homomorphisms $f^\ast_{U'}:S'(U')\ra S(f^{-1}(U'))$, one for each open set $U'\subset X'$, giving commutative diagrams
\bea\bt
 U'\ar[d,hook,"i_{U'V'}"] && S'(U')\ar[rr,"f^\ast_{U'}"]\ar[from=d,"{S(i_{U'V'})}"'] && S(f^{-1}(U'))\ar[from=d,"{S(i_{f^{-1}(U')f^{-1}(V')})}"']\\
 V'                  && S'(V')\ar[rr,"f^\ast_{V'}"]  && S(f^{-1}(V')).\nn
\et
\eea

A ringed space $(X,S)$ is a \ul{locally ringed space} if for each $p\in X$, the stalk $S_p$ is a local ring (i.e., a ring with a unique maximal ideal). We will denote the maximal ideal of $S_p$ by $m_p\subset S_p$. A \ul{morphism of locally ring spaces} {\small$f:(X,S)\ra(X',S')$} is a morphism of ringed spaces that is compatible with the local-ring nature of stalks in the following sense: For each $p\in X$, the map induced on stalks by {\small$f^\ast:S'\ra S\circ f^{-1}$}, i.e.,
\bea
\textstyle f_p^\ast:=\varinjlim_{U'\ni f(p)}f^\ast_{U'}:S'_{f(p)}\ra S_p,~(U',s')\mapsto \left(f^{-1}(U'),f^\ast_{U'}(s')\right),\nn
\eea
when restricted to the maximal ideals ~$m_p\subset S_p$, ~$m'_{f(p)}\subset S'_{f(p)}$ ~satisfies
\bea
(f_p^\ast)^{-1}(m_p)=m'_{f(p)}.\nn
\eea

An \ul{affine scheme} is a ringed space $(X,S)$ for which there exists a commutative ring $R$ such that (i) $X=\Spec R$ and (ii) the structure sheaf $S:\I_X\ra \txt{Rings}$ is given (for any open set $U\subset X$) by the following subring consisting of ``\ul{locally constant sections}'':
{\small\[
\textstyle S(U):=\big\{s:U\ra\bigcup_{q\in U}R_q~\big|~\txt{$s(p)\in R_p$, for each $p\in U$, and $s$ is locally constant}\big\}~\subset~\prod_{p\in U}R_p.
\]}
By Lemma \ref{LocRingSch}, an affine scheme is a locally ringed space. A \ul{scheme} is a locally ringed space {\small $(X,S)$} that locally resembles an affine scheme in that there exists an open cover {\small $\{U_i\}$} of {\small $X$} such that for each $i$,
\bea
(U_i,S|_{U_i})\cong(\Spec R_i,S_i),~~~~\txt{for some commutative ring $R_i$},\nn
\eea
where the isomorphism $\cong$ is specified in the category of locally ringed spaces. A \ul{morphism of schemes} $f:(X,S)\ra(X',S')$ is a morphism of locally ringed spaces.
\end{dfn}

\begin{lmm}[\textcolor{blue}{\cite[Proposition 2.2(a), p.71]{hartshorne1977} }]\label{LocRingSch}
If $X=(X,S)=(\Spec R,S)$ is an affine scheme, then $S_p\cong R_p$ for each $p\in\Spec R$.
\end{lmm}
\begin{proof}
Fix $p\in\Spec R$. Consider the (clearly well defined) map $\vphi:S_p\ra R_p$, $[(U,s)]\mapsto s(p)$.
\bit[leftmargin=0.9cm]
\item[(i)] \ul{$\vphi$ is surjective}: Fix ${a\over b}\in R_p$. Let $X_b:=(\Spec R)_b=\{q\in\Spec R:b\not\in q\}\subset\Spec R$ be the distinguished open set associated with $b$. Let $s:X_b\ra\bigcup_{q\in X_b}R_q,~q\mapsto s(q)\in R_q$ be the ``constant'' section $s\in S(X_b)$ given by $s(q):={a\over b}$ for all $q\in X_b$. Since $p\in X_b$ (recall ${a\over b}\in R_p$ implies $b\not\in p$), we see that ${a\over b}=s(p)=\vphi([(X_b,s)])$.

\item[(ii)] \ul{$\vphi$ is injective}: Let $(U,s),(U',s')\in S_p$ be such that $\vphi([(U,s)])=\vphi([(U',s')])$, i.e., $s(p)=s'(p)$. Since $s,s'$ are each locally constant, there exists a nbd $p\in V\subset U\cap U'$ such that $s|_V={a\over b}\in\bigcap^{\txt{rep}}_{q\in V}R_q$ and $s'|_V={a'\over b'}\in\bigcap^{\txt{rep}}_{q\in V}R_q$. But $s(p)=s'(p)$ implies ${a\over b}={a'\over b'}$ in $R_p$ (i.e., ${a\over b}|_p={a'\over b'}|_p$), which is the case $\iff$ (by the kernel of the localization $h_p:R\ra R_p$) in $R$ there exists $b_p\in R\backslash p$ such that
    \bea
    (ab'-a'b)b_p=0,\nn
    \eea
    which (again by the kernel of the localization $h_q:R\ra R_q$) in turn implies ${a\over b}={a'\over b'}$ in every ring $R_q$ ($q\in V$) such that $b_p\in R\backslash q$, i.e., ${a\over b}={a'\over b'}$ in every ring $R_q$ such that $q\in V\cap(\Spec R)_{b_p}$, and so
    \bea
    &&\textstyle s|_{V\cap(\Spec R)_{b_p}}=s'|_{V\cap(\Spec R)_{b_p}}={a\over b}={a'\over b'}\in\bigcap^{\txt{rep}}_{q\in V\cap(\Spec R)_{b_p}}R_q,\nn\\
    &&\textstyle~~\Ra~~(U,s)\sim(U',s'),~~~~\txt{since}~~~~p\in V\cap(\Spec R)_{b_p}\subset U\cap U'.\nn \hspace{3cm}\qedhere
    \eea
\eit
\end{proof}

By Theorem \ref{ProtShRepThm}, the protosheaf representation of $(\Spec R,S)$ has the form $(E,f,\Spec R)$, where $E:=\bigcup_{p\in\Spec R}S_p\cong\bigcup_{p\in\Spec R}R_p$ (with a suitable topology) and ~$f:E\ra \Spec R,~e\mapsto p$ if $e\in E_p:=R_p$~ is a surjective local homeomorphism.

\begin{dfn}[\textcolor{blue}{\index{Affine variety}{Affine variety}, \index{Variety}{Variety}}]
An \ul{affine variety} is an affine scheme $(\Spec A,S)$ such that $A$ is a (finitely generated) $k$-algebra, for a field $k$. A \ul{variety} $(X,S)$ is a scheme that locally resembles an affine variety in the following sense: There exists an open cover $\{U_i\}$ of $X$ such that for each $i$,
\bea
(U_i,S|_{U_i})\cong(\Spec A_i,S_i),~~~~\txt{for some (finitely generated) $k_i$-algebra $A_i$}.\nn
\eea
\end{dfn}

{\flushleft\hrulefill}

\begin{exercise}
Based on the discussion of this chapter, consider writing a \emph{fully technical} essay (say in the form of a typical section of this chapter) on what is known in the mathematics literature as the \index{Riemann-Roch Theorem}{``Riemann-Roch Theorem''}. If you consider topological analysis to be an essential part of the discussion, then you can of course refer to further discussions on topology from the later chapters of these notes.
\end{exercise}

%% file: parts/AlgebraM/GeomAnaII.tex
\chapter{Geometry and Analysis II: Measure Theory and Topological Analysis}\label{GeomAnaII}

For a further reading on familiar material related to this chapter see for example  \cite{folland-book,wheeden-zygmund} on measure theory and \cite{conway} on functional analysis.

\begin{note}[\textcolor{blue}{Recall: Convention for countability}]
As in the previous chapter, for brevity, ``countable'' will continue to mean ``at most countable'' (i.e., ``finite or countable'') unless it is specified otherwise.
\end{note}

In constructing measurable systems, we will attempt to parallel the construction for topological systems. In particular, (i) ``\emph{finite monic limit-closedness}'' will be replaced with ``\emph{monic cocomplement-closedness}'' (to be defined), while (ii) ``\emph{monic colimit-closedness}'' will be replaced with ``\emph{countable monic colimit-closedness}''. Similarly, (i) ``\emph{finite epic colimit-closedness}'' will be replaced with ``\emph{epic complement-closedness}'' (to be defined), while (ii) ``\emph{epic limit-closedness}'' will be replaced with ``\emph{countable epic limit-closedness}''.

\section{Measurable Systems and Measure Systems}
\subsection{Basic definitions and preliminary remarks}
As mentioned earlier, measurable systems below (as well as topological systems seen earlier) are meant to be examples of $Sets$-representable systems in a category. Also, as done earlier for topology, in order to quickly establish a smooth connection between our new categorical definitions and those of familiar (topological) measure theory, we will begin with an unusually long definition containing several footnotes.

In the following definition, items with the familiar/classical prefix ``\magenta{Borel-}'' are colored simply to isolate them for easy identification.
\begin{dfn}[\textcolor{blue}{
\index{Subsystem}{Subsystem},
\index{Countable! subsystem}{Countable subsystem},
\index{Sub(co)limit}{Sub(co)limit},
\index{Countable! sub(co)limit}{Countable sub(co)limit},
\index{Monic! (co)limit-closed system}{Monic (co)limit - closed system},
\index{Countable! monic (co)limit-closed system}{Countable monic (co)limit-closed system},
\index{Monic! (co)complement}{Monic (co)complement},
\index{Epic! (co)complement}{Epic (co)complement},
\index{Monic! (co)complement-closed system}{Monic (co)complement-closed system},
\index{Epic! (co)complement-closed system}{Epic (co)complement-closed system},
\index{Sigma! (Measurable system) in a category}{Sigma (Measurable system) in a category},
\index{Pseudo-sigma (Pseudo-measurable system)}{Pseudo-sigma (Pseudo-measurable system)},
\index{Generated! sigma}{Generated sigma},
\index{Generator of! (Subbase of) a sigma}{Generator (Subbase) of a sigma},
\index{Base for a sigma}{Base for a sigma},
\index{Sigma! (Sigma algebra) on a set}{Sigma (Sigma algebra) on a set},
\index{Measurable! space}{Measurable space},
\index{Measurable! sets}{Measurable sets},
\index{Sigma! subalgebra}{Sigma subalgebra},
\index{Measurable! subspace}{Measurable subspace},
\index{Borel! algebra of a topological space}{\textcolor{magenta}{Borel algebra of a topological space}},
\index{Borel! space}{\textcolor{magenta}{Borel space}},
\index{Measurable! map}{Measurable map},
\index{Borel! map}{\textcolor{magenta}{Borel map}},
\index{Category of! measurable spaces (Measurable)}{Category of measurable spaces (Measurable)},
\index{Isomorphism of! measurable spaces}{Isomorphism of measurable spaces},
\index{Isomorphic! measurable spaces}{Isomorphic measurable spaces}}]~

Let $\I,\C$ be categories, where $\C$ has an initial objects $U\in\Ob\C$ (and a finial object $U'\in\Ob\C$), and let $S:\I\ra\C$ be a system. If $\J\subset\I$ is a subcategory, then $S|_{\J}:\J\subset\I\ra\C$ is a \ul{subsystem} of $S$. (\blue{footnote}\footnote{Recall that this (map restriction)-based notion of a subsystem (although related) is different from the notion of a subfunctor introduced earlier in a fixed category of functors $\D^\C$. For greater generality (if necessary), we can also say that given systems $S:\I\ra\C$, $S':\I'\ra\C$, $S$ is a subsystem of $S'$, written $S\subset S'$, if (i) $\I\subset\I'$ is a subcategory and (ii) every object of $S$ is a subobject of an object of $S'$, i.e., for each $i\in\Ob\I$, $S_i\subset S_{i'}'$ is a subobject for some $i'\in\Ob\I'$.}). A subsystem $S|_{\J}:\J\subset\I\ra\C$ is a \ul{countable subsystem} if $\Ob\J$ is a countable set. A \ul{sublimit} (resp. \ul{subcolimit}) of $S$ is the limit (resp. colimit) of a subsystem of $S$. A \ul{countable sub(co)limit} of $S$ is the (co)limit of a countable subsystem.

A system $S:\I\ra\C$ is \ul{monic limit-closed} (resp. \ul{monic colimit-closed}) if every monic sublimit $\varprojlim_{mo} S|_{\J}$ (resp. monic subcolimit $\varinjlim_{mo} S|_{\J}$) is an object of $S$, i.e., $\varprojlim_{mo} S|_{\J}=S(i_\J)$ (resp. $\varinjlim_{mo} S|_{\J}=S(i_\J)$) for some $i_\J\in\Ob\I$. Similarly, a system $S:\I\ra\C$ is \ul{countable monic limit-closed} if every countable monic sublimit $\varprojlim_{mo} S|_{\J}$ (resp. \ul{countable monic colimit-closed} if every countable monic subcolimit $\varinjlim_{mo} S|_{\J}$) is an object of $S$.

Let $S:\I\ra\C$ be a system, $i\in\Ob\I$, and $A:=S_i:=S(i)$. The object $A$ is
\bit
\item \ul{monic cocomplemented} in $S$ if there exists a monic subcolimit $B:=\varinjlim_{mo}S|_{\J\subset\I}$ such that (i) $A\cap B:=\varprojlim_{mo}\{A,B\}\cong U$ and (ii) $A\cup B:=\varinjlim_{mo}\{A,B\}\cong\varinjlim_{mo}S$.
    The \ul{monic cocomplement} of $A$ in $S$ is ~$(S-A)_{\underrightarrow{mo}}:=B$.
\item \ul{epic complemented} in $S$ if there exists an epic sublimit $E:=\varprojlim_{ep}S|_{\J\subset\I}$ such that (i) $\varinjlim_{ep}\{A,E\}\cong U'$ and (ii) $\varprojlim_{ep}\{A,E\}\cong\varprojlim_{ep}S$.
    The \ul{epic complement} of $A$ in $S$ is ~$(S-A)_{\underleftarrow{ep}}:=E$.
\eit
The system $S$ is \ul{monic cocomplement-closed} if for every $i\in\I$, the object $S_i$ is monic cocomplemented in $S$, and $(S-S_i)_{\underrightarrow{mo}}=S_{i'}$ for some $i'\in\I$. The system $S$ is \ul{epic complement-closed} if for every $i\in\I$, the object $S_i$ is epic complemented in $S$, and $(S-S_i)_{\underleftarrow{ep}}=S_{i'}$ for some $i'\in\I$.

A system $\Sigma:\I\ra\C$ is an \ul{$\I$-sigma} (or \ul{measurable $\I$-system}) in $\C$ if (i) $\Sigma$ is monic cocomplement-closed and (ii) $\Sigma$ is countable monic colimit-closed. (\blue{footnote}\footnote{Equivalently (up to a contravariant equivalence of categories $\C\ra\C^{op}$), a system $\Sigma:\I\ra\C$ is an \ul{$\I$-sigma} (or \ul{measurable $\I$-system}) in $\C$ if $\Sigma$ is both (i) epic complement-closed and (ii) countable epic limit-closed.}). By removing ``monic'', a system $\Sigma:\I\ra\C$ is an \ul{$\I$-pseudo-sigma} (or \ul{pseudo-measurable $\I$-system}) in $\C$ if (i) $\Sigma$ is cocomplement-closed and (ii) $\Sigma$ is countable colimit-closed.

If $\Sigma:\I\ra\C$ is a sigma, $\I$ is directed, and $\{\Sigma(i)\}_{i\in\Ob\I}$ has a unique maximal object $M=\Sigma(i)$ for some $i$, then we may call $\Sigma$ a \ul{sigma on the object} $M\in\Ob\C$. Given a system $S:\I\ra\C$, the \ul{sigma generated} by $S$ is the smallest sigma $\sigma(S):\J\ra\C$ containing $S$ as a subsystem, i.e., such that $\I\subset\J$ is a subcategory and $S=\sigma(S)|_\I$. Given a sigma $\Sigma:\J\ra\C$, a system $S:\I\ra\C$ is a \ul{generator (or subbase)} of $\Sigma$ if $\Sigma=\sigma(S)$ (where it is also understood that $\I\subset \J$ is a subcategory). Given a sigma $\Sigma:\J\ra\C$, a system $S:\I\ra\C$ is a \ul{base} for $\Sigma$ if (i) $\Sigma=\sigma(S)$ and (ii) every object of $\Sigma$ is a countable monic subcolimit of $S$ (i.e., the monic colimit of a countable collection of objects of $S$), i.e., for each $j\in\Ob\J$ we have $\Sigma_j=\varinjlim_{mo} S_{\al_j(i)}$ for some map $\al_j:\I\ra\J$ with $\al_j(\I)\subset\J$ countable.

Let $X$ be a set and $\X=\Big(\P(X),Hom_\X\big(\P(X),\P(X)\big),\circ\Big)\subset Sets$ the full subcategory whose objects $\Ob\X:=\P(X)$ are the subsets of $X$, and $Hom_\X(A,B):=Hom_{Sets}(A,B)$ maps from $A$ to $B$, for any subsets $A,B\subset X$. (\blue{footnote}\footnote{If necessary, we may also consider $X$ to be an object in any category $\C$ with an initial object $I\in\Ob\C$ (replacing $\emptyset$ in $\C=Sets$), and let $\X=\Big(\P(X),Hom_\X\big(\P(X),\P(X)\big),\circ\Big)\subset \C$ be the full subcategory whose objects $\Ob\X:=\P(X)$ are the subobjects of $X$, and $Hom_\X(A,B):=Hom_\C(A,B)$, for any subobjects $A,B\subset X$.}). Then a sigma $\Sigma:\I\ra\X$ such that $X$ is an object of $\Sigma$ is called a \ul{sigma (or sigma algebra) on $X$}, making $X=(X,\Sigma)$ a \ul{measurable space} with the sets $\{\Sigma(i)\subset X\}_{i\in\Ob\I}$  called \ul{$\Sigma$-measurable subsets} (or just \ul{measurable subsets}) of the measurable space $X$. (\blue{footnote}\footnote{Traditionally, we specify a sigma for a measurable space $X$ in terms of \emph{measurable sets}, namely, a collection of subsets $\M\subset\P(X)$ such that (i) $X\in\M$, (ii) $M\in\M$ $\Ra$ $M^c\in\M$, and (iii) $\{O_k\in\M\}_{k\in C}$ (with $C$ a countable set) $\Ra$ $\bigcup_{k\in C} M_k\in\M$ (i.e., $\M$ contains $X$ and is closed under the complement and countable union operations).}). If ${X'}\subset X$ is $\Sigma$-measurable, i.e., ${X'}=\Sigma(i)$ for some $i\in\I$, the \ul{sigma subalgebra} of $\Sigma_{X'}\subset\Sigma$ on ${X'}$ (making $({X'},\Sigma_{X'})$ a \ul{measurable subspace} of $(X,\Sigma)$) is the sigma on $X'$ given by $\Sigma_{X'}:\I\ra\X'\subset\X,~i\mapsto {X'}\cap\Sigma(i)$, or more explicitly, by
\[
\textstyle\Sigma_{X'}:\I\ra\X'\subset\X,~~i\sr{\kappa_{ij}}{\ral}j~\mapsto~ {X'}\cap\Sigma(i)\sr{\Sigma(\kappa_{ij})|}{\ral}{X'}\cap\Sigma(j).
\]

\textcolor{magenta}{Let $X=(X,\O)$ be a topological space, with topology $\O:\I\ra\X$. The \ul{Borel algebra} of $X$ is $\B:=\sigma(\O):\J\supset\I\ra\X$, i.e., the sigma algebra on $X$ generated by open sets of $X$. The associated measurable space $(X,\B)=(X,\O,\B)$ is the \ul{Borel space} of $X=(X,\O)$.}

Given measurable spaces $X=(X,\Sigma_X:\I\ra\X)$ and $Y=(Y,\Sigma_Y:\J\ra\Y)$, and a map $f:X\ra Y$ (or more explicitly, $f:(X,\Sigma_X)\ra(Y,\Sigma_Y)$), consider any process
\bea
P_f:\Y\ra\X,~B\sr{v}{\ral}B'~\mapsto~f^{-1}(B)\sr{P_f(v)}{\ral}f^{-1}(B').\nn
\eea
We say the map $f:X\ra Y$ (or more explicitly, $f:(X,\Sigma_X)\ra(Y,\Sigma_Y)$) is a \ul{measurable map} if the process $P_f:\Y\ra\X$ takes objects of $\Sigma_Y$ to objects of $\Sigma_X$, i.e., for each $j\in\Ob\J$, ~$P_f(\Sigma_Y(j)):=f^{-1}(\Sigma_Y(j))=\Sigma_X(i)$ for some $i\in\Ob\I$.

Let $X,Y$ be topological spaces, and $\Sigma=\Sigma_X$ a sigma algebra on $X$. We say a map $f:X\ra Y$ is \ul{$\Sigma_X$-measurable} if it is measurable as the map of measurable spaces $f:(X,\Sigma_X)\ra(Y,\B_Y)$, where $\B_Y$ is the Borel algebra of $Y$. \textcolor{magenta}{In particular, we say a map $f:X\ra Y$ is \ul{Borel-measurable} (or a \ul{Borel map}) if it is $\B_X$-measurable, i.e., in the sense it is measurable as a map of Borel spaces $f:(X,\B_X)\ra(Y,\B_Y)$.}

The \ul{category of measurable spaces}, denoted by Measurable, is the category whose objects are measurable spaces and its morphisms $Hom_{Measurable}(X,Y)$ are measurable maps $f:X\ra Y$. An isomorphism $X\sr{f}{\ral}Y$ in Measurable (\blue{footnote}\footnote{That is, $f$ is a measurable bijection and the inverse map (bijection-inverse) of $f$ is itself a measurable map.}) is called an \ul{isomorphism of measurable spaces} (making $X,Y$ \ul{isomorphic measurable spaces}).
\end{dfn}

\begin{rmk}[\textcolor{blue}{Nature of a generated sigma}]
Let $X$ be a set, $\C\subset\P(X)$, $\D:=\C\cup\{X\}$, and
\bea
\textstyle S:=\{id,{}^c,\bigcap,\bigcup\}^\Natural=\big\{s:\Natural\ra\{id,{}^c,\bigcap,\bigcup\}\big\}\nn
\eea
sequences of operations on sets. Then (from the definitions) the sigma algebra ~$\sigma(\C):S\times\D^\Natural\ra\P(X)$~ on $X$ generated by $\C$ is given by
\bea
\sigma(\C)=\{A\subset X:\txt{$A$ is formed from sets in $\C\cup\{X\}$ via a sequence of operations ~$s\in S$}\}.\nn
\eea
Some members of $\sigma(\C)$ can be obtained as follows. Let $\D^c:=\{A^c\subset X:A\in\D\}$, $\A_1:=\D\cup\D^c$, $\A_2:=\bigcap^{\txt{fin}}\A_1^\Natural:=\{\bigcap_{i\in\Natural}f(i)~|~\txt{finite}~f:\Natural\ra\A_1\}$, $\A_3:=\bigcup\A_2^\Natural:=\{\bigcup_{i\in\Natural}g(i)~|~g:\Natural\ra\A_2\}$, and
\[
\textstyle \A_4:=\A_3\cup\A_3^c,~~~~~\txt{where}~~~~\A_3^c:=\{A^c:A\in\A_3\}.
\]
Given ~$g:\Natural\ra\A_2,~i\mapsto g(i)$,~ where ~$g(i):\Natural\ra\A_1,~j\mapsto f_i(j)$,~ we have
\begin{align}
&\textstyle\left(\bigcup_ig(i)\right)^c=\left(\bigcup_i\bigcap_jf_i(j)\right)^c=\bigcap_i\bigcup_jf_i(j)^c=\bigcup_{j_1,j_2,\cdots}f_1(j_1)^c\cap f_2(j_2)^c\cap\cdots=\bigcup_{j_1,j_2,\cdots}~\bigcap_if_i(j_i)^c\nn\\
&\textstyle~~~~=\bigcup_{j_1,j_2,\cdots}\left(\bigcup_if_i(j_i)\right)^c~~\in~~\sigma(\C).\nn
\end{align}
This shows that ~$\A_4\subset\sigma(\C)$. This holds mainly because of the finiteness condition in the definition of $\A_2$.
\end{rmk}

\begin{question}
In the above remark, is it true that ~$\sigma(\C)=\A_4$?
\end{question}

\begin{rmk}[\textcolor{blue}{\index{Sigma! induced}{Sigmas induced via a map}}]
Let $(X,\Sigma_X),(Y,\Sigma_Y)$ be measurable spaces, and $f:X\ra Y$. Then through $f$, (i) $\Sigma_Y$ induces a sigma on $X$ given by
\bea
f^{-1}(\Sigma_Y):=\{f^{-1}(B):B\in\Sigma_Y\},~~~~\txt{(where $f$ is a measurable map iff ~$f^{-1}(\Sigma_Y)\subset\Sigma_X$)}\nn
\eea
and (ii) $\Sigma_X$ induces a sigma on $Y$ given by
\bea
f^\ast(\Sigma_X):=\{B\subset Y: f^{-1}(B)\in\Sigma_X\}.\nn
\eea
Observe that if $f$ is surjective, then ~$f^\ast(\Sigma_X)=f\big(\Sigma_X\cap f^{-1}\big(\P(Y)\big)\big)$, where
\[
f\left(\Sigma_X\cap f^{-1}(\P(Y))\right):=\{f(f^{-1}(B)):B\subset Y,~f^{-1}(B)\in\Sigma_X\}=\{B\cap f(X):B\subset Y,~f^{-1}(B)\in\Sigma_X\}.
\]
\end{rmk}

\begin{dfn}[\textcolor{blue}{
\index{Filtered! system}{Filtered system},
\index{Disjoint! system}{Disjoint system},
\index{Monotone functor}{Monotone functor},
\index{Submeasure!}{Submeasure},
\index{Measure!}{Measure},
\index{Measure! space}{Measure space},
\index{Measure! preserving map}{Measure preserving map},
\index{Category of! measure spaces (Measure)}{Category of measure spaces (Measure)},
\index{Isomorphism of! measure spaces}{Isomorphism of measure spaces},
\index{Isomorphic! measure spaces}{Isomorphic measure spaces}}]\label{SubmeasMeasDefn}
Let $\I,\C,\D$ be categories, where $\C,\D$ have initial objects $U\in\Ob\C$ and $V\in\Ob\D$, and $\D$ is an infinity category (i.e., there is an object $\infty\in\Ob\D$ such that $\coprod\{\infty,D\}\cong\infty$ for all $D\in\Ob\D$).  A system $S:\I\ra\C$ is a \ul{filtered system} if $\Ob\I$ is linearly ordered and for all $i,j\in\Ob\I$, if $i\leq j$ then $S(i)\subset S(j)$. A system $S:\I\ra\C$ is a \ul{disjoint system} if $S_i\cap S_{i'}\cong U$ for all distinct $i,i'\in\I$ (i.e., $i\neq i'$). A functor $F:\C\ra\D$ is a \ul{monotone functor} if for all $C,C'\in\Ob\C$, if $C\subset C'$ then $F(C)\subset F(C')$.

Let $S:\I\ra\C$ be a system in $\C$ and $\Sigma:\I\ra\C$ a sigma in $\C$. A \ul{$\D$-submeasure} on $S$ (a submeasure on $S$ with coefficients in $\D$) is a monotone functor $\mu:\C\ra\D$ such that (i) $\mu(U)\cong V$, and (ii) for every subsystem $S|_\J:\J\subset\I\ra\C$, we have ~$\mu\left(\varinjlim^{mo}_{j\in\J}S_j\right)\subset\varinjlim^{mo}_{j\in\J}\mu(S_j)$. Similarly, a \ul{$\D$-measure} on $\Sigma$ (a measure on $\Sigma$ with coefficients in $\D$) is a monotone functor $\mu:\C\ra\D$ such that (i) $\mu(U)\cong V$, and (ii) for every countable filtered subsystem (or if valid, for every countable disjoint subsystem) $\Sigma|_\J:\J\subset\I\ra\C$, we have ~$\mu\left(\varinjlim^{mo}_{j\in\J}\Sigma_j\right)\cong\varinjlim^{mo}_{j\in\J}\mu(\Sigma_j)$. (\blue{footnote}\footnote{Note that the following system does not have to be a filtered (resp., disjoint) system:
{\tiny\[
\mu\circ\Sigma|_\J:\J\subset\I:\I\sr{\Sigma}{\ral}\C\sr{\mu}{\ral}\D,~i\sr{\kappa}{\ral}j~\mapsto~\mu(\Sigma_i)\sr{\mu\circ\Sigma(\kappa_{ij})}{\ral}\mu(\Sigma_j).
\]}
}).

{\flushleft\ul{Notation:}} For brevity, we will write $\mu:\C\supset\Sigma\ra\D$ to mean ``$\mu$ is a $\D$-(sub)measure on $\Sigma$''.

Let $X=(X,\Sigma)$ be a measurable space, $\A$ an infinity additive category, and $\mu:\X\supset\Sigma\ra\A$ an $\A$-measure on $\Sigma$, i.e., $\mu$ is a monotone functor such that (i) $\mu(\emptyset)\cong0$ and (ii) every countable filtered (resp., disjoint) collection $\{E_i\}_{i\in C}\subset\Sigma$ (\blue{footnote}\footnote{$\{E_i\}_{i\in C}$ is a countable disjoint collection if $C$ is a countable set and $E_i\cap E_j=E_i\cap\delta_{ij}$, where $\delta_{ij}:=X$ if $i=j$, and $\delta_{ij}:=\emptyset$ if $i=j$.}) satisfies $\mu(\bigcup_{i\in C} E_i)\cong \varinjlim^{mo}_{i\in C}\mu(E_i)$. Then the $\A$-measure $\mu$ makes $X=(X,\Sigma)=(X,\Sigma,\mu,\A)$ an \ul{$\A$-measure space}. If $X,Y$ are spaces and $\mu_X,\mu_Y$ are $\A$-measures, a measurable map $f:(X,\Sigma_X,\mu_X,\A)\ra (Y,\Sigma_Y,\mu_Y,\A)$ is \ul{$\A$-measure preserving} if, in $\A$, we have $\mu_X(f^{-1}(B))\cong\mu_Y(B)$ for all $B\in\Sigma_Y$.

The \ul{category of $\A$-measure spaces}, denoted by Measure($\A$), is the category whose objects are $\A$-measure spaces and its morphisms $Hom_{Measure(\A)}(X,Y)$ are $\A$-measure preserving maps $f:X\ra Y$. An isomorphism $X\sr{f}{\ral}Y$ in Measure($\A$), i.e., an \ul{isomorphism of $\A$-measure spaces} (making $X,Y$ \ul{isomorphic $\A$-measure spaces}), is an $\A$-measure preserving bijection whose inverse map, i.e., bijection-inverse, is itself $\A$-measure preserving.
\end{dfn}

\subsection{Classical measurable spaces and classical measure spaces}
This section is mainly concerned with Borel measurability, completion, and Lebesgue measurability.

\begin{dfn}[\textcolor{blue}{Recall:
\index{Borel! algebra}{Borel algebra} (of a topological space),
\index{Borel! space}{Borel space},
\index{Borel! map}{Borel map},
\index{Topologically measurable map}{Topologically measurable map}}]
Let $X=(X,\O)$ be a topological space, with topology $\O:\I\ra\X$. The \ul{Borel algebra} of $X$ is $\B:=\sigma(\O):\J\supset\I\ra\X$, i.e., the sigma algebra on $X$ generated by open sets of $X$. The associated measurable space $(X,\B)$ is the \ul{Borel space} of $X$.

Let $X,Y$ be topological spaces. We say a map $f:X\ra Y$ is \ul{Borel-measurable} (or a \ul{Borel map}) if it is $\B_X$-measurable, i.e., if it is measurable as a map of Borel spaces $f:(X,\B_X)\ra(Y,\B_Y)$.

Given any collection of sets $\C\subset\P(X)$, a map $f:X\ra Y$ is said to be \ul{(topologically) $\C$-measurable} if $f^{-1}(V)\in\C$ for every open set $V\subset Y$. Observe that if $f:X\ra Y$ is topologically $\C$-measurable, then $f:(X,\sigma(\C))\ra (Y,\B_Y)$ is measurable (i.e., $f$ is $\sigma(\C)$-measurable).
\end{dfn}

\begin{question}
Let $X,Y$ be topological spaces and $\C\subset\P(X)$. If a map $f:X\ra Y$ is $\sigma(\C)$-measurable, does it follow that $f$ is topologically $\C$-measurable?
\end{question}

\begin{dfn}[\textcolor{blue}{
\index{Null set}{Null set},
\index{Complete! measure (space)}{Complete measure (space)},
\index{Completion of a measure (space)}{Completion of a measure (space)},
\index{Borel! measure (space)}{Borel measure (space)},
\index{Lebesgue! measure (space)}{Lebesgue measure (space)},
\index{Lebesgue! algebra}{Lebesgue algebra},
\index{Lebesgue! space}{Lebesgue space},
\index{Lebesgue! map}{Lebesgue map}}]
Let $\I$ be a category and $\A$ an infinity additive category.

Let $X$ be a set and $\mu:\X\supset\Sigma\ra\A$ a measure on $\Sigma$, i.e., $X=(X,\Sigma,\mu,\A)$ is a measure space. A $\Sigma$-measurable set $N\in\Sigma$ is called a \ul{null set} if $\mu(N)\cong 0$. The measure $\mu$ is a \ul{complete measure} (making $(X,\Sigma,\mu,\A)$ a \ul{complete measure space}) if subsets of null sets are measurable, i.e., for every null set $N\in\Sigma$, we have $\P(N)\subset\Sigma$. The \ul{completion} $\ol{\mu}:\X\supset\ol{\Sigma}\ra\A$ of $\mu:\X\supset\Sigma\ra\A$ (or equivalently, the completion $\ol{(X,\Sigma,\mu,\A)}:=(X,\ol{\Sigma},\ol{\mu},\A)$ of $(X,\Sigma,\mu,\A)$) is the complete measure (or equivalently, the complete measure space) on the sigma algebra $\ol{\Sigma}\subset\P(X)$ given by
\bea
\ol{\Sigma}:=\{E\cup Y:E\in\Sigma,~Y\subset N\in\Sigma,~\mu(N)\cong 0\},~~~~\ol{\mu}(E\cup Y):=\mu(E).\nn
\eea

Let $X$ be a topological space and $\mu:\X\supset\Sigma\ra\A$ a measure on $\Sigma$, i.e., $X=(X,\Sigma,\mu,\A)$ is a measure space. The measure $\mu:\X\supset\Sigma\ra\A$ is a \ul{Borel measure} (making $X=(X,\Sigma,\mu,\A)$ a \ul{Borel measure space}) if $\Sigma=\B:=\B_X$ (i.e., the Borel algebra of $X$). That is, a Borel measure is a measure on the Borel algebra.

A \ul{Lebesgue measure (space)} is the completion (i.e., a complete) Borel measure (space). That is, a measure $\ld:\X\supset\L\ra\A$ is a \ul{Lebesgue measure} (making $X=(X,\L,\ld,\A)$ a \ul{Lebesgue measure space}) if $(X,\L,\ld,\A)=\ol{(X,\B,\mu,\A)}:=(X,\ol{\B},\ol{\mu},\A)$ for a Borel measure $\mu:\X\supset\B\ra\A$.

Let $(X,\B,\mu,\A)$ be a Borel measure space and $(X,\L,\ld,\A)=\ol{(X,\B,\mu,\A)}:=(X,\ol{\B},\ol{\mu},\A)$ the associated Lebesgue measure space. In analogy with the Borel algebra $\B$ and the Borel space $(X,\B)$, we will call $\L:=\ol{\B}$ the \ul{$\mu$-Lebesgue algebra} of $X$ and $(X,\L):=(X,\ol{\B})$ the \ul{$\mu$-Lebesgue space} of $X$ (i.e., the Lebesgue space of $X$ associated with $\mu$). Accordingly, if $(X,\B_X,\mu,\A)$ is a Borel measure space and $Y$ any topological space, we will say a Borel map $f:(X,\B_X)\ra (Y,\B_Y)$ is \ul{Lebesgue measurable} (or a \ul{Lebesgue map}) if it is $\L_X$-measurable, where $\L_X=\L_{X,\mu}:=\ol{\B}_X$, i.e., it is measurable as the map of measurable spaces $f:(X,\L_X)\ra(Y,\B_Y)$.
\end{dfn}

\subsection{The topology of pointwise convergence}
\begin{dfn}[\textcolor{blue}{\index{Topology! of poinwise convergence}{Topology of poinwise convergence}}] (\blue{footnote}\footnote{Some geometric applications based on variants of the topology of pointwise are considered in \cite{aguilar.etal2002} for example.}).
Let $X$ be a set and $Y$ a space. Recall that a natural topology on $Y^X:=\{\txt{maps}~f:X\ra Y\}=\{\txt{maps}~f:X\ra Y,~x\mapsto f(x)\in Y\}=\prod_{x\in X}Y$ is the product topology (in the present case called the \ul{topology of pointwise convergence}), which has a subbase consisting of sets of the following form (based on the definition (\ref{ProdTopBasEq}), page \pageref{ProdTopBasEq}, with $I:=X$ and $X_i:=Y$ for each $i\in I$):~For all elements $x\in X$ and all open sets $V\subset Y$,
{\small\bea
\textstyle V^x:=V\times\big(\prod_{x'\in X\backslash\{x\}}Y\big)=\{f:X\ra Y~|~f(x)\in V\}=\left\{f\in Y^X~|~f(x)\in V\right\}.\nn
\eea}That is, an open set in $Y^X$ is a union of finite intersections of sets of the form $\{V^x:x\in X,~V\subset Y\}$, i.e.,
\bea
V_{1,\cdots,n}^{x_1,...,x_n}:=V_1^{x_1}\cap\cdots\cap V_n^{x_n}=\left\{f\in Y^X~|~f(x_i)\in V_i\right\},~~\txt{for all}~~n\geq 1,~~x_i\in X,~~\txt{open}~~V_i\subset Y.\nn
\eea
Observe that {\small $V_1^x\cap\cdots\cap V_n^x=(V_1\cap\cdots\cap V_n)^x$}, and so with {\small $U:=V_1\cap\cdots\cap V_n$, $W:=V_1\cup\cdots\cup V_n$}, we have
\bea
U^{x_1}\cap\cdots\cap U^{x_n}\subset V_1^{x_1}\cap\cdots\cap V_n^{x_n}\subset W^{x_1}\cap\cdots\cap W^{x_n}.\nn
\eea
If $Y$ is \ul{Hausdorff} (Definition \ref{SpaceTypeDfn}, page \pageref{SpaceTypeDfn}), then the topology of $Y^X$ (which is then also clearly Hausdorff: \blue{footnote}\footnote{Assume $Y$ is Hausdorff. If $f,g:X\ra Y$ and there is $x_0\in X$ such that $f(x_0)\neq g(x_0)$, let $V,W\subset Y$ be disjoint open sets with $f(x_0)\in V$ and $g(x_0)\in W$. Then $f\in V^{x_0}$, $g\in W^{x_0}$, and $V^{x_0}\cap W^{x_0}=(V\cap W)^{x_0}=\emptyset^{x_0}=\emptyset$.})) has a base with only open sets of the following form (\blue{footnote}\footnote{Let $f\in V_1^{x_1}\cap\cdots\cap V_n^{x_n}$, which is the case $\iff$ $f(x_i)\in V_i$ for $i=1,...,n$. So, if $\{f(x_i):i\in A\subset\{1,...,n\}\}$ are distinct, then with disjoint open sets $f(x_i)\in W_i\subset V_i$ for $i\in A$, we get $f\in \bigcap_{i\in A}W_i^{x_i}\subset V_1^{x_1}\cap\cdots\cap V_n^{x_n}$.}):
{\footnotesize\bea
V_{1,\cdots,n}^{x_1,...,x_n}:=V_1^{x_1}\cap\cdots\cap V_n^{x_n}=\left\{f\in Y^X~|~f(x_i)\in V_i\right\},~~\txt{for all}~~n\geq 1,~~\txt{distinct}~~x_i\in X,~~\txt{disjoint open}~~V_i\subset Y.\nn
\eea}
\end{dfn}

Unless it is stated otherwise, we will assume that subsets of $Y^X$ are given the subspace topology induced by the product topology (i.e., the topology of pointwise convergence). The discussion in the following remark requires a fast-forward to section \ref{GeomPrelimsS4} on page \pageref{GeomPrelimsS4}.
\begin{rmk}[\textcolor{blue}{The topology of $Y^X$ is that of pointwise convergence (see section \ref{GeomPrelimsS4} on page \pageref{GeomPrelimsS4})}] Let $X$ be a set and $Y$ a space. (i) $Y$ can be seen as a subspace of $Y^X$ as follows:
\[
Y\cong\{\txt{constant maps}~X\ra Y\}\subset Y^X.
\]
{\flushleft (ii)} (With say $f=g_\al$) $f\ra g$ in $Y^X$ $\iff$ $f(x)\ra g(x)$ in $Y$ for every $x\in X$. The proof is as follows:

($\Ra$) Assume $f\ra g$ in $Y^X$. Let $x_0\in X$. Then for every open set $V\subset Y$, $f\ra g$ in the subspace $V^{x_0}\subset Y^X$. Since $f\in \txt{some nbd}_{Y^X}(g)\subset V^{x_0}$ $\iff$ $f(x_0)\in \txt{some nbd}_Y(g(x_0))\subset V$, with $V:=\txt{nbd}_Y(g(x_0))$, we see that $f(x_0)\ra g(x_0)$ in $V\subset Y$ (i.e., $f(x_0)\ra g(x_0)$ in $Y$).

($\La$) Assume $f(x)\ra g(x)$ in $Y$ for all $x\in X$. Suppose $f=g_\al\not\ra g$ in $Y^X$, i.e., there exists a base neighborhood $g\in N(g):= V_1^{x_1}\cap\cdots\cap V_n^{x_n}\subset Y^X$ such that $N(g)$ excludes a member of every tail of $g_\al$. Let $g_{\al_\ld}\not\in N(g)$ for all $\ld\in\Ld$. Then some $g_{\al_{\ld_i}}\not\in V_i^{x_i}$ for all $\ld_i\in\Ld_i\subset\Ld$, i.e., $g_{\al_{\ld_i}}(x_i)\not\in V_i$ for all $\ld_i\in\Ld_i\subset\Ld$. Since there are only finitely many $i$, we get a contradiction: Indeed $g_{\al}(x_i)\ra g(x_i)$ in $Y$ implies the nbd $V_i\ni g(x_i)$ contains a tail of $g_{\al}(x_i)$.
\end{rmk}

\begin{question}
Let $X$ be a set and $Y$ a space. (i) Let $x\in X$. Is the evaluation map $e_x:Y^X\ra Y,~f\mapsto f(x)$ continuous? (\blue{footnote}\footnote{Recall that $Y^X=\prod_{x\in X}Y$, and so the canonical projections $p_x=e_x:Y^X\ra Y$ are continuous by definition. Alternatively, we can also proceed as follows: In general, it is of course enough to show $e_x$ is continuous on a base set $V_1^{x_1}\cap\cdots\cap V_n^{x_n}=\left\{f\in Y^X~|~f(x_i)\in V_i\right\}$. However, an answer can be obtained directly using a continuity criterion given later in
Proposition \ref{NetContPrp}. Thus, the answer is yes, due to the pointwise convergence property of the topology of $Y^X$.}). (ii) Let $A\subset Y^X$, and consider the map $A:X\ra \P(Y),~x\mapsto A(x):=e_x(A)=\{f(x):f\in A\}$. Is it true that {\footnotesize $A=\{\txt{maps}~f:X\ra \bigcup_{x\in X} A(x),~x\mapsto f(x)\in A(x)\}=\prod_{x\in X}A(x)$}? When is $A$ compact in $Y^X$? (\textcolor{magenta}{footnote}\footnote{\textcolor{magenta}{A partial answer can be obtained using Tychonoff's theorem (Theorem \ref{TychThm}), namely, if it happens that there exists a continuous surjective map $h:\prod_{x\in X}e_x(A)\ra A$, then $A$ is compact in $Y^X$ iff each $e_x(A)$, $x\in X$, is compact in $Y$ (recall that a continuous map takes a compact set to a compact set).}}).
\end{question}

\begin{question}[\textcolor{blue}{Closed subsets of $Y^X$}]
Let $\{X_i\}_{i\in I}$ be spaces. If $C_i\subset X_i$ are closed in $X_i$ for each $i$, then $\prod_{i\in I}C_i\subset\prod_{i\in X}X_i:\{x:I\ra\bigcup_{i\in I} X_i,~i\mapsto x_i\in X_i\}$ is closed as an intersection of closed sets, since
\bea
\textstyle\prod_{i\in I}C_i=\left\{x\in\prod_{i\in I} X_i:x_i\in C_i,i\in I\right\}=\bigcap_{i\in I}\left\{x\in\prod_{i\in I}X_i:x_i\in C_i\right\}=\bigcap_{i\in I}\wt{C}_i,\nn
\eea
where $\wt{C}_i:=\left\{x\in\prod_{i\in I}X_i:x_i\in C_i\right\}=C_i\times\prod_{j\neq i}X_j$ is closed because of the open complement
\bea
\textstyle(\wt{C}_i)^c=\left\{x\in\prod_{i\in I}X_i:x_i\not\in C_i\right\}=C_i^c\times\prod_{j\neq i}X_j.\nn
\eea
Let $X$ be a set and $Y$ a space. Is it true that a set $A\subset Y^X$ is closed if and only if $e_x(A):=\{f(x):f\in A\}\subset Y$ is closed for each $x\in X$? (\textcolor{magenta}{footnote}\footnote{\textcolor{magenta}{As done for compactness in the previous question, we have a partial answer: If it happens that (i) there exists a homeomorphism $A\cong\prod_{x\in X}e_x(A)$ and (ii) the evaluation maps $e_x:Y^X\ra Y$ are closed, then $A$ is closed in $Y^X$ iff each $e_x(A)$, $x\in X$, is closed in $Y$.}})
\end{question}

\begin{dfn}[\textcolor{blue}{\index{Reflexive space}{Reflexive space}}]
Let $X,Y$ be spaces, $X_Y^\ast:=Cont(X,Y)\subset Y^X$ the space of continuous maps $X\ra Y$, and {\small $X_Y^{\ast\ast}:=(X_Y^\ast)_Y^\ast=Cont(X_Y^\ast,Y)=Cont\big(Cont(X,Y),Y\big)$}. The space $X$ is reflexive wrt $Y$ (or $Y$-reflexive) if $X\cong X_Y^{\ast\ast}$, i.e., if the following evaluation map is an isomorphism:
\bea
e:X\ra X_Y^{\ast\ast},~x\mapsto e_x,~~e_x(f):=f(x),~~\txt{for all}~~f\in X_Y^\ast.\nn
\eea
\end{dfn}

\section{Functional Analysis: A Descriptive Framework}
Functional analysis involves the study of (categories of) topological algebraic objects such as topological multiplicative sets, topological groups, topological rings, topological modules, and topological algebras. The term \emph{functional} will be reserved for a \emph{general map}, i.e., a \emph{type of process}, between \emph{categories}.

\subsection{Topological setup and measuring}
\begin{dfn}[\textcolor{blue}{
\index{Topological! binary operation}{Topological binary operation} (\index{Topological! multiplication}{Topological multiplication}),
\index{Topological! multiplicative set}{Topological multiplicative set},
\index{Topological! group}{Topological group},
\index{Topological! ring}{Topological ring},
\index{Topological! module}{Topological module},
\index{Difference-invariant topological module}{Difference-invariant topological module},
\index{Difference-invariance scale}{Difference-invariance scale (DI-scale)},
\index{Topological! algebra}{Topological algebra}}]
Let $X,Y,Z$ be topological spaces. A multiplication or binary operation (i.e., any map) $X\times Y\ra Z$ is \ul{topological} if it is continuous (where $X\times Y$ is given the product topology by default).

A \ul{topological multiplicative set} is a topological space with a topological (i.e., continuous) binary operation: that is,  a multiplicative set $S=[S,\cdot]$ together with a topology with respect to which multiplication $\cdot:(S\times S)\ra S$ is continuous.

A group $G=[G,\cdot,()^{-1}]$ is a \ul{topological group} if it is (i) a topological space and (ii) multiplication $\cdot:G\times G\ra G$ and inversion $()^{-1}:G\ra G$ are continuous. The category of topological groups is the common subcategory $\txt{Groups}\cap\Top$ of both $\txt{Groups}$ and $\Top$ (whose morphisms are continuous group homomorphisms).

A ring $R=[R,\cdot,+]$ is a \ul{topological ring} if it is (i) a topological space and (ii) addition and multiplication $+,\cdot:R\times R\ra R$ are continuous. The category of topological rings is the common subcategory $\txt{Rings}\cap\Top$ of both $\txt{Rings}$ and $\Top$ (whose morphisms are continuous ring homomorphisms).

An $R$-module $M=[M,+,R]$ is a \ul{topological $R$-module} if (i) $R$ is a topological ring, (ii) $(M,+)$ is a topological group, and (iii) scalar multiplication $R\times M\ra M$ is continuous. A topological $R$-module $M$ is \ul{difference-invariant} if there exists $r_0\in R$ (call it a \ul{difference-invariance scale} or \ul{DI-scale} of $M$) such that for every neighborhood $U_0\subset M$ of $0$, we have the containment
\[
U_0-U_0\subset r_0U_0,~~~~\txt{where}~~U_0-U_0:=\{a-b:a,b\in U_0\}~~\txt{and}~~r_0U_0:=\{r_0u:u\in U_0\}.
\]
The category of topological $R$-modules is the common subcategory $(R\txt{-mod})\cap\Top$ of both $R$-mod and $\Top$ (whose morphisms are continuous $R$-module homomorphisms).

Let $R$ be a commutative ring. An $R$-algebra  $A=[A,\cdot,+,R]$ is a \ul{topological $R$-algebra} if (i) $(A,\cdot,+)$ is a topological ring and (ii) $(A,+,R)$ is a topological $R$-module. The category of topological $R$-algebras is the common subcategory $(R\txt{-Alg})\cap\Top$ of both $R$-Alg and $\Top$ (whose morphisms are continuous $R$-algebra homomorphisms).
\end{dfn}

\begin{dfn}[\textcolor{blue}{
\index{Topological! module of maps}{Topological module of maps},
\index{Topological! algebra of maps}{Topological algebra of maps}}]
Let $X$ be a set. Given a topological $R$-module $M$, the \ul{topological $R$-module of maps} $f:X\ra M$ is the space
\bea
\T(X,M):=M^X,\nn
\eea
as a topological $R$-module with addition $+:\T(X,M)\times\T(X,M)\ra\T(X,M)$, scalar multiplication $\cdot:R\times\T(X,M)\ra\T(X,M)$, and ``zero'' given respectively by
\bea
(f+g)(x):=f(x)+g(x),~~(rg)(x):=rf(x),~~0_{\M(X,M)}(x):=0_R,~~~~\txt{for all}~~~~x\in X,~r\in R.\nn
\eea

Given a topological $R$-algebra $A$ (for a commutative ring $R$), the \ul{topological $R$-algebra of maps} $f:X\ra A$ is the topological $R$-module of such maps $\T(X,A)\subset A^X$ as a ring with multiplication $\cdot:\T(X,A)\times\T(X,A)\ra\T(X,A)$ and ``one'' given respectively by
\bea
(fg)(x):=f(x)g(x),~~1_{\T(X,A)}(x):=1_R,~~~\txt{for all}~~x\in X.\nn
\eea
\end{dfn}

\begin{rmk}[\blue{\magenta{Caution}: The above operations are \ul{not} automatically continuous. Continuity conditions}]
Let $X$ be a set, $M$ a topological $R$-module (any $R$), and $A$ a topological $R$-algebra (commutative $R$).
{\flushleft\bf Continuity of addition}: To show addition $+:\T(X,M)\times \T(X,M)\ra \T(X,M)$ is continuous, it suffices (by properties of the product topology) to show that for each $x\in X$ (with $e_x$ denoting the evaluation at $x$, i.e., the $x$th canonical projection), the associated partial addition below is continuous:
\[
+_x:=e_x\circ +:\T(X,M)\times \T(X,M)\sr{+}{\ral}\T(X,M)\sr{e_x}{\ral}M,~~(f,g)\mapsto f(x)+g(x).
\]
Moreover, by the discussion in Definition \ref{PtwsContDfn} (page \pageref{PtwsContDfn}), we only need to show that given $x\in X$, $f,g\in\T(X,M)$, and an open neighborhood $(e_x\circ +)(f,g)=f(x)+g(x)\in V\subset M$, we can find open neighborhoods $f\in O_1\subset \T(X,M)$ and $g\in O_2\subset \T(X,M)$ such that
\[
(O_1+O_2)(x):=(e_x\circ+)(O_1\times O_2)\subset V,
\]
where $(O_1+O_2)(x)=\{f'(x)+g'(x):f'\in O_1,g'\in O_2\}$. That is, given $x\in X$ and $f,g\in\T(X,M)$,
\[
f(x)+g(x)\in V~~\Ra~~f\in O_1~~\txt{and}~~g\in O_2~~~\txt{that satisfy}~~~(O_1+O_2)(x)\subset V.
\]
The easiest case is to consider $O_1=V_1^x$ and $O_2=V_2^x$ for open sets $V_1,V_2\subset M$, in which case, the continuity statement reads as follows: Given $x\in X$ and $f,g\in\T(X,M)$, it \ul{suffices} (not necessarily) that
\[
f(x)+g(x)\in V~~\Ra~~f(x)\in V_1~~\txt{and}~~g(x)\in V_2~~~\txt{that satisfy}~~~(V_1^x+V_2^x)(x)\subset V,
\]
where $(V_1^x+V_2^x)(x)\subset V$ means that for any $f',g'\in\T(X,M)$,
\[
f'(x)\in V_1~~\txt{and}~~g'(x)\in V_2~~\Ra~~f'(x)+g'(x)\in V.\nn
\]
That is, to show addition on $\T(X,M)$ is continuous, it \ul{suffices} (but is not necessary) to show that for each $x\in X$, $f,g\in \T(X,M)$, and an open nbd $f(x)+g(x)\subset V\subset M$, there exist open nbds $f(x)\in V_1\subset M$ and $g(x)\in V_2\subset M$ such that $V_1+V_2\subset V$, or briefly: It \ul{suffices} that given $x\in X$ and $f,g\in\T(X,M)$,
\[
f(x)+g(x)\in V~~\Ra~~f(x)\in V_1~~\txt{and}~~g(x)\in V_2~~~\txt{that satisfy}~~~V_1+V_2\subset V.
\]

{\flushleft\bf Continuity of multiplication}: It is clear by symmetry between addition and multiplication that the above discussion also holds for multiplication (in place of addition). That is, to show multiplication $\cdot:\T(X,A)\times \T(X,A)\ra\T(X,A)$ is continuous, it suffices to show that for each $x\in X$, the partial multiplication below is continuous:
\[
\cdot_x:=e_x\circ \cdot:\T(X,M)\times \T(X,M)\sr{\cdot}{\ral}\T(X,M)\sr{e_x}{\ral}M,~~(f,g)\mapsto f(x)g(x).
\]
Moreover, by the discussion in Definition \ref{PtwsContDfn} (page \pageref{PtwsContDfn}), we only need to show that given $x\in X$, $f,g\in\T(X,A)$, and an open set $(e_x\circ \cdot)(f,g)=f(x)g(x)\in V\subset A$, we can find open sets $f\in O_1\subset \T(X,A)$ and $g\in O_2\subset \T(X,A)$ such that
\[
(O_1O_2)(x):=(e_x\circ\cdot)(O_1\times O_2)\subset V,
\]
where $(O_1O_2)(x)=\{f'(x)g'(x):f'\in O_1,g'\in O_2\}$. That is, given $x\in X$ and $f,g\in\T(X,M)$,
\[
f(x)g(x)\in V~~\Ra~~f\in O_1~~\txt{and}~~g\in O_2~~~\txt{that satisfy}~~~(O_1O_2)(x)\subset V.
\]
In particular, to show multiplication on $\T(X,A)$ is continuous, it \ul{suffices} (but is not necessary) to show that for each $x\in X$, $f,g\in \T(X,A)$, and an open nbd $f(x)g(x)\subset V\subset A$, there exist open nbds $f(x)\in V_1\subset A$ and $g(x)\in V_2\subset A$ such that $V_1V_2\subset V$, or briefly: It \ul{suffices} that given $x\in X$ and $f,g\in\T(X,M)$,
\[
f(x)g(x)\in V~~\Ra~~f(x)\in V_1~~\txt{and}~~g(x)\in V_2~~~\txt{that satisfy}~~~V_1V_2\subset V.
\]
\end{rmk}

Despite the main point in the above remark, that addition and multiplication in $\M(X,M)$ or $\M(X,A)$ are not automatically continuous in general, we will for simplicity assume the topologies of $M$ and $A$ have been chosen such that addition and multiplication in $\M(X,M)$ or $\M(X,A)$ are continuous as desired.

\begin{dfn}[\textcolor{blue}{
\index{Space of! measurable maps}{Space of measurable maps},
\index{Topological! module of measurable maps}{Topological module of measurable maps},
\index{Topological! algebra of measurable maps}{Topological algebra of measurable maps}}]
Let $X,Y$ be topological spaces and $M$ a topological $R$-module (resp. topological $R$-algebra). Let {\small $(X,\Sigma_X)$, $(Y,\Sigma_Y)$, $(M,\Sigma_M)$} be measurable spaces. The \ul{space of measurable maps} from $(X,\Sigma_X)$ to $(Y,\Sigma_Y)$ is the subspace
\bea
\M(X,Y)=\M(X,Y~|~\Sigma_X,\Sigma_Y):=Hom_{\txt{Measurable}}(X,Y)\subset Y^X.\nn
\eea
Similarly, the \ul{topological $R$-module (resp. $R$-algebra) of measurable maps} from $(X,\Sigma_X)$ to $(M,\Sigma_M)$ is the subspace
\bea
\M(X,M)=\M(X,M~|~\Sigma_X,\Sigma_M):=Hom_{\txt{Measurable}}(X,M)\subset\T(X,M)=M^X.\nn
\eea
\end{dfn}

\begin{dfn}[\textcolor{blue}{
\index{Space of! continuous maps}{Space of continuous maps},
\index{Space of! continuous linear maps}{Space of continuous linear maps},
\index{Dual space of a topological module}{Dual space of a topological module}}]
Let $X,Y$ be topological spaces and $M,N$ topological $R$-modules. The \ul{space of continuous maps} from $X$ to $Y$ is the subspace
\bea
\C(X,Y):=Hom_{\Top}(X,Y)\subset Y^X.\nn
\eea
The \ul{space of continuous linear maps} from $M$ to $N$ is the subspace
\bea
\C\L(M,N):=Hom_{(R\txt{-mod})\cap\Top}(M,N)\subset Hom_{\Top}(M,N)=\C(M,N).\nn
\eea
In particular, the \ul{dual space} of $M$ is the subspace
\bea
M^\ast:=\C\L(M,R)=Hom_{(R\txt{-mod})\cap\Top}(M,R)\subset\C(M,R).\nn
\eea
\end{dfn}

We will now introduce our notion of a functional, which is defined on a category with morphisms being maps that are both continuous and measurable. However, more general variants can also be considered if necessary (say to allow functionals to include/generalize submeasures and measures on sigmas in arbitrary categories as given in Definition \ref{SubmeasMeasDefn}).

\begin{dfn}[\textcolor{blue}{
\index{Topological! measurable space (TMS)}{Topological measurable space (TMS)},
\index{Category of! TMS's}{Category of TMS's},
\index{Functional}{Functional},
\index{Continuous functional}{Continuous functional},
\index{Linear! functional}{Linear functional}}]
A measurable space $X=(X,\Sigma_X)$ is a \ul{TMS} (topological measurable space)  if $X$ is a topological space. The \ul{category of TMS's} $\Top\Sigma:=\Top\cap\txt{Measurable}$ is the category whose objects are TMS's $X=(X,\Sigma_X)$ and whose morphisms are continuous measurable maps, i.e.,
\bea
Hom_{\Top\Sigma}(X,Y):=\C\M(X,Y):=Hom_{\Top\cap\txt{Measurable}}(X,Y),~~~~\txt{for all}~~X,Y\in\Ob\Top\Sigma,\nn
\eea
where $\C\M(X,Y)$ denotes continuous measurable maps from $X$ to $Y$.

A \ul{functional} is any process (not necessarily functorial) of the form
\bea
F:\C\subset\Top\Sigma\ra \D\subset\Top\Sigma,~~(X,\Sigma_X)\sr{f}{\ral}(Y,\Sigma_Y)~~\mapsto~~(FX,\Sigma_{FX})\sr{F(f)}{\ral}(FY,\Sigma_{FY}).\nn
\eea
The functional $F$ is a \ul{continuous functional} if for all $X,Y\in\Ob\Top\Sigma$, the map ~{\small $F:Hom_{\Top\Sigma}(X,Y)\ra Hom_{\Top\Sigma}\big(F(X),F(Y)\big)$}~ is continuous. The functional $F$ is a \ul{linear functional} if $X,Y\in\Ob\Top\Sigma$, whenever $Y$ is a topological $R$-module (in which case, by definition/construction, $Hom_{\Top\Sigma}(X,Y)$ is an $R$-module), the map ~{\small $F:Hom_{\Top\Sigma}(X,Y)\ra Hom_{\Top\Sigma}\big(F(X),F(Y)\big)$}~ is an $R$-module homomorphism, i.e.,
\bea
F(f+g)=F(f)+F(g)~~\txt{and}~~F(rf)=rF(f)~~~\txt{for all}~~f,g\in Hom_{\Top\Sigma}(X,Y),~~r\in R.\nn
\eea
\end{dfn}

\subsection{Integration and differentiation}
Meanwhile our functional above acts on maps that are both continuous and measurable, our integral and derivative below act on maps that need only be measurable.

\begin{dfn}[\textcolor{blue}{
\index{Integral!}{Integral},
\index{Derivative}{Derivative},
\index{Integrable map}{Integrable map},
\index{Differentiable map}{Differentiable map},
\index{Anti-derivative}{Anti-derivative},}]\label{IntDerDef}
Let $X$ be a set, $M$ a topological $R$-module (any $R$), and $A$ a topological $R$-algebra (commutative $R$). Let $(X,\Sigma_X)$, $(M,\Sigma_M)$, $(A,\Sigma_A)$ be measurable spaces and $\M(X,M)$ the topological $R$-module (resp. $\M(X,A)$ the topological $R$-algebra) of measurable maps.

An \ul{(indefinite) integral} on $\M(X,M)$ is an $R$-homomorphism $I:\F\subset\M(X,M)\ra\M(X,M)$. (\blue{footnote}\footnote{If non $R$-linear integrals are desired, then $R$-linearity can be dropped accordingly. Continuity of the integral can be considered as well.}). In which case, the elements of $\F$ are said to be \ul{$I$-integrable maps}, and for each $f\in\F$ the image $I(f)\in\M(X,M)$ is called the \ul{(indefinite) integral} of $f$.

Let $I:\F\subset\M(X,A)\ra\M(X,A)$ be an integral. An \ul{$I$-derivative} on $\M(X,A)$ is an $R$-homomorphism $D:\G\subset\M(X,A)\ra\M(X,A)$ such that (i) $I|_{\F\cap D(\G)}\circ D=id_\G:\G\ra\G$, (ii) $D|_{\G\cap I(\F)}\circ I=id_\F:\F\ra\F$, and (iii) for every $g_1,g_2\in\G$, we have $D(g_1g_2)(x)=(Dg_1)(x)g_2(x)+g_1(x)(Dg_2)(x)$ for all $x\in X$. (\blue{footnote}\footnote{Just as for the integral, if non $R$-linear derivatives are desired, then $R$-linearity can be dropped accordingly. Also, continuity of the derivative can be considered.}). In this case, the elements of $\G\subset\M(X,A)$ are said to be \ul{$(I,D)$-differentiable}, and for each $g\in\G$, $x\in X$, the image $Dg$ is called the \ul{derivative} of $g$ and $(Dg)(x)$ is called the derivative of $g$ at $x$.

Let $I:\F\subset\M(X,A)\ra\M(X,A)$ be an integral and $D:\G\subset\M(X,A)\ra\M(X,A)$ an $I$-derivative. Then the integral $I$ is also called an \ul{anti-derivative} of the derivative $D$.
\end{dfn}

In typical applications, given a set $X$ and an $R$-algebra $A$, special integrals of the form $I:\M(X,A)\ra R\subset\M(X,A)$ with certain ``positivity'' properties or operations (e.g., based on a linear ordering on a subring $R_0\subset R$) are used to introduce new topologies on $\M(X,A)$. Of course, such a new topology in general can contain (i.e., be stronger than), lie in (i.e., be weaker than), intersect, or even avoid our default product topology on $\M(X,A)=A^X$.

\begin{rmk}[\textcolor{blue}{The integral as a functional}]
Let $X$ be a topological space, $M$ a topological $R$-module, and $(X,\Sigma_X)=(X,\Sigma_X:\J\ra\X)$, $(M,\Sigma_M)$ TMS's. In practice, a measure ~$\mu:\X\supset\Sigma_X\ra\txt{SubRings}(R)\subset Rings$~ on $\Sigma_X$ might uniquely define an integral as a \ul{linear selection map}
{\small\begin{align}
\label{IntSelectEq1}&\textstyle I=I_\mu:\F\subset\M(X,M)\ra\M(X,M),~~X\sr{f}{\ral}M~\mapsto~X\sr{I(f)}{\ral}M,\\
&\textstyle ~~~~\txt{given by}~~~I(f):x\mapsto I(f)(x)~\in~\varprojlim^{mo}_{\{j\in \J~|~x\in \Sigma_X(j)\}}\mu\left(\Sigma_X(j)\right)f\left(\Sigma_X(j)\right),\nn
\end{align}}where $\mu\left(\Sigma_X(j)\right)f\left(\Sigma_X(j)\right)$ denotes multiplication $R_jM_j$ of the subring $R_j:=\mu\left(\Sigma_X(j)\right)\subset R$ and the subset $M_j:=f\left(\Sigma_X(j)\right)\subset M$, and the limit is to be specified in the appropriate category (which is either Modules or Sets). (Note that $R_jM_j$ is an $R_j$-module since $R_j$ is a ring.)

In this case, we can also view the integral $I=I_\mu$ as a linear functional
{\small\bea
\textstyle I=I_\mu:\C\subset\Top\Sigma\ra \Top\Sigma,~~(X,\Sigma_X)\sr{f}{\ral}(Y,\Sigma_Y)~~\mapsto~~(X,\Sigma_X)\sr{I(f)}{\ral}(Y,\Sigma_Y).\nn
\eea}
\end{rmk}

\begin{dfn}[\textcolor{blue}{
\index{Fuzzy integration}{Fuzzy integration},
\index{Fuzzy differentiation}{Fuzzy differentiation}}]
Let $X$ be a topological space, $M$ a topological $R$-module, and $(X,\Sigma_X)=(X,\Sigma_X:\J\ra\X)$, $(M,\Sigma_M)$ TMS's. Further assume we also have a preferred topology on $\P(M)$ based on the topology of $M$, and view $\P(M)$ as an $R$-set with the obvious $R$-action
\[
\textstyle R\times\P(M)\ra\P(M),~(C,E)\mapsto CE:=\left\{\txt{finite sums}~\sum c_ie_i:c_i\in C,e_i\in E\right\}.
\]
Consider a measure ~$\mu:\X\supset\Sigma_X\ra\txt{SubRings}(R)\subset Rings$~ on $\Sigma_X$. Let $s:\P(M)\ra M,~P\mapsto s(P)\in P$ be a selection map. A \ul{fuzzy integral} on $\M(X,\P(M))$ is a modified version of the selection map in (\ref{IntSelectEq1}) in the form of a \ul{subset-selection map}
{\small\begin{align}
\label{IntSelectEq2}&\textstyle I=I_{\mu,s}:\F\subset\M(X,\P(M))\ra\M\big(X,\P(M)\big),~~X\sr{f}{\ral}\P(M)~\mapsto~X\sr{I(f)}{\ral}\P(M),\\
&\textstyle ~~~~\txt{given by}~~~I(f):x\mapsto I(f)(x)~\subset~\varprojlim^{mo}_{\{j\in \J~|~x\in \Sigma_X(j)\}}\mu\left(\Sigma_X(j)\right)s\circ f\left(\Sigma_X(j)\right).\nn
\end{align}}Given a topological $R$-algebra $A$ and a fuzzy integral {\small $I=I_{\mu,s}:\M\big(X,\P(A)\big)\ra\M(X,\P(A))$}, a \ul{fuzzy $I$-derivative} on $\M(X,\P(A))$ can similarly be defined as a map
\[
D=D_{\mu,s}:\G\subset\M(X,\P(A))\ra\M(X,\P(A))\nn
\]
satisfying (modified versions of) the required three defining conditions (i),(ii),(iii) in Definition \ref{IntDerDef}.
\end{dfn}

\section{Topological Analysis: Preliminary Tools}\label{GeomPrelimsCh}
\subsection{Closure, interior, and classical space types}\label{GeomPrelimsS1}
\begin{dfn}[\textcolor{blue}{
\index{Closure}{Closure} of a set,
\index{Dense! set}{Dense set},
\index{Separable space}{Separable space}}]
Let $X$ be a space. If $A\subset X$, the \ul{closure} $\ol{A}$ or $\Cl(A)$ of $A$ is the smallest closed set containing $A$, i.e., {\small $\ol{A}=\Cl(A):=\bigcap\{C\supset A:C~\txt{closed}\}$}. (\blue{footnote}\footnote{It is clear that $A$ is closed iff $A=\ol{A}$.}). A set $A\subset X$ is \ul{dense} if $\ol{A}=X$ (i.e., every \emph{nonempty} open set $U\subset X$ contains an element of $A$, since $(\ol{A})^c=\bigcup\{O:O~\txt{open},~O\cap A=\emptyset\}$). The space $X$ is \ul{separable} if it contains a countable dense set, i.e., there exists a countable set $A\subset X$ such that $\ol{A}=X$.
\end{dfn}

\begin{lmm}[\blue{Closure criterion}]\label{ClosureCrit}
Let $X$ be a space and $A\subset X$. The closure of $A$ satisfies
\bea
\ol{A}=\{x\in X:N(x)\cap A\neq\emptyset~\txt{for every nbd $N(x)$ of $x$}\}.\nn
\eea
\end{lmm}
\begin{proof}
Let $B:=\{x\in X:N(x)\cap A\neq\emptyset~\txt{for every nbd $N(x)$ of $x$}\}$. Then $B$ is closed because, $x\in B^c$ iff some $N(x)\cap A=\emptyset$, which in turn implies that for every $y\in N(x)$, some $N(y)\cap A=\emptyset$, which implies $N(x)\subset B^c$. It is also clear that $A\subset B$, and so $\ol{A}\subset\ol{B}=B$.

If $x\in B$, i.e., every $N(x)\cap A\neq\emptyset$, suppose $x\not\in\ol{A}$, i.e., there is a closed set $C\supset A$ with $x\not\in C$. Then $x\in C^c$ = some $N(x)$ and so $C^c\cap A\neq\emptyset$ (a contradiction since $A\subset C$). It follows that $B\subset \ol{A}$.
\end{proof}

\begin{dfn}[\textcolor{blue}{\index{Limit! point}{Limit points} of a set}]
Let $X$ be a space and $A\subset X$. A point $x\in X$ is a limit point of $A$, written $x\in A'$, if $\big(N(x)-\{x\}\big)\cap A\neq\emptyset$ for every nbd $N(x)$ of $x$. (\blue{footnote}\footnote{It follows from Lemma \ref{ClosureCrit} that $\ol{A}=A\cup A'$. Indeed it is clear that $A\cup A'\subset\ol{A}$. Let $x\in\ol{A}$, i.e., every $N(x)\cap A\neq\emptyset$. We have $N(x)\cap A=\big(\{x\}\cap A\big)\cup\big((N(x)-\{x\})\cap A\big)$. Either some $(N(x)-\{x\})\cap A=\emptyset$ (in which case $x\in A$) or every $(N(x)-\{x\})\cap A\neq\emptyset$ (in which case $x\in A'$), i.e., $\ol{A}\subset A\cup A'$.}).
\end{dfn}

\begin{lmm}\label{SecCountSep}
Every second countable space (i.e., space with a countable base) $X$ is separable.
\end{lmm}
\begin{proof}
Consider a countable base $\B:=\{O_n\}_{n\in\Natural}$ for $X$. Then any sequence $x_n\in O_n$ gives a countable dense set $\{x_n\}_{n\in\Natural}\subset X$. Indeed, by the base criterion, for any nonempty open set $U\subset X$, we have $O_{n_U}\subset U$ for some $n_U$, and so $x_{n_U}\in O_{n_U}\subset U$. That is, {\small $U\cap\{x_n\}_{n\in\Natural}\neq\emptyset$}, and so $\ol{\{x_n\}_{n\in\Natural}}=X$.
\end{proof}

\begin{dfn}[\textcolor{blue}{\index{Interior}{Interior} of a set}] Let $X$ be a space and $A\subset X$. The interior $A^o$ or $\Int(A)$ of $A$ is the largest open set contained in $A$, i.e., ~$A^o=\Int(A):=\bigcup\big\{O\subset A:O~\txt{open}\big\}$.  (\blue{footnote}\footnote{It is clear that $A$ is open iff $A=A^o$.})
\end{dfn}

{\flushleft\ul{Note}}: Using DeMorgan's laws, we easily see from the definitions that
\bea
\label{ComplEqs}\left(\ol{A}\right)^c=(A^c)^o~~~~\txt{and}~~~~(A^o)^c=\ol{A^c},
\eea
which simply mean that the complement of the closure (resp. the interior) is the interior (resp. a closure) of the complement.

\begin{lmm}\label{UnionClIntLmm}
Let $X$ be a space and $A,B\subset X$. Then (i) $\ol{A\cup B}=\ol{A}\cup\ol{B}$ and (ii) $(A\cap B)^o=A^o\cap B^o$.
\end{lmm}
\begin{proof}
(i) {\small$\ol{A\cup B}\subset\ol{\ol{A}\cup\ol{B}}=\ol{A}\cup\ol{B}$}. If $x\in X$, let $N(x)$ denote a nbd of $x$. Then $x\in \ol{A}\cup\ol{B}$ implies $x\in\ol{A}$ or $x\in\ol{B}$, which implies every $N(x)\cap A\neq\emptyset$ or every $N(x)\cap B\neq\emptyset$, which implies every {\small $N(x)\cap(A\cup B)\neq\emptyset$ (i.e., $\ol{A}\cup\ol{B}\subset\ol{A\cup B}$)}. (ii) Thus, {\small $(A\cap B)^o=\left[\ol{(A\cap B)^c}\right]^c=\left[^{}\ol{A^c\cup B^c}^{}\right]^c=\left[^{}\ol{A^c}\cup\ol{B^c}^{}\right]^c=A^o\cap B^o$}.
\end{proof}

\begin{crl}\label{UnionClIntCrl}
Let $X$ be a space and $B\subset A\subset X$. The \emph{subspace-closure} $\Cl_A(B)$ and \emph{subspace-interior} $\Int_A(B)$ of $B$ in $A$ (i.e., closure and interior of $B$ in the subspace topology of $A$) satisfy
\bea
\label{SubSpCl}&&\textstyle\Cl_A(B)=\bigcap\{C_1:~B\subset C_1,~\txt{$C_1\subset A$ closed in $A$}\}=\bigcap\{C\cap A:~B\subset C,~\txt{$C\subset X$ closed}\}\nn\\
&&~~~~=\Cl(B)\cap A,\\
\label{SubSpInt}&&\Int_A(B)\sr{(\ref{ComplEqs})}{=} \big[\Cl_A(B^{c_A})\big]^{c_A}=\big[\Cl(B^c\cap A)\cap A\big]^c\cap A=\big[\Cl(B^c\cap A)\big]^c\cap A\nn\\
&&~~~~=\Int(B\cup A^c)\cap A.
\eea
\end{crl}

\begin{dfn}[\textcolor{blue}{\index{Boundary}{Boundary} of a set, \index{Clopen set}{Clopen (or closed open) set}}] Let $X$ be a space. The boundary of a set $A\subset X$ is $\del A:=\ol{A}-A^o=\ol{A}\cap(A^o)^c\sr{(\ref{ComplEqs})}{=}\del A^c$. A set $A\subset X$ is clopen (equivalently, both closed and open) if $\del A=\emptyset$.
\end{dfn}

Note that $\del A^o,\del\ol{A}\subset\del A$ (where equalities hold if $A\subset X$ is a ball in a metric space $X$; see Definition \ref{MetricSp}), but the example $X=\Real$, $A=\Rational\subset\Real$ shows no equality holds in general. Since $\del A=\del A^c$, \ul{$A$ is clopen iff $A^c$ is clopen}. Another relation showing symmetry between $A,A^c$ is $\del A^o\cap\del\ol{A}\sr{(\ref{ComplEqs})}{=}\ol{A^o}\cap\ol{(A^c)^o}$.

\begin{lmm}\label{SubSpBdLmm}
Let $X$ be a space, $B\subset A\subset X$, and $B_1\subset X$. The \emph{subspace-boundary} $\del_A$ in $A$ (i.e., boundary in the subspace topology of $A$) satisfies (i) $\del_AB\subset (\del B)\cap A$ and (ii) $\del_A(B_1\cap A)\subset (\del B_1)\cap A$. (Equality holds in (i) if $A$ is closed in $X$ and $B$ is open in $X$.)
\end{lmm}
\begin{proof}
(i) By direct calculation, we have
\begin{align}
&\del_AB:=\Cl_A(B)-\Int_A(B)=[\ol{B}\cap A]\cap[(B\cup A^c)^o\cap A]^c=\ol{B}\cap A\cap\ol{B^c\cap A}\nn\\
&~~~~\subset \ol{B}\cap A\cap\ol{B^c}\cap\ol{A}=\ol{B}\cap(B^o)^c\cap A=(\del B)\cap A,\nn
\end{align}
where equality holds if $A$ is closed in $X$ and $B$ is open in $X$. (ii) Moreover, if $B:=B_1\cap A$, then
\begin{align}
&\del_A(B_1\cap A)=\del_AB=\ol{B}\cap A\cap\ol{B^c\cap A}=\ol{B_1\cap A}\cap A\cap\ol{B_1^c\cap A}\nn\\
&~~~~\subset\ol{B_1}\cap \ol{A}\cap A\cap\ol{B_1^c}=\ol{B_1}\cap A\cap\ol{B_1^c}=\ol{B_1}\cap A\cap(B_1^o)^c=(\del B_1)\cap A.\nn\qedhere
\end{align}
\end{proof}

For any $A,B\subset X$, we have $\del(A\cup B)=\ol{A\cup B}-(A\cup B)^o\subset \ol{A}\cup\ol{B}-(A^o\cup B^o)=(\del A\cap\ol{B^c})\cup(\ol{A^c}\cap\del B)$, where equality holds if $A,B$ are open. Similarly, $\del(A\cap B)=\del(A^c\cup B^c)\subset(\del A\cap\ol{B})\cup(\ol{A}\cap\del B)$, where equality holds if $A,B$ are closed.

\begin{lmm}\label{ProdClIntLmm}
Let $X,Y$ be spaces, $A\subset X$, and $B\subset Y$. In $X\times Y$, (i) $\ol{A\times B}=\ol{A}\times\ol{B}$ and (ii) $(A\times B)^o=A^o\times B^o$. Hence, the boundary of a product of sets $A\times B\subset X\times Y$ satisfies ~$\del(A\times B)=(\del A\times\ol{B})\cup(\ol{A}\times\del B)$. (\blue{footnote}\footnote{The complement and intersection of products of sets $A,C\subset X$ and $B,D\subset Y$ are given by
\bea
&&(A\times B)^c:=\{(x,y)\in X\times Y:(a,b)\not\in A\times B\}=\{(x,y)\in X\times Y:x\not\in A~\txt{or}~y\not\in B\}\nn\\
&&~~~~=\{(x,y)\in X\times Y:x\not\in A\}\cup\{(x,y)\in X\times Y:y\not\in B\}=(A^c\times Y)\cup(X\times B^c),\nn\\
&& (A\times B)\cap(C\times D)=\{(x,y)\in X\times Y:(x,y)\in A\times B,(x,y)\in C\times D\}=\{(x,y):x\in A,y\in B,x\in C,y\in D\}\nn\\
&&~~~~=\{(x,y)\in X\times Y:x\in A\cap C,~y\in B\cap D\}=(A\cap C)\times(B\cap D).\nn
\eea})
\end{lmm}
\begin{proof}
(i) $\ol{A\times B}\subset \ol{\ol{A}\times\ol{B}}=\ol{A}\times\ol{B}$ (because $\ol{A}\times\ol{B}$ is closed). Also, with $N(x)$ denoting a nbd of $x$,
\bea
&& (a,b)\in \ol{A}\times\ol{B}~~\Ra~~a\in\ol{A},~~b\in\ol{B},\nn\\
&&~~\Ra~~\txt{every}~~[N(a)\times N(b)]\cap[A\times B]=[N(a)\cap B]\times[N(b)\cap B]\neq\emptyset\nn\\
&&~~\Ra~~\txt{every}~~N(a,b)\cap[A\times B]\neq\emptyset,~~\Ra~~(a,b)\in\ol{A\times B}.\nn
\eea
(ii) It follows from (i) that $(A\times B)^o=A^o\times B^o$ in $X\times Y$, since
\begin{align}
&(A\times B)^o=\left[\ol{(A\times B)^c}\right]^c=\left[\ol{(A^c\times Y)\cup(X\times B^c)}\right]^c=\left[\ol{A^c\times Y}\cup\ol{X\times B^c}\right]^c\nn\\
&~~~~=\left(\ol{A^c}\times Y\right)^c\cap\left(X\times\ol{B^c}\right)^c=\left([\ol{A^c}]^c\times Y\right)\cap\left(X\times[\ol{B^c}]^c\right)\nn\\
&~~~~=(A^o\times Y)\cap(X\times B^o)=A^o\times B^o.\nn\qedhere
\end{align}
\end{proof}

\begin{dfns}[\textcolor{blue}{\index{Classical! space types}{Classical space types}}]\label{SpaceTypeDfn}
A space $X$ is called a
{\flushleft (a) \index{Kolmogorov space ($T_0$) space}{\ul{Kolmogorov space} (or \emph{$T_0$ space})}} if for any two distinct points $x,y\in X$, $x$ has a neighborhood that excludes $y$, or $y$ has a neighborhood that excludes $x$.
{\flushleft (b) \index{Fr$\acute{\txt{e}}$chet ($T_1$) space}{\ul{Fr$\acute{\txt{e}}$chet space}} (or \emph{$T_1$ space})} if for any two distinct points $x,y\in X$, $x$ has a neighborhood that excludes $y$, and $y$ also has a neighborhood that excludes $x$ (Equivalently, $X$ is a space in which every point is closed).
{\flushleft (c) \index{Hausdorff ($T_2$) space}{\ul{Hausdorff space}} (or \emph{$T_2$ space})} if  any two distinct points $x,y\subset X$ have disjoint neighborhoods.
{\flushleft (d) \index{Regular! space}{\ul{Regular space}}} if any closed set $C\subset X$ and any point $x\not\in C$ have disjoint neighborhoods (Equivalently, for every $x\in X$, every neighborhood $x\in U\subset X$ contains the closure of another neighborhood $x\in V\subset\ol{V}\subset U\subset X$). A regular Hausdorff space is called a \index{$T_3$ space}{\emph{$T_3$ space}}.

{\flushleft (e) \index{Normal space}{\ul{Normal space}}} if any two disjoint closed sets $C_1,C_2\subset X$, $C_1\cap C_2=\emptyset$, have disjoint neighborhoods. (Equivalently, for every closed set $C\subset X$, every neighborhood $C\subset U\subset X$ contains the closure of another neighborhood $C\subset V\subset\ol{V}\subset U\subset X$). A normal Hausdorff space (i.e., a normal $T_1$ space) is called a \index{$T_4$ space}{\emph{$T_4$ space}}.
\end{dfns}

\subsection{Metric (metrizable) spaces and quotient pseudometric spaces}\label{GeomPrelimsS2}
\begin{dfn}[\textcolor{blue}{
\index{Metric!}{Metric},
\index{Ball}{Ball},
\index{Metric! topology}{{Metric topology}},
\index{Metric! space}{Metric space},
\index{Metrizable space}{Metrizable space},
\index{Pseudometric}{Pseudometric},
\index{Pseudometric space}{Pseudometric space}}]\label{MetricSp}
Let $X$ be a set. A \ul{metric} on $X$ is a function $d:X\times X\ra\Real$ such that for all $x,y,z\in X$,
\bit
\item[] $d(x,y)=d(y,x)\leq d(x,z)+d(z,y)$, and $d(x,y)=0$ iff $x=y$.
\eit
Given a metric $d$ on $X$, the \ul{$d$-ball} of radius $R>0$ centered at a point $x$ in $X$ is the set
\[
B_R(x)=B_R^d(x):=\{x'\in X:d(x,x')<R\}.
\]
Given a metric $d$ on $X$, the \ul{$d$-metric topology} (or just \ul{$d$-topology}) on $X$ (making $X=(X,d)$ a \ul{metric space}) is the topology on $X$ with open sets defined as follows: A (nonempty) set $A\subset X$ is \ul{$d$-open} if for every $x\in A$, there exists a $d$-ball $B_R(x)\subset A$. (As shown in Remark \ref{DtopProofRmk} below, the $d$-topology is indeed a topology in the usual sense, with the $d$-balls as $d$-open sets that form a base for the topology.)

A space $X$ is \ul{metrizable} if it is homeomorphic to a metric space, i.e., if there exists a homeomorphism $h:X\ra Y$, for some metric space $Y=(Y,d)$. A function $d:X\times X\ra\Real$ satisfying $d(x,y)=d(y,x)\leq d(x,z)+d(z,y)$, for all $x,y,z\in X$, is called a \ul{pseudometric} (making $X=(X,d)$ a \ul{pseudometric space} with respect to the $d$-topology above).
\end{dfn}

\begin{rmk}[\blue{A metric space is a topological space}]\label{DtopProofRmk}
Let $X=(X,d)$ be a metric space. Observe that any $d$-ball $B_R(x)\subset X$ is open because for each $y\in B_R(x)$, we have $B_{R-d(x,y)}(y)\subset B_R(x)$. It follows that a set $A\subset X$ is $d$-open $\iff$ it is a union of $d$-balls. Also, the intersection $B_R(x)\cap B_{R'}(x')$ of two $d$-balls is open because for each $y\in B_R(x)\cap B_{R'}(x')$, if ~$R(x,x'):=\min\big(R-d(x,y),R'-d(x',y)\big)$~ then
\[
B_{R(x,x')}(y)\subset B_R(x)\cap B_{R'}(x').
\]
It follows that finite intersections of $d$-open sets are $d$-open sets, since
\[
\textstyle \big(\bigcup_\al B_{R_\al}(x_\al)\big)\cap\big(\bigcup_\beta B_{\beta}(x_\beta)\big)=\bigcup_{\al,\beta}~B_{R_\al}(x_\al)\cap B_{R_\beta}(x_\beta).
\]
Hence, $d$-open sets give a topology on $X$, with $d$-balls forming a base for the topology.
\end{rmk}

\begin{convention}[\textcolor{blue}{
\index{Product metric space}{Product metric space},
\index{Continuity of metrics}{Continuity of metrics}}]
Let $X,Y$ be metric spaces. As a metric space, unless stated otherwise, the product space $X\times Y$ will be given the metric {\footnotesize $d_{X\times Y}\big((x,y),(x',y')\big):=\max\big(d_X(x,x'),d_Y(y,y')\big)$}. (\blue{footnote}\footnote{If $X=(X,d_X)$ and $Y=(Y,d_Y)$ are metric spaces, then the product space $X\times Y$ is metrizable, in the sense that the identity map $id_{X\times Y}:X\times Y\ra\big(X\times Y,d_{X\times Y}\big)$ is a homeomorphism due to ~$B_R^{d_{X\times Y}}\big((x,y)\big)=B_R^{d_X}(x)\times B_R^{d_Y}(y)$.}). Thus, any metric $d_X:X\times X\ra\Real$ on $X$ is a continuous function (see Corollary \ref{ContnyCrit3}, page \pageref{ContnyCrit3}), since for all $(x,y),(x',y')\in X\times X$, we have
\bea
|d_X(x,y)-d_X(x',y')|\leq d_X(x,x')+d_X(y,y')\leq 2d_{X\times X}\big((x,y),(x',y')\big).\nn
\eea
\end{convention}

\begin{dfn}[\textcolor{blue}{
\index{Diameter of a set}{Diameter of a set},
\index{Distance between sets}{Distance between sets}}]
Let $X$ be a metric space and $A,B\subset X$. The \ul{diameter} of $A$ is ~$\diam A:=\sup\{d(a,a'):a,a'\in A\}$. ~The \ul{distance} between $A$ and $B$ is ~$\dist(A,B):=\inf\{d(a,b):a\in A,b\in B\}$. (\blue{footnote}\footnote{Using the infimum criterion, it is easy to see that for any $x\in X$, $A,B\subset X$, (i) $\dist(x,A)=0$ $\iff$ $x\in\ol{A}$, and so (ii) if $A\cap\ol{B}\neq\emptyset$ or $\ol{A}\cap B\neq\emptyset$ then $\dist(A,B)=0$ (but not the reverse as can be seen using $A:=\{1/n\}$, $B:=[-1,0)$ in $(\Real,|\cdot|)$).}).
\end{dfn}

\begin{dfn}[\textcolor{blue}{
\index{Quotient! pseudometric}{Quotient pseudometric},
\index{Quotient! pseudometric topology}{{Quotient pseudometric topology}},
\index{Quotient! pseudometric space}{Quotient pseudometric space}}]\label{QuPseudoMetSp}
Let $(X,d)$ be a metric space and $\sim$ an equivalence relation on $X$. The \ul{quotient pseudometric} $\rho$ on the quotient space ${X\over\sim}=\{[x]=x_\sim:x\in X\}$, making $\big({X\over\sim},\rho\big)$ a \ul{quotient pseudometric space}, is
\bea
&&\rho([x],[y]):=\inf\Big\{l_\delta(c)~\big|~c:\{0,1,...,n\}\ra X,~c_0\in[x],~c_n\in[y],~n\in\Natural\Big\}\leq\dist([x],[y]),\nn\\
&&\textstyle l_\delta(c):=\sum\limits_{i=0}^n\delta(c_{i-1},c_i),
~~~~\delta(u,v):=\left\{
                     \begin{array}{ll}
                       0, & u\sim v  \\
                       d(u,v), & \txt{otherwise}
                     \end{array}
                   \right\},\nn
\eea
which is the smallest $\delta$-length $l_\delta(c)$ of finite chains $c=\{c_0,c_1,...,c_n\}$ of points in $X$ with endpoints in the equivalence classes $[x]$ and $[y]$, i.e., $c_0\in [x]$, $c_n\in[y]$.
\end{dfn}
It is clear that $\rho([x],[y])=\rho([y],[x])$. Also, $\rho([x],[y])\leq\rho([x],[z])+\rho([z],[y])$, since for any $\vep>0$, there exist two finite chains $c_{[x],[z]}$, $c_{[z],[y]}$ such that, with the concatenation $c_{[x],[y]}:=c_{[x],[z]}\cdot c_{[z],[y]}$ (another such finite chain by the transitivity of $\sim$), we have
\bea
\rho([x],[y])\leq l_\delta\left(c_{[x],[y]}\right)\leq l_\delta\left(c_{[x],[z]}\right)+l_\delta\left(c_{[z],[y]}\right)<\rho([x],[z])+\rho([z],[y])+2\vep.\nn
\eea
For the special case where ${X\over\sim}={X\over A}$ for some $A\subset X$, observe that for any $x,y\in X$, we have both $\rho([x],[y])\leq \dist([x],[y])$ and (with any $a\in A$) $\rho([x],[y])\leq \rho([x],[a])+\rho([a],[y])=\rho([x],A)+\rho(A,[y])\leq \dist([x],A)+\dist(A,[y])$. That is,
\[
\rho([x],[y])\leq\min\big\{\dist([x],[y]),\dist([x],A)+\dist(A,[y])\big\}.
\]
Equality holds as the following proposition shows.
\begin{prp}\label{PseudoMetPrp}
Let $X$ be a space and $A\subset X$. The quotient pseudometric on ${X\over\sim}:={X\over A}$ satisfies
\bea
\rho([x],[y]) = \min\big\{d(x,y)~,~\dist(x,A)+\dist(y,A)\big\},~~~~\txt{for all}~~x,y\in X.\nn
\eea
\end{prp}
\begin{proof}
Fix $x,y\in X$. We have three exhaustive special cases as follows:
\bit[leftmargin=0.0cm]
\item[]\ul{Case 1 ($x,y\in A$)}:~~ $\rho([x],[y])=\rho(A,A)=0$.
\item[]\ul{Case 2 ($x\not\in A,y\in A$)}:~~ $\rho([x],[y])=\rho(\{x\},A)=\dist(x,A)$.
\item[]\ul{Case 3 ($x,y\not\in A$)}:~~ $\rho([x],[y])=\rho(\{x\},\{y\})\leq\min\big\{d(x,y),\dist(x,A)+\dist(A,y)\big\}$
\eit
For case 3 above, observe that any finite chain $c_{x,y}:=c_{\{x\},\{y\}}=c_{[x],[y]}$ from $x$ to $y$ has length
\[
l_\delta(c_{x,y})\geq\left\{
                       \begin{array}{ll}
                          d(x,y), & \txt{if ~$c_{x,y}\cap A=\emptyset$} \\
                         \dist(x,A)+\dist(A,y), & \txt{if ~$c_{x,y}\cap A\neq\emptyset$}
                       \end{array}
                     \right\}\geq \min\big\{d(x,y),\dist(x,A)+\dist(A,y)\big\},
\]
and so ~$\rho(\{x\},\{y\})\geq\min\big\{d(x,y),\dist(x,A)+\dist(A,y)\big\}$.
\end{proof}

\begin{dfn}[\blue{Subset visitation path (SV-path), Length of a SV-path}]
Let $X=(X,d)$ be a metric space and $n\geq 0$ an integer. A \ul{subset visitation path} (\ul{SV-path}) of order $n$ is a map of the form $p:\{0,1,2,\cdots,n\}\ra\P(X)$, also written as a sequence $p=(p_0,p_1,...,p_n)$. The \ul{length} of $p$ is
\[
\textstyle l(p)=l(p_0,\cdots,p_n):=\sum_{i=1}^n\dist(p_{i-1},p_i)=\dist(p_0,p_1)+\dist(p_1,p_2)+\cdots+\dist(p_{n-1},p_n).
\]
\end{dfn}

\begin{crl}\label{PseudoMetCrl1}
Let $X=(X,d)$ be a metric space and $A,B\subset X$ disjoint subsets, i.e., $A\cap B=\emptyset$. Let $\sim$ be the equivalence relation on $X$ given by ``$x\sim y$ if $x=y$ or $x,y\in A$ or $x,y\in B$''. The quotient pseudometric on $X[A,B]:={X\over\sim}$ satisfies the following: For all $x,y\in X$,
\[
\rho([x],[y]) = \min\big\{l(x,y),l(x,A,y),l(x,B,y),l(x,A,B,y),l(x,B,A,y)\big\}.\nn
\]
where $l(p)$ is the length of the SV-path $p$ in $\big\{\{x\},\{y\},A,B\big\}\subset\P(X)$, and $x,y$ represent the singleton sets $\{x\},\{y\}$ respectively.
\end{crl}
\begin{proof}
This is a straightforward generalization of the proof of Proposition \ref{PseudoMetPrp}.
\end{proof}

\begin{crl}\label{PseudoMetCrl2}
Let $X=(X,d)$ be a metric space and $\sim$ an equivalence relation on $X$. The quotient pseudometric on ${X\over\sim}$ satisfies the following: For all $x,y\in X$,
\[
\textstyle\rho([x],[y]):=\inf\Big\{l(p)~\big|~p:\{0,1,...,n\}\ra {X\over\sim},~p_0=[x],~p_n\in[y],~n\in\Natural\Big\},
\]
which is the smallest SV-path length $l(p)$ of SV-paths $p=\{p_0,p_1,...,p_n\}$ in ${X\over\sim}\subset\P(X)$ from $[x]$ to $[y]$.
\end{crl}
\begin{proof}
Consider the obvious extension of Corollary \ref{PseudoMetCrl1} by induction (i.e., instead of two disjoint subsets, consider an arbitrary number of disjoint subsets $A_1,A_2,\cdots,A_t\subset X$ and the pseudometric on the associated quotient space $X[A_1,A_2,\cdots,A_t]$). Then merge its proof with the original definition of the pseudometric from Definition \ref{QuPseudoMetSp}.
\end{proof}

\subsection{Pointwise characterization of continuity and openness of maps}\label{ContOpenSec}~\\~
For our discussion in this section and subsequently, if $X=(X,d)$ is a metric space and $x\in X$, the phrase ``\ul{near $x$}'' will mean ``\ul{in a sufficiently small neighborhood of $x$}''.
\begin{dfn}[\textcolor{blue}{\index{Continuous map}{Recall: Continuous map}, Continuity at a point}]\label{PtwsContDfn} Let $X,Y$ be spaces. A map $f:X\ra Y$ is continuous, also written $f\in \C(X,Y):=\{\txt{continuous}~f:X\ra Y\}$, if for any open set $V\subset Y$, the preimage $f^{-1}(V)\subset X$ is open. That is, $f:(X,\T_X)\ra (Y,\T_Y)$ is continuous if $f^{-1}(\T_Y)\subset\T_X$. Equivalently, a map of spaces is continuous if the preimage of each closed set is closed.

A map $f:X\ra Y$ is continuous at $x\in X$ if for every open neighborhood $V$ of $f(x)$, the preimage $f^{-1}(V)$ is an open neighborhood of $x$ (equivalently, for any open set $V\ni f(x)$, there exists an open set $U\ni x$ such that $f(U)\subset V$).
\end{dfn}

\begin{lmm}
Let $X,Y$ be spaces. A map $f:X\ra Y$ is continuous $\iff$ every set $A\subset X$ satisfies {\footnotesize $\ol{A}\subset f^{-1}\left(\ol{f(A)}\right)$}, i.e.,{\footnotesize $f(\ol{A})\subset\ol{f(A)}$}. Equivalently, either directly or, by setting $A:=f^{-1}(B^c)$ in {\footnotesize $\ol{A}\subset f^{-1}\left(\ol{f(A)}\right)$} and applying the result (\ref{ComplEqs}), we see that $f:X\ra Y$ is continuous $\iff$ every set $B\subset Y$ satisfies $f^{-1}(B^o)\subset [f^{-1}(B)]^o$.
\end{lmm}
\begin{proof} If $f$ is continuous then {\footnotesize $f^{-1}\left(\ol{f(A)}\right)$} is closed, and so {\footnotesize $A\subset f^{-1}(f(A))\subset f^{-1}\left(\ol{f(A)}\right)$} implies {\footnotesize $\ol{A}\subset f^{-1}\left(\ol{f(A)}\right)$}. Conversely, if {\footnotesize $\ol{A}\subset f^{-1}\left(\ol{f(A)}\right)$} for any $A\subset X$, then for an open set $V\subset Y$, the preimage $A:=f^{-1}(V^c)$ is closed because it satisfies {\footnotesize $\ol{f^{-1}(V^c)}\subset f^{-1}\left(\ol{f(f^{-1}(V^c))}\right)\subset f^{-1}\left(\ol{V^c}\right)=f^{-1}(V^c)\subset \ol{f^{-1}(V^c)}$}.
\end{proof}

\begin{lmm}[\blue{Continuity criterion}]\label{ContnyCrit1}
Let $X,Y$ be spaces. A map $f:X\ra Y$ is continuous $\iff$ continuous at each point $x\in X$.
\end{lmm}
\begin{proof}
If $f$ is continuous then it is clear that for every open set $V\ni f(x)$, $f^{-1}(V)$ is an open set containing $x$. Conversely, if $f$ is continuous at every point $x\in X$, then for any open set $V\subset Y$ (with $V_y\subset V$ denoting an open set containing $y$) we have the open set
\[
\textstyle f^{-1}(V)=f^{-1}\big(\bigcup_{y\in V}V_y\big)=\bigcup_{y\in V}f^{-1}(V_y). \qedhere
\]
\end{proof}

\begin{crl}[\blue{Metric space continuity}]\label{ContnyCrit2}
Let $X,Y$ be metric spaces. A map $f:X\ra Y$ is continuous $\iff$ for any $\vep>0$ and $x\in X$, there exists $\delta_x(\vep)>0$ such that $f\big(B_{\delta_x(\vep)}(x)\big)\subset B_\vep\big(f(x)\big)$, i.e., for every $x'\in X$, $d(x,x')<\delta_x(\vep)$ implies $d(f(x),f(x'))<\vep$.
\end{crl}
Observe that $f\big(B_{\delta_x(\vep)}(x)\big)\subset B_\vep\big(f(x)\big)$ implies $f\big(B_\delta(x)\big)\subset f\big(B_{\delta_x(\vep)}(x)\big)\subset B_\vep\big(f(x)\big)$ for all $0<\delta<\delta_x(\vep)$, and so by replacing $\delta_x(\vep)$ with a sufficiently small such $\delta$, we can assume wlog that $\delta_x(\vep)\leq\vep$ (or more simply, that $\delta_x(\vep)\ra 0$ as $\vep\ra 0$). This proves the following, which is simply one direction of Proposition \ref{NetContPrp} on page \pageref{NetContPrp}:
\begin{crl}\label{ContnyCrit3n}
Let $X,Y$ be metric spaces and $f:X\ra Y$ a continuous map. If $x\ra x'$ in $X$, then $f(x)\ra f(x')$ in $Y$ (that is, if $d(x,x')\ra 0$, then $d(f(x),f(x'))\ra 0$). (\blue{footnote}\footnote{Note that if $f:X\ra Y$ is continuous, $f(x)\ra f(x')$ does not imply $x\ra x'$ in general, since a constant map is continuous.})
\end{crl}

\begin{crl}\label{ContnyCrit3}
Let $X,Y$ be metric spaces. A map $f:X\ra Y$ is continuous $\iff$ for each $x\in X$ there exists a function $\omega_x:[0,\infty)\ra[0,\infty)$ such that (i) $\omega_x(t)$ is nondecreasing near $t=0$, (ii) $\lim_{t\ra 0}\omega_x(t)=\omega_x(0)=0$, and (iii) $d(f(x),f(x'))\leq\omega_x\big(d(x,x')\big)$ for all $x'\in X$.
\end{crl}
\begin{proof}
($\Ra$): Assume $f$ is continuous, i.e., for any $\vep>0$ and $x\in X$, there exists $\delta_x(\vep)>0$ (where $\delta_x(\vep)\leq\vep$ wlog) such that $d(x,x')<\delta_x(\vep)$ implies $d(f(x),f(x'))<\vep$, for all $x'\in X$. Define the function
\[
\omega_x(t):=\sup\left\{d(f(x),f(x')):d(x,x')\leq t,x'\in X\right\},
\]
which is nondecreasing, $\omega_x(0)=0$, and $d(f(x),f(x'))\leq\omega_x(d(x,x'))$ for all $x'\in X$. Then
\[
d(x,x')<\delta_x(\vep)~~\Ra~~d(f(x),f(x'))\leq\omega_x(d(x,x'))\leq \omega_x(\delta_x(\vep))\leq\vep.
\]
Thus by taking $\vep\ra 0$, and noting $\delta_x(\vep)\leq\vep$, we see that $\omega_x(t)\ra 0$ as $t\ra 0$.

($\La$): Conversely, assume there exists a function $\omega_x:[0,\infty)\ra[0,\infty)$ such that (i) $\omega_x(t)$ is nondecreasing near $t=0$, (ii) $\lim_{t\ra 0}\omega_x(t)=\omega_x(0)=0$, and (iii) $d(f(x),f(x'))\leq\omega_x\big(d(x,x')\big)$ for all $x'\in X$. Then for any $\vep>0$, because $\omega_x(t)\ra 0$ as $t\ra 0$, we can choose a $\delta=\delta_x(\vep)>0$ such that
\[
d(x,x')<\delta~~\Ra~~d(f(x),f(x'))\leq\omega_x\big(d(x,x')\big)\leq \omega_x\big(\delta\big)<\vep. \qedhere
\]
\end{proof}

\begin{dfn}[\textcolor{blue}{
\index{Uniform! continuity}{Uniform continuity},
\index{Modulated uniform continuity}{Modulated uniform continuity},
\index{Lipschitz continuity}{Lipschitz continuity},
\index{H\"older continuity}{H\"older continuity}}]
Let $X,Y$ be metric spaces and $f:X\ra Y$ a continuous map. If the function $\delta_x$ in Corollary \ref{ContnyCrit2} is independent of $x$, we say $f$ is \ul{uniformly continuous}. If (in some neighborhood $[0,\al)$ of $0$ in $[0,\infty)$) the function $\omega_x$ in Corollary \ref{ContnyCrit3} is independent of $x$, we say $f$ is \ul{$\omega$-uniformly continuous} (or \ul{$\omega$-modulated uniformly continuous}). That is, $f$ is \ul{$\omega$-uniformly continuous} if there exists a function $\omega:[0,\infty)\ra[0,\infty)$ such that (i) $\omega(t)$ is nondecreasing near $t=0$, (ii) $\lim_{t\ra 0}\omega(t)=\omega(0)=0$, and (iii) $d(f(x),f(x'))\leq\omega\big(d(x,x')\big)$ for all $x,x'\in X$. (\blue{footnote}\footnote{Equivalently, a map of metric spaces $f:X\ra Y$ is \ul{uniformly continuous} iff for any $\vep>0$, there exists $\delta(\vep)>0$ such that $f\big(B_{\delta(\vep)}(x)\big)\subset B_\vep\big(f(x)\big)$, i.e., for every $x,x'\in X$, $d(x,x')<\delta(\vep)$ implies $d(f(x),f(x'))<\vep$.}).

If $\omega(t)=ct$ for some constant $c$, we say $f$ is \ul{$c$-Lipschitz continuous}. If $\omega(t)=ct^\al$ for some constants $c$ and $0<\al<1$, we say $f$ is \ul{$(c,\al)$-H\"older continuous}. That is, $f$ is $c$-Lipschitz continuous if $d(f(x),f(x'))\leq cd(x,x')$ for all $x,x'\in X$, and $f$ is $(c,\al)$-H\"older continuous if $d(f(x),f(x'))\leq cd(x,x')^\al$ for all $x,x'\in X$.
\end{dfn}

As indicated in \cite{heinonen2}, if $X,Y$ are metric spaces, a Lipschitz map $f:X\ra Y$ can be seen as a generalization of a differentiable real-valued function $f:X\ra\Real$.

\begin{dfn}[\textcolor{blue}{Recall:
\index{Open! map}{Open map},
\index{Closed map}{Closed map},
\index{Openness at a point}{Openness at a point}}]
Let $X,Y$ be spaces. A map $f:X\ra Y$ is \ul{open} if for each open set $O\subset X$, the image $f(O)\subset Y$ is an open set. Similarly, a map $f:X\ra Y$ is \ul{closed} if for each closed set $C\subset X$, the image $f(C)\subset Y$ is closed in $Y$.

A map $f:X\ra Y$ is open at $x\in X$ if for any open neighborhood $U$ of $x$, $f(U)$ is an open neighborhood of $f(x)$ (equivalently, for any open set $U\ni x$, there exists an open set $V\ni f(x)$ such that $V\subset f(U)$).
\end{dfn}

\begin{lmm}[\blue{Open map criterion}]
Let $X,Y$ be spaces. A map $f:X\ra Y$ is open $\iff$ open at each point $x\in X$.
\end{lmm}
\begin{proof}
If $f$ is open, it is clear that $f$ is open at each $x\in X$. Conversely, if $f$ is open at every point $x\in X$, then for any open set $U\subset X$ (with $U_x\subset U$ denoting an open set containing $x$) we have the open set ~$f(U)=f\left(\bigcup_{x\in U}U_x\right)=\bigcup_{x\in U}f(U_x)$.
\end{proof}

\begin{crl}
Let $X,Y$ be metric spaces. A map $f:X\ra Y$ is open $\iff$ for any $\vep>0$ and $x\in X$, there exists $\delta_x(\vep)>0$ such that $B_{\delta_x(\vep)}\big(f(x)\big)\subset f\big(B_\vep(x)\big)$, i.e., for every $x'\in X$ there exists $x''\in f^{-1}(f(x'))$ such that $d\big(f(x),f(x')\big)<\delta_x(\vep)$ implies $d(x,x'')<\vep$.
\end{crl}

\begin{crl}\label{OpnssCrit3}
Let $X,Y$ be metric spaces. A map $f:X\ra Y$ is open $\iff$ for each $x\in X$ there exists a function $\omega_x:[0,\infty)\ra[0,\infty)$ such that (i) $\omega_x(t)$ is nondecreasing near $t=0$, (ii) $\lim_{\vep\ra 0}\omega_x(\vep)=\omega_x(0)=0$, and (iii) for any $x'\in X$ there exists $x''\in f^{-1}(f(x'))$ such that $d(x,x'')\leq \omega_x\Big(d\big(f(x),f(x')\big)\Big)$. In particular,
\[
\dist\Big(x,f^{-1}\big(f(x')\big)\Big)\leq \omega_x\Big(d\big(f(x),f(x')\big)\Big),~~~~\txt{for all}~~~~x'\in X.
\]
\end{crl}
\begin{proof}
($\Ra$): Assume $f$ is open, i.e., for any $\vep>0$ and $x\in X$, there exists $\delta_x(\vep)>0$ (where $\delta_x(\vep)\leq\vep$ wlog) such that for every $x'\in X$ there is $x_\vep'\in f^{-1}(f(x'))$ such that $d(f(x),f(x'))<\delta_x(\vep)$ implies $d(x,x_\vep')<\vep$. Define the function
\[
\omega_x(t):=\left\{
               \begin{array}{ll}
                 0, & t=0 \\
                 \sup\left\{d(x,x_t'):d(f(x),f(x'))\leq t,~x'\in X\right\}, & t>0
               \end{array}
             \right\},
\]
which is nondecreasing, $\omega_x(0)=0$, and for all $x'\in X$ we have $x'':=x'_{d(f(x),f(x'))}\in f^{-1}(f(x'))$ such that $d(x,x'')\leq\omega\big(d(f(x),f(x'))\big)$. Then for all $x'\in X$, there is $x''\in f^{-1}(f(x'))$ such that
\[
d(f(x),f(x'))<\delta_x(\vep)~~\Ra~~d(x,x'')\leq \omega_x(d(f(x),f(x')))\leq \omega_x(\delta_x(\vep))\leq\vep.
\]
Thus by taking $\vep\ra 0$, and noting $\delta_x(\vep)\leq\vep$, we see that $\omega_x(t)\ra 0$ as $t\ra 0$.

($\La$): Conversely, assume there exists a function $\omega_x:[0,\infty)\ra[0,\infty)$ such that (i) $\omega_x(t)$ is nondecreasing near $t=0$, (ii) $\lim_{t\ra 0}\omega_x(t)=\omega_x(0)=0$, and (iii) for any $x'\in X$ there exists $x''\in f^{-1}(f(x'))$ such that $d(x,x'')\leq \omega_x\big(d(f(x),f(x'))\big)$. Then, given $x\in X$ and $\vep>0$, because $\omega_x(t)\ra 0$ as $t\ra 0$, we can choose a $\delta=\delta_x(\vep)>0$ such that for any $x'\in X$ there is $x''\in f^{-1}(f(x'))$ satisfying
\[
d(f(x),f(x'))<\delta~~\Ra~~d(x,x'')\leq\omega_x\big(d(f(x),f(x'))\big)\leq \omega_x\big(\delta\big)<\vep. \qedhere
\]
\end{proof}

\begin{crl}\label{OpnssCrit4}
Let $X,Y$ be metric spaces. A \ul{continuous} map $f:X\ra Y$ is open $\iff$ for each $x\in X$ there exists a function $\omega_x:[0,\infty)\ra[0,\infty)$ such that (i) $\omega_x(t)$ is nondecreasing near $t=0$, (ii) $\lim_{\vep\ra 0}\omega_x(\vep)=\omega_x(0)=0$, and (iii) {\small$\dist\Big(x,f^{-1}\big(f(x')\big)\Big)\leq \omega_x\Big(d\big(f(x),f(x')\big)\Big)$}~ for all $x'\in X$.
(This is immediate from Corollary \ref{OpnssCrit3}, because $f^{-1}(f(x'))$ is a closed set.)
\end{crl}

\begin{lmm}[\blue{Continuity on compact sets}]
Let $X,Y$ be metric spaces. A continuous map $f:X\ra Y$ is uniformly continuous on compact subsets of $X$.
\end{lmm}
\begin{proof}
It suffices to assume $X$ is compact (since restrictions of $f$ are continuous). Fix $\vep>0$. Then for each $x\in X$, there is $\delta_x=\delta_x(\vep)$ such that $d(x,x')<\delta_x$ implies $d(f(x),f(x'))<\vep/2$. Since $X$ is compact, let $\{B_{\delta_{x_i}/2}(x_i)\}_{i=1}^n$ cover $X$. Define $\delta:=\min_i\delta_{x_i}/2$. Pick any $u,v\in X$. Let $u\in B_{\delta_{x_j}/2}(x_j)$ for some $j$. Then $d(u,v)<\delta$ implies $u,v\in B_{\delta_{x_j}}(x_j)$, which in turn implies
\[
d(f(u),f(v))\leq d(f(u),f(x_j))+d(f(x_j),f(v))<\vep/2+\vep/2=\vep.\nn\qedhere
\]
\end{proof}

\begin{crl}
Let $X,Y$ be metric spaces. If $X$ is \ul{locally compact} (in the sense that every point of $X$ has a neighborhood whose closure is compact), then every continuous map $f:X\ra Y$ is \ul{locally uniformly continuous} (in the sense that every point of $X$ has a neighborhood in which $f$ is uniformly continuous).
\end{crl}

\begin{dfn}[\textcolor{blue}{Recall: \index{Homeomorphism}{Homeomorphism}, \index{Imbedding}{Imbedding}}]\label{HomeoImbed}
Let $X,Y$ be spaces. A map $h:X\ra Y$ is a \ul{homeomorphism} if it is bijective and both $h,h^{-1}$ are continuous (in which case, we say $X,Y$ are homeomorphic, or $X\cong Y$). A map $f:X\ra Y$ is an \ul{imbedding} (of $X$ into $Y$) if it is a homeomorphism onto its image (i.e., $f:X\ra f(X)$ is a homeomorphism).
\end{dfn}

\begin{dfn}[\textcolor{blue}{
\index{BiLipschitz map}{BiLipschitz map},
\index{Isometric map}{Isometric map},
\index{BiLipschitz equivalence}{BiLipschitz equivalence},
\index{Isometry}{Isometry}}]
Let $X=(X,d_X)$ and $Y=(Y,d_Y)$ be metric space and $c>0$ a positive constant. A map $f:X\ra Y$ is \ul{$c$-biLipschitz} (or a \ul{$c$-biLipschitz imbedding}) if $d_X(x,x')/c$ $\leq$ $d_Y\big(f(x),f(x')\big)$ $\leq$ $cd_X(x,x')$ for all $x,x'\in X$. A \ul{biLipschitz equivalence} is a surjective biLipschitz map. An \ul{isometric map} is a $1$-biLipschitz map (i.e., a map $f:X\ra Y$ such that $d_Y\big(f(x),f(x')\big)=d_X(x,x')$ for all $x,x'\in X$).  An \ul{isometry} is a surjective isometric map (i.e., an isometric homeomorphism).
\end{dfn}

\begin{lmm}[\blue{Bijection as homeomorphism}]\label{BijVsHomeo0}
Let $X,Y$ be spaces and $f:X\ra Y$ a \ul{bijective continuous} map. Then $f$ is a homeomorphism ($\iff$ open) $\iff$ closed.
\end{lmm}
\begin{proof}
For open $U\subset X$, $f(U^c)=f(U)^c$ is closed $\iff$ $f(U)$ is open $\iff$ $f^{-1}$ is continuous.
\end{proof}
Therefore, for any quotient map $q:X\ra Y$, if $q$ is injective on a closed (resp. open) set $A\subset X$, then $q|_A:A\ra Y$ is an imbedding. This is because, if $B\subset A$ is closed (resp. open), then $B=q^{-1}\big(q(B)\big)$ implies $q(B)\subset q(A)$ is closed (resp. open).

\begin{facts}
(i) A closed subset of a compact space is compact. (ii) A continuous map takes compact sets to compact sets. (iii) A compact subset of a Hausdorff space is closed.
\end{facts}
\begin{proof}
Let $X,Y$ be spaces. (i) Assume $X$ is compact. Let $C\subset X$ be closed. If $\{U_\al\}_{\al\in A}$ is an open cover of $C$, then $\{C^c\}\cup\{U_\al\}_{\al\in A}$ is an open cover of $X$, and so has a finite subscover $\{C^c\}\cup\{U_{\al_i}\}_{i=1}^n$. Hence $\{U_{\al_i}\}_{i=1}^n$ is a finite subcover of $C$.

(ii) Let $f:X\ra Y$ be a continuous map, and $K\subset X$ a compact set. If $\{V_\al\subset Y\}_{\al\in A}$ is an open cover of $f(K)$, then $\{f^{-1}(V_\al)_{\al\in A}\}$, is an open cover of $K$, containing a finite subcover $\{f^{-1}(V_{\al_i})\}_{i=1}^n$. Thus $\{V_{\al_i}\}_{i=1}^n$ is a finite subcover of $f(K)$, i.e., $f(K)$ is compact.

(iii) Assume $X$ is Hausdorff. Let $K\subset X$ be compact. Pick $x\in K^c$. Then for every $k\in K$, there exist disjoint open sets $U_k\ni x$ and $V_k\ni k$. Since $K$ is compact, its open cover $\{V_k\}_{k\in K}$ has a finite subcover $\{V_{k_i}\}_{i=1}^m$. Let $U:=\bigcap_i U_{k_i}$ and $V:=\bigcup_i V_{k_i}$. Then $x\in U\subset K^c$, since $K\subset V$ and $U\cap V=\emptyset$. Hence $K^c$ is open, i.e., $K$ is closed.
\end{proof}

\begin{crl}[\blue{Closed image}]\label{ClosedImLmm}
Let $f :X \ra Y$ be a continuous map. If $X$ is compact and $Y$ is Hausdorff, then $f(X)$ is both compact and closed.
\end{crl}

\begin{lmm}[\blue{Bijection as homeomorphism}]\label{BijVsHomeo}
Let $f:X\ra Y$ be a continuous map. If $X$ is compact and $Y$ is Hausdorff, then $f$ is a homeomorphism $\iff$ a bijection.
\end{lmm}
\begin{proof}
If $f$ is a homeomorphism, it is clear that $f$ is a bijection. Conversely, assume $f$ is a bijection. To show $f^{-1}:Y\ra X$ is continuous (i.e., that $f$ is open), let $O\subset X$ be open. Then the image ~$f(O)=(f^{-1})^{-1}\big((O^c)^c\big)=\big[(f^{-1})^{-1}(O^c)\big]^c=\big[f(O^c)\big]^c$~
is open since $O^c$ is compact (as a closed subset of a compact space) and so $f(O^c)$ is closed (by Corollary \ref{ClosedImLmm}).
\end{proof}

\subsection{Complete normality, local compactness, and compactification}
\begin{dfn}[\textcolor{blue}{\index{Separated sets}{Separated sets}}]
Let $X$ be a space and $A,B\subset X$. The subsets $A,B$ are \ul{separated} if $\ol{A}\cap B=\emptyset=A\cap\ol{B}$ (i.e., $A$ and $B$ are disjoint from the closure of each other).
\end{dfn}

\begin{lmm}
Let $X$ be a space and $A,B\subset X$. The subsets $A,B$ are separated (i.e., $\ol{A}\cap B=\emptyset=A\cap\ol{B}$) $\iff$ there exist open sets $U\supset A$ and $V\supset B$ such that $\ol{A}\cap V=\emptyset=U\cap\ol{B}$. (\blue{footnote}\footnote{That is, $A,B$ each have a neighborhood that excludes (i.e., is disjoint from) the closure of the other.}).
\end{lmm}
\begin{proof}
($\Ra$): If $A,B$ are separated, let $U:=\big(\ol{B}\big)^c$ and $V:=\big(\ol{A}\big)^c$. ($\La$): The converse is clear.
\end{proof}

\begin{dfn}[\textcolor{blue}{\index{Completely normal space}{Completely normal space}}]
A space $X$ is \ul{completely normal} if every two separated sets in $X$ have disjoint neighborhoods.
\end{dfn}

\begin{lmm}
A space is completely normal $\iff$ every subspace is normal. (In particular, a metric space is completely normal, since every subspace of a metric space is a metric space.)
\end{lmm}
\begin{proof}
($\Ra$): Assume $X$ is completely normal, and let $A\subset X$. We need to show $A$ is normal. Let $C,D\subset A$ be disjoint sets that are closed in $A$. Then {\small $\ol{C}\cap D=\emptyset=C\cap\ol{D}$}, since
\bea
\ol{C}\cap D=\ol{C}\cap A\cap D=\Cl_A(C)\cap D=C\cap D=\emptyset~~\txt{and similarly}~~C\cap\ol{D}=\emptyset.\nn
\eea
Thus (by the complete normality of $X$) there are disjoint open sets $U_1\supset C$ and $V_1\supset D$, which give disjoint sets $U:=U_1\cap A\supset C$ and $V:=V_1\cap A\supset D$ that are open in $A$.

($\La$): Assume every subset of a space $X$ is normal. Let $A,B\subset X$ be separated. We need to show $A,B$ have disjoint neighborhoods. Let $Z:=X-(\ol{A}\cap\ol{B})=(\ol{A})^c\cup(\ol{B})^c$. Then $\ol{A}\cap Z$ and $\ol{B}\cap Z$ are closed in $Z$ and so have disjoint containing sets $U\supset \ol{A}\cap Z\supset A$ and $V\supset\ol{B}\cap Z\supset B$ that are open in $Z$, and hence also open in $X$ since $Z$ is open in $X$.
\end{proof}

\begin{dfn}[\textcolor{blue}{
\index{Closure-compact set}{Closure-compact set},
\index{Locally! compact space}{Locally compact space}}]\label{LocCompSp}
Let $X$ be a space and $A\subset X$. The set $A$ is \ul{closure-compact} if its closure $\ol{A}$ is compact. The space $X$ is \ul{locally compact} if each point $x\in X$ has a closure-compact neighborhood.
\end{dfn}

\begin{dfn}[\textcolor{blue}{\index{Compactification}{Compactification} of a space, \index{One-point compactification}{One-point compactification}}]
Let $X$ and $Y$ be spaces. Then $Y$ is a \ul{compactification} of $X$ if (i) $Y$ is a compact Hausdorff space and (ii) there exists an imbedding $f:X\hookrightarrow Y$ with a proper dense image, i.e., $f(X)\subsetneq Y$ and $\ol{f(X)}=Y$. Given $Y$ as a compactification of $X$, if $Y\backslash f(X)$ is a singleton (i.e., consists of a single point), then $Y$ is a \ul{one-point compactification} of $X$ (which is unique up to homeomorphism by Theorem \ref{LCompChar} below).
\end{dfn}

\begin{thm}[\blue{\cite[Theorem 29.1, p.183]{munkres}}]\label{LCompChar}
Let $X$ be a space. Then $X$ is a locally compact Hausdorff space $\iff$ $X$ has a one-point compactification, in the sense there exists a space $Y$ satisfying the following conditions: (i) $X$ is a subspace of $Y$. (ii) The set $Y-X$ consists of a single point. (iii) $Y$ is a compact Hausdorff space. Moreover, if $Y$ and $Z$ are two spaces satisfying these conditions, then there is a homeomorphism $h:Y\ra Z$ such that $h|_X=id_X:X\ra X$.
\end{thm}
\begin{proof}
{\flushleft \ul{Step 1 (Uniqueness of $Y$)}}: Let $Y,Z$ be spaces satisfying (i),(ii),(iii). Then $Y=X\sqcup\{p\}$, $Z=X\sqcup\{q\}$. Define a map $h:Y\ra Z$ by $h(p):=q$, $h|_X:=id_X$. We will show $h$ is open (which implies $h$ is a homeomorphism by symmetry). Let $U\subset Y$ be open. If $p\not\in U$, then $h(U)=U\subset X\subset Z$ is open in $Z$ since $X$ is open in $Z$. So, assume $p\in U$. Then $C:=Y-U\subset X\subset Y$ is compact (as a closed subset of the compact space $Y$). Thus $C=h(C)$ is a compact subspace of $X\subset Z$, and hence a compact subspace of $Z$. Since $Z$ is Hausdorff, $C$ is closed in $Z$, and so $h(U)=h(Y-C)=h([Y-C-p]\cup p)=[Y-C-p]\cup q=Z-C$ is open in $Z$.

{\flushleft \ul{Step 2 ($\Ra$) (Existence of $Y$)}}:
Assume $X$ is a locally compact Hausdorff space. Let $p$ be any object that is not an element of $X$, and define $Y:=X\cup\{p\}$. (Note that $p$ is often denoted by the symbol $\infty$). Let us give $Y$ the following topology (call it $\T$). $A\subset Y$ is open iff (i) $A$ is an open subset of $X$ or (ii) $A=Y-K$ for a compact set $K\subset X$.

(An open set in $Y$ of the type ~$Y-K=\{p\}\cup(X-K)$~ is called a \ul{neighborhood of $\infty:=p$}.)

\bit[leftmargin=0.7cm]
\item[(a)] It is clear that $\emptyset,Y\in\T$ since $\emptyset$ is both compact and open in $X$.
\item[(b)] $\T$ is closed under finite intersections because if $O_1,O_2$ are open in $X$ and $K_1,K_2$ are compact in $X$, then $O_1\cap O_2$, $(Y-K_1)\cap(Y-K_2)=Y-(K_1\cup K_2)$, and $O_1\cap(Y-K_1)=O_1\cap(X-K_1)$ are all in $\T$.
\item[(c)] Similarly, $\T$ is closed under arbitrary unions because if $\{O_\al\}$ are open in $X$ and $\{K_\beta\}$ are compact in $X$, then $\bigcup O_\al$, $\bigcup(Y-K_\beta)=Y-\bigcap K_\beta$, and $\left[\bigcup O_\al\right]\cup \left[\bigcup(Y-K_\beta)\right]=Y-\left(\bigcap K_\beta-\bigcup O_\al\right)$ are all in $\T$.
\eit

\ul{$X$ is a subspace of $Y$}: If $O$ is open in $X$ and $K$ is compact in $X$, then $O\cap X=O$ and $(Y-K)\cap X=X-K$ are both open in $X$. Moreover, $X\cap\T:=\{A\cap X:A\in\T\}$ exhausts all open subsets of $X$. Hence, the topology of $X$ is the subspace topology induced by $\T$.

\ul{$Y$ is compact}: Let $\U$ be an open cover of $Y$. Then $\U$ must contain a set of the type $Y-K$ where $K\subset X$ is compact (otherwise, open subsets of $X$ alone cannot cover the point $p\in Y$). The sets $\V:=\{U\cap X:U\in\U,~U\neq Y-K\}$ form an open cover of $K$ in $X$. Since $K$ is compact, finitely many of these sets $\V':=\{U_1\cap X,...,U_n\cap X\}$ cover $K$. Hence, $\U':=\{U_1,...,U_n\}\cup\{Y-K\}\subset\U$ is a finite subcover.

\ul{$Y$ is Hausdorff}: Let $u,v\in Y$ be distinct points. If $u,v\in X$, it is clear that $u,v$ have disjoint open neighborhoods in $Y$. If $u\in X$ and $v=p$, then because $X$ is locally compact, we can choose a compact set $K\subset X$ containing an open neighborhood $U\subset X$ of $u$. It follows that $U\ni u$ and $Y-K\ni v$  are disjoint open neighborhoods of $u,v$.

{\flushleft \ul{Step 3 ($\La$) (Converse)}}: Assume $X\subset Y$, where $Y=X\cup\{p\}$ satisfies (i),(ii),(iii). Then $X$ is Hausdorff (as a subspace of a Hausdorff space). Let $x\in X$. Let $U\ni x$ and $V\ni p$ be disjoint open sets in $Y$. Then $C:=Y-V$ is closed, hence compact, in $Y$. It follows that $C$ is a compact subspace of $X$ containing the neighborhood $U$ of $x$ (and so $\ol{U}\subset C$ is compact as a closed subset of the compact subspace $C\subset X$).
\end{proof}

\subsection{Length analysis of paths and quasiconvexity}\label{PrlMET3}
\begin{dfn}[\textcolor{blue}{
\index{Path}{Path},
\index{Parametrization}{Parametrization},
\index{Length}{Length} of a path,
\index{Rectifiable path}{Rectifiable path},
\index{ML parametrization}{Minimum length (ML) parametrization},
\index{Natural! parametrization}{Natural parametrization},
\index{Constant speed path}{Constant speed path},
\index{Subconstant speed path}{Subconstant speed path}}]\label{PathLenDfn}
Let $X$ be a space. A \ul{path} in $X$ is a continuous map $\gamma:[0,1]\ra X$. Given paths $\gamma,\eta:[0,1]\ra X$, $\eta$ is a \ul{(re)parametrization} of $\gamma$, written $\eta\in\txt{Para}(\gamma)$, if it has the same image and endpoints as $\gamma$, i.e., if
\bea
\im\eta=\im\gamma,~~\eta(0)=\gamma(0),~~\eta(1)=\gamma(1).\nn
\eea

Let $X=(X,d)$ be a metric space and $\gamma:[0,1]\ra X$ a path (or any map). The \ul{length} of $\gamma$ is
\begin{equation*}
l(\gamma)~:=~\sup\big\{l_P(\gamma):P\subset [0,1]~~\txt{a finite set with $0,1\in P$}\big\},
\end{equation*}
where $l_P(\gamma):=\sum_{i=1}^kd(\gamma(t_{i-1}),\gamma(t_i))$ is the length of $\gamma$ over $P:=\{0=t_0<t_1<\cdots<t_k=1\}$. If $l(\gamma)<\infty$, we say $\gamma$ is \ul{rectifiable}. A parametrization $\eta$ of $\gamma$ is a \ul{minimum length (ML) parametrization} if
\bea
l(\eta)\leq l(\eta')~~~~\txt{for all}~~\eta'\in \txt{Para}(\gamma).\nn
\eea
A continuous map $\gamma:[a,b]\ra X$ (for $a\leq b$ in $\Real$) is a \ul{natural parametrization} if
\bea
l(\gamma|_{[t,t']})=|t-t'|~~~~\txt{for all}~~t,t'\in[a,b].\nn
\eea
A path $\gamma:[0,1]\ra X$ has \ul{constant speed} $c\geq 0$ if
\bea
l(\gamma|_{[t,t']})=c|t-t'|~~~~\txt{for all}~~t,t'\in[0,1].\nn
\eea
A path $\gamma:[0,1]\ra X$ has \ul{subconstant speed} $c\geq0$ if
\bea
l(\gamma|_{[t,t']})\leq c|t-t'|~~~~\txt{for all}~~t,t'\in[0,1].\nn
\eea
\end{dfn}

In the above definition, $l(\gamma)$ depends on the way $\gamma$ is parameterized, i.e., if $\gamma,\eta:[0,1]\ra X$ are parametrizations of the same path, then we can have $l(\gamma)\neq l(\eta)$. Also, $l(\gamma)=0$ if and only if $\gamma$ is constant.
\begin{question}
Let $X=(X,d)$ be a metric space. For any two paths $\gamma,\eta:[0,1]\ra X$, we have
\bea
&&\textstyle|l(\gamma)-l(\eta)|\leq \sup_{P}\inf_{P'}|l_P(\gamma)-l_{P'}(\eta)|\leq \sup_{P}\inf_{P'}\sum_i|d(\gamma(t_{i-1}),\gamma(t_i))-d(\eta(t'_{i-1}),\eta(t'_i))|\nn\\
&&\textstyle~~~~\leq \sup_{P}\inf_{P'}\sum_i\left[d(\gamma(t_{i-1}),\eta(t'_{i-1}))+d(\gamma(t_i),\eta(t'_i))\right]
\leq \sup_{P}\inf_{P'}2\sum_id(\gamma(t_i),\eta(t'_i)),\nn
\eea
where as usual $P,P'\subset[0,1]$ are given by $P=\{0=t_0<t_1<\cdots<t_n=1\}$ and $P'=\{0=t'_0<t'_1<\cdots<t'_n=1\}$. Is this a useful estimate (say for investigating continuity of length as a map $l:\C([0,1],X)\ra\Real$ or otherwise)?
\end{question}

\begin{lmm}[\blue{\cite[Proposition 2.5.9, p.46]{BBI}}]\label{ConstSpeedPar}
Let $X=(X,d)$ be a metric space, $\gamma:[0,1]\ra X$ a rectifiable path, and $l:=l(\gamma)$. Then there exists a nondecreasing (\blue{footnote}\footnote{A function $f:\Real\ra\Real$ is nondecreasing if $t<t'$ implies $f(t)\leq f(t')$.}) continuous map $\vphi:[0,1]\ra [0,l]$ and a natural parametrization $\ol{\gamma}:[0,l]\ra X$ such that $\gamma=\ol{\gamma}\circ\vphi:[0,1]\sr{\vphi}{\ral}[0,l]\sr{\ol{\gamma}}{\ral}X$.
\end{lmm}
\begin{proof}
Define the map ~$\vphi:[0,1]\ra[0,l],~t\mapsto l\big(\gamma|_{[0,t]}\big)$,~ which is clearly nondecreasing. $\vphi$ is continuous
(\magenta{footnote}\footnote{\magenta{Fix $t\in[0,1]$. Let $A(t,t'):=l(\gamma|_{[t,t']})$, $A(t):=\inf_{t'\geq t}A(t,t')=\lim_{|t-t'|\ra 0}A(t,t')$, and $A_P(t,t'):=l_{P|_{[t,t']}}(\gamma)$.\\ (i) Given $n\geq 0$, there is $\delta_{n,t}$ such that
\[
|t-t'|<\delta_{n,t}~~~~\Ra~~~~|A(t,t')-A(t)|<1/n.
\]
(ii) Similarly, given $n\geq 0$, there is $\delta_{n,t,t'}$ such that (with $|t-t'|<\delta_{n,t}$) for any $P=\{0=t_0<t_1<\cdots<t_{|P|}=1\}$,
\[
\mu(P)<\delta_{n,t,t'}~~~~\Ra~~~~|A(t,t')-A_P(t,t')|<1/n.
\]
(iii) Since $\gamma$ is uniformly continuous (as a continuous map on the compact space $[0,1]$), using $\vep_n:={1\over n|P|}$ (for any $P$ with $\mu(P)<\delta_{n,t,t'}$) we can choose $\delta_{\vep_n}$ such that $|a-b|<\delta_{\vep_n}$ implies $d(\gamma(a),\gamma(b))<\vep_n$. Therefore (with $\mu(P)<\delta_{n,t,t'}$),
\begin{align}
&\textstyle |t-t'|<\min(\delta_{n,t},\delta_{\vep_n})~~~~\Ra~~~~|t_{i-1}-t_i|<\delta_{\vep_n}~~~~\Ra~~~~A_P(t,t')=\sum_id(\gamma(t_{i-1}),\gamma(t_i))<|P|\vep_n=1/n\nn\\
&\textstyle~~\txt{while}~~~~|A(t,t')-A(t)|<1/n~~~~\txt{and}~~~~|A(t,t')-A_P(t,t')|<1/n,\nn\\
&~~\txt{and so}~~A(t)\leq|A(t)-A(t,t')|+|A(t,t')-A_P(t,t')|+A_P(t,t')<3/n\ra0.\nn
\end{align}
Hence $A(t)=0$.
}}), since with $P|_{[a,b]}:=\{a,b\}\cup P\cap[a,b]$ and $\mu(P):=\max\{|t_{i-1}-t_i|:t_i\in P\}$,
{\small\begin{align}
&\textstyle \vphi(t)=l(\gamma|_{[0,t]})=\sup_{P\subset[0,1]}l_{P|_{[0,t]}}(\gamma)\sr{\txt{(rf)}}{=}\lim_{\mu(P)\ra 0}l_{P|_{[0,t]}}(\gamma),\nn\\
&\textstyle |\vphi(t)-\vphi(t')|=l(\gamma|_{[t,t']})=\sup_{P\subset[0,1]}l_{P|_{[t,t']}}(\gamma)\sr{\txt{(rf)}}{=}\lim_{\mu(P)\ra 0}l_{P|_{[t,t']}}(\gamma),\nn
\end{align}}where step (rf) holds because (i) it is enough (by the uniqueness of supremum in $\Real$) to take the supremum over a maximal chain $C$ in $\{\txt{finite}~P\subset[0,1]\}$, viewed as a poset under inclusion, and (ii) for any $P,Q\in C$,
\[
P\subset Q ~~\iff~~ \mu(P)\geq\mu(Q) ~~\iff~~ l_{P|_{[0,t]}}(\gamma)\leq l_{Q|_{[0,t]}}(\gamma).
\]
Next, define
\bea
\ol{\gamma}:[0,l]\ra X,~{s}\mapsto\gamma(t)\in\gamma\left(\vphi^{-1}({s})\right),~~\txt{i.e.,}~~t\in\vphi^{-1}({s}),~~\txt{or}~~\vphi(t)={s}.\nn
\eea
$\ol{\gamma}$ is well defined because for any $t,t'\in\vphi^{-1}({s})$, $t\leq t'$,
\bea
&&d(\gamma(t),\gamma(t'))\leq l(\gamma|_{[t,t']})=l(\gamma|_{[0,t']})-l(\gamma|_{[0,t]})=\vphi(t')-\vphi(t)={s}-{s}=0,\nn\\
&&~~\Ra~~d(\gamma(t),\gamma(t'))=0,~~\Ra~~\gamma(t')=\gamma(t).\nn
\eea
Also, for any $t\in[0,1]$, we have $\gamma(t)=\ol{\gamma}(\vphi(t))$ since $t\in\vphi^{-1}(\vphi(t))$. Finally, with $t_s\in\vphi^{-1}(s)$, $t_{s'}\in\vphi^{-1}(s')$,
\[
l(\ol{\gamma}|_{[{s},{s}']})=l(\gamma|_{[t_{s},t_{s}']})=\big|l(\gamma|_{[0,t_{s}]})-l(\gamma|_{[0,t_{{s}'}]})\big|=\big|\vphi(t_{s})-\vphi(t_{{s}'})\big|=|{s}-{s}'|. \qedhere
\]
\end{proof}

\begin{crl}[\blue{Constant speed parametrization of a rectifiable path}]\label{ConstSpeedRP}
Let $X$ be a metric space. If $\gamma:[0,1]\ra X$ is a rectifiable path, there exists a path $\eta:[0,1]\ra X$ such that $\eta([0,1])=\gamma\big([0,1]\big)$, $\eta(0)=\gamma(0)$, $\eta(1)=\gamma(1)$, and ~$d(\eta(t),\eta(t'))\leq l(\eta|_{[t,t']})=l(\gamma)|t-t'|$~ for all ~$t,t'\in[0,1]$.
\end{crl}
\begin{proof}
Let $l:=l(\gamma)$. By Lemma \ref{ConstSpeedPar}, $\gamma=\ol{\gamma}\circ\vphi:[0,1]\sr{\vphi}{\ral}[0,l]\sr{\ol{\gamma}}{\ral}X$, for a nondecreasing continuous map $\vphi:[0,1]\ra [0,l]$ and a natural parametrization $\ol{\gamma}:[0,l]\ra X$. Thus, with $\psi:[0,1]\ra[0,l],~~t\mapsto lt$ we get the new parametrization {\small$\eta:=\ol{\gamma}\circ\psi:[0,1]\sr{\psi}{\ral}[0,l]\sr{\ol{\gamma}}{\ral}X$}, which satisfies $l(\eta|_{[t,t']})=l|t-t'|$ for all $t,t'\in[0,1]$.
\end{proof}

\begin{dfn}[\textcolor{blue}{
\index{Quasigeodesic}{Quasigeodesic},
\index{Standard! parametrization}{Standard parametrization},
\index{Quasiconvex space}{Quasiconvex space},
\index{Geodesic!}{Geodesic},
\index{Geodesic! space}{Geodesic space},
\index{Uniquely geodesic space}{Uniquely geodesic space}}]\label{QsiGeoDfn}
Let $X$ be a metric space and $\ld\geq1$. A path  $\gamma:[0,1]\ra X$ is a \ul{$\ld$-quasigeodesic} if it has subconstant speed $c=\ld d(\gamma(0),\gamma(1))$, i.e.,
\bea
\textstyle l(\gamma|_{[t,t']})\leq\ld d(\gamma(0),\gamma(1))|t-t'|,~~~~\txt{for all}~~~~t,t'\in[0,1].\nn
\eea
We say $X$ is a \ul{$\ld$-quasiconvex space} if for every $x,y\in X$, there exists a $\ld$-quasigeodesic $\gamma:[0,1]\ra X$ from $x$ to $y$, i.e., such that $\gamma(0)=x$, $\gamma(1)=y$. A $1$-quasigeodesic is called a \ul{geodesic}, and similarly, a $1$-quasiconvex space is called a \ul{geodesic space}. A geodesic space $X$ is a \ul{uniquely geodesic space} if for every $x,y\in X$, there exists only one geodesic from $x$ to $y$.
\end{dfn}

Note that a $\ld$-quasigeodesic is also called a \ul{$\ld$-quasiconvex path}, \cite[p.205]{hakobyan-herron2008}. In \cite[p.317]{tyson-wu2005}, a quasigeodesic is differently defined to be a path that is a bi-Lipschitz imbedding. Injectivity of the path is not required in our definition above. An alternative definition of a geodesic in terms of paths that are parameterized by arc length can be found in \cite[Definition 2.2.1, p.56]{papado2014}.

\begin{lmm}[\blue{Characterization and Sufficient condition for quasigeodesics}]\label{QgeodCharLmm}
Let $X$ be a metric space, $\gamma:[0,1]\ra X$ a path, and $\ld,\ld_1,...,\ld_n\geq 1$. Then the following are true:
\bit[leftmargin=0.8cm]
\item[(i)] $\gamma$ is a $\ld$-quasigeodesic $\iff$ $d(\gamma(t),\gamma(t'))\leq\ld d(\gamma(0),\gamma(1))|t-t'|$, for all $t,t'\in[0,1]$.
\item[(ii)] If {\small $[0,1]=\bigcup_{j=1}^n[a_j,b_j]$} such that $b_j=a_{j+1}$, $\gamma_j:=\gamma|_{[a_j,b_j]}$ is a $\ld_j$-quasigeodesic (for all $j=1,...,n$), and {\small $\ld:=\max_j\ld_j{d(\gamma(a_j),\gamma(b_j))\over d(\gamma(0),\gamma(1))}$}, then $\gamma$ is a $\ld$-quasigeodesic.
\eit
\end{lmm}
\begin{proof}
{\flushleft (i)} If $\gamma$ is a $\ld$-quasigeodesic, then $d(\gamma(t),\gamma(t'))\leq l(\gamma|_{[t,t']})\leq\ld d(\gamma(0),\gamma(1))|t-t'|$ for all $t,t'\in[0,1]$. Conversely, if  $d(\gamma(t),\gamma(t'))\leq\ld d(\gamma(0),\gamma(1))|t-t'|$ for all $t,t'\in[0,1]$, then
\bea
&&\textstyle l(\gamma|_{[t,t']})=\sup\{l_P(\gamma)~|~P:=\{t=t_0<t_1<\cdots<t_k=t'\}\subset[t,t']\},\nn\\
&&\textstyle l_P(\gamma)=\sum_{i=1}^kd(\gamma(t_{i-1}),\gamma(t_i))\leq \sum_{i=1}^k\ld d(\gamma(0),\gamma(1))|t_{i-1}-t_i|=\ld d(\gamma(0),\gamma(1))|t-t'|,\nn\\
&&~~\Ra~~l(\gamma|_{[t,t']})\leq \ld d(\gamma(0),\gamma(1))|t-t'|,~~~~\txt{for all}~~t,t'\in[0,1].\nn
\eea
{\flushleft (ii)} Let $t,t'\in[0,1]$. Then by calculating directly, we get
{\small\begin{align}
&\textstyle l(\gamma|_{[t,t']})=\sum_jl(\gamma([t,t']\cap[a_j,b_j]))\leq \sum_j\ld_jd(\gamma(a_j),\gamma(b_j))\big|[t,t']\cap[a_j,b_j]\big|\nn\\
&\textstyle~~~~\leq \max_j\ld_jd(\gamma(a_j),\gamma(b_j))~\sum_j\big|[t,t']\cap[a_j,b_j]\big|=\ld d(\gamma(0),\gamma(1))|t-t'|.\nn \qedhere
\end{align}}
\end{proof}

\begin{lmm}[\blue{Characterization of geodesics: Related to \cite[Sec. 2.2, pp 56-60]{papado2014}}]\label{GeodCharLmm}
Let $X$ be a metric space and $\gamma:[0,1]\ra X$ a path. Then (i) $\gamma$ is a geodesic $\iff$ (ii) $d(\gamma(t),\gamma(t'))\leq d(\gamma(0),\gamma(1))|t-t'|$ for all $t,t'\in[0,1]$, $\iff$ (iii) up to a reparametrization, $d(\gamma(t),\gamma(t'))=d(\gamma(0),\gamma(1))|t-t'|$ for all $t,t'\in[0,1]$.
\end{lmm}
\begin{proof}
{\flushleft \ul{(i)$\Ra$(ii)}:} This is immediate, since ~$d(\gamma(t),\gamma(t'))\leq l(\gamma|_{[t,t']})$, for all $t,t'\in[0,1]$.
{\flushleft \ul{(ii)$\Ra$(i),(iii)}:} (ii) implies $d(\gamma(t),\gamma(t'))\leq d(\gamma(0),\gamma(1))|t-t'|$ for all $t,t'\in[0,1]$. By the definition of $l(\gamma)$, $l(\gamma)\leq d(\gamma(0),\gamma(1))\leq l(\gamma)$, i.e., $l(\gamma)=d(\gamma(0),\gamma(1))$. Also, observe that $d(\gamma(0),\gamma(1))\leq l_P(\gamma)\leq l(\gamma)=d(\gamma(0),\gamma(1))$ for any finite partition $P\subset[0,1]$. That is, $l_P(\gamma)=l_{P'}(\gamma)$ for any two finite partitions $P,P'\subset [0,1]$. In particular, if $P:=\{0,t,t',1\}$ and $P':=P\cup Q=\{0\}\cup Q\cup\{1\}$ for any finite partition $Q\subset[t,t']$, then
\bea
l_P(\gamma)=l_{P'}(\gamma)~~\Ra~~d(\gamma(t),\gamma(t'))=l_Q(\gamma|_{[t,t']})=l(\gamma|_{[t,t']}).\nn
\eea
Hence, up to a reparametrization (by Corollary \ref{ConstSpeedRP}), ~{\small $d(\gamma(t),\gamma(t'))=d(\gamma(0),\gamma(1))|t-t'|$} ~for all ~$t,t'\in[0,1]$.
{\flushleft \ul{(iii)$\Ra$(ii)}:} This is again immediate.
\end{proof}

\begin{dfn}[\textcolor{blue}{\index{Minimum-length path}{Minimum-length path}}]
Let $X$ be a metric space, $x,y\in X$, and $\P_{x,y}(X):=\{\txt{paths}~\gamma:[0,1]\ra X,~\gamma(0)=x,\gamma(1)=y\}\subset \C\big([0,1],X\big)$. Then a path $\gamma:[0,1]\ra X$ is a \ul{minimum-length path} (or path of minimum length) if
\bea
l(\gamma)=\inf\big\{l(\eta)~|~\eta\in\P_{\gamma(0),\gamma(1)}(X)\big\}.\nn
\eea
\end{dfn}
It follows immediately from Lemma \ref{GeodCharLmm} that every geodesic is a minimum-length path.

\begin{lmm}[\blue{\index{Gluing! lemma I}{Gluing lemma I}}]\label{GluLmmI}
Let $X$ be a $\lambda$-quasiconvex space, $\{C_j:j=1,...,n\}$ a finite cover of $X$ by closed sets $C_j$, and $Y$ a metric space. If $f:X\ra Y$ is a map such that each restriction $f_j=f|_{C_j}:C_j\ra Y$ is $c$-Lipschitz, then $f$ is $\lambda c$-Lipschitz.
\end{lmm}
\begin{proof}
We know $f$ is continuous because for any closed set $C'\subset Y$, we have the finite union of closed sets $f^{-1}(C')=\bigcup_j C_j\cap f^{-1}(C')=\bigcup_j f_j^{-1}(C')$, where each $f_j^{-1}(C')$ is closed in $X$ because a relatively closed subset of a closed set is closed.

Given $x,x'\in X$, let $\gamma:[0,1]\ra X$ be a $\ld$-quasigeodesic from $x$ to $x'$ (with only finitely many intersection points $\bigcup_j\gamma\cap \del C_j$ with the boundaries of the cover). Consider any partition $P=\{0=t_0<t_1<\cdots<t_k=1\}$ of $[0,1]$ containing the intersection points $\bigcup_j\gamma\cap \del C_j$ in the sense that $\bigcup_j\gamma\cap\del C_j\subset\gamma(P)$, which is possible (by refinement) since we have only finitely many $C_j$. Then for each $i$, the segment $\gamma|_{[t_{i-1},t_i]}$ of $\gamma$ lies in some $C_j$. Therefore,
\begin{align}
&\textstyle d\big(f(x),f(x')\big)\leq \sum_{i=1}^kd\big(f(\gamma(t_{i-1})),f(\gamma(t_i))\big)=\sum_{i=1}^kd\big(f_{j_{i-1}}(\gamma(t_{i-1})),f_{j_i}(\gamma(t_i))\big)\nn\\
&\textstyle~~~~\sr{(s)}{\leq}c\sum_{i=1}^kd\big(\gamma(t_{i-1}),\gamma(t_i)\big)\leq c~\txt{length}(\gamma)\leq c\ld d(x,x'),\nn
\end{align}
where step (s) holds because of the choice of the partition $P$.
\end{proof}

\begin{dfn}[\textcolor{blue}{\index{Locally! Lipschitz map}{Locally Lipschitz map}}]\label{LocLipMap}
Let $X,Y$ be metric spaces and $c\geq 0$. A map $f:X\ra Y$ is \ul{locally $c$-Lipschitz} if for each $x\in X$, there exists a ball $B_{r_x}(x)$, $r_x=r_{x,f}>0$, such that
\bea
d(f(x),f(z))\leq cd(x,z)~~~~\txt{for all}~~~~z\in B_{r_x}(x).\nn
\eea
\end{dfn}

\begin{lmm}[\blue{\index{Gluing! lemma II}{Gluing lemma II}}]\label{GluLmmII}
Let $X$ be a $\ld$-quasigeodesic space and $Y$ a metric space. If a map $f:X\ra Y$ is locally $c$-Lipschitz, then it is $\ld c$-Lipschitz. (The converse is trivially true as well.)
\end{lmm}
\begin{proof}
Fix $x,x'\in X$. Let $\gamma:[0,1]\ra X$ be a $\ld$-quasigeodesic from $x$ to $x'$. Define $C:=cd(x,x')$. Then $f\circ\gamma:[0,1]\sr{\gamma}{\ral}X\sr{f}{\ral} Y$ is locally $\ld C$-Lipschitz on $[0,1]$, since for any $t\in[0,1]$, there exists a ball $B_{r_{\gamma(t)}}\big(\gamma(t)\big)$, $r_{\gamma(t)}>0$, such that
\bea
d\big(f(\gamma(t)),f(z)\big)\leq cd\big(\gamma(t),z\big)~~~~\txt{for all}~~~~z\in B_{r_{\gamma(t)}}\big(\gamma(t)\big),\nn
\eea
and so we get the neighborhood $U:=\gamma^{-1}\big(B_{r_{\gamma(t)}}\big(\gamma(t)\big)\big)$ of $t$ in $[0,1]$ satisfying
\bea
d\big(f(\gamma(t)),f(\gamma(s))\big)\leq cd\big(\gamma(t),\gamma(s)\big)\leq C\ld|t-s|,~~~~\txt{for all}~~s\in U.\nn
\eea

Since $[0,1]$ is compact, for any $t,t'\in [0,1]$ we can choose a  partition $P=\{t=t_0<t_1<\cdots<t_k=t'\}$ of $[t,t']$ such that for some refinement $Q=\{t=s_0<s_1<\cdots<s_{l}=t'\}\supset P$ of $P$,
\vspace{-0.2cm}
\begin{align}
\textstyle d\big(f(\gamma(t)),f(\gamma(t'))\big)\leq \sum_{i=1}^kd\big(f(\gamma(t_{i-1})),f(\gamma(t_i))\big)\sr{(a)}{\leq} \ld C\sum_{j=1}^l|s_{j-1}-s_j|=\ld C|t-t'|,\nn
\end{align}
where at step (a), we choose $P$ fine enough so that for each $i\in\{1,...,k\}$ there exists $s\in[t_{i-1},t_i]$ satisfying
\bea
d\big(f(\gamma(t_{i-1})),f(\gamma(s))\big)\leq\ld C|t_{i-1}-s|,~~~~d\big(f(\gamma(s)),f(\gamma(t_i))\big)\leq\ld C|s-t_i|.\nn
\eea
This shows $f\circ\gamma$ is $\ld C$-Lipschitz. Thus, $d\big(f(x),f(x')\big)\leq \ld C=\ld cd(x,x')$.
\end{proof}

\subsection{Metric spaces with nonrectifiable paths: Snowflake metric spaces}~\\~
Let $\ld\geq 1$. Recall that a map between metric spaces $f:(X,d)\ra (Z,d_1)$ is $\ld$-biLipschitz if $d(x,y)/\ld\leq d_1(f(x),f(y))\leq \ld d(x,y)$ for all $x,y\in X$.

\begin{dfn}[\textcolor{blue}{\index{Bi-Lipschitz equivalent metrics}{Bi-Lipschitz equivalent metrics}}]
Let $X$ be a set and $1<p\leq\infty$. Two metrics $d,d':X\times X\ra\Real$ are \ul{$\ld$-biLipschitz equivalent} if the identity map $id_X:(X,d)\ra (X,d')$ is $\ld$-biLipschitz, i.e., $d(x,y)/\ld\leq d'(x,y)\leq \ld d(x,y)$ for all $x,y\in X$.
\end{dfn}

\begin{dfn}[\textcolor{blue}{\index{$L^p$-metric}{$L^p$-metric}, \index{$L^p$-metric space}{$L^p$-metric space}}]
Let $(X,d)$ be a metric space and $1<p\leq\infty$. The metric $d$ is an \ul{$L^p$-metric} (making $(X,d)$ an \ul{$L^p$-metric space}) if for all $x,y,z\in X$,
\bea
\label{SnowFlakeEq}d(x,y)\leq\left\{
            \begin{array}{ll}
              \big(d(x,z)^p+d(z,y)^p\big)^{1\over p}, &\txt{if}~~1<p<\infty \\
              \max\big(d(x,z),d(z,y)\big), &\txt{if}~~ p=\infty
            \end{array}
          \right\}.
\eea
\end{dfn}

For example, if $(X,d)$ is a metric space, then for any $1<p<\infty$, the metric $d_p(x,y):=d(x,y)^{1\over p}$ is an $L^p$-metric on $X$. Thus, $(X,d_p)$ is a snowflake (as defined below).

\begin{dfn}[\textcolor{blue}{\index{Snowflake metric space}{Snowflake metric space}}]
Let $1<p\leq\infty$. A metric space $(X,d)$ is a \ul{$p$-snowflake} if the metric $d$ is biLipschitz equivalent to an $L^p$-metric ~$d':X\times X\ra\Real$~ on $X$.
\end{dfn}

\begin{lmm}[\blue{\cite[Proposition 2.3, p.318]{tyson-wu2005}}]
Let $1<p<\infty$. A metric space $(X,d)$ is a $p$-snowflake $\iff$ there exists a constant $c>0$ such that for any finite set of points $x_0,x_1,..,x_N\in X$, we have ~$\sum_{i=1}^Nd(x_{i-1},x_i)^p\geq cd(x_0,x_N)^p$.
\end{lmm}
\begin{proof}
($\Ra$): Assume $(X,d)$ is a $p$-snowflake. Let $d'$ be an $L^p$-metric on $X$ such that $d(x,y)/\ld\leq d'(x,y)\leq\ld d(x,y)$. Then for any finite set of points $x_0,x_1,...,x_N\in X$,
\bea
\textstyle d(x_0,x_N)^p/\ld^p\leq d'(x_0,x_N)^p\sr{(\ref{SnowFlakeEq})}{\leq} \sum_{i=1}^Nd'(x_{i-1},x_i)^p\leq \ld^p\sum_{i=1}^Nd(x_{i-1},x_i)^p.\nn
\eea
($\La$): Conversely, assume there exists $c$ such that $\sum_{i=1}^Nd(x_{i-1},x_i)^p\geq cd(x_0,x_N)^p$ for every finite set of points $x_0,x_1,...,x_N\in X$. Let
{\small\begin{align}
\textstyle d'(x,y):=\inf\big\{l(c)~|~c=\{x=x_0,x_1,...,x_N=y\}\big\},~~\txt{where}~~l(c):=\Big(\sum_{i=1}^Nd(x_{i-1},x_i)^p\Big)^{1\over p}.\nn
\end{align}}
Then given $\vep>0$ and $z\in X$, we can pick two finite chains $c_{xz}$ (a finite chain from $x$ to $z$) and $c_{zy}$ (a finite chain from $z$ to $y$) such that with the concatenation $c_{xy}:=c_{xz}\cdot c_{zy}$,
\bea
d'(x,y)^p\leq l(c_{xy})^p=l(c_{xz})^p+l(c_{zy})^p\leq d'(x,z)^p+d'(z,y)^p+2\vep.\nn
\eea
This shows $d'$ is an $L^p$-metric on $X$ satisfying $c^{1/p} d\leq d'\leq d$, and so $(X,d)$ is a $p$-snowflake.
\end{proof}

\begin{lmm}\label{SnwfNoRectCrv}
Let $1<p<\infty$ and $(X,d)$ a $p$-snowflake. Then every nonconstant path $\gamma:[0,1]\ra X$ is non-rectifiable (i.e., has infinite length).
\end{lmm}
\begin{proof}
Let $\gamma:[0,1]\ra X$ be a nonconstant path, i.e., there exist $t,t'\in[0,1]$, $t<t'$, such that $\gamma(t)\neq\gamma(t')$. Suppose $l(\gamma)<\infty$. Let $t=t_0,t_1,...,t_N=t'\in[t,t']$ be a finite sequence of points such that $l\left(\gamma|_{[t,t']}\right)/N\geq d(\gamma(t_{i-1}),\gamma(t_i))$ for all $i=1,...,N$. Then, upon exponentiating and summing,
\bea
\textstyle l\left(\gamma|_{[t,t']}\right)^p/N^{p-1}\geq\sum\limits_{i=1}^Nd(\gamma(t_{i-1}),\gamma(t_i))^p\geq cd(\gamma(t),\gamma(t'))^p.\nn
\eea
Taking the limit $N\ra\infty$, we get a contradiction.
\end{proof}

\begin{crl}\label{SnwfNoLipCrv}
Let $1<p<\infty$ and $X$ a $p$-snowflake. Every Lipschitz path $\gamma:[0,1]\ra X$ is constant.
\end{crl}

\begin{crl}\label{SnwfNoLipCon}
$Z$ be a quasiconvex space (e.g., geodesic space), $1<p<\infty$, and $X$ a $p$-snowflake. Every Lipschitz map $f:Z\ra X$ is constant.
\end{crl}

\begin{dfn}[\textcolor{blue}{\index{Snowflake curve}{Snowflake curve}}]
Let $X$ be a metric space. A path $\gamma:[0,1]\ra X$ is a \ul{snowflake curve} if the subspace $\gamma\big([0,1]\big)\subset X$ is both (i) a snowflake metric space and (ii) a curve (i.e., a $1$-manifold).
\end{dfn}

\begin{example}[\textcolor{blue}{von Koch snowflake curve, $K$: See \cite{tukia81,koskela94,mckemie1987}}]
In $\Real^2$, let {\small $K_0:=I=[(0,0),(1,0)]$. Next, let ~$K_1:=\left[\left(0,0\right),\left({1\over 3},0\right)\right]\cup\big[\left({1\over 3},0\right),\big({1\over 2},{1\over 3}{\sqrt{3}\over 2}\big)\big]\cup\big[\big({1\over 2},{1\over 3}{\sqrt{3}\over 2}\big),\left({2\over 3},0\right)\big]\cup\left[\left({2\over 3},0\right),\left(1,0\right)\right]$},~
a union of 4 line segments each of length $1/3$ formed from the line segment $K_0$. Similarly, we can apply the process to each of the 4 line segments of $K_1$ to obtain $K_2$ as a union of $4\times 4=4^2$ line segment each of length ${1\over 3}\times{1\over 3}$. Continuing this way, at the $j$th step, we get a set $K_j\subset\Real^2$ consisting of $4^j$ line segments each of length ${1\over 3^j}$. Note that the length of $K_j$ is $4^j\times{1\over 3^j}=\left(4\over 3\right)^j$.

In the limit $j\ra\infty$, we obtain a curve $K=\gamma(I)$, $\gamma:[0,1]\ra\Real^2$, of infinite length between $(0,0)$ and $(1,0)$. Moreover, between any two distinct points $\gamma(t),\gamma(t')\in K$, the segment $\gamma|_{[t,t']}$ also has infinite length.
\end{example}

{\flushleft A further} discussion of snowflake concepts can be found in \cite{tyson-wu2005}.

{\flushleft\hrulefill}
\begin{exercise}
Based on the discussion of this chapter, consider writing a \emph{fully technical} essay (say in the form of a typical chapter of these notes) on what is known in the mathematics literature as ``Measure Theory''.
\end{exercise} 

%% file: parts/AlgebraM/GeomAnaIII.tex
\chapter{Geometry and Analysis III: Topological Analysis with Convergence}\label{GeomAnaIII}

This chapter is a continuation of the previous chapter, and we also retain the following convention.

\begin{note}[\textcolor{blue}{Recall: Convention for countability}]
As in the previous chapter, for brevity, ``\emph{countable}'' will mean ``\emph{at most countable}'' (i.e., ``finite or countable'') unless it is specified otherwise.
\end{note}

\section{Topological Convergence and Characterizations}
\subsection{Convergence, completeness, and continuity}\label{GeomPrelimsS4}
\begin{dfn}[\blue{
\index{Additive! topology}{Additive topology},
\index{Scalar topology}{Scalar topology},
\index{Translation-invariant topology}{Translation-invariant topology},
\index{Scale-invariant topology}{Scale-invariant topology}}]
Let $M$ be a topological $R$-module. The topology of $M$ is \ul{additive} if addition $+:M\times M\ra M$ is open. (We will show below that the topology of $M$ is always additive.). The topology of $M$ is \ul{scalar} scalar multiplication $s:R\times M\ra M$ is open.

The topology of $M$ is \ul{translation-invariant} if for each $m\in M$, the continuous restriction $+|_{\{m\}\times M}:\{m\}\times M\ra M$ (or equivalently, $+_m:M\ra M,~x\mapsto m+x$) is open. (\magenta{Observe that, for any ring $R$, each of the continuous maps $+_m$ has a continuous inverse $+_{-m}$, and so is a homeomorphism. It follows that the topology of $M$ is always (i) translation-invariant, and so (ii) for any subset $A\subset M$ and any open set $U\subset M$, the set $A+U=\bigcup_{a\in A}(a+U)=\bigcup_{a\in A}+_a(U)\subset M$ is open as a union of open sets, and hence in particular (iii) the topology of $M$ is always additive.}).

The topology of $M$ is \ul{scale-invariant} if for each nonzero scalar $r\in R\backslash 0$, the continuous restriction $s|_{\{r\}\times M}:\{r\}\times M\ra M$ (or equivalently, $s_r:M\ra M,~m\mapsto rm$) is open. (\magenta{Observe that if $r\in U(R)$, then the continuous map $s_r$ has a continuous inverse $s_{r^{-1}}$, and so is a homeomorphism, in which case, if $U\subset M$ is open, then so are $rU,r^{-1}U\subset M$. Consequently, for any set of units $A\subset U(R)$ and any open set $U\subset M$, the set $AU=\bigcup\{a_1U+\cdots+a_nU:a_i\in A,~n\geq 1\}\subset M$ is open as a union of open sets. In particular, if $R$ is a division ring, then the topology of $M$ is scale-invariant.}).
\end{dfn}

\begin{crl}\label{InvOfTopLmm}
If $D$ is a division ring and $M$ a topological $D$-module, then the topology of $M$ is translation-invariant, scale-invariant, and additive.
\end{crl}

\begin{dfn}[\textcolor{blue}{\index{Convergent! sequence}{Convergent sequence},
\index{Difference! Cauchy sequence}{Difference Cauchy sequence},
\index{Difference! sequentially complete module}{Difference sequentially complete module}}]
Let $X$ be a space and $M$ a topological $R$-module. A sequence in $X$ is a map $s=(s_n)_{n\in\Natural}:\Natural\ra X,~n\mapsto s_n$. A sequence $(x_n):=(x_n)_{n\in\Natural}\subset X$ \ul{converges} (or is \ul{convergent}) to a point $x$, written $x_n\ra x$, if for any open nbd $U$ of $x$, there exists $N$ such that $\{x_N,x_{N+1},\cdots\}\subset U$ (i.e., $U$ contains all but finitely many $x_n$).

A sequence $(m_n)\subset M$ is \ul{(difference) Cauchy} if $(m_n)\ominus(m_n):=(m_n-m_{n'})_{(n,n')\in\Natural}\ra 0$ in the sense that for any open nbd $U$ of $0$, there exists $N$ such that $m_n-m_{n'} \in U$ for all $n,n'\geq N$ (i.e., $U$ contains all but finitely many $m_n-m_{n'}$). The module $M$ is \ul{(difference) sequentially complete} if every Cauchy sequence in $M$ converges.
\end{dfn}

Recall from the previous chapter that a topological $R$-module $M$ is \ul{difference-invariant} if there exists an element $r_0\in R$ called a \ul{DI-scale} (difference-invariance scale) of $M$ such that for every neighborhood $U_0\subset M$ of $0$, we have the containment
\[
U_0-U_0\subset r_0U_0,~~~~\txt{where}~~U_0-U_0:=\{a-b:a,b\in U_0\}~~\txt{and}~~r_0U_0:=\{r_0u:u\in U_0\}.
\]

\begin{prp}[\blue{Difference Cauchyness of a convergent sequence}]
Let $M$ be a topological $R$-module. If $M$ is difference-invariant with a DI-scale $\al$ that is a unit, i.e., $\al\in U(R)$, then every convergent sequence in $M$ is Cauchy.
\end{prp}
\begin{proof}
If $(m_n)_{n\in\Natural}\subset M$ converges to $m\in M$, then any neighborhood $U_0$ of $0$ contains all but finitely many $\al(m_n-m)$, $n\in\Natural$, and so also contains all but finitely many $m_n-m_{n'}=(m_n-m)-(m_{n'}-m)\in \al^{-1}(U_0-U_0)\subset U_0$ ($n,n'\in\Natural$).
\end{proof}

\begin{crl}
If $D$ is a division ring and $M$ a difference-invariant topological $D$-module (e.g., a difference-invariant topological vector space), then every convergent sequence in $M$ is Cauchy.
\end{crl}

\begin{dfn}[\textcolor{blue}{\index{Net}{Net},
\index{Sequence}{Sequence},
\index{Tail of a net}{Tail of a net},
\index{Limit! of a net}{Limit of a net},
\index{Convergent! net}{Convergent net},
\index{Convergence in a subspace}{Convergence in a subspace},
\index{Cluster point of a net}{Cluster point of a net},
\index{Cofinal subset}{Cofinal subset},
\index{Cofinal map}{Cofinal map},
\index{Subnet}{Subnet}}]
Let $X$ be a space and $I$ a \ul{directed set}. A net in $X$ is a map of the form $f:I\ra X$. If there is an order isomorphism $I\cong\Natural$ (or more generally, if $I$ is countable and well-ordered), then the net $f:I\ra X$ is called a \ul{sequence} in $X$.

A \ul{tail} of a net $I\sr{f}{\ral} X$ is the image $f\big(U_{i_0}\big)=\{f(i)\}_{i\geq i_0}$ of the upper section $U_{i_0}:=[i_0,I]=\{i\in I:i\geq i_0\}$ at some $i_0\in I$. (\blue{footnote}\footnote{In general, we can of course have two tails with none containing the other. However, any finite collection of tails $f([i_1,I])$, ..., $f([i_n,I])$ has a common intersection point, since every finite set in $I$ (hence $\{i_1,...,i_n\}$) has an upperbound.}).

A point $x_0\in X$ is a \ul{limit} of a net $I\sr{f}{\ral}X$ (written $x_0\in\lim f$) if every neighborhood of $x_0$ contains a tail of $f$. If $x_0$  is a limit of $f$, we also say $f$ \ul{converges to} $x_0$ (written $f\ra x_0$, or $x_0\in\lim f$).

If a net $f:I\ra X$ converges to at least one point (i.e., $\lim f\neq\emptyset$) we say $f$ is \ul{convergent}. If every limit of $f:I\ra X$ lies in a subset $A\subset X$ (i.e., $\lim f\subset A\subset X$) we say $f$ \ul{converges in} $A$ (as a subspace).

A point $x_0\in X$ is a \ul{cluster point} of a net $I\sr{f}{\ral}X$ if every neighborhood of $x_0$ intersects each/every tail of $f$. (Note that a limit of a net is clearly a cluster point of the net.)

A subset $K\subset I$ is a \ul{cofinal subset} (as defined before) if every $i\in I$ has an upper bound in $K$, i.e., a $k\in K$ such that $i\leq k$. (\blue{footnote}\footnote{A cofinal set is directed, because if $k_1,k_2\in K$, then there exist $i\in I$ and $k\in K$ such that $k_1,k_2\leq i\leq k$. Also, it may be convenient to note that a cofinal set $K\subset I$ is precisely an exhaustive set of upper bounds in $I$.}). Given any directed sets $K,I$, a map $\phi:K\ra I$ is a \ul{cofinal map} if for each $i\in I$, there exists $k_i\in K$ such that $\phi(k)\geq i$ for all $k\geq k_i$, i.e., $\phi([k_i,K])\subset[i,I]$ (which implies $\phi(K)\subset I$ is a cofinal subset). Given a net $f:I\ra X$, a \ul{subnet} of $f$ is a any right-composition $f\circ\phi:K\sr{\phi}{\ral}I\sr{f}{\ral}X$ of $f$ with a cofinal map $\phi:K\ra I$ (where we write $f\circ\phi\subset f$). (\blue{footnote}\footnote{That is, given nets $I\sr{f}{\ral}X$ and $J\sr{g}{\ral}X$, $g$ is a \ul{subnet} of $f$ (written $g\subset f$) if $g$ factors through $f$ by a cofinal map, in the sense $g=f\circ\phi:J\sr{\phi}{\ral}I\sr{f}{\ral}X$ for some cofinal map $\phi:J\ra I$.}).
\end{dfn}
\begin{rmk}
Let $X$ be a space, $f:I\ra X$ a net, and $\phi:K\ra I$ a cofinal map. For each $i\in I$, let
{\small\[
K_i:=\{k'\in K:\phi(k)\geq i~~\txt{for all}~~k\geq k'\}=\{k'\in K:\phi([k',K])\subset[i,I]\}=\{k'\in K:[k',K]\subset\phi^{-1}([i,I])\}.
\]}
Then associated with $\phi$ are the selection maps $s_\phi:I\ra K,~i\mapsto k_i\in K_i$,
which give the cofinal maps
\[
\phi\circ s_\phi:I\sr{s_\phi}{\ral}K\sr{\phi}{\ral}I,~i\mapsto\phi(k_i)\in[i,I]\cap\phi(K).
\]

In the subnet $f\circ\phi:K\sr{\phi}{\ral}I\sr{f}{\ral}X$, if $\phi$ is an order-imbedding then we can assume wlog that $K\subset I$, because we can replace $\phi:K\ra I$ with the cofinal map given by the inclusion $i_{\phi(K)}:\phi(K)\hookrightarrow I,~\phi(k)\mapsto\phi(k)$. Thus, whenever applicable, we can simply assume $K\subset I$.

For the above situation (where $\phi$ is an order-imbedding), by introducing some redundancy, we can assume wlog that $K=I$, because we can replace $\phi:K\ra I$ with an extension of the inclusion $i_{\phi(K)}:\phi(K)\hookrightarrow I,~\phi(k)\mapsto\phi(k)$ given by any of the cofinal maps $\phi\circ s_\phi$ that satisfies $\phi\circ s_\phi|_{\phi(K)}=id_{\phi(K)}\eqv i_{\phi(K)}$.

Consequently, in the subnet $f\circ\phi:K\sr{\phi}{\ral}I\sr{f}{\ral}X$, if possible (i.e., if no loss of generality is incurred), we will for simplicity assume $K=I$ by replacing $\phi:K\ra I$ with one of the cofinal maps $\phi\circ s_\phi$ that best reflects/describes $\phi$.
\end{rmk}

\begin{convention}
Let $X$ be a space and $I$ a directed set. A net $f:I\ra X,~i\mapsto f_i:=f(i)$ is often given simply as $(f_i)_{i\in I}\subset X$, or as $(f_i)\subset X$ when the role of $I$ is already understood. Similarly, a subnet
\bea
f\circ\phi:K\sr{\phi}{\ral} I\sr{f}{\ral}X,~k\mapsto f_{\phi(k)}=f_{\phi_k}:=f(\phi(k))\nn
\eea
of $f$  is often given simply as ~$(f_{\phi(k)})_{k\in K}\subset X$ ~or as ~$(f_{\phi(k)})\subset X$. Sometimes (as we have been, rather naively, doing so far) we \ul{stack} indices by writing $\phi(k)$ simply as $i_k$, in which case, a subnet of $(f_i)_{i\in I}$ takes the form $(f_{i_k})_{k\in K}$ (in which the old index $i\in I$ is now instead viewed a cofinal map $i:K\ra I,~k\mapsto i_k$).

The above mentioned index stacking is naive because the symbol $i_k$ (for a fixed $k\in K$) also occurs more naturally as a permutation $i_k:I\ra I$ in the following sense: Observe that for each $k\in K$, the inclusion $i_I:I\hookrightarrow K\times I$ factors in the form {\footnotesize $i_I=i_{\{k\}\times I}\circ i_k:I\sr{i_k}{\ral}\{k\}\times I\sr{i_{\{k\}\times I}}{\ral}K\times I$}, where the bijection {\small $i_k:I\ra \{k\}\times I,~i\mapsto(k,i)$} is useful in (re)indexing as follows: Consider an indexing of nets {\small $F:K\ra\{\txt{nets}~I\ra X\},~k\mapsto F_k$} (also written $(F_k)_{k\in K}$ or $\{F_k\}_{k\in K}$) which is equivalently a map $F:K\times I\ra X$. Then as shown in the following commutative diagram of sets, we have the equivalent expressions
\[
F=(F_k)_{k\in K}=(F|_{\{k\}\times I})_{k\in I}=(F_{i_k})_{k\in K}.
\]
\[\adjustbox{scale=0.9}{\bt
I\ar[ddddrr,bend right,"F_{i_k}"']\ar[dr,"i_k"] &   & & & I\ar[ddddll,bend left,"F_{i_{k'}}"]\ar[dl,"i_{k'}"'] \\
 &\{k\}\times I\ar[dddr,bend right,"F_k"]\ar[dr,hook,"i_{\{k\}\times I}"] & & \{k'\}\times I\ar[dddl,bend left,"F_{k'}"']\ar[dl,hook,"i_{\{k'\}\times I}"'] & \\
 &              & K\times I\ar[dd,"F"] &                & \\
 &              &           &                & \\
 &              & X         &                &
\et}\]
Despite this possibility of confusion, the appropriate meaning of $i_k$ is often apparent from the context.
\end{convention}

\begin{lmm}[\textcolor{blue}{Hausdorffness criterion}]\label{NetHausLmm}
A space $X$ is Hausdorff (Definition \ref{SpaceTypeDfn}, page \pageref{SpaceTypeDfn}) $\iff$ every convergent net $f:I\ra X$ has a unique limit.
\end{lmm}
\begin{proof}
($\Ra$): Assume $X$ is Hausdorff. Let a net $f:I\ra X$ converge to a point $x_0\in X$. The limit $x_0$ is unique; otherwise, if $x_1\neq x_0$ is another limit, then $x_0$ and $x_1$ have disjoint neighborhoods which must, however, intersect because the two tails of $f$ contained in these neighborhoods must intersect (a contradiction).

($\La$): Assume $X$ is not Hausdorff. Let $x,x'\in X$ be distinct points with no disjoint neighborhoods, i.e., for every nbd $O$ of $x$ and every nbd $O'$ of $x'$, we have $O\cap O'\neq\emptyset$. Let $I:=\{O\cap O':\txt{for open sets}~O\ni x,O'\ni x'\}$ as a directed set under containment (i.e., given $i,j\in I$, we write $i\leq j$ iff $i\supset j$). Then a selection map $f:I\ra X,~i\mapsto f(i)\in i$ is a convergent net with both $x$ and $x'$ as limits.
\end{proof}

\begin{lmm}[\textcolor{blue}{\index{Closedness criterion}{Closedness criterion}}]\label{NetClosedLmm}
Let $X$ be a space and $A\subset X$. Then $A$ is closed $\iff$ every convergent net $f:I\ra A$ converges in $A$.
\end{lmm}
\begin{proof}
($\Ra$): If $A$ is closed and a net $f:I\ra A$ converges to a point $x_0\in X$, then $x_0\in A=\ol{A}$, otherwise a tail of $f$ will lie outside of $A$. ($\La$): Conversely, assume every convergent net $f:I\ra A$ converges in $A$. Suppose $A$ is not closed. Then there is a point $x_1\in\ol{A}-A$. So, for every neighborhood $N(x_1)$ of $x_1$, we have $(N(x_1)\backslash\{x_1\})\cap A\neq\emptyset$. Let $I$ consist of the neighborhoods $N(x_1)$ ordered by containment $\supset$ (i.e., {\small $N_1(x_1)\leq N_2(x_1)$ iff $N_1(x_1)\supset N_2(x_1)$}). Then we can pick a net {\small $f:I\ra A$, $N\mapsto f(N)\in (N\backslash\{x_1\})\cap A$} that clearly converges to $x_1\not\in A$ (a contradiction).
\end{proof}

\begin{prp}[\textcolor{blue}{\index{Continuity criterion}{Continuity criterion}}]\label{NetContPrp}
Let $X,Y$ be spaces. A function $f:X\ra Y$ is (i) continuous $\iff$ (ii) for every net $s:I\ra X$ and $x\in X$, $s\ra x$ implies $f\circ s\ra f(x)$.
\end{prp}
\begin{proof}
{\flushleft (i)$\Ra$(ii)} Assume $f$ is continuous. Let $x\in X$ and $s:I\ra X$ be a net such that $s\ra x$. Suppose $f\circ s\not\ra f(x)$. Then, with $\{T_i^{f\circ s}\}_{i\in I}$ denoting the collection of all tails $T_i^{f\circ s}:=f\circ s\big([i,I]\big)$ of $f\circ s$, there is an open nbd $V$ of $f(x)$ that misses a point $y_i=f(s(i'))$, $i'\geq i$, from $T_i^{f\circ s}$ (i.e., $y_i\in T_i^{f\circ s}\backslash V$) for all $i\in I$. Thus, the open nbd $f^{-1}(V)$ of $x$ also misses a point $x_i:=s(i')$ from every tail $T_i^s:=s\big([i,I]\big)$ of $s$:
\bea
x_i\in (f^{-1}(y_i)\cap T_i^s)\backslash f^{-1}(V)\subset T_i^s\backslash f^{-1}(V)\subset f^{-1}(T_i^{f\circ s}\backslash V),~~\txt{for all}~~i\in I~~\txt{(a contradiction)}.\nn
\eea

{\flushright (ii)$\Ra$(i)} Assume that for every net $s:I\ra X$ and $x\in X$, $s\ra x$ implies $f\circ s\ra f(x)$. Let $C\subset Y$ be closed. We need to show $f^{-1}(C)$ is closed. By Lemma \ref{NetClosedLmm}, it suffices to show that if $s:I\ra f^{-1}(C)\subset X$ is a net that converges in $X$, then any limit of $s$ lies in $f^{-1}(C)$. Suppose $s\ra x$ but $x\not\in f^{-1}(C)$, i.e., $f(x)\not\in C$. Then, because $f\circ s:I\ra C$ (i.e., $f\circ s$ is a net in $C$), we have $f\circ s\not\ra f(x)$ by Lemma \ref{NetClosedLmm} (since $C$ is closed). This contradicts our starting assumption.
\end{proof}

\begin{dfn}[\textcolor{blue}{
\index{Difference! Cauchy net}{Difference Cauchy net},
\index{Difference! complete module}{Difference complete module}}]
Let $M$ be a topological $R$-module. A net $f:I\ra M,i\mapsto f_i$ is \ul{(difference) Cauchy} if the associated net $g:=f\ominus f:I\times I\ra M,~(i,j)\mapsto f_i-f_j$ converges to $0$. The module $M$ is \ul{(difference) complete} if every Cauchy net in $M$ converges.
\end{dfn}

\begin{prp}[\blue{Difference Cauchyness of a convergent net}]
Let $M$ be a topological $R$-module. If $M$ is difference-invariant with a DI-scale $\al$ that is a unit, i.e., $\al\in U(R)$, then every convergent net in $M$ is Cauchy.
\end{prp}
\begin{proof}
Note that if $I,J$ are posets and $(i,j)\in I\times J$, then $[(i,j),I\times J]=[i,I]\times[j,J]$. If $(m_i)_{i\in I}\subset M$ converges, then any neighborhood $U_0$ of $0$ contains a tail of $\al(m_i-m)$, $i\in I$, and so also contains a tail of $m_i-m_j=(m_i-m)-(m_j-m)\in\al^{-1}(U_0-U_0)\subset U_0$, ($i,j\in I$). Expressed differently,
\[
\al(m_{[i,I]}-m)~\subset~U_0~~\Ra~~m_{[i,I]}-m_{[j,I]}=(m_{[i,I]}-m)-(m_{[j,I]}-m)\subset \al^{-1}(U_0-U_0)\subset U_0. \qedhere
\]
\end{proof}

\begin{crl}
If $D$ is a division ring and $M$ a difference-invariant topological $D$-module (e.g., a difference-invariant topological vector space), then every convergent net in $M$ is Cauchy.
\end{crl}

\begin{dfn}[\textcolor{blue}{
\index{Cauchy net}{Cauchy net},
\index{Cauchy space}{Cauchy space},
\index{Complete! space}{Complete space},
\index{Uniformly continuous map}{Uniformly continuous map}}]
Let $X,Z$ be spaces, and $u:X\times X\ra Z$ a continuous map. A net $f:I\ra X$ in $X$ is \ul{$u$-Cauchy} if the associated net $u\circ(s\times s):I\times I\sr{s\times s}{\ral}X\times X\sr{u}{\ral}Z,~(i,j)\mapsto u(s(i),s(j))$ in $Z$ converges. The space $X=(X,u)$ is \ul{$u$-Cauchy} if every convergent net in $X$ is a $u$-Cauchy net. The space $X$ is \ul{$u$-complete} (equiv., the space $X=(X,u)$ is \ul{complete}) if every $u$-Cauchy net in $X$ converges.

Let $X,Y$ be spaces and $u:X\times X\ra Z$ a fix continuous map as above. A map $f:X\ra Y$ is \ul{$u$-continuous} (or \ul{uniformly continuous wrt $u$}) if it maps $u$-Cauchy nets to $u$-Cauchy nets (i.e., for every $u$-Cauchy net $s:I\ra X$ in $X$, $f\circ s:I\sr{s}{\ral}X\sr{f}{\ral}Y$ is a $u$-Cauchy net in $Y$). (\blue{footnote}\footnote{The definition of uniform continuity here is designed to resemble the continuity criterion in Proposition \ref{NetContPrp}.})
\end{dfn}
In the above definition, a $u$-complete space is not required to be $u$-Cauchy, but $u$-completeness might require $u$-Cauchyness in order to be truly useful (depending on the type of application considered). A concept of \ul{uniform spaces} (not of immediate relevance to us) exists in the literature, but the map $u$ introduced above (although possibly related) is independent of the idea of a \ul{uniform structure} (the structure that underlies a uniform space).
\begin{examples}
Two examples of the continuous map $u:X\times X\ra Z$ in the above definition are the following: (i) If $X=(X,d)$ is a metric space, an example of $u$ is the \ul{metric} $d:X\times X\ra\Real$. (ii) If $X$ is a topological $R$-module, an example of $u$ is the \ul{difference map} (or \ul{subtraction}) $X\times X\ra X,~(x,x')\mapsto x-x'$.
\end{examples}

\begin{question}
Let $X,Z$ be spaces and $u:X\times X\ra Z$ a continuous map. Based on the idea of a DI scale (for the case where $X$ is a topological $R$-module), what possible conditions on $u$ ensure that $X$ is $u$-Cauchy (i.e., that every convergent net in $X$ is $u$-Cauchy)? Note that a metric space is an example of a Cauchy space (with the metric playing the role of $u$).
\end{question}

\begin{lmm}
Let $M$ be a topological $R$-module (with a scale-invariant topology), $m_0,m_0'\in M$, and $s,s':I\ra M$ nets. Let $s+s':I\ra M,~i\mapsto s_i+s'_i$ and for $\al\in R$, let $\al s:I\ra M,~i\mapsto \al s_i$. If $s\ra m_0$ and $s'\ra m_0'$, then $s+s'\ra m_0+m_0'$ and $\al s\ra \al m_0$.
\end{lmm}
\begin{proof}
By hypotheses, addition and the scaling $m_\al:M\ra M,~m\mapsto \al m$, are continuous, while
\[
s+s'=+\circ(s\times s'):I\times I\sr{s\times s'}{\ral} M\times M\sr{+}{\ral}M,~~~~\al s=m_\al\circ s:I\sr{s}{\ral}M\sr{m_\al}{\ral}M,
\]where $s\times s'\ra (m_0,m_0')$. Hence the results follow immediately from Proposition \ref{NetContPrp}.
\end{proof}

\begin{crl}
Let $M$ be a topological $R$-module (with a scale-invariant topology), $m_0\in M$, and $u:M\times M\ra M,~(m,m')\mapsto m-m'$ the difference map. If $M$ is $u$-Cauchy (i.e., every convergent net in $M$ is $u$-Cauchy), then every $u$-Cauchy net $s:I\ra M$ with a convergent subnet $s\circ\phi:K\sr{\phi}{\ral}I\sr{s}{\ral} M$ (such that $s\circ\phi\ra m_0$) is also convergent (such that $s\ra m_0$).
\end{crl}
\begin{proof}
Observe that $s\circ\phi\ra m_0$ $\iff$ $s\circ\phi-m_0\ra 0$. Thus, by \ul{freely} adding nets and subnets by \ul{viewing} them as maps $I\times K\ra M$ (modulo the inclusions $I,K\hookrightarrow I\times K$) in the form
\begin{align}
&s-m_0:I\times K\ra M,~(i,k)\mapsto s(i)-m_0,~~~~s-s\circ\phi:I\times K\ra M,~(i,k)\mapsto s(i)-s\circ\phi(k),\nn\\
&s\circ\phi-m_0:I\times K\ra M,~(i,k)\mapsto s\circ\phi(k)-m_0,\nn
\end{align}
we see that $s-m_0=(s-s\circ\phi)+(s\circ\phi-m_0)\ra 0+0$ (i.e., $s\ra m_0$).

(Note: Earlier, we used the symbol $\ominus$, instead of $-$, for a similar type of ``differences''. Therefore, another alternative is to simply complement the operation $\ominus$ by replacing $+$ with an analogously defined operation $\oplus$, such that $s_1\ominus s_2=s_1\oplus(-s_2)$.)
\end{proof}

\subsection{Compactness and Tychonoff's theorem}
\begin{dfn}[\textcolor{blue}{
\index{Eventual containment}{Eventual containment},
\index{Frequent containment}{Frequent containment}}]
Let $f:I\ra X$ be a net and $A\subset X$. $f$ is \ul{eventually} in $A$ if $A$ contains a tail of $f$. (\blue{footnote}\footnote{Note/recall: If $f$ is eventually in each $A_1,...,A_n\subset X$ then $A_1\cap \cdots\cap A_n\neq\emptyset$, because any finite collection of tails $f([i_1,I])$,...,$f([i_n,I])$ has a common intersection point, namely, $f(i)$ for a common upper bound $i\geq\{i_1,...,i_n\}$ which exists because $I$ is directed.}). $f$ is \ul{frequently} in $A$ if $f$ is not eventually in $X-A$, i.e., if every tail of $f$ intersects $A$  (equivalently, $f([i,I])\cap A\neq\emptyset$ for all $i\in I$, i.e., for each $i\in I$, there is $j\geq i$ such that $f(j)\in A$).
\end{dfn}

\begin{rmks*}
(i) A net $f:I\ra X$ is frequently in $A\subset X$ iff $f^{-1}(A)\subset I$ is cofinal.
{\flushleft (ii)} A net $f:I\ra X$ is frequently in $A\subset X$ iff $f(K)\subset A$ for some cofinal set $K\subset I$ (i.e., $A$ contains the image of a cofinal set). Indeed, if $f$ is frequently in $A$, then we can set $K=f^{-1}(A)$. Conversely, if $f(K)\subset A$ for some cofinal set $K\subset I$, then for any $i\in I$, there is $k\in K$ such that $k\geq i$ and $f(k)\in A$.
{\flushleft (iii)} A net $f:I\ra X$ converges to $x_0\in X$ iff $f$ is eventually in every neighborhood of $x_0$.
{\flushleft (iv)} A point $x_0\in X$ is a cluster point of a net $f:I\ra X$ iff $f$ is frequently in every neighborhood of $x_0$.
\end{rmks*}

\begin{lmm}
Let $X$ be a space and $A\subset X$. (i) A net $f:I\ra X$ is eventually in $A$ $\iff$ every subnet {\footnotesize $f\circ\phi:K\sr{\phi}{\ral}I\sr{f}{\ral}X$} is eventually in $A$. (ii) A net $f:I\ra X$ is frequently in $A\subset X$ $\iff$ some subnet {\footnotesize $f\circ\phi_A:K\sr{\phi_A}{\ral}I\sr{f}{\ral}X$} is eventually in $A$.
\end{lmm}
\begin{proof}
(i) If $f:I\ra X$ is eventually in $A$, i.e., some $f([i,I])\subset A$, and $f\circ\phi\subset f$ is a subnet, then with some $k\in K$ such that $\phi([k,K])\subset[i,I]$, we have $f\circ\phi([k,K])\subset f([i,I])\subset A$). The converse is clear since $f\subset f$.

(ii) If $f:I\ra X$ is frequently in $A$, then $f^{-1}(A)\subset I$ is cofinal, and so the associated subnet $f|_{K:=f^{-1}(A)}:K\hookrightarrow I\sr{f}{\ral} X$ of $f$ is eventually in $A$. Conversely, if some subnet $f\circ\phi\subset f$ is eventually in $A$, i.e., some $f\circ\phi([k,K])\subset A$, then for each $i\in I$, with some $k_i\in K$ such that $\phi(k_i)\geq i,\phi(k)$ and $\phi([k_i,K])\subset[i,I]$, we get $f([i,I])\cap A\neq\emptyset$.
\end{proof}

\begin{lmm}[\blue{\cite[Lemma 5, p.70]{kelley1975}}]\label{ClustSubLmm}
Let $X$ be a space, $f:I\ra X$ a net, and $\A\subset\P(X)$ a family of subsets such that (i) $f$ is frequently in each $A\in\A$, and (ii) $A_1,A_2\in\A$ implies $A_1\cap A_2\in\A$. Then there is a subnet {\footnotesize $g=f\circ\phi:J\sr{\phi}{\ral}I\sr{f}{\ral}X$} that is eventually in each $A\in\A$.
\end{lmm}
\begin{proof}
By hypotheses, $(A,\leq):=(\A,\supset)$ is a directed set. Let $K:=\bigcup_{A\in\A}f^{-1}(A)\times\{A\}=\{(i,A):i\in I,~A\in\A,~f(i)\in A\}\subset I\times\A$ (as a directed set wrt $(i,A)\leq (j,B)$ iff $i\leq j$ and $A\supset B$). (\blue{footnote}\footnote{Recall that a finite product of directed sets is a directed set.}). Then $K\subset I\times\A$ is directed, because if $(i,A),(j,B)\in K$, then $f(i)\in A$, $f(j)\in B$, and so with $C:=A\cap B\in\A$, $t_1\geq\{i,j\}$, and $t\in [t_1,I]\cap f^{-1}(C)\neq\emptyset$, we have $f(t)\in C$, which implies $(t,C)\geq \{(i,A),(j,B)\}$.

Define $\phi:K\ra I$ by $\phi(i,A):=i$. Then $\phi$ is cofinal, since for any $i\in I$, we can pick any $(i,A)\in K$ (which satisfies $i\leq i=\phi(i,A)$), and for any $(i',A')\in K$ satisfying $(i,A)\leq (i',A')$, i.e., $i\leq i'$ and $A\supset A'$, we have $i\leq i'=\phi(i',A')$. Thus, ~{\small $f\circ\phi:K\sr{\phi}{\ral}I\sr{f}{\ral}X$}~ is a subnet of $f$.

To show $f\circ \phi$ is eventually in every member of $\A$, fix $A\in\A$. Pick $i\in I$ such that $f(i)\in A$, i.e., $(i,A)\in K$. Then for any $(j,B)\in K$ satisfying $(i,A)\leq (j,B)$, i.e., $i\leq j$ and $A\supset B$, we have {\small $f\circ\phi(j,B)=f(j)\in B\subset A$}. That is, $f\circ\phi\big([(i,A)\leq K]\big)\subset A$.
\end{proof}

\begin{thm}[\blue{\cite[Theorem 6, p.71]{kelley1975}}]
Let $X$ be a space, $f:I\ra X$ a net, and $x_0\in X$. Then $x_0$ is a cluster point of $f$ $\iff$ a subnet of $f$ converges to $x_0$.
\end{thm}
\begin{proof}
($\Ra$): Assume $x_0\in X$ is a cluster point of $f$. Then the set of neighborhoods $\A$ of $x_0$ satisfies the hypothesis of Lemma \ref{ClustSubLmm}, and so a subnet of $f$ converges to $x_0$.

($\La$): Assume $x_0\in X$ is not a cluster point of $f$. Then there is a neighborhood $N(x_0)$ of $x_0$ such that $f$ is not frequently in $N(x_0)$. This means $f$ is eventually in $X-N(x_0)$. It follows that every subnet of $f$ is eventually in $X-N(x_0)$, and so cannot converge to $x_0$.
\end{proof}

\begin{crl}[\textcolor{blue}{Compactness criterion}]\label{NetCompactCrl}
Let $X$ be a space. Then $X$ is compact $\iff$ every net in $X$ has a convergent subnet (equivalently, every net in $X$ has a cluster point).
\end{crl}
\begin{proof}
($\Ra$): Assume $X$ is compact. Suppose $f:I\ra X$ is a net with no cluster point, i.e., for each $x\in X$, there exists a neighborhood $N(x)$ of $x$ such that $f$ is not frequently in $N(x)$, i.e., $f$ is eventually in $X-N(x)$. Since $X$ is compact, we have a finite cover $X\subset\bigcup_{i=1}^nN_i(x_i)$. This means $\bigcap_{i=1}^n(X-N_i(x_i)\big)=\emptyset$. But, since $f$ is eventually in each $X-N_i(x_i)$, we also have $\bigcap_{i=1}^n(X-N_i(x_i)\big)\neq\emptyset$ (a contradiction).

($\La$): Assume $X$ is not compact. Then $X$ has an open cover $\{U_\al\}_{\al\in A}$ with no finite subcover. Let $I:=\{\txt{finite}~F\subset A\}$, which is a directed set under inclusion $\subset$ (i.e., $F_1\leq F_2$ iff $F_1\subset F_2$). Let $f:I\ra X$, $F\mapsto f(F)\in X\backslash\bigcup_{\al\in F}U_\al$ (which is well defined by hypotheses). Suppose $f$ has a cluster point $x_0\in U_{\al_0}$ (for some $\al_0\in A$). Then for any $F\in I$ there is $F'\geq F$ (i.e., $F'\supset F$) such that $f(F')\in U_{\al_0}\cap\big(X\backslash\bigcup_{\al\in F'}U_\al\big)$. In particular, with $F=\{\al_0\}$, any choice $F'\supset F=\{\al_0\}$ gives $f(F')\in\emptyset$ (a contradiction).
\end{proof}

\begin{thm}[\blue{\index{Tychonoff's theorem}{Tychonoff's theorem}: \cite{chernoff-1992}}]\label{TychThm}
If $\{X_\al\}_{\al\in A}$ are compact spaces, then so is $\prod_\al X_\al$.
\end{thm}
\begin{proof}
Let $\{(X_{\al},\T_{\al})\}_{\al\in A}$ be a family of nonempty compact spaces, and  for any $B\subset A$, let $X_B:=\prod_{\al\in B}X_{\al}=\{(x_{\al})_{\al\in B}:x_{\al}\in X_{\al}\}$ be given the product topology. Recall that a base-neighborhood of a point $x\in X_B$ has the form $N_F(x)$ for some finite set $F\subset B$, with
{\footnotesize\bea
\label{nbd-supp-eqn}\textstyle N_F(x)=\{y\in X_B:y_{\al}\in N(x_{\al})\in\T_{\al}~~\txt{for}~\al\in F\}=\prod_{\al\in F}N(x_{\al})\times\prod_{\beta\not\in F}X_\beta=\bigcap_{\al\in F}N_{\al}(x),
\eea}
where {\small $N_{\al}(x):=\left\{y\in X:y_{\al}\in N(x_{\al})\in\T_{\al}\right\}=N(x_{\al})\times\prod_{\beta\neq \al}X_\beta$} is a ``strip through $N(x_{\al})$''.

Let $f:I\ra X_A,~\al\mapsto f_\al$ be a net. We need to show $f$ has a cluster point. A \ul{partial cluster point} of $f$ is a cluster point of the net {\footnotesize \bt[column sep=small] f_B:I\ar[r,"f"]& X_A\ar[r,two heads,"p^B"]&X_B\et}, for some $B\subset A$ (where $p^B(x)=p^B\big((x_\al)_{\al\in A}\big):=(x_\al)_{\al\in B}:=x|_B$ is the projection). Let $\P:=\{\txt{all partial cluster points of $f$}\}$ as a poset with the ordering
\bea
\textstyle x_B\in X_B\leq x_{B'}\in X_{B'}~~~\txt{if}~~~B\subset B'~~\txt{and}~~x_{B'}\big|_B=x_B~~~~(\txt{recall that}~~ B\sr{x_B}{\ral}\bigcup_{\al\in B}X_{\al}).\nn
\eea
$\P$ is nonempty because for any $\al\in A$, with $B:=\{\al\}$, the net \bt[column sep=small] f_{\{\al\}}:I\ar[r,"f"]& X_A\ar[r,two heads,"p^{\{\al\}}"]&X_{\{\al\}}\et has a cluster point since $X_\al\cong X_{\{\al\}}$ is compact (Corollary \ref{NetCompactCrl}). Let {\footnotesize$L=\big\{x^{(\ld)}_{B_\ld}\in X_{B_\ld}:\ld\in\Ld\big\}$} be a chain (linearly ordered set) in $\P$, where $x^{(\ld)}_{B_\ld}$ is a cluster point of {\footnotesize \bt f_{B_\ld}:I\ar[r,"f"]&X_A\ar[r,two heads,"p^{B_\ld}"]& X_{B_\ld}\et} (i.e., $f_{B_\ld}$ is frequently in every neighborhood of $x^{(\ld)}_{B_\ld}$).

Define {\footnotesize$z:=\bigcup_{\ld\in\Ld}x^{(\ld)}_{B_\ld}:=\bigcup_{\ld\in \Ld}\big(x^{(\ld)}_{\al}:\al\in B_\ld\big)\in X_{\bigcup B_\ld}$} and {\footnotesize\bt f_{\bigcup B_\ld}:I\ar[r,"f"]&X_A\ar[r,two heads,"p^{\bigcup B_\ld}"]&X_{\bigcup B_\ld}\et}, where
\bit[leftmargin=0.8cm]
\item[(i)] {\footnotesize
$f_{\bigcup B_\ld}(i)=p^{\bigcup B_\ld}\circ f(i)=\bigcup_\ld p^{B_\ld}\circ f(i)=\bigcup_\ld f_{B_\ld}(i)~~\Ra~~f_{\bigcup B_\ld}=\bigcup f_{B_\ld}$},
\item[(ii)] $x\in X_{\bigcup B_\ld}$ $\iff$ $x|_{B_\ld}\in X_{B_\ld}$ for all $\ld$, and
\item[(iii)] with the subspace ``imbedding'' $X_{B_\ld}\subset X_{\bigcup B_\ld}$, for any $x\in X_{\bigcup B_\ld}$ each base-neighborhood {\footnotesize $N(x)=\bigcup_\ld\big(N(x)\cap X_{B_\ld}\big)$} takes the form {\footnotesize $N(x)=\prod_{\al\in F}N(x_{\al})\times\prod_{\beta\in\left(\bigcup B_{\ld'}\right)\backslash F}X_\beta= \big(N(x)\cap X_{B_\ld}\big)\times\prod_{\beta\in\left(\bigcup B_{\ld'}\right)\backslash B_\ld}X_\beta$} for some $\ld$ such that $F\subset B_\ld$. It follows that $z$ is a cluster point of $f_{\bigcup B_\ld}$, otherwise, if $f_{\bigcup B_\ld}=\bigcup f_{B_\ld}$ is not frequently in some base-nbd ~{\footnotesize $N\big(\bigcup x^{(\ld)}_{B_\ld}\big)=\prod_{\al\in F}N(z_{\al})\times\prod_{\beta\in\left(\bigcup B_\ld\right)\backslash F}X_\beta$}~ of ~$z=\bigcup x^{(\ld)}_{B_\ld}$,
\eit
\bit[leftmargin=0.5cm]
\item i.e., $f_{\bigcup B_\ld}=\bigcup f_{B_\ld}$ is eventually in ~{\footnotesize $X_{\bigcup B_\ld}-N\big(\bigcup x^{(\ld)}_{B_\ld}\big)=\left[X_F-\prod_{\al\in F}N(z_{\al})\right]\times\prod_{\beta\in\left(\bigcup B_\ld\right)\backslash F}X_\beta$},
\item i.e., {\footnotesize $f_{\bigcup B_\ld}\left([i_0,I]\right)\subset X_{\bigcup B_\ld}-N\big(\bigcup x^{(\ld)}_{B_\ld}\big)=\left[X_F-\prod_{\al\in F}N(z_{\al})\right]\times\prod_{\beta\in\left(\bigcup B_\ld\right)\backslash F}X_\beta$},
\eit
then ~{\footnotesize $f_{B_\ld}\left([i_0,I]\right)\subset \left[X_F-\prod_{\al\in F}N(z_{\al})\right]\times\prod_{\beta\in B_\ld\backslash F}X_\beta=X_{B_\ld}-N\big(x^{(\ld)}_{B_\ld}\big)$}~ for some $\ld$ such that $F\subset B_\ld$. That is, $f_{B_\ld}$ is not frequently in some nbd $N\big(x^{(\ld)}_{B_\ld}\big)$ of $x^{(\ld)}_{B_\ld}$, which is a contradiction.

(Note: It is enough to use base-neighborhoods as we have done, because a net that is not frequently in a given neighborhood $N(x)=\bigcup N_b(x)$ is also not frequently in some base-neighborhood $N_b(x)$.)

Thus, $L$ has an upper bound {\small$\bigcup_{\ld\in\Ld}x^{(\ld)}_{B_\ld}\in \P$}, and so by Zorn's lemma, $\P$ has a maximal element $x_B\in X_B$ for some $B\subset A$. It remains to show that $B=A$. Suppose $B\subsetneq A$. Then some $c\in A\backslash B$. We know that $x_B$ is a cluster point of {\footnotesize \bt[column sep=small]f_B:I\ar[r,"f"]&X_A\ar[r,two heads,"p^B"]& X_B\et}. Also, since $X_c$ is compact and nonempty, the net {\footnotesize \bt[column sep=small]f_{\{c\}}:I\ar[r,"f"]&X_A\ar[r,"p^{\{c\}}"]&X_c\et} has a cluster point $x_c\in X_c$. Consider point $x_{B\cup\{c\}}\in X_{B\cup\{c\}}$ given by
{\small\[
x_{B\cup\{c\}}:=(x_B,x_c):B\cup\{c\}\ra X_{B\cup\{c\}},~~~~x_{B\cup\{c\}}(\al):=\left\{
                                                                    \begin{array}{ll}
                                                                      x_B(\al), & \al\in B \\
                                                                      x_c, & \al=c
                                                                    \end{array}
                                                                  \right\}.
\]}The point $x_{B\cup\{c\}}$ is a cluster point of {\footnotesize\bt f_{B\cup\{c\}}:I\ar[r,"f"]&X_A\ar[r,"{p^{B\cup\{c\}}}"]&X_{B\cup\{c\}}\et}, because if $f_B([i,I])\cap N(x_B)\neq\emptyset$ and $f_{\{c\}}([i,I])\cap N(x_c)\neq\emptyset$, then with $N(x_B\cup\{c\}):=N(x_B)\times N(x_c)$, we also have
\begin{align}
&f_{B\cup\{c\}}([i,I])\cap N(x_{B\cup\{c\}})=\big(f_B([i,I])\times f_{\{c\}}([i,I])\big)\cap\big(N(x_B)\times N(x_c)\big)\nn\\
&~~~~=\big(f_B([i,I])\cap N(x_B)\big)\times\big(f_{\{c\}}([i,I])\cap N(x_c)\big)\neq\emptyset.\nn
\end{align}
That is, $x_B<x_{B\cup\{c\}}\in\P$, which is a contradiction.
\end{proof}

\section{Linear Algebra of Topological Modules (a concretizing digression)}\label{LinAlgModSec2}
Let $M$ be a topological $R$-module. From section \ref{LinAlgModSec1}, we introduced the notions of (i) the \emph{span} of a subset of $M$, (ii) the \emph{linear independence} of a subset of $M$, (iii) a \emph{basis} for $M$, and so on. To further study linear algebra on $M$, we will introduce a closure-based topological notion of a module basis (namely, a \emph{density basis}), in addition to the earlier notion of a module basis (as a linearly independent spanning set).

\subsection{Dense spanning and density-bases}
\begin{dfn}[\textcolor{blue}{
\index{Span-dense set}{Span-dense set},
\index{Density-basis}{Density-basis for a module},
\index{Density-free module}{Density-free module}}]
Let $M$ be a topological $R$-module. A subset $D\subset M$ is \ul{span-dense} if $\Span_RD:=RD$ is dense in $M$, i.e., $\ol{\Span_RD}=M$. Equivalently, $D\subset M$ is span-dense if for any \emph{nonempty} open set $U\subset M$, we have $U\cap\Span_RD\neq\emptyset$ (i.e., there exists a linear combination $\sum_{i=1}^nr_id_i\in U\cap\Span_RD$, where $r_i\in R$ and $d_i\in D$).

A subset $B\subset M$ is a \ul{density-basis} for $M$ if (i) $B$ is linearly independent and (ii) $B$ is span-dense in $M$. The topological $R$-module $M$ is \ul{density-free} if it has (contains, admits) a density-basis.
\end{dfn}

Note that in the literature, our ``base'' for a topology is also called a ``basis'' for the topology. However, we have a strict distinction between the two, because for us, a ``base'' \ul{union-wise} generates the topology of the module while a ``basis'' \ul{linearly} generates, or spans, (a dense subset of) the module. Until now, such a distinction did not seem necessary, but it is now clear that it was (and therefore is) strictly essential owing to the nonspecialized/heterogeneous nature of our discussion throughout.

\begin{rmk}[\textcolor{blue}{Existence of a density-basis}]
Let $M$ be a topological $R$-module. If $M$ has a basis (i.e., if $M$ is free), then it is clear that $M$ has a density-basis (since every basis is clearly a density-basis). In particular, every topological vector space has a density-basis.
\end{rmk}

\begin{lmm}\label{TopLinExtLmm}
Let $M,N$ be topological $R$-modules (each with a scale-invariant topology), $A\subset M$ a linearly independent subset, and $f:A\subset M\ra N$ a continuous map. If the topology of $N$ has a base consisting of sets of the following form (\blue{footnote}\footnote{We are not yet sure if this condition follows from the hypotheses so far (i.e., from scale-invariance of the topologies).}):
\[
\textstyle\big\{\sum_{n\in N} r_nV_n:V_n\subset N~\txt{open},~r_n\in R,~r_n=0~\txt{a.e.f.}\big\},
\]
then $f$ extends to a unique continuous $R$-linear map $h=h_f:\Span_RA\subset M\ra N$.
\end{lmm}
\begin{proof}
Let $h:\Span_RA\subset M\ra N$, $\sum_{a\in A} r_aa\mapsto\sum_{a\in A} r_af(a)$, which is (i) well-defined due to the linear independence of $A$ and (ii) satisfies $h|_A=f$. Since $f$ is continuous, for any $a\in A$ and any open nbd $V_{f(a)}\ni f(a)$, there is an open nbd $U_a\ni a$ (in $M$) such that $f(A\cap U_a)\subset V_{f(a)}$. Given $\sum_a r_aa\in \Span_RA$, let $U_{\sum r_aa}:=\sum_a r_aU_a$ and $V_{\sum r_af(a)}:=\sum_a r_aV_{f(a)}$. Then for any open nbd $V\ni\sum r_af(a)$, we have the open nbd $U_{\sum r_aa}\ni\sum r_aa$ such that
\begin{align}
&\textstyle h(h^{-1}(V)\cap U_{\sum r_aa})\subset h(h^{-1}(V))\cap h(RA\cap U_{\sum r_aa})=h(h^{-1}(V))\cap\sum_a r_af(A\cap U_a)\nn\\
&\textstyle~~~~\subset h(h^{-1}(V))\cap\sum_a r_aV_{f(a)}=h(h^{-1}(V))\cap V_{\sum r_af(a)}\subset V.\nn
\end{align}
It therefore follows from the hypotheses that $h$ is continuous.
\end{proof}

\begin{rmk}[\blue{A question on extendability}]
Let $M$ be a density-free $R$-module with density-basis $B$, and $N$ any topological $R$-module (and assume the topologies of $M,N$ satisfy the hypotheses of Lemma \ref{TopLinExtLmm}). Given a continuous map $f:B\subset M\ra N$, we know it has a unique continuous $R$-linear extension $h_f:\Span_RB\subset M\ra N,~\sum_{b\in B} r_bb\mapsto\sum_{b\in B} r_bf(b)$. Can we extend $f$ to a continuous $R$-homomorphism $F:M\ra N$? (\blue{footnote}\footnote{The same questions can also be asked for the topological $R$-module of maps $\T(M,N)$ whose topology (the topology of pointwise convergence) is based on that of $N$.})

Of course, if $N$ is injective (i.e., if $N$ = the injective envelope of $N$) as an $R$-module, then an $R$-linear extension from $\Span_RB$ is always possible (including when $B\subset M$ is just any linearly independent subset that is not necessarily dense in $M$), in which case we only need to check the possibility of selecting a continuous extension.

Also, if in addition we assume (i) $M$ is Cauchy, i.e., every convergent net in $M$ is Cauchy, (ii) $N$ is complete, and (iii) $h_f$ is uniformly continuous, then a continuous extension might be possible in the following form: Given $m\in M$, pick a net $d^{(m)}:I\ra\Span_RB\subset M,~i\mapsto d_i^{(m)}$, such that $d^{(m)}\ra m$. Then the net $h_f\circ d^{(m)}:I\sr{d^{(m)}}{\ral}\Span_RB\sr{h_f}{\ral} N,~i\mapsto h_f\big(d_i^{(m)}\big)$ is Cauchy in $N$ (since it is Cauchy in $\Span_RB$ and $h_f$ is uniformly continuous), and so converges. We can therefore consider a selection map
\bea
\textstyle F:M\ra N,~m\mapsto F(m)\in\lim h_f\circ d^{(m)}.\nn
\eea
\end{rmk}

\begin{questions}
Let $M,N$ be topological $R$-modules, $A\subset M$ a linearly independent subset, and $f:A\subset M\ra N$ a map. Further assume (i) $f$ is uniformly continuous, (ii) $M$ is Cauchy, and (iii) $N$ is complete.

(a) (When) Does $f$ extend to a uniformly continuous linear map $h_f:\Span_RA\subset M\ra N$?

(b) (When) Does $h_f$ extend to a uniformly continuous map $H_f:\ol{\Span_RA}\subset M\ra N$?

(c) When the answer to (b) above is positive and $N$ is Hausdorff, is the extension $H_f$ unique?
\end{questions}
Compare question (b) above with the case of metric spaces considered later in Lemma \ref{UniExtThm}.

\begin{questions}
 Let $M$ be a topological $R$-module and $X$ a set. (i) If $M$ is density-free, does it follow that $M^X$ is also density-free for all $X$? (ii) Conversely, if $M^X$ is density-free for all $X$, is $M$ also density-free?

As possible hints for these questions, consider the following:
{\flushleft (i)} Is this correct? If $B\subset M$ is a density-basis and $\wt{B}:=\{\txt{finite maps}~f:X\ra B\}$, (\blue{footnote}\footnote{Recall that a \ul{finite map} is a map whose image a finite set.}), then $\Span_R\wt{B}=\{\txt{finite maps}~f:X\ra \Span_RB\}$ is a density-basis for $M^X$, since base sets for $M^X$ are
\bea
V_1^{x_1}\cap\cdots\cap V_n^{x_n}=\{f\in M^X:f(x_i)\in V_i\},~~~~\txt{for open}~~V_i\subset M,\nn
\eea
which each contain a finite map $f:X\ra \Span_RB$.
{\flushleft (ii)} Consider the special case where $X=\{x_0\}$ is a single element set, so that $M^X=M^{\{x_0\}}\cong M$.
\end{questions}

\begin{questions}
Let $M$ be a topological $R$-module and $X$ a set. (i) If $M$ is complete, does it follow that $M^X$ is also complete for all $X$? (ii) Conversely, if $M^X$ is complete for all $X$, is $M$ also complete?

As possible hints for these questions, consider the following:
{\flushleft (i)} Is this correct? By construction a net $f=\{f_i\}_{i\in I}:I\ra M^X,~i\mapsto f_i$ in $M^X$ is \ul{Cauchy} (i.e., $f\ominus f:I\times I\ra M^X,~(i,j)\mapsto f_i-f_j$ converges to $0$) $\iff$ \ul{pointwise Cauchy} (i.e., $f\ominus f:I\times I\ra M^X\sr{e_x}{\ral}M,~(i,j)\mapsto f_i(x)-f_j(x)$ converges to $0$ for each $x\in X$), and so if $M$ is complete then so is $M^X$.
{\flushleft (ii)} Consider the special case where $X=\{x_0\}$ is a single element set, so that $M^X=M^{\{x_0\}}\cong M$.
\end{questions}

\subsection{Internal density-decomposition. Quotient density-splitting}
\begin{dfn}[\textcolor{blue}{
\index{Sum (or span) of submodules}{Sum (or span) of submodules},
\index{Density-spanning! family}{Density-spanning family},
\index{Density-free! family of submodules}{Density-free family of submodules},
\index{Density-internal direct sum}{Density-internal direct sum},
\index{Density-complemented submodule}{Density-complemented submodule}}]
Let $M$ be a topological $R$-module and $\F:=\{M_i:i\in I\}$ a family of submodules of $M$. Recall that the \ul{sum} (or \ul{span}) of $\F$ is the span of the union $\bigcup_{i\in I} M_i$ (i.e., the submodule generated by $\bigcup_{i\in I} M_i$): (\blue{footnote}\footnote{Given subsets $\{A_i\subset M\}_{i\in I}$, we can also consider a generalization of the usual finite sum of subsets as follows:
\[
\textstyle \sum_{i\in I} A_i:=\{a_{i_1}+\cdots+a_{i_n}:a_{i_k}\in A_{i_k},~i_1,...,i_k\in I,~n\in\Natural\}~=~\bigcup_{\txt{finite}~F\subset I}~\sum_{i\in F}A_i,
\]
where if we enumerate $F$ as ~$F=\{f_1,\cdots,f_{|F|}\}$,~ then ~$\sum_{i\in F}A_i=A_{f_1}+\cdots+A_{f_{|F|}}:=\{a_1+\cdots+a_{|F|}:a_t\in A_{f_t}\}$.})
\bea
\textstyle\Span_R\F=\sum_{i\in I}M_i:=\Span_R(\bigcup_{i\in I}M_i)=\left\langle\bigcup_{i\in I}M_i\right\rangle=\left\{m_1+\cdots+m_n~|~m_j\in\bigcup M_i,~n\in\Natural\right\}.\nn
\eea
The family $\F$ is a \ul{density-spanning family} if $\ol{\sum_{i\in I} M_i}=M$. The family $\F$ is a \ul{density-free family} if $\ol{M_i}\cap\ol{\left(\sum_{j\neq i}M_j\right)}=0$ for all $i\in I$ (in which case, the sum's closure $\ol{\sum_{i\in I}M_i}$ is a \ul{density-internal direct sum}, and written as $\ol{\bigoplus}_{i\in I}M_i$).

A submodule $N\subset M$ is \ul{density-complemented} in $M$ if there exists a submodule $N'\subset M$ such that $N\ol{\oplus}N'=M$ (i.e., such that $\ol{N}\cap\ol{N'}=0$ and $\ol{N+N'}=M$).

\end{dfn}
It is clear that a topological $R$-module is a density-internal direct sum iff it has a density-spanning density-free family of submodules.

\begin{rmk}
Let $M$ a topological $R$-module, $\{A_i\subset M\}_{i\in I}$ subsets, and $\{M_i\subset M\}_{i\in I}$ submodules. Since $\sum_iRA_i\subset \sum_iR\ol{A_i}\subset\ol{\sum_iRA_i}$, it follows that
\[
\textstyle \ol{\sum_iRA_i}=\ol{\sum_iR\ol{A_i}}~~~~\txt{and}~~~~\ol{\sum_iM_i}=\ol{\sum_iR\ol{M_i}}.\nn
\]
\end{rmk}

\begin{dfn}[\textcolor{blue}{\index{Topological! quotient group}{Topological quotient group}, \index{Topological! quotient ring}{Topological quotient ring}, \index{Topological! quotient module}{Topological quotient module}, \index{Topological! quotient algebra}{Topological quotient algebra}}]
Let $N\subset M$ be topological $R$-modules. We will give the quotient module $M/N$ the quotient topology induced by the natural $R$-homomorphism $\pi:M\ra M/N,~m\mapsto m+N$, i.e., a set $V\subset M/N$ is open iff the preimage $\pi^{-1}(V)\subset M$ is open.
\end{dfn}

A topological quotient group, topological quotient ring, and topological quotient algebra are defined similarly (i.e., they are also given the quotient topologies induced by their respective natural maps).

\begin{lmm}[\textcolor{blue}{\index{Quotient! map (openness)}{Openness of the quotient map}}]
Let $N\subset M$ be topological $R$-modules. The quotient map $\pi:M\ra M/N,~m\mapsto m+N$ is open.
\end{lmm}
\begin{proof}
By the definitions, $\pi$ is open $\iff$ $\pi(U)\subset M/N$ is open (for every open set $U\subset M$), $\iff$ $\pi^{-1}\pi(U)\subset M$ is open (for every open set $U\subset M$), where
\bea
&&\textstyle\pi^{-1}\pi(U)=\pi^{-1}\left(U+N\over N\right)=\{m:m+N\in (U+N)/N\}=\{m:m+N=u+N,~u\in U\}\nn\\
&&\textstyle~~~~=\{m:m=u+n,~u\in U,n\in N\}=U+N=\bigcup_{n\in N}(n+U).\nn
\eea
Recall: By definition, the topology of $M$ is such that $+:M\times M\ra M$ is continuous, and so the bijective maps $f_n=+|_{M\times\{n\}}:M\ra M,~m\mapsto m+n$ are continuous, with continuous inverses $f_n^{-1}=f_{-n}:M\ra M$, and so are homeomorphisms. Hence, for each $n\in N$, the image $f_n(U)=n+U$ is an open set.
\end{proof}

We can show similarly that if $H\subset G$ are topological quotient groups, with $H$ a normal subgroup of $G$, then the quotient map $\pi:G\ra G/H$ is open.

\begin{lmm}
Let $X,Y$ be spaces and $f:X\ra Y$ a surjective continuous map. If $A\subset X$ is a dense set (i.e., $\ol{A}=X$), then so is $f(A)\subset Y$.
\end{lmm}
\begin{proof}
Let $A\subset X$ be dense. Pick any nonempty open set $V\subset Y$. Since $V$ is nonempty (and $f$ is surjective), so is the open set $f^{-1}(V)$. Therefore $f^{-1}(V)\cap A\neq\emptyset$. Hence $V\cap f(A)\neq\emptyset$, since
\[
f\big(f^{-1}(V)\cap A\big)\subset f\big(f^{-1}(V)\big)\cap f(A)=V\cap f(A). \qedhere
\]
\end{proof}
\begin{lmm}
Let $X$ be a space, $A\subset X$ a dense set, and $O\subset X$ an open set. Then $A\cap O$ is dense in $O$.
\end{lmm}
\begin{proof}
Recall that a relatively open set is an open set: Let $U\subset O$ be open in $O$, i.e., $U=O_U\cap O$ for an open subset $O_U\subset X$. Since $O$ is open in $X$, so is $U=O_U\cap O\subset X$. Therefore, if $U\neq\emptyset$ then $U\cap(A\cap O)=U\cap A\neq\emptyset$.
\end{proof}

\begin{prp}
Let $N\subset M$ be topological $R$-modules and $\pi:M\ra{M\over N}$ the natural map.
\bit[leftmargin=0.9cm]
\item[(i)] If $D\subset M$ is a span-dense set, then so is ${D+N\over N}\subset{M\over N}$.
\item[(ii)] If $D\subset M$ is a span-dense set, then for any open nbd $O_N\supset N$, ~$(\Span_RD)\cap O_N$ is dense in $O_N$.
\item[(iii)] If $N$ is closed in $M$, $A\subset N$ is span-dense, and $C\subset {M\over N}$ is span-dense, then $A\cup(\pi^{-1}(C)\backslash N)\subset M$ is also span-dense. Moreover, if $C$ is linearly independent (i.e., $C$ is a density-basis in ${M\over N}$), then
\[
\textstyle (\Span_RA)\cap\big(\Span_R(\pi^{-1}(C)\backslash N)\big)=\{0\}. \qedhere
\]
\eit
\end{prp}
\begin{proof}
{\flushleft (i)} Recall that $\pi:M\ra {M\over N}$ is a surjective continuous map. Therefore, if $\Span_RD\subset M$ is dense, then (by the preceding lemmas) so is $\pi(\Span_RD)=\Span_R\pi(D)=\Span_R{D+N\over N}\subset {M\over N}$.

{\flushleft (ii)} If $\Span_RD\subset M$ is dense and $O_N\supset N$ is an open nbd in $M$, then by the preceding lemmas,
\[
O_N=\Cl_{O_N}(O_N\cap\Span_RD)\subset \Cl_M(O_N\cap\Span_RD)=\ol{O_N\cap\Span_RD}.
\]

{\flushleft (iii)} Without loss of generality, assume $N\not\in C\subset{M\over N}$. Let {\small $B:=\pi^{-1}(C)\backslash N$}. Then
\begin{align}
\textstyle C=C\backslash\{N\}=\pi\big(\pi^{-1}(C\backslash\{N\})\big)=\pi\big(\pi^{-1}(C)\backslash\pi^{-1}(\{N\})\big)=\pi\big(\pi^{-1}(C)\backslash N\big)=\pi(B)={B+N\over N}.\nn
\end{align}
Since $N\subset\ol{\Span_RA}$ by hypotheses, it suffices to show $M\backslash N\subset\ol{\Span_R(A\cup B)}$.

Let $U\subset M\backslash N$ be a nonempty open set. Then because $\pi$ is open, {\footnotesize $V:=\pi(U)={U+N\over N}\subset M/N$} is also a nonempty open set. Thus, some linear combination {\footnotesize $\sum r_i(b_i+N)\in V\cap\Span_R{B+N\over N}$}, i.e.,
\bea
&&\textstyle\sum r_ib_i+N\in {U+N\over N}\cap\Span_R{B+N\over N},~~\Ra~~\sum r_ib_i+N=u+N,~~\txt{for some}~~b_i\in B,~~u\in U,\nn\\
&&\textstyle~~\Ra~~\sum r_ib_i-u\in N,~~\Ra~~\sum r_ib_i-u=n\in N=\ol{\Span_RA},\nn\\
&&\textstyle~~\Ra~~u=-n+\sum r_ib_i~\in~U\cap\left(N+\Span_RB\right)\subset U\cap\ol{\left(\Span_RA+\Span_RB\right)},\nn
\eea
where $\Span_RA+\Span_RB=\Span_R(A\cup B)$. Thus, every nonempty open set $U\subset M\backslash N$ intersects {\small $\ol{\Span_R(A\cup B)}$}, and so also intersects $\Span_R(A\cup B)$, since ~{\small $\ol{\big(\ol{\Span_R(A\cup B)}\big)}=\ol{\Span_R(A\cup B)}$}.

For the ``moreover'' part, recall that $C=\pi(B)={B+N\over N}$, and so
\[
\textstyle \Span_RC=\Span_R\pi(B)=\pi(\Span_RB)={\Span_RB+N\over N},~~\Ra~~\pi^{-1}(\Span_RC)=\Span_RB+N.\nn
\]
If $C$ is linearly independent, then for any $\{r_i\}_{i=1}^n\subset R$ and $\{c_i=b_i+N\}_{i=1}^n\subset C$ (where $b_i\not\in N$), we have
\[
\textstyle 0_{M\over N}=\sum_i r_ic_i=\sum_i r_ib_i+N=N~~\left(\iff~~\sum_i r_ib_i\in N\right)~~\iff~~r_1=\cdots=r_n=0,
\]
which shows ~$\Span_RB~\cap~\ol{\Span_RA}=\Span_RB~\cap~N=\{0\}$.
\end{proof}

\begin{crl}
Let $N\subset M$ be topological $R$-modules, with $N$ closed in $M$. If $N$ and $M/N$ are density-free, then so is $M$.
\end{crl}

\begin{question}
Let $N\subset M$ be topological $R$-modules, with $N$ closed in $M$. If $N$ and $M/N$ are density-free (so $M$ is also density-free) does it follows that $N$ is density-complemented in $M$, i.e., does there exist a submodule $N'\subset M$ such $N'\cong M/N$ and $M=N\ol{\oplus}N'$, where $\cong$ denotes $R$-homeomorphism (i.e., $R$-linear homeomorphism)?
\end{question}

\begin{rmk}[\blue{Recall}]
Let $N\subset M$ be $R$-modules. If $M/N$ is free (hence projective), then the SES
\bea
0\ra N\ra M\ra M/N\ra 0\nn
\eea
splits, giving an $R$-isomorphism $M\cong N\oplus M/N$. If $N$ is also free (i.e., both $N$ and $M/N$ are free) then $M$ is free as a direct sum of free $R$-modules.
\end{rmk}

\begin{question}
Let $N\subset M$ be topological $R$-modules. If $M$ is complete, does it follow that $M/N$ is complete?
\end{question}

\subsection{Series convergence. Series approximation} ~\\~
Given a density-free $R$-module $M$ with density-basis $B$, because $M=\ol{\Span_RB}$, we can think of $\Span_RB$ as an ``approximation'' of $M$, written $\Span_RB\approx M$. How can we describe such an approximation pointwise, i.e., given $m\in M$, which collection of elements of $\Span_RB$ best approximates $m$, and how? In the following discussion, we indicate how to obtain answers to this question using topological convergence.
\begin{dfn}[\textcolor{blue}{
\index{Series associated with a map}{Series associated with a map},
\index{Convergent! series}{Convergent series},
\index{Sum of a series}{Sum of a series}}]
Let $A$ be a set, $M$ a topological $R$-module, and $x:A\ra M,~a\mapsto x_a$ any map. With $FS(A):=\{\txt{finite subsets}~F\subset A\}$ as a directed set wrt inclusion $\subset$, the \ul{series} in $M$ associated with $x$ is the net $S=S_x:FS(A)\ra M,~F\mapsto\sum_{a\in F}x_a$. If $m\in M$, we say $S$ \ul{converges} to $m$ (making $S$ a \ul{convergent series} and $m$ a \ul{sum} of $S$), written $\sum_{a\in A} x_a=m$, if $S\ra m$ (i.e., $m$ is a limit of the net $S$).
\end{dfn}

\begin{dfn}[\textcolor{blue}{
\index{Fourier! series (expansion/approximation)}{Fourier series (expansion/approximation)},
\index{Fourier! coefficients}{Fourier coefficients},
\index{Parsevalian density-basis}{Parsevalian density-basis},
\index{Orthogonal set}{Orthogonal set},
\index{Orthonormal set}{Orthonormal set},
\index{Orthonormal basis}{Orthonormal basis},
\index{Parseval's identity}{Parseval's identity},
\index{Bessel's inequality}{Bessel's inequality},
\index{Fourier! coefficient function}{Fourier coefficient function},
\index{Fourier! integral}{Fourier integral}}]
Let $M$ be a density-free $R$-module with density-basis $B$, and let $m\in M$. Let $\O_m:=\{\txt{open sets $O\subset M$ containing $m$}\}$ as a directed set wrt inclusion $\subset$. For each $O\in\O_m$, select {\footnotesize $\sum_{F^O\subset B}r^O_bb\in O\cap\Span_RB$} (for a finite set $F^O\subset B$) to obtain a net
\bea
\textstyle S=S_m:\O_m\ra M,~O\mapsto \sum_{b\in F^O\subset B}r^O_bb~\in~ O\cap\Span_RB.\nn
\eea
Any limit of such a net is a \ul{Fourier series (expansion/approximation)} of $m$ with respect to $B$, and written ~$S_m\ra\sum_{b\in B}r_{m,b}b$ ~or just ~$\sum_{b\in B}r_{m,b}b$. The numbers $r_{m,b}$ are called \ul{Fourier coefficients} of $m$ wrt $B$. The density-basis $B$ is \ul{Parsevalian} if $\sum_{b\in B}r_{m,b}b=m$ for all $m\in M$. (See also Example \ref{HilbFourExm} below.)
\end{dfn}

By Lemma \ref{NetHausLmm} on page \pageref{NetHausLmm}, the Fourier series (whenever it exists) is unique if $M$ is a Hausdorff space. 

In practical applications, to see how good the series approximation is, one needs to understand (those optimal net selection) conditions under which either (i) $\sum_{b\in B}r_{m,b}b=m$, i.e., $S_m\ra m$, or (ii) $\sum_{b\in B}r_{m,b}b$ is as close (in a precise sense) to $m$ as possible. Typically, based on special/classical $R$-valued selection/optimization functions, the optimal or best possible ``\ul{Fourier coefficient function}'' $M\times M\ra R,~(m,b)\mapsto r_{m,b}$ is $R$-bilinear (i.e., linear in both $m$ and $b$ separately) up to very slight alterations depending on the type/purpose of application. In particular, if $M=\M(X,W)$ is the topological $R$-module of measurable maps (from a set $X$ to a topological $R$-module $W$), then such an optimal Fourier coefficient function $M\times M\ra R,~(m,b)\mapsto r_{m,b}$ would behave like an $R$-valued integral (called \ul{Fourier integral}) in each variable.

Using further details not of immediate importance to us, it is straightforward to verify the familiar basic facts stated in the following example.

\begin{example}\label{HilbFourExm}
In a Hilbert space (Definition \ref{InnProdSpDfn} on page \pageref{InnProdSpDfn}), $\H=(\H,\langle\rangle)$, a set $A\subset\H$ is called \ul{orthogonal} if $\langle a,a'\rangle=0$ for all $a\neq a'$ in $A$, and a set $A\subset\H$ is called \ul{orthonormal} if it is orthogonal and $\langle a,a\rangle=1$ for all $a\in A$. Here, an example of a Parsevalian density-basis is given by an orthonormal span-dense set (called an \ul{orthonormal basis}), i.e., an orthonormal set $B\subset\H$ such that $\ol{Span_{\mathbb{F}}B}=\H$. In this case, for any $x\in\H$, the relation $\sum_{b\in B}r_{x,b}b=x$ satisfies $r_{x,b}=\langle x,b\rangle$, and gives the equality
\bea
\textstyle \sum_{b\in B}|\langle x,b\rangle|^2=|\langle x,x\rangle|^2,\nn
\eea
which is called \ul{Parseval's identity} or \ul{Parseval's theorem}. Meanwhile, for any orthonormal set $A\subset\H$,
\bea
\textstyle \sum_{a\in A}|\langle x,a\rangle|^2\leq |\langle x,x\rangle|^2,~~\txt{for all}~~x\in\H,\nn
\eea
which is called \ul{Bessel's inequality}, where equality holds for all $x\in\H$ if and only if $A$ is span-dense (besides orthonormal, hence linearly independent).
\end{example}

\section{Classical Spaces in Analysis: Topological Vector Spaces}

\begin{dfn}[\textcolor{blue}{Recall:
\index{Topological! multiplicative set}{Topological multiplicative set},
\index{Topological! group}{Topological group}}]
A multiplicative set $M=(M,\cdot)$ is a \ul{topological multiplicative set} if it is a topological space such that multiplication $\cdot:M\times M\ra M,~(m,m')\mapsto mm'$ is continuous. A group $G=(G,\cdot,()^{-1})$ is a \ul{topological group} if it is a topological space such that multiplication $G\times G\ra G,~(g,h)\mapsto gh$ and inversion $G\ra G,~g\mapsto g^{-1}$ are continuous.
\end{dfn}

\begin{dfn}[\textcolor{blue}{
\index{Classical! field}{Classical field},
\index{Topological! vector space (TVS)}{Topological vector space (TVS)
}}]
Let $\mathbb{F}$ denote a \ul{classical field}, which is either $\Real$ or $\Complex$, with its default topology. (\blue{footnote}\footnote{Here, the topology of $\Complex$ is the quotient topology induced by the natural bijection $\pi:\Real^2\ra\Complex,~(a,b)\mapsto a+ib$. For any open set $O\subset\Real^2$, $\pi^{-1}(\pi(O))=O$ is open, and so $\pi(O)$ is open, i.e., $\pi$ is an open map. It follows that $\pi$ is a homeomorphism, i.e.,
\bea
\Complex\cong\Real^2.\nn
\eea}). An $\mathbb{F}$-vector space $X=(X,+,\mathbb{F})$ is a \ul{topological vector space} (TVS) if it is a topological space such that addition $+:X\times X\ra X,~(x,y)\mapsto x+y$ and scalar multiplication $\cdot:\mathbb{F}\times X\ra X,~(\al,x)\mapsto \al x$ are continuous.
\end{dfn}

Recall from Corollary \ref{InvOfTopLmm} that the topology of a TVS is translation-invariant, scale-invariant, and additive. Examples of topological vector spaces are seminormed spaces in Definition \ref{VecSpExs} below.

\begin{dfn}[\textcolor{blue}{
\index{Associative! action}{Action of a topological associative set on a space},
\index{$M$-space}{$M$-space},
\index{$M$-set}{$M$-set}}]
Let $M$ be a topological associative set and $X$ a space. An \ul{action} of $M$ on $X$ (making $X$ an \ul{$M$-space}) is a \ul{continuous map} $M\times X\ra X$, $(m,x)\mapsto mx$ such that ~$m(m'x)=(mm')x$, ~for all $m,m'\in M$, $x\in X$.

Let $M$ be an associative set and $X$ a set. An \ul{action} of $M$ on $X$ (making $X$ an \ul{$M$-set}) is a map $M\times X\ra X$, $(m,x)\mapsto mx$ such that ~$m(m'x)=(mm')x$, ~for all $m,m'\in M$, $x\in X$.
\end{dfn}

\begin{dfn}[\textcolor{blue}{
\index{Group! action}{Action of a topological group on a space},
\index{$G$-space}{$G$-space},
\index{$G$-set}{$G$-set}}]
Let $G$ be a topological group (with identity $1_G$) and $X$ a space. An \ul{action} of $G$ on $X$ (making $X$ a \ul{$G$-space}) is a \ul{continuous map} $G\times X\ra X$, $(g,x)\mapsto gx$ such that ~(i) $g(g'x)=(gg')x$, and (ii) $1_Gx=x$, ~for all $g,g'\in G$, $x\in X$.

Let $G$ be a group and $X$ a set. An \ul{action} of $G$ on $X$ (making $X$ a \ul{$G$-set}) is a map $G\times X\ra X$, $(g,x)\mapsto gx$ such that ~(i) $g(g'x)=(gg')x$, and (ii) $1_Gx=x$, ~for all $g,g'\in G$, $x\in X$.
\end{dfn}

\begin{dfn}[\textcolor{blue}{
\index{Translation!}{Translation} in a space,
\index{Translation! invariant space}{Translation invariant space}}]
Let $X$ be a space and $V$ a TVS. A \ul{$V$-translation} in $X$ is an action of $V$ on $X$ of the form
\bea
(V,+)\times X\ra X,~~(v,x)\mapsto x+v,~~~~(x+v_1)+v_2=x+(v_1+v_2)\eqv x+v_1+v_2.\nn
\eea
If such a translation is specified (or exists in a canonical way), we say $X$ is a \ul{$V$-translation invariant space}.
\end{dfn}

\begin{dfn}[\textcolor{blue}{
\index{Scaling}{Scaling} in a space,
\index{Scale invariant space}{Scale invariant space}}]
Let $X$ be a space and $\mathbb{F}$ a classical field ($\Real$ or $\Complex$). An \ul{$\mathbb{F}$-scaling} in $X$ is an action of $\mathbb{F}$ on $X$ of the form
\bea
(\mathbb{F},\cdot)\times X\ra X,~~(v,x)\mapsto \al x,~~~~\al_1(\al_2 x)=(\al_1\al_2)x\eqv \al_1\al_2x.\nn
\eea
If such a scaling is specified (or exists in a canonical way), we say $X$ is an \ul{$\mathbb{F}$-scale invariant space}.
\end{dfn}

\begin{dfn}[\textcolor{blue}{
\index{Interval}{Interval},
\index{Convex set}{Convex set},
\index{Convex function}{Convex function},
\index{Concave function}{Concave function}}]\label{ConvConcFnDef}
Let $V$ be an $\mathbb{F}$-vector space. If $u,v\in V$, the \ul{interval} in $V$ from $u$ to $v$ is the set $[u,v]:=\{(1-t)u+tv:0\leq t\leq 1\}$. A set $K\subset V$ is a \ul{convex set} if for all $u,v\in K$, we have $[u,v]\subset V$. A function $f:A\subset V\ra\Real$ is a \ul{convex function} if (i) $A\subset V$ is a convex set and (ii) $f((1-t)u+tv)\leq (1-t)f(u)+tf(v)$ for all $0\leq t\leq 1$. A function $f:A\subset V\ra\Real$ is a \ul{concave function} if $-f$ is a convex function.
\end{dfn}

\begin{dfn}[\textcolor{blue}{
\index{Nondecreasing function}{}\index{Increasing function}{Increasing function},
\index{Nonincreasing function}{}\index{Decreasing function}{Decreasing function},
\index{Strictly! increasing function}{Strictly increasing function},
\index{Strictly! decreasing function}{Strictly decreasing function},
\index{Monotone function}{Monotone function},
\index{Strictly! monotone function}{Strictly monotone function}}]
A function $f:\Real\ra\Real$ is said to be
\bit
\item \ul{increasing} (or \ul{nondecreasing}) if $t<t'$ implies $f(t)\leq f(t')$ for all $t,t'$.
\item \ul{decreasing} (or \ul{nonincreasing}) if $-f$ is increasing.
\item \ul{strictly increasing} if $t<t'$ implies $f(t)<f(t')$ for all $t,t'$. ~~~~ (\blue{footnote}\footnote{Here, a strictly increasing function is just another name for an order-invariant, or order-preserving, map.})
\item \ul{strictly decreasing} if $-f$ is strictly increasing.
\item \ul{monotone} (or \ul{monotonic}) if it is either increasing or decreasing.
\item \ul{strictly monotone} (or \ul{strictly monotonic}) if it is either strictly increasing or strictly decreasing.
\eit
\end{dfn}

\begin{dfn}[\textcolor{blue}{
\index{Seminorm}{Seminorm},
\index{Seminormed space}{Seminormed space},
\index{Complete! set of seminorms}{Complete set of seminorms},
\index{Locally! convex space (LCS)}{Locally convex space (LCS)},
\index{Norm}{Norm},
\index{Normed space}{Normed space},
\index{Banach space}{Banach space}}]\label{VecSpExs}
Let $X$ be an $\mathbb{F}$-vector space ($\mathbb{F}=\Real$ or $\Complex$). A \ul{seminorm} on $X$ is a function $p:X\ra \Real$ such that $p(\al x+\beta y)\leq|\al|p(x)+|\beta|p(y)$ for all $x,y\in X$ and scalars $\al,\beta\in \mathbb{F}$ (or equivalently, $p(\al x)=|\al|p(x)$ and $p(x+y)\leq p(x)+p(y)$ for all $\al,x,y$). Associated with any seminorms $\{p_{a}\}_{{a}\in A}$ on $X$ is the topology (called \ul{seminorm topology}, making $X$ a \ul{seminormed space}) with \emph{subbase} (or generator) consisting of all ``strips'' of the form $U_{\vep,{a}}(x):=\{y\in X~|~p_{a}(x-y)<\vep\}$, $\vep>0$, ${a}\in A$, $x\in X$.

A \ul{locally convex space} is a seminormed space $(X,\{p_{a}\}_{{a}\in A})$ such that for any $x\in X$, if $p_{a}(x)=0$ for all ${a}\in A$ then $x=0$. (\blue{footnote}\footnote{In this case, we say $\{p_{a}\}_{{a}\in A}$ is a \ul{complete set of seminorms}, which is the case iff the seminorm topology is Hausdorff. Indeed, if $\{p_{a}\}_{a\in A}$ is complete, then $x\neq x'$ iff $x-x'\neq 0$ implies $\Ra$ $p_{a}(x-x')\neq 0$ for some ${a}$, and so
\bea
\{y:p_{a}(x-y)<p_{a}(x-x')/2\}\cap\{y:p_{a}(x'-y)<p_{a}(x-x')/2\}=\emptyset.\nn
\eea
Conversely, if $\{p_{a}\}_{a\in A}$ is not complete, then there exists $x\neq 0$ such that $p_{a}(x)=0$ for all ${a}$, and so $x$ and $0$ have no disjoint neighborhoods.}).

A \ul{norm} on $X$ (making $X$ a \ul{normed space}) is a seminorm $p$ such that $p(x)=0$ implies $x=0$, in which case we write $p(x)$ as $\|x\|$. Every normed space $(X,\|\|)$ is automatically considered to be a metric space with respect to the metric $d(x,y):=\|x-y\|$. A complete normed space is called a \ul{Banach space}.
\end{dfn}

Common examples of Banach spaces (due to the completeness of $(\mathbb{F},|~|)$ and Lemma \ref{CompUniMet} later) are given by the following.

\begin{examples}[\blue{
\index{Infinity-normed space}{Infinity-norm space},
\index{Infinity-norm}{Infinity-norm}}]
Let $A$ be a set. The \ul{infinity-normed space} associated with $A$ is the $\mathbb{F}$-subspace (of the $\mathbb{F}$-algebra of maps $\mathbb{F}^A$) given by
{\small\[
\textstyle \ell^\infty(A,\mathbb{F}):=\left\{x\in\mathbb{F}^A~\big|~\sup_{a\in A}|x(a)|<\infty\right\}\subset \mathbb{F}^A:=\{\txt{maps}~~f:A\ra \mathbb{F}\}
\]} as a normed space with respect to the norm ~$\|x\|:=\sup_{a\in A}|x(a)|$, called \ul{infinity-norm} or (\ul{$\infty$-norm}).
\end{examples}

\begin{rmks}
Let $X$ be an $\mathbb{F}$-vector space and $\{p_{a}\}_{{a}\in A}$ seminorms on $X$.
{\flushleft (i)} By the subbase-to-base criterion (i.e., base sets = finite intersections of subbase sets), the seminorm topology on $X$ associated with $\{p_{a}\}_{{a}\in A}$ has a \emph{base} consisting of $\emptyset,X$ and \emph{finite intersections} of the ``strips''
\begin{align}
&\textstyle U_{\vep,{a}_1,...,{a}_n}(x)=\bigcap_{i=1}^nU_{\vep,{a}_i}(x),~~~~~~n\geq 1,~~\vep>0,~~a_i\in A,~~x\in X,\nn\\
&\textstyle~~~~=\{y\in X:p_{{a}_1,...,{a}_n}(x-y)<\vep\},~~~~~~~~p_{{a}_1,...,{a}_n}:=\max\{p_{{a}_1},...,p_{{a}_n}\}.\nn
\end{align}
{\flushleft (ii)} If $p:X\ra\Real$ is any seminorm, the translation of $U_{\vep,p}(0):=\{p<\vep\}$ by a point $x\in X$ is
\bea
x+U_{\vep,p}(0)=x+\{z:p(z)<\vep\}=\{x+z:p(x+z-x)<\vep\}=\{y:p(y-x)<\vep\}=U_{\vep,p}(x).\nn
\eea
\end{rmks}

\begin{lmm}[\blue{\index{Minkowski function}{Minkowski function} of a convex set: \cite[Proposition 1.14, p.102]{conway}}]\label{MinkGauge}
Let $X$ be a locally convex space and $U\subset X$ a convex set that is ``\emph{symmetric}'' in the sense $\ld U=U$ for any $\ld\in\mathbb{F}$ with $|\ld|=1$. Then the map $p:X\ra\Real$, $p(x):=\inf\left\{t> 0:{x\over t}\in U\right\}=\inf\left\{t> 0:x\in tU\right\}$, is a seminorm.

(We call the map $p$ above the \ul{Minkowski function} of the convex set $U$.)
\end{lmm}
\begin{proof}
Let $x,y\in X$. For any scalar $\al\in\mathbb{F}$, $p(\al x)=\inf\left\{t>0:{\al x\over t}\in U\right\}$ = $|\al|\inf\big\{{t\over|\al|}>0:{|\al|e^{i\arg\al}x\over t}\in e^{i\arg\al}U\big\}$ = $|\al|\inf\big\{s>0:{x\over s}\in U\big\}=|\al|p(x)$. Also, for any $t,t'>0$ with ${x\over t},{y\over t'}\in U$, we can choose $\ld\in (0,1)$ such that (by the convexity of $U$) $\ld {x\over t}+(1-\ld){y\over t'}\in U$, along with ${\ld\over t}={1-\ld\over t'}$, imply ${1\over\ld}=1+{t'\over t}$ and ${\ld\over t}(x+y)={1\over t+t'}(x+y)\in U$, which in turn imply
\begin{align}
p(x)&+p(y)\textstyle=\inf\left\{t>0:{x\over t}\in U\right\}+\inf\left\{t>0:{y\over t}\in U\right\}\nn\\
&\textstyle\geq \inf\left\{t+t'>0:\ld{x\over t}+(1-\ld){y\over t'}\in U\right\}
 =  \inf\left\{t+t'>0:{x+y\over t+t'}\in U\right\}=p(x+y).\nn\qedhere
\end{align}
\end{proof}
Note (from the proof above) that if $U\subset X$ is any convex set, then (even without symmetry) the resulting Minkowski function $p:X\ra\Real$ is clearly a \index{Sublinear functional}{\ul{sublinear functional}} (see also \cite[Proposition 3.2,p.108]{conway}) in the sense it satisfies
\bea
p(\al x)=\al p(x)~~\txt{and}~~p(x+y)\leq p(x)+p(y),~~\txt{for all}~~\al\geq 0,~~x,y\in X.\nn
\eea

\begin{crl}\label{LCScrit}
Let $X$ be a TVS. Then $X$ is a seminormed space (e.g., LCS) $\iff$ $X$ has a base consisting of convex sets.
\end{crl}
\begin{proof}
($\Ra$) Seminorms produce convex neighborhoods: Seminorms are convex functions
\bea
\txt{( so seminorm$^{-1}$($[0,\vep)$) is convex )}\nn
\eea
and intersections of convex sets are convex sets.

($\La$) Conversely, convex neighborhoods come from seminorms: Let $U\subset X$ be a convex set. Recall that the topology of a TVS is both translation-invariant and scale-invariant. From Lemma \ref{MinkGauge}, wlog, if $U$ is ``\emph{symmetric}'' (in the sense $\ld U=U$ for any $\ld\in\mathbb{F}$ with $|\ld|=1$), then its Minkowski function $p(x):=\inf\left\{t> 0:{x\over t}\in U\right\}=\inf\left\{t> 0:x\in tU\right\}$ is a seminorm. Again wlog, further choose the convex set $U$ to be ``\emph{absorbing}'' (in the sense $\ld U\subset U$ for any $\ld\in\mathbb{F}$ with $|\ld|<1$). (\blue{footnote}\footnote{Note that if $U\subset X$ is a convex set with $0\in U$, then the absorbing property holds automatically (due to convexity).}). Observe that by definition, for any $t>0$, $x\in tU~\Ra~p(x)<t$. Therefore, $tU\subset \{p<t\}:=\{x\in X:p(x)<t\}$ for all $t>0$. In particular, $U\subset\{p<1\}$. On the other hand, suppose $x\in\{p<1\}$ but $x\not\in U$. Then for all $t$ such that $p(x)<t<1$, we have $x\not\in tU\subset U$. This contradicts the fact that $p(x)$ is the infimum of all $t$ satisfying $x\in tU$. Thus $\{p<1\}\subset U$, and so $U=\{p<1\}$.
\end{proof}

\begin{prp}\label{ExtToLinf}
Let $X$ be a metric space. There is an isometric imbedding {\small $X\hookrightarrow\ell^\infty(A,\mathbb{F})$}, for a set $A$.
\end{prp}
\begin{proof}
Fix $x_0\in X$. Consider the map ~$\vphi:X\ra\ell^\infty(X,\mathbb{F}),~x\mapsto\vphi_x$,~ given for each $x\in X$ by
\[
\vphi_x:X\ra\mathbb{F},~~\vphi_x(a):=d(x,a)-d(x_0,a),~~\txt{for all}~~a\in A:=X~(\txt{as a set}).
\]
The map $\vphi$ is well defined since {\small $\|\vphi_x\|=\sup_a|d(x,a)-d(x_0,a)|\leq d(x,x_0)$}. Also, {\small $\|\vphi_x-\vphi_{x'}\|=d(x,x')$}, since
\[
\textstyle\|\vphi_x-\vphi_{x'}\|=\sup_a|\vphi_x(a)-\vphi_{x'}(a)|\leq d(x,x')=|\vphi_x(x)-\vphi_{x'}(x)|\leq\|\vphi_x-\vphi_{x'}\|. \qedhere
\]
\end{proof}

\begin{dfn}[\textcolor{blue}{
\index{Inner! product}{Inner product},
\index{Inner! product space}{Inner product space},
\index{Hilbert space}{Hilbert space},
\index{Semi-inner product}{Semi-inner product},
\index{Pseudo-inner product space}{Pseudo-inner product space},
\index{Semi-inner product space}{Semi-inner product space}}]\label{InnProdSpDfn}
Let $X$ be an $\mathbb{F}$-vector space ($\mathbb{F}$ = $\Real$ or $\Complex$). An \ul{inner product} on $X$ (making $X$ an \ul{inner product space}) is a function $\langle\rangle:X\times X\ra\mathbb{F}$ satisfying the following: For any $x,y,z\in X$ and $\al\in \mathbb{F}$,
\bit
\item[(i)] $\langle \al x,y\rangle=\al\langle x,y\rangle$ and $\langle x+y,z\rangle=\langle x,z\rangle+\langle y,z\rangle$. ~~~~(Linearity)
\item[(ii)] $\langle x,y\rangle=\overline{\langle y,x\rangle}$. ~~~~ (Symmetry)
\item[(iii)] $\langle x,x\rangle>0$ if $x\neq 0$. ~~~~ (Faithfulness)
\eit
Every inner product space $(X,\langle\rangle)$ is a normed space with norm $\|x\|:=\sqrt{\langle x,x\rangle}$. A complete inner product space is called a \ul{Hilbert space}. A function $\langle\rangle:X\times X\ra\mathbb{F}$ that satisfies (i) and (ii) above is a \ul{pseudo-inner product} on $X$ (making $X=(X,\langle\rangle)$ a \ul{pseudo-inner product space}). A function $\langle\rangle:X\times X\ra\mathbb{F}$ that satisfies (i) and (iii) above is a \ul{semi-inner product} on $X$ (making $X=(X,\langle\rangle)$ a \ul{semi-inner product space}).
\end{dfn}

\section{Topological Convergence and Characterizations in Metric Spaces}

\subsection{Completeness and compactness of metric spaces}
\begin{dfn}[\textcolor{blue}{Recall:
\index{Sequence}{Sequence},
\index{Convergent! sequence}{Convergent sequence},
\index{Sequentially! closed}{Sequentially closed},
\index{Sequentially! continuous}{Sequentially continuous},
\index{Sequentially! compact}{Sequentially compact}}]
Let $X$ be a space. A \ul{sequence} in $X$ is a map $s:\Natural\ra X$, $n\mapsto s_n$, which is often written as $(s_n)\subset X$ or simply as $s_n\in X$. A sequence $(x_n)\subset X$ \ul{converges} to $x$ (written $x_n\ra x$) if every neighborhood of $x$ contains all but finitely many $x_n$ (i.e., for any open set $U\ni x$, we have $\{x_N,x_{N+1},\cdots\}\subset U$ for some $N$). A set $A\subset X$ is \ul{sequentially closed} if for every sequence $x_n\in A$, $x_n\ra x\in X$ implies $x\in A$. A map of spaces $f:X\ra Y$ is \ul{sequentially continuous} if for every sequence $x_n\ra x$ in $X$, we have $f(x_n)\ra f(x)$ in $Y$. A set $K\subset X$ is \ul{sequentially compact} if every sequence $x_n\in K$ has a convergent subsequence $x_{n_k}\ra x\in K$.
\end{dfn}
Note that in a metric space $X=(X,d)$, a net (or sequence in particular) $x_\al\ra x$ in $X$ if and only if the net (or sequence in particular) $d(x,x_\al)\ra 0$ in $\Real$. However, as we show below, sequences (rather than nets in general) suffice for characterizing many topological properties of metric spaces.

\begin{rmk}[\blue{Negating cauchyness and convergence in a metric space}]
Let $X$ be a metric space and $(x_n)\subset X$ a sequence. It is clear that $(x_n)$ is Cauchy iff $\forall \vep>0$, $\exists N_\vep$ s.t. $\forall n,m\geq N_\vep$, $d(x_n,x_m)<\vep$. Therefore, $(x_n)$ is \ul{not} Cauchy iff $\exists\vep_0>0$ s.t. $\forall N$, $\exists n_N,m_N$ s.t. $d(x_{n_N},x_{m_N})\geq\vep_0$. That is, $(x_n)$ is not Cauchy iff there exists $\vep_0>0$ and a subsequence $(x_{n_k})\subset (x_n)$ satisfying $d(x_{n_k},x_{n_{k'}})\geq\vep_0$ for all $k,k'$.

Let $x\in X$. Then similarly, $x_n\ra x$ iff $\forall \vep>0$, $\exists N_{\vep,x}$ s.t. $\forall n\geq N_{\vep,x}$, $d(x_n,x)<\vep$. Therefore, $x_n\not\ra x$ iff $\exists\vep_{0,x}>0$ s.t. $\forall N$, $\exists n_{N,x}$ s.t. $d(x_{n_{N,x}},x)\geq\vep_{0,x}$. That is, $x_n\not\ra x$ iff there exists $\vep_{0,x}>0$ and a subsequence $(x_{n_{k,x}})\subset(x_n)$ such that $d(x_{n_{k,x}},x)\geq\vep_{0,x}$ for all $k$.
\end{rmk}

\begin{facts}[\blue{Metric space peculiarities}]
Let $X,Y$ be metric spaces. The following are true:
\bit
\item[(i)] A set $A\subset X$ is closed $\iff$ sequentially closed.
\item[(ii)] A map $f:X\ra Y$ is continuous $\iff$ sequentially continuous.
\item[(iii)] A set $A\subset X$ is compact $\iff$ sequentially compact.
\eit
\end{facts}
\begin{proof}
{\flushleft (i)} If $A$ is closed then by Lemma \ref{ClosureCrit} (on page \pageref{ClosureCrit}), every $x_n\in A$ satisfies ``$x_n\ra x\in X$ implies $x\in\ol{A}=A$''. If $A$ is not closed, then with $x\in \ol{A}\backslash A$, we can pick a sequence $x_n\in B_{1/n}(x)\cap A$, which satisfies $1/n\geq d(x,x_n)\ra 0$ and so $x_n\ra x\not\in A$.
{\flushleft (ii)} If $f$ is continuous, it follows from Corollary \ref{ContnyCrit3n} that ``$x_n\ra x$ implies $f(x_n)\ra f(x)$''. If $f$ is not continuous, then for some closed set $C\subset Y$, $f^{-1}(C)$ is not closed, i.e., by (i) there exists $x_n\in f^{-1}(C)$ such that $x_n\ra x\not\in f^{-1}(C)$, which implies $f(x_n)\not\!\!\ral f(x)$ since $C$ is closed and $f(x_n)\in C$ but $f(x)\not\in C$.
{\flushleft (iii)} This is Corollary \ref{MetricCompact2} below.
\end{proof}

\begin{dfn}[\textcolor{blue}{
\index{Bounded! set}{Bounded set},
\index{Totally bounded set}{Totally bounded set},
\index{Cauchy sequence}{Cauchy sequence},
\index{Complete! space}{Complete space},
\index{Bounded! map}{bounded map}}] Let $X$ be a metric space and $A\subset X$. The set $A$ is \ul{bounded} if $A\subset B_R(x)$ for some $R>0$, $x\in X$. The set $A$ is \ul{totally bounded} if for any $\vep>0$ there exist points $x_1,...,x_n\in X$ such that $A\subset\bigcup_{i=1}^nB_\vep(x_i)$ (i.e., $A$ can be covered with a finite number of balls of any radius). A sequence $x_n\in X$ is \ul{Cauchy} if $d(x_n,x_m)\ra 0$ as $m,n\ra\infty$. The subspace $A\subset X$ is \ul{complete} if every Cauchy sequence in $A$ converges (to a point) in $A$. If $Z$ is a set, a map $f:Z\ra X$ is a \ul{bounded map} if its image $f(Z)\subset X$ is a bounded set.
\end{dfn}

\begin{rmk}[\blue{\index{Completeness of $\Real$}{Completeness of $\Real$}}]
The lub (least upper bound) property of $\Real$ makes $\Real$ a complete metric space. Indeed, a Cauchy sequence $(x_n)_{n\in\Natural}\subset\Real$ is bounded, and so the decreasing sequence $u_n:=\sup_{m\geq n}x_m$ converges to $x:=\limsup_{n\ra\infty}x_n:=\inf\{u_n:n\in\Natural\}=\inf_{n\in\Natural}\sup_{m\geq n}x_m$. Hence $x_n\ra x$.

This result also follows from Theorem \ref{MetricCompact1} and Corollary \ref{HeinBorelCrl}: Indeed, because a Cauchy sequence in $\Real$ is bounded (hence lies in a compact set), it has a convergent subsequence, and so converges.
\end{rmk}

\begin{thm}[\blue{Compactness in a metric space}]\label{MetricCompact1}
Let $X$ be a metric space and $A\subset X$. The set $A$ is compact $\iff$ complete and totally bounded.
\end{thm}
\begin{proof}
{\flushleft ($\Ra$):} Assume $A$ is compact. It is clear that $A$ is totally bounded since every open cover of $A$ (using balls in particular) has a finite subcover. Also, $A$ is sequentially compact [[i.e., every sequence $(a_i)\subset A$ has a convergent subsequence: Indeed if $a_i\in A$ is an infinite sequence with no convergent subsequence, then $K=\{a_i\}$ is a closed set with no limit point, and so $A\subset(A-K)\cup\bigcup_i B_{\vep_i}(a_i)$ with $B_{\vep_i}(a_i)\cap K=a_i$ for all $i$. By the compactness of $A$, $K$ is finite (a contradiction)]]. Hence $A$ is complete, since a Cauchy sequence converges iff it has a convergent subsequence.
{\flushleft ($\La$):} Assume $A$ is complete and totally bounded. Then $A$ is sequentially compact (i.e., every sequence in $A$ has a subsequence that converges in $X$, hence converges in $A$ by completeness). Otherwise, if an infinite sequence $a_i\in A$ has no subsequence that converges in $X$, then $d(a_i,a_j)\not\ra 0$ since $A$ is complete. Thus, there is a subsequence $a_{i_k}$ and $\vep_0>0$ such that $d(a_{i_k},a_{i_{k'}})\geq 2\vep_0$ for all $k,k'$, and so we need an infinite number of balls $B_{\vep_0}(x_r)$ to cover $A$, i.e., $A\not\subset\bigcup_{r=1}^nB_{\vep_0}(x_r)$ for all $n\geq 1$ (a contradiction).

Let $A\subset\bigcup_{\al\in I}U_\al$, $U_\al\subset A$ be an infinite open cover of $A$. Let $A\subset \bigcup_{j=1}^mB_\vep(x_j^\vep)$, which is possible by total boundedness. Then by choosing $\vep$ small enough, we get $B_\vep(x_j^\vep)\cap A\subset U_{\al_j}$, and so $A\subset \bigcup_{j=1}^mU_{\al_j}$. Otherwise, if for every $\vep>0$, some $B_{\vep}(x_{j_\vep}^\vep)\cap A\not\subset U_\al$ for all $\al$, then for every $n$, some $a_n\in \left(B_{1/n}\left(x_n\right)\cap A\right)\backslash U_\al$ for all $\al$, where $x_n:=x_{j_{1/n}}^{1/n}$. By the sequential compactness of $A$, a subsequence $a_{n_k}\ra a\in A$, and so
\bit
\item[] $d(x_{n_k},a)\leq d(x_{n_k},a_{n_k})+d(a_{n_k},a)\leq 1/n_k+d(a_{n_k},a)\ra 0~~\Ra~~x_{n_k}\ra a\in A.$
\eit
Since $\{U_\al\}$ covers $A$, some $U_{\al_0}\ni a$, i.e., some $B_\delta(a)\cap A\subset U_{\al_0}$. So, for sufficiently large $k$,
\bit
\item[] $a_{n_k}\in B_{1/n_k}\left(x_{n_k}\right)\cap A\subset B_\delta(a)\cap A\subset U_{\al_0}$, ~ (a contradiction).\qedhere
\eit
\end{proof}

\begin{crl}[\blue{Compactness in a metric space}]\label{MetricCompact2}
Let $X$ be a metric space and $A\subset X$. The set $A$ is compact $\iff$ sequentially compact.
\end{crl}
\begin{proof}
($\Ra$): If $A$ is compact, it is clear from the proof of Theorem \ref{MetricCompact1} that $A$ is sequentially compact. ($\La$): Conversely, assume $A$ is sequentially compact. Then $A$ is complete since every cauchy sequence in $A$ has a convergent subsequence, and so converges in $A$. To show $A$ is totally bounded, fix $\vep>0$. Then for any $a\in A$, the sequence $a_1:=a\in A$, $a_2\in A\backslash N_\vep(a_1)$, $a_3\in A\backslash\big(N_\vep(a_1)\cup N_\vep(a_2)\big)$, $\cdots$, $a_i\in A\backslash\big(N_\vep(a_1)\cup\cdots\cup N_\vep(a_{i-1})\big)$, $\cdots$
satisfies $d(a_i,a_j)\geq \vep$ for all $i,j$, and so terminates (otherwise we get an infinite sequence with no convergent subsequence). That is, $A\subset\bigcup_{i=1}^nN_\vep(a_i)$ for some $n$.
\end{proof}

\begin{crl}\label{HeinBorelCrl}
A subset of $\Real$ (and of $\Real^n$ by induction on $n$) is compact $\iff$ closed and bounded.
\end{crl}
\begin{proof}
($\Ra$) A compact subset of $\Real$ is closed (since $\Real$ is Hausdorff) and bounded (since it is totally bounded).
($\La$) Conversely, a bounded closed subset of $\Real$ is compact (as a closed subset of a compact set) due to the following: The closed-ended interval $[0,1]\subset\Real$ (hence every closed-ended interval $[a,b]\subset\Real$) is compact. To see this, suppose $\{O_i\}_{i\in I}$ is an open cover of $[0,1]$ with no finite subcover. Then one of $[0,1/2]$ and $[1/2,1]$, say $[a_1,b_1]$, cannot be covered by a finite number of $O_i$'s. Similarly, one of $[a_1,(a_1+b_1)/2]$ and $[(a_1+b_1)/2,b_1]$, say $[a_2,b_2]$, cannot be covered by a finite number of $O_i$'s. Continuing this way, at the $n$th step with the interval $[a_n,b_n]$, one of $[a_n,(a_n+b_n)/2]$ and $[(a_n+b_n)/2,b_n]$, say $[a_{n+1},b_{n+1}]$, cannot be covered by a finite number of $O_i$'s. We get a strictly decreasing sequence of closed intervals $\{[a_n,b_n]\}_{n\in\Natural}$ with the following properties: For each $n=0,1,2,...$ (where $[a_0,b_0]:=[0,1]$),
\bit
\item[(i)] $[a_n,b_n]\supsetneq[a_{n+1},b_{n+1}]$, and so in particular $a_n\leq a_{n+1}\leq b_{n+1}\leq b_n$.
\item[(ii)] If $x,y\in [a_n,b_n]$, then $|x-y|\leq {1\over 2^n}$.
\item[(iii)] $[a_n,b_n]$ is not covered by a finite number of $O_i$'s.
\eit
Let $a:=\sup\{a_n\}_{n\in\Natural}$. Then $a\in [a_n,b_n]$ for all $n$, i.e., $\bigcap_n[a_n,b_n]\neq\emptyset$. Let $a\in O_{i_a}$. Then $B_r(a)\subset O_{i_a}$ for some $r>0$. Choose $n_r$ such that $2^{-n_r}<r$. Then $|x-y|<r$ for all $x,y\in[a_{n_r},b_{n_r}]$, and so $[a_{n_r},b_{n_r}]\subset B_r(a)\subset O_{i_a}$ (a contradiction).
\end{proof}

\begin{lmm}
A metric space is separable $\iff$ second countable (i.e., it has a countable base).
\end{lmm}
\begin{proof}
($\Ra$): Let $X$ be a separable metric space, with a countable dense set $\{x_n\}_{n\in\Natural}$. Then, by the base criterion, the balls of rational radii centered at the $x_n$'s, $\{B_q(x_n):n\in\Natural,q\in\Rational\}$, form a countable base for $X$. ($\La$): Conversely, recall (from Lemma \ref{SecCountSep}) that any second countable space is separable.
\end{proof}

\subsection{Classical metric function spaces}
\begin{dfn}[\textcolor{blue}{
\index{Retraction}{Retraction},
\index{Retract}{Retract}}]
Let $X$ be a space and $A\subset X$. A continuous map $f:X\ra A$ is a \ul{retraction} of $X$ onto $A$ (making $A$ a \ul{retract} of $X$) if it is an extension of the identity map $1_A:A\ra A$ (i.e., $r|_A=1_A$).
\end{dfn}

\begin{lmm}[\blue{Extension of continuous maps}]
Let $X$ be a space. A subspace $A\subset X$ is a retract of $X$ $\iff$ every continuous map $f:A\ra Y$ extends to a continuous map $F:X\ra Y$.
\end{lmm}
\begin{proof}
($\Ra$): If $A\subset X$ is a retract, then every continuous map $f:A\ra Y$ extends to a continuous map {\small $F=f\circ r:X\sr{r}{\ral}A\sr{f}{\ral}Y$} (where $r$ is a retraction). ($\La$): Conversely, if every continuous map $f:A\ra Y$ extends to a continuous map $F:X\ra Y$, then in particular, $1_A:A\ra A$ extends to a retraction $r:X\ra A$.
\end{proof}

\begin{lmm}\label{ConcaveMetric}
Let $(X,d)$ be a metric space and $\vphi:[0,\infty)\ra[0,\infty)$ a concave increasing map such that (i) $\vphi$ is strictly increasing in $[0,\vep)$ for some $\vep>0$ and (ii) $\vphi(0)\geq0$. Then $d_\vphi:=\vphi\circ d$ is a metric on $X$.
\end{lmm}
\begin{proof}
Since $\vphi$ is strictly increasing near $0$, it suffices to prove the triangle inequality for $d_\vphi$. Observe that for $0\leq \al\leq 1$, we have $\al\vphi(t)\leq \al\vphi(t)+(1-\al)\vphi(0)\leq\vphi(\al t+(1-\al)0)=\vphi(\al t)$. Therefore,
\bea
\textstyle \vphi(t+t')={t\over t+t'}\vphi(t+t')+{t'\over t+t'}\vphi(t+t')\leq\vphi(t)+\vphi(t'),~~~~\txt{for all}~~t,t'\in[0,\infty).\nn
\eea
Given $x,y,z\in X$, let $d_1:=d(x,y)$, $d_2:=d(x,z)$, $d_3:=d(y,z)$. Then because $d_1\leq d_2+d_3$ and $\vphi$ is increasing, we have $\vphi(d_1)\leq\vphi(d_2+d_3)\leq\vphi(d_2)+\vphi(d_3)$, which shows $d_\vphi$ is a metric on $X$.
\end{proof}

\begin{lmm}[\blue{\index{Bounded! metric}{Bounded metric}}]\label{BdedMet}
Let $(X,d)$ be a metric space. Then (i) $d_b(x,y):=\min\big(d(x,y),1\big)$, for $x,y\in X$, is a metric on $X$, and (ii) $(X,d)=(X,d_b)$ topologically, i.e., $(X,d)\cong(X,d_b)$.
\end{lmm}
\begin{proof}
{\flushleft (i)} Set $\vphi(t):=\min(t,1)=
\left\{
  \begin{array}{ll}
    t, & 0\leq t\leq 1 \\
    1, & t>1
  \end{array}
\right\}$ in Lemma \ref{ConcaveMetric}.
{\flushleft (ii)} $d_b$ induces the same topology as $d$, since $d_b=d$ on any ball of radius $<1$, and balls of radii $<1$ form a common base for both topologies.
\end{proof}

\begin{dfn}[\textcolor{blue}{
\index{Uniform! metric}{Uniform metric},
\index{Uniform! metric topology}{{Uniform metric topology}}
\index{Topology! of uniform convergence}{{(Topology of uniform convergence)}},
\index{Metric! function space}{Metric function space},
\index{Product metric}{Product metric}}]
Let $X$ be a set and $Y=(Y,d)$ a metric space. Recall from Lemma \ref{BdedMet} that $d_b(y,y'):=\min\big(d(y,y'),1\big)$ is another metric on $Y$ that induces the same topology as $d$. The \ul{uniform metric} on the set of functions $Y^X:=\{\txt{maps}~f:X\ra Y\}$ (making $Y^X$ a \ul{(uniform) metric function space}) is the map defined as follows:
\[
d_u(f,g):=\sup\big\{d_b(f(x),g(x)):x\in X\big\},~~~~\txt{for}~~~~f,g\in Y^X.
\]
If $Y^X$ is given as a \ul{metric space}, unless it is stated otherwise, we will assume the $d_u$-topology (also called the \ul{topology of uniform convergence}) on $Y^X$. We write $\U(X,Y)$ for the metric space $Y^X$ with the $d_u$-topology.
\end{dfn}
If $Y$ is bounded, then (up to biLipschitz equivalence: \blue{footnote}\footnote{If $\ld\geq 1$, then with $d_b:=\min(d,1)$ and $d_\ld:=\min(d,\ld)$, we have ~$d_b\leq d_\ld\leq \ld d_b$.}) we can simply define $d_u$ by
\[
d_u(f,g):=\sup\big\{d(f(x),g(x)):x\in X\big\},~~~~\txt{for}~~~~f,g\in Y^X.
\]
In subsequent discussions involving $d_u$, we will for convenience not distinguish between $d_b=\min(d,1)$ and $d$. That is, $d$ in $d_u(f,g):=\sup\big\{d(f(x),g(x)):x\in X\big\}$ will be understood to be $d_b$ (especially when $Y$ is unbounded).

More generally, if $\{(X_i,d_i)\}_{i\in I}$ is a family of metric spaces then, as a metric space, $\prod_{i\in I}X_i$ will be given the following metric (the topology of which contains the product topology; see the remark below):
\bea
d_b\big((x_i)_{i\in I},(x'_i)_{i\in I}\big):=\sup_{i\in I}d_{i,b}(x_i,x'_i):=\sup_{i\in I}\min\big(d_i(x_i,x'_i),1\big),\nn
\eea
which becomes the uniform metric $d_u$ whenever there is metric space $(Y,d)$ such that $(X_i,d_i)=(Y,d)$ for all $i\in I$, in which case, $\prod_{i\in I}X_i=\prod_{i\in I}Y=Y^I$. (\blue{footnote}\footnote{(i) When the family $\{X_i\}_{i\in I}$ is uniformly bounded (i.e., $\diam X_i\leq\ld$ for all $i\in I$), then  up to biLipschitz equivalence, we can replace $d_{i,b}$ with $d_i$ in $d_b:=\sup_{i\in I}d_{i,b}$. Also,
(ii) when $I=\{1,2,\cdots,n\}$ is a finite set, then $\sup_{i\in I}d_i(x_i,x'_i)<\infty$ for all $(x_i)_{i\in I},(x'_i)_{i\in I}\in\prod_{i\in I}X_i$. In this case, the natural metric to consider on the product space is the \ul{product metric}
\bea
d\big((x_i)_{i\in I},(x'_i)_{i\in I}\big):=\max_{i\in I}d_i(x_i,x'_i).\nn
\eea}).

\begin{rmk}[\textcolor{blue}{The uniform metric topology contains the topology of pointwise convergence}]~\\~
Let $X$ be a set and $Y=(Y,d)$ a metric space.~\\~
(i) We can see the assertion immediately by recalling that the product topology is the smallest topology on $Y^X$ with respect to which the evaluations $p_x:=e_x:Y^X\ra Y,~f\mapsto f(x)$ are continuous, while
\[
d(e_x(f),e_x(g))=d(f(x),g(x))\leq d_u(f,g),~~~~\txt{(i.e., $e_x$ is $1$-Lipschitz wrt $d_u$)}
\]
\[\bt
Y &  & Y \\
  & Y^X\ar[ul,"p_x"']\ar[ur,"p_{x'}"] & \\
  & (Y^X,d_u)\ar[u,dashed,"id"']\ar[uul,bend left,"p_x"]\ar[uur,bend right,"p_{x'}"'] &
\et~~~~
\begin{minipage}{10cm}
The continuity of projections $p_x$'s, hence of the identity $id:(Y^X,d_u)\ra Y^X$, implies product-open sets are $d_u$-open.
\end{minipage}
\]
We also have the following direct alternatives:\\
(ii) Let $V\subset Y$ be an open set and $z\in X$. Let $f\in V^z:=\{g\in Y^X:g(z)\in V\}$, i.e., $f(z)\in V$. Pick $\vep>0$ such that $B_{2\vep}(f(z))\subset V$. Then
\bea
&&\textstyle B_\vep^{d_u}(f)=\{g\in Y^X:d_u(f,g)<\vep\}\subset \ol{B}_\vep^{d_u}(f)=\{g\in Y^X:d_u(f,g)\leq\vep\}\nn\\
&&\textstyle~~~~=\{g\in Y^X:d(f(x),g(x))\leq\vep,~x\in X\}=\bigcap_{x\in X}\{g\in Y^X:d(f(x),g(x))\leq\vep\}\nn\\
&&\textstyle~~~~=\bigcap_{x\in X}\{g\in Y^X:g(x)\in\ol{B}_\vep(f(x))\}\subset V^z,\nn
\eea
where the last containment holds because for $x=z$,~ $g(z)\in\ol{B}_\vep(f(z))\subset B_{2\vep}(f(z))\subset V$ implies $g\in V^z$. Hence $d_u$-convergence implies pointwise convergence. We can alternatively see this with nets as follows:\\
(iii) Observe that if a net $f_\al\sr{d_u}{\ral} f$ in $Y^X$, i.e., for any $\vep>0$, there exists $\al_\vep$ such that
\bea
d_u(f_\al,f)<\vep~~~~\txt{for all}~~\al\geq\al_\vep,\nn
\eea
then for any $x\in X$, we have $f_\al(x)\ra f(x)$ in $Y$, since
\bea
d(f_\al(x),f(x))=d(e_x(f_\al),e_x(f))\leq d_u(f_\al,f)<\vep~~~~\txt{for all}~~\al\geq\al_\vep,\nn
\eea
and so uniform convergence implies pointwise convergence.
\end{rmk}

\begin{lmm}[\blue{Completeness of the uniform metric}]\label{CompUniMet}
Let $X$ be a set, $Y=(Y,d)$ a metric space, and $\U(X,Y)$ the uniform-metric function space. If $Y$ is complete, then so is $\U(X,Y)$.
\end{lmm}
\begin{proof}
Assume $(Y,d)$ is complete, and let $(f_n)\subset\U(X,Y)$ be a Cauchy sequence. Then for any $\vep>0$, there exists $n_\vep$ such that ~$d_u(f_n,f_m)<\vep$ for all $n,m\geq n_\vep$. For every $x\in X$, $\big(f_n(x)\big)\subset Y$ is Cauchy (and so converges in $Y$) since ~$d\big(f_n(x),f_m(x)\big)\leq d_u(f_n,f_m)\ra 0$.

For each $x\in X$, let $f_n(x)\ra f(x)$ in $Y$. Then $(f_n)\subset\U(X,Y)$ converges to $f$. Indeed, we know ``$d_u(f_n,f_m)<\vep$ for $n,m\geq n_\vep$''. So, for any fixed $n\geq n_\vep$, we have
\begin{align}
&d\left(f_{n+m}(x),f_n(x)\right)\leq d_u\left(f_{n+m},f_n\right)<\vep,~~~~\txt{for all}~~m\geq 1,~~x\in X,\nn\\
&~~\sr{(s)}{\Ra}~~d\left(f(x),f_n(x)\right)=\lim_{m\ra\infty}d\left(f_{n+m}(x),f_n(x)\right)\leq\vep,~~~~\txt{for all}~~x\in X,\nn\\
&~~\Ra~~d_u(f,f_n)\leq \vep,~~\Ra~~f_n\ra f~~\txt{in}~~\U(X,Y),\nn
\end{align}
where step (s) holds by {\footnotesize $\Big|d\big(f_{n+m}(x),f_n(x)\big)-d\big(f(x),f_n(x)\big)\Big|\leq d\big(f_{n+m}(x),f(x)\big)\ra0$} as $m\ra\infty$.
\end{proof}

\begin{lmm}[\blue{\index{Continuous limit theorem}{Continuous limit theorem}}]\label{ContLimThm}
Let $X$ be a space, $Y=(Y,d)$ a metric space, and $\U(X,Y)$ the metric function space. (i) The space of continuous functions $C(X,Y)\subset\U(X,Y)$ is closed. (ii) If $X$ is also a metric space, the space of uniformly continuous functions $UC(X,Y)\subset\U(X,Y)$ is closed.
\end{lmm}
\begin{proof}
(i) Let $(f_n)\subset C(X,Y)$ be such that $f_n\ra f\in\U(X,Y)$. We need to show $f\in C(X,Y)$. That is, for any $x\in X$ and $\vep>0$, there exists an open set $O\ni x$ such that
\bea
f(O)\subset B_\vep(f(x)):=\{y\in Y:d(y,f(x))<\vep\}.\nn
\eea
Since $f_n\ra f$ uniformly, there exists $N=N_\vep$ such that
\bea
d\big(f_n(x),f(x)\big)\leq d_u(f_n,f)<\vep/3,~~~~\txt{for all}~~n\geq N,~~x\in X.\nn
\eea
Fix $x\in X$. Since $f_N$ is continuous at $x$, there exists a neighborhood $O_x$ of $x$ such that
\bea
d\big(f_N(x),f_N(x')\big)<\vep/3,~~~~\txt{for all}~~x'\in O_x.\nn
\eea
Therefore, $f$ is also continuous at $x$, because for all $x'\in O_x$,
\bea
d\big(f(x),f(x')\big)\leq d\big(f(x),f_N(x)\big)+d\big(f_N(x),f_N(x')\big)+d\big(f_N(x'),f(x')\big)<\vep.\nn
\eea
{\flushleft (ii)} Now, assume $X$ is a metric space and, as before, let $\{g_n\}\subset UC(X,Y)$ be such that $g_n\ra g\in\U(X,Y)$. We need to show $g\in UC(X,Y)$. Since $g_n\ra g$ uniformly, there exists $N=N_\vep$ such that
\bea
d\big(g_n(x),g(x)\big)\leq d_u(g_n,g)<\vep/3,~~\txt{for all}~~n\geq N,~~x\in X.\nn
\eea
By the uniform continuity of $g_N$, there exists $\delta=\delta_\vep>0$ such that
\bea
d(x,x')<\delta~~\Ra~~d\big(g_N(x),g_N(x')\big)<\vep/3,~~\txt{for all}~~x,x'\in X.\nn
\eea
Thus, $g$ is also uniformly continuous, because (for all $x,x'\in X$) the bound $d(x,x')<\delta$ implies
\bea
d\big(g(x),g(x')\big)\leq d\big(g(x),g_N(x)\big)+d\big(g_N(x),g_N(x')\big)+d\big(g_N(x'),g(x')\big)<\vep.\nn\qedhere
\eea
\end{proof}

\begin{crl}[\blue{Continuous limit theorem}]\label{ContLimCrl}
Let $X$ be a space, $Y=(Y,d)$ a complete metric space, and $\U(X,Y)$ the uniform-metric function space. (i) The space of continuous functions $C(X,Y)\subset\U(X,Y)$ is complete. (ii) If $X$ is also a metric space, the space of uniformly continuous functions $UC(X,Y)\subset\U(X,Y)$ is complete.
\end{crl}

\begin{lmm}[\blue{\index{Uniform! extension theorem}{Uniform extension theorem}}]\label{UniExtThm}
Let $X,Y$ be metric spaces and $f:E\subset X\ra Y$ a map. If (i) $f$ is uniformly continuous, (ii) $E$ is dense, and (iii) $Y$ is complete, then $f$ extends to a unique uniformly continuous map $F:X\ra Y$.

{\flushleft(Note: By its proof below, the result remains true if "uniformly" is replaced with "Lipschitz".)}
\end{lmm}
\begin{proof}
Define a map $F:X\ra Y$ by $F(x):=\lim f(e_n)$ for any sequence $e_n\in E$ such that $e_n\ra x$. $F$ is \ul{\emph{well defined}} because if $e_n\ra x$, $e_n'\ra x$ and $f(e_n)\ra y$, $f(e_n')\ra y'$, then
\bea
d(y,y')\leq d(y,f(e_n))+d\left(f(e_n),f(e_n')\right)+d(f(e_n'),y')\ra 0,~~\Ra~~y=y',\nn
\eea
where $d\left(f(e_n),f(e_n')\right)\ra 0$ by the uniform continuity of $f$ (because $d(e_n,e_n')\leq d(e_n,x)+d(x,e_n')\ra 0$). It is also clear that $F|_E=f$, i.e., $F$ is an \ul{\emph{extension}} of $f$.

To see that $F$ is \ul{\emph{uniformly continuous}} (hence \ul{\emph{unique}} since continuous maps that agree on a dense set agree everywhere), observe that
\bea
d(F(x),F(x'))=d\left(\lim f(e_n),\lim f(e_n')\right)=\lim_{k,k'} d\left(f(e_k),f(e_{k'}')\right)\sr{(s)}{\ral} 0~~\txt{uniformly as}~~x\ra x',\nn
\eea
where step (s) holds because there always exists $N_{x,x'}=N\big(d(x,x')\big)\geq 1$ such that
\bea
&& d(e_k,x)\leq d(x,x')~~\txt{and}~~d(e'_{k'},x')\leq d(x,x'),~~~~\txt{for all}~~k,k'\geq N_{x,x'},\\
&&~~\Ra~~d(e_k,e_{k'}')\leq d(e_k,x)+d(x,x')+d(x',e_{k'}')\leq 3d(x,x'),~~~~\txt{for all}~~k,k'\geq N_{x,x'}.\nn
\eea
\ul{Step (s)}: By the uniform continuity of $f$, for any $\vep>0$, there exists $\delta=\delta_\vep$ such that for all $k,k'\geq 1$,
\bea
d(e_k,e'_{k'})<3\delta~~\Ra~~d(f(e_k),f(e'_{k'}))<\vep/3.
\eea
Also, there exist $0<\delta_x,\delta_{x'}\leq d(x,x')$ such that for any $k,k'\geq N_{x,x'}$,~ $d(x,e_k)<\delta_x$ $\Ra$ $d(F(x),f(e_k))<\vep/3$, and $d(x',e'_{k'})<\delta_{x'}$ $\Ra$ $d(F(x'),f(e'_{k'}))<\vep/3$. Hence, $d(x,x')<\delta$ implies (with $k,k'\geq N_{x,x'}$)
\[
d(F(x),F(x'))\leq d(F(x),f(e_k))+d(f(e_k),f(e'_{k'}))+d(f(e'_{k'}),F(x'))<\vep. \qedhere
\]
\end{proof}

\subsection{Baire spaces and completion of a metric space}
\begin{dfn}[\textcolor{blue}{
\index{$G$-delta set}{$G$-delta set},
\index{$F$-sigma set}{$F$-sigma set}}]
Let $X$ be a space. A set $G\subset X$ is a \ul{$G$-delta set} if it is a countable intersection of open sets, i.e., $G=\bigcap_{i=1}^\infty G_i$, for open sets $G_i\subset X$. A set $F\subset X$ is an \ul{$F$-sigma set} if it is a countable union of closed sets (i.e., the complement of a $G$-delta set).
\end{dfn}

\begin{thm}[\blue{\index{Baire category theorem}{Baire category theorem}}]\label{BairCatThm1}
Let $X=(X,d)$ be a complete metric space. If $\{D_n\}_{n=1}^\infty$ is a countable collection of dense open (or dense $G$-delta) sets in $X$, then $\bigcap_{n=1}^\infty D_n$ is dense in $X$.
\end{thm}
\begin{proof}
Let $O\subset X$ be a nonempty open set. We need to show $O\cap\bigcap_{n=1}^\infty D_n\neq\emptyset$. Observe that  $O\cap\bigcap_{n=1}^\infty D_n=\bigcap_{n=1}^\infty O\cap D_n=\bigcap_{n=1}^\infty E_n$,  where $E_n:=O\cap D_n$ are dense open sets in the complete metric space $Y:=\ol{O}=(\ol{O},d)$.

Let $x_1\in E_1$. Since $E_1$ is open in $Y$, some $B_{2r_1}(x_1)\subset E_1$, and so $\ol{B_{r_1}(x_1)}\subset E_1$. Since $E_2$ is dense in $Y$, some $x_2\in E_2\cap B_{r_1}(x_1)\backslash\{x_1\}$, and so (as above) some $\ol{B_{r_2}(x_2)}\subset E_2\cap B_{r_1}(x_1)\backslash\{x_1\}$. Similarly, some $\ol{B_{r_3}(x_3)}\subset E_3\cap B_{r_2}(x_2)\backslash\{x_2\}$. Continuing this way, at the $n$th step, we get $\ol{B_{r_n}(x_n)}\subset E_n\cap B_{r_{n-1}}(x_{n-1})\backslash\{x_{n-1}\}$, resulting in a decreasing sequence of closed sets $\ol{B_{r_n}(x_n)}\supset\ol{B_{r_{n+1}}(x_{n+1})}$ of sizes $2^{-n}\geq r_n>r_{n+1}\ra 0$. This means the sequence $\{x_n\}\subset Y$ is Cauchy, since for $n<m$,
\bea
\textstyle d(x_n,x_m)\leq\sum_{i=n}^{m-1}d(x_i,x_{i+1})\leq \sum_{i=n}^{m-1} 2^{-i}\leq \sum_{i=n}^\infty 2^{-i}=2^{-n+1},\nn
\eea
and so converges to a point $z\in Y$, since $Y$ is complete. Because $\bigcap_n\ol{B_{r_n}(x_n)}$ is a closed set containing every $x_n$, we have $z\in\bigcap_n\ol{B_{r_n}(x_n)}\subset \bigcap_n E_n$. Hence $\bigcap_n E_n\neq\emptyset$. (\blue{footnote}\footnote{Notice that, prior to convergence in terms of the metric, the first part of this proof relies mainly on the fact that a metric space is a regular Hausdorff space.}).
\end{proof}

\begin{dfn}[\textcolor{blue}{\index{Baire space}{Baire space}}]
A space $X$ is a Baire space if every countable intersection of dense open (or dense $G$-delta) sets is dense. (\blue{footnote}\footnote{By the Baire category theorem (BCT), every complete metric space (hence every complete metrizable space) is a Baire space.}).
\end{dfn}

\begin{lmm}
Let $X$ be a space and $\{K_a\}_{a\in A}$ a collection of compact subsets. If for every finite subcollection $\{K_{a_i}\}_{i=1}^n\subset\{K_a\}_{a\in A}$ the intersection $\bigcap_{i=1}^nK_{a_i}\neq\emptyset$, then we have $\bigcap_a K_a\neq\emptyset$.
\end{lmm}
\begin{proof}
Suppose on the contrary that $\bigcap_a K_a=\emptyset$. Fix $K_b\in \{K_a\}_{a\in A}$. Since no point of $K_b$ lies in every other $K_a$, it follows that $\{K_a^c:a\neq b\}$ is an open cover of $K_b$. Since $K_b$ is compact, we have a finite subcover $\{K^c_{a_1},...,K^c_{a_n}\}$ of $K_b$ (with each $a_i\neq b$), and so $K_b\cap K_{a_1}\cap...\cap K_{a_n}=\emptyset$ (a contradiction).
\end{proof}

\begin{crl}
Let $X$ be a space and $\{K_n\}_{n\in\Natural}$ compact subsets. If $K_n\neq\emptyset$ and $K_n\supset K_{n+1}$ (for all $n$), then $\bigcap_{n\in\Natural}K_n\neq\emptyset$.
\end{crl}

\begin{lmm}
A locally compact Hausdorff space is regular and normal.
\end{lmm}
\begin{proof}
Let $X$ be a locally compact Hausdorff space, $x\in X$, and $C\subset X$ a closed set not containing $x$. Let $X\subset Y$, where $Y=\{p\}\cup X$ is the one-point compactification of $X$ (hence $Y$ is a compact Hausdorff space). Then $C$ is a compact subset of $Y$. For each $c\in C$, let $N_c(x)$ and $N(c)$ be disjoint open nbds of $x$ and $c$ respectively in $Y$, i.e., $N_c(x)\cap N(c)=\emptyset$ for all $c\in C$. Since $C$ is compact in $Y$, $C\subset \bigcup_{i=1}^nN(c_i)$ for some $c_1,...,c_n\in C$. Let $U:=\bigcap_i N_{c_i}(x)$ and $V:=\bigcup_iN(c_i)$. Then $U$ and $V$ are disjoint open nbds of $x$ and $C$ in $Y$. Hence $U\cap X$ and $V\cap X$ are disjoint open nbds of $x$ and $C$ in $X$. That is, $X$ is a regular space.

Given disjoint closed sets $C,D\subset X$, by repeating the same steps above with $D$ replacing $x$ (and then using normality to get disjoint open nbds $N_c(D)\supset D$ and $N(c)\ni c$ for each $c\in C$), we also see that $X$ is a normal space.
\end{proof}

\begin{thm}\label{BairCatThm2}
A locally compact Hausdorff space is a Baire space.
\end{thm}
\begin{proof}
Using the preceding lemmas, follow the same steps as in the proof of Theorem \ref{BairCatThm1}.
\end{proof}

\begin{dfn}[\textcolor{blue}{Recall: \index{Dense! set}{Dense set}, \index{Nowhere dense set}{Nowhere dense set}}]
Let $X$ be a space and $A\subset X$. Then $A$ is \ul{dense} in $X$ if $\left(\ol{A}\right)^o=X$ (or equivalently, $\ol{A}=X$). The set $A$ is \ul{nowhere dense} in $X$ if
\bea
\big(\ol{A}\big)^o=\emptyset~~\txt{or equivalently},~~\ol{\left(\ol{A}\right)^c}=X~~\txt{or equivalently},~~\ol{(A^c)^o}=X,\nn
\eea
where the equivalences are due to the identities
\bea
\left(\ol{A}\right)^c=(A^c)^o~~~~\txt{and}~~~~(A^o)^c=\ol{A^c}.\nn
\eea

If $A$ is nowhere dense in $X$, it is clear that $\ol{A}^c$ (hence $A^c$ also) is dense in $X$. However, denseness of $A^c$ in $X$ (i.e., $\ol{A^c}=X$) is equivalent to $A^o=\emptyset$, and so insufficient for nowhere denseness of $A$.
\end{dfn}

\begin{lmm}\label{NoWhDnsUn}
Let $X$ be a Baire space. If $\{C_n\}_{n=1}^\infty\subset\P(X)$ is a countable collection of nowhere dense closed sets, then $\bigcup_n C_n$ has empty interior (hence $\bigcup_n C_n\neq X$).
\end{lmm}
\begin{proof}
By hypotheses (i.e., $C_n$ is closed and $\big(\ol{C}_n\big)^o=\emptyset$), $C_n^c$ is open, and dense because
{\footnotesize\[
\textstyle \Big(\ol{C_n^c}\Big)^o=\Big(\big(C_n^o\big)^c\Big)^o=\Big(\big(\emptyset\big)^c\Big)^o =X^o=X.
\]}Since $X$ is Baire, it follows that $\bigcap C_n^c$ is dense, i.e., {\footnotesize $\Big(\ol{\bigcap C_n^c}\Big)^o=X$}, which implies
{\footnotesize\[
\textstyle \emptyset=\left[\Big(\ol{\bigcap C_n^c}\Big)^o\right]^c=\ol{\Big(\ol{\bigcap C_n^c}\Big)^c}=\ol{\Big(\big(\bigcap C_n^c\big)^c\Big)^o}=\ol{\Big(\bigcup C_n\Big)^o}. \qedhere
\]}
\end{proof}

\begin{dfn}[\blue{\index{Completely metrizable space}{Completely metrizable space}}]
A space is \ul{completely metrizable} if it is homeomorphic to a complete metric space.
\end{dfn}

\begin{crl}
Let $X_n\subset X_{n+1}$ be an increasing sequence of spaces in a subcategory $\C\subset Top$ of spaces such that $X_n$ is closed and nowhere dense in $X_{n+1}$. The \ul{colimit} (\ul{inductive limit}) $X:=\varinjlim_n X_n$ (i.e., the set-union $\bigcup_nX_n$ structured as an object) in the category $\C$ is not completely metrizable. (\blue{footnote}\footnote{The topology of $X$ is given by the following: A subset $A\subset X$ is open (in $X$) $\iff$ for each $n$, $A\cap X_n$ is open in $X_n$.}).
\end{crl}
\begin{proof}
Assume $X$ is metrizable (otherwise there is nothing to prove), i.e., $X\cong Y$ for a metric space $Y$. Suppose $Y$ is complete. Then $Y$ is Baire. Since $X$ (hence $Y$) is a union of nowhere dense closed sets, we get a contradiction of Lemma \ref{NoWhDnsUn}.
\end{proof}

\begin{dfn}[\textcolor{blue}{\index{Dense! map}{Dense map}}]
Let {\small $X,Y$} be spaces. A map {\small $f:X\ra Y$} is \ul{dense} if its image is dense in $Y$.
\end{dfn}
\begin{dfn}[\textcolor{blue}{\index{Completion}{Completion} of a metric space}]
Let $X$ be a metric space. A metric space $Z$ is a \ul{completion} of $X$ if (i) $Z$ is complete and (ii) there exists a dense isometric imbedding $\vphi:X\hookrightarrow Z$.
\end{dfn}
By the following theorem, every metric space has a completion that is unique up to isometry. The existence part of the theorem (with an alternative proof below) is an immediate consequence of the isometric imbedding $X\hookrightarrow\ell^\infty(X,\mathbb{F})$ from Proposition \ref{ExtToLinf}.

\begin{thm}[\blue{Existence and uniqueness of a completion}]\label{MetComThm}
Let $X$ be metric space. There is a complete metric space $\wt{X}$ and an isometric imbedding $f:X\hookrightarrow\wt{X}$ such that $f(X)$ is dense in $\wt{X}$. Moreover, $\wt{X}$ is unique up to isometry.
\end{thm}
\begin{proof}
{\flushleft\bf Uniqueness}: Let $f_1:X\hookrightarrow\wt{X}_1$ and $f_2:X\hookrightarrow\wt{X}_2$ be isometric imbeddings such that $\wt{X}_1,\wt{X}_2$ are complete, $f_1(X)$ is dense in $\wt{X}_1$, and $f_2(X)$ is dense in $\wt{X}_2$. Then we get the map
\bea
h:f_1(X)\subset\wt{X}_1\ra \wt{X}_2,~f_1(x)\mapsto f_2(x)\nn
\eea
which is uniformly continuous on $f_1(X)$ as an isometric map, because
\bea
d(h(f_1(x)),h(f_1(y)))=d(f_2(x),f_2(y))=d(x,y)=d(f_1(x),f_1(y)),~~~~\txt{for all}~~x,y\in X.\nn
\eea
Since $f_1(X)$ is dense in $\wt{X}_1$, and $\wt{X}_2$ is complete, it follows by Theorem \ref{UniExtThm} that $h$ extends to a unique continuous map $\wt{X}_1\ra\wt{X}_2$. Hence, by symmetry, we get an isometry $h:\wt{X}_1\ra \wt{X}_2$.

{\flushleft\bf Existence}: Let $C:=\{\txt{all Cauchy sequences in $X$}\}$. Define a relation $\sim$ on $C$ as follows:
\bea
(x_j)\sim(y_j)~~~~\txt{if}~~~~d(x_j,y_j)\ra 0,~~~~~~~~\txt{for all}~~~~(x_j),(y_j)\in C.\nn
\eea
Then $\sim$ is an equivalence relation due to the following.
\bit
\item[(a)] $d(x_j,x_j)=0\ra 0$ (i.e., $\sim$ is reflexive).
\item[(b)] $d(x_j,y_j)\ra 0$ implies $d(y_j,x_j)=d(x_j,y_j)\ra 0$ (i.e., $\sim$ is symmetric).
\item[(c)] $d(x_j,y_j)\ra 0$, $d(y_j,z_j)\ra 0$ implies $d(x_j,z_j)\leq d(x_j,y_j)+d(y_j,z_j)\ra 0$ (i.e., $\sim$ is transitive).
\eit
Let $\big[(x_j)\big]:=\big\{(y_j)\in C:(y_j)\sim(x_j)\big\}$ be the equivalence class of $(x_j)$ in $C$. Consider the set of equivalence classes ~$\wt{X}:={C\over\sim}:=\big\{\wt{x}=\big[(x_j)\big]:(x_j)\in C\big\}$~ as a metric space with respect to the metric
\bea
\wt{d}\left(\wt{x},\wt{y}\right):=\lim_{j\ra\infty} d(x_j,y_j),\nn
\eea
where $\wt{d}$ is well-defined (i.e., exists and is independent of the choice of representatives of the equivalence classes) because $\big(d(x_j,y_j)\big)$ converges as a Cauchy sequence in $\Real$, and, if $(x'_j)\sim(x_j)$ and $(y'_j)\sim(y_j)$, then $|d(x'_j,y'_j)-d(x_j,y_j)|\leq d(x'_j,x_j)+d(y_j,y'_j)\ra 0$ implies
\bea
\lim_{j\ra\infty} d(x'_j,y'_j)=\lim_{j\ra\infty} d(x_j,y_j).\nn
\eea

\ul{Define the inclusion} $f:X\hookrightarrow \wt{X}$ by $f(x):=\big[(x,x,\cdots)\big]$, where $(x,x,\cdots)$ denotes the constant sequence. Then $f$ is isometric since~ $\wt{d}(f(x),f(y))=\lim d(x,y)=d(x,y)$. Note that every sequence $(x_j)$ such that $x_j\ra x$ belongs to the equivalence class $\big[(x,x,\cdots)\big]$, and so
\bea
f(x):=\big[(x,x,\cdots)\big]=\big\{(x_j)\in C:~x_j\ra x\big\}.\nn
\eea

\ul{To show that $f(X)$ is dense} in $\wt{X}$ we will show that $\wt{X}\subset \overline{f(X)}=f(X)\cup f(X)'$. Let $\wt{x}=\big[(x_j)\big]\in\wt{X}$. If $(x_j)$ converges, then $\wt{x}\in f(X)\subset\overline{f(X)}$. If $(x_j)$ does not converge, let $\vep>0$ be given. We need to show that $B_\vep(\wt{x})\backslash\{\wt{x}\}$ contains an element of $f(X)$. Recall that $\wt{y}=\big[(y_j)\big]\in B_\vep(\wt{x})$
\bea
~~\iff~~\wt{d}(\wt{x},\wt{y})=\lim_j d(x_j,y_j)<\vep.\nn
\eea
Since $(x_j)$ is Cauchy, there is $N=N(\vep)$ such that $d(x_j,x_{j'})<\vep$ for $j,j'\geq N$, and so
for any $k\geq N$,
\bea
&&\wt{d}\left(\wt{x},\big[f(x_k)\big]\right)=\lim_j d(x_j,x_k)<\vep,~~~~\Ra~~~~\big[f(x_k)\big]=\big[(x_k,x_k,\cdots)\big]\in B_\vep(\wt{x})\backslash\{\wt{x}\}.\nn
\eea
Therefore $f(X)$ is dense in $\wt{X}$.

\ul{Next we will show that $\wt{X}$ is complete}. Let $(\wt{x}_k)$, where $\wt{x}_k=\big[(x_{jk})\big]$, be a Cauchy sequence in $\wt{X}$. We need to show that $(\wt{x}_k)$ converges. Let $\vep>0$ be given. Then there is $N_1=N_1(\vep)$ such that
\bea
&&\wt{d}(\wt{x}_k,\wt{x}_{k'})=\lim_j~d(x_{jk},x_{jk'})<\vep~~~~\txt{for}~~~~k,k'\geq N_1.\nn
\eea
Since for each pair $k,k'\geq N_1$, $d(x_{jk},x_{jk'})\ra \wt{d}(\wt{x}_k,\wt{x}_{k'})$, it follows that there is $N_2(k,k')=N_2(\vep,k,k')$ with
\bea
&& |d(x_{jk},x_{jk'})-\wt{d}(\wt{x}_k,\wt{x}_{k'})|<\vep~~~~for~~~~j\geq N_2(k,k'),~~~~k,k'\geq N_1,\nn\\
&&~~\Ra~~d(x_{jk},x_{jk'})< \wt{d}(\wt{x}_k,\wt{x}_{k'})+\vep<2\vep~~~~for~~~~j\geq N_2(k,k'),~~~~k,k'\geq N_1.\nn
\eea
Consider the diagonal sequence $(x_{kk})$. Then because each representative of $\wt{x}_k$ is Cauchy,
\bea
d(x_{kk},x_{k'k'})\leq d(x_{kk},x_{jk})+d(x_{jk},x_{jk'})+d(x_{jk'},x_{kk'})<3(2\vep)=6\vep~~~~~\txt{for}~~k,k'\geq N_1,\nn
\eea
where $j\geq N_2(k,k')$. Thus $(x_{kk})$ is Cauchy, i.e., $\big[(x_{kk})\big]\in \wt{X}$. Therefore, $(\wt{x}_k)$ converges in $\wt{X}$, since
\[
\wt{d}\left(\wt{x}_k,\big[(x_{jj})\big]\right)=\lim_j~d(x_{jk},x_{jj})<6\vep~~~~\txt{for}~~k\geq N_1,~~~~\Ra~~~~\wt{x}_k\ra\big[(x_{jj})\big]\in\wt{X}. \qedhere
\]
\end{proof}

{\flushleft\hrulefill}

\begin{exercise}
Based on the discussion of this chapter and the previous chapter, consider writing a \emph{fully technical} essay (say in the form of a typical chapter of these notes) on what is known in the mathematics literature as ``Functional Analysis''.
\end{exercise}

%% file: parts/front/Vita.tex
\hrulefill
\doublespacing
\vspace{1cm}
\begin{center}
\bit
\item[] \hspace{0cm}{\Large About the author}
\vspace{0.5cm}
\item[] \textbf{Name of author:}~~ Earnest Akofor
\vspace{0.3cm}

%

\item[] \textbf{Degrees awarded:}~
\vspace{0.2cm}
\bit
\item[] Ph.D. Mathematics, Syracuse University, 2020
\vspace{0.2cm}

\item[] Ph.D. Electrical and Computer Engineering, Syracuse University, 2016
\vspace{0.2cm}

\item[] Ph.D. Physics, Syracuse University, 2010
\vspace{0.2cm}

\item[] Diploma. Mathematical Sciences, AIMS, 2004
\vspace{0.2cm}

\item[] Diploma. High Energy Physics, ICTP, 2003
\vspace{0.2cm}

\item[] B.Sc. Physics, University of Buea, 2001
\eit
\vspace{0.5cm}

\item[] \textbf{Professional experience:}~
\vspace{0.2cm}
\bit
\item[] Review, MathRev and zbMATH, 2020
\vspace{0.2cm}

\item[] Instructor \& Teaching Assistant, Syracuse University, 2016-2020
\vspace{0.2cm}

\item[] Peer Review, IJDSN (2015) and IEEE Transactions (2015-2016)
\vspace{0.2cm}

\item[] Research Assistant, Syracuse University, 2011-2016
\vspace{0.2cm}

\item[] Adjunct Instructor, ITT Technical Institute, 2011
\vspace{0.2cm}

\item[] Teaching \& Research Assistant, Syracuse University, 2004-2009
\eit
\eit
\end{center}